# $\theta(p_c, \mathbb{Z}^d) = 0$?

Raphaël Cerf*

August 22, 2026

## 1 $\theta(p_c, \mathbb{Z}^d) = 0$?

We consider the Bernoulli site percolation model on $\mathbb{Z}^d$ with $d \geq 3$. Each site is declared open with probability $p$ and closed with probability $1-p$, and the sites are independent. The percolation probability is

$$\theta(p, \mathbb{Z}^d) \,=\, P_p\big(0 \longleftrightarrow \infty\big)\,.$$

For a parameter $p$ such that $\theta(p, \mathbb{Z}^d) = 0$, with probability one, all the open clusters are finite. On the contrary, if $\theta(p, \mathbb{Z}^d) > 0$, then, with probability one, there exists an infinite open cluster. The critical point of the model is

$$p_c \,=\, \sup\,\big\{\, p \in [0,1] : \theta(p, \mathbb{Z}^d) = 0 \,\big\}\,.$$

It follows from the classical Peierls argument that $p_c$ is strictly between 0 and 1. The definition of $p_c$ readily implies that $\theta(p, \mathbb{Z}^d) = 0$ for $p < p_c$ and $\theta(p, \mathbb{Z}^d) > 0$ for $p > p_c$. Now, what is the value of $\theta(p_c, \mathbb{Z}^d)$?

In dimension 2, the case of bond percolation has been settled in 1980 by Kesten [82] (using also a former result of Harris [69] from 1960) and the case of site percolation by Russo [117] in 1981. In 1990, Hara and Slade [66, 67] proved that $\theta(p_c, \mathbb{Z}^d) = 0$ for $d \geq 19$. In 1991, Barsky, Grimmett and Newman [11, 12] proved the analogous result for bond percolation in a $d$-dimensional half-space. After these fundamental progresses, the situation has been stalled for 25 years. In 2015, Damron, Newman and Sidoravicius [35] handled the case of a two-dimensional sandwich of thickness two, that is $\mathbb{Z}^2 \times \{0,1\}$. This result was extended one year later to finite slabs by Duminil-Copin, Sidoravicius and Tassion [45]. In 2017, the mean-field behaviour was proved to hold for $d \geq 11$ by Fitzner and van der Hofstad [53, 54]. Very recently, Kozma and Nitzan [90] exhibited an inequality, which, if true, would settle the problem.

The main open question in percolation theory is to prove that there is no percolation at the critical point in dimension $d = 3$, i.e., that $\theta(p_c, \mathbb{Z}^3) = 0$.

Let us state formally this conjecture.

**Conjecture 1.1.** *For any $d \geq 3$, we have $\theta(p_c, \mathbb{Z}^d) = 0$.*

*Université Paris-Saclay, CNRS, Laboratoire de mathématiques d'Orsay, 91405, Orsay.

The search for the proof has been a long endeavour, and unfortunately, it was not crowned with success. Time and time again, we have lost our way on side roads that were not conclusive, but which have nonetheless led to interesting partial results. The best result achieved so far is the following theorem.

**Theorem 1.2.** *Let $d \geq 3$ and let $p_c$ be the critical point for the Bernoulli site percolation in $\mathbb{Z}^d$. Denoting by $C(0)$ the cluster of the origin, we have either $\theta(p_c, \mathbb{Z}^d) = 0$ or*

$$\liminf_{n \to \infty} \frac{1}{\ln n} \ln \ln \Big( \frac{1}{P_{p_c}\big(n \leq |C(0)| < +\infty\big)} \Big) \;=\; 0\,. \tag{1.1}$$

This result can be rephrased as follows: for the Bernoulli percolation model in $\mathbb{Z}^d$ with $d \geq 3$, at the critical point $p_c$, either there is no infinite cluster or the tail of the finite cluster distribution is not a stretched exponential.

A proper understanding of the percolation model requires an analysis of the tail of the finite cluster size distribution. Of course, for any value of $p$ in $[0,1]$, we have

$$\lim_{n \to \infty} P_p\big(n \leq |C(0)| < +\infty\big) \;=\; 0\,, \tag{1.2}$$

but it is a major challenge to obtain the speed of convergence towards 0 in (1.2). In the subcritical regime, it is known that

$$\forall p < p_c \qquad \lim_{n \to \infty} \frac{1}{n} \ln P_p\big(n \leq |C(0)| < +\infty\big) \;<\; 0\,. \tag{1.3}$$

The existence of the limit follows from a classical subadditive argument (see [63], chapter 6.3), but the more important fact that it is negative is a consequence of the theorem of Menshikov [105], Aizenman and Barsky [1], recently simplified by Duminil-Copin and Tassion [46], together with tricky tree estimates due to Aizenman and Newman (see theorem 6.75 in [63]). The limit (1.3) readily implies that

$$\forall p < p_c \qquad \lim_{n \to \infty} \frac{1}{\ln n} \ln \ln \Big( \frac{1}{P_p\big(n \leq |C(0)| < +\infty\big)} \Big) \;=\; 1\,. \tag{1.4}$$

In the supercritical regime, Kesten and Zhang [84] proved that

$$\forall p > p_c \qquad \lim_{n \to \infty} \frac{1}{\ln n} \ln \ln \Big( \frac{1}{P_p\big(n \leq |C(0)| < +\infty\big)} \Big) \;=\; 1 - \frac{1}{d}\,. \tag{1.5}$$

In fact, Kesten and Zhang proved the limit (1.5) under the stronger condition that $p > \widehat{p}_c$, where $\widehat{p}_c$ is the slab percolation threshold, which was eventually shown to coincide with the critical point $p_c$ by Grimmett and Marstrand [62]. So we see that the existence and the determination of the two limits (1.4) and (1.5) are the outcome of several seminal works on the Bernoulli percolation model. At criticality, the tail of the finite cluster size distribution is expected to decay as an inverse power of $n$, and this would imply the result (1.1). The limit (1.1) tells us that the tail of the finite cluster size distribution is not governed by a stretched exponential, of the form $\exp(-n^\alpha)$ for some positive $\alpha$. Of course this result is everything but spectacular, and one might wonder why it is so difficult to obtain. To be honest, the goal of our enterprise was to prove the conjecture $\theta(p_c, \mathbb{Z}^d) = 0$. To this end, we implemented the classical one step renormalization scheme, and we tried to push it as far as possible. The best result we could achieve was the following.

**Theorem 1.3.** *Let $d \geq 3$. For any $p \in ]0,1[$, we have the implication*

$$\theta(p, \mathbb{Z}^d) > 0 \quad \text{and} \quad \liminf_{n\to\infty} \frac{1}{\ln n} \ln\ln\Big(\frac{1}{P_p\big(n \leq |C(0)| < +\infty\big)}\Big) > 0 \qquad \Longrightarrow \qquad p > p_c\,.$$

Obviously, theorem 1.3 implies theorem 1.2. It is also very likely that the following stronger implication holds: for $d \geq 3$ and for any $p \in ]0,1[$,

$$\liminf_{n\to\infty} \frac{1}{\ln n} \ln\ln\Big(\frac{1}{P_p\big(n \leq |C(0)| < +\infty\big)}\Big) > 0 \qquad \Longrightarrow \qquad p > p_c\,. \tag{1.6}$$

However, we did not try very intensively to prove the implication (1.6), because this implication would not be enough to prove the conjecture $\theta(p_c, \mathbb{Z}^d) = 0$. Our main effort has been to try to get rid of the second condition in theorem 1.3 in order to prove the following implication:

$$\forall p \in ]0,1[ \qquad \theta(p, \mathbb{Z}^d) > 0 \quad \Longrightarrow \quad p > p_c\,.$$

This would readily imply that $\theta(p_c, \mathbb{Z}^d) = 0$. Unfortunately, there is still something missing to achieve this goal. At this point, we believe that the proof of theorem 1.3 is more interesting than the result itself, as it is very promising, and it seems close to succeeding. The crux of the matter is to obtain an adequate control on the negative higher pivots for an event of the disconnection type. The second condition in theorem 1.3 rules out the existence of higher pivots which are far apart. So it remains to control the scenarii involving what we call local higher pivots, a task that we achieve with the help of a process of multiple intertwined explorations. If only we could obtain the same control for the higher pivots without the locality condition, the whole scheme of the proof would go through. In brief, to complete the proof of the conjecture $\theta(p_c, \mathbb{Z}^d) = 0$, one hope is to understand the missing ingredient to make the scheme of proof of theorem 1.3 work. Another hope would be to prove the following conjecture.

**Conjecture 1.4.** *Let $d \geq 3$. For any $p \in ]0,1[$, we have the implication*

$$\theta(p, \mathbb{Z}^d) > 0 \quad \Longrightarrow \quad \liminf_{n\to\infty} \frac{1}{\ln n} \ln\ln\Big(\frac{1}{P_p\big(n \leq |C(0)| < +\infty\big)}\Big) > 0\,.$$

Needless to say, this still seems to be a formidable task, but it could be that the correct way to prove that $\theta(p, \mathbb{Z}^d) = 0$ is to prove separately theorem 1.3 and conjecture 1.4.

**Thanks.** This question took hold on my mind in 1995, when, as a young post-doc, I attended Geoffrey's Grimmett beautiful lectures in Cambridge. Over the last thirty years, I have been discussing this question with every percolationist I met, and it has always been a great pleasure to do so. I am deeply grateful to Jeff Steif, who gave me the most precious piece of advice, during my visit to Göteborg in 2008. He told me to go and look at the half-forgotten uniqueness proof of Aizenman, Kesten and Newman [4]. This gave me a new impulse at a time when I was essentially stuck on the problem. It is also fair to say that [4] contains a watermark of the premises of the proof presented here.

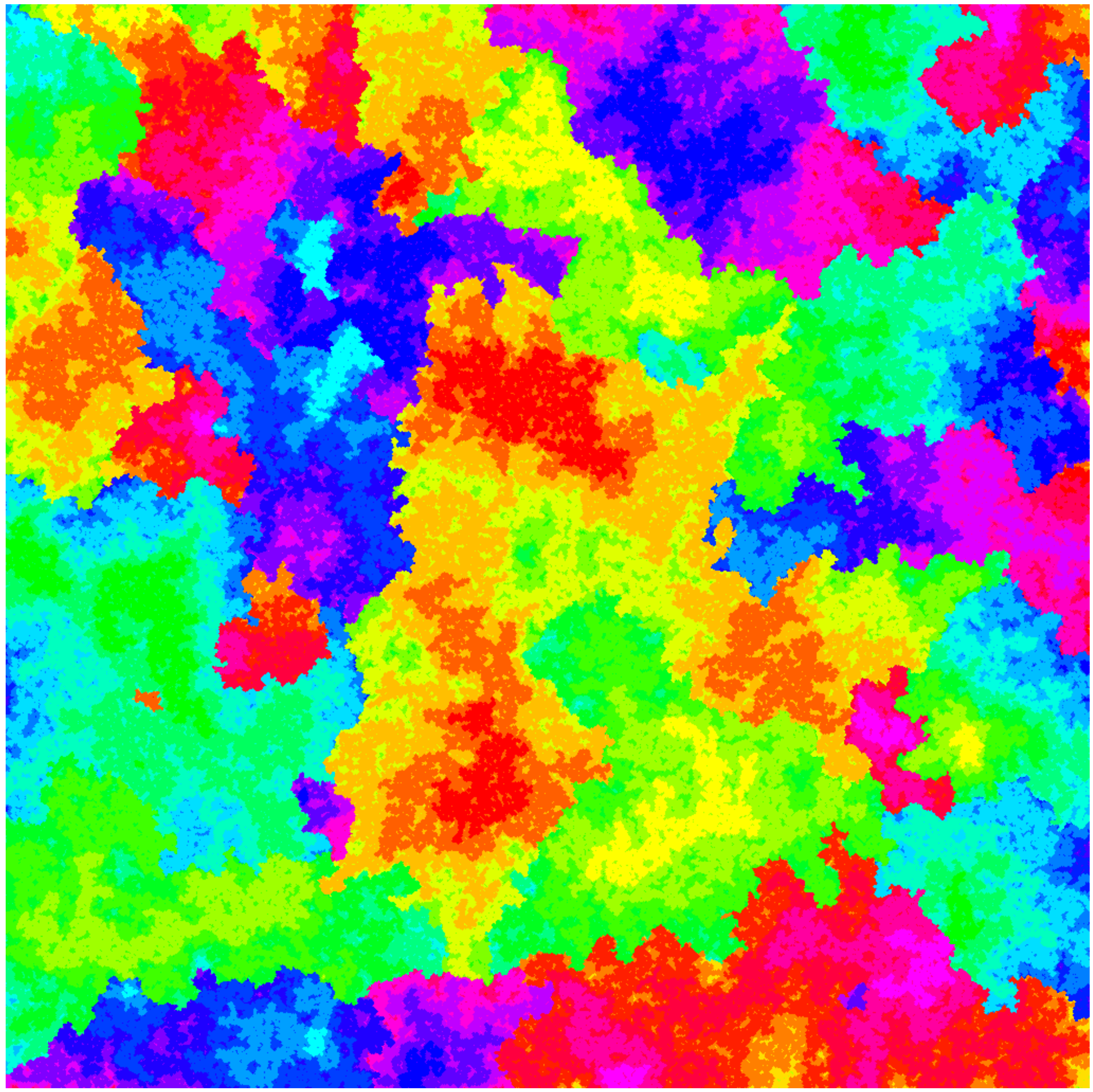

Figure 1: Fifty intertwined explorations, $\Lambda(1024)$, site percolation, $p = 0.57$.

**Regrets.** I regret that Harry Kesten and Vladas Sidoravicius, percolation masters, if ever there were, could not read the proof of theorem 1.3: they would certainly find what is missing to remove the question mark from the title.

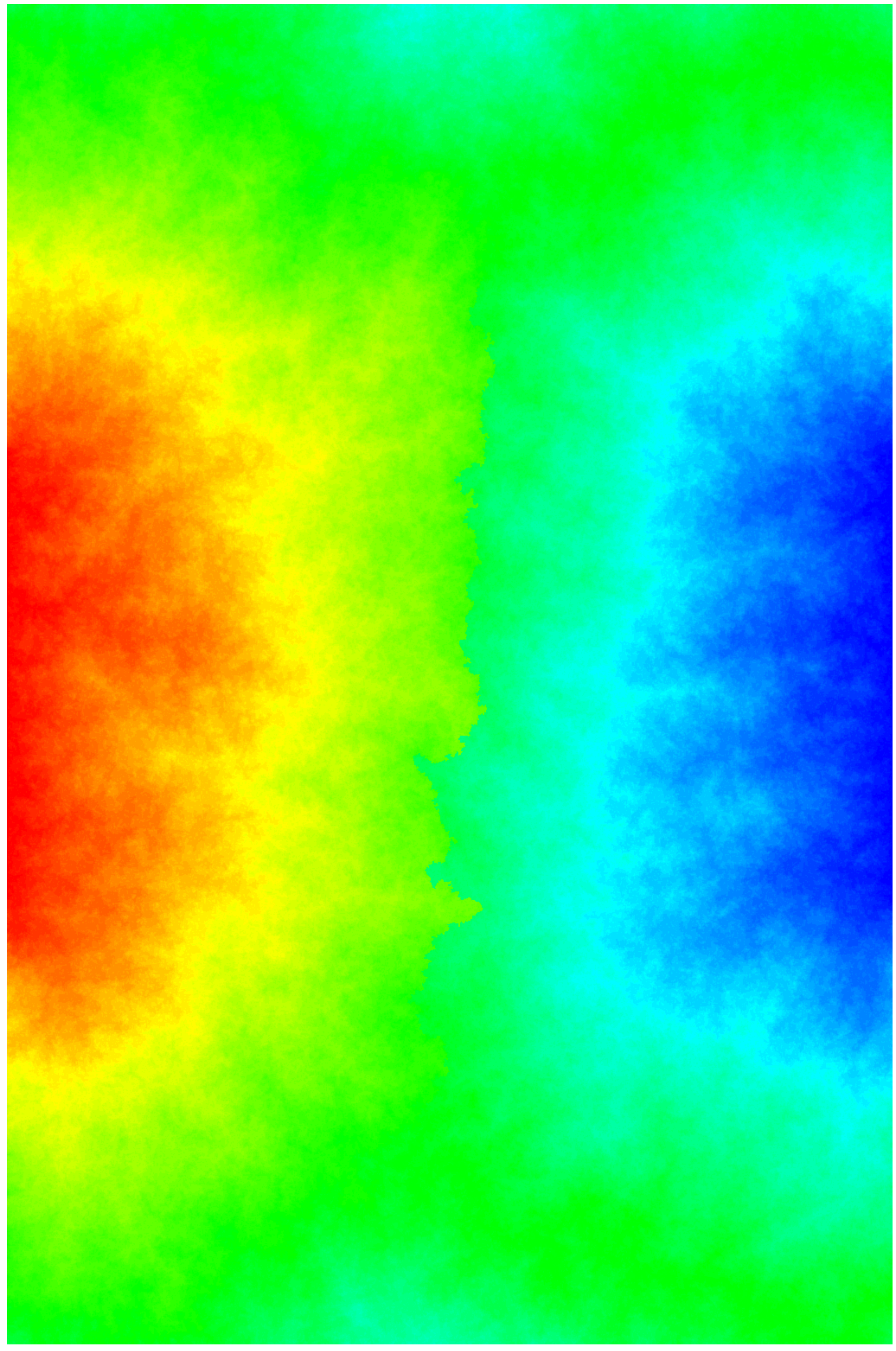

Figure 2: The heart of the matter: disconnection between the middle thirds of left and right sides in the rectangle $1550 \times 1024$, bond percolation, $p = 0.44$.

# Contents

## 2 A guide to the text

Usually, mathematical texts contain perfect proofs and the many unsuccessful attempts to obtain them remain hidden. The $\theta(p_c, \mathbb{Z}^d) = 0$ problem is so fascinating that unsuccessful attempts have themselves generated some extremely interesting problems that do not deserve to be forgotten. It was also by trying to write them down precisely that we were able to understand the gaps that were to be filled and to make slow progress. Numerous simulations have also been implemented, in order to check some partial conjectures, or to improve the understanding of the obstacles. So the text contains a selection of the material developed over the years around the main line of argument (simulations, comments, failed attempts, related mathematical results). The potential reader should not be frightened by the length of the text, which is not optimized to furnish the proof of the main results in the shortest possible space. To facilitate the reading and make a reasonably self-contained text, the relevant tools are introduced properly and a lot of pictures (171 in total!) are included. As a matter of fact, we also took the time and liberty to discuss various concepts pertaining to mathematical fields tangential or even quite distant from percolation theory. One instance is monotone boolean functions, another one is extremal set theory and the Kruskal-Katona theorem. So this text is not intended to be a mathematical article providing the most effective proof of a selection of a few results. The reader familiar with the Bernoulli percolation model who is only interested in the actual proof of theorem 1.3 might skip the sections marked with a circled asterisk $\circledast$ and focus on the sections marked with a sun $\odot$. Indeed, the whole proof of theorem 1.3 for the site model is contained in the sections marked with $\odot$ and does not make appeal to the material presented in the sections marked with $\circledast$. We list next the numbers of the sections and subsections containing the proof, along with the corresponding number of pages:

| sections | 8.1-8.4 | 9-11 | 12-14 | 16-17 | 39-40 |
|---|---|---|---|---|---|
| page numbers | 48-53 | 60-66 | 68-88 | 100-107 | 268-277 |
| number of pages | 6 | 7 | 21 | 8 | 10 |

| sections | 44.1 | 49 | 80-81 | 84-93 |
|---|---|---|---|---|
| page numbers | 291 | 373 | 490-523 | 527-695 |
| number of pages | 1 | 1 | 34 | 169 |

In total, the number of pages needed for proving theorem 1.3 is approximately

$$257 \;=\; 6+7+21+8+10+1+1+34+169\,.$$

More than one third of the text is devoted to this proof. If one wants also to include the proof for the bond model, one should at least add the 24 pages for the construction of the intertwined bond exploration presented in section 29. The rest of the text presents various failed attempts, and also partial results.

We have made every effort to provide detailed proofs, which makes some passages difficult to read. Some might say that the text is too detailed. In fact, we tried to ensure that no calculations were left out, and the formulas are presented in such a way that they can be followed mentally without the need for additional calculations. In several places, an argument had to be repeated twice, because it required several modifications. Confronted with this situation, there are typically three possible choices:

- to find a slightly more general setting in which both arguments fit;
- to say that it suffices to repeat the argument with minor modifications;
- to rewrite the argument with the required modifications.

If we had opted for the shortest possible presentation, we would have always opted for the second choice. However, for all matters related to the main proof, we felt that it was necessary to give all the important details. So we opted for either the first or the third choice. The third choice has the drawback of resulting in some repetition. However, these passages are well identified and can easily be skipped.

The initial plan was to work exclusively in the context of Bernoulli site percolation, but this turned out to be far too difficult. The site model is undeniably more complicated to comprehend than the bond model. Several arguments had to be written in the first place for the bond model, before getting a chance to understand how to proceed for the site model. Some inconclusive results have been written only for the bond model. Unfortunately, this state of affairs does not facilitate the reading of the text, which switches constantly between the two models. Here is a little sign to help the reader. When a section is marked with a bond sign ×—× , it means that it has been written for the bond model. When a section is marked with a site sign $\circ$, it means that it has been written for the site model. In the absence of these two signs to the right of the title, the content of the section is valid for both models, either directly as it is, or with minor straightforward adaptations. When a part (respectively a section) is marked with a sign $\circledast$, $\odot$, ×—× or $\circ$, the same sign applies to all the sections of the part (respectively the subsection of the section), unless an antagonist sign pops up. The heart of the matter is to obtain a quantitative estimate for a certain disconnection event. Throughout the years, we focused on the disconnection between the left and right faces of a cubic box. The pages which are traditionally left blank, for page layout purposes are used to present a simulation of this event, colored according to an intertwined exploration process. Towards the end, this event is replaced by pictures of 7 intertwined explorations.

**AI.** No artificial intelligence was used to generate any mathematical content nor to perform heavy numerical simulations. We used the AIs ChatGPT, Claude and Le Chat during the last stages of the writing to configure Neovim, to generate Python code for five illustrations, and to do some proofreading. For the proofreading, AI was used to identify potential issues in selected passages, which were then verified and corrected manually.

**English.** DeepL was used occasionally to better the quality of the English text. The text was spell-checked by TeXtidote, a program written by Sylvain Hallé.

# List of Algorithms

The above list contains only the algorithms which are relevant for the mathematical proofs and which are presented in pseudocode format. The details of the algorithms involved in the numerical simulations and to create the pictures are not provided. The programs for simulating bond and site percolation configurations are written in C and in C++. The programs for intensive simulations, typically for estimating the curve of the function $\theta(p, \mathbb{Z}^3)$, are written in C and they rely on the GNU Scientific Library. They were run on the computer cluster of the mathematics department of Orsay's University. The programs for illustrative simulations, typically for obtaining two-dimensional pictures in squares $1024 \times 1024$, are written in C++20 and they rely on the Standard Template Library. Modern computers can easily process percolation configurations in larger squares, say $10000 \times 10000$; however there is not much gain in the visual quality of the final image. So, in order to limit the size of the final pdf file, we restrict ourselves to simulations in squares of size $1024 \times 1024$ or rectangles $1550 \times 1024$. To plot the experimental curves, we use the Gnuplot program. Some of the simulation pictures have been enhanced with the help of the GNU Image Manipulation program, typically to enhance the colors. A small number of pictures had to undergo more sophisticated treatment, because some important details were no longer visible when the percolation picture was included and slightly reduced in the pdf file. To compensate for this problem, we applied a standard operator of mathematical morphology, namely the dilation by a small square. For this, and for other simple format conversions, we used the ImageMagick program. We reworked Python code initially generated with the help of the AIs Claude, ChatGPT, Le Chat to obtain pictures of the hypercube graphs and the representations of the examples of monotone subsets and their boundaries (the complicated examples themselves were obtained with the help of C++ programs developed manually). We relied heavily on the computer to illustrate the text, and to test several conjectures on increasing subsets. To some extent, the proof presented here has thus been partially computer-assisted, but this was only through simulations and numerical computations. Neither AI nor machine calculation has been directly involved in the mathematical part of the proof.

# List of Tables

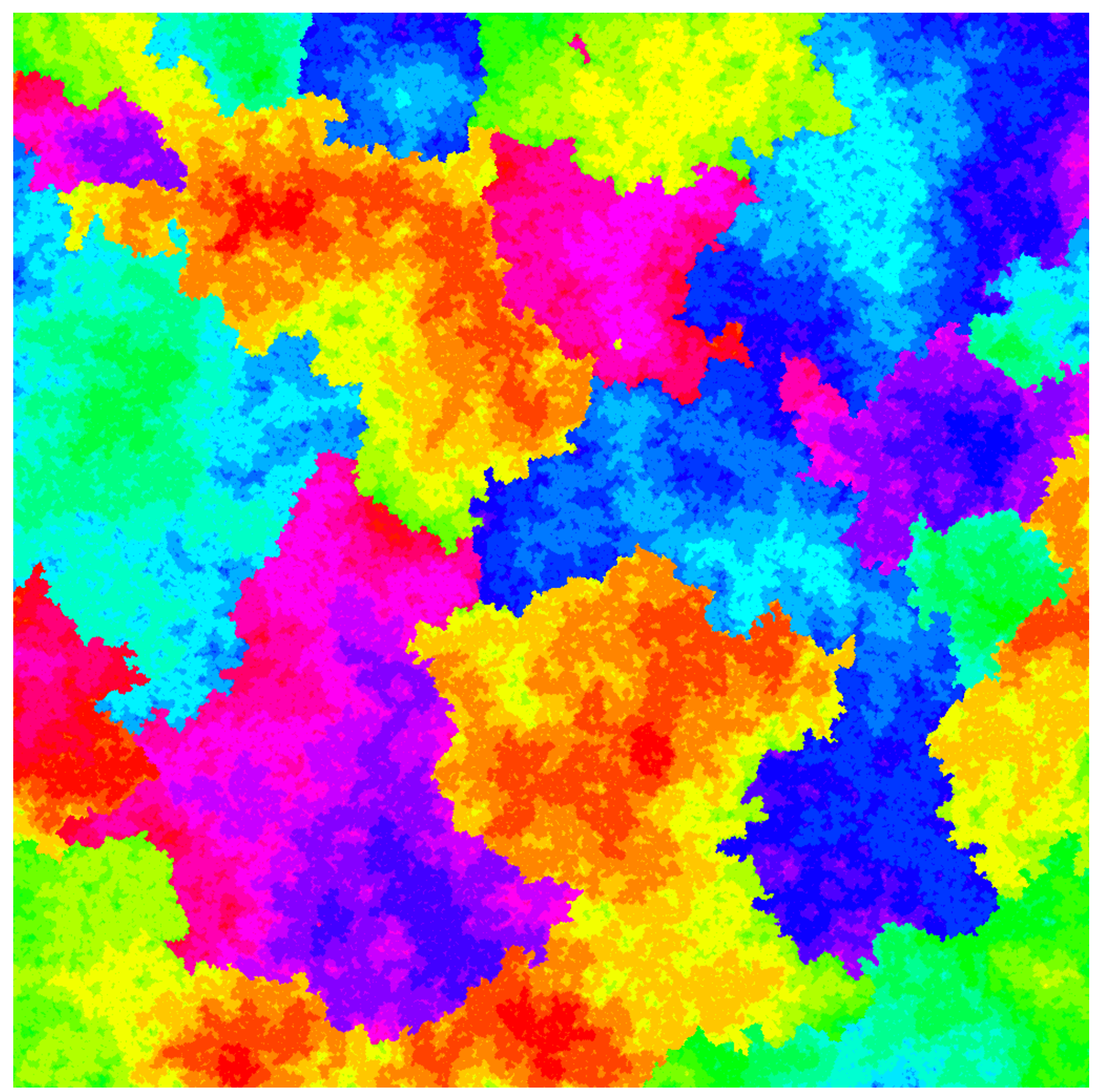

Figure 3: 25 intertwined explorations, $\Lambda(1024)$, site percolation, $p = 0.57$.

# List of Figures

Most of the figures were done in the old-fashioned way, except for a few. For the simulation pictures, a png file was generated by a handwritten C++ program. For the diagrams, we used mainly the pstricks package. For the graphs, we used the Gnuplot program to plot the data created by the C++ programs. All of this was done before the advent of LLMs and AI. We used AI to create the last pictures. Namely, we used Claude, ChatGPT and Le Chat to write the initial version of Python programs to generate the configurations of the Erdős-Rényi model appearing in part VI, as well as the pictures representing the hypercube and the examples of monotone sets appearing in part VIII.

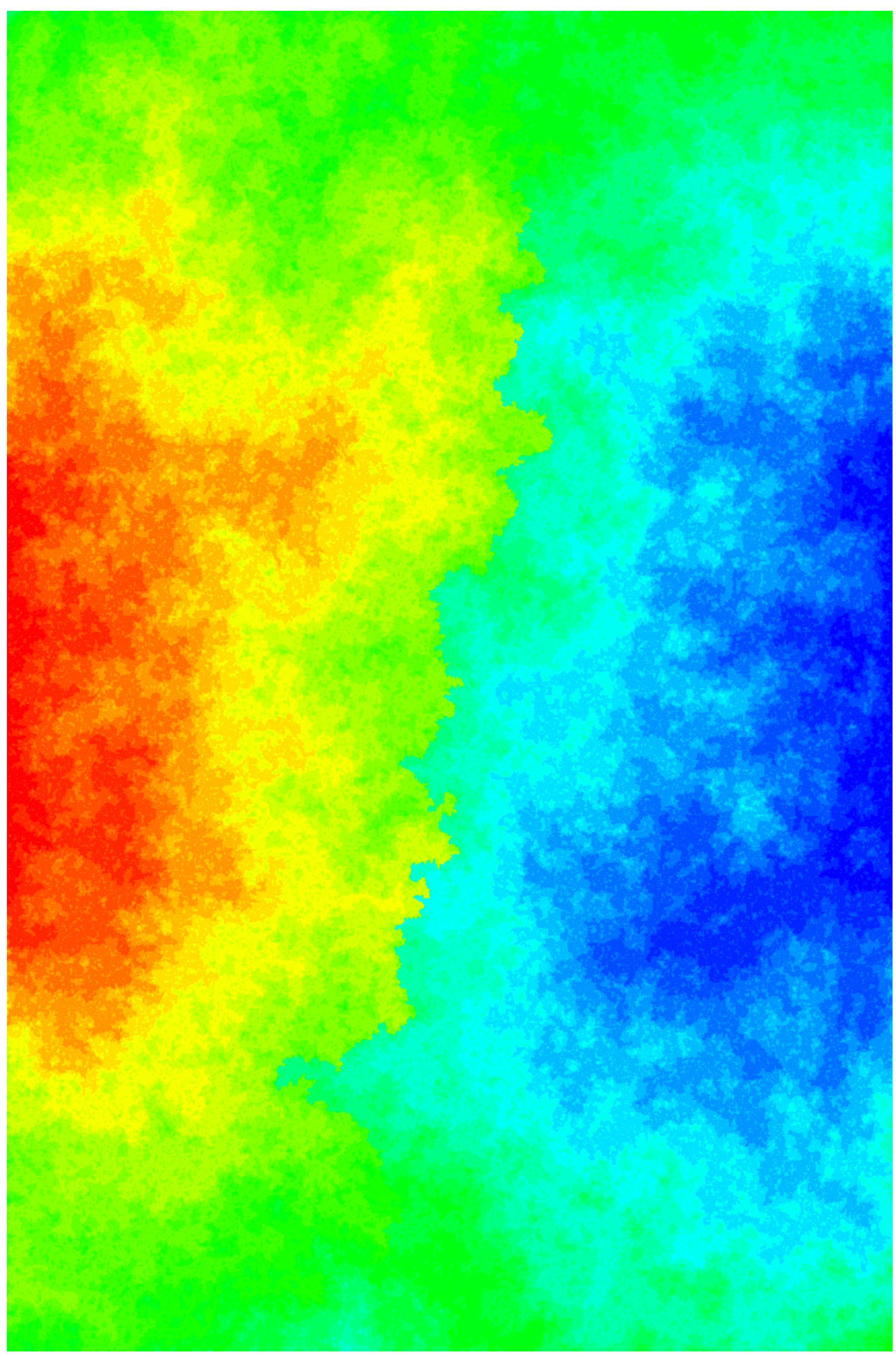

Figure 4: The heart of the matter: disconnection between the middle thirds of left and right sides in the rectangle $1550 \times 1024$, bond percolation, $p = 0.48$.

# Part I
# Preliminaries

In this preliminary part, we introduce the site and the bond percolation models in sections 3 and 4. We discuss various ways to simulate accurately the percolation probability $\theta(p,\mathbb{Z}^d)$ and we show numerical simulations supporting the conjecture $\theta(p_c,\mathbb{Z}^d)=0$ in section 5. We present two attempts for proving the conjecture. The first one, in section 6, is through an extension in three dimensions of a classical two-dimensional argument due to Zhang. The second one, in section 7, is a natural try to control the variance of the quantity $\rho(\Lambda)$, which plays a central role in the simulations. The most promising strategy for proving the conjecture $\theta(p_c,\mathbb{Z}^d)=0$ has been known for quite a long time, it is based on a static renormalization procedure. We present it in section 8.

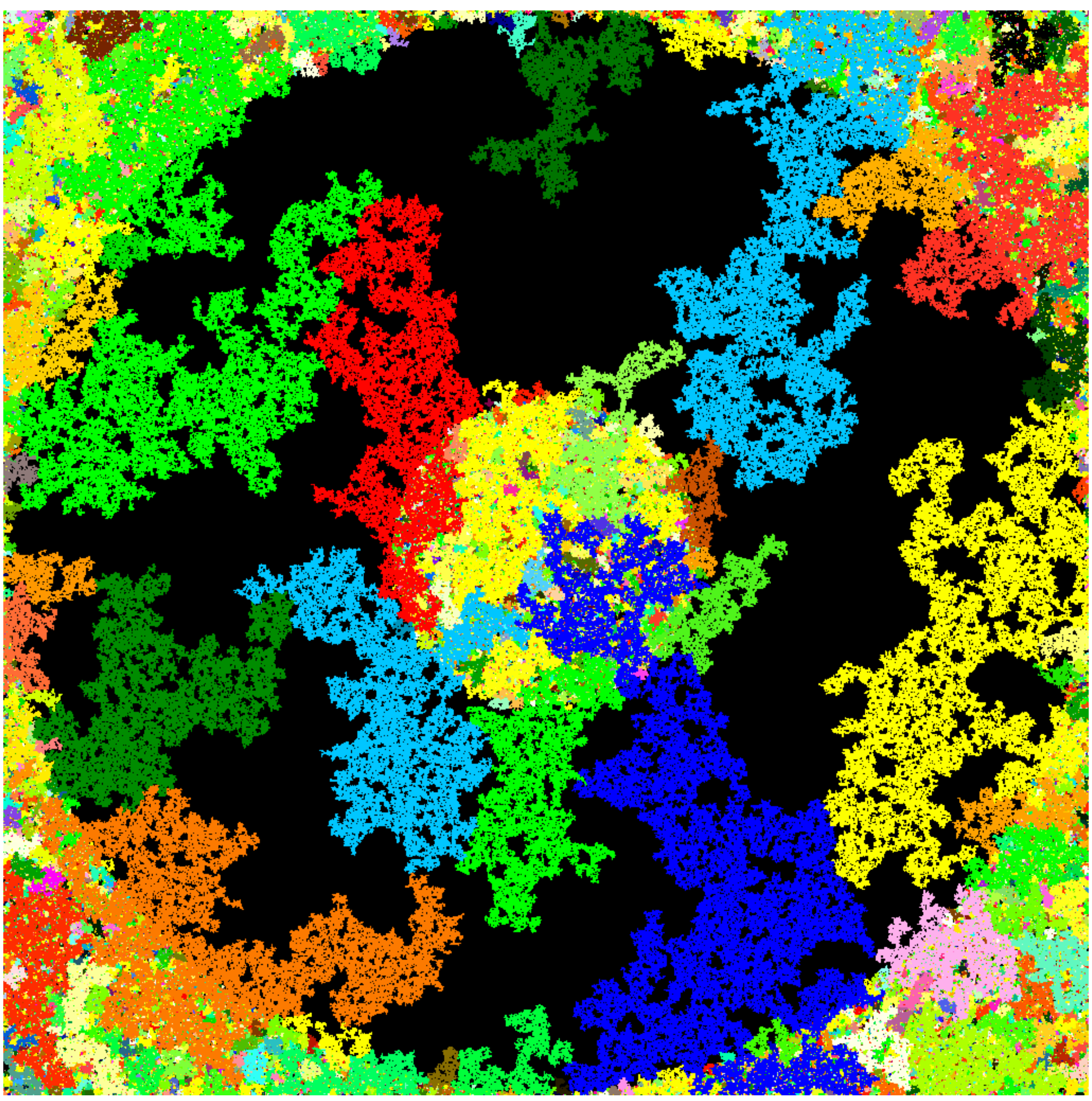

Figure 5: Bond percolation, p=0.495, 1024x1024 box.
The clusters touching the inner ball or the outside of the large ball are colored.

## 3 Bernoulli site percolation in $\mathbb{Z}^d$ ⊛ ○

We gather here the basic standard notation describing the objects involved in the Bernoulli site percolation model in $\mathbb{Z}^d$. More advanced detailed notations will be provided in section 8.

**The $d$-dimensional cubic lattice**. We turn $\mathbb{Z}^d$ into a graph with vertex set $\mathbb{Z}^d$ and edge set

$$\mathbb{E}^d = \big\{ \{x,y\} : x \in \mathbb{Z}^d,\, y \in \mathbb{Z}^d,\, |x-y| = 1 \big\}\,,$$

where $|\cdot|$ is the usual Euclidean norm. This graph is called the $d$-dimensional cubic lattice. A path $\gamma$ in $(\mathbb{Z}^d, \mathbb{E}^d)$ is an alternating sequence

$$x_0, e_0, x_1, e_1, \cdots, e_{n-1}, x_n, \cdots$$

of distinct vertices $x_i$ and edges $e_i$, where $e_i$ is the edge between $x_i$ and $x_{i+1}$. The path is said to connect every pair of its vertices. If the path terminates at some vertex $x_n$ it is said to have length $n$, otherwise it is infinite.

**The probability space.** Let $p$ be a parameter in $[0,1]$. The nearest neighbour Bernoulli site percolation model on the cubic lattice at density $p$ is defined by independently choosing each site of $\mathbb{Z}^d$ to be either open with probability $p$ or closed with probability $1-p$. The configuration space $\Omega$ consists of all the functions from $\mathbb{Z}^d$ to $\{0,1\}$ (1 stands for open and 0 for closed). We equip $\Omega$ with the product $\sigma$-field $\mathcal{P}(\{0,1\})^{\otimes \mathbb{Z}^d}$ (where $\mathcal{P}(\{0,1\}) = \{\varnothing, \{0\}, \{1\}, \{0,1\}\}$ is the collection of the subsets of $\{0,1\}$). We denote by $P_p$, or simply $P$, the Bernoulli product probability measure on $\Omega$ with parameter $p$.

**The open clusters.** Let $\omega$ be a configuration in $\Omega$. We consider the random graph having as vertices the sites of $\mathbb{Z}^d$ which are open in $\omega$ only and edge set the edges whose both endpoints are open in $\omega$ only; equivalently, we remove the sites which are closed in $\omega$ as well as the edges emanating from them. The open clusters in $\omega$ are the connected components of this random graph. We denote by $C(0)$ the open cluster containing the origin (if 0 is closed, then $C(0) = \varnothing$). Two subsets $A, B$ of $\mathbb{Z}^d$ are said to be connected in the configuration $\omega$ if there is a path connecting a site of $A$ to a site of $B$ whose sites are all open in $\omega$; we denote this event by $\{\, A \longleftrightarrow B \,\}$.

**Order on $\Omega$.** There is a natural order on $\Omega$ defined by the following relation: $\omega_1 \leq \omega_2$ if and only if all the sites open in $\omega_1$ are open in $\omega_2$. An event $E$ is said to be increasing (respectively decreasing) if its characteristic function $1_E$ is non-decreasing (respectively non-increasing) with respect to this partial order.

**Cubic boxes.** The box of side length $n \in \mathbb{N}$ centered at the origin is the set

$$\Lambda(n) \,=\, \big\{ (x_1, \cdots, x_d) \in \mathbb{R}^d : \forall i \in \{1, \dots, d\} \;\; -n/2 \leq x_i \leq n/2 \big\}\,.$$

For $A$ a subset of $\mathbb{Z}^d$, we denote by $\partial^{\,in} A$ the sites of $A$ which have a nearest neighbour in $A^c$, i.e.,

$$\partial^{\,in} A \,=\, \big\{ x \in A : \exists y \in A^c \quad |x-y| = 1 \big\}\,.$$

# 4 Bernoulli bond percolation in $\mathbb{Z}^d$ ⊛ $\times\!\!\longrightarrow\!\!\times$

This model is defined on the $d$-dimensional cubic lattice, that is the graph with vertex set $\mathbb{Z}^d$ and edge set $\mathbb{E}^d$ introduced in section 3. We use the same definition for paths as for the site model. The structural difference between the site and the bond models is that open or closed states are assigned to the vertices in the site model and to the edges in the bond model. As a consequence, the mathematical constructions of the two models differ in the definitions of the probability space and of the open clusters. For the bond percolation model, they are as follows.

**The probability space.** The nearest neighbour Bernoulli bond percolation model on the cubic lattice at density $p$ is defined by independently choosing each bond of $\mathbb{E}^d$ to be either open with probability $p$ or closed with probability $1-p$. The configuration space $\Omega$ consists of all the functions from $\mathbb{E}^d$ to $\{0,1\}$ (1 stands for open and 0 for closed). We equip $\Omega$ with the product $\sigma$-field $\mathcal{P}(\{0,1\})^{\otimes \mathbb{E}^d}$ (where $\mathcal{P}(\{0,1\}) = \{\varnothing, \{0\}, \{1\}, \{0,1\}\}$). The space $\Omega$ is endowed with the usual product order. We denote by $P_p$ the Bernoulli product probability measure on $\Omega$ with parameter $p$.

**The open clusters.** Let $\omega$ be a configuration. We consider the random graph having as vertices the set $\mathbb{Z}^d$ and edge set the bonds of $\mathbb{E}^d$ which are open in $\omega$ only; equivalently, we remove the bonds which are closed in $\omega$. The open clusters in $\omega$ are the connected components of this random graph. We denote by $C(0)$ the open cluster containing the origin. Often we identify an open cluster with its vertex set, so that $C(0)$ can be seen as a subset of $\mathbb{Z}^d$. For the site model, the cluster $C(0)$ is empty in case 0 is closed, while for the bond model the cluster $C(0)$ always contains 0. Two subsets $A, B$ of $\mathbb{Z}^d$ are said to be connected in the configuration $\omega$ if there is a path connecting a vertex of $A$ to a vertex of $B$ whose edges are all open in $\omega$; we denote this event by $\{ A \longleftrightarrow B \}$.

**Terminology.** When dealing with the cubic lattice $(\mathbb{Z}^d, \mathbb{E}^d)$ without making reference to a specific model, the elements of $\mathbb{Z}^d$ are called the vertices and the elements of $\mathbb{E}^d$ are called the edges. In percolation, it is customary to call the vertices of $\mathbb{Z}^d$ sites when working with the site percolation model, and to call the edges of $\mathbb{E}^d$ bonds when working with the bond percolation model. This convention will be adhered to as much as possible, although some degree of inconsistency is inevitable.

**Outer boundaries.** The geometric definitions associated to the cubic lattice $(\mathbb{Z}^d, \mathbb{E}^d)$ can interchangeably be used for both the site model and the bond model. For $A$ a subset of $\mathbb{Z}^d$, we define its outer vertex boundary $\partial^{\,out} A$ as the set of the vertices of $A^c$ which have a nearest neighbour in $A$, i.e.,

$$\partial^{\,out} A \,=\, \big\{\, x \in A^c : \exists\, y \in A \quad |x-y| = 1 \,\big\} \,. \tag{4.1}$$

We define also the edge boundary $\Delta A$ of $A$ as the set of the edges which have one endpoint in $A$ and the other in $A^c$, i.e.,

$$\Delta A \,=\, \big\{\, e = \langle x, y \rangle \in \mathbb{E}^d : x \in A, y \in A^c \,\big\} \,.$$

# 5 Numerical simulations ⊛

In order to support various attempts to prove that $\theta(p_c, \mathbb{Z}^3) = 0$, we implemented several simulations of the Bernoulli percolation model in three dimensions. One major goal of these simulations was to obtain an experimental curve of the function $p \in [0,1] \mapsto \theta(p, \mathbb{Z}^3)$ with the best possible precision. Yet the computation of an approximate value of $\theta(p, \mathbb{Z}^3)$ for a given value of $p$ is already a challenging problem. Since, by definition,

$$\theta(p, \mathbb{Z}^3) \,=\, P_p\big(0 \longleftrightarrow \infty\big)\,,$$

the most obvious strategy consists in generating a percolation configuration with parameter $p$ and then exploring the cluster of the origin in order to check whether it is infinite or not. If this cluster is infinite, then the exploration will never stop, so this approach will involve some stopping rule to decide whether the cluster is likely to be infinite or not. However, a standard exploration process stores the whole set of explored sites and this set becomes huge very quickly. In practice, the algorithm would have to stop as soon as the memory of the computer is filled, thereby limiting the size of the region that can be explored, and the stopping rule would produce a poor estimate. Another problem is that huge memory structures lead to the deterioration of the performances, because it takes time to access a specific element in a gigantic array. As is customary in intensive numerical computations, we have to struggle with two limiting factors: the memory and the computation time. So we should try to design an exploration algorithm for three-dimensional clusters which requires a reduced storage. The most natural algorithm exploring a percolation cluster is the following.

**Standard exploration algorithm**. The algorithm uses a sequence of i.i.d. random variables $(U_k)_{k\in\mathbb{N}}$ which follow the Bernoulli distribution with parameter $p$. At each iteration, the algorithm updates three sets of sites:

- The set $A_k$: these are the active sites, which are to be explored.
- The set $O_k$: these are the open sites which have been explored.
- The set $C_k$: these are the closed sites which have been explored.

Initially, we set $O_0 = C_0 = \varnothing$ and $A_0 = \{\,0\,\}$. Suppose that the sets $A_k, O_k, C_k$ are built and let us explain how to build the sets $A_{k+1}, O_{k+1}, C_{k+1}$. If $A_k = \varnothing$, the algorithm terminates and the cluster of the origin is $C(0) = O_k$ and its exterior boundary is $C_k$. If $A_k$ is not empty, we pick up an element $X_k$ of $A_k$. The site $X_k$ has not been explored previously, and its state will be decided by the random variable $U_k$. We consider two cases, according to the value of $U_k$:

- $U_k = 0$. The site $X_k$ is declared closed, and we set

$$A_{k+1} \,=\, A_k \setminus \{\, X_k \,\}\,, \quad O_{k+1} \,=\, O_k\,, \quad C_{k+1} \,=\, C_k \cup \{\, X_k \,\}\,.$$

- $U_k = 1$. The site $X_k$ is declared open, and we set

$$O_{k+1} \,=\, O_k \cup \{\, X_k \,\}\,, \quad C_{k+1} \,=\, C_k\,,$$
$$A_{k+1} \,=\, A_k \cup \mathcal{N}(X_k) \setminus \big(\{\, X_k \,\} \cup O_k \cup C_k\big)\,,$$

where $\mathcal{N}(X_k)$ is the set of the sites of $\mathbb{Z}^3$ which are neighbours of $X_k$.

**Algorithm 5.1** Standard exploration algorithm

```
A(0) ← { 0 }, O(0) ← ∅, C(0) ← ∅, k ← 0
repeat
    Pick up X_k ∈ A(k)
    Draw U_k a random Bernoulli variable with parameter p
    if U_k = 0 then
        X_k is closed
        A(k+1) ← A(k) \ {X_k}
        O(k+1) ← O(k)
        C(k+1) ← C(k) ∪ {X_k}
    else if U_k = 1 then
        X_k is open
        A(k+1) ← A(k) ∪ N(X_k) \ ({X_k} ∪ O(k) ∪ C(k))
        O(k+1) ← O(k) ∪ {X_k}
        C(k+1) ← C(k)
    end if
    k ← k+1
until A(k) = ∅
return O(k), C(k)
```

The pseudocode of this exploration algorithm is given in Algorithm 5.1. This algorithm is quite convenient to obtain nice pictures of percolation. For instance, one can pick up a starting point at random, explore its cluster, and repeat indefinitely the operation. This is the procedure used to generate figure 6 (a critical two-dimensional configuration) and figure 7 (the nasty clusters which cross an annulus). The problem with this algorithm is that it stores the full set of the already explored sites. However, at a given step $k \geq 1$, the relevant information to proceed through the next step is the set of the active sites $A_k$ and the inner boundary of the set of the already explored sites $\partial^{\,in}(O_k \cup C_k)$. Indeed, when we discover a new open site $X_k$, the new sites which become active are the neighbors of $X_k$ which have not been explored previously and these are the sites in $\mathcal{N}(X_k) \setminus \big(\{X_k\} \cup O(k) \cup C(k)\big)$. Yet the site $X_k$ belongs to $A_k$ and

$$A_k \subset \Big\{\, x \in \mathbb{Z}^d \setminus (O_k \cup C_k) : \exists y \in O_k \cup C_k \quad |x-y| = 1 \,\Big\}\,.$$

This yields the inclusion

$$\mathcal{N}(X_k) \setminus \big(\{X_k\} \cup O(k) \cup C(k)\big) \,\subset\, \mathcal{N}(X_k) \setminus \big(\{X_k\} \cup \partial^{\,in}(O_k \cup C_k)\big)\,.$$

The key point is to prevent the exploration process from reentering the set of already explored sites. However the entrance into this set can occur solely through a boundary site, so we need only to store the inner boundary of the already explored region. Unfortunately, it turns out that, for slightly supercritical clusters, the size of this boundary is still of volume order (see for instance theorem 8.99 in [63]). So we should seek a better strategy to estimate $\theta(p)$.

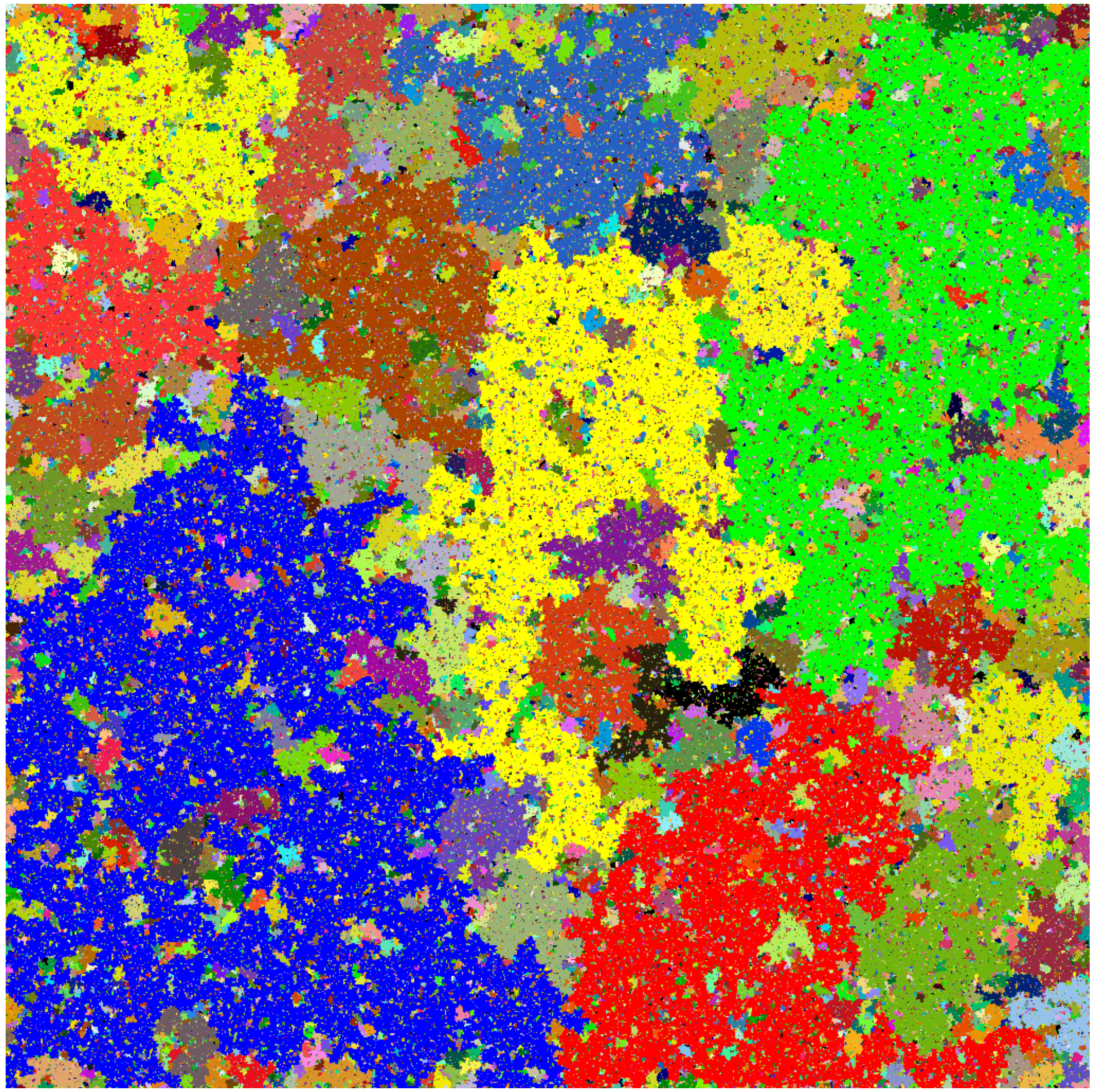

Figure 6: Critical bond percolation in a box of size $1024 \times 1024$

Here is another possibility. Recall that

$$\theta(p,\mathbb{Z}^3) \,=\, P_p\big(0 \longleftrightarrow \infty\big)\,.$$

We look for a finite volume approximation of this function. The natural candidate which comes to mind is the probability that the origin is connected to the boundary in a large finite box $\Lambda(n)$, that is,

$$\theta_n(p,\mathbb{Z}^3) \,=\, P_p\big(0 \longleftrightarrow \partial^{\,in}\Lambda(n)\big)\,.$$

Indeed, we have

$$\theta(p,\mathbb{Z}^3) \,=\, \lim_{n\to\infty} \theta_n(p,\mathbb{Z}^3)\,.$$

Notice that the map $p\in[0,1]\mapsto \theta_n(p,\mathbb{Z}^3)$ is a polynomial function of $p$, which is in addition non-decreasing. Moreover, the sequence of functions $\big(\theta_n(p,\mathbb{Z}^3)\big)_{n\geq 1}$ is also non-increasing in the sense that

$$\forall n\geq 1 \quad \forall p\in[0,1] \qquad \theta_{n+1}(p,\mathbb{Z}^3) \,\leq\, \theta_n(p,\mathbb{Z}^3)\,.$$

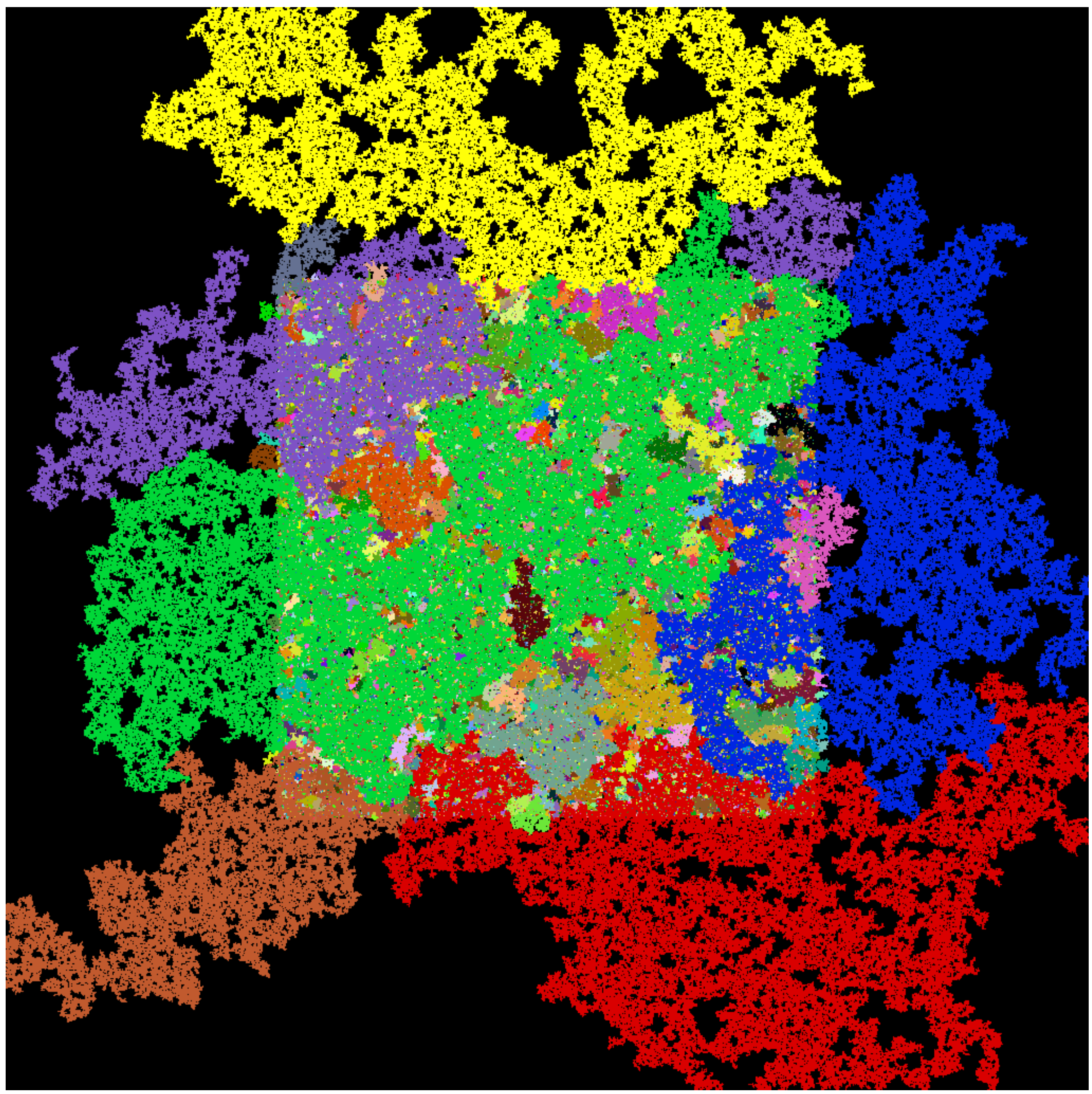

Figure 7: Critical bond percolation in a box of size $1024 \times 1024$
Only the clusters intersecting the central square are colored

It is quite remarkable that the sequence of functions $\big(\theta_n(p,\mathbb{Z}^3)\big)_{n\geq 1}$ satisfies the hypothesis of both Dini's uniform convergence theorems, which we recall next.

**Theorem 5.1.** *If $(f_n)_{n\geq 1}$ is a non-decreasing sequence of continuous functions defined on $[0,1]$ with values in $\mathbb{R}$ which converges pointwise towards a continuous function $f$, then the convergence is uniform.*

**Theorem 5.2.** *If $(f_n)_{n\geq 1}$ is a sequence of non-decreasing functions defined on $[0,1]$ with values in $\mathbb{R}$ which converges pointwise towards a continuous function $f$, then the convergence is uniform.*

To be precise, the sequence $-\theta_n(p,\mathbb{Z}^3)$ satisfies the hypothesis of theorem 5.1, the sequence $\theta_n(p,\mathbb{Z}^3)$ satisfies the hypothesis of theorem 5.2, but our huge problem is to decide whether their pointwise limit $\theta(p,\mathbb{Z}^3)$ is continuous. Each of Dini's theorems leads to the following proposition.

**Proposition 5.3.** *The following two statements are equivalent:*
*i)* $\theta(p_c, \mathbb{Z}^3) = 0$ *,*
*ii) The sequence* $\big(\theta_n(p, \mathbb{Z}^3)\big)_{n\geq 1}$ *converges uniformly towards* $\theta(p, \mathbb{Z}^3)$ *on* $[0, 1]$*.*

*Proof.* To be in position to apply one of the two Dini's theorems, we need to know that the limiting function is continuous. It has been known for a long time that the limiting function $\theta(p, \mathbb{Z}^3)$ is continuous on $]p_c, 1]$ (this is a nice consequence of the uniqueness of the infinite cluster; see theorem 8.1 in [63]). The only troublesome point is the critical point $p_c$. Using the monotonicity properties and the continuity of the functions $\theta_n$, we have

$$\begin{aligned}\theta(p_c, \mathbb{Z}^3) \;=\; \inf_{n\geq 1} \theta_n(p_c, \mathbb{Z}^3) \;&=\; \inf_{n\geq 1} \inf_{p>p_c} \theta_n(p, \mathbb{Z}^3) \\ &=\; \inf_{p>p_c} \inf_{n\geq 1} \theta_n(p, \mathbb{Z}^3) \;=\; \lim_{\substack{p\to p_c \\ p>p_c}} \theta(p, \mathbb{Z}^3)\,.\end{aligned}$$

Obviously, we have $\theta(p, \mathbb{Z}^3) = 0$ for $p < p_c$. Once we know that $\theta(p_c, \mathbb{Z}^3) = 0$, we know indeed that $\theta(p, \mathbb{Z}^3)$ is continuous on $[0, 1]$. □

Unfortunately, it seems very hard to prove the uniform convergence of the sequence of functions $\theta_n(p, \mathbb{Z}^3)$. In any case, this sequence converges pointwise towards $\theta(p, \mathbb{Z}^3)$. So, a natural way to compute numerically $\theta(p, \mathbb{Z}^3)$ consists in simulating a percolation configuration in a box of side length $n$, checking whether 0 is connected to the boundary $\partial^{in}\Lambda(n)$ and repeating this operation a large number of times to estimate the value of $\theta_n(p, \mathbb{Z}^3)$, which itself approximates the value of $\theta(p, \mathbb{Z}^3)$. The advantage of this strategy is that the simulations always take place in a finite box having a fixed side length. Let us discuss which side length might be reachable in practice. For instance, with a memory of 8 gigabytes, that is $8.10^9$ bytes, we can potentially explore a cubic box of side length 2000. This would also require some ingenious encoding of the percolation configuration, because modern computers typically work with variables encoded with 64 bits. Now ingenious encoding implies a slowing down of the algorithm. If we use only simple memory structures to store the percolation configuration, then with 8 gigabytes, we might be able to explore a cubic box of side length at most 1000. This is still a bit optimistic, because, in addition to a cubic box, an efficient labeling algorithm requires an additional memory structure to process the labels, whose number is of the same order as the size of the box. The standard algorithm for labelling the connected components of an array is the two pass algorithm (also called the Hoshen-Kopelman algorithm [72], it was originally designed to study percolation!), it involves a union-find data structure to record equivalent classes of labels. The pseudocode of the two pass algorithm is given in Algorithm 5.2. The complexity of this algorithm is $\Theta(n^3)$, where $n$ is the diameter of the box (the symbol $\Theta(n^3)$ means that there are both upper and lower bounds of order $n^3$). In the end, we might reach cubic boxes of size at most 500. In order to avoid the implementation of the union-find structure and to save some memory, a simple alternative is to use a region growing algorithm, which consists in growing individually each cluster with the help of the standard exploration algorithm. This algorithm uses

**Algorithm 5.2** Labelling connected components in a cube

```
label is a 3D array. A label ℓ has three attributes:
  ℓ.cardinality is the number of sites in the cluster with label ℓ
  ℓ.boundary is true if the cluster with label ℓ intersects the boundary
  ℓ.active is true if the label ℓ is active
closed sites are labelled with the null label 0

for [row,column,height] in label do                    ▷ Initialization of labels
    if uniform random number > p then
        label[row,column,height]=0
    else
        ℓ=new label                                    ▷ Creation of a new label
        ℓ.cardinality=1
        ℓ.active=true
        if [row,column,height] is a boundary site then ℓ.boundary=true
        else  ℓ.boundary=false
        end if
        label[row,column,height]=ℓ
    end if
end for
for [row,column,height] in label do                    ▷ First pass
    ℓ=label[row,column,height]
    if ℓ ≠ 0 then
        L = the non zero labels of the neighbors of [row,column,height]
        for ℓ in L do
            union({ r : r ∼ ℓ },L)                     ▷ Merging of equivalent labels
        end for
    end if
end for
for [row,column,height] in label do                    ▷ Second pass
    ℓ=label[row,column,height]
    r=find(ℓ)              ▷ r is the representative of the labels equivalent to ℓ
    if r ≠ ℓ then
        label[row,column,height]=r
        ℓ.active=false                                 ▷ Label ℓ becomes inactive
    end if
end for
for ℓ active label do                                  ▷ Update label information
    ℓ.cardinality=∑_{r∼ℓ} r.cardinality
    ℓ.boundary=true if there exists r ∼ ℓ such that r.boundary=true
end for
release inactive labels
return the sum of the cardinals of all boundary labels
```

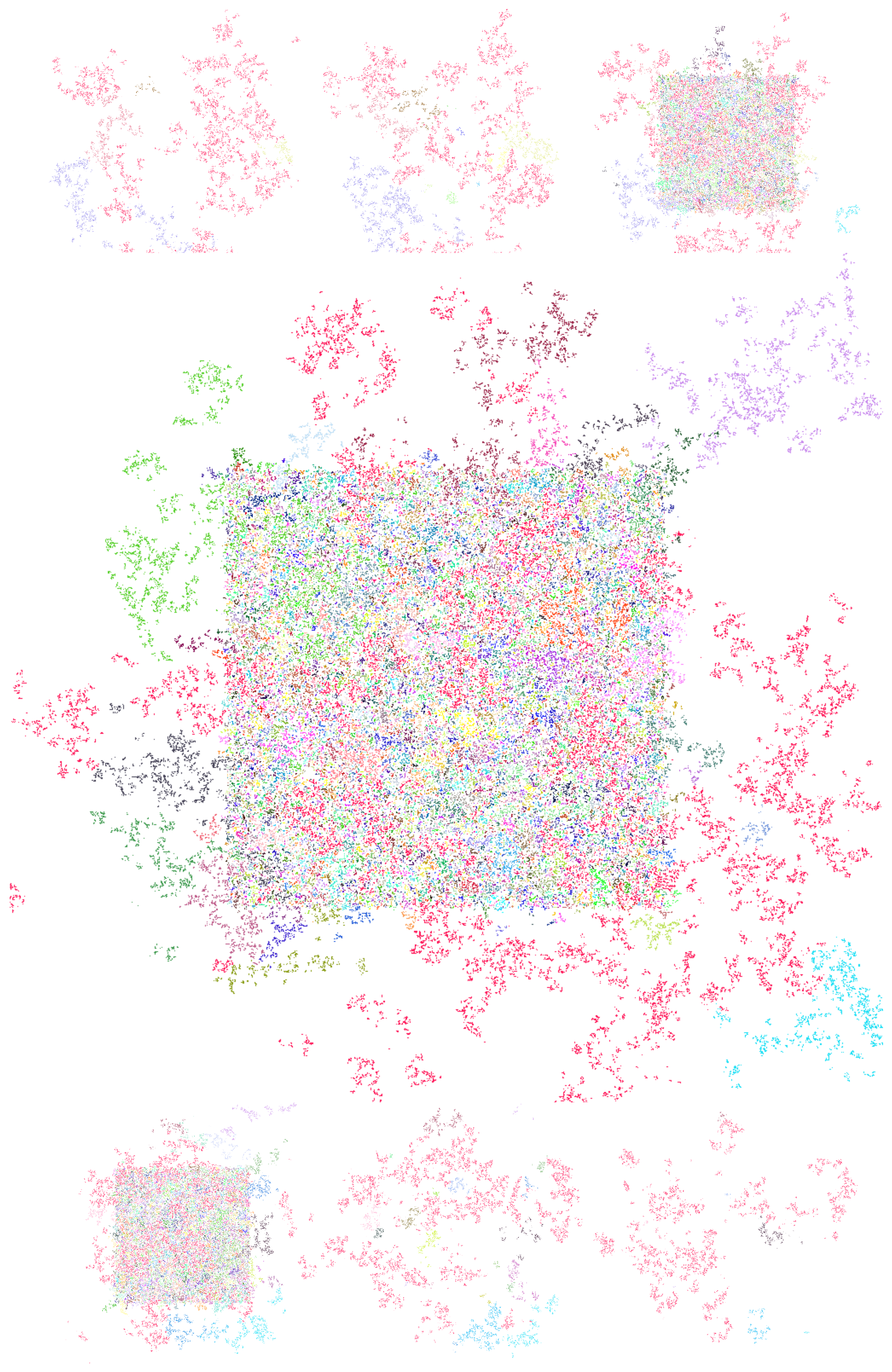

Figure 8: Site percolation in a box of size $1024 \times 1024 \times 1024$, $p = 0.3117$: traces at heights $128, 256, 384, 512, 640, 768, 896$ of the clusters touching the middle box

less memory, but it is computationally less efficient. Yet this is a significant advantage, because the limiting factor for these algorithms in three dimensions is the memory rather than the computation time. We implemented this algorithm in order to visualize directly the structure of a critical three-dimensional percolation configuration, and we could handle easily cubes of side length 1024. In figure 8, we see all the clusters which intersect the middle box of side length 512 in a box of side length 1024, more precisely, we see the intersection of these clusters with the horizontal planes of heights $128, 256, 384, 512, 640, 768, 896$.

There exists another very natural, yet less known estimator for $\theta(p, \mathbb{Z}^3)$. For $A \subset \mathbb{Z}^d$, we denote by $\partial^{\,in} A$ the sites of $A$ which have a nearest neighbour in $\mathbb{Z}^d \setminus A$ and we define

$$\varrho(A) \,=\, \big|\{\, x \in A : x \longleftrightarrow \partial^{\,in} A \,\}\big|\,. \tag{5.1}$$

The map $\varrho$ is subadditive, i.e.,

$$\forall A, B \subset \mathbb{Z}^d, \quad A \cap B = \emptyset \quad \Longrightarrow \quad \varrho(A \cup B) \,\leq\, \varrho(A) + \varrho(B)\,.$$

Indeed, if $x \in A$ is connected by an open path to $\partial^{\,in}(A \cup B)$, then this path must cross $\partial^{\,in} A$, so that $\{\, x \longleftrightarrow \partial^{\,in} A \,\}$ occurs. Given a realization of the percolation configuration on the infinite lattice $\mathbb{Z}^d$, we define the sequence of random variables $\varrho(\Lambda(n)), n \geq 1$, so that we can consider its almost sure convergence after a suitable renormalization.

**Proposition 5.4.** *Let $d \geq 2$. For any $p \in [0,1]$, we have with probability one*

$$\lim_{n \to \infty} \frac{\varrho(\Lambda(n))}{|\Lambda(n)|} \,=\, \theta(p, \mathbb{Z}^d)\,.$$

The multiparameter subadditive ergodic theorem (see [5], theorem 2.4 or [122], theorem 2.1) implies that the above limit exists almost surely. We reproduce below the proof given in [26], Proposition 8.1, which relies only on the spatial ergodic theorem (see for instance [92], section 6.1) to prove the existence of the limit and to identify its value.

*Proof.* Let $m \geq 2n$ be two integers. For $x \in \Lambda(n)$, we have $\Lambda(n) \subset x + \Lambda(m)$, whence

$$\varrho(\Lambda(n)) \,\geq\, \sum_{x \in \Lambda(n)} 1_{\{\, x \longleftrightarrow x + \partial^{\,in} \Lambda(m) \,\}}\,,$$

and, sending $m$ to $\infty$,

$$\frac{\varrho(\Lambda(n))}{|\Lambda(n)|} \,\geq\, \frac{1}{|\Lambda(n)|} \sum_{x \in \Lambda(n)} 1_{\{\, x \longleftrightarrow \infty \,\}} \,=\, \frac{|C_\infty \cap \Lambda(n)|}{|\Lambda(n)|}\,.$$

Sending now $n$ to $\infty$ and applying the ergodic theorem, we obtain

$$\liminf_{n \to \infty} \frac{\varrho(\Lambda(n))}{|\Lambda(n)|} \,\geq\, P(0 \longleftrightarrow \infty)\,.$$

Conversely, for $m \leq n$,

$$\varrho(\Lambda(n)) \,\leq\, |\Lambda(n) \setminus \Lambda(n-m)| + \sum_{x \in \Lambda(n-m)} 1_{\{\, x \longleftrightarrow x + \partial^{\,in}\Lambda(m)\,\}}\,,$$

whence

$$\frac{\varrho(\Lambda(n))}{|\Lambda(n)|} \,\leq\, \frac{|\Lambda(n) \setminus \Lambda(n-m)|}{|\Lambda(n)|} + \frac{1}{|\Lambda(n-m)|} \sum_{x \in \Lambda(n-m)} 1_{\{\, x \longleftrightarrow x + \partial^{\,in}\Lambda(m)\,\}}.$$

We keep $m$ fixed and we send $n$ to $\infty$. The ergodic theorem yields

$$\lim_{n \to \infty} \frac{1}{|\Lambda(n-m)|} \sum_{x \in \Lambda(n-m)} 1_{\{\, x \longleftrightarrow x + \partial^{\,in}\Lambda(m)\,\}} \,=\, P\big(0 \longleftrightarrow \partial^{\,in}\Lambda(m)\big)\,,$$

so that

$$\limsup_{n \to \infty} \frac{\varrho(\Lambda(n))}{|\Lambda(n)|} \,\leq\, P\big(0 \longleftrightarrow \partial^{\,in}\Lambda(m)\big)\,.$$

We conclude by sending $m$ to $\infty$. □

To compute numerically the quantity $\varrho(\Lambda(n))$, the straightforward approach consists in labelling all the clusters in the box $\Lambda(n)$ and adding up the cardinality of the clusters which intersect the boundary of $\Lambda(n)$, as in Algorithm 5.2. This still requires a large memory, and the memory is the limiting factor, so we shall consider a simpler variant of $\varrho(\Lambda(n))$, which leads to a much more efficient algorithm. Let us denote by $L(n)$ the left face of the cubic box $\Lambda(n)$ and let us define

$$\varrho_L(\Lambda(n)) \,=\, \big|\{\, x \in \Lambda(n) : x \longleftrightarrow L(n)\,\}\big|\,.$$

**Proposition 5.5.** *For any $p \in [0,1] \setminus \{\, p_c\,\}$, the quantity $\varrho_L(\Lambda(n))/|\Lambda(n)|$ converges in probability towards $\theta(p, \mathbb{Z}^3)$, that is,*

$$\forall \varepsilon > 0 \qquad \lim_{n \to \infty} P\Big(\Big|\frac{\varrho_L(\Lambda(n))}{|\Lambda(n)|} - \theta(p, \mathbb{Z}^3)\Big| > \varepsilon\Big) \,=\, 0\,.$$

*Proof.* Obviously, we have $\varrho_L(\Lambda(n)) \leq \varrho(\Lambda(n))$. The case $p < p_c$ follows from proposition 5.4. To deal with the case $p > p_c$, we make appeal to a block estimate, for instance Lemma 7.97 in [63]. □

Moreover, speed estimates associated to the convergence result of proposition 5.5 are also available. In the subcritical regime $p < p_c$, the convergence occurs at a speed of order $\exp(-O(n^3))$. This is a consequence of the exponential decay of the cluster size (see lemma 5.1 in [30]). In the supercritical regime $p > p_c$, the convergence occurs at a speed of order $\exp(-O(n^2))$. This is a consequence of the renormalization scheme of Pisztora [114]. As is customary, the situation for $p = p_c$ is more involved and requires a specific treatment. If we suppose that conjecture 1.1 is true, proposition 5.4 readily implies that the convergence occurs as well for $p = p_c$. However, at that point, we have no proof of this result for $p = p_c$.

Without the proof of conjecture 1.1, we have only the following partial result:

$$\theta(p_c,\mathbb{Z}^3)=0 \qquad \Longleftrightarrow \qquad P_{p_c}\Big(\lim_{n\to\infty}\frac{\varrho_L(\Lambda(n))}{|\Lambda(n)|}=0\Big)=1\,. \tag{5.2}$$

The left to right implication is a consequence of proposition 5.4 and the inequality $\varrho_L(\Lambda(n))\leq\varrho(\Lambda(n))$. Let us prove the converse implication. Denoting by $F_i$, $1\leq i\leq 6$, the six faces of $\Lambda(n)$, we have

$$\varrho(\Lambda(n))\leq\sum_{1\leq i\leq 6}\varrho_{F_i}(\Lambda(n))\,,$$

whence, taking the limit, using proposition 5.4 and the right assertion of (5.2),

$$\theta(p_c,\mathbb{Z}^3)=\lim_{n\to\infty}\frac{\varrho(\Lambda(n))}{|\Lambda(n)|}\leq\lim_{n\to\infty}\sum_{1\leq i\leq 6}\frac{\varrho_{F_i}(\Lambda(n))}{|\Lambda(n)|}=0\,.$$

In the end, our preferred method for computing numerical estimates of the function $\theta(p,\mathbb{Z}^3)$ consists in simulating the quantity $\varrho_L(\Lambda(n))$. The main advantage is that there exists an algorithm allowing to compute $\varrho_L(\Lambda(n))$ in a three-dimensional box $\Lambda(n)$ which requires an amount of memory of order $\Theta(n^2)$. Since this algorithm is used only for the numerical simulations and not for the mathematical proofs, we do not give its precise mathematical formulation.

**Simulation of $\varrho_L(\Lambda(n))$.** The algorithm works with two two-dimensional arrays of labels, named square and nextsquare. Its pseudocode is given in Algorithm 5.3. The algorithm starts from the left face $L(n)$ and performs one left to right traversal of the cubic box $\Lambda(n)$. We conducted extensive simulations of the function $\theta(p)$ in three dimensions with the help of this algorithm. We ran it on the computer cluster of the mathematics department of Orsay's University.[1] With a lot of time and numerous simulations, we could explore three-dimensional boxes of side length 10000. The results of these simulations are presented in several graphics. Figure 9 is a synthesis of the whole curve of $\theta$. Looking at this curve, who would believe that $\theta(p)$ does not vanish at $p_c$? In the subsequent figures, we zoom around the critical point. In figure 10, we see monotone simulations, while in figure 11, we see independent simulations. Simulations corresponding to disjoint lines in figure 10 are stochastically coupled in the sense that the random number generator (in this case, Mersenne twister) is initialized with the same seed, only the percolation parameter $p$ is changed from one simulation to another. As a little thought reveals, this does not imply that the portion of the configuration explored by the algorithm is monotone with respect to $p$. This is why some of the lines in figure 10 are not monotone. Finally, figure 12 presents the evolution of interesting dynamical quantities computed during the simulations, namely the numbers of crossing clusters and of sites connected to the opposite face, plotted as functions of the distance to the left face. This amounts to explore the configuration in a cylinder, whose basis is the left face and whose length is the distance to the right face. It would be of crucial importance to obtain a mathematical control of these quantities.

[1] I thank Yves Misiti and Sylvain Faure for helping me to use the computer cluster.

**Algorithm 5.3** Computation of the number of sites connected to the left face

square is a 2D array isomorphic to the left face of the box
$\ell$.cardinal is the number of sites in the cluster with label $\ell$
$\ell$.boundary is true if the cluster with label $\ell$ intersects the left face
$\ell$.active is true if the label $\ell$ is active
closed sites are labelled with the null label 0, square is initialized with 0

no=0 ▷ no is the numero of the square which is explored
**repeat**
no=no+1
**for** [row,col] in square **do** ▷ Transmission of labels to nextsquare
**if** uniform random number $> p$ **then** nextsquare[row,col]=0
**else**
$\ell$=square[row,col]
**if** $\ell \neq 0$ **then**
$\ell$.cardinality+=1 ▷ Transmission of the previous label
**else**
$\ell$=new label ▷ Creation of a new label
$\ell$.cardinality=1, $\ell$.active=true
**if** no=1 **then** $\ell$.boundary=true **else** $\ell$.boundary=false
**end if**
**end if**
nextsquare[row,col]=$\ell$
**end for**
**for** [row,col] in nextsquare **do** ▷ First pass
$\ell$=nextsquare[row,col]
**if** $\ell \neq 0$ **then**
L= the non zero labels of neighbors of [row,col] in nextsquare
**for** $\ell$ in L **do**
union($\{ r : r \sim \ell \}$,L) ▷ Merging of equivalent labels
**end for**
**end if**
**end for**
**for** [row,col] in nextsquare **do** ▷ Second pass
$\ell$=nextsquare[row,col]
r=find($\ell$) ▷ r is the representative of the labels equivalent to $\ell$
**if** $r \neq \ell$ **then** nextsquare[row,col]=r, $\ell$.active=false
**end for**
**for** $\ell$ active label **do** ▷ Update label information
$\ell$.cardinal=$\sum_{r \sim \ell}$ r.cardinal
$\ell$.boundary=true if there exists $r \sim \ell$ such that r.boundary=true
**end for**
release inactive labels
square $\leftarrow$ nextsquare
**until** the box is crossed from left to right or there is no boundary label
**return** the sum of the cardinals of all boundary labels

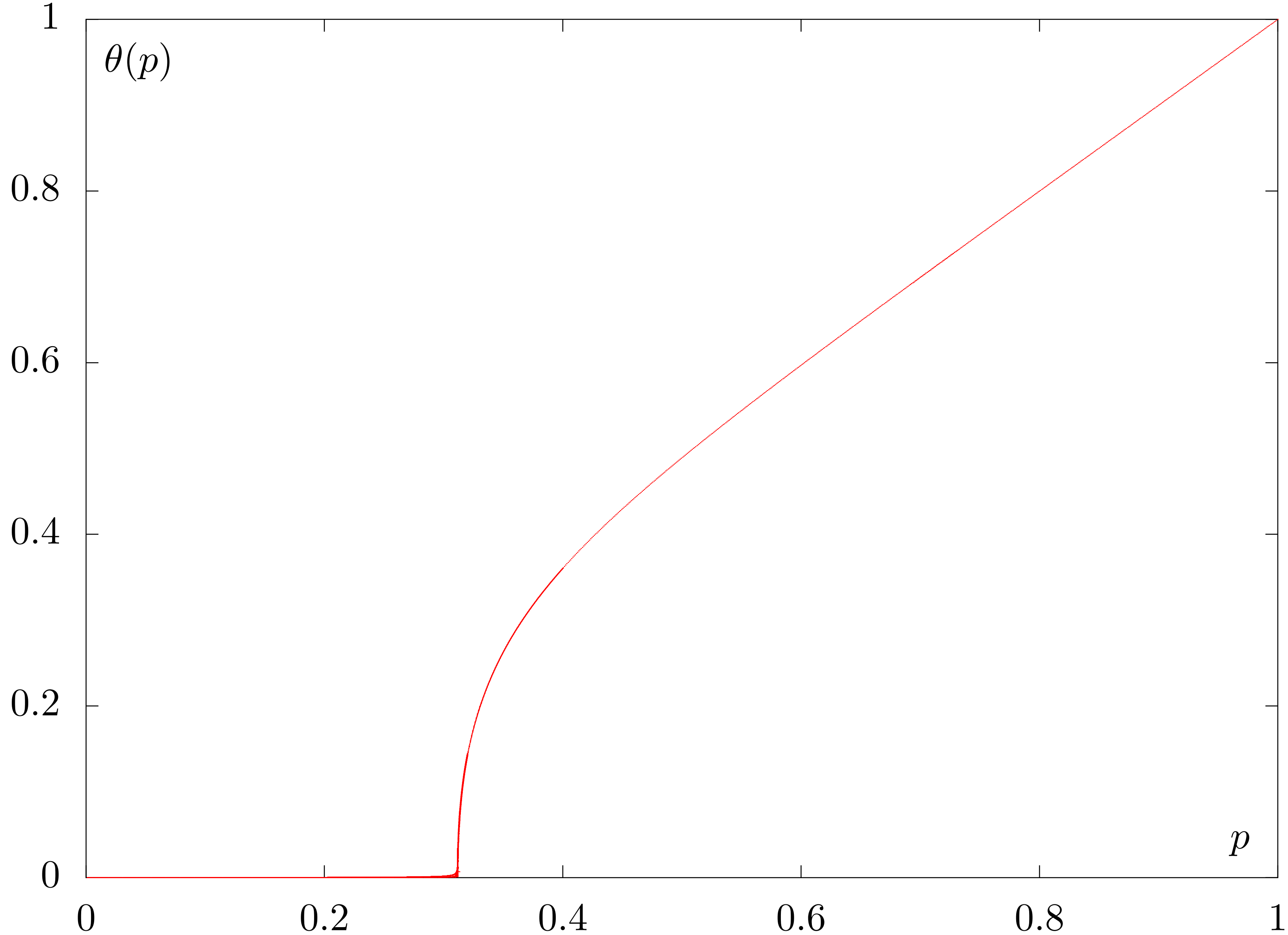

Figure 9: 20000 simulations inside cubic boxes of side length between 1000 and 10000. They were realized between the years 2000 and 2020 on the computer cluster of Orsay's mathematics department. The total computation time is equal to 875 days.

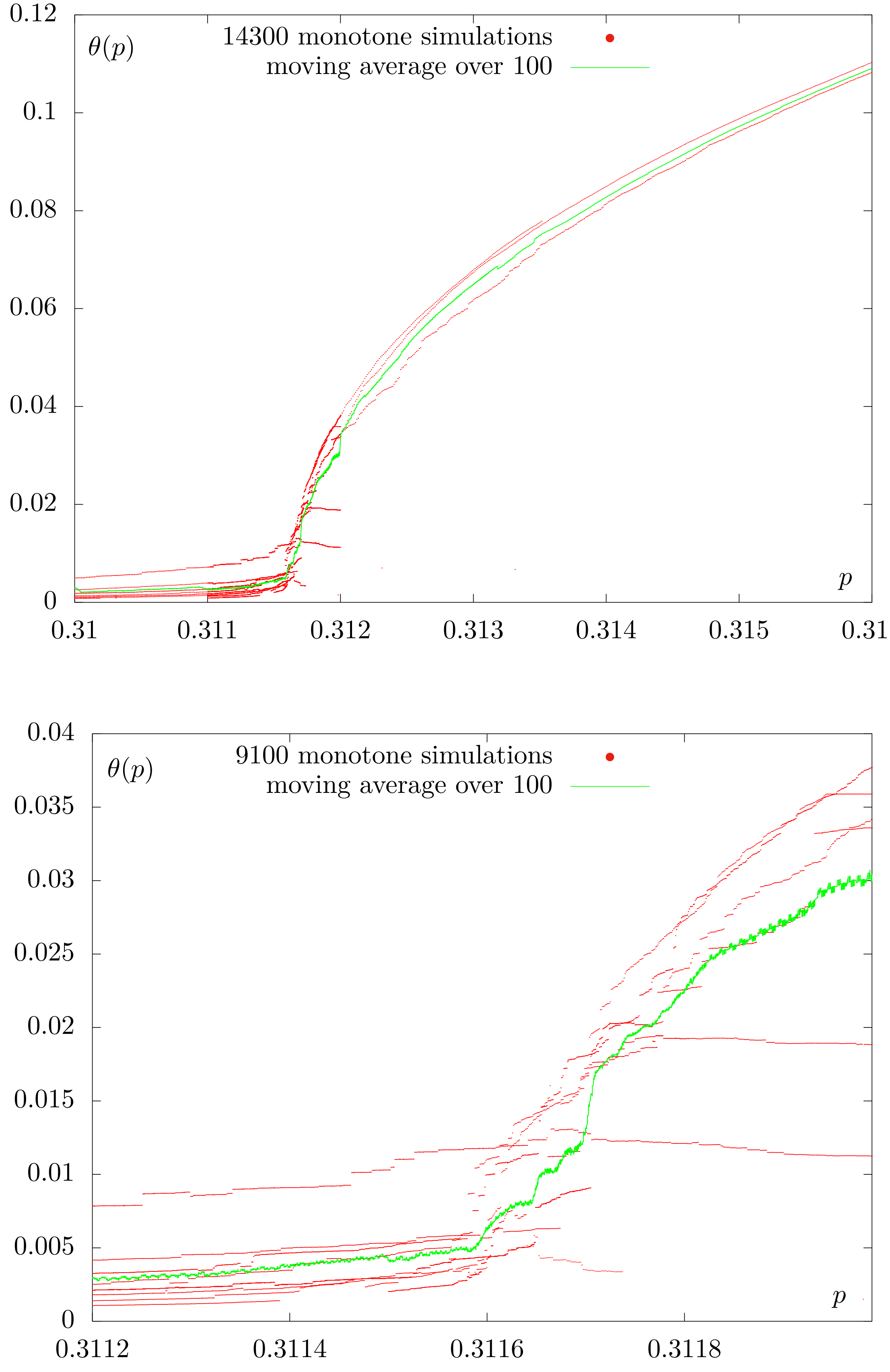

Figure 10: Monotone simulations around the critical point

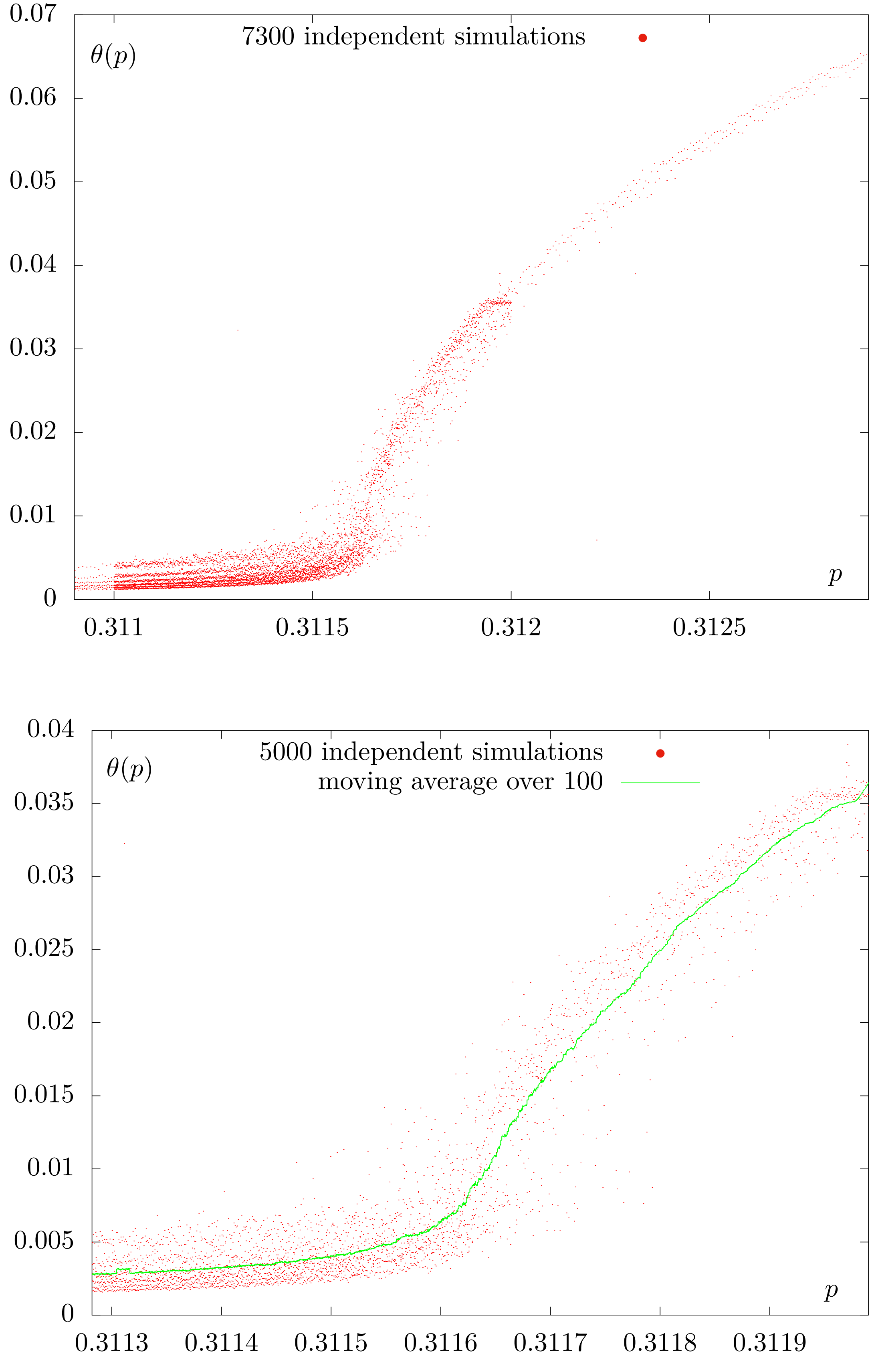

Figure 11: Independent simulations around the critical point

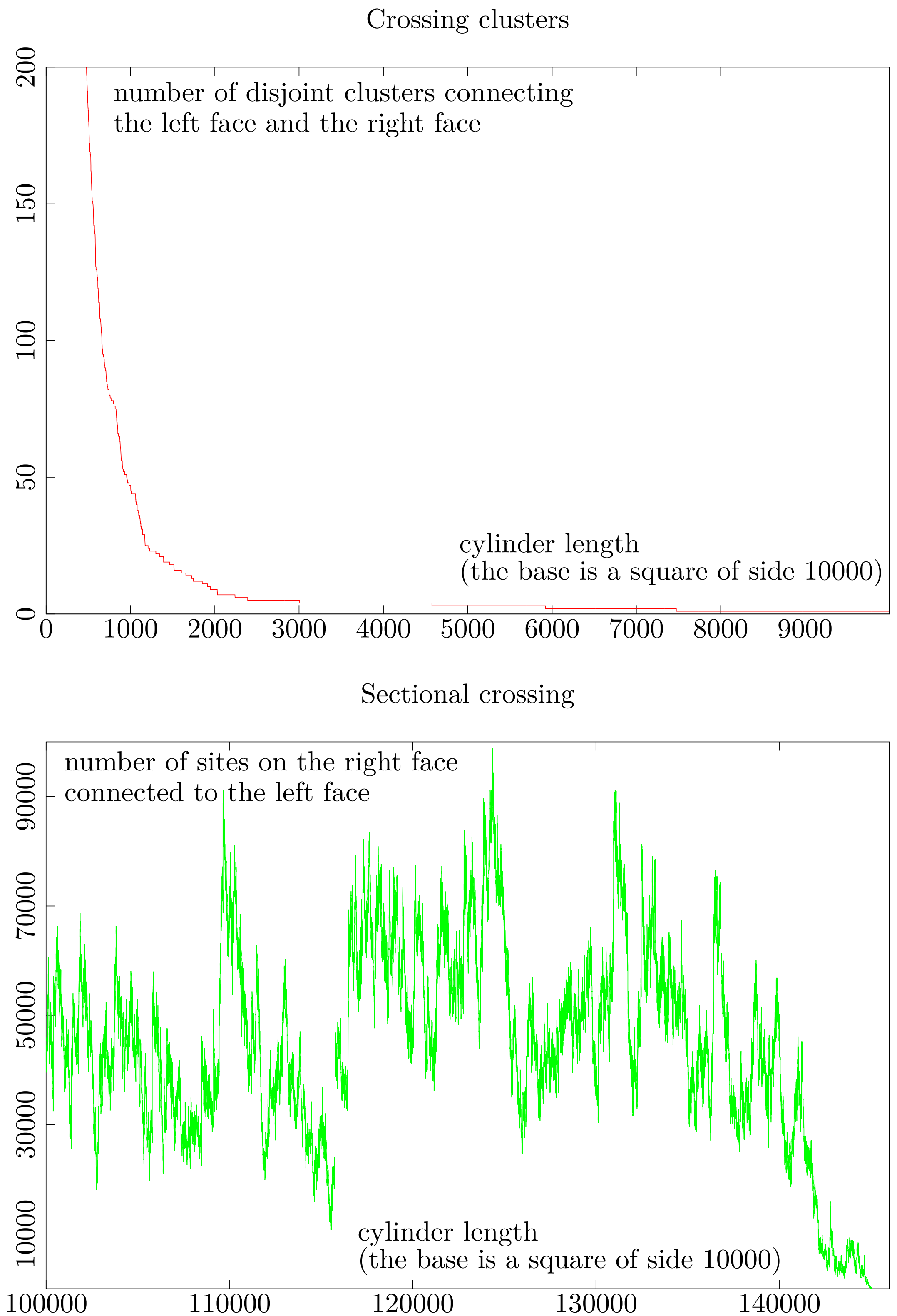

Figure 12: Numbers of crossing clusters and of sites connected to the opposite face as a function of the cylinder length in a cylinder of base 10000x10000 at parameter p=0.311651. All the clusters have merged at length 7474 and they die at length 145040. The simulation took 30 full days of computation time on the computer cluster of Orsay's mathematics department.

# 6 Extension of Zhang's argument ⊛ ○

An important step in Kesten's proof that the critical point of two-dimensional bond percolation is 1/2 is to prove that

$$\theta\big(1/2, \mathbb{Z}^2\big) \,=\, 0\,. \tag{6.1}$$

Zhang found a beautiful argument to prove (6.1), which rests on duality (as quoted in [63], proof of lemma 11.12, this argument does not seem to appear in a published paper). Let us try to extend this argument to higher dimensions. Of course, the percolation model is not self-dual in dimensions $d \geq 3$. Yet it is well-known that there is a form of duality between two different notions of connectedness, namely:
• The usual connectedness, which involves nearest neighbours vertices.
• The $*$-connectedness, which involves second nearest neighbours vertices.
Let us define more precisely these two notions. We denote by $|\cdot|$ the usual Euclidean norm. Two vertices $x, y$ of $\mathbb{Z}^d$ are neighbours if $|x-y| = 1$, and they are $*$-neighbours if $|x-y|_\infty = 1$. These two definitions of neighbours give rise to two different graph structures on $\mathbb{Z}^d$, hence we can consider subsets of $\mathbb{Z}^d$ which are either connected or $*$-connected. Obviously, connectedness in the usual sense implies $*$-connectedness. We define the external $*$-boundary $\partial_*^{\,ext} A$ of a connected subset $A$ of $\mathbb{Z}^d$ as

$$\partial_*^{\,ext} A \,=\, \left\{ \begin{matrix} x \in \mathbb{Z}^d \setminus A : \text{there exists a } *\text{-neighbour } y \text{ of } x \text{ in A and there} \\ \text{exists a connected path of sites in } \mathbb{Z}^d \setminus A \text{ going from } x \text{ to } \infty \end{matrix} \right\}.$$

The set $\partial_*^{\,ext} A$ consists of the vertices in the outer $*$-boundary of $A$ which are visible from $\infty$. Symmetrically, we define the external boundary $\partial^{\,ext} A$ of a $*$-connected subset $A$ of $\mathbb{Z}^d$ as

$$\partial^{\,ext} A \,=\, \left\{ \begin{matrix} x \in \mathbb{Z}^d \setminus A : \text{there exists a neighbour } y \text{ of } x \text{ in A and there} \\ \text{exists a } *\text{-connected path of sites in } \mathbb{Z}^d \setminus A \text{ going from } x \text{ to } \infty \end{matrix} \right\}.$$

Notice that, whenever the set $A$ is infinite, its external boundaries $\partial_*^{\,ext} A$ and $\partial^{\,ext} A$ might be void. However, we will only consider external boundaries of finite sets in the sequel. With a topological argument, Kesten has proved that, if $A$ is finite and $*$-connected, then its external boundary $\partial^{\,ext} A$ is connected (see [83], lemma 2.23). Symmetrically, if $A$ is finite and connected, then its external $*$-boundary $\partial_*^{\,ext} A$ is $*$-connected. This is further proved by Deuschel and Pisztora (see [37], lemma 2.1). Timár [130] has given a simpler proof of these facts using arguments of graph theory. We define now two percolation functions, the usual one associated to the infinite open cluster, given by

$$\theta(p, \mathbb{Z}^d) \,=\, P_p\big(0 \text{ belongs to an infinite open cluster}\big)\,,$$

and a new one associated to the $*$-infinite closed cluster, given by

$$\theta^*(p, \mathbb{Z}^d) \,=\, P_p\big(0 \text{ belongs to an infinite closed } *\text{-cluster}\big)\,.$$

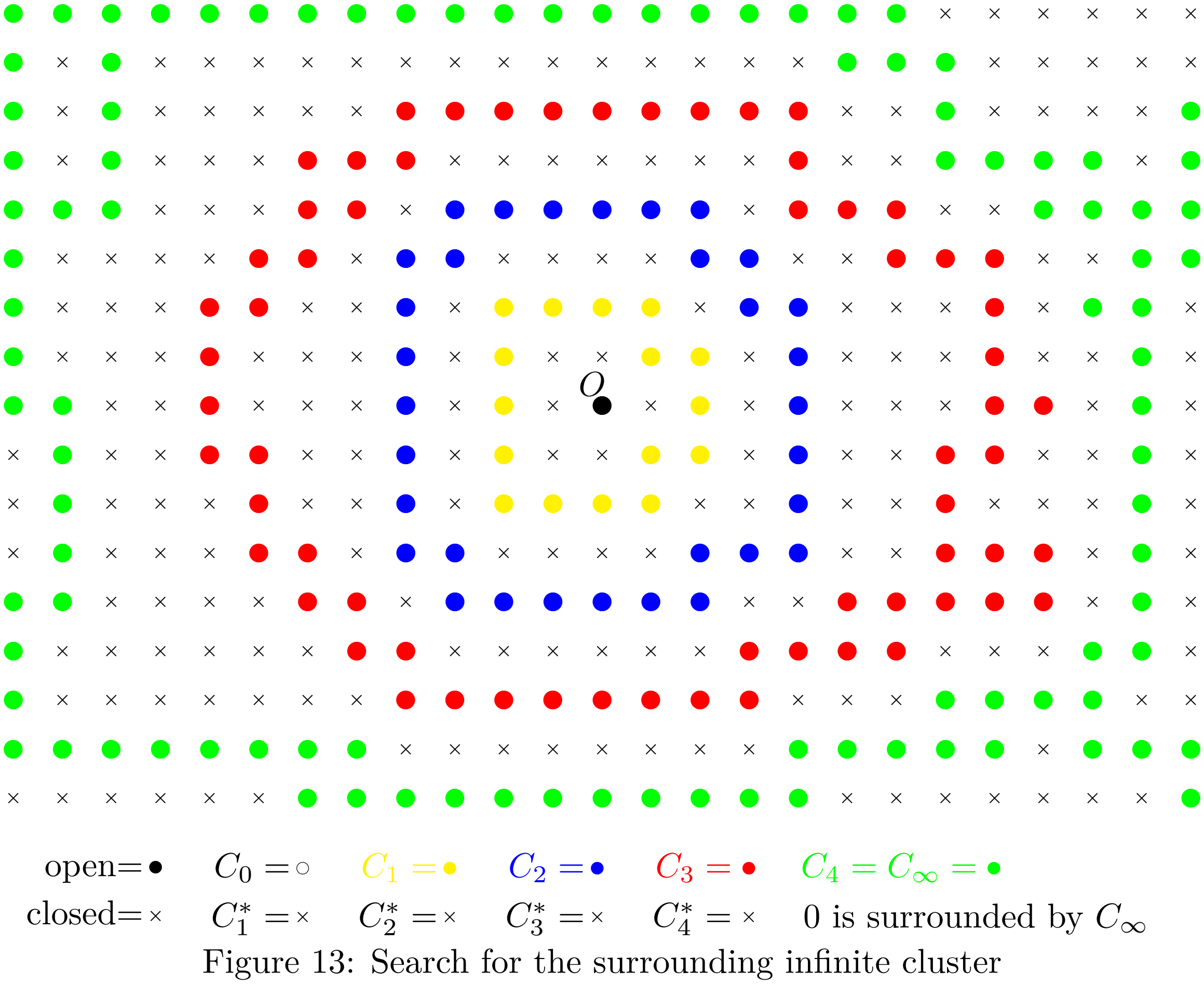

Figure 13: Search for the surrounding infinite cluster

By the classical argument of Burton-Keane (see for instance [63], section 8.2), we know that there exists at most one infinite open cluster. We denote it by $C_\infty$ (if no such cluster exists, we set $C_\infty = \varnothing$). With the same argument, one can prove that, with probability 1, there exists at most one infinite closed $*$-cluster. Let us denote it by $C_\infty^*$. We consider now the following algorithm.

**Algorithm for the surrounding infinite cluster**. We look at the state of 0. If 0 is open, then we build the open cluster $C(0)$ and we set $C_0 = C(0)$. If 0 is closed, then we set $C_0 = \varnothing$, we build the closed $*$-cluster $C^*(0)$ containing 0, and we set $C_1^* = C^*(0)$. For the sake of the discussion, suppose for instance that we are in the first case, that is, the site 0 is open. If $C_0 = C_\infty$, then the algorithm terminates. Otherwise, if $C(0)$ is finite and non-empty, then $\partial_*^{ext} C_0$ is a finite $*$-connected set. By construction, all the sites in $\partial_*^{ext} C_0$ are closed, so we look at the closed $*$-cluster $C_1^*$ containing it. If $C_1^* = C_\infty^*$, then the algorithm terminates. Otherwise, if $C_1^*$ is finite, then it is non-empty and $\partial^{ext} C_1^*$ is a finite connected set. By construction, all the sites in $\partial^{ext} C_1^*$ are open, so we look at the open cluster $C_1$ containing it. The algorithm continues this way until it hits one of the two infinite clusters $C_\infty$ or $C_\infty^*$.

This algorithm terminates after a finite number of steps if and only if at least one of the two clusters $C_\infty$ or $C_\infty^*$ is present in the configuration. Indeed, at step $k$, the algorithm has uncovered a $*$-connected set which surrounds the ball centered at 0 of radius $k$ for the norm $|\cdot|_1$. Necessarily, this sequence of growing regions will eventually meet $C_\infty$ or $C^*_\infty$, if one of them exists, leading to the termination of the algorithm. In two dimensions, the two clusters $C_\infty$ and $C^*_\infty$ never coexist, and none of them exists at criticality, and this fact is the essence of Zhang's argument. In particular, the previous algorithm would never terminate if it is started on a two-dimensional critical configuration. From dimensions 3 onwards, a very different phenomenon occurs. Indeed, there is always one of these two clusters present in the configuration, and moreover there is an intermediate phase where they are both present. More precisely, let us define the two critical points

$$\begin{aligned} p_c &= \sup \left\{ p \in [0,1] : \theta(p, \mathbb{Z}^d) = 0 \right\}, \\ p_c^* &= \inf \left\{ p \in [0,1] : \theta^*(p, \mathbb{Z}^d) = 0 \right\}. \end{aligned}$$

Campanino and Russo [23] have proved the following result.

**Theorem 6.1.** *In dimensions $d \geq 3$, we have* $0 < p_c < 1/2$.

Moreover, the critical point of the $*$-connected open infinite cluster is less or equal than $p_c$. So, theorem 6.1 implies the further inequality

$$0 < p_c < 1/2 < p_c^* < 1$$

(we recall that $p_c^*$ is the critical point of the $*$-connected closed infinite cluster). Therefore, for any value of $p$, with probability 1, at least one of the two clusters $C_\infty$ or $C^*_\infty$ is present in the configuration. This ensures that the previous algorithm terminates with probability 1 for any value of $p$ in dimensions $d \geq 3$. We say that the origin is surrounded by the infinite open cluster $C_\infty$ (respectively by the infinite closed $*$-cluster) if the algorithm terminates with $C_\infty$ (respectively $C^*_\infty$). This could have been defined more formally in the following way. We say that 0 is surrounded by $C_\infty$ (respectively $C^*_\infty$) if it belongs to $C_\infty$ (respectively $C^*_\infty$) or to a bounded $*$-component of $\mathbb{Z}^d \setminus C_\infty$ (respectively a bounded component of $\mathbb{Z}^d \setminus C^*_\infty$). Let us now define the functions

$$\begin{aligned} \rho(p, \mathbb{Z}^d) &= P_p\big(0 \text{ is surrounded by } C_\infty\big), \\ \rho^*(p, \mathbb{Z}^d) &= P_p\big(0 \text{ is surrounded by } C^*_\infty\big). \end{aligned}$$

The termination of the algorithm, or a direct geometrical argument, shows that 0 is surrounded either by $C_\infty$ or by $C^*_\infty$, and the two possibilities are mutually exclusive. Thus these two events form a partition of $\Omega$, and we can conclude safely that

$$\forall p \in [0,1] \qquad \rho(p, \mathbb{Z}^d) + \rho^*(p, \mathbb{Z}^d) = 1. \tag{6.3}$$

Let us next focus on the critical point $p_c$. By theorem 6.1, we know that $p_c < p_c^*$, therefore the point $p_c$ belongs to the supercritical regime $[0, p_c^*[$ of the infinite closed $*$-cluster and thus $\rho^*$ is continuous at $p_c$. The above identity (6.3) implies then that $\rho$ is also continuous at $p_c$, hence $\rho(p_c, \mathbb{Z}^d) = 0$. Needless to say, we have that $\rho(p, \mathbb{Z}^d) \geq \theta(p, \mathbb{Z}^d)$, hence $\theta(p_c, \mathbb{Z}^d) = 0$!
This is my most efficient proof [2] [3] of $\theta(p_c, \mathbb{Z}^d) = 0$.

[2] Unfortunately, it is wrong. Can you find the problem?

[3] Over the years, I have tested this argument with my fellow percolationists. Many of them, experts in the field, remained confused for quite a long time.

## 7 The variance of $\rho(\Lambda)$ $\circledast$ $\times\!\!\longleftrightarrow\!\!\times$

We introduced the quantity $\varrho$ in (5.1) in order to simulate efficiently the function $\theta$. Recall that, for $\Lambda$ a box included in $\mathbb{Z}^d$, we have

$$\varrho(\Lambda) \,=\, \big|\{\, x\in\Lambda : x \longleftrightarrow \partial^{\,in}\Lambda \,\}\big|\,.$$

In proposition 5.4, we proved a law of large numbers for $\varrho(\Lambda)$ renormalized by $|\Lambda|$, the limit being precisely $\theta(p,\mathbb{Z}^d)$. A major problem is to precise this convergence. In fact, the control of this convergence is a key point for proving that $\theta(p_c,\mathbb{Z}^d)=0$. So, let us try to estimate the variance of $\varrho(\Lambda)$. We write

$$\varrho(\Lambda) \,=\, \sum_{x\in\Lambda} 1_{\{\, x \longleftrightarrow \partial^{\,in}\Lambda \,\}}\,,$$

and we compute successively

$$\begin{aligned}
E(\varrho(\Lambda)) \,&=\, \sum_{x\in\Lambda} P\big(\, x \longleftrightarrow \partial^{\,in}\Lambda \,\big)\,,\\
\big(E(\varrho(\Lambda))\big)^2 \,&=\, \sum_{x,y\in\Lambda} P\big(\, x \longleftrightarrow \partial^{\,in}\Lambda \,\big) P\big(\, y \longleftrightarrow \partial^{\,in}\Lambda \,\big)\,,\\
E\big(\varrho(\Lambda)^2\big) \,&=\, \sum_{x,y\in\Lambda} P\big(\, x \longleftrightarrow \partial^{\,in}\Lambda,\ y \longleftrightarrow \partial^{\,in}\Lambda \,\big)\,.
\end{aligned}$$

### 7.1 Non-disjoint occurrence of connections

Subtracting the last two equalities, we obtain

$$\begin{aligned}
\mathrm{Var}(\varrho(\Lambda)) \,=\, \sum_{x,y\in\Lambda} \Big( & P\big(\, x \longleftrightarrow \partial^{\,in}\Lambda,\ y \longleftrightarrow \partial^{\,in}\Lambda \,\big)\\
& - P\big(\, x \longleftrightarrow \partial^{\,in}\Lambda \,\big) P\big(\, y \longleftrightarrow \partial^{\,in}\Lambda \,\big)\Big)\,. \qquad (7.1)
\end{aligned}$$

Applying the BK inequality (see for instance [63], section 2.3), we have

$$P\big(\{\, x \longleftrightarrow \partial^{\,in}\Lambda \,\} \circ \{\, y \longleftrightarrow \partial^{\,in}\Lambda \,\}\big) \,\leq\, P\big(\, x \longleftrightarrow \partial^{\,in}\Lambda \,\big) P\big(\, y \longleftrightarrow \partial^{\,in}\Lambda \,\big)\,. \qquad (7.2)$$

Substituting inequality (7.2) into formula (7.1), we get

$$\begin{aligned}
\mathrm{Var}(\varrho(\Lambda)) \,\leq\, \sum_{x,y\in\Lambda} \Big( & P\big(\, x \longleftrightarrow \partial^{\,in}\Lambda,\ y \longleftrightarrow \partial^{\,in}\Lambda \,\big)\\
& - P\big(\{\, x \longleftrightarrow \partial^{\,in}\Lambda \,\} \circ \{\, y \longleftrightarrow \partial^{\,in}\Lambda \,\}\big)\Big)\\
\leq\, & \sum_{x,y\in\Lambda} P\big(\{\, x \longleftrightarrow \partial^{\,in}\Lambda \,\} \not\circ \{\, y \longleftrightarrow \partial^{\,in}\Lambda \,\}\big)\,, \qquad (7.3)
\end{aligned}$$

where $\not\circ$ means that the events occur, but not disjointly (which is not the same as saying that they do not occur disjointly!). From now onwards, we continue the computation in the framework of the bond model. It can certainly be conducted as well in the site model, but there is an extra difficulty. In our specific case where the two events are connection events, the non-disjoint occurrence of the two events $\{ x \longleftrightarrow \partial^{in}\Lambda \}$ and $\{ y \longleftrightarrow \partial^{in}\Lambda \}$ implies the two following facts:

• First fact: the open clusters of $x$ and of $y$ in $\Lambda$ are not disjoint, so that the vertices $x$ and $y$ are connected by an open path inside $\Lambda$.

• Second fact: there exists a bond $e$ which is pivotal for both events.

We present next two very different upper bounds on $\mathrm{Var}(\varrho(\Lambda))$, based on the above two different facts.

## 7.2 First upper bound on Var($\varrho(\Lambda)$)

Using the first of the two facts, we conclude that

$$P\big(\{ x \longleftrightarrow \partial^{in}\Lambda \} \not\circ \{ y \longleftrightarrow \partial^{in}\Lambda \}\big) \,\leq\, P\big( x \longleftrightarrow y \text{ in } \Lambda,\, x \longleftrightarrow \partial^{in}\Lambda \big)\,. \quad (7.4)$$

Substituting inequality (7.4) into formula (7.3), we have

$$\begin{aligned} \mathrm{Var}(\varrho(\Lambda)) \,&\leq\, \sum_{x,y\in\Lambda} P\big( x \longleftrightarrow y \text{ in } \Lambda,\, x \longleftrightarrow \partial^{in}\Lambda \big) \\ &= \sum_{x,y\in\Lambda} E\Big( 1_{\{ x \longleftrightarrow y \text{ in } \Lambda \}} 1_{\{ x \longleftrightarrow \partial^{in}\Lambda \}} \Big) \\ &= \sum_{x\in\Lambda} E\Big( \sum_{y\in\Lambda} 1_{\{ x \longleftrightarrow y \text{ in } \Lambda \}} 1_{\{ x \longleftrightarrow \partial^{in}\Lambda \}} \Big) \\ &\leq \sum_{x\in\Lambda} E\Big( \big|C(x,\Lambda)\big| 1_{\{ x \longleftrightarrow \partial^{in}\Lambda \}} \Big)\,, \end{aligned} \quad (7.5)$$

where $C(x,\Lambda)$ is the open cluster of $x$ in $\Lambda$. Inequality (7.5) implies that

$$\mathrm{Var}(\varrho(\Lambda)) \,\leq\, |\Lambda| E\Big( \max\big\{ \big|C(x,\Lambda)\big| : x\in\partial^{in}\Lambda \big\} \Big)\,.$$

Thus, in order to control the variance of $\varrho(\Lambda)$, we have to control the expected cardinality of the largest open cluster in $\Lambda$ which is connected to $\partial^{in}\Lambda$. This upper bound can be exploited in the regime where this expected cardinality is negligible compared to $|\Lambda|$. Of course this is not the case in the supercritical regime, yet we can start an argument by supposing that it is the case at criticality and try to prove that it stands in contradiction with the hypothesis $\theta(p_c) > 0$. Indeed, as we will explain in subsection 8.3, the crucial problem to solve the conjecture $\theta(p_c, \mathbb{Z}^d) = 0$ is to prove that, whenever $\theta(p) > 0$, there exists a cluster in $\Lambda$ with cardinality close to $\theta|\Lambda|$. We can also extend the previous computation to subdomains of $\Lambda$ as follows. For $\Gamma \subset \Lambda$, we define

$$\varrho(\Gamma,\Lambda) \,=\, \big|\{ x\in\Gamma : x \longleftrightarrow \partial^{in}\Lambda \}\big|\,. \quad (7.6)$$

The same computation as above yields that

$$\mathrm{Var}(\varrho(\Gamma,\Lambda)) \,\leq\, |\Lambda| E\Big( \max\big\{ \big|C(x,\Lambda)\cap\Gamma\big| : x\in\partial^{in}\Lambda \big\} \Big)\,. \quad (7.7)$$

### 7.3 Second upper bound on Var($\varrho(\Lambda)$)

We derive next our second upper bound on $\varrho(\Lambda)$. Using the second of the two facts, we conclude that

$$P\big(\{\,x \longleftrightarrow \partial^{\,in}\Lambda\,\} \not\circ \{\,y \longleftrightarrow \partial^{\,in}\Lambda\,\}\big) \;\le\; P\begin{pmatrix} \exists e\in\mathbb{E}^d(\Lambda) \quad e \text{ is open pivotal for} \\ \text{both } x \longleftrightarrow \partial^{\,in}\Lambda \text{ and } y \longleftrightarrow \partial^{\,in}\Lambda \end{pmatrix}. \tag{7.8}$$

Substituting (7.8) into (7.3) and using a standard union bound, we obtain

$$\begin{aligned}
\mathrm{Var}(\varrho(\Lambda)) \;&\le\; \sum_{x,y\in\Lambda}\sum_{e\in\mathbb{E}^d(\Lambda)} P\begin{pmatrix} e \text{ is pivotal for both events} \\ x \longleftrightarrow \partial^{\,in}\Lambda \text{ and } y \longleftrightarrow \partial^{\,in}\Lambda \end{pmatrix} \\
&= \sum_{e\in\mathbb{E}^d(\Lambda)}\sum_{x,y\in\Lambda} \frac{1}{1-p} P\begin{pmatrix} e \text{ is closed pivotal for both} \\ x \longleftrightarrow \partial^{\,in}\Lambda \text{ and } y \longleftrightarrow \partial^{\,in}\Lambda \end{pmatrix} \\
&= \sum_{e\in\mathbb{E}^d(\Lambda)}\sum_{x,y\in\Lambda} \frac{1}{1-p} E\left(1_{\{\,e \text{ is closed}\,\}} 1_{\left\{\substack{e \text{ is pivotal for both events} \\ x\longleftrightarrow\partial^{\,in}\Lambda \text{ and } y\longleftrightarrow\partial^{\,in}\Lambda}\right\}}\right) \\
&= \sum_{e\in\mathbb{E}^d(\Lambda)} \frac{1}{1-p} E\left(1_{\{\,e \text{ is closed}\,\}} \sum_{x,y\in\Lambda} 1_{\left\{\substack{e \text{ is pivotal for both events} \\ x\longleftrightarrow\partial^{\,in}\Lambda \text{ and } y\longleftrightarrow\partial^{\,in}\Lambda}\right\}}\right)
\end{aligned} \tag{7.9}$$

Let $e = \langle v, z\rangle$ be a closed bond of the configuration restricted to $\Lambda$ and let us evaluate the sum

$$\sum_{x,y\in\Lambda} 1_{\left\{\substack{e \text{ is pivotal for both events} \\ x\longleftrightarrow\partial^{\,in}\Lambda \text{ and } y\longleftrightarrow\partial^{\,in}\Lambda}\right\}}. \tag{7.10}$$

For this sum to be non-zero, it must be the case that one endpoint of $e$ is connected by an open path to $\partial^{\,in}\Lambda$, say $v$, and the other, say $z$, is not connected to $v$ inside $\Lambda$. In this case, the vertices $x, y$ which contribute to the sum are the vertices of the open cluster of $z$, and the sum is simply equal to $|C(z,\Lambda)|^2$. Since the bond $e$ is closed, then the open cluster $C(z,\Lambda)$ does not intersect $\partial^{\,in}\Lambda$ and in fact $C(z,\Lambda) = C(z)$. Therefore the sum (7.10) is equal to

$$1_{\{\,v\longleftrightarrow\partial^{\,in}\Lambda\,\}} 1_{\{\,z\not\longleftrightarrow\partial^{\,in}\Lambda\,\}} |C(z)|^2\,. \tag{7.11}$$

In fact, we should also add the symmetric term obtained by exchanging the roles of $v$ and $z$. In any case, the term (7.11) is less or equal than

$$\max\,\big\{\,|C(z)|^2 : z\in\Lambda, z \not\longleftrightarrow \partial^{\,in}\Lambda\,\big\}\,.$$

Substituting this upper bound in (7.9), we obtain

$$\mathrm{Var}(\varrho(\Lambda)) \;\le\; \frac{2}{1-p}\big|\mathbb{E}^d(\Lambda)\big|\, E\Big(\max\,\big\{\,|C(z)|^2 : z\in\Lambda, z \not\longleftrightarrow \partial^{\,in}\Lambda\,\big\}\Big)\,.$$

As in the first case, we can extend this computation to a subdomain $\Gamma$ of $\Lambda$, and we would get

$$\mathrm{Var}(\varrho(\Gamma,\Lambda)) \,\leq\, \frac{2}{1-p}\big|\mathbb{E}^d(\Lambda)\big|\,E\Big(\max\,\big\{\,\big|C(z)|\cap\Gamma\big|^2 : z\in\Lambda, z\not\longleftrightarrow\partial^{\,in}\Lambda\,\big\}\Big)\,, \tag{7.12}$$

where the quantity $\varrho(\Gamma,\Lambda)$ was defined in (7.6).

## 7.4 Upper bound on the disconnection between $\Gamma$ and $\partial^{\,in}\Lambda$

Once we have an upper bound on $\mathrm{Var}(\varrho(\Lambda))$, we can obtain an upper bound on the probability of a disconnection between $\Gamma$ and $\partial^{\,in}\Lambda$, as follows:

$$\begin{aligned} P\big(\,\Gamma\not\longleftrightarrow\partial^{\,in}\Lambda\,\big) \,&=\, P\big(\,\varrho(\Gamma,\Lambda)=0\,\big)\\ &\leq\, P\Big(\,\big|\varrho(\Gamma,\Lambda)-E\big(\varrho(\Gamma,\Lambda)\big)\big|\geq\frac{1}{2}E\big(\varrho(\Gamma,\Lambda)\big)\Big)\\ &\leq\,\frac{4\mathrm{Var}(\varrho(\Gamma,\Lambda))}{\Big(E\big(\varrho(\Gamma,\Lambda)\big)\Big)^2}\,. \end{aligned} \tag{7.13}$$

Furthermore, we have

$$E\big(\varrho(\Gamma,\Lambda)\big)\,\geq\,E\big(\big|C_\infty\cap\Gamma\big|\big)\,=\,\theta(p,\mathbb{Z}^d)|\Gamma|\,. \tag{7.14}$$

Combining the inequalities (7.13) and (7.14), we obtain that

$$P\big(\,\Gamma\not\longleftrightarrow\partial^{\,in}\Lambda\,\big)\,\leq\,\frac{4\mathrm{Var}(\varrho(\Gamma,\Lambda))}{\big(\theta(p,\mathbb{Z}^d)\big)^2|\Gamma|^2}\,. \tag{7.15}$$

On the one hand, if we use the inequality (7.7) on $\mathrm{Var}(\varrho(\Gamma,\Lambda))$, we obtain

$$P\big(\,\Gamma\not\longleftrightarrow\partial^{\,in}\Lambda\,\big)\,\leq\,\frac{4|\Lambda|E\Big(\max\big\{\,\big|C(x,\Lambda)\cap\Gamma\big| : x\in\partial^{\,in}\Lambda\,\big\}\Big)}{\big(\theta(p,\mathbb{Z}^d)\big)^2|\Gamma|^2}\,. \tag{7.16}$$

On the other hand, if we use the inequality (7.12) on $\mathrm{Var}(\varrho(\Gamma,\Lambda))$, we obtain

$$P\big(\,\Gamma\not\longleftrightarrow\partial^{\,in}\Lambda\,\big)\,\leq\,\frac{8\big|\mathbb{E}^d(\Lambda)\big|E\Big(\max\,\big\{\,\big|C(z)\cap\Gamma\big|^2 : z\in\Lambda, z\not\longleftrightarrow\partial^{\,in}\Lambda\,\big\}\Big)}{(1-p)\big(\theta(p,\mathbb{Z}^d)\big)^2|\Gamma|^2}\,. \tag{7.17}$$

In the end, we have obtained two rather different inequalities on the probability of disconnection between $\Gamma$ and $\partial^{\,in}\Lambda$. Surprisingly, the first one involves a maximum with respect to the open clusters in $\Lambda$ which are connected to $\partial^{\,in}\Lambda$, while the second one involves a maximum with respect to the open clusters in $\Lambda$ which are not connected to $\partial^{\,in}\Lambda$. Suppose that $\theta(p)>0$ and let us take $\Lambda=\Lambda(n)$, $\Gamma=\Lambda(n^\alpha)$ for some $\alpha\in]0,1[$. With these specific choices, we expect that the right-hand quantity in (7.16) does not go to 0 as $n$ goes to $\infty$, while we do expect that the right-hand quantity in (7.17) goes to 0, at least if the exponent $\alpha$ is not small. Yet we are still quite far from being able to prove these claims.

# 8 Static renormalization

The general strategy to attack conjecture 1.1 has been known for quite a long time, it is based on a static renormalization procedure. This technique has been successfully used to prove the absence of percolation at criticality in a half-space by Barsky, Grimmett and Newman [11, 12], in $\mathbb{Z}^2 \times \{0,1\}$ by Damron, Newman and Sidoravicius [35], and in finite slabs $\mathbb{Z}^2 \times \{0,\dots,k\}$ by Duminil-Copin, Sidoravicius and Tassion [45]. The general method of static renormalization is perfectly exposed in the seventh chapter of Grimmett's book [63]. It involves a dependent site percolation process defined on a renormalized lattice associated to a specific good event. There exist many variations on this scheme, so in the next section we define precisely the one we shall use.

## 8.1 Renormalized percolation $\odot$

In order to renormalize the percolation configuration, we define successively a renormalized lattice, block events and a good event.

**The renormalized lattice.** Let $K$ be a fixed positive integer. The renormalized lattice is the sublattice $K\mathbb{Z}^d$ of $\mathbb{Z}^d$. It is isomorphic to $\mathbb{Z}^d$ and we denote it by $\underline{\mathbb{Z}}^d$. To $x \in \mathbb{Z}^d$, we associate the site $\underline{x} = Kx$ of $\underline{\mathbb{Z}}^d$. This defines a one to one correspondence between $\mathbb{Z}^d$ and $\underline{\mathbb{Z}}^d$. In general, we use underline in the notation to emphasize that we are dealing with renormalized objects. When working on the renormalized lattice, it will be convenient to use the supremum norm $|\cdot|_\infty$ defined by

$$\forall (x_1,\dots,x_d) \in \mathbb{R}^d \qquad |(x_1,\dots,x_d)|_\infty \,=\, \max_{1\le i\le d} |x_i|\,.$$

**Block events.** We recall that the cubic box $\Lambda(K)$ of side length $K$ is the set

$$\Lambda(K) \,=\, \big\{\, x \in \mathbb{R}^d : |x|_\infty \le K/2 \,\big\}\,.$$

For $\underline{x} = Kx \in \underline{\mathbb{Z}}^d$, we define the block event indexed by $\underline{x}$ as

$$B(\underline{x}) \,=\, Kx + \Lambda(4K) \,=\, \big\{\, y \in \mathbb{R}^d : |Kx - y|_\infty \le 2K \,\big\}\,.$$

Notice that the blocks associated to two neighbors in $\underline{\mathbb{Z}}^d$ overlap (see figure 14).

**The good event.** Given a cubic box $\Lambda$, we define an event $\mathcal{E}(\Lambda)$ which depends only on the percolation configuration restricted to $\Lambda$. The events are defined in a translation-invariant way: for any $x \in \mathbb{Z}^d$, any box $\Lambda$, a percolation configuration $\omega$ realizes the event $\mathcal{E}(\Lambda)$ if and only if the configuration $\omega(\cdot - x)$ translated by $-x$ realizes the event $\mathcal{E}(x + \Lambda)$.

**Renormalized configuration.** We define a percolation configuration $\underline{\omega}$ on the renormalized lattice $\underline{\mathbb{Z}}^d$ by setting

$$\forall \underline{x} \in \underline{\mathbb{Z}}^d \qquad \underline{\omega}(\underline{x}) \,=\, 1_{\mathcal{E}(B(\underline{x}))}\,.$$

In other words, a site $\underline{x} = Kx$ of $\underline{\mathbb{Z}}^d$ is occupied if the event $\mathcal{E}(B(\underline{x}))$ occurs on the initial lattice $\mathbb{Z}^d$.

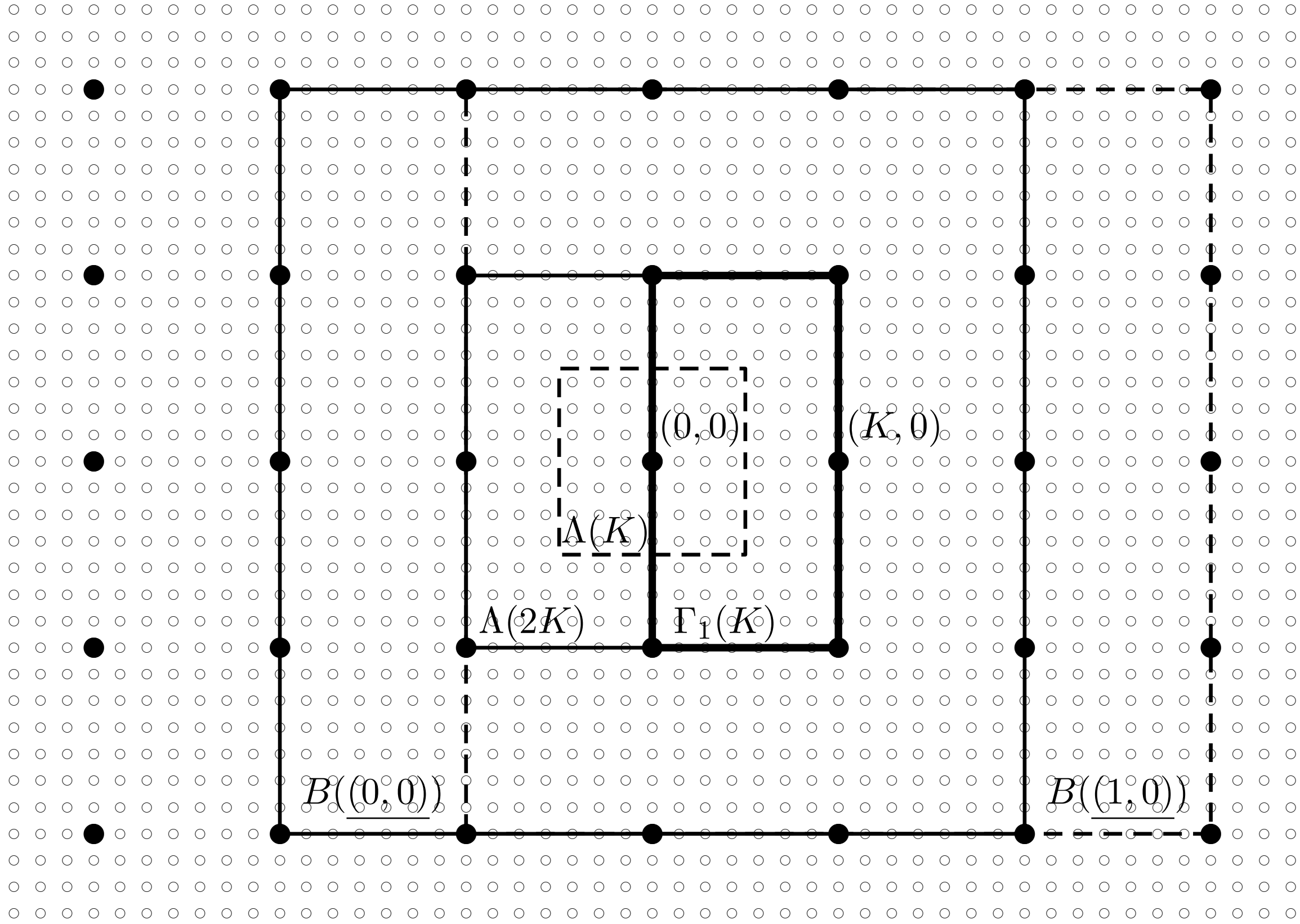

Figure 14: The two block events $B(\underline{(0,0)})$, $B(\underline{(1,0)})$ and the region $\Gamma_1(K)$

## 8.2 $k$-dependent percolation ⊙

Because the boxes $Kx + \Lambda(4K), x \in \mathbb{Z}^d$, are overlapping, the states of the sites of the renormalized lattice are not independent. However, they are $k$-dependent (with the choice $k = 5$), meaning that, for $\underline{A}, \underline{B}$ two subsets of $\underline{\mathbb{Z}}^d$, the states of the sites in $\underline{A}$ and the states of the sites in $\underline{B}$ are independent as soon as

$$\forall \underline{x} \in \underline{A} \quad \forall \underline{y} \in \underline{B} \qquad |\underline{x} - \underline{y}|_\infty > k\,.$$

More generally, if $\underline{x}_1, \cdots, \underline{x}_r$ are $r$ sites of $\underline{\mathbb{Z}}^d$ such that

$$\forall i, j \in \{\,1, \dots, r\,\} \qquad i \neq j \quad \Longrightarrow \quad |\underline{x}_i - \underline{x}_j|_\infty > k\,,$$

then the states of the sites $\underline{x}_1, \cdots, \underline{x}_r$ are independent. The translation invariance of the definition of the events $\mathcal{E}(\Lambda)$ and of the probability measure $P_p$ readily implies that

$$\begin{aligned} \forall \underline{x} \in \underline{\mathbb{Z}}^d \qquad P_p\big(\underline{x} \text{ is occupied}\big) \,&=\, P_p\big(\mathcal{E}(Kx + \Lambda(4K))\big) \\ &=\, P_p\big(\mathcal{E}(\Lambda(4K))\big) \,=\, P_p\big(\underline{0} \text{ is occupied}\big)\,. \end{aligned}$$

From the two previous remarks, we deduce the following crude estimate.

**Lemma 8.1.** *For any finite subset $\underline{A}$ of $\underline{\mathbb{Z}}^d$, we have*

$$P_p\big(\forall\ \underline{x}\in\underline{A}\quad \underline{\omega}(\underline{x})=0\big)\ \leq\ \Big(1-P_p\big(\underline{0}\ \textit{is occupied}\big)\Big)^{(k+1)^{-d}|\underline{A}|}\,.$$

*Proof.* We introduce an equivalence relation on $\underline{\mathbb{Z}}^d$: $\underline{x}\approx\underline{y}$ if and only if $k+1$ divides each component of $\underline{y}-\underline{x}$. There are $(k+1)^d$ distinct classes in $\underline{\mathbb{Z}}^d$ for the relation $\approx$, hence there exists $\underline{x}^*\in\underline{\mathbb{Z}}^d$ such that the intersection of $\underline{A}$ and the equivalence class of $\underline{x}^*$ has cardinality at least $(k+1)^{-d}|\underline{A}|$. Since the configuration $\underline{\omega}$ is $k$-dependent and $P_p$ is translation invariant, then the variables $\big(\underline{\omega}(\underline{x}),\underline{x}\in\underline{\mathbb{Z}}^d,\underline{x}\approx\underline{x}^*\big)$ are i.i.d. Bernoulli with parameter $P_p\big(\underline{0}\text{ is occupied}\big)$. Thus

$$\begin{aligned}P_p\big(\forall\ \underline{x}\in\underline{A}\quad \underline{\omega}(\underline{x})=0\big)\ &\leq\ P_p\Big(\forall\ \underline{x}\in\big\{\,\underline{y}\in\underline{A}:\underline{y}\approx\underline{x}^*\big\}\quad \underline{\omega}(\underline{x})=0\Big)\\ &\leq\ \prod_{\underline{y}\in\underline{A},\underline{y}\approx\underline{x}^*}P_p\big(\underline{\omega}(\underline{y})=0\big)\ \leq\ \Big(1-P_p\big(\underline{0}\text{ is occupied}\big)\Big)^{(k+1)^{-d}|\underline{A}|}\,,\end{aligned}$$

as claimed. □

This estimate will be used in the next section to show that, when the density of the occupied sites on the renormalized lattice is sufficiently high, the probability of percolation on the renormalized lattice becomes positive. We will not make appeal to a more sophisticated result when working on the renormalized lattice (like for instance the general stochastic domination result developed in [114, 95]).

## 8.3 The crucial problem ⊙

The crucial problem to implement successfully the classical strategy is to design a good event $\mathcal{E}$ such that:

i) If there is percolation, then the event $\mathcal{E}$ is typical in a large box:

$$\theta(p)>0\qquad\Longrightarrow\qquad \lim_{\Lambda\uparrow\mathbb{Z}^d}P_p\big(\mathcal{E}(\Lambda)\big)\ =\ 1\,.$$

ii) The existence of an infinite occupied cluster on the renormalized lattice implies the existence of an infinite open cluster on the initial lattice.

The condition ii) is implied by the following condition:

ii') If $\underline{x}_0,\dots,\underline{x}_r$ is a path of occupied sites in $\underline{\mathbb{Z}}^d$, then there exists an open path in $\mathbb{Z}^d$ starting in $B(\underline{x}_0)$, ending in $B(\underline{x}_r)$ and staying inside $B(\underline{x}_0)\cup\cdots\cup B(\underline{x}_r)$.

Suppose that this program has been completed and let us explain how this would yield that $\theta(p_c)=0$. A finite set of sites $\underline{A}$ of $\underline{\mathbb{Z}}^d$ is said to separate $\underline{0}$ from $\underline{\infty}$ if there is no connected path from $\underline{0}$ to $\underline{\infty}$ in the graph having for vertices the sites in $\underline{\mathbb{Z}}^d\setminus\underline{A}$. Recall that for the notion of $*$-connectedness, the second-nearest neighbours are taken into account, not only the nearest neighbours (see section 6 for the precise definitions). If the occupied cluster $\underline{C}(\underline{0})$ of $\underline{0}$ in the percolation configuration on the renormalized lattice $\underline{\mathbb{Z}}^d$ is finite, then its external boundary $\partial^{ext}\underline{C}(\underline{0})$ is a finite $*$-connected set, which separates $\underline{0}$ from $\underline{\infty}$ in $\underline{\mathbb{Z}}^d$. Let $q$ be a parameter in $]0,1[$. We have

$$\begin{aligned}
P_q&\begin{pmatrix}\text{there is no infinite occupied}\\ \text{path in } \underline{\mathbb{Z}}^d \text{ starting at } \underline{0}\end{pmatrix}\\
&= P_q\big(\text{the occupied cluster } \underline{C}(\underline{0}) \text{ of } \underline{0} \text{ in } \underline{\mathbb{Z}}^d \text{ is finite}\big)\\
&= P_q\begin{pmatrix}\text{there exists a finite $*$-connected set } \underline{A} \text{ of non}\\ \text{occupied sites in } \underline{\mathbb{Z}}^d \text{ which separates } \underline{0} \text{ from } \underline{\infty}\end{pmatrix}\\
&\leq \sum_{\underline{A}} P_q\big(\text{the sites of } \underline{A} \text{ are not occupied}\big)\,,
\end{aligned}$$

where the last sum runs over the finite $*$-connected subsets $\underline{A}$ of $\underline{\mathbb{Z}}^d$ which separates $\underline{0}$ from $\underline{\infty}$. We decompose the above sum according to the cardinality of $\underline{A}$. There exists a constant $\rho = \rho(d)$ such that, for any $i \geq 4$, the number of such sets having cardinality $i$ is less than $i\rho^i$. Using lemma 8.1, we obtain

$$\begin{aligned}
P_q&\begin{pmatrix}\text{there is no infinite occupied}\\ \text{path in } \underline{\mathbb{Z}}^d \text{ starting at } \underline{0}\end{pmatrix}\\
&\leq P_q\big(\underline{0} \text{ is not occupied}\big) + \sum_{i\geq 4}\sum_{\underline{A}:|\underline{A}|=i} P_q\begin{pmatrix}\text{the sites of } \underline{A}\\ \text{are not occupied}\end{pmatrix}\\
&\leq P_q\big(\underline{0} \text{ is not occupied}\big) + \sum_{i\geq 4} i\rho^i\Big(1 - P_q\big(\underline{0} \text{ is occupied}\big)\Big)^{6^{-d}i}\,. \qquad (8.1)
\end{aligned}$$

We are ready to implement the core argument of the proof. We start with a parameter $p$ such that $\theta(p) > 0$. Let $\delta_0$ in $]0,1/4[$ be such that

$$\sum_{i\geq 4} i\rho^i(\delta_0)^{6^{-d}i} < \frac{1}{4}\,. \qquad (8.2)$$

The condition $i)$ above implies that

$$\exists K \geq 1 \qquad P_p\big(\underline{0} \text{ is occupied}\big) = P_p\big(\mathcal{E}(\Lambda(4K))\big) > 1 - \delta_0\,.$$

Yet the probability $P_p\big(\mathcal{E}(\Lambda(4K))\big)$ is a continuous function of $p$, thus

$$\exists q < p \qquad P_q\big(\underline{0} \text{ is occupied}\big) > 1 - \delta_0\,. \qquad (8.3)$$

Plugging (8.3) into (8.1) and using (8.2), we conclude that

$$P_q\begin{pmatrix}\text{there is an infinite occupied}\\ \text{path in } \underline{\mathbb{Z}}^d \text{ starting at } \underline{0}\end{pmatrix} \geq \frac{1}{2}\,.$$

By condition ii), this entails the existence of an infinite open cluster in the percolation configuration on $\mathbb{Z}^d$ with parameter $q$, whence $\theta(q) > 0$. We have shown that

$$\forall p \in ]0,1] \qquad \Big[\quad \theta(p) > 0 \quad\Longrightarrow\quad \exists\, q < p \quad \theta(q) > 0 \quad\Big]\,.$$

Since $\theta(p) = 0$ for $p < p_c$, we can conclude tranquilly that $\theta(p_c) = 0$.

In fact, the realization of the above program amounts to obtain a kind of finite size criterion for the existence of an infinite open cluster. Indeed, the previous argument shows that

$$\exists\delta_0>0\quad\forall p\in]0,1[$$
$$\Big[\quad\exists K\geq 1\quad P_p\big(\mathcal{E}(\Lambda(4K))\big)\;>\;1-\delta_0\quad\Longrightarrow\quad\theta(p)>0\quad\Big]\,.$$

Putting this together with condition i), we see that

$$\exists\delta_0>0\quad\forall p\in]0,1[$$
$$\Big[\quad\theta(p)>0\quad\Longleftrightarrow\quad\exists K\geq 1\quad P_p\big(\mathcal{E}(\Lambda(4K))\big)\;>\;1-\delta_0\quad\Big]\,.$$

The probabilities of events occurring in a finite box are continuous functions of the parameter $p$, this is the reason why the finite size criterion implies the continuity of the function $\theta$. The difficulty in the above program is to design a good event which satisfies simultaneously the conditions i) and ii). If we put too much demand on the good event, then the condition ii) becomes very easy to prove, but the condition i) is out of reach. For example, we could say that $\mathcal{E}(\Lambda)$ is the event that all the sites of the box $\Lambda$ are open. On the contrary, if the requirements for a good event are too weak, then the condition i) can possibly be proved, but the condition ii) does not work. We describe next the prototype candidate for the good event.

## 8.4 Crossing clusters $\odot$

An open cluster of the percolation configuration restricted to $\Lambda$ is called crossing for $\Lambda$ if it intersects each of the $2d$ faces of $\partial^{\,in}\Lambda$. Let us consider the event $U(\Lambda)$ defined by

$$U(\Lambda)=\{\,\text{there exists an open crossing cluster } C^*\text{ in }\Lambda\,\}\,. \tag{8.4}$$

We can prove that the event $U(\Lambda)$ satisfies the condition i).

**Lemma 8.2.** *Let $d\geq 3$ and let $p$ be such that $\theta(p)>0$. We have*

$$\lim_{n\to\infty}P_p\big(U(\Lambda(n))\big)\;=\;1\,.$$

*Proof.* Let $p$ be such that $\theta(p)>0$. Let $\varepsilon>0$. Since there exists an infinite cluster with $P_p$ probability one, then there exists $m$ large enough so that

$$P_p\big(\Lambda(m)\longleftrightarrow\infty\big)\;>\;1-\varepsilon\,. \tag{8.5}$$

Let us denote by $F_i(n)$, $1\leq i\leq 2d$, the $2d$ faces of $\Lambda(n)$. For any $n\geq m$, we have

$$\big\{\,\Lambda(m)\longleftrightarrow\infty\,\big\}\;\subset\;\big\{\,\Lambda(m)\longleftrightarrow\partial^{\,in}\Lambda(n)\,\big\}$$
$$\subset\;\bigcup_{1\leq i\leq 2d}\big\{\,\Lambda(m)\longleftrightarrow F_i(n)\text{ in }\Lambda(n)\,\big\}\,.$$

Passing to the complementary event, we get

$$\bigcap_{1\leq i\leq 2d} \left\{ \Lambda(m) \not\longleftrightarrow F_i(n) \text{ in } \Lambda(n) \right\} \subset \left\{ \Lambda(m) \not\longleftrightarrow \partial^{\,in}\Lambda(n) \right\},$$

whence, by the FKG inequality,

$$\prod_{1\leq i\leq 2d} P\big( \Lambda(m) \not\longleftrightarrow F_i(n) \text{ in } \Lambda(n) \big) \;\leq\; P\big( \Lambda(m) \not\longleftrightarrow \partial^{\,in}\Lambda(n) \big)\,.$$

Using (8.5) and the symmetries of the lattice, we conclude that

$$P\big( \Lambda(m) \not\longleftrightarrow F_1(n) \text{ in } \Lambda(n) \big)^{2d} \;\leq\; \varepsilon\,.$$

Another application of the FKG inequality yields that, for any $n \geq m$,

$$P\big( \forall i \in \{\,1,\dots,2d\,\} \quad \Lambda(m) \longleftrightarrow F_i(n) \text{ in } \Lambda(n) \big) \;\geq\; \big(1-\varepsilon^{1/2d}\big)^{2d}\,. \tag{8.6}$$

Since the infinite cluster is unique, then, for $m$ fixed,

$$P\left( \bigcap_{n=m}^{+\infty} \left\{ \begin{matrix} \text{two sites of } \Lambda(m) \text{ are connected to the boundary of} \\ \text{the box } \Lambda(n) \text{ by two disjoint open clusters of } \Lambda(n) \end{matrix} \right\} \right) = 0\,,$$

and there exists $N \geq m$ large enough so that

$$\forall n \geq N \qquad P_p \begin{pmatrix} \text{two sites of } \Lambda(m) \text{ are connected to the boundary of} \\ \text{the box } \Lambda(n) \text{ by two disjoint open clusters of } \Lambda(n) \end{pmatrix} < \varepsilon\,. \tag{8.7}$$

Combining the inequalities (8.6) and (8.7), we obtain that

$$\forall n \geq N \qquad P_p\big(U(\Lambda(n))\big) \;\geq\; \big(1-\varepsilon^{1/2d}\big)^{2d} - \varepsilon\,,$$

which is the desired result. □

Now the trouble with the event $U(\Lambda)$ is that the two crossing clusters of two neighboring event blocks do not necessarily intersect, thus the condition ii) is not satisfied. So the success of the strategy relies ultimately on the choice of the good event. How should we proceed to get the right event? A sensible approach consists in examining all the finite block events that are typical in large blocks under the sole hypothesis that $\theta(p)$ is positive, and check whether they satisfy the condition ii), or, if not, try to understand what is missing. We would then try to strengthen slightly the requirements on the good event. A recent and up-to-date discussion of this argument can be found in [44]. The problem is that most of the existing results describing the percolation configuration in a finite box are derived under the stronger hypothesis that $p > p_c$. A notable exception is the result on the travel times [27], we shall come back to it later. What we can do for the time being is to examine the existing results proved in the supercritical regime and try to see what remains under the weaker hypothesis that $\theta(p) > 0$. We start this program in the next subsection.

## 8.5 Refinements ⊛

In fact, if it is the case that $\theta(p_c) = 0$, then the hypothesis $\theta(p) > 0$ implies that $p > p_c$, and any block event which is typical in the supercritical regime is eligible for the fulfillment of the program. For instance, the block events $U, R, V, W$ employed by Pisztora [114] should work. The event $U$ has already been defined in (8.4), we recall next the definitions of the events $R, V, W$. For $m \in \mathbb{N}$ and $\varepsilon \in ]0,1[$, we define the following events:

$$\begin{aligned} R(\Lambda, m) &= U(\Lambda) \cap \{\, C^* \text{ crosses every sub-box of } \Lambda \text{ with diameter } \geq m \,\} \\ &\qquad \cap \{\, \text{there exists a unique open cluster in } \Lambda \text{ with diameter } \geq m \,\}, \\ V(\Lambda, \varepsilon) &= U(\Lambda) \cap \{\, \theta(1-\varepsilon)\, |\Lambda| \leq |C^*| \,\}, \\ W(\Lambda, \varepsilon) &= \{\, \big|\{\, x \in \Lambda : x \longleftrightarrow \partial^{in}\Lambda \,\}\big| \,\leq\, (1+\varepsilon)\, \theta\, |\Lambda| \,\}. \end{aligned}$$

The slab percolation threshold $\widehat{p}_c$ is the limit of the critical points associated to percolation in finite slabs $S(L) = \mathbb{Z}^{d-1} \times [1, L]$ as $L$ goes to $\infty$. Pisztora [114] has proved that, for $p > \widehat{p}_c$, all the above events are typical, i.e., their probabilities converge to 1 as the size of the box goes to $\infty$, for a suitable choice of the length $m$. Grimmett and Marstrand [62] have proved that the slab percolation threshold $\widehat{p}_c$ coincides with the critical point $p_c$. Combining the results of these two fundamental works, we have the following proposition.

**Proposition 8.3.** *Let $d \geq 3$. For any $p > p_c$, any $\varepsilon > 0$, there exists a positive constant $\kappa$ such that*

$$\begin{aligned} \lim_{n\to\infty} P_p\big(U(\Lambda(n))\big) &= 1, \\ \lim_{n\to\infty} P_p\big(R(\Lambda(n), \kappa \ln n)\big) &= 1, \\ \lim_{n\to\infty} P_p\big(V(\Lambda(n), \varepsilon)\big) &= 1, \\ \lim_{n\to\infty} P_p\big(W(\Lambda(n), \varepsilon)\big) &= 1. \end{aligned}$$

Of course the big difficulty is to reach the same conclusions starting only with the hypothesis that $\theta(p) > 0$. Unfortunately it is believed and supported by numerical simulations that multiple crossing clusters do exist at the critical point in dimensions 3 and higher, and it seems very difficult to get rid of them with the sole hypothesis $\theta(p) > 0$. Thus we should definitely weaken the conditions involved in the event blocks of Pisztora. Because multiple crossing clusters seem to exist at $p_c$, we should certainly drop the condition that there exists a unique crossing cluster. What we really need is that there exists a big macroscopic cluster in a good block, together with a mechanism to ensure that the macroscopic clusters of two adjacent good blocks are connected. Let us consider the box $\Lambda(n)$ of side length $n$. By the spatial ergodic theorem, we have, with probability 1,

$$\lim_{n\to\infty} \frac{1}{|\Lambda(n)|} \sum_{x \in \Lambda(n)} 1_{\{\, x \longleftrightarrow \infty \,\}} \;=\; P\big(0 \longleftrightarrow \infty\big) \;=\; \theta(p). \qquad (8.8)$$

Let $p \in [0,1]$ be such that $\theta(p) > 0$ and let $\varepsilon$ be such that $0 < \varepsilon < \theta(p)$. The limit (8.8) implies that, for $n$ large enough,

$$\frac{1}{|\Lambda(n)|} \sum_{x \in \Lambda(n)} 1_{\{\, x \longleftrightarrow \infty \,\}} \ge \theta(p) - \varepsilon \,. \tag{8.9}$$

Let $\ell \ge n$ and let us consider now the box of side length $\ell$. For any $x \in \Lambda(n)$, we have the inclusion $\{\, x \longleftrightarrow \infty \,\} \subset \{\, x \longleftrightarrow \partial^{in}\Lambda(\ell) \,\}$, whence

$$\sum_{x \in \Lambda(n)} 1_{\{\, x \longleftrightarrow \infty \,\}} \le \sum_{x \in \Lambda(n)} 1_{\{\, x \longleftrightarrow \partial^{in}\Lambda(\ell) \,\}} \,. \tag{8.10}$$

Let us reformulate this inequality in terms of clusters. Let $\mathcal{C}(n,\ell)$ be the collection of the clusters of the configuration restricted to $\Lambda(\ell)$ which intersect both $\Lambda(n)$ and $\partial^{in}\Lambda(\ell)$. When we look at the intersection of the infinite cluster $C_\infty$ with the box $\Lambda(n)$, we might obtain a multitude of different connected components. Yet each of these components has to be included in a cluster of the collection $\mathcal{C}(n,\ell)$, so that we have the following cluster version of (8.10):

$$C_\infty \cap \Lambda(n) \subset \bigcup_{C \in \mathcal{C}(n,\ell)} C \cap \Lambda(n) \,. \tag{8.11}$$

Inequalities (8.9) and (8.11) imply that

$$\sum_{C \in \mathcal{C}(n,\ell)} \big| C \cap \Lambda(n) \big| \ge \big(\theta(p) - \varepsilon\big) |\Lambda(n)| \,. \tag{8.12}$$

## 8.6 Control of the number of clusters ⊛

We would like to prove that there are not too many clusters in $\mathcal{C}(n,\ell)$. Let us denote by $\mathcal{N}(n,\ell)$ the number of clusters in the collection $\mathcal{C}(n,\ell)$. For fixed $n$, the random number $\mathcal{N}(n,\ell)$ is a non-increasing function of $\ell$: indeed, when we increase $\ell$ to $\ell+1$, a cluster in $\mathcal{C}(n,\ell)$ either stays as such, or it disappears, or it coalesces with another cluster of $\mathcal{C}(n,\ell)$. Moreover any cluster joining $\Lambda(n)$ to $\partial^{in}\Lambda(\ell+1)$ has to be formed from clusters that were already present in $\mathcal{C}(n,\ell)$. In addition, we know that the infinite cluster, when it exists, is unique. Therefore, with probability 1,

$$\lim_{\ell \to \infty} \mathcal{N}(n,\ell) = 1 \,.$$

The crucial problem is to control the speed of this convergence. The following theorem was proved in [28].

**Theorem 8.4.** *Let $d \ge 2$ and let $p_c$ be the critical probability of the site percolation model in $d$ dimensions. Let $\alpha$ be such that*

$$\alpha > \frac{2d^2 + 2d - 2}{2d^2 + 3d - 3}(4d^2 + 5d - 5) \,. \tag{8.13}$$

*We have*

$$\lim_{n \to \infty} P_{p_c}\big(\mathcal{N}(n, n^\alpha) \ge 2\big) = 0 \,. \tag{8.14}$$

For $d = 3$, this gives

$$\lim_{n\to\infty} P_{p_c}\big(\mathcal{N}(n, n^{43}) \geq 2\big) \;=\; 0\,.$$

This result was obtained by revisiting the beautiful argument of Aizenman, Kesten and Newman [4] for the uniqueness of the infinite cluster. Their argument yields a control of the two arms event, in the following form: for $d \geq 2$,

$$\exists\,\kappa > 0 \quad \forall n \geq 2 \qquad P_{p_c}\big(\text{two-arms}(0, n)\big) \;\leq\; \frac{\kappa \ln n}{\sqrt{n}}\,. \tag{8.15}$$

The event "two-arms$(0, n)$" is the event that two neighbours of 0 are connected to the boundary of the box $\Lambda(n)$ by two disjoint open clusters. The deep uniqueness proof of [4] was originally written for a quite general percolation model. A simplified and illuminating version has been worked out by Gandolfi, Grimmett and Russo [57]. Nowadays, the uniqueness of the infinite cluster in percolation is proved with the help of the more robust Burton-Keane argument: see for instance [63], section 8.2 or [64], second half of page 292. The inequality (8.15) is absolutely remarkable, because it holds for any value of $p$ in $]0,1[$, even at the critical point $p_c$ (the constant $\kappa$ might depend on $p$). It is further discussed in [131, 132]. Combining (8.12) and (8.14), we see that, with probability going to 1 as $n$ goes to $\infty$, for an exponent $\alpha$ satisfying (8.13), inside the box $\Lambda(n^\alpha)$, there will be a cluster $C$ whose trace on the box $\Lambda(n)$ satisfies

$$\big|C\cap\Lambda(n)\big|\,\geq\,\big(\theta(p)-\varepsilon\big)|\Lambda(n)|\,.$$

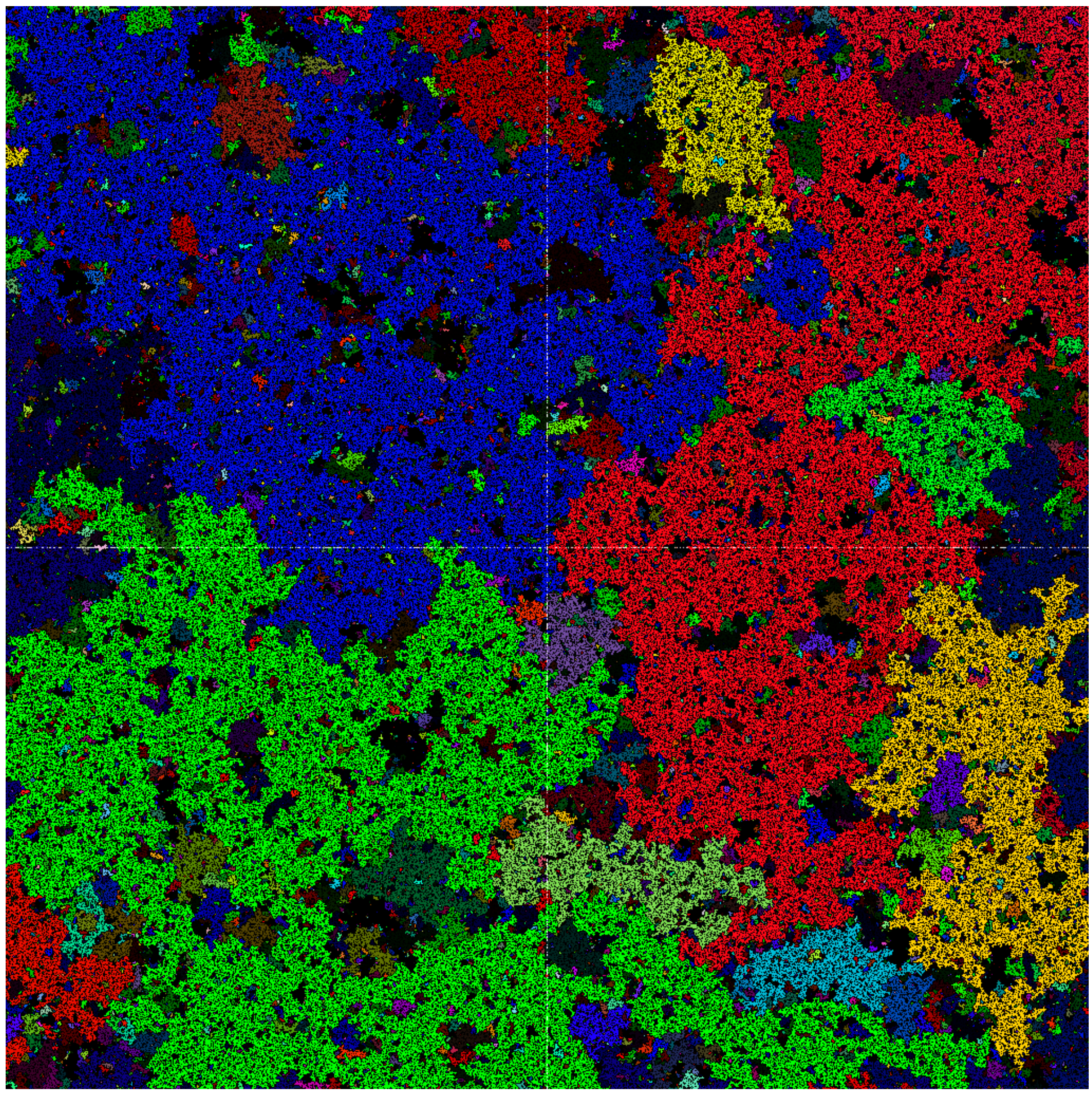

Figure 15: two-arms$(0, 1024)$. Window centered around a pivotal site for left to right connection in a box of size $2000 \times 2000$, site percolation at $p = 0.5905$.

Unfortunately, this is not enough to fulfill our program. The point which does not work is the Peierls type argument. Indeed, the range of dependence of the block event is $n^\alpha$, and this ruins completely the computation done in (8.1). Naturally, we should try to improve this range $n^\alpha$. The exponent $\alpha$ in (8.13) is obtained through a bootstrapping procedure initialized with the inequality (8.15). Although this argument looked promising, it could only improve slightly the exponent 1/2. If we use the additional hypothesis that $\theta(p)$ is positive, we get a further serious improvement on the exponent governing the establishment of long-range order in a box (see theorem 1.3 of [28]), but this is not enough. What we really need is to show that $\ell$ can be taken of the same order as $n$, i.e., that

$$\exists\, c>0\qquad \lim_{n\to\infty}\,\mathcal{N}(n,cn)\,=\,1\quad\text{almost surely}\,. \tag{8.16}$$

Notice that, in the supercritical regime $p>p_c$, Kesten and Zhang [84] have proved with the help of the slab technology that

$$P\big(\mathcal{N}(n,2n)\geq 2\big)\,=\,\exp(-O(n))\,, \tag{8.17}$$

so that the estimate (8.16) we seek should indeed be true. However, since multiple crossing clusters seem to exist at the critical point $p_c$, it is certainly too ambitious to hope to prove something like (8.17). So, the question is, how can we improve the control of $\mathcal{N}(n,\ell)$?

## 8.7 Some failed attempts ⊛

The trick to prove the estimate (8.14) was to obtain a control of the two arms event, of the following form (see lemma 5.1 of [28], but beware that the variable $\ell$ used here corresponds to $n+\ell$ in [28]). For any $p$ in $]0,1[$, any $n\geq 1$, $\ell\geq n$, we have the inequality (where $c=c(p)=p(1-p)$)

$$P(\text{two-arms}(0,n+\ell))\leq 2d(\ln n)E\Big(\sqrt{\frac{\mathcal{N}(n,\ell)}{\big|\Lambda(n)\big|}}\Big)+\frac{4d}{c}\big|\Lambda(n)\big|^2\exp\Big(-(c\ln n)^2\Big).$$

The main problem is then to obtain a control on the number $\mathcal{N}(n,\ell)$ of clusters in the collection $\mathcal{C}(n,\ell)$. The uniqueness result of [4] rests ultimately on the observation that $\mathcal{N}(n,\ell)\,\leq\,\big|\partial\Lambda(n)\big|\,\leq\,2dn^{d-1}$. Plugging the latter inequality in the former, we obtain the estimate on the two arms event stated in (8.15). Theorem 1.1 of [28] yields the following inequality:

$$\exists\,\gamma>\frac{1}{2}\quad\exists\,c>0\quad\forall n\geq 1\qquad P_{p_c}\big(\text{two-arms}(0,n)\big)\;\leq\;\frac{c}{n^\gamma}\,. \tag{8.18}$$

We would like to improve this estimate, or to get some better control in order to define an adequate good event. Suppose that we have a lot of clusters realizing the connection between $\Lambda(n)$ and $\partial^{in}\Lambda(\ell)$. We will then see many pairs of clusters inside $\Lambda(\ell)$ which realize the two arms event. To obtain an upper bound on the probability of this situation, we naturally think of applying the BK inequality, or rather Reimer's inequality (because the two arms event is not monotone). What we would get is something of the following sort:

$$P\begin{pmatrix}\exists\,C_1,\dots,C_k\in\mathcal{C}(n,n^\alpha)\qquad\forall i,j\in\{\,1,\dots,k\,\}\\ i\neq j\quad\Rightarrow\quad\forall x\in C_i\quad\forall y\in C_j\quad|x-y|_\infty>1\end{pmatrix}\;\leq\;\Big(\frac{c}{n^\gamma}\Big)^k\,. \tag{8.19}$$

This looks very promising, because if we take $k$ of the form $n^\beta$, we obtain a very fast convergence to 0, and we could hope that this would compensate for the long range dependence $n^\alpha$ of the renormalized process. Unfortunately, the inequality (8.19) does not yield a direct control on the random number $\mathcal{N}(n,\ell)$. The problem is that, to apply Reimer's inequality, we must consider events occurring disjointly, and two clusters which have a common boundary point do not occur disjointly. In any case, the inequality (8.19) tells us that, in a group of clusters emanating from $\Lambda(n)$ and reaching the boundary of $\Lambda(\ell)$, with very high probability, there exists at least a closed site belonging to the boundary of two distinct clusters. This indicates that large clusters have to meet quickly, but only along closed sites, and this does not preclude that they remain disconnected from each other. Let us mention another argument indicating that large clusters have to meet quickly. Imagine an exploration algorithm which starts from the box $\Lambda(n)$ and which explores successively each shell

$$\Lambda(n+1)\setminus\Lambda(n)\,,\quad \Lambda(n+2)\setminus\Lambda(n+1)\,,\quad\cdots$$

At step $k$, we will have uncovered the collection $\mathcal{C}(n,n+k)$. Let us look at the trace of $\mathcal{C}(n,n+k)$ on the boundary $\partial^{in}\Lambda(n+k)$. If two clusters $C_1,C_2$ of $\mathcal{C}(n,n+k)$ have a common boundary site belonging to $\partial^{in}\Lambda(n+k)$, then, at the next step, they will coalesce with a probability larger or equal than $p^5$ ($p^3$ works for sites having a common neighbour and belonging to the same face of $\Lambda(n+k)$, we put $p^5$ to handle also the case of sites belonging to different faces). If two clusters $C_1,C_2$ of $\mathcal{C}(n,n+k)$ do not have a common boundary site belonging to $\partial^{in}\Lambda(n+k)$, then, this means that they don't feel each other, and their next exploration steps will be independent (this is a bit over-simplistic, because two clusters might already have a common boundary site deep inside $\Lambda(n+k)$, and this would create a dependence). Thus the clusters of $\mathcal{C}(n,n+k)$ which stay apart perform independent explorations of $\partial^{in}\Lambda(n+k)$. Now, in dimension 3 at least, the set $\partial^{in}\Lambda(n+k)$ is a two-dimensional set, so we have a group of clusters which explore independently this two-dimensional set. If we could embed a simple exploration mechanism inside each cluster, ideally like a random walk, then we would certainly get an interesting quantitative control on $\mathcal{N}(n,\ell)$. Unfortunately this approach has not succeeded until now.

## Part II

# The ultimate goal ⊙

At the end of the previous part, we presented the general strategy for proving the conjecture $\theta(p_c, \mathbb{Z}^d) = 0$ and we discussed some failed attempts to implement it. These failures were due to the use of an inappropriate good event. In section 9, we propose a choice for the good event, which should be the right one to complete the proof. We prove that, if this event is typical, then the strategy involving the static renormalization will be successful. Of course, that is the easy part. The ultimate goal is to show that the event in question is typical. In section 10, we offer a vague rationale to support the validity of this approach. Finally, in section 11, we state the discrete isoperimetric inequality in $\mathbb{Z}^d$, that will undoubtedly play a central role in any proof.

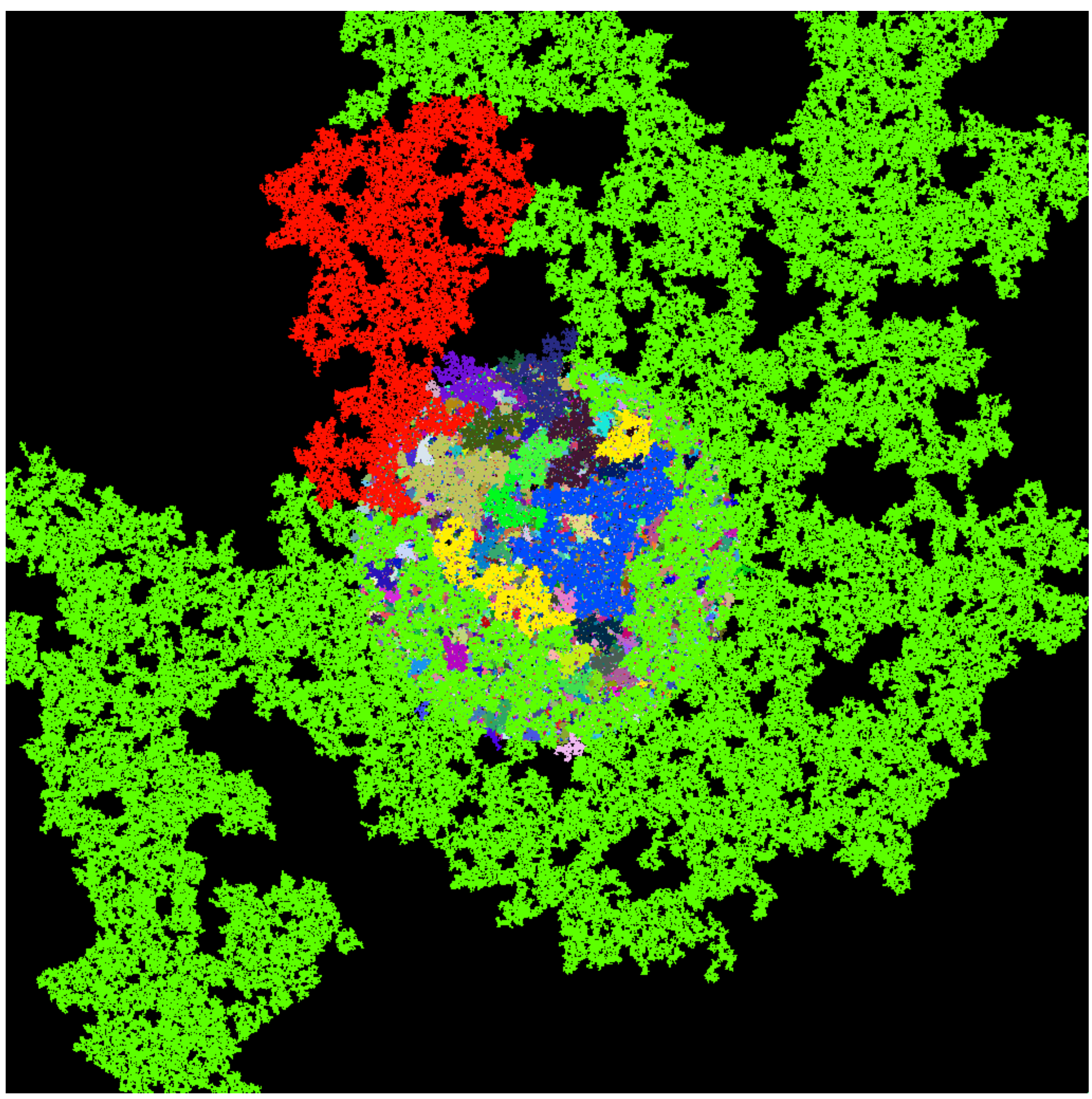

Figure 16: Bond percolation, p=0.5001, 1024x1024 box.
The good event is realized by the green cluster.

# 9 The good event

We proceed now to define the good event, which, in our view, is the right one to complete the proof that $\theta(p_c,\mathbb{Z}^d)=0$. For $i\in\{1,\dots,d\}$, we denote by $e_i=(0,\dots,0,1,0,\dots,0)$ the $i$-th vector of the canonical basis of $\mathbb{R}^d$. Let $K\geq 1$. We define the parallelepipeds

$$\begin{aligned}\Gamma_i(K) &= \Lambda(2K)\cap\big(Ke_i+\Lambda(2K)\big)\,,\\ \Gamma_{-i}(K) &= \Lambda(2K)\cap\big(-Ke_i+\Lambda(2K)\big)\,.\end{aligned}$$

The set $\Gamma_1(K)$ is depicted in figure 14. We recall that, for $\Gamma$ a subset of $\mathbb{Z}^d$, we had defined in (5.1) the quantity

$$\varrho(\Gamma) = \big|\{\,x\in\Gamma : x\longleftrightarrow\partial^{\,in}\Gamma\,\}\big|\,.$$

Our choice for the block event $\mathcal{E}$ is the following. Given $K\geq 1$ and $\varepsilon>0$, we define

$$\mathcal{E}(\Lambda(4K))=\left\{\begin{matrix}\text{there exists an open cluster } C^* \text{ in } \Lambda(4K)\\ \text{which touches the boundary } \partial^{\,in}\Lambda(4K) \text{ and}\\ \text{which satisfies } \big|C^*\cap\Lambda(2K)\big|\,\geq\,(\theta-\varepsilon)\,|\Lambda(2K)|\,,\\ \forall i\in\{\,-d,\dots,d\,\}\setminus\{\,0\,\}\quad \varrho(\Gamma_i(K))\,\leq\,(\theta+\varepsilon)\,|\Gamma_i(K)|\end{matrix}\right\}. \tag{9.1}$$

The ultimate goal is to prove that the event $\mathcal{E}$ defined above satisfies the two conditions $i)$ and $ii)$ stated at the beginning of section 8.3. We first show that the condition $ii')$ holds. Let $\underline{x}_0,\dots,\underline{x}_r$ be a path of vertices in $\underline{\mathbb{Z}}^d$, and let us suppose that all the sites of this path are occupied, meaning that the events $\mathcal{E}(B(\underline{x}_0)),\dots,\mathcal{E}(B(\underline{x}_r))$ occur simultaneously. For each $m\in\{\,1,\dots,r\,\}$, we denote by $C^*_m$ an open cluster in $B(\underline{x}_m)$ which realizes the event $\mathcal{E}(B(\underline{x}_m))$ defined in (9.1). Let us fix $m\in\{\,0,\dots,r-1\,\}$. We claim that the two clusters $C^*_m$ and $C^*_{m+1}$ are connected within $B(\underline{x}_m)\cup B(\underline{x}_{m+1})$. Indeed, since $\underline{x}_m$ and $\underline{x}_{m+1}$ are neighbours, there exists $i_0\in\{\,-d,\dots,d\,\}\setminus\{\,0\,\}$ such that $\underline{x}_{m+1}=\underline{x}_m+\text{sign}(i_0)e_{|i_0|}$. We will look at the sets

$$\begin{aligned}\Gamma^*=\Gamma(m,i_0) &= K\underline{x}_m+\Gamma_{i_0}(K) = K\underline{x}_{m+1}+\Gamma_{-i_0}(K) = \Gamma(m+1,-i_0)\,,\\ \Gamma(m,-i_0) &= K\underline{x}_m+\Gamma_{-i_0}(K)\,,\quad \Gamma(m+1,i_0) = K\underline{x}_{m+1}+\Gamma_{i_0}(K)\,.\end{aligned}$$

The simultaneous occurrence of the events $\mathcal{E}(B(\underline{x}_m))$, $\mathcal{E}(B(\underline{x}_{m+1}))$ implies that

$$\big|C^*_m\cap\big(K\underline{x}_m+\Lambda(2K)\big)\big|\,\geq\,(\theta-\varepsilon)\,|\Lambda(2K)|\,, \tag{9.2}$$

$$\big|C^*_{m+1}\cap\big(K\underline{x}_{m+1}+\Lambda(2K)\big)\big|\,\geq\,(\theta-\varepsilon)\,|\Lambda(2K)|\,, \tag{9.3}$$

$$\big|C^*_m\cap\Gamma(m,-i_0)\big|\,\leq\,\varrho(\Gamma(m,-i_0))\,\leq\,(\theta+\varepsilon)\,|\Gamma(m,-i_0)|\,, \tag{9.4}$$

$$\big|C^*_{m+1}\cap\Gamma(m+1,i_0)\big|\,\leq\,\varrho(\Gamma(m+1,i_0))\,\leq\,(\theta+\varepsilon)\,|\Gamma(m+1,i_0)|\,, \tag{9.5}$$

$$\varrho(\Gamma^*)\,\leq\,(\theta+\varepsilon)\,|\Gamma^*|\,. \tag{9.6}$$

Notice that we used the fact that the clusters $C^*_m$ and $C^*_{m+1}$ are connected to the inner boundaries of their respective event blocks to obtain the inequalities (9.4) and (9.5): indeed, all the vertices of $C^*_m \cap \Gamma(m,-i_0)$ are connected to the boundary of $K\underline{x}_m + \Lambda(4K)$, therefore they are also connected to $\partial^{\,in}\Gamma(m,-i_0)$ and they contribute to $\varrho(\Gamma(m,-i_0))$ (the same argument works for $C^*_{m+1} \cap \Gamma(m+1,i_0)$). Since we have in addition

$$K\underline{x}_m + \Lambda(2K) \;=\; \Gamma(m,-i_0) \cup \Gamma^*\,,$$

then, using (9.2) and (9.4), we get

$$\begin{aligned} \left|C^*_m \cap \Gamma^*\right| \;\geq\; & \left|C^*_m \cap \big(K\underline{x}_m + \Lambda(2K)\big)\right| - \left|C^*_m \cap \Gamma(m,-i_0)\right| \\ \geq\; & (\theta-\varepsilon)\,|\Lambda(2K)| - (\theta+\varepsilon)\,|\Gamma(m,-i_0)|\,. \end{aligned} \tag{9.7}$$

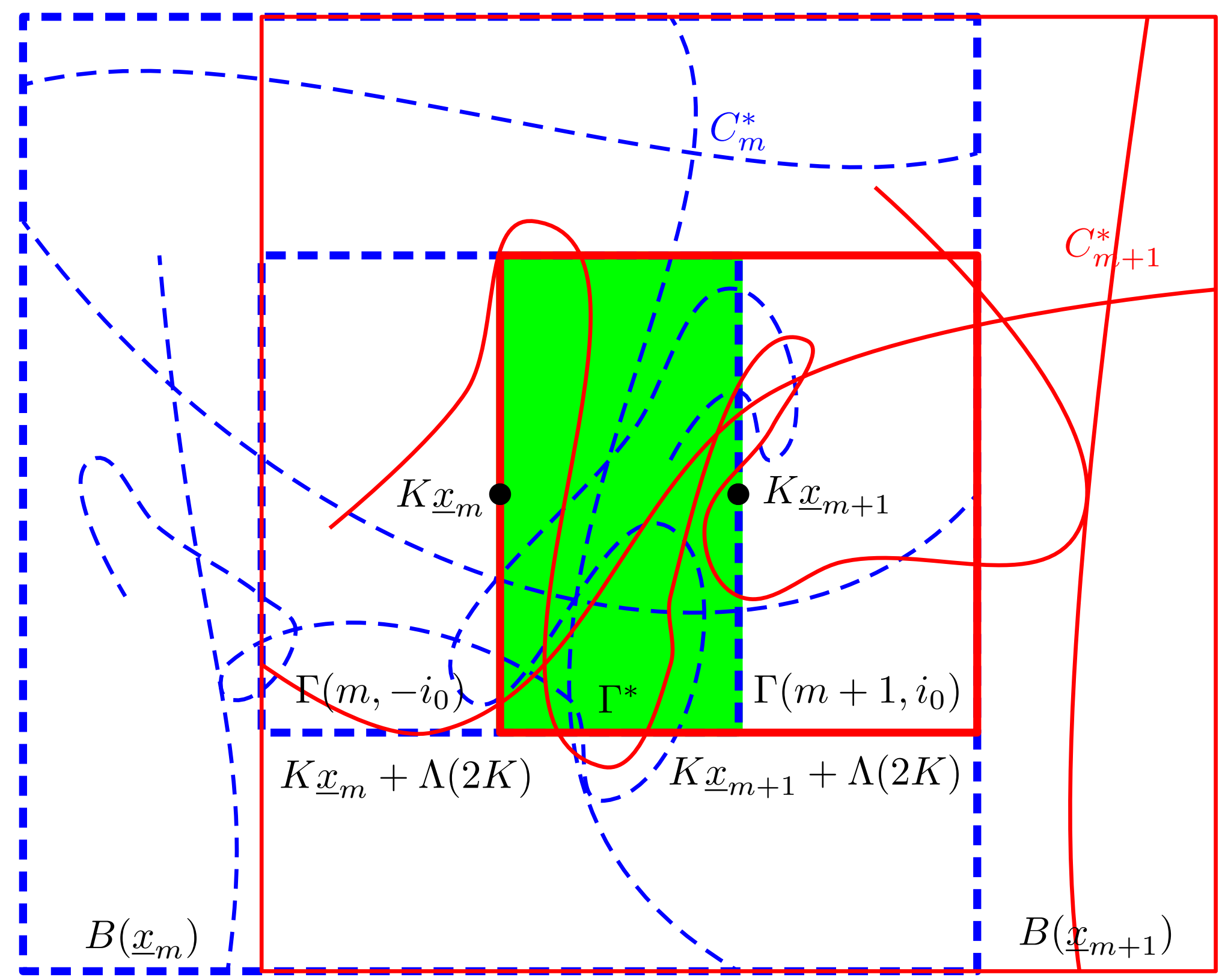

Figure 17: The two clusters $C^*_m$, $C^*_{m+1}$ must meet within the region $\Gamma^*$

Similarly, the inequalities (9.3) and (9.5) yield that

$$\left|C^*_{m+1} \cap \Gamma^*\right| \;\geq\; (\theta-\varepsilon)\,|\Lambda(2K)| - (\theta+\varepsilon)\,|\Gamma(m+1,i_0)|\,. \tag{9.8}$$

Using once more the fact that the clusters $C^*_m$ and $C^*_{m+1}$ are connected to the inner boundaries of their respective event blocks, we see that all the sites belonging to $C^*_m \cap \Gamma^*$ or $C^*_{m+1} \cap \Gamma^*$ are also connected to $\partial^{\,in}\Gamma^*$, hence they contribute to $\varrho(\Gamma^*)$ and

$$\begin{aligned} \varrho(\Gamma^*) \;\geq\; & \left|\big(C^*_m \cap \Gamma^*\big) \cup \big(C^*_{m+1} \cap \Gamma^*\big)\right| \\ \geq\; & \left|C^*_m \cap \Gamma^*\right| + \left|C^*_{m+1} \cap \Gamma^*\right| - \left|C^*_m \cap C^*_{m+1} \cap \Gamma^*\right|\,. \end{aligned} \tag{9.9}$$

We finally use the inequalities (9.6), (9.7) and (9.8) to estimate the two sides of (9.9), together with the fact that the three sets $\Gamma(m,-i_0)$, $\Gamma(m+1,i_0)$, $\Gamma^*$ have the same cardinality, and we get

$$\begin{aligned}\big|C^*_m \cap C^*_{m+1} \cap \Gamma^*\big| \,&\geq\, 2(\theta-\varepsilon)\,|\Lambda(2K)| - 2(\theta+\varepsilon)\,|\Gamma^*| - (\theta+\varepsilon)\,|\Gamma^*| \\ &\geq\, \Big(\frac{\theta}{2}-4\varepsilon\Big)\,|\Lambda(2K)| - 4(2K+1)^{d-1}\,,\end{aligned}$$

where we have used the inequality $2|\Gamma^*| \,\leq\, |\Lambda(2K)| + (2K+1)^{d-1}$ for the last step. We suppose that $K$ is large enough so that $4(2K+1)^{d-1} \leq \varepsilon\big|\Lambda(2K)\big|$ and that $\varepsilon < \theta/20$. We can then conclude that

$$\big|C^*_m \cap C^*_{m+1} \cap \Gamma^*\big| \,\geq\, \Big(\frac{\theta}{2}-5\varepsilon\Big)\,|\Lambda(2K)| \,\geq\, \frac{\theta}{4}\,|\Lambda(2K)| \,>\, 0\,.$$

Therefore the clusters $C^*_m$ and $C^*_{m+1}$ share a common vertex and they belong to the same cluster of the percolation configuration restricted to $B(\underline{x}_m)\cup B(\underline{x}_{m+1})$. This is true for any $m\in\{\,0,\dots,r-1\,\}$, therefore all the clusters $C^*_0,\dots,C^*_r$ are connected together, and they contain an open path starting in $B(\underline{x}_0)$, ending in $B(\underline{x}_r)$ and staying inside $B(\underline{x}_0)\cup\cdots\cup B(\underline{x}_r)$. Thus condition ii') is proved. We will next show that the event $\mathcal{E}(\Lambda(4K))$ is typical, i.e., its probability converges to 1 as the size of the box goes to $\infty$. The easy part is the inequalities on $\varrho(\Gamma_i(K))$ for $i\in\{\,-d,\dots,d\,\}\setminus\{\,0\,\}$, it is the object of the next proposition.

**Proposition 9.1.** *For any* $i\in\{\,-d,\dots,d\,\}\setminus\{\,0\,\}$, *we have*

$$\forall\varepsilon>0\qquad \lim_{K\to\infty} P\Big(\varrho\big(\Gamma_i(K)\big) \,\leq\, (\theta+\varepsilon)|\Gamma_i(K)|\Big) \,=\, 1\,.$$

*Proof.* We do the proof for the set $\Gamma_1$. We introduce the notation

$$\forall\varepsilon_1,\dots,\varepsilon_d\in\,\{\,-1,1\,\}\qquad e_K(\varepsilon_1,\dots,\varepsilon_d) \,=\, \frac{K}{2}\sum_{1\leq i\leq d}\varepsilon_i e_i\,.$$

We decompose $\Gamma_1$ as the union of cubic boxes of side length $K$ as follows:

$$\Gamma_1 \,=\, \bigcup_{\varepsilon_2,\dots,\varepsilon_d\in\{\,-1,1\,\}} \Big(e_K(1,\varepsilon_2,\dots,\varepsilon_d)+\Lambda(K)\Big)\,.$$

These boxes are not disjoint, however two boxes can intersect only along their faces, so that

$$|\Gamma_1(K)| \,\geq \sum_{\varepsilon_2,\dots,\varepsilon_d\in\{\,-1,1\,\}} \Big|e_K(1,\varepsilon_2,\dots,\varepsilon_d)+\Lambda(K)\Big| - 2^{d-1}\Big|\partial^{\,in}\Lambda(K)\Big|\,. \tag{9.10}$$

By the subadditivity of $\varrho$, we have

$$\varrho(\Gamma_1(K)) \,\leq \sum_{\varepsilon_2,\dots,\varepsilon_d\in\{\,-1,1\,\}} \varrho\Big(e_K(1,\varepsilon_2,\dots,\varepsilon_d)+\Lambda(K)\Big)\,. \tag{9.11}$$

Let $\varepsilon$ be such that $\theta+\varepsilon<1$ (otherwise there is nothing to prove) and suppose

$$\varrho(\Gamma_1(K)) \,>\, (\theta+\varepsilon)|\Gamma_1(K)|\,. \tag{9.12}$$

Inequalities (9.10), (9.11) and (9.12) together yield that

$$\sum_{\varepsilon_2,\dots,\varepsilon_d\in\{\,-1,1\,\}} \varrho\Big(e_K(1,\varepsilon_2,\dots,\varepsilon_d)+\Lambda(K)\Big) \;>\; (\theta+\varepsilon)\sum_{\varepsilon_2,\dots,\varepsilon_d\in\{\,-1,1\,\}}\Big|e_K(1,\varepsilon_2,\dots,\varepsilon_d)+\Lambda(K)\Big|-2^{d-1}\Big|\partial^{\,in}\Lambda(K)\Big|\,. \tag{9.13}$$

We suppose that $K$ is large enough so that

$$\big|\partial^{\,in}\Lambda(K)\big| \;\leq\; \frac{\varepsilon}{2}\Big|e_K(1,\varepsilon_2,\dots,\varepsilon_d)+\Lambda(K)\Big|\,. \tag{9.14}$$

Inequalities (9.13), (9.14) imply that there exists $\varepsilon_2,\dots,\varepsilon_d\in\{\,-1,1\,\}$ such that

$$\varrho\Big(e_K(1,\varepsilon_2,\dots,\varepsilon_d)+\Lambda(K)\Big) \;>\; (\theta+\frac{\varepsilon}{2})\Big|e_K(1,\varepsilon_2,\dots,\varepsilon_d)+\Lambda(K)\Big|\,.$$

Thus

$$\begin{aligned}
&P\big(\varrho(\Gamma_1(K)) \;>\; (\theta+\varepsilon)|\Gamma_1(K)|\big) \;\leq\\
&\sum_{\varepsilon_2,\dots,\varepsilon_d\in\{\,-1,1\,\}} P\bigg(\varrho\Big(e_K(1,\varepsilon_2,\dots,\varepsilon_d)+\Lambda(K)\Big)>(\theta+\frac{\varepsilon}{2})\Big|e_K(1,\varepsilon_2,\dots,\varepsilon_d)+\Lambda(K)\Big|\bigg)\\
&\quad\leq 2^{d-1}P\bigg(\varrho\Big(e_K(1,\dots,1)+\Lambda(K)\Big) \;>\; (\theta+\frac{\varepsilon}{2})\Big|e_K(1,\dots,1)+\Lambda(K)\Big|\bigg)\\
&\qquad\qquad\leq 2^{d-1}P\Big(\varrho\big(\Lambda(K)\big) \;>\; (\theta+\frac{\varepsilon}{2})\big|\Lambda(K)\big|\Big)\,.
\end{aligned} \tag{9.15}$$

We have used the symmetries of the model to say that the deviation probability appearing in the sum is the same for all the boxes and the translation invariance to come back to the box $\Lambda(K)$ (honestly, if $K$ is odd, then we rather end up with the box $e_1(1,\dots,1)+\Lambda(K)$). By proposition 5.4, the right-hand side of (9.15) goes to 0 as $K$ goes to $\infty$. □

It remains to prove that the other conditions appearing in the definition (9.1) of the good event $\mathcal{E}$ are typical. Unfortunately, we were only able to do so with the help of an additional assumption.

**Theorem 9.2.** *Suppose that $p$ is such that $\theta(p,\mathbb{Z}^d)>0$ and*

$$\liminf_{n\to\infty}\frac{1}{\ln n}\ln\ln\Big(\frac{1}{P_p\big(n\leq|C(0)|<+\infty\big)}\Big) \;>\; 0\,. \tag{9.16}$$

*For any $\varepsilon>0$, we have*

$$\limsup_{K\to\infty} P\left(\begin{array}{c}\text{there exists an open cluster } C^* \text{ in } \Lambda(4K)\\ \text{which touches the boundary } \partial^{\,in}\Lambda(4K) \text{ and}\\ \text{which satisfies } \big|C^*\cap\Lambda(2K)\big| \;\geq\; (\theta-\varepsilon)\,|\Lambda(2K)|\end{array}\right) \;=\; 1\,.$$

To complete the proof that $\theta(p_c,\mathbb{Z}^d)=0$, we should get rid of condition (9.16). The remaining part of the text is devoted to the proof of theorem 9.2.

# 10 Isoperimetry and concentration

We started the work in the subsection 8.5. Let us sum up this important initial step, with a slightly different notation (which is more adapted to the sequel). Let $\mathcal{C}_{\text{bd}}(4n)$ be the collection of the clusters of the percolation configuration restricted to $\Lambda(4n)$ which intersect both $\Lambda(2n)$ and $\partial^{\,in}\Lambda(4n)$. We have

$$\sum_{C\in\mathcal{C}_{\text{bd}}(4n)} \big|C\cap\Lambda(2n)\big| \,\geq\, \sum_{x\in\Lambda(2n)} 1_{\{\,x\longleftrightarrow\partial^{\,in}\Lambda(4n)\,\}} \,\geq\, \sum_{x\in\Lambda(2n)} 1_{\{\,x\longleftrightarrow\infty\,\}}\,. \tag{10.1}$$

By the spatial ergodic theorem, we have, with probability 1,

$$\lim_{n\to\infty}\frac{1}{|\Lambda(2n)|}\sum_{x\in\Lambda(2n)} 1_{\{\,x\longleftrightarrow\infty\,\}} \,=\, P\big(0\longleftrightarrow\infty\big) \,=\, \theta(p)\,. \tag{10.2}$$

The almost sure convergence implies the convergence in probability, therefore it follows from the inequalities (10.1) and (10.2) that

$$\forall\varepsilon>0 \qquad \lim_{n\to\infty} P\big(\mathcal{F}(n,\varepsilon)\big) \,=\, 1\,,$$

where $\mathcal{F}(n,\varepsilon)$ is the event defined by

$$\mathcal{F}(n,\varepsilon) \,=\, \Big\{\sum_{C\in\mathcal{C}_{\text{bd}}(4n)} \big|C\cap\Lambda(2n)\big| \,\geq\, \big(\theta(p)-\varepsilon\big)|\Lambda(2n)|\Big\}\,.$$

In other words, the event $\mathcal{F}(n,\varepsilon)$ is typical in the box $\Lambda(4n)$. This is the starting point from which we have to prove that the much stronger event defined in (9.1) is also typical. The challenge is to prove that, with high probability, the trace of the infinite cluster $C_\infty$ in the finite box $\Lambda(4n)$ contains one big cluster which itself contains a large fraction of the vertices of $\Lambda(2n)$ which are connected to $\infty$.

Let us consider a typical configuration in the box $\Lambda(4n)$ in the supercritical regime $p>p_c$ and let us consider the largest cluster $C^*$ in this box. When $n$ is large, with probability close to 1, the cardinality of $C^*$ is close to $\theta|\Lambda(4n)|$, and $C^*$ has density approximately $\theta$ in every intermediate box included in $\Lambda(4n)$, as demonstrated by the coarse-graining results of Pisztora (see proposition 8.3). In fact, the cluster $C^*$ is essentially the trace on $\Lambda(4n)$ of the infinite cluster. Pete [113] has proved that the isoperimetric profile of the infinite cluster basically coincides with the one of the original lattice. It is therefore to be expected that the cluster $C^*$ enjoys an isoperimetric inequality similar to the isoperimetric inequality of the lattice $\mathbb{Z}^d$, in the following sense. For any $\varepsilon>0$, in order to remove $\varepsilon n^d$ vertices from $C^*\cap\Lambda(2n)$, it is necessary to close a number of sites of $\Lambda(4n)$ which is of order $\varepsilon n^{d-1}$. From the concentration inequality for the Bernoulli product measure in $\Lambda(4n)$, we know that, up to the modification of a number of sites of order $n^{d/2}$, the configuration has to be typical. In a nutshell, the random fluctuations might change the states of at most $n^{d/2}$ sites, but to create a hole of order $\varepsilon n^d$ should require the closure of $\varepsilon n^{d-1}$ sites. This is unlikely to happen if $d/2<d-1$. Seductive as this argument may be, the road to make it rigorous is perilous and full of pitfalls. Of course, we cannot rely on the results of Pete and Pisztora, which require $p>p_c$.

# 11 The isoperimetric inequality

A fundamental tool destined to play an important role in our story is the discrete isoperimetric inequality in the lattice $\mathbb{Z}^d$. This is a whole subject in itself. However, we will need only a weak form, which we state and explain next. An isoperimetric inequality is an inequality between the cardinality of a subset of $\mathbb{Z}^d$ and some quantitative measure of its boundary. For a discrete set, there exist several possible definitions of the boundary. Let $A$ be a finite subset of $\mathbb{Z}^d$. We have already defined the outer boundary $\partial^{\,out} A$ of $A$ in (4.1). Of particular relevance is also the edge boundary $\partial^{\,edge} A$ of $A$, defined as

$$\partial^{\,edge} A \,=\, \big\{\, e = \langle x, y\rangle \in \mathbb{E}^d : x \in A,\, y \in \mathbb{Z}^d \setminus A \,\big\}\,.$$

We state next two classical isoperimetric inequalities in $\mathbb{Z}^d$.

**Theorem 11.1.** *There exist positive constants $c_{iso}^{out}(d)$, $c_{iso}^{edge}(d)$ which depend on the dimension $d$ only such that, for any finite subset $A$ of $\mathbb{Z}^d$,*

$$|A| \,\leq\, c_{iso}^{out}(d)\,\big|\partial^{\,out} A\big|^{\frac{d}{d-1}}\,,\qquad |A| \,\leq\, c_{iso}^{edge}(d)\,\big|\partial^{\,edge} A\big|^{\frac{d}{d-1}}\,. \tag{11.1}$$

Difficult related questions are to find the optimal isoperimetric constants, and to describe the sets realizing the equality. Bollobás and Leader [18] have studied edge-isoperimetric inequalities in the box $\Lambda(n)$, with free or toric boundary conditions, and they give essentially sharp bounds. A notable difficulty is that the sets realizing the equality in the isoperimetric inequality are not unique, and the structure of the optimal sets depends greatly on the value of $|A|$ (see [8] for a specific analysis in $\mathbb{Z}^3$). At least, if we look at optimal sets of increasing cardinality and if we rescale them properly, then they converge to a limiting shape, called the Wulff shape, which in this case turns out to be a cube. The recent work [102] studies the deviations of large optimal sets from the limiting Wulff shape. If one does not care about the precise values of the constants, nor on the optimal sets, then things become easier. Indeed, we have

$$\forall A \subset \mathbb{Z}^d \qquad |\partial^{\,out} A| \,\leq\, |\partial^{\,edge} A| \,\leq\, 2d\,|\partial^{\,out} A|\,.$$

Therefore each of the two isoperimetric inequalities stated in (11.1) imply the other, albeit with a different constant. So, from now onwards, we will simply denote by $c_{\text{iso}}$ a suitable isoperimetric constant. The weak isoperimetry takes the following form: there exists a positive constant $c_{\text{iso}}$ which depends on the dimension $d$ only such that, for any finite subset $A$ of $\mathbb{Z}^d$,

$$|A| \,\leq\, c_{\text{iso}}\,\big|\partial A\big|^{\frac{d}{d-1}}\,, \tag{11.2}$$

where $\partial A$ is any of the discrete boundaries of $A$. This inequality is ubiquitous in the study of random models defined on the lattice $\mathbb{Z}^d$. Let us mention two prominent instances where it played a crucial role:

- Agoston Pisztora's work on surface large deviations, see formula (2.1) in [37];
- Roberto Schonmann's work on the Ising model, see formula (2) of [119].

Agoston Pisztora provides a direct proof of (11.2) which is based on the Loomis-Whitney inequality [97] and he obtains the following constant:

$$c_{\text{iso}} = \left(\frac{2d}{2^d}\right)^{\frac{1}{d-1}}.$$

The Loomis-Whitney inequality rests on an induction and the classical Hölder inequality for a finite sum. The whole proof is thus completely elementary and it seems to be the most efficient strategy to prove (11.2) (see also [138] for a recent alternative proof of the Loomis-Whitney inequality). Another strategy to prove (11.2) is to derive it from the continuous isoperimetric inequality in $\mathbb{R}^d$. Starting with the finite discrete set $A$, we put at each of its vertices a unit cube, thereby obtaining a continuous subset $\overline{A}$ of $\mathbb{R}^d$. The volume of $\overline{A}$ is equal to the cardinality of $A$, and its perimeter, or the surface measure of its boundary, is equal to the cardinality of $\partial^{edge} A$. Applying the isoperimetric inequality in $\mathbb{R}^d$ to $\overline{A}$, we recover the discrete isoperimetric inequality (11.2). Now, the question is, which isoperimetric inequality in $\mathbb{R}^d$ are we appealing to? Indeed, the world of continuous isoperimetric inequalities is vast and complex. An important problem is to choose properly the class of the sets which are admissible in the inequality. Usually, the volume of the set is given by its Lebesgue measure (hence the set has to be Lebesgue measurable at least). The more delicate question is to define the perimeter of the set, which is a generalization of the $d-1$ surface measure of the boundary for a smooth set (for instance if the boundary is a $C^1$ hypersurface in $\mathbb{R}^d$). Roberto Schonmann [119] made appeal to the continuous Wulff theorem of Taylor [129] and he obtained the inequality (11.2) with the explicit constant

$$c_{\text{iso}} = \left(\frac{1}{2d}\right)^{\frac{d}{d-1}}.$$

However, the Wulff theorem of Taylor is quite sophisticated, it rests on geometric measure theory, and although this approach is completely valid, it is certainly not the simplest one. In fact, the simplest isoperimetric inequality in $\mathbb{R}^d$ is the one obtained by Loomis and Whitney [97]. So, if our aim is to get the discrete version, then the best way is to proceed as Pisztora did. Another classical approach for proving isoperimetric inequalities is to use symmetrization procedures. This technique was introduced by Steiner in 1838 [123]. Bollobás and Leader [18] rely on symmetrization techniques and partly on a continuous embedding. In fact, they prove an isoperimetric inequality for rectilinear bodies, these are sets that are finite unions of rectangular cuboids (also called right rectangular prisms). The advantage of this more sophisticated approach is that they obtain essentially sharp constants in the inequality, and also they describe true solutions of the isoperimetric problem. Another classical and powerful approach for proving isoperimetric inequalities in $\mathbb{R}^d$ is to rely on the famous Brunn-Minkowski inequality: see for instance [21]. Amazingly, the classical proof of the Brunn-Minkowski inequality is done first for rectilinear bodies. When dealing with rectilinear bodies, there is no difficulty in defining the perimeter.

## Part III

# The travel time in $\mathbb{Z}^d$ $\odot$

A first key tool for the proofs that are to come is the travel time, defined in section 12, where we give also an algorithm for its computation. A second key tool is a probabilistic control of a class of exploration algorithms, that we call genuine. This class and the control are presented in section 13. In section 14, we show how we can control the travel time with the help of the percolation function $\theta(p)$. In order to improve this control, in section 15, we try to obtain an estimate on the tail of the distribution of the sizes of the finite clusters, but unsuccessfully. We introduce the notion of $\gamma(n)$-large clusters in section 16. Finally, we revisit the travel time algorithm in order to design a shell exploration algorithm in section 17. We start to discuss the remote goal of constructing intertwined explorations.

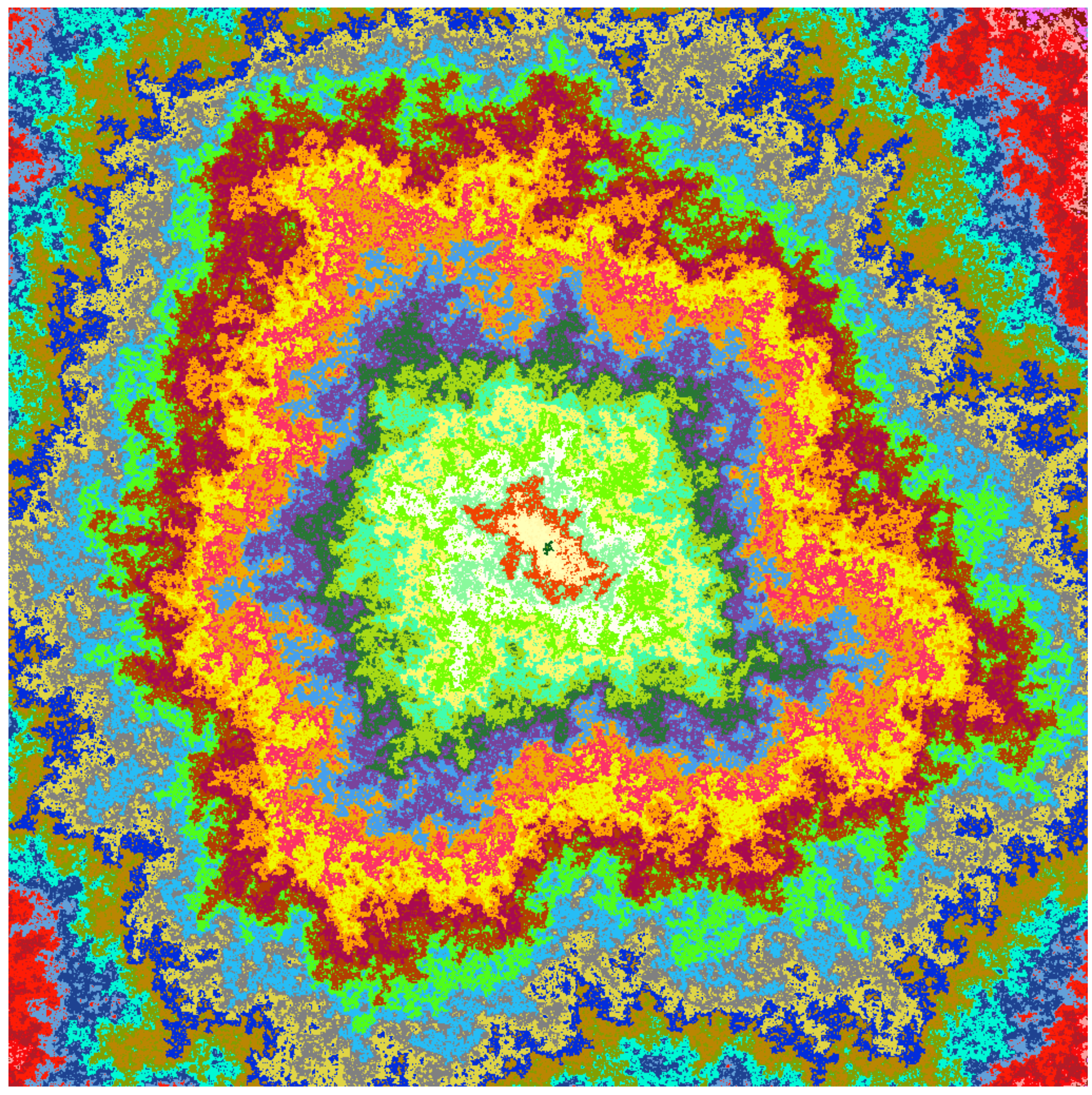

Figure 18: Bond percolation, shells around 0, p=0.46, 1024x1024 box.

# 12 The travel time to $\infty$ $\circ$

We continue here to investigate the relevant facts that we can prove starting with the hypothesis that $\theta(p)$ is positive. A fundamental fact, which has to be exploited, is that, when the density of the infinite cluster is positive, the infinite cluster is not far from any given point of the lattice, in the sense that the travel time to the infinite cluster is small with very high probability. We shall try to quantify this statement, but before that, we explain what we mean by the travel time. The travel time is similar to the passage time in first passage percolation, yet we will use a slight variant, more adapted to the proofs that are to follow. This requires the introduction of some formal definitions.

## 12.1 Definition of the travel time

Let $D$ be a subset of $\mathbb{Z}^d$ and let $A$ be a subset of $D$. We define its internal boundary $\partial_D^{in} A$ in $D$ and its external boundary $\partial_D^{out} A$ in $D$ by

$$\partial_D^{in} A = \big\{ x \in A : \exists y \in D \setminus A \quad |x-y| = 1 \big\},$$
$$\partial_D^{out} A = \big\{ x \in D \setminus A : \exists y \in A \quad |x-y| = 1 \big\},$$

where $|\cdot|$ is the Euclidean norm. Recall that a path $z_0, \dots, z_r$ in $D$ is a sequence of sites of $D$ such that each site is a neighbour of its predecessor:

$$\forall i \in \{0, \dots, r-1\} \qquad |z_{i+1} - z_i| = 1\,.$$

For $x, y$ in $D$, we define the graph distance $d_D(x,y)$ between $x$ and $y$ in $D$ as

$$d_D(x,y) = \inf \Big\{ r \geq 0 : z_0, \dots, z_r \text{ path in } D \text{ from } z_0 = x \text{ to } z_r = y \Big\}, \quad (12.1)$$

and the travel time $T_D(x,y)$ between $x$ and $y$ in $D$ by

$$T_D(x,y) = \inf \Big\{ \sum_{i=0}^{r-1} 1_{z_i \text{ closed}} : z_0, \dots, z_r \text{ path in } D \text{ from } z_0 = x \text{ to } z_r = y \Big\}. \quad (12.2)$$

We make the convention that $T_D(x,x) = 0$ for any $x \in D$ (this is compatible with the previous definition, if we consider that $z_0 = x$ is a path from $x$ to $x$; the length $r$ of this path is 0 and the sum in (12.2) is empty). With our definitions, we have $T_D(x,y) \leq d_D(x,y)$ for any $x, y \in D$. For $x$ in $D$ and $B$ a subset of $D$, we define the travel time between $x$ and $B$ in $D$ by

$$T_D(x,B) = \inf \big\{ T_D(x,y) : y \in B \big\},$$

and symmetrically the travel time between $B$ and $x$ in $D$ by

$$T_D(B,x) = \inf \big\{ T_D(y,x) : y \in B \big\}.$$

For $A, B$ two subsets of $D$, we define the travel time between $A$ and $B$ in $D$ by

$$T_D(A,B) = \inf \big\{ T_D(x,y) : x \in A,\ y \in B \big\}. \quad (12.3)$$

The next problem is, how do we compute the travel time $T_D(A,B)$?

## 12.2 The travel time algorithm

The mathematical definition of $T_D(A, B)$ given in (12.3) looks very simple, but this is rather fallacious. A fundamental goal of the first passage percolation model is to understand the behavior of these travel times between distant sets for more general distributions. For the time being, we are interested in the actual computation of $T_D(A, B)$. When $D$ is finite, the most simple algorithm consists in enumerating all the paths from $A$ to $B$ in $D$, computing their travel times and sorting out those which realize the infimum. However, this is catastrophic, not only in terms of the computational cost, but also from a probabilistic point of view, because this algorithm would reveal all the sites of $D$, even those whose state is completely irrelevant for the value of $T_D(A, B)$. We present next a natural algorithm which is already much better. This algorithm implements the most naive approach to compute the travel time. It starts from the set $A$, and it explores all the clusters of the sites of $A$. Exploring a cluster means exploring all the open sites of the cluster, but also the closed sites at the boundary of the cluster. After this first phase, the algorithm has computed the set of the sites which are at travel distance 0 from $A$. In the second phase, the algorithm explores all the clusters of the unexplored sites which are neighbours of a site belonging to the set of the sites which are at travel distance 0 from $A$. After this second phase, the algorithm has computed the set of the sites which are at travel distance 1 from $A$. The algorithm goes on this way until $B$ is hit or $D$ is exhausted (in case the sets $A$ and $B$ are not in the same connected component of $D$). The algorithm updates iteratively four sets of sites, denoted by $A, O, C$ and $W$. The sets $A(k)$, $O(k)$, $C(k)$ and $W(k)$ are the sets when exactly $k$ sites have been explored and their contents are the following:

- The set $A(k)$: the set of the active sites, which are to be explored.
- The set $O(k)$: the set of the open sites, which have already been explored.
- The set $C(k)$: the set of the closed sites, which have already been explored.
- The set $W(k)$: the set of the waiting sites, which wait to become active.

The sets $A(k)$, $O(k)$, $C(k)$, $W(k)$ are all included in $D$. The algorithm computes also an increasing sequence of integers

$$K_0 = 0 \,<\, K_1 \,<\, \cdots \,<\, K_T$$

corresponding to the regeneration times of the set $A(k)$ of the active sites. For $t \in \{\, 0, \dots, T-1 \,\}$, the set of the sites which are at travel distance exactly equal to $t$ from $A$ inside $D$ is the set

$$\big(C(K_{t+1}) \cup O(K_{t+1})\big) \setminus \big(C(K_t) \cup O(K_t)\big)\,.$$

The algorithm terminates at step $K$ and returns the travel time $T = T_D(A, B)$. Its pseudocode is presented as Algorithm 12.1 and it uses the following notation: for $x \in D$, we define $\mathcal{N}(x, D)$ as the set of the neighbours of $x$ in $D$, i.e.,

$$\mathcal{N}(x, D) \,=\, \big\{\, y \in D : |x - y| = 1 \,\big\}\,. \tag{12.4}$$

**Algorithm 12.1** Travel time algorithm

**Require:** a non-empty starting set $A$, a target set $B$, an authorized domain $D$
**Ensure:** $T$ is the travel time from $A$ to $B$ inside $D$

```
  A(0) ← ∅, O(0) ← ∅, C(0) ← ∅, W(0) ← A
  K_0 ← 0, t ← 0, k ← 0
                                  ▷ t is the number of the shell currently explored
                                      ▷ k is the number of sites already explored
  repeat
      A(k) ← W(k), W(k) ← ∅          ▷ Regeneration of the set of active sites
      repeat
          Pick up x ∈ A(k)
          if x is closed then
              A(k+1) ← A(k) \ {x}
              O(k+1) ← O(k)
              C(k+1) ← C(k) ∪ {x}
              W(k+1) ← W(k) ∪ N(x, D)
          else if x is open then
              A(k+1) ← A(k) ∪ N(x, D) \ ({x} ∪ O(k) ∪ C(k))
              O(k+1) ← O(k) ∪ {x}
              C(k+1) ← C(k)
              W(k+1) ← W(k)
          end if
          k ← k + 1
      until A(k) = ∅ or x ∈ B
      t ← t + 1, K_t ← k
      W(k) ← W(k) \ (O(k) ∪ C(k))  ▷ Remove explored sites from waiting set
  until W(k) = ∅ or x ∈ B
  T ← t − 1, K ← K_t
  return T, O(K), C(K)
```

Each algorithm which explores the percolation configuration performs the construction of a specific random process. Naturally, there exists a mathematical definition which is the counterpart to each such algorithm. In fact, the travel time algorithm 12.1 does much more than the mere computation of $T_D(A, B)$. Indeed, it explores a specific random subset of $D$, and it builds four relevant sequences of random sets. We give next the corresponding mathematical construction of the sequence of sets $A(k)$, $O(k)$, $C(k)$, $W(k)$ and the sequence of integers $K_t$ computed by this algorithm. The index $k$ is the number of sites explored by the algorithm. The integer $t$ is the travel time between $A$ and the sites which are currently explored. The integer $K_t$ is the number of sites which are at travel time from $A$ strictly less than $t$. Our goal is to derive interesting probabilistic estimates on the travel times. In its mathematical version, the algorithm is driven by a sequence of i.i.d. Bernoulli random variables $(X_k)_{k\geq 1}$ with parameter $p$. The random variable $X_k$ will be used to decide the state of the $k$-th site visited by the algorithm. We will build four sequences of random sets $A(k)$, $O(k)$, $C(k)$, $W(k)$ and a random sequence of integers $K_t$, which are deterministic functions of the sequence $(X_k)_{k\geq 1}$. Moreover, for any $k \geq 1$, the random sets $A(k)$, $O(k)$, $C(k)$, $W(k)$ are measurable with respect to $X_1, \dots, X_k$. Also, for any $k \geq 1$ and $t \geq 0$, the event $\{ K_t = k \}$ is measurable with respect to $X_1, \dots, X_k$. We give next the precise mathematical construction.

• **Initialization:** We set $W(0) = O(0) = C(0) = \varnothing$ and $A(0) = A$, $K_0 = 0$.

• **Iteration:** Suppose that the sets $A(k), O(k), C(k), W(k)$ and the integer $K_t$ are built and let us explain how to build the sets $A(k+1)$, $O(k+1)$, $C(k+1)$, $W(k+1)$ and the integer $K_{t+1}$. We pick up an element $x_{k+1}$ of $A(k)$ according to some deterministic procedure. The site $x_{k+1}$ has not been explored previously, and its state will be decided by the random variable $X_{k+1}$.
Site exploration. We consider two cases, according to the value of $X_{k+1}$:
*First case:* $X_{k+1} = 0$. The site $x_{k+1}$ is declared closed, and we set

$$O(k+1) \,=\, O(k)\,, \quad C(k+1) \,=\, C(k) \cup \{ x_{k+1} \}\,,$$

$$A_* \,=\, A(k) \setminus \{ x_{k+1} \}\,, \quad W_* \,=\, W(k) \cup \mathcal{N}(x_{k+1}, D)\,,$$

where the notation $\mathcal{N}(x, D)$ was introduced in (12.4).
*Second case:* $X_{k+1} = 1$. The site $x_{k+1}$ is declared open, and we set

$$O(k+1) \,=\, O(k) \cup \{ x_{k+1} \}\,, \quad C(k+1) \,=\, C(k)\,,$$

$$A_* \,=\, A(k) \cup \mathcal{N}(x_{k+1}, D) \setminus \Big( \{ x_{k+1} \} \cup O(k) \cup C(k) \Big)\,, \quad W_* \,=\, W(k)\,.$$

Active and waiting sites update. We consider two cases, according to whether there are still some unexplored sites at travel distance $t$ from $A$ or not:
⋆ $A_* \neq \varnothing$. We set $A(k+1) = A_*$, $W(k+1) = W_*$.
⋆ $A_* = \varnothing$. We set $A(k+1) = W_* \setminus \big(O(k) \cup C(k)\big)$, $W(k+1) = \varnothing$, $K_{t+1} = k+1$.
The algorithm terminates if either $x_{k+1}$ belongs to the target set $B$ or the set $A(k+1)$ is empty. Otherwise, we start a new iteration.

• **Termination:** The algorithm terminates with the time value $T$ at step $K_{T+1}$. It returns the sequence of sets $O(k)$, $C(k)$, $0 \leq k \leq K_{T+1}$ and the sequence of integers $K_t$, $0 \leq t \leq T+1$.

The algorithm terminates either when both $A(k)$ and $W(k)$ are empty, i.e., when there are no more sites to explore in $D$, or when a site $x_k$ belonging to the target set $B$ is explored. Since $O(k)$, $C(k)$, $A(k)$ are disjoint subsets of $D$ and the sequence of sets $O(k) \cup C(k), k \geq 0$, is increasing, necessarily $A(k)$ is empty after at most $|D|$ steps and the termination is guaranteed.

As we see, the precise mathematical description is a bit cumbersome, so we will not repeat it for the other more complicated algorithms that are to follow. However, the mathematical description of the algorithm has two great benefits, compared to the mere definition of $T_D(A, B)$ given in (12.3). First, it leads to extremely useful probabilistic estimates, and second it helps to introduce very relevant random subsets.

## 12.3 Definition of the shells

The travel time algorithm proceeds by exploring in order the sites which are at travel time 0 from $A$, then those which are at travel time 1 from $A$, and so on, until the set $B$ is hit. These successive sets are random subsets of $D$, and they will play a key role in the sequel. So we call them the shells around $A$, and we define them precisely in this subsection. For $A$ a subset of $D$, we define the shell around $A$ in $D$ by

$$\text{Shell}\,(A,D)\,=\,\big\{\,x\in D: T_D(A,x)=0\,\big\}\,, \tag{12.5}$$

and for $t\geq 0$ an integer, we define the shell number $t$ around $A$ in $D$ by

$$\text{Shell}\,(A,D,t)\,=\,\big\{\,x\in D: T_D(A,x)=t\,\big\}\,, \tag{12.6}$$

and the closed ball of radius $t$ for the travel time $T_D$ by

$$\mathcal{B}(A,D,t)\,=\,\big\{\,x\in D: T_D(A,x)\leq t\,\big\}\,. \tag{12.7}$$

Obviously, for any set $A$, we have $\text{Shell}\,(A,D)=\text{Shell}\,(A,D,0)$ and the successive shells around $A$ are pairwise disjoint. With our definitions, for any $x$ in $D$, the site $x$ belongs to $\text{Shell}\,(\{\,x\,\},D)$, and more generally

$$\forall x\in D\quad\forall t\geq 0\qquad\big\{\,y\in D: d_D(x,y)\leq t\,\big\}\,\subset\,\mathcal{B}(\{\,x\,\},D,t) \tag{12.8}$$

(the graph distance $d_D$ was defined in (12.1)). Let $A$ be any subset of $D$. We have an identity for balls, which reads

$$\forall s\geq 0\quad\forall t\geq 0\qquad\mathcal{B}(A,D,s+t)\,=\,\mathcal{B}\big(\mathcal{B}(A,D,s),D,t\big)\,, \tag{12.9}$$

and an identity for shells, which is

$$\forall s\geq 0\quad\forall t\geq 0\qquad\text{Shell}\,(A,D,s+t)\,=\,\text{Shell}\,\big(\mathcal{B}(A,D,s),D,t\big)\,. \tag{12.10}$$

Unfortunately, the functional $T_D(x,y)$ is not symmetric: indeed the starting point $x$ and the ending point $y$ do not play the same roles in the definition (12.2). The advantage of our definition is that the triangle inequality holds:

$$\forall x,y,z\in D\qquad T_D(x,y)\,\leq\,T_D(x,z)\,+\,T_D(z,y)\,,$$

and furthermore, we have

$$\begin{aligned}&\forall A,B\subset D\quad\forall t\in\{\,0,\dots,T_D(A,B)\,\}\\ &\qquad\qquad\qquad\qquad T_D(A,B)\,=\,t+T_D(\mathcal{B}(A,D,t),B)\,.\end{aligned} \tag{12.11}$$

If we had included both extremities of the path when defining the travel time between two sites, then we would have the symmetry of $T_D$ and the triangle inequality, but not the identity (12.11). The set $\text{Shell}\,(A,D)$ is the union of the clusters of the sites of $A$, together with their outer boundaries, in the percolation configuration restricted to $D$. We define the outer boundary $\partial_D^{out}A$ of $A$ as

$$\partial_D^{out}A\,=\,D\cap\partial^{out}A\,.$$

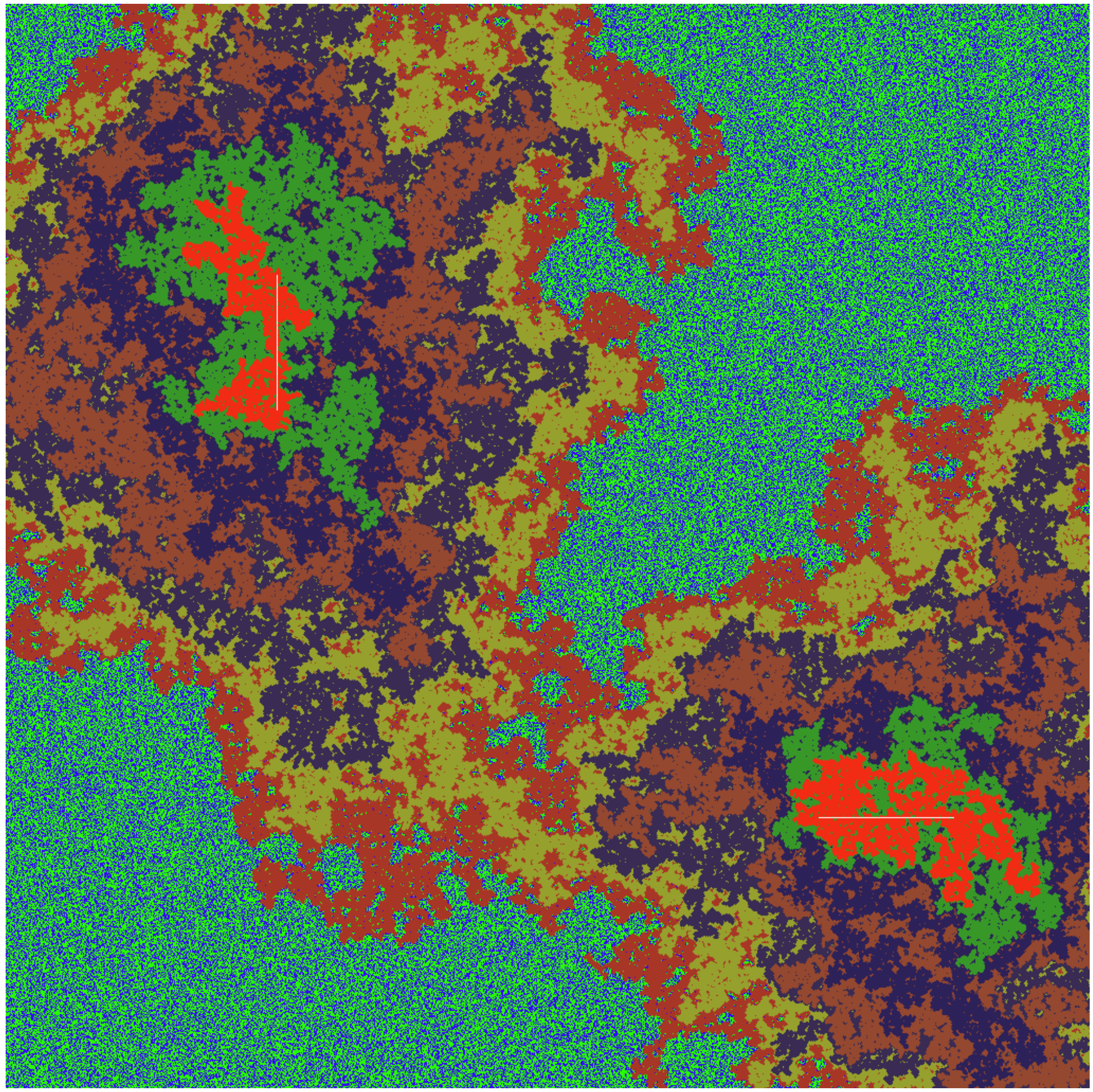

Figure 19: Site percolation in a box of size $1024 \times 1024$, $p = 0.569$. The first 7 shells around the 2 white segments. Open sites are green, closed sites are blue.

For $x \in D$, we define $C(x, D)$ as the open cluster of $x$ in $D$, it consists of all the sites $y$ which are connected to $x$ by an open path whose sites all belong to $D$. We make also the convention that, for a closed site $x$ of $D$, we have $C(x, D) = \varnothing$ and $\partial_D^{out} C(x, D) = \{ x \}$. This is coherent with the previous definitions. Indeed, we have $T_D(x, x) = 0$ because $x$ is a path of length 0 which is admissible in the definition of $T_D(x, x)$, however it does not constitute an open path whenever $x$ is closed. Furthermore, we define

$$\forall x \in D \qquad \overline{C}(x, D) \,=\, C(x, D) \cup \partial_D^{out} C(x, D) \,. \tag{12.12}$$

With this notation, the set Shell $(A, D)$ can be written as

$$\text{Shell}\,(A, D) \,=\, \bigcup_{x \in A} \overline{C}(x, D) \,.$$

# 13 Genuine site exploration algorithms

As we have seen in the subsection 12.2, the travel time algorithm is driven by a sequence of i.i.d. Bernoulli random variables $(X_k)_{k\geq 1}$ with parameter $p$. Let us introduce the associated partial sums $(S_k)_{k\geq 1}$, defined as

$$\forall k \geq 1 \qquad S_k = X_1 + \cdots + X_k \,.$$

From the mathematical description of the algorithm, for any $k \geq 1$, we have

$$\big|O(k)\big| \,=\, S_k \,, \quad \big|C(k)\big| \,=\, k - S_k \,. \tag{13.1}$$

Thus, for a fixed (deterministic) $k \geq 1$, the number of the open sites discovered during the first $k$ steps of the algorithm follows the binomial distribution with parameters $k$, $p$. Yet the algorithm terminates after a random number of iterations $K$. Upon termination, the random variable $\big|O(K)\big| = S_K$ has a much more complicated distribution. However, it is still the case that $\big|O(K)\big|$ is equal to some value found along the sequence $(S_k)_{k\geq 1}$. We have also an easy upper bound on $K$, given by the number of sites in $D$. As a consequence, any probabilistic estimate on the finite sequence $(S_k)_{1\leq k\leq |D|}$ will yield some information on the behaviour of the exploration algorithm. The most useful deviation inequality for the binomial distribution is the famous Hoeffding inequality [70]:

$$\forall k \geq 1 \quad \forall t \geq 0 \qquad P\Big(\big|S_k - kp\big| \geq t\Big) \,\leq\, 2\exp\Big(-\frac{2t^2}{k}\Big)\,. \tag{13.2}$$

The identities (13.1) yield that

$$\begin{aligned} S_k - kp \,&=\, \big|O(k)\big| - \Big(\big|O(k)\big| + \big|C(k)\big|\Big)p \\ &=\, \big|O(k)\big|(1-p) - \big|C(k)\big|p \,=\, p(1-p)\Big(\frac{\big|O(k)\big|}{p} - \frac{\big|C(k)\big|}{1-p}\Big)\,. \end{aligned} \tag{13.3}$$

Substituting (13.3) in (13.2), and making the change of variable $t \to p(1-p)t$, we obtain

$$\forall k \geq 1 \quad \forall t \geq 0 \qquad P\Big(\Big|\frac{\big|O(k)\big|}{p} - \frac{\big|C(k)\big|}{1-p}\Big| \geq t\Big) \,\leq\, 2\exp\Big(-\frac{2}{k}t^2(p(1-p))^2\Big)\,.$$

We convert this in a control over the sequence until the index $|D|$ with the help of a simple union bound:

$$\begin{aligned} \forall t \geq 0 \qquad P\Big(\exists k \in \{\,1,\dots,|D|\,\} \quad \Big|\frac{\big|O(k)\big|}{p} - \frac{\big|C(k)\big|}{1-p}\Big| \geq t\Big) \\ \leq\, 2|D|\exp\Big(-\frac{2}{|D|}t^2(p(1-p))^2\Big)\,. \end{aligned}$$

Let us generalize slightly this inequality. For $A$ a finite subset of $\mathbb{Z}^d$, we define

$$S(A) \,=\, \frac{1}{p}\Big|\big\{\, x \in A : x \text{ is open}\,\big\}\Big| - \frac{1}{1-p}\Big|\big\{\, x \in A : x \text{ is closed}\,\big\}\Big|\,. \tag{13.4}$$

Let us denote by $\mathcal{E}$ the random set of the sites explored by the travel time algorithm upon termination. We have

$$\forall \lambda \geq 0 \qquad P\big(\big|S(\mathcal{E})\big| \geq \lambda\big) \;\leq\; 2|D| \exp\Big(-\frac{2}{|D|}(p(1-p))^2\lambda^2\Big)\,. \tag{13.5}$$

This estimate is in fact available for any genuine exploration algorithm. By a genuine algorithm, we mean an algorithm which exploits only the information discovered by himself: at each step, the choice of the next site to be explored has to be a deterministic function of the set of the sites already explored (including the states of these sites). Notice that $\mathcal{E}$ is a random subset of $D$, and of course the inequality (13.5) has no chance to be true for an arbitrary random set. The key point is that, in a genuine exploration algorithm, the state of the next site to be explored is always independent of the past history of the exploration. Let us take the time to define formally a genuine exploration algorithm. A typical exploration algorithm updates iteratively three sets:

- The set $E(k)$: these are the sites which have been explored.
- The set $O(k)$: these are the open sites which have been explored.
- The set $C(k)$: these are the closed sites which have been explored.

Initially, we have $E(0) = O(0) = C(0) = \varnothing$. The definition of the algorithm rests on a decision rule for choosing the next site to be explored at step $k$. This choice depends on the history of the algorithm until step $k$, plus possibly some additional source of randomness. We model this source of randomness by a sequence $(U_k)_{k\geq 0}$ of i.i.d. random variables, which are uniformly distributed over $[0,1]$, and which is independent of the percolation process. So, in its most general form, the decision rule at step $k$ is given by a map

$$\phi_k : D^k \times \{\,0,1\,\}^k \times [0,1] \to D\,.$$

For the first step, the map $\phi_0$ is simply a map from $[0,1]$ to $D$ and $x_1 = \phi_0(U_0)$ is the first site that the algorithm will explore. After completion of the step $k$, the algorithm will have explored the sites $x_1,\dots,x_k$ and $E(k) = \{\,x_1,\dots,x_k\,\}$. The states of these sites is given by the sequence of random variables $X_1,\dots,X_k$, so that

$$O(k) \;=\; \bigcup_{1\leq i\leq k, X_i=1} \{\,x_i\,\}\,, \qquad C(k) \;=\; \bigcup_{1\leq i\leq k, X_i=0} \{\,x_i\,\}\,.$$

The site explored at step $k+1$ is then given by

$$x_{k+1} \;=\; \phi_k\big(x_1,\dots,x_k,X_1,\dots,X_k,U_k\big)\,.$$

For a non-random algorithm, the function $\phi_k$ does not depend on the random variable $U_k$. An interesting subcategory of exploration algorithms corresponds to decision rules $\phi_k$ which depend only on the sets $E(k), O(k), C(k)$ (and not on the order in which the sites of $E(k)$ were revealed). We consider also exploration algorithms which have a random termination time $T$. The termination time $T$ will always be a stopping time with respect to the sequence of random variables $(X_k)_{k\geq 1}$, i.e., for any $t\geq 0$, the event $\{\,T=t\,\}$ is in the $\sigma$-field generated by the variables $X_1,\dots,X_t$. For a cluster exploration algorithm, we will naturally take for $T$ the time when the exploration of the desired clusters is complete. The total region explored by the algorithm is $\mathcal{E} = E(T)$. A typical example is the following. We consider an initial set $W\subset D$, and we explore all the clusters of the sites of $W$. The algorithm terminates when all these clusters have been explored. In this example, the order in which the clusters of the sites of $W$ are explored is encoded in the decision rule $\phi_k$, $k\geq 0$. For instance, we can order the sites of $W$ according to some deterministic order. We start the algorithm by picking up the first site of $W$, and we explore its cluster. Once this first cluster has been explored, we pick up the second site of $W$. If this site has not been explored, we explore its cluster, otherwise we continue with the third site, until all the sites of $W$ are exhausted. This exploration algorithm is genuine, as long as the set $W$ is a fixed deterministic set. This is not any more the case if we consider an algorithm which explores the clusters of a random set of sites $\mathcal{W}$ (as we will do later). In this case, it is much more difficult to develop a control on the region explored by the algorithm, because the starting set might incorporate some information on the configuration that is to be explored, and this information will alter in an irreversible way the distribution of the sites of the clusters of $\mathcal{W}$. For instance, say that $\mathcal{W}$ is a site belonging to the largest cluster in a box $\Lambda$ (or one them if there are several). It is to be expected that, due to this positive information, the proportion of open sites visited when the exploration starts from $\mathcal{W}$ will be much higher than for a genuine exploration.

Now that we have defined precisely what we mean by a genuine exploration algorithm, we can control the tail distribution of the random variable $S(\mathcal{E})$ as follows. For any $\lambda>0$, we decompose the probability

$$P\big(|S(\mathcal{E})|\geq\lambda\big) \;=\; \sum_{1\leq t\leq |D|} P\big(|S(\mathcal{E})|\geq\lambda, T=t\big)\,.$$

The sum starts at $t=1$ because if $T=0$, then $\mathcal{E}=\varnothing$ and the probability inside the sum is 0 for $t=0$. On the event $\{\,T=t\,\}$, we have $\mathcal{E}=E(t)$ and

$$S_t - tp \;=\; p(1-p)S(E(t))\,,$$

whence

$$P\big(|S(\mathcal{E})|\geq\lambda\big) \;\leq\; \sum_{1\leq t\leq |D|} P\big(|S_t-tp|\geq p(1-p)\lambda, T=t\big)\,.$$

Applying Hoeffding's inequality (13.2), we get

$$P\big(|S(\mathcal{E})|\geq\lambda\big) \;\leq\; \sum_{1\leq t\leq |D|} 2\exp\Big(-\frac{2}{t}(p(1-p))^2\lambda^2\Big)\,,$$

which implies inequality (13.5). Of course, the factor $|D|$ in the denominator inside the exponential in (13.5) is quite large and undesirable, and there is certainly room for improvement. Suppose that we have in addition a quantitative control on the termination time. We could then bound the probability as follows:

$$
\begin{aligned}
P\big(|S(\mathcal{E})| \geq \lambda\big) \;\leq\; & \sum_{1\leq v\leq t} P\big(|S_v - vp| \geq p(1-p)\lambda, T=v\big) \;+\; P(T>t) \\
& \leq\; 2t\exp\Big(-\frac{2}{t}(p(1-p))^2\lambda^2\Big) \;+\; P(T>t)\,.
\end{aligned}
$$

Unfortunately, to control the tail of the termination time for a cluster exploration amounts to control the tail of the cluster size distribution, and we did not succeed in obtaining any quantitative control at the critical point $p_c$ (see the disappointing section 15).

Another category of tools at our disposal are the classical martingale results. Indeed, the sequence $(S_k - kp)_{k\geq 1}$ is a martingale, the stopping time $T$ is bounded by $|D|$ and $S_T - Tp = p(1-p)S(\mathcal{E}(T))$. Applying the optional sampling theorem, we get that $E\big(S(\mathcal{E}(T))\big) = 0$.

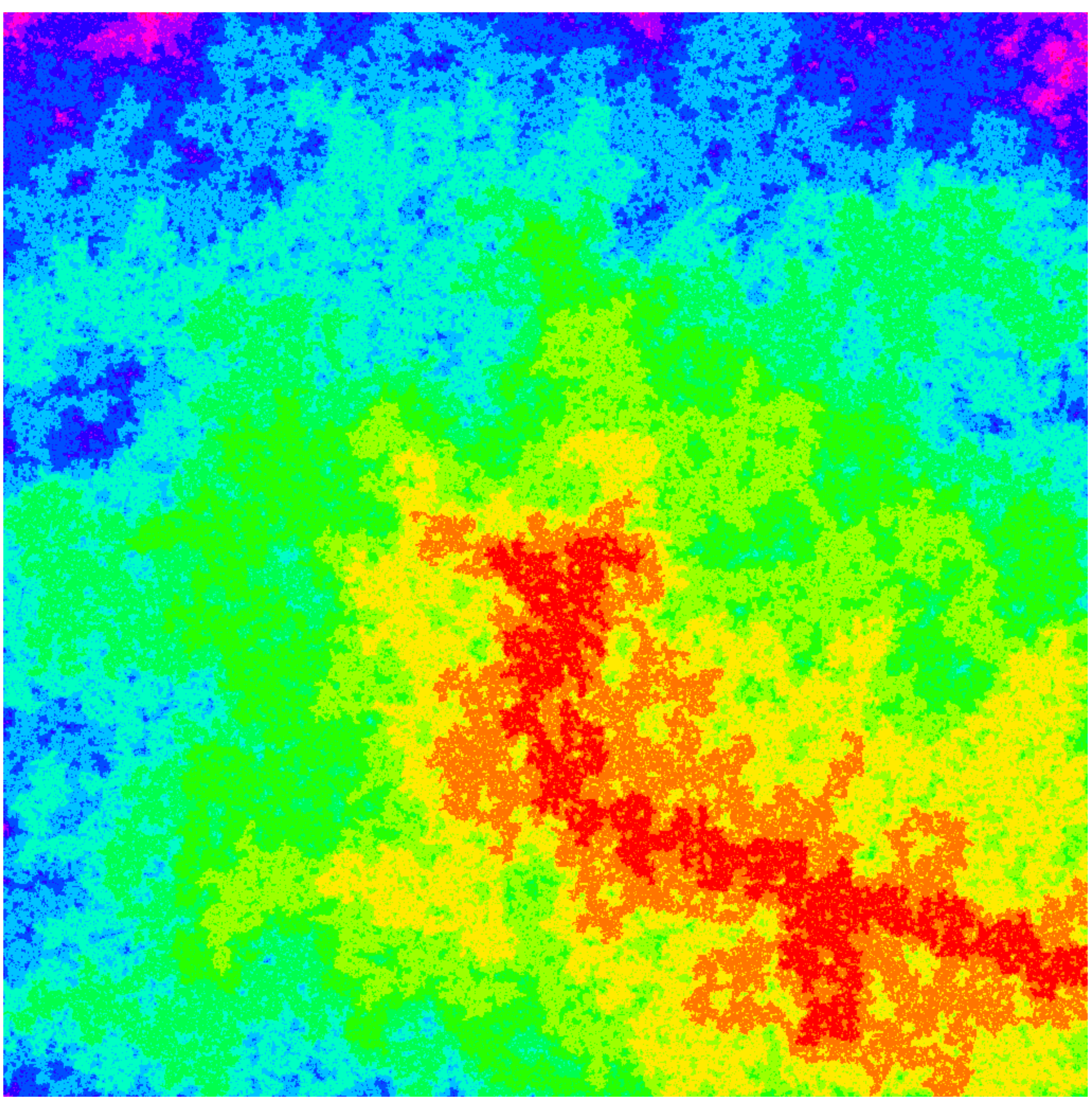

Figure 20: The shell exploration algorithm starting from the center of the box

# 14 Control of the travel time

In this section, we develop various quantitative estimates on the travel time from a site to the infinite cluster, or to the boundary of a domain, under the hypothesis that $\theta(p)$ is positive.

## 14.1 An application of the BK inequality

A routine application of the BK inequality shows that the travel time to exit a finite domain is dominated by a geometric distribution with parameter $1-\theta(p)$, and this yields also a control on the travel time until the infinite cluster.

**Lemma 14.1.** *Let $D$ be a finite subset of $\mathbb{Z}^d$ and let $x \in D$. We have*

$$\forall t \geq 1 \qquad P\big(\,T_D(x, \partial^{\,in} D)\;\geq\;t\,\big)\;\leq\;\big(1-\theta(p)\big)^t\,.$$

*Proof.* Let $D$ be a finite subset of $\mathbb{Z}^d$ and let $x \in D$. The event $\{\,x \longleftrightarrow \infty\,\}$ is included in the event $\{\,x \longleftrightarrow \partial^{\,in} D\,\}$, thus

$$P\big(T_D(x, \partial^{\,in} D) = 0\big)\;\geq\;P\big(x \longleftrightarrow \infty\big)\;=\;\theta(p)\,,$$

or, passing to the complementary event,

$$P\big(T_D(x, \partial^{\,in} D) \geq 1\big)\;\leq\;1-\theta(p)\,.$$

To deal with the general case where $t \geq 1$, let us run the travel time algorithm 12.1 in the domain $D$, starting from the set $A = \{\,x\,\}$ and with target set $B = \partial^{\,in} D$. This algorithm explores successively the shells around $x$ in $D$. Suppose that $T_D(x, \partial^{\,in} D) \geq t$. In this case, the first $t$ shells, $\text{Shell}\,(x, D, s)$, $0 \leq s < t$, do not intersect $\partial^{\,in} D$. Moreover, for $0 \leq s < t$, the inner boundary of the closed ball $\mathcal{B}(x, D, s)$ for the travel time $T_D$ (defined in (12.7)), namely the set $\partial^{\,in}\mathcal{B}(x, D, s)$, contains a subset of closed sites which disconnects $x$ from $\partial^{\,in} D$. Let us denote by $\mathcal{D}(s)$ one such set of sites. The set $\mathcal{D}(s)$ realizes the event $\{\,x \not\longleftrightarrow \partial^{\,in} D\,\}$. By construction, the sets $\mathcal{D}(s)$, $0 \leq s < t$, are pairwise disjoint. Thus we have

$$\Big\{\,T_D(x, \partial^{\,in} D)\;\geq\;t\,\Big\}\;\subset\;\Big\{\,x \not\longleftrightarrow \partial^{\,in} D\;\text{ occurs disjointly } t \text{ times}\,\Big\}\,.$$

Applying the BK inequality (see for instance [63]), we conclude that

$$\begin{aligned} P\big(\,T_D(x, \partial^{\,in} D)\;\geq\;t\,\big)\;&\leq\;P\big(\,x \not\longleftrightarrow \partial^{\,in} D\,\big)^t \\ &\leq\;P\big(\,x \not\longleftrightarrow \infty\,\big)^t\;\leq\;\big(1-\theta(p)\big)^t\,, \end{aligned}$$

as required. ∎

Our next task will be to improve this estimate, and for this we introduce in the next subsection the function $\theta(n, p)$, which coincides with $\theta(p)$ when $n$=1.

## 14.2 The function $\theta(n,p)$

For $n \geq 1$, we define

$$\theta(n,p) \;=\; \inf\big\{\, P_p(E \leftrightarrow \infty) : E \text{ subset of } \mathbb{Z}^d, |E| = n \,\big\}\,.$$

In particular, $\theta(1,p)$ is equal to the percolation probability $\theta(p)$. We shall show below that

$$\theta(p) > 0 \qquad \Longrightarrow \qquad \lim_{n\to\infty} \theta(n,p) \;=\; 1\,. \tag{14.1}$$

A slightly weaker (but almost equivalent) form of this statement is proved by Martineau and Tassion, see lemma 3.6 in [103]. Their goal was to replace the use of the seeds in the Grimmett-Marstrand construction by long self-avoiding paths. Their argument is a compactness argument in the space of self-avoiding paths. Here we use a more refined packing argument and we provide a control on the speed of convergence in the limit (14.1).

**Proposition 14.2.** *Let $p$ be such that $0 < \theta(p) < 1$. Let $\phi$ be the function defined on $[0,1]$ by $\phi(t) = t\ln t$. We have*

$$\forall n \geq \max\left(\Big( - 2\frac{13^d}{\theta}\ln\frac{\theta}{2}\Big)^{d+1}, \Big( - 2\frac{98^d}{\theta^2}\ln\big(p(1-p)\big)\Big)^{d(d+1)}\right) \tag{14.2}$$

$$1-\theta(n,p) \;\leq\; \frac{2}{\ln(1-\theta)}\phi\Big(P_p\big(n^{1/(d+1)} \leq |C(0)| < +\infty\big)\Big)\,. \tag{14.3}$$

*Proof.* Let $n \geq 1$ and let $E$ be a subset of $\mathbb{Z}^d$ of cardinality $n$. Let $m$ be another integer such that $1 \leq m \leq n$. An $m$-packing of $E$ is a subset $F$ of $E$ satisfying

$$\forall y,z \in F \qquad y \neq z \quad \Longrightarrow \quad |y-z|_\infty > m\,.$$

The $m$-packing number of $E$ is the maximal cardinality of an $m$-packing of $E$, we denote it by $N^{\text{pack}}(E,m)$. Let $F$ be an $m$-packing of $E$ having maximal cardinality. Necessarily, each point of $E$ is at distance less than or equal to $2m$ from a point of $F$:

$$\forall x \in E \quad \exists\, y \in F \qquad |x-y|_\infty \leq 2m\,.$$

Indeed, if there existed a point of $E$ which didn't satisfy the above property, we could add it to $F$ and we would obtain an $m$-packing of $E$ having a strictly larger cardinality. Therefore, we have the following inequality:

$$n = |E| \;\leq\; (4m+1)^d N^{\text{pack}}(E,m)\,. \tag{14.4}$$

Let $m \geq 1$ and let $N$ be an integer smaller than or equal to $N^{\text{pack}}(E,3m)$. Let $F$ be a $3m$-packing of $E$ of cardinality $N$. We write

$$\begin{aligned} P_p\big(E \not\leftrightarrow \infty\big) \;&\leq\; P_p\big(\forall x \in F \quad x \not\leftrightarrow \infty\big) \\ &\leq\; P_p\big(\forall x \in F \quad |C(x)| \leq m\big) + P_p\big(\exists\, x \in F \quad |C(x)| > m, x \not\leftrightarrow \infty\big)\,. \end{aligned} \tag{14.5}$$

Since $F$ is a $3m$-packing, the events $\{\,|C(x)| \leq m\,\}$, $x \in F$, are independent. Indeed, the event $\{\,|C(x)| \leq m\,\}$ belongs to the $\sigma$-field generated by the sites of the box $x + \Lambda(2m)$, and the boxes $x + \Lambda(2m)$, $x \in F$, are pairwise disjoint. We have thus

$$P_p\big(\forall x \in F \quad |C(x)| \leq m\big) \;=\; \prod_{x \in F} P_p\big(|C(x)| \leq m\big)$$
$$= P_p\big(|C(0)| \leq m\big)^N \;\leq\; P_p\big(0 \not\longleftrightarrow \infty\big)^N \;\leq\; \big(1 - \theta(p)\big)^N . \quad (14.6)$$

Moreover, we have

$$P_p\big(\exists\, x \in F \quad |C(x)| > m, x \not\longleftrightarrow \infty\big) \;\leq\; N P_p\big(m < |C(0)| < +\infty\big) . \qquad (14.7)$$

Reporting the inequalities (14.6), (14.7) into (14.5), we get

$$P_p\big(E \not\longleftrightarrow \infty\big) \;\leq\; \big(1 - \theta(p)\big)^N + N P_p\big(m < |C(x)| < +\infty\big) . \qquad (14.8)$$

From the inequality (14.4), we see that the above inequality (14.8) holds for any integers $m, N$ such that

$$(12m + 1)^d N \;\leq\; n .$$

We take the infimum of the right-hand side of (14.8) and we obtain

$$1 - \theta(n, p) \;\leq\; \inf\,\Big\{\, \big(1 - \theta(p)\big)^N + N P_p\big(m < |C(x)| < +\infty\big) : (12m+1)^d N \;\leq\; n \,\Big\} . \qquad (14.9)$$

To warm up, let us prove the statement (14.1). Let $\varepsilon > 0$. Let $N$ be an integer such that $(1 - \theta(p))^N < \varepsilon$. Let $m$ be an integer such that

$$P_p\big(m < |C(x)| < +\infty\big) \;\leq\; \frac{\varepsilon}{N} .$$

For any $n > (12m + 1)^d N$, we have $1 - \theta(n, p) < 2\varepsilon$, and this concludes the proof of (14.1). We try next to estimate the infimum in (14.9). To this end, we will need the following lemma.

**Lemma 14.3.** *Suppose that $p$ is such that $\theta(p) > 0$. We have the lower bound*

$$\forall m \geq \frac{3}{\theta(p)} \qquad P_p\Big(m \leq \big|C(0)\big| < +\infty\Big) \;\geq\; \frac{\theta(p)}{2}\,\big(p(1-p)\big)^{\frac{4d2^d}{\theta(p)}(m)^{1-1/d}} . \qquad (14.10)$$

*Proof.* We follow the argument of Aizenman, Delyon, and Souillard [2] (which can also be found in Theorem 8.61 in [63]). In fact, although their result is presented for a parameter $p > p_c$, it holds as well with the sole hypothesis that $\theta(p) > 0$. We reproduce the argument here, because we need to keep track of the constants in the successive inequalities in order to relate them with the percolation probability $\theta(p)$. So, let $p$ be a parameter such that $\theta(p) > 0$. To alleviate the notation, we write simply $\theta$ instead of $\theta(p)$. For any box $\Lambda$, we have

$$E\big(|\Lambda \setminus C_\infty|\big) \;=\; E\Big(\sum_{x \in \Lambda} 1_{\{\,x \not\leftrightarrow \infty\,\}}\Big) \;=\; |\Lambda|(1 - \theta) .$$

By the Markov inequality, we have then, for any $\alpha \geq 1$,

$$P_p\big(\big|\Lambda \setminus C_\infty\big| \geq \alpha|\Lambda|(1-\theta)\big) \;\leq\; \frac{1}{\alpha}\,,$$

which we rewrite as

$$P_p\big(\big|\Lambda \setminus C_\infty\big| < \alpha|\Lambda|(1-\theta)\big) \;\geq\; 1-\frac{1}{\alpha}\,,$$

or equivalently

$$P_p\Big(\big|\Lambda \cap C_\infty\big| \geq |\Lambda|\big(1-\alpha(1-\theta)\big)\Big) \;\geq\; 1-\frac{1}{\alpha}\,,$$

Taking $\alpha = (1-\theta/2)/(1-\theta)$, we get

$$P_p\Big(\big|\Lambda \cap C_\infty\big| \geq \frac{\theta}{2}|\Lambda|\Big) \;\geq\; 1-\frac{2(1-\theta)}{2-\theta} \;=\; \frac{\theta}{2-\theta}\,. \tag{14.11}$$

Let $\mathcal{E}$ be the event

$$\mathcal{E} \;=\; \Big\{\,\big|\{\,x \in \Lambda : x \longleftrightarrow \partial^{\,in}\Lambda\,\}\big| \geq \frac{\theta}{2}|\Lambda|\,\Big\}\,.$$

We have the inclusion

$$\Big\{\,\big|\Lambda \cap C_\infty\big| \geq \frac{\theta}{2}|\Lambda|\,\Big\} \;\subset\; \mathcal{E}\,. \tag{14.12}$$

With the notation $e_1 = (1,0,\dots,0)$, let $\mathcal{F}$ be the event

$$\mathcal{F} \;=\; \left\{ \begin{array}{l} \text{the sites of the inner boundary } \partial^{\,in}\Lambda \text{ are open} \\ \text{and the sites of the segment } \mathbb{R}e_1 \cap \Lambda \text{ are open} \end{array} \right\},$$

and let $\mathcal{G}$ be the event

$$\mathcal{G} \;=\; \big\{\text{the sites of the outer boundary } \partial^{\,out}\Lambda \text{ of } \Lambda \text{ are closed}\big\}\,.$$

If both events $\mathcal{F}$ and $\mathcal{G}$ occur, then the cluster of the origin contains all the sites of $\Lambda$ connected to the boundary $\partial^{\,in}\Lambda$, and it is finite. We have thus

$$P_p\Big(\frac{\theta}{2}|\Lambda| \leq \big|C(0)\big| < +\infty\Big) \;\geq\; P_p\big(\mathcal{E} \cap \mathcal{F} \cap \mathcal{G}\big)\,. \tag{14.13}$$

The events $\mathcal{E}$, $\mathcal{F}$ belong to the $\sigma$-field generated by the sites of $\Lambda$, while the event $\mathcal{G}$ depends on the sites of $\partial^{\,out}\Lambda$. Thus $\mathcal{G}$ is independent of $\mathcal{E}$ and $\mathcal{F}$. Moreover the events $\mathcal{E}$ and $\mathcal{F}$ are increasing. By the FKG inequality, we have

$$P_p\big(\mathcal{E} \cap \mathcal{F} \cap \mathcal{G}\big) \;=\; P_p(\mathcal{E} \cap \mathcal{F})P_p(\mathcal{G}) \;\geq\; P_p(\mathcal{E})P_p(\mathcal{F})P_p(\mathcal{G})\,. \tag{14.14}$$

Thus, using the lower bound (14.11), the inclusion (14.12), and the inequalities (14.13), (14.14), we have

$$P_p\Big(\frac{\theta}{2}|\Lambda| \leq \big|C(0)\big| < +\infty\Big) \;\geq\; \frac{\theta}{2-\theta}\, p^{|\partial^{\,in}\Lambda| + \text{diameter}\,|\Lambda|}(1-p)^{|\partial^{\,out}\Lambda|}\,. \tag{14.15}$$

Let us apply the inequality (14.15) to the box $\Lambda=\Lambda(m)=[-m/2,m/2]^d$. For this box, we have

$$|\Lambda(m)|\,\geq\,m^d\,,\quad\big|\mathbb{Z}e_1\cap\Lambda(m)\big|\,\leq\,m+1\,,$$
$$|\partial^{\,in}\Lambda(m)|\,\leq\,2d(m+1)^{d-1}\,,\quad|\partial^{\,out}\Lambda(m)|\,\leq\,2d(m+3)^{d-1}\,,$$

thus we get from inequality (14.15) that

$$P_p\Big(\frac{\theta}{2}m^d\leq\big|C(0)\big|<+\infty\Big)\,\geq\,\frac{\theta}{2}\,p^{2d(m+1)^{d-1}+m+1}(1-p)^{2d(m+3)^{d-1}}\,.\qquad(14.16)$$

Let us set $m'\,=\,\big\lfloor\theta m^d/2\big\rfloor$, so that

$$\frac{\theta}{2}m^d-1\,\leq\,m'\,\leq\,\frac{\theta}{2}m^d\,.$$

For $m\geq 3/\theta$, we have $m\geq 3$, $m'\geq 1$ and

$$(m+3)^{d-1}\,\leq\,2^{d-1}\Big(\frac{2}{\theta}(m'+1)\Big)^{1-1/d}\,\leq\,\frac{2^{d+1}}{\theta}(m')^{1-1/d}\,.\qquad(14.17)$$

Inequalities (14.16) and (14.17) yield

$$P_p\Big(m'\leq\big|C(0)\big|<+\infty\Big)\,\geq\,\frac{\theta}{2}\,\big(p(1-p)\big)^{\frac{4d2^d}{\theta}(m')^{1-1/d}}\,,$$

which is precisely the inequality (14.10). □

Let us come back to the infimum in (14.9). We take, for $n\geq 1$ and $\alpha>0$,

$$m=n^\alpha\,,\qquad N\,=\,\frac{1}{\ln(1-\theta)}\ln P_p\Big(m\leq\big|C(0)\big|<+\infty\Big)\,.\qquad(14.18)$$

We first check the constraint $(12m+1)^dN\,\leq\,n$. By inequality (14.10), we have

$$\forall n\geq\Big(\frac{3}{\theta(p)}\Big)^{1/\alpha}\qquad N\,\leq\,\frac{1}{\ln(1-\theta)}\Big(\ln\frac{\theta}{2}+\frac{4d2^d}{\theta}\big(n^\alpha\big)^{1-1/d}\ln\big(p(1-p)\big)\Big)\,.$$

Using the bound $\ln(1-\theta)\leq-\theta$, we obtain, for $n\geq(3/\theta(p))^{1/\alpha}$,

$$(12m+1)^dN\,\leq\,-\frac{13^d}{\theta}n^{\alpha d}\ln\frac{\theta}{2}-\frac{98^d}{\theta^2}n^{\alpha(d+1-1/d)}\ln\big(p(1-p)\big)\,.$$

We choose $\alpha=1/(d+1)$. The constraint $(12m+1)^dN\,\leq\,n$ is satisfied as soon as

$$-2\frac{13^d}{\theta}\ln\frac{\theta}{2}\leq n^{1/d+1}\,,\qquad-\frac{98^d}{\theta^2}\ln\big(p(1-p)\big)\leq\frac{1}{2}n^{1/(d(d+1))}\,.$$

For $n$ satisfying the constraint (14.2), the inequality (14.9) and the equality (14.18) together yield that

$$\begin{aligned}1-\theta(n,p)\,&\leq\,\big(1-\theta(p)\big)^N+NP_p\big(n^{1/(d+1)}<|C(x)|<+\infty\big)\\&\leq\,(1+N)P_p\big(n^{1/(d+1)}\leq|C(x)|<+\infty\big)\\&\leq\,\frac{2}{\ln(1-\theta)}\phi\Big(P_p\Big(n^{1/(d+1)}\leq\big|C(0)\big|<+\infty\Big)\Big)\,,\end{aligned}$$

where $\phi(t)=t\ln t$. This is the desired inequality (14.3). □

## 14.3 The travel time to exit a finite domain

We improve here the estimate of lemma 14.1 with the help of the function $\theta(n,p)$. Indeed, we shall show that the travel time to exit a finite domain is in fact dominated by a super-geometric distribution. For $x=(x_1,\dots,x_d)\in\mathbb{Z}^d$, we define the norm $|x|_1$ of $x$ by setting

$$|x|_1 \;=\; |x_1|+\cdots+|x_d|\,,$$

and we denote by $B_1(x,t)$ the closed ball centered at $x$ of radius $t$ for the norm $|\cdot|_1$, i.e.,

$$B_1(x,t) \;=\; \big\{\, y\in\mathbb{Z}^d : |x-y|_1\leq t \,\big\}\,. \tag{14.19}$$

**Proposition 14.4.** *Let $D$ be a finite subset of $\mathbb{Z}^d$ and let $x\in D$. Let $t\geq 1$ be such that $B_1(x,t)\subset D$. We have*

$$\forall t\geq 1 \qquad P\big(\,T_D(x,\partial^{\,in}D)\,>\,t\,\big)\;\leq\;\prod_{1\leq s\leq t}\Big(1-\theta\big(\big|B_1(x,s-1)\big|,p\big)\Big)\,. \tag{14.20}$$

*Proof.* Let $D$ be a finite subset of $\mathbb{Z}^d$, let $x$ belong to $D$ and let $t\geq 1$. Using (12.11), we have (we write simply $\mathcal{B}(x,D,t-1)$ instead of $\mathcal{B}(\{\,x\,\},D,t-1)$)

$$\forall t<T_D(x,\partial^{\,in}D)\qquad T_D(x,\partial^{\,in}D)\;=\;t-1+T_D\big(\mathcal{B}(x,D,t-1),\partial^{\,in}D\big)\,.$$

Therefore

$$\begin{aligned}P\big(T_D(x,\partial^{\,in}D)>t\big)\;&=\;P\Big(T_D\big(\mathcal{B}(x,D,t-1),\partial^{\,in}D\big)>1\Big)\\ &=\sum_{B\subset D,B\cap\partial^{\,in}D=\varnothing}P\Big(\mathcal{B}(x,D,t-1)=B,\;T_D(B,\partial^{\,in}D)>1\Big)\,.\end{aligned} \tag{14.21}$$

Let $B\subset D$ be such that $B\cap\partial^{\,in}D=\varnothing$ and suppose that $\mathcal{B}(x,D,t-1)=B$. Necessarily, all the sites in $\partial^{\,in}B$ must be closed and therefore we have

$$T_D(B,\partial^{\,in}D)\;=\;1+T_D(\partial^{\,out}B,\partial^{\,in}D)\,.$$

Furthermore

$$T_D(\partial^{\,out}B,\partial^{\,in}D)>0\quad\Longrightarrow\quad B\cup\partial^{\,out}B\not\longleftrightarrow\partial^{\,in}D\,,$$

whence, for any $t\geq 1$,

$$\begin{aligned}P\big(\mathcal{B}(x,D,t-1)=B,T_D(B,\partial^{\,in}D)>1\big)&\\ \leq\;P\big(\mathcal{B}(x,D,t-1)=B,T_D(\partial^{\,out}B,\partial^{\,in}D)>0\big)&\\ \leq\;P\big(\mathcal{B}(x,D,t-1)=B,B\cup\partial^{\,out}B\not\longleftrightarrow\partial^{\,in}D\big)&\,.\end{aligned} \tag{14.22}$$

The event $\{\,\mathcal{B}(x,D,t-1)=B\,\}$ is measurable with respect to the sites of $B$, while the event $\{\,B\cup\partial^{\,out}B\not\longleftrightarrow\partial^{\,in}D\,\}$ is measurable with respect to the sites of $D\setminus B$, thus these two events are independent and

$$\begin{aligned}P\big(\mathcal{B}(x,D,t-1)=B,T_D(B,\partial^{\,in}D)>1\big)&\\ \leq\;P\big(\mathcal{B}(x,D,t-1)=B\big)P\big(B\cup\partial^{\,out}B\not\longleftrightarrow\partial^{\,in}D\big)&\,.\end{aligned} \tag{14.23}$$

The event $\{ B \cup \partial^{out} B \not\longleftrightarrow \partial^{in} D \}$ is included in $\{ B \cup \partial^{out} B \not\longleftrightarrow \infty \}$, thus

$$P\big(B \cup \partial^{out} B \not\longleftrightarrow \partial^{in} D\ \big) \ \leq\ 1 - \theta\big(\big|B \cup \partial^{out} B\big|, p\big)\,. \tag{14.24}$$

From (12.8), the hypothesis $B_1(x,t) \subset D$, and the fact that $\mathcal{B}(x,D,t-1) = B$, it follows that

$$B_1(x,t-1) \ =\ \Big\{\, y \in D : d_D(x,y) \leq t-1 \,\Big\} \ \subset\ B\,. \tag{14.25}$$

Reporting (14.25) into (14.24), (14.23) and (14.22), we obtain

$$\begin{aligned} P\big(\mathcal{B}(x,D,t-1) = B, T_D(B,\partial^{in} D) > 1\big) \ \leq& \\ & P\big(\mathcal{B}(x,D,t-1) = B\big)\Big(1 - \theta\big(\big|B_1(x,t-1)\big|, p\big)\Big)\,. \end{aligned}$$

We plug this into (14.21) and we sum over $B$ to get

$$P\big(T_D(x,\partial^{in} D) > t\big) \ \leq\ P\big(T_D(x,\partial^{in} D) > t-1\big)\Big(1 - \theta\big(\big|B_1(x,t-1)\big|, p\big)\Big)\,.$$

Iterating this inequality, we obtain the inequality (14.20). □

**Corollary 14.5.** *Suppose that $p$ is such that $\theta(p) > 0$. Let $D$ be a finite subset of $\mathbb{Z}^d$ and let $x \in D$. There exists a non-decreasing function $\alpha$ defined on $\mathbb{N}$ such that $\alpha(1) \geq \theta(p)$, $\lim_{t\to+\infty} \alpha(t) = +\infty$ and*

$$\forall t \geq 0 \qquad B_1(x,t) \subset D \quad \Longrightarrow \quad P\big(\,T_D(x,\partial^{in} D)\ >\ t\,\big) \ \leq\ \exp\big(-t\alpha(t)\big)\,. \tag{14.26}$$

*Recalling that $\phi(t) = t \ln t$, we have in addition*

$$\forall t \ \geq\ 4d \max\Big(1, \Big(-2\frac{13^d}{\theta(p)} \ln \frac{\theta(p)}{2}\Big)^{1+1/d}, \Big(-2\frac{98^d}{\theta(p)^2} \ln\big(p(1-p)\big)\Big)^{d+1}\Big) \tag{14.27}$$

$$\alpha(t) \ \geq\ \ln\Big(\frac{\theta(p)}{2}\Big) - \frac{1}{2}\ln\bigg(-\phi\Big(P_p\Big(\Big(\frac{t}{4d}\Big)^{d/(d+1)} < |C(0)| < +\infty\Big)\Big)\bigg)\,. \tag{14.28}$$

The condition $B_1(x,t) \subset D$ in (14.26) could be removed. In fact, without this condition, the probability $P\big(T_D(x,\partial^{in} D) > t\big)$ is equal to 0. Why do we keep the ugly expression (14.27)? We could write that the inequality (14.28) holds as soon as $t$ is larger than some value depending on $\theta(p), p, d$. However, the inequality will play a key role, so we feel it is important to present an explicit expression that clearly shows how this value depends on the other parameters.

The common sense dictates that the probability $P_p\big(n < |C(0)| < +\infty\big)$ decays at least like an inverse power of $n$. If that is the case, then it would follow from the inequality (14.28) that

$$\alpha(t) \ \geq\ \frac{1}{5} \ln t$$

for $t$ large enough. Unfortunately, we did not manage to obtain a bound of this kind (see the disappointing section 15).

*Proof.* We define the function $\alpha(t)$ by

$$\alpha(t) \,=\, -\frac{1}{t}\sum_{1\leq s\leq t}\ln\Big(1-\theta\big(\big|B_1(x,s-1)\big|,p\big)\Big)\,.$$

Since the map $s\mapsto\theta\big(\big|B_1(x,s-1)\big|,p\big)$ is non-decreasing, then the map $\alpha(t)$ is also non-decreasing. Moreover

$$\alpha(1) \,=\, -\ln(1-\theta(p)) \,\geq\, \theta(p)\,.$$

With this function $\alpha$, the inequality (14.20) can be rewritten as

$$P\big(\,T_D(x,\partial^{\,in}D)\,>\,t\,\big)\,\leq\,\exp\big(-t\alpha(t)\big)\,. \tag{14.29}$$

By (14.1) or proposition 14.2, we have that $\lim_{n\to\infty}\theta(n,p)\,=\,1$. It follows from Cesaro's lemma that

$$\lim_{t\to+\infty}\alpha(t)=+\infty\,.$$

Here is a more direct argument. Using the monotonicity in $n$ of $\theta(n,p)$, we have

$$\alpha(t)\,\geq\,-\frac{1}{t}\Big(t-\Big\lfloor\frac{t}{2}\Big\rfloor+1\Big)\ln\Big(1-\theta\big(\big|B_1(x,\lfloor t/2\rfloor-1)\big|,p\big)\Big)\,, \tag{14.30}$$

and the right-hand quantity goes to $+\infty$ as $t$ goes to $+\infty$. Let us exploit further the inequality (14.29). From now onwards, we suppose that $t\geq 8$. We have

$$\Big(t-\Big\lfloor\frac{t}{2}\Big\rfloor+1\Big)\,\geq\,\frac{t}{2}\,,\quad\Big|B_1\Big(x,\Big\lfloor\frac{t}{2}\Big\rfloor-1\Big)\Big|\,\geq\,\Big|B_1\Big(x,\frac{t}{4}\Big)\Big|\,\geq\,\Big(\frac{t}{4d}\Big)^d\,. \tag{14.31}$$

Since $\theta(n,p)$ is non-decreasing in $n$, the previous inequalities (14.31) together with (14.30) yield

$$\alpha(t)\,\geq\,-\frac{1}{2}\ln\Big(1-\theta\Big(\Big(\frac{t}{4d}\Big)^d,p\Big)\Big)\,.$$

It follows from proposition 14.2 that

$$\Big(\frac{t}{4d}\Big)^d\,\geq\,\max\bigg(\Big(-2\frac{13^d}{\theta}\ln\frac{\theta}{2}\Big)^{d+1},\Big(-2\frac{98^d}{\theta^2}\ln\big(p(1-p)\big)\Big)^{d(d+1)}\bigg) \tag{14.32}$$

$$\implies\quad\alpha(t)\,\geq\,-\frac{1}{2}\ln\bigg(\frac{2}{\ln(1-\theta)}\phi\Big(P_p\Big(\Big(\frac{t}{4d}\Big)^{d/(d+1)}<|C(0)|<+\infty\Big)\Big)\bigg)$$

We use the inequality

$$-\frac{1}{2}\ln\Big(-\frac{2}{\ln(1-\theta)}\Big)\,\geq\,\ln\Big(\frac{\theta}{2}\Big)$$

and we rewrite the implication (14.32) as

$$t\,\geq\,4d\max\bigg(\Big(-2\frac{13^d}{\theta}\ln\frac{\theta}{2}\Big)^{1+1/d},\Big(-2\frac{98^d}{\theta^2}\ln\big(p(1-p)\big)\Big)^{d+1}\bigg)$$

$$\implies\quad\alpha(t)\,\geq\,\ln\Big(\frac{\theta}{2}\Big)-\frac{1}{2}\ln\bigg(-\phi\Big(P_p\Big(\Big(\frac{t}{4d}\Big)^{d/(d+1)}<|C(0)|<+\infty\Big)\Big)\bigg)\,.$$

This is the inequality (14.28) we were seeking. □

## 14.4 The classical square root trick

We focus here on the specific situation where $D = \Lambda$ is a symmetric cubic box and $x = 0$ is the center of the box. Let us denote by $F_i$, $1 \leq i \leq 2d$, the $2d$ faces of $\Lambda$. Each face $F_i$ is a $(d-1)$-dimensional square, which is itself the union of $2^{d-1}$ $(d-1)$-dimensional symmetric squares $F_i^j$, $1 \leq j \leq 2^{d-1}$. Each of these squares shares a vertex with a vertex of the boundary of $\Lambda$ and admits the center of $F_i$ as another vertex. By the definition of the travel time, we have

$$T_\Lambda(0, \partial^{in}\Lambda) = \min_{1 \leq i \leq 2d} \min_{1 \leq j \leq 2^{d-1}} T_\Lambda(0, F_i^j),$$

whence, for any $t \geq 1$,

$$\{ T_\Lambda(0, \partial^{in}\Lambda) \geq t \} = \bigcap_{1 \leq i \leq 2d} \bigcap_{1 \leq j \leq 2^{d-1}} \{ T_\Lambda(0, F_i^j) \geq t \}.$$

All these events are decreasing. Applying the FKG inequality, we get

$$P\big( T_\Lambda(0, \partial^{in}\Lambda) \geq t \big) \geq \prod_{1 \leq i \leq 2d} \prod_{1 \leq j \leq 2^{d-1}} P\big( T_\Lambda(0, F_i^j) \geq t \big). \tag{14.33}$$

From the symmetry of the model, we see that the random variables $T_\Lambda(0, F_i^j)$, $1 \leq i \leq 2d$, $1 \leq j \leq 2^{d-1}$, are identically distributed. It follows from the inequality (14.33) that

$$\forall i \in \{ 1, \dots, 2d \} \quad \forall j \in \{ 1, \dots, 2^{d-1} \}$$
$$P\big( T_\Lambda(0, F_i^j) \geq t \big) \leq \Big( P\big( T_\Lambda(0, \partial^{in}\Lambda) \geq t \big) \Big)^{1/(d2^d)}. \tag{14.34}$$

The inequality (14.34) and corollary 14.5 imply the following estimate.

**Corollary 14.6.** *Suppose that $p$ is such that $\theta(p) > 0$. Let $\alpha(t)$ be the same function as in corollary 14.5. We have*

$$\forall i \in \{ 1, \dots, 2d \} \quad \forall j \in \{ 1, \dots, 2^{d-1} \}$$
$$\forall t \geq 0 \qquad P\big( T_\Lambda(0, F_i^j) \geq t \big) \leq \exp\big( -t\alpha(t)/(d2^d) \big).$$

Let $Q$ be a rectangular box in $\mathbb{Z}^d$ with sides parallel to the axes of $\mathbb{Z}^d$. Let $x$ be any vertex belonging to $Q$. Let $n(x) = d(x, \partial^{in}Q)$ be the graph distance between $x$ and $\partial^{in}Q$. The box $\Lambda(x, 2n(x))$ is the largest cubic box centered at $x$ contained in $Q$. Since $Q$ is a rectangular box, then at least one face of $\Lambda(x, 2n(x))$ is included in $\partial^{in}Q$, and therefore at least one of the squares $F_i^j$ which tile the faces of $\Lambda(x, 2n(x))$ is included in $\partial^{in}Q$. Let $F_{i*}^{j^*}$ be one such square. Using the translation invariance of the model, we can apply the estimate of corollary 14.6 to control the travel time between $x$ and $F_{i*}^{j^*}$ in $Q$.

## 14.5 The travel time to the infinite cluster

A noteworthy point is that proposition 14.4 does not rely explicitly on the BK inequality, but rather on the exploration algorithm which discovers the successive shells around $x$. As a consequence, we do not need to restrict ourselves to a finite domain, and the same argument yields an estimate on the travel time $T_\infty$ to $+\infty$, defined as follows: for $A$ a finite subset of $\mathbb{Z}^d$, we set

$$T_\infty(A) \,=\, \inf\Big\{\,\sum_{i=0}^{+\infty} 1_{z_i \text{ closed}} : z_0, z_1, \dots \text{ infinite path starting from } z_0 \in A\Big\}\,.$$

In case $A = \{\,x\,\}$ is a singleton, we write simply $T_\infty(x)$ instead of $T_\infty(\{\,x\,\})$. By the way, the definitions involving a domain $D$ are written for $D$ being a finite subset of $\mathbb{Z}^d$. However, some of the definitions are perfectly valid for infinite subsets of $\mathbb{Z}^d$. When this is the case, we can even take $D = \mathbb{Z}^d$, and we recover the usual definitions for configurations on the infinite lattice. When we wish to take $D = \mathbb{Z}^d$ and there is no risk of confusion, we shall simply remove $D$ from the notation. For example, the open cluster of a site $x$ is $C(x) = C(x, \mathbb{Z}^d)$, and we define $\overline{C}(x)$ as

$$\overline{C}(x) \,=\, \overline{C}(x, \mathbb{Z}^d) \,=\, C(x) \cup \partial^{\,out} C(x)\,.$$

Similarly, we denote by $T$ the travel time on the whole lattice, i.e., $T = T_{\mathbb{Z}^d}$. For $x \in \mathbb{Z}^d$ and $t \geq 0$ an integer, we define the shell number $t$ around $x$ by

$$\text{Shell}\,(x,t) \,=\, \big\{\, y \in \mathbb{Z}^d : T(x,y) = t\,\big\}\,,$$

and the closed ball of radius $t$ for the travel time $T_{\mathbb{Z}^d}$ by

$$\mathcal{B}(x,t) \,=\, \big\{\, y \in \mathbb{Z}^d : T(x,y) \leq t\,\big\}\,.$$

Obviously, we have $\text{Shell}\,(x) = \text{Shell}\,(x,0) = \overline{C}(x)$ and the successive shells around $x$ are pairwise disjoint.

**Proposition 14.7.** *Let $x$ be a point of $\mathbb{Z}^d$. We have*

$$\forall t \geq 1 \qquad P\big(\,T_\infty(x) \,\geq\, t\,\big) \,\leq\, \prod_{1 \leq s \leq t} \Big(1 - \theta\big(\big|B_1(x, s-1)\big|, p\big)\Big)\,. \qquad (14.35)$$

*Proof.* We adapt the proof of proposition 14.4 to the case $D = \mathbb{Z}^d$. Let $x \in \mathbb{Z}^d$ and let $t \geq 1$. From the identity (12.11), we have

$$\forall t \leq T_\infty(x) \qquad T_\infty(x) \,=\, t - 1 + T_\infty(\mathcal{B}(x, t-1))\,.$$

If the travel time between $\text{Shell}\,(x, t-1)$ and $+\infty$ is positive, then it must be the case that $\text{Shell}\,(x, t-1)$ is a finite set, otherwise it would already contain an infinite open path. We have therefore

$$\begin{aligned} P\big(T_\infty(x) \geq t\big) \,&=\, P\big(T_\infty(\mathcal{B}(x,t-1)) \geq 1\big) \\ &=\, \sum_{B \subset \mathbb{Z}^d, |B| < \infty} P\big(\mathcal{B}(x, t-1) = B, T_\infty(B) \geq 1\big)\,. \qquad (14.36) \end{aligned}$$

On the event $\{\,\mathcal{B}(x,t-1)=B\,\}$, all the sites in $\partial^{\,in}B$ are closed, and we have $T_\infty(B)=1+T_\infty(\partial^{\,out}B)$. Furthermore

$$T_\infty(\partial^{\,out}B)>0 \quad\Longleftrightarrow\quad B\cup\partial^{\,out}B \not\longleftrightarrow C_\infty\,,$$

whence, for any $t\geq 1$,

$$\begin{aligned}P\big(\mathcal{B}(x,t-1)=B,T_\infty(B)>1\big) \;&=\; P\big(\mathcal{B}(x,t-1)=B,T_\infty(\partial^{\,out}B)>0\big)\\ &\leq\; P\big(\mathcal{B}(x,t-1)=B,B\cup\partial^{\,out}B\not\longleftrightarrow C_\infty\big)\,. \end{aligned}\tag{14.37}$$

The event $\{\,\mathcal{B}(x,t-1)=B\,\}$ is measurable with respect to the sites of $B$, while the event $\{\,B\cup\partial^{\,out}B\not\longleftrightarrow C_\infty\,\}$ is measurable with respect to the sites of $\mathbb{Z}^d\setminus B$, thus these two events are independent and

$$\begin{aligned}&P\big(\mathcal{B}(x,t-1)=B,T_\infty(B)>1\big)\\ &\qquad\leq\; P\big(\mathcal{B}(x,t-1)=B\big)P\big(B\cup\partial^{\,out}B\not\longleftrightarrow C_\infty\big)\,.\end{aligned}\tag{14.38}$$

The event $\{\,B\cup\partial^{\,out}B\not\longleftrightarrow C_\infty\,\}$ is included in the event $\{\,B\cup\partial^{\,out}B\not\longleftrightarrow \infty\,\}$, thus

$$P\big(B\cup\partial^{\,out}B\not\longleftrightarrow C_\infty\,\big)\;\leq\;1-\theta\big(\big|B\cup\partial^{\,out}B\big|,p\big)\,.\tag{14.39}$$

From the inclusion (12.8) and the fact that $\mathcal{B}(x,t-1)=B$, we have

$$B_1(x,t-1)\;\subset\;\Big\{\,y\in\mathbb{Z}^d:d(x,y)\leq t-1\,\Big\}\;\subset\;B\,.\tag{14.40}$$

Reporting (14.40) into (14.39), (14.38) and (14.37), we obtain

$$\begin{aligned}&P\big(\mathcal{B}(x,t-1)=B,T_\infty(B)>1\big)\;\leq\\ &\qquad\qquad P\big(\mathcal{B}(x,t-1)=B\big)\Big(1-\theta\big(\big|B_1(x,t-1)\big|,p\big)\Big)\,.\end{aligned}$$

We plug this into (14.36), and we sum over $B$ to get

$$P\big(T_\infty(x)\geq t\big)\;\leq\;P\big(T_\infty(x)\geq t-1\big)\Big(1-\theta\big(\big|B_1(x,t-1)\big|,p\big)\Big)\,.$$

Iterating this inequality, we obtain the inequality (14.35). □

**Corollary 14.8.** *Suppose that $p$ is such that $\theta(p)>0$. The travel time from 0 to the infinite cluster $C_\infty$ admits exponential moments of any order, i.e.,*

$$\forall\lambda>0\qquad E\big(\exp(\lambda T_\infty)\big)\;<\;+\infty\,.$$

*Proof.* We simply apply proposition 14.7, and we use the upper bound (14.29) on the product appearing in (14.35) that we developed in the course of the proof of corollary 14.5. This yields that

$$\forall t\geq 0\qquad P\big(\,T_\infty(x)\;>\;t\,\big)\;\leq\;\exp\big(-t\alpha(t)\big)\,,$$

where $\alpha(t)$ goes to $\infty$ as $t$ goes to $\infty$. This readily implies that $T_\infty$ admits exponential moments of any order. □

## 15 Two disappointing results on clusters at $p_c$ ⊛

The probability $P_p\big(n < |C(0)| < +\infty\big)$ plays a key role in the estimates on the travel time developed throughout section 12. We would very much like to have an explicit quantitative upper bound on the asymptotic decay of this probability as $n$ goes to $\infty$, however weak it may be. Unfortunately, we did not manage to obtain such a result. We prove here the following two non-obvious (but disappointing) facts on the distribution of the clusters at the critical point $p_c$:

• For any $\alpha > 1 - 1/d$, we have

$$E_{p_c}(|C(0)|^\alpha) \,=\, +\infty\,.$$

• There exists a positive constant $\kappa$ such that

$$\forall n \geq 1 \qquad P_{p_c}(|C(0)| = n) \,\leq\, \kappa\sqrt{\frac{\ln n}{n}}\,.$$

### 15.1 On the critical exponent $\delta$

In this subsection, we prove the first of the two above facts. A classical theorem of Hammersley [65] states that, if for some value $p$ we have

$$E_p(|C(0)|) \,<\, +\infty\,, \tag{15.1}$$

then the radius of $C(0)$ decays exponentially fast, so that $p < p_c$. We shall improve the condition (15.1) as follows.

**Theorem 15.1.** *Let $p \in [0,1]$. Suppose that, for some $\alpha > 1 - 1/d$, we have*

$$E_p(|C(0)|^\alpha) \,<\, +\infty\,. \tag{15.2}$$

*Then there exists $\sigma(p) > 0$ such that*

$$\forall n \geq 1 \qquad P_p\big(0 \longleftrightarrow \partial^{\,in}\Lambda(n)\big) \,\leq\, \exp(-\sigma(p)n)\,.$$

Any value of the percolation parameter $p$ for which there is exponential decay of the radius of the cluster at the origin is strictly smaller than $p_c$. This yields immediately the following corollary.

**Corollary 15.2.** *For any $\alpha > 1 - 1/d$, we have*

$$E_{p_c}(|C(0)|^\alpha) \,=\, +\infty\,.$$

There are numerous conjectures for the percolation model concerning the existence of the critical exponents. For instance, it is believed that there exists an exponent $\delta$ such that

$$P_{p_c}(|C(0)| = n) \,\approx\, n^{-1-\frac{1}{\delta}} \quad \text{as} \quad n \to \infty\,, \tag{15.3}$$

where $\approx$ denotes asymptotic equivalence of the logarithms, i.e., formula (15.3) means

$$\lim_{n\to\infty} \frac{1}{\ln n} \ln P_{p_c}(|C(0)| = n) \,=\, -1 - \frac{1}{\delta}\,.$$

**Corollary 15.3.** *Suppose that the critical exponent $\delta$ exists, i.e., that* (15.3) *holds. Then this exponent satisfies $\delta \geq d/(d-1)$.*

In fact, there exist already two rigorous proofs that $\delta \geq 2$, provided that $\delta$ exists: they are due to Aizenman and Barsky [1] and to Newman [108] (these proofs are also presented in Proposition 10.29 in [63]). Their arguments are more sophisticated than the one presented here, and their results are better in dimensions $d \geq 3$. In contrast, the result of corollary 15.2 is completely proved, meaning that it is not a conditional statement relying on the existence of an exponent, like corollary 15.3.

*Proof of theorem 15.1.* The end of the argument is the same as in the proof of Hammersley's theorem. The main difficulty here consists in exploiting the weaker condition (15.2) in order to control the number $N_n$ of the vertices in $\partial\Lambda(n)$ which are connected to the origin, i.e.,

$$N_n \,=\, \big|\{\, x \in \partial\Lambda(n) : 0 \longleftrightarrow x \,\}\big|$$

(throughout this section, the inner boundary $\partial^{\,in}\Lambda$ of the box $\Lambda$ is denoted simply by $\partial\Lambda$). This shall be done with the help of a classical convexity inequality. Let $0 < \alpha < 1$ and $n \geq 1$. We have

$$|C(0)|^\alpha \,\geq\, \Big(\sum_{k=0}^{+\infty} N_{2k}\Big)^\alpha \,\geq\, \Big(\sum_{k=0}^{n-1} N_{2k}\Big)^\alpha .$$

By convexity, or rather concavity of the map $t \in \mathbb{R}^+ \mapsto t^\alpha$, we have

$$\Big(\frac{1}{n}\sum_{k=0}^{n-1} N_{2k}\Big)^\alpha \,\geq\, \frac{1}{n}\sum_{k=0}^{n-1} \big(N_{2k}\big)^\alpha ,$$

whence

$$\frac{1}{n^\alpha}|C(0)|^\alpha \,\geq\, \frac{1}{n}\sum_{k=0}^{n-1} \big(N_{2k}\big)^\alpha . \tag{15.4}$$

Taking the expectation in inequality (15.4), we get

$$\frac{1}{n^\alpha}E_p\big(|C(0)|^\alpha\big) \,\geq\, \frac{1}{n}\sum_{k=0}^{n-1} E_p\big((N_{2k})^\alpha\big) ,$$

Now, there exists an index $\ell \in \{\, 0, \dots, n-1 \,\}$ such that

$$E_p\Big(\big(N_{2\ell}\big)^\alpha\Big) \,\leq\, \frac{1}{n^\alpha}E_p\big(|C(0)|^\alpha\big) . \tag{15.5}$$

Next, since

$$|N_{2\ell}| \,\leq\, \big|\partial\Lambda(2\ell)\big| \,\leq\, 2d(2\ell+1)^{d-1} ,$$

we have, using (15.5),

$$\begin{aligned} E_p\big(N_{2\ell}\big) \;&\leq\; E_p\big((N_{2\ell})^\alpha\big)\Big(2d(2\ell+1)^{d-1}\Big)^{1-\alpha} \\ &\leq\; \frac{1}{n^\alpha}E_p\big(|C(0)|^\alpha\big)\Big(2d(2n+1)^{d-1}\Big)^{1-\alpha} \\ &\leq\; E_p\big(|C(0)|^\alpha\big)(3^d d)n^{d-1-d\alpha}\,. \end{aligned} \tag{15.6}$$

Suppose that $\alpha > 1-1/d$ and that $E_p\big(|C(0)|^\alpha\big) < +\infty$. Then the right-hand quantity in inequality (15.6) goes to 0 as $n$ goes to $+\infty$, and we conclude that there exists $\ell \geq 1$ such that $E_p(N_{2\ell}) < 1$. From this point onwards, we can reuse the original argument of Hammersley [65]. $\square$

## 15.2 An upper bound on $P_{p_c}(|C(0)| = n)$

At criticality, it is believed that $P_{p_c}(|C(0)| = n)$ decays like an inverse power of $n$, as expressed in the conjecture (15.3). However, it does not seem to be easy to obtain a quantitative upper bound on $P_{p_c}(|C(0)| = n)$. We provide here such an upper bound.

**Theorem 15.4.** *Let $p \in [0,1]$. There exists a positive constant $\kappa = \kappa(p,d)$, which depends on $p$ and the dimension $d$ only, such that*

$$\forall n \geq 1 \qquad P_p(|C(0)| = n) \;\leq\; \kappa\sqrt{\frac{\ln n}{n}}\,.$$

Of course this theorem is interesting only at $p = p_c$. Indeed, it is known that, in the subcritical regime, the tail of the cluster distribution decreases exponentially fast, i.e.,

$$\forall p < p_c \quad \exists \sigma > 0 \quad \exists c > 0 \quad \forall n \geq 1 \qquad P_p(|C(0)| = n) \;\leq\; c\exp(-\sigma\, n)\,,$$

while in the supercritical regime, it decays like a stretched exponential, i.e.,

$$\forall p > p_c \quad \exists \gamma > 0 \quad \exists c > 0 \quad \forall n \geq 1 \qquad P_p(|C(0)| = n) \;\leq\; c\exp\big(-\gamma\, n^{(d-1)/d}\big)\,.$$

The bound given in theorem 15.4 is a bit disappointing, since we know that the sum of the series $\sum_{n\geq 0} P_p(|C(0)| = n)$ is at most 1. Still, it does not seem obvious at all to obtain a polynomial decay of these probabilities. We shall rely here on a technique due to van der Hofstad and Kager.

### 15.2.1 The 2-pattern theorem of van der Hofstad and Kager

We import literally the notation of the paper [134], and we restate the two-pattern theorem of van der Hofstad and Kager. We introduce a minor modification with respect to the notation of [134], namely, the patterns are denoted by the calligraphic letters $\mathcal{P}, \mathcal{P}'$ in order to keep the usual letter $P$ for the probability. We denote by $Q$ the closed cube of side length 2, whose smallest corner for the lexicographic order is the origin, i.e.,

$$Q \;=\; \big\{\, x = (x_1, \dots, x_d) \in \mathbb{Z}^d : 0 \leq x_i \leq 2, \quad 1 \leq i \leq d \,\big\}\,.$$

The extended cube $\overline{Q}$ is obtained by extending $Q$ by one unit in all directions, i.e.,

$$\overline{Q} = \{ x = (x_1, \dots, x_d) \in \mathbb{Z}^d : -1 \le x_i \le 3, \quad 1 \le i \le d \}.$$

The boundary $\partial Q$ of the cube $Q$ is defined as

$$\partial Q = \partial^{in} \overline{Q} = \overline{Q} \setminus Q$$

and the boundary $\partial \overline{Q}$ of $\overline{Q}$ is

$$\partial \overline{Q} = \partial^{in} \overline{\overline{Q}} = \overline{\overline{Q}} \setminus \overline{Q}.$$

We denote by $Q_x$ the cube $Q$ translated by $x$, i.e., $Q_x = x + Q$. Naturally, we have

$$\partial Q_x = x + \partial Q, \quad \overline{Q}_x = x + \overline{Q}, \quad \partial \overline{Q}_x = x + \partial \overline{Q}.$$

A pattern $\mathcal{P}$ is a prescribed configuration of the states of the sites in the cube $Q$, i.e., $\mathcal{P} = (\mathcal{P}(x), x \in Q)$ is an element of $\{0, 1\}^Q$. We say that the pattern $\mathcal{P}$ occurs at the site $x$ in the configuration $\omega$ if

$$\forall y \in Q \qquad \omega(x + y) = \mathcal{P}(y).$$

We say that the pattern $\mathcal{P}$ occurs at $x$ on the occupied cluster of the origin $C(0)$ if the pattern $\mathcal{P}$ occurs at $x$ and all the sites of $\partial Q_x$ belong to $C(0)$. We introduce next the grid $V = 5\mathbb{Z}^d + (1, 1, \cdots, 1)$. We denote by $N_{\mathcal{P}}$ the total number of occurrences of the pattern $\mathcal{P}$ on $C(0)$ at distinct sites of $V$, i.e.,

$$N_{\mathcal{P}} = \big|\{ x \in V : \text{the pattern } \mathcal{P} \text{ occurs at } x \text{ on } C(0)\}\big|. \tag{15.7}$$

There are two important reasons that led van der Hofstad and Kager to introduce these definitions. First, different patterns $\mathcal{P}$ and $\mathcal{P}'$ occurring on $V$ are well separated and cannot overlap. Second, whenever a pattern $\mathcal{P}$ occurs on $C(0)$ at a site $x$ of $V$, it always contributes to the same number of sites to $C(0)$. The framework of the paper [134] is much more general than the Bernoulli percolation model, and the following two theorems require some additional hypothesis, namely, the finite energy property and the existence of the limit

$$\lim_{n \to \infty} \big(P(|C(0)| = n)\big)^{1/n}$$

(see the condition (1.3) of [134] on the cluster size distribution). However, these two hypotheses are automatically satisfied in our context (see for instance theorem 6.78 of [63]), so we state directly the relevant results for the Bernoulli percolation model without mentioning further the hypothesis appearing in [134]. Van der Hofstad and Kager prove first a variant of the pattern theorem for clusters previously proved in [101].

**Theorem 15.5** (Pattern theorem, particular case of theorem 1.8 of [134])**.** *Let $\mathcal{P}$ be a pattern. There exist positive constants $a, \alpha, \beta > 0$ such that*

$$\forall n \ge 1 \qquad P_p\big(N_{\mathcal{P}} \le an \,\big|\, |C(0)| = n\big) \le \alpha \exp(-\beta n).$$

More importantly, van der Hofstad and Kager prove the following two-pattern theorem, which is the main ingredient used in the proof of theorem 15.4.

**Theorem 15.6** (Two-pattern theorem, particular case of theorem 1.9 of [134])**.** *Let $\mathcal{P}, \mathcal{P}'$ be two distinct patterns such that*

$$P_p\big(\forall x \in Q \quad \omega(x) = \mathcal{P}(x)\big) \;=\; P_p\big(\forall x \in Q \quad \omega(x) = \mathcal{P}'(x)\big)\,. \tag{15.8}$$

*For any $\varepsilon > 0$, there exist positive constants $\alpha, \beta > 0$ such that*

$$\forall n \geq 1 \qquad P_p\big(\big|N_{\mathcal{P}} - N'_{\mathcal{P}}\big| \geq \varepsilon n \,\big|\, |C(0)| = n\big) \;\leq\; \alpha \exp(-\beta n)\,.$$

With the help of the two-pattern theorem, van der Hofstad and Kager obtained a general ratio limit theorem for percolation clusters. We restate this result for the Bernoulli percolation model at the critical point (hence $\mu = 1$ in corollary 1.11 of [134]).

**Corollary 15.7** (Ratio limit theorem)**.** *The following limits exist and satisfy*

$$\lim_{n\to\infty} \frac{P_{p_c}(|C(0)| = n+1)}{P_{p_c}(|C(0)| = n)} \;=\; \lim_{n\to\infty} \frac{P_{p_c}(|C(0)| \geq n+1)}{P_{p_c}(|C(0)| \geq n)} \;=\; 1\,.$$

In fact, this ratio limit theorem for percolation clusters had been previously derived by Madras [101], and the strategy of using patterns goes back to Kesten, who obtained an analogous result for self-avoiding walks in [81]. Yet the use of the two-pattern theorem in [134] yields a simpler and more elegant proof.

### 15.2.2 A strengthened ratio limit theorem

Our contribution here is a slight strengthening of corollary 15.7, in the following form.

**Theorem 15.8** (Strengthened ratio limit theorem)**.** *Let $\gamma > 0$ be fixed. There exists a positive constant $\kappa$ such that, for each $n \geq 1$, either*

$$P_{p_c}(|C(0)| = n) \;\leq\; \frac{1}{n^\gamma}$$

*or*

$$\left|\frac{P_{p_c}(|C(0)| = n+1)}{P_{p_c}(|C(0)| = n)} \;-\; 1\right| \;\leq\; \kappa\sqrt{\frac{\ln n}{n}}\,.$$

Notice that corollary 1.12 of [134] is a similar strengthening of the ratio limit theorem, however it requires that the cluster size distribution has a stretched exponential tail, and this is believed not to be true at the critical point.

*Proof.* We define $\mathcal{N}_{\mathcal{P}}$ as the set of the sites $x$ in $V$ where the pattern $\mathcal{P}$ occurs at $x$ on $C(0)$, i.e.,

$$\mathcal{N}_{\mathcal{P}} \;=\; \big\{\, x \in V : \text{the pattern } \mathcal{P} \text{ occurs at } x \text{ on } C(0) \big\}\,.$$

Notice that $\mathcal{N}_\mathcal{P}$ is a random set of sites, and that the random number $N_\mathcal{P}$ defined in (15.7) is precisely the cardinality of $\mathcal{N}_\mathcal{P}$. We shall use the same patterns $\mathcal{P}$, $\mathcal{P}'$ as in [134], we recall their definitions here for convenience. The pattern $\mathcal{P}$ is such that the site $(1,\dots,1)$ is occupied and the other sites of $Q$ are vacant, the pattern $\mathcal{P}'$ is such that the origin $(0,\dots,0)$ is occupied and the other sites of $Q$ are vacant, see figure 21. These two patterns $\mathcal{P},\mathcal{P}'$ satisfy the condition (15.8).

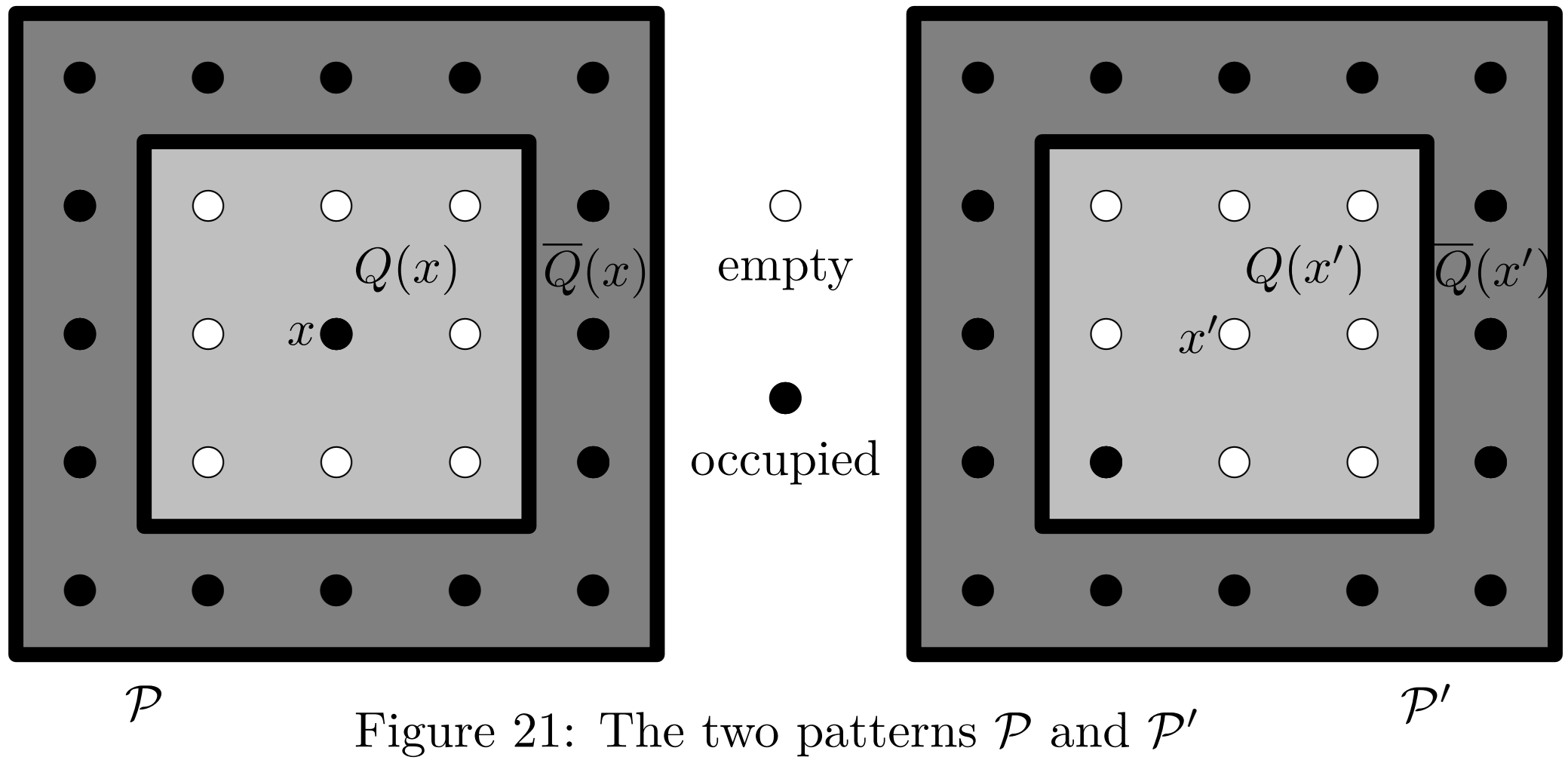

Figure 21: The two patterns $\mathcal{P}$ and $\mathcal{P}'$

Throughout the proof, we fix the parameter $p\in]0,1[$, and we write simply $P$ instead of $P_p$. We introduce the notation

$$\forall n\in\mathbb{N}\qquad c_n\,=\,P(|C(0)|=n)\,,$$

and also

$$\forall n\in\mathbb{N}\quad\forall i,j\in\mathbb{N}\qquad c_n(i,j)\,=\,P\big(|C(0)|=n,N_\mathcal{P}=i,N_{\mathcal{P}'}=j\big)\,.$$

Naturally, for $i>n$ or $j>n$, we have $c_n(i,j)=0$.

**Lemma 15.9.** *For any $n,j\in\mathbb{N}\setminus\{\,0\,\}$ and $i\in\mathbb{N}$, we have*

$$c_{n+1}(i+1,j-1)\,=\,\frac{j}{i+1}c_n(i,j)\,,\tag{15.9}$$

$$c_n(i,j)\,\leq\,2\exp\Big(-2\frac{(i-j)^2}{i+j}\Big)\sum_{k=n-i}^{n+j}c_k\,.\tag{15.10}$$

*Proof.* Given a percolation configuration, we introduce the random set $\mathcal{L}$ of the sites $x\in V$ such that $\partial\overline{Q}_x\subset C(0)$. We fix $n,j\in\mathbb{N}\setminus\{\,0\,\}$ and $i\in\mathbb{N}$ and we decompose $c_n(i,j)$ as follows:

$$\begin{aligned}c_n(i,j)\,&=\,P\big(|C(0)|=n,N_\mathcal{P}=i,N_{\mathcal{P}'}=j\big)\\&=\sum_{L\subset V}\sum_{\substack{K\subset L\\|K|=i+j}}P\begin{pmatrix}|C(0)|=n,N_\mathcal{P}=i,N_{\mathcal{P}'}=j\\\mathcal{L}=L,\mathcal{N}_\mathcal{P}\cup\mathcal{N}_{\mathcal{P}'}=K\end{pmatrix}\end{aligned}$$

$$= \sum_{L\subset V} \sum_{\substack{K\subset L\\|K|=i+j}} \sum_{A} P\begin{pmatrix} C(0)=A, N_{\mathcal{P}}=i, N_{\mathcal{P}'}=j \\ \mathcal{L}=L, \mathcal{N}_{\mathcal{P}}\cup\mathcal{N}_{\mathcal{P}'}=K \end{pmatrix}, \tag{15.11}$$

where the inner sum runs over the connected sets $A$ containing 0 of cardinality $n$. Of course, for the probability in the sum to be non-zero, the set $A$ must satisfy some additional constraints, so that the conditions $\mathcal{L}=L$, $\mathcal{N}_{\mathcal{P}}\cup\mathcal{N}_{\mathcal{P}'}=K$ can be realized; we call such sets $A$ admissible, and we denote by $\mathcal{A}(L,K)$ the collection of all these sets. To the subset $K$ of $V$, we associate the region

$$Q(K) = \bigcup_{x\in K} Q_x .$$

If $A$ is an admissible set belonging to $\mathcal{A}(L,K)$, then its trace on $Q(K)$ has cardinality exactly $i$ (because there exists a configuration such that $C(0)=A$, $\mathcal{N}_{\mathcal{P}}\cup\mathcal{N}_{\mathcal{P}'}=K$, $N_{\mathcal{P}}=i$, $N_{\mathcal{P}'}=j$), therefore its trace $A|_{Q(K)^c}$ on $Q(K)^c$ has cardinality exactly $n-i$. Moreover, all the sites belonging to the boundary of $Q(K)$ must be included in $A$, thus $A|_{Q(K)^c}$ is still connected. In order to condition on $C(0)$ outside $Q(K)$, we introduce the set

$$\widetilde{\mathcal{A}}(L,K) = \left\{ A|_{Q(K)^c} : A\in\mathcal{A}(L,K) \right\},$$

and we rewrite the expression of $c_n(i,j)$ obtained in (15.11) as

$$c_n(i,j) = \sum_{L\subset V} \sum_{\substack{K\subset L\\|K|=i+j}} \sum_{\widetilde{A}\in\widetilde{\mathcal{A}}(L,K)} P\begin{pmatrix} C(0)|_{Q(K)^c}=\widetilde{A}, N_{\mathcal{P}}=i, N_{\mathcal{P}'}=j \\ \mathcal{L}=L, \mathcal{N}_{\mathcal{P}}\cup\mathcal{N}_{\mathcal{P}'}=K \end{pmatrix}$$

$$= \sum_{L\subset V} \sum_{\substack{K\subset L\\|K|=i+j}} \sum_{\widetilde{A}\in\widetilde{\mathcal{A}}(L,K)} P\begin{pmatrix} C(0)|_{Q(K)^c}=\widetilde{A}, N_{\mathcal{P}}=i, N_{\mathcal{P}'}=j \\ \mathcal{N}_{\mathcal{P}}\cup\mathcal{N}_{\mathcal{P}'}=K \end{pmatrix}. \tag{15.12}$$

The last equality comes from the fact that, for any $\widetilde{A}\in\widetilde{\mathcal{A}}(L,K)$, the event $\{C(0)|_{Q(K)^c}=\widetilde{A}, \mathcal{N}_{\mathcal{P}}\cup\mathcal{N}_{\mathcal{P}'}=K\}$ is automatically included in the event $\{\mathcal{L}=L\}$. In the next step, we estimate the probability appearing inside (15.12), with the help of the following conditioning:

$$P\begin{pmatrix} C(0)|_{Q(K)^c}=\widetilde{A}, N_{\mathcal{P}}=i \\ N_{\mathcal{P}'}=j, \mathcal{N}_{\mathcal{P}}\cup\mathcal{N}_{\mathcal{P}'}=K \end{pmatrix} = \sum_{\substack{I,J\subset K, I\cup J=K\\|I|=i,|J|=j}} P\begin{pmatrix} C(0)|_{Q(K)^c}=\widetilde{A}, \mathcal{N}_{\mathcal{P}}=I \\ \mathcal{N}_{\mathcal{P}'}=J, \mathcal{N}_{\mathcal{P}}\cup\mathcal{N}_{\mathcal{P}'}=K \end{pmatrix}$$

$$= \sum_{\substack{I,J\subset K, I\cup J=K\\|I|=i,|J|=j}} P\left(\begin{matrix} \mathcal{N}_{\mathcal{P}}=I \\ \mathcal{N}_{\mathcal{P}'}=J \end{matrix}\,\middle|\, \begin{matrix} C(0)|_{Q(K)^c}=\widetilde{A} \\ \mathcal{N}_{\mathcal{P}}\cup\mathcal{N}_{\mathcal{P}'}=K \end{matrix}\right) P\begin{pmatrix} C(0)|_{Q(K)^c}=\widetilde{A} \\ \mathcal{N}_{\mathcal{P}}\cup\mathcal{N}_{\mathcal{P}'}=K \end{pmatrix}. \tag{15.13}$$

Now the conditional probability in the last sum is equal to

$$P\left(\begin{matrix} \mathcal{N}_{\mathcal{P}}=I \\ \mathcal{N}_{\mathcal{P}'}=J \end{matrix}\,\middle|\, \begin{matrix} C(0)|_{Q(K)^c}=\widetilde{A} \\ \mathcal{N}_{\mathcal{P}}\cup\mathcal{N}_{\mathcal{P}'}=K \end{matrix}\right)$$

$$= P\left(\begin{matrix} \forall x\in I \quad \omega|_{Q_x}=\mathcal{P} \\ \forall x\in J \quad \omega|_{Q_x}=\mathcal{P}' \end{matrix}\,\middle|\, \begin{matrix} C(0)|_{Q(K)^c}=\widetilde{A} \\ \forall x\in I\cup J \quad \omega|_{Q_x}=\mathcal{P} \text{ or } \omega|_{Q_x}=\mathcal{P}' \end{matrix}\right)$$

$$= \ P\begin{pmatrix} \forall x\in I \quad \omega|_{Q_x}=\mathcal{P} \\ \forall x\in J \quad \omega|_{Q_x}=\mathcal{P}' \end{pmatrix} \begin{matrix} \forall x\in I\cup J \\ \omega|_{Q_x}=\mathcal{P}\text{ or }\mathcal{P}' \end{matrix} \Bigg) \ = \ \frac{1}{2^{|I|+|J|}}\,, \tag{15.14}$$

because the two patterns $\mathcal{P}$, $\mathcal{P}'$ have the same probability of occurring. Substituting (15.14) in (15.13), we obtain

$$P\begin{pmatrix} C(0)|_{Q(K)^c}=\widetilde{A}, N_{\mathcal{P}}=i \\ N_{\mathcal{P}'}=j, \mathcal{N}_{\mathcal{P}}\cup\mathcal{N}_{\mathcal{P}'}=K \end{pmatrix} \ = \ \binom{i}{i+j}\frac{1}{2^{i+j}}P\begin{pmatrix} C(0)|_{Q(K)^c}=\widetilde{A} \\ \mathcal{N}_{\mathcal{P}}\cup\mathcal{N}_{\mathcal{P}'}=K \end{pmatrix}. \tag{15.15}$$

Reporting now (15.15) in (15.12), we have the formula

$$c_n(i,j) \ = \ \binom{i}{i+j}\frac{1}{2^{i+j}}\sum_{L\subset V}\sum_{\substack{K\subset L\\|K|=i+j}}\sum_{\widetilde{A}\in\widetilde{\mathcal{A}}(L,K)} P\begin{pmatrix} C(0)|_{Q(K)^c}=\widetilde{A} \\ \mathcal{N}_{\mathcal{P}}\cup\mathcal{N}_{\mathcal{P}'}=K \end{pmatrix}. \tag{15.16}$$

Notice that the set $\widetilde{\mathcal{A}}(L,K)$ associated to the subsets $K,L$ and the integers $n,i,j$ is the same as the set $\widetilde{\mathcal{A}}(L,K)$ associated to the subsets $K,L$ and the integers $n+1,i+1,j-1$. In fact, the set $\widetilde{\mathcal{A}}(L,K)$ depends on the integers $n,i,j$ only through $n-i$ (the number of occupied sites in $A$ outside $Q(K)$) and $i+j$ (the cardinality of $L$). Formula (15.16) implies therefore that there exists a function $f$ defined on $\mathbb{N}\times\mathbb{N}$ such that

$$c_n(i,j) \ = \ \binom{i}{i+j}\frac{1}{2^{i+j}}f(n-i,i+j)\,.$$

This implies furthermore that, for any $n,j\in\mathbb{N}\setminus\{\,0\,\}$ and $i\in\mathbb{N}$, we have

$$\forall i\geq 1 \qquad \binom{i+1}{i+j}c_n(i,j) \ = \ \binom{i}{i+j}c_{n+1}(i+1,j-1)\,,$$

which yields formula (15.9). To get the inequality (15.10), we bound the sums appearing in (15.16) as follows (recall that the event in the first line of next formula is included in $\{\,\mathcal{L}=L\,\}$):

$$\begin{aligned}
&\sum_{L\subset V}\sum_{\substack{K\subset L\\|K|=i+j}}\sum_{\widetilde{A}\in\widetilde{\mathcal{A}}(L,K)} P\begin{pmatrix} C(0)|_{Q(K)^c}=\widetilde{A} \\ \mathcal{N}_{\mathcal{P}}\cup\mathcal{N}_{\mathcal{P}'}=K \end{pmatrix} \\
&\qquad\qquad \leq \sum_{L\subset V}\sum_{\substack{K\subset L\\|K|=i+j}} P\begin{pmatrix} \big|C(0)|_{Q(K)^c}\big|=n-i \\ \mathcal{L}=L, \mathcal{N}_{\mathcal{P}}\cup\mathcal{N}_{\mathcal{P}'}=K \end{pmatrix} \\
&\qquad\qquad \leq \sum_{L\subset V}\sum_{\substack{K\subset L\\|K|=i+j}} P\begin{pmatrix} n-i\leq\big|C(0)\big|\leq n+j \\ \mathcal{L}=L, \mathcal{N}_{\mathcal{P}}\cup\mathcal{N}_{\mathcal{P}'}=K \end{pmatrix} \\
&\qquad\qquad\qquad\qquad \leq \ P\big(n-i\leq\big|C(0)\big|\leq n+j\big)\,.
\end{aligned} \tag{15.17}$$

Substituting the bound (15.17) in (15.16), we obtain

$$c_n(i,j) \ \leq \ \binom{i}{i+j}\frac{1}{2^{i+j}}\sum_{k=n-i}^{n+j}c_k\,. \tag{15.18}$$

We use next a classical bound on the binomial coefficients:

$$\binom{i}{i+j}\frac{1}{2^{i+j}}\;\leq\;2\exp\Big(-\frac{(i-j)^2}{2(i+j)}\Big)\,. \tag{15.19}$$

This inequality is for instance a consequence of Hoeffding's inequality applied to a Binomial random variable $S_{i+j}$ with parameters $(i+j,1/2)$:

$$\begin{aligned}\binom{i}{i+j}\frac{1}{2^{i+j}}\;&=\;P(S_{i+j}=i)\\ &\leq\;P\Big(\Big|S_{i+j}-\frac{i+j}{2}\Big|\;\geq\;\frac{|i-j|}{2}\Big)\;\leq\;2\exp\Big(-\frac{(i-j)^2}{2(i+j)}\Big)\,.\end{aligned}$$

Plugging (15.19) into (15.18), we get the inequality (15.10). □

We come back to the proof of theorem 15.8. It follows from theorem 15.5 that there exist positive constants $a,\alpha,\beta>0$ such that

$$\forall n\geq 0\qquad \sum_{i\leq an}P\big(N_{\mathcal{P}}=i\,\big|\,|C(0)|=n\big)\;\leq\;\alpha\exp(-\beta n)\,. \tag{15.20}$$

Using (15.20), we have, for any $n\geq 0$,

$$\sum_{i\leq an}\sum_{j\geq 0}c_n(i,j)\;\leq\;P\big(|C(0)|=n,N_{\mathcal{P}}\leq an\big)\;\leq\;\alpha\exp(-\beta n)\,. \tag{15.21}$$

Symmetrically, we have

$$\forall n\geq 0\qquad \sum_{i\geq 0}\sum_{j\leq an}c_n(i,j)\;\leq\;\alpha\exp(-\beta n)\,. \tag{15.22}$$

We shall compare $c_{n+1}$ and $c_n$. Let $n\geq 0$ be fixed and let us decompose $c_{n+1}$ as follows:

$$c_{n+1}\;=\;\sum_{i\geq 0}\sum_{j\geq 0}c_{n+1}(i,j)\,. \tag{15.23}$$

Of course the quantity $c_n(i,j)$ is equal to 0 if $i$ or $j$ is strictly larger than $n$, so the above sum is in fact finite, as well as all the sums appearing below. With the help of the inequalities (15.21) and (15.22), we control all the terms in the decomposition (15.23) for which $i\leq a(n+1)$ or $j\leq a(n+1)$, and in fact all the terms for which $i\leq a'(n+1)$ or $j\leq a'(n+1)$ where $a'$ is any number smaller than $a$. For a technical reason, we work with $a'=a/2$ and we get

$$c_{n+1}\;=\;\sum_{i>a'(n+1)}\;\sum_{j>a'(n+1)}c_{n+1}(i,j)\;+\;O\big(\exp(-\beta n)\big)\,. \tag{15.24}$$

We use next the identity (15.9) of lemma 15.9, and we obtain

$$c_{n+1}\;=\;\sum_{i>a'(n+1)}\;\sum_{j>a'(n+1)}\;\frac{j+1}{i}c_n(i-1,j+1)\;+\;O\big(\exp(-\beta n)\big)\,.$$

We reindex the sums with $i' = i - 1$ and $j' = j + 1$:

$$c_{n+1} \,=\, \sum_{i' > a'(n+1)-1} \; \sum_{j' > a'(n+1)+1} \frac{j'}{i'+1} c_n(i',j') \,+\, O\big(\exp(-\beta n)\big)\,. \qquad (15.25)$$

Proceeding in the same way that led to (15.24), we have for $c_n$ that

$$c_n \,=\, \sum_{i > a'n} \sum_{j > a'n} c_n(i,j) \,+\, O\big(\exp(-\beta n)\big)\,. \qquad (15.26)$$

Subtracting (15.26) to (15.25), we obtain

$$\begin{aligned} c_{n+1} - c_n \,=\, & \sum_{i > a'n} \; \sum_{j > a'(n+1)+1} \Big(\frac{j}{i+1} - 1\Big) c_n(i,j) \\ & + \sum_{\substack{i > a'(n+1)-1 \\ i \leq a'n}} \; \sum_{j > a'(n+1)+1} \frac{j}{i+1} c_n(i,j) \\ & - \sum_{i > a'n} \; \sum_{\substack{j > a'n \\ j \leq a'(n+1)+1}} c_n(i,j) \,+\, O\big(\exp(-\beta n)\big)\,. \qquad (15.27) \end{aligned}$$

We control separately each sum. We suppose that $n$ is large enough to ensure that

$$\frac{an}{3} \,\leq\, a'(n+1) - 1\,, \qquad a'(n+1) + 1 \,\leq\, an\,.$$

For the second sum in (15.27), we have

$$\sum_{\substack{i > a'(n+1)-1 \\ i \leq a'n}} \; \sum_{j > a'(n+1)+1} \frac{j}{i+1} c_n(i,j) \,\leq\, \frac{n}{an/3} \sum_{i \leq a'n} \sum_{j \geq 0} c_n(i,j) \,\leq\, \frac{3\alpha}{a} \exp(-\beta n)\,. \qquad (15.28)$$

For the third sum in (15.27), we have

$$\sum_{i > a'n} \; \sum_{\substack{j > a'n \\ j \leq a'(n+1)+1}} c_n(i,j) \,\leq\, \sum_{i \geq 0} \sum_{j \leq an} c_n(i,j) \,\leq\, \alpha \exp(-\beta n)\,. \qquad (15.29)$$

Plugging (15.28) and (15.29) in (15.27) and taking the absolute values, we obtain

$$\big|c_{n+1} - c_n\big| \,\leq\, \sum_{i > a'n} \sum_{j > a'n} \Big|\frac{j}{i+1} - 1\Big| c_n(i,j) \,+\, O\big(\exp(-\beta n)\big)\,.$$

It remains to control this last sum. Let $m \geq 1$ be a parameter, to be chosen later. We split the sum in two and we write, with the help of inequality (15.10) of lemma 15.9,

$$\big|c_{n+1}-c_n\big| \le \sum_{i>a'n}\sum_{j:|j-i-1|\le m}\cdots + \sum_{i>a'n}\sum_{j:|j-i-1|>m}\cdots + O\big(\exp(-\beta n)\big)$$
$$\le \frac{m}{a'n}\sum_{i,j}c_n(i,j) + \sum_{\substack{i,j\le n\\ |j-i-1|>m}}\frac{4}{a'}\exp\Big(-2\frac{(j-i)^2}{i+j}\Big) + O\big(\exp(-\beta n)\big)$$
$$\le \frac{m}{a'n}c_n + \frac{4n^2}{a'}\exp\Big(-\frac{m^2}{n}\Big) + O\big(\exp(-\beta n)\big)\,.$$

The time to choose $m$ has come. We take $m$ of the form $m=\sqrt{(\gamma+3)n\ln n}$, and we obtain

$$\big|c_{n+1}-c_n\big| \le \frac{1}{a'}\sqrt{(\gamma+3)\frac{\ln n}{n}}c_n + \frac{4}{a'n^{\gamma+1}} + O\big(\exp(-\beta n)\big)\,.$$

Thus we have proved the following:

$$\exists b>0\quad \forall n\ge 1\qquad \big|c_{n+1}-c_n\big| \le b\sqrt{\frac{\ln n}{n}}c_n + \frac{b}{n^{\gamma+1}}\,. \tag{15.30}$$

Now, if $n$ is such that $c_n\ge n^{-\gamma}$, then

$$b\sqrt{\frac{\ln n}{n}}c_n + \frac{b}{n^{\gamma+1}} \le 2b\sqrt{\frac{\ln n}{n}}c_n\,,$$

and this concludes the proof of theorem 15.8. □

We complete now the proof of theorem 15.4.

*Proof.* Let $n,j\ge 1$ and let us sum the inequality (15.30) with $\gamma=1$ from $n$ to $n+j-1$:

$$\big|c_{n+j}-c_n\big| \le \sum_{i=0}^{j-1}\left(b\sqrt{\frac{\ln(n+i)}{n+i}}c_{n+i} + \frac{b}{(n+i)^2}\right)$$
$$\le b\sqrt{\frac{\ln(n+j)}{n}}\sum_{i=0}^{j-1}c_{n+i} + \frac{bj}{n^2}\,.$$

This inequality implies that

$$c_{n+j} \ge c_n - b\sqrt{\frac{\ln(n+j)}{n}}\sum_{i=0}^{j-1}c_{n+i} - \frac{bj}{n^2}\,. \tag{15.31}$$

Let $k\ge 1$. We sum again (15.31) over $j$ between $j=0$ and $j=k-1$, and we obtain

$$\sum_{j=0}^{k-1}c_{n+j} \ge kc_n - bk\sqrt{\frac{\ln(n+k)}{n}}\sum_{i=0}^{k-1}c_{n+i} - \frac{bk^2}{n^2}\,. \tag{15.32}$$

Recalling that $\sum_{n\ge 0}c_n \le 1$, and taking $k=\sqrt{n/\ln n}$ in (15.32), we get the desired conclusion. □

## 16 The large clusters of $\Lambda(4n)$

We did not succeed in obtaining an explicit control on the function $\theta(n,p)$, yet we can still use the bound on the travel time to describe typical configurations in large finite boxes. We pursue this program here.

**Proposition 16.1.** *Let $p\in[0,1]$ be such that $\theta(p)>0$. There exists a non-decreasing function $\beta(n)$ defined for $n\geq 3$ such that*

$$\forall n\geq 3\qquad \beta(n)\leq\sqrt{\ln n}\,,$$
$$\beta(3)\geq\theta(p)\,,\qquad \lim_{n\to+\infty}\beta(n)=+\infty\,,$$

*and moreover*

$$\forall n\geq 3\qquad P\Big(\exists x\in\Lambda(4n)\quad T_{\Lambda(4n)}\big(x,\partial^{\,in}\Lambda(4n)\big)\,>\,\frac{\ln n}{\beta(n)}\Big)\,\leq\,n^{3d-\beta(n)}\,. \tag{16.1}$$

*Proof.* We apply corollary 14.5 to the box $\Lambda(4n)$ and with a starting point $x$ in $\Lambda(4n)$:

$$\forall x\in\Lambda(4n)\quad\forall t\in\{\,0,\dots,2n\,\}$$
$$P\big(T_{\Lambda(4n)}\big(x,\partial^{\,in}\Lambda(4n)\big)\,>\,t\,\big)\,\leq\,\exp\big(-t\alpha(t)\big)\,.$$

From the standard union bound, it follows that

$$\forall t\in\{\,0,\dots,2n\,\}$$
$$P\Big(\exists x\in\Lambda(4n)\quad T_{\Lambda(4n)}\big(x,\partial^{\,in}\Lambda(4n)\big)\,>\,t\Big)\,\leq\,|\Lambda(4n)|\exp\big(-t\alpha(t)\big)\,. \tag{16.2}$$

Let $\beta(n)$ be the function defined by

$$\forall n\geq 3\qquad \beta(n)\,=\,\min\Big(\sqrt{\ln n},\sqrt{\alpha\big(\sqrt{\ln n}\big)}\Big)\,. \tag{16.3}$$

Thanks to the fact that $\lim_{t\to\infty}\alpha(t)=+\infty$, we have also

$$\lim_{n\to\infty}\beta(n)\,=\,+\infty\,.$$

Recalling that the function $\alpha(t)$ is non–decreasing, we have furthermore

$$\frac{\ln n}{\beta(n)}\,\geq\,\sqrt{\ln n}\,,\qquad \alpha\Big(\frac{\ln n}{\beta(n)}\Big)\,\geq\,\alpha\big(\sqrt{\ln n}\big)\,\geq\,\beta(n)^2\,, \tag{16.4}$$

whence

$$\forall n\geq 3\qquad \frac{1}{\beta(n)}\alpha\Big(\frac{\ln n}{\beta(n)}\Big)\,\geq\,\beta(n)\,. \tag{16.5}$$

Since $\alpha(1)\geq\theta(p)$ by corollary 14.5, we have also

$$\beta(3)\,=\,\min\Big(\sqrt{\ln 3},\sqrt{\alpha\big(\sqrt{\ln 3}\big)}\Big)\,\geq\,\min\big(1,\sqrt{\alpha(1)}\big)\,\geq\,\sqrt{\theta(p)}\,\geq\,\theta(p)\,.$$

We take $t=(\ln n)/\beta(n)$ in the inequality (16.2) and, using (16.5) and the simple bound $|\Lambda(4n)|\leq(6n)^d\leq n^{3d}$, we get the desired inequality (16.1). □

In proposition 16.1, we saw that, with very high probability, the travel time between any site of $\Lambda(2n)$ and $\partial^{in}\Lambda(4n)$ is negligible compared to $\ln n$. Yet a path which travels a distance of order $n$ in a time less than $\ln n$ has to go through at least one cluster of diameter larger than $n/\ln n$. So we conclude that, typically, any site of $\Lambda(2n)$ is within travel distance at most $\ln n$ of a cluster of $\Lambda(4n)$ having diameter larger than $n/\ln n$. We provide next a precise formulation for this important result.

We define the diameter $\mathrm{diam}_\infty A$ of a subset $A$ of $\mathbb{Z}^d$ as its diameter for the norm $|\cdot|_\infty$, i.e.,

$$\mathrm{diam}_\infty A \;=\; \sup\big\{\,|x-y|_\infty : x,y\in A\,\big\}\,.$$

Let $L\geq 1$ be a positive integer. We say that an open cluster $C$ is $L$-large if $\mathrm{diam}_\infty C\geq L$. We denote by $\mathcal{C}(n,L)$ the collection of all the $L$-large open clusters of the percolation configuration restricted to the box $\Lambda(n)$.

**Proposition 16.2.** *Let $p\in[0,1]$ be such that $\theta(p)>0$. There exists a non-decreasing function $\beta(n)$ defined for $n\geq 3$ such that*

$$\begin{gathered}\forall n\geq 3\qquad \beta(n)\leq\sqrt{\ln n}\,,\\ \beta(3)\geq\theta(p)\,,\qquad \lim_{n\to+\infty}\beta(n)=+\infty\,,\end{gathered}$$

*and if we set*

$$\forall n\geq 3\qquad \gamma(n)\;=\;\Big\lfloor\frac{\beta(n)n}{2\ln n}\Big\rfloor\,,\tag{16.6}$$

*then we have*

$$\forall n\geq 3\qquad P\begin{pmatrix}\forall x\in\Lambda(2n)\quad \exists C\in\mathcal{C}(4n,\gamma(n))\\ T_{\Lambda(4n)}(x,C)\;\leq\;(\ln n)/\beta(n)\end{pmatrix}\;\geq\;1-n^{3d-\beta(n)}\,.$$

*Proof.* Let $\beta(n)$ be the function given by proposition 16.1 and let $\mathcal{E}_n$ be the event

$$\mathcal{E}_n\;=\;\Big\{\forall x\in\Lambda(2n)\quad T_{\Lambda(4n)}\big(x,\partial^{in}\Lambda(4n)\big)\;\leq\;\frac{\ln n}{\beta(n)}\Big\}\,.$$

It follows from the inequality (16.1) of proposition 16.1 that

$$\forall n\geq 3\qquad P(\mathcal{E}_n)\;\geq\;1-n^{3d-\beta(n)}\,.$$

Suppose that the event $\mathcal{E}_n$ occurs and let $x$ belong to $\Lambda(2n)$. There exists a path $z_0,\dots,z_r$ from $z_0=x$ to a point $z_r=y$ in $\partial^{in}\Lambda(4n)$ such that

$$T_{\Lambda(4n)}(x,y)\;=\;\sum_{i=0}^{r-1}1_{z_i\text{ closed}}\;\leq\;\frac{\ln n}{\beta(n)}\,.$$

The diameter of this path is larger than or equal to

$$|x-y|_\infty\;=\;|z_0-z_r|_\infty\;\geq\;\inf\{\,|a-b|_\infty : a\in\Lambda(2n),\,b\in\partial^{in}\Lambda(4n)\,\}\;\geq\;n\,.\tag{16.7}$$

Let $\gamma(n)$ be defined as in (16.6). To the path $z_0,\dots,z_r$, we associate its $\gamma(n)$-skeleton (a classical construction which was instrumental in the proof of the Wulff construction in two dimensions, see [6]), defined by the following algorithm. We set $i_0=0$. Suppose that for some $k\geq 0$, the indices $i_0,\dots,i_k$ have been defined and that $i_k<r$. If

$$\forall i\in\{\,i_k,\dots,r\,\}\qquad |z_i-z_{i_k}|_\infty<\gamma(n)\,,$$

then we set $i_{k+1}=r$ and the algorithm terminates. Otherwise, we define

$$i_{k+1}\;=\;\min\,\big\{\,i>i_k:|z_i-z_{i_k}|_\infty\geq\gamma(n)\,\big\}\,,$$

and we iterate the construction. This algorithm terminates at some step $\ell\leq r$ with the index $i_\ell=r$. The $\gamma(n)$-skeleton of the path $z_0,\dots,z_r$ is the sequence

$$s_0=z_0,\,s_1=z_{i_1},\dots,s_{\ell-1}=z_{i_{\ell-1}},\,s_\ell=z_r\,.$$

The very construction of the $\gamma(n)$-skeleton ensures that

$$\forall j\in\{\,0,\dots,\ell-2\,\}\qquad \big|s_j-s_{j+1}\big|_\infty\;=\;\gamma(n)\,,$$

and $\big|s_{\ell-1}-s_\ell\big|_\infty\;\leq\;\gamma(n)$. Summing the previous inequalities, we get

$$\big|z_0-z_r\big|_\infty\;=\;\big|s_0-s_\ell\big|_\infty\;\leq\;\ell\gamma(n)\,,$$

which, together with (16.7), yields

$$\ell\;\geq\;\frac{n}{\gamma(n)}\,. \tag{16.8}$$

Let us consider the disjoint sub-intervals

$$\{\,0,\dots,i_1-1\,\},\,\{\,i_1,\dots,i_2-1\,\},\,\dots,\,\{\,i_{\ell-2},\dots,i_{\ell-1}-1\,\}\,.$$

There are $\ell-1$ of them. Now, from (16.8), the choice (16.6) of $\gamma(n)$ and the first inequality of formula (16.4), we have, for $n\geq 3$,

$$\ell-1\;\geq\;\frac{n}{\gamma(n)}-1\;\geq\;\frac{2\ln n}{\beta(n)}-1\;\geq\;\frac{\ln n}{\beta(n)}+\sqrt{\ln n}-1\;>\;\frac{\ln n}{\beta(n)}\,.$$

The path $z_0,\dots,z_{r-1}$ contains at most $(\ln n)/\beta(n)$ closed sites, therefore at least one of the sub-intervals, say $\{\,i_k,\dots,i_{k+1}-1\,\}$, is such that the sites of the associated sub-path $\{\,z_{i_k},z_{i_k+1},\dots,z_{i_{k+1}-1}\,\}$ are all open. Let $C$ be the open cluster of $z_{i_k}$ in $\Lambda(4n)$. This cluster $C$ contains both sites $z_{i_k}$ and $z_{i_{k+1}-1}$, therefore

$$\mathrm{diam}_\infty C\;\geq\;\big|z_{i_k}-z_{i_{k+1}-1}\big|_\infty\;=\;\gamma(n)\,,$$

and moreover $T_{\Lambda(4n)}(x,C)\;\leq\;(\ln n)/\beta(n)\,.$ □

Pictorially, the $\gamma(n)$-large clusters of the box $\Lambda(4n)$ are black holes that rapidly absorb all the matter in the box $\Lambda(2n)$ at lightning speed. They will play a central role in the final steps of the proof.

# 17 Shell exploration ◦

In the previous section 16, we have seen that, with high probability, each site of the box $\Lambda(2n)$ is within sub-logarithmic travel distance from the $\gamma(n)$-large clusters of the box $\Lambda(4n)$. A central problem for proving the conjecture $\theta(p_c, \mathbb{Z}^d) = 0$ is to get a control on the number of $\gamma(n)$-large clusters. Suppose that many $\gamma(n)$-large clusters coexist in the box $\Lambda(4n)$. On the one hand, this requires the existence of many long open connections in the box, but this is not problematic when $\theta(p)$ is positive. On the other hand, for these connections to belong to disjoint clusters, it must be the case that they are surrounded by closed sites, and this requires the presence of many separating sets of closed sites. This is much more problematic when $\theta(p)$ is positive, or at least this is what we hope to prove. To accomplish this daunting task, we will design a specific algorithm for exploring the configuration in order to reveal the presence of the large interfaces separating the $\gamma(n)$-large clusters: the intertwined explorations.

## 17.1 Concurrent exploration algorithms

The standard exploration algorithm described in section 5 is sequential, and it explores at most one site at each iteration. This algorithm has been used in computer science since the early days of digital image analysis (see the pioneering book [116]). It belongs to the class of "growing regions algorithms". This algorithm can be speeded up considerably, notably a lot of site explorations can be performed concurrently or in parallel. The design of massively parallel algorithms to find the connected components of a graph is still an active area of research (see [22, 33, 96] for very recent works on this question). In these algorithms, several explorers are launched simultaneously, and they explore the configuration in parallel. The main issue to tackle is the communication between the explorers, as well as their synchronization. Indeed, the explorers have to share some information during the exploration process, and this has to be orchestrated in an efficient way. Neither do we wish that an explorer is stalled because he is waiting for the others to complete their jobs, nor do we wish that some sites get explored several times by different explorers. These would lead to a total waste of the computing resources. We will not delve into these delicate algorithmic questions, however in the construction of our random processes we will face subtle problems which are reminiscent of these. Indeed, the intertwined explorations will exchange information through taboo sets and this will create a form of stochastic dependence. Yet we will make it sure that the mathematical constructions always correspond to sequential algorithms which reveal at most one site at each step, so we will avoid the more complex issues raised by concurrent algorithms. In fact, we will put more emphasis on the synchronization aspect: the intertwined explorations are the mathematical counterpart of a basic mutex process (for mutual exclusion) in computer science.

Before starting this program, we reformulate the basic exploration algorithm into a shell exploration algorithm for which we do not specify the order of the exploration of the set of the active sites.

## 17.2 The shell exploration algorithm

The shell exploration algorithm takes as input a starting set $A$ and an authorized domain $D$. The set $D$ is a finite subset of $\mathbb{Z}^d$ and $A$ is a subset of $D$. Upon termination, the shell exploration algorithm will have explored the shell around $A$ in $D$, defined in subsection 12.3 as

$$\text{Shell}\,(A,D)\;=\;\big\{\,x\in D: T_D(A,x)=0\,\big\}\,.$$

The set $\text{Shell}\,(A,D)$ is the union of the clusters of the sites of the initial set $A$, together with their outer boundaries. Moreover, these clusters are the clusters for the percolation configuration restricted to $D$. The algorithm returns two sets, the set $S=\text{Shell}\,(A,D)$ and the set

$$W\;=\;\big\{\,x\in \text{Shell}\,(A,D): x \text{ is closed}\,\big\}\,. \tag{17.1}$$

The set $W$ consists of the closed sites of $\text{Shell}\,(A,D)$. These are the sites waiting for the next exploration round, once the exploration of $\text{Shell}\,(A,D)$ is completed. The pseudocode of the shell exploration algorithm is given in Algorithm 17.1.

**Algorithm 17.1** Explore_Shell $(A,D)$

**Require:** an authorized domain $D$, a non-empty starting set $A$ included in $D$
**Ensure:** $S$ is the set of the sites of $\text{Shell}\,(A,D)$, $W$ is the set of waiting sites

```
S ← ∅, W ← ∅
repeat
    Pick up x ∈ A          ▷ This step can take advantage of the parallelism
    S ← S ∪ {x}
    if x is closed then
        A ← A \ {x}
        W ← W ∪ {x}
    else if x is open then
        A ← A ∪ (N(x) ∩ D) \ ({x} ∪ S)
    end if
until A = ∅
return S, W
```

Among the waiting sites, we can distinguish the subset of the sites having a neighbour which have not been visited. These are the interesting waiting sites, because the other waiting sites will just be cancelled during the next round of exploration. An alternative possibility would therefore be to take the set of the waiting sites as the set

$$\partial_D^{in}\text{Shell}\,(A,D)\;=\;\big\{\,x\in \text{Shell}\,(A,D): x \text{ is closed}, \mathcal{N}(x)\cap D\setminus A\neq\varnothing\,\big\}\,.$$

The advantage is that it avoids storing sites that do not lead to an unvisited site. We will adopt this point of view when defining the waiting site aggregates in the subsection 80.1. The drawback is that it might render some of the formulas more complicated. Upon reflection after completing the writing, if we were to do it again, we would probably stick to the choice (17.1).

## 17.3 Revisiting the travel time algorithm

Once we have at our disposition the shell exploration algorithm, we can easily build an algorithm to compute the travel time $T = T_D(A, B)$ in a finite set $D$ between two subsets $A, B$ of $D$. The pseudocode of the travel time algorithm is given in Algorithm 17.2.

---

**Algorithm 17.2** Travel_Time $T_D(A, B)$

---

**Require:** a domain $D$, a starting set $A$, a target set $B$ disjoint from $A$
**Ensure:** $T$ is the travel time from $A$ to $B$ inside $D$
  $A(0) \leftarrow A$, $t \leftarrow 0$
  $S, W =$ Explore_Shell $(A(0), D)$
  $E(0) \leftarrow S$, $W(0) \leftarrow W$
  **repeat** ▷ $t$ is the number of the last shell which has been explored
    $A(t+1) \leftarrow \varnothing$
    **while** $W(t) \neq \varnothing$ **do** ▷ Building the set of active sites
      Pick up $x \in W(t)$
      $A(t+1) \leftarrow A(t+1) \cup \mathcal{N}(x, D) \setminus E(t)$
      $W(t) \leftarrow W(t) \setminus \{x\}$
    **end while**
    $S, W =$ Explore_Shell $(A(t+1), D \setminus E(t))$
    $E(t+1) \leftarrow E(t) \cup S$, $W(t+1) \leftarrow W$
    $t \leftarrow t + 1$
  **until** $W(t) = \varnothing$ or $E(t) \cap B \neq \varnothing$
  $T \leftarrow t$
  **return** $T$

---

This algorithm updates iteratively three sets of sites, denoted by $E$, $A$, $W$ and it computes three sequences $\big(E(t),\ A(t),\ W(t),\ 0 \leq t \leq T_D(A, B)\big)$ of subsets of $D$, whose contents are the following:

- $E(t)$: the sites which have been explored during the first $t$ iterations;
- $A(t)$: the active sites, from where the round $t$ of the exploration starts;
- $W(t)$: the waiting sites, where the round $t$ of the exploration is stuck.

For $t$ in $\{0, \dots, T\}$, the set $E(t)$ is the set of the sites which are at travel distance less than or equal to $t$ from $A$ inside $D$, i.e.,

$$E(t) \,=\, \bigcup_{0 \leq s \leq t} \text{Shell}\,(A, D, s) \,=\, \mathcal{B}(A, D, t)\,,$$

while the set $W(t)$ of the waiting sites at time $t$ is

$$W(t) \,=\, \big\{\, x \in \text{Shell}\,(A, D, t) : x \text{ is closed} \,\big\}\,.$$

We set $A(0) = A$ and, for $t \geq 1$, the set $A(t)$ consists of the sites of $D \setminus E(t-1)$ which have a neighbour in $W(t-1)$, i.e.,

$$A(t) \,=\, \big\{\, x \in D \setminus E(t-1) : \exists\, y \in W(t-1) \quad x \in \mathcal{N}(y) \,\big\}\,.$$

## 17.4 Intertwined explorations

Let $\Lambda$ be a box and let $W$, $W'$ be two disjoint subsets of $\Lambda$. We consider a percolation configuration in $\Lambda$ in which the sets $W$ and $W'$ are disconnected. This disconnection event entails the presence of a kind of interface between $W$ and $W'$, that is a set of closed sites separating them. If the sets $W$ and $W'$ are quite large, this event is unlikely in the supercritical regime, and we wish to control the probability of such a scenario starting with the sole hypothesis that $\theta(p) > 0$. Our strategy consists in exploring the configuration in an adequate way in order to reveal an interface between $W$ and $W'$. This is accomplished by a combination of two exploration processes, which we call the intertwined explorations. The basic component for building the intertwined explorations is the shell exploration algorithm. We recall next its relevant features.

**The shell exploration algorithm.** This algorithm is initialized with a starting set and it explores the successive shells around this starting set. These shells are naturally explored in increasing order, starting with the 0-th shell. Let $k \geq 0$ be fixed. We call the exploration of the $k$-th shell the $k$-th round of the algorithm. The $k$-th shell is fully explored before the exploration of the $(k+1)$-th shell. After the completion of the $k$-th round, the algorithm has built a set of sites called the waiting set, which consists of the vertices of the $(k+1)$-th shell reachable from the $k$-th shell. To explore the $(k+1)$-th shell, the algorithm will explore all the clusters of the vertices of the waiting set. We explain next the general principle underlying the construction of the intertwined explorations.

**Principle of the construction.** To each starting set $W$ and $W'$, we associate an explorer. The explorer starting from $W$ (respectively $W'$) is called the first (respectively second) explorer. Both explorers will perform a shell exploration algorithm, but with an essential modification. Indeed, we will create an interaction between the two explorers, as follows. We start by calling the first explorer and we let him complete his 0-th round. We then call the second explorer, but we inform him of the region explored by the first explorer, and we let him complete his 0-th round. Let $k \geq 0$ and suppose that each explorer has completed his $k$-th round. We inform the first explorer of the region explored by the second explorer, and the first explorer is asked to perform his $(k+1)$-th round, but he is forbidden to explore a cluster already visited by the second explorer. Practically, the vertices of the waiting set of the first explorer after his $k$-th round which belong to a cluster visited by the second explorer are ignored. Finally, we inform the second explorer of the region explored by the first explorer after his $(k+1)$-th round, and the second explorer is asked to perform his $(k+1)$-th round, but he is forbidden to explore a cluster already visited by the first explorer. As was the case for the first explorer, the vertices of the waiting set of the second explorer after his $k$-th round which belong to a cluster visited by the first explorer are ignored. Since the sets $W$ and $W'$ are disconnected, the 0-th rounds of the two explorers are not affected by the new rule. Although this mechanism sounds relatively simple, it is delicate to implement. In order to enforce the rule that each explorer is forbidden to enter the territory explored by the other, we will define a taboo shell exploration algorithm. The taboo set is a forbidden region, that the explorer cannot enter. Later on, the taboo set will be replaced by the region explored by the inactive explorer. This mechanism induces an exchange of information between the two explorers. A crucial point is to decide which information is to be transmitted, as we discuss next.

**Information transmission.** The mathematical definitions of the random objects can be encapsulated at various levels and the same is true for the corresponding algorithms. For instance, the travel time algorithm 17.2 is built upon the shell exploration algorithm 17.1. This might look innocent, however this has consequences for the information transmitted upon termination of the algorithm. In the presentation of algorithms 17.1 and 17.2, the final output of the algorithm is simply the travel time between the input sets $A$ and $B$. We could ask the travel time algorithm 17.2 to transmit in addition the successive shells that are explored. However, the states of the sites in the explored shells would still stay partially hidden, because this information is buried deep down in the shell exploration algorithm and it is not transmitted to the higher level. When we build explorers which interact and transmit information to each other, it is essential to decide precisely what information are to be transmitted: this changes fundamentally the probabilistic estimates that can be derived. For instance, if we know that a site has been explored (without knowing its state), there is still a chance that it is open. The same delicate questions arise when building random objects which are not independent. Usually, algorithms are more precise than the corresponding mathematical definitions, but this needs not to be the case. Anyway, in the sequel, we shall put more emphasis on the mathematical aspects of the constructions. Our aim is to build interacting explorations that will reveal the presence of large interfaces. These interactive explorations will communicate by sending to each other the territory they have explored. From the point of view of one explorer, the territory explored by the other explorers will be treated as a taboo set. So, in the next stage of our program, we shall add a layer of complexity to the previous shell exploration algorithm by introducing a taboo set $\mathcal{T}$. This taboo set plays the role of a forbidden zone, and we will design an exploration process which never explores the open clusters of $\mathcal{T}$.

**Taboo exploration.** The taboo exploration algorithm takes as input two subsets $A$, $\mathcal{T}$ of $D$, and it operates as the shell exploration algorithm, except that the open sites of the taboo set $\mathcal{T}$ are ignored. Notice here that a taboo set is not merely a forbidden set. Indeed, the taboo exploration will still look at the sites on the boundary of the taboo set. Typically, these sites will be closed, because they belong to the outer boundary of an open cluster which has been explored by another explorer. However, the mechanism driving the taboo exploration will prevent an explorer from going beyond this outer boundary, so that he will not visit open sites that are inside the taboo set.

As the reader might have guessed from the previous brief description, the design of the interacting explorations is quite delicate. The fine details of the intertwined explorations are different for the bond and the site models, and they are simpler for the former than for the latter. So we will first carry out this program in the bond model, and we will later on come back to the site model.

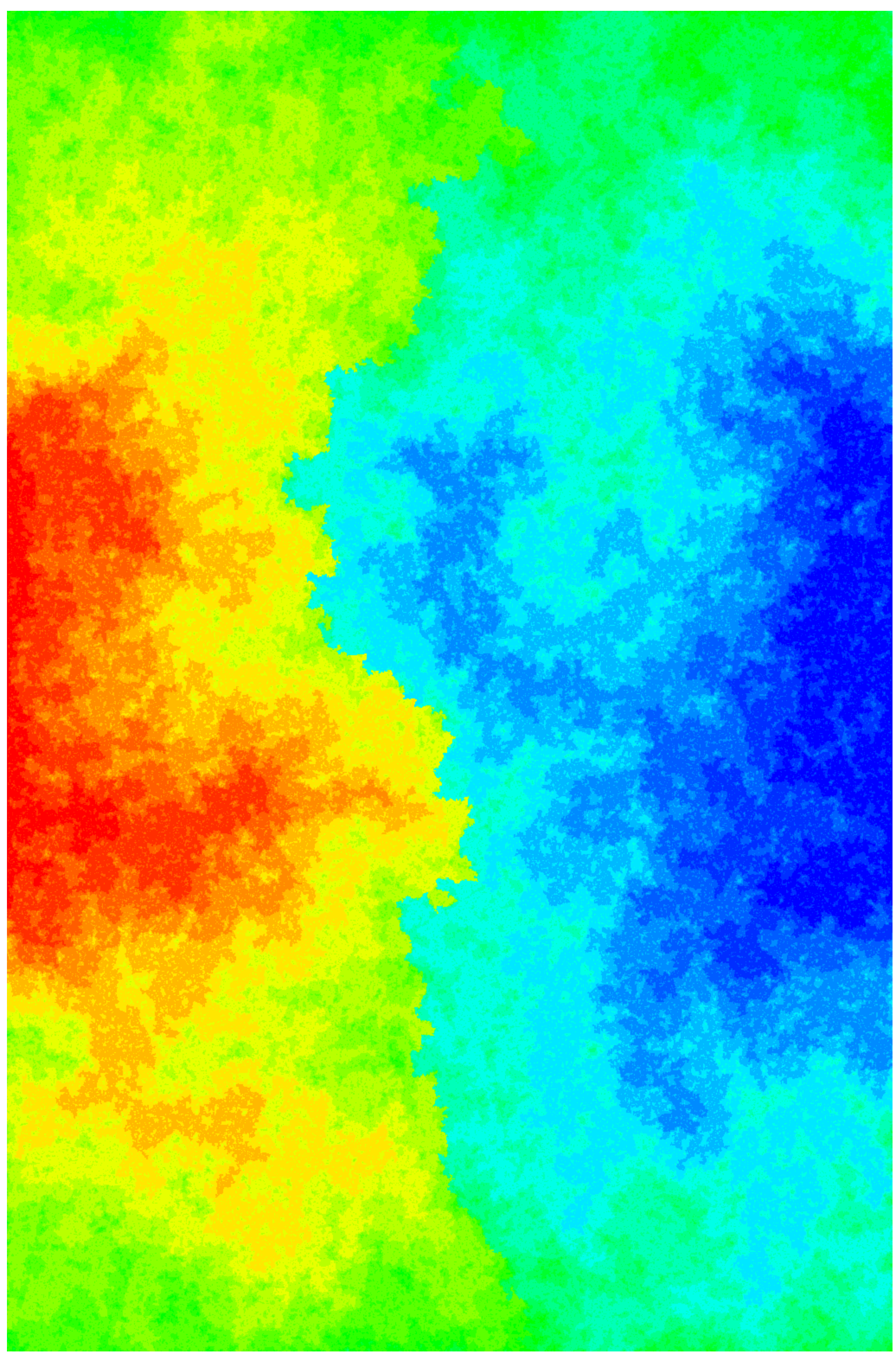

Figure 22: The heart of the matter: disconnection between the middle thirds of left and right sides in the rectangle $1550 \times 1024$, bond percolation, $p = 0.485$.

## Part IV

# The travel time in a finite box

We try here to obtain estimates for the travel times between two sites inside a finite box. The case of sites which are close to the boundary is handled separately in section 18. For sites which are deep inside the box, we present in section 19 a method which works only in dimension three! To handle the dimensions $d \geq 3$, we have to resort to a multiscale construction, explained in section 25.2. The basic step involved on each scale is presented in section 21, it is a variation on a technique borrowed from the Grimmett-Marstrand construction. The synthesis of the geometric construction and the probabilistic estimates are the object of sections 22, 23, 24, 25.1 and 25.2. For a long time, we thought that these results would play an important role to prove theorem 9.2, but ultimately it turned out that this was not the case. Nevertheless, they are interesting in their own right and may be useful for other purposes.

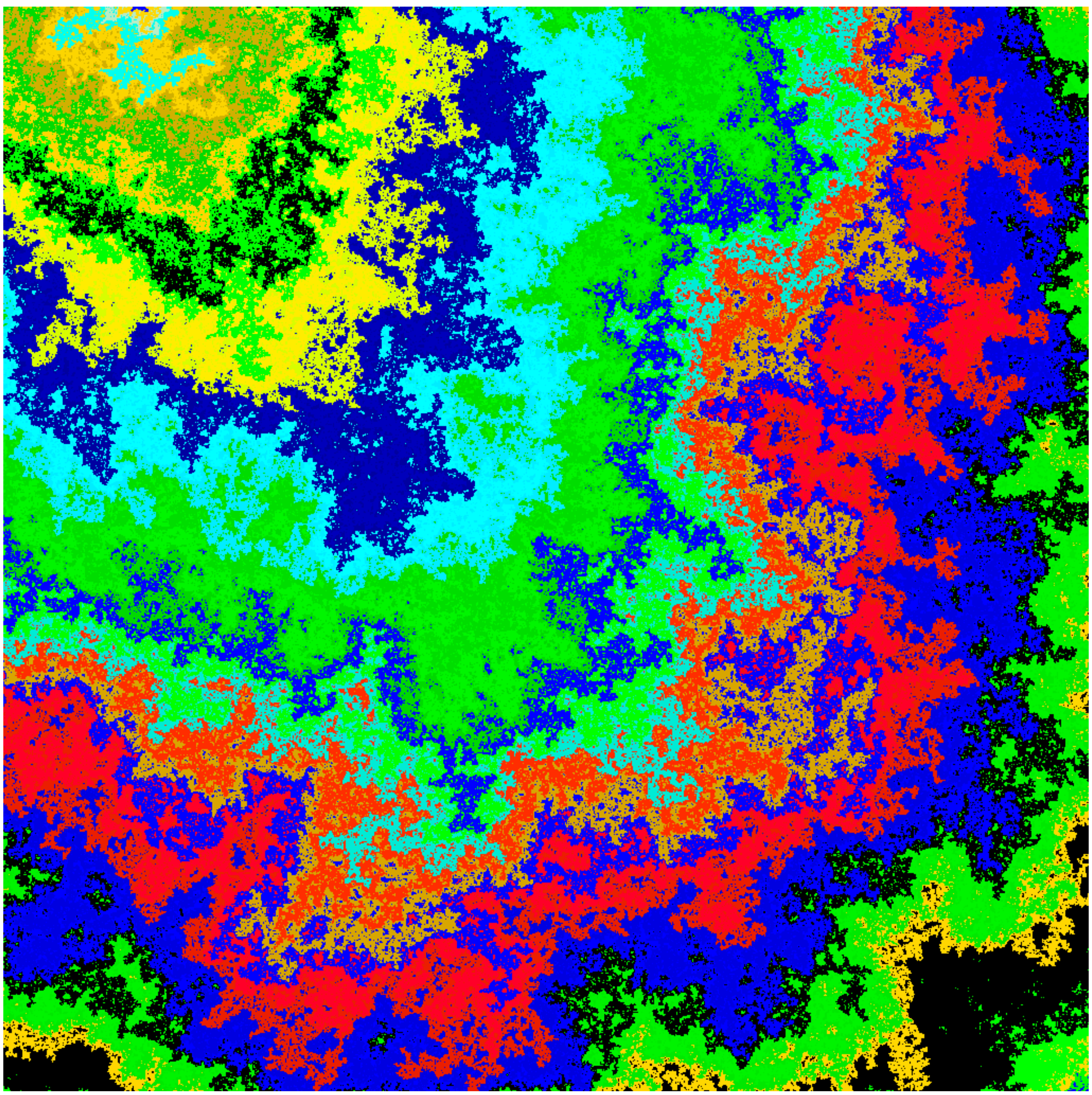

Figure 23: Bond percolation, p=0.47, shells starting at $(-384, 512)$ in $\Lambda(1024)$.

# 18 Stepping away from the boundary

The results presented here will be used later to estimate the travel time to go from the boundary of a box $\Lambda(n)$ to the inner half-box $\Lambda(n/2)$. The technique consists in building an adequate path, along which the travel time can be controlled. This path is built iteratively. For the elementary step of the construction, we associate to each site $y$ of $\Lambda(n)\setminus\Lambda(n/2)$ a target set which gets closer to the box $\Lambda(n/2)$. We must move away from the boundary of $\Lambda(n)$ at an adequate speed, and later on, we will need to compute a quantitative estimate on the travel time between a site $y$ and its target set. To fulfill these goals, we have to choose carefully target$(y)$. This can be done thanks to the symmetries of the lattice and the translation invariance of the model. We introduce first some additional geometric notation before defining target$(y)$.

**Scalar product, distance.** We denote by $\cdot$ the usual scalar product in $\mathbb{R}^d$. Let $u_1,\dots,u_d$ be the vectors of the canonical basis of $\mathbb{R}^d$, i.e.,

$$u_1=(1,0,\dots,0)\,,u_2=(0,1,0,\dots,0)\,,\dots\,,u_d=(0,\dots,0,1)\,.$$

It will be convenient to work with the distance $d_2$ associated to the Euclidean norm $|\cdot|_2$. For $x\in\mathbb{R}^d$ and a subset $A$ of $\mathbb{R}^d$, the distance $d_2(x,A)$ between $x$ and $A$ is defined as

$$d_2(x,A)\,=\,\inf\big\{\,|x-a|_2:a\in A\,\big\}\,.$$

The distance between two subsets $A,B$ of $\mathbb{R}^d$ is defined as

$$d_2(A,B)\,=\,\inf\big\{\,|a-b|_2:a\in A,\,b\in B\,\big\}\,.$$

**The faces of a box $\Lambda$.** Let $\Lambda$ be a box centered at 0 and let $r$ be its side length. We denote by $\partial\Lambda$ the topological boundary of $\Lambda$, defined as

$$\partial\Lambda\,=\,\big\{\,x\in\mathbb{R}^d:\exists i\in\{\,1,\dots,d\,\}\quad x\cdot u_i=\pm r/2\,\big\}\,.$$

Let $v$ be a vector in the list $-u_1,\dots,-u_d,u_1,\dots,u_d$. We define the face $F(\Lambda,v)$ associated to $v$ as

$$F(\Lambda,v)\,=\,\big\{\,x\in\Lambda:x\cdot v=r/2\,\big\}\,.\tag{18.1}$$

Furthermore, we define, for $\varepsilon=(\varepsilon_1,\cdots,\varepsilon_d)$ in $\{-1,+1\}^d$,

$$F(\Lambda,v,\varepsilon)\,=\,\big\{\,x\in F(\Lambda,v):\forall j\in\{\,1,\dots,d\,\}\setminus\{\,i\,\}\quad\varepsilon_j(x\cdot u_j)\geq 0\,\big\}\,.\tag{18.2}$$

We denote by $\mathcal{F}(\Lambda)$ the collection of these $d2^d$ subsets of $\partial\Lambda$, i.e,

$$\mathcal{F}(\Lambda)\,=\,\big\{\,F(\Lambda,v,\varepsilon):v=\pm u_i,1\leq i\leq d,\varepsilon\in\{-1,+1\}^d\,\big\}\,.\tag{18.3}$$

We extend these definitions to boxes $\Lambda$ centered at any point of $\mathbb{R}^d$ by setting

$$\forall x\in\mathbb{R}^d\quad\forall r>0\qquad\mathcal{F}\big(x+\Lambda(r)\big)\,=\,x+\mathcal{F}\big(\Lambda(r)\big)\,.\tag{18.4}$$

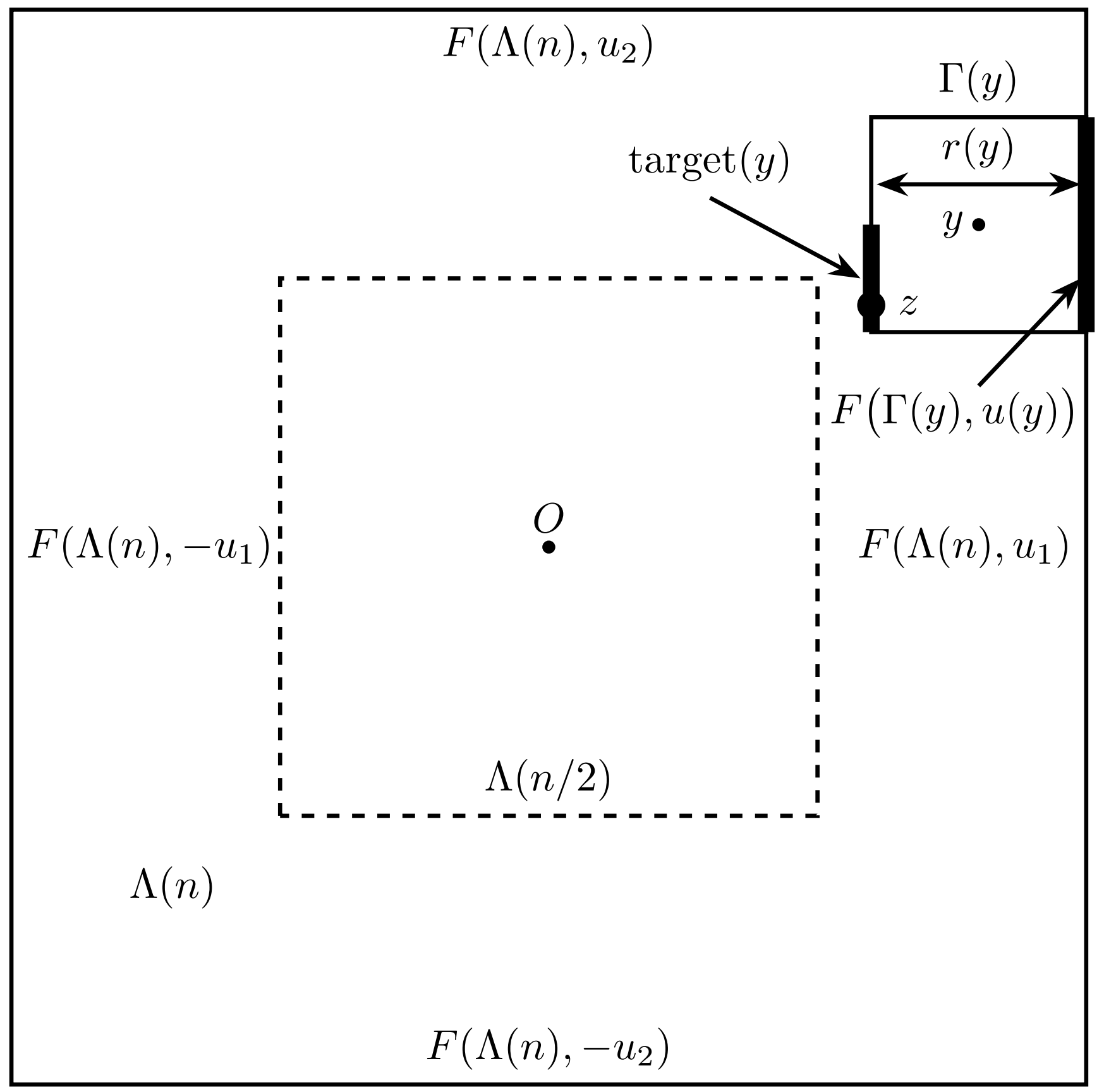

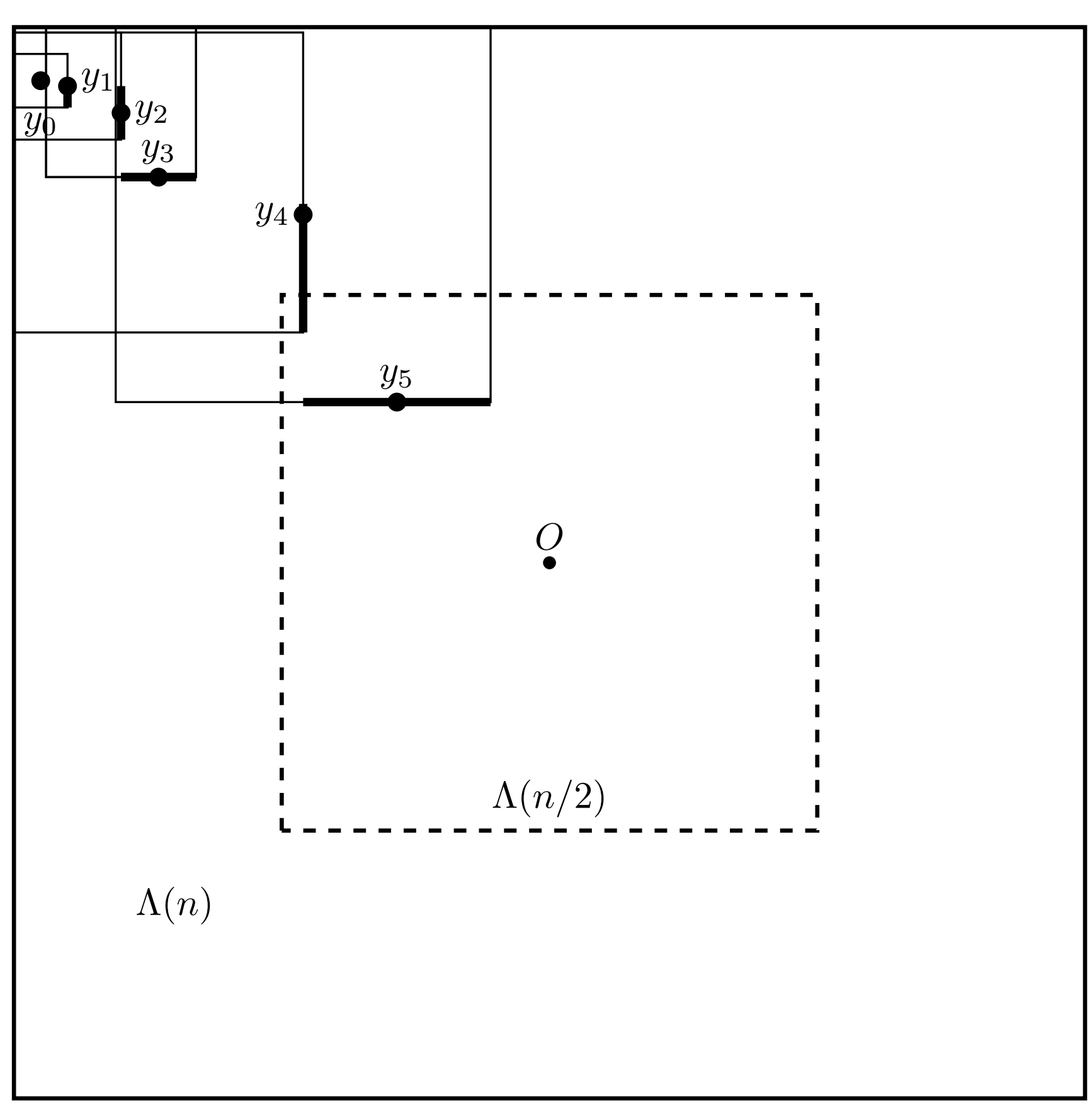

Figure 24: Top: the faces $F\big(\Gamma(y), u(y)\big)$, $F(\Lambda(n), \pm u_i), i = 1, 2$ and target$(y)$.
Bottom: an example of a sequence $(y_k, 0 \leq k \leq K)$.

**The set target$(y)$.** We proceed next to the precise definition of the set target$(y)$. Let us fix $y\in\Lambda(n)\setminus\big(\Lambda(n/2)\cup\partial^{\,in}\Lambda(n)\big)$. Let $\Gamma(y)$ be the largest cubic box centered at $y$ and included in $\Lambda$, i.e., $\Gamma(y)=y+\Lambda\big(r(y)\big)$ where

$$r(y)\,=\,\max\big\{\,r\geq 1: y+\Lambda(r)\subset\Lambda(n)\,\big\}\,=\,2d_2\big(y,\partial\Lambda(n)\big)\,.$$

Since $y$ is not in $\Lambda(n/2)$, then $r(y)<n/2$. We denote also by $u(y)$ the first vector in the list

$$-u_1,\dots,-u_d,u_1,\dots,u_d$$

which satisfies

$$d_2\Big(y,F\big(\Lambda(n),u(y)\big)\Big)\,=\,\frac{r(y)}{2}\,=\,d_2\big(y,\partial\Lambda(n)\big)\,.$$

We have therefore

$$\begin{aligned}d_2\Big(F\big(\Gamma(y),-u(y)\big),F\big(\Gamma(y),u(y)\big)\Big)\,&=\,d_2\Big(F\big(\Gamma(y),-u(y)\big),F\big(\Lambda(n),u(y)\big)\Big)\\&=\,r(y)\,=\,2d_2\Big(y,F\big(\Lambda(n),u(y)\big)\Big)\,.\end{aligned}$$

We set

$$\forall\ell\in\{\,1,\dots,d\,\}\qquad e_\ell(y)\,=\,\begin{cases}-1&\text{if } d_2\Big(y,F\big(\Lambda(n),u_\ell\big)\Big)\,\leq\,\dfrac{n}{2}\,,\\[2ex]+1&\text{if } d_2\Big(y,F\big(\Lambda(n),u_\ell\big)\Big)\,>\,\dfrac{n}{2}\,,\end{cases}$$

and

$$e(y)\,=\,\big(e_1(y),\dots,e_d(y)\big)\in\{-1,+1\}^d\,.$$

The target set associated to $y$ is finally defined as

$$\text{target}(y)\,=\,F\big(\Gamma(y),-u(y),e(y)\big)\,.$$

**Geometric inequalities**. When moving from $y$ to target$(y)$, we get closer to the box $\Lambda(n/2)$. We provide next several inequalities controlling the distance between an arbitrary site of target$(y)$ and $\Lambda(n/2)$. Let $y,z$ be two sites such that

$$y\in\Lambda(n)\setminus\big(\Lambda(n/2)\cup\partial^{\,in}\Lambda(n)\big)\,,\quad z\in\text{target}(y)\,.$$

Since $z$ belongs to the face $F\big(\Gamma(y),-u(y)\big)$, then

$$d_2\Big(z,F\big(\Lambda(n),u(y)\big)\Big)\,\geq\,2d_2\Big(y,F\big(\Lambda(n),u(y)\big)\Big)\,.$$

We compute next the distance between $z$ and the other faces of $\Lambda(n)$. Let $i$ be the unique index such that $u(y)=\pm u_i$ and let $\ell\in\{\,1,\dots,d\,\}\setminus\{\,i\,\}$. We have

$$d_2\Big(z,F\big(\Lambda(n),-u_\ell\big)\Big)+d_2\Big(z,F\big(\Lambda(n),u_\ell\big)\Big)\,=\,n\,.$$

We distinguish two cases, according to the value of $e_\ell(y)$:

- $e_\ell(y) = -1$. In this case, we have

$$d_2\Big(z, F\big(\Lambda(n), u_\ell\big)\Big) \;\geq\; d_2\Big(y, F\big(\Lambda(n), u_\ell\big)\Big)\,,$$
$$d_2\Big(z, F\big(\Lambda(n), -u_\ell\big)\Big) \;\geq\; d_2\Big(y, F\big(\Lambda(n), -u_\ell\big)\Big) - \frac{r(y)}{2} \;>\; \frac{n}{2} - \frac{n}{4} \;=\; \frac{n}{4}\,.$$

- $e_\ell(y) = +1$. In this case, we have

$$d_2\Big(z, F\big(\Lambda(n), -u_\ell\big)\Big) \;\geq\; d_2\Big(y, F\big(\Lambda(n), -u_\ell\big)\Big)\,,$$
$$d_2\Big(z, F\big(\Lambda(n), u_\ell\big)\Big) \;\geq\; d_2\Big(y, F\big(\Lambda(n), u_\ell\big)\Big) - \frac{r(y)}{2} \;>\; \frac{n}{2} - \frac{n}{4} \;=\; \frac{n}{4}\,.$$

In the next proposition, we sum up what we have learnt from the previous inequalities.

**Proposition 18.1.** *Let $n \geq 1$ and let $y$ be a site in $\Lambda(n) \setminus \big(\Lambda(n/2) \cup \partial^{\,in}\Lambda(n)\big)$. For any $z$ in target$(y)$, there is a vector $u(y)$ in the list $-u_1, \dots, -u_d, u_1, \dots, u_d$ such that*

$$d_2\Big(z, F\big(\Lambda(n), u(y)\big)\Big) \;\geq\; 2d_2\Big(y, F\big(\Lambda(n), u(y)\big)\Big) \;=\; 2d_2\big(y, \partial\Lambda(n)\big)\,,$$

*and moreover*

$$\forall \ell \in \{\,1, \dots, d\,\} \quad \forall \varepsilon \in \{-1, +1\}$$
$$d_2\Big(z, F\big(\Lambda(n), \varepsilon u_\ell\big)\Big) \;\geq\; \min\left(d_2\Big(y, F\big(\Lambda(n), \varepsilon u_\ell\big)\Big), \frac{n}{4}\right).$$

The point is that, when we move from a site $y$ to a site belonging to target$(y)$, the distance to the face of $\Lambda(n)$ which is closest to $y$ is doubled (or one such face if there are several), while the distance to the other faces does not decrease, unless it is larger than or equal to $n/4$. By iterating the inequalities of proposition 18.1, we obtain the following corollary.

**Corollary 18.2.** *Let $n, K \geq 1$ and let $(y_k, 0 \leq k \leq K)$ be a sequence of sites in $\Lambda(n)$ such that*

$$\begin{aligned} &y_0 \in \Lambda(n) \setminus \big(\Lambda(n/2) \cup \partial^{\,in}\Lambda(n)\big)\,,\\ \forall k \in \{\,1, \dots, K\,\} \qquad &y_k \in \mathit{target}(y_{k-1})\,. \end{aligned} \tag{18.5}$$

*We have*

$$\forall k \in \{\,0, \dots, K\,\} \qquad d_2\big(y_k, \partial\Lambda(n)\big) \;\geq\; \min\left(2^{\lfloor \frac{k}{d} \rfloor - 1}, \frac{n}{4}\right). \tag{18.6}$$

*Proof.* A straightforward induction yields the inequality (18.6). Indeed, as long as the sequence is outside of the box $\Lambda(n/2)$, the distance to one of the closest face of $\Lambda(n)$ doubles at each step, while the distances to the other faces do not decrease or stay larger than or equal to $n/4$. In particular, the distance to $\partial\Lambda(n)$ doubles every $d$ steps, until the sequence enters the box $\Lambda(n/2)$. □

A simple consequence of corollary 18.2 is that a sequence satisfying (18.5) enters the box $\Lambda(n/2)$ after at most $d(\ln n/\ln 2) + 1$ steps.

# 19 Travelling in the box

Thanks to the constructions performed in the previous section, we will be able to travel from the boundary of $\Lambda(n)$ until we reach the box $\Lambda(n/2)$. Our next goal is to travel between two arbitrary sites $x, y$ belonging to $\Lambda(n/2)$. We will rely on the same technique as in the previous section. However, the problem is more complicated because the ultimate target is much smaller. So we need to make the most of the available symmetries, and we describe in the next subsection the relevant group of transformations.

## 19.1 The hyperoctahedral group

The hyperoctahedral group is the group of the symmetries of a $d$-dimensional hypercube centered at the origin. We consider the $d$ hyperplanes of $\mathbb{R}^d$ having for equations

$$x_i = 0\,, \quad 1 \leq i \leq d\,.$$

These $d$ hyperplanes tile $\mathbb{R}^d$ into $2^d$ hyperoctants. We consider in addition the $d(d-1)$ hyperplanes of $\mathbb{R}^d$ having for equations

$$\begin{aligned} x_i &= x_j\,, \quad 1 \leq i < j \leq d\,,\\ x_i &= -x_j\,, \quad 1 \leq i < j \leq d\,. \end{aligned}$$

We denote by $\mathcal{H}$ the collection of all these hyperplanes. In total, the collection $\mathcal{H}$ contains $d + d(d-1) = d^2$ hyperplanes. The collection $\mathcal{H}$ tiles the space $\mathbb{R}^d$ into a finite number of closed regions. We denote by $\mathcal{R}$ the collection of these regions. Since all the hyperplanes of $\mathcal{H}$ contain the origin, the regions in the collection $\mathcal{R}$ are positive cones emanating from the origin. These cones are delimited by lines in $\mathbb{R}^d$, which are all the lines obtained by intersecting some of the hyperplanes of $\mathcal{H}$. The equations of these lines are of the form:

$$\begin{gathered} x_{i_1} = x_{i_2} = \cdots = x_{i_k} = 0\,, \quad 1 \leq i_1 < \cdots < i_k \leq d\,,\\ \varepsilon_{k+1} x_{i_{k+1}} = \cdots = \varepsilon_{d-1} x_{i_{d-1}} = x_{i_d}\,, \quad 1 \leq i_{k+1} < \cdots < i_d \leq d\,,\\ \varepsilon_{k+1}, \dots, \varepsilon_{d-1} \in \{\,-1, +1\,\}\,,\\ 0 \leq k \leq d-1\,, \quad \{\,1, \dots, d\,\} \,=\, \{\,i_1, \dots, i_k\,\} \cup \{\,i_{k+1}, \dots, i_d\,\}\,. \end{gathered}$$

Let $G(\mathcal{H})$ be the group of the transformations of $\mathbb{R}^d$ generated by the orthogonal symmetries with respect to the hyperplanes of the collection $\mathcal{H}$. Since the cubic lattice $\mathbb{Z}^d$ is globally invariant with respect to these orthogonal symmetries, it is still invariant under the action of $G(\mathcal{H})$. Moreover the group $G(\mathcal{H})$ acts transitively on the collection of regions $\mathcal{R}$:

$$\forall R_1, R_2 \in \mathcal{R} \quad \exists\, g \in G(\mathcal{H}) \quad g(R_1) = R_2\,.$$

What is more, the Bernoulli site percolation model is invariant under the action of $G(\mathcal{H})$. As a consequence, for any subset $D$ of $\mathbb{Z}^d$, for any subsets $A, B$ of $D$, the distribution of the random variable $T_D(A, B)$ is the same as the distribution of the random variable $T_{g(D)}(g(A), g(B))$, for any transformation $g$ in $G(\mathcal{H})$.

## 19.2 The natural attempt with boxes

Let $r \geq 1$. The collection of regions $\mathcal{R}$ induces a tiling of the boundary of the cubic box $\Lambda(r)$ into a finite number of tiles. These tiles are obtained as the intersection between a region $R$ of $\mathcal{R}$ and $\partial\Lambda(r)$. A generic tile will be denoted by the letter $T$, the collection $\mathcal{T}$ of all the tiles is given by

$$\mathcal{T} \,=\, \big\{\, \partial\Lambda(r) \cap R : R \in \mathcal{R} \,\big\}\,.$$

The group $G(\mathcal{H})$ acts transitively on the collection $\mathcal{T}$ of the tiles, any tile can be transformed into any tile by an isometry of $G(\mathcal{H})$. Let us focus on the tile $T_0$ associated to the region $R_0$ delimited by the lines of equations

$$x_1 = \cdots = x_k = 0\,, \quad x_{k+1} = \cdots = x_d\,, \quad 0 \leq k \leq d-1\,.$$

These lines intersect the boundary $\partial\Lambda(r)$ at the following points:

$$s_1 \,=\, \Big(0,\dots,0,\frac{r}{2}\Big)\,, \quad s_2 \,=\, \Big(0,\dots,0,\frac{r}{2},\frac{r}{2}\Big)\,,$$
$$s_3 \,=\, \Big(0,\dots,\frac{r}{2},\frac{r}{2},\frac{r}{2}\Big)\,, \quad \dots\,, \quad s_d \,=\, \Big(\frac{r}{2},\dots,\frac{r}{2}\Big)\,.$$

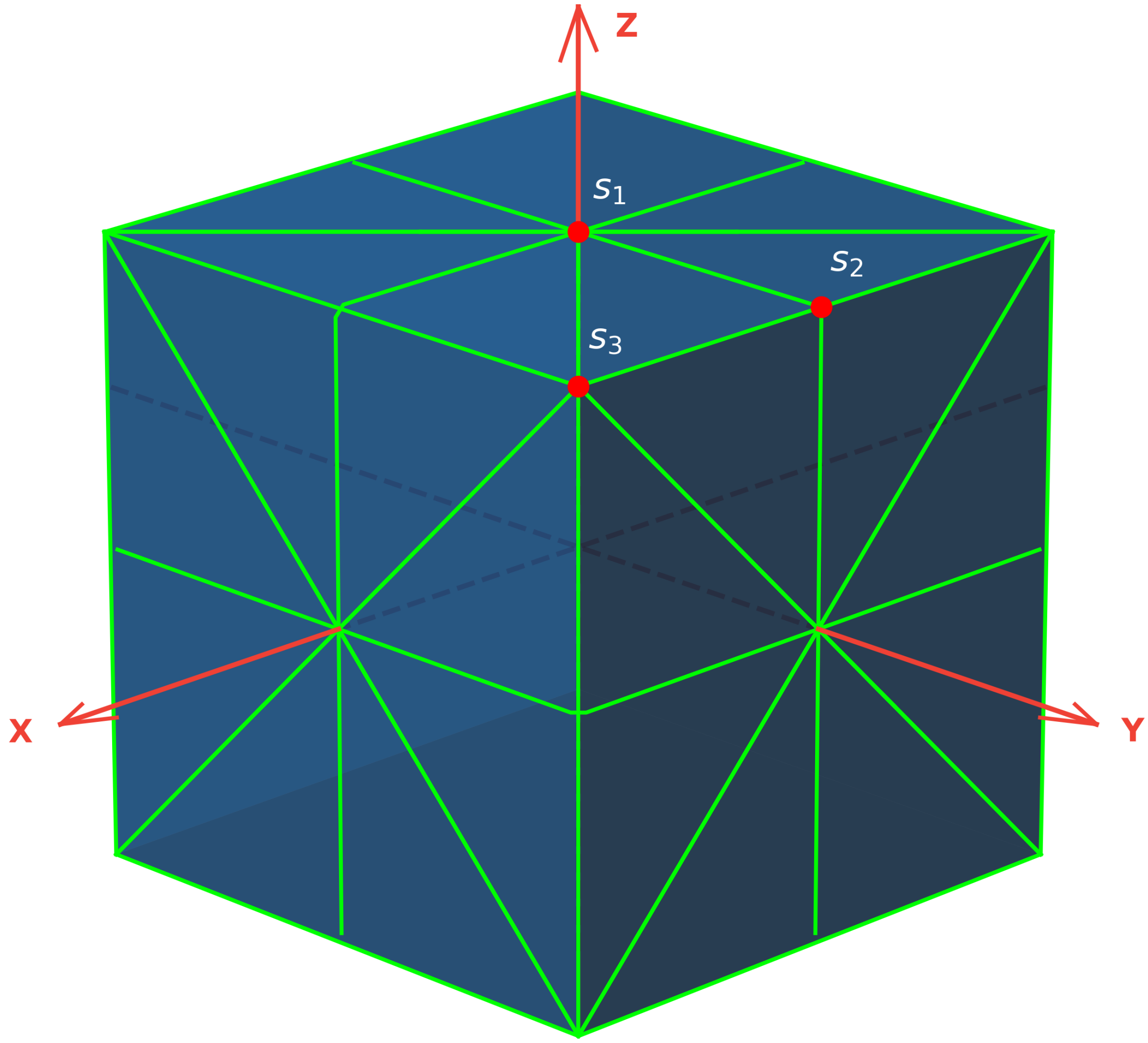

Figure 25: The tiling of the 3D cube

Let $n \geq 1$ be fixed. Let $x, y$ be two arbitrary sites belonging to $\Lambda(n/2)$. The site $x$ is the starting site and the site $y$ is the target. Our goal is to control the travel time between $x$ and $y$ inside the box $\Lambda(n)$. A natural attempt is to consider the smallest box $\Gamma(x,y)$ centered at $x$ containing $y$ in its boundary, i.e.,

$$\Gamma(x,y) \,=\, x + \Lambda\big(2|x-y|_\infty\big)\,.$$

We want the box $\Gamma(x,y)$ to be included in $\Lambda(n)$, so we add the additional constraint that $|x-y|_\infty \leq n/4$. We would then define the set $\mathrm{target}(x,y)$ as the tile of the boundary of $\Gamma(x,y)$ which contains $y$. All the sites in this tile have one of their $d$ coordinates which coincides with the one of $y$. Thus we hope that, after travelling from $x$ to this tile , the distance to the target site $y$ is seriously reduced. Unfortunately, all the tiles of $\partial\Gamma(x,y)$ are isometric right-angled isosceles triangles, whose diameter is larger than or equal to $|x-y|_\infty$, whether using the Euclidean or the supremum norm. In the end, we cannot improve the distance to the target with such a move. The reason is that the tiles of the boundary of the cube have a too large diameter. This suggests that the cube is not the right geometric shape to work with. In the next section, we try to implement this strategy with an Euclidean ball instead of a cube.

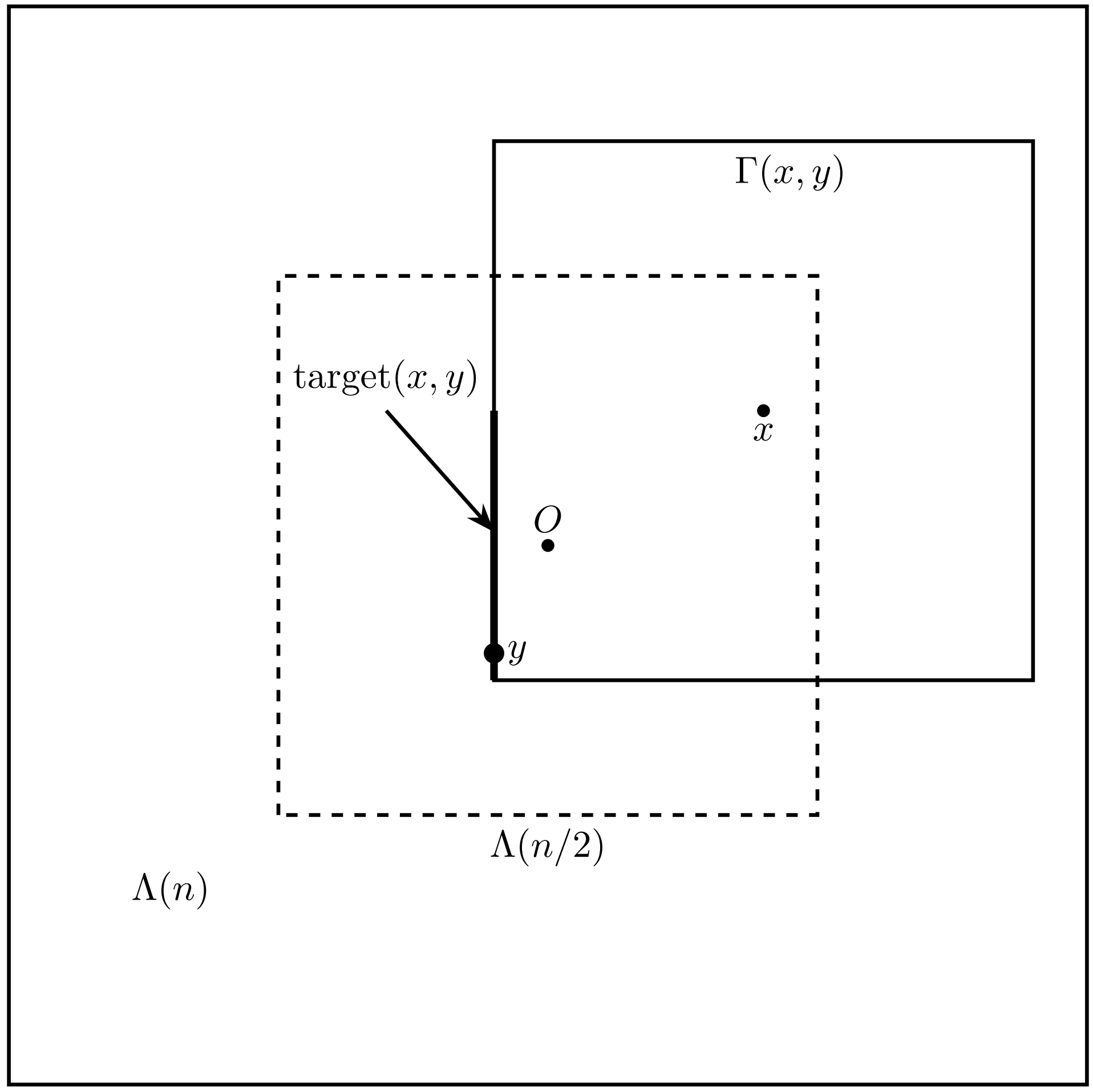

Figure 26: An attempt for $\mathrm{target}(x,y)$, $x, y \in \Lambda(n/2)$ such that $|x-y|_\infty \leq n/4$

## 19.3 The tiling of the sphere

We denote by $S^{d-1}$ the $(d-1)$-dimensional unit sphere of $\mathbb{R}^d$. The collection of regions $\mathcal{R}$ induces a tiling of the unit sphere $S^{d-1}$ into a finite number of tiles. These tiles are obtained as the intersection between a region $R$ of $\mathcal{R}$ and the unit sphere $S^{d-1}$. A generic tile will be denoted by the letter $T$, the collection $\mathcal{T}$ of all the tiles is given by

$$\mathcal{T} = \left\{ S^{d-1} \cap R : R \in \mathcal{R} \right\}. \tag{19.1}$$

Once more, the group $G(\mathcal{H})$ acts transitively on the collection $\mathcal{T}$ of the tiles. Let us look more carefully at the geometry of these tiles. Since any tile can be transformed into any tile by an isometry of $G(\mathcal{H})$, we can focus ourselves on the tile $T_0$ associated to the region $R_0$ delimited by the lines of equations

$$x_1 = \cdots = x_k = 0\,, \quad x_{k+1} = \cdots = x_d\,, \quad 0 \leq k \leq d-1\,.$$

These lines intersect the sphere $S^{d-1}$ at the following points:

$$s_1 = (0, \dots, 0, 1)\,, \quad s_2 = \Big(0, \dots, 0, \frac{\sqrt{2}}{2}, \frac{\sqrt{2}}{2}\Big)\,,$$
$$s_3 = \Big(0, \dots, 0, \frac{\sqrt{3}}{3}, \frac{\sqrt{3}}{3}, \frac{\sqrt{3}}{3}\Big)\,, \quad \dots\,, \quad s_d = \Big(\frac{\sqrt{d}}{d}, \dots, \frac{\sqrt{d}}{d}\Big)\,.$$

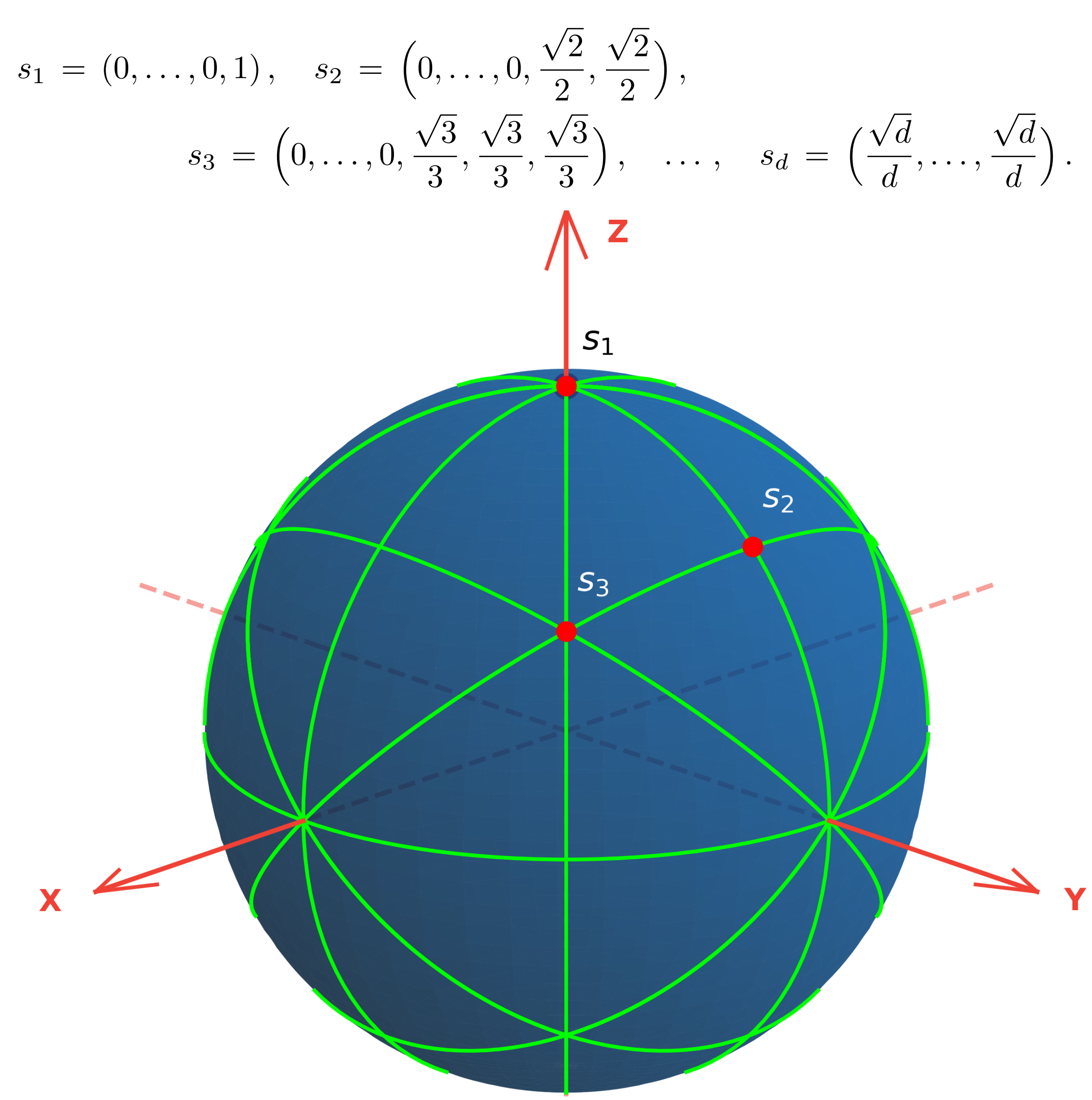

Figure 27: The tiling of the 3D sphere

These $d$ points delimit a tile $T_0$ on $S^{d-1}$. The other tiles of $\mathcal{T}$ are obtained as the orbit of $T_0$ under the action of $G(\mathcal{H})$. In fact, the points of the sphere $S^{d-1}$ which are obtained from the $d$ points above by permuting two coordinates or adding a minus sign constitute the vertices of a convex polytope. In dimension 2, this convex polytope is a regular octogon. In dimension 3, this convex polytope is a polyhedron with 48 triangular faces, it is a Catalan solid called the disdyakis dodecahedron or the hexakis octahedron [71, 137].

**Lemma 19.1.** *The diameter of the tiles of $\mathcal{T}$ is strictly less than* 0.92 *in dimension $d = 3$, it is equal to* 1 *in dimension $d = 4$ and strictly larger than* 1 *in dimensions $d \geq 5$.*

*Proof.* Let $x, y$ be two points of $T_0$. We denote by $\angle(x, y)$ the angle between the two vectors $\overrightarrow{Ox}$, $\overrightarrow{Oy}$. The Euclidean distance $|x - y|_2$ is equal to

$$|x - y|_2 \;=\; \sqrt{2 - 2\cos\angle(x, y)}\,.$$

Taking the supremum over $x, y$ in $T_0$, we conclude that

$$\text{diameter}\, T_0 \;=\; \sup\,\Big\{\,\sqrt{2 - 2\cos\angle(x, y)} : x, y \in T_0\,\Big\}\,. \tag{19.2}$$

The tile $T_0$ is included in the positive hyperoctant, thus the angles $\angle(x, y)$ are all in the interval $[0, \pi/2]$ and the corresponding cosines are non-negative. The supremum in (19.2) can be attained only when $x, y$ are on the lines delimiting the tile $T_0$, therefore

$$\begin{aligned} \text{diameter}\, T_0 \;&=\; \max\,\Big\{\,\sqrt{2 - 2\cos\angle(s_i, s_j)} : 1 \leq i < j \leq d\,\Big\} \\ &=\; \max\,\Big\{\,|s_i - s_j|_2 : 1 \leq i < j \leq d\,\Big\}\,. \end{aligned} \tag{19.3}$$

The maximum in (19.3) is reached when the cosine is minimal, or when the distance between the vertices is maximal. We compute, for $1 \leq i < j \leq d$,

$$|s_i - s_j|_2 \;=\; \sqrt{i\Big(\frac{\sqrt{i}}{i} - \frac{\sqrt{j}}{j}\Big)^2 + (j - i)\Big(\frac{\sqrt{j}}{j}\Big)^2} \;=\; \sqrt{2 - 2\sqrt{\frac{i}{j}}}\,.$$

The distance is maximal for $i = 1$ and $j = d$, and we conclude that

$$\text{diameter}\, T_0 \;=\; \sqrt{2 - \frac{2}{\sqrt{d}}}\,.$$

In dimension 3, we find that the diameter is strictly less than 0.92. In dimension 4, the diameter is exactly equal to 1. In dimensions 5 and higher, the diameter is strictly larger than 1. ☐

The good news is that, at least in dimension 3, the diameter of a tile is strictly less than 1. So there is hope that the attempt of subsection 19.2, which was based on cubes, might work with balls. We examine this possibility in the next subsection.

## 19.4 The target set in three dimensions

In this subsection, we work only in dimension 3. In order to define the target set, we will rely on adequate dilations of the unit ball. Let $r > 0$. We define a discrete version of the ball of radius $r$ by setting

$$B_r \,=\, \big\{\, (x,y,z) \in \mathbb{Z}^3 : x^2+y^2+z^2 \leq r^2 \,\big\}\,.$$

Let $T$ belong to the collection $\mathcal{T}$ defined in (19.1). We define

$$T_r \,=\, \big\{\, x \in B_r : d(x, rT) \leq 3 \,\big\}\,.$$

Lemma 19.1 yields that

$$\forall y,z \in T_r \qquad |y-z|_2 \,\leq\, 6 + r\,\text{diameter}(T) \,\leq\, 6+0.92r\,. \tag{19.4}$$

The inequality (19.4) readily implies the following result.

**Corollary 19.2.** *For any $r \geq 600$, we have*

$$\forall y,z \in T_r \qquad |y-z|_2 \,\leq\, 0.93\,r\,.$$

Let $n \geq 1$ be fixed. Let $y$ belong to $\Lambda(n/2)$ and let $x$ in $\mathbb{Z}^d$ be such that $|x-y|_2 \leq n/8$. This ensures that $x+\Lambda(n/4) \subset \Lambda(n)$. The site $x$ is the starting site and the site $y$ is the target. We consider the smallest ball $\Gamma(x,y)$ centered at $x$ containing $y$ in its boundary, i.e.,

$$\Gamma(x,y) \,=\, x + B_r\,, \qquad \text{where} \quad r = |x-y|_2\,.$$

We define the set $\text{target}(x,y)$ as follows:

• If $|x-y|_2 \leq 600$, then we set $\text{target}(x,y) = \{\, x\,\}$.

• If $|x-y|_2 > 600$, then we define $\text{target}(x,y)$ to be the tile of the boundary of $\Gamma(x,y)$ which contains $y$. If there are several such tiles, we pick one of them according to some deterministic procedure.

This definition and corollary 19.2 imply that either $\text{target}(x,y) = \{\, x\,\}$ or

$$\forall z \in \text{target}(x,y) \qquad |y-z|_2 \,\leq\, 0.93\,|x-y|_2\,, \tag{19.5}$$

therefore $\text{target}(x,y)$ is also included in $\Lambda(n)$. By iterating inequality (19.5), we obtain the following corollary.

**Corollary 19.3.** *Let $n \geq 1$, let $y$ belong to $\Lambda(n/2)$ and let $x$ in $\mathbb{Z}^d$ be such that $|x-y|_2 \leq n/8$. Let $(y_k, 0 \leq k \leq K)$ be a sequence such that $y_0 = x$ and*

$$\forall k \in \{\,1,\dots,K\,\} \qquad y_k \in \mathit{target}(y_{k-1}, y)\,.$$

*We have*

$$\forall k \in \{\,0,\dots,K\,\} \qquad |y-y_k|_2 \,\leq\, \max\Big(0.93^k|x-y|_2, 600\Big)\,. \tag{19.6}$$

The inequality (19.6) and the definition of $\text{target}(x,y)$ imply that the sequence $(y_k, 0 \leq k \leq K)$ becomes stationary after at most

$$\frac{\ln n - \ln 600}{-\ln 0.93} \,\leq\, 14 \ln n$$

steps. It is then equal to a site $y^*$ such that $|y-y^*|_2 \,\leq\, 600$.

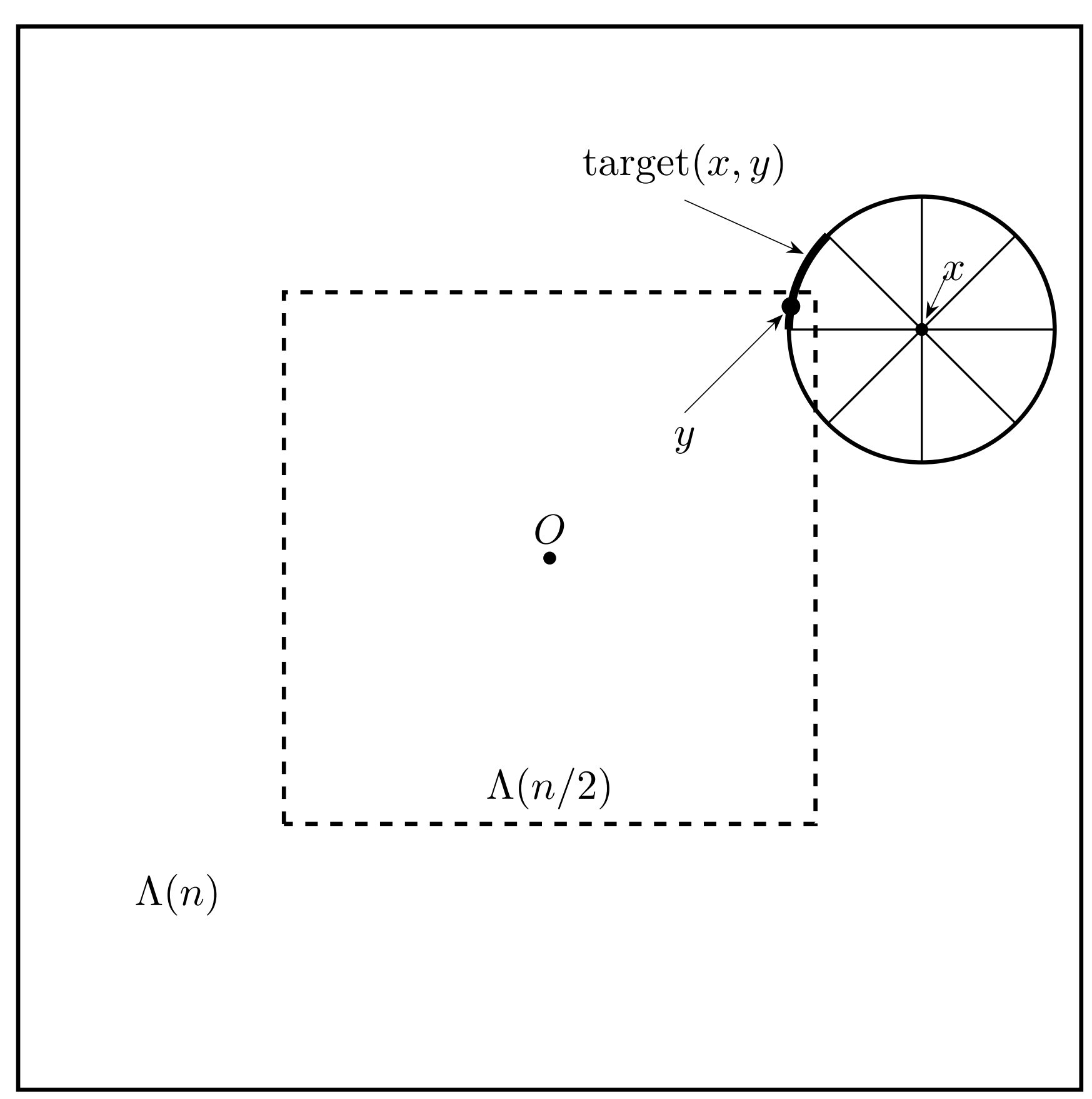

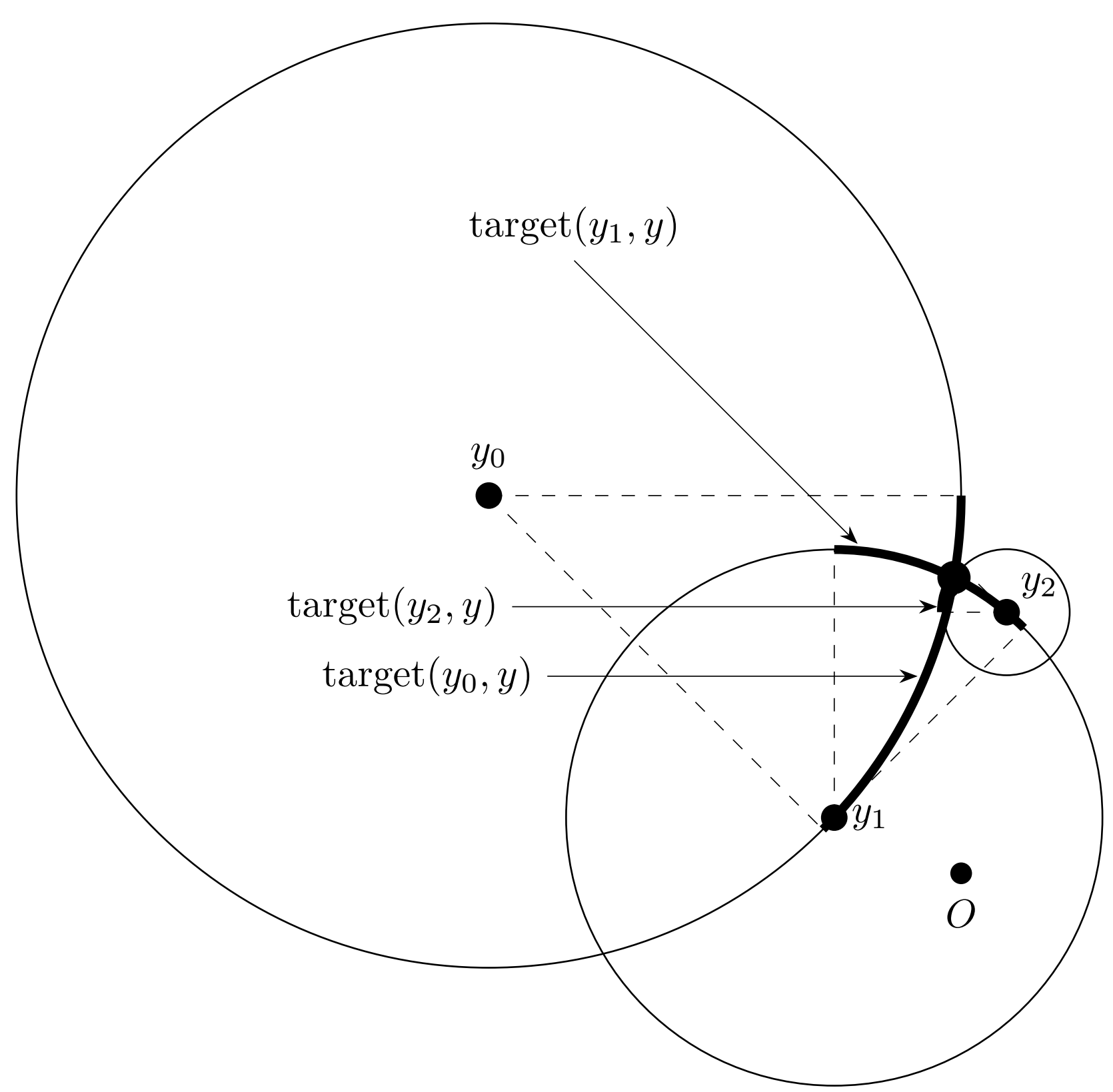

Figure 28: The set target$(x, y)$ and an example of a sequence $(y_k, 0 \leq k \leq K)$

# 20 Travelling towards a box

The method employed in the previous section works only in dimension three. We shall develop a different strategy to deal with the higher dimensions. There is a cost to pay, this strategy is more complicated because it involves different scales. In this section, we present the construction that will be used on each intermediate scale in order to reach a specific orthant. It will also be used to control the travel time between a cubical shell and its inner cube in proposition 25.5.

Throughout the section, we fix a site $x$ in $\mathbb{Z}^d$ and a cubic box $\Lambda^*$ in $\mathbb{R}^d$. The site $x$ is the starting site and the box $\Lambda^*$ is the target. Our goal is to control the travel time between $x$ and $\Lambda^*$. To that end, we will construct iteratively a path between $x$ and $\Lambda^*$ along which the travel times can be controlled. Here we perform the basic step of the construction.

**The set target$(x,\Lambda^*)$.** We define next the set $\text{target}(x,\Lambda^*)$ associated to $x$ and $\Lambda^*$. We denote by $r>0$ the side length of $\Lambda^*$ and by $v^*$ its center, so that

$$\Lambda^* \,=\, v^* + \Lambda(r)\,.$$

If $x$ already belongs to the box $\Lambda^*$, then we set $\text{target}(x,\Lambda^*)=\{\,x\,\}$. Otherwise, let $\Gamma$ be the box

$$\Gamma \,=\, x+\Lambda(r/2)\,.$$

Let $u^*$ be the first vector in the list $-u_1,\dots,-u_d,u_1,\dots,u_d$ which satisfies

$$d_2\big(v^*,F\big(\Gamma,u^*\big)\big) \,=\, d_2\big(v^*,\Gamma\big)\,.$$

Recall that $F\big(\Gamma,u^*\big)$ is the face of $\Gamma$ associated to $u^*$, see (18.1). The vector $u^*$ satisfies also

$$\big|u^*\cdot(v^*-x)\big| \,=\, \max\,\big\{\,\big|u_\ell\cdot(v^*-x)\big|:\ell\in\{\,1,\dots,d\,\}\,\big\}\,,$$

where $\cdot$ is the usual scalar product in $\mathbb{R}^d$. We set

$$\forall \ell\in\{\,1,\dots,d\,\}\qquad e^*_\ell \,=\, \begin{cases}-1 & \text{if } u_\ell\cdot(v^*-x)\leq 0\,,\\ +1 & \text{if } u_\ell\cdot(v^*-x)> 0\,,\end{cases}$$

and $e^*=\big(e^*_1,\dots,e^*_d\big)$. The signs of the components of $e^*$ are chosen so that the segment $[x,v^*]$ intersects the set $F\big(\Gamma,u^*,e^*\big)$ (defined in (18.2)). The target set associated to $x$ and $\Lambda^*$ is finally defined as

$$\text{target}(x,\Lambda^*) \,=\, F\big(\Gamma,u^*,e^*\big)\,. \tag{20.1}$$

**Geometric inequalities**. When moving from $x$ to $\text{target}(x,\Lambda^*)$, we get closer to the center $v^*$ of the box $\Lambda^*$. We derive next an inequality to control how the distance to $v^*$ changes. Suppose that $x$ is not in the box $\Lambda^*$ and let $z$ belong to $\text{target}(x,\Lambda^*)$. For the direction associated to $u^*$, we have

$$\begin{aligned}\big(v^*-z\big)\cdot u^* \,&\geq\, 0\,,\\ \big(v^*-x\big)\cdot u^* \,&=\, \big(v^*-z\big)\cdot u^* \,+\, \frac{r}{4}\,.\end{aligned} \tag{20.2}$$

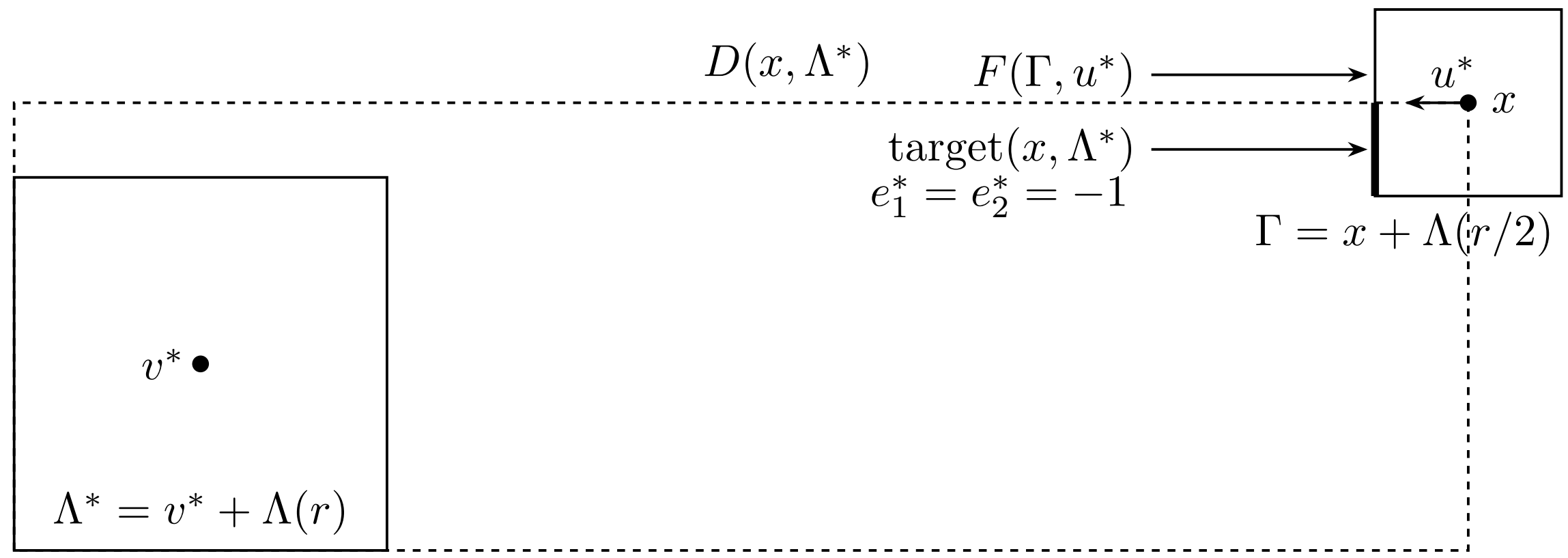

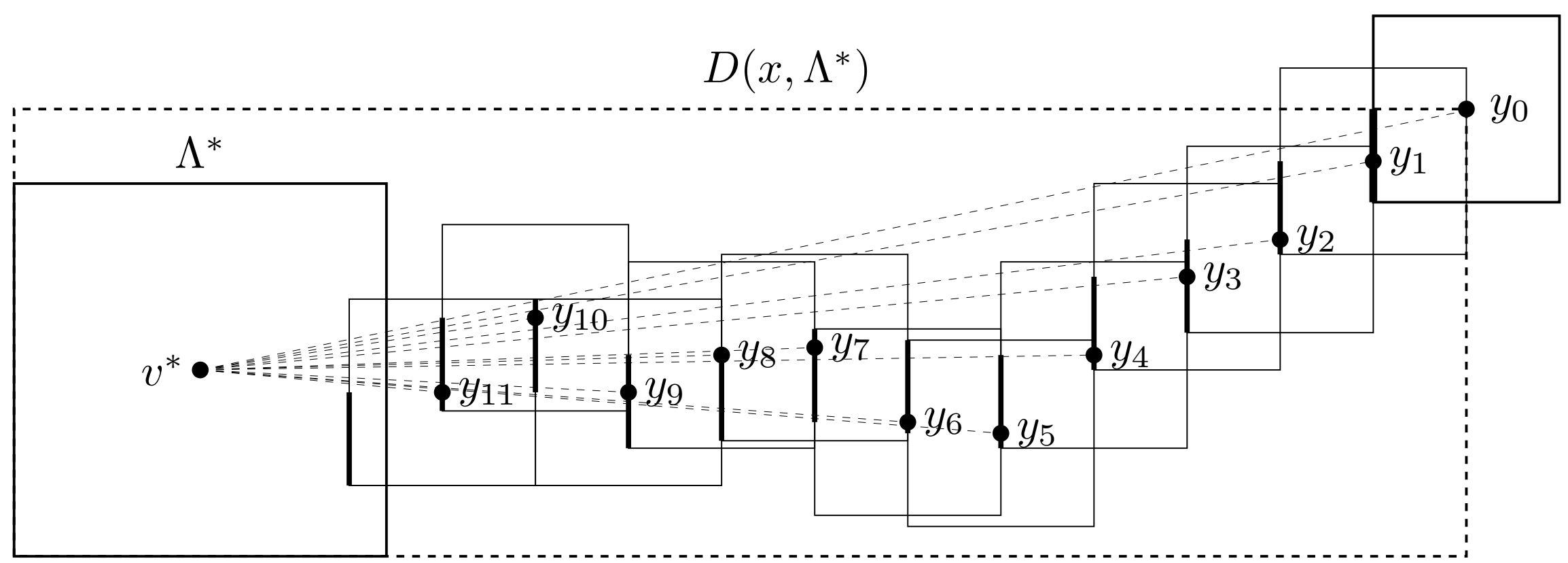

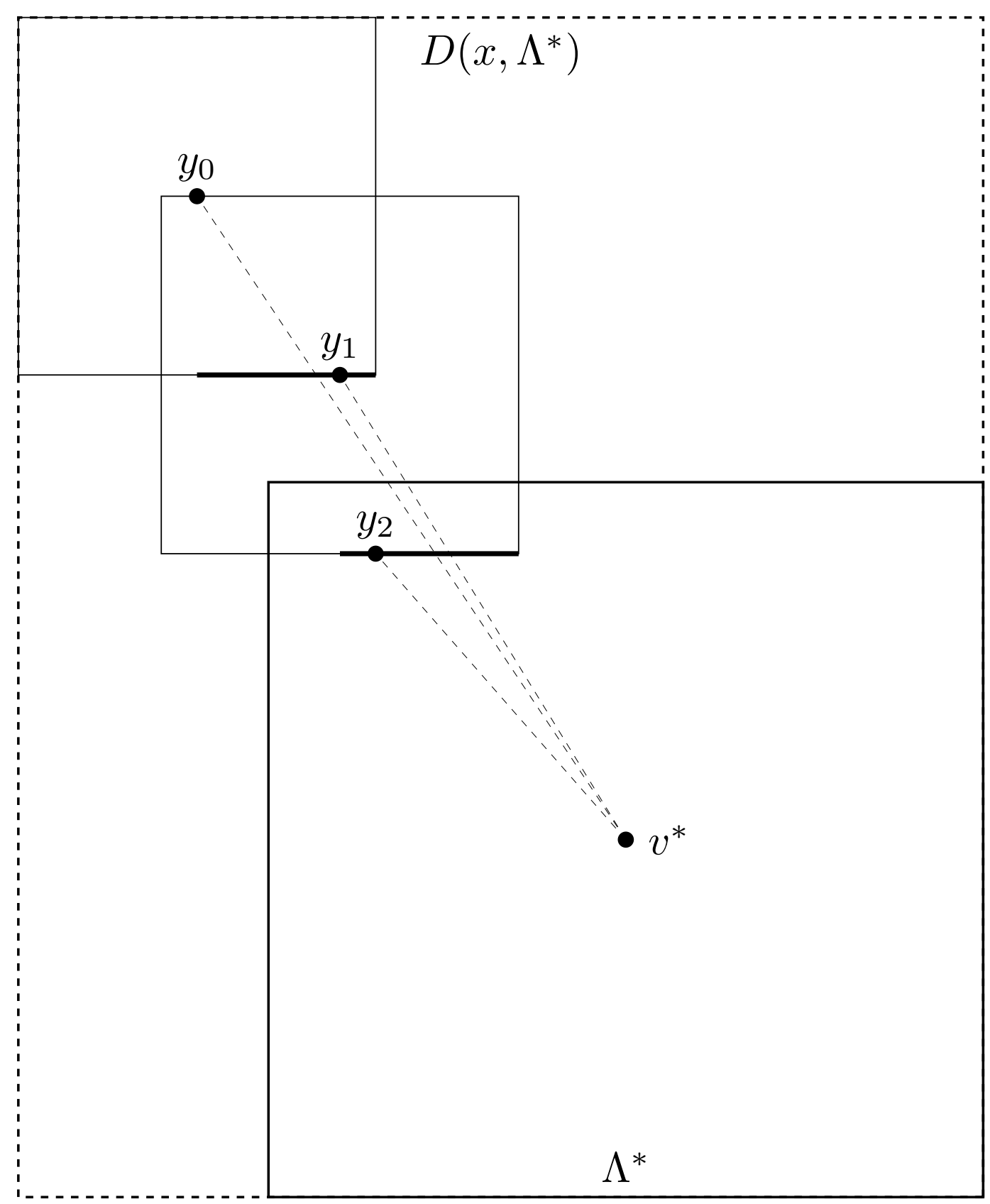

Figure 29: The set target$(x, \Lambda^*)$ and two sequences $(y_k, 0 \leq k \leq K)$.

Let next $\ell$ in $\{1, \ldots, d\}$ be such that $u_\ell \neq u^*$ and $u_\ell \neq -u^*$. We write

$$\big(v^* - x\big) \cdot u_\ell \,=\, \big(v^* - z\big) \cdot u_\ell \,+\, (z - x) \cdot u_\ell \,.$$

It follows from the definition of $e^*$, target$(x, \Lambda^*)$ that $\big(v^* - x\big) \cdot u_\ell$ and $(z - x) \cdot u_\ell$ are of the same sign. Moreover, we have $\big|(z - x) \cdot u_\ell\big| \leq r/4$. Let us discuss each case and state the relevant inequalities:

- $e^*_\ell = -1$. In this case, we have $u_\ell \cdot (v^* - x) \leq 0$ and

$$\big(v^* - x\big) \cdot u_\ell \,\leq\, \big(v^* - z\big) \cdot u_\ell \,\leq\, \big(v^* - x\big) \cdot u_\ell \,+\, \frac{r}{4} \,. \tag{20.3}$$

- $e^*_\ell = +1$. In this case, we have $u_\ell \cdot (v^* - x) > 0$ and

$$\big(v^* - x\big) \cdot u_\ell \,-\, \frac{r}{4} \,\leq\, \big(v^* - z\big) \cdot u_\ell \,\leq\, \big(v^* - x\big) \cdot u_\ell \,. \tag{20.4}$$

Inequalities (20.3) and (20.4) together yield that

$$\Big|\big(v^* - z\big) \cdot u_\ell\Big| \,\leq\, \max\Big(\Big|\big(v^* - x\big) \cdot u_\ell\Big|, \frac{r}{4}\Big) \,.$$

Let $D(x, \Lambda^*)$ be the smallest axis-parallel hyperrectangle containing $x$ and $\Lambda^*$. The inequalities (20.2), (20.3) and (20.4) imply that $z$ belongs to $D(x, \Lambda^*)$, hence target$(x, \Lambda^*)$ is included in $D(x, \Lambda^*)$. Unless $x$ is already in $\Lambda^*$, the box $\Gamma$ is not included in $D(x, \Lambda^*)$, but in any case it is included in the enlargement

$$D(x, \Lambda^*) + \Lambda\Big(\frac{r}{2}\Big) \,=\, \Big\{\, y \in \mathbb{R}^d : \exists\, z \in D(x, \Lambda^*) \quad |y - z|_\infty \leq \frac{r}{2} \,\Big\} \,.$$

The point is that target$(x, \Lambda^*)$ is a tile of the boundary of $\Gamma$ which is included in $D(x, \Lambda^*)$. However, to control the travel time between $x$ and target$(x, \Lambda^*)$, we will consider paths inside the above enlargement. In the next proposition, we sum up what we have learnt from the previous inequalities.

**Proposition 20.1.** *Let $\Lambda^*$ be a cubic box in $\mathbb{R}^d$, and let $x$ be a site in $\mathbb{Z}^d \setminus \Lambda^*$. Let $r > 0$ be the side length of $\Lambda^*$ and let $v^*$ be its center. For any $z$ in target$(x, \Lambda^*)$, there exists $\ell^*$ in $\{1, \ldots, d\}$ such that*

$$\big|u_{\ell^*} \cdot (v^* - x)\big| \,=\, \max\,\Big\{\, \big|u_\ell \cdot (v^* - x)\big| : \ell \in \{1, \ldots, d\}\,\Big\} \,,$$

$$\big|\big(v^* - z\big) \cdot u_{\ell^*}\big| \,=\, \big|\big(v^* - x\big) \cdot u_{\ell^*}\big| \,-\, \frac{r}{4} \,, \tag{20.5}$$

$$\forall \ell \in \{1, \ldots, d\} \setminus \{\ell^*\} \quad \Big|\big(v^* - z\big) \cdot u_\ell\Big| \,\leq\, \max\Big(\Big|\big(v^* - x\big) \cdot u_\ell\Big|, \frac{r}{4}\Big) \,. \tag{20.6}$$

To sum up, when we move from the site $x$ to a site $z$ of target$(x, \Lambda^*)$, and we look at the differences $v^* - x$, $v - z$, the absolute values of the components do not increase, unless they are smaller than $r/4$. Moreover, the absolute value of one component decreases by $r/4$.

**Corollary 20.2.** *Let $\Lambda^*$ be a cubic box in $\mathbb{R}^d$, and let $x$ be a site in $\mathbb{Z}^d \setminus \Lambda^*$. Let $r > 0$ be the side length of $\Lambda^*$ and let $v^*$ be its center. Let $K \geq 1$ and let $(y_k, 0 \leq k \leq K)$ be a sequence of sites such that $y_0 = x$ and*

$$\forall k \in \{\,1,\dots,K\,\} \qquad y_k \in \mathit{target}(y_{k-1},\Lambda^*)\,. \tag{20.7}$$

*Then the whole sequence $(y_k, 0 \leq k \leq K)$ is included in $D(x,\Lambda^*)$ and it satisfies*

$$\forall k \in \{\,0,\dots,K\,\} \qquad \big|y_k - v^*\big|_\infty \leq \max\Big(\big|x-v^*\big|_\infty - \Big\lfloor\frac{k}{d}\Big\rfloor\frac{r}{4},\frac{r}{4}\Big)\,. \tag{20.8}$$

*Proof.* Recall that the supremum norm is given by

$$\forall z \in \mathbb{R}^d \qquad |z|_\infty \,=\, \max\,\big\{\,\big|u_\ell\cdot z\big| : \ell \in \{\,1,\dots,d\,\}\,\big\}\,.$$

It follows from the inequalities (20.5) and (20.6) that

$$\forall k \in \{\,1,\dots,K\,\} \quad \forall \ell \in \{\,1,\dots,d\,\}$$
$$\Big|\big(v^*-y_k\big)\cdot u_\ell\Big| \,\leq\, \max\Big(\Big|\big(v^*-y_{k-1}\big)\cdot u_\ell\Big|,\frac{r}{4}\Big)\,,$$

whence

$$\forall k \in \{\,1,\dots,K\,\} \quad \big|v^*-y_k\big|_\infty \,\leq\, \max\Big(\big|v^*-y_{k-1}\big|_\infty,\frac{r}{4}\Big)\,. \tag{20.9}$$

The inequality (20.9) readily implies the inequality (20.8) for $0 \leq k < d$. We proceed next by induction. Let $m \geq 1$ and suppose that the inequality (20.8) has been proved for $k$ such that $0 \leq k \leq \max(md-1,K)$. For $k = (m-1)d$, we have therefore

$$\big|y_{(m-1)d}-v^*\big|_\infty \,\leq\, \max\Big(\big|x-v^*\big|_\infty - (m-1)\frac{r}{4},\frac{r}{4}\Big)\,. \tag{20.10}$$

If $\big|v^*-y_{(m-1)d}\big|_\infty \leq r/4$, it follows from (20.9) that

$$\forall k \in \{\,(m-1)d,\dots,K\,\} \qquad \big|v^*-y_k\big|_\infty \,\leq\, \frac{r}{4}\,.$$

Otherwise, suppose that $\big|v^*-y_{(m-1)d}\big|_\infty > r/4$. By the inequality (20.5), there exists an index $\ell^*$ such that

$$\big|u_{\ell^*}\cdot(v^*-y_{(m-1)d})\big| \,=\, \max\,\big\{\,\big|u_\ell\cdot(v^*-y_{(m-1)d})\big| : \ell \in \{\,1,\dots,d\,\}\,\big\}\,,$$
$$\big|\big(v^*-y_{(m-1)d+1}\big)\cdot u_{\ell^*}\big| \,=\, \big|\big(v^*-y_{(m-1)d}\big)\cdot u_{\ell^*}\big| - \frac{r}{4}\,,$$

and the same will happen for the subsequent steps, as long as the sequence is outside of $v^*+\Lambda(r/2)$. This implies that, after $d$ steps, we have

$$\big|y_{md}-v^*\big|_\infty \,\leq\, \max\Big(\big|y_{(m-1)d}-v^*\big|_\infty - \frac{r}{4},\frac{r}{4}\Big)\,. \tag{20.11}$$

Combining (20.10) and (20.11), we obtain (20.8) for $k = md$. The inequality (20.9) then yields (20.8) for $k$ such that $md \leq k \leq \max((m+1)d-1,K)$, and the induction is completed. □

A consequence of corollary 20.2 is that a sequence starting at $x$ in $\mathbb{Z}^d \setminus \Lambda^*$ and satisfying (20.7) enters $\Lambda^*$ after a number of steps smaller than $(4d/r)\big|x-v^*\big|_\infty$.

# 21 Travelling towards an orthant

Throughout the section, we fix a cubic box $\Lambda$ included in $\mathbb{R}^d$ and two sites $x, y$ in $\Lambda$. The site $x$ is the starting site and the site $y$ is the target. Our goal is to control the travel time between $x$ and $y$ inside a slight enlargement of the box $\Lambda$. To that end, we will construct iteratively a path between $x$ and $y$ along which the travel times can be controlled. Here we perform the basic step of the construction.

**The set target$(x, y, \Lambda)$.** We define next the set $\text{target}(x, y, \Lambda)$ associated to $x$, $y$ and the box $\Lambda$. For simplicity, we suppose that the box $\Lambda$ is centered at the origin, although we will later apply the construction to an arbitrary box. We denote by $r > 0$ the side length of $\Lambda$, so that $\Lambda = \Lambda(r)$ and we define

$$\forall \varepsilon = (\varepsilon_1, \cdots, \varepsilon_d) \in \{-1, +1\}^d \qquad v(\varepsilon) \,=\, \frac{r}{4}\varepsilon\,.$$

The sites $v(\varepsilon), \varepsilon \in \{-1, +1\}^d$, are the $2^d$ vertices of the box $\Lambda(r/2)$ and the $2^d$ boxes $v(\varepsilon) + \Lambda(r/2)$, $\varepsilon \in \{-1, +1\}^d$, cover $\Lambda$. Let $\varepsilon(y, \Lambda)$ be an element of $\{-1, +1\}^d$ such that the box $v(\varepsilon(y, \Lambda)) + \Lambda(r/2)$ contains the target site $y$. If there exist several choices for $\varepsilon(y, \Lambda)$, then we select one of them according to a deterministic procedure, for instance we can select the smallest one for the lexicographic order. To alleviate the notation, we set

$$v^* \,=\, v\big(\varepsilon(y, \Lambda)\big)\,. \tag{21.1}$$

If $x$ already belongs to the box $v^* + \Lambda(r/2)$, then we set $\text{target}(x, y, \Lambda) = \{\, x \,\}$. Otherwise, we set

$$\text{target}(x, y, \Lambda) \,=\, \text{target}\Big(x, v^* + \Lambda\Big(\frac{r}{2}\Big)\Big)\,,$$

where the right-hand member was defined in (20.1). Thus the box $v^* + \Lambda(r/2)$ plays here the role of the target box considered in section 20. We restate next the result of corollary 20.2 specialized to our current context.

**Corollary 21.1.** *Let $r > 0$ and let $x, y$ be two sites in $\Lambda(r)$. Let $K \geq 1$ and let $(y_k, 0 \leq k \leq K)$ be a sequence of sites in $\Lambda(r)$ such that $y_0 = x$ and*

$$\forall k \in \{\, 1, \ldots, K \,\} \qquad y_k \in \mathit{target}\big(y_{k-1}, y, \Lambda(r)\big)\,. \tag{21.2}$$

*Let $v^* = v\big(\varepsilon(y, \Lambda)\big)$ be the site associated to $y$ and $\Lambda$ defined in (21.1). We have*

$$\forall k \in \{\, 0, \ldots, K \,\} \qquad \big|y_k - v^*\big|_\infty \,\leq\, \max\Big(\big|x - v^*\big|_\infty - \Big\lfloor\frac{k}{d}\Big\rfloor\frac{r}{8}, \frac{r}{8}\Big)\,.$$

A consequence of corollary 21.1 is that a sequence starting in $\Lambda(r)$ and satisfying (21.2) enters the box $v^* + \Lambda(r/2)$ after at most $7d$ steps.

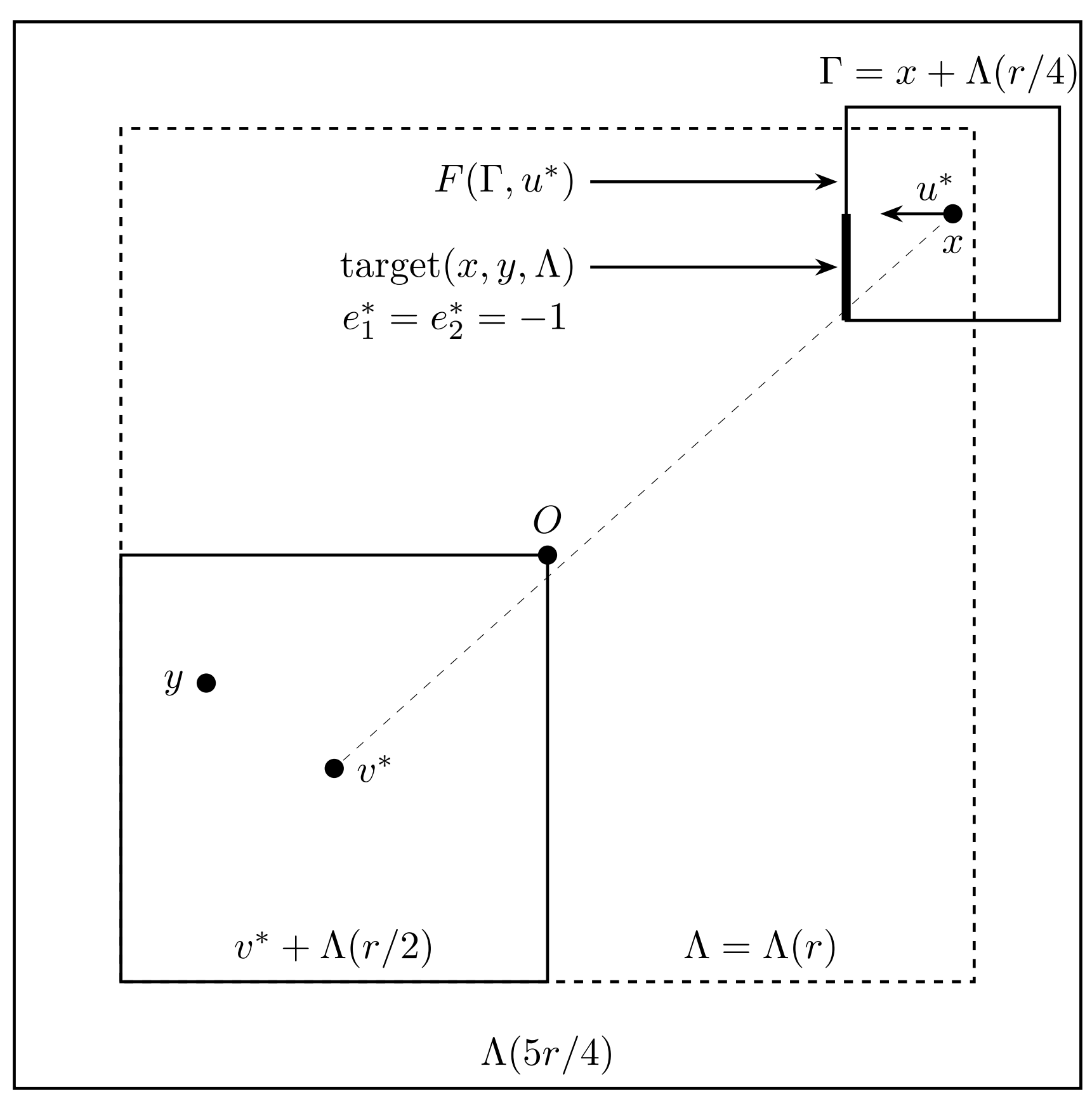

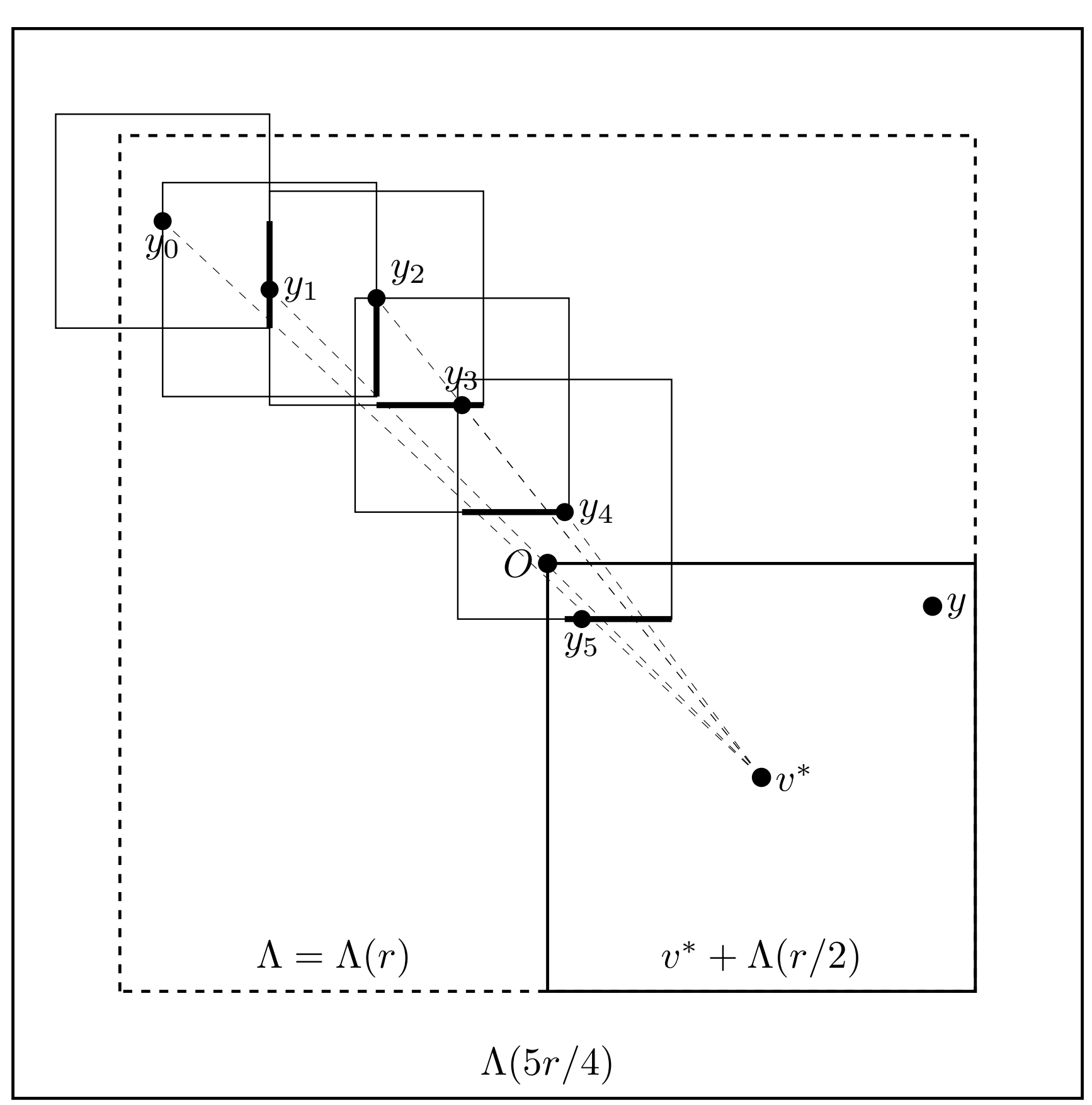

Figure 30: The set target$(x, y, \Lambda)$ and an example of a sequence $(y_k, 0 \leq k \leq K)$.

# 22 Stochastic comparison

Throughout this section, we fix a cubic box $\Lambda$ in dimension $d \geq 3$, and a subbox $\Lambda'$ of $\Lambda$ (possibly equal to $\Lambda$). To each site $y$ of $\Lambda'$, we associate a target set, denoted by $\text{target}(y)$, as in the previous sections. This target set is a deterministic subset of $\Lambda'$, chosen according to the desired goal. Typically, we aim at reaching a target region of the box $\Lambda'$ with a quantitative control on the travel time. We will carry out different constructions, and in each case the following hypothesis will be verified, albeit with specific choices for $\Lambda'$ and $\mu$.

**Hypothesis 22.1.** *There exists a probability measure $\mu$ on $\mathbb{N}$ such that, for any $y$ in $\Lambda'$, the travel time in $\Lambda$ between $y$ and target$(y)$ is stochastically dominated by the distribution $\mu$, i.e.,*

$$\forall y \in \Lambda' \quad \forall t \geq 0 \qquad P\big(T_\Lambda(y, \mathit{target}(y)) \geq t\big) \;\leq\; \mu\big([t, +\infty[\big)\,.$$

Let $x$ be a starting site in $\Lambda'$. We build iteratively a finite sequence

$$(Y_k, \mathcal{D}_k, \tau_k)\,, \quad 0 \leq k \leq K\,,$$

of sites $Y_k$, domains $\mathcal{D}_k$ and times $\tau_k$ as follows. For $k = 0$, we set

$$Y_0 = x\,, \quad \mathcal{D}_0 \;=\; \big\{\, y \in \Lambda : T_\Lambda(x, y) \;=\; 0 \,\big\}\,, \quad \tau_0 = 0\,.$$

Let $k \geq 0$ and suppose that $Y_k$, $\mathcal{D}_k$, $\tau_k$ have been built. We select a site $Y_{k+1}$ in $\text{target}(Y_k)$ such that

$$T_\Lambda(\mathcal{D}_k, \text{target}(Y_k)) \;=\; T_\Lambda(\mathcal{D}_k, Y_{k+1})\,.$$

Since we assumed that $\text{target}(Y_k) \subset \Lambda'$, then $Y_{k+1}$ is still in $\Lambda'$. We set also

$$\begin{aligned}
\tau_{k+1} \;&=\; \tau_k + T_\Lambda(\mathcal{D}_k, \text{target}(Y_k))\,,\\
\mathcal{D}_{k+1} \;&=\; \big\{\, y \in \Lambda : T_\Lambda(x, y) \;\leq\; \tau_{k+1} \,\big\}\,.
\end{aligned}$$

The construction terminates at step $K$ when the desired stopping criterion is fulfilled, typically when the site $Y_k$ has reached the desired goal. The quantitative control of the travel time between the starting site $x$ and the goal will be achieved by controlling the termination step $K$ and using the stochastic comparison result stated in the next proposition. The symbol $\preceq$ stands for stochastic domination, i.e., for two non-negative random variables $U, V$, we have

$$U \preceq V \quad \iff \quad \forall t \geq 0 \qquad P\big(U \geq t\big) \;\leq\; P\big(V \geq t\big)\,.$$

**Proposition 22.2.** *Let $(Y_k,\, \mathcal{D}_k,\, \tau_k)_{0 \leq k \leq K}$ be a sequence constructed as above. Suppose that $K \geq 1$ and that the hypothesis 22.1 holds. We have then*

$$T_\Lambda(Y_0, Y_K) \;\preceq\; K - 1 + T_1 + \cdots + T_K\,,$$

*where $T_1, \ldots, T_K$ are i.i.d. with common distribution $\mu$.*

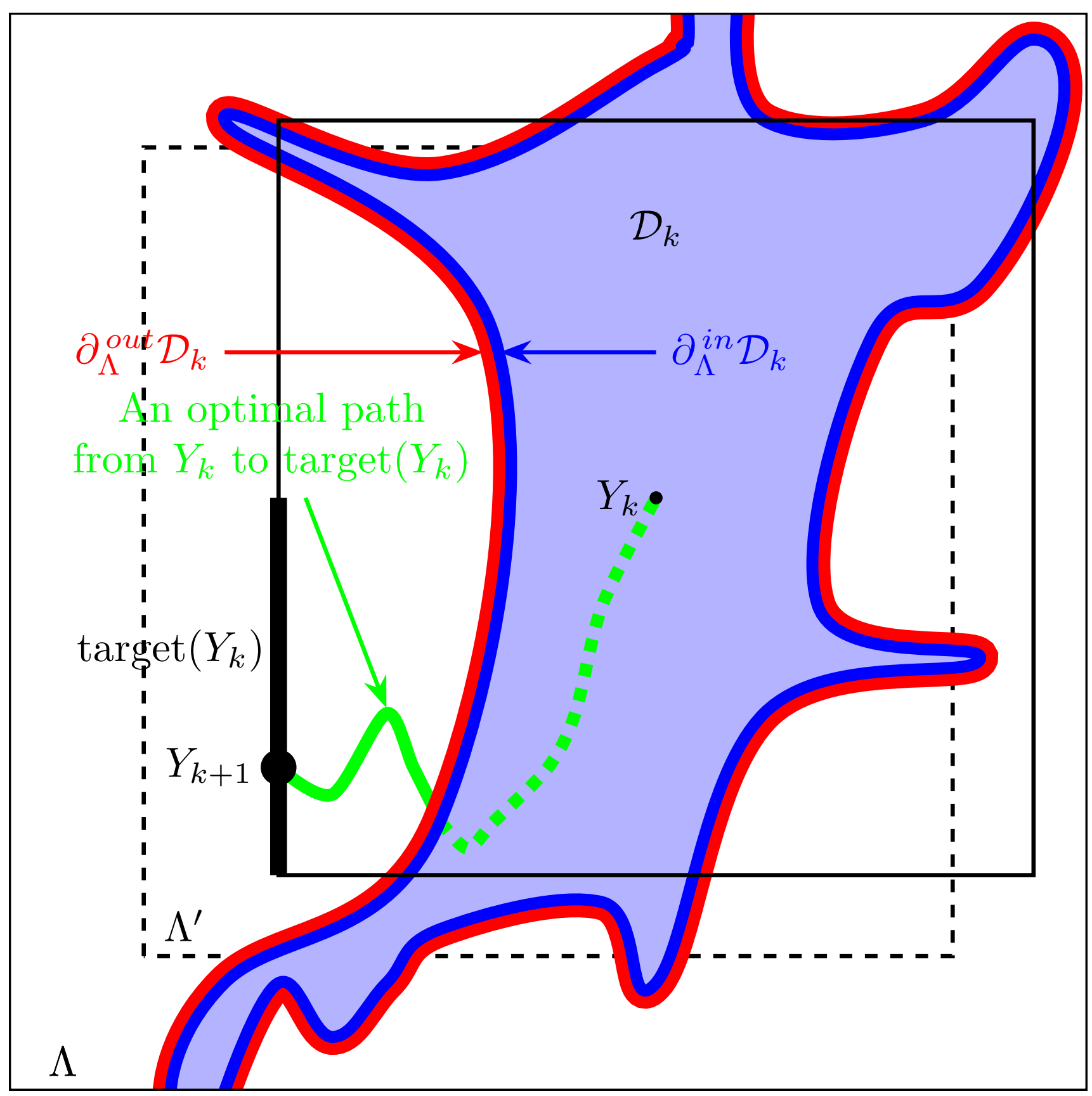

Figure 31: The sites $Y_k$, $Y_{k+1}$ and the sets $\mathcal{D}_k$, target$(Y_k)$

*Proof.* We do the proof by induction on $K$. For $K = 1$, we have

$$T_\Lambda(Y_0, Y_1) \;=\; T_\Lambda\big(\mathcal{D}_0, \mathrm{target}(Y_0)\big) \;=\; \tau_1 \;=\; T_\Lambda\big(x, \mathrm{target}(x)\big)\,,$$

and the hypothesis 22.1 yields readily that $T_\Lambda(Y_0, Y_1) \preceq \mu$. Let next $K \geq 1$ and suppose that the result has been proved at rank $K$. By definition, we have

$$T_\Lambda(Y_0, Y_{K+1}) \;=\; \tau_{K+1} \;=\; T_\Lambda(Y_0, Y_K) + T_\Lambda\big(\mathcal{D}_K, \mathrm{target}(Y_K)\big)\,. \tag{22.1}$$

Let $t \geq 0$. Using (22.1), we decompose the probability $P\big(T_\Lambda(Y_0, Y_{K+1}) \geq t\big)$ as

$$\begin{aligned} P\big(T_\Lambda(Y_0, Y_{K+1}) \geq t\big) \;&=\; P\big(T_\Lambda(Y_0, Y_K) \geq t\big) \\ &+ \sum_{0 \leq s < t} P\Big(T_\Lambda(Y_0, Y_K) = s, T_\Lambda\big(\mathcal{D}_K, \mathrm{target}(Y_K)\big) \geq t - s\Big)\,. \end{aligned} \tag{22.2}$$

We fix next $s$ such that $0 \leq s < t$ and we take care of the probability appearing in the previous sum. We perform a further conditioning on the region $\mathcal{D}_K$ and the site $Y_K$:

$$\begin{aligned} &P\Big(T_\Lambda(Y_0, Y_K) = s, T_\Lambda\big(\mathcal{D}_K, \mathrm{target}(Y_K)\big) \geq t - s\Big) \;= \\ &\sum_{D \subset \Lambda, y \in \Lambda} P\Big(T_\Lambda(x, y) = s, T_\Lambda\big(D, \mathrm{target}(y)\big) \geq t - s, \mathcal{D}_K = D, Y_K = y\Big)\,. \end{aligned} \tag{22.3}$$

For the probability in the sum to be non-zero, the set $D$ and the site $y$ have to satisfy several constraints. Typically, the region $D$ has to be a connected subset of $\Lambda$ containing $y$. Let us fix one such region $D$ and one such site $y$ and suppose that the event

$$\Big\{\, T_\Lambda(x,y)=s, T_\Lambda\big(D,\text{target}(y)\big)\geq t-s, \mathcal{D}_K=D, Y_K=y \,\Big\} \tag{22.4}$$

occurs. Since $T_\Lambda\big(D,\text{target}(y)\big)$ is positive, then necessarily the set target$(y)$ does not intersect $D$. We define the inner boundary $\partial_\Lambda^{in} D$ of $D$ relative to $\Lambda$ as

$$\partial_\Lambda^{in} D \,=\, \big\{\, x\in D : \exists y\in\Lambda\setminus D \quad |x-y|_2=1 \,\big\}\,,$$

and the outer boundary $\partial_\Lambda^{out} D$ of $D$ relative to $\Lambda$ as

$$\partial_\Lambda^{out} D \,=\, \big\{\, x\in \Lambda\setminus D : \exists y\in D \quad |x-y|_2=1 \,\big\}\,.$$

On the event (22.4), we have also that $\tau_K=s$ and

$$\mathcal{D}_K \,=\, D \,=\, \big\{\, z\in\Lambda : T_\Lambda(x,z)\,\leq\, s \,\big\}\,,$$

therefore the sites in $\partial_\Lambda^{in}\mathcal{D}_K=\partial_\Lambda^{in}D$ must be closed. A path which realizes the minimum travel time between $D$ and $y$ has to go through a site of $\partial_\Lambda^{in}D$. By considering the portion of such a path after its last visit to $\partial_\Lambda^{in}D$, we see that

$$T_\Lambda\big(D,\text{target}(y)\big) \,=\, 1+T_{\Lambda\setminus D}\big(\partial_\Lambda^{out}D,\text{target}(y)\big)\,. \tag{22.5}$$

Now the random variable $T_{\Lambda\setminus D}\big(\partial_\Lambda^{out}D,\text{target}(y)\big)$ is measurable with respect to the sites belonging to $\Lambda\setminus D$, while the event

$$\Big\{\, T_\Lambda(x,y)=s, \mathcal{D}_K=D, Y_K=y \,\Big\}$$

is measurable with respect to the sites belonging to $D$. Using in addition the equality (22.5), we have therefore

$$\begin{aligned} &P\Big(\, T_\Lambda(x,y)=s, T_\Lambda\big(D,\text{target}(y)\big)\geq t-s, \mathcal{D}_K=D, Y_K=y\Big) \qquad (22.6)\\ &=\, P\Big(\, T_\Lambda(x,y)=s, \mathcal{D}_K=D, Y_K=y, T_{\Lambda\setminus D}\big(\partial_\Lambda^{out}D,\text{target}(y)\big)\geq t-s-1\,\Big)\\ &=\, P\Big(\, T_\Lambda(x,y)=s, \mathcal{D}_K=D, Y_K=y\Big)P\Big(T_{\Lambda\setminus D}\big(\partial_\Lambda^{out}D,\text{target}(y)\big)\geq t-s-1\Big). \end{aligned}$$

The identity (22.5) and the fact that $y$ is in $D$, which is a subset of $\Lambda$, yield that

$$T_{\Lambda\setminus D}\big(\partial_\Lambda^{out}D,\text{target}(y)\big) \,\leq\, T_\Lambda\big(D,\text{target}(y)\big) \,\leq\, T_\Lambda\big(y,\text{target}(y)\big)\,. \tag{22.7}$$

The inequality (22.7) implies furthermore that

$$P\Big(T_{\Lambda\setminus D}\big(\partial_\Lambda^{out}D,\text{target}(y)\big)\geq t-s-1\Big) \,\leq\, P\Big(T_\Lambda\big(y,\text{target}(y)\big)\geq t-s-1\Big)\,. \tag{22.8}$$

It follows from (22.6), the inequality (22.8) and the hypothesis 22.1 that

$$\begin{aligned}
P\Big( T_\Lambda(x,y) = s, T_\Lambda\big(D, \text{target}(y)\big) \geq t-s, \mathcal{D}_K = D, Y_K = y \Big) \\
\leq P\Big( T_\Lambda(x,y) = s, \mathcal{D}_K = D, Y_K = y \Big) P\Big(T_\Lambda\big(y, \text{target}(y)\big) \geq t-s-1\Big) \\
\leq P\Big( T_\Lambda(x,y) = s, \mathcal{D}_K = D, Y_K = y \Big)\, \mu\big([\,t-s-1, +\infty[\big)\,. \qquad (22.9)
\end{aligned}$$

Substituting the inequality (22.9) into the sum appearing in (22.3), we obtain

$$\begin{aligned}
P\Big(T_\Lambda(Y_0, Y_K) = s, T_\Lambda\big(\mathcal{D}_K, \text{target}(Y_K)\big) \geq t-s\Big) \\
\leq \sum_{D \subset \Lambda, y \in \Lambda} P\Big( T_\Lambda(x,y) = s, \mathcal{D}_K = D, Y_K = y \Big)\, \mu\big([\,t-s-1, +\infty[\big) \\
= P\big(T_\Lambda(Y_0, Y_K) = s\big)\, \mu\big([\,t-s-1, +\infty[\big)\,. \qquad (22.10)
\end{aligned}$$

We substitute next inequality (22.10) into (22.2) and we get

$$\begin{aligned}
P\big(T_\Lambda(Y_0, Y_{K+1}) \geq t\big) \;\leq\; P\big(T_\Lambda(Y_0, Y_K) \geq t\big) + \\
\sum_{0 \leq s < t} P\big(T_\Lambda(Y_0, Y_K) = s\big)\, \mu\big([\,t-s-1, +\infty[\big)\,. \qquad (22.11)
\end{aligned}$$

Let $X$ be a random variable distributed according to $\mu$ and which is independent of $T_\Lambda(Y_0, Y_K)$. We rewrite (22.11) as

$$\begin{aligned}
P\big(T_\Lambda(Y_0, Y_{K+1}) \geq t\big) \\
\leq P\big(T_\Lambda(Y_0, Y_K) \geq t\big) + \sum_{0 \leq s < t} P\big(T_\Lambda(Y_0, Y_K) = s\big)\, P(X + 1 \geq t - s) \\
= P\big(T_\Lambda(Y_0, Y_K) + X + 1 \geq t\big)\,. \qquad (22.12)
\end{aligned}$$

By the induction hypothesis, we have

$$T_\Lambda(Y_0, Y_K) \;\preceq\; K - 1 + T_1 + \cdots + T_K\,, \qquad (22.13)$$

where $T_1, \dots, T_K$ are i.i.d. with common distribution $\mu$. It follows from (22.12) and (22.13) that

$$P\big(T_\Lambda(Y_0, Y_{K+1}) \geq t\big) \;\leq\; P\big(K + T_1 + \cdots + T_K + X \geq t\big)\,,$$

where $X$ is independent of $T_1, \dots, T_K$. This is true for any $t \geq 0$, therefore

$$T_\Lambda(Y_0, Y_{K+1}) \;\preceq\; K + T_1 + \cdots + T_K + X\,,$$

and the induction step is completed. □

# 23 Geometric domination

In this section, we consider a fixed cubic box $\Lambda$ in dimension $d \geq 3$. We check that the hypothesis 22.1 on the domination of the travel time between a site and its target is satisfied with $\Lambda' = \Lambda$ and $\mu$ an adequate geometric distribution.

Let $\Gamma$ be any symmetric box included in $\Lambda$ centered at a site $x$ of $\mathbb{Z}^d \cap \Lambda$. Denoting by $r$ the side length of $\Gamma$, we have thus $\Gamma = x + \Lambda(r)$. We work with the collection $\mathcal{F}(\Gamma)$ of the tiles of the boundary of $\Gamma$ defined in (18.3) and (18.4). Since

$$\partial^{in}\Gamma \subset \bigcup_{F \in \mathcal{F}(\Gamma)} F \,,$$

we have

$$T_\Gamma(x, \partial^{in}\Gamma) \;=\; \min\big\{\, T_\Gamma(x,F) : F \in \mathcal{F}(\Gamma) \,\big\} \,,$$

and by the FKG inequality, for any $k \geq 0$,

$$\begin{aligned} P\big(T_\Gamma(x, \partial^{in}\Gamma) \geq k\big) \;&=\; P\big(\forall F \in \mathcal{F}(\Gamma) \quad T_\Gamma(x,F) \geq k\big) \\ &\geq\; \prod_{F \in \mathcal{F}(\Gamma)} P\big(T_\Gamma(x,F) \geq k\big) \,. \end{aligned} \tag{23.1}$$

The percolation model is invariant under the action of the group $G(\mathcal{H})$ and also invariant by translation. Therefore all the probabilities appearing in the last product are equal, and we deduce from (23.1) that

$$\forall F \in \mathcal{F}(\Gamma) \qquad P\big(T_\Gamma(x,F) \geq k\big) \;\leq\; \Big(P\big(T_\Gamma(x, \partial^{in}\Gamma) \geq k\big)\Big)^{1/|\mathcal{F}(\Gamma)|} \,. \tag{23.2}$$

The number of tiles is $|\mathcal{F}(\Gamma)| = d2^d$. It follows from lemma 14.1 that

$$\forall k \geq 0 \qquad P\big(\, T_\Gamma(x, \partial^{in}\Gamma) \;\geq\; k \,\big) \;\leq\; \big(1 - \theta(p)\big)^k \,. \tag{23.3}$$

Combining the inequalities (23.2) and (23.3), we obtain

$$\forall F \in \mathcal{F}(\Gamma) \quad \forall k \geq 0 \qquad P\big(T_\Gamma(x,F) \geq k\big) \;\leq\; \big(1 - \theta(p)\big)^{k/(d2^d)} \,.$$

As $\Gamma$ is included in $\Lambda$, we have also

$$\forall F \in \mathcal{F}(\Gamma) \qquad T_\Gamma(x,F) \geq T_\Lambda(x,F) \,,$$

and we conclude finally that

$$\forall F \in \mathcal{F}(\Gamma) \quad \forall k \geq 0 \qquad P\big(T_\Lambda(x,F) \geq k\big) \;\leq\; \big(1 - \theta(p)\big)^{k/(d2^d)} \,. \tag{23.4}$$

The good news here is that the upper bound (23.4) is uniform with respect to the side length of the box $\Gamma$. In the next constructions, the target set associated to a site $x$ of $\Lambda \setminus \partial^{in}\Lambda$ will always be a tile of the boundary of some box centered at $x$ and included in $\Lambda$. The choice of the tile will depend in a deterministic way on $x$ and on the goal we are trying to reach. In any case, the hypothesis 22.1 is satisfied with $\mu$ the geometric distribution of parameter $1 - \big(1 - \theta(p)\big)^{1/(d2^d)}$. In addition, the same is true if we work with discrete balls instead of cubes, and we use the tiling of their boundaries described in subsection 19.3. The only difference is that the number of tiles would be $d2^{d+1}$ instead of $d2^d$.

# 24 Large deviations estimate

Another ingredient that we will need to derive our estimates on the travel times is a large deviations inequality for a sum of i.i.d. geometric variables. Naturally, we could make appeal to the classical Chernoff inequality. However, we prefer to use a simpler elementary inequality, which relies on the following little trick that can be found in a work of Svante Janson (see [75], the very beginning of the proof of theorem 2.1).

**Lemma 24.1.** *Let* $\lambda$ *in* $]0,1[$ *and let* $T$ *be a random variable distributed according to the geometric law Geo*$(\lambda)$. *We have*

$$\forall t<\lambda \qquad E\big(e^{tT}\big)\,\leq\,\frac{\lambda}{\lambda-t}\,. \tag{24.1}$$

*Proof.* Let $t<\lambda$. We have $0\leq e^t(1-\lambda)\leq e^{t-\lambda}<1$, hence we can sum the geometric progression of common ratio $e^t(1-\lambda)$ to get

$$E\big(e^{tT}\big)\,=\,\sum_{k\geq 1}e^{tk}(1-\lambda)^{k-1}\lambda\,=\,\lambda e^t\frac{1}{1-e^t(1-\lambda)}\,=\,\frac{\lambda}{e^{-t}-1+\lambda}\,.$$

We then use the inequality $e^{-t}-1\geq -t$ to obtain the desired result (24.1). □

Once we have a bound on the Laplace transform of the geometric distribution, we obtain automatically a large deviations inequality. The advantage of the bound (24.1) is that it is very simple and allows to derive an explicit upper bound with minimal computations.

**Proposition 24.2.** *Let* $\lambda\in]0,1[$, $n\geq 1$, *and let* $T_1,\dots,T_n$ *be* $n$ *i.i.d. random variables distributed according to the geometric law Geo*$(\lambda)$ *of parameter* $\lambda$. *We have*

$$\forall t\geq 1 \qquad P\Big(T_1+\cdots+T_n\geq t\frac{n}{\lambda}\Big)\,\leq\,\big(t\,e^{1-t}\big)^n\,. \tag{24.2}$$

*Proof.* Let $t\geq 1$ be fixed and let $u>0$ be another parameter, to be chosen later. We write

$$\begin{aligned} P\Big(T_1+\cdots+T_n\geq t\frac{n}{\lambda}\Big)\,&=\,P\Big(e^{u(T_1+\cdots+T_n)}\geq e^{utn/\lambda}\Big)\\ &\leq\,e^{-utn/\lambda}E\Big(e^{u(T_1+\cdots+T_n)}\Big)\\ &=\,e^{-utn/\lambda}\Big(E\big(e^{uT_1}\big)\Big)^n\,. \end{aligned} \tag{24.3}$$

Using the inequality (24.1), we deduce from (24.3) that, for $u<\lambda$,

$$P\Big(T_1+\cdots+T_n\geq t\frac{n}{\lambda}\Big)\,\leq\,e^{-utn/\lambda}\Big(\frac{\lambda}{\lambda-u}\Big)^n\,=\,\Big(e^{-ut/\lambda}\frac{1}{1-u/\lambda}\Big)^n\,.$$

We choose now $u=\big(1-1/t\big)\lambda$, so that $u/\lambda=1-1/t$, and we get the desired inequality (24.2). □

# 25 Synthesis

In this section, we use the geometric constructions of sections 18-21 and the estimates of sections 22-24 to obtain crucial results on the travel times.

## 25.1 The three-dimensional case

We consider a fixed three-dimensional cubic box $\Lambda$ and, under the sole hypothesis that $\theta(p,\mathbb{Z}^3)>0$, we obtain the following estimate on the travel time between two arbitrary points of $\Lambda$.

**Theorem 25.1.** *We consider the site percolation model in dimension* 3 *and a parameter $p$ such that $\theta(p,\mathbb{Z}^3)>0$. For any $\alpha>0$, setting*

$$\kappa(\alpha,p)\,=\,\frac{6\alpha+14000}{1-\big(1-\theta(p,\mathbb{Z}^3)\big)^{1/48}}\,,$$

*we have*

$$\forall n\geq 100\qquad P\big(\forall x,y\in\Lambda(n)\quad T_{\Lambda(n)}(x,y)\,\leq\,\kappa(\alpha,p)\ln n\big)\,\geq\,1-\frac{1}{n^\alpha}\,.\quad(25.1)$$

The following result had previously been proved in [27].

**Theorem 25.2** (Theorem 1.1 of [27])**.** *Let $p$ be such that $\theta(p)>0$ and let $\alpha>0$. There exists a constant $\kappa$, depending on $\alpha$ and $p$, such that*

$$\forall n\geq 2\qquad P\begin{pmatrix}\text{every pair of sites of the box }\Lambda(n)\\ \text{are joined by a path in }\Lambda(n)\text{ having}\\ \text{at most }\kappa(\ln n)^2\text{ closed sites}\end{pmatrix}\,\geq\,1-\frac{1}{n^\alpha}\,.$$

Theorem 25.1 constitutes a major improvement over theorem 25.2, indeed the travel times are shown to be of order $\ln n$ instead of $(\ln n)^2$. The geometric constructions involved in the proof of theorem 25.1 are the same as the ones used in theorem 25.2. The improvement of the result comes from the use of the stochastic comparison developed in section 22.

*Proof of theorem 25.1.* Let $n\geq 100$ be fixed and let $x$ be an arbitrary point of $\Lambda(n)$. We will use the constructions of the previous sections in order to build a random path in $\Lambda(n)$ starting at $x$ and terminating at a site close to the center 0 of $\Lambda(n)$. In a first step, we ensure that the starting point is not on the boundary $\partial^{\,in}\Lambda(n)$. If $x$ is in $\Lambda(n)\setminus\partial^{\,in}\Lambda(n)$, then we set $Y_0=x$. Otherwise, if $x$ belongs to $\partial^{\,in}\Lambda(n)$, we pick up according to a deterministic procedure a site $y_0$ in $\Lambda(n)\setminus\partial^{\,in}\Lambda(n)$ such that $|x-y_0|_\infty\leq 1$ and we set $Y_0=y_0$. This way, we have selected a site $Y_0$ which is not on the boundary $\partial^{\,in}\Lambda(n)$. Moreover, with our choice for $Y_0$, we have $T_{\Lambda(n)}(x,Y_0)\,\leq\,3$. We use next the construction of a random finite sequence $(Y_k,\mathcal{D}_k,\tau_k)$, $0\leq k\leq K$, described at the beginning of section 22. Recall that this construction was iterative: at each step $k$, we selected a site $Y_{k+1}$ in $\text{target}(Y_k)$ such that

$$T_\Lambda\big(\mathcal{D}_k,\text{target}(Y_k)\big)\,=\,T_\Lambda(\mathcal{D}_k,Y_{k+1})\,.$$

The point is that, at each step, we are free to choose the target set $\text{target}(Y_k)$, depending on the goal we wish to achieve. We use the beginning of the path to step away from the boundary, and we employ the target sets defined in the section 18. We do so until the sequence enters the half-box $\Lambda(n/2)$. Corollary 18.2 ensures that this occurs at some random index $L$ satisfying

$$L \,\leq\, 3\frac{\ln n}{\ln 2} + 1\,. \tag{25.2}$$

If the starting site $Y_0$ was already in $\Lambda(n/2)$, we set $L = 0$. From the index $L$ onwards, the goal is to reach the center 0 of the box. Thus, for any step $k \geq L$, we would like to employ the target set $\text{target}(Y_k, 0)$ associated to $Y_k$ and 0, as defined in the subsection 19.4. However, we cannot do so directly, because the set $\text{target}(y, z)$ is defined only for pairs of sites $y, z$ satisfying $|y-z|_2 \leq n/4$. We have therefore to design a finite number of intermediate targets, which we denote by $Z_i$, $0 \leq i \leq 3$. We define also random indices $L_i$, $0 \leq i \leq 3$, corresponding to the activation time of these targets, as follows. First we set $Z_0 = Y_L$ and $L_0 = L$. Suppose next that $Z_i$ and $L_i$ have been defined for some $i \geq 0$, and let us explain how we define $Z_{i+1}$ and $L_{i+1}$. For $Z_{i+1}$, we choose a site such that

$$|Y_{L_i} - Z_{i+1}|_2 \,\leq\, \frac{n}{4}\,, \qquad |Z_{i+1}|_2 \,\leq\, \max\Big(|Y_{L_i}|_2 - \frac{n}{4} + 6, 600\Big)\,. \tag{25.3}$$

From the index $L_i$ onwards, the site $Z_{i+1}$ becomes the new intermediate goal, and we employ the target set $\text{target}(Y_k, Z_{i+1})$ associated to $Y_k$ and $Z_{i+1}$, until we reach a site at Euclidean distance less than or equal to 600 from $Z_{i+1}$. More precisely, we define

$$L_{i+1} \,=\, \min\big\{\, k \geq L_i : |Y_k - Z_{i+1}|_2 \leq 600\,\big\}\,. \tag{25.4}$$

Corollary 19.3 ensures that

$$\forall k \geq L_i \qquad |Z_{i+1} - Y_k|_2 \,\leq\, \max\Big(0.93^{k-L_i}|Z_{i+1} - Y_{L_i}|_2, 600\Big)\,,$$

therefore

$$L_{i+1} - L_i \,\leq\, \frac{\ln n - \ln 600}{-\ln 0.93} \,\leq\, 14\ln n\,. \tag{25.5}$$

It follows from (25.3) and (25.4) that

$$|Y_{L_{i+1}}|_2 \,\leq\, \max\Big(|Y_{L_i}|_2 - \frac{n}{4} + 6, 600\Big) + 600 \,\leq\, \max\Big(|Y_{L_i}|_2 - \frac{n}{4} + 606, 1200\Big)\,. \tag{25.6}$$

Iterating the inequality (25.6), we obtain

$$|Y_{L_i}|_2 \,\leq\, \max\Big(|Y_{L_0}|_2 - i\Big(\frac{n}{4} - 606\Big), 1200\Big)\,.$$

Since $|Y_{L_0}|_2 \leq 3n/4$, then $|Y_{L_3}|_2 \leq 2000$. We end the path at the index $K = L_3$. The inequalities (25.2) and (25.5) imply that, for $n \geq 100$,

$$K \,\leq\, 3\frac{\ln n}{\ln 2} + 1 \,+\, 42\ln n \,\leq\, 47\ln n\,. \tag{25.7}$$

We have completed the construction of the sequence $(Y_k, \mathcal{D}_k, \tau_k)$, $0 \leq k \leq K$. It was checked in section 23 that the target sets employed for the construction do satisfy the hypothesis 22.1 with $\Lambda' = \Lambda$ and $\mu$ being the geometric distribution of parameter $1 - \rho_3$, with $\rho_3 = \big(1-\theta(p)\big)^{1/24}$. For the initial part from $Y_0$ to $Y_L$, this is a consequence of the inequality (23.4). For the subsequent parts, for $L_0 \leq k \leq L_3 - 1$, we worked with balls instead of cubes, but the very same reasoning as the one done in section 23 yields the analog of (23.4) for balls. The only difference is that the constant $\rho_3$ should be replaced by $\rho_3^* = \big(1-\theta(p)\big)^{1/48}$. This constant is worse because there are twice as many tiles on the balls. We did not write formally the details, because the argument is the same, and also because a method valid for all dimensions $d \geq 3$ is exposed in subsection 25.2. In order to have a uniform bound, we work with $\rho_3^*$. Applying proposition 22.2, we conclude that

$$\begin{aligned} T_\Lambda(x,0) \;&\leq\; T_\Lambda(x,Y_0) + T_\Lambda(Y_0,Y_K) + 3|Y_{L_3}|_2 \\ &\preceq\; 3 + K - 1 + T_1 + \cdots + T_K + 6000\,, \end{aligned}$$

where $T_1, \dots, T_K$ are i.i.d. with geometric distribution of parameter $1 - \rho_3^*$. Let us set $M = \lfloor 47 \ln n \rfloor$. Since we know from (25.7) that $K \leq M$, we conclude that

$$T_\Lambda(x,0) \;\preceq\; M + 6002 + T_1 + \cdots + T_M\,. \tag{25.8}$$

For $n \geq 100$, we have $M \geq 46 \ln n$. We apply next proposition 24.2 with $\lambda = 1 - \rho_3^*$, and we get

$$\forall t \geq 1 \qquad P\Big(T_1 + \cdots + T_M \geq \frac{47t}{1-\rho_3^*} \ln n\Big) \;\leq\; \big(t\, e^{1-t}\big)^{46 \ln n}\,. \tag{25.9}$$

Let $\alpha > 0$. We take $t = \alpha + 2$. Since

$$46(\ln(\alpha+2) - \alpha - 1) \;=\; 46\Big(\ln 2 + \ln\Big(1 + \frac{\alpha}{2}\Big) - \alpha - 1\Big) \;\leq\; -14 - 23\alpha\,,$$

we deduce from (25.8) and (25.9) that

$$P\Big(T_\Lambda(x,0) \geq 47 \ln n + 6002 + \frac{47\alpha + 94}{1-\rho_3^*} \ln n\Big) \;\leq\; \frac{1}{n^{14+23\alpha}}\,. \tag{25.10}$$

This is true for any $x$ in $\Lambda = \Lambda(n)$, thus by a simple union bound, we have

$$P\Big(\forall x \in \Lambda \quad T_\Lambda(x,0) < 47 \ln n + 6002 + \frac{3\alpha + 94}{1-\rho_3^*} \ln n\Big) \;\geq\; 1 - \frac{1}{n^{10+\alpha}} \tag{25.11}$$

(we have used (25.10) with $\alpha/23$ instead of $\alpha$). We notice finally that

$$\forall x,y \in \Lambda \quad T_\Lambda(x,y) \;\leq\; T_\Lambda(x,0) + T_\Lambda(0,y) \;\leq\; T_\Lambda(x,0) + 1 + T_\Lambda(y,0)\,, \tag{25.12}$$

and the desired inequality (25.1) follows from (25.11) and (25.12). □

Unfortunately, the previous proof does not work in dimensions $d \geq 4$. That is the height of irony: we have now a method to get an essential estimate on the travel time in a finite box which works in three dimensions, but not in higher dimensions!

## 25.2 The case of higher dimensions

Let $d \geq 3$. In this subsection, we consider a fixed $d$-dimensional symmetric box $\Lambda$ and, under the sole hypothesis that $\theta(p, \mathbb{Z}^d) > 0$, we obtain the following estimate on the travel time between two arbitrary points of $\Lambda$.

**Theorem 25.3.** *We consider the site percolation model in dimension $d$ and a parameter $p$ such that $\theta(p, \mathbb{Z}^d) > 0$. For any $\alpha > 0$, setting*

$$\kappa(\alpha, p) \;=\; \frac{6\alpha + 60d^2}{1 - \big(1 - \theta(p, \mathbb{Z}^d)\big)^{1/(d2^d)}}\,, \tag{25.13}$$

*we have*

$$\forall n \geq 2 \qquad P\Big(\forall x, y \in \Lambda(n) \quad T_{\Lambda(n)}(x, y) \;\leq\; \kappa(\alpha, p) \ln n\Big) \;\geq\; 1 - \frac{1}{n^\alpha}\,. \tag{25.14}$$

*Proof.* Let $n \geq 2$ be fixed and let $x$ be an arbitrary point of $\Lambda(n)$. The strategy of the proof is similar to the proof of theorem 25.1 in dimension 3. We will use the constructions of sections 18-21 in order to build a random path in $\Lambda(n)$ starting at $x$ and terminating at a site close to the center 0 of $\Lambda(n)$. However, the construction of the path inside $\Lambda(n/2)$ is more involved, it requires a multiscale argument. To ensure that this proof can be read independently of the proof of theorem 25.1, we reproduce the beginning of the construction, although it is the same. In a first step, we choose a starting point outside the boundary $\partial^{\,in}\Lambda(n)$. If $x$ is in $\Lambda(n) \setminus \partial^{\,in}\Lambda(n)$, then we set $Y_0 = x$. Otherwise, if $x$ belongs to $\partial^{\,in}\Lambda(n)$, we pick up according to a deterministic procedure a site $y_0$ in $\Lambda(n) \setminus \partial^{\,in}\Lambda(n)$ such that $|x - y_0|_\infty \leq 1$ and we set $Y_0 = y_0$. This way we have selected a site $Y_0$ which is not on the boundary $\partial^{\,in}\Lambda(n)$. Moreover, we have $T_{\Lambda(n)}(x, Y_0) \leq d$. We use next the construction of a random finite sequence $(Y_k, \mathcal{D}_k, \tau_k)$, $0 \leq k \leq K$, described at the beginning of section 22. Recall that this construction was iterative: at each step $k$, we selected a site $Y_{k+1}$ in $\text{target}(Y_k)$ such that $T_\Lambda\big(\mathcal{D}_k, \text{target}(Y_k)\big) = T_\Lambda(\mathcal{D}_k, Y_{k+1})$. The point is that, at each step, we are free to choose the target set $\text{target}(Y_k)$, depending on the goal we wish to achieve. We use the beginning of the path to step away from the boundary of $\Lambda(n)$, and we employ the target sets defined in the section 18. We do so until the sequence enters the half-box $\Lambda(n/2)$. Corollary 18.2 ensures that this occurs at some random index $L$ satisfying

$$L \;\leq\; d\frac{\ln n}{\ln 2} + 1\,. \tag{25.15}$$

If the starting site $Y_0$ was already in $\Lambda(n/2)$, we set $L = 0$. From the index $L$ onwards, the goal is to reach the center 0 of the box. We cannot proceed as in the proof of theorem 25.1, because we do not have at our disposal a method to perform a move which divides the distance to the origin by a constant factor larger than 1. The problem is that the tiles of the boundary of the boxes or of the $(d-1)$-dimensional sphere induced by the symmetries of the lattice have a too large diameter. We will therefore use a multiscale argument. On each intermediate scale, we will use the construction done in section 21, which allows to travel between two boxes in a finite number of controllable moves. The target site will always be the origin, but the starting site and the associated box will be redefined on each scale. Let us denote by $S$ the random number of changes of scale involved in the construction. For each $i$ in $\{1,\dots,S\}$, we define the following objects, all of them random:

• the beginning $L_{i-1}$ and the end $L_i$ of the use of the scale number $i$ in the construction of the sequence $(Y_k,\mathcal{D}_k,\tau_k)$, $0\leq k\leq K$;

• the starting site $X_i$ employed on the scale $i$;

• the box $\Lambda_i$ used to define the target during the scale $i$. The center of $\Lambda_i$ is denoted by $C_i$ and its side length by $R_i$.

We give now the precise details of the construction. For the first scale, we set

$$L_0=L\,,\quad X_0=Y_L\,,\quad C_0=0\,,\quad R_0=\frac{n}{2}\,,\quad \Lambda_0=C_0+\Lambda(R_0)\,. \tag{25.16}$$

Suppose next that $L_i$, $X_i$, $C_i$, $R_i$ and $\Lambda_i$ have been defined for some $i\geq 0$, and let us explain how we build $L_{i+1}$, $X_{i+1}$, $C_{i+1}$, $R_{i+1}$ and $\Lambda_{i+1}$. As the origin is always the target site, throughout the scale $i$, we employ the target set $\mathrm{target}(Y_k,0,\Lambda_i)$ defined in section 21 associated to $Y_k$ and 0 in the box $\Lambda_i$. In section 21, we worked with boxes centered at 0, but of course the geometric results obtained there are valid for translated boxes as well. We denote by $V_i$ the center of the sub-box of $\Lambda_i$ which dictates the choice of the tile $\mathrm{target}(Y_{L_i},0,\Lambda_i)$. Thus the sub-box corresponding to $v^*+\Lambda(r/2)$ in section 21 is denoted here $V_i+\Lambda(R_i/2)$. We start from $Y_{L_i}$ and we iterate the construction until we reach a site belonging to the sub-box $V_i+\Lambda(R_i/2)$. More precisely, we define

$$L_{i+1}\;=\;\min\big\{\,k\geq L_i:Y_k\in V_i+\Lambda(R_i/2)\,\big\}\,.$$

A consequence of corollary 20.2 is that at most $7d$ steps are required to enter the box $V_i+\Lambda(R_i/2)$, therefore

$$L_{i+1}-L_i\;\leq\;7d\,. \tag{25.17}$$

We define $X_{i+1}=Y_{L_{i+1}}$. For the box $\Lambda_{i+1}$, we choose a cubic box which contains both $X_{i+1}$ and the origin $O$, whose center is a site of $\Lambda(n/2)\cap\mathbb{Z}^d$ and which has minimum side length. Typically, there exist a lot of such boxes, we can for instance choose a box such that $X_{i+1}$ or $O$ is one of the corners of the box. We denote by $C_{i+1}$ the center of $\Lambda_{i+1}$ and its side length by $R_{i+1}$. Since $X_{i+1}$ is in the box $V_i+\Lambda(R_i/2)$, which contains also 0, then $|X_{i+1}|_\infty\leq R_i/2$ and

$$R_{i+1}\;\leq\;\frac{R_i}{2}+1\,. \tag{25.18}$$

The presence of $+1$ is due to the constraint that the center of the box has to belong to $\mathbb{Z}^d$. Iterating the inequalities (25.17), (25.18), and using the bounds on the initial terms given in (25.15) and (25.16), we obtain

$$L_i\;\leq\;7di+d\frac{\ln n}{\ln 2}+1\,,\qquad R_i\;\leq\;\frac{n}{2^i}+2\,. \tag{25.19}$$

The construction is properly defined as long as we deal with boxes whose faces can be tiled into $2^{d-1}$ isometric pieces. For safety, we require that the side length of the smallest box involved is larger than $2d$, and we stop the construction as soon as $R_i$ becomes smaller than $8d$. Thus we define

$$I \;=\; \min\big\{\, i \geq 0 : R_i \;\leq\; 8d \,\big\}\,,$$

and we end the path at the index $K = L_I$. We have thus $|Y_K|_\infty = |X_I|_\infty \leq 8d$. The inequalities (25.19) imply that, for $n \geq 2$,

$$I \;\leq\; \frac{\ln n}{\ln 2}\,, \qquad K \;\leq\; 7d\frac{\ln n}{\ln 2} + d\frac{\ln n}{\ln 2} + 1 \;\leq\; 13d \ln n\,. \tag{25.20}$$

We have completed the construction of the sequence $(Y_k, \mathcal{D}_k, \tau_k)$, $0 \leq k \leq K$. The target sets used in the construction are tiles of the boundary of cubic boxes. It was checked in section 23 that these target sets satisfy the hypothesis 22.1 with $\Lambda' = \Lambda$ and $\mu$ the geometric law of parameter

$$\lambda \;=\; 1 - \big(1 - \theta(p)\big)^{1/(d2^d)}\,.$$

Applying proposition 22.2, we conclude that

$$\begin{aligned} T_\Lambda(x, O) \;\leq\; & T_\Lambda(x, Y_0) + T_\Lambda(Y_0, Y_K) + T_\Lambda(Y_K, O) \\ & \qquad\qquad \preceq\; d + K - 1 + T_1 + \cdots + T_K + 8d^2\,, \end{aligned}$$

where $T_1, \dots, T_K$ are i.i.d. with geometric distribution of parameter $\lambda$. Let us set

$$M \;=\; \lfloor 13d \ln n \rfloor\,.$$

We know from (25.20) that $K \leq M$, therefore

$$T_\Lambda(x, O) \;\preceq\; M + 9d^2 + T_1 + \cdots + T_M\,. \tag{25.21}$$

When specialized to $d = 3$, the inequality (25.21) is better than (25.8), which was obtained with balls instead of cubes in dimension three. For $n \geq 1$ and $d \geq 3$, we have

$$12d \ln n \;\leq\; M \;\leq\; 14d \ln n\,.$$

Applying proposition 24.2, we get

$$\forall t \geq 1 \qquad P\Big(T_1 + \cdots + T_M \geq \frac{14dt}{\lambda} \ln n\Big) \;\leq\; \big(t\, e^{1-t}\big)^{12d \ln n}\,. \tag{25.22}$$

Let $\alpha > 0$. We take $t = \alpha + 2$. Since

$$12d\big(\ln(\alpha+2) - \alpha - 1\big) \;=\; 12d\Big(\ln 2 + \ln\Big(1 + \frac{\alpha}{2}\Big) - \alpha - 1\Big) \;\leq\; -3d - 6d\alpha\,,$$

we deduce from (25.21) and (25.22) that

$$P\Big(T_\Lambda(x, O) \geq 14d \ln n + 9d^2 + \frac{14d(\alpha+2)}{\lambda} \ln n\Big) \;\leq\; \frac{1}{n^{3d+6d\alpha}}\,.$$

We replace $\alpha$ by $\alpha/(6d)$, we rework a bit the inequality, and we get

$$P\Big(T_\Lambda(x,O) \geq \frac{3\alpha+28d^2}{\lambda}\ln n\Big) \;\leq\; \frac{1}{n^{3d+\alpha}}\,.$$

This is true for any $x$ in $\Lambda = \Lambda(n)$, thus by a simple union bound, we have

$$P\Big(\forall x\in\Lambda \quad T_\Lambda(x,O) < \frac{3\alpha+28d^2}{\lambda}\ln n\Big) \;\geq\; 1-\frac{1}{n^{d+\alpha}}\,. \tag{25.23}$$

We notice next that

$$\forall x,y\in\Lambda \quad T_\Lambda(x,y) \;\leq\; T_\Lambda(x,O)+T_\Lambda(O,y) \;\leq\; T_\Lambda(x,O)+1+T_\Lambda(y,O)\,. \tag{25.24}$$

It follows from (25.23) and (25.24) that

$$P\Big(\forall x,y\in\Lambda(n) \quad T_{\Lambda(n)}(x,y) < \frac{6\alpha+56d^2}{\lambda}\ln n+1\Big) \;\geq\; 1-\frac{2}{n^{d+\alpha}}\,.$$

This yields the desired inequality (25.14) □

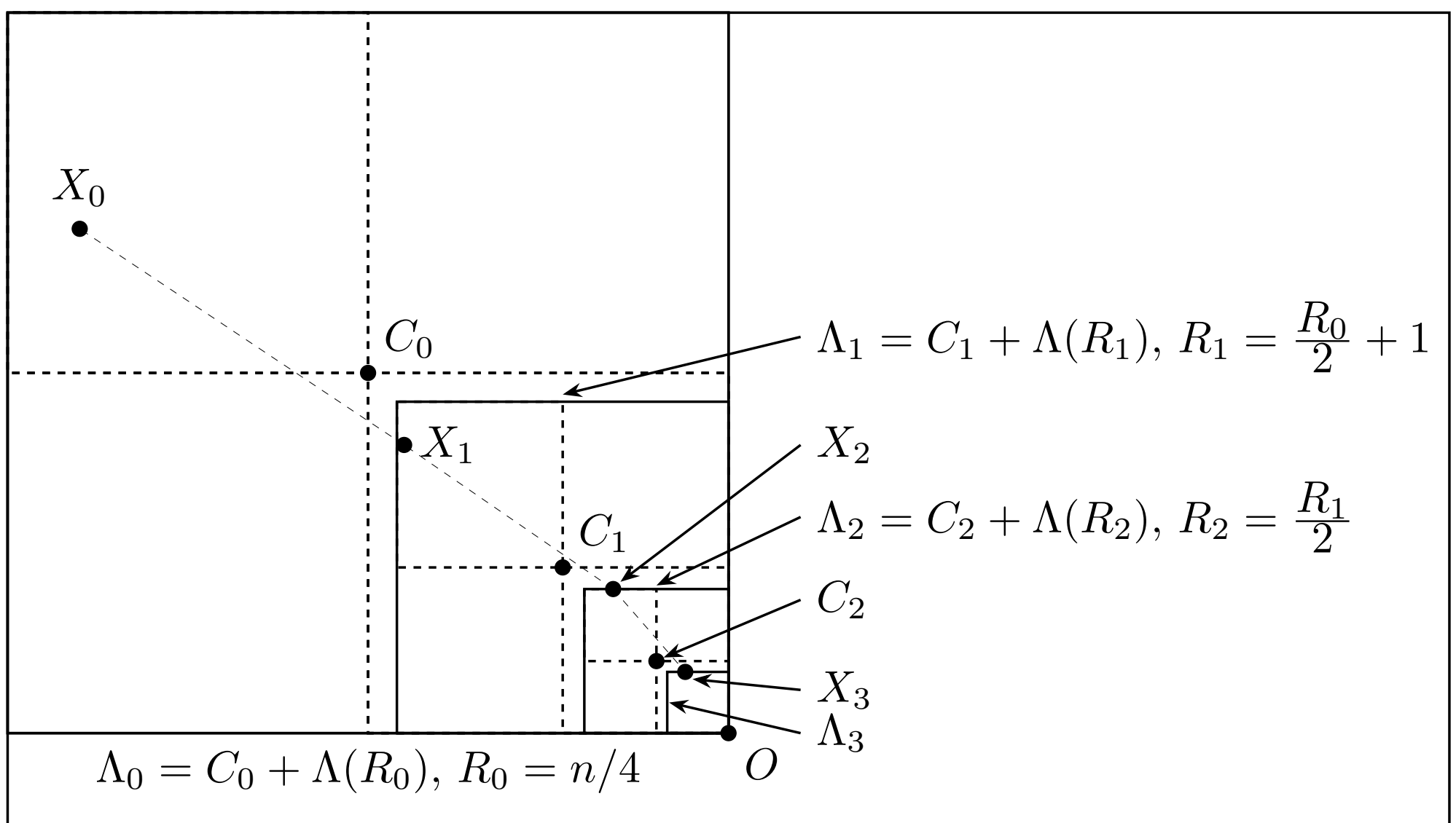

Figure 32: The sites $X_i$ and the boxes $\Lambda_i$ for $0 \leq i \leq 3$

## 25.3 Grimmett-Marstrand

It will not have escaped to the reader familiar with the Grimmett-Marstrand theorem that the construction performed on each scale is the very same as the driving process used by Grimmett and Marstrand to travel between two boxes of their renormalized lattice. To some extent, theorem 25.3 can be seen as a finite volume version of the Grimmett-Marstrand result. Still, this is not exact, because an essential ingredient in Grimmett-Marstrand is the famous sprinkling procedure, which is absent here. We absolutely want to avoid resorting to sprinkling, because we need a quantitative estimate valid at the parameter $p$, and we cannot afford to increase it. The advantage is that the proof is considerably simpler than the full proof of the Grimmett-Marstrand theorem. The price we have to pay is a weakening of the final result. Instead of getting a path which is fully open, we get a path with some closed sites, but we have a control on their number. In the reverse direction, we could try to apply the multiscale construction of the previous section to the percolation configuration in infinite volume, instead of working inside a finite box. We believe that it would lead to the following result.

**Potential theorem 25.4.** *We consider the site percolation model in dimension $d$ and a parameter $p$ such that $\theta(p,\mathbb{Z}^d)>0$. For any $\alpha>0$, there exist two positive constants $\kappa_\infty(\alpha,p)$, $c_\infty(d)$ such that*

$$\forall x,y\in\mathbb{Z}^d \qquad P\Big(T(x,y)\,\geq\,\kappa_\infty(\alpha,p)\ln\big(|x-y|_2\big)+c_\infty(d)\Big)\,\leq\,\frac{1}{\big(|x-y|_2\big)^\alpha}\,.$$

An important issue is whether the constant $\kappa(\alpha,p)$ defined in (25.13) can be taken arbitrarily small. The inequality (25.14) was obtained with the help of the stochastic domination result of proposition 22.2, applied with the geometric law of parameter $\lambda=1-\big(1-\theta(p)\big)^{1/(d2^d)}$. The result would be improved if we could take $\lambda$ close to 1. A first attempt would be to use an idea of Martineau and Tassion. In their work [103], they implement a version of the Grimmett-Marstrand construction in which the cubic seeds are replaced by a portion of the path discovered by the exploration mechanism. However, this would cause trouble, because we would lose the symmetry used at each step to apply the square root trick, and it does not seem plausible that each tile of a cube has the same probability of being connected to a given initial path. A second attempt would be to use large finite seeds as in the Grimmett-Marstrand construction. It seems likely that the previous construction can be carried over with seeds, and that the parameter $\lambda$ of the geometric variables appearing in the proof of theorem 25.3 can be taken as close to 1 as desired. Yet there is also a geometric content in (25.14), which corresponds to the number of steps involved in the multiscale construction. This number is proportional to the number of intermediate scales needed to travel from one site to another in the box $\Lambda(n)$, and it seems very difficult to perform a geometric construction which would give a number negligible compared to $\ln n$. Although the value $\kappa(\alpha,p)$ defined in (25.13) is not optimal, there is little hope that it can be improved further than a barrier $\kappa'(d)>0$ which depends on the dimension $d$ only.

## 25.4 Travelling from $\Lambda(4n)$ to $\Lambda(2n)$ inside $\Lambda(6n)$

Let $d \geq 3$. In this subsection, we consider a fixed $d$-dimensional symmetric box $\Lambda$ and, under the sole hypothesis that $\theta(p, \mathbb{Z}^d) > 0$, we obtain an estimate on the travel time between an arbitrary point of a cubical shell and its inner cube.

**Proposition 25.5.** *Let $p \in [0,1]$ be such that $\theta(p, \mathbb{Z}^d) > 0$. There exists a non-decreasing function $\beta(n)$ defined for $n \geq 3$ such that*

$$\forall n \geq 3 \qquad \beta(n) \leq \sqrt{\ln n}\,,$$
$$\beta(3) \geq \frac{\theta(p)}{d2^d}\,, \qquad \lim_{n \to +\infty} \beta(n) = +\infty\,,$$

*and moreover*

$$\forall n \geq 3 \qquad P\Big( \exists x \in \Lambda(4n) \quad T_{\Lambda(6n)}\big(x, \Lambda(2n)\big) > \frac{\ln n}{\beta(n)} \Big) \leq n^{5d - \beta(n)}. \tag{25.25}$$

*Proof.* Let $n \geq 1$ be fixed and let $x$ be an arbitrary point in $\Lambda(4n) \setminus \Lambda(2n)$. We are going to build a random path in $\Lambda(6n)$ starting at $x$ and terminating at a site belonging to $\Lambda(2n)$. As before, we use the construction of a random finite sequence $(Y_k, \mathcal{D}_k, \tau_k)$, $0 \leq k \leq K$, described at the beginning of section 22. The goal here is to enter the box $\Lambda(2n)$. We rely on the construction performed in section 20. The box $\Lambda(6n)$ plays the role of the fixed box of section 20 throughout the construction. We start from $Y_0 = x$. At each step, we employ the target set target$(Y_k, \Lambda(2n))$ defined in section 20 associated to $Y_k$ and the box $\Lambda^* = \Lambda(2n)$. We stop the sequence at the first index $K$ such that $Y_K$ is in $\Lambda(2n)$. This completes the construction of the sequence $(Y_k, \mathcal{D}_k, \tau_k)$, $0 \leq k \leq K$. We apply next corollary 20.2 to the site $x$ and the box $\Lambda^* = \Lambda(2n)$. We have $v^* = 0$, $r = 2n$, thus

$$K \leq \frac{4d}{r} \big| x - v^* \big|_\infty \leq 4d\,, \tag{25.26}$$

and moreover the sequence $(Y_k, 0 \leq k \leq K)$ is included in $D(x, \Lambda(2n))$, which is itself included in $\Lambda(4n)$ (recall that $D(x, \Lambda(2n))$ is the smallest axis-parallel hyperrectangle containing $x$ and $\Lambda(2n)$). In addition, all the intermediate boxes $\Gamma$ intervening in the construction of the sequence $(Y_k, 0 \leq k \leq K)$ are translates of $\Lambda(n)$ centered at a site of $D\big(x, \Lambda(2n)\big)$, these boxes are included in the enlargement

$$D(x, \Lambda^*) + \Lambda\Big(\frac{r}{2}\Big) = D\big(x, \Lambda(2n)\big) + \Lambda(n) \subset \Lambda(4n) + \Lambda(n) = \Lambda(6n)\,,$$

and this legitimates the use of $\Lambda(6n)$ as the fixed box into which the travel times are computed during the construction of the sequence $(Y_k, 0 \leq k \leq K)$. Furthermore, for any $y$ in $\Lambda(4n)$, the target set associated to $y$ is chosen deterministically among the collection of tiles $\mathcal{F}(y + \Lambda(n))$ defined in (18.3) and (18.4). Let us set

$$\mathbf{1} = (1, \cdots, 1) \in \{-1, +1\}^d\,.$$

It follows from the symmetries of the model and the translation invariance that

$$\begin{aligned}
&\forall y\in\Lambda(4n)\quad \forall i\in\{\,1,\dots,d\,\}\quad \forall v\in\{\,\pm u_i\,\}\quad \forall\varepsilon\in\{-1,+1\}^d\quad \forall t\geq 0\\
&P\Big(T_{\Lambda(6n)}\big(y,y+F(\Lambda(n),v,\varepsilon)\big)\geq t\Big)\;\leq\; P\Big(T_{y+\Lambda(n)}\big(y,y+F(\Lambda(n),v,\varepsilon)\big)\geq t\Big)\\
&\qquad\qquad\qquad\qquad = P\Big(T_{\Lambda(n)}\big(0,F(\Lambda(n),u_1,\mathbf{1})\big)\geq t\Big)\,. \qquad (25.27)
\end{aligned}$$

Let $\nu_n$ be the distribution of the random variable $T_{\Lambda(n)}\big(0,F(\Lambda(n),u_1,\mathbf{1})\big)$. The inequality (25.27) yields that

$$\forall y\in\Lambda(4n)\quad \forall t\geq 0\qquad P\big(T_{\Lambda(6n)}(y,\text{target}(y))\geq t\big)\;\leq\;\nu_n\big([t,+\infty[\big)\,.$$

Thus the hypothesis 22.1 on the domination of the travel times is satisfied with $\Lambda=\Lambda(6n)$, $\Lambda'=\Lambda(4n)$ and the distribution $\nu_n$. Applying proposition 22.2 together with the bound (25.26), we conclude that

$$\begin{aligned}
T_{\Lambda(6n)}\big(x,\Lambda(2n)\big)\;\leq\; T_\Lambda(Y_0,Y_K)\;&\preceq\; K-1+T_1+\cdots+T_K\\
&\preceq\; 4d+T_1+\cdots+T_{4d}\,, \qquad (25.28)
\end{aligned}$$

where $T_1,\dots,T_{4d}$ are i.i.d. with common distribution $\nu_n$. If we had only used the box $\Lambda(4n)$ to compute the travel times, then the best we could have done is to bound stochastically the common distribution of the intermediate travel times by the geometric distribution of parameter $1-\big(1-\theta(p)\big)^{1/(d2^d)}$, as we did in the previous section. Indeed, the boxes whose centers are close to the boundary of $\Lambda(4n)$ have a small side length, and there is no room to improve upon the geometric estimate of (23.4). The main advantage of using the larger box $\Lambda(6n)$ to compute the travel times is that, for sites in $\Lambda(4n)$, we can always use boxes of side length $n$ for the intermediate steps. This allows us to reach $\Lambda(2n)$ with only $4d$ steps, and the estimate on the tail of the travel times given by corollary 14.6 becomes particularly relevant. It yields that

$$\forall t\geq 0\qquad \nu_n\big([t,+\infty]\big)\;\leq\;\exp\Big(-\frac{t\alpha(t)}{d2^d}\Big)\,, \qquad (25.29)$$

where $\alpha(t)$ is the function appearing in corollary 14.5. Thus the tail $\nu_n$ decays super-exponentially fast. Now the inequalities (25.28) and (25.29) give that

$$\forall x\in\Lambda(4n)\quad \forall t\geq 0\qquad P\Big(T_{\Lambda(6n)}\big(x,\Lambda(2n)\big)\geq 4d(1+t)\Big)\;\leq\;4d\exp\Big(-\frac{t\alpha(t)}{d2^d}\Big)\,.$$

From the standard union bound, it follows that, for any $t\geq 1$,

$$P\Big(\exists x\in\Lambda(4n)\quad T_{\Lambda(6n)}\big(x,\Lambda(2n)\big)\geq 8dt\Big)\;\leq\;|\Lambda(4n)|\,4d\exp\Big(-\frac{t\alpha(t)}{d2^d}\Big)\,. \qquad (25.30)$$

Let $\beta(n)$ be the function defined by

$$\forall n\geq 3\qquad \beta(n)\;=\;\frac{1}{d2^d}\min\Big(\sqrt{\ln n},\sqrt{\alpha\big(\sqrt{\ln n}\big)}\Big)\,.$$

Thanks to the fact that $\lim_{t\to\infty}\alpha(t)=+\infty$, we have also $\lim_{n\to\infty}\beta(n)=+\infty$. Recalling that the function $\alpha(t)$ is non–decreasing, we have furthermore

$$\frac{\ln n}{d2^d\beta(n)}\,\geq\,\sqrt{\ln n}\,,\qquad \alpha\Big(\frac{\ln n}{d2^d\beta(n)}\Big)\,\geq\,\alpha\big(\sqrt{\ln n}\big)\,\geq\,\Big(d2^d\beta(n)\Big)^2\,,$$

whence

$$\forall n\geq 3\qquad \frac{1}{d2^d\beta(n)}\alpha\Big(\frac{\ln n}{d2^d\beta(n)}\Big)\,\geq\,d2^d\beta(n)\,. \tag{25.31}$$

Since $\alpha(1)\geq\theta(p)$ by corollary 14.5, we have also

$$d2^d\beta(3)\,=\,\min\Big(\sqrt{\ln 3},\sqrt{\alpha\big(\sqrt{\ln 3}\big)}\Big)\,\geq\,\min\big(1,\sqrt{\alpha(1)}\big)\,\geq\,\sqrt{\theta(p)}\,\geq\,\theta(p)\,.$$

We take $t=\ln n/\big(d2^d\beta(n)\big)$ in the inequality (25.30) and, using (25.31) and the simple bound $4d|\Lambda(4n)|\leq(5n)^d\leq n^{5d}$, we get the desired inequality (25.25). □

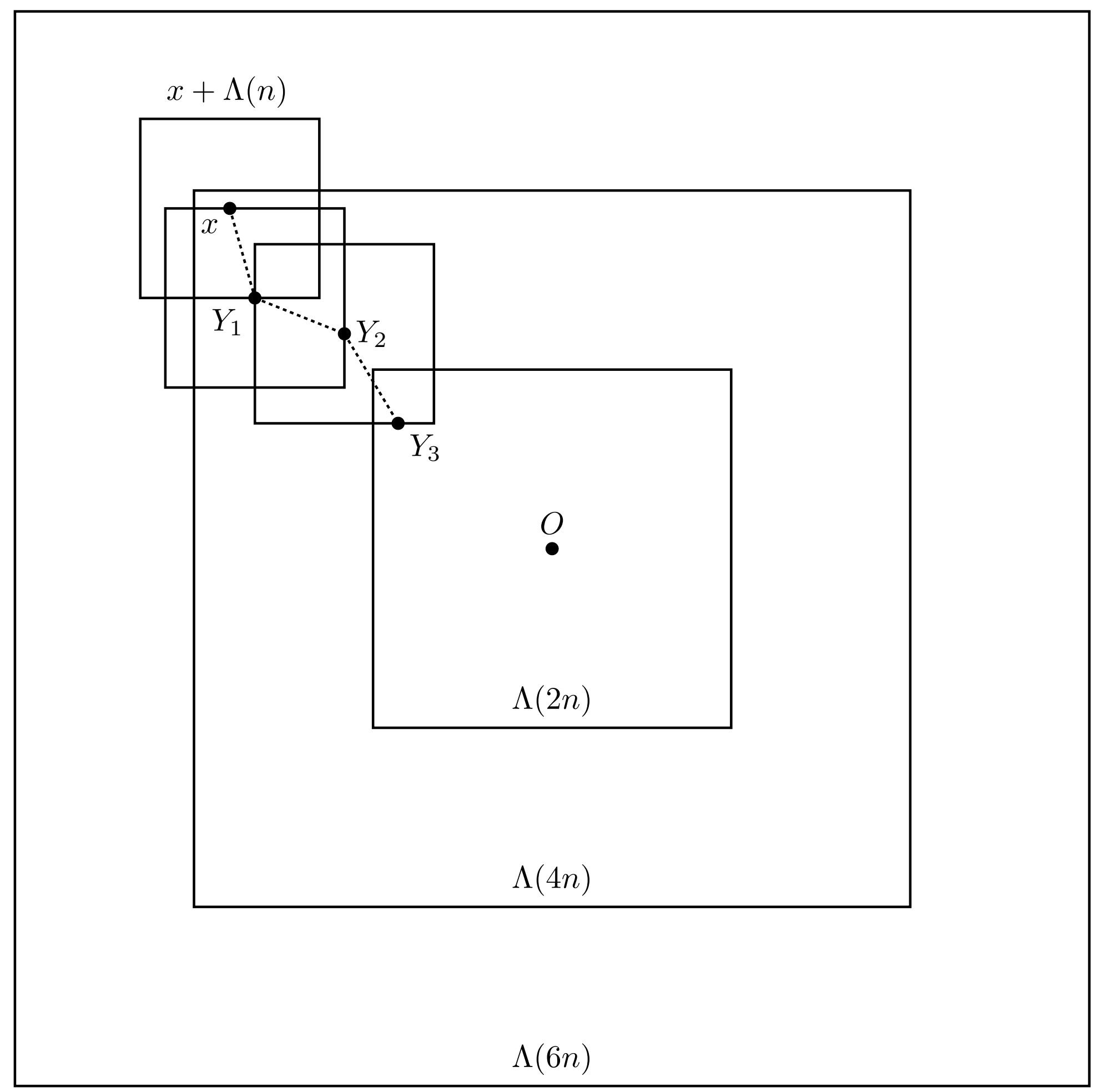

Figure 33: Travelling from $\Lambda(4n)\setminus\Lambda(2n)$ to $\Lambda(2n)$ within $\Lambda(6n)$

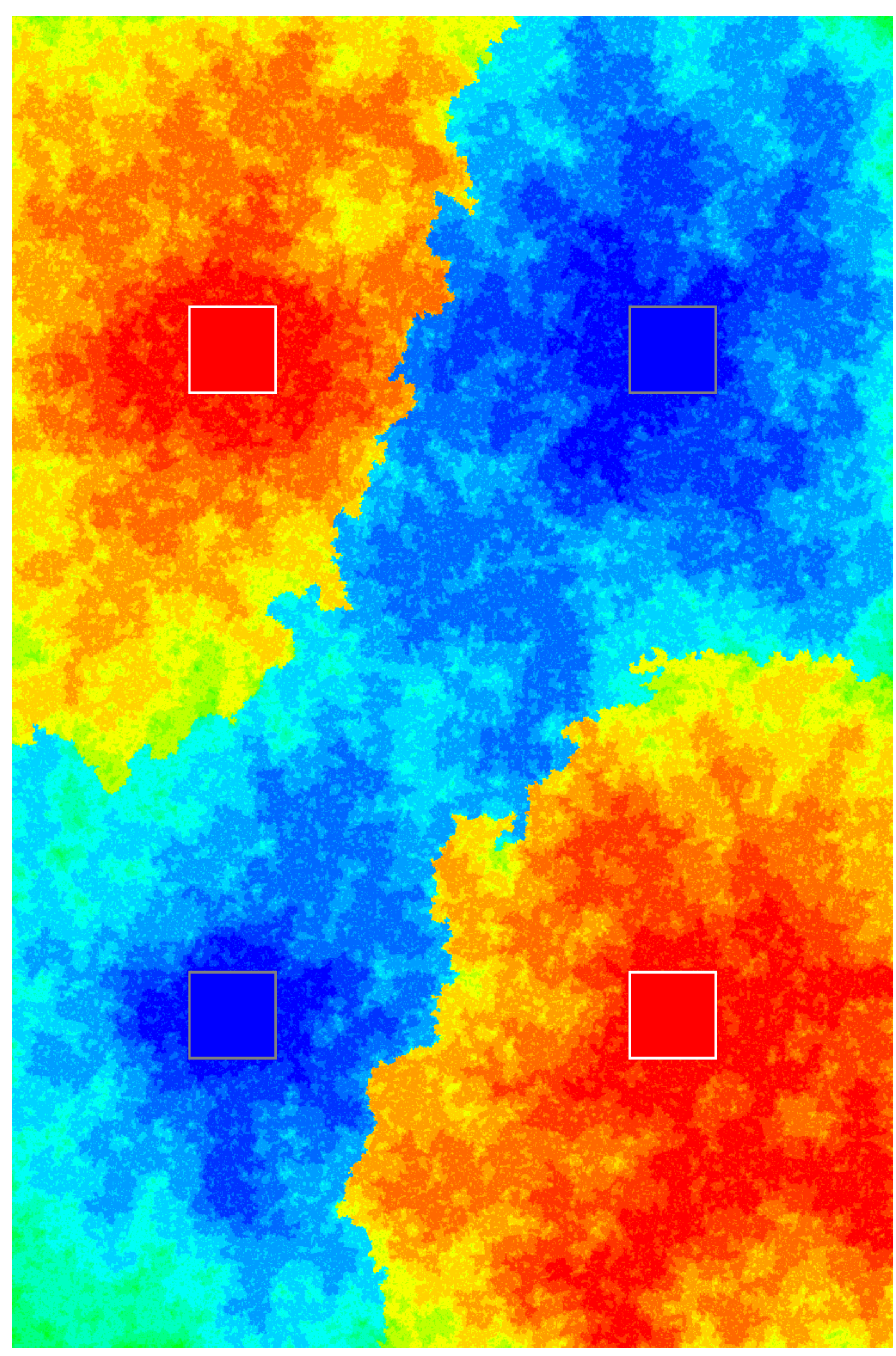

Figure 34: Intertwined bond exploration in the rectangle $1024 \times 1550$, $p = 0.485$ $W, W' =$ white, gray squares, lighter colors correspond to higher rounds.

**Part V** ⊛ ⤫

# Warming up in the bond model

In this part, we focus on the bond model. In section 26, we compare the site and the bond models. In section 27, we present two different ways of visualizing a configuration of bond percolation. All the results seen so far were written for the site model. We explain in section 28 how to adapt them to the bond model. The intertwined bond explorations are constructed in section 29. In section 30, we describe three challenging questions, that should reasonably be attacked before trying the conjecture $\theta(p_c, \mathbb{Z}^d) = 0$. The main problem is to control the bond intersections revealed by the intertwined explorations. We present a brutal way to do that in section 31. Although it gives a non-obvious result, it is not enough to reach an interesting conclusion for the three challenging questions. In section 33, we present new attempts to control the second intersections.

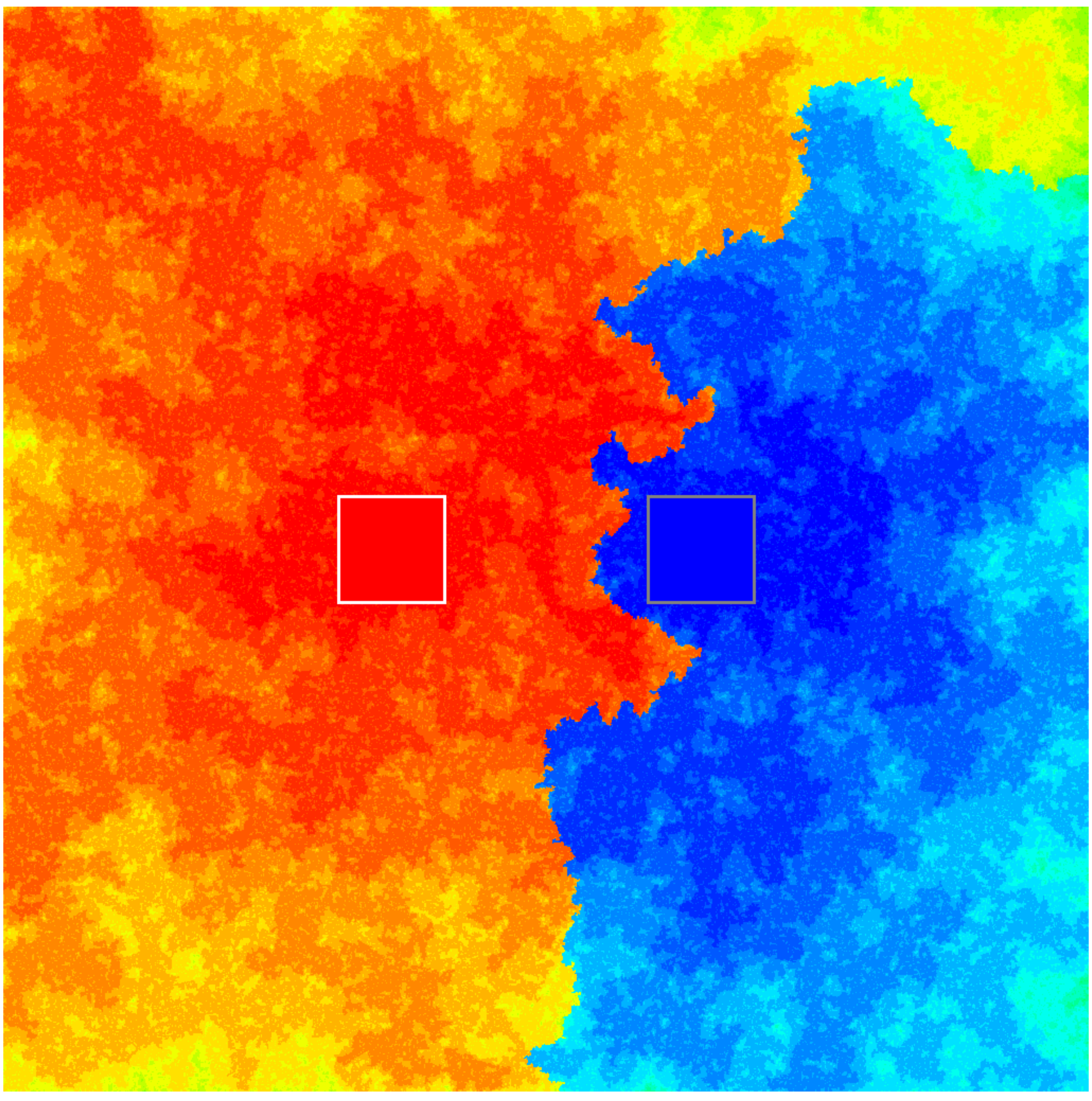

Figure 35: Intertwined exploration in $\Lambda(1024)$ for bond percolation at $p = 0.49$, $W$ = white square, $W'$ = gray square, lighter colors correspond to higher rounds.

# 26 The bond model versus the site model

Before delving into the details of the constructions, we discuss the differences between the site and the bond models. The fundamental difference between the two models is that the random blocking mechanism is built on the vertices in the site model and on the edges in the bond model. At first sight, one might think that this innocent-looking variation should have little effect on the behavior of the model. Indeed, many results and proofs are literally the same for both models. However, this is not always the case. For instance, the critical point for bond percolation on $\mathbb{Z}^2$ is $1/2$, while the exact value of the critical point for site percolation on $\mathbb{Z}^2$ is not even known. It is also known that the site model is in principle more general than the bond model, in the sense that any bond percolation model is isomorphic to a site percolation model, albeit on a different graph, while the converse is not true (see [63], section 1.6). Concerning the simulations, the site model is very convenient because we can use very simple memory structures, like $2d$ or $3d$ arrays to store all the relevant information. The bond model is more involved to implement, because it requires a specific structure to store the states of the bonds. When working on the lattice, we can still use $2d$ or $3d$ arrays (essentially by transforming the bond model into a site model). A more natural approach is to implement fully the graph structure, by maintaining an adjacency list associated to each vertex of the graph. This requires to build additional memory structures, but it is more efficient, and that is indeed how we proceed to simulate the bond model.

With regard to our problem, our initial plan was to write everything for the site model, as it is more general than the bond model. It turns out that the mathematical proofs involving interacting explorations are far more intricate for site percolation than for bond percolation. After several unsuccessful attempts to write these proofs directly for the site model, we had to come back to the bond model. So we present first the construction of the intertwined explorations for the bond model, and we write in the framework of the bond model several interesting auxiliary results, even if they are not conclusive. When we will be ready to write the final proof, we will come back to the site percolation model. Surprisingly, a similar phenomenon occurs when we try to implement certain specific simulations. A notable example is Tarjan's algorithm for the second positive pivots in subsection 45.7. Some heavy numerical computations have been conducted for the site model, typically to simulate the function $\theta(p,\mathbb{Z}^3)$. The programs are written in C++ and they make appeal to the Gnu Scientific Library [60]. The simulations whose aim is to illustrate the text are done for the bond model, unless they concern a question specific to the site model. The corresponding programs are written in C++ with the help of the Standard Template Library [78]. From an aesthetic point of view, the pictures of bond percolation are much prettier, because there is always a positive density of uninteresting black pixels in the pictures of site percolation (when a site is closed, its cluster is declared empty and the site is colored in black). These black pixels behave like an independent Bernoulli field, and they create a background noise that is detrimental to the quality of the pictures.

## 27 Visualization of the bonds

Until now, the pictures of the percolation configurations presented only the vertices of the graph. This is fine for the site model, but it is a bit annoying for the bond model, because the pictures do not yield complete information on the states of the bonds. Typically, we color the open clusters of the configuration. If two neighbouring vertices have different colors, then the bond joining them must be closed, but if they have the same color, the state of the bond joining them remains indeterminate. To remedy this inconvenience (which, by the way, is reminiscent of the question on information transmission raised in subsection 17.4), when we wish to have access to the full information on the bonds, we shall present the percolation configuration in an array of pixels twice the size of the real configuration. The vertices correspond to the pixels whose coordinates are both even. The pixels whose coordinates are both odd are just holes corresponding to nothing. The remaining pixels are in one to one correspondence with the bonds (this is how the bond model is mapped onto a site model). If a bond is closed, its associated pixel stays black. If a bond is open, its associated pixel is colored in the same way as its two endpoints (which obviously belong to the same open cluster). This representation will be very useful to depict the pivotal bonds in future sections. It will also help us to understand the dynamics of the intertwined explorations. In figures 36 and 37, we see two bond percolation configurations, and for each of them, the visualization with or without the bonds. For the configuration in $\Lambda(32)$ in figure 36, one cannot decide with the help of the left image which bonds inside the big red cluster are open. Note that the pictures without the bonds have been rescaled to force them to have the same size as the pictures with the bonds. For configurations in larger boxes, as in figure 37, the overall effect of visualizing the bonds is to blur the picture. This visual degradation is even worse than the blurring induced by the presence of the black pixels in the site percolation pictures.

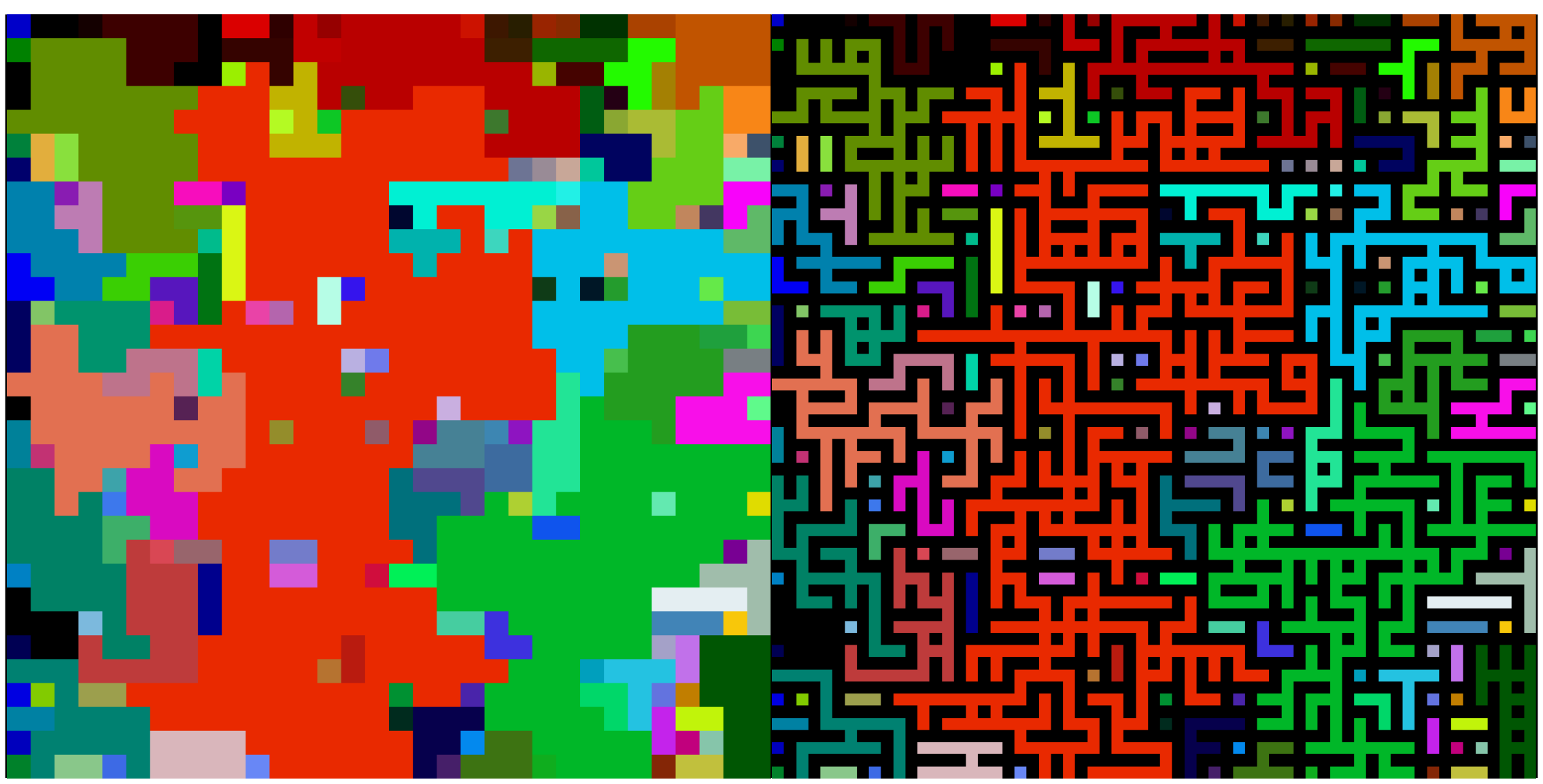

Figure 36: Two visualizations of the same bond percolation in $\Lambda(32)$, $p = 0.45$. Left: colored vertices only. Right: colored vertices and bonds.

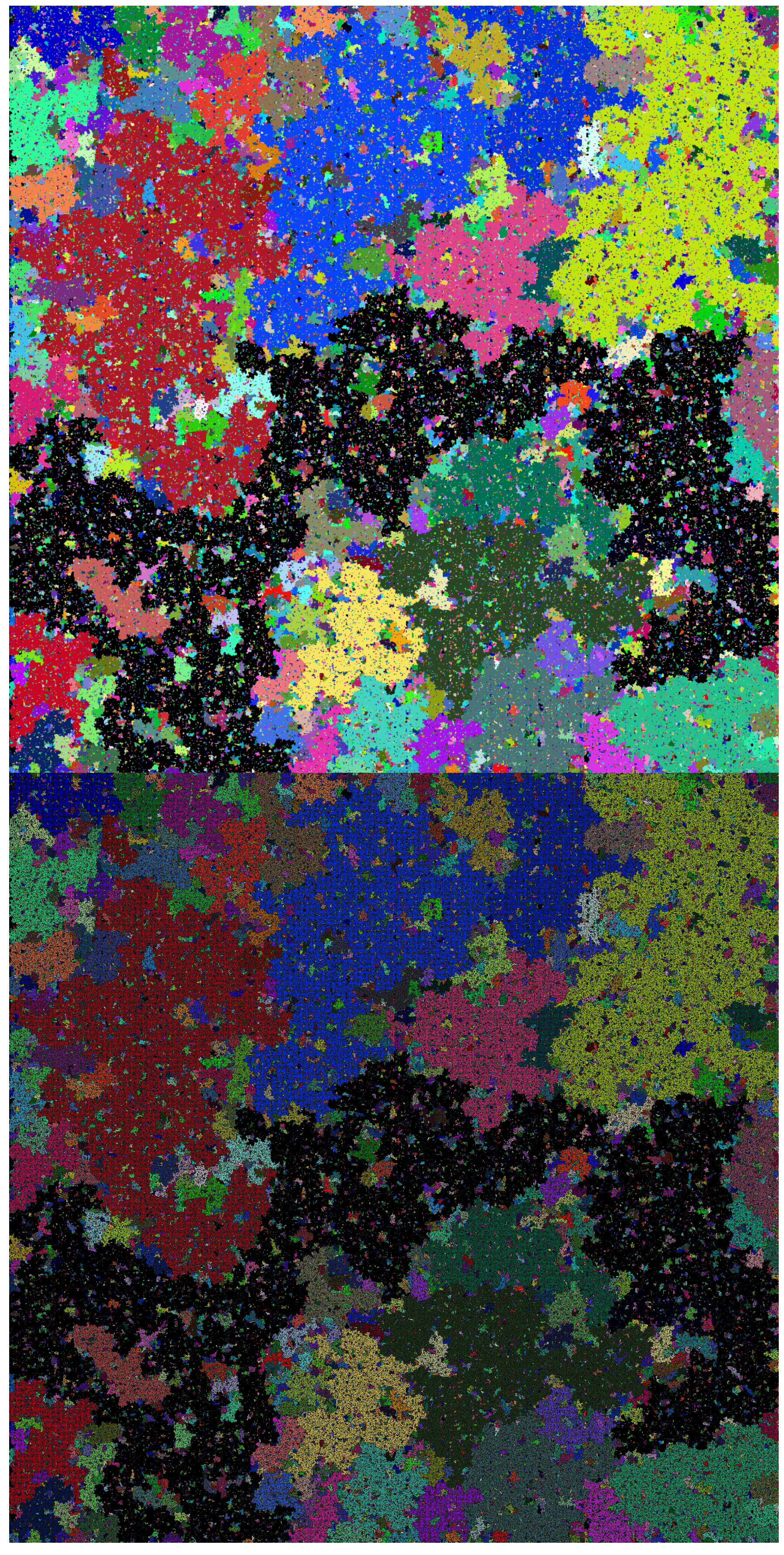

Figure 37: Top: colored vertices only. Bottom: colored vertices and bonds. Bond percolation, $\Lambda(1024)$, $p = 0.5$.

# 28 Adaptation to bond percolation

As a preliminary step, we shall work in the framework of the bond model to build the intertwined explorations and to prove several interesting results on these processes. Before starting this program, we explain how the results presented so far for the site model can be adapted to the bond model. The adaptation of the estimates for genuine exploration algorithms is quite straightforward, they are presented in subsection 28.1. The definitions of the travel times and the shells are more natural in the bond model, they are presented in subsection 28.2.

## 28.1 Genuine bond exploration algorithms

Let us explain briefly how the probabilistic estimates of section 13 can be adapted to the bond model. We still consider a sequence of i.i.d. Bernoulli random variables $(X_k)_{k\geq 1}$ with parameter $p$. This sequence will now be used to decide the states of the bonds. We introduce the associated partial sums $(S_k)_{k\geq 1}$, defined as

$$\forall k\geq 1 \qquad S_k = X_1+\cdots+X_k\,.$$

We consider a finite subset $D$ of $\mathbb{Z}^d$ and a genuine exploration algorithm in $D$. The main difference with the site model is that the three sets maintained by the algorithm will be sets of bonds instead of sets of sites. So, at step $k\geq 0$, we will have:

- The set $E(k)$: these are the bonds which have been explored.
- The set $O(k)$: these are the open bonds which have been explored.
- The set $C(k)$: these are the closed bonds which have been explored.

Initially, we have $E(0)=O(0)=C(0)=\varnothing$. As for the site model, we will have

$$\forall k\geq 1 \qquad \big|E(k)\big| \,=\, k\,,\quad \big|O(k)\big| \,=\, S_k\,,\quad \big|C(k)\big| \,=\, k-S_k\,.$$

The definition of the algorithm rests on a decision rule for choosing the next bond to be explored at step $k$. We use the same formalism as for the site model. Let $(U_k)_{k\geq 0}$ be a sequence of i.i.d. random variables, which are uniformly distributed over $[0,1]$, and which is independent of the percolation process. The decision rule at step $k$ is given by a map

$$\phi_k : \mathbb{E}^d(D)^k\times\{\,0,1\,\}^k\times[0,1]\to\mathbb{E}^d(D)\,.$$

For the first step, the map $\phi_0$ is simply a map from $[0,1]$ to $D$ and $e_1=\phi_0(U_0)$ is the first bond that the algorithm will explore. After completion of the step $k$, the algorithm will have explored the set of bonds $E(k)=\{\,e_1,\dots,e_k\,\}$. The states of these bonds is given by the random variables $X_1,\dots,X_k$, so that

$$O(k) \,=\, \bigcup_{1\leq i\leq k, X_i=1}\{\,e_i\,\}\,,\qquad C(k) \,=\, \bigcup_{1\leq i\leq k, X_i=0}\{\,e_i\,\}\,.$$

The bond explored at step $k+1$ is then given by

$$e_{k+1} \,=\, \phi_k\big(e_1,\dots,e_k,X_1,\dots,X_k,U_k\big)\,.$$

For $E$ a finite subset of $\mathbb{E}^d$, we define

$$S(E) \;=\; \frac{1}{p}\Big|\big\{\, e\in E : e \text{ is open}\,\big\}\Big| - \frac{1}{1-p}\Big|\big\{\, e\in E : e \text{ is closed}\,\big\}\Big| \,. \tag{28.1}$$

The map $S$ is an additive function on the finite sets of edges, that is, if $E, F$ are two disjoint finite subsets of $\mathbb{E}^d$, then

$$S(E\cup F) \;=\; S(E)+S(F)\,.$$

We denote by $T$ the random termination time of the algorithm. The termination time $T$ will be a stopping time with respect to the sequence of random variables $(X_k)_{k\geq 1}$, i.e., for any $t\geq 0$, the event $\{\,T=t\,\}$ is in the $\sigma$-field generated by the variables $X_1,\dots,X_t$. The total set of bonds explored by the algorithm is $\mathcal{E} = E(T)$. Proceeding as in the case of the site model, we obtain the following control on the tail distribution of the random variable $S(\mathcal{E})$:

$$\forall\lambda\geq 0 \qquad P\big(\big|S(\mathcal{E})\big|\geq\lambda\big) \;\leq\; 2|\mathbb{E}^d(D)|\exp\Big(-\frac{2(p(1-p))^2}{|\mathbb{E}^d(D)|}\lambda^2\Big)\,. \tag{28.2}$$

The proof is exactly the same as for the site model. The only noticeable difference with the inequality (13.5) for the site model is that the value $|D|$ has been replaced by $|\mathbb{E}^d(D)|$. Of course, if we had a quantitative control on the termination time $T$, then we could improve considerably the estimate (28.2).

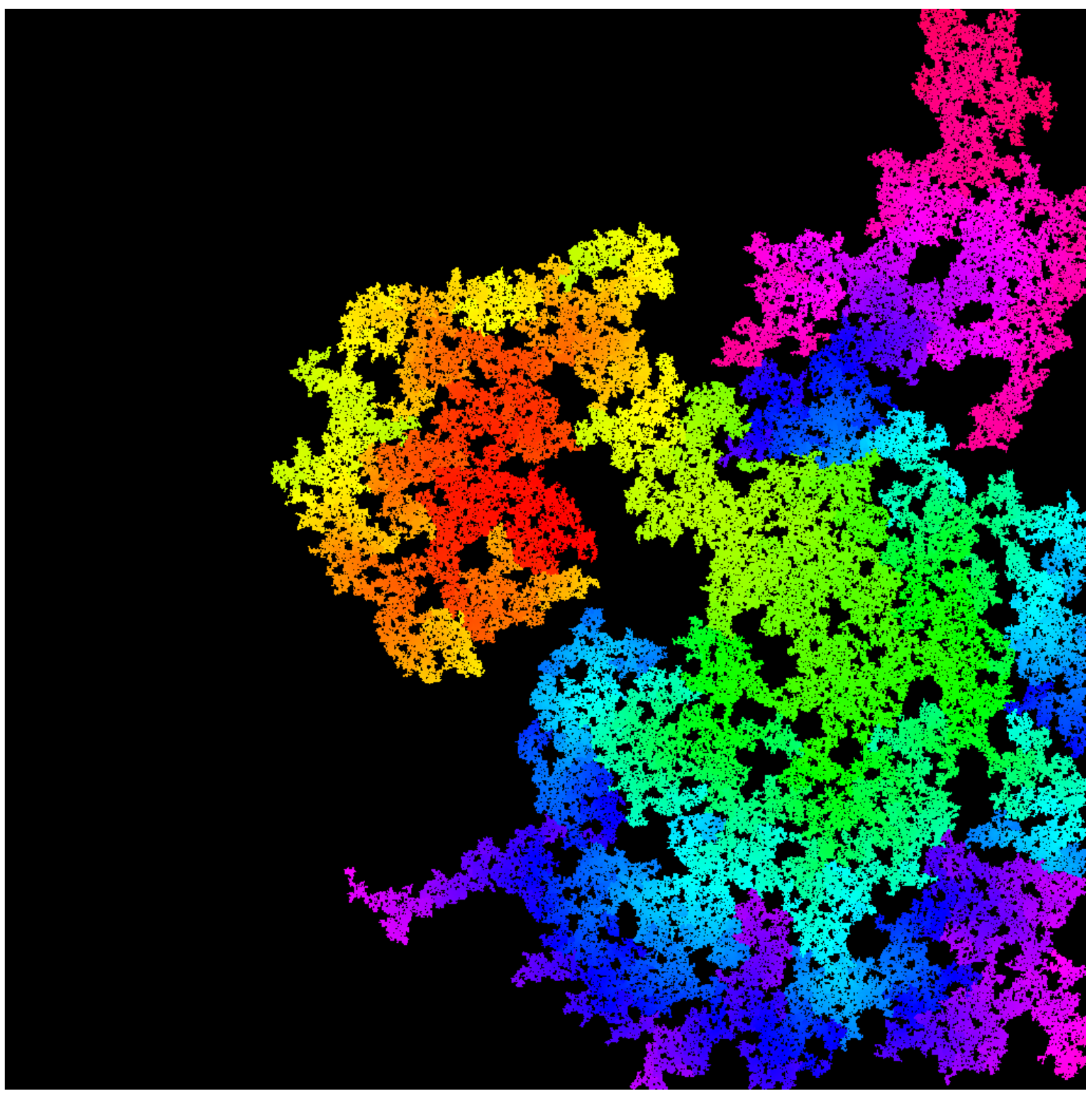

Figure 38: Bond percolation, $\Lambda(1024)$, $p=0.5$, a genuine exploration of $C(0)$. The vertex visited at step $k$ is attributed the rainbow color number $k$, 0 is red.

## 28.2 Travel time, balls and shells for bond percolation

We adapt here the definitions of the travel time and the shells for the bond percolation model. Let $D$ be a subset of $\mathbb{Z}^d$. The travel time $T_D(x,y)$ between $x$ and $y$ in $D$ is defined by

$$T_D(x,y) \,=\, \inf \Big\{ \sum_{i=0}^{r-1} 1_{e_i \text{ closed}} : \\ z_0, e_0, z_1, e_1, \cdots, e_{r-1}, z_r \text{ path in } D \text{ from } z_0 = x \text{ to } z_r = y \Big\}. \quad (28.3)$$

We make the convention that $T_D(x,x) = 0$ for any $x \in D$ (this is compatible with the previous definition, if we consider that $z_0 = x$ is a path from $x$ to $x$; the length $r$ of this path is 0 and the sum in (28.3) is empty). As in the case of site percolation, we have

$$\forall x,y \in D \qquad T_D(x,y) \,\leq\, d_D(x,y)\,,$$

where $d_D$ is the graph distance defined by

$$d_D(x,y) \,=\, \inf \Big\{ r \geq 0 : z_0, \dots, z_r \text{ path in } D \text{ from } z_0 = x \text{ to } z_r = y \Big\}. \quad (28.4)$$

For $x$ in $D$ and $B$ a subset of $D$, we define the travel time between $x$ and $B$ in $D$ by

$$T_D(x,B) \,=\, \inf \big\{ T_D(x,y) : y \in B \big\}.$$

For $A$, $B$ two subsets of $D$, we define the travel time between $A$ and $B$ in $D$ by

$$T_D(A,B) \,=\, \inf \big\{ T_D(x,y) : x \in A,\, y \in B \big\}. \quad (28.5)$$

Let us remark that the travel time $T_D(x,y)$ for the bond model is symmetric in $x$ and $y$. In particular, the travel time $T_D(A,B)$ defined in (28.5) is also symmetric. This is one of the reasons why the bond model is more pleasant to analyze than the site model. We proceed now to define the shells for the bond model. For $A, D$ two finite subsets of $\mathbb{Z}^d$, we define the shell around $A$ in $D$ by

$$\text{Shell}\,(A,D) \,=\, \big\{ x \in D : T_D(A,x) = 0 \big\}.$$

The set $\text{Shell}\,(A,D)$ is the union of the clusters of the percolation configuration restricted to $D$ which intersect $A$. In case $A$ is not a subset of $D$, only the vertices of $A \cap D$ are taken into account, thus we have

$$\text{Shell}\,(A,D) \,=\, \text{Shell}\,(A \cap D, D)\,. \quad (28.6)$$

For $t \geq 0$ an integer, we define the closed ball of radius $t$ around $A$ for the travel time $T_D$ by

$$\mathcal{B}(A,D,t) \,=\, \big\{ x \in D : T_D(A,x) \leq t \big\},$$

and the shell number $t$ around $A$ in $D$ by

$$\text{Shell}\,(A,D,t) \,=\, \big\{ x \in D : T_D(A,x) = t \big\}. \quad (28.7)$$

For any set $A$, we have $\mathcal{B}(A,D,0) \,=\, \text{Shell}\,(A,D) \,=\, \text{Shell}\,(A,D,0)$ and

$$\text{Shell}\,(A,D,1) \,=\, \text{Shell}\,\big(\partial_D^{out}\text{Shell}\,(A,D),D\big)\,. \tag{28.8}$$

Another direct consequence of the definition (28.7) is that the successive shells around $A$ are pairwise disjoint and

$$\forall t\geq 0 \qquad \mathcal{B}(A,D,t) \,=\, \bigcup_{0\leq s\leq t} \text{Shell}\,(A,D,s)\,.$$

Moreover, the balls for the graph distance $d_D$ (defined in (28.4)) are contained in the balls for the travel time $T_D$, in the following sense:

$$\forall x\in D \quad \forall t\geq 0 \qquad \big\{\,y\in D: d_D(x,y)\leq t\,\big\}\,\subset\,\mathcal{B}(\{\,x\,\},D,t)\,.$$

We have also the semigroup identity. Let $A$ be any subset of $D$. We have an identity for balls, which reads

$$\forall s\geq 0 \quad \forall t\geq 0 \qquad \mathcal{B}(A,D,s+t) \,=\, \mathcal{B}\big(\mathcal{B}(A,D,s),D,t\big)\,,$$

and an identity for shells, which is

$$\forall s\geq 0 \quad \forall t\geq 1 \qquad \text{Shell}\,(A,D,s+t) \,=\, \text{Shell}\,\big(\mathcal{B}(A,D,s),D,t\big)\,. \tag{28.9}$$

Notice that the identity (28.9) for shells holds for $t = s = 0$, whereas the corresponding identity (12.10) for the site model was not valid with $t = s = 0$. The semigroup identity yields the following equation for the travel time:

$$\begin{aligned} &\forall A,B\subset D \quad \forall t\in\{\,0,\dots,T_D(A,B)-1\,\} \\ &\qquad\qquad\qquad T_D(A,B) \,=\, t+T_D(\text{Shell}\,(A,D,t),B)\,. \end{aligned}$$

Finally, if we combine the semigroup identity (28.9) and formula (28.8), we obtain that, for any $t\geq 0$,

$$\begin{aligned} \text{Shell}\,(A,D,t+1) \,&=\, \text{Shell}\,\big(\mathcal{B}(A,D,t),D,1\big) \,=\, \text{Shell}\,\big(\partial_D^{out}\mathcal{B}(A,D,t),D\big) \\ &=\, \text{Shell}\,\big(\partial_D^{out}\text{Shell}\,(A,D,t),D\big)\setminus\mathcal{B}(A,D,t)\,. \end{aligned} \tag{28.10}$$

Yet there is something quite different for the shells in the site and the bond models. Let us take $D=\mathbb{Z}^d$ and let $x$ be a vertex in $\mathbb{Z}^d$. For the bond model, the open cluster of $x$ coincides with the 0-th shell of $\{\,x\,\}$, i.e., we have

$$\text{Shell}\,(\{\,x\,\},\mathbb{Z}^d,0) \,=\, C(x)\,,$$

but this is not true for the site model! Indeed, in the site model, the set $\text{Shell}\,(\{\,x\,\},\mathbb{Z}^d,0)$ contains also the closed sites belonging to the outer boundary of $C(x)$. The same problem arises when we start from a finite set of sites $A$. We could have fixed this annoying detail by taking into account the states of the two extremities of the path when defining the travel time $T_D(x,y)$ in the site model. Yet in doing so, we would have lost the nice identity (12.11). We conclude this section with a technical lemma that will play an important role for analyzing the intertwined explorations.

**Lemma 28.1.** *Let $A, E, F$ be finite subsets of $\mathbb{Z}^d$ such that $A \cap E = A \cap F$ and $E \subset F$. If $t$ is such that $T_F(A, F \setminus E) > t \geq 0$, then $\mathcal{B}(A,E,t) = \mathcal{B}(A,F,t)$. Moreover, we have*

$$\forall s \in \{0, \dots, t\} \qquad \mathit{Shell}(A,E,s) = \mathit{Shell}(A,F,s). \tag{28.11}$$

*Proof.* We proceed by induction over $t$ to prove that $\mathcal{B}(A,E,t) = \mathcal{B}(A,F,t)$. We start with $t = 0$. We have

$$\mathcal{B}(A,E,0) = \text{Shell}(A,E), \qquad \mathcal{B}(A,F,0) = \text{Shell}(A,F).$$

The inclusion $\text{Shell}(A,E) \subset \text{Shell}(A,F)$ holds as soon as $E \subset F$. Conversely, let $x$ be a vertex in $\text{Shell}(A,F)$. If $x$ belongs to $A$, then $x$ is $A \cap F = A \cap E$, so that $x$ is in $\text{Shell}(A,E)$ as well. If $x$ does not belong to $A$, then there exists an open path in $F$ connecting $x$ to a vertex $y$ in $A \cap F = A \cap E$. It follows from the hypothesis that all the vertices of this path must belong to $E$, so that $x$ is also in $\text{Shell}(A,E)$. Suppose now that the result has been proved at rank $t - 1 \geq 0$. The inclusion $E \subset F$ implies that $T_E \geq T_F$, which in turn yields the inclusion $\mathcal{B}(A,E,t) \subset \mathcal{B}(A,F,t)$. Conversely, let $x$ belong to $\mathcal{B}(A,F,t)$. If $x$ is in $\mathcal{B}(A,F,t-1)$, then, by the induction hypothesis, it is also in $\mathcal{B}(A,E,t-1)$, which is included in $\mathcal{B}(A,E,t)$. Suppose that $x$ is in

$$\mathcal{B}(A,F,t) \setminus \mathcal{B}(A,F,t-1) = \text{Shell}\big(\mathcal{B}(A,F,t-1), F, 1\big).$$

Using the induction hypothesis and formula (28.8), we have

$$\text{Shell}\big(\mathcal{B}(A,F,t-1), F, 1\big) = \text{Shell}\big(\partial_F^{out}\mathcal{B}(A,E,t-1), F\big).$$

Since $x$ belongs to this last set, then there exists a path $z_0, e_0, z_1, \dots, e_{r-1}, z_r$ in $F$ joining a vertex $z_0$ of $\partial_F^{out}\mathcal{B}(A,E,t-1)$ to $z_r = x$ such that the edges $e_0, \dots, e_{r-1}$ are open. This implies that $T_F(A, z_i) \leq t$ for $0 \leq i \leq r$, and it follows from the hypothesis that the vertices $z_0, \dots, z_r$ all belong to $E$. Therefore $x$ is in $\text{Shell}\big(\partial_E^{out}\mathcal{B}(A,E,t-1), E\big)$, which is included in $\mathcal{B}(A,E,t)$. This concludes the induction step. We prove finally the second assertion (28.11). We remark first that any integer $s$ such that $0 \leq s \leq t$ satisfies the hypothesis of the lemma, therefore we have also that

$$\forall s \in \{0, \dots, t\} \qquad \mathcal{B}(A,E,s) = \mathcal{B}(A,F,s). \tag{28.12}$$

The claim (28.11) on the equality of the shells is a consequence of the equalities of the balls stated in (28.12). $\square$

Once the travel times and the shells are properly defined, all the results proved within the framework of the site model in section 14 also apply to the bond model, with only minor modifications. The proofs are almost identical: it is simply a matter of adjusting certain geometric constants that appear in several places.

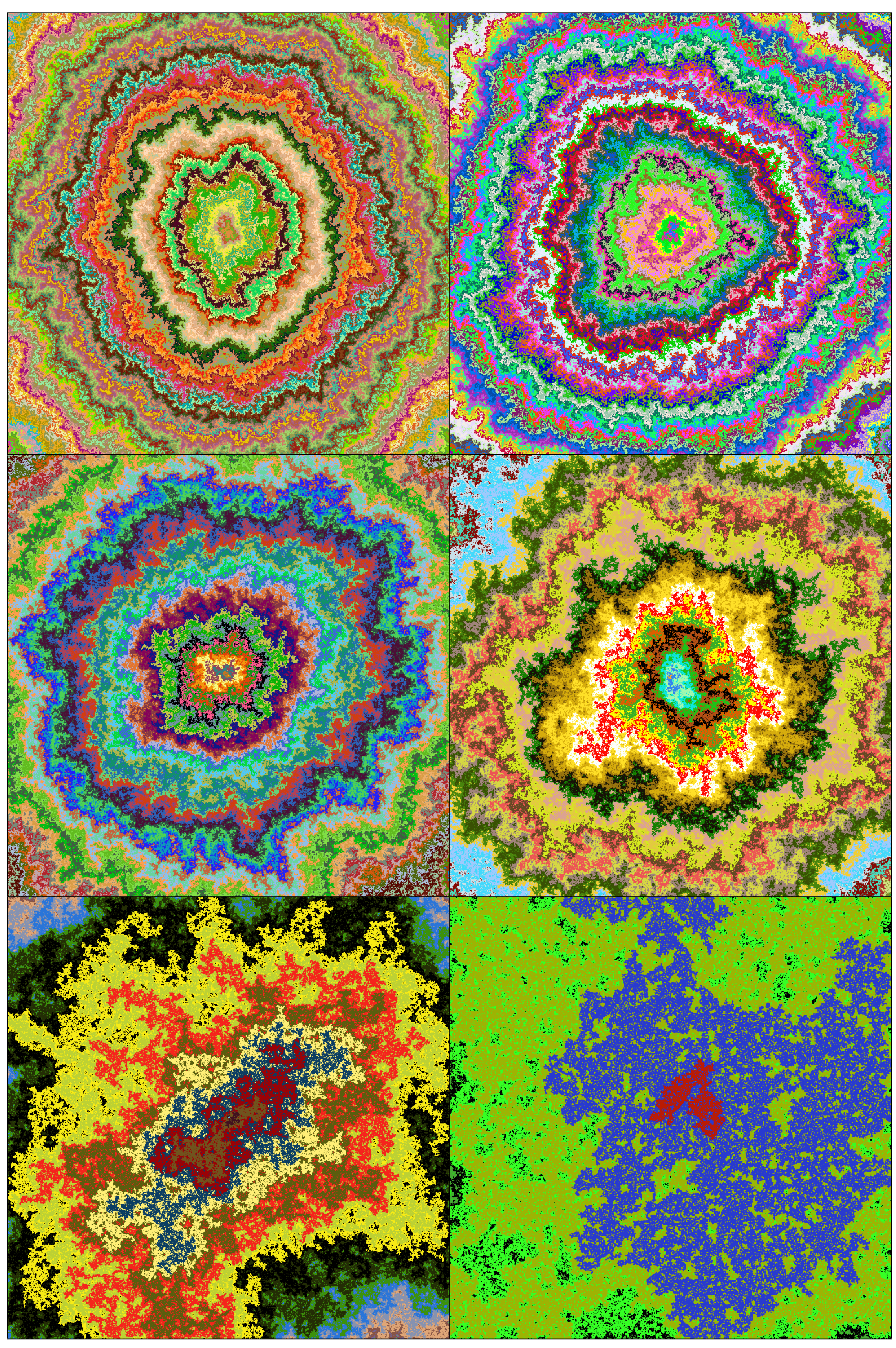

Figure 39: Shells around 0 for bond percolation in a box of size $1024 \times 1024$ with parameters $p = 0.4, 0.42, 0.44, 0.46, 0.48, 0.5$ (left to right, top to bottom).

# 29 Intertwined bond explorations

In this section, we shall at last define the intertwined bond explorations. These are interacting explorations involving two explorers. Each explorer maintains an object called a waiting bond aggregate, defined in subsection 29.1. The two explorers are called in turn to explore the configuration. During each round of exploration, they run a taboo exploration algorithm, as defined in subsection 29.2, which is encoded in the map Coherent_Bond_Explore defined in subsection 29.3. Later on, we will develop mathematical proofs based on these interacting explorations. So we will take the luxury to furnish a full algorithmic version, together with its mathematical counterpart (see subsection 29.4). These two versions transmit all the information acquired during the exploration process. We define also the intersections associated to the intertwined explorations in subsection 29.5 and we furnish a proof of termination for the intertwined exploration process in subsection 29.6.

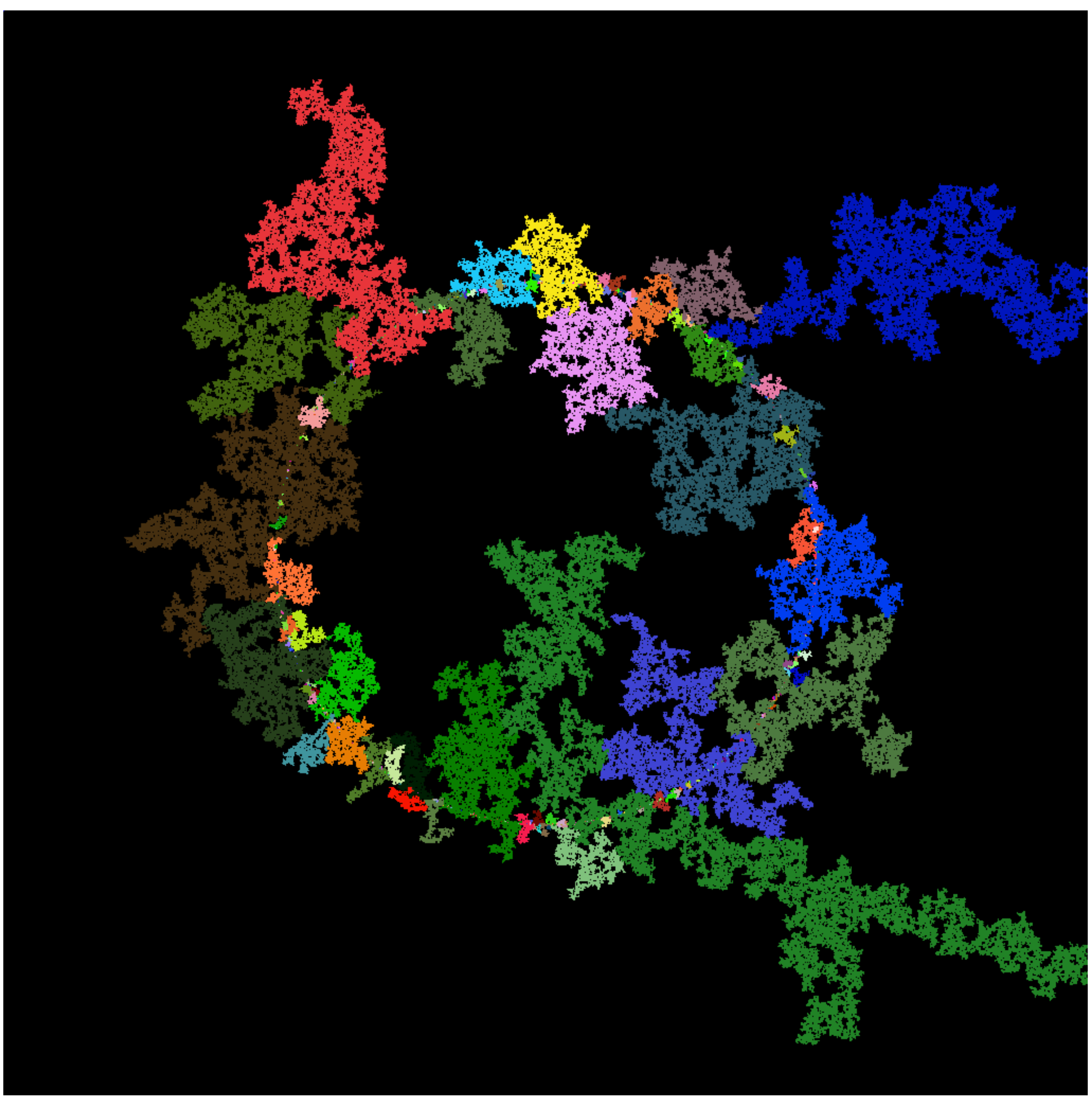

Figure 40: Clusters of $C\big(\mathrm{circle}(0,256),\Lambda(1024)\big)$, bond percolation at $p=0.49$.

## 29.1 Waiting bond aggregates

The taboo exploration algorithm will involve objects called waiting aggregates, that we define in this subsection. Let $D$ be a finite subset of $\mathbb{Z}^d$. For $x \in D$, the open cluster of $x$ in $D$ is the set

$$C(x, D) = \big\{ y \in D : x \longleftrightarrow y \text{ in } D \big\}.$$

A subset $\mathcal{A}$ of $D$ is an aggregate of clusters in $D$, or simply an aggregate, if

$$\forall x \in \mathcal{A} \qquad C(x, D) \subset \mathcal{A}.$$

In words, an aggregate is the union of some clusters of $D$. Obviously, the union of two aggregates is still an aggregate. Also, if $\mathcal{A}$ is an aggregate and $x$ is not in $\mathcal{A}$, then $C(x, D) \subset D \setminus \mathcal{A}$. For $D$ a finite subset of $\mathbb{Z}^d$ and $A$ a subset of $D$, we define

$$C(A, D) = \bigcup_{x \in A} C(x, D).$$

Thus $C(A, D)$ is the union of all the clusters of the vertices of $A$, for the percolation configuration restricted to $D$. In particular, this set is an aggregate of clusters in $D$. A set $\mathcal{W}$ is a waiting set associated to an aggregate $\mathcal{A}$ in a domain $D$ if it is a subset of its outer boundary in $D$, i.e.,

$$\mathcal{W} \subset \partial_D^{out} \mathcal{A}.$$

An aggregate $\mathcal{A}$ endowed with a waiting set $\mathcal{W}$ is called a waiting aggregate. An aggregate corresponds to a region which has been explored. Its associated waiting set corresponds to the vertices which are to be explored when a next round of exploration starts.

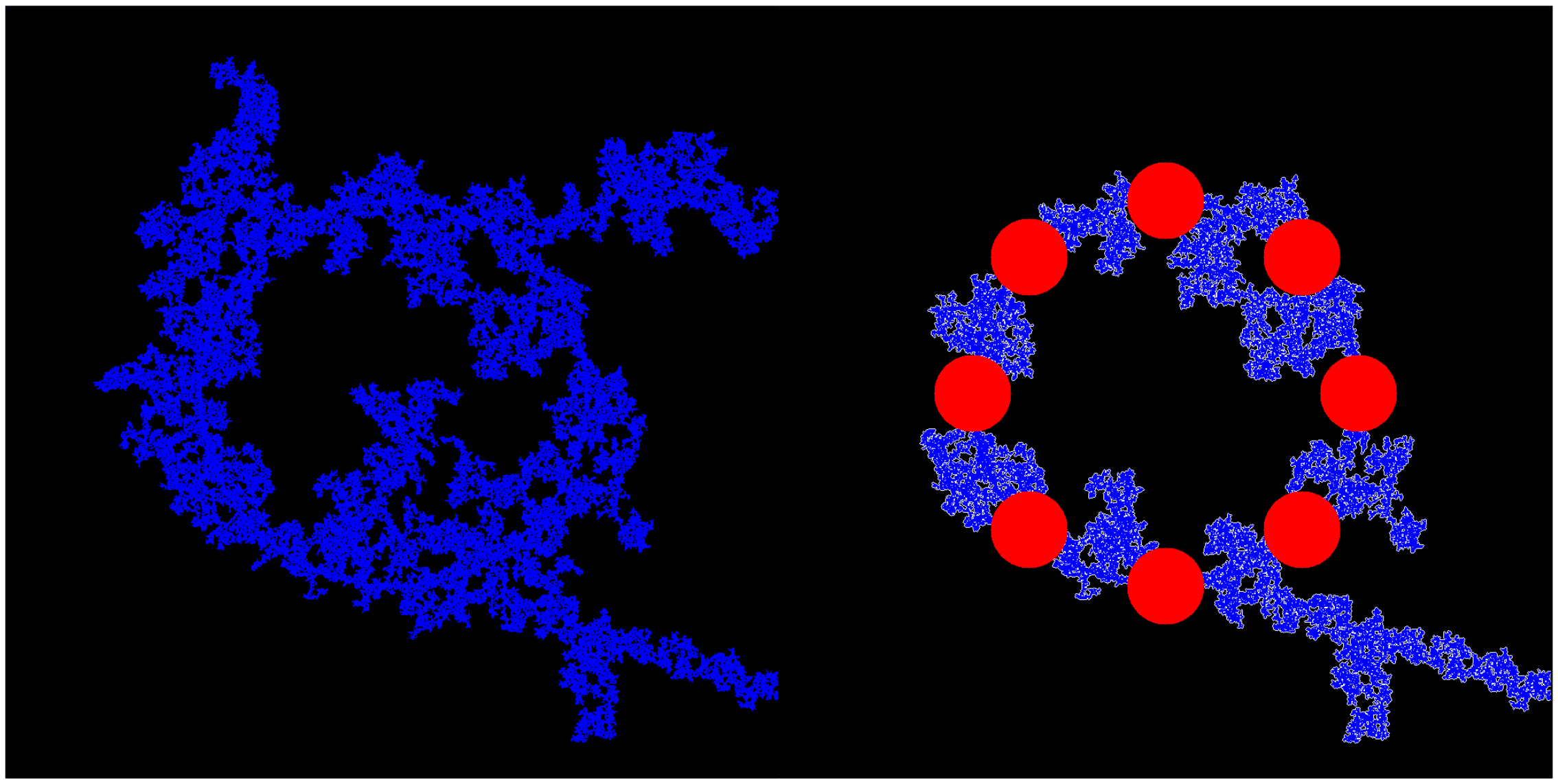

Figure 41: Bond percolation in $D = \Lambda(1024)$ with $p = 0.49$, $A = \text{circle}(0, 256)$, $\mathcal{T}$ = the 8 red discs. Left: the set $C(A, D)$ (the clusters of $C(A, D)$ are visible in figure 40). Right: Clusters $(A, D, \mathcal{T})$ is in blue, Waiting $(A, D, \mathcal{T})$ is in white.

## 29.2 Taboo bond exploration algorithm

We start from the standard exploration algorithm, and we introduce a taboo set $\mathcal{T}$, which the algorithm is not permitted to enter during the exploration process. More precisely, the algorithm will discard from the set of active vertices the vertices belonging to $\mathcal{T}$. However, the algorithm might examine a bond emanating from the set of active vertices and ending into $\mathcal{T}$. The pseudocode of the taboo exploration algorithm is presented in algorithm 29.1.

**Algorithm 29.1** Taboo_Bond_Explore_Clusters $(A, D, \mathcal{T})$

```
Require: a non-empty starting set A, an authorized domain D, a taboo set T
Ensure: C is Clusters (A, D, T)
        W is the set of waiting vertices
        E is the set of edges which have been explored
  C ← ∅, W ← ∅, E ← ∅
  repeat
    Pick up x ∈ A
    A ← A \ {x}                          ▷ The vertex x is discarded from A
    if x is in D \ T then
      C ← C ∪ {x}
      for y ∈ N(x) ∩ D do
        e ← ⟨x, y⟩
        if e ∉ E then        ▷ Check that e has not been previously examined
          E ← E ∪ {e}               ▷ The bond e is examined even if y ∈ T
          if y ∉ T then
            if e is closed then
              W ← W ∪ {y}   ▷ The vertex y is put on the waiting list
            else if e is open then
              if y ∉ C then        ▷ Check that y has not been explored
                A ← A ∪ {y}
              end if
            end if
          end if
        end if
      end for
    end if
  until A = ∅
  W ← W \ C              ▷ Remove the explored vertices from the waiting set
  return C, W, E
```

We give now the corresponding mathematical definitions. Let $A$,$\mathcal{T}$ be two subsets of $D$. Upon termination, the taboo exploration algorithm will have explored the edges of $D$ having one endpoint in the set

$$\text{Clusters}\,(A, D, \mathcal{T}) \;=\; C\big(A \cap D \setminus \mathcal{T}, D \setminus \mathcal{T}\big)\,. \tag{29.1}$$

This set is precisely the set of the vertices which are connected to $A \cap D \setminus \mathcal{T}$ by an open path in $D$ staying outside $\mathcal{T}$. It is also the set of the vertices which can be reached from $A \setminus \mathcal{T}$ with a null travel time in $D$. Notice that $\text{Clusters}\,(A, D, \mathcal{T})$ is non-increasing with respect to $\mathcal{T}$, i.e.,

$$\mathcal{T}_1 \subset \mathcal{T}_2 \qquad \Longrightarrow \qquad \text{Clusters}\,(A, D, \mathcal{T}_2) \subset \text{Clusters}\,(A, D, \mathcal{T}_1)\,.$$

If we give an empty taboo set to the algorithm, it will simply explore all the clusters in $D$ which intersect $A$, so that $\text{Clusters}\,(A, D, \varnothing) = C(A \cap D, D)$. During the exploration, the algorithm will also record a set of waiting vertices, which is

$$\text{Waiting}\,(A, D, \mathcal{T}) \;=\; \partial^{out}_{D\setminus\mathcal{T}} C\big(A \cap D \setminus \mathcal{T}, D \setminus \mathcal{T}\big)\,. \tag{29.2}$$

In principle, to build the sets $\text{Clusters}\,(A, D, \mathcal{T})$ and $\text{Waiting}\,(A, D, \mathcal{T})$, it is not necessary to examine the states of the bonds emanating from $\text{Clusters}\,(A, D, \mathcal{T})$ which lead to a vertex in the taboo set $\mathcal{T}$, because $\text{Clusters}\,(A, D, \mathcal{T})$ and $\text{Waiting}\,(A, D, \mathcal{T})$ are disjoint from $\mathcal{T}$ (notice that we used the outer boundary in $D \setminus \mathcal{T}$ to define $\text{Waiting}\,(A, D, \mathcal{T})$). However, for a reason that will appear later, the taboo exploration will examine these bonds, and they are incorporated in the set of bonds explored by the algorithm. The next round of exploration will start from the set $\text{Waiting}\,(A, D, \mathcal{T})$. More precisely, the algorithm will look at all the bonds which have one endpoint belonging to $\text{Waiting}\,(A, D, \mathcal{T})$ and which have not been examined previously. The algorithm 29.1 returns in addition the set $E$ of the edges it has explored, that is

$$\text{Explored}\,(A, D, \mathcal{T}) \;=\; \mathbb{E}^d\big(\text{Clusters}\,(A, D, \mathcal{T})\big) \cup \Delta_D \text{Clusters}\,(A, D, \mathcal{T})\,.$$

The fact that $\text{Explored}\,(A, D, \mathcal{T})$ contains $\Delta_D \text{Clusters}\,(A, D, \mathcal{T}) \cap \Delta_D \mathcal{T}$ results from a deliberate choice in the design of the taboo exploration algorithm. Indeed, we could define an algorithm which does not examine the bonds having an endpoint in $\mathcal{T}$, in which case the explored bonds would be the smaller set

$$\mathbb{E}^d\big(\text{Clusters}\,(A, D, \mathcal{T}) \cup \text{Waiting}\,(A, D, \mathcal{T})\big)\,.$$

The reason for our choice will become clear later on. For the time being, let us say loosely that we wish to design genuine exploration algorithms (or at least partly genuine) so that the intersection of their respective sets of explored bonds is a set of specific closed bonds that plays a crucial role in some disconnection events.

Notice also that the definitions of $\text{Clusters}\,(A, D, \mathcal{T})$ and $\text{Waiting}\,(A, D, \mathcal{T})$ depend on $D$ and $\mathcal{T}$ only through their difference $D \setminus \mathcal{T}$. Thus we could in principle give a simpler definition involving only the set $D \setminus \mathcal{T}$ instead of the two sets $D, \mathcal{T}$. There are two reasons for not doing so. The first reason is that the set $\text{Explored}\,(A, D, \mathcal{T})$ cannot be defined in this way. The second reason is that it would obscure the formulas, because the sets $D$ and $\mathcal{T}$ do not play the same role.

The pair of sets $\big(\text{Clusters}\,(A, D, \mathcal{T}), \text{Waiting}\,(A, D, \mathcal{T})\big)$ is a waiting aggregate in $D \setminus \mathcal{T}$. An example of the sets $\text{Clusters}\,(A, D, \mathcal{T})$ and $\text{Waiting}\,(A, D, \mathcal{T})$ is presented in figure 41. Notice that the taboo set has the effect of disconnecting some open clusters included in the set $C(A, D)$ (the open clusters included in $C(A, D)$ can be seen in figure 40).

## 29.3 Coherence bond exploration

The intertwined exploration processes which are to be defined shortly will involve pairs of aggregates satisfying certain conditions. We call such pairs coherent, and we define them next.

**Definition 29.1.** *Two waiting aggregates* $(\mathcal{A},\mathcal{W}),(\mathcal{A}',\mathcal{W}')$ *are coherent if*
- *the sets* $\mathcal{A}$ *and* $\mathcal{A}'$ *are disjoint;*
- *we have the inclusion* $\partial_D^{out}\big(\mathcal{A}\cup\mathcal{A}'\big)\,\subset\,\mathcal{W}\cup\mathcal{W}'$ *;*
- *either* $\mathcal{W}'\cap\mathcal{A}=\varnothing$ *or* $\mathcal{W}\cap\mathcal{A}'=\varnothing$*.*

The state of the intertwined explorations after each round of exploration will be encoded in a pair of coherent waiting aggregates. The aggregates $\mathcal{A}$ and $\mathcal{A}'$ will correspond to the regions explored by each explorer, this is why they are disjoint. The sets $\mathcal{W}$ and $\mathcal{W}'$ will correspond to the waiting sets of each explorer. The second condition will ensure that the successive shells around the starting sets are explored in their natural increasing order. The third condition is more tricky, it is a consequence of the construction of the intertwined explorations.

We proceed now to the definition of the map Coherent_Bond_Explore, which is the mathematical object corresponding to one round of exploration by one of the two explorers.

**Definition 29.2.** *Suppose that we are given two coherent waiting aggregates* $(\mathcal{A},\mathcal{W})$, $(\mathcal{A}',\mathcal{W}')$ *in* $D$*. We define a map Coherent_Bond_Explore on the pairs of coherent waiting aggregates by setting*

$$\begin{multline}\textit{Coherent_Bond_Explore}\Big((\mathcal{A},\mathcal{W}),(\mathcal{A}',\mathcal{W}')\Big)\,=\\ \Big(\mathcal{A}\cup \textit{Clusters}\,(\mathcal{W},D\setminus\mathcal{A},\mathcal{A}'),\,\textit{Waiting}\,(\mathcal{W},D\setminus\mathcal{A},\mathcal{A}')\Big)\,. \tag{29.3}\end{multline}$$

Let us introduce the notation

$$(\mathcal{A}'',\mathcal{W}'')\;=\;\text{Coherent_Bond_Explore}\,\big((\mathcal{A},\mathcal{W}),(\mathcal{A}',\mathcal{W}')\big)\,.$$

The definition of the sets $\mathcal{A}''$, $\mathcal{W}''$ can be rewritten as follows:

$$\mathcal{A}''\;=\;\mathcal{A}\cup C\big(\mathcal{W}\setminus\mathcal{A}',D\setminus(\mathcal{A}\cup\mathcal{A}')\big)\,, \tag{29.4}$$
$$\mathcal{W}''\;=\;\partial^{out}_{D\setminus(\mathcal{A}\cup\mathcal{A}')}C\big(\mathcal{W}\setminus\mathcal{A}',D\setminus(\mathcal{A}\cup\mathcal{A}')\big)\,. \tag{29.5}$$

Three examples of the map Coherent_Bond_Explore are presented in picture 42, with the following code of colors:

$$\mathcal{A}=C(\text{white square},D),\quad \mathcal{W}=\partial^{\,out}\mathcal{A},\quad \mathcal{A}'=C(\text{gray square},D),$$
$$\mathcal{W}'=\partial^{\,out}\mathcal{A}'\setminus\mathcal{A},\quad \mathcal{W}\cap\mathcal{A},\quad \mathcal{W}'\cap\mathcal{A}'',$$
$$\mathcal{A}''\setminus\mathcal{A}=\text{Clusters}\,(\mathcal{W},D\setminus\mathcal{A},\mathcal{A}'),\quad \mathcal{W}''=\text{Waiting}\,(\mathcal{W},D\setminus\mathcal{A},\mathcal{A}')\,.$$

Formula (29.4) results from a direct application of the definitions (29.1) and (29.3), while formula (29.5) comes directly from the definitions (29.2) and (29.3). We have also used the fact that $\mathcal{W}\setminus\mathcal{A}'\subset D\setminus(\mathcal{A}\cup\mathcal{A}')$ to simplify the above formulas.

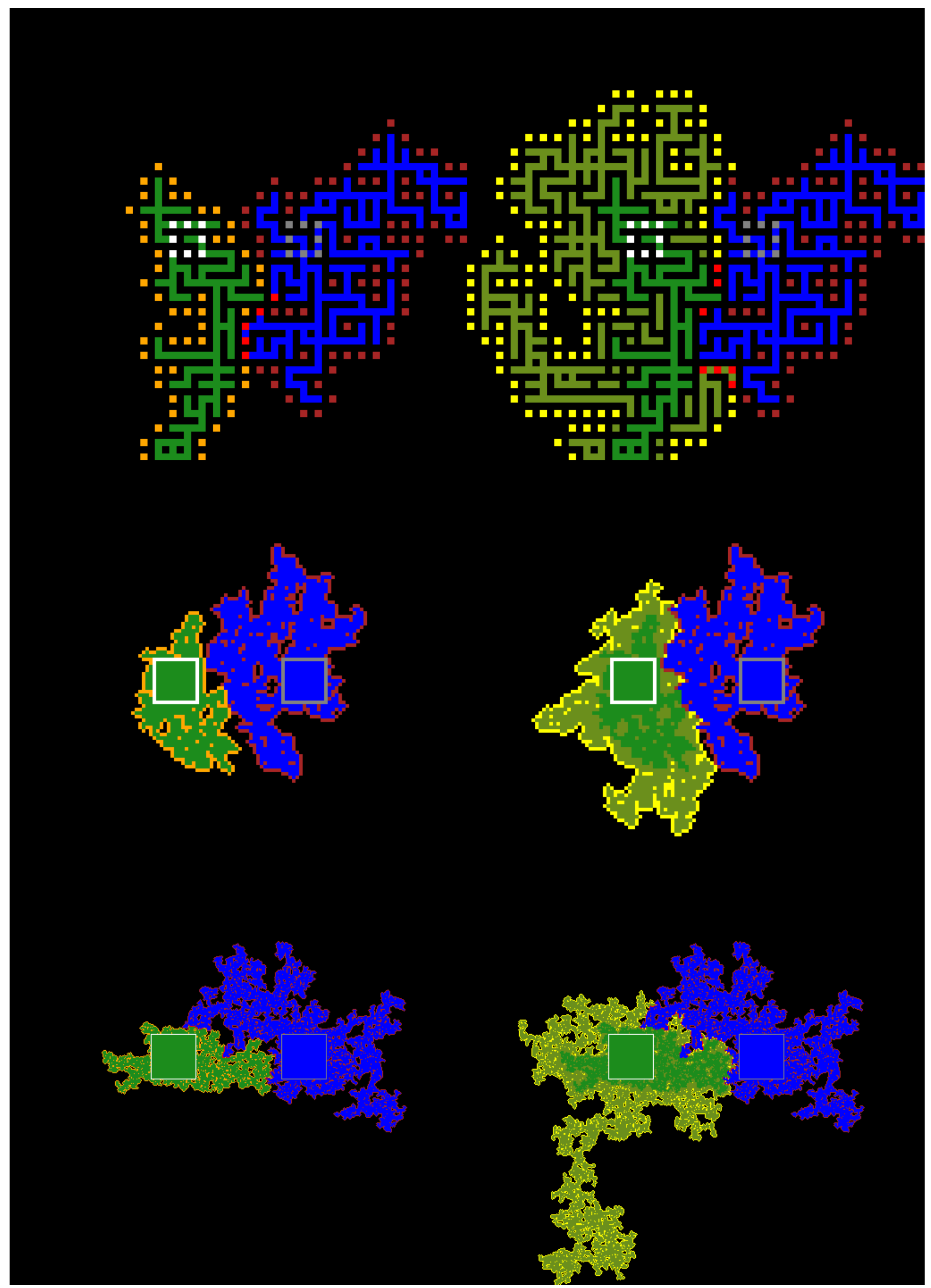

Figure 42: Three examples of Coherent_Bond_Explore, one example per row. The input is on the left and the output on the right, and $D, p$ are equal to: $\Lambda(32), 0.4042$ (top row), $\Lambda(128), 0.45$ (middle row), $\Lambda(512), 0.4812$ (bottom row). $\mathcal{A} = C(\text{white } \square, D)$, $\mathcal{W} = \partial^{out}\mathcal{A}$, $\mathcal{A}' = C(\square, D)$, $\mathcal{W}' = \partial^{out}\mathcal{A}' \setminus \mathcal{A}$, $\mathcal{W} \cap \mathcal{A}'$ $\mathcal{W}' \cap \mathcal{A}''$, $\mathcal{A}'' \setminus \mathcal{A} = \text{Clusters}\,(\mathcal{W}, D \setminus \mathcal{A}, \mathcal{A}')$, $\mathcal{W}'' = \text{Waiting}\,(\mathcal{W}, D \setminus \mathcal{A}, \mathcal{A}')$.

**Proposition 29.3.** *Suppose that* $(\mathcal{A},\mathcal{W})$, $(\mathcal{A}',\mathcal{W}')$ *are two coherent waiting aggregates in* $D$. *Then*

$$(\mathcal{A}'',\mathcal{W}'') \;=\; \textit{Coherent_Bond_Explore}\big((\mathcal{A},\mathcal{W}),(\mathcal{A}',\mathcal{W}')\big)$$

*is a waiting aggregate in* $D$ *and moreover* $(\mathcal{A}',\mathcal{W}')$ *and* $(\mathcal{A}'',\mathcal{W}'')$ *are coherent.*

We invite the reader to read the proof with the help of the colors in figure 42.

*Proof.* Let us prove first that $\mathcal{A}''$ is an aggregate in $D$. Let $x$ be a vertex belonging to $\mathcal{A}''$. If $x$ belongs to $\mathcal{A}$, then $C(x,D)$ is included in $\mathcal{A}$ because $\mathcal{A}$ is an aggregate in $D$, and since $\mathcal{A}\subset\mathcal{A}''$, we conclude that $C(x,D)$ is included in $\mathcal{A}''$. Suppose that $x$ is in $\mathcal{A}''\setminus\mathcal{A}=\text{Clusters}\big(\mathcal{W},D\setminus\mathcal{A},\mathcal{A}'\big)$. The set $\mathcal{A}$ is an aggregate in $D$ and $x$ is not in $\mathcal{A}$, therefore $C(x,D)\subset D\setminus\mathcal{A}$. By the definition of Clusters$\big(\mathcal{W},D\setminus\mathcal{A},\mathcal{A}'\big)$ (see formulas (29.1) or (29.4)), there exists $y\in\mathcal{W}\setminus\mathcal{A}'$ such that $y\longleftrightarrow x$ in $D\setminus(\mathcal{A}\cup\mathcal{A}')$ (recall also that $\mathcal{W}\setminus\mathcal{A}'\subset D\setminus(\mathcal{A}\cup\mathcal{A}')$). It follows that

$$C\big(x,D\setminus(\mathcal{A}\cup\mathcal{A}')\big) \;=\; C\big(y,D\setminus(\mathcal{A}\cup\mathcal{A}')\big)\,. \tag{29.6}$$

Since $\mathcal{A}$ and $\mathcal{A}'$ are two aggregates in $D$ and $y\notin\mathcal{A}\cup\mathcal{A}'$, then it follows that $C(y,D)\;\subset\;D\setminus(\mathcal{A}\cup\mathcal{A}')$, and this implies furthermore that

$$C(y,D) \;=\; C\big(y,D\setminus(\mathcal{A}\cup\mathcal{A}')\big)\,. \tag{29.7}$$

The equalities (29.6) and (29.7) imply that $C(x,D)=C(y,D)\subset\mathcal{A}''$, and we conclude that $\mathcal{A}''$ is an aggregate in $D$. From the definitions (29.4) and (29.5), we see that $\mathcal{W}''\subset\partial_D^{out}\mathcal{A}''$, so that $(\mathcal{A}'',\mathcal{W}'')$ is a waiting aggregate in $D$. We prove finally that $(\mathcal{A}',\mathcal{W}')$ and $(\mathcal{A}'',\mathcal{W}'')$ are coherent. Since $(\mathcal{A},\mathcal{W})$ and $(\mathcal{A}',\mathcal{W}')$ are coherent, then $\mathcal{A}\cap\mathcal{A}'=\varnothing$. Moreover, the set Clusters$\big(\mathcal{W},D\setminus\mathcal{A},\mathcal{A}'\big)$ does not intersect $\mathcal{A}'$, therefore $\mathcal{A}''\cap\mathcal{A}'=\varnothing$. It follows directly from the definition (29.3) (or formula (29.5)) that $\mathcal{W}''\cap\mathcal{A}'=\varnothing$. It remains to prove that

$$\partial_D^{out}\big(\mathcal{A}'\cup\mathcal{A}''\big) \;\subset\; \mathcal{W}'\cup\mathcal{W}''\,. \tag{29.8}$$

Let $x$ be a vertex in $\partial_D^{out}\big(\mathcal{A}'\cup\mathcal{A}''\big)$. By definition, there exists $y$ in $\mathcal{A}'\cup\mathcal{A}''$ such that $|x-y|=1$. We consider two cases.

• First case: $y\in\mathcal{A}\cup\mathcal{A}'$. Since $x\notin\mathcal{A}'\cup\mathcal{A}''$, then $x\notin\mathcal{A}\cup\mathcal{A}'$, and we see that $x$ belongs to $\partial_D^{out}\big(\mathcal{A}\cup\mathcal{A}'\big)$, which is itself included in $\mathcal{W}\cup\mathcal{W}'$ because $(\mathcal{A},\mathcal{W})$ and $(\mathcal{A}',\mathcal{W}')$ are coherent. Suppose that $x$ is in $\mathcal{W}$. We would then have $x\in\mathcal{W}\setminus\mathcal{A}'$ and $x\in D\setminus\mathcal{A}$. Yet

$$\big(\mathcal{W}\setminus\mathcal{A}'\big)\cap\big(D\setminus\mathcal{A}\big) \;\subset\; \text{Clusters}\,\big(\mathcal{W},D\setminus\mathcal{A},\mathcal{A}'\big) \;\subset\; \mathcal{A}''\,,$$

and $x$ would belong to $\mathcal{A}''$, which is absurd. Thus $x$ is in $\mathcal{W}'$.

• Second case: $y\in\mathcal{A}''\setminus(\mathcal{A}\cup\mathcal{A}')$. We see from formula (29.4) that $y$ is in

$$C(\mathcal{W}\setminus\mathcal{A}',D\setminus(\mathcal{A}\cup\mathcal{A}'))\,.$$

Since $x$ is not in $\mathcal{A}\cup\mathcal{A}'$, we conclude from formula (29.5) that $x$ belongs to $\mathcal{W}''$. This completes the proof of the inclusion (29.8). □

## 29.4 Two intertwined bond explorations

We now combine two taboo explorations to create the intertwined explorations. Suppose that we are given two sets of vertices $W$ and $W'$ and a percolation configuration in the finite domain $D$ such that $W$ and $W'$ are not connected. We will launch two taboo exploration algorithms, or rather two explorers. The first explorer starts from the set $W$, while the second explorer starts from the set $W'$. Each explorer possesses a personal waiting set and a personal taboo set. The two explorers will be run alternatively and they will interact in a complicated way. In turn, each explorer will complete the exploration of the vertices which are connected to his waiting set by an open path staying outside his taboo set. After each round of exploration of one explorer, the taboo set of the other explorer is updated with some information coming from the exploration just completed. There are various ways to describe formally this process. We would like to have a rigorous formulation which involves only the relevant objects. Mathematically, we construct two sequences of random sets $(\mathcal{A}(t), \mathcal{W}(t), t \geq 0)$ whose meaning is the following:

- $\mathcal{A}(t)$ is the set of vertices explored at round $t$;
- $\mathcal{W}(t)$ is the set of the waiting vertices at round $t$.

These sequences of sets are defined as follows. The initial sets for $t = 0, 1$ are

$$\begin{aligned} \mathcal{A}(0) &= \text{Clusters}\big(W, D, \varnothing\big),\ \mathcal{W}(0) = \text{Waiting}\big(W, D, \varnothing\big)\,, \\ \mathcal{A}(1) &= \text{Clusters}\big(W', D, \mathcal{A}(0)\big),\ \mathcal{W}(1) = \text{Waiting}\big(W', D, \mathcal{A}(0)\big)\,. \end{aligned}$$

Recall that the sets $\text{Clusters}(\cdot,\cdot,\cdot)$ and $\text{Waiting}(\cdot,\cdot,\cdot)$ are defined in formulas (29.1) and (29.2). For $t \geq 2$, we define inductively the sequences by setting

$$\begin{aligned} &\big(\mathcal{A}(t), \mathcal{W}(t)\big) = \\ &\text{Coherent_Bond_Explore}\Big(\big(\mathcal{A}(t-2), \mathcal{W}(t-2)\big), \big(\mathcal{A}(t-1), \mathcal{W}(t-1)\big)\Big)\,. \end{aligned} \tag{29.9}$$

In particular, the increment between $\mathcal{A}(t-2)$ and $\mathcal{A}(t)$ is given by

$$\mathcal{A}(t) \setminus \mathcal{A}(t-2) = \text{Clusters}\big(\mathcal{W}(t-2), D \setminus \mathcal{A}(t-2), \mathcal{A}(t-1)\big)\,. \tag{29.10}$$

The sequence of the sets with even indices $(\mathcal{A}(2t), \mathcal{W}(2t), t \geq 0)$ correspond to the set of vertices discovered by the first explorer and his set of waiting vertices. The sequence of the sets with odd indices $(\mathcal{A}(2t+1), \mathcal{W}(2t+1), t \geq 0)$ are their counterparts for the second explorer. We define the final iteration number $T$ as the smallest index $T$ such that the domain $D$ has been completely visited by the two explorers, i.e.,

$$T = \min\big\{\, t \geq 1 : \mathcal{A}(t-1) \cup \mathcal{A}(t) = D \,\big\}\,. \tag{29.11}$$

The sequence $\mathcal{A}(2t)$ (respectively $\mathcal{A}(2t+1)$) is stationary for $2t \geq T-1$ (resp. $2t+1 \geq T-1$). Indeed, if we apply the map Coherent_Bond_Explore to a pair of aggregates whose union is $D$, the resulting aggregate is equal to the first aggregate. Beware also that the final iteration number $T$ is not a stopping time with respect to the sequence of Bernoulli random variables driving the algorithm. It is simply a random number associated to the intertwined explorations. In fact, each iteration of the intertwined explorations corresponds to a round of exploration of a taboo shell exploration, which is itself driven by the sequence of Bernoulli random variables. Notice also that the waiting set of one explorer might become void well before the other. If that is the case, the explorer is stuck and his sequence of explored vertices becomes stationary, so that he simply waits until the other explorer has an empty set of waiting vertices. The pseudocode of the algorithm is presented in algorithm 29.2. However, for the subsequent

**Algorithm 29.2** Intertwined_Bond_Exploration $(W, W', D)$

**Require:** two non-empty starting sets $W, W'$, an authorized domain $D$
**Ensure:** $T$ is the number of iterations of the intertwined explorations

```
A(0), W(0) ← Taboo_Bond_Explore_Clusters (W, D, ∅)
A(1), W(1) ← Taboo_Bond_Explore_Clusters (W', D, A(0))
t ← 1
repeat
    t ← t + 1
    T = A(t − 1)
    C, W = Taboo_Bond_Explore_Clusters (W(t − 2), D \ A(t − 2), T)
    A(t) = A(t − 2) ∪ C
    W(t) = W
until A(t − 1) ∪ A(t) = D
T ← t
return T
```

analysis, we will rely on the previous mathematical construction. We start with the following simple result.

**Lemma 29.4.** *For any $t \geq 0$, the aggregates $\mathcal{A}(t)$ and $\mathcal{A}(t+1)$ are disjoint and*

$$\mathcal{A}(2t) = \mathcal{B}\big(W, \mathcal{A}(2t), t\big)\,, \qquad \mathcal{A}(2t+1) = \mathcal{B}\big(W', \mathcal{A}(2t+1), t\big)\,. \tag{29.12}$$

*Proof.* We prove the result by induction on $t$. Let us examine first the case $t = 0$. By definition, we have

$$\mathcal{A}(0) \;=\; \text{Clusters}\,(W, D)\,, \qquad \mathcal{A}(1) \;=\; \text{Clusters}\,\big(W', D, \mathcal{A}(0)\big)\,. \tag{29.13}$$

Since the sets $W$ and $W'$ are not connected, then $\mathcal{A}(0) \cap \mathcal{A}(1) = \varnothing$. Moreover, it follows from (29.13) that $\mathcal{A}(0) = \mathcal{B}(W, \mathcal{A}(0), 0)$ and $\mathcal{A}(1) = \mathcal{B}(W', \mathcal{A}(1), 0)$, hence the result holds for $t = 0$. We perform next the induction step for the claim on the disjointness. Let $t \geq 1$ and suppose that the result has been proved at rank $t-1$, so that $\mathcal{A}(t-1) \cap \mathcal{A}(t) = \varnothing$. This implies that

$$\mathcal{A}(t) \cap \mathcal{A}(t+1) \;=\; \mathcal{A}(t) \cap \big(\mathcal{A}(t+1) \setminus \mathcal{A}(t-1)\big)\,. \tag{29.14}$$

However the increment $\mathcal{A}(t+1) \setminus \mathcal{A}(t-1)$ is disjoint from $\mathcal{A}(t)$, because the latter set is used as a taboo set when building the former (see formula (29.10)).

So we conclude from (29.14) that $\mathcal{A}(t) \cap \mathcal{A}(t+1) = \varnothing$. We perform finally the induction step for the inclusions (29.12). Let $t \geq 0$ and suppose that the result has been proved at rank $t$. It follows from formula (29.10) at rank $2t+2$ that

$$\begin{aligned}\mathcal{A}(2t+2) \setminus \mathcal{A}(2t) \,\subset\, & C\big(\partial_D^{out}\mathcal{A}(2t), \mathcal{A}(2t+2)\big) \\ & \subset \text{Shell}\,\big(\mathcal{B}(W, \mathcal{A}(2t), t), \mathcal{A}(2t+2), 1\big) \,\subset\, \mathcal{B}(W, \mathcal{A}(2t+2), t+1)\,.\end{aligned}$$

This implies that $\mathcal{A}(2t+2) \subset \mathcal{B}(W, \mathcal{A}(2t+2), t+1)$. The converse inclusion is obvious. Similarly, we have

$$\mathcal{A}(2t+3) \setminus \mathcal{A}(2t+1) \,\subset\, C\big(\partial_D^{out}\mathcal{A}(2t+1), \mathcal{A}(2t+3)\big) \,\subset\, \mathcal{B}(W', \mathcal{A}(2t+3), t+1)\,,$$

whence $\mathcal{A}(2t+3) \subset \mathcal{B}(W', \mathcal{A}(2t+3), t+1)$, and the result holds at rank $t+1$. $\square$

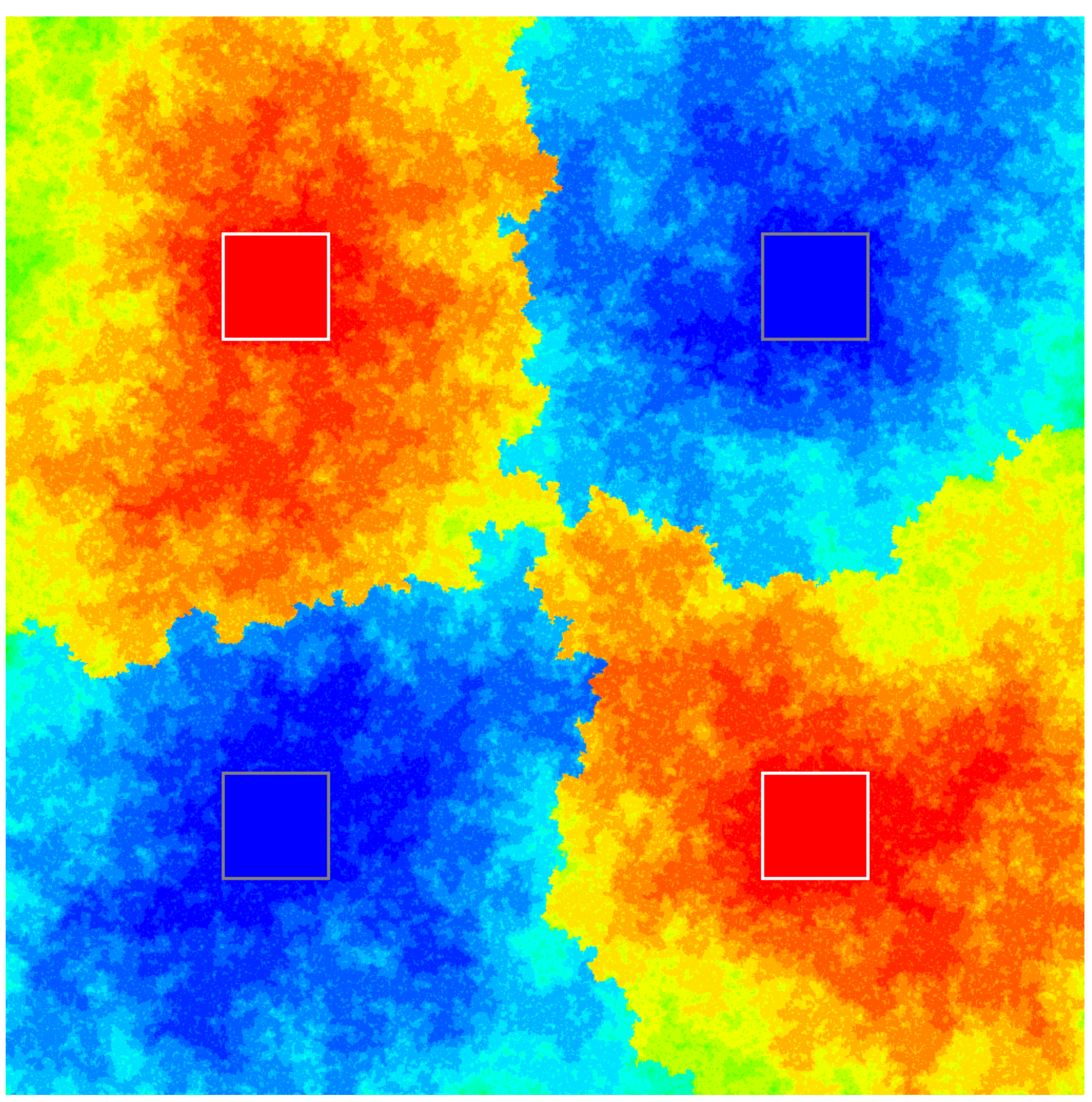

Figure 43: Intertwined explorations in $\Lambda(1024)$ for bond percolation at $p = 0.48$ $W, W' =$white,gray squares, lighter colors correspond to higher rounds.

## 29.5 Intersections in the bond model

Suppose that we are given two sets of vertices $W$ and $W'$ and a percolation configuration in the finite domain $D$ such that $W$ and $W'$ are not connected. Let $(\mathcal{A}(t), \mathcal{W}(t), t \geq 0)$ be intertwined explorations starting from $W$ and $W'$. We define the first intersection set $\mathcal{I}(1)$ as

$$\mathcal{I}(1) \,=\, \Delta_D\mathcal{A}(0) \cap \Delta_D\mathcal{A}(1)\,.$$

Of course this set might be empty. In fact, the bonds of $\mathcal{I}(1)$ are the pivotal bonds for the connection between $W$ and $W'$. On the event that $W$ and $W'$ are not connected, these bonds must be closed. We define the successive intersection sets $(\mathcal{I}(t), t \geq 0)$ by setting

$$\forall t \geq 1 \qquad \mathcal{I}(t) \,=\, \Delta_D\mathcal{A}(t-1) \cap \Delta_D\mathcal{A}(t) \setminus \bigcup_{1 \leq s \leq t-1} \mathcal{I}(s)\,. \tag{29.15}$$

It follows from the definition of the final iteration number $T$ that

$$\forall t > T \qquad \mathcal{I}(t) \,=\, \varnothing\,. \tag{29.16}$$

If a bond $e$ belongs to $\mathcal{I}(2)$, then $e$ is non-pivotal, closed and there exists another closed bond $f$ such that, if we open $f$, then the bond $e$ becomes pivotal for the connection between $W$ and $W'$. We call such a bond a pivot of second order. However, not all such bonds are included in $\mathcal{I}(2)$. For instance, if we exchange the roles of the sets $W$ and $W'$, we will end up with the same set $\mathcal{I}(1)$, but with a completely different set $\mathcal{I}(2)$ (disjoint from the previous one). In order to analyze the sets $\mathcal{I}(t)$, $t \geq 1$, we introduce a new important definition.

**The pivotal bonds of higher order.** Let $D$ be a finite subset of $\mathbb{Z}^d$. For $\omega$ in $\{0,1\}^{\mathbb{E}^d(D)}$, $E$ a subset of $\mathbb{E}^d(D)$, and a configuration $\eta$ in $\{0,1\}^E$, we denote by $\omega(E \leftarrow \eta)$ the configuration obtained by keeping unchanged the configuration $\omega|_{\mathbb{E}^d(D)\setminus E}$ and by setting the configuration $\omega|_E$ to be equal to $\eta$, i.e.,

$$\forall e \in \mathbb{E}^d(D) \qquad \omega(E \leftarrow \eta)(e) \,=\, \begin{cases} \omega(e) & \text{if } e \notin E \\ \eta(e) & \text{if } e \in E \end{cases}.$$

Let $\mathcal{D}$ be an event which is measurable with respect to the bonds in $\mathbb{E}^d(D)$ and let $k \geq 1$. We say that a bond $e$ is pivotal of order $k$ for the event $\mathcal{D}$ in the configuration $\omega$, or $k$-pivotal, if there exists a set $E$ of bonds of cardinality $k-1$ and a configuration $\eta$ of the bonds of $E$ such that, if we transform the configuration $\omega$ by forcing the edges of $E$ to be in the same states as in $\eta$, then the edge $e$ becomes pivotal for $\mathcal{D}$. We denote by $\mathcal{P}_k(\mathcal{D}, \omega)$ the set of the $k$-pivotal bonds for $\mathcal{D}$ in the configuration $\omega$. For $k = 1$, the set $\mathcal{P}_1(\mathcal{D}, \omega)$ coincides with the classical set $\mathcal{P}(\mathcal{D}, \omega)$ of the bonds which are pivotal for $\mathcal{D}$ in $\omega$. With the previous notation, we have, for $k \geq 1$,

$$\mathcal{P}_k(\mathcal{D}, \omega) \,=\, \Big\{ e \in \mathbb{E}^d(D) : \exists\, E \subset \mathbb{E}^d(D) \quad |E| \leq k-1 \\ \exists\, \eta \in \{0,1\}^E \quad e \in \mathcal{P}\big(\mathcal{D}, \omega(E \leftarrow \eta)\big) \Big\}\,.$$

With this definition, we see that the sets $(\mathcal{P}_k)_{k\geq 1}$ form a non–decreasing sequence of subsets of $\mathbb{E}^d(D)$. After some index $k \leq |\mathbb{E}^d(D)|$, the sequence becomes stationary.

We apply next the previous definitions to the disconnection event

$$\mathcal{D} = \{\, W \not\longleftrightarrow W' \,\}\,.$$

The bonds in the successive intersections $(\mathcal{I}(t), t \geq 1)$ have a particular significance for the event $\mathcal{D}$. The first intersection set $\mathcal{I}(1)$ is exactly the set of the pivotal bonds $\mathcal{P}(\mathcal{D})$. This is not true for the next sets $\mathcal{I}(k)$, $k \geq 2$. However, we have the following weaker but nonetheless important result.

**Lemma 29.5.** *For $k \geq 1$, the intersection set $\mathcal{I}(k)$ is included in $\mathcal{P}_k(\mathcal{D})$, i.e., the bonds in $\mathcal{I}(k)$ are $k$-pivotal.*

*Proof.* We do the proof for the even indices, the odd indices can be handled similarly. Let $\omega$ be a configuration in $\mathcal{D}$, let $k \geq 1$ and let us consider a bond $e$ belonging to $\mathcal{I}(2k)$. By the definition of $\mathcal{I}(2k)$, we can write $e = \langle x, x' \rangle$, where $x$ is in $\mathcal{A}(2k)$ and $x'$ in $\mathcal{A}(2k-1)$. It follows from lemma 29.4 that in addition $T_{\mathcal{A}(2k)}(W, x) \leq k$ and $T_{\mathcal{A}(2k-1)}(W', x') \leq k-1$. In particular, there exists a path

$$x_0, e_0, x_1, e_1, \cdots, e_{\ell-1}, x_\ell$$

in $\mathcal{A}(2k)$ joining a vertex $x_0$ of $W$ to $x_\ell = x$ such that at most $k$ bonds among $e_0, \cdots, e_{\ell-1}$ are closed. We denote by $E$ this set of bonds. Similarly, there exists a path

$$x'_0, e'_0, x'_1, e'_1, \cdots, e'_{\ell'-1}, x'_{\ell'}$$

in $\mathcal{A}(2k-1)$ joining a vertex $x'_0$ of $W'$ to $x'_{\ell'} = x'$ such that at most $k-1$ bonds among $e'_0, \cdots, e'_{\ell'-1}$ are closed. We denote by $E'$ this set of bonds. Let $\widetilde{\omega}$ be the configuration obtained from $\omega$ by opening all the bonds of $E \cup E'$. We claim that this configuration $\widetilde{\omega}$ is still in $\mathcal{D}$. Indeed, the edges of $E$ (respectively $E'$) are all included in $\mathbb{E}^d(\mathcal{A}(2k))$ (respectively $\mathbb{E}^d(\mathcal{A}(2k-1))$), so their opening cannot create a connection between $W$ (respectively $W'$) and a vertex of $\mathcal{A}(2k-1)$ (respectively $\mathcal{A}(2k)$). Therefore the configuration $\widetilde{\omega}$ is in $\mathcal{D}$. Moreover, the opening of the bond $e$ in $\widetilde{\omega}$ creates a connection between $x_0$ and $x'_0$, thus $e$ is pivotal for $\mathcal{D}$ in $\widetilde{\omega}$. We conclude that $e$ is in $\mathcal{P}_{2k}(\mathcal{D}, \omega)$. □

Notice that $\mathcal{I}(1) \cup \mathcal{I}(2)$ is in general strictly included in $\mathcal{P}_2(\mathcal{D})$. However, all the bonds in $\mathcal{P}_2(\mathcal{D})$ can be identified with the help of several intertwined explorations. In fact, we should run two intertwined explorations, a first one starting from $W$ and $W'$, and a second one starting from $W'$ and $W$. This is not our goal for the time being, but we might come back to this point in the future. Similarly, for $k \geq 2$, the intersection set $\mathcal{I}(k)$ does not exhaust the pivots of order $k$. Later, we will describe more precisely the structure of the set $\mathcal{P}_k(\mathcal{D})$, and we will design a system of intertwined explorations in order to reveal the whole set $\mathcal{P}_k(\mathcal{D})$.

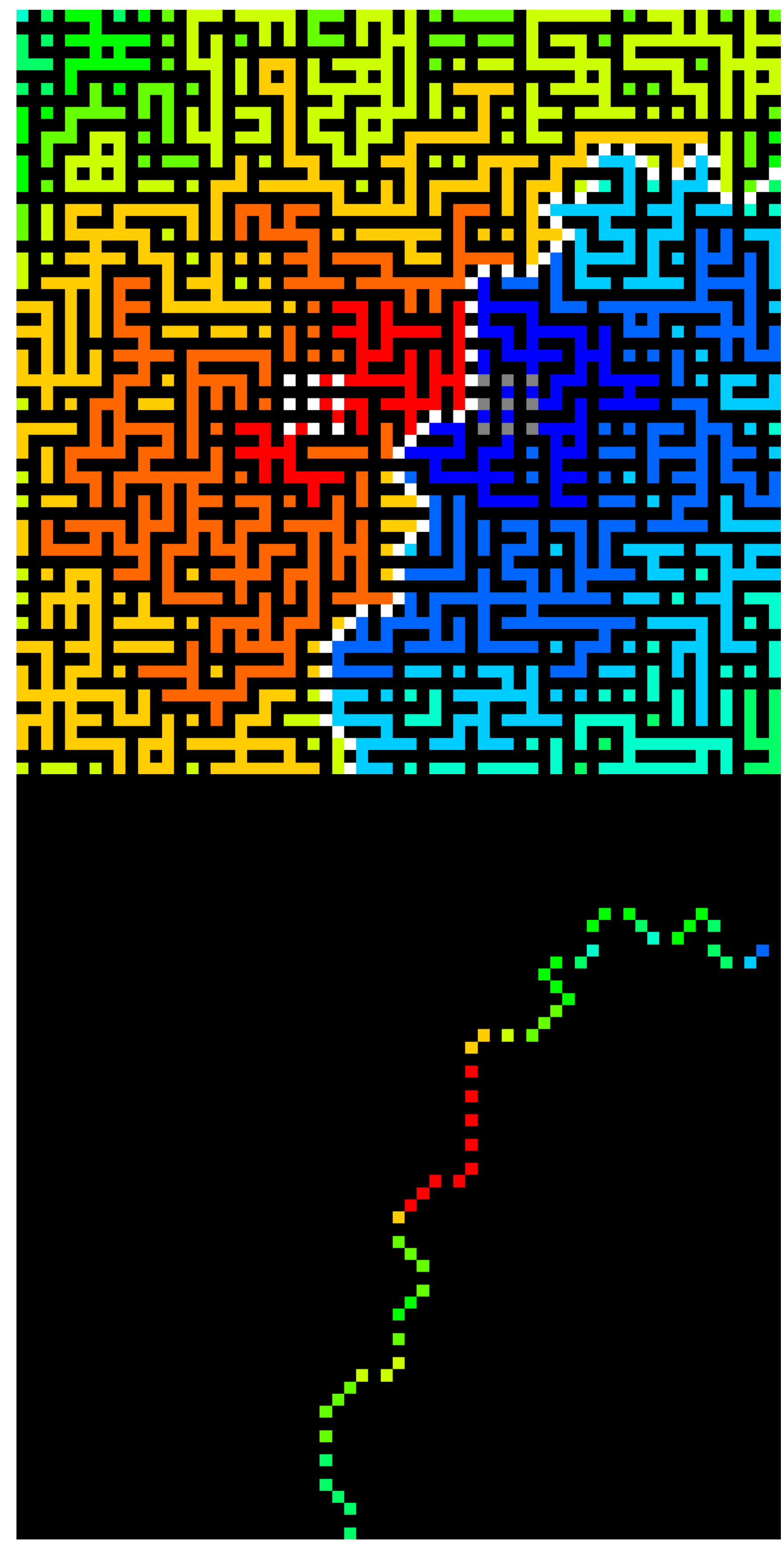

Figure 44: Top: Intertwined explorations, intersections are colored in white. Bottom: Intersections only, colored from red to blue according to the order. 1st, 2nd explorers: red to yellow, blue to green. $W, W'$ = white,gray squares. Bond percolation, $\Lambda(32)$, $p = 0.406$.

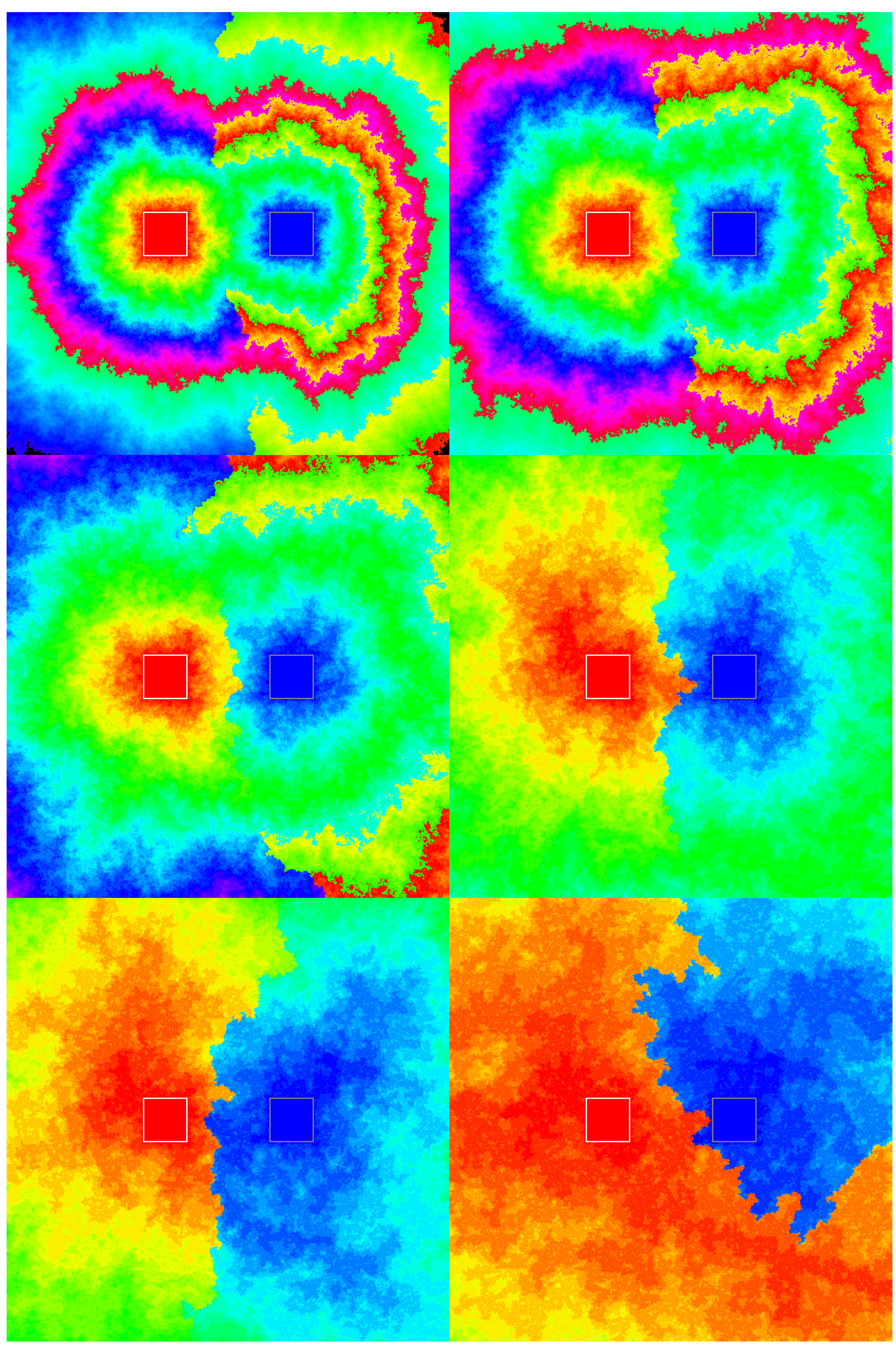

Figure 45: Intertwined bond explorations, $W, W'$ = white,gray squares, $\Lambda(1024)$. Rainbow colors starting from red for the first explorer, and blue for the second. From left to right and top to bottom: $p = 0.4, 0.425, 0.45, 0.47, 0.48, 0.49$.

## 29.6 Termination

The role of the taboo sets used in the intertwined explorations is to prevent each explorer from entering into the set of the vertices already discovered by the other explorer. Although each explorer forbids the other to penetrate his own territory, the union of the two sets discovered by the two explorers coincides with the set that a single non-taboo exploration would visit if it were started from the union of $W$ and $W'$. This implies in particular that the intertwined explorations terminate after a finite number of iterations. We give a precise formulation of this important result in the next proposition.

**Proposition 29.6.** *Let $W$, $W'$ be two subsets of a finite domain $D$ and let us consider a percolation configuration in $D$ such that $W$ and $W'$ are not connected. Let $(\mathcal{A}(t), \mathcal{W}(t), t \geq 0)$ be the intertwined explorations starting from $W$ and $W'$. For any $t \geq 0$, we have*

$$\mathcal{A}(2t) \cup \mathcal{A}(2t+1) \,=\, \mathcal{B}\big(W \cup W', D, t\big)\,, \tag{29.17}$$

$$\mathcal{A}(2t+1) \cup \mathcal{A}(2t+2) \,=\, \mathcal{B}\big(W', D, t\big) \cup \mathcal{B}\big(W, D, t+1\big)\,. \tag{29.18}$$

*Proof.* Since the sets $\text{Shell}\,(W \cup W', D, s)$, $0 \leq s \leq t$, are pairwise disjoint and

$$\mathcal{B}\big(W \cup W', D, t\big) \,=\, \bigcup_{0 \leq s \leq t} \text{Shell}\,(W \cup W', D, s)\,,$$

the formula (29.17) is equivalent to saying that

$$\text{Shell}\,(W \cup W', D, 0) \,=\, \mathcal{A}(0) \cup \mathcal{A}(1) \tag{29.19}$$

and for $t \geq 1$,

$$\text{Shell}\,(W \cup W', D, t) \,=\, \big(\mathcal{A}(2t) \cup \mathcal{A}(2t+1)\big) \setminus \big(\mathcal{A}(2t-2) \cup \mathcal{A}(2t-1)\big)\,. \tag{29.20}$$

We prove first the case $t = 0$ of formula (29.17), which is equivalent to (29.19). We have $\mathcal{A}(0) \,=\, \text{Clusters}\,(W, D)$ and since $W$ and $W'$ are not connected inside the domain $D$, then

$$\mathcal{A}(1) \,=\, \text{Clusters}\,\big(W', D, \mathcal{A}(0)\big) \,=\, \text{Clusters}\,\big(W', D, \varnothing\big) \,=\, C(W', D)\,,$$

whence

$$\mathcal{A}(0) \cup \mathcal{A}(1) \,=\, \text{Shell}\,(W \cup W', D)\,.$$

We proceed next by induction on $t$ to prove the formula (29.17). Suppose that the result has been proved at rank $t$ for some $t \geq 0$. We will prove that the equality (29.20) holds at rank $t+1$. By the definition of $\mathcal{A}(2t+2)$ (see formulas (29.3) and (29.9)), we have the inclusions

$$\begin{aligned}\mathcal{A}(2t+2) \setminus \mathcal{A}(2t) \,&\subset\, \text{Clusters}\,\big(\mathcal{W}(2t), D \setminus \mathcal{A}(2t), \mathcal{A}(2t+1)\big)\\ &\subset\, \text{Clusters}\,\big(\mathcal{W}(2t), D \setminus \mathcal{A}(2t)\big) \,\subset\, \text{Shell}\,\big(\mathcal{A}(2t), D, 1\big)\,.\end{aligned} \tag{29.21}$$

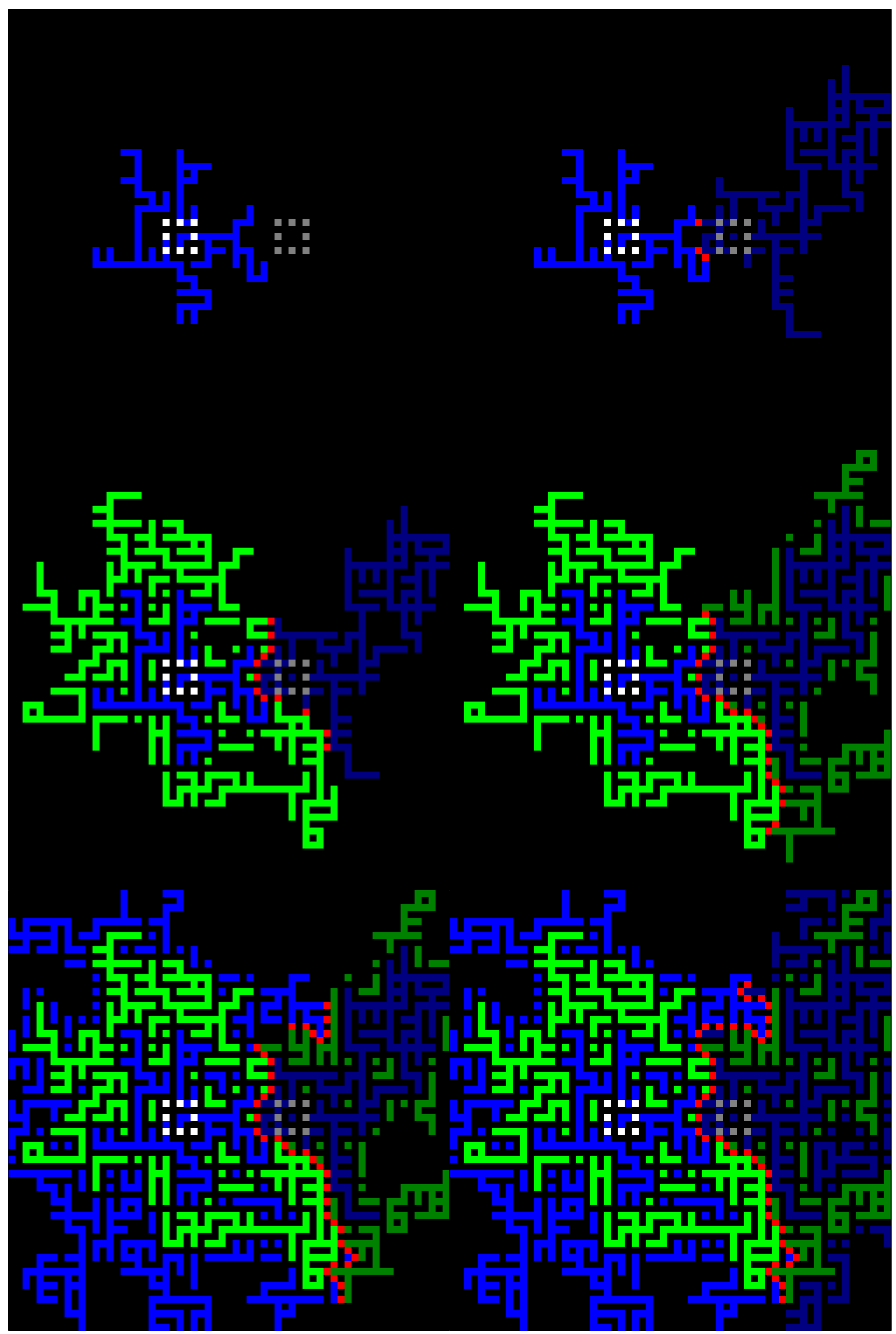

Figure 46: Intertwined explorations in $\Lambda(32)$ for bond percolation at $p = 0.39$. $W$ = white square, $W'$ = gray square. From left to right and top to bottom: $\mathcal{A}(0), \mathcal{A}(0) \cup \mathcal{A}'(0), \mathcal{A}(1) \cup \mathcal{A}'(0), \mathcal{A}(1) \cup \mathcal{A}'(1), \mathcal{A}(2) \cup \mathcal{A}'(1), \mathcal{A}(2) \cup \mathcal{A}'(2)$. The right explorer has darker colors than the left explorer.

Similarly, we have

$$\begin{aligned}\mathcal{A}(2t+3)\setminus\mathcal{A}(2t+1)\ &\subset\ \text{Clusters}\big(\mathcal{W}(2t+1),D\setminus\mathcal{A}(2t+1),\mathcal{A}(2t+2)\big)\\ &\subset\ \text{Clusters}\big(\mathcal{W}(2t+1),D\setminus\mathcal{A}(2t+1)\big)\ \subset\ \text{Shell}\big(\mathcal{A}(2t+1),D,1\big)\,.\end{aligned}\tag{29.22}$$

Using the induction hypothesis, we have furthermore

$$\mathcal{A}(2t)\cup\mathcal{A}(2t+1)\ \subset\ \bigcup_{0\le s\le t}\text{Shell}\,(W\cup W',D,s)\,,\tag{29.23}$$

whence

$$\begin{aligned}\text{Shell}\big(\mathcal{A}(2t)\cup\mathcal{A}(2t+1),D,1\big)\ &\subset\ \text{Shell}\Big(\bigcup_{0\le s\le t}\text{Shell}\,(W\cup W',D,s),D,1\Big)\\ &\subset\ \text{Shell}\,(W\cup W',D,t+1)\,.\end{aligned}\tag{29.24}$$

Putting together (29.21), (29.22), (29.23) and (29.24), we conclude that

$$\big(\mathcal{A}(2t+2)\cup\mathcal{A}(2t+3)\big)\setminus\big(\mathcal{A}(2t)\cup\mathcal{A}(2t+1)\big)\ \subset\ \text{Shell}\,(W\cup W',D,t+1)\,.\tag{29.25}$$

We prove next the converse inclusion. Using the semigroup identity (12.10), the induction hypothesis, and formula (28.8), we have

$$\begin{aligned}\text{Shell}\,(W\cup W',D,t+1)\ &=\ \text{Shell}\Big(\bigcup_{0\le s\le t}\text{Shell}\,(W\cup W',D,s),D,1\Big)\\ &=\ \text{Shell}\big(\mathcal{A}(2t)\cup\mathcal{A}(2t+1),D,1\big)\ =\ C\big(\partial_D^{out}\big(\mathcal{A}(2t)\cup\mathcal{A}(2t+1)\big),D\big)\,.\end{aligned}\tag{29.26}$$

Thanks to proposition 29.3, two consecutive waiting aggregates of the sequence are coherent, whence

$$\partial_D^{out}\big(\mathcal{A}(2t)\cup\mathcal{A}(2t+1)\big)\ \subset\ \mathcal{W}(2t)\cup\mathcal{W}(2t+1)\,.\tag{29.27}$$

The inclusions (29.26) and (29.27) imply furthermore that

$$\begin{aligned}\text{Shell}\,(W\cup W',D,t+1)\ &\subset\ C\big(\mathcal{W}(2t)\cup\mathcal{W}(2t+1),D\big)\\ &\subset\ C\big(\mathcal{W}(2t),D\big)\cup C\big(\mathcal{W}(2t+1),D\big)\,.\end{aligned}\tag{29.28}$$

Suppose that $x$ is in $\text{Shell}\,(W\cup W',D,t+1)$. Using the inclusion (29.28), two cases can occur:

• $x\in C\big(\mathcal{W}(2t),D\big)$. In this case, there exists a path $z_0,e_0,z_1,\dots,e_{r-1},z_r$ in $D$ joining a vertex $z_0=y$ of $\mathcal{W}(2t)$ to $z_r=x$ such that the bonds $e_0,\dots,e_{r-1}$ are open. From the inclusion (29.23), we know that $x$ is neither in $\mathcal{A}(2t)$ nor in $\mathcal{A}(2t+1)$. If one of the vertices $z_0,\dots,z_{r-1}$ were to belong to $\mathcal{A}(2t)\cup\mathcal{A}(2t+1)$, then $x$ would also be in this set, so this is not the case. Therefore, the connection between $x$ and $\mathcal{W}(2t)$ is realized outside $\mathcal{A}(2t)\cup\mathcal{A}(2t+1)$ and thus $x$ is in $\text{Clusters}\big(\mathcal{W}(2t),D\setminus\mathcal{A}(2t),\mathcal{A}(2t+1)\big)$, which is itself included in $\mathcal{A}(2t+2)$.

• $x \in C\big(\mathcal{W}(2t+1), D\big)$. In this case, there exists a path $z_0, e_0, z_1, \dots, e_{r-1}, z_r$ in $D$ joining a vertex $z_0 = y$ of $\mathcal{W}(2t+1)$ to $z_r = x$ such that the bonds $e_0, \dots, e_{r-1}$ are open. From the inclusion (29.23), we know that $x$ is not in $\mathcal{A}(2t+1)$, and this implies that none of the vertices $z_0, \dots, z_{r-1}$ is in $\mathcal{A}(2t+1)$. Furthermore, if one of these vertices belongs to $\mathcal{A}(2t+2)$, then so does $x$. Otherwise, the connection between $x$ and $\mathcal{W}(2t+1)$ is realized outside $\mathcal{A}(2t+1) \cup \mathcal{A}(2t+2)$ and thus $x$ belongs to $\text{Clusters}\big(\mathcal{W}(2t+1), D \setminus \mathcal{A}(2t+1), \mathcal{A}(2t+2)\big)$, which is itself included in $\mathcal{A}(2t+3)$.

In both cases, we see that $x$ belongs to $\mathcal{A}(2t+2) \cup \mathcal{A}(2t+3)$. This yields the converse inclusion of (29.25), and the equality (29.20) at rank $t+1$, which implies the formula (29.17) at rank $t+1$. This concludes the induction step.

We prove now formula (29.18). We have by definition

$$\begin{aligned} \mathcal{A}(1) \;&=\; \mathcal{B}\big(W', D, 0\big)\,, \\ \mathcal{A}(2) \;&=\; \text{Shell}\,(W, D, 0) \cup \text{Shell}\,(W, D \setminus \mathcal{A}(1), 1)\,, \end{aligned}$$

from which we deduce that

$$\mathcal{A}(1) \cup \mathcal{A}(2) \;=\; \mathcal{B}\big(W', D, 0\big) \cup \mathcal{B}\big(W, D, 1\big)\,. \tag{29.29}$$

This is the case $t = 0$ of formula (29.18). We show finally how the identity (29.18) for $t \geq 1$ can be deduced from the identity (29.17). Let us set

$$\widetilde{W} = \mathcal{A}(1)\,, \qquad \widetilde{W}' = \mathcal{A}(2)\,,$$

and let $(\widetilde{\mathcal{A}}(t), \widetilde{\mathcal{W}}(t), t \geq 0)$ be the intertwined explorations starting from $\widetilde{W}, \widetilde{W}'$. It can be checked that the sequence $(\widetilde{\mathcal{A}}(t), t \geq 0)$ is the left shift of the sequence $(\mathcal{A}(t), t \geq 0)$, i.e.,

$$\forall t \geq 0 \qquad \widetilde{\mathcal{A}}(t) \;=\; \mathcal{A}(t+1)\,, \tag{29.30}$$

Applying the identity (29.17) and using (29.29), we obtain that, for $t \geq 0$,

$$\begin{aligned} \widetilde{\mathcal{A}}(2t) \cup \widetilde{\mathcal{A}}(2t+1) \;&=\; \mathcal{B}\big(\widetilde{W} \cup \widetilde{W}', D, t\big) \;=\; \mathcal{B}\big(\mathcal{A}(1) \cup \mathcal{A}(2), D, t\big)\,, \\ &=\; \mathcal{B}\Big(\mathcal{B}\big(W', D, 0\big) \cup \mathcal{B}\big(W, D, 1\big), D, t\Big)\,, \\ &=\; \mathcal{B}\Big(\mathcal{B}\big(W', D, 0\big), D, t\Big) \cup \mathcal{B}\Big(\mathcal{B}\big(W, D, 1\big), D, t\Big) \\ &=\; \mathcal{B}\big(W', D, t\big) \cup \mathcal{B}\big(W, D, t+1\big)\,. \end{aligned} \tag{29.31}$$

Formulas (29.30) and (29.31) together yield formula (29.18) for $t \geq 1$. □

We recall the definition of the final iteration time $T$, introduced in (29.11):

$$T \;=\; \min\,\big\{\, t \geq 1 : \mathcal{A}(t-1) \cup \mathcal{A}(t) = D \,\big\}\,.$$

With the help of proposition 29.6, we can obtain a control on the final iteration time $T$ in a finite connected domain $D$. Indeed, the identity (29.17) implies that, for any $t \geq 0$, the union $\mathcal{A}(2t) \cup \mathcal{A}(2t+1)$ contains all the vertices that are within travel distance less or equal than $t$ from $W \cup W'$. So the final iteration number $T$ is less than or equal to one plus twice the maximal travel distance between a vertex of $D$ to the set $W \cup W'$. We formalize this in the next corollary.

**Corollary 29.7.** *Let $D$ be a finite connected subset of $\mathbb{Z}^d$. Let $W$, $W'$ be two subsets of $D$ and let us consider a percolation configuration in $D$ such that $W$ and $W'$ are not connected. The final iteration time $T$ is finite and satisfies*

$$T \;\leq\; 2\max\big\{\,T_D(W\cup W',x) : x\in D\,\big\}+1\,. \tag{29.32}$$

*Proof.* Let us define $T^*$ as

$$T^* \;=\; \max\big\{\,T_D(W\cup W',x) : x\in D\,\big\}\,.$$

It follows from proposition 29.6 that

$$\forall t\geq T^* \qquad \mathcal{A}(2t)\cup\mathcal{A}(2t+1) \;=\; D\,,$$

therefore the final iteration number $T$ satisfies $T\leq 2T^*+1$. □

A major goal is to control the intersection sets. In the next corollary, with the help of proposition 29.6, we compute the travel time between the sets $W, W'$ and the endpoints of the bonds in $\mathcal{I}(k)$.

**Corollary 29.8.** *Let $k\geq 1$ and let $e=\langle x,x'\rangle$ be a bond belonging to $\mathcal{I}(k)$, where $x$ is in $\mathcal{A}(k)$ and $x'$ in $\mathcal{A}(k-1)$. If $k=2\ell$ is even, we have*

$$T_{\mathcal{A}(2\ell)}(W,x)=T_D(W\cup W',x)=\ell\,, \tag{29.33}$$

$$T_{\mathcal{A}(2\ell-1)}(W',x')=T_D(W\cup W',x')=\ell-1\,. \tag{29.34}$$

*If $k=2\ell+1$ is odd, we have*

$$T_{\mathcal{A}(2\ell+1)}(W,x)=T_D(W\cup W',x)=\ell\,,$$
$$T_{\mathcal{A}(2\ell)}(W',x')=T_D(W\cup W',x')=\ell\,.$$

*Proof.* We do the proof for the even indices, the odd indices can be handled similarly. Let $\omega$ be a configuration in $\mathcal{D}$, let $\ell\geq 1$ and let $e=\langle x,x'\rangle$ be a bond belonging to $\mathcal{I}(2\ell)$, where $x$ is in $\mathcal{A}(2\ell)$ and $x'$ in $\mathcal{A}(2\ell-1)$. Since $e$ is not in $\mathcal{I}(2\ell-1)$, then $x$ is not in $\mathcal{A}(2\ell-2)$. Suppose that $x'$ is in $\mathcal{A}(2\ell-3)$. By lemma 29.4, we would have $T_{\mathcal{A}(2\ell-3)}(W',x')\leq\ell-2$, whence

$$T_D(W',x)\;\leq\; T_{\mathcal{A}(2\ell-3)}(W',x')+1\;\leq\;\ell-1\,.$$

It would follow from the identity (29.17) that $x$ is in $\mathcal{A}(2\ell-2)\cup\mathcal{A}(2\ell-1)$, therefore $x$ would be in $\mathcal{A}(2\ell-1)$ as well, which is absurd. Thus $x$ is not in $\mathcal{A}(2\ell-2)\cup\mathcal{A}(2\ell-1)$ and $x'$ is not in $\mathcal{A}(2\ell-3)\cup\mathcal{A}(2\ell-2)$. Using again proposition 29.6, we conclude that

$$x\notin\mathcal{B}\big(W\cup W',D,\ell-1\big)\,,\qquad x'\notin\mathcal{B}\big(W',D,\ell-2\big)\cup\mathcal{B}\big(W,D,\ell-1\big)\,. \tag{29.35}$$

From lemma 29.4, we have also

$$T_{\mathcal{A}(2\ell)}(W,x)\leq\ell\,,\qquad T_{\mathcal{A}(2\ell-1)}(W',x')\leq\ell-1\,. \tag{29.36}$$

The equalities (29.33) and (29.34) follow from (29.35) and (29.36). □

Finally, we compute the order of pivotality of the bonds in $\mathcal{I}(k)$ for $k \geq 1$.

**Corollary 29.9.** *For $k \geq 1$, the intersection set $\mathcal{I}(k)$ is disjoint from $\mathcal{P}_{k-1}(\mathcal{D})$.*

*Proof.* We do the proof for the even indices, the odd indices can be handled similarly. Let $\omega$ be a configuration in $\mathcal{D}$, let $k \geq 1$ and let us consider a bond $e$ belonging to $\mathcal{I}(2k)$. By the definition of $\mathcal{I}(2k)$, we can write $e = \langle x, x' \rangle$, where $x$ is in $\mathcal{A}(2k)$ and $x'$ in $\mathcal{A}(2k-1)$. Suppose that $e$ is in $\mathcal{P}_{2k-1}(\mathcal{D}, \omega)$. By definition, there exists a set $E$ of bonds of cardinality $2k-2$ and a configuration $\eta$ of the bonds of $E$ such that $e$ is pivotal for $\mathcal{D}$ in $\omega(E \leftarrow \eta)$. In particular, there exists a path $x_0, e_0, x_1, e_1, \cdots, e_{\ell-1}, x_\ell$ joining a vertex $x_0$ of $W$ to a vertex $x_\ell$ of $W'$, which goes through the bond $e$, and such that all the bonds of this path are open in $\omega(E \leftarrow \eta)$, except $e$. Let $m$ be the index of $e$ in the path, so that

$$e_m = e\,, \qquad \{\, x_{m-1}, x_m \,\} \;=\; \{\, x, x' \,\}\,. \tag{29.37}$$

The configurations $\omega(E \leftarrow \eta)$ and $\omega$ differ in at most $2k-2$ bonds, therefore at most $2k-2$ bonds of this path, in addition to $e_m$, are closed in $\omega$. This implies that, in $\omega$,

$$T_D(x_0, x_{m-1}) + T_D(x_m, x_\ell) \;\leq\; 2k - 2\,. \tag{29.38}$$

According to (29.37), two cases can occur:

- $x_{m-1} = x, x_m = x'$. We have then

$$T_D(x_0, x_{m-1}) \;\geq\; T_D(W \cup W', x)\,, \quad T_D(x_m, x_\ell) \;\geq\; T_D(W \cup W', x')\,.$$

- $x_{m-1} = x', x_m = x$. We have then

$$T_D(x_0, x_{m-1}) \;\geq\; T_D(W \cup W', x')\,, \quad T_D(x_m, x_\ell) \;\geq\; T_D(W \cup W', x)\,.$$

In both cases, we have

$$T_D(x_0, x_{m-1}) + T_D(x_m, x_\ell) \;\geq\; T_D(W \cup W', x) + T_D(W \cup W', x')\,. \tag{29.39}$$

Using the inequalities (29.33) and (29.34) of corollary 29.8, together with inequality (29.39), we obtain that

$$T_D(x_0, x_{m-1}) + T_D(x_m, x_\ell) \;\geq\; 2k - 1\,,$$

but this stands in contradiction with inequality (29.38). We conclude that the bond $e$ is not in $\mathcal{P}_{2k-1}(\mathcal{D}, \omega)$. ☐

Putting together lemma 29.5 and corollary 29.9, we conclude that, for $k \geq 1$, the bonds in $\mathcal{I}(k)$ are $k$-pivotal but not $(k-1)$-pivotal.

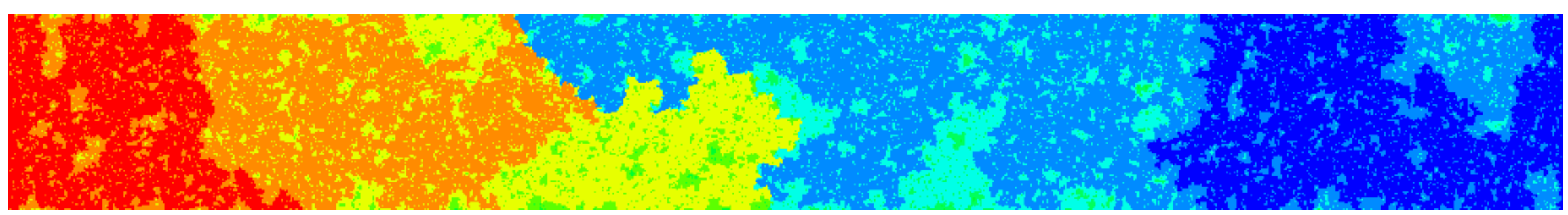

Figure 47: Intertwined explorations in the rectangle $1024 \times 128$, $p = 0.515$.

## 29.7 The envisaged strategy

Here we outline the strategy that we will try to implement in the rest of the text. Let us consider a finite subset $D$ of $\mathbb{Z}^d$ and an event $\mathcal{D}$ occurring in the percolation configuration restricted to $D$. Typically, we have in mind an event $\mathcal{D}$ which is unlikely to occur, and we aim at proving that $P(\mathcal{D})$ goes to 0 as the size of $D$ goes to $\infty$ (the set $D$ is usually a cubic box and we send its side length to $\infty$). Ideally, we also wish to obtain a quantitative estimate on the speed of convergence to 0 of $P(\mathcal{D})$. The prototype of the event $\mathcal{D}$ is a disconnection event: there exist two fixed subsets $W$ and $W'$ of $D$ and $\mathcal{D} = \big\{ W \not\leftrightarrow W' \big\}$ is the event that $W$ and $W'$ are not connected by an open path in the percolation configuration restricted to $D$. However, we might consider more complicated events. This will be the case in subsection 30.2, when we study the largest finite cluster in a box. The first essential ingredient of the strategy is an adequate inequality for a specific random min-cut problem between $W$ and $W'$ in $D$. We explain this in the next paragraph.

**Random min-cut problem.** We start by defining properly the relevant random min-cut problem between $W$ and $W'$ in $D$. The central notion is that of a separating set or a cut set, we recall it next.

**Definition 29.10.** *A set of edges $E \subset \mathbb{E}^d(D)$ is said to separate $W$ and $W'$ in $D$ if there is no path in the graph $(D, \mathbb{E}^d(D) \setminus E)$ connecting $W$ and $W'$.*

This was the notion used to define surface tension in the percolation model [25]. After an adequate extension to the FK model, it led to an operational definition of the surface tension in the Ising model, which was instrumental to prove the Wulff theorem [31, 26]. Here is another equivalent and more classical definition.

**Definition 29.11.** *A set of edges $E \subset \mathbb{E}^d(D)$ is a cut between $W$ and $W'$ in $D$ if every path in $D$ from $W$ to $W'$ goes through an edge of $E$.*

Among cut sets between $W$ and $W'$, of particular interest are the cuts which are minimal for the inclusion: they are called the $(W, W')$-cutsets in $D$. The random min-cut problem between the sets $W$ and $W'$ for a percolation configuration in $D$ realizing the event $\mathcal{D}$ is to find a set $E$ of closed edges in $\mathbb{E}^d(D)$ that separates $W$ and $W'$ and that has minimum cardinality. More precisely, we define a random variable $\text{Min-Cut}\,(W, W', D, \mathcal{D})$ by setting

$$\text{Min-Cut}\,(W, W', D, \mathcal{D}) \,=\, +\infty \quad \text{on the event} \quad \mathcal{D}^c\,, \tag{29.40}$$

and for a percolation configuration in $D$ which realizes the event $\mathcal{D}$,

$$\text{Min-Cut}\,(W, W', D, \mathcal{D}) \,=\, \min \Big\{ |E| : E \subset \mathbb{E}^d(D),\, \forall e \in E \quad e \text{ is closed}, \\ E \text{ separates } W \text{ and } W' \text{ in } D \Big\}. \tag{29.41}$$

Typically, the event $\mathcal{D}$ is the event that the sets $W$ and $W'$ are disconnected in the percolation configuration restricted to $D$ (or a subset of this event).

If this disconnection event does not occur, then there is no admissible candidate for the minimum in formula (29.41) and its value is infinite. However, we might consider events $\mathcal{D}$ that involve additional conditions besides the mere disconnection between $W$ and $W'$, hence it is safer to impose automatically the definition (29.40). If the sets $W$ and $W'$ are not connected in the percolation configuration restricted to $D$, as long as we deal with a finite set $D$, the above min-cut problem (29.41) is well posed and admits at least a solution $E^*$, i.e., there exists a cut $E^*$ between $W$ and $W'$ in $D$ such that all the edges of $E^*$ are closed and

$$|E^*| \,=\, \text{Min-Cut}\,(W, W', D, \mathcal{D})\,.$$

Obviously, any solution $E^*$ has to be a $(W, W')$-cutset in $D$. This kind of problem has been studied a lot in the more general context of the first passage percolation model (see [32] and the references therein, or [36] for a more recent work). In this model, each edge is endowed with a capacity, and these capacities are i.i.d. random variables. The Bernoulli percolation model corresponds to the specific case of the first passage percolation model where the common distribution of the capacities (also called the passage times) is the Bernoulli law with parameter $p$. We define then the capacity of a cut $E$ as the sum of the capacities of its edges, and the aim is to find a cut having minimum capacity. To obtain exactly the formulation given in (29.41), we consider the Bernoulli percolation model, and we define the capacity of a bond to be 1 if it is closed and infinite if it is open. The famous max-flow min-cut theorem asserts that this kind of minimum cut problem is dual to a maximal flow problem. For a maximal flow problem, the famous Ford-Fulkerson algorithm provides an efficient method to compute a solution (see the classical book [55]). By the way, it does not seem plausible that there is a method to compute the minimum (29.41) without actually exhibiting at least one solution. In general, for a given min-cut problem, there is no simple formula expressing the value of Min-Cut $(W, W', D, \mathcal{D})$, especially when the underlying graph is random, as in our case. In fact, the quantity Min-Cut $(W, W', D, \mathcal{D})$ is a random variable, whose distribution is extremely complex. So, we will contend ourselves with an inequality on Min-Cut $(W, W', D, \mathcal{D})$, of the form:

$$\text{Min-Cut}\,(W, W', D, \mathcal{D}) \,\geq\, \underline{\text{min-cut}}\,(W, W', D, \mathcal{D})\,, \tag{29.42}$$

where $\underline{\text{min-cut}}\,(W, W', D, \mathcal{D})$ is a deterministic number depending on the event $\mathcal{D}$ and the geometry of $W, W', D$. The lower bound depends on the event $\mathcal{D}$ because, thanks to the definition (29.40), we need only to consider percolation configurations realizing the event $\mathcal{D}$ for the min-cut problem (29.41). For the value $\underline{\text{min-cut}}\,(W, W', D, \mathcal{D})$, we can take the minimum cardinality for a cut between $W$ and $W'$ in the graph with vertex set $D$ and edge set $\mathbb{E}^d(D)$. Indeed, the admissible sets of edges in the min-cut problem (29.41) are sets of closed bonds which constitute a cut between $W$ and $W'$ in $D$. So if we forget about the constraint that the edges have to be closed, we obtain a larger set of candidates and this provides a deterministic lower bound on Min-Cut $(W, W', D, \mathcal{D})$.

**Intertwined explorations.** The second ingredient of the story is the intertwined explorations $(\mathcal{A}(t), \mathcal{W}(t), t \geq 0)$ inside $D$ starting from $W$ and $W'$. Whenever the event $\mathcal{D}$ occurs, the sets $W$ and $W'$ are disconnected in the percolation configuration restricted to $D$, therefore there exists a cut set $E^*$ between $W$ and $W'$ in $D$ consisting exclusively of closed bonds. The purpose of the intertwined explorations is to reveal such a cut set $E^*$ through a procedure that allows to develop useful probabilistic estimates. We suppose that $D$ is a finite connected subset of $\mathbb{Z}^d$. It follows from corollary 29.7 that the final iteration number $T$ of the intertwined explorations, defined in formula (29.11), is finite and moreover the whole set $D$ is visited after the completion of the iteration number $T$. The two aggregates $\mathcal{A}(T-1), \mathcal{A}(T)$ are disjoint, and they are such that $D = \mathcal{A}(T-1) \cup \mathcal{A}(T)$ and

$$\text{either} \quad W \subset \mathcal{A}(T-1),\, W' \subset \mathcal{A}(T) \quad \text{or} \quad W' \subset \mathcal{A}(T-1),\, W \subset \mathcal{A}(T)\,.$$

Therefore the set of bonds $\Delta_D \mathcal{A}(T-1) \cap \Delta_D \mathcal{A}(T)$ is a set of closed bonds which separates the sets $W$ and $W'$ in $D$. In particular, it is an admissible candidate for the min-cut problem (29.41) and this implies that

$$\text{Min-Cut}\,(W, W', D, \mathcal{D}) \;\leq\; \big|\Delta_D \mathcal{A}(T-1) \cap \Delta_D \mathcal{A}(T)\big|\,. \tag{29.43}$$

**Control on the final iteration number.** The next step consists in getting a control on the tail of the distribution of the final iteration number $T$. This is typically achieved with the help of the inequality (29.32) of corollary 29.7, which tells us that

$$T \;\leq\; 2 \max \big\{\, T_D(x, W \cup W') : x \in D \,\big\} + 1\,. \tag{29.44}$$

To bound the right-hand side of (29.44), we will rely on estimates of the travel time like those developed in section 12, or more sophisticated ones.

**Synthesis.** We shall combine the previous ingredients to show that, whenever the event $\mathcal{D}$ occurs, there exists an index $t$ which is rather small such that the intersection set $\mathcal{I}(t)$ is rather large. The whole crux of the matter is in the meaning of the two "rather" in the previous sentence. So let us precise a bit this mechanism. We first choose $t^*$ so that the probability $P(T \geq t^*)$ is very small and we write

$$P(\mathcal{D}) \;\leq\; P(\mathcal{D}, T \leq t^*) \,+\, P(T \geq t^*)\,. \tag{29.45}$$

We need now to control the probability $P(\mathcal{D}, T \leq t^*)$. Suppose that the event $\mathcal{D}$ occurs and that $T \leq t^*$. On the one hand, combining the inequalities (29.42) and (29.43), we have

$$\big|\Delta_D \mathcal{A}(T-1) \cap \Delta_D \mathcal{A}(T)\big| \;\geq\; \underline{\text{min-cut}}\,(W, W', D, \mathcal{D})\,. \tag{29.46}$$

On the other hand, from the very definition (29.15) of the successive intersection sets $(\mathcal{I}(t), 1 \leq t \leq T)$, we have that

$$\Delta_D \mathcal{A}(T-1) \cap \Delta_D \mathcal{A}(T) \;=\; \bigcup_{1 \leq t \leq T} \mathcal{I}(t)\,. \tag{29.47}$$

Combining (29.46), (29.47) and the fact that $T \leq t^*$, we conclude that

$$\exists\, t \in \{\, 1, \dots, t^* \,\} \qquad \big|\mathcal{I}(t)\big| \,\geq\, \frac{1}{t^*}\underline{\text{min-cut}}\,(W, W', D, \mathcal{D})\,. \tag{29.48}$$

Coming back to (29.45) and using a union bound to control the probability of the event (29.48), we obtain that

$$P(\mathcal{D}) \,\leq\, t^* \max_{1\leq t\leq t^*} P\Big(\big|\mathcal{I}(t)\big| \,\geq\, \frac{1}{t^*}\underline{\text{min-cut}}\,(W, W', D, \mathcal{D})\Big) \,+\, P(T \geq t^*)\,. \tag{29.49}$$

It all comes down to showing that the probability

$$P\big(\big|\mathcal{I}(t)\big| \,\geq\, i\big) \tag{29.50}$$

goes to 0 for $t$ sufficiently small and $i$ sufficiently large. The meaning of "small" and "large" in the previous sentence depends on the parameter governing the size of the problem, typically the side length of $D$ if $D$ is a box, which is sent to $\infty$. Ideally, we wish to obtain a quantitative estimate on the probability (29.50), this would yield a quantitative upper bound on $P(\mathcal{D})$. This last step has remained the main obstacle for proving that $\theta(p_c, \mathbb{Z}^d) = 0$ for approximately fifteen years.

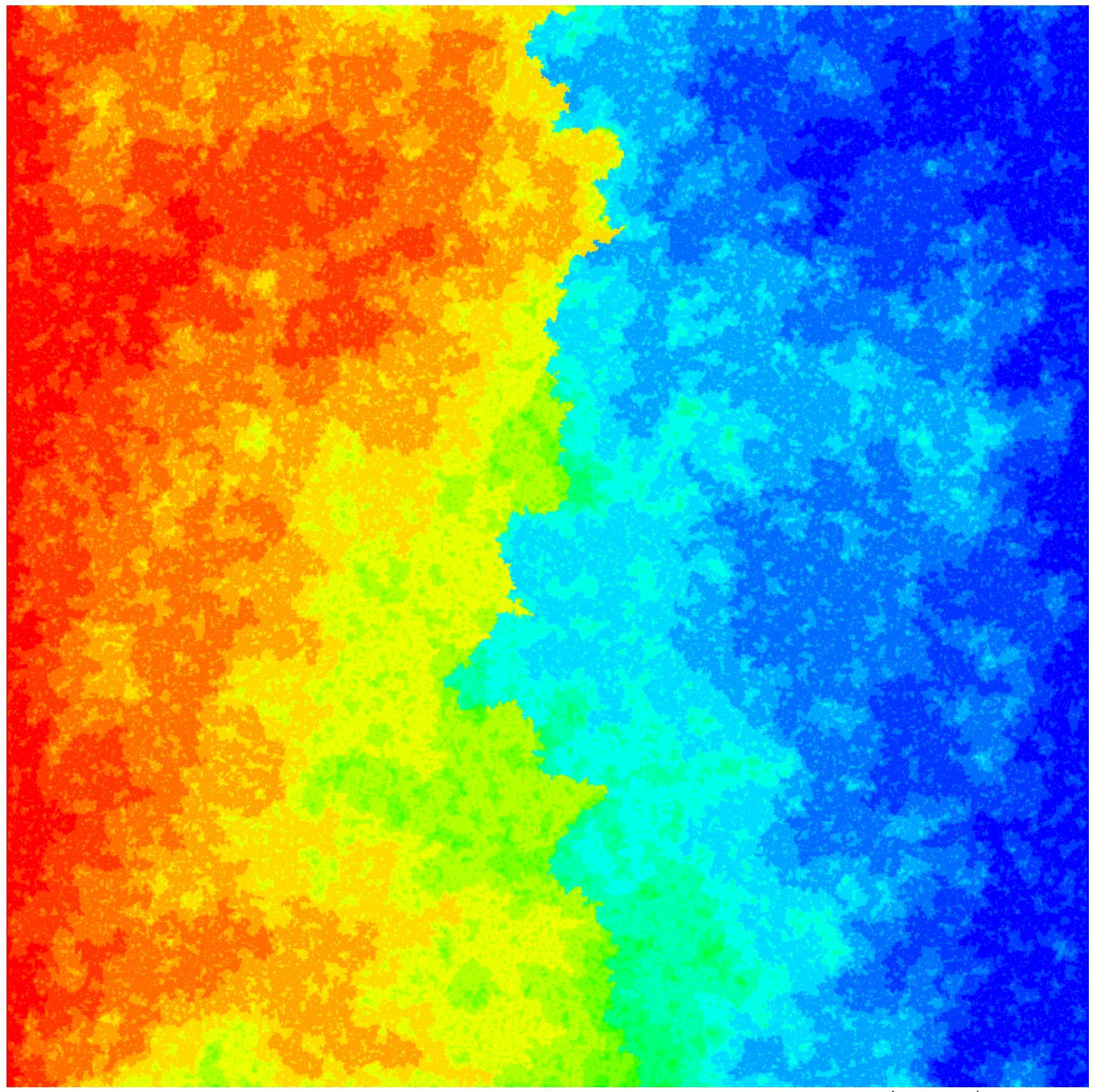

Figure 48: Intertwined explorations, bond percolation in $\Lambda(1024)$, $p = 0.488$, $W, W' =$ left,right boundaries. Lighter colors correspond to higher rounds.

# 30 Three challenging questions

Now the main characters of the story have been introduced. Let us try to put them into action. Ultimately, our goal is to understand better the typical configurations in a finite large box under the sole hypothesis that $\theta(p) > 0$. To solve the $\theta(p_c, \mathbb{Z}^d) = 0$ conjecture, we have to prove that, with probability going to 1 as $n$ goes to $\infty$, inside the box $\Lambda(n)$, there exists a large cluster of cardinality larger than $(\theta - \varepsilon)|\Lambda(n)|$. Maybe it is too ambitious to attack directly this question, and it might be wise to address first more humble issues. In the next three subsections, we describe three reasonable playing fields, each of them is associated with a delicate question that we believe can play a key role in solving the problem. The first question is to control the probability that a sub-box of intermediate size is not connected to the inner boundary of a larger surrounding box. The second question is to obtain some quantitative estimate on the size of the largest cluster in a finite box. The third question is to control the probability that, in a rectangular prism, two squares centered on opposite faces are not connected. In each case, we show how the question can be reduced to the problem of estimating the cardinality of the intersection sets of an adequate intertwined exploration algorithm. These three questions are simpler than the ultimate goal presented in section 9, yet they seem to be very difficult to tackle. Among the three questions, the third one seem to be the less difficult one, so it is sensible to concentrate our efforts on it in the first instance.

## 30.1 Disconnection of a sub-box

Let $\Lambda$ be a cubic box and let $\Gamma$ be a cubic sub-box of $\Lambda$. We consider the disconnection event

$$\mathcal{D}(\Gamma, \Lambda) = \big\{ \Gamma \not\longleftrightarrow \partial^{\,in}\Lambda \big\} ,$$

that we started to discuss at the end of section 7. We would like to have an upper bound on the probability of $\mathcal{D}(\Gamma, \Lambda)$. A crude upper bound on $P_p\big(\mathcal{D}(\Gamma, \Lambda)\big)$ involving $\varrho(\Gamma, \Lambda)$ was derived in (7.15), and it was rewritten in (7.16) with the help of the mean of the maximal cardinality of a cluster in $\Lambda$. Unfortunately, we have no quantitative control on either of these quantities at the critical point $p_c$. So we will try another method to control $P_p\big(\mathcal{D}(\Gamma, \Lambda)\big)$. Suppose that $p$ is such that $\theta(p) > 0$. Let $m \leq n$ be two integers. For a vertex $x$ of $\Lambda(n)$, the box $x + \Lambda(m)$ is included in $\Lambda(2n)$, and we would like to bound from above the probability of the event $\mathcal{D}\big(x + \Lambda(m), \Lambda(2n)\big)$. So we fix $x$ in $\Lambda(n)$, and we start as follows:

$$\begin{aligned} P_p\big(\, x + \Lambda(m) \not\longleftrightarrow \partial^{\,in}\Lambda(2n) \,\big) \,&\leq\, P_p\big(\, x + \Lambda(m) \not\longleftrightarrow \infty \,\big) \\ &\leq\, P_p\big(\, \Lambda(m) \not\longleftrightarrow \infty \,\big) \,. \end{aligned}$$

The good point is that this last upper bound does not depend any more or $n$ nor on $x$. Taking the supremum of the left-hand side with respect to $x$, we have

$$\sup_{x\in\Lambda(n)} P_p\big(\,x+\Lambda(m)\not\longleftrightarrow \partial^{\,in}\Lambda(2n)\,\big)\;\leq\; P_p\big(\,\Lambda(m)\not\longleftrightarrow \infty\,\big)\,. \tag{30.1}$$

We know that

$$\theta(p)>0\qquad\Longrightarrow\qquad \lim_{m\to\infty} P_p\big(\,\Lambda(m)\not\longleftrightarrow\infty\,\big)\;=\;0\,.$$

In the supercritical regime $p > p_c$, we can prove with the help of block estimates and Pisztora's renormalization scheme [114] that there exists a positive constant $c = c(p, d)$ depending on $p$ and the dimension $d$ only such that

$$\exists m_0\geq 1\quad \forall m\geq m_0\qquad P_p\big(\,\Lambda(m)\not\longleftrightarrow\infty\,\big)\;\leq\;\exp\big(-cm^{d-1}\big)\,. \tag{30.2}$$

The rationale behind this bound is that the disconnection between $\Lambda(m)$ and $\infty$ is typically realized by a set of closed bonds separating $\Lambda(m)$ and $\infty$ which is concentrated near the boundary of $\Lambda(m)$. Substituting inequality (30.2) into (30.1) and using a standard union bound, we conclude that, for $m \geq m_0$,

$$P_p\Big(\,\exists\, x\in\Lambda(n)\quad x+\Lambda(m)\not\longleftrightarrow \partial^{\,in}\Lambda(2n)\,\Big)\;\leq\;|\Lambda(n)|\exp\big(-cm^{d-1}\big)\,.$$

This estimate is very informative. We can take $m$ to be equal to

$$m(n)\;=\;\Big(\frac{3d}{c}\ln n\Big)^{1/(d-1)}$$

and we obtain an upper bound going to 0:

$$P_p\Big(\,\exists\, x\in\Lambda(n)\quad x+\Lambda\big(m(n)\big)\not\longleftrightarrow \partial^{\,in}\Lambda(2n)\,\Big)\;\leq\;\frac{1}{n^d}\,. \tag{30.3}$$

This implies that, typically, every sub-box of $\Lambda(n)$ of side length larger than $m(n)$ is connected to $\partial^{\,in}\Lambda(2n)$. Moreover, if we increased slightly the length $m(n)$, the upper bound in (30.3) would go to 0 much faster. The trouble is that we relied on the hypothesis $p > p_c$ to derive the inequality (30.3). If the conjecture $\theta(p_c, \mathbb{Z}^d) = 0$ is true and $p$ is such that $\theta(p, \mathbb{Z}^d) > 0$, then we should in principle be able to reach the same conclusion starting with the sole hypothesis that $\theta(p, \mathbb{Z}^d) > 0$. This would be a very precious piece of information to prove the conjecture $\theta(p_c, \mathbb{Z}^d) = 0$. This type of result could be used to control efficiently the holes in the infinite cluster. So, the challenging question is the following. Suppose that $p$ is such that $\theta(p) > 0$. Can we obtain an explicit sequence $m(n)$ and a quantitative upper bound on

$$P_p\Big(\,\exists\, x\in\Lambda(n)\quad x+\Lambda\big(m(n)\big)\not\longleftrightarrow \partial^{\,in}\Lambda(2n)\,\Big)\,?$$

We already tried to attack this question in subsection 7.4, but we did not obtain any quantitative results. Let us see now how we can implement the strategy outlined in subsection 29.7. To this end, we specify next the adequate choices for the sets $D, W, W'$, the event $\mathcal{D}$, the formulation of the Min-cut problem, and how to control the final iteration number.

**Min-cut problem.** We take $D = \Lambda(2n)$ and $W = \partial^{\,in}\Lambda(2n)$. We introduce another parameter $m \leq n$ (which will depend on $n$ in the end) and we fix a sub-box $\Gamma$ of $\Lambda(2n)$ which is a translate of $\Lambda(m)$ centered at a vertex of $\Lambda(n)$. Finally we take $W' = \Gamma$. The disconnection event $\mathcal{D}$ is as usual $\mathcal{D} = \{\, W \not\longleftrightarrow W' \,\}$. We consider the associated min-cut problem (29.41). For the lower bound, we take

$$\underline{\text{min-cut}}\,(\partial^{\,in}\Lambda(2n), \Gamma, \Lambda(2n), \mathcal{D}) \;=\; 2dm^{d-1}\,.$$

Indeed, a cut $E$ between $\Gamma$ and $\partial^{\,in}\Lambda(2n)$ must contain at least one edge along each segment joining $\Gamma$ and $\partial^{\,in}\Lambda(2n)$ and parallel to the axes, so its cardinality must satisfy $|E| \geq 2dm^{d-1}$.

**Control on the final iteration number.** We consider the intertwined explorations $(\mathcal{A}(t), \mathcal{W}(t), t \geq 0)$ inside $\Lambda(2n)$ starting from $\partial^{\,in}\Lambda(2n)$ and $\Gamma$. Our next goal is to obtain a control on the final iteration number $T$ defined in (29.11). This is done in the following proposition.

**Proposition 30.1.** *Let $p \in [0,1]$ be such that $\theta(p) > 0$. There exists a non-decreasing function $\beta(n)$ defined for $n \geq 6$ such that $\lim_{n\to+\infty} \beta(n) = +\infty$, $\beta(6) \geq \theta(p)$ and the final iteration number $T$ satisfies*

$$\forall n \geq 6 \qquad P\Big(\, T \,>\, 3\frac{\ln n}{\beta(n)} \,\Big) \;\leq\; n^{3d-\beta(n)}\,. \tag{30.4}$$

*Proof.* We know from (29.44) that

$$\begin{aligned} T \;&\leq\; 2\max\big\{\, T_{\Lambda(2n)}\big(x, \partial^{\,in}\Lambda(2n) \cup \Gamma\big) : x \in \Lambda(n) \,\big\} + 1 \\ &\leq\; 2\max\big\{\, T_{\Lambda(2n)}\big(x, \partial^{\,in}\Lambda(2n)\big) : x \in \Lambda(n) \,\big\} + 1\,. \end{aligned} \tag{30.5}$$

The estimate (16.1) obtained in proposition 16.1 provides the required control on the right-hand side of (30.5) and yields the inequality (30.4) of the proposition (we have replaced the function $\beta(n)$ of proposition 16.1 by $\beta(n/2)$). ☐

**Synthesis.** We see from proposition 30.1 that, with very high probability, the final iteration number $T$ is less or equal than $3(\ln n)/\beta(n)$. We proceed exactly as in the synthesis of subsection 29.7. By applying the inequality (29.49) in conjunction with proposition 30.1, we obtain the following result.

**Corollary 30.2.** *Let $p \in [0,1]$ be such that $\theta(p) > 0$. There exists a non-decreasing function $\beta(n)$ defined for $n \geq 6$ such that $\lim_{n\to+\infty} \beta(n) = +\infty$, $\beta(6) \geq \theta(p)$ and, for any $n \geq 6$,*

$$P(\mathcal{D}) \;\leq\; (3\ln n) \max_{1 \,\leq\, t \,\leq\, 3\frac{\ln n}{\beta(n)}} P\Big(\, |\mathcal{I}(t)| \,\geq\, \frac{dm^{d-1}\beta(n)}{\ln n} \,\Big) + n^{3d-\beta(n)}\,. \tag{30.6}$$

So it remains to estimate

$$P\Big(\, |\mathcal{I}(t)| \,\geq\, \frac{dm^{d-1}\beta(n)}{\ln n} \,\Big)$$

for $t \leq 3(\ln n)/\beta(n)$.

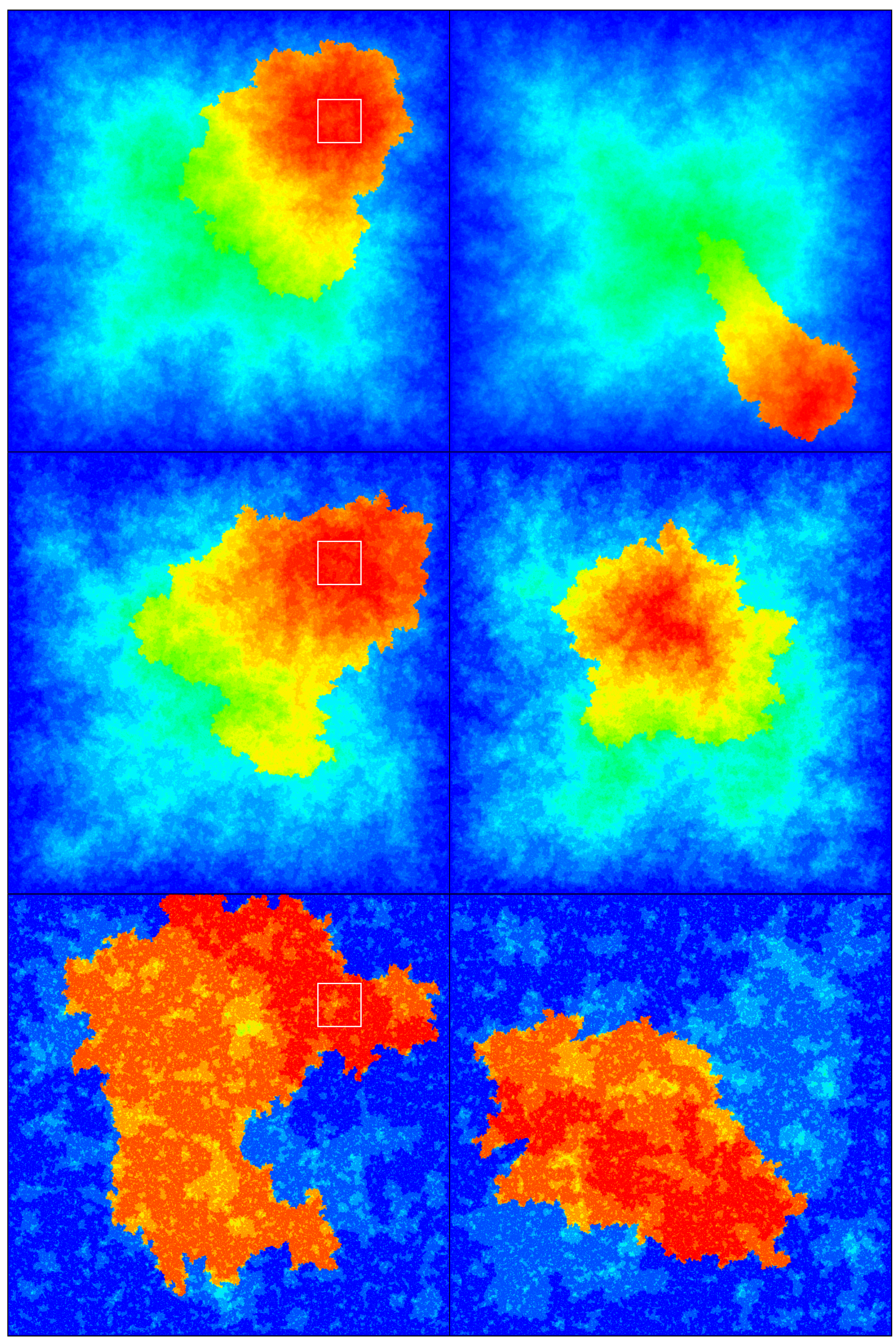

Figure 49: Disconnection from the boundary of the white sub-box (left column) and of the maximal cluster (right column). Bond percolation in a box of size $1024 \times 1024$ with parameters $p = 0.45$, $0.47$, $0.497$ (top row to bottom row).

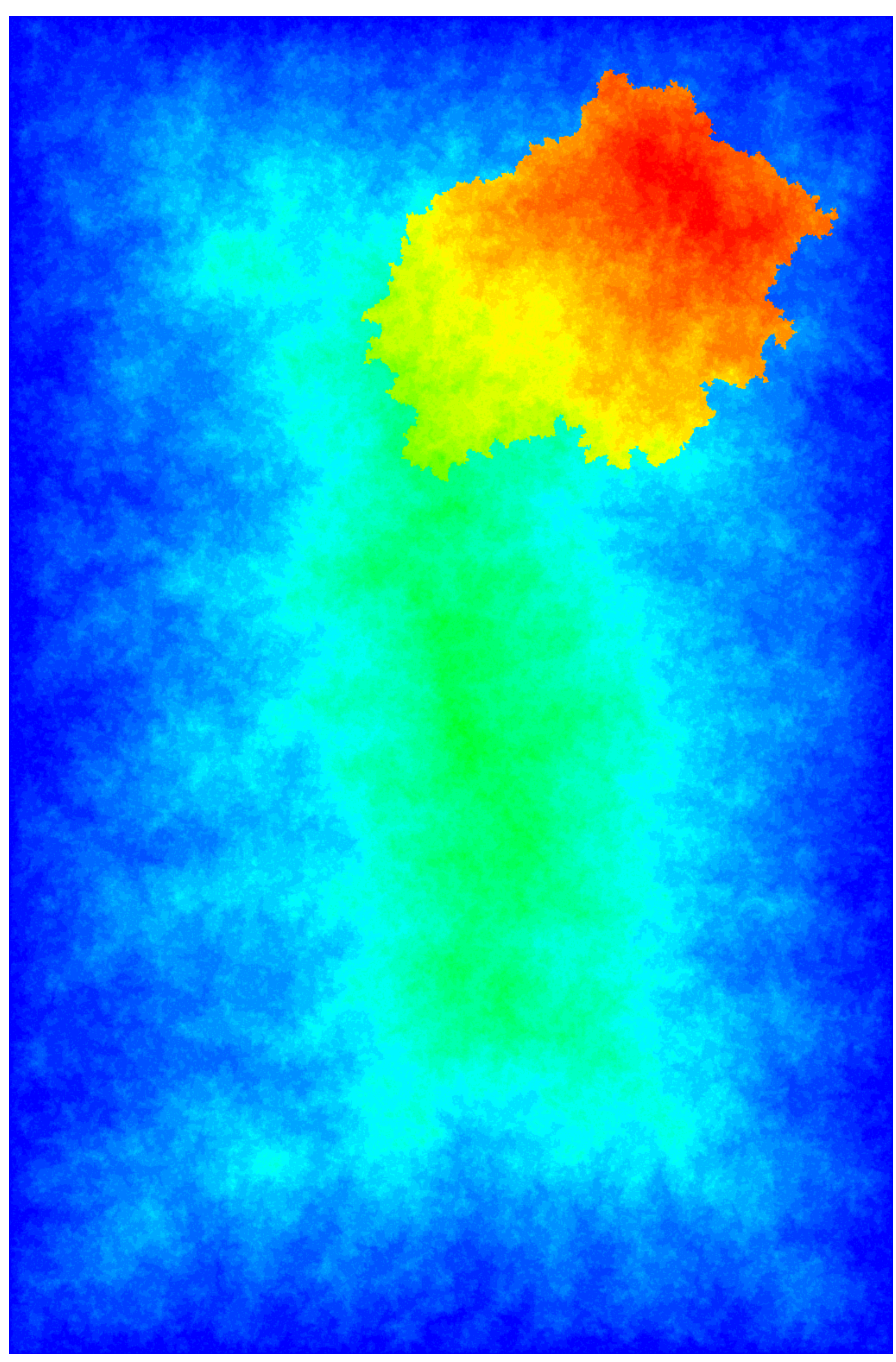

Figure 50: Intertwined exploration in the rectangle $1024 \times 1550$, $p = 0.46$. $W, W' =$ boundary, largest cluster, lighter colors correspond to higher rounds.

## 30.2 The largest cluster in a box

In section 7, we tried to control the variance of the quantity $\varrho(\Gamma, \Lambda)$ and in (7.17) we obtained an inequality involving the maximal cardinality of a cluster in $\Lambda$ which does not intersect $\partial^{in}\Lambda$, i.e.,

$$\max \left\{ |C(x)| : x \in \Lambda, x \not\longleftrightarrow \partial^{in}\Lambda \right\}$$

(there was a square in inequality (7.17), but we can take the square out of the maximum). It is a real challenge to obtain an upper bound on the above quantity. Various interesting results have been obtained by van der Hofstad and Redig [135] on this quantity away from the critical point. To make progress towards a proof of the conjecture $\theta(p_c, \mathbb{Z}^d) = 0$, we would like to have some results in the same direction, but at the critical point $p_c$. Here we shall address the following simpler question. Suppose that $p$ is such that $\theta(p) > 0$. Is it possible to find some exponent $\alpha < d$ and to prove that

$$\lim_{n\to\infty} P_p\Big( \max \left\{ |C(x)| : x \in \Lambda(n), x \not\longleftrightarrow \partial^{in}\Lambda(n) \right\} \geq n^\alpha \Big) = 0\,? \qquad (30.7)$$

In addition, can we obtain a quantitative estimate controlling the speed of convergence to 0 in (30.7)? With this goal in mind, let us see how we can implement the strategy outlined in subsection 29.7.

We take $D = \Lambda(n)$ and $W = \partial^{in}\Lambda(n)$. We introduce another parameter $m$ (which will depend on $n$ in the end) and we fix a vertex $x$ in $\Lambda(n)$. Finally we take $W' = \{ x \}$. The event $\mathcal{D}$ is not taken as a mere disconnection event, here we consider the event $\mathcal{D}(x, m)$ defined by

$$\mathcal{D}(x, m) = \left\{ x \not\longleftrightarrow \partial^{in}\Lambda(n)\,, |C(x)| \geq m^d \right\}.$$

To alleviate the notation, we write simply $\mathcal{D}$ instead of $\mathcal{D}(x, m)$ in the sequel. We consider the min-cut problem (29.41). For the lower bound, we take

$$\underline{\text{min-cut}}\,(\partial^{in}\Lambda(n), \{ x \}, \Lambda(n), \mathcal{D}) = c_{\text{iso}}(d)\, m^{d-1}\,,$$

where $c_{\text{iso}}(d)$ is the constant of the discrete isoperimetric inequality (11.2) associated to the edge boundary. Indeed, thanks to the definition (29.40), we need only to take into account the configurations realizing the event $\mathcal{D}$ when computing the deterministic lower bound (29.42). When $\mathcal{D}$ occurs, a set of closed edges $E$ separating $x$ and $\partial^{in}\Lambda(n)$ must although separate all the vertices of $C(x)$ from $\partial^{in}\Lambda(n)$. Therefore the set $E$ surrounds at least $m^d$ vertices. By the discrete isoperimetric inequality (11.2) on the lattice $\mathbb{Z}^d$, the cardinality of $E$ has to be larger than $c_{\text{iso}}(d)m^{d-1}$. The next steps of the strategy are the same as in subsection 30.1. We obtain the same control of the termination time as in proposition 30.1 and this yields an inequality very similar to (30.6), the only little difference being that the factor $2d$ in front of $m^{d-1}$ has to be replaced by the isoperimetric constant $c_{\text{iso}}(d)$. So, the main difference with the question of subsection 30.1 is that we had to make appeal to the discrete isoperimetric inequality in order to provide a useful lower bound in the min-cut problem.

### 30.3 Disconnection in a cuboid

We consider a right rectangular prism in $\mathbb{Z}^d$ (also called a rectangular cuboid in $\mathbb{Z}^3$) of the form

$$\Pi(n) \;=\; \big\{\, (x_1,\cdots,x_d)\in\mathbb{R}^d : -n/2\leq x_1\leq n/2, \\ \forall i\in\{\,2,\dots,d\,\}\quad -n\leq x_i\leq n\,\big\}\,.$$

We call an hypersquare a $(d-1)$-dimensional symmetric box embedded in $\mathbb{R}^d$. We denote by $L$ (respectively $R$) the hypersquare of side length $n$ centered at $(-n/2,0,\cdots,0)$ (respectively at $(+n/2,0,\cdots,0)$), see figure 51 below.

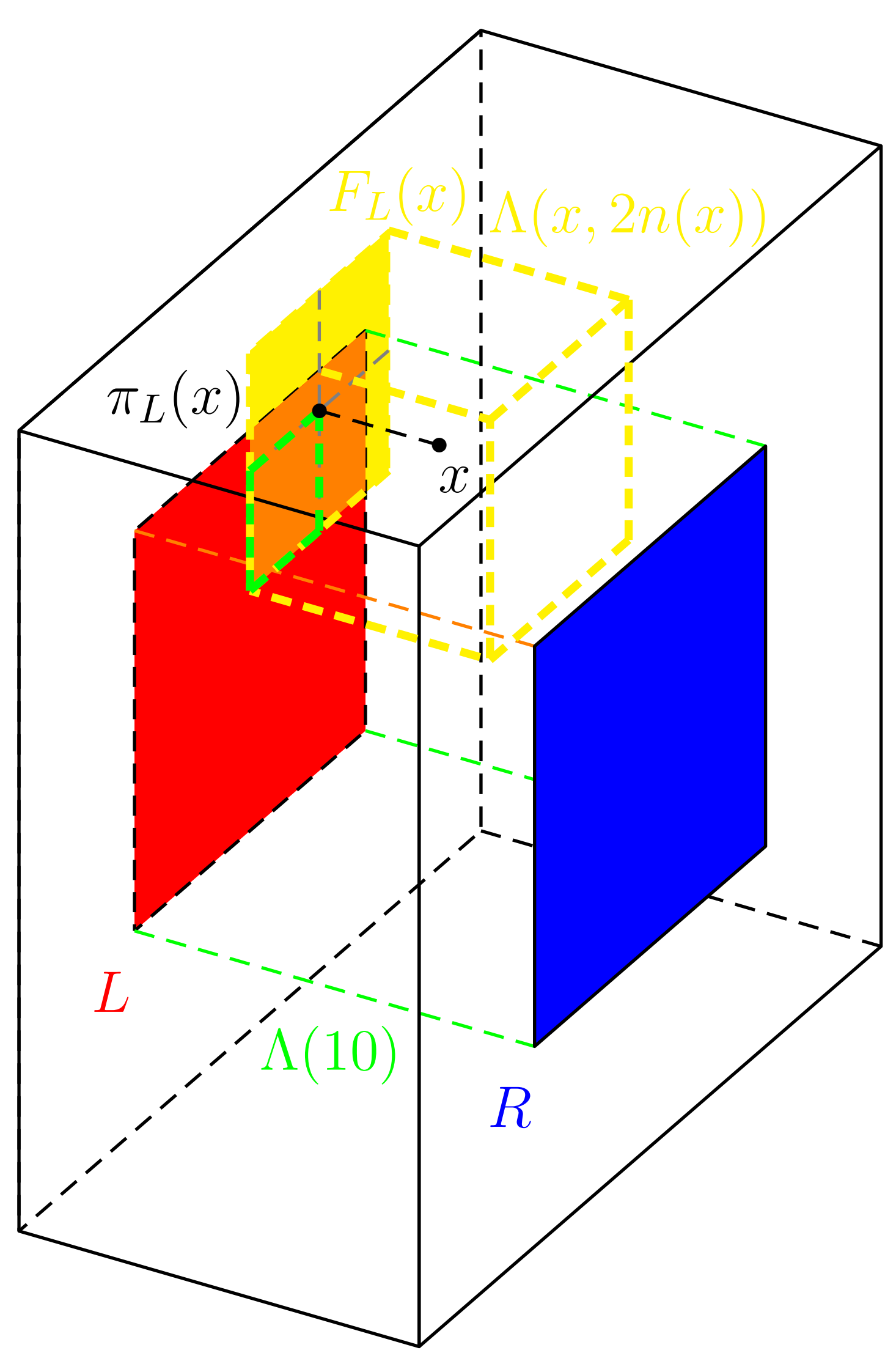

Figure 51: The prism $\Pi(20)$, the box $\Lambda(10)$, left square $L$ and right square $R$. The intersection $L\cap F_L(x)$ contains the hypersquare with the green boundary.

We consider the disconnection event

$$\mathcal{D}(n) \,=\, \big\{\, L \not\longleftrightarrow R \text{ in } \Pi(n) \,\big\}\,.$$

Suppose that $p$ is such that $\theta(p) > 0$. In lemma 8.2, we proved that, with probability going to 1 as $n$ goes to $\infty$, there is a crossing cluster in the box $\Lambda(n)$. The presence of a crossing cluster ensures the existence of a connection between the left and right faces of the box. If we apply this result to the box $\Lambda(n)$, which is included in the prism $\Pi(n)$, we obtain that the event $\mathcal{D}(n)$ becomes very unlikely when $n$ is large, i.e.,

$$\lim_{n\to\infty} P_p\big(\mathcal{D}(n)\big) \,=\, 1\,.$$

In fact, in the supercritical regime $p > p_c$, we can prove with the help of block estimates and Pisztora's renormalization scheme [114] that there exists a positive constant $c = c(p, d)$ depending on $p$ and the dimension $d$ only such that

$$\exists n_0 \geq 1 \quad \forall n \geq n_0 \qquad P_p\big(\mathcal{D}(n)\big) \,\leq\, \exp\big(-cn^{d-1}\big)\,.$$

If the conjecture $\theta(p_c, \mathbb{Z}^d) = 0$ is true and $p$ is such that $\theta(p, \mathbb{Z}^d) > 0$, then we should in principle be able to reach the same conclusion starting with the sole hypothesis that $\theta(p, \mathbb{Z}^d) > 0$. This would be an important step for proving that $\theta(p_c, \mathbb{Z}^d) = 0$, because this type of result could be used to localize efficiently the long connections in the infinite configuration. So, the challenging question is the following. Suppose that $p$ is such that $\theta(p) > 0$. Can we obtain a quantitative estimate controlling the speed of convergence to 0 of $P_p\big(\mathcal{D}(n)\big)$?

Of course, we can ask the same question inside the symmetric box $\Lambda(n)$, but it is more difficult than for the prism. Indeed, the extra space in the prism around the left and right squares opens the possibility to use the square root trick to obtain quantitative estimates on the travel time in $\Pi(n)$ between the vertices of $\Lambda(n)$ and $L, R$, as we shall explain shortly.

**Min-cut problem.** We take $D = \Pi(n)$, $W = L$ and $W' = R$. We consider the min-cut problem (29.41). For the lower bound, we take

$$\underline{\text{min-cut}}\,\big(L, R, \Pi(n), \mathcal{D}(n)\big) \,=\, n^{d-1}\,.$$

Indeed, a cut $E$ between $L$ and $R$ in $\Pi(n)$ contains at least one edge along each segment parallel to the first axis joining $L$ and $R$, so its cardinality $|E|$ has to be larger or equal than the number of these segments, which is at least $n^{d-1}$.

**Intertwined explorations.** We consider next the intertwined explorations $(\mathcal{A}(t), \mathcal{W}(t), t \geq 0)$ inside the prism $\Pi(n)$ starting from $W = L$ and $W' = R$. Instead of the final iteration number $T$, we shall consider the iteration number $T'$ of the intertwined explorations when all the sites of $\Lambda(n)$ have been visited:

$$T' \,=\, \min\big\{\, t \geq 1 : \Lambda(n) \subset \mathcal{A}(t-1) \cup \mathcal{A}(t) \,\big\}\,. \tag{30.8}$$

The introduction of $T'$ is the major difference with the strategy outlined in subsection 29.7. We do this because we do not have yet a control on the final iteration number $T$ defined in (29.11). Since $\Pi(n)$ is finite and connected, we have, thanks to corollary 29.7,

$$\Pi(n) \,=\, \mathcal{A}(T-1) \cup \mathcal{A}(T)\,,$$

and this readily implies that $T' \leq T$. We develop next a specific technique to control $T'$, which relies on the symmetries of the model and the square root trick. However, this technique does not work for the larger number $T$.

**Control on the iteration number $T'$.** Let $x = (x_1, \dots, x_d)$ be a vertex in the box $\Lambda(n)$. The distance between $x$ and $\partial^{in}\Pi(n)$ depends only on the first coordinate $x_1$ and is equal to

$$n(x) \,=\, d\big(x, \partial^{in}\Pi(n)\big) \,=\, \begin{cases} \dfrac{n}{2} + x_1 & \text{if } x_1 \leq 0\,, \\ \dfrac{n}{2} - x_1 & \text{if } x_1 > \dfrac{n}{2}\,. \end{cases}$$

We consider the box

$$\Lambda(x, 2n(x)) \,=\, x + \Lambda(2n(x))\,,$$

which is the largest symmetric box centered at $x$ and contained in $\Pi(n)$. Let us denote by $F_L(x)$ the left face of $\Lambda(x, 2n(x))$. Suppose for instance that we are in the case $x_1 \leq 0$. We have then $n(x) = n/2 + x_1$ and $F_L(x)$ is a subset of the left face of $\Pi(n)$, although it is not fully included in $L$. Yet this left face $F_L(x)$ is a hypersquare of side length $2n(x)$, and it is the union of $2^{d-1}$ hypersquares of side length $n(x)$. Moreover, at least one of these hypersquares is fully included in $L$ (see the example on figure 51). Indeed, let $\Lambda_{d-1}(n)$ be the hypersquare defined by

$$\Lambda_{d-1}(n) \,=\, \big\{\, (x_1, \cdots, x_d) \in \mathbb{R}^d : x_1 = 0,\, \forall i \in \{\, 2, \dots, d \,\} \;\; -n/2 \leq x_i \leq n/2 \,\big\}\,,$$

and let $\pi_L(x)$ be the orthogonal projection of $x$ on $L$. We have

$$F_L(x) \,=\, \pi_L(x) + \Lambda_{d-1}(n(x))\,, \qquad L \,=\, \pi_L(0) + \Lambda_{d-1}(n)\,,$$

whence

$$F_L(x) \cap L \,=\, \big(\pi_L(0) + \Lambda_{d-1}(n)\big) \cap \big(\pi_L(x) + \Lambda_{d-1}(n(x))\big)\,.$$

Suppose for instance that $x$ is such that $x_2 \geq 0, \dots, x_d \geq 0$. In this case, the set $F_L(x) \cap L$ contains

$$\Big\{\, y = (y_1, \cdots, y_d) \in \mathbb{R}^d : y_1 = -\frac{n}{2},\, \forall i \in \{\, 2, \dots, d \,\} \quad x_i - \frac{n}{2} \leq y_i \leq x_i \,\Big\}\,.$$

In the general case, where the signs of $x_2, \dots, x_d$ are arbitrary, we would define

$$\begin{aligned} I^- \,&=\, \big\{\, i \in \{\, 2, \dots, d \,\} : x_i < 0 \,\big\}\,, \\ I^+ \,&=\, \big\{\, i \in \{\, 2, \dots, d \,\} : x_i \geq 0 \,\big\}\,, \end{aligned}$$

and $F_L(x) \cap L$ would contain the hypersquare

$$\Big\{ y = (y_1, \cdots, y_d) \in \mathbb{R}^d : y_1 = -\frac{n}{2},\\ \forall i \in I^- \quad x_i \leq y_i \leq x_i + \frac{n}{2}, \quad \forall i \in I^+ \quad x_i - \frac{n}{2} \leq y_i \leq x_i \Big\}.$$

So, let us denote by $F_L^*(x)$ one of the $2^{d-1}$ hypersquares of $F_L(x)$ which is fully included in $L$. We apply the square root trick result of corollary 14.6 in the box $\Lambda(n(x))$ to get

$$\forall t \geq 0 \qquad P\big( T_{\Lambda(x,n(x))}(x, F_L^*(x)) \,\geq\, t \big) \,\leq\, \exp\big( -t\alpha(t)/(d2^d) \big) . \tag{30.9}$$

Since $\Lambda(x, n(x))$ is included in the prism $\Pi(n)$ and $F_L^*(x)$ is included in $L$, then (30.9) implies that

$$\forall t \geq 0 \qquad P\big( T_{\Pi(n)}(x, L) \,\geq\, t \big) \,\leq\, \exp\big( -t\alpha(t)/(d2^d) \big) .$$

This inequality holds whenever $x_1 \leq 0$. Of course, when $x_1 \geq 0$, we can derive a similar inequality but for the travel time between $x$ and $R$. Putting these two inequalities together, we conclude that

$$\forall x \in \Lambda(n) \quad \forall t \geq 0 \qquad P\big( T_{\Pi(n)}(x, L \cup R) \,\geq\, t \big) \,\leq\, \exp\big( -t\alpha(t)/(d2^d) \big) . \tag{30.10}$$

The good news is that the estimate (30.10) holds uniformly for all the vertices $x$ of $\Lambda(n)$. It is more difficult to obtain such a uniform estimate if we restrict ourselves to the box $\Lambda(n)$. Indeed, the extra space around $\Lambda(n)$ in the prism $\Pi(n)$ has allowed us to use the square root trick through corollary 14.6. We will come back to that point later. For the time being, we use (30.10) and a simple union bound to get

$$\forall t \geq 0 \qquad P\big( \exists\, x \in \Lambda(n) \quad T_{\Pi(n)}(x, L \cup R) \,\geq\, t \big) \,\leq\, \big|\Lambda(n)\big| \exp\big( -t\alpha(t)/(d2^d) \big) .$$

Proceeding as in the proof of proposition 16.1, we obtain the following result.

**Corollary 30.3.** *Let $p \in [0,1]$ be such that $\theta(p) > 0$. There exists a non-decreasing function $\beta(n)$ defined for $n \geq 3$ such that*

$$\forall n \geq 3 \qquad \beta(n) \leq \sqrt{\ln n} \,,$$
$$\beta(3) \geq \theta(p) \,, \qquad \lim_{n \to +\infty} \beta(n) = +\infty \,,$$

*and moreover*

$$\forall n \geq 3 \qquad P\Big( \exists\, x \in \Lambda(n) \quad T_{\Pi(n)}(x, L \cup R) \,>\, \frac{\ln n}{\beta(n)} \Big) \,\leq\, n^{2d - \beta(n)} .$$

Corollary 30.3 in turn yields a control on the random number $T'$.

**Corollary 30.4.** *Let $p \in [0,1]$ be such that $\theta(p) > 0$. There exists a non-decreasing function $\beta(n)$ defined for $n \geq 3$ such that $\lim_{n\to+\infty} \beta(n) = +\infty$, $\beta(3) \geq \theta(p)$ and the iteration number $T'$ defined in* (30.8) *satisfies*

$$\forall n \geq 3 \qquad P\Big(\,T' \,>\, 3\frac{\ln n}{\beta(n)}\,\Big) \,\leq\, n^{2d-\beta(n)}\,. \tag{30.11}$$

*Proof.* It follows from proposition 29.6 that

$$\mathcal{A}(2t) \cup \mathcal{A}(2t+1) \,=\, \mathcal{B}_{\Pi(n)}(R \cup L, t)\,.$$

Therefore, if $t$ is such that

$$t \,\geq\, \max\big\{\, T_{\Pi(n)}(x, L \cup R) : x \in \Lambda(n)\,\big\}\,,$$

then certainly $\Lambda(n)$ is included in $\mathcal{A}(2t) \cup \mathcal{A}(2t+1)$. We conclude that

$$T' \,\leq\, 2\max\big\{\, T_{\Pi(n)}(x, L \cup R) : x \in \Lambda(n)\,\big\} + 1\,.$$

The estimate obtained in corollary 30.3 implies the inequality (30.11) and the statement of the corollary. ☐

**Synthesis.** Let us consider the aggregates $\mathcal{A}(T'-1), \mathcal{A}(T')$ obtained after the iteration $T'$ of the intertwined explorations. These two aggregates are disjoint, and they are such that $\Pi(n) \subset \mathcal{A}(T'-1) \cup \mathcal{A}(T')$ and

$$\text{either} \quad L \subset \mathcal{A}(T'-1),\ R \subset \mathcal{A}(T') \quad \text{or} \quad R \subset \mathcal{A}(T'-1),\ L \subset \mathcal{A}(T')\,.$$

Therefore the set of bonds $\mathbb{E}^d(\Lambda(n)) \cap \Delta_D\mathcal{A}(T'-1) \cap \Delta_D\mathcal{A}(T')$ is a set of closed bonds which separates the sets $R$ and $L$ in $\Lambda(n)$. In particular, it is an admissible candidate for the min-cut problem (29.41) and this implies that

$$\text{Min-Cut}\,\big(R, L, \Lambda(n), \mathcal{D}(n)\big) \,\leq\, \big|\mathbb{E}^d(\Lambda(n)) \cap \Delta_D\mathcal{A}(T'-1) \cap \Delta_D\mathcal{A}(T')\big|\,.$$

The final steps are the same as in the previous subsections, and in the end we obtain the following result.

**Corollary 30.5.** *Let $p \in [0,1]$ be such that $\theta(p) > 0$. There exists a non-decreasing function $\beta(n)$ defined for $n \geq 3$ such that $\lim_{n\to+\infty} \beta(n) = +\infty$, $\beta(3) \geq \theta(p)$, for any $n \geq 3$,*

$$P(\mathcal{D}(n)) \,\leq\, 3(\ln n) \max_{1 \,\leq\, t \,\leq\, 3\frac{\ln n}{\beta(n)}} P\Big(\,|\mathcal{I}(t)| \,\geq\, \frac{n^{d-1}\beta(n)}{2\ln n}\,\Big) + n^{2d-\beta(n)}\,.$$

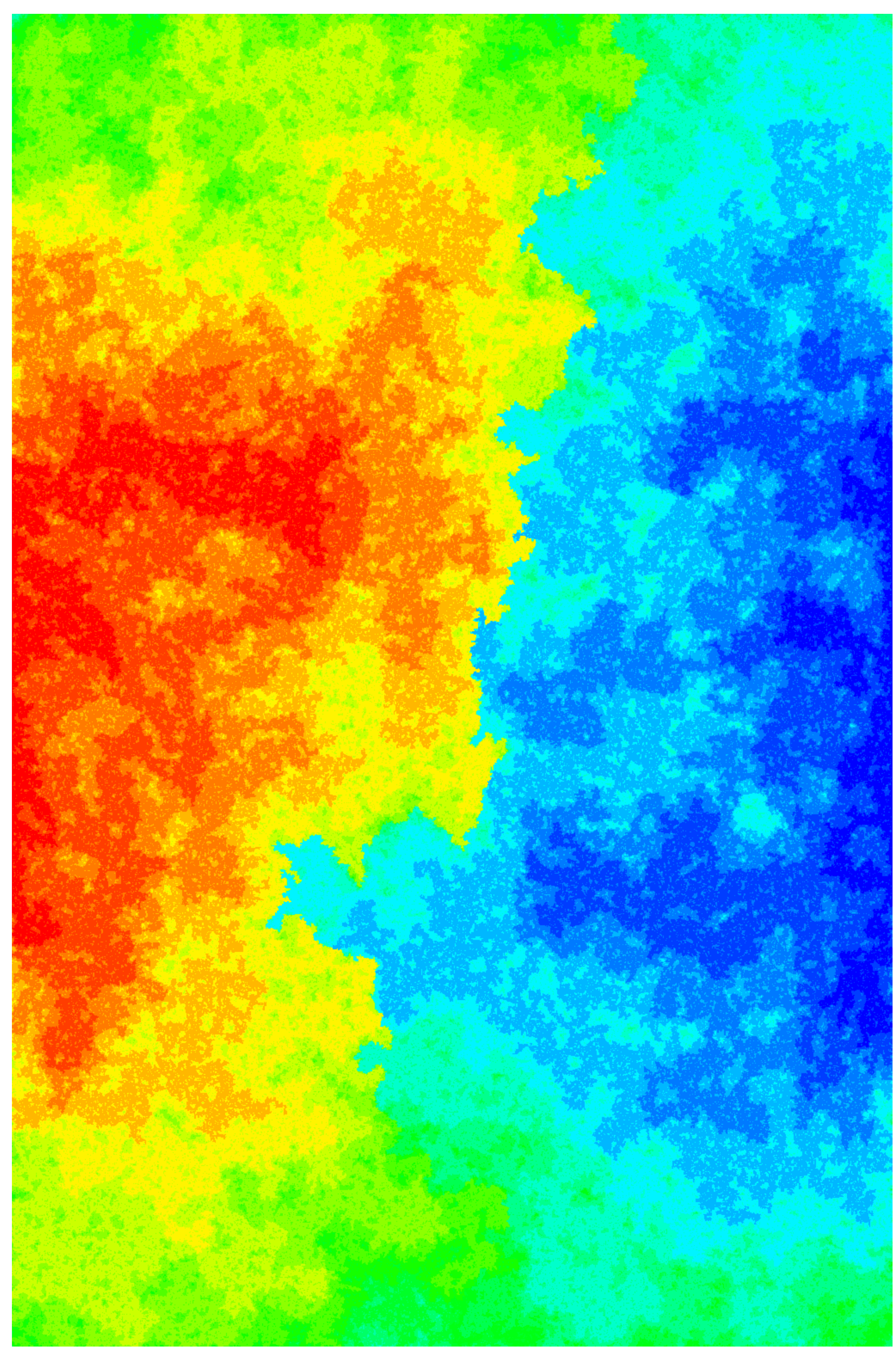

Figure 52: Disconnection between the middle thirds of left and right sides in the rectangle $1550 \times 1024$, $p = 0.49$, intertwined explorations of bond percolation.

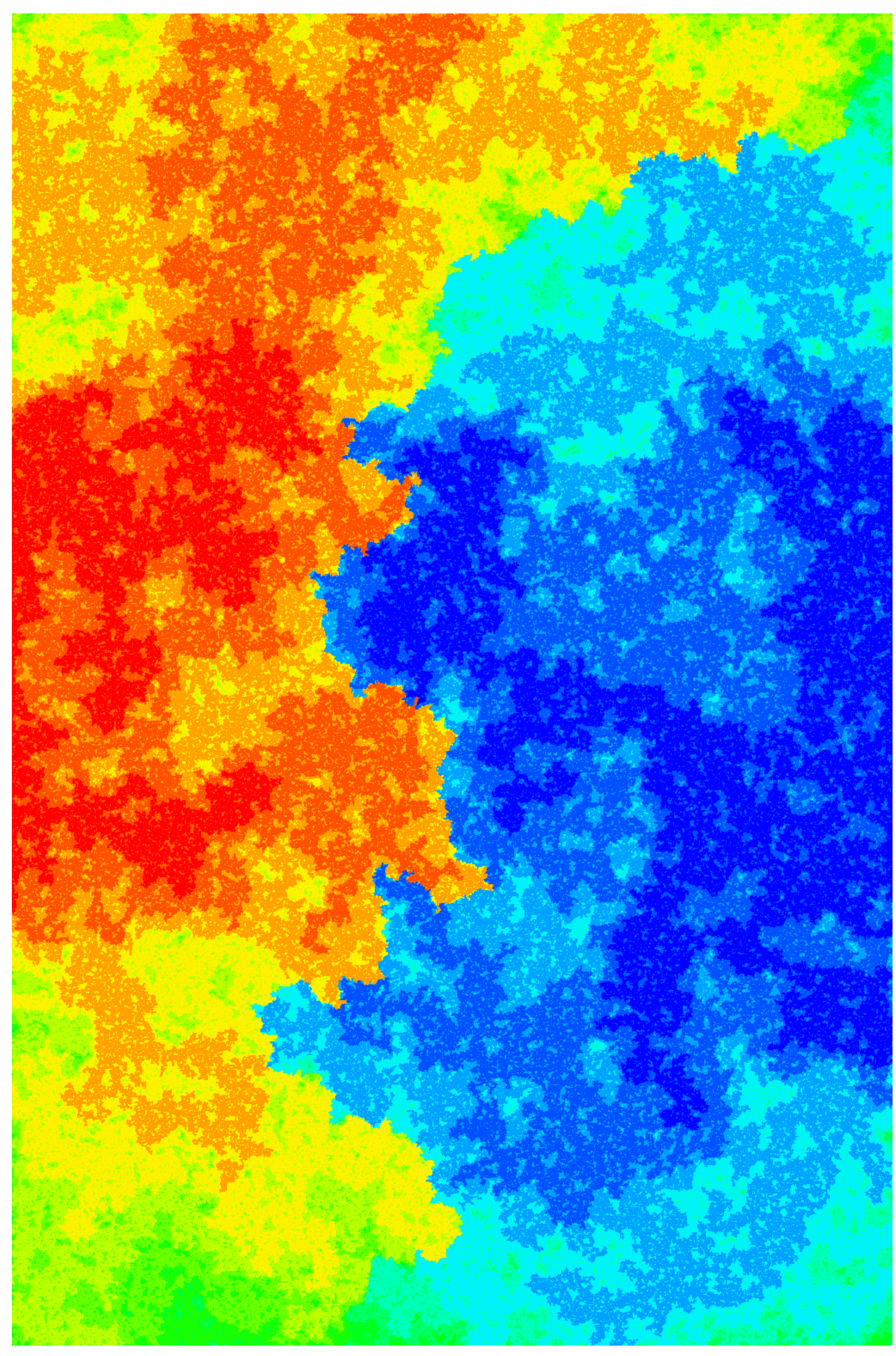

Figure 53: Same as left picture, but with $p = 0.494$.

# 31 Brutal control of bond intersections

In this section, we start to study the intersections, and we try to derive quantitative estimates on their sizes. We start with the first intersections in subsection 31.1, then we attack the second intersections in subsection 31.2. Afterwards, we implement a strategy to control the intersection sets of higher order. To achieve this goal, we provide an alternative construction of the aggregates built by the intertwined explorations in subsection 31.3. This construction involves a taboo growth exploration algorithm, which is presented in subsection 31.5. In subsection 31.4, we show that the successive taboo sets can be restricted with the help of the intersection sets. This provides a mean to derive a general inequality for controlling the successive intersections, which is presented in subsection 31.6.

## 31.1 Control of the first intersection

Let $W$ and $W'$ be two subsets of a finite domain $D$ in $\mathbb{Z}^d$ and let us consider a percolation configuration in $D$ in which $W$ and $W'$ are not connected. Let $(\mathcal{A}(t),\mathcal{W}(t),t\geq 0)$ be the intertwined explorations starting from $W$ and $W'$. We wish to control the cardinality of the first intersection set $\mathcal{I}(1)$. Let $i_1$ be a positive integer. By definition, we have

$$P\big(|\mathcal{I}(1)|\geq i_1, W\not\longleftrightarrow W'\big)\,=\,P\big(\big|\Delta_D\mathcal{A}(0)\cap\Delta_D\mathcal{A}(1)\big|\geq i_1, W\not\longleftrightarrow W'\big)\,.$$

On the event $\{\,W\not\longleftrightarrow W'\,\}$, the bonds of $\Delta_D\mathcal{A}(0)\cap\Delta_D\mathcal{A}(1)$ are the pivotal bonds for the event $\{\,W\longleftrightarrow W'\,\}$. Moreover, these bonds are precisely the bonds of $D$ that are visited by both cluster explorations starting from $W$ and from $W'$. Let us be more precise. We denote by $\mathcal{E}(W,D)$ the set of the bonds that are explored by the shell exploration algorithm which builds $\text{Shell}\,(W,D,0)$, that is

$$\mathcal{E}(W,D)\,=\,\mathbb{E}^d\big(\text{Shell}\,(W,D,0)\cup\partial_D^{out}\text{Shell}\,(W,D,0)\big)\,.\tag{31.1}$$

In words, the set $\mathcal{E}(W,D)$ is the set of the bonds of $D$ which have at least one endpoint in $\text{Shell}\,(W,D,0)$. We consider similarly the set $\mathcal{E}(W',D)$. An example of the sets $\mathcal{E}(W,D)$, $\mathcal{E}(W',D)$ and $\mathcal{I}(1)$ can be seen on figure 54. Next, we have

$$P\big(|\mathcal{I}(1)|\geq i_1, W\not\longleftrightarrow W'\big)\,=\,P\big(\big|\mathcal{E}(W,D)\cap\mathcal{E}(W',D)\big|\geq i_1, W\not\longleftrightarrow W'\big)\,.\tag{31.2}$$

Naturally, we have

$$\mathcal{E}(W,D)\cup\mathcal{E}(W',D)\,=\,\mathcal{E}(W\cup W',D)\,,\tag{31.3}$$

and therefore, using the additivity of the function $S$ defined in (28.1), we have

$$S\big(\mathcal{E}(W\cup W',D)\big)+S\big(\mathcal{E}(W,D)\cap\mathcal{E}(W',D)\big)\,=\,S\big(\mathcal{E}(W,D)\big)+S\big(\mathcal{E}(W',D)\big)\,.\tag{31.4}$$

On the event $\{\,W\not\longleftrightarrow W'\,\}$, the bonds of $\mathcal{E}(W,D)\cap\mathcal{E}(W',D)$ are all closed, therefore

$$S\big(\mathcal{E}(W,D)\cap\mathcal{E}(W',D)\big)\,=\,-\frac{1}{1-p}\big|\mathcal{E}(W,D)\cap\mathcal{E}(W',D)\big|\,.\tag{31.5}$$

Using the probabilistic control (28.2) on the three explorations starting from $W$, $W'$ or $W \cup W'$, the equalities (31.2), (31.3), (31.4) and (31.5), we obtain that

$$P\big(|\mathcal{I}(1)| \geq i_1, W \not\leftrightarrow W'\big) \,\leq\, 6|\mathbb{E}^d(D)| \exp\Big(-\frac{2p^2 i_1^2}{9|\mathbb{E}^d(D)|}\Big)\,. \tag{31.6}$$

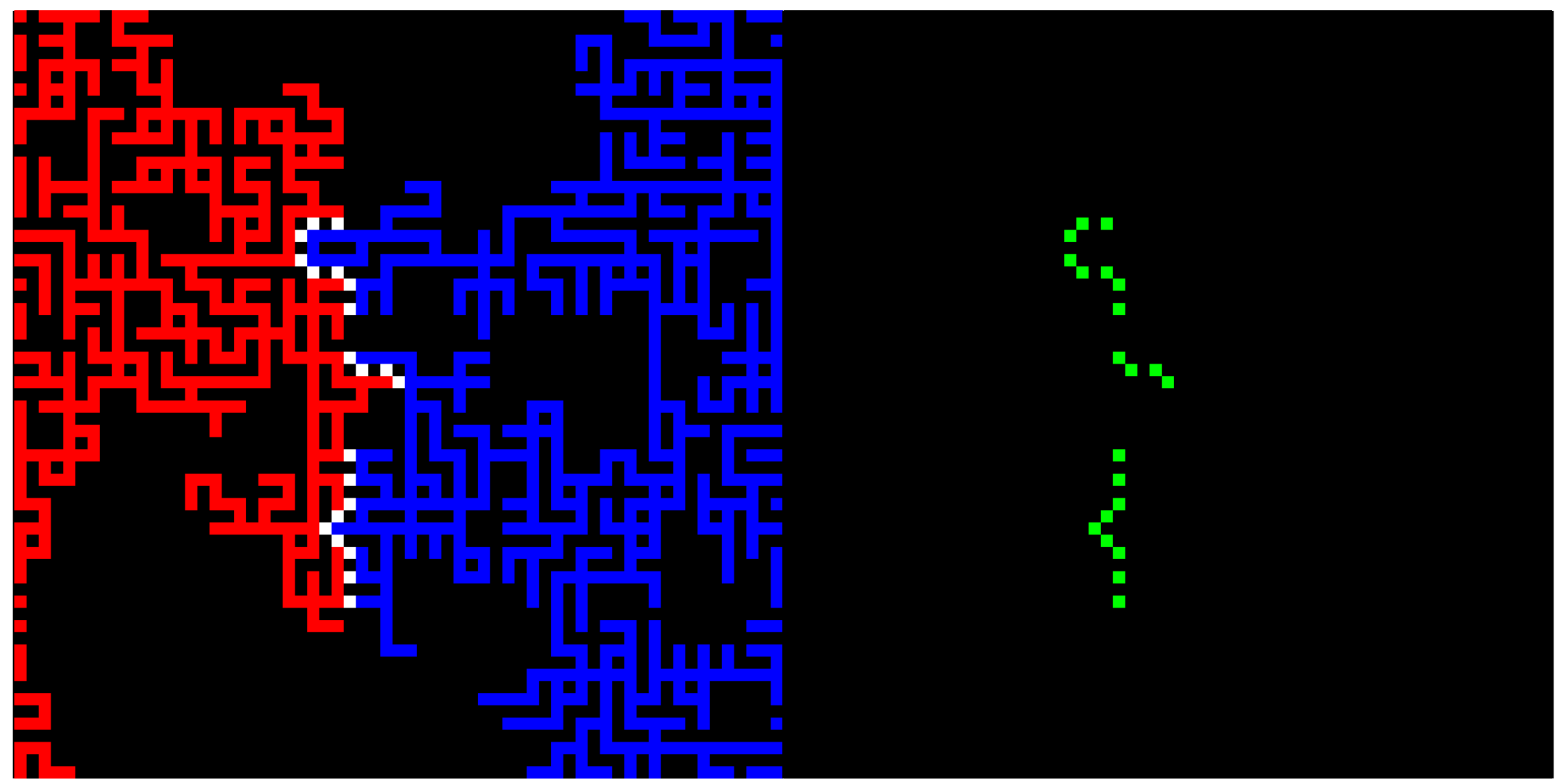

Figure 54: Two rounds of intertwined bond explorations, $D = \Lambda(32)$, $p = 0.5$. The sets of bonds $\mathcal{E}(W, D)$, $\mathcal{E}(W', D)$ are displayed on left, $\mathcal{I}(1)$ on right (only the open bonds belonging to these sets are colored, as well as their endpoints). The bonds of $\mathcal{I}(1)$ are colored in white on the left picture. The initial sets $W$ and $W'$ are the left and the right boundaries of the square.

Let us apply this result with $D = \Lambda(n)$ for some $n \geq 2d$. Using the inequalities

$$|\mathbb{E}^d(\Lambda(n))| \,\leq\, d|\Lambda(n)| \,\leq\, d(n+1)^d \,\leq\, d2^d n^d\,, \tag{31.7}$$

we obtain

$$P\big(|\mathcal{I}(1)| \geq i_1, W \not\leftrightarrow W'\big) \,\leq\, 6d2^d n^d \exp\Big(-\frac{2p^2 i_1^2}{9d2^d n^d}\Big)\,. \tag{31.8}$$

With the choice

$$i_1 \,=\, \frac{3}{p} d2^d \sqrt{n^d \ln n}\,,$$

and using the inequalities $4d \leq d2^d \leq n^d$, we conclude that

$$P\Big(|\mathcal{I}(1)| \geq \frac{3}{p} d2^d \sqrt{n^d \ln n}, W \not\leftrightarrow W'\Big) \,\leq\, 6n^{2d} \exp\big(-4d \ln n\big) \,\leq\, \frac{6}{n^{2d}}\,. \tag{31.9}$$

In the next step, we would like to have a similar control over all possible starting sets $W$ and $W'$. Of course this is too optimistic, so we suppose that there is a control on the cardinalities of $W$ and $W'$. Let $N \geq 1$ be a fixed integer. By the union bound, we have

$$P\begin{pmatrix} \exists W, W' \subset \mathbb{E}^d(D) \quad |W| \leq N\,, \\ |W'| \leq N\,, W \not\leftrightarrow W'\,, |\mathcal{I}(1)| \geq i_1 \end{pmatrix} \leq \sum_{\substack{W,W' \subset \mathbb{E}^d(D) \\ |W| \leq N\,, |W'| \leq N}} P\big(W \not\leftrightarrow W', |\mathcal{I}(1)| \geq i_1\big). \tag{31.10}$$

We would like to find a regime in which $n, N$ and $i_1$ go simultaneously to $\infty$ and such that the right-hand quantity in (31.10) goes to 0. We use the probabilistic control (31.6), and we bound crudely the number of subsets of $\mathbb{E}^d(D)$ of cardinality at most $N$ by $|\mathbb{E}^d(D)|^{N+1}$ to conclude that

$$P\begin{pmatrix} \exists W, W' \subset \mathbb{E}^d(D) \quad |W| \leq N\,, \\ |W'| \leq N\,, W \not\leftrightarrow W'\,, |\mathcal{I}(1)| \geq i_1 \end{pmatrix} \leq 6|\mathbb{E}^d(D)|^{2N+3} \exp\Big(-\frac{2p^2 i_1^2}{9|\mathbb{E}^d(D)|}\Big)\,. \tag{31.11}$$

This inequality is interesting only when $i_1$ is sufficiently large to ensure that the prefactor $6|\mathbb{E}^d(D)|^{2N+1}$ is absorbed in the exponential. More precisely, for $i_1$ such that

$$i_1 \geq \frac{6}{p}\sqrt{N|\mathbb{E}^d(D)| \ln(|\mathbb{E}^d(D)|)}\,, \tag{31.12}$$

we have

$$P\begin{pmatrix} \exists W, W' \subset \mathbb{E}^d(D) \quad |W| \leq N\,, \\ |W'| \leq N\,, W \not\leftrightarrow W'\,, |\mathcal{I}(1)| \geq i_1 \end{pmatrix} \leq 6\exp\Big(-\frac{p^2 i_1^2}{9|\mathbb{E}^d(D)|}\Big)\,.$$

Let us apply this result with $D = \Lambda(n)$ for some $n \geq 2d$. Thanks to the inequalities (31.7), the inequality (31.12) is satisfied with the choice

$$i_1 = \frac{6}{p} d 2^d \sqrt{N n^d \ln n}$$

(we bound $\ln(|\mathbb{E}^d(D)|)$ by $2d \ln n$, indeed $|\mathbb{E}^d(D)| \leq d2^d n^d \leq n^{2d}$ for $n \geq 2d$). This yields that, for $n \geq 2d$, we have

$$P\begin{pmatrix} \exists W, W' \subset \mathbb{E}^d(D) \quad |W| \leq N\,, |W'| \leq N\,, \\ W \not\leftrightarrow W'\,, |\mathcal{I}(1)| \geq \frac{6}{p} d2^d \sqrt{N n^d \ln n} \end{pmatrix} \leq \frac{6}{n^{dN}}\,. \tag{31.13}$$

This inequality tells us that, typically, in the box $\Lambda(n)$, the number of pivotal bonds for the connection event between two subsets of vertices of cardinality less than $N$ is at most of the order $\sqrt{Nn^d \ln n}$. Although we are still very far from what we need, this is certainly a valuable piece of information. An unfortunate nasty fact is the presence of $N$ at the exponential level in the inequality (31.11), because it deteriorates seriously the estimates. To avoid this prefactor which comes from a crude combinatorial bound, in a first approach, we will restrict the possible choices for the sets $W$ and $W'$, and we will deal later with this prefactor. So, for the time being, let us not worry on the presence of $N$.

## 31.2 Control of the second intersections

We try next to control the cardinality of the second intersections, defined by

$$\mathcal{I}(2) \,=\, \Delta_D\mathcal{A}(1)\cap\Delta_D\mathcal{A}(2)\setminus\mathcal{I}(1)\,. \tag{31.14}$$

We start as we did for the first intersections $\mathcal{I}(1)$. Let $W$, $W'$ be two subsets of a finite domain $D$ and let us consider a percolation configuration in $D$ in which $W$ and $W'$ are not connected. Let $(\mathcal{A}(t),\mathcal{W}(t),t\geq 0)$ be the intertwined explorations starting from $W$ and $W'$. We would like to control, for $i_2\geq 1$,

$$P\big(\big|\Delta_D\mathcal{A}(1)\cap\Delta_D\mathcal{A}(2)\setminus\mathcal{I}(1)\big|\geq i_2, W\not\longleftrightarrow W'\big)\,.$$

By construction, on the event $\{\,W\not\longleftrightarrow W'\,\}$, the sets $\mathcal{A}(1)$ and $\mathcal{A}(2)$ are two disjoint aggregates in $D$, therefore all the bonds in $\Delta_D\mathcal{A}(1)\cup\Delta_D\mathcal{A}(2)$ are closed. We would like to realize the set $\Delta_D\mathcal{A}(1)\cup\Delta_D\mathcal{A}(2)\setminus\mathcal{I}(1)$ as the intersection of the outputs of two exploration algorithms. By definition, we have

$$\mathcal{A}(1) \,=\, \text{Clusters}\big(W',D,\mathcal{A}(0)\big)\,.$$

On the event $\{\,W\not\longleftrightarrow W'\,\}$, the taboo set $\mathcal{A}(0)$ plays no role at all in the previous formula, so that in fact

$$\mathcal{A}(1) \,=\, \text{Clusters}\big(W',D,\varnothing\big) \,=\, \text{Shell}\,(W',D)\,,$$

and the set $\mathcal{A}(1)$ is simply obtained by running an exploration algorithm starting from $W'$. We denote by $\mathcal{E}(W',D)$ the set of the bonds explored by this algorithm. The set $\mathcal{A}(2)$ is built in two stages. First we build the sets

$$\begin{aligned}\mathcal{A}(0) \,&=\, \text{Clusters}\big(W,D,\varnothing\big) \,=\, \text{Shell}\,(W,D)\,,\\ \mathcal{W}(0) \,&=\, \text{Waiting}\big(W,D,\varnothing\big) \,=\, \partial_D^{out}\text{Shell}\,\big(W,D\big)\,.\end{aligned}$$

At this point, we have explored the set of bonds $\mathcal{E}(W,D)$ defined in (31.1). In a second step, we run an exploration algorithm from $\mathcal{W}(0)$ with taboo set $\mathcal{A}(1)$. This exploration returns the set

$$\text{Clusters}\big(\mathcal{W}(0),D\setminus\mathcal{A}(0),\mathcal{A}(1)\big) \,=\, \text{Shell}\,(W,D\setminus\mathcal{A}(1),1)\,, \tag{31.15}$$

and we have

$$\mathcal{A}(2) \,=\, \text{Shell}\,(W,D,0)\cup\text{Shell}\,(W,D\setminus\mathcal{A}(1),1)\,. \tag{31.16}$$

Let $\mathcal{E}(W,D,\mathcal{A}(1),1)$ be the set of the bonds that are explored by the shell exploration algorithm which builds $\text{Shell}\,(W,D\setminus\mathcal{A}(1),1)$. In fact, this set is the set of the bonds of $D\setminus\mathcal{A}(0)$ which have at least one endpoint in $\mathcal{A}(2)$, i.e.,

$$\mathcal{E}(W,D,\mathcal{A}(1),1) \,=\, \mathbb{E}^d\big(\mathcal{A}(2)\cup\partial^{out}_{D\setminus\mathcal{A}(0)}\mathcal{A}(2)\big)\,.$$

From the previous discussion, we see that

$$\Delta_D\mathcal{A}(1)\cap\Delta_D\mathcal{A}(2)\setminus\mathcal{I}(1) \,=\, \mathcal{E}(W',D)\,\cap\,\mathcal{E}(W,D,\mathcal{A}(1),1)\,. \tag{31.17}$$

We know that $\mathcal{E}(W', D)$ is the output of a genuine exploration algorithm in $D$, hence we have the usual control on $S(\mathcal{E}(W', D))$. Of course, if $\mathcal{E}(W, D, \mathcal{A}(1), 1)$ was also the output of a genuine exploration algorithm in $D$, then we could proceed as for $\mathcal{I}(1)$, and we would get the same estimate for $\mathcal{I}(2)$ as for $\mathcal{I}(1)$. The big trouble is that, to build $\mathcal{E}(W, D, \mathcal{A}(1), 1)$, we have used $\mathcal{A}(1)$ as taboo set, and it is no longer the case that $\mathcal{E}(W, D, \mathcal{A}(1), 1)$ is the output of a genuine exploration algorithm. Indeed, the exploration algorithm is confined in $D \setminus \mathcal{A}(1)$ and knows that $\mathcal{A}(1)$ is an aggregate which is not reachable from $\mathcal{W}(0) \setminus \mathcal{A}(1)$. As a consequence, any bond which leads inside $\mathcal{A}(1)$ has to be closed. If only the algorithm could forget about this unwanted information, then we could derive a probabilistic estimate like (31.6) and we could proceed. Unfortunately, we can't hide from the truth, and it doesn't seem possible to pretend not to have this information. We have to face the fact that the algorithm knows that certain bonds are closed.

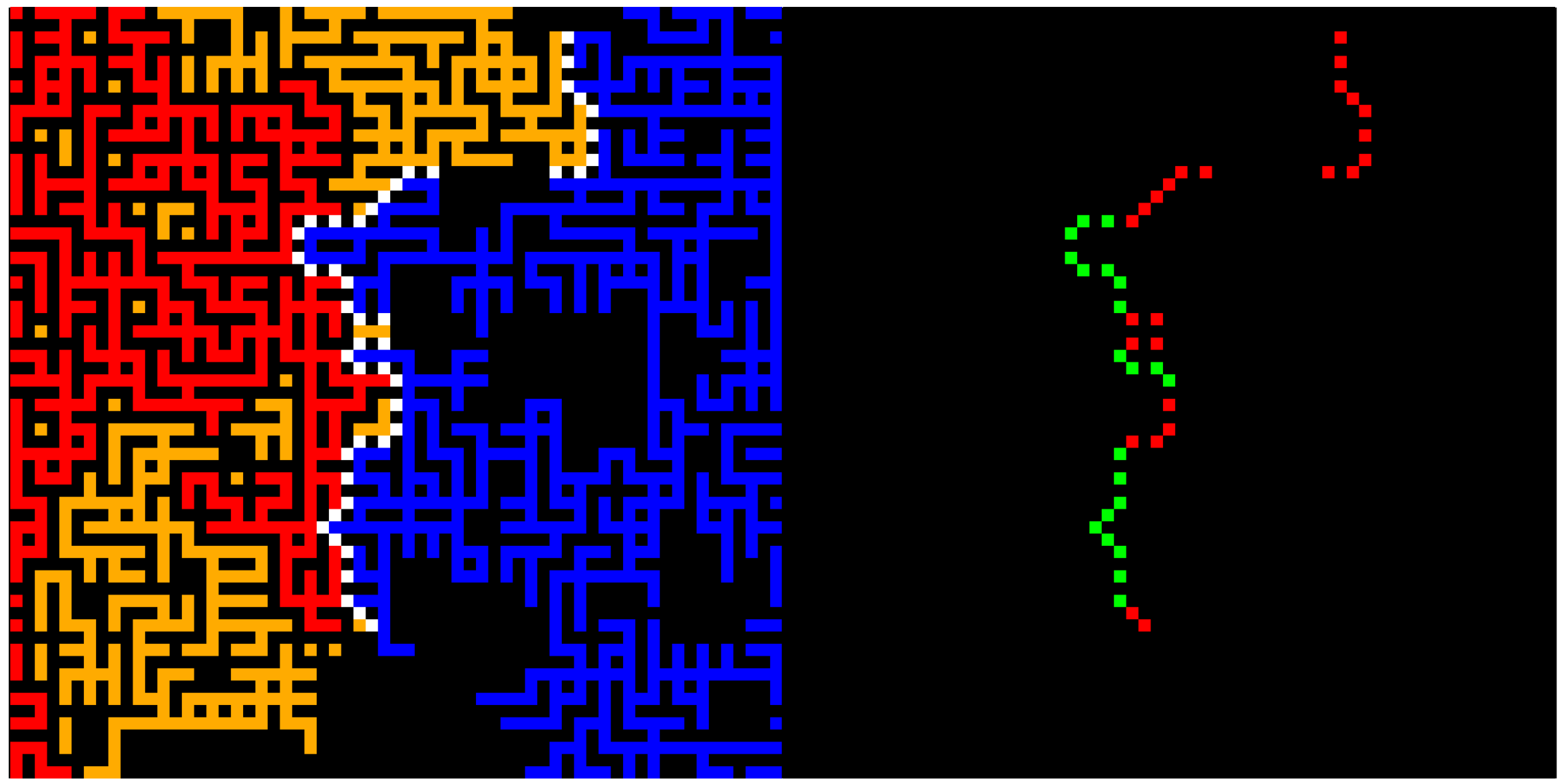

Figure 55: Three rounds of intertwined bond explorations and the set $\mathcal{T}(1)$. Left: $\mathcal{E}(W, D)$, $\mathcal{E}(W', D)$, $\mathcal{E}(W, D, \mathcal{A}(1), 1)$ (only the open bonds belonging to these sets are colored, as well as their endpoints), $\mathcal{I}(1) \cup \mathcal{I}(2)$ in white, $\mathcal{T}(1)$. Right: $\mathcal{I}(1)$, $\mathcal{I}(2)$, $\mathcal{T}(1)$. $D = \Lambda(32)$, $p = 0.5$, $W, W' =$ left,right boundaries.

A natural idea to go around this problem is to perform a conditioning on the taboo set sent by the second explorer to the first. However, this is a bad idea, because at some point we will have to sum over all possible sets to get rid of the conditioning, and this would create a too large combinatorial factor. So we will try to condition on the smallest relevant set.

When the first explorer starts its second round of exploration, he knows that $\mathcal{A}(1)$ is a taboo set, yet only a small fraction of this taboo set will matter for this second round of exploration, namely the vertices of $\mathcal{A}(1)$ which have one endpoint in $\mathcal{I}(1)$. Indeed, these vertices are the only possible entry points into $\mathcal{A}(1)$ when starting an exploration from $\mathcal{W}(0)$. So, if we denote by $\mathcal{T}(1)$ these vertices, i.e.,

$$\mathcal{T}(1) \,=\, \big\{\, x \in \mathcal{A}(1) : \exists y \in D \setminus \mathcal{A}(1) \quad \langle x, y\rangle \in \mathcal{I}(1) \,\big\}\,,$$

then we have that

$$\begin{aligned} \text{Shell}\,(W, D \setminus \mathcal{A}(1), 1) \,&=\, \text{Shell}\,(W, D \setminus \mathcal{T}(1), 1)\,,\\ \mathcal{E}(W, D, \mathcal{A}(1), 1) \,&=\, \mathcal{E}(W, D, \mathcal{T}(1), 1)\,. \end{aligned} \tag{31.18}$$

The situation is already much better when working with the taboo set $\mathcal{T}(1)$, because $|\mathcal{T}(1)| \leq |\mathcal{I}(1)|$ (a bond of $\mathcal{I}(1)$ can give rise to at most one vertex in $\mathcal{T}(1)$), and we have already a control on $|\mathcal{I}(1)|$. So we will perform the conditioning on $\mathcal{I}(1)$ and we write

$$\begin{aligned} &P\big(|\mathcal{I}(2)| \geq i_2, W \not\longleftrightarrow W'\big) \,\leq\\ &P\big(|\mathcal{I}(2)| \geq i_2, |\mathcal{I}(1)| \leq i_1, W \not\longleftrightarrow W'\big) \,+\, P\big(|\mathcal{I}(1)| \geq i_1, W \not\longleftrightarrow W'\big)\,. \end{aligned} \tag{31.19}$$

The last term in (31.19) has already been dealt with, so we focus now on the term before, which we decompose according to $\mathcal{T}(1)$:

$$\begin{aligned} &P\big(|\mathcal{I}(2)| \geq i_2, |\mathcal{I}(1)| \leq i_1, W \not\longleftrightarrow W'\big) \,=\\ &\qquad \sum_{T_1 \subset D, |T_1| \leq i_1} P\big(|\mathcal{I}(2)| \geq i_2, |\mathcal{I}(1)| \leq i_1, \mathcal{T}(1) = T_1, W \not\longleftrightarrow W'\big)\,. \end{aligned} \tag{31.20}$$

We take care of the term inside the sum. Let us fix a set $T_1$. Using the equalities (31.14), (31.17) and (31.18), we obtain

$$\begin{aligned} &P\big(|\mathcal{I}(2)| \geq i_2, |\mathcal{I}(1)| \leq i_1, \mathcal{T}(1) = T_1, W \not\longleftrightarrow W'\big) \,\leq\\ &\qquad P\big(\big|\mathcal{E}(W', D) \,\cap\, \mathcal{E}(W, D, T_1, 1)\big| \geq i_2, \mathcal{T}(1) = T_1, W \not\longleftrightarrow W'\big)\\ &\qquad\qquad \leq\, P\big(\big|\mathcal{E}(W', D) \,\cap\, \mathcal{E}(W, D, T_1, 1)\big| \geq i_2, W \not\longleftrightarrow W'\big)\,. \end{aligned} \tag{31.21}$$

By forgetting the condition $\mathcal{T}(1) = T_1$ and taking the upper bound, we are left with the problem of controlling $\mathcal{E}(W, D, T_1, 1)$. This is now possible, because the set $T_1$ is fixed, and it does not carry with him any information on the percolation configuration in $D \setminus T_1$. As a consequence, the set $\mathcal{E}(W, D, T_1, 1)$ is again the output of a genuine exploration algorithm and the probabilistic control (28.2) is in force. Using the additivity of the function $S$, we have

$$\begin{aligned} S\big(\mathcal{E}(W, D, T_1, 1) \cap \mathcal{E}(W', D)\big) \,&=\, S\big(\mathcal{E}(W, D, T_1, 1)\big) + S\big(\mathcal{E}(W', D)\big)\\ &\qquad - S\big(\mathcal{E}(W, D, T_1, 1) \cup \mathcal{E}(W', D)\big)\,. \end{aligned} \tag{31.22}$$

Yet $\mathcal{E}(W, D, T_1, 1) \cup \mathcal{E}(W', D)$ is also the output of a genuine exploration algorithm. So the probabilistic control (28.2) holds for the three terms of the right-hand side of (31.22), from which we deduce that

$$\begin{aligned} &P\big(|S\big(\mathcal{E}(W, D, T_1, 1) \cap \mathcal{E}(W', D)\big)| \geq i_2, W \not\longleftrightarrow W'\big)\\ &\qquad\qquad \leq\, 6|\mathbb{E}^d(D)| \exp\Big(-\frac{2p^2(1-p)^2 i_2^2}{9|\mathbb{E}^d(D)|}\Big)\,. \end{aligned} \tag{31.23}$$

In addition, all the bonds in $\mathcal{E}(W,D,T_1,1)\cap\mathcal{E}(W',D)$ must be closed, because the aggregates $\mathcal{A}(1)$ and $\mathcal{A}(2)$ are disjoint, therefore

$$S\big(\mathcal{E}(W,D,T_1,1)\cap\mathcal{E}(W',D)\big) \;=\; -\frac{1}{1-p}\big|\mathcal{E}(W,D,T_1,1)\cap\mathcal{E}(W',D)\big|\,. \quad (31.24)$$

Combining the inequalities (31.21), (31.23) and (31.24), we obtain that

$$P\big(|\mathcal{I}(2)|\geq i_2, |\mathcal{I}(1)|\leq i_1, \mathcal{T}(1)=T_1, W \not\longleftrightarrow W'\big) \leq 6|\mathbb{E}^d(D)|\exp\Big(-\frac{2p^2 i_2^2}{9|\mathbb{E}^d(D)|}\Big). \quad (31.25)$$

Substituting (31.25) into (31.20) and bounding the number of terms in the sum by $|D|^{i_1}$, we obtain

$$P\big(|\mathcal{I}(2)|\geq i_2, |\mathcal{I}(1)|\leq i_1, W \not\longleftrightarrow W'\big) \;\leq\; |D|^{i_1}6|\mathbb{E}^d(D)|\exp\Big(-\frac{2p^2 i_2^2}{9|\mathbb{E}^d(D)|}\Big)\,.$$

We report this estimate (31.19) and we use the control on $|\mathcal{I}(1)|$ obtained in (31.6) to get finally that, for any $i_1,i_2\geq 1$,

$$P\big(|\mathcal{I}(2)|\geq i_2, W \not\longleftrightarrow W'\big) \;\leq$$
$$|D|^{i_1}6|\mathbb{E}^d(D)|\exp\Big(-\frac{2p^2 i_2^2}{9|\mathbb{E}^d(D)|}\Big) + 6|\mathbb{E}^d(D)|\exp\Big(-\frac{2p^2 i_1^2}{9|\mathbb{E}^d(D)|}\Big)\,. \quad (31.26)$$

We apply this result with $D=\Lambda(n)$ for some $n\geq 2d$ and we use the inequalities (31.7). The inequality (31.26) becomes

$$P\big(|\mathcal{I}(2)|\geq i_2, W \not\longleftrightarrow W'\big) \;\leq$$
$$6(d2^d n^d)^{i_1+1}\exp\Big(-\frac{2p^2 i_2^2}{9d2^d n^d}\Big) + 6d2^d n^d\exp\Big(-\frac{2p^2 i_1^2}{9d2^d n^d}\Big)\,. \quad (31.27)$$

For $i_1, i_2$, we choose

$$i_1 \;=\; \frac{3}{p}d2^d\sqrt{n^d\ln n}\,,\quad i_2 \;=\; \frac{8}{p^2}d^2 2^{2d}n^{3d/4}\big(\ln n\big)^{3/4}\,.$$

These values ensure that $i_2\geq\sqrt{2}i_1^{3/2}$ and

$$\frac{2p^2 i_1^2}{9d2^d n^d} \;\geq\; 4d\ln n\,, \quad (31.28)$$

$$\frac{2p^2 i_2^2}{9d2^d n^d} \;\geq\; 2i_1\frac{2p^2 i_1^2}{9d2^d n^d} \;\geq\; (i_1+1)4d\ln n\,. \quad (31.29)$$

Substituting the inequalities (31.28) and (31.29) in (31.27), we obtain that

$$P\Big(|\mathcal{I}(2)|\geq\frac{8}{p^2}d^2 2^{2d}n^{3d/4}(\ln n)^{3/4}, W \not\longleftrightarrow W'\Big) \;\leq\; \frac{12}{n^{2d}}\,. \quad (31.30)$$

This is to be compared with the control on $|\mathcal{I}(1)|$ obtained in (31.9). The main lesson so far is that the bound on $|\mathcal{I}(1)|$ is of order $\sqrt{n^d\ln n}$, while the bound on $|\mathcal{I}(2)|$ is of order $(n^d\ln n)^{3/4}$. Unfortunately, there is a deterioration in the exponent.

## 31.3 Reconstruction through shells

In order to control the higher intersections, we try to use the same strategy as for the second intersections. To that end, we wish to express the intersections as the output of two exploration algorithms. Since the explorations involved in the construction of the intersections are intertwined, the dynamics of each explorer depends heavily on the other. Let us focus for instance on the first explorer. When the first explorer is activated and performs a round of exploration, it is given a taboo set which is the set of vertices already visited by the second explorer. In fact, the entire behavior of the first explorer is somehow determined by the sequence of these taboo sets. Using this idea, we shall disentangle the intertwined explorations and we shall reconstruct the territory associated to each explorer with the help of an adequate taboo exploration. The key result to realize this program is the next proposition. It describes the increments between the successive aggregates as shells around $W$ and $W'$ with an adequate taboo set. Before stating the results, let us recall the definition of the balls and the shells associated to the travel time in a domain $D$. For $A$ a subset of $D$ and for $t \geq 0$ an integer, we define the closed ball of radius $t$ for the travel time $T_D$ by

$$\mathcal{B}(A, D, t) \,=\, \big\{\, x \in D : T_D(A, x) \leq t \,\big\}\,, \tag{31.31}$$

and the shell number $t$ around $A$ in $D$ by

$$\mathrm{Shell}\,(A, D, t) \,=\, \big\{\, x \in D : T_D(A, x) = t \,\big\}\,.$$

**Proposition 31.1.** *Let $W$, $W'$ be two subsets of a finite domain $D$ and let us consider a percolation configuration in $D$ in which $W$ and $W'$ are not connected. Let $(\mathcal{A}(t), \mathcal{W}(t), t \geq 0)$ be the intertwined explorations starting from $W$ and $W'$. For any $t \geq 0$, we have (with the convention that $\mathcal{A}(-1) = \varnothing$)*

$$\mathcal{A}(2t) \,=\, \mathcal{B}\big(W, D \setminus \mathcal{A}(2t-1), t\big)\,, \tag{31.32}$$

$$\begin{aligned} \mathit{Clusters}\,\big(\mathcal{W}(2t), D \setminus \mathcal{A}(2t), \mathcal{A}(2t+1)\big) \,=&\, \\ &\mathit{Shell}\,\big(W, D \setminus \mathcal{A}(2t+1), t+1\big)\,, \end{aligned} \tag{31.33}$$

$$\mathcal{A}(2t+1) \,=\, \mathcal{B}\big(W', D \setminus \mathcal{A}(2t), t\big)\,,$$

$$\begin{aligned} \mathit{Clusters}\,\big(\mathcal{W}(2t+1), D \setminus \mathcal{A}(2t+1), \mathcal{A}(2t+2)\big) \,=&\, \\ &\mathit{Shell}\,\big(W', D \setminus \mathcal{A}(2t+2), t+1\big)\,. \end{aligned}$$

*Proof.* We do the proof for the even indices, the argument for the odd indices is essentially the same. Let us consider the case $t = 0$. By definition, we have

$$\mathcal{A}(0) \,=\, \mathrm{Clusters}\,\big(W, D, \varnothing\big) \,=\, \mathrm{Shell}\,(W, D) \,=\, \mathcal{B}(W, D, 0)\,.$$

Moreover, since $W$ and $W'$ are not connected, then $\mathcal{A}(1) = \mathrm{Shell}\,(W', D)$ is disjoint from $\mathcal{A}(0)$ and

$$\mathcal{B}(W, D, 0) \,=\, \mathcal{B}\big(W, D \setminus \mathcal{A}(1), 0\big)\,.$$

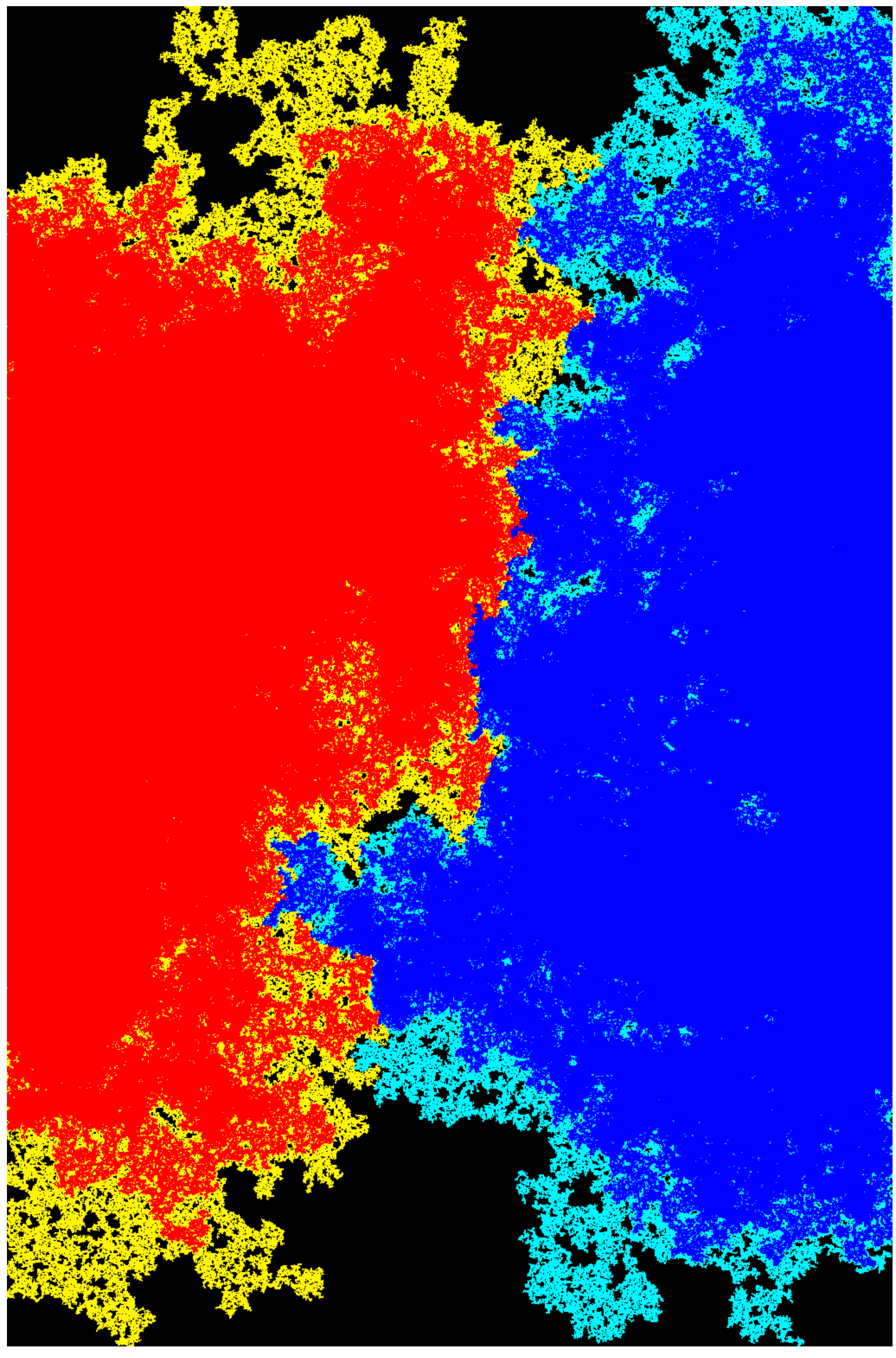

Figure 56: Same simulation as in figure 52, left right disconnection, $p = 0.49$, $\mathcal{A}(6) = \mathcal{B}\big(W, D \setminus \mathcal{A}(5), 3\big)$ , Clusters $\big(\mathcal{W}(6), D \setminus \mathcal{A}(6), \mathcal{A}(7)\big)$,
$\mathcal{A}(7) = \mathcal{B}\big(W', D \setminus \mathcal{A}(6), 3\big)$ , Clusters $\big(\mathcal{W}(6), D \setminus \mathcal{A}(6), \mathcal{A}(7)\big)$ .

This proves formula (31.32) for $t = 0$. We already proved in (31.15) the formula (31.33) for $t = 0$ when we studied the second intersections. We proceed next by induction over $t$. Let $t \geq 1$ and suppose that the formulas (31.32) and (31.33) have been proved at rank $t-1 \geq 0$. We first prove formula (31.32) at rank $t$. By the very definition of the intertwined explorations, we have

$$\mathcal{A}(2t) \;=\; \mathcal{A}(2t-2) \cup \text{Clusters}\big(\mathcal{W}(2t-2), D \setminus \mathcal{A}(2t-2), \mathcal{A}(2t-1)\big)\,. \quad (31.34)$$

Using the induction hypothesis and substituting formulas (31.32) and (31.33) at rank $t-1$ into (31.34), we obtain

$$\begin{aligned}\mathcal{A}(2t) \;&=\; \mathcal{B}\big(W, D \setminus \mathcal{A}(2t-3), t-1\big) \cup \text{Shell}\big(W, D \setminus \mathcal{A}(2t-1), t\big)\\ &=\; \mathcal{B}\big(W, D \setminus \mathcal{A}(2t-1), t-1\big) \cup \text{Shell}\big(W, D \setminus \mathcal{A}(2t-1), t\big)\\ &=\; \mathcal{B}\big(W, D \setminus \mathcal{A}(2t-1), t\big)\,.\end{aligned}$$

The second equality holds because $\mathcal{A}(2t)$ and $\mathcal{A}(2t-1)$ are disjoint. For $t = 1$, we use the convention that $\mathcal{A}(-1) = \varnothing$. This completes the proof of formula (31.32) at rank $t$. We attack now the more delicate proof of formula (31.33) at rank $t$. By the very definition of the set Clusters $(\cdot,\cdot,\cdot)$ (see formula (29.1)), we have

$$\begin{aligned}&\text{Clusters}\big(\mathcal{W}(2t), D \setminus \mathcal{A}(2t), \mathcal{A}(2t+1)\big)\\ &\qquad= C\Big(\mathcal{W}(2t) \setminus \mathcal{A}(2t+1), D \setminus \big(\mathcal{A}(2t) \cup \mathcal{A}(2t+1)\big)\Big)\\ &\qquad= \text{Shell}\Big(\mathcal{W}(2t) \setminus \mathcal{A}(2t+1), D \setminus \big(\mathcal{A}(2t) \cup \mathcal{A}(2t+1)\big)\Big)\\ &\qquad= \text{Shell}\Big(\mathcal{W}(2t), D \setminus \big(\mathcal{A}(2t) \cup \mathcal{A}(2t+1)\big)\Big)\,. \quad (31.35)\end{aligned}$$

The last equality holds thanks to the property (28.6). The set $\mathcal{W}(2t)$ is in turn defined as (see formulas (29.1) and (29.2))

$$\begin{aligned}\mathcal{W}(2t) \;&=\; \text{Waiting}\big(\mathcal{W}(2t-2), D \setminus \mathcal{A}(2t-2), \mathcal{A}(2t-1)\big) \quad (31.36)\\ &=\; \partial^{out}_{D\setminus(\mathcal{A}(2t-2)\cup\mathcal{A}(2t-1))}\Big(\text{Clusters}\big(\mathcal{W}(2t-2), D \setminus \mathcal{A}(2t-2), \mathcal{A}(2t-1)\big)\Big)\,.\end{aligned}$$

Using the induction hypothesis at rank $t-1$, we have

$$\text{Clusters}\big(\mathcal{W}(2t-2), D \setminus \mathcal{A}(2t-2), \mathcal{A}(2t-1)\big) \;=\; \text{Shell}\big(W, D \setminus \mathcal{A}(2t-1), t\big)\,. \quad (31.37)$$

Substituting successively (31.37) in (31.36), and (31.36) in (31.35), we have

$$\begin{aligned}&\text{Clusters}\big(\mathcal{W}(2t), D \setminus \mathcal{A}(2t), \mathcal{A}(2t+1)\big) \;= \quad (31.38)\\ &\text{Shell}\Big(\partial^{out}_{D\setminus(\mathcal{A}(2t-2)\cup\mathcal{A}(2t-1))}\text{Shell}\big(W, D\setminus\mathcal{A}(2t-1), t\big), D\setminus(\mathcal{A}(2t) \cup \mathcal{A}(2t+1))\Big).\end{aligned}$$

We apply next lemma 28.1 with $A = W$, $E = D \setminus \mathcal{A}(2t+1)$, $F = D \setminus \mathcal{A}(2t-1)$ and we check first the hypothesis of the lemma. We have $A \cap E = A \cap F = W$ because $W \subset \mathcal{A}(0)$ and neither $\mathcal{A}(2t-1)$ nor $\mathcal{A}(2t+1)$ intersects $W$ by coherence. Next,

$$\begin{aligned} T_F(A, F \setminus E) \,&=\, T_{D\setminus\mathcal{A}(2t-1)}\big(W, \mathcal{A}(2t+1) \setminus \mathcal{A}(2t-1)\big) \\ &\geq\, T_{D\setminus\mathcal{A}(2t-1)}\big(W, D \setminus \mathcal{A}(2t)\big)\,, \end{aligned} \tag{31.39}$$

where we have used that $\mathcal{A}(2t) \cap \mathcal{A}(2t+1) = \varnothing$ by coherence. We have already proved at the beginning of the induction step that

$$\mathcal{A}(2t) \,=\, \mathcal{B}\big(W, D \setminus \mathcal{A}(2t-1), t\big)\,, \tag{31.40}$$

and this implies that $T_{D\setminus\mathcal{A}(2t-1)}\big(W, D\setminus\mathcal{A}(2t)\big) > t$, which, together with (31.39), shows that the hypothesis $T_F(A, F \setminus E) > t$ is satisfied. Applying lemma 28.1, we obtain that

$$\mathcal{B}\big(W, D \setminus \mathcal{A}(2t-1), t\big) \,=\, \mathcal{B}\big(W, D \setminus \mathcal{A}(2t+1), t\big)\,.$$

It is also a consequence of the equalities (28.11) proved in lemma 28.1 that

$$\text{Shell}\,\big(W, D \setminus \mathcal{A}(2t-1), t\big) \,=\, \text{Shell}\,\big(W, D \setminus \mathcal{A}(2t+1), t\big)\,. \tag{31.41}$$

To alleviate the formulas, we introduce the notation

$$\mathcal{S} \,=\, \text{Shell}\,\big(W, D \setminus \mathcal{A}(2t+1), t\big)\,. \tag{31.42}$$

It follows from formula (31.37), (31.41) and (31.42) that we have the inclusion

$$\mathcal{S} \,\subset\, \mathcal{A}(2t) \setminus \big(\mathcal{A}(2t-2) \cup \mathcal{A}(2t+1)\big)\,. \tag{31.43}$$

With the help of (31.41) and (31.42), we can rewrite (31.38) as

$$\begin{aligned} &\text{Clusters}\,\big(\mathcal{W}(2t), D \setminus \mathcal{A}(2t), \mathcal{A}(2t+1)\big) \\ &\qquad =\, \text{Shell}\,\Big(\partial^{\,out}_{D\setminus\big(\mathcal{A}(2t-2)\cup\mathcal{A}(2t-1)\big)}\mathcal{S}, D \setminus \big(\mathcal{A}(2t) \cup \mathcal{A}(2t+1)\big)\Big)\,. \end{aligned} \tag{31.44}$$

We claim that

$$\begin{aligned} &\Big(\partial^{\,out}_{D\setminus\big(\mathcal{A}(2t-2)\cup\mathcal{A}(2t-1)\big)}\mathcal{S}\Big) \setminus \big(\mathcal{A}(2t) \cup \mathcal{A}(2t+1)\big) \,= \\ &\qquad\qquad \Big(\partial^{\,out}_{D\setminus\mathcal{A}(2t+1)}\mathcal{S}\Big) \setminus \big(\mathcal{A}(2t) \cup \mathcal{A}(2t+1)\big)\,. \end{aligned} \tag{31.45}$$

This equality holds because, thanks to the inclusion (31.43), we have

$$\mathcal{S} \setminus \big(\mathcal{A}(2t-2) \cup \mathcal{A}(2t-1)\big) \,=\, \mathcal{S} \setminus \mathcal{A}(2t+1) \,=\, \mathcal{S}\,. \tag{31.46}$$

Let us provide the full details. We prove (31.45) by double inclusion. Let $x$ belong to

$$\Big(\partial^{\,out}_{D\setminus\big(\mathcal{A}(2t-2)\cup\mathcal{A}(2t-1)\big)}\mathcal{S}\Big) \setminus \big(\mathcal{A}(2t) \cup \mathcal{A}(2t+1)\big)\,. \tag{31.47}$$

By definition, the site $x$ is in $D \setminus \big(\mathcal{S} \cup \mathcal{A}(2t) \cup \mathcal{A}(2t+1)\big)$ and there exists $y$ in $\mathcal{S} \setminus \big(\mathcal{A}(2t-2) \cup \mathcal{A}(2t-1)\big)$ such that $|x-y| = 1$. By (31.46), the site $y$ is also in $D \setminus \mathcal{A}(2t+1)$ and we see that $x$ belongs also to the set

$$\Big(\partial^{out}_{D\setminus\mathcal{A}(2t+1)}\mathcal{S}\Big) \setminus \big(\mathcal{A}(2t) \cup \mathcal{A}(2t+1)\big)\,. \tag{31.48}$$

Conversely, let $x$ belong to the set (31.48). By definition, the site $x$ is in the set $D \setminus \big(\mathcal{S} \cup \mathcal{A}(2t) \cup \mathcal{A}(2t+1)\big)$ and there exists $y$ in $\mathcal{S} \setminus \mathcal{A}(2t+1)$ such that $|x-y| = 1$. By (31.46), the site $y$ is also in $\mathcal{S} \setminus \big(\mathcal{A}(2t-2) \cup \mathcal{A}(2t-1)\big)$ and we see that $x$ belongs also to the set (31.47).

Now that (31.45) is proved, we are close to the conclusion of the proof. Using (31.45), we can rewrite (31.44) as

$$\begin{aligned}&\text{Clusters}\big(\mathcal{W}(2t), D \setminus \mathcal{A}(2t), \mathcal{A}(2t+1)\big)\\ &\qquad\qquad = \text{Shell}\Big(\partial^{out}_{D\setminus\mathcal{A}(2t+1)}\mathcal{S}, D \setminus (\mathcal{A}(2t) \cup \mathcal{A}(2t+1))\Big)\,.\end{aligned} \tag{31.49}$$

We claim that

$$\begin{aligned}&\text{Shell}\Big(\partial^{out}_{D\setminus\mathcal{A}(2t+1)}\mathcal{S}, D \setminus (\mathcal{A}(2t) \cup \mathcal{A}(2t+1))\Big)\\ &\qquad\qquad = \text{Shell}\Big(\partial^{out}_{D\setminus\mathcal{A}(2t+1)}\mathcal{S}, D \setminus \mathcal{A}(2t+1)\Big) \setminus \mathcal{A}(2t)\,.\end{aligned} \tag{31.50}$$

The set of the left-hand side of (31.50) is obviously included in the set of the right-hand side. Let next $x$ be a vertex in the set of the right-hand side of (31.50). Then $x$ is not in $\mathcal{A}(2t)$ and it is connected by an open path in $D \setminus \mathcal{A}(2t+1)$ to a site $y$ in $\partial^{out}_{D\setminus\mathcal{A}(2t+1)}\mathcal{S}$. Yet $\mathcal{A}(2t)$ is an aggregate, therefore the whole path must be included in $D \setminus \mathcal{A}(2t)$ and we see that $x$ belongs to the set of the left-hand side of (31.50). Substituting (31.50) into (31.49), we get

$$\begin{aligned}&\text{Clusters}\big(\mathcal{W}(2t), D \setminus \mathcal{A}(2t), \mathcal{A}(2t+1)\big)\\ &\qquad\qquad = \text{Shell}\Big(\partial^{out}_{D\setminus\mathcal{A}(2t+1)}\mathcal{S}, D \setminus \mathcal{A}(2t+1)\Big) \setminus \mathcal{A}(2t)\,.\end{aligned} \tag{31.51}$$

Furthermore, using (31.40) and the fact that $\mathcal{A}(2t)$ and $\mathcal{A}(2t+1)$ are disjoint, we have

$$\mathcal{A}(2t) \,=\, \mathcal{B}\big(W, D \setminus \mathcal{A}(2t+1), t\big)\,. \tag{31.52}$$

From (31.51) and (31.52), we conclude that

$$\begin{aligned}&\text{Clusters}\big(\mathcal{W}(2t), D \setminus \mathcal{A}(2t), \mathcal{A}(2t+1)\big)\\ &\quad = \text{Shell}\Big(\partial^{out}_{D\setminus\mathcal{A}(2t+1)}\mathcal{S}, D \setminus \mathcal{A}(2t+1)\Big) \setminus \mathcal{B}\big(W, D \setminus \mathcal{A}(2t+1), t\big)\\ &\qquad\qquad\qquad = \text{Shell}\big(W, D \setminus \mathcal{A}(2t+1), t+1\big)\,.\end{aligned}$$

The last equality is a consequence of the identity (28.10) applied to $W$ and $D \setminus \mathcal{A}(2t+1)$. This completes the proof of formula (31.33) at rank $t$. □

As painful as it was, we had to provide a precise and detailed proof, because this result is crucial. Moreover, at some point in the future, we will need an analogous result for the site model, and with $N$ explorers instead of just 2!

## 31.4 Restricting the taboo set

In proposition 31.1, we provided a way to disentangle the intertwined explorations. The aggregate of the first explorer at time $2t$ is completely defined as the output of a taboo exploration algorithm, the taboo set being the aggregate of the second explorer at time $2t-1$. When we perform a taboo exploration with a deterministic taboo set, we have again a genuine exploration algorithm, and we are free to use the probabilistic control (28.2). Unfortunately, this is not any more true when the taboo set is random, because it can incorporate information on the percolation configuration. This is notably the case if we take as taboo set the aggregate $\mathcal{A}(2t-1)$. For instance, if we tell the first explorer the current value of $\mathcal{A}(2t-1)$, then he will know that this set is an aggregate, and that each bond leading from the outer boundary $\partial^{out}\mathcal{A}(2t-1)$ to a vertex of $\mathcal{A}(2t-1)$ has to be closed. As a consequence, the distribution of the bonds discovered by the first explorer is not any more an i.i.d. Bernoulli field. Unfortunately this is in the nature of things. A possible strategy to obtain a control on the aggregate of the first explorer is to condition on the aggregate of the second explorer, and then try to say something about the ensuing conditional distribution of the bonds explored by the first explorer. Of course this conditioning is very nasty to handle. Our first try will consist in making a simple union bound to get rid of the conditioning. Now, if we implement this strategy and we condition on $\mathcal{A}(2t-1)$, nothing interesting comes out, because the combinatorial factor induced by the union bound is far too large. Therefore, it is desirable to condition on a smaller set than the full aggregate of the second explorer. By the way, only a subset of the taboo set $\mathcal{A}(2t-1)$ matters for the first explorer, namely the vertices which are an endpoint of a bond belonging to an intersection of order at most $2t-1$. Indeed these vertices are precisely the vertices at which the first explorer would enter the territory of the second explorer, if he were to explore the configuration without the taboo constraint. Let us be a bit more precise. When the first explorer starts its $t$-th round of exploration, he will explore the clusters of the vertices of $\mathcal{W}(2t-2)$, and he will use the set $\mathcal{A}(2t-1)$ as a taboo set. Since $\mathcal{A}(2t-1)$ is an aggregate, the only vertices of $\mathcal{W}(2t-2)$ which can be connected by an open path in $D$ to a vertex of $\mathcal{A}(2t-1)$ are those vertices of $\mathcal{W}(2t-2)$ which belong to $\mathcal{A}(2t-1)$. Therefore the vertices of the taboo set which really matter during this round of exploration are the vertices belonging to $\mathcal{A}(2t-1)$ which are the endpoint of an edge belonging to

$$\Delta_D\mathcal{A}(2t-2)\cap\Delta_D\mathcal{A}(2t-1) \;=\; \bigcup_{1\leq s\leq 2t-1}\mathcal{I}(s)\,.$$

We formalize this in the next proposition. To this end, we introduce an adequate restricted taboo set $\mathcal{T}(2t)$ and we show that the aggregate $\mathcal{A}(2t)$ of the first explorer is indeed the result of a taboo exploration using $\mathcal{T}(2t)$ for taboo set. Naturally, we state a symmetric result for the second explorer.

**Proposition 31.2.** *For $t \geq 1$, we define*

$$\mathcal{T}(2t-1) = \Big\{ x \in \mathcal{A}(2t-1) : \exists\, y \in \mathcal{A}(2t-2) \quad \langle x,y \rangle \in \bigcup_{1\leq s\leq 2t-1} \mathcal{I}(s) \Big\} \tag{31.53}$$

$$\mathcal{T}(2t) = \Big\{ x \in \mathcal{A}(2t) : \exists\, y \in \mathcal{A}(2t-1) \quad e = \langle x,y \rangle \in \bigcup_{1\leq s\leq 2t} \mathcal{I}(s) \Big\} \tag{31.54}$$

*We have then*

$$\begin{aligned} \mathcal{A}(2t) &= \mathcal{B}\big(W, D \setminus \mathcal{T}(2t-1), t\big)\,, \\ \mathcal{A}(2t+1) &= \mathcal{B}\big(W, D \setminus \mathcal{T}(2t), t\big)\,. \end{aligned} \tag{31.55}$$

*Proof.* We do the proof for the even indices, the argument for the odd indices is essentially the same. The case $t = 1$ has been treated when we studied the second intersections in subsection 31.2. For $t = 1$, we have

$$\begin{aligned} \mathcal{T}(1) &= \big\{ x \in \mathcal{A}(1) : \exists y \in D \setminus \mathcal{A}(1) \quad \langle x,y \rangle \in \mathcal{I}(1) \big\} \\ &= \big\{ x \in \mathcal{A}(1) : \exists y \in \mathcal{A}(0) \quad \langle x,y \rangle \in \mathcal{I}(1) \big\}\,, \end{aligned}$$

and we have shown in (31.16) that

$$\mathcal{A}(2) = \text{Shell}\,(W, D, 0) \cup \text{Shell}\,(W, D \setminus \mathcal{A}(1), 1)\,. \tag{31.56}$$

Next, we proved in (31.18) that

$$\text{Shell}\,(W, D \setminus \mathcal{A}(1), 1) = \text{Shell}\,(W, D \setminus \mathcal{T}(1), 1)\,. \tag{31.57}$$

From (31.56) and (31.57), we conclude that

$$\mathcal{A}(2) = \mathcal{B}\big(W, D \setminus \mathcal{T}(1), 1\big)\,.$$

Let $t \geq 1$ and suppose that formula (31.55) has been proved at rank $t$. From proposition 31.1, we know that

$$\mathcal{A}(2t+2) = \mathcal{B}\big(W, D \setminus \mathcal{A}(2t+1), t+1\big)\,.$$

Since the restricted taboo set $\mathcal{T}(2t+1)$ is a subset of $\mathcal{A}(2t+1)$, then we have

$$\mathcal{B}\big(W, D \setminus \mathcal{A}(2t+1), t+1\big) \subset \mathcal{B}\big(W, D \setminus \mathcal{T}(2t+1), t+1\big)\,.$$

To complete the induction step, we have to prove the converse inclusion. Let $x$ be a vertex in $\mathcal{B}\big(W, D \setminus \mathcal{T}(2t+1), t+1\big)$. We decompose this last set as

$$\begin{aligned} \mathcal{B}\big(W, D \setminus \mathcal{T}(2t+1), t+1\big) &= \mathcal{B}\big(W, D \setminus \mathcal{T}(2t+1), t\big) \\ &\quad \cup \text{Shell}\,\Big(\mathcal{B}\big(W, D \setminus \mathcal{T}(2t+1), t\big), D \setminus \mathcal{T}(2t+1), 1\Big)\,, \end{aligned} \tag{31.58}$$

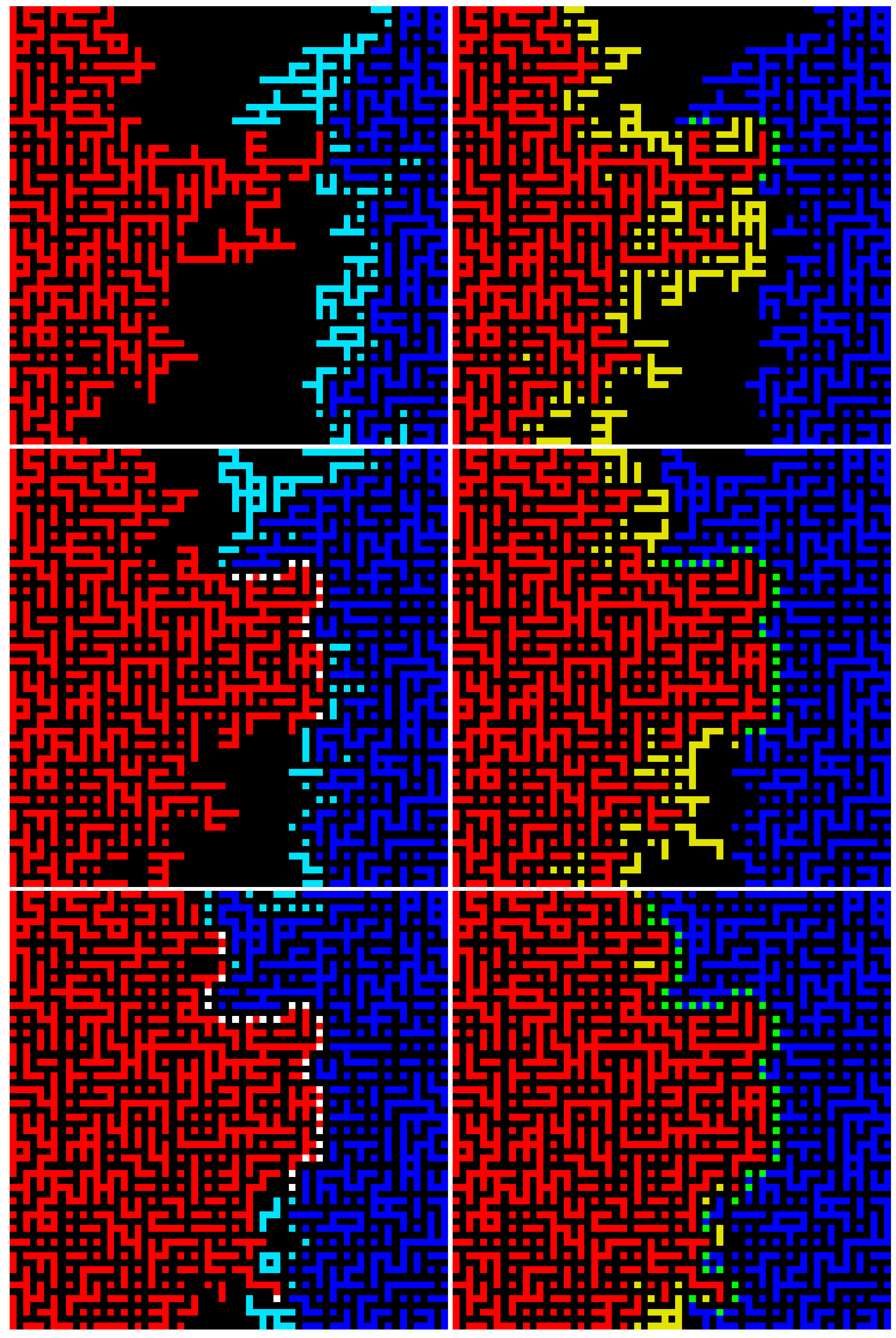

Figure 57: Taboo restrictions $W$ and $W' =$ left and right sides, $\Lambda(32)$, $p = 0.35$. $\mathcal{A}(2), \mathcal{A}(4), \mathcal{A}(6)$, increments of left explorer, $\mathcal{T}(1) = \varnothing, \mathcal{T}(3), \mathcal{T}(5), \mathcal{T}(7)$, $\mathcal{A}(3), \mathcal{A}(5), \mathcal{A}(7)$, increments of right explorer, $\mathcal{T}(2), \mathcal{T}(4)$ in white, 6 rounds.

and we study separately the two sets appearing in the union. Since $\mathcal{T}(2t-1) \subset \mathcal{T}(2t+1)$, then

$$\mathcal{B}\big(W, D \setminus \mathcal{T}(2t+1), t\big) \subset \mathcal{B}\big(W, D \setminus \mathcal{T}(2t-1), t\big) = \mathcal{A}(2t)\,, \tag{31.59}$$

where the last equality comes from the induction hypothesis. Using now (31.59), we have

$$\text{Shell}\left(\mathcal{B}\big(W, D\setminus\mathcal{T}(2t+1), t\big), D\setminus\mathcal{T}(2t+1), 1\right) \subset \text{Shell}\left(\mathcal{A}(2t), D\setminus\mathcal{T}(2t+1), 1\right). \tag{31.60}$$

Let $x$ be a vertex in Shell $\big(\mathcal{A}(2t), D \setminus \mathcal{T}(2t+1), 1\big)$. By definition, we have

$$T_{D\setminus\mathcal{T}(2t+1)}(\mathcal{A}(2t), x) = 1\,, \tag{31.61}$$

hence there exists a path

$$z_0, e_0, z_1, \dots, e_{r-1}, z_r$$

in $D \setminus \mathcal{T}(2t+1)$ joining a vertex $z_0 = y$ of $\mathcal{A}(2t)$ to $z_r = x$ such that, among the bonds $e_0, \dots, e_{r-1}$, exactly one is closed and the other are opened. Let $i$ be the index of the closed bond. The bonds $e_0, \dots, e_{i-1}$ are opened and $z_0 \in \mathcal{A}(2t)$. Since $\mathcal{A}(2t)$ is an aggregate in $D$, then all the vertices $z_0, \dots, z_i$ must belong to $\mathcal{A}(2t)$ as well. Let us now take a closer look at the bond $e_i = \langle z_i, z_{i+1}\rangle$. The vertex $z_i$ is in $\mathcal{A}(2t)$, the bond $e_i$ is closed, and $z_{i+1}$ is not in $\mathcal{A}(2t)$ (otherwise the vertices $z_{i+1}, \dots, z_r = x$ would also belong to $\mathcal{A}(2t)$, and this would contradict (31.61)). Therefore the bond $e_i$ is in the edge boundary $\Delta_D\mathcal{A}(2t)$ of $\mathcal{A}(2t)$. However, the vertex $z_{i+1}$ does not belong to $\mathcal{T}(2t+1)$. This prevents $e_i$ from being in the union of the intersections

$$\bigcup_{1\le s\le 2t+1} \mathcal{I}(s)\,.$$

However the above set is equal to $\Delta_D\mathcal{A}(2t)\cap\Delta_D\mathcal{A}(2t+1)$. In particular, $e_i$ does not belong to this last set and therefore $z_{i+1}$ is not in $\mathcal{A}(2t+1)$. Since $\mathcal{A}(2t+1)$ is an aggregate in $D$, we conclude that none of the vertices $z_{i+1}, \dots, z_r = x$ is in $\mathcal{A}(2t+1)$, so that $x$ is in Shell $\big(\mathcal{A}(2t), D \setminus \mathcal{A}(2t+1), 1\big)$ and we have proved that

$$\text{Shell}\left(\mathcal{A}(2t), D \setminus \mathcal{T}(2t+1), 1\right) \subset \text{Shell}\left(\mathcal{A}(2t), D \setminus \mathcal{A}(2t+1), 1\right). \tag{31.62}$$

Substituting the set inclusions (31.59), (31.60) and (31.62) into (31.58), we obtain

$$\mathcal{B}\big(W, D\setminus\mathcal{T}(2t+1), t+1\big) \subset \mathcal{A}(2t) \cup \text{Shell}\left(\mathcal{A}(2t), D\setminus\mathcal{A}(2t+1), 1\right). \tag{31.63}$$

Using formula (31.32), we have also $\mathcal{A}(2t) = \mathcal{B}\big(W, D \setminus \mathcal{A}(2t+1), t\big)$, which together with (31.63), yields that

$$\mathcal{B}\big(W, D \setminus \mathcal{T}(2t+1), t+1\big) \subset \mathcal{B}\big(W, D \setminus \mathcal{A}(2t+1), t+1\big)\,,$$

and the proof of the induction step is complete. □

## 31.5 The taboo growth algorithm

In the previous subsection, we have seen that the sequence of the aggregates $(\mathcal{A}(t), t \geq 0)$ can be built as closed balls around $W, W'$ for the travel time $T_D$, but with adequate random taboo sets. Our goal here is to realize the aggregates $(\mathcal{A}(t), t \geq 0)$ as the output of an exploration algorithm, and ideally a genuine one. To that end, we introduce here a taboo growth algorithm. This taboo growth algorithm is simply built upon the taboo exploration algorithm presented in subsection 29.2. It consists in iterating the taboo exploration algorithm to explore the successive shells around the starting set. Let us introduce a useful notation before describing the algorithm. For $C$ a subset of $D$, we define $\mathbb{E}^d(C, D)$ as the set of the bonds in $D$ having one or two endpoints in $C$, i.e.,

$$\begin{aligned} \mathbb{E}^d(C, D) \;&=\; \bigcup_{x \in C} \big\{ e = \langle x, y \rangle \in \mathbb{E}^d(D) : y \in \mathcal{N}(x, D) \,\big\} \\ &=\; \mathbb{E}^d\big(C \cup \partial_D^{out} C\big)\,. \end{aligned} \tag{31.64}$$

The pseudocode of the taboo growth algorithm is given in Algorithm 31.1.

**Algorithm 31.1** Taboo_Bond_Growth$\,(A, D, \mathcal{T}, t)$

**Require:** a starting set $A$, a domain $D$, a taboo set $\mathcal{T}$, a radius $t$
**Ensure:** $C \,=\, \mathcal{B}(A, D \setminus \mathcal{T}, t)\,, E \,=\, \mathbb{E}^d\big(C \cup \partial_D^{out} C\big)$

  $C(0) \leftarrow \varnothing$, $E(0) \leftarrow \varnothing$, $A(0) \leftarrow A$, $s \leftarrow 0$
    ▷ $s$ is the number of the shell around $A$ which is currently explored
  **repeat**
    $C, W, E = $ Taboo_Bond_Explore_Clusters$\,(A(s), D, \mathcal{T})$
    $C(s+1) \leftarrow C(s) \cup C$
    $E(s+1) \leftarrow E(s) \cup E$
    $A(s+1) \leftarrow W$
    $s \leftarrow s + 1$
  **until** $A(s) = \varnothing$ or $s > t$
  **return** $C(s), E(s)$

This algorithm updates iteratively four sets, denoted by $C$, $E$, $A$ and $W$ and it computes two sequences

$$\big(C(s), E(s), 0 \leq s \leq t\big)$$

of subsets of vertices and bonds of $D$, whose contents are the following:

- $E(s)$: the bonds which have been explored during the first $s$ iterations.
- $C(s)$: the vertices which have been explored during the first $s$ iterations.

Initially, we set $C(0) = \varnothing$, $E(0) = \varnothing$. For $s$ in $\{\,0, \dots, t\,\}$, the set $C(s)$ is the set of the vertices which are at travel distance less than or equal to $s$ from $A$ inside $D \setminus \mathcal{T}$, i.e.,

$$C(s) \;=\; \mathcal{B}(A, D \setminus \mathcal{T}, s)\,,$$

while the set $E(s)$ is the set of the bonds in $D$ having one or two endpoints in $C(s)$, i.e.,

$$E(s) \;=\; \mathbb{E}^d(C(s), D)\,.$$

The sets $A, W$ are two sets maintained internally for the algorithm. The set $A$ is the set of the active vertices for the current iteration, the set $W$ is the set of the waiting vertices, which will become active at the next iteration. We stress that the taboo growth algorithm in $D$ with taboo set $\mathcal{T}$ is different from a simple growth algorithm in $D \setminus \mathcal{T}$. Indeed, the taboo growth algorithm can explore bonds having one endpoint in $\mathcal{T}$ (as the taboo exploration algorithm does), whereas such bonds would not be taken into account for an exploration algorithm in $D \setminus \mathcal{T}$.

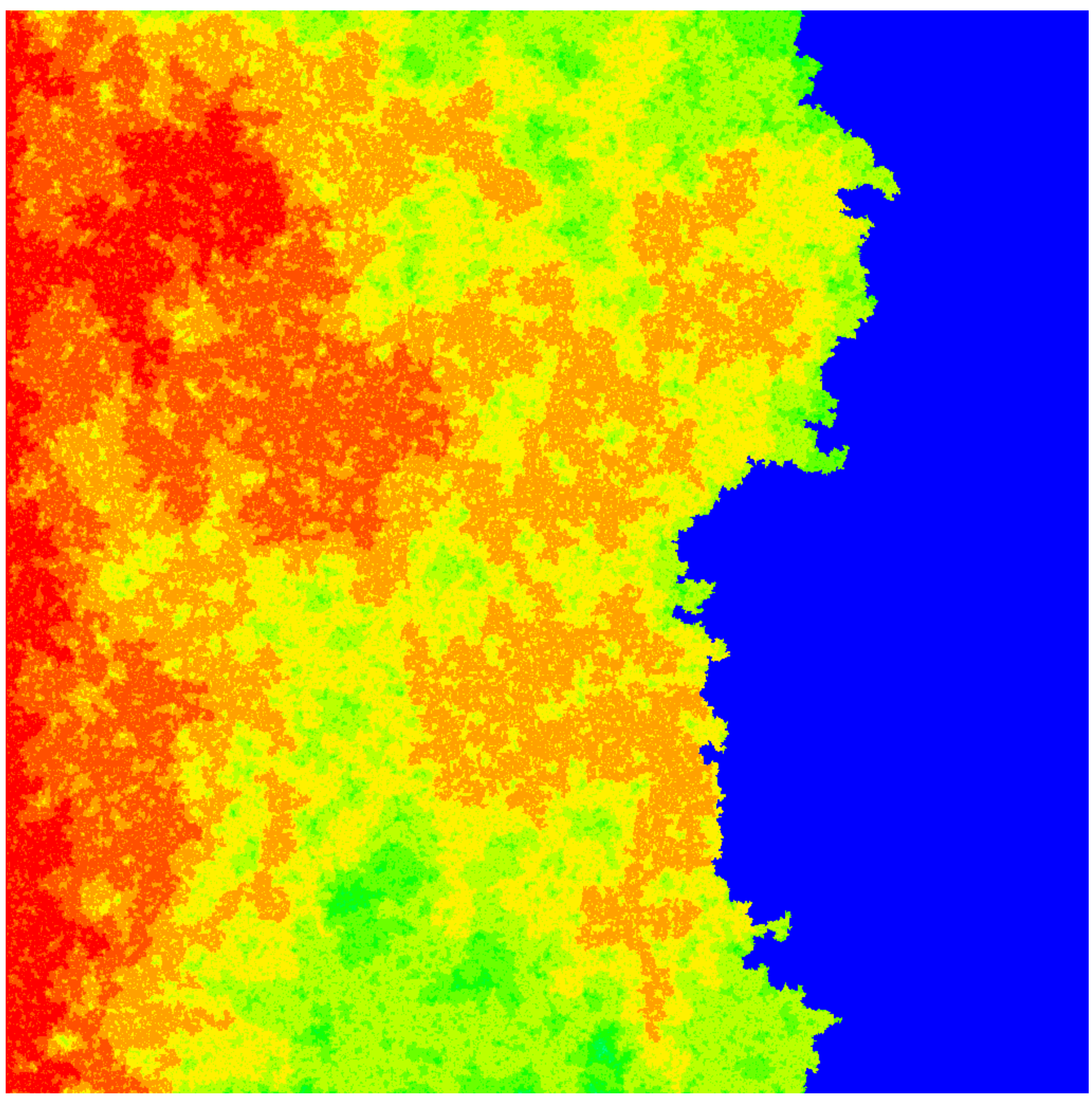

Figure 58: Taboo bond growth from the left side of the box in $\Lambda(1024)$, $p = 0.49$, $A$ is the left side, $\mathcal{T}$ is the blue zone, the result would be the same if $\mathcal{T}$ was the left boundary of the blue zone. The shells are in rainbow colors from red.

## 31.6 Brutal control of the higher intersections

We would like to propagate the estimates (31.13) and (31.26) on the intersections and to get some control on the successive sets $\mathcal{I}(t)$, $t \geq 2$. Let $t \geq 1$ and let $i_1, \dots, i_t$ be $t$ integers. We will try to control the probability

$$P\big(\exists s \in \{\,1, \dots, t\,\} \quad |\mathcal{I}(s)| > i_s,\, W \not\longleftrightarrow W'\big)\,.$$

In the first step, we decompose it as follows:

$$\begin{multline}
P\big(\exists s \in \{\,1, \dots, t\,\} \quad |\mathcal{I}(s)| > i_s,\, W \not\longleftrightarrow W'\big) \;=\; P\big(|\mathcal{I}(1)| > i_1, W \not\longleftrightarrow W'\big) \\
+ \sum_{1 \leq s \leq t} P\big(|\mathcal{I}(s)| > i_s, |\mathcal{I}(s-1)| \leq i_{s-1}, \cdots, |\mathcal{I}(1)| \leq i_1, W \not\longleftrightarrow W'\big)\,.
\end{multline} \tag{31.65}$$

The control on $|\mathcal{I}(1)|$ has already been obtained in (31.6). The next step is to get a control on the terms inside the sum. This is the purpose of the following proposition.

**Proposition 31.3.** *Let $t \geq 1$ and let $i_1, \dots, i_t$ be $t$ integers. We have*

$$\begin{multline}
P\big(|\mathcal{I}(t)| > i_t, |\mathcal{I}(t-1)| \leq i_{t-1}, \cdots, |\mathcal{I}(1)| \leq i_1, W \not\longleftrightarrow W'\big) \;\leq \\
6|\mathbb{E}^d(D)|^{i_1 + \cdots + i_{t-1}} \exp\Big(-\frac{2p^2 i_t^2}{9|\mathbb{E}^d(D)|}\Big)\,.
\end{multline} \tag{31.66}$$

*Proof.* We do the proof for the even indices, the argument for the odd indices is essentially the same. Let $t \geq 1$ and let $i_1, \dots, i_{2t}$ be $2t$ integers. We wish to control the probability of the event

$$\mathcal{E} \;=\; \left\{ \begin{matrix} |\mathcal{I}(2t)| > i_{2t}, |\mathcal{I}(2t-1)| \leq i_{2t-1}, \\ |\mathcal{I}(2t-2)| \leq i_{2t-2}, \cdots, |\mathcal{I}(2)| \leq i_2, \\ |\mathcal{I}(1)| \leq i_1, W \not\longleftrightarrow W' \end{matrix} \right\}\,.$$

Recalling the definition (29.15) of $\mathcal{I}(2t)$, we see that

$$\mathcal{I}(2t) \;\subset\; \Delta_D \mathcal{A}(2t-1) \cap \Delta_D \mathcal{A}(2t)\,,$$

and therefore

$$|\mathcal{I}(2t)| > i_{2t} \qquad \Longrightarrow \qquad \Big|\Delta_D \mathcal{A}(2t-1) \cap \Delta_D \mathcal{A}(2t)\Big| \;>\; i_{2t}\,. \tag{31.67}$$

The two aggregates $\mathcal{A}(2t-1)$ and $\mathcal{A}(2t)$ are disjoint, therefore

$$\Delta_D \mathcal{A}(2t-1) \cap \Delta_D \mathcal{A}(2t) \;=\; \mathbb{E}^d(\mathcal{A}(2t-1), D) \cap \mathbb{E}^d(\mathcal{A}(2t), D)$$

(the notation $\mathbb{E}^d(\cdot, D)$ was introduced in (31.64)). Moreover, all the bonds belonging to the above intersection are closed, so that

$$S\Big(\mathbb{E}^d(\mathcal{A}(2t-1), D) \cap \mathbb{E}^d(\mathcal{A}(2t), D)\Big) \;=\; -\frac{1}{1-p}\Big|\Delta_D \mathcal{A}(2t-1) \cap \Delta_D \mathcal{A}(2t)\Big|\,. \tag{31.68}$$

Using the additivity of the map $S$, we conclude from (31.67) and (31.68) that, if $|\mathcal{I}(2t)| > i_{2t}$, then at least one among the three quantities

$$S\big(\mathbb{E}^d(\mathcal{A}(2t-1),D)\big),\quad S\big(\mathbb{E}^d(\mathcal{A}(2t),D)\big),\quad S\big(\mathbb{E}^d(\mathcal{A}(2t-1)\cup\mathcal{A}(2t),D)\big)$$

must have an absolute value larger than $i_{2t}/(3(1-p))$. Therefore we can bound $P(\mathcal{E})$ as follows:

$$\begin{aligned} P(\mathcal{E}) \,\leq\, P\Big(\mathcal{E}, \big|S\big(\mathbb{E}^d(\mathcal{A}(2t-1),D)\big)\big| > \frac{i_{2t}}{3(1-p)}\Big) \\ + P\Big(\mathcal{E}, \big|S\big(\mathbb{E}^d(\mathcal{A}(2t),D)\big)\big| > \frac{i_{2t}}{3(1-p)}\Big) \\ + P\Big(\mathcal{E}, \big|S\big(\mathbb{E}^d(\mathcal{A}(2t-1)\cup\mathcal{A}(2t),D)\big)\big| > \frac{i_{2t}}{3(1-p)}\Big)\,. \end{aligned} \tag{31.69}$$

We try next to control the three quantities of the right-hand side of (31.69). We have shown in proposition 31.2 that $\mathcal{A}(2t)=\mathcal{B}\big(W,D\setminus\mathcal{T}(2t-1),t\big)$, where the taboo set $\mathcal{T}(2t-1)$ is defined in (31.53). In particular, the aggregate $\mathcal{A}(2t)$ is the output of the taboo growth algorithm presented in subsection 31.5 and we have

$$\Big(\mathcal{A}(2t),\mathbb{E}^d\big(\mathcal{A}(2t),D\big)\Big) \,=\, \text{Taboo_Bond_Growth}\,\big(W,D,\mathcal{T}(2t-1),t\big)\,. \tag{31.70}$$

Similarly, we have

$$\Big(\mathcal{A}(2t-1),\mathbb{E}^d\big(\mathcal{A}(2t-1),D\big)\Big) \,=\, \text{Taboo_Bond_Growth}\,\big(W',D,\mathcal{T}(2t-2),t\big)\,.$$

Finally, we have, thanks to formula (29.18),

$$\mathcal{A}(2t-1)\cup\mathcal{A}(2t) \,=\, \mathcal{B}\big(W',D,t-1\big)\cup\mathcal{B}\big(W,D,t\big)\,,$$

whence

$$\begin{gathered}\Big(\mathcal{A}(2t-1)\cup\mathcal{A}(2t),\mathbb{E}^d\big(\mathcal{A}(2t-1)\cup\mathcal{A}(2t),D\big)\Big) \,=\\ \text{Taboo_Bond_Growth}\,(W',D,\varnothing,t-1)\cup\text{Taboo_Bond_Growth}\,(W,D,\varnothing,t)\end{gathered}$$

(the union of the taboo growth algorithms means that we form the union of their outputs, for the aggregates on the one hand and for the set of explored bonds on the other hand). By the way, the union of the outputs of two genuine exploration algorithms is also the output of a genuine exploration algorithm (to see it, we can for instance run the two algorithms successively in such a way that the second algorithm has full access to the data collected by the first algorithm). So we have the usual control on $S\big(\mathbb{E}^d(\mathcal{A}(2t-1)\cup\mathcal{A}(2t),D)\big)$. The cases of the taboo growth algorithms working with the taboo sets $\mathcal{T}(2t-1)$ and $\mathcal{T}(2t)$ are more delicate. As these sets are random, the resulting algorithm is not any more genuine. Far worse, these algorithms incorporate from the start some important information on the percolation configuration. To deal with this issue, we shall try to condition on these taboo sets. Indeed, the taboo growth algorithm working with a fixed deterministic set is a genuine exploration algorithm. Looking at the definitions (31.53), (31.54) of $\mathcal{T}(2t-1)$ and $\mathcal{T}(2t)$, we see that, on the event $\mathcal{E}$, we have

$$\big|\mathcal{T}(2t-1)\big| \leq i_1+\cdots+i_{2t-1}\,, \qquad \big|\mathcal{T}(2t)\big| \leq i_1+\cdots+i_{2t}\,.$$

We decompose the second term of the right-hand side of (31.69) as follows:

$$P\Big(\mathcal{E}, \big|S\big(\mathbb{E}^d(\mathcal{A}(2t),D)\big)\big| > \frac{i_{2t}}{3(1-p)}\Big) = \sum_{T\subset D, |T|\leq i_1+\cdots+i_{2t-1}} P\Big(\mathcal{E}, \mathcal{T}(2t-1)=T, \big|S\big(\mathbb{E}^d(\mathcal{A}(2t),D)\big)\big| > \frac{i_{2t}}{3(1-p)}\Big)\,. \tag{31.71}$$

Next, we take care of the term inside the sum. Let us fix a set $T$. The map Taboo_Bond_Growth $(\cdot,\cdot,\cdot,\cdot)$ returns a pair $(C,E)$ where $C$ is a set of vertices and $E$ a set of bonds. We denote by $\pi_2$ the projection on the second element of such a pair, so that $\pi_2(C,E)=E$. Formula (31.70) yields that

$$\mathbb{E}^d\big(\mathcal{A}(2t),D\big) = \pi_2\big(\text{Taboo_Bond_Growth}\,(W,D,\mathcal{T}(2t-1),t)\big)\,.$$

With this notation, the probability inside the sum of (31.71) can be rewritten and bounded as follows:

$$\begin{aligned}
&P\Big(\mathcal{E}, \mathcal{T}(2t-1)=T, \big|S\big(\mathbb{E}^d(\mathcal{A}(2t),D)\big)\big| > \frac{i_{2t}}{3(1-p)}\Big) = \\
&P\Big(\mathcal{E}, \mathcal{T}(2t-1)=T, \Big|S\Big(\pi_2\big(\text{Taboo_Bond_Growth}\,(W,D,T,t)\big)\Big)\Big| > \frac{i_{2t}}{3(1-p)}\Big) \\
&\qquad\leq P\Big(\Big|S\Big(\pi_2\big(\text{Taboo_Bond_Growth}\,(W,D,T,t)\big)\Big)\Big| > \frac{i_{2t}}{3(1-p)}\Big)\,.
\end{aligned}$$

In the last inequality, we have forgotten the event $\mathcal{E}$ and the condition $\mathcal{T}(2t-1)=T$. This way we are left with the problem of controlling

$$S\Big(\pi_2\big(\text{Taboo_Bond_Growth}\,(W,D,T,t)\big)\Big)\,.$$

This is now possible, because the set $T$ is fixed, and it does not carry with him any information on the percolation configuration. As a consequence, the taboo growth algorithm working with the taboo set $T$ is again a genuine exploration algorithm. So the probabilistic control (28.2) can be used to control the term inside the sum of formula (31.71) and we obtain

$$P\Big(\mathcal{E}, \mathcal{T}(2t-1)=T, \big|S\big(\mathbb{E}^d(\mathcal{A}(2t),D)\big)\big| > \frac{i_{2t}}{3(1-p)}\Big) \leq 2|\mathbb{E}^d(D)|\exp\Big(-\frac{2p^2 i_{2t}^2}{9|\mathbb{E}^d(D)|}\Big)\,. \tag{31.72}$$

Substituting (31.72) in (31.71) and bounding the number of terms in the sum by $|\mathbb{E}^d(D)|^{i_1+\cdots+i_{2t-1}}$, we obtain

$$P\Big(\mathcal{E}, \big|S\big(\mathbb{E}^d(\mathcal{A}(2t),D)\big)\big| > \frac{i_{2t}}{3(1-p)}\Big) \;\leq\; 2|\mathbb{E}^d(D)|^{i_1+\cdots+i_{2t-1}} \exp\Big(-\frac{2p^2 i_{2t}^2}{9|\mathbb{E}^d(D)|}\Big)\,.$$

The first term of the right-hand side of (31.69) enjoys a similar inequality. Substituting these inequalities in (31.69), we conclude finally that

$$P(\mathcal{E}) \;\leq\; 6|\mathbb{E}^d(D)|^{i_1+\cdots+i_{2t-1}} \exp\Big(-\frac{2p^2 i_{2t}^2}{9|\mathbb{E}^d(D)|}\Big)\,.$$

This is the inequality (31.66) for even indices. □

We come back to the identity (31.65). We use the control on $|\mathcal{I}(1)|$ obtained in (31.6) and inequality (31.66) to control the terms inside the sum and we obtain

$$\begin{gathered} P\big(\exists s\in\{\,1,\dots,t\,\}\quad |\mathcal{I}(s)|>i_s,\, W \not\longleftrightarrow W'\big) \;=\\ 6|\mathbb{E}^d(D)| \exp\Big(-\frac{2p^2 i_1^2}{9|\mathbb{E}^d(D)|}\Big) + \sum_{2\leq s\leq t} 6|\mathbb{E}^d(D)|^{i_1+\cdots+i_{s-1}} \exp\Big(-\frac{2p^2 i_s^2}{9|\mathbb{E}^d(D)|}\Big)\,. \end{gathered} \tag{31.73}$$

We apply this result with $D=\Lambda(n)$ for some $n\geq 2d$ and we use the inequalities (31.7). The inequality (31.73) becomes

$$\begin{gathered} P\big(\exists s\in\{\,1,\dots,t\,\}\quad |\mathcal{I}(s)|>i_s,\, W \not\longleftrightarrow W'\big) \;\leq\\ 6d2^dn^d \exp\Big(-\frac{2p^2 i_1^2}{9d2^dn^d}\Big) + \sum_{2\leq s\leq t} 6(d2^dn^d)^{i_1+\cdots+i_{s-1}} \exp\Big(-\frac{2p^2 i_s^2}{9d2^dn^d}\Big)\,. \end{gathered} \tag{31.74}$$

For the sequence $i_1,\dots,i_t$, we choose

$$\forall s\in\{\,1,\dots,t\,\}\qquad i_s \;=\; \Big(\frac{3}{p}d2^d\Big)^s s!\big(n^d\ln n\big)^{1-2^{-s}}\,. \tag{31.75}$$

These values ensure that $1\leq i_1\leq i_2\leq\cdots\leq i_t$ and

$$\frac{2p^2 i_1^2}{9d2^dn^d} \;=\; 2d2^d\ln n \;\geq\; 3d\ln n\,. \tag{31.76}$$

Moreover, for $s\geq 2$, we have on the one hand

$$(i_1+\cdots+i_{s-1})\ln(d2^dn^d) \;\leq\; 2d(s-1)i_{s-1}\ln n\,, \tag{31.77}$$

and on the other hand

$$\begin{aligned} \frac{2p^2 i_s^2}{9d2^dn^d} \;&=\; \frac{2p^2}{9d2^dn^d}\Big(\frac{3}{p}d2^d\Big)^{2s}(s!)^2\big(n^d\ln n\big)^{2-2^{-s+1}}\\ &\geq\; 2d2^d\Big(\frac{3}{p}d2^d\Big)^{2s-2} s!\big(n^d\ln n\big)^{1-2^{-s+1}}\ln n \;=\; 2d2^d s i_{s-1}\ln n\,. \end{aligned} \tag{31.78}$$

Inequalities (31.77) and (31.78) together yield that

$$(d2^d n^d)^{i_1+\cdots+i_{s-1}} \exp\Big(-\frac{2p^2 i_s^2}{9d2^d n^d}\Big) \,\leq\, \frac{1}{n^d}\,. \tag{31.79}$$

Substituting the inequalities (31.76) and (31.79) in (31.74), we obtain that

$$P\big(\exists s\in\{\,1,\dots,t\,\}\quad |\mathcal{I}(s)|>i_s,\, W\not\longleftrightarrow W'\big) \,\leq\, \frac{6t}{n^d}\,.$$

Of course, these inequalities are interesting only when the values $i_s$ are significantly smaller than $n^d$. Unfortunately, the exponent $1-2^{-s}$ goes to 1 very fast. For a fixed value of $t$, we see that, with probability going fast to 1 (the corresponding speed could be considerably improved by taking slightly larger values for the sequence $i_1,\dots,i_t$), on the event $\big\{\,W\not\longleftrightarrow W'\,\big\}$, we will have

$$\forall s\in\{\,1,\dots,t\,\}\qquad |\mathcal{I}(s)|\,\leq\,\Big(\frac{3}{p}d2^d\Big)^s s!\big(n^d\ln n\big)^{1-2^{-s}}\,. \tag{31.80}$$

What is important here is the exponent of $n$, the factor $\ln n$ can be seen as a correcting factor. The first values in the sequence are

$$n^{d/2},\, n^{3d/4},\, n^{7d/8},\, n^{15d/16},\,\dots$$

Unfortunately, these bounds are too weak and does not allow us to conclude something interesting in any of the three challenging questions of section 30.

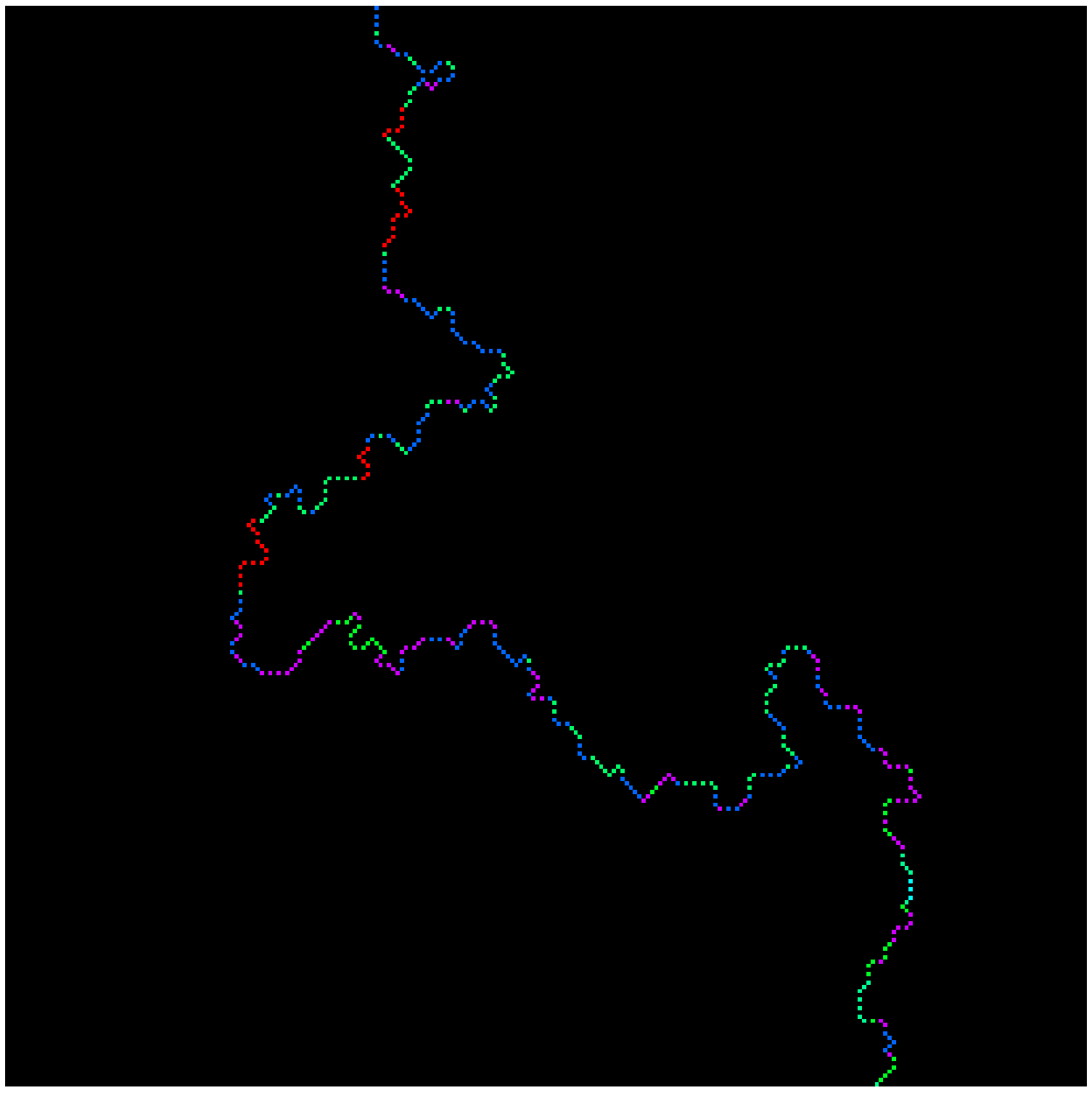

Figure 59: Intertwined bond explorations in the box $\Lambda(128)$ with $p=0.49$. $W,W'=$ left,right boundaries. Intersection sets $\mathcal{I}(1)$, $\mathcal{I}(2)$, $\mathcal{I}(3)$, $\mathcal{I}(4)$, $\mathcal{I}(5)$.

# 32 Simulation of the intersection sets

This section is devoted to the simulation of the intersection sets. The next subsection discusses the implementation of the graph algorithms required to simulate and study the disconnection problem in percolation. Mathematicians not interested in this aspect can skip this subsection and jump to the next one, which presents and comments the simulations.

## 32.1 The free C and C++ graph libraries

When it comes to simulate the intersection sets, we are confronted to a classical computing problem: what is the best way to find the connected components in a graph? What matters is not so much the algorithm itself, but the choice of the computer representation of the graph. We already started discussing this question in section 5, where the experimental curve of the function $\theta(p,\mathbb{Z}^3)$ was presented. For the simulations of $\theta(p,\mathbb{Z}^3)$, we opted for a graph representation closely linked to the geometry of the lattice, because the memory was the limiting factor. However, we will very soon face algorithmic problems of increasing complexity when we wish to simulate the higher intersection sets. Indeed, we will have to implement the intertwined explorations of the graph starting from the sets $W$ and $W'$. For this, we will use repeatedly a breadth-first search of the graph. Later on, when we will try to simulate the second positive pivots, we will use the beautiful and classical algorithm of Tarjan (see subsection 45.7), which involves a depth-first search of the graph. For these problems, it is much better to implement the relevant algorithms in a generic way, which is not dependent on the specific geometry of the 2d or 3d lattices. Obviously, we will lose out on memory, but the algorithms will be much more reliable and reusable, for instance to go from dimension two to three. Needless to be said, these classical graph algorithms have been thoroughly studied in computer science, and we would like to use efficient existing implementations, designed by professional computer scientists. So we looked at the graph libraries available in the free software category, and mainly in C and C++ (we dismissed the libraries in Python because we would like to maximize the computational speed). The free libraries are available under various licenses, namely: the Boost Software License (BSL), the Berkeley Software Distribution license (BSD), the Gnu General Public License (GPL) and its Lesser version (LGPL), the Common Public License (CPL). The licenses are indicated below in parentheses in abbreviated form.

Let us start with C libraries. A first possibility is libcgraph, which is part of the Graphviz package (CPL) [49], a powerful software for visualizing graphs. Graphviz is certainly of interest for studying percolation on graphs more exotic than the lattice $\mathbb{Z}^d$. A second possibility is the igraph library (GPL) [34]. A general drawback of the C libraries is that the naming conventions are quite heavy (long names for functions and a lot of arguments for each function call). However igraph does also have C++ bindings. Moreover it is still being maintained and developed. We come now to C++ libraries, which should provide a more readable and transparent interface. The first very interesting candidate is the Boost Graph Library (BSL) [121]. This is precisely a generic library which implements a generic graph structure, together with generic graph algorithms. The drawback is that the level of abstraction of this library is rather high, which makes it harder to use. It seems also that the high level of abstraction of this library leads to a loss of efficiency and speed. A second candidate is the LEDA library [104], which is free for academic usage. This is a quite general library implementing a large range of data structures and types. A third candidate is the Lemon library (BSL). It implements various data structures and algorithms, with a focus on combinatorial optimization, so its purpose is a bit different from what we seek. Moreover the last release is from 2014. A fourth candidate is the CGAL Library (GPL v3+), which implements sophisticated geometric algorithms of various types. Although it is very interesting, it seems too high-level for our percolation problem. A fifth candidate is the Stanford Network Analysis Platform, SNAP (BSD) [94]. This is a comprehensive library, which is able to handle very large graphs and has been used with success for various large scale projects. It is well documented and still actively developed. A sixth candidate is the graphtool library (LGPL v3). The interface is very transparent and easy to use. There exist also other smaller interesting graph libraries, e.g., the graph template library, that we did not investigate because of the lack of time. At the end of the day, there is no obvious choice, each library has its own advantages and disadvantages. So we tried a bit some of the most promising libraries, with the following (certainly superficial and biased) outcome. The graphtool library was by far the easiest to use and the quickest to lead to nice pictures. However, things were getting more complicated when we had to implement specific algorithms, because naturally we had to go beyond the basic interface of the library. The SNAP library seems to be the most suited for our purpose. Unfortunately, it is rather heavy to put at work and it would require also further implementation of classical algorithms (because we always need some specific variant of what is available). We decided finally to implement the simulation programs from scratch. This might sound quite absurd in the 21st century. One can wonder why the C++ Standard Template Library (STL in short) does not contain basic routines to deal with graphs. Digging deeper into this issue, we became aware of a Graph Library proposal for the STL. The proposal is not to build a new independent library, like the Boost Graph Library, but rather to extend the existing STL library with an adequate graph container interface, equipped with the classical graph algorithms. This project evolved from the library NWGraph [99]. It is currently under development, with experimental versions available on the Github site under the name stdgraph-v2, along with several documents explaining the goals of the project. It might be integrated to the STL library in the future, and it would certainly have been our preferred choice if it was already available. In any case, the current conclusion of these investigations is that a sensible choice is to build the adequate structure for our graph problem with the help of the existing tools of the STL library, and typically to represent the graphs with the help of STL vectors.

## 32.2 The simulations.

The graph algorithms are implemented with the help of C++ Standard Template Library. The base graph is represented as a STL vector of STL vectors, each entry corresponding to a vertex and its adjacency list (consisting either of its neighbouring vertices, or its incident edges, or both). This is the way we chose to implement various simulations revolving around the disconnection problem in percolation. As explained in subsection 32.1, the simulations of the intersection sets and of the pivotal bonds of higher orders involve various interesting graph algorithms. However, since they are not relevant for the mathematical proofs, we will not present them in detail in the text.

We present three different groups of simulations. In each figure, the upper graphics are for dimension 2 and the lower graphics for dimension 3. The first group is a series of simulations of the intersections sets in which the box $\Lambda$ is fixed and $p$ varies from 0 to 1. We consider the four boxes $\Lambda(8)$, $\Lambda(16)$, $\Lambda(32)$, $\Lambda(64)$ in dimensions two and three (figures 60, 61, 62, 63, 64, 65). The second group is a series of simulations of the first seven intersections sets in which the parameter $p$ is fixed and the side length $n$ of the box $\Lambda$ varies from 0 to 1000. We consider the following values for the parameter $p$: 0.49, 0.5, 0.501, 0.51 in dimension two and 0.23, 0.2488, 0.26, 0.27 in dimension three (figures 66, 67, 68, 69). The third group is a series of simulations of the "maximal" intersection sets (figures 70 and 71). In each simulation, we tried to achieve the maximal value for the first pivotal set, through the following procedure. We start from the configuration with all bonds closed. At each step, we choose randomly and uniformly a closed bond and we open it, if it is not pivotal. We stop when the chosen bond is pivotal. The configuration created this way has a complicated distribution, but in some sense, it yields configurations realizing the disconnection event which have the largest pivotal set. Naturally the outcome of the simulations are very unstable when the parameter $p$ is around the critical point (which is the focus of our interest), so we performed a huge averaging to stabilize the curves.

The numerical results seem to indicate that the brutal controls on the intersection sets obtained in section 31 are widely overestimated. Furthermore, it looks like that the intersections sets of order higher than two, except the last one, have comparable cardinalities! Another good point is that the behavior of the various curves are similar in dimensions two and three. The main lesson we can draw from these empirical results is that we have not managed to grasp the mathematical structure of the intersection sets. So we should try a different approach to obtain an adequate control on their cardinalities. The good news is that the curves do not show the dreaded explosion of the second pivotal bonds around the critical point, when the size of the box increases. If we manage to show that, in a $d$-dimensional box, the cardinalities of the very first intersection sets grows slower than $n^{d/2}$, then we will be in good position to fill in the missing step of the strategy outlined in subsection 29.7. Even in the case of the third group of simulations, for an algorithm which tries to create in a reasonable way the largest pivotal sets, the rate of increase of the cardinalities as functions of $n$ seems to stay within controllable bounds.

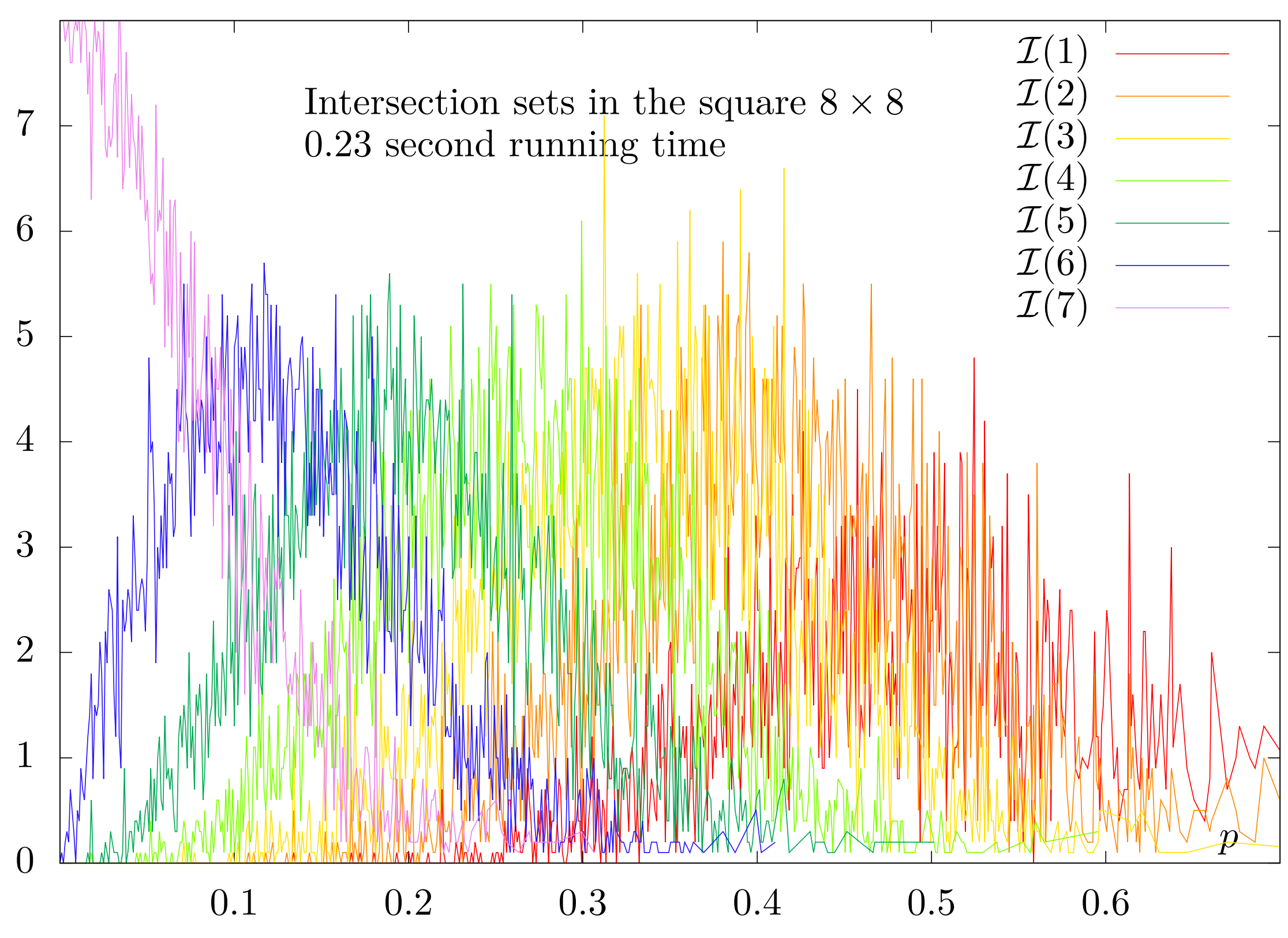

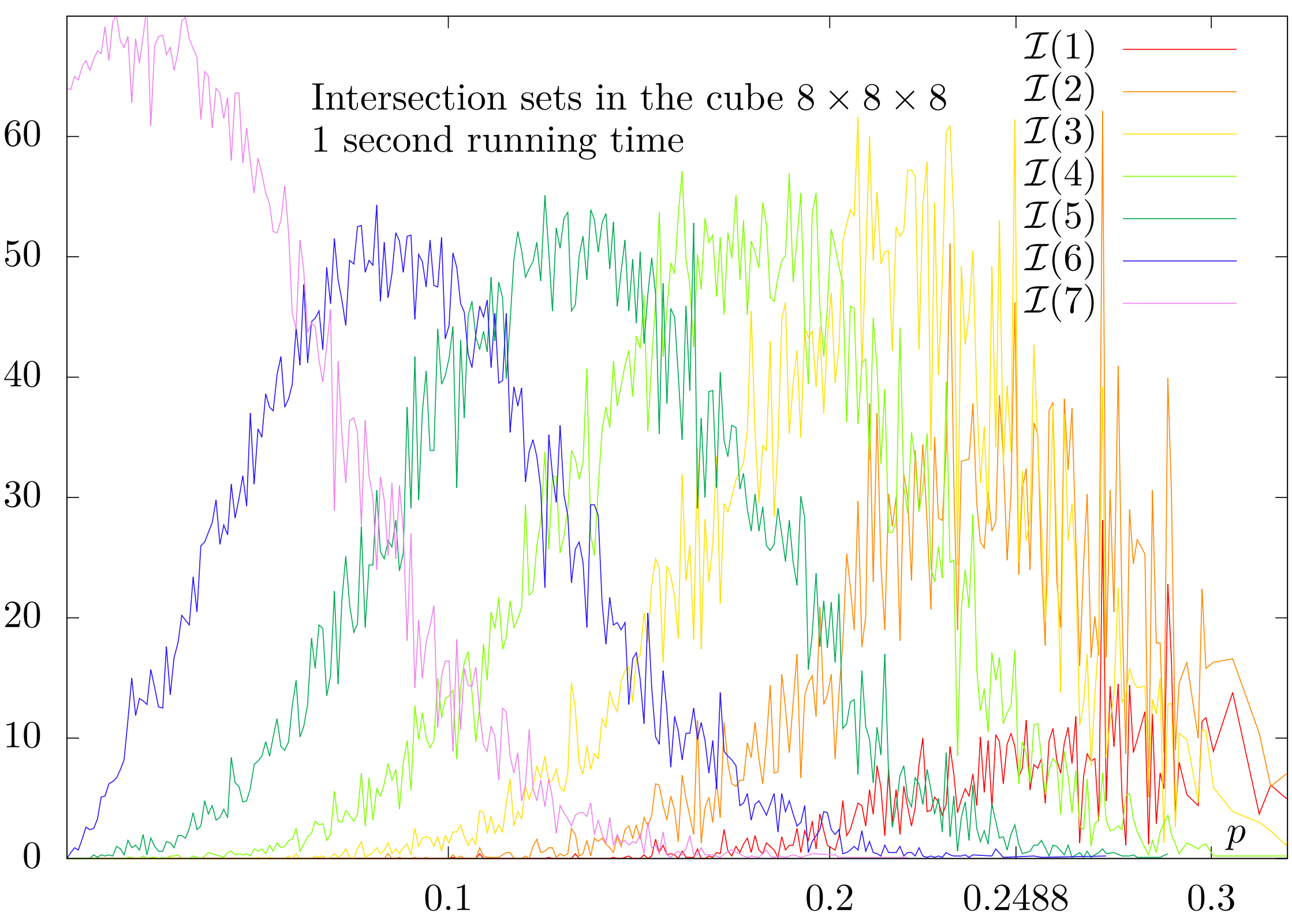

Figure 60: Intersection sets in the square and cube of side 8, as function of $p$. Each point is the average of 10 simulations. The increment is $\Delta p = 0.001$.

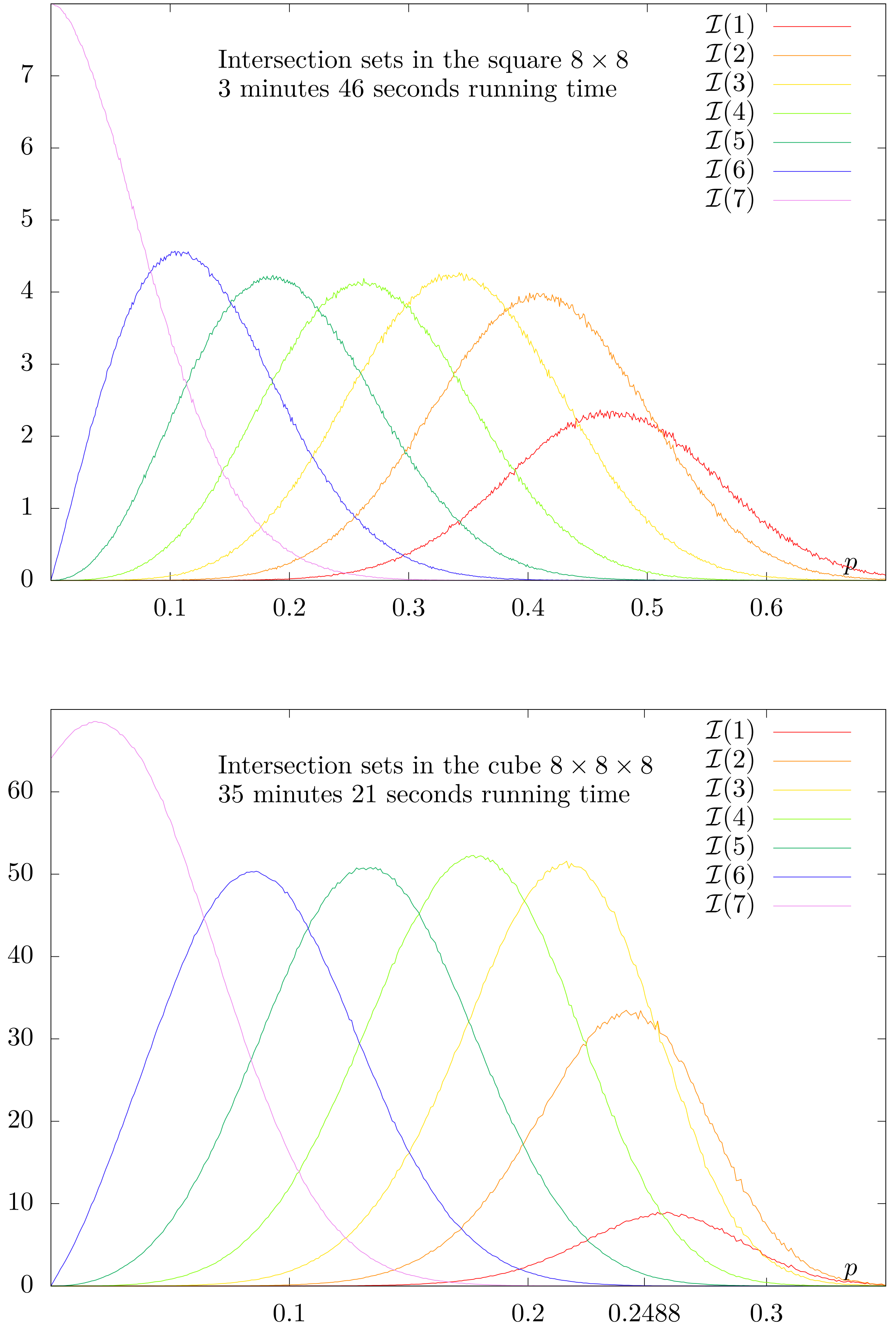

Figure 61: Intersection sets in the square and cube of side 8, as function of $p$. Each point is the average of 10000 simulations. The increment is $\Delta p = 0.001$.

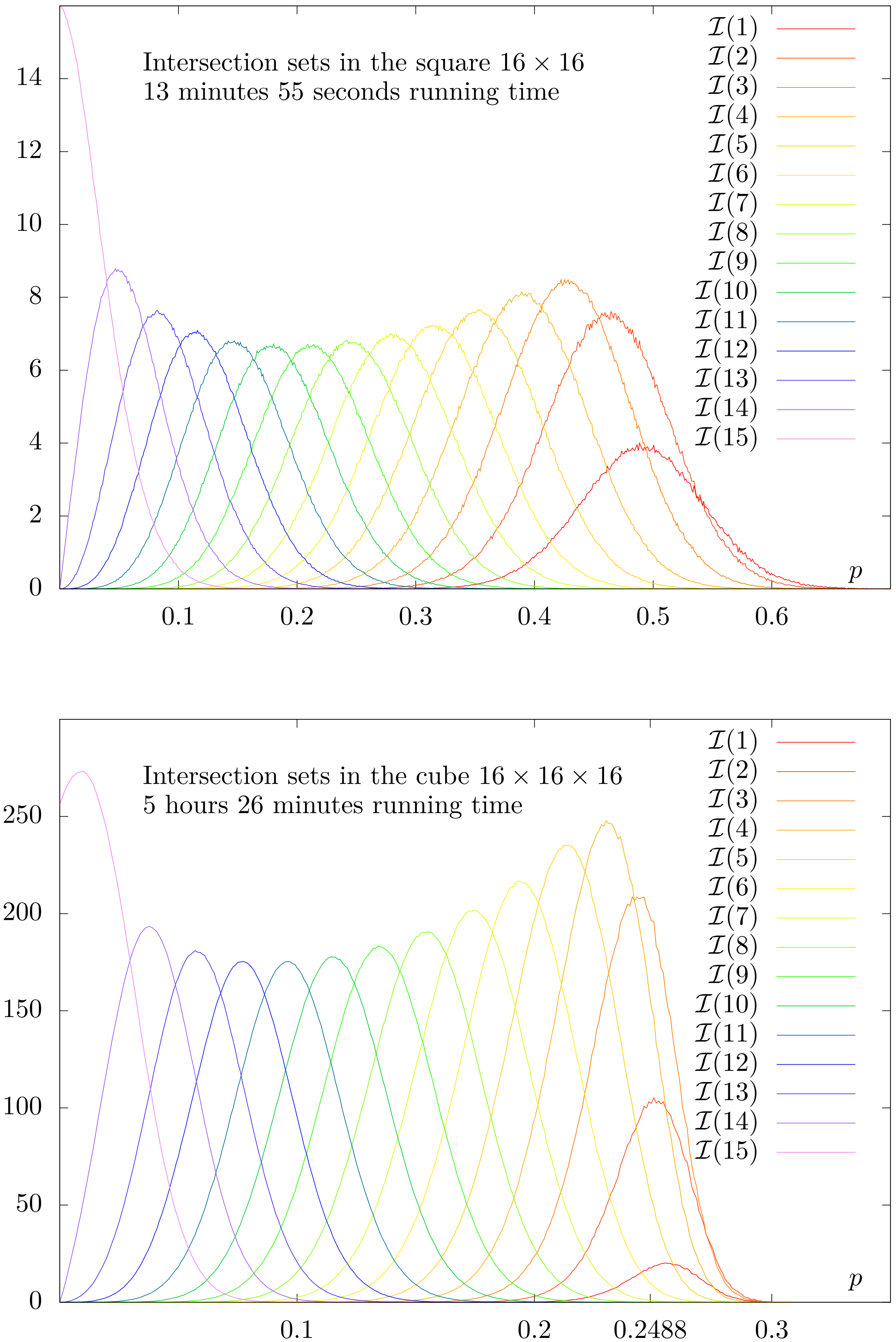

Figure 62: Intersection sets in the square and cube of side 16, as function of $p$. Each point is the average of 10000 simulations. The increment is $\Delta p = 0.001$.

Intersection sets in the square 32 × 32
52 minutes 14 seconds running time

15
10
5
0
0.1 0.2 0.3 0.4 0.5
$p$

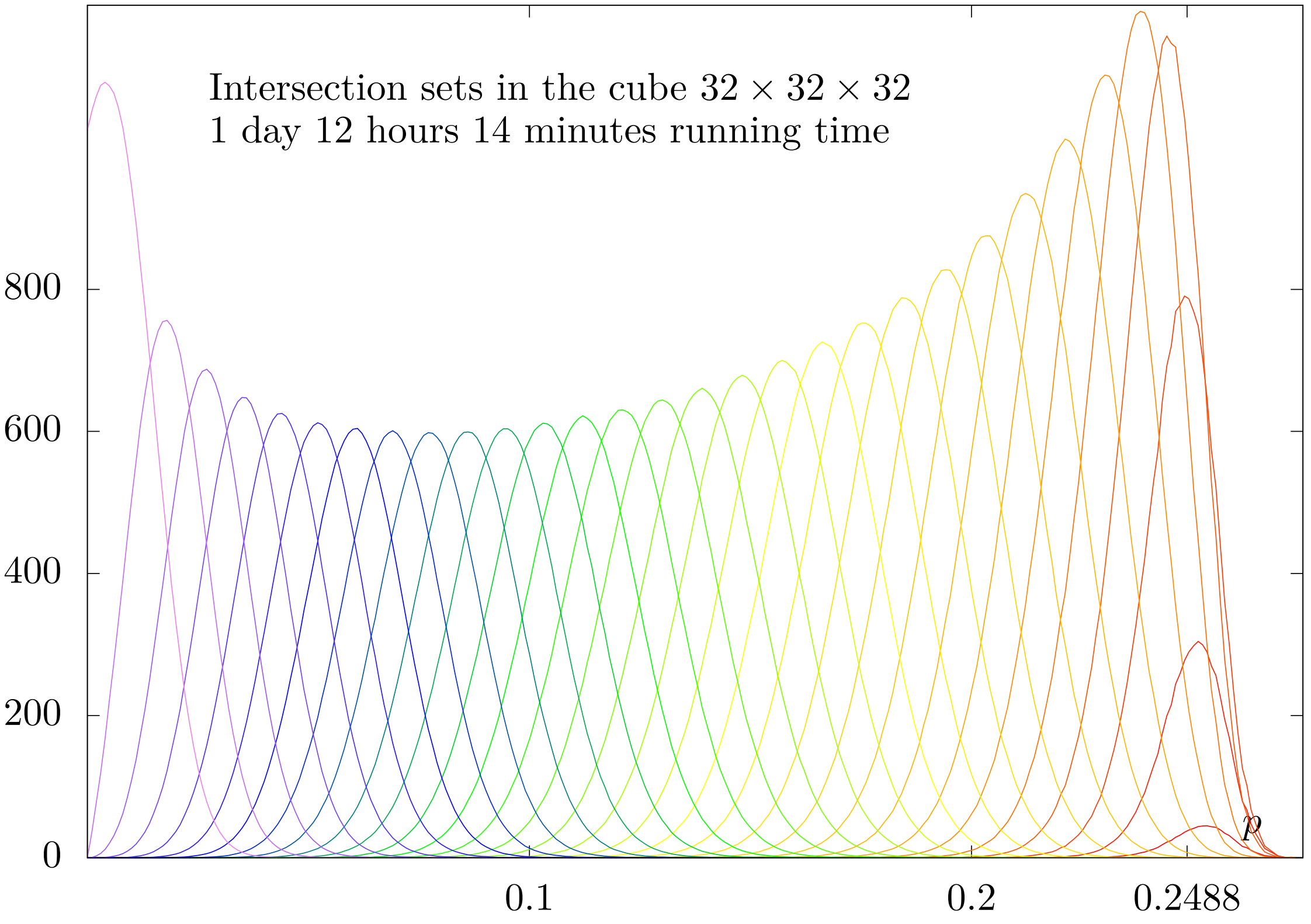

Figure 63: Intersection sets in the square and cube of side 32, as function of $p$. Each point is the average of 10000 simulations. The increment is $\Delta p = 0.001$.

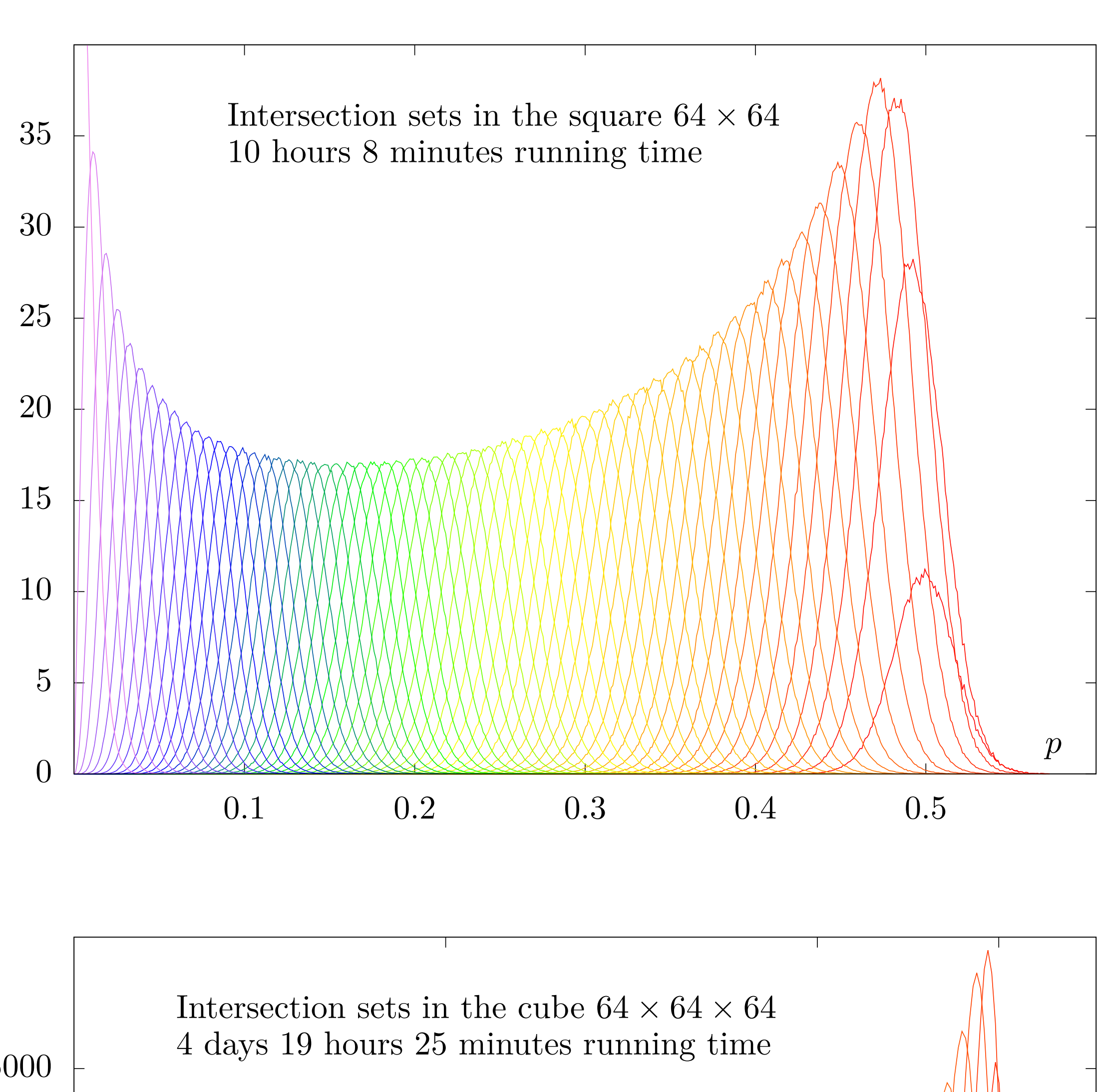

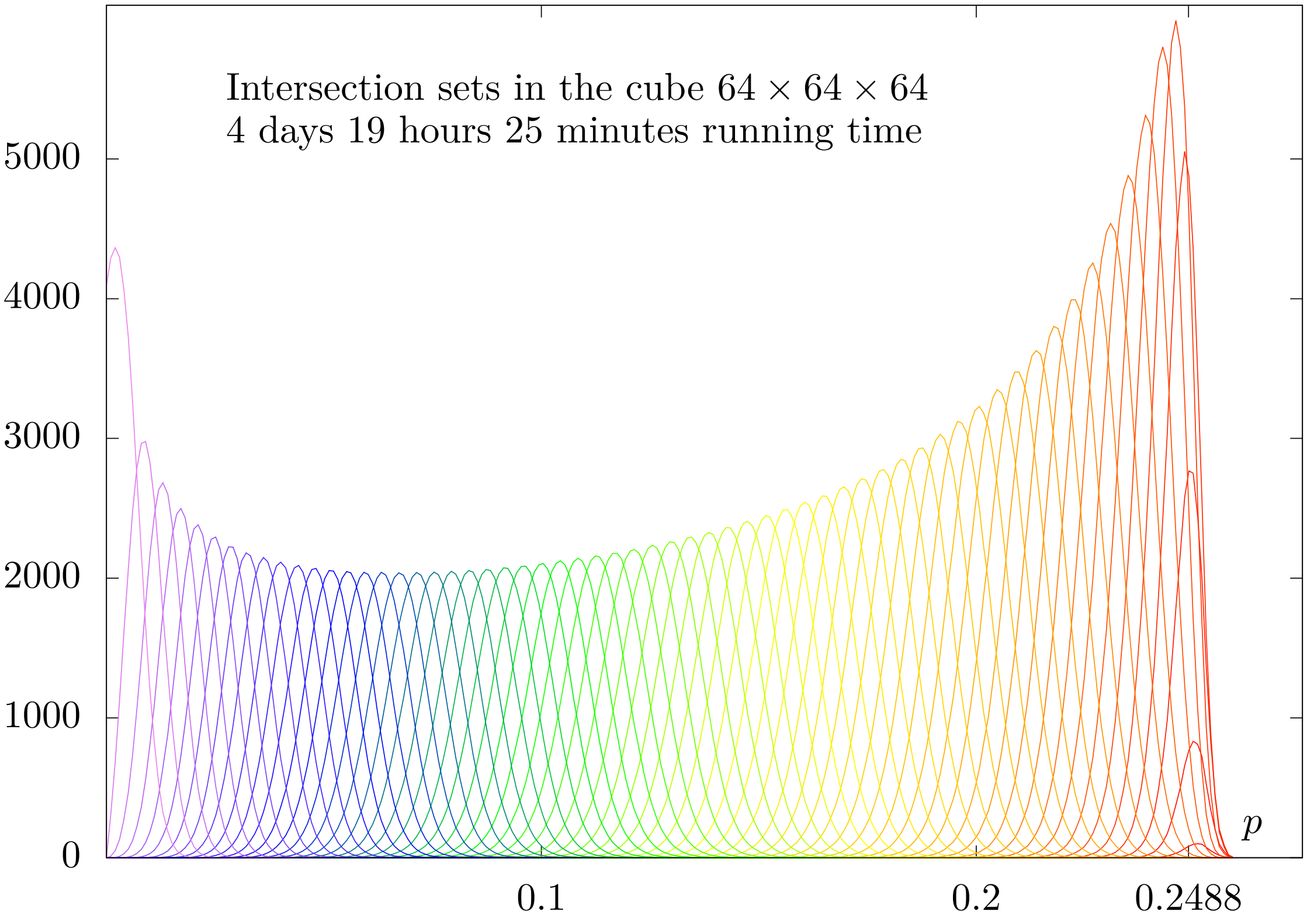

Figure 64: Intersection sets in the square and cube of side 64, as function of $p$. Each point is the average of 10000 simulations. The increment is $\Delta p = 0.001$.

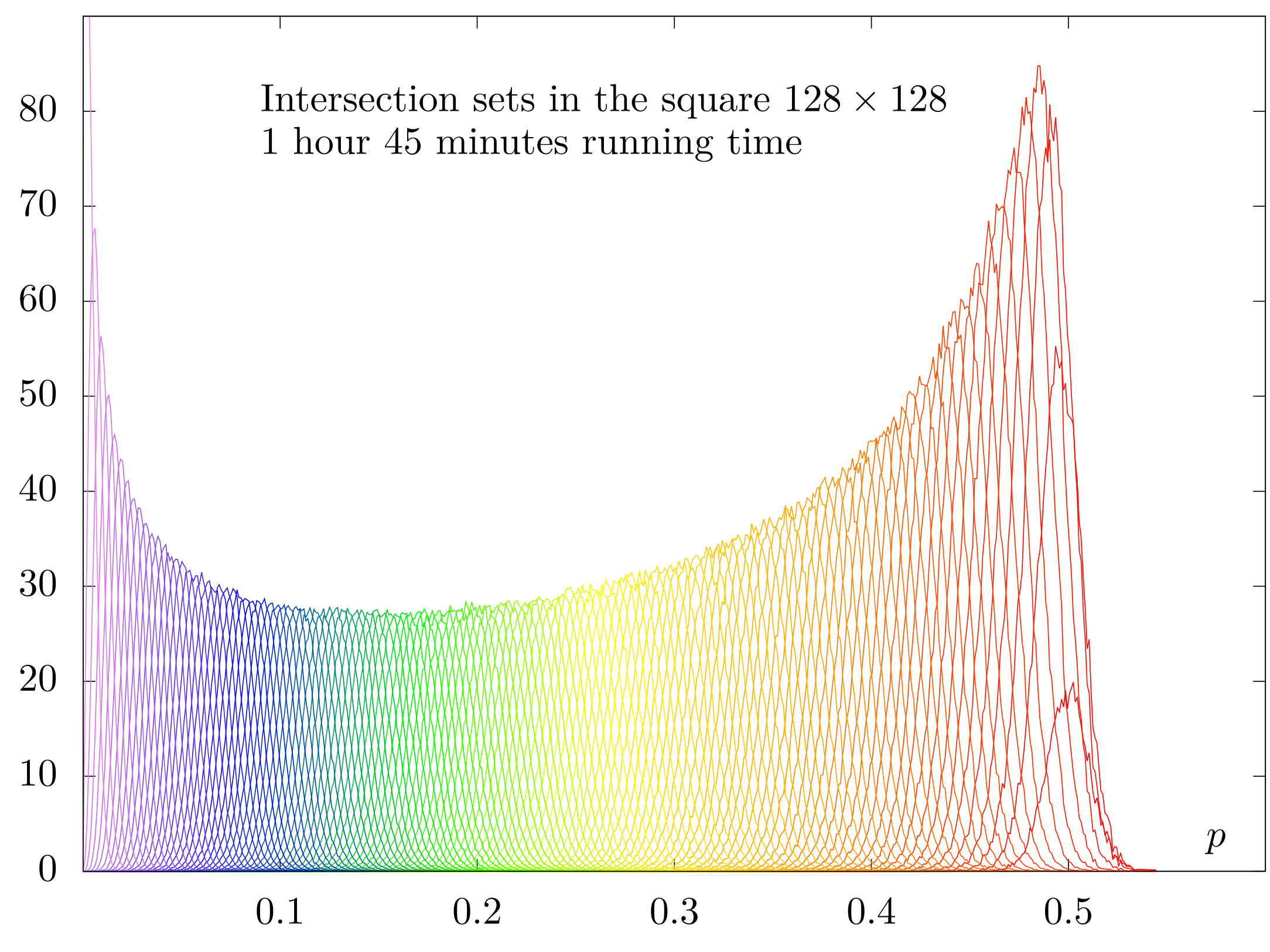

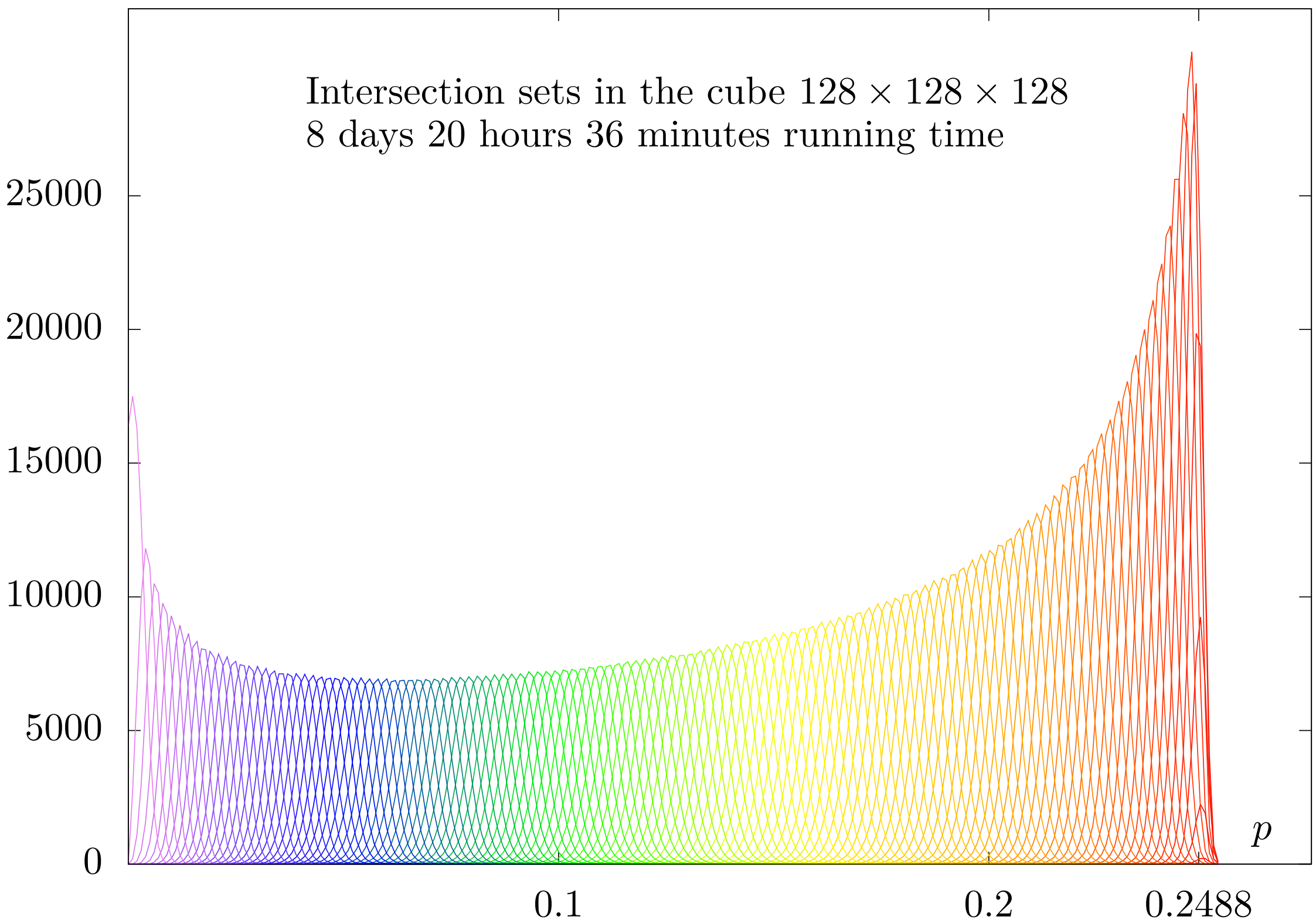

Figure 65: Intersection sets in the square and cube of side 128, as function of $p$. Each point is the average of 1000 simulations. The increment is $\Delta p = 0.001$.

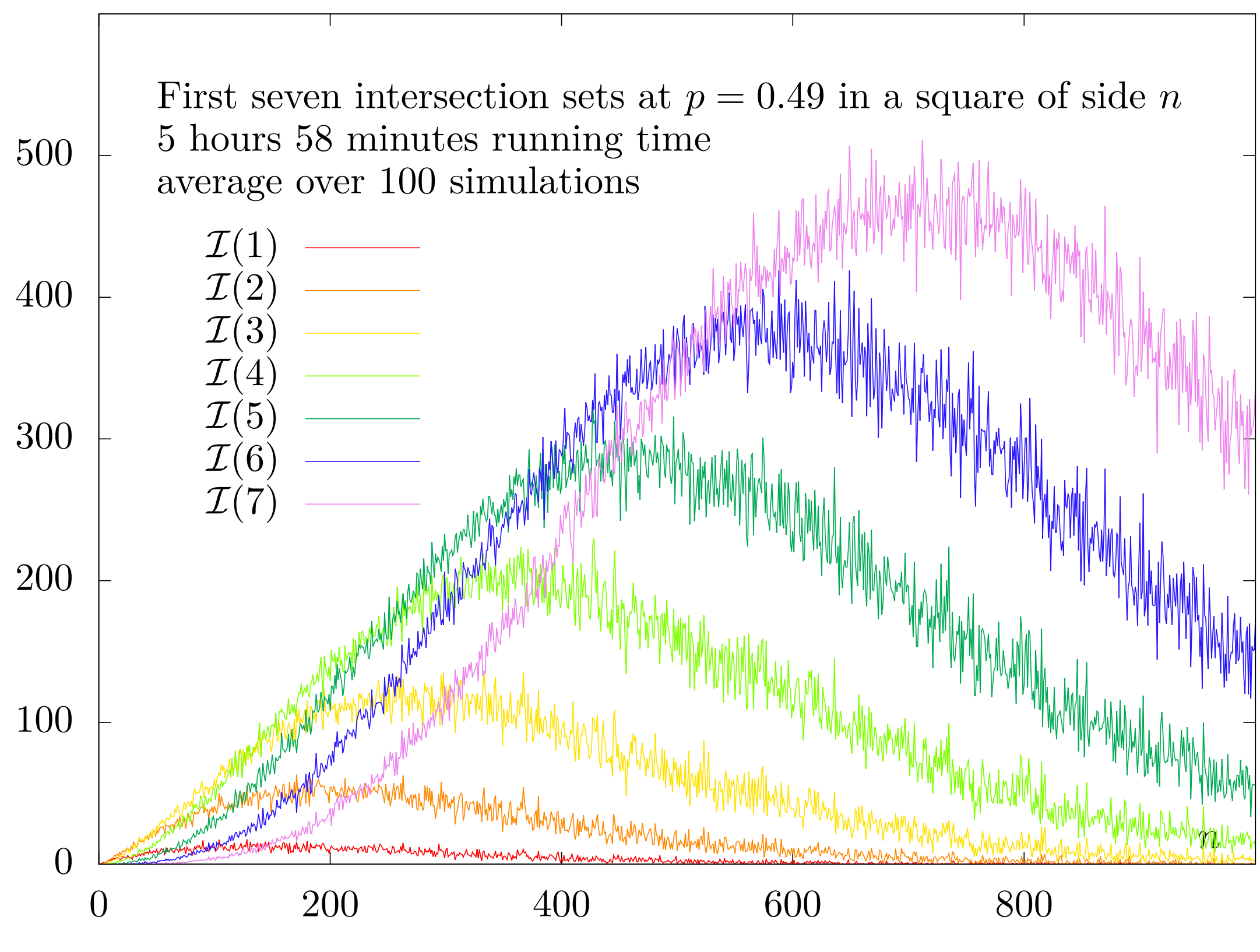

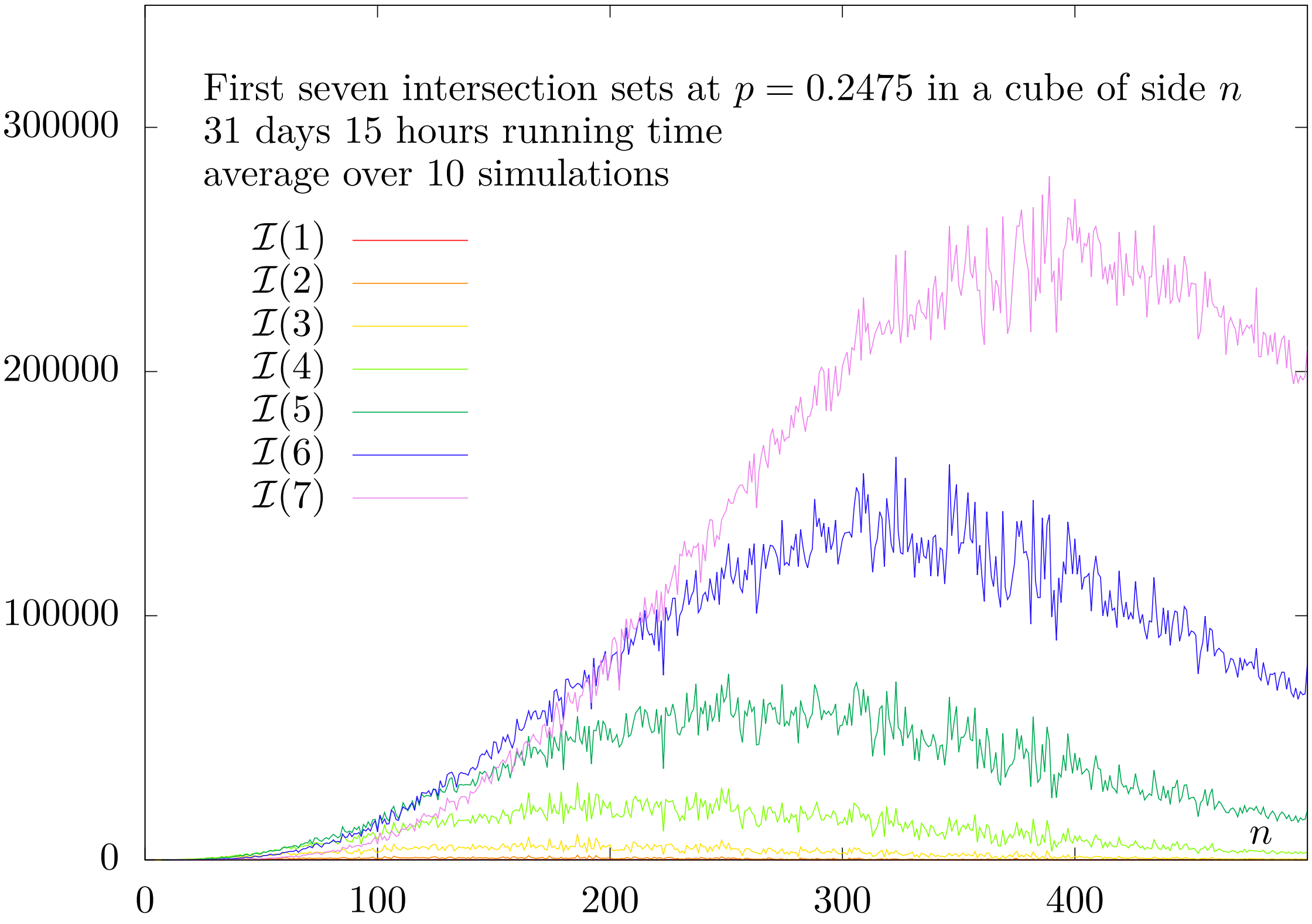

Figure 66: The first seven intersection sets, as functions of the side length $n$, in a square at $p = 0.49$ and in a cube at $p = 0.2475$. The increment for $n$ is $\Delta n = 1$.

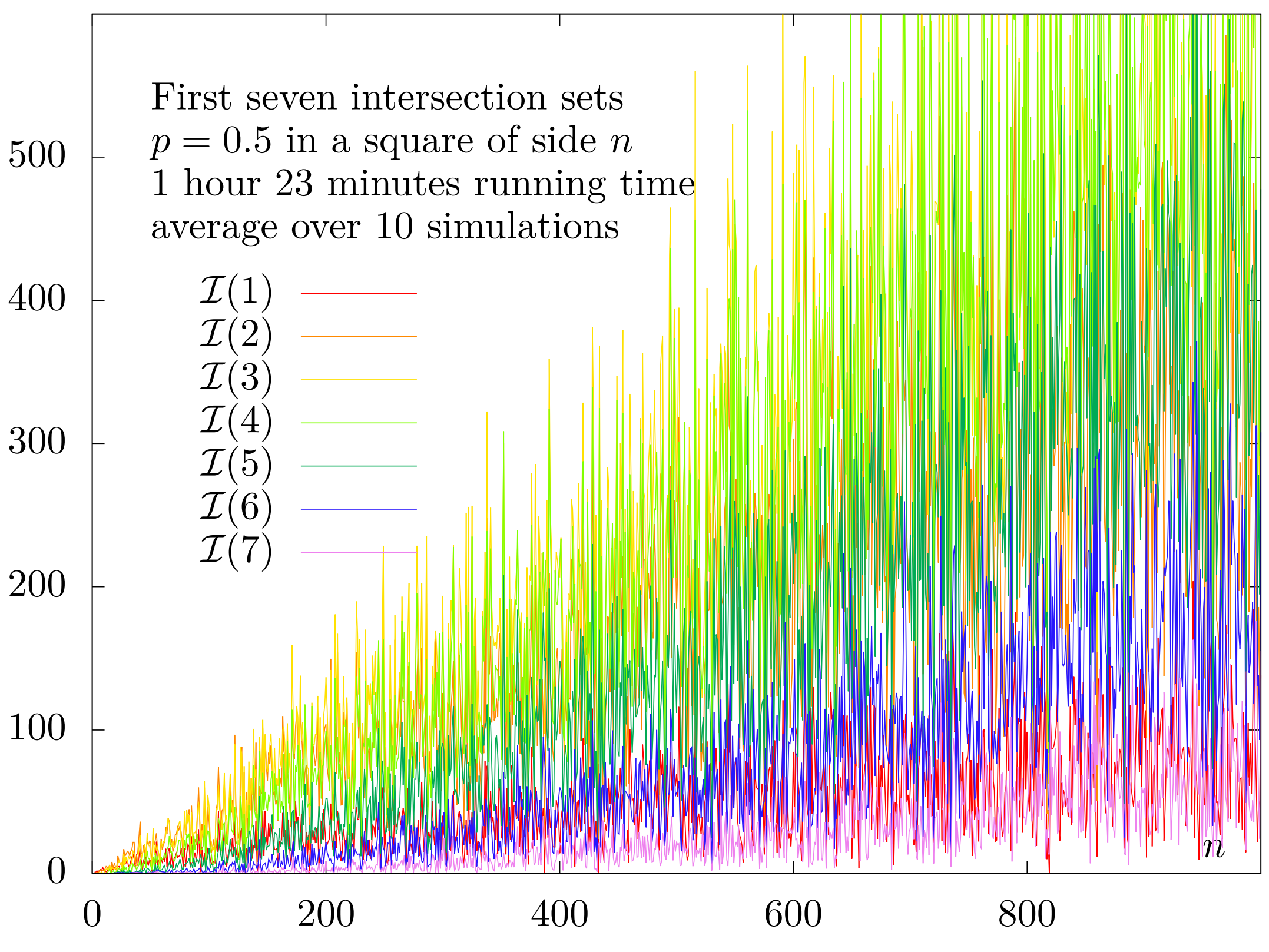

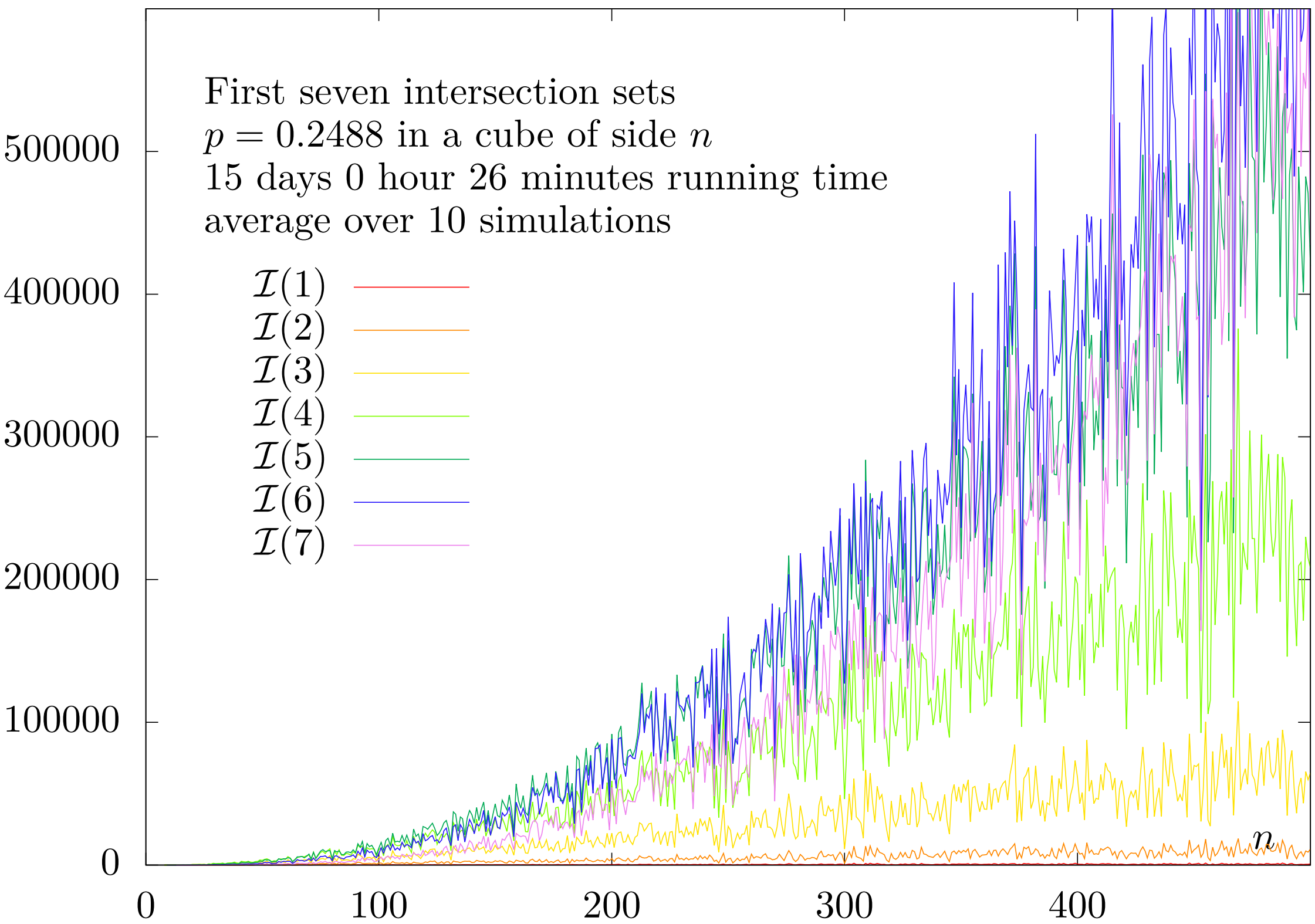

Figure 67: The first seven intersection sets, as functions of the side length $n$, in a square at $p = 0.5$ and in a cube at $p = 0.2488$. The increment for $n$ is $\Delta n = 1$.

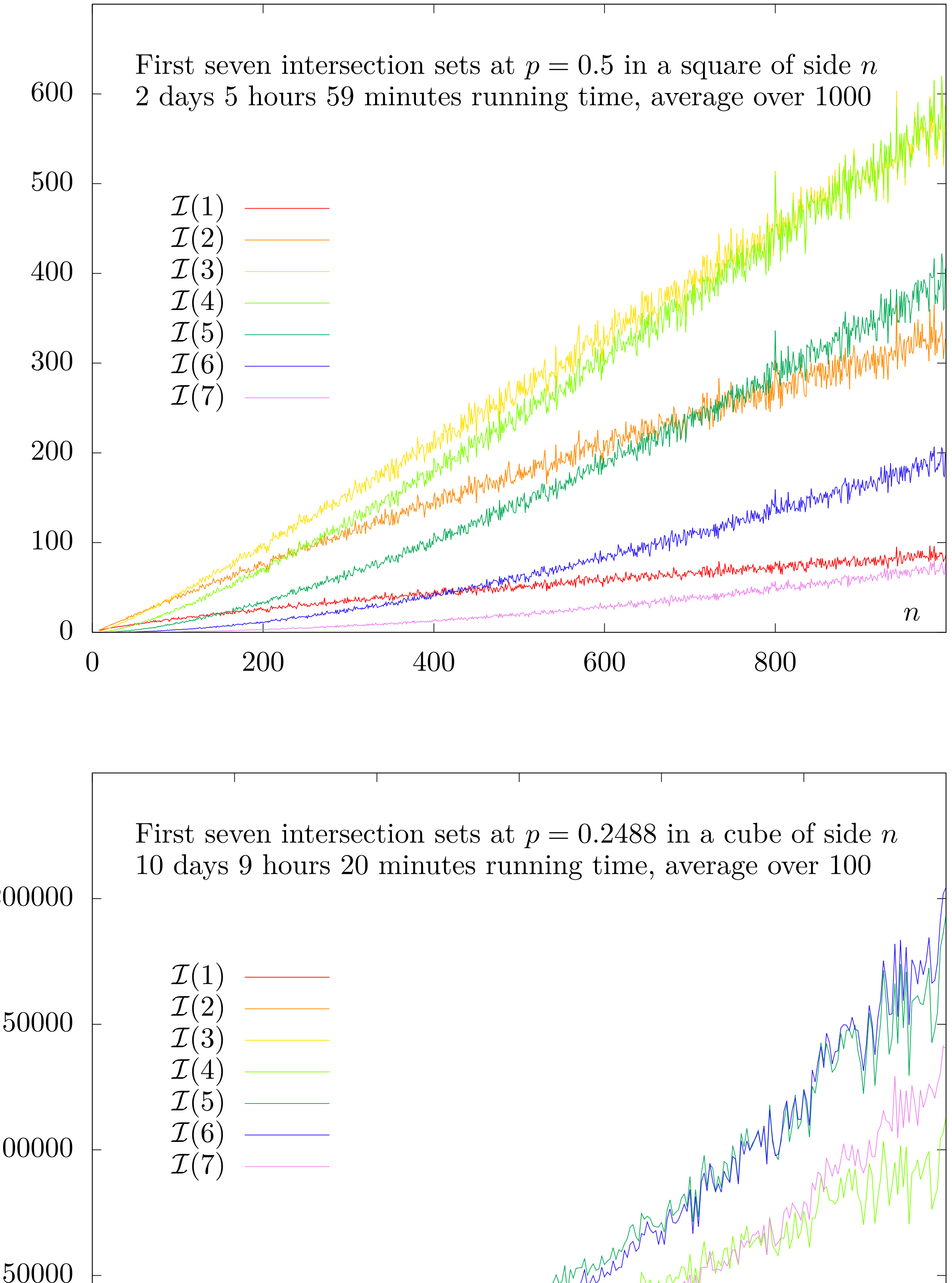

Figure 68: The first seven intersection sets, as functions of the side length $n$, in a square at $p = 0.5$ and in a cube at $p = 0.2488$. The increment for $n$ is $\Delta n = 1$.

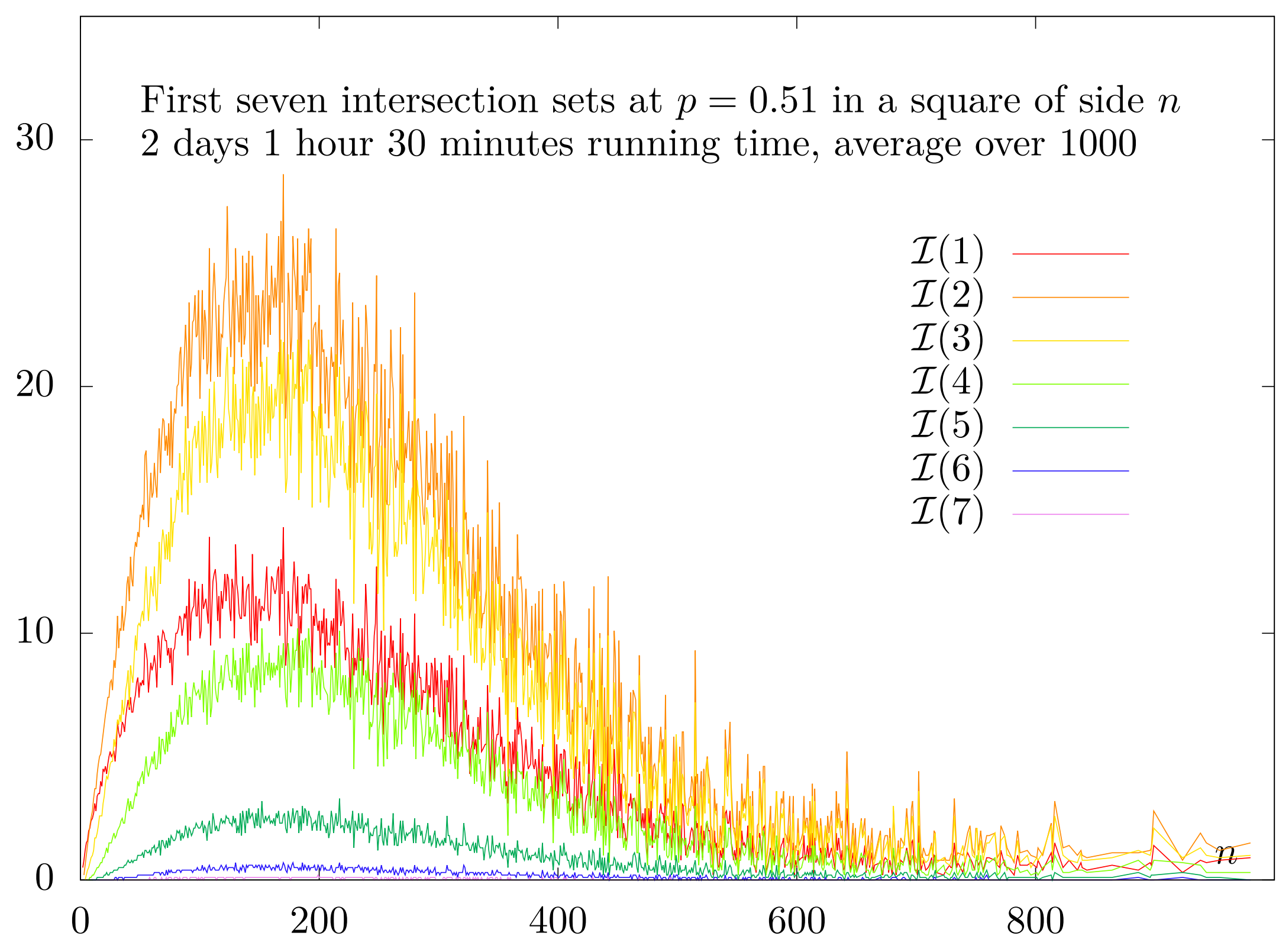

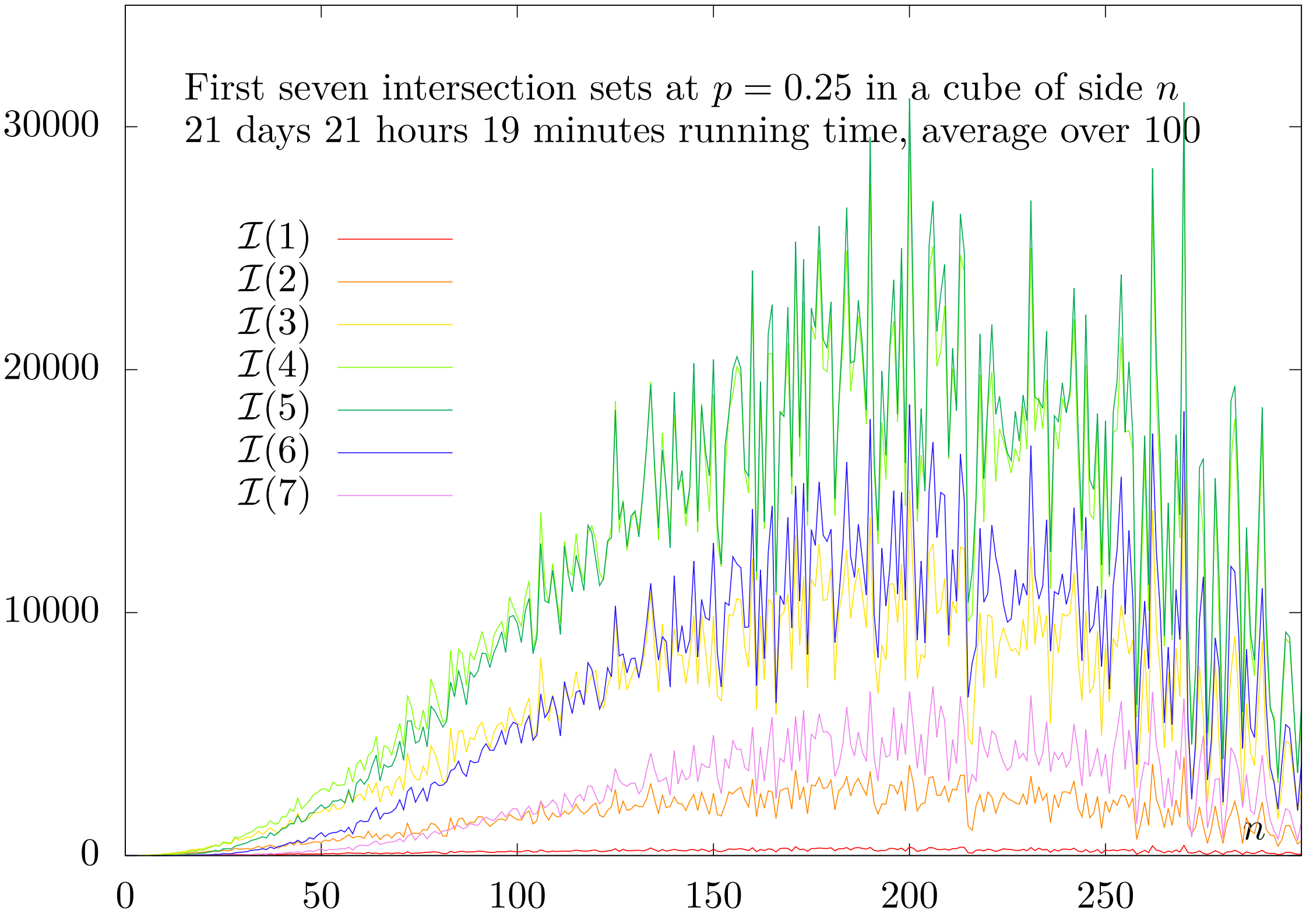

Figure 69: The first seven intersection sets, as functions of the side length $n$, in a square at $p = 0.51$ and in a cube at $p = 0.25$. The increment for $n$ is $\Delta n = 1$.

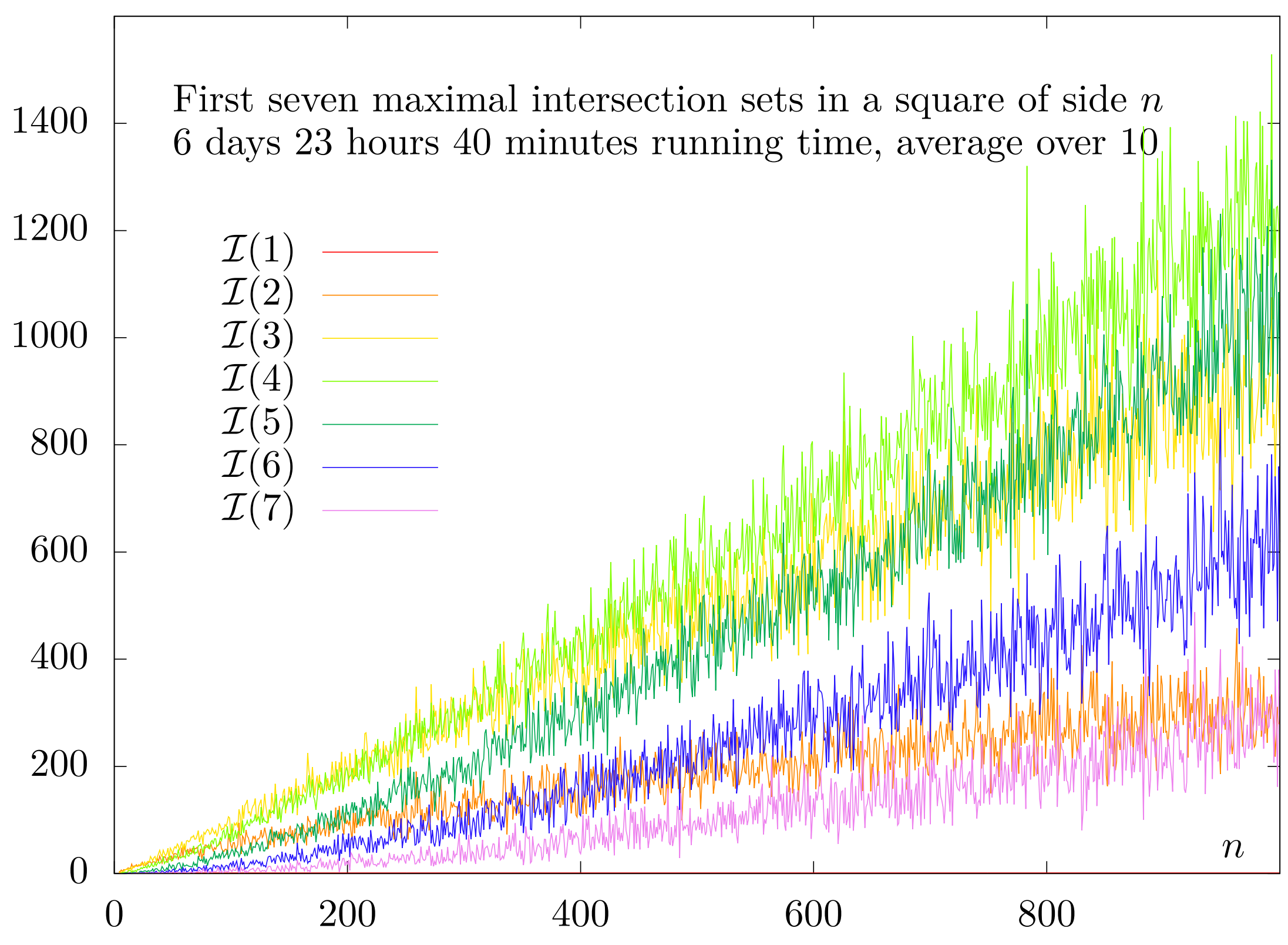

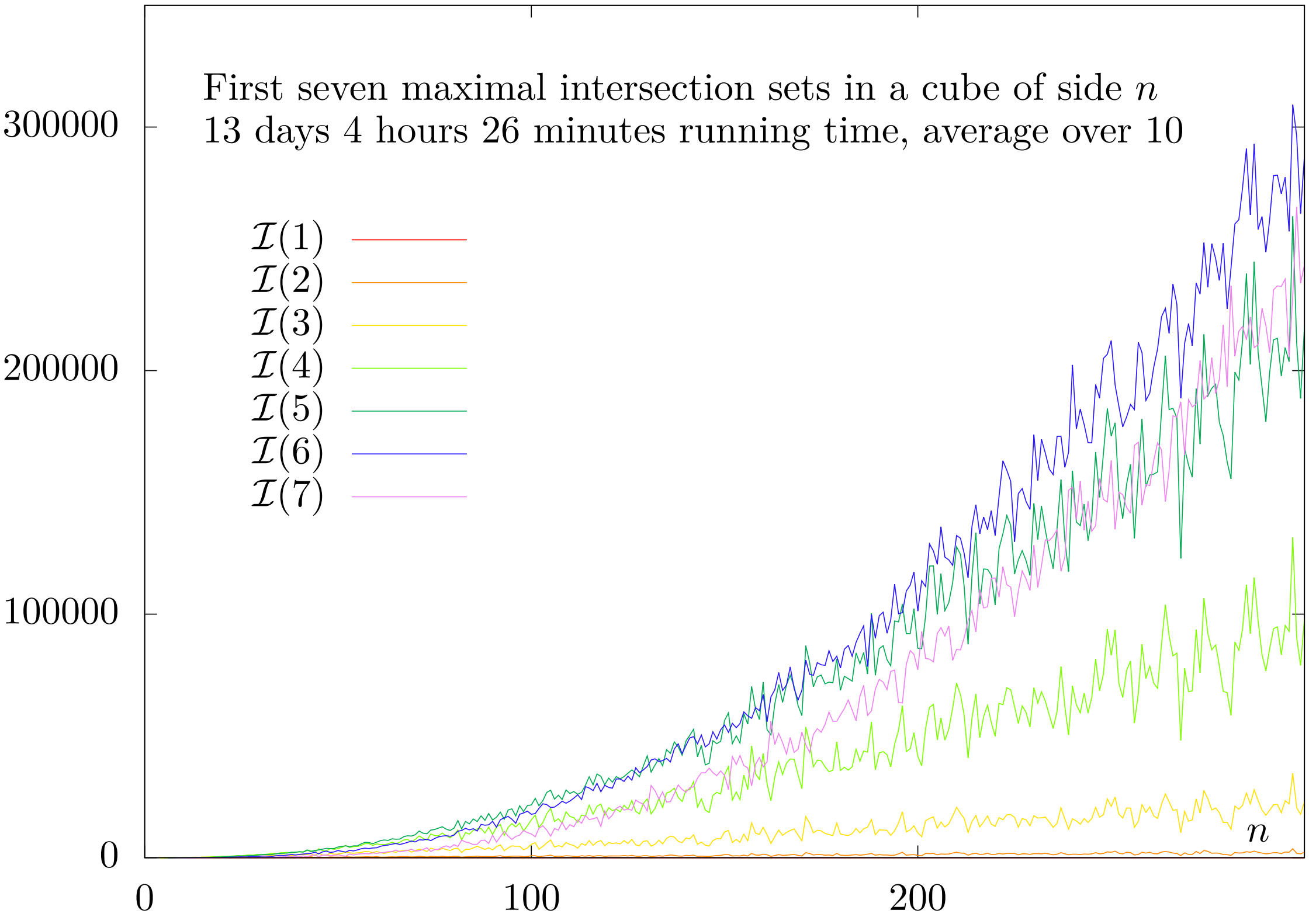

Figure 70: The first seven maximal intersection sets in a square and a cube: ultimate configuration before the connexion. The increment for $n$ is $\Delta n = 1$.

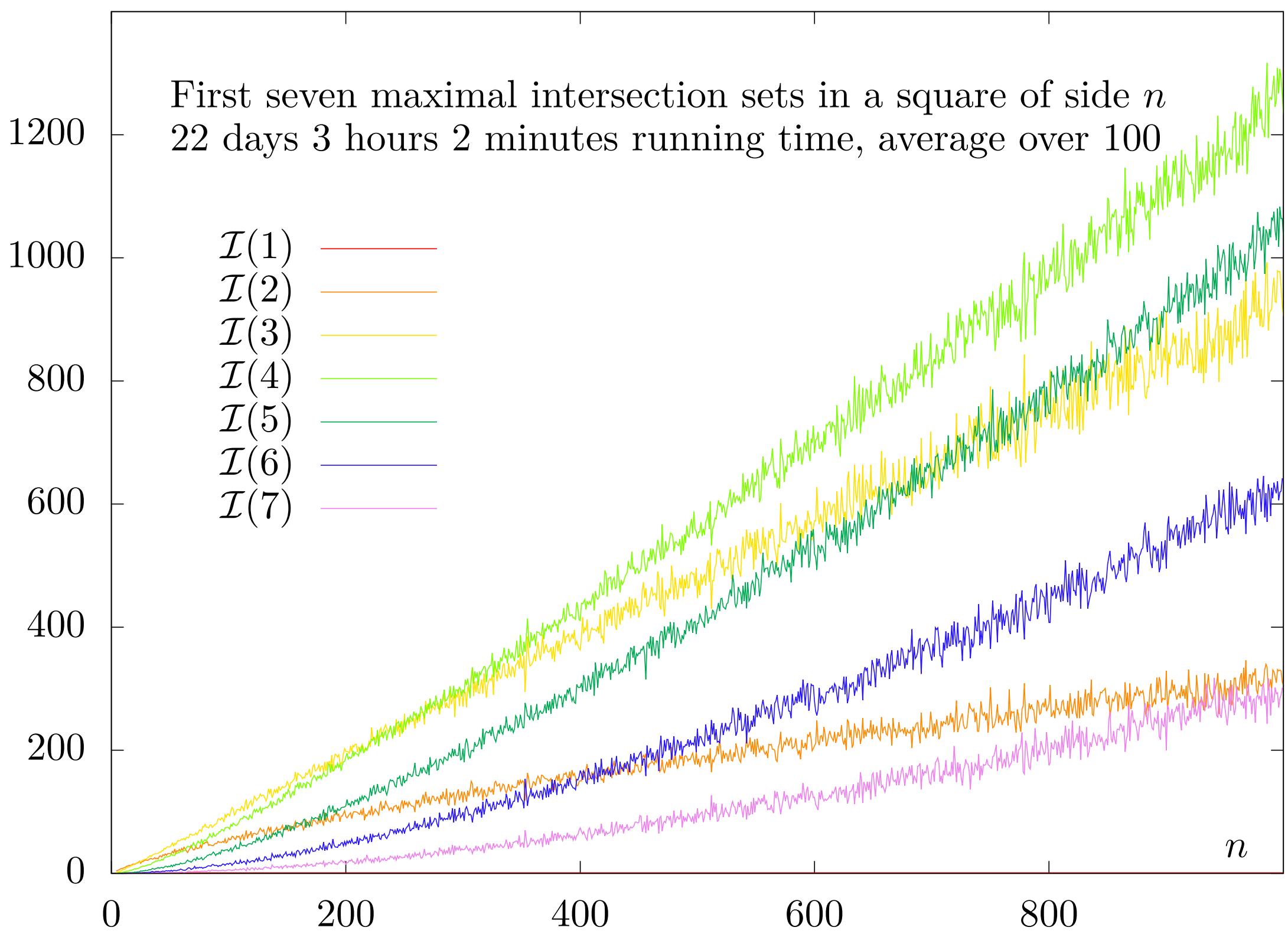

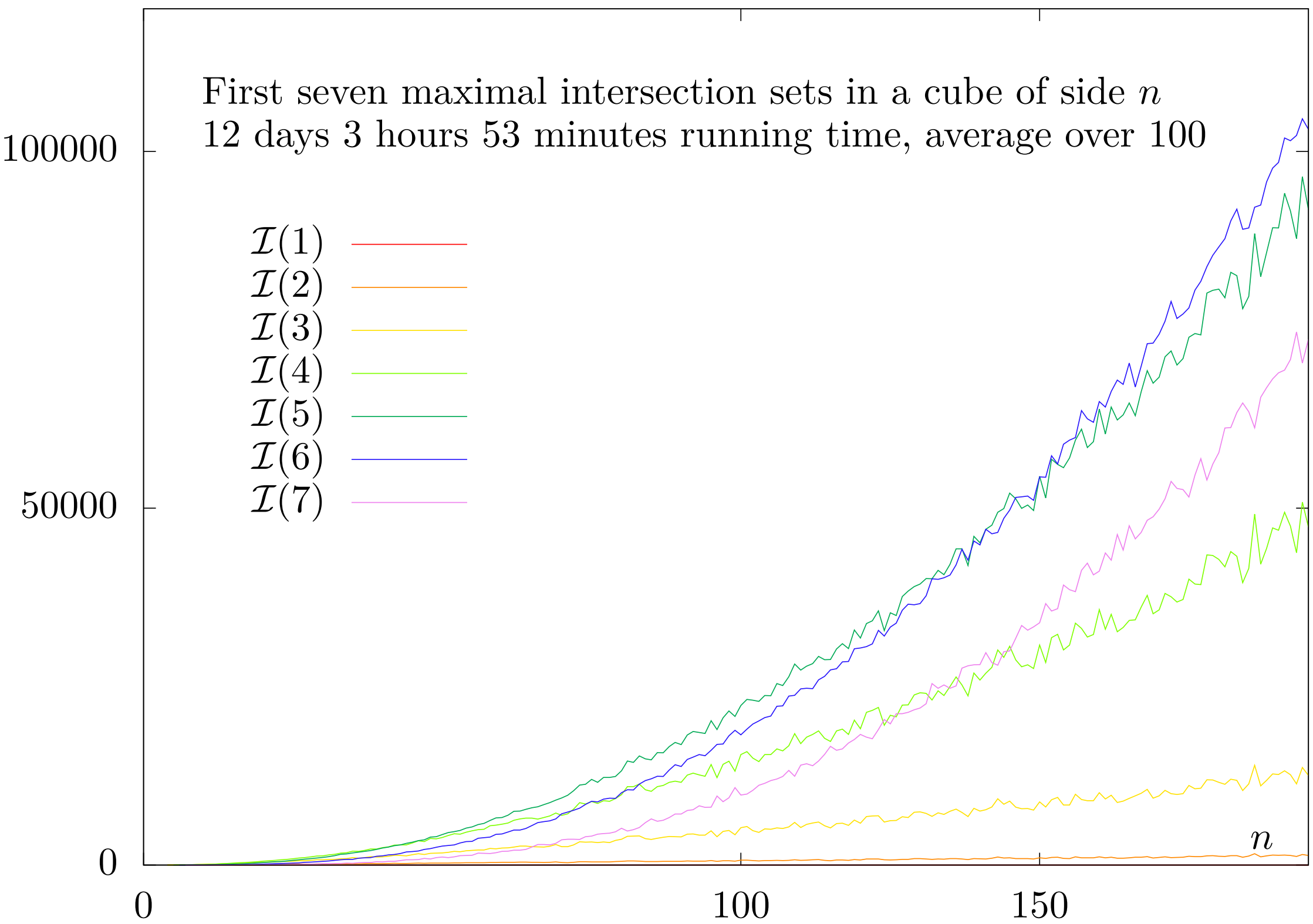

Figure 71: The first seven maximal intersection sets in a square and a cube: ultimate configuration before the connexion. The increment for $n$ is $\Delta n = 1$.

# 33 New attempts on the second intersections

The brutal control on the higher intersections developed in subsection 31.6 fails to provide a useful result for any of the three challenging questions raised in section 30. However, the first estimate concerning the pivotal bonds looked promising. Indeed, in the inequality (31.9) of subsection 31.1, we proved that the cardinality of the pivotal bonds is at most $O(\sqrt{n^d \ln n})$ with very high probability. The failure of the brutal control is due to the rapid degradation of the power of $n$ when it comes to the control of the higher intersections. We should not lose hope at this point, because there is still a lot of room for potential improvement in the brutal control of subsection 31.6. In this section, we present some interesting possibilities. Yet we warn the reader that, despite years of effort, none of these attempts was successful in the end.

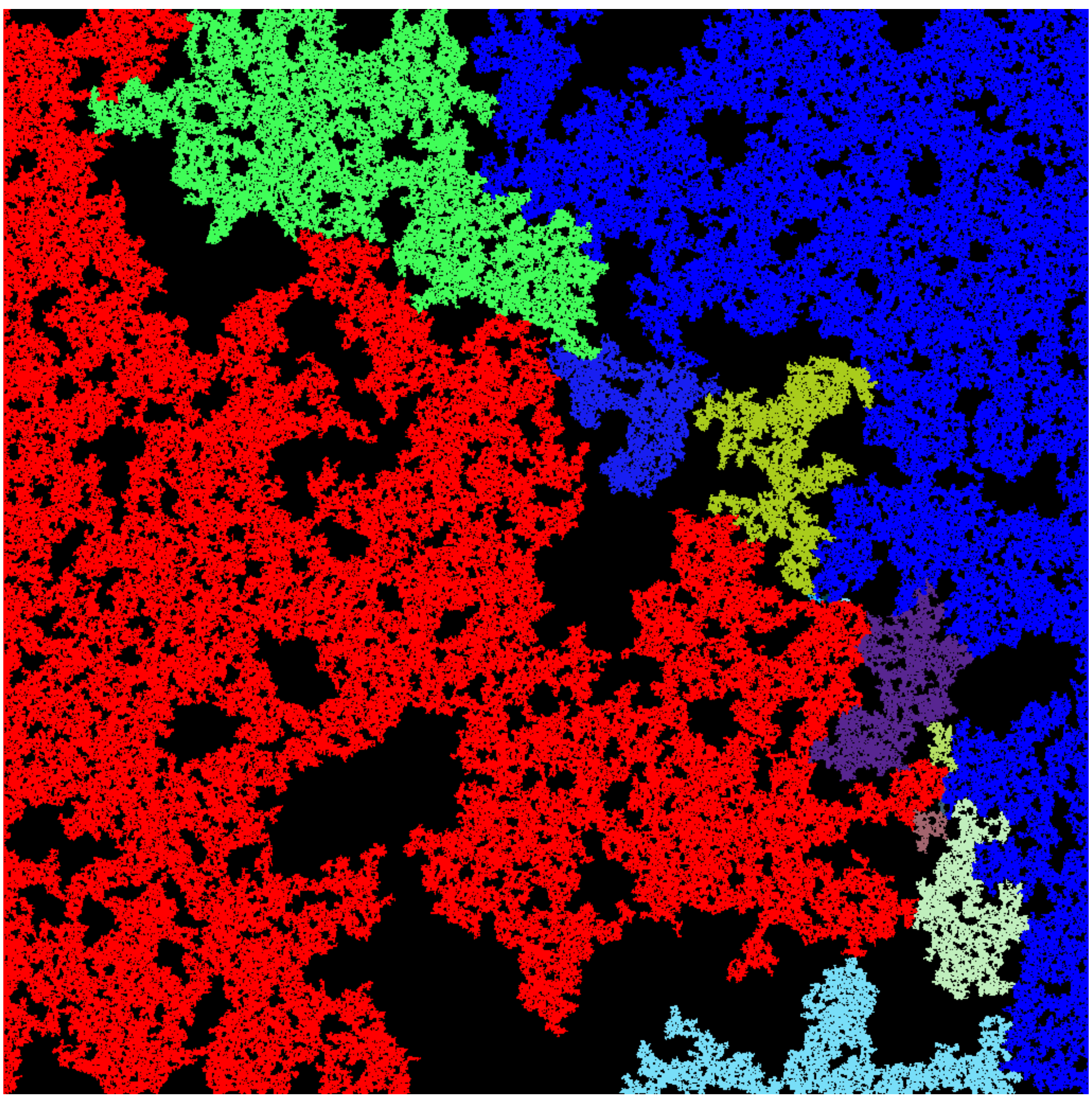

Figure 72: Bond percolation, p=0.5, 1024x1024 box. Clusters touching the left side are in red, clusters touching the right side are in blue. The remaining colored clusters touch simultaneously the red and the blue clusters. The second intersections are the bonds between these clusters and the red or blue clusters.

## 33.1 Incorporation of the volumes

Ultimately, the brutal control rests on the inequality (31.66) of proposition 31.3. In order to obtain a meaningful control inside the box $\Lambda(n)$, the sequence $i_1,\dots,i_t$ has to be adjusted so that

$$|\mathbb{E}^d(\Lambda(n))|^{i_1+\cdots+i_{t-1}}\exp\Big(-\frac{i_t^2}{|\mathbb{E}^d(\Lambda(n))|}\Big) \tag{33.1}$$

goes to 0 as $n$ goes to $\infty$ (to simplify the discussion, the additional constants appearing in the inequality have been removed). If we suppose in addition that the sequence $i_1,\dots,i_t$ is non-decreasing, then the quantity (33.1) is less or equal than

$$|\mathbb{E}^d(\Lambda(n))|^{(t-1)i_{t-1}}\exp\Big(-\frac{i_t^2}{|\mathbb{E}^d(\Lambda(n))|}\Big)$$

and this last quantity goes to 0 as soon as

$$i_t^2 \,\geq\, |\mathbb{E}^d(\Lambda(n))|2(t-1)i_{t-1}\ln\big(|\mathbb{E}^d(\Lambda(n))|\big)\,. \tag{33.2}$$

Moreover, this inequality has to be satisfied for all the indices $t\geq 1$. So, up to the constants that have been removed for the sake of the discussion, the best result that we can achieve this way corresponds essentially to the sequence chosen in (31.75), given by

$$\forall t\geq 1\qquad i_t \,=\, 2^t t!\big(n^d\ln n\big)^{1-2^{-t}}\,. \tag{33.3}$$

Solving the equality case of (33.2) gives a smaller sequence, but with the same power of $n$. The presence of the factorial $t!$ is not too worrisome, because, for the time being, we try only to control the pivotal bonds of very small order. The main trouble comes rather from the presence of the factor $|\mathbb{E}^d(\Lambda(n))|$ in (33.2) (the first one, not the one in the logarithm). Because of this factor, the best value we can pick up for $i_t$ is essentially the geometric mean between $n^d$ and $i_{t-1}$, and the resulting sequence (33.3) converges far too fast towards $n^d$. Indeed, the exponent of $n$ in the sequence $i_t$ converges towards $d$ at geometric speed! This is catastrophic, and this type of estimate is bound to fail, because the limiting exponent will be $d$ as $t$ goes to $\infty$. For our strategy to be successful, we must always have an exponent smaller than the isoperimetric exponent of the lattice, that is $d-1$. So we should try to improve this factor $|\mathbb{E}^d(\Lambda(n))|$. This nasty factor, the cardinality of the number of bonds in $\Lambda(n)$, appeared in the inequality (31.72), because we bounded crudely by $|\mathbb{E}^d(\Lambda(n))|$ the number of bonds that were explored to build the aggregate $\mathcal{A}(2t)$. Let us try to enhance this factor. We start by decomposing more carefully $S\big(\mathbb{E}^d(\mathcal{A}(2t),D)\big)$ as the successive sums corresponding to each round of the intertwined exploration algorithm. Suppose that, during the round $t$, the number of bonds explored by the active explorer is $V_t$. Suppose also that we have at our disposal a control on the probability $P(|V_t|\geq v)$. We could then replace the term (33.1) by

$$|\mathbb{E}^d(\Lambda(n))|^{i_1+\cdots+i_{t-1}}\exp\Big(-\frac{i_t^2}{v}\Big)\,+\,P(|V_t|\geq v)\,,$$

and we would optimize the value of $v$ so that the two terms are of the same order. This could lead to a serious improvement of the sequence obtained with the brutal control. Unfortunately, after years of unfruitful attempts, this direction did not lead to anything satisfactory, because we could never obtain a meaningful control on $P(|V_t|\geq v)$. In fact, this question is related to the tail distribution of the cluster size at the critical point, for which no quantitative estimate of any kind is available (see the disappointing results of section 15).

It should also be noted that the trouble really started with the second intersections. Indeed, in subsection 31.2, we obtained an upper bound of the form (see formula (31.30)):

$$P\big(|\mathcal{I}(2)|\geq c(d,p)(n^d\ln n)^{3/4}, W\not\longleftrightarrow W'\big)\ \leq\ \frac{c(d)}{n^{2d}}\,,$$

where $c(d,p)$, $c(d)$ are positive constants depending only on the dimension $d$ and the parameter $p$, as the parentheses indicate. Thus we could only prove that $|\mathcal{I}(2)|$ is $O\big((n^d\ln n)^{3/4}\big)$, which is already very bad. Let us discuss what we might obtain if we use this estimate to answer the challenging question of subsection 30.3, namely the disconnection in a cuboid. In corollary 30.5, we showed that, whenever the disconnection event $\mathcal{D}(n)=\big\{\,L\not\longleftrightarrow R\quad\text{in }\Pi(n)\,\big\}$ occurs, there exists with high probability an intersection set whose cardinality is larger than $n^{d-1}\beta(n)/\ln n$, or, if we forget the factor $\beta(n)$ on which we don't have any control beyond the mere fact that it goes to $\infty$, larger than $n^{d-1}/\ln n$. Let us compare the exponents of $n$ in these two antagonistic bounds. We have that

$$d-1\ >\ \frac{3d}{4}\qquad\Longleftrightarrow\qquad d\ >\ 4\,,$$

so the control of the second intersections could possibly work only in dimensions $d\geq 5$, which is quite disappointing. What is worse, with the further degradation of the exponent of $n$ for the higher intersection sets (see inequality (31.80)), the argument inevitably breaks down. Therefore, we need to concentrate our efforts on improving the control of the second intersections. As in the case of the higher intersections, the cause of the degradation of the exponent of $n$ for the second intersections is the presence of the factor $|D|^{i_1}$ in inequality (31.26) (which is responsible for the factor $n^{di_1}$ in (31.27)). The origin of this factor can be traced back to the control of the sum appearing in (31.20). Indeed, to bound the probability inside the sum

$$\sum_{T_1\subset D,|T_1|\leq i_1} P\big(|\mathcal{I}(2)|\geq i_2, |\mathcal{I}(1)|\leq i_1, \mathcal{T}(1)=T_1, W\not\longleftrightarrow W'\big)$$

we essentially forgot the event $\{\,\mathcal{T}(1)=T_1\,\}$ in order to obtain an estimate which is uniform over all the terms in the sum (see inequality (31.21)). As a consequence, we had then to sum this uniform estimate over all possible taboo sets $T_1$ in $D$ having cardinality less or equal than $i_1$, and this created the nasty factor $|D|^{i_1}$.

## 33.2 Conditioning on the taboo set $\mathcal{T}(1)$

To get a smarter upper bound, we should condition on the event $\{\, \mathcal{T}(1) = T_1 \,\}$ instead of merely forgetting it. We start with the equality

$$P\big(|\mathcal{I}(2)| \geq i_2, |\mathcal{I}(1)| \leq i_1, \mathcal{T}(1) = T_1, W \not\leftrightarrow W'\big) = \\ P\big(|\mathcal{I}(2)| \geq i_2, |\mathcal{I}(1)| \leq i_1,\, W \not\leftrightarrow W' \,\big|\, \mathcal{T}(1) = T_1\big) P\big(\mathcal{T}(1) = T_1\big)\,. \tag{33.4}$$

The next big problem is to estimate the conditional probability

$$P\big(|\mathcal{I}(2)| \geq i_2, |\mathcal{I}(1)| \leq i_1,\, W \not\leftrightarrow W' \,\big|\, \mathcal{T}(1) = T_1\big)\,.$$

If we obtain a valuable upper bound on this probability, then we could substitute it back in (33.4) and in (31.20). We would then perform the summation over the possible taboo sets $T_1$. This time, the presence of the factor $P\big(\mathcal{T}(1) = T_1\big)$ would save us from the appearance of the nasty factor $|D|^{i_1}$. Of course, a new challenging problem has arisen. When conditioning on $\{\, \mathcal{T}(1) = T_1 \,\}$, the distribution of the remaining bonds is not any more a Bernoulli product, because we have incorporated the unwanted information that $T_1$ is the set of the pivotal bonds for the connection between $W$ and $W'$. However, there are still some further steps of the previous brutal control of subsection 31.6 that can be followed. In fact, once we condition on the taboo set $\mathcal{T}_1$, we might as well condition on the aggregate $\mathcal{A}(0)$. So we write

$$P\big(|\mathcal{I}(2)| \geq i_2, |\mathcal{I}(1)| \leq i_1, W \not\leftrightarrow W'\big) = \\ \sum_{A_0, T_1 \subset D, |T_1| \leq i_1} P\big(|\mathcal{I}(2)| \geq i_2, |\mathcal{I}(1)| \leq i_1, \mathcal{T}(1) = T_1, \mathcal{A}(0) = A_0, W \not\leftrightarrow W'\big)\,. \tag{33.5}$$

Let us fix two subsets $A_0$, $T_1$ of $D$. We have

$$P\big(|\mathcal{I}(2)| \geq i_2, |\mathcal{I}(1)| \leq i_1, \mathcal{T}(1) = T_1, \mathcal{A}(0) = A_0, W \not\leftrightarrow W'\big) = \\ P\left(|\mathcal{I}(2)| \geq i_2,\, |\mathcal{I}(1)| \leq i_1,\, W \not\leftrightarrow W' \,\middle|\, \begin{matrix} \mathcal{T}(1) = T_1 \\ \mathcal{A}(0) = A_0 \end{matrix}\right) P\begin{pmatrix} \mathcal{T}(1) = T_1 \\ \mathcal{A}(0) = A_0 \end{pmatrix}\,.$$

We try now to control the conditional probability. By definition (see for instance (31.14)), we have

$$\mathcal{I}(1) \cup \mathcal{I}(2) \,=\, \Delta_D \mathcal{A}(1) \cap \Delta_D \mathcal{A}(2)\,. \tag{33.6}$$

On the event $\{\, W \not\leftrightarrow W' \,\}$, the two aggregates $\mathcal{A}(1)$, $\mathcal{A}(2)$ are disjoint and

$$\Delta_D \mathcal{A}(1) \cap \Delta_D \mathcal{A}(2) \,=\, \mathbb{E}^d\big(\mathcal{A}(1), D\big) \cap \mathbb{E}^d\big(\mathcal{A}(2), D\big)\,. \tag{33.7}$$

Moreover all the bonds belonging to the above intersection are closed, so that

$$S\Big(\mathbb{E}^d\big(\mathcal{A}(1), D\big) \cap \mathbb{E}^d\big(\mathcal{A}(2), D\big)\Big) \,=\, -\frac{1}{1-p}\Big|\Delta_D \mathcal{A}(1) \cap \Delta_D \mathcal{A}(2)\Big|\,. \tag{33.8}$$

Using the additivity of the map $S$, we conclude from the three identities (33.6), (33.7) and (33.8) that, if $|\mathcal{I}(2)| > i_2$ and $|\mathcal{I}(1)| \leq i_1$, then at least one among the three quantities

$$S\big(\mathbb{E}^d(\mathcal{A}(1), D)\big), \quad S\big(\mathbb{E}^d(\mathcal{A}(2), D)\big), \quad S\big(\mathbb{E}^d(\mathcal{A}(1) \cup \mathcal{A}(2), D)\big)\,, \tag{33.9}$$

must have an absolute value larger than

$$j \;=\; \frac{i_2}{3(1-p)}\,.$$

Now the new problem is that we have to control the conditional probability of these events. The easy terms are the first and the third, the second one is very difficult to apprehend. So, let us start with the first one. In fact, there is no way to control directly the conditional probability

$$P\left(\begin{matrix}|\mathcal{I}(2)| \geq i_2, |\mathcal{I}(1)| \leq i_1,\, W \not\longleftrightarrow W' \\ \big|S\big(\mathbb{E}^d(\mathcal{A}(1), D)\big)\big| \geq j\end{matrix}\,\middle|\, \begin{matrix}\mathcal{T}(1) = T_1 \\ \mathcal{A}(0) = A_0\end{matrix}\right)$$

so we re-inject it in the sum (33.5) and we do a deconditioning:

$$\begin{aligned}
\sum_{\substack{A_0, T_1 \subset D \\ |T_1| \leq i_1}} P\left(\begin{matrix}|\mathcal{I}(2)| \geq i_2, |\mathcal{I}(1)| \leq i_1,\, W \not\longleftrightarrow W' \\ \big|S\big(\mathbb{E}^d(\mathcal{A}(1), D)\big)\big| \geq j\end{matrix}\,\middle|\, \begin{matrix}\mathcal{T}(1) = T_1 \\ \mathcal{A}(0) = A_0\end{matrix}\right) P\left(\begin{matrix}\mathcal{T}(1) = T_1 \\ \mathcal{A}(0) = A_0\end{matrix}\right) \\
= P\left(\begin{matrix}|\mathcal{I}(2)| \geq i_2, |\mathcal{I}(1)| \leq i_1,\, W \not\longleftrightarrow W' \\ \big|S\big(\mathbb{E}^d(\mathcal{A}(1), D)\big)\big| \geq j\end{matrix}\right) \\
\leq P\Big(\big|S\big(\mathbb{E}^d(\mathcal{A}(1), D)\big)\big| \geq j\Big)\,.
\end{aligned} \tag{33.10}$$

The set $\mathcal{A}(1)$ is simply obtained by running an exploration algorithm starting from $W'$. Therefore the set of bonds $\mathbb{E}^d(\mathcal{A}(1), D)$ is the output of a genuine exploration algorithm in $D$, hence we have the usual control on $S\big(\mathbb{E}^d(\mathcal{A}(1), D)\big)$. We examine now the third term in (33.9). We have shown in proposition 31.2 that $\mathcal{A}(2) = \mathcal{B}\big(W, D \setminus \mathcal{T}(1), 1\big)$. The clusters of the vertices of $\mathcal{T}(1)$ are already included in $\mathcal{A}(1)$. So, the union $\mathcal{A}(1) \cup \mathcal{A}(2)$ will contain the whole second shell $\text{Shell}\,(W, D, 1)$ and we have

$$\mathcal{A}(1) \cup \mathcal{A}(2) \;=\; \mathcal{B}\big(W', D, 0\big) \cup \mathcal{B}\big(W, D, 1\big)\,.$$

As a consequence, the set of bonds $\mathbb{E}^d(\mathcal{A}(1) \cup \mathcal{A}(2), D)$ is the output of a genuine exploration algorithm in $D$. To be precise, it is the union of the outputs of the two taboo growth algorithms (defined in subsection 31.5)

$$\text{Taboo_Bond_Growth}\,(W', D, \varnothing, 0)\,, \quad \text{Taboo_Bond_Growth}\,(W, D, \varnothing, 1)\,,$$

and the taboo growth algorithm working with a fixed deterministic taboo set is a genuine exploration algorithm (honestly, when the taboo set is empty, we can use directly a genuine bond growth algorithm, however we did not take the time to define it separately in the foregoing). So we have the usual control on

$$S\big(\mathbb{E}^d(\mathcal{A}(1)\cup\mathcal{A}(2),D)\big)\,.$$

We conclude that the third term in (33.9) can be controlled in a way similar to the first one, through a deconditioning as in formula (33.10). In fact, if we had known in advance that we would do a deconditioning, we could have applied the control from the start. The conditioning is relevant only for the second term in (33.9). We come now to the heart of the matter, namely the control of this second term $S\big(\mathbb{E}^d(\mathcal{A}(2),D)\big)$ and of the conditional probability

$$P\left(\begin{matrix}|\mathcal{I}(2)|\geq i_2, |\mathcal{I}(1)|\leq i_1,\, W\not\longleftrightarrow W' \\ \big|S\big(\mathbb{E}^d(\mathcal{A}(2),D)\big)\big|\geq j\end{matrix}\,\middle|\,\begin{matrix}\mathcal{T}(1)=T_1\\ \mathcal{A}(0)=A_0\end{matrix}\right)\,.$$

The problem is that the construction of the set of bonds $\mathbb{E}^d(\mathcal{A}(2),D)$ is much more complex than for $\mathbb{E}^d(\mathcal{A}(1),D)$ or $\mathbb{E}^d(\mathcal{A}(1)\cup\mathcal{A}(2),D)$, and it cannot be realized as the output of a genuine exploration algorithm. Let us decompose its construction. The first explorer begins by building $\mathcal{A}(0)$ and for this, it explores $\mathbb{E}^d(\mathcal{A}(0),D)$. For his second round of exploration, he will restart from the vertices in $\partial_D^{out}\mathcal{A}(0)\setminus\mathcal{T}(1)$ and he will return the set

$$C\big(\partial_D^{out}\mathcal{A}(0)\setminus\mathcal{T}(1),D\setminus\mathcal{A}(0)\big)\;=\;C\big(\partial_D^{out}\mathcal{A}(0)\setminus\mathcal{T}(1),D\big)\,.$$

The two sets above are equal because all the bonds in $\Delta_D\mathcal{A}(0)$ are closed. So, during his second round, the first explorer will explore the set of bonds

$$\mathbb{E}^d\Big(C\big(\partial_D^{out}\mathcal{A}(0)\setminus\mathcal{T}(1),D\big),D\setminus\mathcal{A}(0)\Big)\,.$$

Indeed, all the bonds emanating from $\mathcal{A}(0)$ have already been explored during the first round, hence the exploration looks only at the bonds in $D\setminus\mathcal{A}(0)$. Moreover, the first explorer is prevented from entering into the territory $\mathcal{A}(1)$ already visited by the second explorer, and this is exactly the goal achieved with the taboo set $\mathcal{T}(1)$. The trouble comes precisely from this taboo set $\mathcal{T}(1)$. It is a random set, which carries with him some information on the configuration of the bonds which do not have an endpoint in $\mathcal{A}(0)$, and the distribution of the bonds explored during the second round is not any more i.i.d.. At least, using the additivity of $S$, we can write

$$S\big(\mathbb{E}^d(\mathcal{A}(2),D)\big)\;=\;S\big(\mathbb{E}^d(\mathcal{A}(0),D)\big)+S\Big(\mathbb{E}^d\Big(C\big(\partial_D^{out}\mathcal{A}(0)\setminus\mathcal{T}(1),D\big),D\setminus\mathcal{A}(0)\Big)\Big).\tag{33.11}$$

Since $\mathbb{E}^d(\mathcal{A}(0),D)$ is the output of a genuine exploration algorithm, then the term $S\big(\mathbb{E}^d(\mathcal{A}(0),D)\big)$ can be taken care of through a deconditioning procedure as the term $S\big(\mathbb{E}^d(\mathcal{A}(1),D)\big)$. So we are left with the second term in (33.11). Conditioning on $\mathcal{A}(0)$ and $\mathcal{T}(1)$, we have to consider the conditional probability

$$P\left(\begin{matrix}|\mathcal{I}(2)|\geq i_2, |\mathcal{I}(1)|\leq i_1,\, W\not\longleftrightarrow W' \\ S\Big(\mathbb{E}^d\Big(C\big(\partial_D^{out}A_0\setminus T_1,D\big),D\setminus A_0\Big)\Big)\geq \frac{j}{2}\end{matrix}\,\middle|\,\begin{matrix}\mathcal{T}(1)=T_1\\ \mathcal{A}(0)=A_0\end{matrix}\right)\,.\tag{33.12}$$

Notice that we have replaced $j$ by $j/2$, because we imagine that the term associated to $S\big(\mathbb{E}^d(\mathcal{A}(0),D)\big)$ has been controlled (we refrain from writing the full details here, our goal is to describe some possible strategies). Let us explain again what we are trying to achieve. We suppose that during his first round, the first explorer has built the aggregate $A_0$. In the absence of a taboo set, the second round of exploration of the first explorer would restart from $\partial_D^{out} A_0$ and the explorer would explore $C(\partial_D^{out} A_0, D)$. Now, the first explorer is told that his taboo set is $T_1$, hence he will only explore the clusters of the vertices of $\partial_D^{out} A_0 \setminus T_1$. As a consequence, when he decides to start the exploration from a vertex of $\partial_D^{out} A_0$, he knows that this vertex will not lead to a vertex of $W'$. With this information in hand, we wish to control

$$S\Big(\mathbb{E}^d\Big(C\big(\partial_D^{out} A_0 \setminus T_1, D\big), D \setminus A_0\Big)\Big)\,. \tag{33.13}$$

At this point, we can certainly forget the events $\{\,|\mathcal{I}(2)| \geq i_2, |\mathcal{I}(1)| \leq i_1\,\}$ (this might create a polynomial factor in the end, and we seek to improve the exponential bounds), as well as the event $\{\,W \not\longleftrightarrow W'\,\}$ (which is realized by the conditioning $\mathcal{A}(0) = A_0$), so we bound the probability (33.12) by

$$P\left(\left|S\Big(\mathbb{E}^d\Big(C\big(\partial_D^{out} A_0 \setminus T_1, D\big), D \setminus A_0\Big)\Big)\right| \geq \frac{j}{2}\,\middle|\, \begin{matrix}\mathcal{T}(1) = T_1\\ \mathcal{A}(0) = A_0\end{matrix}\right)\,. \tag{33.14}$$

Unfortunately, conditioning on the event $\{\,\mathcal{T}(1) = T_1,\, \mathcal{A}(0) = A_0\,\}$ changes radically the distribution of the bonds in $\mathbb{E}^d(D \setminus A_0)$. Because the information conveyed by the conditioning is of negative nature (it tells us that some connections do not exist in the domain that is to be explored), the value of the quantity (33.13) might be much smaller than for a genuine exploration algorithm. So, how on earth can we control it?

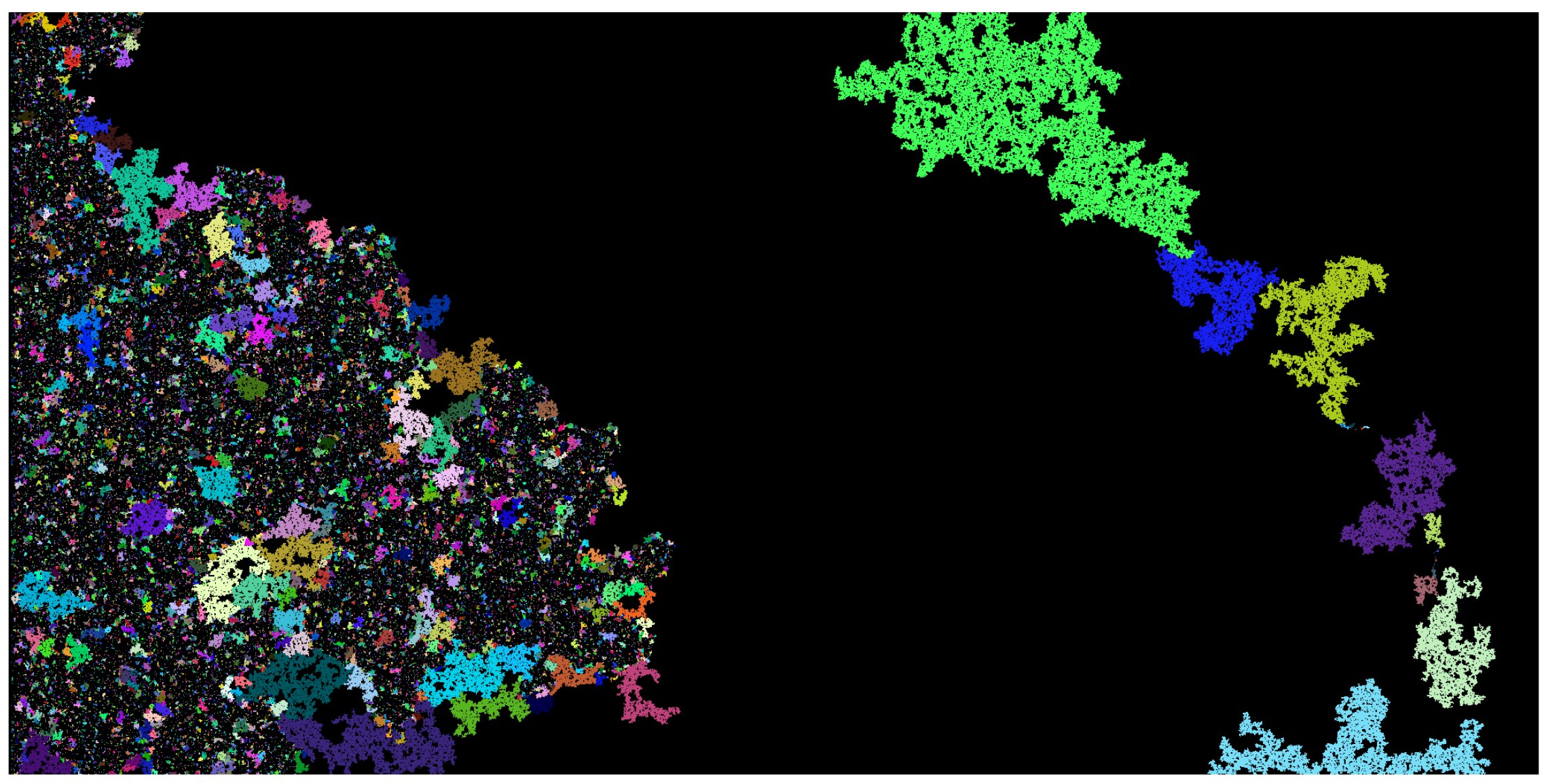

Figure 73: Further analysis of the bond percolation configuration of figure 72. Left: the clusters touching the red clusters of figure 72 but not the blue ones. Right: the clusters touching both the red and the blue clusters of figure 72.

## 33.3 Hiding some finite clusters

We present here a strategy which is applicable only when $D = \Lambda(n)$ and one of the two sets $W$, $W'$, say for instance $W$, is equal to $\partial^{\,in}\Lambda(n)$ (so this would work for the two first questions of section 30, but not for the third one, the disconnection in a cuboid). Imagine that we use a standard exploration algorithm to build

$$\mathbb{E}^d\Big(C\big(\partial_D^{out}A_0\setminus T_1, D\big), D\setminus A_0\Big)\,.$$

Such an algorithm would pick up a vertex in $\partial_D^{out}A_0\setminus T_1$, and then explore its cluster in $D\setminus A_0$ (that is, he will explore all the bonds of $\mathbb{E}^d(D\setminus A_0)$ that have one endpoint in this cluster). Once this cluster has been examined, the algorithm would pick up another vertex in $\partial_D^{out}A_0\setminus T_1$ that has not been visited previously, and he will examine its cluster. This process will go on until all the vertices of $\partial_D^{out}A_0\setminus T_1$ are exhausted. If we were to perform the same task without a taboo set, that is, if we build $\mathbb{E}^d\big(C\big(\partial_D^{out}A_0, D\big), D\setminus A_0\big)$, then we can use also a standard exploration algorithm, but this time, it is a genuine exploration algorithm, and we have the usual control on

$$S\Big(\mathbb{E}^d\Big(C\big(\partial_D^{out}A_0, D\big), D\setminus A_0\Big)\Big)\,. \tag{33.15}$$

The difference between the two algorithms is merely that the vertices of $T_1$ are not accepted as starting vertices when we pick up a new vertex in order to explore its cluster. So we shift our point of view and we try to understand the worst damage that a taboo set can inflict on the quantity (33.15). To this end, we introduce a family of random variables indexed by the subsets of $D$ by setting

$$\forall T\subset D\qquad \Phi(T)\,=\,S\Big(\mathbb{E}^d\Big(C\big(\partial_D^{out}A_0\setminus T, D\big), D\setminus A_0\Big)\Big)\,.$$

Until now, we have been trying to control directly the conditional tail distribution of $\Phi(T_1)$:

$$P\big(|\Phi(T_1)|\geq i\,\big|\,\mathcal{T}(1)=T_1,\mathcal{A}(0)=A_0\big)\,.$$

It might be sensible to try first to control the conditional expectation

$$E\big(|\Phi(T_1)|\,\big|\,\mathcal{T}(1)=T_1,\mathcal{A}(0)=A_0\big)\,.$$

Using the control on the pivotal bonds for the connection between $W$ and $W'$, we know that, with high probability, the cardinality of $\mathcal{I}(1)$ is $O(\sqrt{n^d\ln n})$, and the same is true for $\mathcal{T}(1)$. Thus, in the conditioning, we need only to consider taboo sets $T_1$ such that

$$|T_1|\,=\,O(\sqrt{n^d\ln n})\,.$$

For such a set $T_1$, we have the obvious bound

$$|\Phi(T_1)|\,\leq\,\sup\,\Big\{\,|\Phi(T)|:T\subset D,\,|T|\,=\,O(\sqrt{n^d\ln n})\,\Big\}\,. \tag{33.16}$$

The good point is that the right-hand quantity in (33.16) does not depend any more on $T_1$. Taking the conditional expectation with respect to the event

$\{\mathcal{T}(1)=T_1,\,\mathcal{A}(0)=A_0\,\}$, we get

$$E\Big(|\Phi(T_1)|\,\Big|\,\mathcal{T}(1)=T_1,\mathcal{A}(0)=A_0\Big)\,\leq\\ E\Big(\sup\,\Big\{\,|\Phi(T)|:T\subset D,\,|T|\,=\,O(\sqrt{n^d\ln n})\,\Big\}\,\Big|\,\mathcal{A}(0)=A_0\Big)\,. \quad (33.17)$$

So we should now control the supremum appearing inside the expectation. This is still far from obvious, and it is not clear whether this upper bound would be enough for our purpose in the end. Let us consider more closely this supremum and the role played by the taboo set. The quantity $\Phi(T)$ is a random variable, which is a martingale-like sum over the bonds explored during the taboo exploration of the clusters in $A_0\setminus T$ intersecting the set $\partial_D^{out}A_0\setminus T$. Recall that the exploration of all the clusters intersecting $\partial_D^{out}A_0$ can be done with the help of a genuine exploration algorithm driven by a sequence $(X_k)_{k\geq 1}$ of i.i.d. Bernoulli random variables. An effect of the taboo set intervening in the quantity $\Phi(T)$ is that several intervals of the sequence $(X_k)_{k\geq 1}$ will be completely discarded. An optimal strategy for choosing the taboo set consists in preventing the exploration process from entering into the clusters which create compensations and stabilize the value of $\Phi(T)$ around its mean. For instance, if we wish the value of $\Phi(T)$ to be very negative, then we would try to forbid the entrance into the clusters having a lot of open bonds and comparatively few closed bonds in their boundaries. With regard to our problem, for which the taboo set was taken to be the endpoints of the pivotal bonds, it is to be expected that the vertices of $\mathcal{T}(1)$ correspond to clusters having a positive contribution to $\Phi(T_1)$, because a vertex of $\partial_D^{out}A_0$ is in $\mathcal{T}(1)$ if and only if it is connected to $W'$. More precisely, let $x$ be a fixed vertex in $D$ and let us examine the conditional expectation

$$E\Big(S\big(\mathbb{E}^d(C(x,D\setminus\mathcal{A}(0)))\big)\,1_{x\in\mathcal{T}(1)}\,\Big|\,\mathcal{A}(0)=A_0\Big)\,. \quad (33.18)$$

The exploration building $C(x,D\setminus\mathcal{A}(0))$ takes place in $D\setminus\mathcal{A}(0)$, so that no bond of $\Delta_D\mathcal{A}(0)$ is incorporated in the sum $S\big(\mathbb{E}^d(C(x,D\setminus\mathcal{A}(0)))\big)$. Conditionally on the event $\mathcal{A}(0)=A_0$, the bonds in $\mathbb{E}^d(D\setminus A_0)$ are Bernoulli i.i.d., so the only additional information in (33.18) is that the vertex $x$ belongs to $\mathcal{T}(1)$, meaning that the event

$$x\longleftrightarrow W'\quad\text{in}\quad\mathbb{E}^d(D\setminus A_0) \quad (33.19)$$

occurs, which is of positive nature. By the FKG inequality, we have

$$\forall x\in D\qquad E\Big(S\big(\mathbb{E}^d(C(x,D\setminus\mathcal{A}(0)))\big)\,1_{x\in\mathcal{T}(1)}\,\Big|\,\mathcal{A}(0)=A_0\Big)\,\geq\,0\,.$$

Yet the fact that the expectation (33.19) is non-negative does not mean that the value of $\Phi(T_1)$ will increase dramatically. In a nutshell, what we try to achieve in the supremum problem (33.16) is to hide from the exploration process the clusters which contribute positively and substantially to $\Phi(T)$. A fundamental question is to get some valuable control on the supremum in (33.17). We consider next a toy version of this problem.

## 33.4 A toy version of the supremum problem.

We consider a sequence of i.i.d. Bernoulli random variables $(Y_k)_{k\geq 1}$ such that

$$\forall k\geq 1 \qquad P\Big(Y_k=-\frac{1}{1-p}\Big)\,=\,1-p\,,\quad P\Big(Y_k=\frac{1}{p}\Big)\,=\,p\,, \tag{33.20}$$

and we form the associated partial sums $(S_n)_{n\geq 1}$, defined as

$$\forall n\geq 1 \qquad S_n=Y_1+\cdots+Y_n\,.$$

We have the classical control given by Hoeffding's inequality:

$$\forall u>0\quad \forall n\geq 1 \qquad P\big(|S_n|\geq u\big)\,\leq\,2\exp\Big(-\frac{2p^2(1-p)^2u^2}{n}\Big)\,.$$

A specific proof for this particular case, with the Bernoulli variables $(Y_k)_{k\geq 1}$ defined in (33.20), can be found in proposition 8 of [29]. Let us fix $n\geq 1$ and let us take $u=c\sqrt{n\ln n}$. We obtain

$$\forall c>0 \qquad P\big(|S_n|\geq c\sqrt{n\ln n}\big)\,\leq\,2n^{-2c^2p^2(1-p)^2}\,.$$

Suppose now that we have the right to declare that some specific sub-intervals of $\{\,1,\dots,n\,\}$ are taboo. More precisely, the rule is that we can choose $r$ sub-intervals $I_1,\dots,I_r$ of $\{\,1,\dots,n\,\}$ and remove them when computing the sum $S_n$. So, instead of $S_n$, we compute

$$S(n,I_1,\dots,I_r)\,=\,\sum_{\substack{k\in\{\,1,\dots,n\,\}\\ k\notin I_1\cup\cdots\cup I_r}} Y_k\,.$$

The toy version of the supremum problem amounts to control

$$E\Big(\sup\,\Big\{\,\big|S(n,I_1,\dots,I_r)\big|:I_1,\dots,I_r\subset\{\,1,\dots,n\,\}\,\Big\}\Big)\,. \tag{33.21}$$

This supremum depends also on $r$. Obviously, if we take $r$ of order $n$, then the supremum will also be of order $n$. To have a situation comparable to our problem, we should take $r=\sqrt{n\ln n}$, or say $r=\sqrt{n}$ to make things simpler in a first approach. Unfortunately, it looks like that, for $r$ of order $\sqrt{n}$, the quantity (33.21) is typically of order $n^{3/4}$! Here is a possible strategy to obtain a relevant lower bound. We subdivide $\{\,1,\dots,n\,\}$ into $\sqrt{n}$ intervals of length $\sqrt{n}$. A partial sum over an interval of length $\sqrt{n}$ is typically of order $n^{1/4}$, but its sign can be either negative or positive. We remove those intervals for which the associated partial sum is negative, and we will typically obtain a sum of order $\sqrt{n}\times n^{1/4}=n^{3/4}$.

Thus, transforming the supremum problem (33.16) into a mean-field question like (33.21) will not help to improve our brutal upper bound $O((n\ln n)^{3d/4})$ on the second intersections. So we should use an approach that stays closer to the reality of the percolation model.

## 33.5 Can we get help from Talagrand's book?

Looking again at (33.17), the problem amounts to control the supremum of a complicated stochastic process. The major general reference for this kind of problem is Talagrand's book [126]. One could hope that the various techniques exposed by Talagrand could be helpful, and chiefly, the generic chaining, or some variant of it which would be more suited to the structure of the clusters. Very interesting estimates along these lines have been derived by Ding, Wirth and Xia in a series of papers on the random field Ising model [38, 39, 40]. They implemented the generic chaining scheme of Talagrand to derive estimates on quantities associated to greedy lattice animals, which are very reminiscent of our functional $S$. Thanks to these estimates, they could solve a long-standing difficult question on the random field Ising model. Regrettably, we did not manage to reuse these estimates profitably to make any progress around the conjecture $\theta(p_c,\mathbb{Z}^d)=0$. There is one marginal progress that came out from the general results exposed at the beginning of Talagrand's book. Indeed, by using the Gaussian bound on the distribution of the pivotal set, we can prove a variant of the two-arms inequality (8.15), which we present next. Theorem 33.1 and its proof are written for the site percolation model.

**Theorem 33.1.** *Let $p\in]0,1[$. There exists a constant $\kappa=\kappa(p,d)$ which depends only on the dimension $d$ and the percolation parameter $p$ such that*

$$\forall n\geq 1\qquad P_p\begin{pmatrix}\text{there exist two disjoint clusters } C_1,C_2 \text{ such}\\ \text{that } |C_1|\geq n, |C_2|\geq n, 0\in\partial^{\,out}C_1\cap\partial^{\,out}C_2\end{pmatrix}\;\leq\;\kappa\sqrt{\frac{\ln n}{n}}\,. \tag{33.22}$$

Tom Hutchcroft has proved a stronger inequality, which he called the two-ghost inequality [73]. Not only is his inequality more general than (33.22), but the bound is better because there is no logarithmic factor. The proof of Tom Hutchcroft is based on the mass-transport principle. This inequality plays an important role in the recent resolution of the Schramm locality conjecture by Philip Easo and Tom Hutchcroft [48]. Although the inequality (33.22) is weaker than the two-ghost inequality of Hutchcroft, we shall present its proof, because it relies on a different kind of argument, less elegant than Hutchcroft's proof, but which might prove useful for other purposes.

Of course, the result of theorem 33.1 is only interesting at the critical point $p_c$. Away from the critical point, the probability of the event in (33.22) decays much faster. As discussed in subsection 8.6, it follows from the argument of Aizenman, Kesten and Newman [4] for the uniqueness of the infinite cluster that

$$\exists\,\kappa>0\quad\forall n\geq 1\qquad P_p\big(\text{two-arms}(0,n)\big)\;\leq\;\frac{\kappa\ln n}{\sqrt{n}}\,, \tag{33.23}$$

where "two-arms$(0,n)$" is the event that two neighbours of 0 are connected to the boundary of the box $\Lambda(n)$ by two disjoint open clusters. The relationship between the two inequalities (33.22) and (33.23) is subtle. We stress the fact that the inequality (33.22) deals with an event defined for the infinite configuration in $\mathbb{Z}^d$, while the inequality (33.23) deals with an event defined in a finite box. As a consequence, it is not the case that the inequality (33.23) is implied by the inequality (33.22). Indeed, it is possible that the event two-arms$(0,n)$ occurs in $\Lambda(n)$, and that the two clusters realizing it are connected outside of $\Lambda(n)$, so that 0 is not in the outer boundaries of two distinct large clusters of the configuration in $\mathbb{Z}^d$. In the other direction, if the event appearing in (33.22) occurs, then the two clusters realizing it must have a diameter larger than $c\,n^{1/d}$ for some fixed constant $c>0$, and therefore the event two-arms$(0,cn^{1/d})$ occurs as well. Thus we have, for any $n\geq 1$,

$$P_p\begin{pmatrix}\text{there exist two disjoint clusters } C_1,C_2 \text{ such}\\ \text{that } |C_1|\geq n, |C_2|\geq n, 0\in\partial^{\,out}C_1\cap\partial^{\,out}C_2\end{pmatrix}\,\leq\,P_p\big(\text{two-arms}(0,c\,n^{1/d})\big)\,.$$

Therefore inequality (33.23) implies an inequality which is weaker than (33.22). The improvement from the inequality (33.23) to (33.22) is not spectacular, yet the methods of the proofs are quite different, and any little progress in the understanding of the percolation model at criticality should be duly appreciated.

*Proof.* Let $n$ be a positive integer. We denote by $\mathcal{C}_n$ the collection of the open clusters of the configuration in $\mathbb{Z}^d$ having cardinality larger than or equal to $n$ (so $\mathcal{C}_n$ contains also the infinite cluster $C_\infty$ if present), and we define the set

$$\mathcal{I}(n)\,=\,\bigcup_{\substack{C,D\in\mathcal{C}_n\\ C\neq D}}\partial^{\,out}C\cap\partial^{\,out}D\,.$$

In words, the set $\mathcal{I}(n)$ is the set of the sites which belong to the outer boundaries of two distinct clusters of $\mathcal{C}_n$. We define next the finite volume counterparts of $\mathcal{C}_n$ and $\mathcal{I}(n)$. For $\Lambda$ a finite box of $\mathbb{Z}^d$, we denote by $\mathcal{C}_n(\Lambda)$ the open clusters of the configuration restricted to $\Lambda$ whose cardinality is larger than or equal to $n$ and we define

$$\mathcal{I}(n,\Lambda)\,=\,\bigcup_{\substack{C,D\in\mathcal{C}_n(\Lambda)\\ C\neq D}}\partial^{\,out}C\cap\partial^{\,out}D\,.$$

**Lemma 33.2.** *Let $n,N$ be two integers such that $n\geq d$ and $N\geq 2n$. We have*

$$\mathcal{I}(n)\cap\Lambda(N)\,\subset\,\mathcal{I}\big(n,\Lambda(3N)\big)\,.$$

*Proof.* Let $x$ belong to $\mathcal{I}(n)\cap\Lambda(N)$. By definition, there exist $C,D$ two distinct clusters in $\mathcal{C}_n$ such that $\partial^{\,out}C$ and $\partial^{\,out}D$ contain $x$. If the cluster $C$ does not meet the internal boundary $\partial^{\,in}\Lambda(3N)$ of $\Lambda(3N)$, then $C$ is also a cluster of the configuration restricted to $\Lambda(3N)$. Otherwise, the trace of $C$ on $\Lambda(3N)$ splits up into a finite collection of clusters of the configuration restricted to $\Lambda(3N)$. One of these clusters, say $C'$, is such that

$$x\in\partial^{\,out}C'\,,\quad C'\cap\partial^{\,in}\Lambda(3N)\neq\varnothing\,.$$

This implies furthermore that

$$\text{diameter}\,C'\,\geq\,\inf\big\{\,|a-b|_\infty:a\in\partial^{\,out}\Lambda(N),\,b\in\partial^{\,in}\Lambda(3N)\,\big\}\,\geq\,N-1\,\geq\,n\,,$$

so that $|C'| \geq n$ and $C'$ belongs to $\mathcal{C}_n(\Lambda(3N))$. We apply the same argument to $D$, and we obtain that there exists a cluster $D'$ in $\mathcal{C}_n(\Lambda(3N))$ such that $D' \subset D$ and $x \in \partial^{\,out} D'$. Since $C' \subset C$ and $D' \subset D$, then the clusters $C'$ and $D'$ are distinct. This proves that the vertex $x$ is in $\mathcal{I}(n, \Lambda(3N))$. $\square$

Let $u \geq 0$ and suppose that $\big|\mathcal{I}(n) \cap \Lambda(N)\big| \geq u$. It follows from lemma 33.2 that $\big|\mathcal{I}(n, \Lambda(3N))\big| \geq u$. Let us denote by $(C_i, i \in I)$ the collection of the clusters of $\mathcal{C}_n(\Lambda(3N))$. We have therefore

$$\Big|\bigcup_{k,\ell \in I,\, k \neq \ell} \partial^{\,out} C_k \cap \partial^{\,out} C_\ell\Big| \,\geq\, u\,. \tag{33.24}$$

We use inequality (33.24) and we apply lemma 84.3: there exist two disjoint subsets $K, L$ of $I$ such that

$$\Big|\Big(\bigcup_{k \in K} \partial^{\,out} C_k\Big) \cap \Big(\bigcup_{\ell \in L} \partial^{\,out} C_\ell\Big)\Big| \,\geq\, \frac{u}{8d^2}\,.$$

We pick up a vertex of each cluster of the collection $(C_k, k \in K)$, thereby getting a subset $A$ of $\Lambda(3N)$ such that the collection $(C_k, k \in K)$ is equal to $(C(x), x \in A)$. We do the same to represent the collection $(C_\ell, \ell \in L)$ by means of a subset $B$ of $\Lambda(3N)$, which contains exactly one vertex of each cluster $C_\ell$ for $\ell \in L$. We have then

$$\bigcup_{k \in K} \partial^{\,out} C_k \,=\, \partial^{\,out} C(A, \Lambda(3N))\,, \qquad \bigcup_{\ell \in L} \partial^{\,out} C_\ell \,=\, \partial^{\,out} C(B, \Lambda(3N))\,.$$

Since the clusters $C(x)$, $x \in A \cup B$, all belong to $\mathcal{C}_n(\Lambda(3N))$ and are pairwise disjoint, then the sets $A, B$ must satisfy

$$\big(|A| + |B|\big)n \,\leq\, \big|\Lambda(3N)\big|\,. \tag{33.25}$$

As a consequence, the total number of possible choices for the sets $A, B$ is bounded by

$$\big|\Lambda(3N)\big|^{|A|+|B|} \,\leq\, \big|\Lambda(3N)\big|^{|\Lambda(3N)|/n}\,. \tag{33.26}$$

Remember that we supposed that the event $\big\{\, \big|\mathcal{I}(n) \cap \Lambda(N)\big| \geq u \,\big\}$ occurred. We deduce from the previous considerations that

$$P\big(\big|\mathcal{I}(n) \cap \Lambda(N)\big| \geq u\big) \leq P\begin{pmatrix} \exists A, B \subset \Lambda(3N) \quad \big(|A| + |B|\big)n \,\leq\, \big|\Lambda(3N)\big| \\ A \not\longleftrightarrow B \text{ in } \Lambda(3N) \\ \big|\partial^{\,out} C(A, \Lambda(3N)) \cap \partial^{\,out} C(B, \Lambda(3N))\big| \geq \frac{u}{8d^2} \end{pmatrix}$$
$$\leq \sum_{\substack{A, B \subset \Lambda(3N) \\ (|A|+|B|)n \,\leq\, |\Lambda(3N)|}} P\begin{pmatrix} A \not\longleftrightarrow B \text{ in } \Lambda(3N) \\ \big|\partial^{\,out} C(A, \Lambda(3N)) \cap \partial^{\,out} C(B, \Lambda(3N))\big| \,\geq\, \dfrac{u}{8d^2} \end{pmatrix}\,, \tag{33.27}$$

where we have used a standard union bound in the last step. We fix now two subsets $A, B$ of $\Lambda(3N)$ and we evaluate the probability appearing in the last sum. To this end, we use the first step of the construction of the intertwined exploration process starting from the sets $W = A$, $W' = B$ in the box $\Lambda(3N)$, and we apply the inequality (31.11). The $N$ in inequality (31.11) has not the same meaning as here, it is an upper bound on the cardinality of the sets $A, B$ and it is controlled here by inequality (33.25). We get

$$P\left(\begin{matrix} A \not\longleftrightarrow B \text{ in } \Lambda(3N) \\ \big|\partial^{out} C(A,\Lambda(3N)) \cap \partial^{out} C(B,\Lambda(3N))\big| \,\geq\, \dfrac{u}{8d^2} \end{matrix}\right) \leq \\ 6|\mathbb{E}^d(\Lambda(3N))|^{2|\Lambda(3N)|/n+3} \exp\Big(-\frac{p^2u^2}{288d^4|\mathbb{E}^d(\Lambda(3N))|}\Big)\,. \quad (33.28)$$

This bound is uniform with respect to the sets $A, B$ intervening in the sum (33.27). Substituting (33.28) in (33.27), and using the crude bound (33.26) to control the number of terms in the sum, we obtain

$$P\big(\big|\mathcal{I}(n) \cap \Lambda(N)\big| \geq u\big) \,\leq \quad (33.29)$$

$$\big|\Lambda(3N)\big|^{|\Lambda(3N)|/n} 6|\mathbb{E}^d(\Lambda(3N))|^{2|\Lambda(3N)|/n+3} \exp\Big(-\frac{p^2u^2}{288d^4|\mathbb{E}^d(\Lambda(3N))|}\Big)\,.$$

Recall that $N \geq 2n \geq 2d$ by hypothesis. We use the bounds

$$|\Lambda(3N)| \leq (3N+1)^d \leq (4N)^d\,,$$
$$|\mathbb{E}^d(\Lambda(3N))| \leq 2d|\Lambda(3N)| \leq 2d(4N)^d\,,$$

to deduce from (33.29) that

$$P\big(\big|\mathcal{I}(n) \cap \Lambda(N)\big| \geq u\big) \,\leq\, 6\big(2d(4N)^d\big)^{3(4N)^d/n+3} \exp\Big(-\frac{p^2u^2}{576d^5(4N)^d}\Big)\,. \quad (33.30)$$

Inequality (33.30) is nothing else than a Gaussian bound on the tail of the random variable $\big|\mathcal{I}(n) \cap \Lambda(N)\big|$. We bound the prefactor in (33.30) by

$$6\big(2d(4N)^d\big)^{3(4N)^d/n+3} \,\leq\, (48dN)^{4d(4N)^d/n}\,,$$

and we use lemma A.1 to conclude that

$$E\Big(\big|\mathcal{I}(n) \cap \Lambda(N)\big|\Big) \,\leq\, 2\sqrt{\frac{576d^5(4N)^d}{p^2}\ln\Big((48dN)^{4d(4N)^d/n}\Big)}\,. \quad (33.31)$$

By translation invariance, we have

$$E\Big(\big|\mathcal{I}(n) \cap \Lambda(N)\big|\Big) \,=\, \sum_{x \in \Lambda(N)} P\big(x \in \mathcal{I}(n)\big) \,=\, \big|\Lambda(N)\big| P\big(0 \in \mathcal{I}(n)\big)\,. \quad (33.32)$$

Combining (33.31) and (33.32), we conclude that

$$P\big(0 \in \mathcal{I}(n)\big) \,\leq\, c(d,p)\sqrt{\frac{\ln(48dN)}{n}}\,,$$

where $c(d,p)$ is a positive constant depending only on $d$ and $p$. Taking finally $N = 2n$, we obtain the inequality (33.22) of theorem 33.1. $\square$

## 33.6 The phase of the finite clusters

We present here another strategy to control the conditional probability (33.14). We still consider the situation where $D = \Lambda(n)$ and $W$ is equal to $\partial^{in}\Lambda(n)$. When $W$ is equal to $\partial^{in}\Lambda(n)$, the clusters explored during the second round of the first explorer are clusters of $\Lambda(n)$ which do not intersect $\partial^{in}\Lambda(n)$. These clusters are left unchanged if we look at the whole percolation configuration in $\mathbb{Z}^d$ instead of the configuration restricted to $\Lambda(n)$, because they are not connected to any vertex outside $\Lambda(n)$. We call them simply the finite clusters. In the supercritical regime $p > p_c$, it is well-known that the finite clusters have an exponentially decaying tail distribution. More precisely, we have the following result:

$$\forall p > p_c \quad \exists c > 0 \quad \forall n \geq 1 \qquad P_p\big(\text{diameter}\, C(0) \geq n,\, 0 \not\longleftrightarrow \infty\big) \;\leq\; e^{-cn}\,. \tag{33.33}$$

It would be extremely helpful to prove an estimate of this sort, or even a weaker one, starting from the sole hypothesis that $\theta(p) > 0$. We start by discussing how this kind of estimate could help control the supremum in (33.17). Afterwards, we discuss a strategy for deriving an estimate like (33.33).

Let us return to the exploration algorithm which builds $C\big(\partial_D^{out}\mathcal{A}(0), D\big)$. We consider first the standard exploration algorithm. This algorithm picks up a vertex in $\partial_D^{out}\mathcal{A}(0)$, explores its cluster, and reiterates the operation until all the vertices of $\partial_D^{out}\mathcal{A}(0)$ are explored. This is certainly a genuine exploration algorithm, for which the usual control is in force. If we choose randomly a subset $\mathcal{T}$ of $\partial_D^{out}\mathcal{A}(0)$ in a way that does not use any information on the percolation configuration outside the bonds of $\mathbb{E}^d(\mathcal{A}(0), D)$, and decide to make it a taboo set, then the exploration algorithm building $C\big(\partial_D^{out}\mathcal{A}(0) \setminus \mathcal{T}, D\big)$ is still genuine. The serious trouble begins when we choose a random subset $\mathcal{T}$ of $\partial_D^{out}\mathcal{A}(0)$ which carries some important information on the percolation configuration outside the bonds of $\mathbb{E}^d(\mathcal{A}(0), D)$. In this case, once we know the set $\mathcal{T}$, the distribution of the bonds outside $\mathbb{E}^d(\mathcal{A}(0), D)$ is not any more an i.i.d. Bernoulli field and the usual probabilistic control on the algorithm is not available. This is precisely what happens when we take for the taboo set $\mathcal{T}(1)$ the set of the vertices which are connected by an open path to $W'$. The problem we face is the following: what is the worst damage that such a taboo set can inflict on the usual probabilistic control?

Let us look once more at the algorithm which builds $C\big(\partial_D^{out}\mathcal{A}(0), D\big)$ and imagine, for the sake of the discussion, that $W'$ is a singleton $W' = \{\, x' \,\}$. In this case, all the vertices of $\mathcal{T}(1)$ belong to one single cluster $C'$, namely the cluster of $x'$. At some point, the algorithm will start exploring this cluster $C'$, which contains all the vertices of the set $\mathcal{T}(1)$. Afterwards, the algorithm will never again pick up a vertex of $\mathcal{T}(1)$. Let us compare this standard algorithm with the taboo exploration algorithm building $C\big(\partial_D^{out}\mathcal{A}(0) \setminus \mathcal{T}(1), D\big)$. The difference between the two algorithms is that the taboo exploration will not visit any vertex of $C'$. Therefore it seems that, as far as the sequence of the successive bonds explored by the algorithm is concerned, the net effect of the taboo set is the suppression of a sub-interval in the i.i.d. sequence of Bernoulli random variables $(X_k)_{k\geq 1}$ driving the exploration. If that were the case, then we could bound the supremum in (33.17) by the toy problem (33.21) with $r=1$, and we would indeed get an upper bound of order $\sqrt{n^d \ln n}$. We would then be in a good position to continue the proof. Unfortunately, the effect of the taboo set is far worse. Indeed, after the standard exploration has visited the cluster $C'$, the algorithm knows which bonds belong to $\Delta_D C'$. These bonds are closed, and they might be hit again during a future exploration of another cluster starting from $\partial_D^{out}\mathcal{A}(0)\setminus\mathcal{T}(1)$. So the edge boundary $\Delta_D C'$ is a further random set that influences the states of the bonds that are explored after the cluster $C'$ has been revealed. This completely hinders to recover the usual probabilistic control. If, by a fantastic chance, the cluster $C'$ is the last cluster to be explored when building $C\big(\partial_D^{out}\mathcal{A}(0)\setminus\mathcal{T},D\big)$, then the situation would be good, because this would annihilate the previous effect, according to which $C'$ creates a complex influence on the posterior exploration. The next question is, is it likely that $C'$ is the last cluster to be explored?

In fact, if the taboo set $\mathcal{T}(1)$ has a very small size compared to $\partial_D^{out}\mathcal{A}(0)$, indeed it will take quite some time until $C'$ is discovered. So there is some hope to get an interesting estimate in the regime where

$$|\mathcal{T}(1)| \ll |\partial_D^{out}\mathcal{A}(0)|\,. \tag{33.34}$$

However, when we iterate the argument as in the subsection 31.6, we cannot avoid having to control the sizes of the successive taboo sets, which typically increase too fast and violate very quickly any condition like (33.34). By symmetry, we could also hope that if $C'$ is the very first cluster to be explored, then we could reverse the sequence of the explored bonds in order to get back to the previous case. Unfortunately, we run into a similar problem, and we would have to work in a regime where

$$|\partial_D^{out}\mathcal{A}(0)\setminus\mathcal{T}(1)| \ll |\partial_D^{out}\mathcal{A}(0)|\,.$$

Another possibility is to try to design a procedure of exploration of $\partial_D^{out}\mathcal{A}(0)$ which avoids the taboo set as long as possible. Unfortunately, any such procedure will necessarily rely on some partial knowledge of the set $\mathcal{T}(1)$, and this inevitably ruins the estimates. The next hope is to develop some control on the clusters explored during the taboo exploration. To this end, we can take advantage of the fact that all the clusters which are explored after the first round are finite clusters (because we have supposed that $W=\partial^{\,in}\Lambda(n)$). Imagine that the finite clusters obey an estimate like (33.33). This would change a lot the rules of the exploration game. Indeed, inside the box $\Lambda(n)$, typically, all the finite clusters (that is, the clusters not intersecting the boundary $\partial^{\,in}\Lambda(n)$) would have a diameter which is $O(\ln n)$, hence a cardinality which is $O\big((\ln n)^d\big)$. As a consequence, the negative contribution of a finite cluster to the quantity (33.13) would also be $O\big((\ln n)^d\big)$. In turn, by a standard martingale argument, this yields that the mean contribution to (33.13) of a cluster intersecting $\partial^{\,in}\Lambda(n)$ is also $O\big((\ln n)^d\big)$. Let us precise this claim.

**Proposition 33.3.** *For any $n \geq 1$ and any vertex $x$ of $\Lambda(n)$, we have*

$$E\Big(S\big(\mathbb{E}^d(\overline{C}(x), D)\big)\Big) \,=\, 0\,. \tag{33.35}$$

*Proof.* We perform a standard exploration of the cluster of $x$ in $\Lambda(n)$. Let $(X_k)_{k\geq 1}$ be the sequence of i.i.d. Bernoulli random variables with parameter $p$ which drives the exploration algorithm. At step $k$, the algorithm has explored the set of bonds $O(k) \cup C(k)$ and by formula (13.1), we have

$$\big|O(k)\big| \,=\, S_k\,,\quad \big|C(k)\big| \,=\, k - S_k\,,$$

where, as usual,

$$\forall k \geq 1 \qquad S_k \,=\, X_1 + \cdots + X_k\,.$$

Let $T$ be the termination time of the algorithm. The random time $T$ is a stopping time with respect to the sequence of random variables $(X_k)_{k\geq 1}$, and it is bounded by $\big|\mathbb{E}^d(\Lambda(n))\big|$. Moreover, the sequence $(S_k - kp)_{k\geq 1}$ is a martingale. Applying the optional stopping theorem, we have

$$\forall m \geq 1 \qquad E\big(S_{T\wedge m} - (T \wedge m)p\big) \,=\, 0\,.$$

We send $m$ to $\infty$ and, with the help of the dominated convergence theorem, we get

$$E\big(S_T - Tp\big) \,=\, 0\,. \tag{33.36}$$

Moreover, we have

$$S\big(\mathbb{E}^d(\overline{C}(x), D)\big) \,=\, \frac{\big|O(T)\big|}{p} - \frac{\big|C(T)\big|}{1-p}\,,$$

whence

$$S_T - Tp \,=\, p(1-p)S\big(\mathbb{E}^d(\overline{C}(x), D)\big)\,. \tag{33.37}$$

Taking the expectation in formula (33.37) and using (33.36), we obtain the statement (33.35) of the proposition. ☐

Let $x$ be a vertex of $D = \Lambda(n)$. We decompose the expectation in formula (33.35) as follows:

$$\begin{aligned} E\Big(S\big(\mathbb{E}^d(\overline{C}(x), D)\big)\Big) \,=\, & E\Big(S\big(\mathbb{E}^d(\overline{C}(x), D)\big)1_{\{\,x\longleftrightarrow\partial^{\,in}D\,\}}\Big) \\ & + E\Big(S\big(\mathbb{E}^d(\overline{C}(x), D)\big)1_{\{\,x\not\longleftrightarrow\partial^{\,in}D\,\}}\Big)\,. \end{aligned} \tag{33.38}$$

On the event $\{\,x \not\longleftrightarrow \partial^{\,in}D\,\}$, the cluster $C(x)$ is a finite cluster, and from the previous discussion, its cardinality is $O\big((\ln n)^d\big)$ with overwhelming probability. Using proposition 33.3 together with formula (33.38), from the fact that the finite clusters of $D$ all have cardinality $O\big((\ln n)^d\big)$, we would conclude that

$$E\Big(S\big(\mathbb{E}^d(\overline{C}(x), D)\big)1_{\{\,x\longleftrightarrow\partial^{\,in}D\,\}}\Big) \,=\, O\big((\ln n)^d\big)\,. \tag{33.39}$$

Using the estimate (33.39), we see that the quantity (33.13) would satisfy

$$E\left(S\Big(\mathbb{E}^d\Big(C\big(\partial_D^{out}A_0\setminus T_1,D\big),D\setminus A_0\Big)\Big)\,\middle|\,\begin{matrix}\mathcal{T}(1)=T_1\\ \mathcal{A}(0)=A_0\end{matrix}\right)\;=\;O\big(|T_1|(\ln n)^d\big)\,.$$

Taking again the expectation, up to some additional steps which should be made rigorous, we could hope to conclude that

$$E\big(|\mathcal{I}(2)|\big)\;=\;E\big(|\mathcal{I}(1)|\big)\,O\big((\ln n)^d\big)\,. \tag{33.40}$$

This would be a great improvement compared to the brutal control of subsection 31.6. Regrettably, we are still very far from a rigorous proof of an estimate like (33.40). The previous argument was a bit sketchy, but it could probably be transformed into a rigorous proof. The main missing ingredient is of course an initial estimate like (33.33). We discuss in the next subsection an attempted strategy to obtain such an estimate.

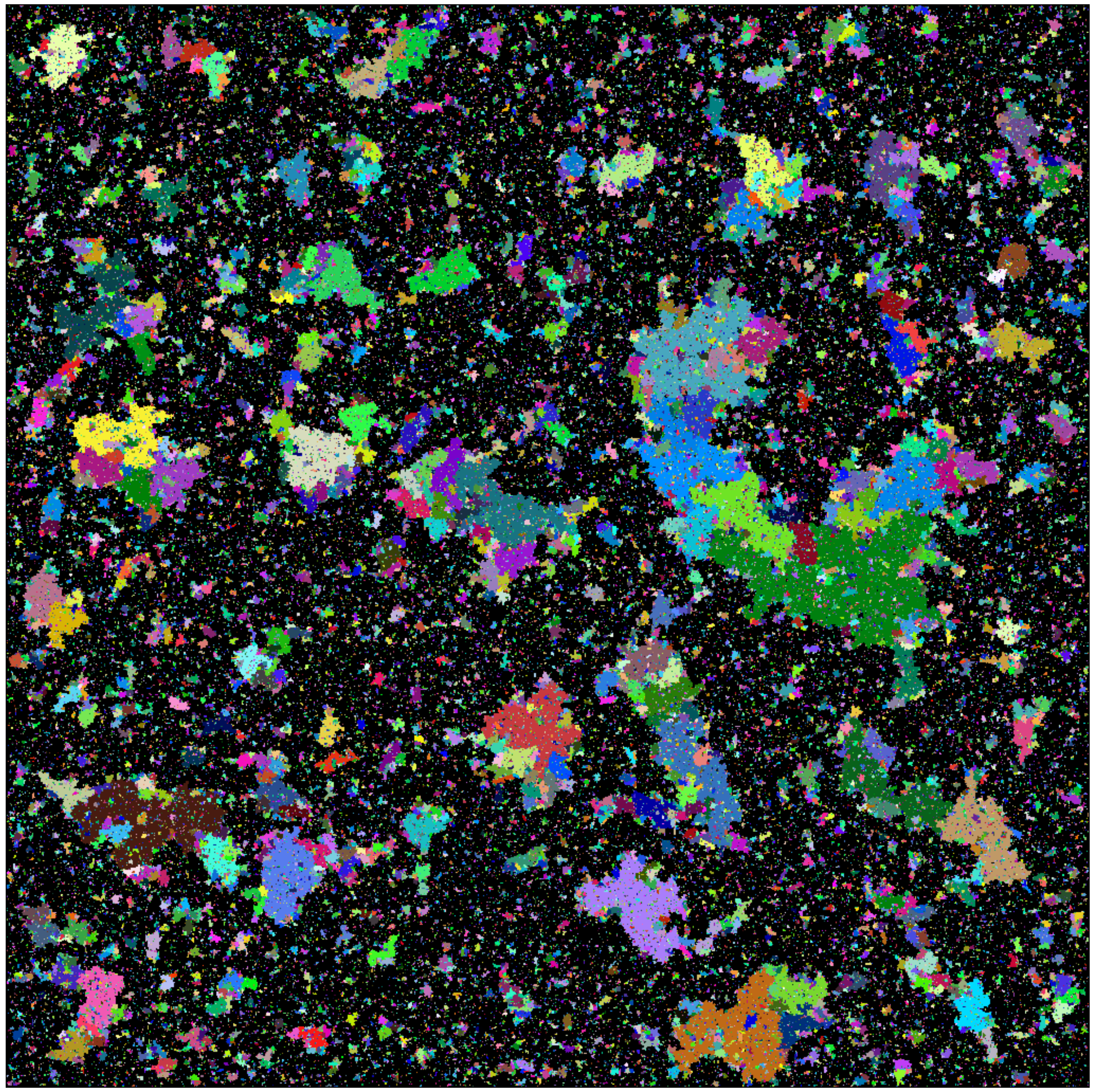

Figure 74: The finite clusters for bond percolation in $\Lambda(1024)$ at $p=0.5$

## 33.7 Removing a translation invariant set $\circ$

As indicated by the $\circ$ sign, this subsection is written for the site model, which, for once, is more convenient than the bond model to state and discuss the problem. Let $H$ be a random subset of $\mathbb{Z}^d$ which is translation invariant. The random field $(1_H(x), x \in \mathbb{Z}^d)$ is then translation invariant. We denote by $\mathcal{I}$ the $\sigma$-field of the translation invariant events. The spatial ergodic theorem (see for instance [92], section 6.1) yields that

$$\frac{1}{|\Lambda|}|H \cap \Lambda| \;=\; \frac{1}{|\Lambda|}\sum_{x\in\Lambda} 1_H(x) \quad\longrightarrow\quad P(0\in H\,|\,\mathcal{I}) \qquad \text{as } \Lambda\uparrow\mathbb{Z}^d\,,$$

where the convergence occurs almost surely and in the $L^1$ sense. Given a realization of the random set $H$, we define an associated percolation process on $\mathbb{Z}^d$ by declaring independently each site in $\mathbb{Z}^d\setminus H$ to be either open with probability $p$ or closed with probability $1-p$, while each site in $H$ is closed with probability 1. We denote by $P_p^H$ the corresponding probability measure on $\{0,1\}^{\mathbb{Z}^d}$. Denoting by $\delta_{\omega(x)=\varepsilon}$ the Dirac mass on the state $\omega(x)=\varepsilon$ and by $\otimes$ the tensor product of measures, we have, for any $\omega\in\{0,1\}^{\mathbb{Z}^d}$,

$$P_p^H(\omega) \;=\; \bigg(\bigotimes_{x\in H}\delta_{\omega(x)=0}\bigg)\otimes\bigg(\bigotimes_{x\in\mathbb{Z}^d\setminus H}\Big((1-p)\delta_{\omega(x)=0}+p\delta_{\omega(x)=1}\Big)\bigg)\,.$$

Alternatively, this percolation process can be described as Bernoulli site percolation on the random set $\mathbb{Z}^d\setminus H$. Notice that, although the random set $H$ is translation invariant, the probability measure $P_p^H$ is not translation invariant. We define next the percolation function $\theta^H$ associated to $P_p^H$. Since $P_p^H$ is not translation invariant, we must be careful and consider arbitrary starting points. The probability $\theta^H(x,p)$ that a given site $x$ of $\mathbb{Z}^d$ is connected to $\infty$ is defined as

$$\theta^H(x,p) \;=\; P_p^H\big(x\longleftrightarrow\infty\big)\,.$$

If $x$ belongs to $H$ or to a bounded component of $\mathbb{Z}^d\setminus H$, then the probability $\theta^H(x,p)$ is null. The probabilities $(\theta^H(x,p), x\in\mathbb{Z}^d)$ constitute a translation invariant random field. By the spatial ergodic theorem, the spatial averages of this random field converge towards the conditional expectation of $\theta^H(0,p)$ with respect to the $\sigma$-field $\mathcal{I}$, i.e.,

$$\frac{1}{|\Lambda|}\sum_{x\in\Lambda}\theta^H(x,p) \quad\longrightarrow\quad E\big(\theta^H(0,p)\,|\,\mathcal{I}\big) \qquad \text{as } \Lambda\uparrow\mathbb{Z}^d\,,$$

where the convergence occurs almost surely and in the $L^1$ sense. The quantity $E\big(\theta^H(0,p)\,|\,\mathcal{I}\big)$ is the spatial density of the infinite clusters (there might be several of them in general) present in the percolation configuration on $\mathbb{Z}^d\setminus H$ at level $p$. It is a consequence of the spatial ergodic theorem that this spatial average does not depend on the specific realization of the random set $H$, but only on the invariant $\sigma$-field $\mathcal{I}$. We define the critical point $p_c(\mathbb{Z}^d\setminus H)$ as the

parameter $p$ for which infinite clusters begin to appear on $\mathbb{Z}^d\setminus H$. More precisely, we set

$$p_c(\mathbb{Z}^d \setminus H) \,=\, \sup\big\{\, p\in[0,1] : \forall x\in\mathbb{Z}^d \quad \theta^H(x,p)=0\,\big\}\,.$$

Obviously, we have the inequality

$$p_c(\mathbb{Z}^d) \,\leq\, p_c(\mathbb{Z}^d \setminus H)\,. \tag{33.41}$$

Notice that the value $p_c(\mathbb{Z}^d \setminus H)$ depends on $H$, hence it is a random variable. However, it is defined in a translation invariant way, so that it is in fact measurable with respect to the invariant $\sigma$-field $\mathcal{I}$. We make the natural conjecture that, whenever the spatial density of the random set $H$ is not null with positive probability, the inequality (33.41) is strict.

**Conjecture**. If $P(O\in H)>0$, then $p_c(\mathbb{Z}^d) < p_c(\mathbb{Z}^d \setminus H)$.

To prove strict inequalities between critical points associated to different graphs is a delicate matter. We try to follow the general strategy presented in chapter 3 of [63].

Given a realization of the random set $H$, we define a percolation process depending on two parameters $p,s$ by declaring independently each site in $\mathbb{Z}^d\setminus H$ to be either open with probability $p$ or closed with probability $1-p$, and each site in $H$ to be either open with probability $s$ or closed with probability $1-s$. We denote by $P^H_{p,s}$ the corresponding probability measure on $\{0,1\}^{\mathbb{Z}^d}$, i.e., for any $\omega\in\{0,1\}^{\mathbb{Z}^d}$, we have

$$P^H_{p,s}(\omega) \,=\, \Bigg(\bigotimes_{x\in H}\Big((1-s)\delta_{\omega(x)=0}+s\delta_{\omega(x)=1}\Big)\Bigg) \otimes \Bigg(\bigotimes_{x\in \mathbb{Z}^d\setminus H}\Big((1-p)\delta_{\omega(x)=0}+p\delta_{\omega(x)=1}\Big)\Bigg)\,.$$

We define also

$$\forall x\in\mathbb{Z}^d \qquad \theta^H(x,p,s) \,=\, P^H_{p,s}\big(x\longleftrightarrow\infty\big)\,.$$

The subcritical region for the parameters $p,s$ is

$$\big\{\,(p,s)\in[0,1]^2 : \forall x\in\mathbb{Z}^d \quad \theta^H(x,p,s) \,=\, 0\,\big\}\,.$$

The strategy designed by Aizenman and Grimmett [3] consists in proving that the boundary of this region contains no vertical segment. To achieve this goal, they compute estimates on the gradient of the percolation function, which depends on two parameters. Here our percolation function $\theta^H(x,p,s)$ depends in addition on the site $x$. So we perform a spatial average to get rid of this dependence. The probability $\theta^H(x,p,s)$ is a random variable, but it is a deterministic function of $H$. Since $H$ is translation invariant, then the random field $(\theta^H(x,p,s), x\in\mathbb{Z}^d)$ is also translation invariant. The function $\phi(p,s)$ defined by

$$\phi(p,s) \,=\, E\big(\theta^H(0,p,s)\,|\,\mathcal{I}\big)$$

is the counterpart in our framework of the two parameters percolation function $\theta(p,s)$ of Aizenman and Grimmett. So the next goal would be to compute the partial derivatives of $\phi$. In order to estimate the variations of $\phi$ with respect to $p$ and $s$, we introduce finite volume approximations of $\phi$, for which we can compute the partial derivatives. From the spatial ergodic theorem, we know that

$$\phi(p,s) \;=\; \lim_{\Lambda\uparrow\mathbb{Z}^d} \frac{1}{|\Lambda|}\sum_{x\in\Lambda} P^H_{p,s}\big(x\longleftrightarrow\infty\big)\,, \tag{33.42}$$

but this is not enough, because the probability $P^H_{p,s}\big(x\longleftrightarrow\infty\big)$ still depends on infinitely many sites, so it is not clear whether it is differentiable. A better quantity to work with is in fact the counterpart of the quantity $\varrho$ used to perform the numerical simulations (see formula (5.1)). So we define, for $\Lambda$ a finite box,

$$\varrho^H(\Lambda,p,s) \;=\; \sum_{x\in\Lambda} P^H_{p,s}\big(x\longleftrightarrow\partial^{\,in}\Lambda\big)\,.$$

We show first that $\frac{1}{|\Lambda|}\varrho(\Lambda,p,s)$ is a finite volume approximation of $\phi(\Lambda,p,s)$.

**Lemma 33.4.** *For any $p,s$ in $[0,1]$, we have*

$$\phi(p,s) \;=\; \lim_{\Lambda\uparrow\mathbb{Z}^d} \frac{1}{|\Lambda|}\varrho^H(\Lambda,p,s) \quad \textit{almost surely and in } L^1\,.$$

*Proof.* For any box $\Lambda$, we have

$$\forall x\in\Lambda\qquad P^H_{p,s}\big(x\longleftrightarrow\infty\big)\;\leq\;P^H_{p,s}\big(x\longleftrightarrow\partial^{\,in}\Lambda\big)\,.$$

Summing over $x$ in $\Lambda$, we get

$$\sum_{x\in\Lambda} P^H_{p,s}\big(x\longleftrightarrow\infty\big)\;\leq\;\varrho^H(\Lambda,p,s)\,.$$

Dividing by $|\Lambda|$ and sending $\Lambda$ to $\mathbb{Z}^d$, we obtain, with the help of the limit (33.42),

$$\phi(p,s)\;\leq\;\liminf_{\Lambda\uparrow\mathbb{Z}^d}\frac{1}{|\Lambda|}\varrho^H(\Lambda,p,s)\,. \tag{33.43}$$

Conversely, we fix $m\leq n$. For any given percolation configuration, we have

$$\sum_{x\in\Lambda(n)} 1_{\{\,x\longleftrightarrow\partial^{\,in}\Lambda(n)\,\}}\;\leq\;|\Lambda(n)\setminus\Lambda(n-m)| \\ +\sum_{x\in\Lambda(n-m)} 1_{\{\,x\longleftrightarrow x+\partial^{\,in}\Lambda(m)\,\}}\,.$$

We take the expectation with respect to $P^H_{p,s}$, and we get

$$\varrho^H(\Lambda(n),p,s)\;\leq\;|\Lambda(n)\setminus\Lambda(n-m)|+\sum_{x\in\Lambda(n-m)} P^H_{p,s}\big(x\longleftrightarrow x+\partial^{\,in}\Lambda(m)\big)\,.$$

We divide by $\Lambda(n)$, and we obtain the inequality

$$\frac{\varrho^H(\Lambda(n),p,s)}{|\Lambda(n)|} \leq \frac{|\Lambda(n)\setminus\Lambda(n-m)|}{|\Lambda(n)|} + \frac{1}{|\Lambda(n-m)|}\sum_{x\in\Lambda(n-m)} P^H_{p,s}\big(x\longleftrightarrow x+\partial^{\,in}\Lambda(m)\big)\,. \quad (33.44)$$

We keep $m$ fixed and we send $n$ to $\infty$. The ergodic theorem yields

$$\lim_{n\to\infty}\frac{1}{|\Lambda(n-m)|}\sum_{x\in\Lambda(n-m)} P^H_{p,s}\big(x\longleftrightarrow x+\partial^{\,in}\Lambda(m)\big) = E\Big(P^H_{p,s}\big(0\longleftrightarrow \partial^{\,in}\Lambda(m)\big)\,\Big|\,\mathcal{I}\Big)\,,$$

so that, passing to the limit in (33.44), we have

$$\limsup_{n\to\infty}\frac{\varrho^H(\Lambda(n),p,s)}{|\Lambda(n)|} \leq E\Big(P^H_{p,s}\big(0\longleftrightarrow \partial^{\,in}\Lambda(m)\big)\,\Big|\,\mathcal{I}\Big)\,. \quad (33.45)$$

Applying the monotone convergence theorem, we have

$$\lim_{m\to\infty} E\Big(P^H_{p,s}\big(0\longleftrightarrow \partial^{\,in}\Lambda(m)\big)\,\Big|\,\mathcal{I}\Big) = E\big(\theta^H(0,p,s)\,|\,\mathcal{I}\big) = \phi(p,s)\,. \quad (33.46)$$

Sending $m$ to $\infty$ in (33.45) and using (33.46), we conclude that

$$\limsup_{n\to\infty}\frac{\varrho^H(\Lambda(n),p,s)}{|\Lambda(n)|} \leq \phi(p,s)\,. \quad (33.47)$$

The two inequalities (33.43) and (33.47) yield the desired conclusion. □

The point is that the quantity $\varrho^H(\Lambda,p,s)$ depends only on the states of the sites in the box $\Lambda$, so we can safely compute its partial derivatives with respect to $p$ and $s$. The states of the sites in $\mathbb{Z}^d\setminus H$ (respectively $H$) are governed by the parameter $p$ (respectively $s$). With the help of Russo's formula, we obtain that for $x\in\Lambda$,

$$\frac{\partial}{\partial p}P^H_{p,s}\big(x\longleftrightarrow\partial^{\,in}\Lambda\big) = \sum_{y\in\Lambda\setminus H} P^H_{p,s}\big(y \text{ is pivotal for } \{\,x\longleftrightarrow\partial^{\,in}\Lambda\,\}\big)\,,$$

$$\frac{\partial}{\partial s}P^H_{p,s}\big(x\longleftrightarrow\partial^{\,in}\Lambda\big) = \sum_{y\in\Lambda\cap H} P^H_{p,s}\big(y \text{ is pivotal for } \{\,x\longleftrightarrow\partial^{\,in}\Lambda\,\}\big)\,,$$

whence, by summing over $x$ in $\Lambda$,

$$\frac{\partial}{\partial p}\varrho^H(\Lambda,p,s) = \sum_{x\in\Lambda}\sum_{y\in\Lambda\setminus H} P^H_{p,s}\big(y \text{ is pivotal for } \{\,x\longleftrightarrow\partial^{\,in}\Lambda\,\}\big)\,,$$

$$\frac{\partial}{\partial s}\varrho^H(\Lambda,p,s) = \sum_{x\in\Lambda}\sum_{y\in\Lambda\cap H} P^H_{p,s}\big(y \text{ is pivotal for } \{\,x\longleftrightarrow\partial^{\,in}\Lambda\,\}\big)\,. \quad (33.48)$$

The next step of the strategy of Aizenman and Grimmett is to obtain inequalities between these partial derivatives which are uniform with respect to the box $\Lambda$. The uniformity ensures that the very same inequalities still hold for the limit when $\Lambda$ grows to $\mathbb{Z}^d$. Regrettably, there is no hope of fulfilling this program here, because of the presence of the random set $H$. The modifications of the graph structure induced by $H$ might be highly non-local, and this prevents us from comparing term by term the probabilities appearing in the sums in (33.48): there is no hope that, by a local modification of the configuration, we can transform a pivotal site $y$ in $H$ into a pivotal site $y$ in $\mathbb{Z}^d \setminus H$. Although there is not much hope to prove that the critical curve does not contain a vertical segment, a weaker result would be enough to prove the strict inequality between the critical points $p_c(\mathbb{Z}^d)$ and $p_c(\mathbb{Z}^d \setminus H)$. We first rewrite the definitions of these critical points with the help of the function $\phi(p,s)$. We have

$$\begin{aligned} p_c(\mathbb{Z}^d) \;&=\; \sup\left\{\, p \in [0,1] : \phi(p,p) = 0 \,\right\}, \\ p_c(\mathbb{Z}^d \setminus H) \;&=\; \sup\left\{\, p \in [0,1] : \phi(p,0) = 0 \,\right\}. \end{aligned}$$

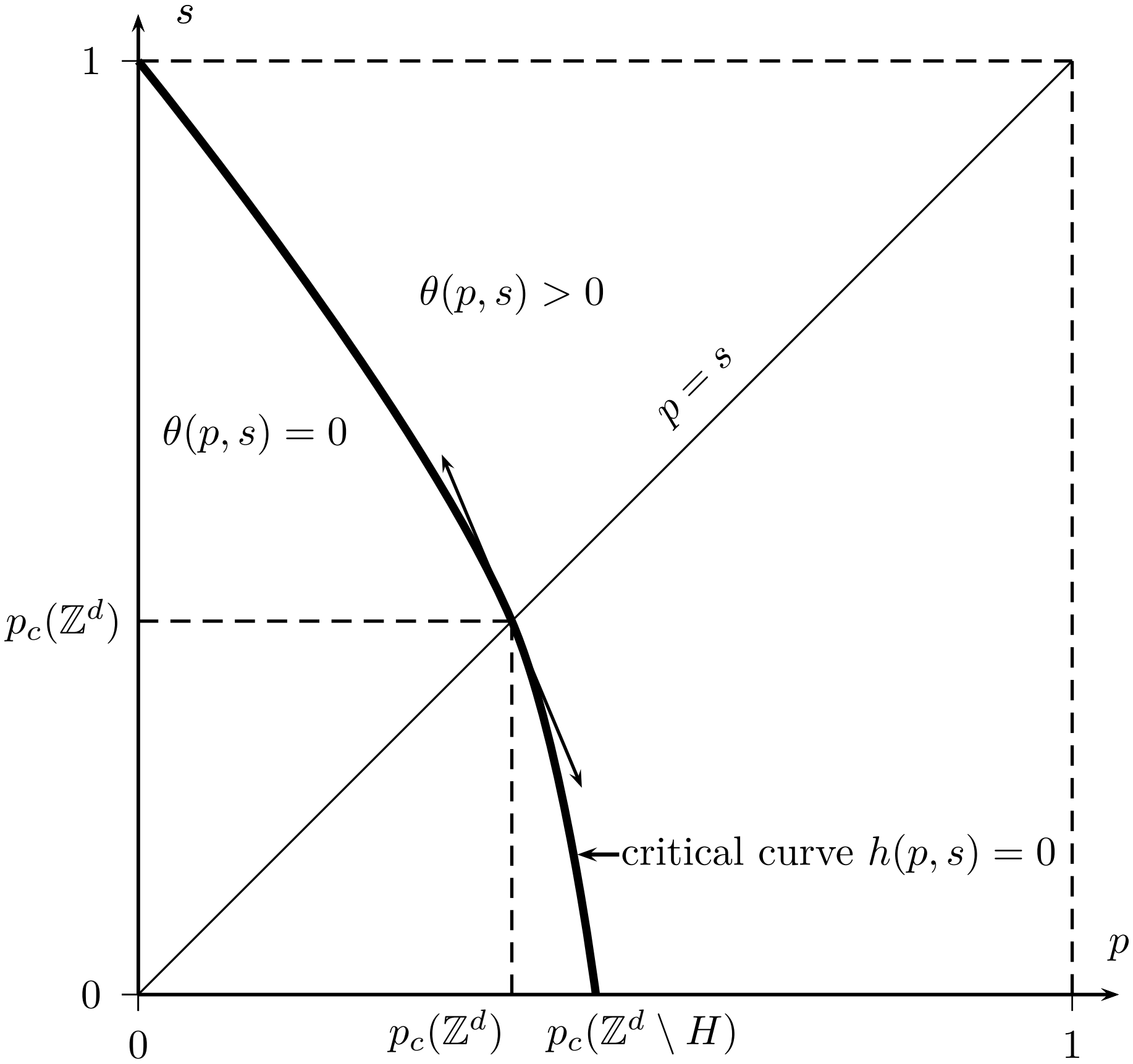

Figure 75: The critical region in the $(p,s)$ plane (Aizenman–Grimmett method)

Suppose, as in chapter 3 of [63], that the curve in the $(p,s)$-plane separating the regime where $\phi(p,s) > 0$ from the regime where $\phi(p,s) = 0$ may be written in the form of an equation $h(p,s) = 0$. In their situation, where the modifications of the lattice are local, Aizenman and Grimmett proved a uniform lower bound on the gradient $\nabla h$ of $h$, of the form

$$\forall p,s\in ]0,1[^2\qquad \frac{\partial h}{\partial s}(p,s)\,\geq\,\alpha(p,s)\frac{\partial h}{\partial p}(p,s)\,.$$

Moreover, they have a positive uniform lower bound on $\alpha(p,s)$ when $p,s$ vary in a compact subset of $]0,1[^2$. In order to prove the strict inequality between $p_c(\mathbb{Z}^d)$ and $p_c(\mathbb{Z}^d\setminus H)$, it would be enough to have a lower bound on $\nabla h$ at the point $(p_c(\mathbb{Z}^d),p_c(\mathbb{Z}^d))$, like

$$\frac{\partial h}{\partial s}\big(p_c(\mathbb{Z}^d),p_c(\mathbb{Z}^d)\big)\,\geq\,\alpha\frac{\partial h}{\partial p}\big(p_c(\mathbb{Z}^d),p_c(\mathbb{Z}^d)\big)\,,$$

where $\alpha>0$. There is some hope to achieve this goal, because, on the diagonal $p=s$, the modified percolation process $P^H_{p,s}$ coincides with the original percolation process $P_p$, namely, we have

$$\forall p\in[0,1]\qquad P^H_{p,s}\,=\,P_p\,. \tag{33.49}$$

The notation becomes a bit confusing, because we wish to estimate the partial derivative of $\varrho^H(\Lambda,p,s)$ with respect to $p$, when evaluated on the diagonal $p=s$. So, we fix a point $p_0\in ]0,1[$, and we apply the formula (33.48) at the point $(p,s)=(p_0,p_0)$ to get

$$\frac{\partial}{\partial p}\varrho^H(\Lambda,p,s)|_{(p,s)=(p_0,p_0)}\;=\sum_{x\in\Lambda}\sum_{y\in\Lambda\setminus H}P^H_{p_0,p_0}\big(y\text{ is pivotal for }\{\,x\longleftrightarrow\partial^{\,in}\Lambda\,\}\big)\,,$$

$$\frac{\partial}{\partial s}\varrho^H(\Lambda,p,s)|_{(p,s)=(p_0,p_0)}\;=\sum_{x\in\Lambda}\sum_{y\in\Lambda\cap H}P^H_{p_0,p_0}\big(y\text{ is pivotal for }\{\,x\longleftrightarrow\partial^{\,in}\Lambda\,\}\big)\,.$$

Thanks to the magical equality (33.49), we see that

$$\frac{\partial}{\partial p}\varrho^H(\Lambda,p,s)|_{(p,s)=(p_0,p_0)}\;=\sum_{x\in\Lambda}\sum_{y\in\Lambda\setminus H}P_{p_0}\big(y\text{ is pivotal for }\{\,x\longleftrightarrow\partial^{\,in}\Lambda\,\}\big)\,,$$

$$\frac{\partial}{\partial s}\varrho^H(\Lambda,p,s)|_{(p,s)=(p_0,p_0)}\;=\sum_{x\in\Lambda}\sum_{y\in\Lambda\cap H}P_{p_0}\big(y\text{ is pivotal for }\{\,x\longleftrightarrow\partial^{\,in}\Lambda\,\}\big)\,.$$

In the next step, we take the conditional expectation with respect to $\mathcal{I}$ (denoted by $E_{\mathcal{I}}$). It follows from the translation invariance of $H$ that $P(y\in\mathbb{Z}^d\setminus H\,|\,\mathcal{I})$ and $P(y\in H\,|\,\mathcal{I})$ do not depend on $y$ in $\Lambda$, thus we obtain

$$\begin{aligned}E_{\mathcal{I}}\Big(\frac{\partial}{\partial p}\varrho^H(\Lambda,p,s)|_{(p,s)=(p_0,p_0)}\Big)\;&=\sum_{x\in\Lambda}\sum_{y\in\Lambda}P(0\in\mathbb{Z}^d\setminus H\,|\,\mathcal{I})\\&\qquad\times P_{p_0}\big(y\text{ is pivotal for }\{\,x\longleftrightarrow\partial^{\,in}\Lambda\,\}\big)\,,\\E_{\mathcal{I}}\Big(\frac{\partial}{\partial s}\varrho^H(\Lambda,p,s)|_{(p,s)=(p_0,p_0)}\Big)\;&=\sum_{x\in\Lambda}\sum_{y\in\Lambda}P(0\in H\,|\,\mathcal{I})\\&\qquad\times P_{p_0}\big(y\text{ is pivotal for }\{\,x\longleftrightarrow\partial^{\,in}\Lambda\,\}\big)\,.\end{aligned}$$

From the two previous equalities, we conclude that

$$E_{\mathcal{I}}\Big(\frac{\partial}{\partial p}\varrho^H(\Lambda,p,s)|_{(p,s)=(p_0,p_0)}\Big)P(0\in H\,|\,\mathcal{I})\,=\,E_{\mathcal{I}}\Big(\frac{\partial}{\partial s}\varrho^H(\Lambda,p,s)|_{(p,s)=(p_0,p_0)}\Big)P(0\in\mathbb{Z}^d\setminus H\,|\,\mathcal{I})\,. \tag{33.50}$$

For any realization of the set $H$ and for any finite box $\Lambda$, the function $(p,s)\mapsto\varrho^H(\Lambda,p,s)$ is a polynomial function in the variables $(p,s)$, so that we can interchange the differentiation and the expectation in (33.50) to conclude that

$$\frac{\partial}{\partial p}E_{\mathcal{I}}\Big(\varrho^H(\Lambda,p,s)|_{(p,s)=(p_0,p_0)}\Big)P(0\in H\,|\,\mathcal{I})\,=\,\frac{\partial}{\partial s}E_{\mathcal{I}}\Big(\varrho^H(\Lambda,p,s)|_{(p,s)=(p_0,p_0)}\Big)P(0\in\mathbb{Z}^d\setminus H\,|\,\mathcal{I})\,. \tag{33.51}$$

Dividing by $\Lambda$, sending $\Lambda$ to $\mathbb{Z}^d$, and using lemma 33.4, we should obtain

$$\frac{\partial\phi}{\partial p}(p_0,p_0)P(0\in H\,|\,\mathcal{I})\,=\,\frac{\partial\phi}{\partial s}(p_0,p_0)P(0\in\mathbb{Z}^d\setminus H\,|\,\mathcal{I})\,. \tag{33.52}$$

Unfortunately, this last step cannot be made rigorous easily. We do not know that the limiting function $\phi$ admits partial derivatives, and we had to interchange a limit and a derivative to go from (33.51) to (33.52). If the issues raised in the previous sentence could be addressed, then we would be close to proving rigorously that $p_c(\mathbb{Z}^d)<p_c(\mathbb{Z}^d\setminus H)$. Now, how would this result help us for the conjecture $\theta(p_c,\mathbb{Z}^d)=0$?

The faint hope is to proceed as follows. Let $p\in]0,1[$ be such that $\theta(p)>0$. Let $\omega$ be a percolation configuration drawn according to $P_p$. Let $C_\infty(\omega)$ be the unique infinite cluster present in $\omega$. We take for $H$ the set

$$H\,=\,\overline{C}_\infty(\omega)\,=\,C_\infty(\omega)\cup\partial^{\,out}C_\infty(\omega)\,=\,\big\{\,x\in\mathbb{Z}^d:\exists\,y\in C_\infty(\omega)\quad|x-y|\leq 1\,\big\}\,.$$

The set $H$ is a random set and it is translation invariant. The associated $\sigma$-field $\mathcal{I}$ is trivial, and $H$ has positive density $\theta(p)$. So, using the previous conjecture, we would have $p_c(\mathbb{Z}^d)<p_c(\mathbb{Z}^d\setminus\overline{C}_\infty(\omega))$. This means that, if we work at the critical point $p_c(\mathbb{Z}^d)$, we are in the subcritical regime for the percolation restricted to $\mathbb{Z}^d\setminus\overline{C}_\infty(\omega)$, and we could hope to get an estimate pertaining to the subcritical phase like (33.33), or a weaker version, maybe for an average over a large box $\Lambda$. Notice that, although we have not been able to prove that $p_c(\mathbb{Z}^d)<p_c(\mathbb{Z}^d\setminus\overline{C}_\infty(\omega))$, we know that

$$\forall x\in\mathbb{Z}^d\setminus\overline{C}_\infty(\omega)\qquad\theta^{\overline{C}_\infty(\omega)}(x,p,p)\,=\,P^{\overline{C}_\infty(\omega)}_{p,p}\big(x\longleftrightarrow\infty\big)\,=\,0\,.$$

Indeed, this is a simple consequence of the uniqueness of the infinite cluster. Anyway, in view of the numerous gaps that need to be filled to make the whole previous argument rigorous, this is all rather science fiction.

Part VI
# The Erdős-Rényi model ⊛

Before delving further into the main argument, we take a pause to examine how it could be put into action in the context of another famous classical model, the Erdős-Rényi model. In section 34, we recall first how the phase transition manifests itself in this model. In section 35, we discuss the classical proofs for the existence and uniqueness of the giant component. In section 36, we reformulate the standard genuine exploration algorithm in the framework of the Erdős-Rényi model, and we derive a specific associated inequality. In section 37, we develop a new inequality to control the probability that two large clusters coexist. In section 38, we discuss two applications of this inequality.

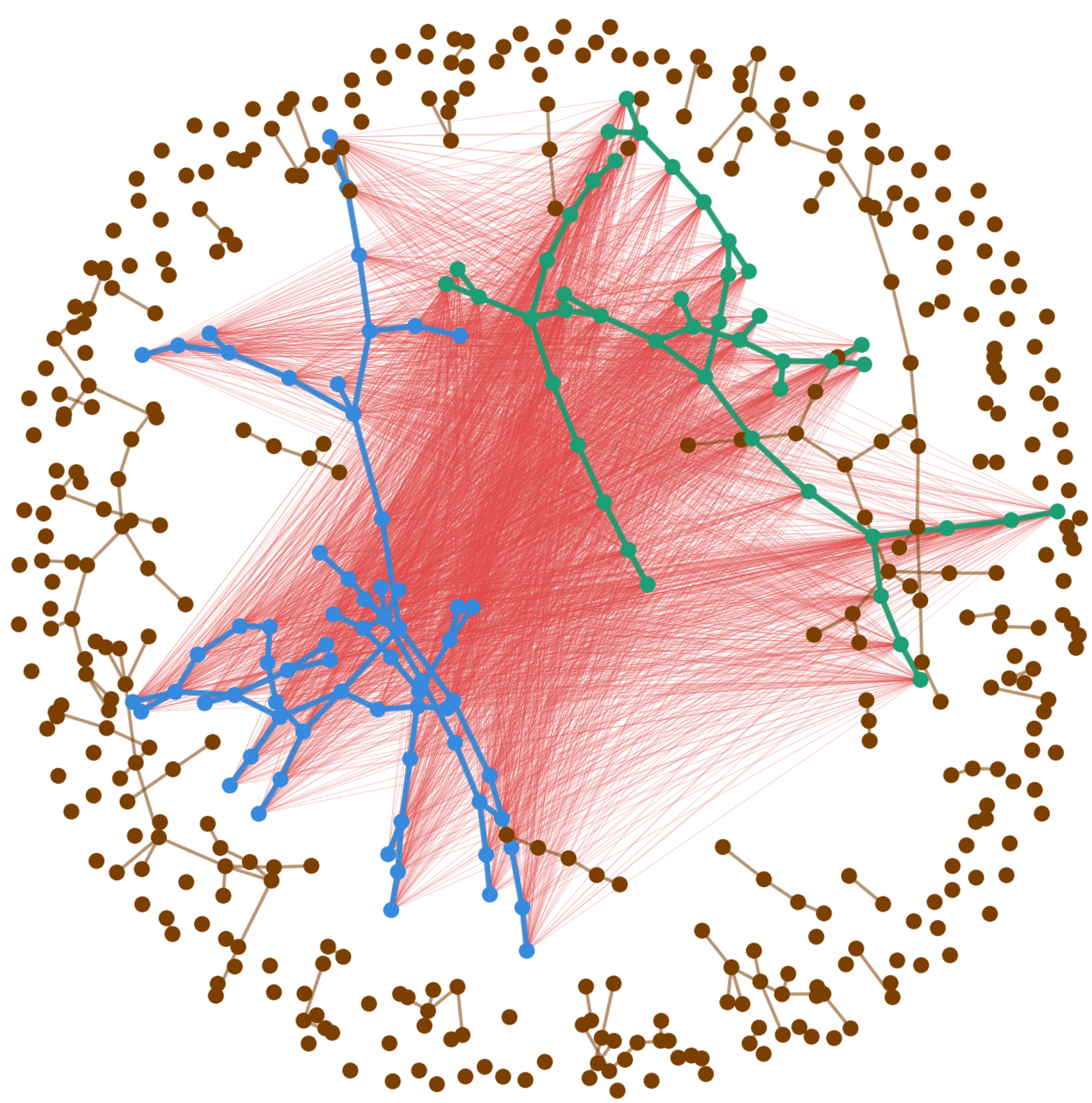

Figure 76: The Erdős-Rényi model with $n = 500$, $p = 0.002$
The thin red lines are the closed bonds between the two largest components.
The largest component has size 67, the second largest 45.

# 34 The phase transition

Let us recall briefly the definition of the model. The Erdős-Rényi model with parameters $n \geq 1$ and $p$ in $[0,1]$ is the bond percolation model with parameter $p$ on the complete graph with $n$ vertices. More precisely, each edge $\{ i, j \}$ between two distinct vertices $i, j$ of $\{ 1, \dots, n \}$ is open with probability $p$ or closed with probability $1-p$, and all the edges are independent. We denote by $G(n,p)$ the random graph with vertex set $\{ 1, \dots, n \}$ and edge set the open edges only.

**Theorem 34.1.** *We reparametrize the model by setting $p = c/n$ where $c > 0$. The model presents a phase transition as $n$ goes to $\infty$, described as follows:*
• *Subcritical regime. For $c < 1$, there exists a positive constant $\kappa$, which depends on $c$, such that, with probability going to 1 as $n$ goes to $\infty$, all the clusters in $G(n, c/n)$ have cardinality less than $\kappa \ln n$.*
• *Supercritical regime. For $c > 1$, there exist two positive constants $\kappa$, $\zeta$, which depend on $c$, such that, with probability going to 1 as $n$ goes to $\infty$, there exists one cluster in $G(n, c/n)$ which has cardinality of order $\zeta n$, and all the other clusters in $G(n, c/n)$ have cardinality less than $\kappa \ln n$.*

This result was proved by Erdős and Rényi in 1960 [51]. The historical approach of Erdős and Rényi [51] is combinatorial in nature, it involves the enumeration of various structures that appear when the percolation parameter is raised. Nowadays, the model and numerous results are presented in various textbooks. Let us first mention two books called random graphs, one by Svante Janson, Tomasz Luczak and Andrzej Ruciński [77], and the other one by Béla Bollobás [17]. These books present very detailed and technical results, mainly based on combinatorial arguments. Another classical presentation can be found in chapter 11 of the book of Noga Alon and Tom Spencer [7]. More recently, Remco van der Hofstad wrote a beautiful monograph [133], which outlines the approach based on the exploration processes. Finally, Nicolas Curien [76] presents several other possible techniques to apprehend this phase transition.

Needless to say, mathematicians have conducted a thorough analysis of the phenomena occurring at and around the critical point $p = 1/n$, and much more is known than the initial facts stated in theorem 34.1. The birth of the giant component is described with a surprising wealth of details in the monumental paper [76], which can leave us longing for the percolation model in $\mathbb{Z}^d$. The strength of the Erdős-Rényi model is that it is a mean-field model, making it amenable to deep mathematical analysis. However, this is also its weakness, as it belongs to a different universality class than the percolation model in $\mathbb{Z}^d$. We certainly won't be presenting any new results on the Erdős-Rényi model here. We simply wish to use this model as a testing ground for our strategy for proving the conjecture $\theta(p_c, \mathbb{Z}^d) = 0$. To begin with, what would be the equivalent question in the context of the Erdős-Rényi model?

As the giant component plays the role of the infinite cluster, and its density is equal to $\zeta = \zeta(c)$, the natural counterpart of the $\theta(p_c, \mathbb{Z}^d) = 0$ conjecture is the statement that $\zeta(1) = 0$. This is a well-known fact that is proven using various methods in the texts mentioned above. However, none of them corresponds to the strategy we intend to use for attacking the conjecture $\theta(p_c, \mathbb{Z}^d) = 0$. The most illuminating text, and the one which is closest to our point of view is the one of Remco van der Hofstad [133]. Chapter 4 of [133] contains a full proof of theorem 34.1, and much more. The proof is based on couplings between the exploration process of a cluster and various branching processes. The offspring distribution of these branching processes is the binomial distribution with parameters $n, c/n$, which can in turn be approximated by the Poisson distribution of parameter $c$. One of the main result of chapter 4 of [133] is the theorem 4.8, which proves that, in the supercritical regime, the cardinality of the largest cluster is polynomially close to $\zeta(c)\, n$, where $\zeta(c)$ is the survival probability of the Poisson branching process with parameter $c$. It is a classical fact of the theory of branching processes that this probability of survival is positive if and only if the mean offspring distribution (which is equal to $c$ in our case) is strictly larger than 1, thus

$$\zeta(c) > 0 \quad \iff \quad c > 1\,. \tag{34.1}$$

This is exactly the result we want to prove. The result (34.1) is an easy piece of cake once $\zeta(c)$ is identified as the survival probability of the Poisson branching process with parameter $c$. Indeed, the extinction probability $1 - \zeta(c)$ is characterized as the smallest fixed point of the generating function of the offspring distribution, hence it satisfies the equation

$$1 - \zeta(c) \;=\; \exp\big( - c\,\zeta(c)\big)\,.$$

No one believes that the percolation probability $\theta(p, \mathbb{Z}^d)$ can be characterized as the fixed point of some nice simple function. An interesting question is to design a more robust approach for proving that $\zeta(1) = 0$. While this can indeed be done, it does not lead to any new insights for either the Erdős-Rényi model or the percolation model in $\mathbb{Z}^d$, so we refrain from presenting it. Instead, we reformulate the technique employed in our so far unsuccessful attempts to attack the conjecture $\theta(p_c, \mathbb{Z}^d) = 0$ within the framework of the Erdős-Rényi model. Although it will not lead to new results, we obtain an interesting inequality controlling the probability that two large clusters coexist.

**Notation.** The Erdős-Rényi model being a bond model, we shall reuse here our notation for the bond percolation model in $\mathbb{Z}^d$. For instance, the open cluster of the vertex $v$ in $\{1, \dots, n\}$ is denoted by $C(v)$. As regards the specific details of the Erdős-Rényi model, we will adhere as closely as possible to the book [133]. Nevertheless, some variations remain. For instance, the clusters are denoted with the calligraphic letter $\mathcal{C}$ in [133], whereas we still use $C$ here. More annoying, we use what is perhaps the more traditional convention to parametrize the probability $p$ involved in the definition of the random graph $G(n, p)$, namely we set

$$p = \frac{c}{n}\,, \quad c > 0\,, \qquad p = \frac{1}{n} + \frac{\lambda}{n^{4/3}}\,, \quad \lambda \in \mathbb{R}\,.$$

We use $c$ to discuss the phase transition and $\lambda$ to describe the critical window.

## 35 Outline of the proof

Our main interest here is the giant component, and its emergence. Let us consider first the supercritical regime and let us discuss how to prove the second part of theorem 34.1. Currently, the most efficient proof consists in coupling the exploration process of a cluster with a branching process, firstly one with binomial offspring distribution, and secondly one with Poisson offspring distribution. This is the strategy employed in the three references [7], [76] and [133]. One needs to obtain matching lower and upper bounds, this requires to develop adequate couplings in both directions. Let us describe the three main steps common to the proofs presented in [7], [76] and [133]:

• Control of the small clusters. The goal of this step is to control the density of the clusters having cardinality smaller than $a \ln n$, where $a$ is a suitable positive constant. In a nutshell, the couplings show that the cardinality of the small clusters has a distribution close to the total progeny of a Poisson branching process (see for instance theorems 4.2 and 4.3 of [133]). An immediate consequence of these first inequalities is that the graph $G(n, c/n)$ typically contains clusters of cardinality of order $a \ln n$, where $a$ is a positive constant which depends on $c$. The determination of the optimal choice of the constant $a$ requires some additional work and is also done in [133]. It is linked with the large deviations rate function of the Poisson distribution of parameter $c$.

• No middle ground. The goal of this step is to show that it is extremely unlikely for the graph $G(n, c/n)$ to contain a cluster having cardinality larger than $a \ln n$ but smaller than $bn$, where $a, b$ are suitable positive constants. This step is more technical than the previous one, as it requires more refined deviations inequalities. Both references [7] and [133] rely on the Chernoff inequality for the binomial distribution. The corresponding paragraph in section 11.9 of [7] explains the required computations. Proposition 4.12 of [133] presents a nice inequality involving the rate function of the Poisson distribution which leads to more precise estimates. Another outcome of these estimates is that the density of the small clusters is close to $1 - \zeta(c)$, where $\zeta(c)$ is the survival probability of the Poisson branching process with parameter $c$.

• Existence and uniqueness of the giant. We call a cluster small if it has cardinality less than $a \ln n$ and giant if it has cardinality larger than $bn$, where $a, b$ are the two positive constants given by the first two steps. Taken together, the first two steps readily imply the existence of at least one giant. Indeed, we know from the first step that the density of the small clusters is close to $1 - \zeta(c)$, which is strictly less than one when $c > 1$. This leaves us with a positive density $\zeta(c)$ of sites belonging to clusters which are not small. But the second step tells us that, with high probability, the clusters are either small or giant. Therefore there is a density close to $\zeta(c)$ of sites belonging to giant clusters. To conclude, it remains to prove that it is very unlikely to observe two giants. This is the most tricky part, that we discuss in the section 38.

In the presentation of [76], proof of theorem 6.2, the last two steps are performed simultaneously with the help of a specific sprinkling procedure.

# 36 The genuine exploration algorithm

We reuse the standard exploration algorithm described for the bond percolation model in subsection 28.1. Let $v$ be any vertex in $\{1,\dots,n\}$. We denote by $\mathcal{E}(v)$ the set of all the bonds having at least one endpoint into the open cluster $C(v)$ of $v$, i.e.,

$$\mathcal{E}(v) \,=\, \big\{\,\{\,i,j\,\} : i\in C(v)\,, j\in\{\,1,\dots,n\,\}\,, i\neq j\,\big\}\,.$$

We develop next a specific variant of the deviations inequality (28.2), well-suited to the Erdős-Rényi model. The main point is that, whenever a cluster in $G(n,p)$ has cardinality $i$, then the number $b(i,n)$ of bonds examined during its exploration is equal to

$$b(i,n) \,=\, \frac{1}{2}i(i-1)+i(n-i) \,=\, in-\frac{1}{2}i(i+1)\,.$$

This is quite different from the situation in $\mathbb{Z}^d$, that is why we write an inequality taking explicitly into account this specificity. For $E$ a set of bonds, we define

$$S(E) \,=\, \frac{1}{p}\Big|\big\{\,e\in E: e\text{ is open}\,\big\}\Big| - \frac{1}{1-p}\Big|\big\{\,e\in E: e\text{ is closed}\,\big\}\Big|\,.$$

Let $(X_k)_{k\geq 1}$ be the sequence of i.i.d. Bernoulli random variables with parameter $p$ which is driving the exploration algorithm and let $(S_k)_{k\geq 1}$ be the associated partial sums:

$$\forall k\geq 1\qquad S_k=X_1+\cdots+X_k\,.$$

Let us fix a vertex $v$ and $i\geq 1$. Suppose that $|C(v)|=i$ and let $k$ be such that $1\leq k\leq i$. Denoting by $O(k)$ (respectively $C(k)$) the set of open (respectively closed) bonds which have been explored during the first $k$ steps of the algorithm, we have the identity (see (13.3) in section 13 and subsection 28.1)

$$S_k-kp \,=\, p(1-p)\Big(\frac{\big|O(k)\big|}{p}-\frac{\big|C(k)\big|}{1-p}\Big)\,. \tag{36.1}$$

Applying the identity (36.1) at step $b(i,n)$, we conclude that

$$S_{b(i,n)}-b(i,n)p \,=\, p(1-p)S(\mathcal{E}(v))\,.$$

Therefore, for any $t\geq 0$, we have

$$P\Big(\big|C(v)\big|=i\,,\big|S(\mathcal{E}(v))\big|\geq t\Big) \,\leq\, P\Big(\big|S_{b(i,n)}-b(i,n)p\big|\geq p(1-p)t\Big)\,. \tag{36.2}$$

For the percolation model in $\mathbb{Z}^d$, we relied extensively on Hoeffding's inequality (13.2). For the Erdős-Rényi model, as the parameter $p$ goes to 0, we need a more refined inequality to control the deviations of the Binomial distribution. A possibility would be to start from the classical Chernoff inequality, and to compute an adequate lower bound on the rate functions in the vicinity of the mean. We choose rather to use a ready-made inequality of the literature, namely the Bernstein inequality, that is recalled next.

**Bernstein's inequality.** For any $n \geq 1$ and any $p \in [0,1]$, we have

$$\forall t > 0 \qquad P_p\big(|S_n - np| \geq t\big) \;\leq\; 2\exp\left(-\frac{t^2}{2np(1-p)+\frac{2}{3}t}\right)\,. \tag{36.3}$$

This inequality holds for more general distributions, the general statement can be found in section 2.8 of [19]. The factor 2 in (36.3) is due to the fact that we have stated the two-sided form of the inequality, which controls $|S_n - np|$ and not only $S_n - np$.

We apply Bernstein's inequality to the right-hand quantity in (36.2), and we obtain the following control on the functional $S$ for the clusters in the Erdős-Rényi model:

$$\forall n \geq 1 \quad \forall p \in [0,1] \quad \forall v \in \{\,1,\dots,n\,\} \quad \forall i \in \{\,1,\dots,n\,\} \quad \forall t \geq 0$$
$$P\big(\big|C(v)\big| = i\,,\big|S(\mathcal{E}(v))\big| \geq t\big) \;\leq\; 2\exp\left(-\frac{p(1-p)t^2}{2b(i,n)+\frac{2}{3}t}\right)\,. \tag{36.4}$$

The essential gain compared with Hoeffding's inequality is that the factor $p$ could be simplified in the fraction in the exponential. Thus only one factor $p$ remains in the numerator. As $p$ is of order $c/n$ in the critical regime, this will make all the difference.

The difference between (36.4) and the inequality (28.2) is that we consider the event $|C(v)| = i$. On this event, the number of bonds visited during the exploration of $C(v)$ is equal to $b(i,n)$, hence the denominator in the exponential.

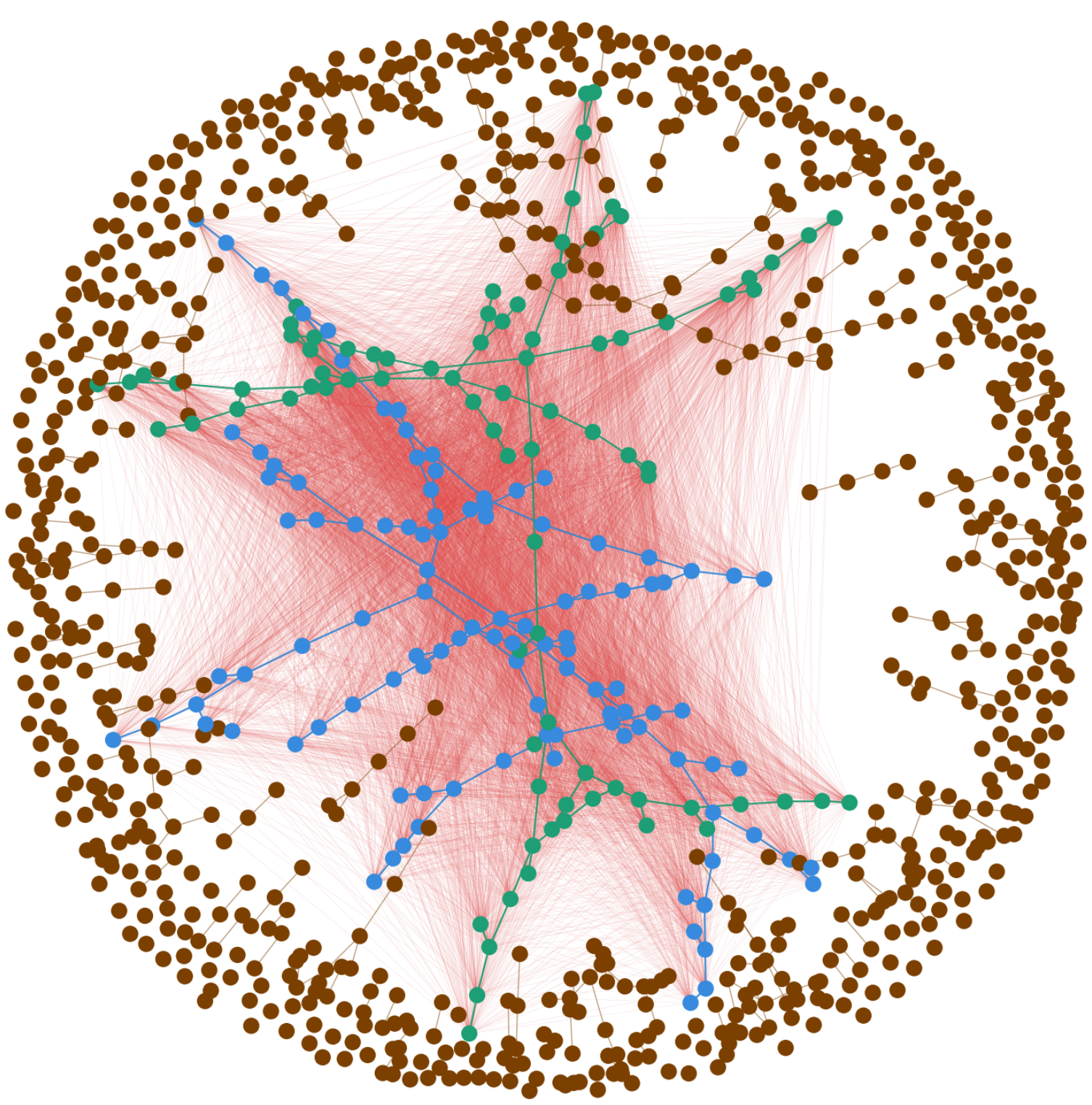

Figure 77: The Erdős-Rényi model with $n = 1000$, $p = 0.001$
The thin red lines are the closed bonds between the two largest components.
The largest component has size 109, the second largest 86.

# 37 Coexistence of two clusters

We would like to understand under which conditions two clusters can coexist in the random graph $G(n, c/n)$, in the critical $c = 1$ or supercritical $c > 1$ regimes. We already know that small clusters can coexist without trouble, and we would like to focus on non-small clusters, either of intermediate size or giant. The main obstacle that prevents the coexistence of two non-small clusters is that it requires the closure of many bonds, namely all the bonds linking the two clusters. Let us begin with two fixed clusters. More precisely, we fix two disjoint sets of vertices $V_1$ and $V_2$ in $\{1, \dots, n\}$, and we bound the probability that these two sets correspond to open clusters in $G(n, c/n)$ as follows:

$$P\begin{pmatrix}V_1 \text{ and } V_2 \text{ are two}\\ \text{two open clusters}\end{pmatrix} \leq P\begin{pmatrix}\text{all the bonds between}\\ V_1 \text{ and } V_2 \text{ are closed}\end{pmatrix}$$
$$\leq \Big(1 - \frac{c}{n}\Big)^{|V_1|\,|V_2|} \leq \exp\Big(-\frac{c}{n}|V_1|\,|V_2|\Big)\,.$$

This upper bound is very good and it goes very fast to 0 whenever we consider subsets $V_1$, $V_2$ of cardinality much larger than $\sqrt{n}$. Of course, there is an antagonistic combinatorial factor, corresponding to the number of possible choices for the sets $V_1$, $V_2$. Let us fix two integers $i_1, i_2$ and let us examine the probability of observing two clusters of sizes $i_1$ and $i_2$. By a standard union bound, we get

$$P\begin{pmatrix}\text{there exist two clusters}\\ \text{of sizes } i_1 \text{ and } i_2\end{pmatrix} \leq \sum_{\substack{V_1, V_2 \subset \{1,\dots,n\}\\ V_1\cap V_2=\varnothing, |V_1|=i_1, |V_2|=i_2}} P\begin{pmatrix}V_1 \text{ and } V_2 \text{ are two}\\ \text{two open clusters}\end{pmatrix}$$
$$\leq \binom{n}{i_1}\binom{n}{i_2}\exp\Big(-\frac{c}{n}i_1 i_2\Big)\,. \qquad (37.1)$$

This inequality is quite crude, but it can be used for $i_1, i_2$ close to $n$. To this end, we shall bound the binomial coefficients. We start with the inequality

$$\forall i \in \{0, \dots, n\} \quad \forall x \in [0, 1] \qquad \binom{n}{i} x^i (1-x)^{n-i} \leq 1\,.$$

Taking $x = i/n$ in the above formula, we obtain the following classical inequality:

$$\forall i \in \{0, \dots, n\} \qquad \binom{n}{i} \leq \exp\Big(-i\ln\frac{i}{n} - (n-i)\ln\frac{n-i}{n}\Big)\,. \qquad (37.2)$$

Let us fix $\alpha$ in $[1/2, 1]$. The map $t \in [0,1] \mapsto t\ln t + (1-t)\ln(1-t)$ is increasing on $[1/2, 1]$, thus

$$\forall t \in [\alpha, 1] \qquad -t\ln t - (1-t)\ln(1-t) \leq -\alpha\ln\alpha - (1-\alpha)\ln(1-\alpha)\,. \quad (37.3)$$

Using (37.2) and (37.3) together, we have

$$\forall i \in \{0, \dots, n\} \qquad i \geq \alpha n \implies \binom{n}{i} \leq \exp\Big(-n\big(\alpha\ln\alpha + (1-\alpha)\ln(1-\alpha)\big)\Big)\,. \qquad (37.4)$$

We consider two integers $i_1, i_2$ that are larger than $\alpha n$. Substituting (37.4) into (37.1), we get

$$P\begin{pmatrix}\text{there exist two clusters}\\ \text{of sizes larger than } \alpha n\end{pmatrix} \leq \exp\Big(-n\big(2\alpha\ln\alpha+2(1-\alpha)\ln(1-\alpha)+c\alpha^2\big)\Big)\,.$$

The quantity $2\alpha\ln\alpha + 2(1-\alpha)\ln(1-\alpha) + c\alpha^2$ converges towards $c$ as $\alpha$ goes to 1. Thus it becomes positive for $\alpha$ close to 1, and the probability on the left decreases exponentially fast with $n$:

$$\forall c > 0 \quad \exists\, \alpha < 1 \quad \exists\, \beta > 0 \qquad P\begin{pmatrix}\text{there exist two clusters}\\ \text{of sizes larger than } \alpha n\end{pmatrix} \leq \exp\big(-\beta n\big)\,.$$

This result is already interesting, but it is not strong enough to prove the uniqueness of the giant: it only rules out the coexistence of two giants of size larger than $\alpha n$, for $\alpha$ close to 1. The powerful technique employed in [7] and [133] is the sprinkling. The idea is to consider first a value of $p = c/n$ for which at least one giant exists with high probability. Supposing then that two giants are present, one increases slightly the probability $p$ to $p'$. The value $p'$ is chosen to ensure that, with overwhelming probability, at least one of the bonds between the two giants, which was closed at level $p$, becomes open at level $p'$. This forces all the giants which were present at the level $p$ to merge into a single giant at level $p'$. There is an additional subtlety there: the values $p$ and $p'$ must be close enough to ensure that no new giant is created when going from $p$ to $p'$, or alternatively that all the giants present at level $p'$ were already present at level $p$. This is the technique used in [7], section 11.9 (see the paragraph entitled 'Sprinkling'). A tricky variation is presented in [76], section 6.3: the sprinkling is performed from a configuration which contains a lot of clusters of sizes larger than a suitably chosen constant $A$. In contrast, the approach in [133] does not rely on sprinkling, but rather on a subtle second moment estimate (see proposition 4.10 in [133]).

**The strategy.** Our aim here is to present an alternative way to prove the uniqueness of the giant in the supercritical regime. Before doing the computations, let us explain briefly the idea. The coexistence of two giants necessitates the closure of all the bonds between them. We wish to avoid conditioning on the whole set of the vertices of each giant, because we have then to deal with huge combinatorial factors, as in (37.1). Instead, we pick up a vertex from each giant, and we perform a genuine exploration of the clusters. On one hand, each closed bond between the giants will be visited during the exploration of each giant, and these bonds are the only bonds that are visited by both these explorations. On the other hand, if we perform a single genuine exploration starting from the two vertices, we will explore the same set of bonds as the two previous explorations together, with the crucial difference that the bonds between the two giants will be visited only once. We compare then the value of the function $S$ associated to each of these three explorations. With the help of the quantitative control (36.4) on the functional $S$, we derive a valuable bound on the size of the interface between the two giants. In fact, this is exactly what we did to control the first intersection set for the bond model in $\mathbb{Z}^d$ in subsection 31.1.

**The computations.** Let us fix two integers $i_1, i_2$ and let us examine again the probability of observing two clusters of sizes $i_1$ and $i_2$. This time, however, we do not condition on the clusters themselves, rather we pick one vertex in each of them. This leads to the following bound:

$$P\begin{pmatrix}\text{there exist two clusters}\\ \text{of sizes } i_1 \text{ and } i_2\end{pmatrix} = P\begin{pmatrix}\exists v_1, v_2 \in \{1,\dots,n\} \quad v_1 \not\longleftrightarrow v_2\\ |C(v_1)| = i_1\,, |C(v_2)| = i_2\end{pmatrix}$$
$$\leq \sum_{\substack{1\leq v_1,v_2\leq n\\ v_1\neq v2}} P\Big(|C(v_1)| = i_1,\ |C(v_2)| = i_2,\ v_1 \not\longleftrightarrow v_2\Big)\,.$$

Since the Erdős-Rényi model is a mean-field model, the probability inside the sum above does not depend on the choices of the vertices $v_1, v_2$. Thus we obtain

$$P\begin{pmatrix}\text{there exist two clusters}\\ \text{of sizes } i_1 \text{ and } i_2\end{pmatrix} \leq n^2 P\Big(|C(1)| = i_1,\ |C(2)| = i_2,\ v_1 \not\longleftrightarrow v_2\Big)\,. \tag{37.5}$$

To control this last probability, we make appeal to the genuine exploration algorithms of the clusters. Let us denote by $\mathcal{D}(i_1, i_2)$ the event

$$\mathcal{D}(i_1, i_2) = \Big\{\, |C(1)| = i_1,\ |C(2)| = i_2,\ v_1 \not\longleftrightarrow v_2 \,\Big\}\,.$$

Suppose that the event $\mathcal{D}(i_1, i_2)$ occurs. The genuine exploration of the cluster of 1 (respectively 2) visits the set $\mathcal{E}(1)$ (respectively $\mathcal{E}(2)$) of all the bonds having at least one endpoint in $C(1)$ (respectively $C(2)$). The intersection between $\mathcal{E}(1)$ and $\mathcal{E}(2)$ consists of all the bonds between a vertex of $C(1)$ and a vertex of $C(2)$. Since 1 and 2 are not connected, all these bonds must be closed. Let us evaluate the functional $S$ on $\mathcal{E}(1) \cap \mathcal{E}(2)$. We get

$$S\big(\mathcal{E}(1) \cap \mathcal{E}(2)\big) = -\frac{1}{1-p}\big|\mathcal{E}(1) \cap \mathcal{E}(2)\big| = -\frac{i_1 i_2}{1-p}\,. \tag{37.6}$$

The map $S$ is additive, therefore we have also

$$S\big(\mathcal{E}(1) \cap \mathcal{E}(2)\big) = S\big(\mathcal{E}(1)\big) + S\big(\mathcal{E}(2)\big) - S\big(\mathcal{E}(1) \cup \mathcal{E}(2)\big)\,. \tag{37.7}$$

By the equality (37.6), the left-hand quantity in (37.7) is equal to $-i_1 i_2/(1-p)$, therefore at least one of the three terms of the right-hand side must be smaller than or equal to $-i_1 i_2/(3(1-p))$. From the preceding discussion, we conclude that

$$P\big(\mathcal{D}(i_1, i_2)\big) \leq P\Big(|C(1)| = i_1\,, S\big(\mathcal{E}(1)\big) \leq -\frac{i_1 i_2}{3(1-p)}\Big)$$
$$+ P\Big(|C(2)| = i_2\,, S\big(\mathcal{E}(2)\big) \leq -\frac{i_1 i_2}{3(1-p)}\Big)$$
$$+ P\Big(\mathcal{D}(i_1, i_2)\,, S\big(\mathcal{E}(1) \cup \mathcal{E}(2)\big) \geq \frac{i_1 i_2}{3(1-p)}\Big)\,. \tag{37.8}$$

To control the first two terms of the right-hand side of (37.8), we use the inequality (36.4). For the third term, we need to develop a specific variant of this inequality. To this end, we run an exploration algorithm which explores simultaneously the clusters of 1 and 2. This is essentially a standard exploration algorithm whose starting set is the pair $\{1,2\}$. On the event $\mathcal{D}(i_1,i_2)$, upon termination, this algorithm will have explored the set $\mathcal{E}(1)\cup\mathcal{E}(2)$ of all the bonds incident to the vertices of $C(1)$ and $C(2)$, of which there are

$$(i_1+i_2)\big(n-(i_1+i_2)\big)+\frac{1}{2}(i_1+i_2)(i_1+i_2-1)\,=\,b(i_1+i_2,n)\,.$$

The counterpart of (36.2) for the exploration starting from the pair $\{1,2\}$ is

$$\forall t\geq 0\qquad P\big(\mathcal{D}(i_1,i_2)\,,\big|S\big(\mathcal{E}(1)\cup\mathcal{E}(2)\big)\big|\geq t\big)\,\leq\, P\big(\big|S_{b(i_1+i_2,n)}-b(i_1+i_2,n)p\big|\geq p(1-p)t\big)\,,$$

and the application of Bernstein's inequality yields

$$\forall t\geq 0\quad P\big(\mathcal{D}(i_1,i_2)\,,\big|S\big(\mathcal{E}(1)\cup\mathcal{E}(2)\big)\big|\geq t\big)\,\leq\,2\exp\left(-\frac{p(1-p)t^2}{2b(i_1+i_2,n)+\frac{2}{3}t}\right)\,.\tag{37.9}$$

We apply now the inequalities (36.4) and (37.9) with $t=i_1i_2/(3(1-p))$, and we deduce from (37.8) that

$$P\big(\mathcal{D}(i_1,i_2)\big)\,\leq\,2\exp\left(-\frac{pi_1^2i_2^2}{18b(i_1,n)(1-p)+2i_1i_2}\right)+$$
$$2\exp\left(-\frac{pi_1^2i_2^2}{18b(i_2,n)(1-p)+2i_1i_2}\right)+2\exp\left(-\frac{pi_1^2i_2^2}{18b(i_1+i_2,n)(1-p)2i_1i_2}\right)\,.\tag{37.10}$$

We can get rid of the factor $1-p$ in the denominator of the fraction in the three exponentials. Furthermore, we have

$$\begin{aligned}\forall i\in\{1,\dots,n\}\qquad & b(i,n)\leq in\,,\\ \forall i_1,i_2\in\{1,\dots,n\}\qquad & i_1i_2\,\leq\,(i_1+i_2)n\,.\end{aligned}\tag{37.11}$$

We deduce from (37.10) and (37.11) that

$$P\big(\mathcal{D}(i_1,i_2)\big)\,\leq\,6\exp\left(-\frac{pi_1^2i_2^2}{20(i_1+i_2)n}\right)\,.\tag{37.12}$$

Substituting (37.12) into (37.5), we obtain the next theorem.

**Theorem 37.1.** *In the Erdős-Rényi model, we have the following upper bound for the simultaneous presence of two clusters of prescribed sizes:*

$$\forall n\geq 1\quad\forall p\in[0,1]\quad\forall i_1,i_2\in\{1,\dots,n\}$$
$$P\begin{pmatrix}\text{there exist two clusters}\\ \text{of sizes } i_1 \text{ and } i_2\end{pmatrix}\,\leq\,6n^2\exp\Big(-\frac{pi_1^2i_2^2}{20(i_1+i_2)n}\Big)\,.\tag{37.13}$$

It is remarkable that this simple bound is valid for all the values of $n,p,i_1,i_2$.

# 38 Applications

To conclude this part, we present here two applications of the inequality (37.13). The first one concerns the proof of the uniqueness of the giant. We show how the sprinkling argument can be bypassed with the help of (37.13). In the second one, we present a little variant of the heuristic of Noga Alon and Tom Spencer [7] explaining the parametrization of the critical window.

**The uniqueness of the giant.** As we saw in section 35, the proof of the second part of theorem 34.1 was done in three steps:

- Control of the small clusters;
- No middle ground;
- Existence and uniqueness of the giant.

We show that our inequality provides an alternative way of proving the last part, namely the uniqueness of the giant. So, suppose that the first two steps have been accomplished. At this point, we know that there is a density close to $\zeta(c)$ of sites belonging to giant clusters, and moreover $\zeta(c) > 0$ whenever $c > 1$. Recall that a giant cluster is a cluster having cardinality larger than $bn$, where $b$ is a positive constant. By a simple union bound, we have

$$P\big(\text{there exist two giants}\big) \;\leq\; \sum_{i_1,i_2\geq bn} P\begin{pmatrix}\text{there exist two clusters}\\ \text{of sizes } i_1 \text{ and } i_2\end{pmatrix}.$$

We apply next the inequality (37.13) and we obtain

$$P\big(\text{there exist two giants}\big) \;\leq\; 6n^4 \exp\Big(-\frac{1}{40}pb^3n^2\Big).$$

For $p$ of the form $p = c/n$, with $c > 0$, this probability goes exponentially fast to 0. Thus, with overwhelming probability, the giant is unique. Notice that the proofs of uniqueness presented in [7] and [133] relied on a sprinkling technique. Therefore, we see that the argument above bypasses the sprinkling. We do not claim that the above argument is simpler than a sprinkling procedure, because the derivation of the inequality itself requires some additional work. However, what we specifically want to do within the framework of the percolation model in $\mathbb{Z}^d$ is to avoid resorting to a sprinkling method. That we can do so in the context of the Erdős-Rényi model is a good sign that we are on the right track! We remark also that the approach in [133] does not rely on sprinkling, however the argument there seems to be very specific to the Erdős-Rényi model. Indeed, it rests on an estimate for the variance of the number of sites belonging to giants (see proposition 4.10 in [133]). The computation of this estimate involves a conditioning on the size $\ell$ of one cluster, and it can be conducted because the conditional distribution of the remaining graph in $G(n, c/n)$ is the same as $G(n-\ell, c/n)$. Unfortunately, there is little hope of adapting this computation to the percolation model in $\mathbb{Z}^d$. By the way, this question is similar to the one we have struggled with time and again throughout part V. We were unable to obtain adequate estimates on the number of pivots of higher order because we could not grasp the conditional distribution of the percolation configuration given the set of the first pivots.

**The exponent** $1/3$. To say that the Erdős-Rényi model presents a critical point is a bit oversimplistic. Indeed, the finiteness of the graph $G(n,p)$ makes the phase transition harder to apprehend, because, as long as $n$ is finite, every observable of the system is analytic. In fact, the giant component appears progressively, and rather than a single critical point, there is a critical window, which is the range of parameters $p$ of the form

$$p = \frac{1}{n}\Big(1 + \frac{\lambda}{n^{1/3}}\Big), \qquad \lambda \in \mathbb{R}. \tag{38.1}$$

It is known that, within the critical window, with probability going to 1 as $n$ goes to $\infty$, the largest cluster has size of order $n^{2/3}$. In the section 11.7 of [7], Noga Alon and Tom Spencer write that the question most frequently asked about the Erdős-Rényi model is: why should the exponent be $1/3$ in the parametrization (38.1)?

They provide also a nice heuristic explaining this exponent, based on the Poisson branching process. They argue as follows. Let us take $p = (1+\varepsilon)/n$, where $\varepsilon$ is positive and small. For this parameter, the largest component in $G(n,p)$ should typically have a size of order $2\varepsilon n$. The small components should be well approximated by the Poisson branching process of parameter $1+\varepsilon$. For this process, the expected size of the progeny is of order $1/\varepsilon^2$. This picture can remain valid as long as these small components are much smaller than the largest one, that is as long as $1/\varepsilon^2 << 2\varepsilon n$. This readily gives that the limiting regime occurs when $\varepsilon$ is of order $1/n^{1/3}$. In this heuristic, the exponent $1/3$ corresponds to a threshold beyond which the approximation of a cluster in $G(n,p)$ by the progeny of a Poisson branching process is not any more valid.

We can propose a variant of this heuristic, based on our inequality (37.13). Here is how it works. For $p$ small, the graph $G(n,p)$ is constituted of a large number of small components. This collection of small components is well approximated by sampling a branching Poisson process of parameter $np$, until the total sum of the progenies of the sample is roughly equal to $n$. As long as this picture holds, the sizes of the two largest components are of the same order. Denoting these sizes by $i_1$ and $i_2$, the inequality (37.13) implies that $i_1^2 i_2 << n$. The maximum sizes which are compatible with this picture are of order $n^{2/3}$. Let $p = (1+\varepsilon)/n$, where $\varepsilon$ is positive and small. The expected size of the progeny of the Poisson branching process of parameter $1+\varepsilon$ is of order $1/\varepsilon^2$. Thus the branching Poisson approximation breaks down when $\varepsilon$ is of order $1/n^{1/3}$.

The difference between the two heuristics is that the second one does not require the estimate on the size of the largest component in $G(n,p)$. We could certainly analyze other aspects of the Erdős-Rényi model with the help of the inequality (37.13). For instance, we could try to study rigorously some parts of the critical window. Indeed, Noga Alon and Tom Spencer [7] show that the giant remains unique in what they call the barely supercritical phase. They use again a sprinkling argument, albeit more technical, and most likely this sprinkling can be bypassed thanks to (37.13). We could also work out an extension of the inequality (37.13) involving $N > 2$ disjoint clusters, but that would take us too far off track and distract us from our main goal.

# Part VII
# Higher pivots

Our major problem is to control the intersection set $\mathcal{I}(k)$ for $k \geq 2$. In lemma 29.5, we saw that the intersection set $\mathcal{I}(k)$ is included in the set $\mathcal{P}_k(\mathcal{D})$ of the $k$-pivotal bonds for the disconnection event $\mathcal{D}$. Therefore, a natural attempt to control the cardinality of $\mathcal{I}(k)$ would be to try to control the cardinality of $\mathcal{P}_k(\mathcal{D})$. This is a quite promising strategy, at least the first step is successful. Indeed, there is a well-known upper bound on the expected number of pivots for a monotone event, which says that, for any monotone event $\mathcal{D}$ depending on the bonds of $\mathbb{E}^d(D)$, we have

$$E\big(|\mathcal{P}(\mathcal{D})|\big) \;\leq\; \sqrt{\frac{\big|\mathbb{E}^d(D)\big|}{p(1-p)}}\,. \tag{V.2}$$

So we will investigate whether similar bounds can be derived for the pivotal edges of higher orders. This is the purpose of this part. This question is extremely interesting in itself, and can conveniently be addressed in a more general framework, without reference to the percolation model.

The classical proofs of the bound (V.2) rest either on Parseval's formula or on the Margulis-Russo formula. This last approach is more amenable to generalizations, so we start by revisiting the well-known Margulis-Russo formula in section 39. In section 40, we prove an exponential inequality linking the probability of an event $A$ and the expected number of its pivots. Afterwards, we generalize the Margulis Russo formula to higher derivatives in section 41. The ensuing formula involves the Bernoulli chaos, which are studied in section 43. By computing the variance of this chaos in subsection 43.1, we obtain the desired generalization of (V.2): for any event $A$, any $\ell \leq n$, we have

$$\Big|\Big(\frac{d}{dp}\Big)^{\ell} P_p(A)\Big| \;\leq\; \ell!\sqrt{P_p(A)\binom{n}{\ell}\frac{1}{p^{\ell}(1-p)^{\ell}}}\,.$$

We then investigate further the Bernoulli chaos and their relationship with the Krawtchouk polynomials in subsections 43.2 and 43.3. After these general computations, we come back to the main question. In section 44, we give the definition of the higher pivots in the general framework and we explore how they are related to the general inequality. In section 45, we examine the specific case of the second pivots. We apply the general results on higher pivots and we discuss an important inequality of Talagrand, Keller and Kindler related to them. We extend the control of the pivots to a random domain, with the hope that it can help to control the second pivots. Unfortunately, these inequalities do not seem powerful enough, so we try to formulate the kind of inequalities that we would need. We present some pictures of the second positive pivots which are computed with the help of Tarjan's algorithm. In the final section 46 of this part, we prove a generalization of two inequalities of Talagrand, Benjamini, Kalai, Schramm, Keller and Kindler, which involve the higher pivots.

# 39 The classical Margulis-Russo formula $\odot$

Let $n \geq 1$ be an integer and let us consider the configuration space $\Omega = \{0,1\}^n$. For $\omega \in \Omega$, we denote by $\omega(i)$ its $i$-th component, so that $\omega = \big(\omega(1), \dots, \omega(n)\big)$. Let $p \in ]0,1[$. We endow $\Omega$ with the product probability measure $P_p$ defined by

$$\forall \omega \in \Omega \qquad P_p(\omega) \,=\, \prod_{i=1}^n p^{\omega(i)}(1-p)^{1-\omega(i)} \,.$$

For $\omega$ a configuration in $\Omega$ and $i \in \{1, \dots, n\}$, we denote by $\omega^i$ (respectively $\omega_i$) the configuration $\omega$ where the $i$-th component is set to 1 (respectively to 0), that is,

$$\forall j \in \{1, \dots, n\}$$

$$\omega^i(j) \,=\, \begin{cases} \omega(j) & \text{if } j \neq i \\ 1 & \text{if } j = i \end{cases}, \qquad \omega_i(j) \,=\, \begin{cases} \omega(j) & \text{if } j \neq i \\ 0 & \text{if } j = i \end{cases}.$$

A component $i$ is said to be pivotal for the event $A$ in the configuration $\omega$ if

$$\text{either} \quad \omega^i \in A,\, \omega_i \notin A \quad \text{or} \quad \omega^i \notin A,\, \omega_i \in A \,.$$

We denote by $\mathcal{P}(A, \omega)$ the set of the pivotal components for $A$ in the configuration $\omega$. The set $\mathcal{P}(A, \omega)$ is a random subset of $\{1, \dots, n\}$. When $\omega$ is hidden in the formula, we write simply $\mathcal{P}(A)$ instead of $\mathcal{P}(A, \omega)$.

**Theorem 39.1** (Margulis-Russo formula)**.** *For $A$ an increasing event, we have*

$$\frac{d}{dp} P_p(A) \,=\, E\big(|\mathcal{P}(A)|\big) \,. \tag{39.1}$$

This formula has played a central role in percolation theory (see for instance [63], section 2.4) and it has numerous applications. There exists also a version for a general event $A$ (not necessarily monotone), which reads

$$\frac{d}{dp} P_p(A) \,=\, \sum_{i=1}^n E\big(\delta_i 1_A\big) \,, \tag{39.2}$$

where for $i \in \{1, \dots, n\}$, the random variable $\delta_i 1_A$ is defined by

$$\forall \omega \in \Omega \qquad \delta_i 1_A(\omega) \,=\, 1_A(\omega^i) - 1_A(\omega_i) \,.$$

The above versions (39.1) and (39.2) are very convenient to study the response of the probability $P_p(A)$ to a slight variation of the parameter $p$, for instance to analyze threshold phenomena. Yet we wish to understand the typical configurations at parameter $p$. We will therefore rewrite the Margulis-Russo formula in a non-differential form, which, as we will see, has a considerable interest on its own. We want also to deal with non-monotone events, so we distinguish between negative and positive pivotal components.

**Definition 39.2.** *A component $i$ is said be to positive pivotal (respectively negative pivotal) for the event $A$ in the configuration $\omega$ if $\omega^i \in A$ and $\omega_i \notin A$ (respectively $\omega^i \notin A$ and $\omega_i \in A$).*

We denote by $\mathcal{P}^+(A,\omega)$ (respectively $\mathcal{P}^-(A,\omega)$) the set of the positive pivotal (respectively negative pivotal) components for $A$ in the configuration $\omega$. Finally, we set, for $i \in \{1,\dots,n\}$,

$$\forall \omega \in \Omega \qquad X_i(\omega) \,=\, \frac{\omega(i)}{p} - \frac{1-\omega(i)}{1-p}\,,$$

and we denote by $S_n(\omega)$ the random sum

$$S_n(\omega) \,=\, \sum_{i=1}^n X_i(\omega)\,.$$

With this notation, we can restate the Margulis-Russo formula as follows.

**Proposition 39.3.** *For any event $A \subset \Omega$, we have*

$$\frac{d}{dp} P_p(A) \,=\, E\big(|\mathcal{P}^+(A)|\big) - E\big(|\mathcal{P}^-(A)|\big) \,=\, E\big(1_A S_n\big)\,. \tag{39.3}$$

*Proof.* We take the derivative of $P_p(A)$ with respect to the parameter $p$:

$$\frac{d}{dp} P_p(A) \,=\, \frac{d}{dp} \sum_{\omega\in\Omega} 1_A(\omega)\, P_p(\omega) \,=\, \sum_{\omega\in\Omega} 1_A(\omega)\, \frac{d}{dp} P_p(\omega)\,.$$

We differentiate $P_p(\omega)$ either as an $n$-fold product, or through its logarithmic derivative. Indeed, we have

$$\forall \omega \in \Omega \qquad \ln P_p(\omega) \,=\, \sum_{i=1}^n \Big(\omega(i)\ln p + \big(1-\omega(i)\big)\ln(1-p)\Big)\,.$$

Whichever way, we get, for any $\omega \in \Omega$,

$$\frac{d}{dp} P_p(\omega) \,=\, \sum_{i=1}^n X_i(\omega)\, P_p(\omega) \,=\, S_n(\omega)\, P_p(\omega)\,, \tag{39.4}$$

whence

$$\frac{d}{dp} P_p(A) \,=\, E\big(1_A S_n\big)\,,$$

and the final part of formula (39.3) is proved. In the next step, we take the sum involved in $S_n$ outside the expectation:

$$E\big(1_A S_n\big) \,=\, E\Big(1_A \sum_{i=1}^n X_i(\omega)\Big) \,=\, \sum_{i=1}^n E\big(1_A X_i\big)\,. \tag{39.5}$$

We compute next $E\big(1_A X_i\big)$, in fact we will compute $E\big(f X_i\big)$ for an arbitrary function $f$. For $f$ a function from $\Omega$ to $\mathbb{R}$ and $i \in \{1,\dots,n\}$, we define the discrete derivative along the $i$-th component as

$$\forall \omega \in \Omega \qquad \delta_i f(\omega) \,=\, f(\omega^i) - f(\omega_i)\,.$$

**Lemma 39.4.** *For any function $f$ from $\Omega$ to $\mathbb{R}$ and for any $i \in \{1, \dots, n\}$, we have*

$$E\big(fX_i\big) \;=\; E\big(\delta_i f\big)\,. \tag{39.6}$$

*Proof.* To compute $E\big(fX_i\big)$, we take the conditional expectation knowing $\omega(1)$, $\dots$, $\omega(i-1)$, $\omega(i+1), \dots, \omega(n)$. Equivalently, we perform the average with respect to $\omega(i)$ and we get

$$E\big(fX_i \,|\, \omega(1), \dots, \omega(i-1), \omega(i+1), \dots, \omega(n)\big) \;= \\ \delta_i f\big(\omega(1), \dots, \omega(i-1), \cdot, \omega(i+1), \dots, \omega(n)\big)\,. \tag{39.7}$$

Naturally, the resulting function does not depend any more on the variable $\omega(i)$. Taking the expectation of formula (39.7), we obtain formula (39.6). □

We conclude now the proof of proposition 39.3. Coming back to formula (39.5), we have, thanks to formula (39.6),

$$\begin{aligned} E\big(1_A S_n\big) \;&=\; \sum_{i=1}^n E\big(1_A X_i\big) \;=\; \sum_{i=1}^n E\big(\delta_i 1_A\big) \;=\; E\Big(\sum_{i=1}^n \delta_i 1_A\Big) \\ &=\; E\big(|\mathcal{P}^+(A)|\big) - E\big(|\mathcal{P}^-(A)|\big) \end{aligned}$$

as requested. □

Of course, there is nothing fundamentally new in the previous computations. We will next forget about the derivative in Russo's formula, and we will try to exploit the term $E\big(1_A S_n\big)$. A straightforward application of Cauchy-Schwarz inequality gives

$$\big|E\big(1_A S_n\big)\big| \;\leq\; \sqrt{P_p(A)\, E\big((S_n)^2\big)}\,. \tag{39.8}$$

Under $P_p$, the random variable $S_n$ is the sum of $n$ i.i.d. centered variables, thus

$$E\big((S_n)^2\big) \;=\; \sum_{i=1}^n E\big((X_i)^2\big) \;=\; \frac{n}{p(1-p)}\,. \tag{39.9}$$

Putting together (39.8) and (39.9), we obtain the inequality stated in the next proposition.

**Proposition 39.5.** *For any event $A \subset \Omega$, we have*

$$\Big|E\big(|\mathcal{P}^+(A)|\big) - E\big(|\mathcal{P}^-(A)|\big)\Big| \;\leq\; \sqrt{\frac{nP_p(A)}{p(1-p)}}\,. \tag{39.10}$$

In the case of an increasing event $A$, the pivotal sites are automatically positive and formula (39.10) reduces to

$$E\big(|\mathcal{P}^+(A)|\big) \;\leq\; \sqrt{\frac{nP_p(A)}{p(1-p)}}\,. \tag{39.11}$$

This bound is well known. For the case $p = 1/2$, it is proved in Proposition IV.6 of [58] with the help of Parseval's identity. It is stated in Corollary 9.16 of [19] with the same proof as above.

# 40 Exponential control $\odot$

We will be interested in controlling the cardinality of the pivotal set whenever the event $A$ occurs. Our first step in this direction consists in developing a formula involving the conditional expectations $E\big(|\mathcal{P}^+(A)|\,\big|A\big)$ and $E\big(|\mathcal{P}^-(A)|\,\big|A\big)$. Whenever a coordinate $i$ belongs to $\mathcal{P}(A)$, its state decides the occurrence of $A$. Moreover the state of the coordinate $i$ itself is independent of the fact that $i$ is pivotal or not, therefore

$$\begin{aligned} P\big(\{\, i \in \mathcal{P}^+(A)\,\} \cap A\big) \;&=\; P\big(\{\, i \in \mathcal{P}^+(A)\,\} \cap \{\, i \text{ is open}\,\}\big) \\ &=\; P\big(\{\, i \in \mathcal{P}^+(A)\,\}\big)\, p\,, \end{aligned}$$

which we rewrite as

$$P\big(\, i \in \mathcal{P}^+(A)\,\big) \;=\; \frac{1}{p} P\big(\{\, i \in \mathcal{P}^+(A)\,\} \cap A\big)\,. \tag{40.1}$$

Similarly, we have

$$P\big(\, i \in \mathcal{P}^-(A)\,\big) \;=\; \frac{1}{1-p} P\big(\{\, i \in \mathcal{P}^-(A)\,\} \cap A\big)\,. \tag{40.2}$$

Summing equations (40.1) and (40.2) over $i$, and putting the summation inside the expectation, we obtain

$$E\big(|\mathcal{P}^+(A)|\big) - E\big(|\mathcal{P}^-(A)|\big) \;=\; \frac{1}{p} E\big(|\mathcal{P}^+(A)|\, 1_A\big) - \frac{1}{1-p} E\big(|\mathcal{P}^-(A)|\, 1_A\big)\,. \tag{40.3}$$

The identities (39.3) and (40.3) yield

$$E\big(1_A S_n\big) \;=\; \frac{1}{p} E\big(|\mathcal{P}^+(A)|\, 1_A\big) - \frac{1}{1-p} E\big(|\mathcal{P}^-(A)|\, 1_A\big)\,.$$

Dividing finally by $P_p(A)$, we obtain the formula stated in the next lemma.

**Lemma 40.1.** *For any event $A \subset \Omega$, we have*

$$E\big(S_n \,\big|\, A\big) \;=\; \frac{1}{p} E\big(|\mathcal{P}^+(A)|\,\big|\, A\big) - \frac{1}{1-p} E\big(|\mathcal{P}^-(A)|\,\big|\, A\big)\,. \tag{40.4}$$

## 40.1 An exponential inequality

Ultimately, we would like to link the size of the pivotal sets $\mathcal{P}^-(A)$ and $\mathcal{P}^+(A)$ to the probability of $A$. Our strategy to do so is the following. Suppose that the difference appearing in the right-hand side of (40.4) is too large in absolute value. This forces that the conditional law of $S_n$ knowing that $A$ occurs is concentrated on subsets of $\Omega$ where $S_n$ is large, but these sets have small probability, so the probability of $A$ itself has to be small. The next proposition presents an exponential inequality whose proof relies on the previous idea. Notice that this inequality for the symmetric case $p = 1/2$ appears at the beginning of the paper [125] of Talagrand (the inequality of Talagrand is restated in proposition 40.7). The proof presented below appears in [29], but the upper bound of the erfc function is missing there.

**Proposition 40.2.** *Let $n$ be a positive integer. For any event $A \subset \Omega$, we have*

$$P_p(A) \,\leq\, 2\exp\Big(-\frac{1}{2n}\Big((1-p)E\big(|\mathcal{P}^+(A)|\,\big|\,A\big)-pE\big(|\mathcal{P}^-(A)|\,\big|\,A\big)\Big)^2\Big)\,. \quad (40.5)$$

*Proof.* We start with

$$E\big(|S_n|\,\big|\,A\big) \,=\, \frac{1}{P_p(A)}\int_A |S_n|\,dP_p\,. \quad (40.6)$$

Using Fubini's theorem, we rewrite the integral in (40.6) as follows:

$$\int_A |S_n|\,dP_p \,=\, \int_A\Big(\int_0^{+\infty} 1_{u\leq|S_n|}\,du\Big)\,dP_p \,=\, \int_0^{+\infty} P_p\big(A,|S_n|\geq u\big)\,du\,.$$

Let $t>0$ and let us split this integral in two:

$$\int_0^{+\infty} P_p\big(A,|S_n|\geq u\big)\,du \,=\, \int_0^t\cdots+\int_t^{+\infty}\cdots \,\leq\, tP_p(A)+\int_t^{+\infty} P_p\big(|S_n|\geq u\big)\,du\,. \quad (40.7)$$

Reporting (40.7) in (40.6), we have

$$E\big(|S_n|\,\big|\,A\big) \,\leq\, \frac{1}{P_p(A)}\int_t^{+\infty} P_p\big(|S_n|\geq u\big)\,du\,+\,t\,. \quad (40.8)$$

Recall that $S_n$ is the random sum $S_n = X_1+\cdots+X_n$, where the variables $X_1,\dots,X_n$ are i.i.d. with distribution given by

$$\forall i\in\{\,1,\dots,n\,\}\qquad P\Big(X_i=\frac{1}{p}\Big) \,=\, p\,,\qquad P\Big(X_i=-\frac{1}{1-p}\Big) \,=\, 1-p\,.$$

The key tool to complete the previous program is the classical Hoeffding inequality, that we recall next.

**Hoeffding inequality.** For any $n\geq 1$ and any $p\in[0,1]$, we have

$$\forall u>0\qquad P_p\big(|S_n|\geq u\big) \,\leq\, 2\exp\Big(-\frac{2p^2(1-p)^2u^2}{n}\Big)\,. \quad (40.9)$$

In fact, a more general inequality is proved by Hoeffding in his paper [70]. A shorter proof of the specific inequality (40.9) can be found in proposition 8 of [29]. Reporting (40.9) into (40.8), we obtain

$$E\big(|S_n|\,\big|\,A\big) \,\leq\, \frac{2}{P_p(A)}\int_t^{+\infty}\exp\Big(-\frac{2p^2(1-p)^2u^2}{n}\Big)\,du\,+\,t\,. \quad (40.10)$$

We make the change of variable $2p(1-p)u=\sqrt{n}v$, so that (40.10) becomes

$$E\big(|S_n|\,\big|\,A\big) \,\leq\, \frac{\sqrt{n}}{p(1-p)P_p(A)}\int_{\frac{2p(1-p)t}{\sqrt{n}}}^{+\infty}\exp\Big(-\frac{v^2}{2}\Big)\,dv\,+\,t\,.$$

Recalling that this inequality holds for any $t > 0$, we replace $t$ by $\frac{\sqrt{n}t}{2p(1-p)}$, and we get

$$E\big(|S_n| \,\big|\, A\big) \,\leq\, \frac{\sqrt{n}}{p(1-p)} \bigg( \frac{1}{P_p(A)} \int_t^{+\infty} \exp\Big(-\frac{v^2}{2}\Big)\,dv \,+\, \frac{t}{2}\bigg). \tag{40.11}$$

The right-hand quantity admits a unique global minimum at

$$t^* \,=\, \sqrt{2\ln\frac{2}{P_p(A)}}\,.$$

We perform an integration by parts to obtain a simple upper bound on the erfc function: for any $t > 0$,

$$\begin{aligned}\int_t^{+\infty} \exp\Big(-\frac{v^2}{2}\Big)\,dv \,&=\, \int_t^{+\infty} \frac{1}{v} v \exp\Big(-\frac{v^2}{2}\Big)\,dv \\ &= \Big[-\frac{1}{v}\exp\Big(-\frac{v^2}{2}\Big)\Big]_t^{+\infty} - \int_t^{+\infty} \frac{1}{v^2}\exp\Big(-\frac{v^2}{2}\Big)\,dv \\ &\leq\, \frac{1}{t}\exp\Big(-\frac{t^2}{2}\Big). \end{aligned} \tag{40.12}$$

Applying inequality (40.12) at $t^*$, we obtain

$$\int_{t^*}^{+\infty} \exp\Big(-\frac{v^2}{2}\Big)\,dv \,\leq\, \frac{1}{t^*}\frac{P(A)}{2}\,. \tag{40.13}$$

Choosing $t = t^*$ in inequality (40.11) and using (40.13), we get

$$E\big(|S_n| \,\big|\, A\big) \,\leq\, \frac{\sqrt{n}}{2p(1-p)}\Big(\frac{1}{t^*} \,+\, t^*\Big) \,\leq\, \frac{\sqrt{n}}{p(1-p)} t^*\,,$$

where the last inequality holds thanks to the fact that $t^* \geq 1$. We conclude that

$$E\big(|S_n| \,\big|\, A\big) \,\leq\, \frac{1}{p(1-p)}\sqrt{2n\ln\frac{2}{P_p(A)}}\,.$$

We rewrite this inequality in exponential form as

$$P_p(A) \,\leq\, 2\exp\Big(-\frac{p^2(1-p)^2}{2n}\Big(E\big(|S_n| \,\big|\, A\big)\Big)^2\Big). \tag{40.14}$$

We have the standard inequality

$$\Big|E\big(S_n \,\big|\, A\big)\Big| \,\leq\, E\big(|S_n| \,\big|\, A\big)\,. \tag{40.15}$$

Furthermore, by lemma 40.1, we have

$$E\big(S_n \,\big|\, A\big) \,=\, \frac{1}{p}E\big(|\mathcal{P}^+(A)| \,\big|\, A\big) - \frac{1}{1-p}E\big(|\mathcal{P}^-(A)| \,\big|\, A\big)\,. \tag{40.16}$$

The inequality (40.5) of the proposition follows from the inequalities (40.14), (40.15) and (40.16). □

Before commenting further on proposition 40.2, we state an immediate corollary.

**Corollary 40.3.** *Let $n$ be a positive integer. For any increasing event $A \subset \Omega$, we have*

$$P_p(A) \le 2\exp\Big(-\frac{1}{2n}\Big((1-p)E\big(|\mathcal{P}^+(A)|\,\big|\,A\big)\Big)^2\Big). \tag{40.17}$$

In fact, the technique employed to bound the expectation $E\big(|S_n|\,\big|\,A\big)$ has been used routinely in various contexts by Talagrand. It is presented in his most recent book (see lemma 2.3.2 in [127]). Notably, Talagrand used this technique to derive his inequality on the correlation of increasing sets in his influential work [125]. However, only the very beginning of the proof of the correlation inequality involves this technique, i.e., the first step of the proof of Proposition 2.2 in [125]. This kind of upper bound is such a routine for Talagrand that the detailed proof of proposition 40.2 corresponds essentially to 6 lines in [125]! So, what is the point of the previous computations? A difference with [125] is that we deal here with the Bernoulli product measure of parameter $p$, whereas Talagrand was considering the uniform case $p = 1/2$. So the sub-Gaussian inequality of the uniform case is replaced here by Hoeffding's inequality. Talagrand did not care about the values of the constants in his inequalities, and it might be desirable to understand their dependence on the parameter $p$.

The inequality (40.17) tells us that, when $E\big(|\mathcal{P}^+(A)|\,\big|\,A\big)$ is very large, then $P_p(A)$ is very small. A natural question which arises is the following: if $A$ is an increasing event having a small probability, can we say something on its pivotal set $\mathcal{P}(A)$? In the fundamental paper [118], Russo has proved the following beautiful inequality (see [118], lemma 3): for any event $A \subset \Omega$,

$$E\big(|\mathcal{P}(A)|\,\big|\,A\big) \ge \frac{1}{\min(\ln p, \ln(1-p))}\ln P_p(A)\,. \tag{40.18}$$

This inequality gives a lower bound on the probability $P_p(A)$ in terms of the average size of the pivotal set $\mathcal{P}(A)$. It tells also that, whenever $P_p(A)$ is small, then $E\big(|\mathcal{P}(A)|\,\big|\,A\big)$ has to be big. However, the forms of the two inequalities (40.17) and (40.18) are quite different. Furthermore, Russo has also proved an inequality going in the reverse direction, that we restate next.

**Lemma 40.4** (Lemma 5 of [118])**.** *Let $S$ be an event of $\{0,1\}^n$ which is of the form $S = A \setminus B$, where $A$ and $B$ are two increasing subsets of $\{0,1\}^n$. If $0 < \alpha < \beta < 1$, then*

$$\int_\alpha^\beta P_p(S)\,dp \le 2\Big[\inf_{p\in[\alpha,\beta]} E_p\Big(\big|\mathcal{P}(S)\big|\,\Big|\,S\Big)\Big]^{-1}. \tag{40.19}$$

The inequality (40.19) yields an upper bound on the average size of the pivotal set $\mathcal{P}(A)$, but its drawback is that it involves an integral over the parameter $p$. Russo says also that the converse inequality of (40.18) does not hold. The inequalities (40.5) and (40.17) give upper bounds on $P_p(A)$, but they become effective only when the pivotal sets have a size larger than $\sqrt{n}$.

## 40.2 Implementing the second step of Talagrand

Let us proceed a bit further along the computations of Talagrand. We implement now the second step of Proposition 2.2 in [125] in our context, and we see where it leads. So, for this second step, we work with the random variable $\widetilde{S}_n$ given by

$$\widetilde{S}_n \,=\, \sum_{i=1}^n \alpha_i X_i \,, \tag{40.20}$$

where the $\alpha_i$'s are real numbers that will be chosen later. Let $A$ be any event. On the one hand, we have

$$E(\widetilde{S}_n 1_A) \,=\, \sum_{i=1}^n \alpha_i E(\delta_i 1_A) \,=\, \sum_{i=1}^n \alpha_i \Big( P\big( i \in \mathcal{P}^+(A) \big) - P\big( i \in \mathcal{P}^-(A) \big) \Big) \,. \tag{40.21}$$

On the other hand, we have, by the Cauchy-Schwarz inequality,

$$\big| E\big(\widetilde{S}_n 1_A\big)\big| \,\leq\, \sqrt{P_p(A)\, E\big((\widetilde{S}_n)^2\big)} \,. \tag{40.22}$$

Under $P_p$, the random variable $\widetilde{S}_n$ is the sum of $n$ independent centered variables, so that

$$E\big((\widetilde{S}_n)^2\big) \,=\, \sum_{i=1}^n (\alpha_i)^2 E\big((X_i)^2\big) \,=\, \frac{1}{p(1-p)} \sum_{i=1}^n (\alpha_i)^2 \,. \tag{40.23}$$

Putting together (40.22) and (40.23), we obtain the following inequality:

$$\Big| \sum_{i=1}^n \alpha_i \Big( P\big( i \in \mathcal{P}^+(A) \big) - P\big( i \in \mathcal{P}^-(A) \big) \Big) \Big| \,\leq\, \sqrt{\frac{P_p(A)}{p(1-p)} \sum_{i=1}^n (\alpha_i)^2} \,.$$

The optimal choice for the $\alpha_i$'s is to take

$$\forall i \in \{\, 1, \dots, n \,\} \qquad \alpha_i \,=\, P\big( i \in \mathcal{P}^+(A) \big) - P\big( i \in \mathcal{P}^-(A) \big) \,, \tag{40.24}$$

this yields the inequality stated in the next proposition.

**Proposition 40.5.** *For any event $A \subset \Omega$, we have*

$$\sum_{i=1}^n \Big( P\big( i \in \mathcal{P}^+(A) \big) - P\big( i \in \mathcal{P}^-(A) \big) \Big)^2 \,\leq\, \frac{P_p(A)}{p(1-p)} \,. \tag{40.25}$$

This inequality is considered immediate in the literature dealing with Boolean functions. Indeed, if we consider the Fourier decomposition of the indicator function $1_A$ with respect to the Walsh functions, then the coefficient $\alpha_i$ given in (40.24) is equal to the Fourier coefficient of $X_i$, and the inequality (40.25) becomes a simple consequence of Parseval's identity. Notice that

$$E((X_i)^2) \,=\, \frac{1}{p(1-p)} \,,$$

hence the random variable $X_i$ should be multiplied by $\sqrt{p(1-p)}$ for its $L^2$ norm with respect to $P_p$ to become 1, this explains the presence of the denominator $p(1-p)$ in inequality (40.25). We remark finally that proposition 40.5 is stronger than proposition 39.5, indeed, by the Cauchy-Schwarz inequality, we have

$$\begin{aligned}\Big(E\big(|\mathcal{P}^+(A)|\big)-E\big(|\mathcal{P}^-(A)|\big)\Big)^2 &= \bigg(\sum_{i=1}^n\Big(P\big(\,i\in\mathcal{P}^+(A)\,\big)-P\big(\,i\in\mathcal{P}^-(A)\,\big)\Big)\bigg)^2\\ &\leq n\sum_{i=1}^n\Big(P\big(\,i\in\mathcal{P}^+(A)\,\big)-P\big(\,i\in\mathcal{P}^-(A)\,\big)\Big)^2,\end{aligned}$$

so that inequality (39.10) is implied by inequality (40.25). Using the same strategy as above, we can also improve proposition 40.2 to get the following deviation inequality.

**Proposition 40.6.** *Let $n$ be a positive integer. For any event $A\subset\Omega$, we have*

$$P_p(A)\,\leq\,2\exp\bigg(-\frac{1}{2}\sum_{i=1}^n\Big((1-p)P\big(i\in\mathcal{P}^+(A)\,\big|\,A\big)-pP\big(i\in\mathcal{P}^-(A)\,\big|\,A\big)\Big)^2\bigg).\tag{40.26}$$

*Proof.* The proof follows the same strategy as the proof of proposition 40.2, the difference being that we work with the random variable $\widetilde{S}_n$ defined in (40.20) instead of $S_n$. Let $t>0$. Formula (40.8) holds for $\widetilde{S}_n$ as well:

$$E\big(|\widetilde{S}_n|\,\big|\,A\big)\,\leq\,\frac{1}{P_p(A)}\int_t^{+\infty}P\big(|\widetilde{S}_n|\geq u\big)\,du\,+\,t\,.$$

Now $\widetilde{S}_n$ is the sum of $n$ independent centered random variables which satisfy, for $i$ in $\{\,1,\dots,n\,\}$,

$$\frac{\alpha_i}{1-p}\,\leq\,\alpha_iX_i\,\leq\,\frac{\alpha_i}{p}\quad\text{if}\quad\alpha_i\geq 0\,,\qquad\frac{\alpha_i}{p}\,\leq\,\alpha_iX_i\,\leq\,-\frac{\alpha_i}{1-p}\quad\text{if}\quad\alpha_i<0\,.$$

Applying Hoeffding's inequality [70] to $\widetilde{S}_n$, we get

$$\forall u>0\qquad P\big(|\widetilde{S}_n|\geq u\big)\,\leq\,2\exp\Big(-\frac{2p^2(1-p)^2u^2}{\sum_{1\leq i\leq n}(\alpha_i)^2}\Big)\,.$$

We perform then exactly the same step as in the proof of proposition 40.2, except that the factor $n$ has to be replaced by $\sum_{1\leq i\leq n}(\alpha_i)^2$. This way we obtain the following inequality:

$$E\big(|\widetilde{S}_n|\,\big|\,A\big)\,\leq\,\frac{1}{p(1-p)}\sqrt{2\sum_{1\leq i\leq n}(\alpha_i)^2\ln\frac{2}{P_p(A)}}\,.\tag{40.27}$$

We use the inequality $\big|E\big(\widetilde{S}_n\,\big|\,A\big)\big|\,\leq\,E\big(|\widetilde{S}_n|\,\big|\,A\big)$ and we rewrite the inequality (40.27) in exponential form to get

$$P_p(A)\,\leq\,2\exp\Big(-\frac{p^2(1-p)^2}{2\sum_{1\leq i\leq n}(\alpha_i)^2}\Big(E\big(\widetilde{S}_n\,\big|\,A\big)\Big)^2\Big)\,.\tag{40.28}$$

Formulas (40.1), (40.2) and (40.21) yield

$$E(\widetilde{S}_n 1_A) \,=\, \sum_{i=1}^n \alpha_i\Big(\frac{1}{p}P\big(\{\, i\in\mathcal{P}^+(A)\,\}\cap A\big) - \frac{1}{1-p}P\big(\{\, i\in\mathcal{P}^-(A)\,\}\cap A\big)\Big)\,,$$

whence, dividing by $P_p(A)$,

$$E\big(\widetilde{S}_n \,\big|\, A\big) \,=\, \sum_{i=1}^n \alpha_i\Big(\frac{1}{p}P\big(i\in\mathcal{P}^+(A)\,\big|\,A\big) - \frac{1}{1-p}P\big(i\in\mathcal{P}^-(A)\,\big|\,A\big)\Big)\,. \quad (40.29)$$

The time has come to choose the $\alpha_i$'s. We take

$$\forall i\in\{\,1,\dots,n\,\}\qquad \alpha_i \,=\, \frac{1}{p}P\big(i\in\mathcal{P}^+(A)\,\big|\,A\big) - \frac{1}{1-p}P\big(i\in\mathcal{P}^-(A)\,\big|\,A\big)\,.$$

With this choice, the equation (40.29) together with the inequality (40.28) yield the inequality (40.26) of the proposition. □

Naturally, proposition 40.6 is stronger than proposition 40.2: a routine application of the Cauchy-Schwarz inequality shows that the inequality (40.5) is implied by the inequality (40.26). At this point, we have obtained the analog for the Bernoulli product measure of parameter $p$ of an inequality proved by Talagrand for the uniform probability measure $\mu = P_{1/2}$ on $\{0,1\}^n$. We restate next Talagrand's inequality, adapted to our notation.

**Proposition 40.7** (Proposition 2.2 of [125])**.** *There exists a universal constant $K$ such that, for any subset $A$ of $\{0,1\}^n$, we have*

$$\sum_{1\le i\le n}\Big(\int_A X_i\,d\mu\Big)^2 \,\le\, K\mu(A)^2\ln\frac{e}{\mu(A)}\,. \quad (40.30)$$

Remarking that

$$\frac{1}{\mu(A)}\int_A X_i\,d\mu \,=\, \frac{1}{2}\mu\big(i\in\mathcal{P}^+(A)\,\big|\,A\big) - \frac{1}{2}\mu\big(i\in\mathcal{P}^-(A)\,\big|\,A\big)\,,$$

we can rewrite the inequality (40.30) in an exponential form to get an inequality similar to (40.26), albeit with different constants (but Talagrand did not care about the values of the constants). By the way, the main purpose of the work of Talagrand [125] was to improve the inequality (40.30). This inequality was just a little warm-up before starting the real work!

Keller and Kindler [80] have extended to the biased case the results of Talagrand, in particular they prove an inequality which is much stronger than the deviations inequality (40.26). So, why did we bother to rederive inequality (40.26) from scratch? It should be noted that the proof of Keller and Kindler is quite complex and difficult (although not as mysterious as Talagrand's proof), and it involves also other deep results, typically hypercontractivity estimates. In the end, we will only need the simplest of the previous inequalities, that is the ones obtained in subsection 40.1, which do not involve the second step performed in subsection 40.2. This is the reason why we kept the computations of subsection 40.1. Not only are they more transparent, but they are also more natural from the probabilistic point of view.

# 41 The successive derivatives of $P_p(\omega)$ ⊛

Our goal here is to obtain useful formulas for the successive derivatives of $P_p(\omega)$, and ultimately to control their expectations. Our first approach is based on formula (39.2). Let $\ell \geq 1$. Applying $\ell$ times formula (39.2), we obtain

$$\begin{aligned}\Big(\frac{d}{dp}\Big)^{\ell} P_p(A) \,&=\, \sum_{\substack{1\leq i_1,\dots,i_\ell\leq n\\ \text{pairwise distinct}}} E\big(\delta_{i_1}\cdots\delta_{i_\ell}1_A\big)\\ &=\, \sum_{\substack{I\subset\{\,1,\dots,n\,\}\\ |I|=\ell}} E\Big(\prod_{i\in I}\delta_i\,1_A\Big)\,.\end{aligned} \tag{41.1}$$

We introduce some notation in order to expand the product in the expectation. Let $I$ be a subset of $\{\,1,\dots,n\,\}$. We denote by $\Omega(I)$ the set of the configurations on $I$, or equivalently the set of the maps from $I$ to $\{\,0,1\,\}$. For $\omega$ a configuration in $\Omega$, we denote by $\omega|_I$ its restriction to the components belonging to $I$. Given a configuration $\omega$ in $\Omega$ and a configuration $\sigma$ in $\Omega(I)$, we denote by $\omega(I\leftarrow\sigma)$ the configuration obtained by keeping unchanged the configuration $\omega|_{\{\,1,\dots,n\,\}\setminus I}$ and by setting the configuration $\omega|_I$ to be equal to $\sigma$, i.e.,

$$\forall i\in\{\,1,\dots,n\,\}\qquad \omega(I\leftarrow\sigma)(i)\,=\,\begin{cases}\omega(i) & \text{if } i\notin I\\ \sigma(i) & \text{if } i\in I\end{cases}. \tag{41.2}$$

We denote by $o(\sigma)$ the number of ones in the configuration $\sigma$, i.e.,

$$o(\sigma)\,=\,\sum_{i\in I}\sigma(i)\,.$$

With this notation, we have

$$\prod_{i\in I}\delta_i\,1_A(\omega)\,=\,\sum_{\sigma\in\Omega(I)}(-1)^{|I|-o(\sigma)}1_A\big(\omega(I\leftarrow\sigma)\big)\,. \tag{41.3}$$

This quantity will play an important role to generalize the Margulis-Russo formula, so we introduce a specific definition for it.

**Definition 41.1.** *Let $I$ be a subset of $\{\,1,\dots,n\,\}$ and let $f$ be a Boolean function defined on $\Omega(I)$. The signature of $f$ is*

$$\mathcal{S}(f)\,=\,\sum_{\sigma\in\Omega(I)}(-1)^{|I|-o(\sigma)}f(\sigma)\,. \tag{41.4}$$

The right-hand side of (41.4) is a sum of $2^{|I|}$ terms, half of them with the plus sign, and half of them with the minus sign. Thus an obvious upper bound is

$$\big|\mathcal{S}(f)\big|\,\leq\,2^{|I|-1}\,. \tag{41.5}$$

Denoting by $1_A\big(\omega(I\leftarrow\cdot)\big)$ the Boolean function $\sigma\in\Omega(I)\mapsto 1_A\big(\omega(I\leftarrow\sigma)\big)$, we can rewrite (41.3) with the help of the previous notation as

$$\prod_{i\in I}\delta_i\,1_A(\omega)\,=\,\mathcal{S}\big(1_A(\omega(I\leftarrow\cdot))\big)\,. \tag{41.6}$$

Substituting this identity into (41.1), we have

$$\Big(\frac{d}{dp}\Big)^\ell P_p(A) \;=\; \sum_{\substack{I\subset\{\,1,\dots,n\,\}\\ |I|=\ell}} E\Big(\mathcal{S}\big(1_A(\omega(I\leftarrow\cdot))\big)\Big)\,. \tag{41.7}$$

From the inequality (41.5), we see that the signature of a Boolean function of $\ell$ variables is between $-2^{\ell-1}$ and $2^{\ell-1}$. In the sum in (41.7), we group together the terms which have the same signature, and we get the generalized Margulis-Russo formula stated in the next theorem.

**Theorem 41.2** (Generalized Margulis-Russo formula)**.** *Let $A$ be an event on $\Omega$. For $\ell\geq 1$ and an integer $s$ such that $|s|\leq 2^{\ell-1}$, we define*

$$\forall\omega\in\Omega\qquad \mathcal{S}^{\ell,s}(A,\omega)\;=\;\Big\{\,I\subset\{\,1,\dots,n\,\}:|I|=\ell,\,\mathcal{S}\big(1_A(\omega(I\leftarrow\cdot))\big)=s\,\Big\}\,. \tag{41.8}$$

*We have then*

$$\Big(\frac{d}{dp}\Big)^\ell P_p(A) \;=\; \sum_{-2^{\ell-1}\leq s\leq 2^{\ell-1}} s\,E\big(\big|\mathcal{S}^{\ell,s}(A)\big|\big)\,. \tag{41.9}$$

Formula (41.9) looks awkward. Why do we call it the generalized Margulis-Russo formula? The case $\ell=1$ reads

$$\frac{d}{dp}P_p(A)\;=\;E\big(\big|\mathcal{S}^{1,1}(A)\big|\big)-E\big(\big|\mathcal{S}^{1,-1}(A)\big|\big)\,, \tag{41.10}$$

and it follows from the definition (41.8) of $\mathcal{S}^{\ell,s}(A)$ that

$$\mathcal{S}^{1,1}(A)=\mathcal{P}^+(A)\,,\quad \mathcal{S}^{1,-1}(A)=\mathcal{P}^-(A)\,.$$

The Margulis-Russo formula is usually stated for monotone events. In case the event $A$ is increasing, the Boolean function $1_A\big(\omega(I\leftarrow\cdot)\big)$ on $\Omega(I)$ obtained by fixing the configuration outside a set $I$ is increasing. If $I$ is a singleton, then it is an increasing Boolean function of one variable. The only non-zero signature for such a function is 1, so that $\mathcal{S}^{1,-1}$ is empty, the last term in (41.10) disappears, and we recover the classical Margulis-Russo formula for an increasing event.

Let us consider the case $\ell=2$. There exist 16 Boolean functions of two variables and there are 5 possible values for the signature $\mathcal{S}$, hence there are 5 terms in the sum of (41.9). The usefulness of the Margulis-Russo formula is due to the simplification occurring when considering increasing events. Indeed, there exist only 3 non-constant increasing Boolean functions in two variables, whose signatures are $-1,0,1$. Thus, for an increasing event $A$, the case $\ell=2$ of formula (41.9) reads

$$\frac{d^2}{dp^2}P_p(A)\;=\;E\big(\big|\mathcal{S}^{2,1}(A)\big|\big)-E\big(\big|\mathcal{S}^{2,-1}(A)\big|\big)\,.$$

This formula is stated and discussed after theorem 2.32 in [63]. We will return later to the case of the higher derivatives.

# 42 A Green formula ⊛

We would like to improve the upper bound obtained in (41.5). Our strategy is to rewrite the formula (41.3) with the help of a Green formula. This will lead to a better control on the derivatives for monotone events. We turn the configuration space $\Omega(I)$ into a graph by declaring that two configurations $\sigma_1$ and $\sigma_2$ are neighbors if they differ on one component exactly, i.e., if

$$\exists j \in \{\,1,\dots,n\,\} \quad \sigma_1(j) \neq \sigma_2(j)\,, \qquad \forall i \in \{\,1,\dots,n\,\} \setminus \{\,j\,\} \qquad \sigma_1(i) = \sigma_2(i)\,.$$

We recall next the Green formula for a finite graph (see [42], Proposition 1.3 or [115], Proposition 1.1).

**Proposition 42.1** (Green formula without boundary term)**.** *Let $N$ be a finite graph and let $u, v$ be two functions defined on $N$ with values in $\mathbb{R}$. We have*

$$\frac{1}{2} \sum_{\substack{x,y\in N \\ x,y \text{ neighbors}}} \big(u(y)-u(x)\big)\big(v(y)-v(x)\big) \;=\; -\sum_{x\in N} m(x)\, v(x)\, \Delta u(x)\,, \tag{42.1}$$

*where $m(x)$ is the degree of the vertex $x$, defined by*

$$m(x) \;=\; \big|\,\{\, y\in N : y \text{ is a neighbor of } x\,\}\,\big|\,,$$

*and $\Delta u(x)$ is the Laplacian of $u$, defined by*

$$\Delta u(x) \;=\; \frac{1}{m(x)} \sum_{y \text{ neighbor of } x} u(y)-u(x)\,.$$

*Proof.* The proof is a done by a direct computation. We notice first that, by symmetry, we have

$$\sum_{\substack{x,y\in N \\ x,y \text{ neighbors}}} \big(u(y)-u(x)\big)v(y) \;=\; \sum_{\substack{x,y\in N \\ x,y \text{ neighbors}}} \big(u(y)-u(x)\big)\big(-v(x)\big)\,. \tag{42.2}$$

We start from the left-hand member of (42.1) and we deduce from (42.2) that

$$\begin{aligned} \frac{1}{2} \sum_{\substack{x,y\in N \\ x,y \text{ neighbors}}} \big(u(y)-u(x)\big)\big(v(y)-v(x)\big) \;&=\; \sum_{x\in N} \sum_{\substack{y\in N \\ x,y \text{ neighbors}}} \big(u(y)-u(x)\big)\big(-v(x)\big) \\ &=\; -\sum_{x\in N} m(x)\,\Delta u(x) v(x)\;, \end{aligned}$$

and this is precisely the Green formula (42.1). ☐

In fact, a more general formula is given in Proposition 1.3 of [42]. In the more general version, the set $N$ is only a subset of the vertices of the graph, and there is an extra sum corresponding to boundary terms. Apart from this point, we have kept the original notation of [42] to state the Green formula (42.1).

Let us come back to our problem and let us fix $\omega$ in $\Omega$. We apply the Green formula (recalled in proposition 42.1) on $\Omega(I)$ to the functions $u, v$ defined by

$$\forall \sigma \in \Omega(I) \qquad u(\sigma) \,=\, (-1)^{o(\sigma)}\,, \quad v(\sigma) \,=\, 1_A(\omega(I \leftarrow \sigma))\,,$$

and we get

$$\frac{1}{2} \sum_{\substack{\sigma,\eta \in \Omega(I) \\ \sigma,\eta \text{ neighbors}}} \big(u(\eta) - u(\sigma)\big)\big(v(\eta) - v(\sigma)\big) \,=\, - \sum_{\sigma \in \Omega(I)} |I|\, v(\sigma)\, \Delta u(\sigma)\,. \tag{42.3}$$

For two neighboring configurations $\sigma, \eta$, we have in fact

$$u(\eta) - u(\sigma) \,=\, (-1)^{o(\eta)} - (-1)^{o(\sigma)} \,=\, 2(-1)^{o(\eta)}\,. \tag{42.4}$$

We compute next the Laplacian of $u$:

$$\Delta u(\sigma) \,=\, \frac{1}{|I|} \sum_{\eta \text{ neighbor of } \sigma} (-1)^{o(\eta)} - (-1)^{o(\sigma)} \,=\, -2(-1)^{o(\sigma)}\,. \tag{42.5}$$

With the help of (42.4) and (42.5), we can rewrite (42.3) as

$$\sum_{\substack{\sigma,\eta \in \Omega(I) \\ \sigma,\eta \text{ neighbors}}} (-1)^{o(\eta)} \big(v(\eta) - v(\sigma)\big) \,=\, 2|I| \sum_{\sigma \in \Omega(I)} 1_A(\omega(I \leftarrow \sigma))(-1)^{o(\sigma)}\,.$$

Comparing with formula (41.3), we see that

$$\begin{aligned}
\prod_{i \in I} \delta_i \, 1_A(\omega) \,&=\, \frac{(-1)^{|I|}}{2|I|} \sum_{\substack{\sigma,\eta \in \Omega(I) \\ \sigma,\eta \text{ neighbors}}} (-1)^{o(\eta)} \big(v(\eta) - v(\sigma)\big) \\
&=\, \frac{(-1)^{|I|}}{|I|} \sum_{\substack{\sigma,\eta \in \Omega(I),\, \sigma < \eta \\ \sigma,\eta \text{ neighbors}}} (-1)^{o(\eta)} \big(v(\eta) - v(\sigma)\big)\,.
\end{aligned} \tag{42.6}$$

We try next to get a better upper bound for monotone events on the product $\prod_{i \in I} \delta_i \, 1_A(\omega)$ than the one obtained in (41.5) for a general event. From now onwards, we consider the case where the event $A$ is increasing. Taking the absolute value in the equality (42.6), we obtain the following upper bound:

$$\Big| \prod_{i \in I} \delta_i \, 1_A(\omega) \Big| \,\leq\, \frac{1}{|I|} \sum_{\substack{\sigma,\eta \in \Omega(I),\, \sigma < \eta \\ \sigma,\eta \text{ neighbors}}} \big(v(\eta) - v(\sigma)\big)\,. \tag{42.7}$$

We split the sum appearing in (42.7) and we compute:

$$\begin{aligned}
\sum_{\substack{\sigma,\eta \in \Omega(I),\, \sigma < \eta \\ \sigma,\eta \text{ neighbors}}} \big(v(\eta) - v(\sigma)\big) \,&=\, \sum_{\eta \in \Omega(I)} \sum_{\substack{\sigma \in \Omega(I),\, \sigma < \eta \\ \sigma,\eta \text{ neighbors}}} v(\eta) - \sum_{\sigma \in \Omega(I)} \sum_{\substack{\eta \in \Omega(I),\, \sigma < \eta \\ \sigma,\eta \text{ neighbors}}} v(\sigma) \\
&=\, \sum_{\eta \in \Omega(I)} v(\eta) o(\eta) - \sum_{\sigma \in \Omega(I)} v(\sigma) \big(|I| - o(\sigma)\big) \\
&=\, \sum_{\sigma \in \Omega(I)} v(\sigma) \big(2o(\sigma) - |I|\big)\,.
\end{aligned} \tag{42.8}$$

Obviously, we have furthermore

$$\sum_{\sigma\in\Omega(I)} v(\sigma)\big(2o(\sigma)-|I|\big) \;\leq\; \sum_{|I|/2<k\leq|I|}\;\sum_{\substack{\sigma\in\Omega(I)\\ o(\sigma)=k}} \big(2k-|I|\big)$$

$$= \sum_{|I|/2<k\leq|I|} \binom{|I|}{k}\big(2k-|I|\big)\,. \quad (42.9)$$

We compute next the last sum appearing in (42.9), which we denote by $S(|I|)$, i.e.,

$$S(|I|) \;=\; \sum_{|I|/2<k\leq|I|} \binom{|I|}{k}\big(2k-|I|\big)\,. \quad (42.10)$$

This quantity has a natural probabilistic interpretation. We consider the configuration space $\Omega\;=\;\{\,0,1\,\}^n$, equipped with the uniform probability measure, and we define the random sum $S_n$ as

$$S_n(\omega) \;=\; \sum_{i=1}^n \big(2\omega(i)-1\big) \;=\; 2\Big(\sum_{i=1}^n \omega(i)-\frac{n}{2}\Big)\,.$$

Alternatively, the random integer $S_n$ is the position of the symmetric random walk on the integers after $n$ steps. Taking $n\;=\;|I|$ and using the symmetry with respect to the sign reversal, the quantity $S(|I|)$ of formula (42.10) can be rewritten as

$$S(n) \;=\; 2^n E\big(1_{\{\,S_n>0\,\}}S_n\big) \;=\; 2^{n-1}E\big(|S_n|\big)\,.$$

We compute next this quantity. We distinguish two cases, depending on the parity of $|I|$.

**First case:** $|I|=2p$ **for some** $p\geq 1$. We write

$$S(2p) \;=\; \sum_{k=p+1}^{2p}\left(2k\binom{2p}{k}-2p\binom{2p}{k}\right)$$

$$= \sum_{k=p+1}^{2p}\left(4p\binom{2p-1}{k-1}-2p\binom{2p}{k}\right)$$

$$= \sum_{k=p}^{2p-1} 4p\binom{2p-1}{k} - \sum_{k=p+1}^{2p} 2p\binom{2p}{k}\,. \quad (42.11)$$

Moreover we have that

$$\sum_{k=p+1}^{2p}\binom{2p}{k} \;=\; \sum_{k=0}^{p-1}\binom{2p}{k} \;=\; \frac{1}{2}\left(2^{2p}-\binom{2p}{p}\right), \quad (42.12)$$

while

$$\sum_{k=p}^{2p-1}\binom{2p-1}{k} \;=\; \sum_{k=0}^{p-1}\binom{2p-1}{k} \;=\; 2^{2p-2}\,. \quad (42.13)$$

Reporting (42.12) and (42.13) in (42.11), we obtain

$$S(2p) \;=\; p\binom{2p}{p}\,.$$

We can check the first values: $S(2) = 2$, $S(4) = 12$, $S(6) = 60$.
**Second case:** $|I| = 2p+1$ **for some** $p \geq 0$. We write

$$\begin{aligned} S(2p+1) \;&=\; \sum_{k=p+1}^{2p+1}\left(2k\binom{2p+1}{k} - (2p+1)\binom{2p+1}{k}\right) \\ &=\; \sum_{k=p+1}^{2p+1}\left(2(2p+1)\binom{2p}{k-1} - (2p+1)\binom{2p+1}{k}\right) \\ &=\; \sum_{k=p}^{2p} 2(2p+1)\binom{2p}{k} - \sum_{k=p+1}^{2p+1}(2p+1)\binom{2p+1}{k}\,. \qquad (42.14) \end{aligned}$$

Moreover we have that

$$\sum_{k=p}^{2p}\binom{2p}{k} \;=\; \sum_{k=0}^{p}\binom{2p}{k} \;=\; \frac{1}{2}\left(2^{2p} + \binom{2p}{p}\right), \qquad (42.15)$$

while

$$\sum_{k=p+1}^{2p+1}\binom{2p+1}{k} \;=\; \sum_{k=0}^{p}\binom{2p+1}{k} \;=\; 2^{2p}\,. \qquad (42.16)$$

Reporting (42.15) and (42.16) in (42.14), we obtain

$$S(2p+1) \;=\; (2p+1)\binom{2p}{p} \;=\; (p+1)\binom{2p+1}{p+1}\,.$$

We can check the first values: $S(1) = 1$, $S(3) = 6$, $S(5) = 30$, $S(7) = 140$. Putting together the results of the two cases, we conclude that

$$S(|I|) \;=\; \left\lfloor \frac{|I|+1}{2} \right\rfloor \binom{|I|}{\left\lfloor \frac{|I|+1}{2} \right\rfloor}\,. \qquad (42.17)$$

Reporting successively (42.17) in (42.10), (42.9), (42.8) and (42.7), we arrive at the following upper bound:

$$\Big|\prod_{i\in I}\delta_i\, 1_A(\omega)\Big| \;\leq\; \frac{1}{|I|}\left\lfloor \frac{|I|+1}{2} \right\rfloor \binom{|I|}{\left\lfloor \frac{|I|+1}{2} \right\rfloor}\,.$$

As $|I|$ goes to $\infty$, the right-hand quantity behaves like $2^{|I|-1/2}/\sqrt{\pi|I|}$, thus this upper bound is a serious improvement of the previous upper bound $2^{|I|-1}$.

# 43 The Bernoulli chaos ⊛

Our next goal is to generalize the classical bound (39.11) on the expected number of pivotal sites. To this end, we will generalize propositions 39.3 and 39.5. More precisely, we would like to obtain a probabilistic representation of the $\ell$-th derivative of $P_p(A)$ which generalizes the nice formula $P_p'(A) = E\big(1_A S_n\big)$. In formula (39.4), we had computed the first derivative of $P_p(\omega)$:

$$\frac{d}{dp}P_p(\omega) \,=\, \sum_{i=1}^n X_i(\omega)\, P_p(\omega)\,.$$

Recall that $X_i$ is defined as

$$\forall \omega\in\Omega \qquad X_i(\omega) \,=\, \frac{\omega(i)}{p} - \frac{1-\omega(i)}{1-p}\,, \tag{43.1}$$

so that it depends also on the parameter $p$. Differentiating again, we get

$$\frac{d^2}{dp^2}P_p(\omega) \,=\, \bigg(\Big(\sum_{i=1}^n X_i(\omega)\Big)^2 + \sum_{i=1}^n \frac{\partial X_i}{\partial p}(\omega)\bigg)\, P_p(\omega)\,.$$

It turns out that

$$\frac{\partial X_i}{\partial p} \,=\, -\frac{\omega(i)}{p^2} - \frac{1-\omega(i)}{(1-p)^2} \,=\, -\big(X_i(\omega)\big)^2\,, \tag{43.2}$$

whence

$$\frac{d^2}{dp^2}P_p(\omega) \,=\, \bigg(\sum_{\substack{1\leq i,j\leq n\\ i\neq j}} X_i(\omega)\, X_j(\omega)\bigg)\, P_p(\omega)\,.$$

We are naturally lead to the following result.

**Proposition 43.1.** *Let $\ell$ be such that $1\leq \ell\leq n$ and let $\omega$ belong to $\Omega$. We define the Bernoulli chaos of order $\ell$ as*

$$Q_\ell(\omega) \,=\, \sum_{\substack{1\leq i_1,\dots,i_\ell\leq n\\ \text{pairwise distinct}}} X_{i_1}(\omega)\,\cdots\, X_{i_\ell}(\omega)\,. \tag{43.3}$$

*The $\ell$-th derivative of $P_p(\omega)$ is equal to*

$$\Big(\frac{d}{dp}\Big)^\ell P_p(\omega) \,=\, Q_\ell(\omega)\, P_p(\omega)\,. \tag{43.4}$$

*Proof.* We prove the formula by induction on $\ell$. The formula has already been proved for $\ell=1$ and $\ell=2$. Suppose that it holds for some $\ell\geq 1$. Differentiating once more formula (43.4), we have

$$\Big(\frac{d}{dp}\Big)^{\ell+1} P_p(\omega) \,=\, \bigg(Q_\ell(\omega)\sum_{1\leq i_{\ell+1}\leq n} X_{i_{\ell+1}}(\omega) + \frac{d}{dp}Q_\ell(\omega)\bigg)\, P_p(\omega)\,. \tag{43.5}$$

Using again (43.2) to compute the derivative of the $\ell$-fold product appearing in $Q_\ell(\omega)$, we see that the products containing a square cancels out, and we obtain the formula (43.4) at rank $\ell+1$. □

## 43.1 Variance and large deviations of the chaos

The Bernoulli chaos of order $\ell$ given by $Q_\ell$ is centered. We wish to control its deviations from its mean. Using the expression (43.3) and the fact that the random variables $X_i$ are i.i.d. and centered, we compute

$$E\big((Q_\ell)^2\big) \,=\, \binom{n}{\ell}\frac{\ell!^2}{p^\ell(1-p)^\ell}\,. \tag{43.6}$$

We are ready to generalize proposition 39.3 to the higher derivatives.

**Proposition 43.2.** *For any event $A\subset\Omega$ and for any $\ell\in\{\,1,\dots,n\,\}$, we have*

$$\Big|\Big(\frac{d}{dp}\Big)^\ell P_p(A)\Big| \,\leq\, \ell!\sqrt{P_p(A)\binom{n}{\ell}\frac{1}{p^\ell(1-p)^\ell}}\,. \tag{43.7}$$

*Proof.* It follows from proposition 43.1 that

$$\begin{aligned}\Big(\frac{d}{dp}\Big)^\ell P_p(A) \,=\, \Big(\frac{d}{dp}\Big)^\ell\sum_{\omega\in A}P_p(\omega) \,&=\, \sum_{\omega\in A}\Big(\frac{d}{dp}\Big)^\ell P_p(\omega)\\ &=\,\sum_{\omega\in A}Q_\ell(\omega)\,P_p(\omega)\,=\,E\big(1_AQ_\ell\big)\,.\end{aligned} \tag{43.8}$$

A straightforward application of Cauchy-Schwarz inequality gives

$$\big|E\big(1_AQ_\ell\big)\big| \,\leq\, \sqrt{P_p(A)\,E\big((Q_\ell)^2\big)}\,. \tag{43.9}$$

Putting together formulas (43.8), (43.9) and (43.6), we obtain the desired inequality (43.7). ☐

In fact, the chaos defined by $Q_\ell$ is a specific instance of an $U$ statistics. There are numerous papers dealing with the large deviations behaviour of general $U$ statistics under various hypotheses. Yet $Q_\ell$ is one of the simplest $U$ statistics, so will rely on the exponential inequality which was obtained by Hoeffding in his seminal paper [70]. The inequality (5.7) of [70], specialized to the case where the function $g$ is taken to be $g(x_1,\dots,x_\ell)\,=\,x_1\cdots x_\ell$, with the bounds $0\,\leq\,g(X_1,\dots,X_\ell)\,\leq\,q^{-\ell}$, where we have set for convenience $q=\min(p,1-p)$, yields the following:

$$\forall t>0\qquad P\big(|Q_\ell|\geq n(n-1)\cdots(n-\ell+1)t\big)\,\leq\,2\exp\Big(-\frac{1}{2}\big\lfloor\frac{n}{\ell}\big\rfloor q^{2\ell}t^2\Big)\,. \tag{43.10}$$

Let us suppose that $\ell\geq 2$ and $n\geq 2\ell$. We have then

$$\big\lfloor\frac{n}{\ell}\big\rfloor\,\geq\,\frac{n}{\ell}-1\,\geq\,\frac{n}{2\ell}\,,\qquad n(n-1)\cdots(n-\ell+1)\,\leq\,n^\ell\,,$$

so that, after setting $s=n(n-1)\cdots(n-\ell+1)t$, formula (43.10) becomes

$$\forall s>0\qquad P\big(|Q_\ell|\geq s\big)\,\leq\,2\exp\Big(-\frac{q^{2\ell}s^2}{4\ell n^{2\ell-1}}\Big)\,.$$

## 43.2 Computation of the chaos

Rewriting the expression (43.4) of $Q_\ell(\omega)$ as

$$\frac{1}{\ell!}\,Q_\ell(\omega)\;=\;\sum_{1\le i_1<\cdots<i_\ell\le n} X_{i_1}(\omega)\,\cdots\,X_{i_\ell}(\omega)\,,$$

we see that $\frac{1}{\ell!}Q_\ell$ is the $\ell$-th elementary symmetric polynomial of $X_1,\dots,X_n$. Using Newton's identities, we know that $Q_\ell$ can be expressed as a polynomial in $\ell$ variables of degree $\ell$ of the first $\ell$ power sums

$$\Sigma_k\;=\;\sum_{1\le i\le n}\big(X_i\big)^k\,,\qquad 1\le k\le \ell\,.$$

Recalling the expression (43.1) of $X_i$, we compute

$$\Sigma_k\;=\;\sum_{1\le i\le n}\frac{\omega(i)}{p^k}+(-1)^k\frac{1-\omega(i)}{(1-p)^k}\;=\;\Big(\frac{1}{p^k}-\frac{(-1)^k}{(1-p)^k}\Big)o(\omega)+\frac{n(-1)^k}{(1-p)^k}\,,$$

where $o(\omega)\;=\;\sum_{i=1}^n\omega(i)$ is the number of ones in the configuration $\omega$. So we see that the power sums $\Sigma_k$ are in fact affine functions of $o(\omega)$, and we conclude that there exists a polynomial $\widetilde{Q}_\ell$ in one variable of degree $\ell$ such that

$$\forall\omega\in\Omega\qquad Q_\ell(\omega)\;=\;\widetilde{Q}_\ell\big(o(\omega)\big)\,.\tag{43.11}$$

In order to determine the polynomial $\widetilde{Q}_\ell$, we shall take another road to compute the derivatives of $P_p(\omega)$. We start from the expression

$$\forall\omega\in\Omega\qquad P_p(\omega)\;=\;p^{o(\omega)}(1-p)^{n-o(\omega)}\,,$$

and we shall apply the classical Leibniz rule for the $\ell$-th derivative of a product. To alleviate the formulas, we set $x=o(\omega)$ and

$$\phi(p,x)\;=\;p^x(1-p)^{n-x}\,.$$

We will use also the Pochhammer symbol given by

$$(x)_n\;=\;x(x+1)\cdots(x+n-1)\,.$$

From Leibniz rule, we have

$$\begin{aligned}\Big(\frac{d}{dp}\Big)^\ell\phi(p,x)\;&=\;\sum_{m=0}^{\ell}\binom{\ell}{m}\Big(\frac{d}{dp}\Big)^m p^x\Big(\frac{d}{dp}\Big)^{\ell-m}(1-p)^{n-x}\\
&=\;\sum_{m=0}^{\ell}\binom{\ell}{m}(x-m+1)_m p^{x-m}(x-n)_{\ell-m}(1-p)^{n-x-(\ell-m)}\\
&=\;\phi(p,x)\frac{\ell!}{(p(1-p))^\ell}\sum_{m=0}^{\ell}\frac{(x-m+1)_m}{m!}p^{\ell-m}\frac{(x-n)_{\ell-m}}{(\ell-m)!}(1-p)^m\,,\end{aligned}\tag{43.12}$$

so we conclude that the polynomial $\widetilde{Q}_\ell$ is given by

$$\widetilde{Q}_\ell(x)\;=\;\frac{\ell!}{(p(1-p))^\ell}\sum_{m=0}^{\ell}\frac{(x-m+1)_m}{m!}p^{\ell-m}\frac{(x-n)_{\ell-m}}{(\ell-m)!}(1-p)^m\,.$$

## 43.3 The Krawtchouk polynomials

The polynomial function corresponding to the sum in formula (43.12) is in fact a Krawtchouk polynomial. Krawtchouk polynomials are to the binomial distribution what Hermite polynomials are to the Brownian motion. The $n$-th Krawtchouk polynomial with parameters $N, p$ is

$$k_n^{(p)}(x, N) = \sum_{k=0}^{n} \frac{(x-k+1)_k}{k!} p^{n-k} \frac{(x-N)_{n-k}}{(n-k)!} (1-p)^k \,, \tag{43.13}$$

so that we have indeed

$$\left(\frac{d}{dp}\right)^{\ell} \phi(p,x) = \frac{\ell!}{(p(1-p))^{\ell}} \, k_\ell^{(p)}(x,n)\, \phi(p,x) \,.$$

From the expression (43.13), we see that the highest degree term of $k_n^{(p)}(x, N)$ is

$$\sum_{k=0}^{n} \frac{x^k}{k!} p^{n-k} \frac{x^{n-k}}{(n-k)!} (1-p)^k = \frac{x^n}{n!} \,. \tag{43.14}$$

Where do the Krawtchouk polynomials come from? They were originally introduced by Mikhail Kravchuk in 1929 [91]. The Krawtchouk polynomials are the orthogonal polynomials associated to the binomial distribution, i.e., they satisfy the orthogonality conditions:

$$\sum_{x=0}^{n} k_\ell^{(p)}(x,n)\, k_m^{(p)}(x,n) \binom{n}{x} p^x (1-p)^{n-x} = \begin{cases} n & \text{if } \ell = m \\ 0 & \text{if } \ell \neq m \end{cases} .$$

There exist mainly two versions of them, which depend on the normalization condition. With the choice of formula (43.13), we have

$$\sum_{x=0}^{n} \left(k_\ell^{(p)}(x,n)\right)^2 \binom{n}{x} p^x (1-p)^{n-x} = \binom{n}{\ell} \left(p(1-p)\right)^{\ell} .$$

These polynomials can be studied within the framework of the general theory of orthogonal polynomials. In the continuous setting, orthogonal polynomials correspond to eigenfunctions of second order differential operators and there exist a complete theory classifying them for the case of symmetric operators. Most of the results of the continuous theory carry over to the discrete setting, and naturally differential operators have to be replaced by difference operators. The Krawtchouk polynomials are discrete orthogonal polynomials (meaning that the weight function defining the scalar product is a discrete measure supported on the set of the integers). There is a huge mathematical literature devoted to the theory of orthogonal polynomials, both in the continuous and in the discrete settings. Here are some excellent book references: [10, 13, 14, 74, 88, 109]. Throughout the literature, there are some notational variations for the Krawtchouk polynomials concerning the indices or the normalization constants. We will follow the notation and convention of the book of Beals and Wong [13]. To make a long story short, there are four main families of discrete orthogonal polynomials which are the eigenfunctions of a symmetric second order difference operator, associated with the names Charlier, Krawtchouk, Meixner, Hahn.

Of course, it is not a hazard that the Krawtchouk polynomials pop up when we differentiate $P_p(\omega)$. Here is a little argument explaining that the formula for the derivatives of $P_p(\omega)$ must involve these orthogonal polynomials. We start again with

$$\phi(p,x) \;=\; p^x(1-p)^{n-x}\,,$$

and we observe the first derivative with respect to $p$:

$$\frac{\partial\phi}{\partial p}(p,x) \;=\; \widetilde{Q}_1(p,x)\phi(p,x)\,,$$

where

$$\widetilde{Q}_1(p,x) \;=\; \frac{x}{p}-\frac{n-x}{1-p} \;=\; \frac{x-np}{p(1-p)}\,. \tag{43.15}$$

We keep the previous notation $\widetilde{Q}_1$, so in fact $\widetilde{Q}_1$ is the polynomial introduced in formula (43.11), the difference being that we write now $\widetilde{Q}_1(p,x)$ instead of simply $\widetilde{Q}_1(x)$. Next, for the second derivative, we have

$$\frac{\partial^2\phi}{\partial p^2}(p,x) \;=\; \Big(\frac{\partial\widetilde{Q}_1}{\partial p}(p,x)+\big(\widetilde{Q}_1(p,x)\big)^2\Big)\phi(p,x) \;=\; \widetilde{Q}_2(p,x)\phi(p,x)\,,$$

where

$$\widetilde{Q}_2(p,x) \;=\; \Big(\frac{x}{p}-\frac{n-x}{1-p}\Big)^2-\frac{x}{p^2}-\frac{n-x}{(1-p)^2}\,.$$

So we guess, or we prove by induction, that

$$\frac{\partial^\ell\phi}{\partial p^\ell}(p,x) \;=\; \widetilde{Q}_\ell(p,x)\phi(p,x)\,, \tag{43.16}$$

where $\widetilde{Q}_\ell(p,x)$ is a polynomial function in $x$ of degree $\ell$, whose coefficients are rational fractions in the variable $p$. In fact, differentiating the relation (43.16), we obtain the recursion formula

$$\widetilde{Q}_{\ell+1}(p,x) \;=\; \frac{\partial\widetilde{Q}_\ell}{\partial p}(p,x)+\widetilde{Q}_\ell(p,x)\widetilde{Q}_1(p,x)\,, \tag{43.17}$$

and the previous statements on $\widetilde{Q}_\ell$ are direct consequences of this formula. Of course, formula (43.17) is the counterpart for $\widetilde{Q}_\ell$ of formula (43.5) for $Q_\ell$. Another noteworthy fact, which can be proved by induction with the help of formula (43.17), is that the highest degree term of $\widetilde{Q}_\ell(p,x)$ is

$$\frac{x^\ell}{(p(1-p))^\ell}\,. \tag{43.18}$$

We had also computed the variance of $Q_\ell$ in (43.6). The transfer formula yields then that

$$\sum_{x=0}^{n}\binom{n}{x}\big(\widetilde{Q}_\ell(p,x)\big)^2\phi(p,x) \;=\; \binom{n}{\ell}\frac{\ell!^2}{p^\ell(1-p)^\ell}\,. \tag{43.19}$$

We will next deduce from the previous results that the family of the polynomials $1,\widetilde{Q}_1,\dots,\widetilde{Q}_\ell$ is an orthogonal family. To start, let us consider the binomial identity:

$$\sum_{x=0}^{n}\binom{n}{x}\phi(p,x) \;=\; 1\,,$$

and let us differentiate with respect to $p$. It comes

$$\sum_{x=0}^{n}\binom{n}{x}\widetilde{Q}_1(p,x)\phi(p,x) \;=\; 0\,,$$

hence $\widetilde{Q}_1$ is orthogonal to the constant polynomial 1 for the binomial weights. Furthermore, we have, using formula (43.15),

$$\sum_{x=0}^{n}\binom{n}{x}\big(\widetilde{Q}_1(p,x)\big)^2\phi(p,x) \;=\;$$

$$\frac{1}{(p(1-p))^2}\sum_{x=0}^{n}\binom{n}{x}(x-np)^2\phi(p,x) \;=\; \frac{n}{p(1-p)}\,,$$

where we have recognized the variance of the binomial distribution in the last line. We rewrite the previous equality as

$$\sum_{x=0}^{n}\binom{n}{x}p(1-p)\big(\widetilde{Q}_1(p,x)\big)^2\phi(p,x) \;=\; n\,. \tag{43.20}$$

Differentiating identity (43.20) with respect to $p$, we obtain

$$\begin{aligned}0 \;=\; &\sum_{x=0}^{n}\binom{n}{x}\frac{\partial}{\partial p}\big(p(1-p)\widetilde{Q}_1(p,x)\big)\widetilde{Q}_1(p,x)\phi(p,x)\\ &+\sum_{x=0}^{n}\binom{n}{x}p(1-p)\widetilde{Q}_1(p,x)\frac{\partial}{\partial p}\big(\widetilde{Q}_1(p,x)\phi(p,x)\big)\\ &\qquad=\sum_{x=0}^{n}\binom{n}{x}p(1-p)\widetilde{Q}_1(p,x)\widetilde{Q}_2(p,x)\phi(p,x)\,,\end{aligned}$$

hence the polynomials $\widetilde{Q}_1(p,x)$ and $\widetilde{Q}_2(p,x)$ are orthogonal. At step $\ell\geq 1$, suppose that we have proved that $1,\widetilde{Q}_1,\dots,\widetilde{Q}_\ell$ is an orthogonal system and let us rewrite (43.19) as

$$\sum_{x=0}^{n}\binom{n}{x}\big(p(1-p)\big)^\ell\big(\widetilde{Q}_\ell(p,x)\big)^2\phi(p,x) \;=\; \binom{n}{\ell}\ell!^2\,. \tag{43.21}$$

We differentiate the relation (43.21) with respect to $p$:

$$0 \;=\; \sum_{x=0}^{n} \binom{n}{x} \frac{\partial}{\partial p}\Big(\big(p(1-p)\big)^{\ell}\widetilde{Q}_\ell(p,x)\Big)\widetilde{Q}_\ell(p,x)\phi(p,x)$$
$$+ \sum_{x=0}^{n} \binom{n}{x} \big(p(1-p)\big)^{\ell}\widetilde{Q}_\ell(p,x)\frac{\partial}{\partial p}\big(\widetilde{Q}_\ell(p,x)\phi(p,x)\big)\,. \quad (43.22)$$

We notice that

$$\frac{\partial}{\partial p}\Big(\big(p(1-p)\big)^{\ell}\widetilde{Q}_\ell(p,x)\Big)$$

is a polynomial in $x$ of degree at most $\ell-1$, hence it follows from the induction hypothesis that it is orthogonal to $\widetilde{Q}_\ell(p,x)$ and the first sum above vanishes. Moreover, we have

$$\frac{\partial}{\partial p}\big(\widetilde{Q}_\ell(p,x)\phi(p,x)\big) \;=\; \widetilde{Q}_{\ell+1}(p,x)\phi(p,x)\,, \qquad (43.23)$$

and we deduce from equations (43.22) and (43.23) that $\widetilde{Q}_\ell$ and $\widetilde{Q}_{\ell+1}$ are orthogonal. For $m<\ell$, we already know that

$$\sum_{x=0}^{n} \binom{n}{x} \widetilde{Q}_m(p,x)\widetilde{Q}_\ell(p,x)\phi(p,x) \;=\; 0\,.$$

We multiply by $\big(p(1-p)\big)^m$ and we differentiate with respect to $p$:

$$0 \;=\; \sum_{x=0}^{n} \binom{n}{x} \frac{\partial}{\partial p}\Big(\big(p(1-p)\big)^{m}\widetilde{Q}_m(p,x)\Big)\widetilde{Q}_\ell(p,x)\phi(p,x)$$
$$+ \sum_{x=0}^{n} \binom{n}{x} \big(p(1-p)\big)^{m}\widetilde{Q}_m(p,x)\frac{\partial}{\partial p}\big(\widetilde{Q}_\ell(p,x)\phi(p,x)\big)$$
$$= \sum_{x=0}^{n} \binom{n}{x} \big(p(1-p)\big)^{m}\widetilde{Q}_m(p,x)\,\widetilde{Q}_{\ell+1}(p,x)\phi(p,x)\,,$$

from which we conclude that $\widetilde{Q}_{\ell+1}$ is orthogonal to $\widetilde{Q}_m$. This completes the induction.

The family of polynomials $1,\widetilde{Q}_1,\dots,\widetilde{Q}_\ell$ is an orthogonal system such that

$$\forall \ell\in\{\,1,\dots,n\,\}\qquad \text{degree}\,\widetilde{Q}_\ell \;=\; \ell\,,$$

therefore they are proportional to the Krawtchouk polynomials. By comparing the terms of highest degree in both polynomials $\widetilde{Q}_\ell(p,x)$ and $k_\ell^{(p)}(x,n)$, as given in (43.14) and (43.18), we conclude that

$$\widetilde{Q}_\ell(p,x) \;=\; \frac{\ell!}{(p(1-p))^{\ell}}k_\ell^{(p)}(x,n)\,.$$

# 44 The higher pivots

In subsection 29.5, we introduced the notion of pivotal bonds of order $k$ for the disconnection event $\mathcal{D}$. This notion is crucial for understanding the structure of the intersections generated by the intertwined explorations. When applied to the bond percolation model, the classical Margulis-Russo formula gives the expected cardinality of the pivotal bonds of order 1, and the generalized Margulis-Russo formula developed in section 41 involves the pivotal bonds of order $k$, but in a more complicated way. In order to understand better the relationship between the generalized Margulis-Russo formula and the pivotal bonds of order $k$, we first define the higher pivots in the general framework of the product space $\Omega = \{0, 1\}^n$, without any reference to the geometry of the lattice $\mathbb{Z}^d$.

## 44.1 Definition of the higher pivots ⊙

Let us recall an important notation (already introduced in (41.2)) before defining the pivots of higher order. For $\omega$ in $\Omega = \{0,1\}^n$, $I$ a subset of $\{1, \dots, n\}$, and a configuration $\sigma$ in $\Omega(I) = \{0,1\}^I$, we denote by $\omega(I \leftarrow \sigma)$ the configuration obtained by keeping unchanged the configuration $\omega|_{\{1,\dots,n\}\setminus I}$ and by setting the configuration $\omega|_I$ to be equal to $\sigma$, i.e.,

$$\forall i \in \{1, \dots, n\} \qquad \omega(I \leftarrow \sigma)(i) = \begin{cases} \omega(i) & \text{if } i \notin I\,, \\ \sigma(i) & \text{if } i \in I\,. \end{cases}$$

Let $A \subset \Omega$ be an event. Let $\omega \in \Omega$ be a configuration and let $k \geq 1$. We say that a component $i$ in $\{1, \dots, n\}$ is pivotal of order $k$ for the event $A$ in the configuration $\omega$, or $k$-pivotal, if there exists a subset $I$ of $\{1, \dots, n\}$ of cardinality strictly less than $k$ and a configuration $\sigma$ in $\Omega(I)$ such that, if we transform the configuration $\omega$ by forcing the components $I$ to be in the same states as in $\sigma$, then the component $i$ becomes pivotal for $A$. We denote by $\mathcal{P}_k(A, \omega)$ the set of the $k$-pivotal components for $A$ in the configuration $\omega$. We make the convention that $\mathcal{P}_0(A, \omega) = \varnothing$.

As for the usual pivots, the set of the pivots of order $k$ is a random subset of $\{1, \dots, n\}$. For $k = 1$, the set $\mathcal{P}_1(A, \omega)$ coincides with the classical set $\mathcal{P}(A, \omega)$ of the components which are pivotal for $A$ in $\omega$. With the previous notation, we have, for $k \geq 1$,

$$\mathcal{P}_k(A, \omega) = \Big\{ i \in \{1, \dots, n\} : \exists\, I \subset \{1, \dots, n\} \quad |I| \leq k - 1 \\ \exists\, \sigma \in \Omega(I) \quad i \in \mathcal{P}\big(A, \omega(I \leftarrow \sigma)\big) \Big\}. \quad (44.1)$$

With this definition, we see that the sets $(\mathcal{P}_k)_{k \geq 1}$ form a non-decreasing sequence of subsets of $\{1, \dots, n\}$. After some index $k \leq n$, the sequence becomes stationary.

## 44.2 Signature and pivotality ⊛

The classical Margulis-Russo formula involves the pivots of first order. The generalized Margulis-Russo formula involves the higher pivots, but in a more complicated way, through the signature of a Boolean function that was defined in (41.4). Let $A$ be an event on $\Omega$. Let $\ell \geq 1$ and let $s$ be an integer such that $|s| \leq 2^{\ell-1}$. The set $\mathcal{S}^{\ell,s}(A)$ was defined in (41.8) as

$$\forall \omega \in \Omega \qquad \mathcal{S}^{\ell,s}(A,\omega) \,=\, \Big\{\, I \subset \{\,1,\dots,n\,\} : |I| = \ell,\, \mathcal{S}\big(1_A(\omega(I \leftarrow \cdot))\big) = s \,\Big\}\,. \tag{44.2}$$

The set $\mathcal{S}^{\ell,s}(A)$ is a subset of the subsets of $\{\,1,\dots,n\,\}$ having exactly $\ell$ elements. We define a projection $\pi$ by setting

$$\pi\big(\mathcal{S}^{\ell,s}(A,\omega)\big) \,=\, \Big\{\, i \in \{\,1,\dots,n\,\} : \exists I \in \mathcal{S}^{\ell,s}(A,\omega)\,, \quad i \in I \,\Big\}\,.$$

We make the convention that $\mathcal{P}_0(A) = \varnothing$ and $\mathcal{S}^{0,s}(A) = \varnothing$ for $s \neq 0$. The relationship between the sets $\mathcal{S}^{\ell,s}$ and $\mathcal{P}_\ell$ is clarified in the next proposition.

**Proposition 44.1.** *Let $A$ be an event on $\Omega$. For any $\ell \in \{\,1,\dots,n\,\}$, we have the inclusions*

$$\forall s \neq 0 \qquad \pi\big(\mathcal{S}^{\ell,s}(A)\big) \,\subset\, \mathcal{P}_\ell(A)\,, \tag{44.3}$$

$$\mathcal{P}_\ell(A) \setminus \mathcal{P}_{\ell-1}(A) \,\subset\, \pi\big(\mathcal{S}^{\ell,-1}(A) \,\cup\, \mathcal{S}^{\ell,+1}(A)\big)\,. \tag{44.4}$$

*Proof.* For $\ell = 0$ and $s \neq 0$, both sets $\mathcal{S}^{0,s}(A)$ and $\mathcal{P}_0(A)$ are empty by convention, therefore we have indeed $\pi\big(\mathcal{S}^{0,s}(A)\big) \subset \mathcal{P}_0(A)$. We consider next the case where $\ell \in \{\,1,\dots,n\,\}$ and $s \neq 0$. Let $\omega \in \Omega$, and let $i^*$ belong to $\pi\big(\mathcal{S}^{\ell,s}(A,\omega)\big)$. By definition, there exists $I$ in $\mathcal{S}^{\ell,s}(A,\omega)$ which contains $i^*$. We have, using (41.3) and (41.4),

$$s \,=\, \mathcal{S}\big(1_A(\omega(I \leftarrow \cdot))\big) \,=\, \prod_{i \in I} \delta_i\, 1_A(\omega)\,.$$

The operators $\delta_i$ commute. We put $\delta_{i^*}$ at the end of the last product and we develop the remaining product over $I \setminus \{i^*\}$ with the help of formula (41.3):

$$\begin{aligned}\prod_{i \in I} \delta_i\, 1_A(\omega) \,=\, &\prod_{i \in I \setminus \{i^*\}} \delta_i\big(\delta_{i^*}\, 1_A\big)(\omega) \\ &\qquad\qquad = \sum_{\sigma \in \Omega(I \setminus \{i^*\})} (-1)^{|I| - o(\sigma)} \big(\delta_{i^*}\, 1_A\big)\big(\omega(I \leftarrow \sigma)\big)\,.\end{aligned}$$

For this sum to be non-zero, it must be the case that at least one term in the sum is non-zero. Therefore

$$\exists\, \sigma \in \Omega(I \setminus \{i^*\}) \qquad \big(\delta_{i^*}\, 1_A\big)\big(\omega(I \leftarrow \sigma)\big) \neq 0\,.$$

But this is equivalent to

$$\exists\, \sigma \in \Omega(I \setminus \{i^*\}) \qquad i^* \in \mathcal{P}\big(A, \omega(I \leftarrow \sigma)\big)\,,$$

hence the component $i^*$ belongs to $\mathcal{P}_\ell(A,\omega)$. This completes the proof of the inclusion (44.3). We turn now to the proof of the inclusion (44.4). Let $i^*$ belong to $\mathcal{P}_\ell(A,\omega)\setminus\mathcal{P}_{\ell-1}(A,\omega)$. For $\ell=1$, we have

$$\mathcal{P}_1(A,\omega) \,=\, \mathcal{S}^{1,-1}(A,\omega)\,\cup\,\mathcal{S}^{1,1}(A,\omega)\,,$$

and there is even equality in formula (44.4). Suppose next that $\ell\geq 2$. By definition, there exist $I\subset\{\,1,\dots,n\,\}$ and $\sigma\in\Omega(I)$ such that $|I|\leq\ell-1$ and $i^*\in\mathcal{P}\big(A,\omega(I\leftarrow\sigma)\big)$. This implies that $\delta_{i^*}\,1_A\big(\omega(I\leftarrow\sigma)\big)$ is equal to $-1$ or $+1$. We compute next the product $\prod_{i\in I\cup\{i^*\}}\delta_i\,1_A(\omega)$. The operators $\delta_i$ commute. We put $\delta_{i^*}$ at the beginning of the product and we develop the remaining product over $I\setminus\{i^*\}$ with the help of formula (41.3):

$$\begin{aligned}\prod_{i\in I\cup\{i^*\}}\delta_i\,1_A(\omega) \,&=\, \delta_{i^*}\Big(\prod_{i\in I}\delta_i\Big)\,1_A(\omega)\\ &=\,\delta_{i^*}\Big(\sum_{\eta\in\Omega(I)}(-1)^{|I|-o(\eta)}1_A\big(\omega(I\leftarrow\eta)\big)\Big)\\ &=\,\sum_{\eta\in\Omega(I)}(-1)^{|I|-o(\eta)}\delta_{i^*}\Big(1_A\big(\omega(I\leftarrow\eta)\big)\Big)\,.\end{aligned}\tag{44.5}$$

Let $\eta$ belong to $\Omega(I)$ and suppose that there exists $i_0$ in $I$ such that $\eta(i_0)=\omega(i_0)$. Letting $J=I\setminus\{i_0\}$, we have $1_A\big(\omega(I\leftarrow\eta)\big)=1_A\big(\omega(J\leftarrow\eta)\big)$, and since $i^*$ does not belong to $\mathcal{P}_{\ell-1}(A,\omega)$ and $|J|\leq\ell-2$, then $i^*$ does not belong to $\mathcal{P}\big(A,\omega(J\leftarrow\eta)\big)$, whence

$$\delta_{i^*}\Big(1_A\big(\omega(I\leftarrow\eta)\big)\Big)\,=\,0\,.$$

As a consequence, in the last sum of formula (44.5), all the terms associated to a configuration $\eta$ which coincides with $\omega$ on one component in $I$ vanish. Only one term remains, corresponding to the configuration $\eta_0$ obtained by modifying all the components of $\omega$ in $I$, i.e., $\eta_0=1-\omega|_I$. By the way, we also knew from the start that $\delta_{i^*}\,1_A\big(\omega(I\leftarrow\sigma)\big)$ does not vanish, so it must be the case that $\sigma=\eta_0$. Anyway, we see that

$$\prod_{i\in I\cup\{i^*\}}\delta_i\,1_A(\omega)\,=\,(-1)^{|I|-o(\sigma)}\delta_{i^*}\Big(1_A\big(\omega(I\leftarrow\sigma)\big)\Big)\,\in\,\{\,-1,+1\,\}\,.$$

Therefore

$$I\cup\{i^*\}\in\mathcal{S}^{\ell,-1}(A,\omega)\,\cup\,\mathcal{S}^{\ell,+1}(A,\omega)\,,$$

and thus $i^*$ belongs to the image of this last set by the projection $\pi$. This concludes the proof of the inclusion (44.4). □

Proposition 44.1 yields the following interesting corollary.

**Corollary 44.2.** *Let $A$ be an event on $\Omega$. For any $k\in\{\,1,\dots,n\,\}$, we have*

$$\mathcal{P}_k(A)\,=\,\bigcup_{1\leq\ell\leq k}\pi\Big(\mathcal{S}^{\ell,-1}(A)\,\cup\,\mathcal{S}^{\ell,+1}(A)\Big)\,.$$

## 44.3 What is missing? ⊛

The computations conducted in the general framework of an arbitrary product space give interesting results on the higher pivots. Although they seem to go in the right direction, they are not strong enough for our purpose. A natural possibility is to resort to more specific aspects of our problem. We consider the case where $\mathcal{D}$ is a disconnection event, as in the first and the third questions of section 30. Our main problem is to control the probability $P\big(\big|\mathcal{I}(t)\big| \geq i\big)$ for $t, i \geq 1$. By Markov's inequality, we have

$$P\big(\big|\mathcal{I}(t)\big| \,\geq\, i\big) \,\leq\, \frac{1}{i} E\big(\big|\mathcal{I}(t)\big|\big)\,. \tag{44.6}$$

For $t = 1$, we know that $\mathcal{I}(1) = \mathcal{P}_1(\mathcal{D})$, and we have an adequate bound on $E\big(\big|\mathcal{P}_1(\mathcal{D})\big|\big)$, which is the counterpart of inequality (39.11) for decreasing events. Suppose now that $t \geq 2$. We write

$$E\big(\big|\mathcal{I}(t)\big|\big) \,=\, \sum_{e\in\mathbb{E}^d(D)} P\big(e \in \mathcal{I}(t)\big)\,. \tag{44.7}$$

Suppose we could prove a bound like the following:

$$\exists\, c \geq 1 \quad \forall t \geq 1 \quad \forall e \in \mathbb{E}^d(D) \qquad P\big(e \in \mathcal{I}(t+1)\big) \,\leq\, c\, P\big(e \in \mathcal{I}(t)\big)\,, \tag{44.8}$$

where the constant $c$ depends on the dimension $d$ and the parameter $p$ only. Summing inequality (44.8) over $e \in \mathbb{E}^d(D)$ and using (44.7), this would imply

$$\forall t \geq 1 \qquad E\big(\big|\mathcal{I}(t+1)\big|\big) \,\leq\, c\, E\big(\big|\mathcal{I}(t)\big|\big)\,. \tag{44.9}$$

Iterating inequality (44.9), we would obtain

$$\forall t \geq 1 \qquad E\big(\big|\mathcal{I}(t)\big|\big) \,\leq\, c^{t-1}\, E\big(\big|\mathcal{P}_1(\mathcal{D})\big|\big)\,. \tag{44.10}$$

Let us focus on the first and the third questions of section 30, namely the disconnection of a sub-box of $\Lambda(n)$ and the disconnection in the prism $\Pi(n)$. We denote by $D(n)$ the domain where the disconnection event $\mathcal{D}$ occurs (hence $D(n)$ is either $\Lambda(n)$ or $\Pi(n)$). In corollaries 30.2 and 30.4, we have reduced the problem of bounding the probability $P\big(\mathcal{D}\big)$ to that of bounding the probability $P\big(\big|\mathcal{I}(t)\big| \geq i\big)$ for $i$ larger than a polynomial in $n$ and for $t$ smaller than $2(\ln n)/\beta(n)$, where $\beta(n)$ is a function tending to $+\infty$ as $n$ goes to $+\infty$. Using inequalities (44.10) and (39.11), we would have

$$\forall t \in \Big\{\, 1, \dots, 2\frac{\ln n}{\beta(n)} \,\Big\} \qquad E\big(\big|\mathcal{I}(t)\big|\big) \,\leq\, c^{2\frac{\ln n}{\beta(n)}} \sqrt{\frac{\big|\mathbb{E}^d(D(n))\big|}{p(1-p)}}\,. \tag{44.11}$$

Substituting (44.11) in inequality (44.6), we would get

$$\max_{1\,\leq\, t\,\leq\, 2\frac{\ln n}{\beta(n)}} P\big(\,|\mathcal{I}(t)| \,\geq\, i\,\big) \,\leq\, \frac{1}{i} c^{2\frac{\ln n}{\beta(n)}} \sqrt{\frac{\big|\mathbb{E}^d(D(n))\big|}{p(1-p)}}\,. \tag{44.12}$$

The point is that the factor $c^{2\frac{\ln n}{\beta(n)}}$ is negligible compared to $n$, so the exponent of $n$ given by this upper bound remains under control.

**Application to the disconnection of a sub-box of $\Lambda(n)$.** In this case, we have $D(n)=\Lambda(n)$ and $|\mathbb{E}^d(D(n))|=O(n^d)$. Using inequality (44.12), together with corollary 30.2, we would get the following control:

$$P_p\big(\,\Lambda(m)\not\leftrightarrow\partial^{\,in}\Lambda(n)\,\big)\;=\;\frac{2(\ln n)^2}{dm^{d-1}\beta(n)}c^{2\frac{\ln n}{\beta(n)}}\,O(n^{d/2})\,.$$

We would conclude that, for any $\alpha$ such that $\alpha>d/(2(d-1))$,

$$P_p\big(\,\Lambda(n^\alpha)\not\leftrightarrow\partial^{\,in}\Lambda(n)\,\big)\;=\;n^{d/2-(d-1)\alpha+o(1)}\,. \tag{44.13}$$

Not only do we see that a sub-box of side length $n^\alpha$, with $\alpha>d/(2(d-1))$, is unlikely to be disconnected from $\partial^{\,in}\Lambda(n)$, but we have also some polynomial control on the probability of this disconnection.

**Application to the disconnection in the prism $\Pi(n)$.** In this case, we have $D(n)=\Pi(n)$ and again $|\mathbb{E}^d(D(n))|=O(n^d)$. Using inequality (44.12), together with corollary 30.4, we would get the following control:

$$P\big(\,L\not\leftrightarrow R\quad\text{in }\Pi(n)\,\big)\;=\;\frac{1}{n^{d/2-1+o(1)}}\,. \tag{44.14}$$

In both cases, we would have obtained a quantitative bound on the disconnection probability, and this would be a major step forward. Remember that these results rest on the unproven inequality (44.8). So let us discuss further this inequality. For simplicity, we shall consider the case of odd indices. Let $t\geq 0$ be fixed. Let $e$ be a bond belonging to $D$ (the domain where the disconnection event $\mathcal{D}$ occurs), let us denote by $\mathcal{E}(e)$ the event

$$\mathcal{E}(e)\;=\;\big\{\,e\in\mathcal{I}(2t+1)\,\big\}\,,$$

and let us consider a configuration $\omega$ in $D$ realizing the event $\mathcal{E}(e)$. We have shown in corollary 29.8 that the bond $e$ can be written $e=\langle a,b\rangle$, where $a$ is in $\mathcal{A}(2t+1)$ and $b$ in $\mathcal{A}(2t)$. Moreover, we have

$$\begin{gathered}T_{\mathcal{A}(2t)}(W,a)=T_D(W\cup W',a)=t\,,\\ T_{\mathcal{A}(2t+1)}(W',b)=T_D(W\cup W',b)=t\,.\end{gathered}$$

We consider the three following scenarii:

• Scenario 1: there exists a closed bond $e_1$ emanating from $b$ which is pivotal for the event $\big\{\,T_{\mathcal{A}(2t+1)}(W',b)=t\,\big\}$.

• Scenario 2: any path in $\mathcal{A}(2t)$ between $W$ and $a$ having travel time $t$ terminates with an open bond emanating from $a$.

• Scenario 3: neither the scenario 1 nor the scenario 2 occurs.

Let us set, for $i=1,2,3$,

$$\mathcal{E}_i(e)\;=\;\big\{\,\omega\in\mathcal{E}:\text{scenario }i\text{ occurs in }\omega\,\big\}\,.$$

We bound the probability of $\mathcal{E}(e)$ as follows:

$$P\big(\mathcal{E}(e)\big) \leq P\big(\mathcal{E}_1(e)\big) + P\big(\mathcal{E}_2(e)\big) + P\big(\mathcal{E}_3(e)\big). \tag{44.15}$$

Suppose that the scenario 1 occurs in the configuration $\omega$. Let $\Phi_1(\omega)$ be the configuration in which the bond $e_1$ has been opened. In this configuration, we have

$$T_D(W',b) = T_{\mathcal{A}(2t-1)}(W',b) = t-1,$$

and $e$ belongs now to the intersection set $\mathcal{I}(2t)$, i.e., we have $e \in \mathcal{I}(2t, \Phi_1(\omega))$. Moreover we changed only the state of one bond to transform $\omega$ into $\Phi_1(\omega)$, thus

$$P(\omega) = \frac{1-p}{p} P(\Phi_1(\omega)) \leq \frac{1}{p} P(\Phi_1(\omega)).$$

Summing over the configurations $\omega$ of $\mathcal{E}_1(e)$, we conclude that

$$P\big(\mathcal{E}_1(e)\big) = \sum_{\omega \in \mathcal{E}_1(e)} P(\omega) \leq \frac{1}{p} \sum_{\omega \in \mathcal{E}_1(e)} P\big(\Phi_1(\omega)\big). \tag{44.16}$$

Reindexing the last sum with $\widetilde{\omega} = \Phi_1(\omega)$, we get

$$\begin{aligned} \sum_{\omega \in \mathcal{E}_1(e)} P(\Phi_1(\omega)) &= \sum_{\widetilde{\omega} \in \Phi_1(\mathcal{E}_1(e))} P(\widetilde{\omega}) \big|\Phi_1^{-1}(\widetilde{\omega})\big| \\ &\leq 2d \sum_{\widetilde{\omega} \in \Phi_1(\mathcal{E}_1(e))} P(\widetilde{\omega}) = 2dP\big(\Phi_1(\mathcal{E}_1(e))\big), \end{aligned} \tag{44.17}$$

because the map $\Phi_1$ modifies only one bond emanating from $b$, hence the preimage $\Phi_1^{-1}(\widetilde{\omega})$ of any configuration $\widetilde{\omega}$ has cardinality at most $2d$. As we discussed previously, the set $\Phi_1(\mathcal{E}_1(e))$ is included in the event $\{\, e \in \mathcal{I}(2t) \,\}$. Thus, it follows from (44.16) and (44.17) that

$$P\big(\mathcal{E}_1(e)\big) \leq \frac{2d}{p} P\big(e \in \mathcal{I}(2t)\big). \tag{44.18}$$

Suppose that the scenario 2 occurs in the configuration $\omega$. Let $\Phi_2(\omega)$ be the configuration in which all the bonds emanating from $a$ and having their second endpoint in $\mathcal{A}(2t)$ are closed. In this configuration, we have

$$T_D(W,a) = T_{\mathcal{A}(2t)}(W,a) = t+1,$$

and the bond $e$ belongs now to the intersection set $\mathcal{I}(2t+2)$, i.e., we have $e \in \mathcal{I}(2t+2, \Phi_2(\omega))$. Moreover we changed the states of at most $2d-1$ bonds to transform $\omega$ into $\Phi_2(\omega)$, thus

$$P(\omega) \leq \Big(\frac{p}{1-p}\Big)^{2d-1} P(\Phi_2(\omega)) \leq \frac{1}{(1-p)^{2d}} P(\Phi_2(\omega)).$$

The map $\Phi_2$ modifies at most $2d$ bonds emanating from $a$, hence the preimage $\Phi_2^{-1}(\widetilde{\omega})$ of any configuration $\widetilde{\omega}$ has cardinality at most $2^{2d}$. Summing over the configurations $\omega$ of $\mathcal{E}_2(e)$, and proceeding as for the event $\mathcal{E}_1(e)$, we obtain

$$\begin{aligned} P\big(\mathcal{E}_2(e)\big) \;&\leq\; \frac{1}{(1-p)^{2d}} \sum_{\omega\in\mathcal{E}_2(e)} P(\Phi_2(\omega)) \;=\; \frac{1}{(1-p)^{2d}} \sum_{\widetilde{\omega}\in\Phi_2(\mathcal{E}_2(e))} P(\widetilde{\omega})\big|\Phi_2^{-1}(\widetilde{\omega})\big| \\ &\leq\; \Big(\frac{2}{1-p}\Big)^{2d} P\big(\Phi_2(\mathcal{E}_2(e))\big) \;\leq\; \Big(\frac{2}{1-p}\Big)^{2d} P(e\in\mathcal{I}(2t+2))\,. \end{aligned} \tag{44.19}$$

The computations leading to the inequalities (44.18) and (44.19) are standard. They are consequences of more general results on random modifications of the percolation configuration that have been stated and used many times in the literature, see for instance lemma 6.3 of [31] or lemma 7 of [45]. However, in simple instances like ours, it seems better to perform the computations from scratch rather than to make appeal to a more general result involving complicated notation. Substituting the inequalities (44.18) and (44.19) into (44.15), we obtain

$$\begin{aligned} &P\big(e\in\mathcal{I}(2t+1)\big) \;\leq \\ &\qquad \frac{2d}{p} P\big(e\in\mathcal{I}(2t)\big) \;+\; \Big(\frac{2}{1-p}\Big)^{2d} P\big(e\in\mathcal{I}(2t+2)\big) \;+\; P\big(\mathcal{E}_3(e)\big)\,. \end{aligned}$$

Of course a similar inequality can be derived for the even indices, so in the end we have an inequality of the following form: there exists a constant $c=c(d,p)$ depending on the dimension $d$ and $p$ only such that

$$\forall t\geq 2 \qquad E\big(\big|\mathcal{I}(t)\big|\big) \;\leq\; c\,E\big(\big|\mathcal{I}(t-1)\big|\big) + c\,E\big(\big|\mathcal{I}(t+1)\big|\big) + \sum_{e\in\mathbb{E}^d(D)} P\big(\mathcal{E}_3(e)\big)\,. \tag{44.20}$$

Unfortunately, this is not quite the desired inequality (44.9). Notice however that an inequality of the form

$$\exists c\geq 1 \quad \forall t\geq 2 \qquad E\big(\big|\mathcal{I}(t)\big|\big) \;\leq\; c\,E\big(\big|\mathcal{I}(t-1)\big|\big) \tag{44.21}$$

would also be enough to obtain very interesting conclusions, similar to (44.13) and (44.14). In fact, the inequality (44.21) readily implies that

$$\forall t\in\Big\{1,\dots,2\frac{\ln n}{\beta(n)}\Big\} \qquad E\big(\big|\mathcal{I}(t)\big|\big) \;\leq\; c^{2\frac{\ln n}{\beta(n)}}\, E\Big(\Big|\mathcal{I}\Big(\Big\lfloor 2\frac{\ln n}{\beta(n)}\Big\rfloor\Big)\Big|\Big)\,. \tag{44.22}$$

We know from (29.16) that the intersection set $\mathcal{I}(t)$ becomes empty for $t$ larger than the final iteration number $T$. Moreover, we have at our disposal a quantitative control on the final iteration number, in proposition 30.1 for the disconnection of a sub-box of $\Lambda(n)$ and in corollary 30.4 for the disconnection in the prism $\Pi(n)$. Using these controls, we bound the last expectation in (44.22) as follows:

$$E\Big(\Big|\mathcal{I}\Big(\Big\lfloor 2\frac{\ln n}{\beta(n)}\Big\rfloor\Big)\Big|\Big) \;\leq\; |\mathbb{E}^d(D(n))|P\Big(T\;>\;2\frac{\ln n}{\beta(n)}\Big) \;\leq\; O(n^d)\,n^{3d-\beta(n)}\,. \tag{44.23}$$

The inequalities (44.22) and (44.23) yield a result stronger than the inequality (44.11), and their application to the disconnection of a sub-box of $\Lambda(n)$ and the disconnection in the prism $\Pi(n)$ would yield even better results than (44.13) and (44.14). So, the question is, can we get rid of the last term appearing in the right-hand side of inequality (44.20), and possibly remove one of the first two?

Let us discuss the last term, involving the third scenario. We would like to show that this scenario is unlikely. A possibility is to put more constraints on it, so we consider instead the three following scenarii, which have also the advantage of being more symmetric with respect to the endpoints $a, b$ of $e$:

• Scenario 1: there exists a closed bond $e_1$ emanating from $a$ which is pivotal for the event $\{ T_{\mathcal{A}(2t)}(W, a) = t \}$ or a closed bond $e_2$ emanating from $b$ which is pivotal for the event $\{ T_{\mathcal{A}(2t+1)}(W', b) = t \}$.
• Scenario 2: any path in $\mathcal{A}(2t)$ between $W$ and $a$ having travel time $t$ terminates with an open bond emanating from $a$, and any path in $\mathcal{A}(2t+1)$ between $W'$ and $b$ having travel time $t$ terminates with an open bond emanating from $b$.
• Scenario 3: neither the scenario 1 nor the scenario 2 occurs.

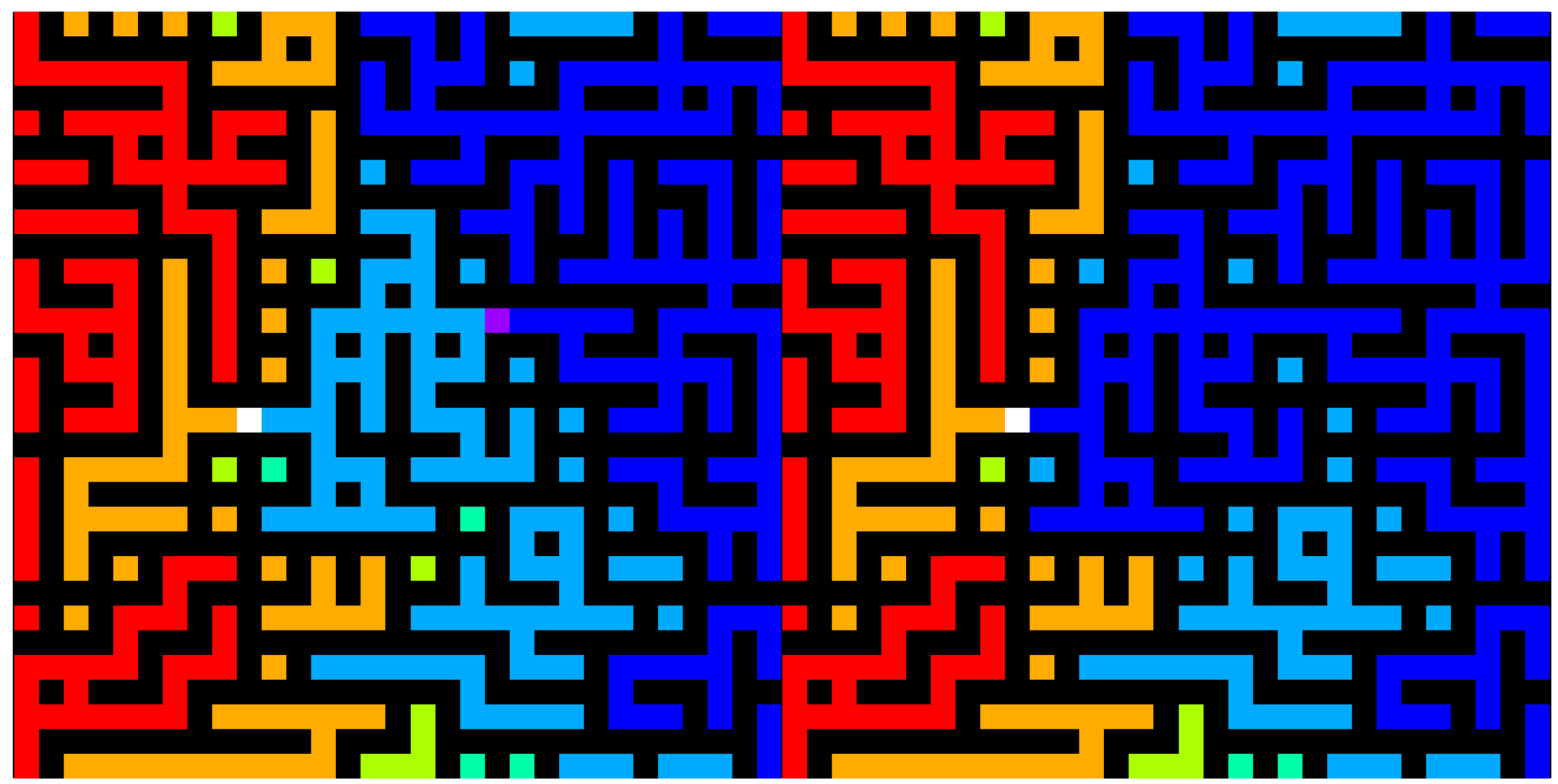

Figure 78: The scenario 4 in $D = \Lambda(16)$, $p = 0.47$. The white bond is initially in $\mathcal{I}(3)$. The opening of the purple bond makes the white bond belong to $\mathcal{I}(2)$.

It is still not clear that the probability of the third scenario is small. So, instead of considering only the bonds emanating from $e$, we allow ourselves to modify the configuration in a neighborhood of $e$. Let $\Gamma(e, m)$ be a box of side length $m$ centered at the middle of $e$, and let us consider the three following scenarii:

• Scenario 4: there exists a closed bond $e_1$ in $\Gamma(e, m)$ which is pivotal for the event $\{ T_{\mathcal{A}(2t)}(W, a) = t \}$ or a closed bond $e_2$ in $\Gamma(e, m)$ which is pivotal for the event $\{ T_{\mathcal{A}(2t+1)}(W', b) = t \}$.
• Scenario 5: any path in $\mathcal{A}(2t)$ between $W$ and $a$ having travel time $t$ visits only open bonds of $\Gamma(e, m)$, and any path in $\mathcal{A}(2t+1)$ between $W'$ and $b$ having travel time $t$ visits only open bonds of $\Gamma(e, m)$.
• Scenario 6: neither the scenario 4 nor the scenario 5 occurs.

Now, whenever the sixth scenario occurs, it is the case that the two endpoints $a, b$ of $e$ are connected by open paths to the boundary of $\Gamma(e,m)$, but they are not connected together inside $\Gamma(e,m)$. Therefore the two arms event occurs and we can use the inequality (8.18) to control the probability of the scenario 6. However, this control is too weak to reach a useful conclusion. Another strategy consists in modifying the configuration inside $\Gamma(e,m)$ as follows. When the sixth scenario occurs, and we have the two arms event, we close simultaneously all the bonds emanating from $a$ and $b$. In this new configuration, the bond $e$, which was initially in the intersection set $\mathcal{I}(2t+1)$, is reached two rounds later by the intertwined explorations, and it becomes an element of $\mathcal{I}(2t+3)$. More generally, if in the box $\Gamma(e,m)$ there exist an open bond $f_a$ which is pivotal for the travel time between $a$ and $L$, and an open bond $f_b$ which is pivotal for the travel time between $b$ and $R$, we close simultaneously these two bonds. This has again the effect of transforming the bond $e$ into an element of $\mathcal{I}(2t+3)$. We can even allow ourselves to close a finite number of bonds to increase the travel times of $a$ and $b$ to the boundary. When the box $\Gamma(e,m)$ is large, it is very likely that such a modification is possible. So we could hope to get an inequality of the form

$$\exists c\geq 1\quad \forall t\geq 2\qquad E\big(\big|\mathcal{I}(t)\big|\big)\,\leq\, c\,E\big(\big|\mathcal{I}(t+2)\big|\big)\,.$$

Such an inequality would again perfectly serve our purpose. Unfortunately, we can prove this inequality with a constant $c$ independent of $n$ (the parameter controlling the size of the domain $D$), but not independent of $t$. Indeed, given an integer $t\geq 1$ and a bond $e$ in $\mathcal{I}(2t+1)$, we can certainly perform the required modification in the box $\Gamma(e,3t)$, but we end up with a constant $c$ depending on $t$, and which grows very fast with $t$. Furthermore, it is not a technical limitation and it does seem to be avoidable. The reason is that, whenever a bond $e$ belongs to an intersection set $\mathcal{I}(2t+1)$ with $t$ large, it will be the case that it is surrounded by an abnormally large number of closed bonds. In this situation, it becomes increasingly difficult to find open bonds in the neighborhood of $e$ whose closure might increase the travel times between the endpoints of $e$ and the boundary of $D$. Finally, to make things worse, even if we work with these new scenarii $4,5,6$, we encounter problems to deal with the first two scenarii. For instance, in scenario 4, whenever we find in $\Gamma(e,m)$ a closed bond $e_1$ which is pivotal for the travel time between $a$ and $L$, it might be the case that the opening of this bond makes the bond $e$ disappear from the intersection sets. In fact, the effect of opening $e_1$ will be to change the order of pivotality of $e$. So, we could at least hope to obtain an inequality involving the higher pivotal bonds instead of the intersection sets, of the form

$$\exists c=c(d,p)>0\quad \forall t\geq 2\qquad E\big(\big|\mathcal{P}(t)\big|\big)\,\leq\, c\,E\big(\big|\mathcal{P}(t-1)\big|\big)+c\,E\big(\big|\mathcal{P}(t+1)\big|\big)\,.$$

Pivotal bonds of order $t$ are created (respectively destroyed) by opening (respectively closing) a pivotal bond of order $t+1$, or by closing (respectively opening) a pivotal bond of order $t-1$. Thus there must exist some relations involving the number of pivotal bonds of different orders.

# 45 New attempts on the second pivots ⊛

In this section, we focus on the second pivots and we try to obtain a useful control on their cardinality.

## 45.1 Application of the general results

In subsection 44.2, we have established a rigorous link between the sets $(\mathcal{S}^{\ell,s}, 1 \leq \ell \leq n)$ and $(\mathcal{P}_\ell, 1 \leq \ell \leq n)$. Moreover, the generalized Margulis-Russo formula and the control of the variance of the Bernoulli chaos yield an inequality controlling the expected cardinality of the sets $\mathcal{S}^{\ell,s}$ for $1 \leq \ell \leq n$ and $|s| \leq 2^{\ell-1}$. In this subsection, we shall examine what can be inferred from this inequality on the second pivots. Let us first rewrite the definition of the set $\mathcal{S}^{\ell,s}$ for $\ell = 2$. Let $A$ be an event on $\Omega$ and let $\omega$ be an element of $\Omega$. For an integer $s$ such that $|s| \leq 2$, the set $\mathcal{S}^{2,s}(A,\omega)$ is given by

$$\mathcal{S}^{2,s}(A,\omega) \,=\, \Big\{ \{ i,j \} \subset \{ 1,\dots,n \} : i \neq j,\, \mathcal{S}\big(1_A(\omega(\{ i,j \} \leftarrow \cdot))\big) = s \Big\},$$

where $\mathcal{S}\big(1_A(\omega(\{ i,j \} \leftarrow \cdot))\big)$ is the signature given by

$$\begin{aligned} \mathcal{S}\big(1_A(\omega(\{ i,j \} \leftarrow \cdot))\big) \,&=\, \delta_i\delta_j 1_A(\omega) \\ &=\, 1_A(\omega^{ij}) - 1_A(\omega^i_j) - 1_A(\omega^j_i) + 1_A(\omega_{ij})\,. \end{aligned} \quad (45.1)$$

For a general event $A$, the signature (45.1) can take the five values $-2$, $-1$, $0$, $1$, $2$. However, we are mainly concerned with a disconnection event in percolation, which is decreasing. So, from now onwards, we suppose that the event $A$ is decreasing. For such an event $A$, the signature can take only the three values $-1, 0, 1$, and, as we already discussed after theorem 41.2, the generalized Margulis formula reads in this case

$$\frac{d^2}{dp^2} P_p(A) \,=\, E\big(\big|\mathcal{S}^{2,+1}(A)\big|\big) - E\big(\big|\mathcal{S}^{2,-1}(A)\big|\big)\,. \quad (45.2)$$

Let us take a closer look at the sets $\mathcal{S}^{2,-1}$ and $\mathcal{S}^{2,+1}$.
**The set $\mathcal{S}^{2,-1}$.** The pair $\{ i,j \}$ belongs to $\mathcal{S}^{2,-1}$ if and only if

$$1_A(\omega^{ij}) = 0, \quad 1_A(\omega^i_j) = 1, \quad 1_A(\omega^j_i) = 1, \quad 1_A(\omega_{ij}) = 1\,.$$

**The set $\mathcal{S}^{2,+1}$.** The pair $\{ i,j \}$ belongs to $\mathcal{S}^{2,+1}$ if and only if

$$1_A(\omega^{ij}) = 0, \quad 1_A(\omega^i_j) = 0, \quad 1_A(\omega^j_i) = 0, \quad 1_A(\omega_{ij}) = 1\,.$$

In particular, for $\omega$ in $\Omega$ such that $\{ i,j \} \in \mathcal{S}^{2,+1}(A,\omega)$, we have the equivalence

$$\omega \in A \quad\Longleftrightarrow\quad \omega(i) = \omega(j) = 0\,.$$

It turns out that the set $\mathcal{S}^{2,+1}(A)$ can be completely described with the help of the set $\mathcal{P}_1(A)$. Indeed, for any $i,j \in \{ 1,\dots,n \}$, $i \neq j$, we have the equivalence:

$$\omega \in A\,, \quad \{ i,j \} \in \mathcal{S}^{2,+1} \qquad\Longleftrightarrow\qquad \{ i,j \} \subset \mathcal{P}_1(A,\omega)\,. \quad (45.3)$$

We see also from the expression (45.1) that the signature $\mathcal{S}\big(1_A(\omega(\{\,i,j\,\}\leftarrow\cdot))\big)$ is a function of the configuration $\omega$ outside the components $i,j$. Therefore, the membership of a pair $\{\,i,j\,\}$ to a set $\mathcal{S}^{2,s}(A,\omega)$ is independent of the specific values of $\omega(i),\omega(j)$. Using the previous facts, we can decompose the first expectation appearing in (45.2) as follows:

$$
\begin{gathered}
E\big(\big|\mathcal{S}^{2,+1}(A)\big|\big) \,=\, \sum_{1\leq i<j\leq n} E\big(1_{\{\,i,j\,\}\in\mathcal{S}^{2,+1}(A)}\big) \\
=\sum_{1\leq i<j\leq n}\frac{1}{(1-p)^2}E\big(1_A 1_{\{\,i,j\,\}\in\mathcal{S}^{2,+1}(A)}\big) \,=\, \frac{1}{(1-p)^2}\sum_{1\leq i<j\leq n} E\big(1_{\{\,i,j\,\}\subset\mathcal{P}_1(A)}\big) \\
=\frac{1}{2(1-p)^2}E\Big(\big|\mathcal{P}_1(A)\big|\big(\big|\mathcal{P}_1(A)\big|-1\big)\Big)\,. \qquad (45.4)
\end{gathered}
$$

We apply next proposition 43.2 for $\ell=2$, and we get the inequality

$$
\Big|\Big(\frac{d}{dp}\Big)^2 P_p(A)\Big| \,\leq\, 2!\sqrt{P_p(A)\frac{n(n-1)}{2p^2(1-p)^2}} \,\leq\, \frac{2n}{p(1-p)}\,. \qquad (45.5)
$$

Combining formulas (45.2), (45.4) and (45.5), we obtain

$$
E\big(\big|\mathcal{S}^{2,-1}(A)\big|\big) \,\leq\, \frac{1}{2(1-p)^2}E\Big(\big|\mathcal{P}_1(A)\big|\big(\big|\mathcal{P}_1(A)\big|-1\big)\Big) \,+\, \frac{2n}{p(1-p)}\,. \qquad (45.6)
$$

Ideally, we would like to obtain an inequality involving $E(\big|\mathcal{P}_2(A)\big|)$, but the inequality (45.6) does not seem to lend itself to further simplification, probably because it is too general.

To move forward, we consider the specific case where the event $A$ is a disconnection event $\mathcal{D}$ in percolation. In order to control the right-hand side of (45.6), we would like to control the second moment of $|\mathcal{P}_1(\mathcal{D})|$. We can achieve this with the help of the control of the first intersection set that was obtained in subsection 31.1. We consider two disjoint subsets $W,W'$ of the box $\Lambda(n)$ and the specific case where $\mathcal{D}$ is the disconnection event $\mathcal{D}=\{\,W\not\longleftrightarrow W'\,\}$. In this situation, the first intersection set $\mathcal{I}(1)$ coincides with the pivotal set $\mathcal{P}_1(\mathcal{D})$ and the inequality (31.8) gives

$$
\forall i_1\geq 1 \qquad P\big(|\mathcal{P}_1(\mathcal{D})|\geq i_1,\mathcal{D}\big) \,\leq\, 6d2^dn^d\exp\Big(-\frac{2p^2i_1^2}{9d2^dn^d}\Big)\,. \qquad (45.7)
$$

Thus the integer-valued random variable $|\mathcal{P}_1(\mathcal{D})|1_{\mathcal{D}}$ satisfies a Gaussian bound. We can apply lemma A.1 with

$$
A\,=\,6d2^dn^d\,\geq 3\,,\qquad B^2\,=\,\frac{9d2^dn^d}{2p^2}\,\geq 1\,,
$$

and we get

$$
E\big(|\mathcal{P}_1(\mathcal{D})|1_{\mathcal{D}}\big) \,\leq\, \frac{6d2^d}{p}n^{d/2}\sqrt{\ln(6d2^dn^d)}\,. \qquad (45.8)
$$

In fact, we would rather need to control the second moment of $|\mathcal{P}_1(\mathcal{D})|1_{\mathcal{D}}$. It follows from the exercise 2.3.8, question b) of [126] that the Gaussian bound yields an adequate control on all the moments, so that the inequality (45.7) implies that

$$E\big(|\mathcal{P}_1(\mathcal{D})|^2 1_{\mathcal{D}}\big) \,=\, O\Big(n^d\sqrt{\ln n}\Big)\,,$$

with the understanding that the constant associated with the notation $O(\cdot)$ depends on $d, p$ only. It follows then from inequality (45.6) (where $n$ has to be replaced by $n^d$) that

$$E\big(\big|\mathcal{S}^{2,-1}(\mathcal{D})\big|\big) \,=\, O\Big(n^d\sqrt{\ln n}\Big)\,. \tag{45.9}$$

What can we infer from (45.9) concerning the set of the second pivots $\mathcal{P}_2(\mathcal{D})$? We know from corollary 44.2 applied with $\ell = 2$ that

$$\pi\big(\mathcal{S}^{2,-1}(\mathcal{D}) \,\cup\, \mathcal{S}^{2,+1}(\mathcal{D})\big) \,=\, \mathcal{P}_1(\mathcal{D}) \cup \mathcal{P}_2(\mathcal{D})\,. \tag{45.10}$$

We have already a control on the expectation of $|\mathcal{P}_1(\mathcal{D})|1_{\mathcal{D}}$ in inequality (45.8). It follows from the equivalence (45.3) that

$$\forall \omega \in \mathcal{D} \qquad \pi\big(\mathcal{S}^{2,+1}(\mathcal{D},\omega)\big) \,=\, \mathcal{P}_1(\mathcal{D},\omega)\,. \tag{45.11}$$

Combining the formulas (45.8), (45.10) and (45.11), we obtain that

$$E\big(|\mathcal{P}_2(\mathcal{D})|1_{\mathcal{D}}\big) \,=\, E\Big(|\pi\big(\mathcal{S}^{2,-1}(\mathcal{D})\big)|1_{\mathcal{D}}\Big) + O\Big(n^{d/2}\sqrt{\ln n}\Big)\,.$$

At this point, we would like to use the inequality (45.9) on $E\big(\big|\mathcal{S}^{2,-1}(\mathcal{D})\big|\big)$ in order to control $E\big(|\pi\big(\mathcal{S}^{2,-1}(\mathcal{D})\big)\big|\big)$. Yet this is far from obvious. There are some situations in which the inequality (45.9) gives a useful control. Let us focus on the third question of section 30, namely the disconnection in a cuboid. Since we work here with the bond percolation model, the previous results on the pivots and second pivots should be translated in terms of bonds, but this can be done smoothly. We consider a percolation configuration realizing the disconnection event $\mathcal{D}$. In this context, a pair of bonds $\{\, e = \langle a, b\rangle, f = \langle c, d\rangle \,\}$ belongs to $\mathcal{S}^{2,-1}(\mathcal{D})$ if, in the graph with vertex set $D$ and edge set $\mathbb{E}^d(D)\setminus\{\, e, f \,\}$, up to a change of order in the names of the vertices $a, b, c, d$, the following events occur:

$$L \not\longleftrightarrow R\,,\quad L \longleftrightarrow a\,,\quad L \not\longleftrightarrow c\,,\quad b \longleftrightarrow c\,,\quad c \not\longleftrightarrow R\,,\quad d \longleftrightarrow R\,.$$

In particular, in the pair $\{\, e, f \,\}$, there is exactly one bond connected to the left side $L$ and exactly one bond connected to the right side $R$. The bonds of $\pi\big(\mathcal{S}^{2,-1}(\mathcal{D})\big)$ which are connected to $L$ (respectively $R$) are called the left (respectively right) second pivots, and we denote the set of these bonds by $\mathcal{P}_2^L(\mathcal{D})$ (respectively $\mathcal{P}_2^R(\mathcal{D})$). With this notation, we have

$$\pi\big(\mathcal{S}^{2,-1}(\mathcal{D})\big) \,=\, \mathcal{P}_2^L(\mathcal{D}) \cup \mathcal{P}_2^R(\mathcal{D})\,,$$

and everything boils down to controlling the cardinality of $\mathcal{P}_2^L(\mathcal{D})$ and $\mathcal{P}_2^R(\mathcal{D})$. Consider now a percolation configuration $\omega$ realizing the event $\mathcal{D}$ in which all the bonds of $\mathcal{S}^{2,-1}(\mathcal{D})$ belong to the edge boundary of one single cluster $C^*$. In this specific case, every bond in $\mathcal{P}_2^L(\mathcal{D},\omega)$ has one endpoint which is connected to one endpoint of every bond in $\mathcal{P}_2^R(\mathcal{D},\omega)$, and this implies that $\mathcal{S}^{2,-1}(\mathcal{D},\omega)$ has a simple product structure:

$$\mathcal{S}^{2,-1}(\mathcal{D},\omega) \,=\, \Big\{ \{ e,f \} : e\in \mathcal{P}_2^L(\mathcal{D},\omega),\, f\in \mathcal{P}_2^R(\mathcal{D},\omega) \Big\}\,,$$

whence

$$\big|\mathcal{S}^{2,-1}(\mathcal{D},\omega)\big| \,=\, \big|\mathcal{P}_2^L(\mathcal{D},\omega)\big| \times \big|\mathcal{P}_2^R(\mathcal{D},\omega)\big|\,. \tag{45.12}$$

Imagine that, in addition, we had a control guaranteeing that a positive proportion, say 1/4, of the second pivots is of the left type, and a similar control for the right type, and imagine that we could use these controls when taking the expectation of identity (45.12), then we would have

$$E\Big(\big|\mathcal{S}^{2,-1}(\mathcal{D},\omega)\big|1_{\mathcal{D}}\Big) \,\geq\, \frac{1}{16}E\Big(\big|\mathcal{P}_2(\mathcal{D},\omega)\big|^2 1_{\mathcal{D}}\Big)\,. \tag{45.13}$$

Using this inequality (45.13) in conjunction with the inequality (45.9), we would indeed conclude that $E\big(\big|\mathcal{P}_2(\mathcal{D},\omega)\big|1_{\mathcal{D}}\big)$ is $O\big(n^{d/2}\sqrt{\ln n}\big)$. We would like to obtain a similar control for the higher pivots. For any $\ell\in\{1,\dots,n\}$, we have, by proposition 43.2,

$$\Big|\Big(\frac{d}{dp}\Big)^{\ell}P_p(\mathcal{D})\Big| \,\leq\, \sqrt{P_p(\mathcal{D})\binom{n}{\ell}\frac{1}{p^{\ell}(1-p)^{\ell}}} \,=\, O\big(n^{\ell d/2}\big)\,. \tag{45.14}$$

With some magical induction on $\ell$, we could hope to get an inequality generalizing inequality (45.13) for the pivots of order $\ell$. Together with the generalized Margulis formula (41.8), and some control on $E\big(\big|\mathcal{S}^{\ell,s}(\mathcal{D})\big|1_{\mathcal{D}}\big)$, this might yield

$$E\Big(\big|\mathcal{P}_{\ell}(\mathcal{D})\big|^{\ell}1_{\mathcal{D}}\Big) \,=\, O\Big(\Big|\Big(\frac{d}{dp}\Big)^{\ell}P_p(\mathcal{D})\Big|\Big)\,. \tag{45.15}$$

Of course the constant associated with $O(\cdot)$ might depend on $\ell$. Finally, combining the estimates (45.14) and (45.15), we would conclude that

$$E\big(\big|\mathcal{P}_{\ell}(\mathcal{D})\big|1_{\mathcal{D}}\big) \,=\, O\big(n^{d/2}\sqrt{\ln n}\big)\,.$$

Unfortunately, it took a lot of imagination to get to inequality (45.13), let alone inequality (45.15). To begin with, in a typical configuration $\omega$ realizing the event $\mathcal{D}$, the bonds of $\mathcal{S}^{2,-1}(\mathcal{D})$ belong to the edge boundary of several clusters. Denoting by $(C_i, i\in I)$ the collection of these clusters, the inequality (45.12) should be replaced by the less advantageous inequality

$$\big|\mathcal{S}^{2,-1}(\mathcal{D},\omega)\big| \,=\, \sum_{i\in I}\big|\mathcal{P}_2^L(\mathcal{D},\omega)\cap\Delta_D C_i\big| \times \big|\mathcal{P}_2^R(\mathcal{D},\omega)\cap\Delta_D C_i\big|\,. \tag{45.16}$$

We could still hope that the major contribution to the sum in (45.16) comes from one or at least a few clusters. Unfortunately, we did not manage to obtain a useful quantitative estimate with this strategy.

## 45.2 Talagrand and the second pivots

In his influential work [125], Talagrand derived a beautiful inequality on the correlation on increasing sets for the uniform probability measure. In the course of the proof, Talagrand has obtained several inequalities which involves the second pivots. Let us present briefly these important inequalities. We adapt partly Talagrand's notation to our framework. We denote by $\mu = P_{1/2}$ the uniform measure on $\{0,1\}^n$ and we use the random variables $X_i$, $1 \leq i \leq n$, defined by

$$\forall i \in \{\,1,\dots,n\,\} \quad \forall \omega \in \{0,1\}^n \qquad X_i(\omega) \,=\, 4\omega(i) - 2\,.$$

So the variables $X_i$ are the double of the variables $r_i$ employed by Talagrand. The first result proved by Talagrand in [125] is the following inequality: there exists a universal constant $K$ such that, for any subset $A$ of $\{0,1\}^n$, we have

$$\sum_{1\leq i\leq n} \Big(\int_A X_i\, d\mu\Big)^2 \,\leq\, K\mu(A)^2 \ln \frac{e}{\mu(A)}\,. \tag{45.17}$$

In subsection 40.2, we already restated and discussed this inequality (see proposition 40.7), and we proved also the corresponding inequality for the biased measure $P_p$ (see proposition 40.6). Yet Talagrand went much further. Noticing that inequality (45.17) holds for any set $A$, he applied it to the event $\{\,k \in \mathcal{P}(A)\,\}$, where $k$ is an index in $\{\,1,\dots,n\,\}$ and $A$ is an increasing set. The variable $X_k$ and the event $\{\,k \in \mathcal{P}(A)\,\}$ are independent, therefore

$$E\big(X_k 1_{\{\,k\in\mathcal{P}(A)\,\}}\big) \,=\, 0\,.$$

Using formula (39.6) in the reverse direction, we can rewrite the expectation for $i \neq k$ as follows:

$$\forall i \neq k \qquad E\big(X_i 1_{\{\,k\in\mathcal{P}(A)\,\}}\big) \,=\, E\big(X_i X_k 1_A\,\big)\,,$$

and the previous inequality becomes

$$\sum_{\substack{1\leq i\leq n\\ i\neq k}} \Big(\int_A X_i X_k\, d\mu\Big)^2 \,\leq\, K\Big(\mu\big(k \in \mathcal{P}(A)\big)\Big)^2 \ln\Big(\frac{e}{\mu\big(\,k \in \mathcal{P}(A)\,\big)}\Big)\,.$$

Talagrand next summed this inequality over $k$ to get

$$\sum_{\substack{1\leq i,k\leq n\\ i\neq k}} \Big(\int_A X_i X_k\, d\mu\Big)^2 \,\leq\, K \sum_{1\leq k\leq n} \Big(\mu\big(k \in \mathcal{P}(A)\big)\Big)^2 \ln\Big(\frac{e}{\mu\big(\,k \in \mathcal{P}(A)\,\big)}\Big)\,.$$

This inequality is certainly relevant to our problem, because it involves the second pivots, albeit in a complicated way. Finally, in one of those feats of strength for which only he has the secret, Talagrand proved the following much stronger inequality.

**Proposition 45.1** (Proposition 2.3 of [125]). *There exists a universal constant $K$ such that, for any increasing subset $A$ of $\{0,1\}^n$, we have*

$$\sum_{\substack{1\leq i,k\leq n\\ i\neq k}} \Big(\int_A X_i X_k\,d\mu\Big)^2 \,\leq\, K \sum_{1\leq k\leq n} \mu\big(k\in\mathcal{P}(A)\big)^2 \ln \frac{K}{\sum_{1\leq k\leq n}\mu\big(k\in\mathcal{P}(A)\big)^2}\,. \tag{45.18}$$

The considerable improvement lies in the presence of the sum inside the logarithm. Talagrand's proof is beautiful and mysterious. It rests on a complicated intermediate result which involves a sort of bootstrap mechanism. The inequality (45.18) attracted the interest of several researchers and led to numerous unexpected developments. For instance, the whole theory of noise sensitivity of Benjamini, Kalai and Schramm [15] was born with a generalization of this inequality. An important question was whether inequality (45.18) holds for the product measure $P$ for $p\neq 1/2$. It is known that the inequalities involving the uniform measure can in principle be generalized to the Bernoulli product measure. Keller [79] has designed a procedure to deduce automatically some inequalities for the biased measure from inequalities involving the uniform measure. Anyway, the natural approach consists in adapting the proofs written for the uniform measure, as we did in this paper for the first inequality of Talagrand. However, this might turn out to be very difficult in some cases. After more than twenty years, Talagrand's inequality has been extended to the biased case by Keller and Kindler [80], with a new proof which is simpler than the original proof of Talagrand. We need to introduce some notation to state their result. For $i\in\{\,1,\dots,n\,\}$, we define the normalized random variable

$$Y_i \,=\, \sqrt{p(1-p)}X_i \,=\, \begin{cases} -\sqrt{\frac{1-p}{p}} & \text{if } \omega(i)=0\,,\\ \sqrt{\frac{p}{1-p}} & \text{if } \omega(i)=1\,,\end{cases}$$

and for $S$ a subset of $\{\,1,\dots,n\,\}$, we define the Walsh product

$$Y_S \,=\, \prod_{i\in S} Y_i\,.$$

For $f:\{0,1\}^n\to\{0,1\}$ a Boolean function and $i\in\{1,\dots,n\}$, we define the influence of the $i$-th coordinate on $f$ as

$$I_i(f) \,=\, P_p\big(\{\,\omega\in\Omega : f(\omega^i)\neq f(\omega_i)\,\}\big)\,,$$

and the Walsh coefficient of $f$ associated to a subset $S$ of $\{\,1,\dots,n\,\}$ as

$$\widehat{f}(S) \,=\, E\big(fY_S\big) \,=\, \sum_{\omega\in\Omega}\Big(f(\omega)Y_S(\omega)P_p(\omega)\Big)\,.$$

Finally, the hypercontractivity constant $B(p)$ is given by

$$B(p) \,=\, \frac{\dfrac{1-p}{p}-\dfrac{p}{1-p}}{2\ln\dfrac{1-p}{p}}\,.$$

We restate next the inequality proved by Keller and Kindler (see [80], lemma 6).

**Lemma 45.2.** *Let* $f:\{0,1\}^n\to\{0,1\}$*, and denote*

$$\mathcal{W}(f) \,=\, p(1-p)\sum_{1\leq i\leq n} I_i(f)^2\,.$$

*For all* $d\geq 2$*, if*

$$\mathcal{W}(f)\leq \exp(-2(d-1)),$$

*then we have*

$$\sum_{|S|=d}\widehat{f}(S)^2\leq \frac{5e}{d}\left(\frac{2B(p)e}{d-1}\right)^{d-1}\mathcal{W}(f)\left(\log\left(\frac{d}{\mathcal{W}(f)}\right)\right)^{d-1}.$$

Unfortunately, the direct application of proposition 45.1 or lemma 45.2 to our disconnection problem does not lead to a serious progress, as far as the conjecture $\theta(p_c,\mathbb{Z}^d)=0$ is concerned (although these inequalities were instrumental to prove that the disconnection between two opposite sides of a square in two-dimensional critical percolation is noise–sensitive [15]).

Since our primary goal is to develop a quantitative control on the second pivotal bonds, we must examine in detail the proofs of proposition 45.1 and lemma 45.2 in order to see if these proofs can be improved in the specific case of a disconnection event in percolation. Notice that Keller and Kindler managed to prove the inequality for any boolean function, not necessarily non–decreasing. However, upon inspection, it seems that the proof of Talagrand would work also for any boolean function! In fact, the proof of Talagrand rests on a beautiful but very technical estimate. Keller and Kindler managed to bypass this estimate, thereby simplifying the overall argument. However, there is some hope to adapt this estimate in our context in order to get a better quantitative result. So let us state the key technical estimate of Talagrand (rewritten with our notation).

**Lemma 45.3** (Lemma 3.1 of [125])**.** *Consider a partition* $I,J$ *of* $\{\,1,\dots,n\,\}$*. Consider an increasing set* $A$ *and* $S>0$*. Consider*

$$L\,=\,\left\{k\in J:\left(\sum_{i\in I}\Big(\int_{\{\,k\in\mathcal{P}(A)\,\}}X_i\,d\mu\Big)^2\right)^{1/2}\,\geq\,S\mu\big(k\in\mathcal{P}(A)\big)\right\}.$$

*There exists a universal constant* $C>0$ *such that*

$$\sum_{k\in L}\mu\big(k\in\mathcal{P}(A)\big)^2\,\leq\,C\exp\left(-\frac{S^2}{C}\right)\,. \tag{45.19}$$

To obtain his refinement of the classical FKG inequality, Talagrand applies lemma 45.3 and performs an average over all possible pairs of sets $I,J$. The hope, in our situation, would be to apply it to a more specific choice of $I,J$, which is relevant for the disconnection event. For instance, a tantalising choice would be to take for one of the two sets the domain discovered when exploring all the clusters touching the first set $W$. We should therefore try to prove a version of lemma 45.3 for some specific random sets. We start this program in the next subsection.

## 45.3 The pivots in a random domain

We try here to extend the inequalities controlling the pivotal bonds to the case of a random subset of the index set. We hope that this might be useful to control the second pivotal bonds.

### 45.3.1 Extension of Margulis-Russo formula

We will consider random subsets $\mathcal{D} = \mathcal{D}(\omega)$ of $\{1, \dots, n\}$ which satisfy the following hypothesis.

**Hypothesis 45.4.** *The random subset $\mathcal{D} = \mathcal{D}(\omega)$ of $\{1, \dots, n\}$ is not empty and moreover*

$$\forall i \in \{1, \dots, n\} \qquad \{i \in \mathcal{D}\} \in \sigma\big(X_j, j \neq i\big) . \tag{45.20}$$

It is not convenient at all to take the derivative with respect to the parameters $p$ associated to the sites belonging to a random subset, so we rather generalize the second part of the Margulis-Russo formula (39.3).

**Proposition 45.5.** *Let $\mathcal{D} = \mathcal{D}(\omega)$ be a random subset of $\{1, \dots, n\}$ satisfying hypothesis 45.4. We define*

$$S(\mathcal{D}) = \sum_{i \in \mathcal{D}} X_i .$$

*For any event $A \subset \Omega$, we have*

$$E\big(1_A S(\mathcal{D})\big) = E\big(|\mathcal{P}^+(A) \cap \mathcal{D}|\big) - E\big(|\mathcal{P}^-(A) \cap \mathcal{D}|\big) . \tag{45.21}$$

$$E\big(S(\mathcal{D}) \,|\, A\big) = \frac{1}{p} E\big(|\mathcal{P}^+(A) \cap \mathcal{D}| \,\big|\, A\big) - \frac{1}{1-p} E\big(|\mathcal{P}^-(A) \cap \mathcal{D}| \,\big|\, A\big) . \tag{45.22}$$

*Proof.* We write

$$E\big(1_A S(\mathcal{D})\big) = \sum_{1 \leq i \leq n} E\big(1_A 1_{\{i \in \mathcal{D}\}} X_i\big) , \tag{45.23}$$

and we apply formula (39.6) to get

$$E\big(1_A 1_{\{i \in \mathcal{D}\}} X_i\big) = E\big(\delta_i 1_{A \cap \{i \in \mathcal{D}\}}\big) . \tag{45.24}$$

Now, for any $\omega \in \Omega$,

$$\delta_i 1_{A \cap \{i \in \mathcal{D}\}}(\omega) = 1_A(\omega^i) 1_{\{i \in \mathcal{D}\}}(\omega^i) - 1_A(\omega_i) 1_{\{i \in \mathcal{D}\}}(\omega_i) .$$

By the hypothesis (45.20) on the random set $\mathcal{D}$, we have, for any $\omega \in \Omega$,

$$1_{\{i \in \mathcal{D}\}}(\omega^i) = 1_{\{i \in \mathcal{D}\}}(\omega_i) = 1_{\{i \in \mathcal{D}\}}(\omega) ,$$

whence

$$\begin{aligned} \delta_i 1_{A \cap \{i \in \mathcal{D}\}}(\omega) &= 1_{\{i \in \mathcal{D}\}}(\omega)\big(1_A(\omega^i) - 1_A(\omega_i)\big) = 1_{\{i \in \mathcal{D}\}}(\omega)\, \delta_i 1_A(\omega) \\ &= 1_{\{i \in \mathcal{D}\}}(\omega)\big(1_{\{i \in \mathcal{P}^+(A)\}} - 1_{\{i \in \mathcal{P}^-(A)\}}\big) . \end{aligned} \tag{45.25}$$

Reporting the two identities (45.24) and (45.25) in (45.23), we obtain

$$E\big(1_A S(\mathcal{D})\big) \;=\; \sum_{1\le i\le n} E\big(1_{\{i\in\mathcal{D}\cap\mathcal{P}^+(A)\}} - 1_{\{i\in\mathcal{D}\cap\mathcal{P}^-(A)\}}\big)$$
$$=\; E\big(|\mathcal{P}^+(A)\cap\mathcal{D}|\big) - E\big(|\mathcal{P}^-(A)\cap\mathcal{D}|\big)\,, \quad (45.26)$$

and this yields formula (45.21). In addition, for any $i\in\{\,1,\dots,n\,\}$, we have

$$E\big(1_{\{i\in\mathcal{D}\cap\mathcal{P}^+(A)\}}\big) \;=\; P\big(i\in\mathcal{D}\cap\mathcal{P}^+(A)\big)\,. \quad (45.27)$$

Thanks to the hypothesis (45.20) and the definition of a pivotal site, the event $\{\,i\in\mathcal{D}\cap\mathcal{P}^+(A)\,\}$ belongs to the $\sigma$–field $\sigma\big(\omega(j), j\neq i\big)$, therefore

$$P\big(i\in\mathcal{D}\cap\mathcal{P}^+(A), A\big) \;=\; P\big(i\in\mathcal{D}\cap\mathcal{P}^+(A), \omega(i)=1\big)$$
$$=\; p\,P\big(i\in\mathcal{D}\cap\mathcal{P}^+(A)\big)\,. \quad (45.28)$$

An analogous computation can be done for the negative pivotal sites. Plugging (45.28) and (45.27) into (45.26), we obtain

$$E\big(1_A S(\mathcal{D})\big) \;=\; \sum_{1\le i\le n} \Big(\frac{1}{p}P\big(i\in\mathcal{D}\cap\mathcal{P}^+(A), A\big) - \frac{1}{1-p}P\big(i\in\mathcal{D}\cap\mathcal{P}^-(A), A\big)\Big)\,.$$

We simply divide by $P(A)$ to obtain formula (45.22). □

We shall need additional assumptions on the random set $\mathcal{D}$ in order to derive some useful probabilistic estimates. We make the following hypothesis.

**Hypothesis 45.6.** *There exists a sequence $(Y_i)_{1\le i\le n}$ of i.i.d. random variables and a random stopping time $\mathcal{T}$ associated to the sequence $(Y_i)_{1\le i\le n}$ such that*

$$\forall i\in\{\,1,\dots,n\,\}\qquad P(Y_i=1)\;=\;p\,,\quad P(Y_i=0)\;=\;1-p\,, \quad (45.29)$$

*and moreover the two random vectors*

$$\Big(|\mathcal{D}|, \sum_{i\in\mathcal{D}} X_i\Big) \quad \textit{and} \quad \Big(\mathcal{T}, \sum_{1\le i\le\mathcal{T}} \Big(\frac{Y_i}{p} - \frac{1-Y_i}{1-p}\Big)\Big)$$

*have the same joint distribution.*

The typical situation that we have in mind is an exploratory process. Suppose that, for each $k\in\{\,1,\dots,n\,\}$, we are given a deterministic function $f_k$ from $\big(\{\,1,\dots,n\,\}\times\{0,1\}\big)^k$ to $\{\,0,\dots,n\,\}$ satisfying

$$\forall (i_0,\dots,i_{k-1},\varepsilon_0,\dots,\varepsilon_{k-1})\in\big(\{\,1,\dots,n\,\}\times\{0,1\}\big)^k$$
$$f_k(i_0,\dots,i_{k-1},\varepsilon_0,\dots,\varepsilon_{k-1})\notin\{\,i_0,\dots,i_{k-1}\,\}\,. \quad (45.30)$$

We define next a sequence of random indices as follows. Let $i_0$ be an index in $\{\,1,\dots,n\,\}$. We set $I_0=i_0$ and for $k\ge 1$, we set

$$I_k \;=\; f_k\big(I_0,\dots,I_{k-1}, p(1-p)X_{I_0}+p,\dots,p(1-p)X_{I_{k-1}}+p\big)$$

and $\mathcal{T} = k$ if $I_k = 0$. The sequence terminates at the random index $\mathcal{T}$, which corresponds to the first index $k$ such that $I_k = 0$. Since the indices of the sequence $I_0, I_1, \cdots$ are distinct by construction, we must have $1 \leq \mathcal{T} \leq n$. Finally, we define the random set $\mathcal{D}$ as

$$\mathcal{D} = \{ I_0, \dots, I_{\mathcal{T}-1} \}.$$

The next proposition is the generalization of proposition 39.5 to a random set.

**Proposition 45.7.** *Let $\mathcal{D}$ be a random subset of $\{ 1, \dots, n \}$ satisfying the hypotheses 45.4 and 45.6. For any event $A \subset \Omega$, we have*

$$\Big| E\big(|\mathcal{P}^+(A) \cap \mathcal{D}|\big) - E\big(|\mathcal{P}^-(A) \cap \mathcal{D}|\big) \Big| \leq \sqrt{\frac{E(|\mathcal{D}|)}{p(1-p)}}. \tag{45.31}$$

*Proof.* The Cauchy–Schwarz inequality gives

$$\big| E\big(S(\mathcal{D}) 1_A\big) \big| \leq \sqrt{E\big((S(\mathcal{D}))^2\big) P_p(A)}. \tag{45.32}$$

Let us define

$$\forall k \in \{ 1, \dots, n \} \qquad S_k = \sum_{i=1}^{k} \Big( \frac{Y_i}{p} - \frac{1 - Y_i}{1-p} \Big).$$

The sequence $(S_k)_{1 \leq k \leq n}$ is a martingale with respect to the variables $(Y_k)_{1 \leq k \leq n}$. From hypothesis 45.6, we have

$$E\big((S(\mathcal{D}))^2\big) = E\big((S_{\mathcal{T}})^2\big). \tag{45.33}$$

We use the classical martingale method to compute $E\big((S_{\mathcal{T}})^2\big)$. We have

$$\forall i \in \{ 1, \dots, n \} \qquad E\left( \Big( \frac{Y_i}{p} - \frac{1 - Y_i}{1-p} \Big)^2 \right) = \frac{1}{p(1-p)},$$

and we introduce the sequence defined by

$$\forall k \in \{ 1, \dots, n \} \qquad \widetilde{S}_k = (S_k)^2 - \frac{k}{p(1-p)}.$$

We check that the sequence $(\widetilde{S}_k)_{1 \leq k \leq n}$ is also a martingale with respect to the variables $(Y_k)_{1 \leq k \leq n}$. Since $\mathcal{T}$ is a bounded stopping time, by Doob's optional stopping theorem, we have $E\big(\widetilde{S}_{\mathcal{T}}\big) = 0$, therefore

$$E\big((S_{\mathcal{T}})^2\big) = \frac{1}{p(1-p)} E(\mathcal{T}). \tag{45.34}$$

The hypothesis 45.6 implies that $E(\mathcal{T}) = E(|\mathcal{D}|)$, hence we conclude from the formulas (45.32), (45.33) and (45.34) that

$$\big| E\big(S(\mathcal{D}) 1_A\big) \big| \leq \sqrt{\frac{E(|\mathcal{D}|) P_p(A)}{p(1-p)}}. \tag{45.35}$$

Inequalities (45.21) and (45.35) together yield inequality (45.31). □

### 45.3.2 Extension of the exponential control

We will next try to generalize the deviations inequality of proposition 40.2 to the case of a random set.

**Theorem 45.8.** *Let $\mathcal{D}=\mathcal{D}(\omega)$ be a random subset of $\{1,\dots,n\}$ satisfying the hypotheses 45.4 and 45.6. For any event $A\subset\Omega$, we have*

$$P_p(A)\ \le\ 4n^2\exp\left(-\left(\frac{(1-p)E\big(|\mathcal{P}^+(A)\cap\mathcal{D}|\,\big|\,A\big)-pE\big(|\mathcal{P}^-(A)\cap\mathcal{D}|\,\big|\,A\big)}{2E\big(\sqrt{|\mathcal{D}|}\,\big|\,A\big)}\right)^2\right). \tag{45.36}$$

*Proof.* We use essentially the same strategy as for the proof of proposition 40.2. There is an extra difficulty due to the fact that the number of terms involved in $S(\mathcal{D})$ is random, and the computations become more complicated. The typical fluctuations of $S(\mathcal{D})$ are of order $\sqrt{|\mathcal{D}|}$, so we shall split the expectation according to whether $\big|S(\mathcal{D})\big|$ is larger than $t\sqrt{|\mathcal{D}|}$ or not, where $t$ is an adequately chosen parameter. So, let $t>0$ and let us write

$$\begin{aligned}
E\big(|S(\mathcal{D})|\,\big|\,A\big)\ &=\ \frac{1}{P_p(A)}\int_A |S(\mathcal{D})|\,dP_p\\
&= \frac{1}{P_p(A)}\int_{A\cap\{\,|S(\mathcal{D})|\geq t\sqrt{|\mathcal{D}|}\,\}}|S(\mathcal{D})|\,dP_p+\frac{1}{P_p(A)}\int_{A\cap\{\,|S(\mathcal{D})|< t\sqrt{|\mathcal{D}|}\,\}}|S(\mathcal{D})|\,dP_p\\
&\leq\ \frac{1}{P_p(A)}\int_{\{\,|S(\mathcal{D})|\geq t\sqrt{|\mathcal{D}|}\,\}}|S(\mathcal{D})|\,dP_p+tE\big(\sqrt{|\mathcal{D}|}\,\big|\,A\big)\,.
\end{aligned} \tag{45.37}$$

Using Fubini's theorem, we rewrite the first integral in (45.37) as follows:

$$\begin{aligned}
&\int_{\{\,|S(\mathcal{D})|\geq t\sqrt{|\mathcal{D}|}\,\}}|S(\mathcal{D})|\,dP_p\\
&\qquad=\int_{\{\,|S(\mathcal{D})|\geq t\sqrt{|\mathcal{D}|}\,\}}\left(\int_0^{+\infty}\sqrt{|\mathcal{D}|}1_{\{\,u\sqrt{|\mathcal{D}|}\leq|S(\mathcal{D})|\,\}}\,du\right)dP_p\\
&\qquad\qquad\leq\ \sqrt{n}\int_0^{+\infty}P\big(|S(\mathcal{D})|\geq\max(u,t)\sqrt{|\mathcal{D}|}\big)\,du\,.
\end{aligned} \tag{45.38}$$

To control the probability inside the integral, we decompose it according to the value of $|\mathcal{D}|$ and we use the hypothesis (45.29). For any $v>0$, we have

$$\begin{aligned}
P\big(|S(\mathcal{D})|\geq v\sqrt{|\mathcal{D}|}\big)\ &=\ \sum_{m=1}^n P\big(|S(\mathcal{D})|\geq v\sqrt{|\mathcal{D}|},\,|\mathcal{D}|=m\big)\\
&=\ \sum_{m=1}^n P\left(\Big|\sum_{i=1}^m\Big(\frac{Y_i}{p}-\frac{1-Y_i}{1-p}\Big)\Big|\geq v\sqrt{m},\,\mathcal{T}=m\right)\\
&\leq\ \sum_{m=1}^n P\left(\Big|\sum_{i=1}^m\Big(\frac{Y_i}{p}-\frac{1-Y_i}{1-p}\Big)\Big|\geq v\sqrt{m}\right).
\end{aligned} \tag{45.39}$$

Applying Hoeffding's inequality [70], we get

$$\forall v > 0 \quad \forall m \in \{ 1, \dots, n \}$$
$$P\left( \Big| \sum_{i=1}^m \Big( \frac{Y_i}{p} - \frac{1-Y_i}{1-p} \Big) \Big| \geq v\sqrt{m} \right) \leq 2\exp\big( -2p^2(1-p)^2v^2 \big) . \quad (45.40)$$

Reporting (45.40) into (45.39), we obtain

$$\forall v > 0 \qquad P\big( |S(\mathcal{D})| \geq v\sqrt{|\mathcal{D}|} \big) \leq 2n\exp\big( -2p^2(1-p)^2v^2 \big) . \qquad (45.41)$$

Substituting (45.41) in (45.38) with $v = \max(u,t)$, we get

$$\int_{\{ |S(\mathcal{D})| \geq t\sqrt{|\mathcal{D}|} \}} |S(\mathcal{D})| \, dP_p \leq 2n^{3/2} \int_0^{+\infty} \exp\big( -2p^2(1-p)^2 \max(u,t)^2 \big) \, du \, .$$

To bound the integral, we use the simple inequality $2\max(u,t)^2 \geq u^2 + t^2$ and the Gauss integral, and we obtain

$$\int_{\{ |S(\mathcal{D})| \geq t\sqrt{|\mathcal{D}|} \}} |S(\mathcal{D})| \, dP_p \leq \frac{n^{3/2}\sqrt{\pi}}{p(1-p)} \exp\big( -p^2(1-p)^2t^2 \big) . \qquad (45.42)$$

Reporting (45.42) in (45.37), we conclude that

$$E\big( |S(\mathcal{D})| \,\big|\, A \big) \leq \frac{n^{3/2}\sqrt{\pi}}{P_p(A)p(1-p)} \exp\big( -p^2(1-p)^2t^2 \big) + tE\big( \sqrt{|\mathcal{D}|} \,\big|\, A \big) .$$

We choose now

$$t = \frac{1}{p(1-p)} \sqrt{ \ln \frac{2n^{3/2}\sqrt{\pi}}{P_p(A)E\big( \sqrt{|\mathcal{D}|} \,\big|\, A \big)} } \geq \frac{1}{p(1-p)} \sqrt{\ln 2n\sqrt{\pi}} \geq \frac{1}{p(1-p)} ,$$

and we conclude that

$$E\big( |S(\mathcal{D})| \,\big|\, A \big) \leq \frac{2}{p(1-p)} \sqrt{ \ln \frac{2n^{3/2}\sqrt{\pi}}{P_p(A)E\big( \sqrt{|\mathcal{D}|} \,\big|\, A \big)} } E\big( \sqrt{|\mathcal{D}|} \,\big|\, A \big) . \qquad (45.43)$$

Using the facts that $|\mathcal{D}| \geq 1$ and $\big| E\big( S(\mathcal{D}) \,\big|\, A \big) \big| \leq E\big( |S(\mathcal{D})| \,\big|\, A \big)$, we rewrite inequality (45.43) in exponential form as

$$P_p(A) \leq 2n^{3/2}\sqrt{\pi} \exp\left( -\frac{1}{4} \left( \frac{p(1-p)E\big( S(\mathcal{D}) \,\big|\, A \big)}{E\big( \sqrt{|\mathcal{D}|} \,\big|\, A \big)} \right)^2 \right) .$$

We use finally the identity (45.22) of proposition 45.5 to obtain the inequality (45.36) of the theorem. □

## 45.4 Talagrand's lemma for the disconnection event?

The inequalities of Talagrand are very general, which is their strength, but they are too weak for our percolation problem. So we should try to develop specific variants of these inequalities for the disconnection event in percolation. The hope is that, by taking advantage of the geometric structure of the problem, we will end up with a stronger control on the second pivotal bonds. The first step consists in extending Talagrand's lemma 45.3 to the biased case. Although we do not provide the details, this step does not present any unexpected difficulty. The more challenging part is to find a sensible version of the inequality (45.19), without even addressing the question of the related proof. The main idea would be to work with two adequate random sets $I, J$ (remember that Talagrand proved the lemma 45.3 with two fixed deterministic sets, and subsequently performed an average over all possible choices of sets $I, J$ to derive his main inequality on the correlation of increasing events). The natural candidates for $I, J$ would be the regions revealed by the intertwined explorations. Now these two regions overlap precisely on the set of the pivotal bonds, thus we must find an adequate rule for setting the pivotal bonds in one of these two regions. Let us try to precise this scheme. For simplicity, we focus on the event which is the disconnection between two opposite faces in the cubic box $\Lambda(n)$. For the sake of the discussion, we consider a three-dimensional cubic box, and we take for $W$ its left face and for $W'$ its right face. We first run an exploration starting from the left face $W$. This exploration will find the random set $\mathcal{L}$ of the vertices of the box that are connected to $W$, and it will have revealed all the bonds emanating from $\mathcal{L}$, that is the set of bonds

$$\mathbb{L} \,=\, \mathbb{E}^d\big(\mathcal{L}\cup\partial^{out}_{\Lambda(n)}\mathcal{L}\big)\,. \tag{45.44}$$

At this point, the remaining bonds of $\Lambda(n)$, i.e., $\mathbb{E}^d(\Lambda(n))\setminus\mathbb{L}$, which have not been examined, are still distributed according to the Bernoulli product measure. Once we condition on the region $\mathcal{L}$, and we fix a bond $e$ in $\Delta\mathcal{L}$, the event $\big\{\, e\in\mathcal{P}(\mathcal{D})\,\big\}$ is an increasing event of the configuration of the bonds belonging to $\mathbb{E}^d(\Lambda(n))\setminus\mathbb{L}$. In particular, the expected number of pivotal bonds for the event $\big\{\, e\in\mathcal{P}(\mathcal{D})\,\big\}$ which belong to $\mathbb{E}^d(\Lambda(n))\setminus\mathbb{L}$ is $O(n^{d/2})$. Unfortunately, this bound is not good enough for our purpose, because in the end, we would have to sum it over the pivotal bonds of $\mathcal{I}(1)$ and we would get something like

$$\big|\mathcal{I}(2)\big| \,\leq\, \big|\mathcal{I}(1)\big|\,O(n^{d/2})\,. \tag{45.45}$$

We would like to prove a specific version of Talagrand's inequality (45.19), optimized for our problem. For the set $J$, we would take the random set $\mathbb{L}$ of all the bonds revealed by the exploration starting from the left face, and $I$ would be the set of the remaining bonds. So these sets would now be random, and we denote them by calligraphic letters: $\mathcal{J}=\mathbb{L}$, $\mathcal{I}=\mathbb{E}^d(\Lambda(n))\setminus\mathcal{J}$. The set $L$ would also be random, but since the letter $\mathcal{L}$ is already used, we denote it by $\mathcal{H}$, and it would be given by

$$\mathcal{H} \,=\, \left\{ f\in\mathcal{J} : \left(\sum_{e\in\mathcal{I}}\Big(E\big(X(e)1_{\{\,f\in\mathcal{P}(\mathcal{D})\,\}}\big)\Big)^2\right)^{1/2} \geq\, SP\big(f\in\mathcal{P}(\mathcal{D})\big)\right\}, \tag{45.46}$$

and the inequality (45.19) would become

$$\sum_{f\in\mathcal{H}} P\big(f\in\mathcal{P}(\mathcal{D})\big)^2\,\leq\,C\exp\left(-\frac{S^2}{C}\right)\,. \tag{45.47}$$

Yet the formulas (45.46) and (45.47) above do not make much sense when $\mathcal{I}$, $\mathcal{J}$ and $\mathcal{H}$ are random sets! So, in a first step, we will look for a weaker version of Talagrand's lemma, which has the advantage of being potentially extended to random sets. To obtain his inequalities, Talagrand made a systematic use of the following geometric trick. Whenever he deals with a sum of i.i.d. random variables, say $\sum_{i\in I}X_i$, he considered the weighted sum $\sum_{i\in I}\alpha_iX_i$, where the $\alpha_i$'s are non–negative real numbers such that $\sum_{i\in I}(\alpha_i)^2=1$. He performs the computations and in the end he takes a specific choice for the $\alpha_i$'s to get the best possible inequality. This is the method used in subsection 40.2, where we implemented the second step of Talagrand in order to enhance the initial exponential inequality. The benefit of this trick is typically to replace a linear sum of i.i.d. random variables by its Euclidean norm, and this yields the sums over the squared expectations in (45.46) or probabilities in (45.47). We start by stating the inequality that would have resulted from Talagrand's proof without the use of this geometric trick. Suppose for the time being that the sets $\mathcal{I},\mathcal{J},\mathcal{H}$ are non-random. Let us begin with inequality (45.47). By the Cauchy–Schwarz inequality, we have

$$\sum_{f\in\mathcal{H}} P\big(f\in\mathcal{P}(\mathcal{D})\big)\,\leq\,\sqrt{\big|\mathcal{H}\big|}\sqrt{\sum_{f\in\mathcal{H}} P\big(f\in\mathcal{P}(\mathcal{D})\big)^2}\,. \tag{45.48}$$

We can bound $\big|\mathcal{H}\big|$ by $\big|\mathbb{E}^d(\Lambda(n))\big|$, and the good point is that we can now incorporate the set $\mathcal{H}$ in the expectation:

$$\sum_{f\in\mathcal{H}} P\big(f\in\mathcal{P}(\mathcal{D})\big)\,=\,E\big(\big|\mathcal{P}(\mathcal{D})\cap\mathcal{H}\big|\big)\,.$$

We continue with the more complicated formula (45.46), and we rework the expectation it involves, as follows. For any bond $f$, we compute the expectation with the help of a conditioning on the random set $\mathbb{L}$ defined in (45.44):

$$E\Big(X(e)1_{\{\,f\in\mathcal{P}(\mathcal{D})\,\}}1_{e\notin\mathbb{L}}\Big)\,=\,E\bigg(1_{e\notin\mathbb{L}}E\Big(X(e)1_{\{\,f\in\mathcal{P}(\mathcal{D})\,\}}\,\Big|\,\mathbb{L}\Big)\bigg)\,. \tag{45.49}$$

The conditional distribution knowing $\mathbb{L}$ of an edge $e$ not belonging to $\mathbb{L}$ is Bernoulli with parameter $p$, therefore, using lemma 39.4, we have

$$\begin{aligned}1_{e\notin\mathbb{L}}E\Big(X(e)1_{\{\,f\in\mathcal{P}(\mathcal{D})\,\}}\,\Big|\,\mathbb{L}\Big)\,&=\,1_{e\notin\mathbb{L}}P\Big(e\in\mathcal{P}^+\big(f\in\mathcal{P}(\mathcal{D})\big)\,\Big|\,\mathbb{L}\Big)\\&=\,\frac{1}{p}1_{e\notin\mathbb{L}}P\Big(e\in\mathcal{P}^+\big(f\in\mathcal{P}(\mathcal{D})\big),f\in\mathcal{P}(\mathcal{D})\,\Big|\,\mathbb{L}\Big)\,.\end{aligned} \tag{45.50}$$

To get this equality, we have also used the fact that the event $\{ e \in \mathcal{P}(\mathcal{D}) \}$ is an increasing event of the configuration of the bonds in $\mathbb{E}^d(\Lambda(n)) \setminus \mathbb{L} = \mathcal{I}$. So we see that the conditional expectation is non–negative when $e$ is in $\mathcal{I}$. This allows us to put the sums inside the expectations and we use (45.49) and (45.50) to obtain

$$\left( \sum_{e \in \mathcal{I}} \Big( E\big( X(e) 1_{\{ f \in \mathcal{P}(\mathcal{D}) \}} \big) \Big)^2 \right)^{1/2} \leq \frac{1}{p} E\Big( \big| \mathcal{P}^+\big( f \in \mathcal{P}(\mathcal{D}) \big) \setminus \mathbb{L} \big| ; f \in \mathcal{P}(\mathcal{D}) \Big) . \tag{45.51}$$

In view of (45.51), our tentative definition for the set $\mathcal{H}$ is

$$\mathcal{H} = \Big\{ f \in \mathcal{J} : E\Big( \big| \mathcal{P}^+\big( f \in \mathcal{P}(\mathcal{D}) \big) \setminus \mathbb{L} \big| \,\Big|\, f \in \mathcal{P}(\mathcal{D}) \Big) \geq pS \Big\} . \tag{45.52}$$

At this point, we would like to prove a result similar to Talagrand's inequality: there exists a universal constant $C > 0$ such that

$$E\big( \big| \mathcal{P}(\mathcal{D}) \cap \mathcal{H} \big| \big) \leq C \sqrt{ \big| \mathbb{E}^d(\Lambda(n)) \big| } \exp\left( -\frac{S^2}{C} \right) . \tag{45.53}$$

So far, we did not manage to prove an estimate like (45.53). To arrive at this inequality, on the one hand we have weakened Talagrand's inequality by considering a specific event and using inequality (45.48), and on the other hand we have strengthened it by using inequality (45.51). What we can do rigorously is to prove a conditional version of Talagrand's inequality, by conditioning from the start with respect to the set $\mathbb{L}$, and working with the set $\mathcal{H}$ defined by

$$\mathcal{H} = \left\{ f \in \mathbb{L} : \left( \sum_{e \in \mathcal{I}} \Big( E\Big( X(e) 1_{\{ f \in \mathcal{P}(\mathcal{D}) \}} \,\Big|\, \mathbb{L} \Big) \Big)^2 \right)^{1/2} \geq SP\big( f \in \mathcal{P}(\mathcal{D}) \,\big|\, \mathbb{L} \big) \right\} .$$

However, it seems very difficult to improve this result or to obtain an inequality involving the set $\mathcal{H}$ defined in (45.52). Another interesting possibility is to use variants for the sets $\mathcal{I}$ and $\mathcal{J}$, in which the set of the pivotal bonds is assigned to the region $\mathcal{I}$ instead of $\mathcal{J}$. This creates an extra difficulty, because the conditional expectation (45.49) becomes negative for the bonds $e$ belonging to $\mathcal{P}(\mathcal{D})$, since these bonds are necessarily closed. As a consequence, the inequality (45.51) does not hold any more. We can still group together the terms in the sum according to the sign of the conditional expectation, and instead of (45.51), we would obtain a weaker inequality with the additional term

$$\frac{1}{1-p} \big| \mathcal{P}(\mathcal{D}) \big| \, P\big( f \in \mathcal{P}(\mathcal{D}) \big) .$$

We explain next the expected benefits of an inequality like (45.53). For a bond $f$ in $\mathcal{P}(\mathcal{D}) \setminus \mathcal{H}$, we have

$$E\Big( \big| \mathcal{P}^+\big( f \in \mathcal{P}(\mathcal{D}) \big) \setminus \mathbb{L} \big| ; f \in \mathcal{P}(\mathcal{D}) \setminus \mathcal{H} \Big) \leq p\, S\, P\big( f \in \mathcal{P}(\mathcal{D}) \big) .$$

We would then sum this inequality and we would get on one hand

$$E\Big(\sum_{f\in\mathcal{P}(\mathcal{D})\setminus\mathcal{H}} \big|\mathcal{P}^+\big(f\in\mathcal{P}(\mathcal{D})\big)\setminus\mathbb{L}\big|\Big) \,\leq\, p\,S\,E\big(\big|\mathcal{P}(\mathcal{D})\setminus\mathcal{H}\big|\big)\,. \tag{45.54}$$

On the other hand, we have, using (45.53),

$$\begin{aligned} E\Big(\sum_{f\in\mathcal{P}(\mathcal{D})\cap\mathcal{H}} \big|\mathcal{P}^+\big(f\in\mathcal{P}(\mathcal{D})\big)\setminus\mathbb{L}\big|\Big) \,&\leq\, \big|\mathbb{E}^d(\Lambda(n))\big|\,E\big(\big|\mathcal{P}(\mathcal{D})\cap\mathcal{H}\big|\big) \\ &\leq\, C\Big(\big|\mathbb{E}^d(\Lambda(n))\big|\Big)^{3/2}\exp\Big(-\frac{S^2}{C}\Big)\,. \end{aligned} \tag{45.55}$$

Adding together inequalities (45.54) and (45.55), we obtain

$$\begin{aligned} E\Big(\sum_{f\in\mathcal{P}(\mathcal{D})} \big|\mathcal{P}^+\big(f\in\mathcal{P}(\mathcal{D})\big)\setminus\mathbb{L}\big|\Big) \,&\leq \\ p\,S\,E\big(\big|\mathcal{P}(\mathcal{D})\big|\big) &+ C\Big(\big|\mathbb{E}^d(\Lambda(n))\big|\Big)^{3/2}\exp\Big(-\frac{S^2}{C}\Big)\,. \end{aligned} \tag{45.56}$$

By taking $S=\kappa\sqrt{\ln n}$, with a sufficiently large constant $\kappa$, we would conclude from (45.56) and the classical bound (39.11) on the pivotal bonds that

$$E\Big(\Big|\bigcup_{f\in\mathcal{P}(\mathcal{D})} \mathcal{P}^+\big(f\in\mathcal{P}(\mathcal{D})\big)\setminus\mathbb{L}\Big|\Big) \,\leq\, 2\kappa\sqrt{\frac{\ln n}{p(1-p)}}\Big(\big|\mathbb{E}^d(\Lambda(n))\big|\Big)^{1/2}\,.$$

By exchanging the roles of the left and the right faces of $\Lambda(n)$, we would get a similar inequality for the positive pivotal bonds for the set $\mathcal{P}(\mathcal{D})$ in the complement of the union of the clusters touching the right face. Adding these two inequalities, we would have

$$E\big(\big|\mathcal{P}^+\big(\mathcal{P}(\mathcal{D})\big)\big|\big) \,\leq\, 4\kappa\sqrt{\frac{\ln n}{p(1-p)}}\Big(\big|\mathbb{E}^d(\Lambda(n))\big|\Big)^{1/2}\,, \tag{45.57}$$

where the set $\mathcal{P}^+\big(\mathcal{P}(\mathcal{D})\big)$ is defined as

$$\mathcal{P}^+\big(\mathcal{P}(\mathcal{D})\big) \,=\, \bigcup_{f\in\mathcal{P}(\mathcal{D})} \mathcal{P}^+\big(f\in\mathcal{P}(\mathcal{D})\big)\,. \tag{45.58}$$

The set $\mathcal{P}^+\big(\mathcal{P}(\mathcal{D})\big)$ is the set of the positive pivotal bonds for the set $\mathcal{P}(\mathcal{D})$, that is the set of the bonds which are open, and the closure of which would modify the set $\mathcal{P}(\mathcal{D})$. This set is a subset of $\mathcal{P}_2(\mathcal{D})$, we call it the second positive pivots. As we will explain in the next subsection, a correct control on the second positive pivots might be enough to obtain an adequate control on the full set $\mathcal{P}_2(\mathcal{D})$ of the second pivots. Unfortunately, even the inequality (45.57) doesn't meet our needs. Indeed, the logarithm factor is quite bad, because, upon iteration, after $\ln n$ steps, we would have a super-polynomial factor. Anyway, the not so useful inequality (45.57) is based on the unproven inequality (45.53). In conclusion, this subsection is yet another piece of science-fiction!

## 45.5 The second pivots, again and again

In this subsection, we consider again the third question of section 30, namely the disconnection in a cuboid, and we focus on the second pivotal bonds. As always, our goal is to obtain a quantitative control either on $\mathcal{I}(2)$ or on $\mathcal{P}_2(\mathcal{D})$. Instead of trying to control directly the probability $P\big(|\mathcal{I}(2)|>i_2,\, W\not\longleftrightarrow W'\big)$, as we did in subsection 31.2, we try to implement a new approach involving the expectation $E\big(|\mathcal{I}(2)|;\, W\not\longleftrightarrow W'\big)$. Here and in the sequel, we use the following notation: for $A$ an event and $X$ a random variable, we define

$$E(X;A)\,=\,E(X1_A)\,.$$

The new approach relies on the fact that the bonds of $\mathcal{I}(2)$ or $\mathcal{P}_2(\mathcal{D})$ are in fact pivotal for the random variable $|\mathcal{I}(1)|$. In our first attempt conducted in subsection 31.2, we used the fact the bonds of $\mathcal{I}(2)$ are pivotal for the connection between $\partial^{\,out}\mathcal{A}(0)\setminus\mathcal{A}(1)$ and $W'$. Yet we had to condition on the set $\mathcal{I}(1)$ and since we could not control the ensuing conditional probability in a satisfactory way, we had to sum brutally over the possible values of $\mathcal{I}(1)$. This created a catastrophic super-exponential factor $n^{dn^d/2}$, which was the source of the failure of the method. The advantage of working with the variable $|\mathcal{I}(1)|$ is that the number of its possible values is $O(n^{d/2})$, which is polynomial in $n$, instead of super-exponential. Naturally, what we gain here has to be paid elsewhere, because the conditioning event $\{\,|\mathcal{I}(1)|=i_1\,\}$ is not monotone, and there is no straightforward bound on the number of its pivotal bonds. Let us start the precise computation. We decompose and we bound the expectation of $|\mathcal{I}(2)|$ as follows. We fix an adequate large constant $\kappa>0$ and we write

$$E\big(|\mathcal{I}(2)|;\, W\not\longleftrightarrow W'\big)\,\leq\sum_{0\leq i_1\leq\kappa\sqrt{n^d\ln n}}E\big(|\mathcal{I}(2)|;\, W\not\longleftrightarrow W',\,|\mathcal{I}(1)|=i_1\big)$$
$$+\,P\Big(W\not\longleftrightarrow W',\,|\mathcal{I}(1)|>\kappa\sqrt{n^d\ln n}\Big)\,.$$

We already know how to control the last term above (see inequality (31.9)). The essential quantity we are now aiming for is $E\big(|\mathcal{I}(2)|;\, W\not\longleftrightarrow W',\,|\mathcal{I}(1)|=i_1\big)$. The bonds of $\mathcal{I}(2)$ are negative pivotal bonds for the event $\{\,|\mathcal{I}(1)|=i_1\,\}$, and they are included in $\mathbb{E}^d(D\setminus\mathcal{A}(0))$. Notice that the bonds of $\mathcal{I}(1)$ are also negative pivotal bonds for the event $\{\,|\mathcal{I}(1)|=i_1\,\}$, but they have one endpoint in $\mathcal{A}(0)$, so they do not belong to $\mathbb{E}^d(D\setminus\mathcal{A}(0))$. The good point is that the distribution of the bonds of $\mathbb{E}^d(D\setminus\mathcal{A}(0))$ is unaffected by the value of $\mathcal{A}(0)$, so our strategy consists in conditioning on $\mathcal{A}(0)$. We write

$$E\big(|\mathcal{I}(2)|;\, W\not\longleftrightarrow W',\,|\mathcal{I}(1)|=i_1\big)\,=\,\sum_{A_0}E\big(|\mathcal{I}(2)|;\,\mathcal{A}(0)=A_0,\,|\mathcal{I}(1)|=i_1\big)\,,\tag{45.59}$$

where the sum runs over the connected sets $A_0$ satisfying $W\subset A_0\subset D$ and $W'\cap A_0=\varnothing$. For these sets $A_0$, the event $\{\,\mathcal{A}(0)=A_0\,\}$ is automatically included in the event $\{\,W\not\longleftrightarrow W'\,\}$. Let $A_0$ be such a set and let $\omega$ be a configuration belonging to the event $\{\,\mathcal{A}(0)=A_0,\,|\mathcal{I}(1)|=i_1\,\}$. Let $\omega_0$ be the percolation configuration $\omega$ restricted to the edges $\mathbb{E}^d(D\setminus\mathcal{A}(0))$. We have then

$$\mathcal{I}(2)\;=\;\mathcal{P}^-\big(\{\,|\mathcal{I}(1)|=i_1\,\},\omega_0\big)\,. \tag{45.60}$$

Unfortunately, the classical bound (39.11) on the pivotal bonds cannot be used, because the event $\{\,|\mathcal{I}(1)|=i_1\,\}$ is not monotone. So we use instead the inequality (39.10), applied to the configuration restricted to the edges $\mathbb{E}^d(D\setminus\mathcal{A}(0))$. We introduce the set of the right positive pivots for $\{\,|\mathcal{I}(1)|=i_1\,\}$, defined as

$$\mathcal{P}^+_R\big(\{\,|\mathcal{I}(1)|=i_1\,\}\big)\;=\;\mathcal{P}^+\big(\{\,|\mathcal{I}(1)|=i_1\,\},\omega_0\big)\,.$$

We have then

$$\Big|E\Big(\big|\mathcal{P}^+_R\big(\{\,|\mathcal{I}(1)|=i_1\}\big)\big|-\big|\mathcal{P}^-_R\big(\{\,|\mathcal{I}(1)|=i_1\}\big)\big|\,\Big|\,\mathcal{A}(0)=A_0\Big)\Big|\leq\sqrt{\frac{\big|\mathbb{E}^d\big(\Pi(n)\big)\big|}{p(1-p)}}\,.$$

We rewrite the left-hand side with the help of the identity (40.1), we multiply by $P(\mathcal{A}(0)=A_0)$ and we obtain

$$\begin{aligned}\Big|E\Big(\frac{1}{p}\big|\mathcal{P}^+_R\big(\{\,|\mathcal{I}(1)|=i_1\}\big)\big|-\frac{1}{1-p}\big|\mathcal{P}^-_R\big(\{\,|\mathcal{I}(1)|=i_1\}\big)\big|;|\mathcal{I}(1)|=i_1,\mathcal{A}(0)=A_0\Big)\Big|&\\ \leq\sqrt{\frac{\big|\mathbb{E}^d\big(\Pi(n)\big)\big|}{p(1-p)}}\,P(\mathcal{A}(0)=A_0)\,.&\end{aligned} \tag{45.61}$$

We use next the identity (45.60) to replace $\mathcal{P}^-_R\big(\{\,|\mathcal{I}(1)|=i_1\,\}\big)$ by $\mathcal{I}(2)$ in the second term of the inequality (45.61). We sum the resulting inequality over the possible values of the set $A_0$, and we use the identity (45.59) to get finally

$$\begin{aligned}\Big|E\Big(\frac{1}{p}\big|\mathcal{P}^+_R\big(\{\,|\mathcal{I}(1)|=i_1\,\}\big)\big|-\frac{1}{1-p}\big|\mathcal{I}(2)\big|;W\not\longleftrightarrow W',\,|\mathcal{I}(1)|=i_1\Big)\Big|&\\ \leq\sqrt{\frac{\big|\mathbb{E}^d\big(\Pi(n)\big)\big|}{p(1-p)}}\,P(W\not\longleftrightarrow W')\,.&\end{aligned} \tag{45.62}$$

Naturally, we can define the symmetric set $\mathcal{P}^+_L\big(\{\,|\mathcal{I}(1)|=i_1\,\}\big)$ of the left positive pivotal bonds for the variable $|\mathcal{I}(1)|$, whose cardinality satisfies the same inequality (45.62). As a matter of fact, the sets of the left and right positive pivots for $|\mathcal{I}(1)|$ can be defined directly, without a conditioning, as follows.

**Definition 45.9.** *A bond $e=\langle a,b\rangle$ belongs to the set of the right second positive pivots, which we denote by $\mathcal{P}^+_R(\mathcal{P}(\mathcal{D}))$, if and only if there exists a bond $f=\langle c,d\rangle$ belonging to $\mathcal{I}(1)$ such that the following events occur:*

$$W\longleftrightarrow c\,,\quad d\longleftrightarrow a\,,\quad d\not\longleftrightarrow W'\quad \textit{in}\quad \mathbb{E}^d\big(\Pi(n)\big)\setminus\{\,e\,\}\,,\quad b\longleftrightarrow W'\,.$$

*Symmetrically, a bond $e=\langle a,b\rangle$ belongs to the set of the left second positive pivots, which we denote by $\mathcal{P}^+_L(\mathcal{P}(\mathcal{D}))$, if the same event as above occurs, but with $W$ and $W'$ exchanged.*

In other words, the bonds of $\mathcal{P}^+_R(\mathcal{P}(\mathcal{D}))$ are the open bonds of $\mathbb{E}^d\big(\Pi(n)\big)$ which are pivotal for the connection between $W'$ and a vertex in the outer vertex boundary of $\mathcal{A}(0)$. It can be checked that

$$\begin{aligned}\mathcal{P}^+_L(\mathcal{P}(\mathcal{D})) \,&=\, \bigcup_{i_1\geq 0} \mathcal{P}^+_L\big(\{\,|\mathcal{I}(1)|=i_1\,\}\big)\,,\\ \mathcal{P}^+_R(\mathcal{P}(\mathcal{D})) \,&=\, \bigcup_{i_1\geq 0} \mathcal{P}^+_R\big(\{\,|\mathcal{I}(1)|=i_1\,\}\big)\,.\end{aligned}$$

Moreover the union of the two sets $\mathcal{P}^+_L(\mathcal{P}(\mathcal{D}))$ and $\mathcal{P}^+_R(\mathcal{P}(\mathcal{D}))$ is precisely the set $\mathcal{P}^+(\mathcal{P}(\mathcal{D}))$ defined in (45.58), which popped up when we tried to adapt Talagrand's lemma to the disconnection event. Some examples of these sets in simulations are presented in figures 79 and 80. Inequality (45.62) shows that the problem of controlling the expected cardinality of $\mathcal{I}(2)$ can be reduced to the problem of controlling the expected cardinality of the right positive pivotal bonds for the variable $|\mathcal{I}(1)|$. Of course, the left positive pivotal bonds for the variable $|\mathcal{I}(1)|$ satisfy the same inequality (45.62). However, there is no obvious way of getting a control on the cardinality of the sets $\mathcal{P}^-(\mathcal{P}(\mathcal{D}))$ or $\mathcal{P}^+(\mathcal{P}(\mathcal{D}))$.

Let us present an attempt of reducing the problem to monotone events. We write the event $\{\,|\mathcal{I}(1)|=i_1\,\}$ as the intersection of two monotone events:

$$\{\,|\mathcal{I}(1)|=i_1\,\} \,=\, \{\,|\mathcal{I}(1)|\leq i_1\,\}\cap\{\,|\mathcal{I}(1)|\geq i_1\,\}\,. \tag{45.63}$$

We will try to take advantage of the following little lemma.

**Lemma 45.10.** *Let $A$ be an increasing event and let $B$ be a decreasing event. For any configuration $\omega$, we have*

$$\mathcal{P}^-(A\cap B,\omega)\,\subset\,\mathcal{P}^-(B,\omega)\,,\qquad \mathcal{P}^+(A\cap B,\omega)\,\subset\,\mathcal{P}^+(A,\omega)\,.$$

*Proof.* This result holds for a general product space, but for coherence we do it here for the bond percolation model. Let us fix a configuration $\omega$. Let $e$ be a bond belonging to $\mathcal{P}^-(A\cap B,\omega)$. By definition, we have

$$\omega_e\in A\cap B\,,\qquad \omega^e\not\in A\cap B\,.$$

Now the event $A$ is increasing, thus

$$\omega_e\in A\quad\Longrightarrow\quad\omega^e\in A\,.$$

Therefore

$$\omega^e\not\in A\cap B\quad\Longrightarrow\quad\omega^e\not\in B\,,$$

and we conclude that $e$ belongs to $\mathcal{P}^-(B,\omega)$. Similarly, suppose that the bond $e$ belongs to $\mathcal{P}^+(A\cap B,\omega)$. By definition, we have

$$\omega_e\not\in A\cap B\,,\qquad \omega^e\in A\cap B\,.$$

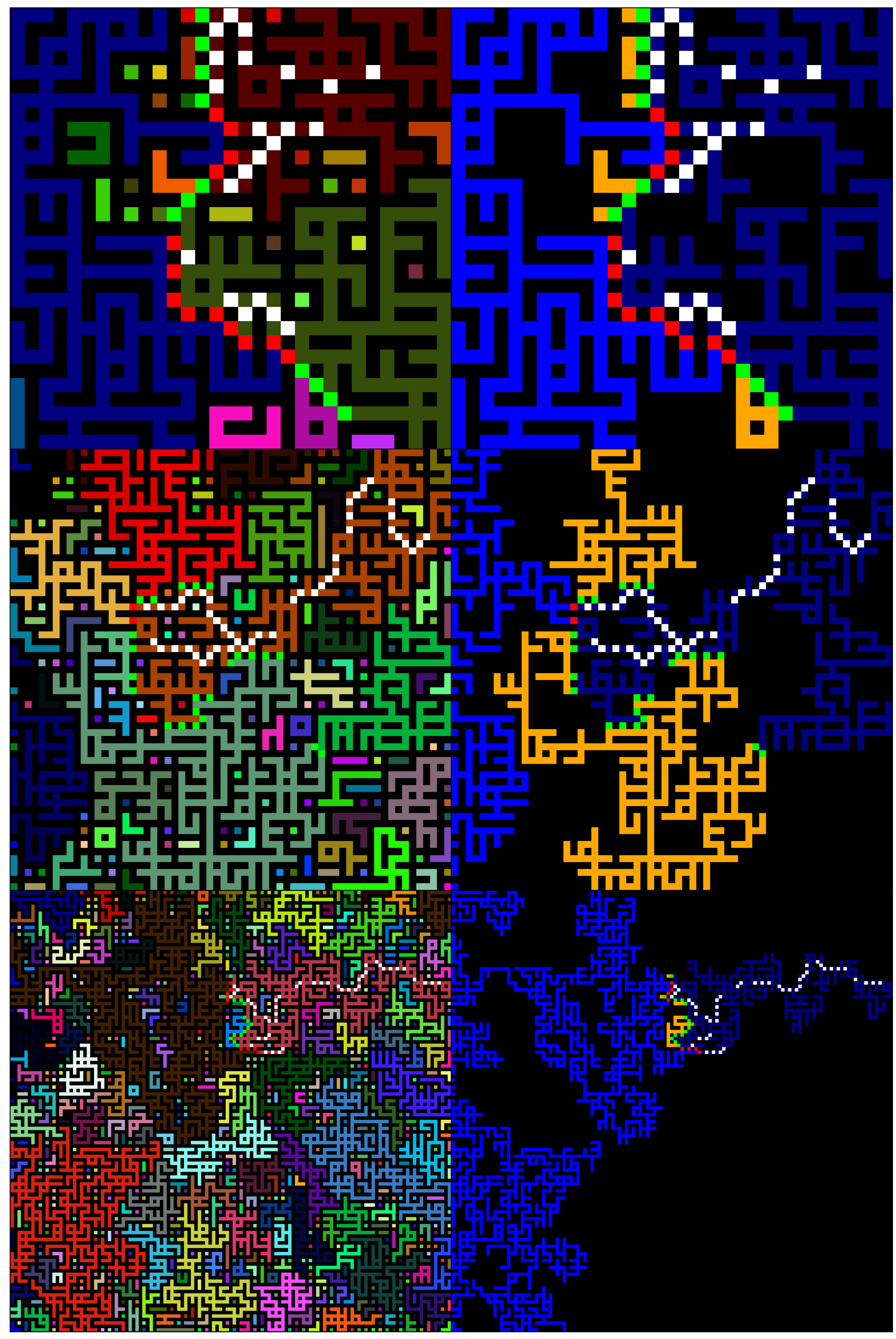

Figure 79: Disconnection between left and right in a square: pivotal bonds, negative right second pivotal bonds, positive right second pivotal bonds. Left column: the colored clusters, the pivots and the second right pivots. Right column: $\mathcal{A}(0)$, clusters meeting $\Delta\mathcal{A}(0)$ and negative right second pivots, clusters meeting the right face and pivots (first or second), and the pivots. Boxes $\Lambda(16), \Lambda(32), \Lambda(64)$, $p = 0.443$, $0.361$, $0.340$ (top row to bottom row).

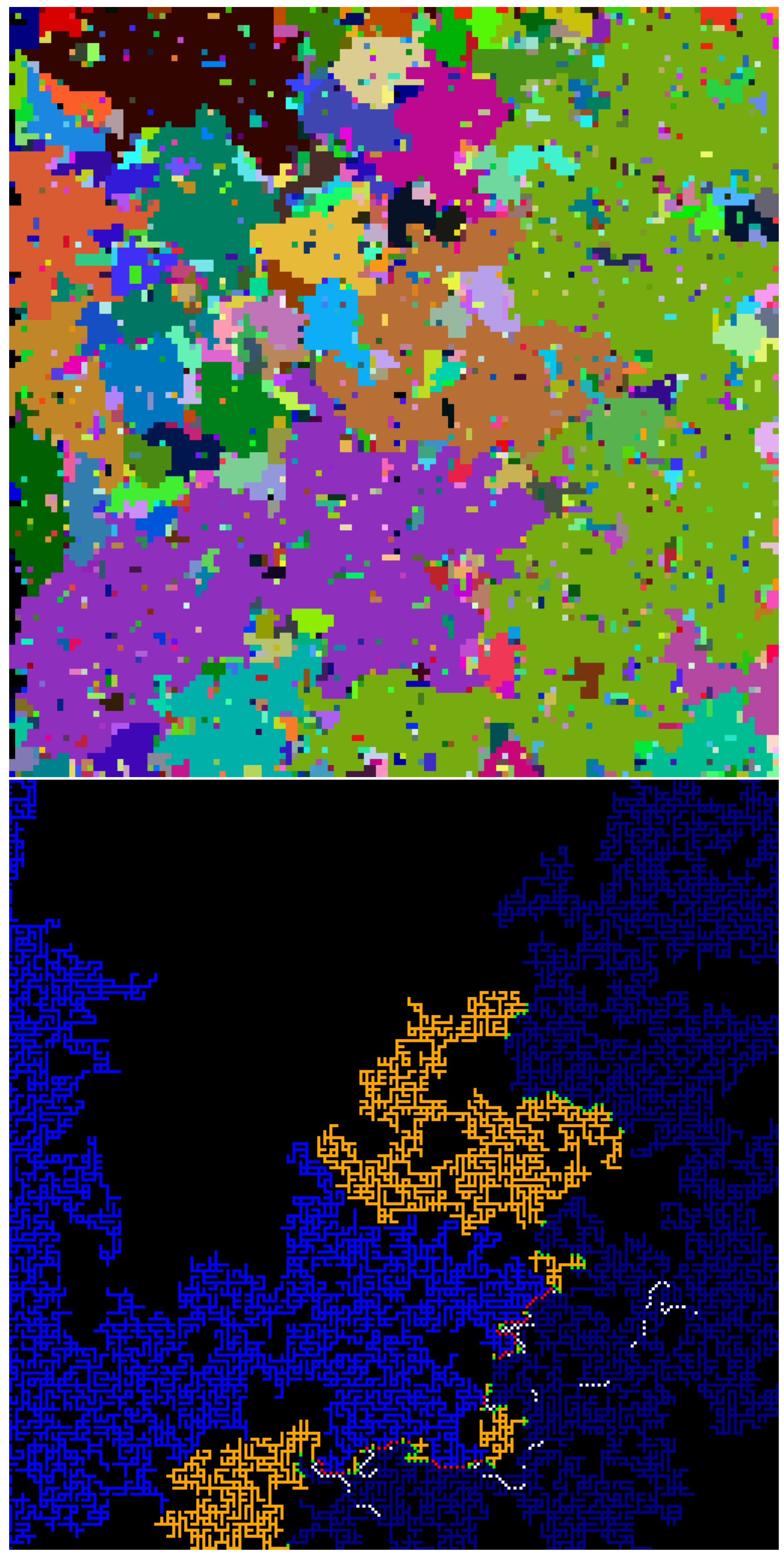

Figure 80: Continuation of figure 79. Bond percolation in $\Lambda(128)$, $p = 0.391$.

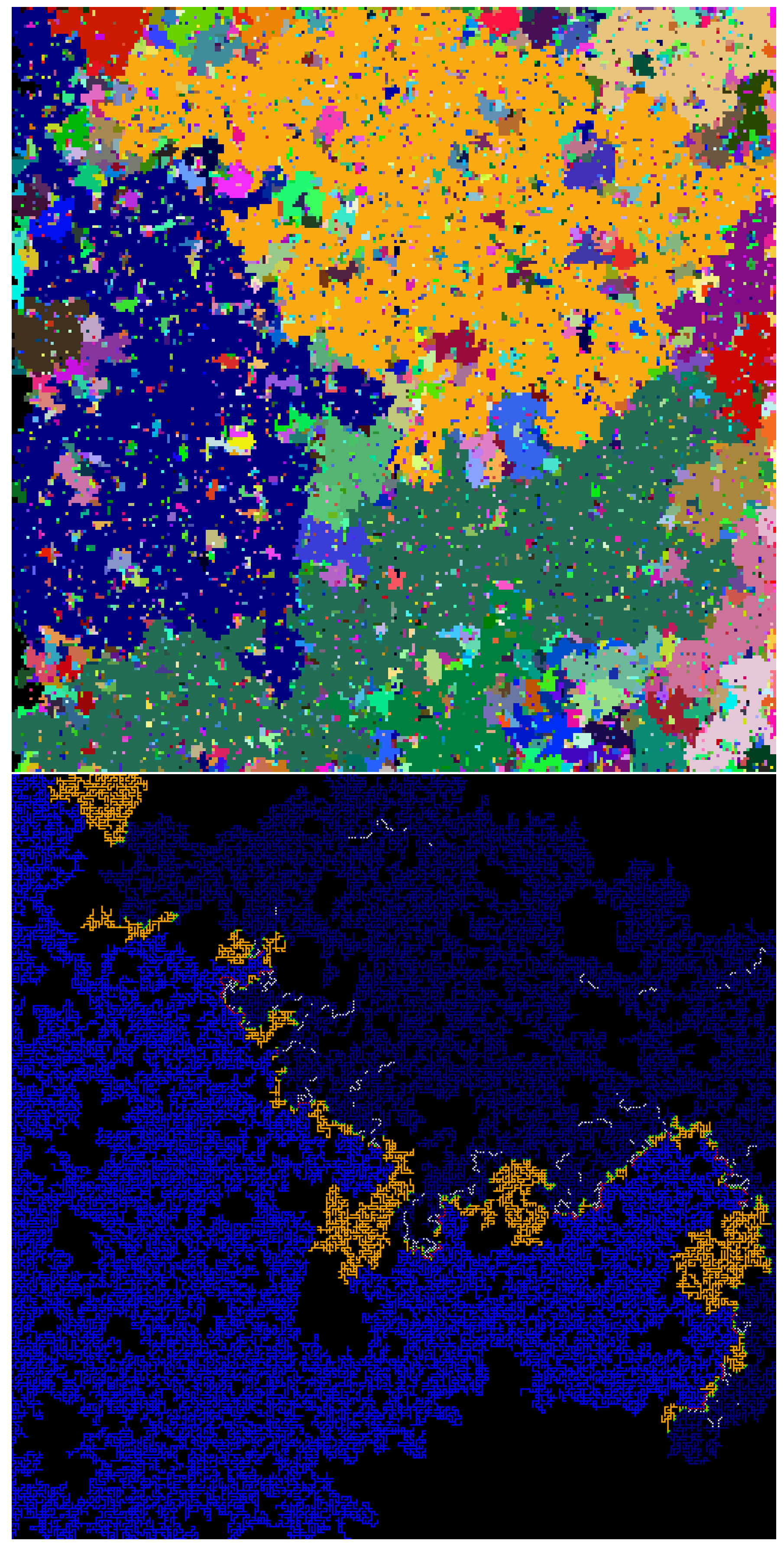

Figure 81: Continuation of figure 80. Bond percolation in $\Lambda(256)$, $p = 0.396$.

Now the event $B$ is decreasing, thus

$$\omega^e \in B \quad \Longrightarrow \quad \omega_e \in B\,.$$

Therefore

$$\omega_e \notin A \cap B \quad \Longrightarrow \quad \omega_e \notin A\,,$$

and we conclude that $e$ belongs to $\mathcal{P}^+(A, \omega)$. $\square$

We restart the computations from formula (45.59). Let us fix a set $A_0$ among those arising in the sum of (45.59). On the event $\{\, \mathcal{A}(0) = A_0 \,\}$, the values of $|\mathcal{I}(1)|$ and $|\mathcal{I}(2)|$ are completely determined by the configuration $\omega_0$ restricted to $\mathbb{E}^d(D \setminus \mathcal{A}(0))$, whence

$$\begin{aligned} E\big(|\mathcal{I}(2)|;\, \mathcal{A}(0) = A_0,\, |\mathcal{I}(1)| = i_1\big) \;&=\; E\big(|\mathcal{I}(2)| 1_{\{\, |\mathcal{I}(1)| = i_1 \,\}};\, \mathcal{A}(0) = A_0\big) \\ &=\; E_{A_0}\big(|\mathcal{I}(2)| 1_{\{\, |\mathcal{I}(1)| = i_1 \,\}}\big) P\big(\mathcal{A}(0) = A_0\big)\,, \qquad (45.64) \end{aligned}$$

where $E_{A_0}$ means that we have conditioned on $\{\, \mathcal{A}(0) = A_0 \,\}$ and we take the expectation with respect to the edges in $\mathbb{E}^d(D \setminus A_0)$, that is we average over the configuration that was previously denoted by $\omega_0$. Using the identity (45.60), we write next

$$E_{A_0}\big(|\mathcal{I}(2)| 1_{\{\, |\mathcal{I}(1)| = i_1 \,\}}\big) \;=\; E_{A_0}\Big(\big|\mathcal{P}^-\big(\{\, |\mathcal{I}(1)| = i_1 \,\}, \omega_0\big)\big| 1_{\{\, |\mathcal{I}(1)| = i_1 \,\}}(\omega_0)\Big)\,. \qquad (45.65)$$

Using the representation (45.63) and lemma 45.10, we obtain furthermore

$$\begin{aligned} &E_{A_0}\Big(\big|\mathcal{P}^-\big(\{\, |\mathcal{I}(1)| = i_1 \,\}, \omega_0\big)\big| 1_{\{\, |\mathcal{I}(1)| = i_1 \,\}}(\omega_0)\Big) \;\leq \\ &\qquad\qquad E_{A_0}\Big(\big|\mathcal{P}^-\big(\{\, |\mathcal{I}(1)| \leq i_1 \,\}, \omega_0\big)\big| 1_{\{\, |\mathcal{I}(1)| = i_1 \,\}}(\omega_0)\Big)\,. \qquad (45.66) \end{aligned}$$

The point is that the event $\{\, |\mathcal{I}(1)| \leq i_1 \,\}$ is now decreasing with respect to $\omega_0$, hence we can use the classical bound (39.11) to obtain

$$E_{A_0}\Big(\big|\mathcal{P}^-\big(\{\, |\mathcal{I}(1)| \leq i_1 \,\}, \omega_0\big)\big| 1_{\{\, |\mathcal{I}(1)| = i_1 \,\}}(\omega_0)\Big) \;\leq\; \sqrt{\frac{\big|\mathbb{E}^d(\Lambda(n))\big|}{p(1-p)}}\,. \qquad (45.67)$$

Substituting (45.67) successively in (45.66), (45.65), (45.64), and (45.59), we get

$$E\big(|\mathcal{I}(2)|;\, W \not\longleftrightarrow W',\, |\mathcal{I}(1)| = i_1\big) \;\leq\; \sqrt{\frac{\big|\mathbb{E}^d(\Lambda(n))\big|}{p(1-p)}}\,. \qquad (45.68)$$

The inequality (45.68) is rather nice and goes in the right direction, but it is still not sufficient to reach a decisive conclusion. The remaining problem is the conditioning with respect to the event $\{\, |\mathcal{I}(1)| = i_1 \,\}$. To control the expectation of $|\mathcal{I}(2)|$, we would need to sum (45.68) over the possible values of $i_1$, and this would create a factor of order $O(n^{d/2})$, which would once more ruin the ensuing estimate. Indeed, we would end up with a useless estimate like the one predicted in (45.45). Anyway, one benefit of conditioning with respect to the variable $|\mathcal{I}(1)|$ is to highlight the role played by the set $\mathcal{P}^+(\mathcal{P}(\mathcal{D}))$ of the second positive pivots, which coincides with the positive pivots for the random variable $|\mathcal{I}(1)|$.

## 45.6 Extension to higher pivots

We would like to derive further inequalities like (45.68) involving higher pivots, or rather intersection sets. The idea is the following. We will condition on the cardinalities of the intersection sets. So we fix $k \geq 1$, $k$ integers $i_1, \dots, i_k$ and we try to estimate the probability of the event $\{ |\mathcal{I}(1)| = i_1, \dots, |\mathcal{I}(k)| = i_k \}$. Naturally, we have at our disposal an estimate like

$$\forall k \geq 1 \qquad E\Big( \big|\mathcal{P}^+\big(|\mathcal{I}(k)| = i_k\big)\big| - \big|\mathcal{P}^-\big(|\mathcal{I}(k)| = i_k\big)\big| \Big) \,=\, O(n^{d/2})\,. \quad (45.69)$$

Here the constant associated to the $O(\cdot)$ should not depend on $k$, only on $d$ and $p$. For the third question of section 30, we can prove that, when the disconnection event occurs, some intersection set $\mathcal{I}(k)$ has cardinality larger than $n^{d-1}$ with overwhelming probability, for an index $k$ which is negligible compared to $\ln n$. So, what we are missing is another control like (45.69) which would link the expected cardinalities of $\mathcal{P}^+\big(|\mathcal{I}(k)| = i_k\big)$ and $\mathcal{P}^-\big(|\mathcal{I}(k+1)| = i_{k+1}\big)$ for $k \geq 1$.

Our great hope is to obtain additional inequalities similar to (45.62), using either successive conditioning or the formula for the pivots in a random domain. A natural attempt consists in designing a more sophisticated exploration algorithm to find the intersection sets. Ideally, at step $k$, when building $\mathcal{I}(k+1)$, we would like to avoid visiting the bonds belonging to $\mathcal{P}^+\big(|\mathcal{I}(k)| = i_k\big)$. However, this seems to be too demanding. During the last decade, most of our efforts were spent trying to make progress along these directions, but in vain. Maybe it was a deadly trap. Nevertheless, we describe next quickly a valuable attempt.

We remark first that the set $\mathcal{I}(k)$ is included in the set $\mathcal{P}^-\big(|\mathcal{I}(k)| = i_k\big)$, but also in the set $\mathcal{P}_k^-\big(|\mathcal{I}(1)| = i_1\big)$ (this set was defined in section 44). Yet the set $\mathcal{P}_k^-\big(|\mathcal{I}(1)| = i_1\big)$ contains many more bonds. Among them, we single out the bonds which can be detected with the help of adequate intertwined explorations, as we define next. Let us consider a partition $\pi$ of $k$

$$\pi = (k_1, \dots, k_\ell)\,, \quad k \,=\, k_1 + k_2 + \cdots + k_\ell\,,$$

into $\ell$ positive integers. To this partition $\pi$, we associate a sequence of random sets $(\mathcal{A}^\pi(t), \mathcal{W}^\pi(t), t \geq 0)$ which is built as in subsection 29.4 with the difference that, at each step, instead of applying the map Coherent_Bond_Explore only once, we iterate it a certain number of times, as specified by the partition, while keeping the same second argument. We make the convention that $\mathcal{A}(-2) = \mathcal{A}(-1) = \mathcal{W}(-2) = \mathcal{W}(-1) = \varnothing$. At step $t \geq 0$, the increment between $\mathcal{A}^\pi(t-2)$ and $\mathcal{A}^\pi(t)$ will consist of the vertices which are at travel distance $k_t - 1$ from $\mathcal{W}(t-2)$ in the domain $D \setminus \big(\mathcal{A}(t-2) \cup \mathcal{A}(t-1)\big)$. The associated waiting set $\mathcal{W}^\pi(t)$ consists of the vertices from which the round of exploration would restart after these $k_t$ iterations. So we obtain a sequence of intersection sets $(\mathcal{I}^\pi(t), t \geq 0)$ associated to this mechanism and the partition $\pi$ of $k$. Finally, we consider the union of these sets over all the partitions $\pi$ of $k$. We refrain from providing here a heavy formal definition, because for the time being, we are still missing an essential estimate to make this construction successful. The counterpart of the curves in figures 60 to 71 is presented in figures 86 to 96. Various simulations of these intersection sets are presented in figures 82, 83 and 84.

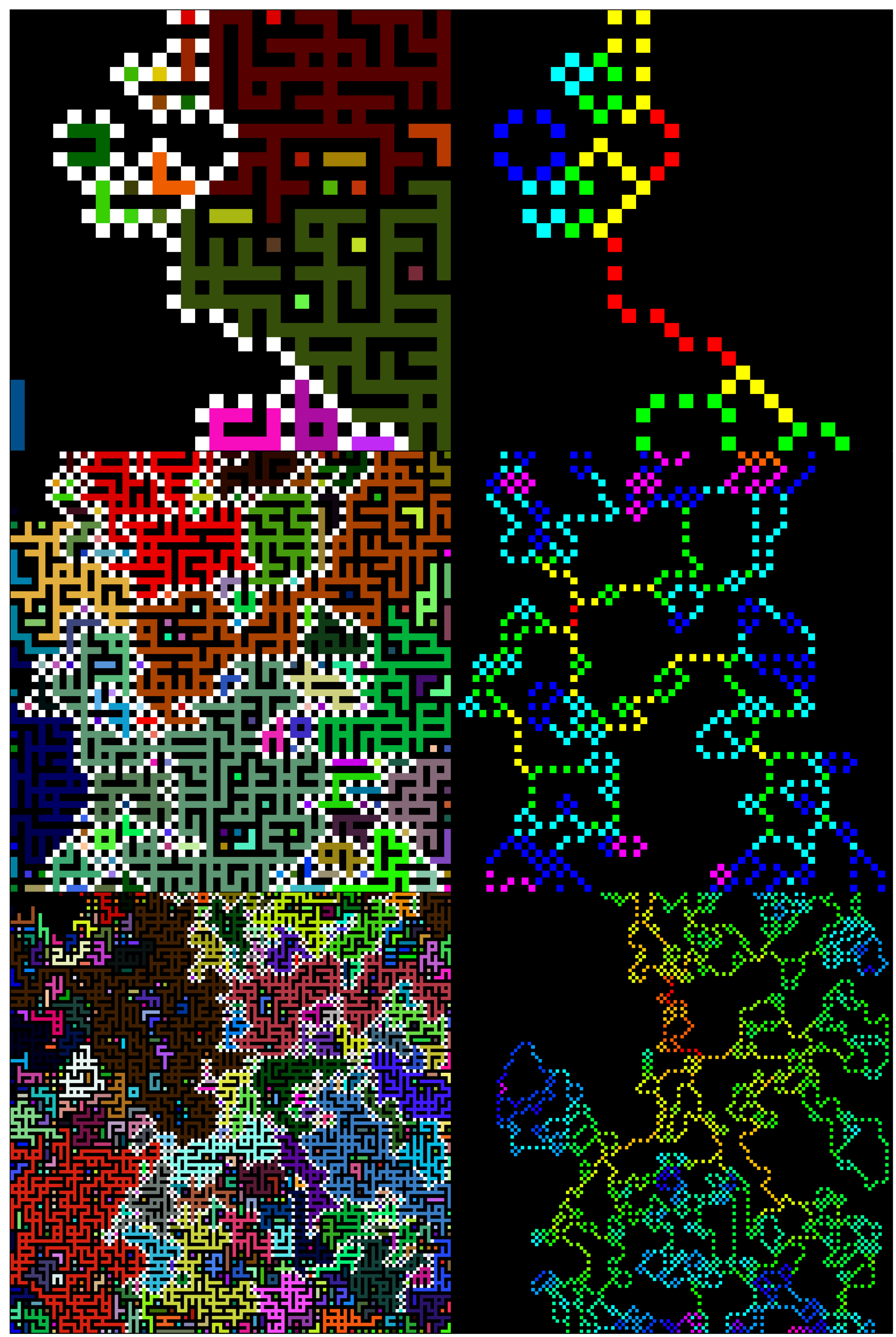

Figure 82: Higher pivots, disconnection between left and right faces of a square:
Left column: the colored clusters and the negative pivots of all orders in white.
Right column: the pivots of all orders with rainbow colors, first, second, third,...
Top row to bottom row: boxes $\Lambda(16), \Lambda(32), \Lambda(64)$, $p = 0.443$, $0.361$, $0.340$.

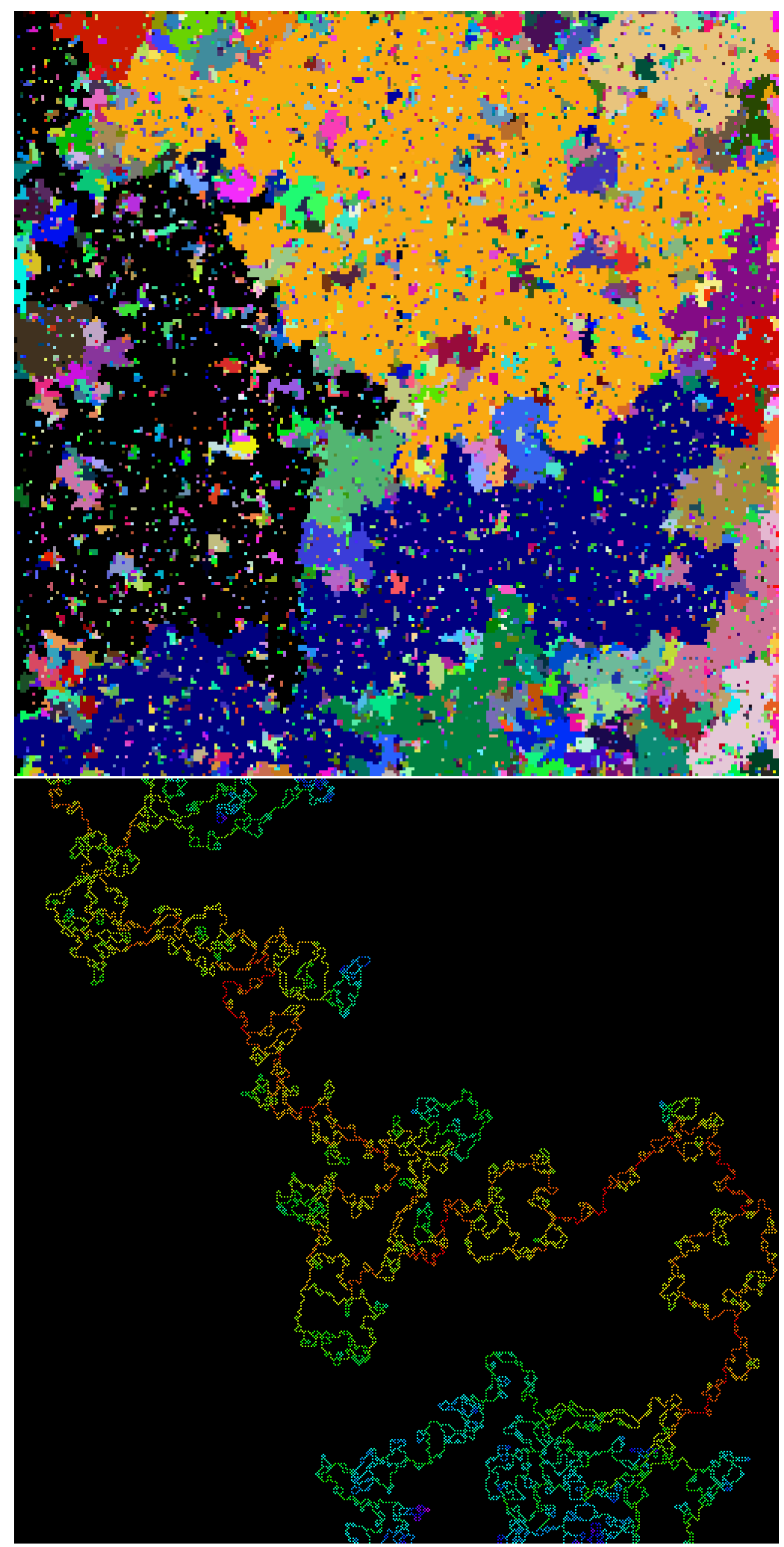

Figure 83: Continuation of figure 82. Bond percolation in $\Lambda(256)$, $p = 0.396$. With a larger box, the pivotal bonds are no longer discernible on the picture.

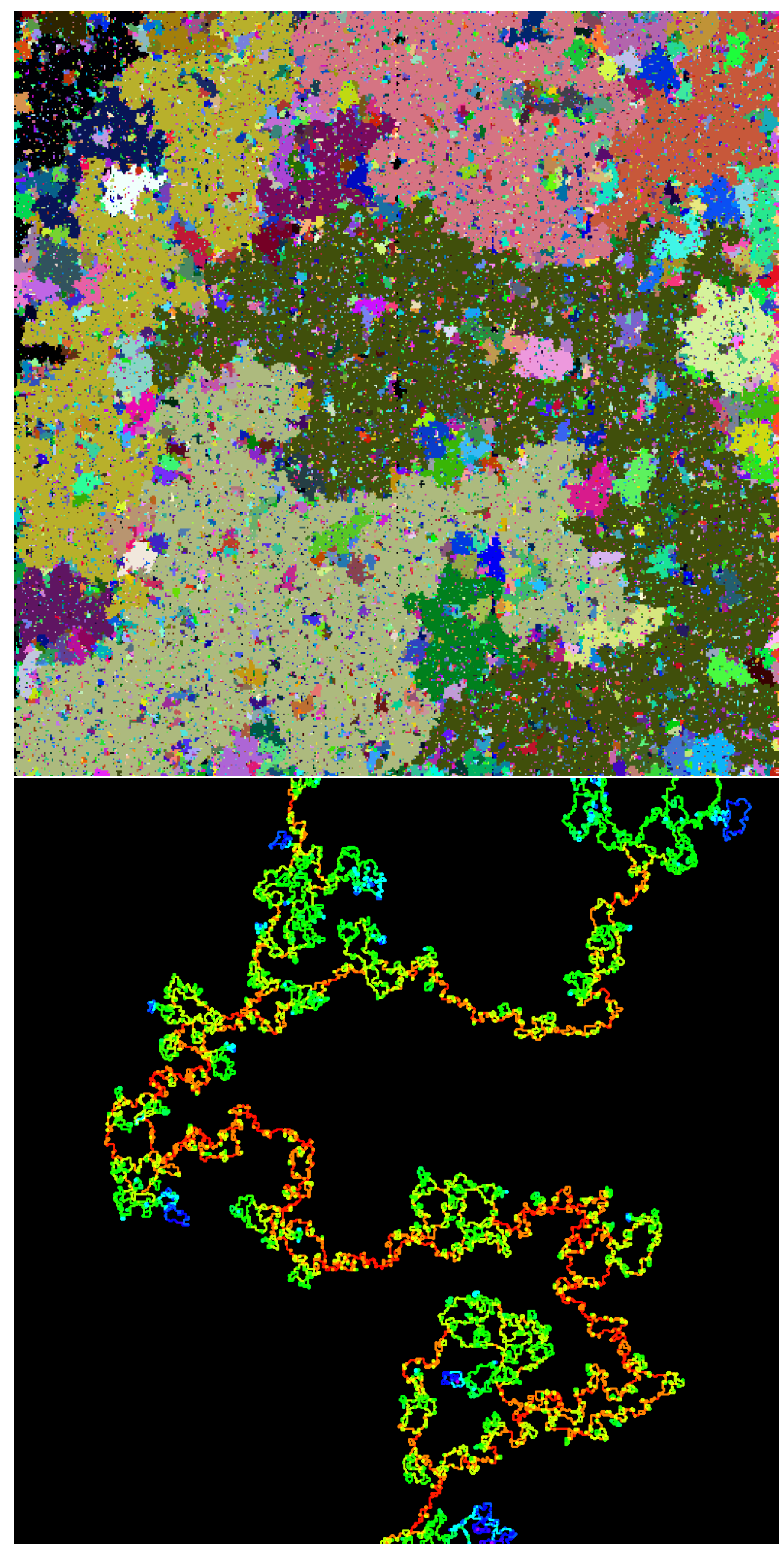

Figure 84: Continuation of figure 83. Bond percolation in $\Lambda(512)$, $p = 0.400$. The morphology operator "dilation by a square 1" was applied to the pivots.

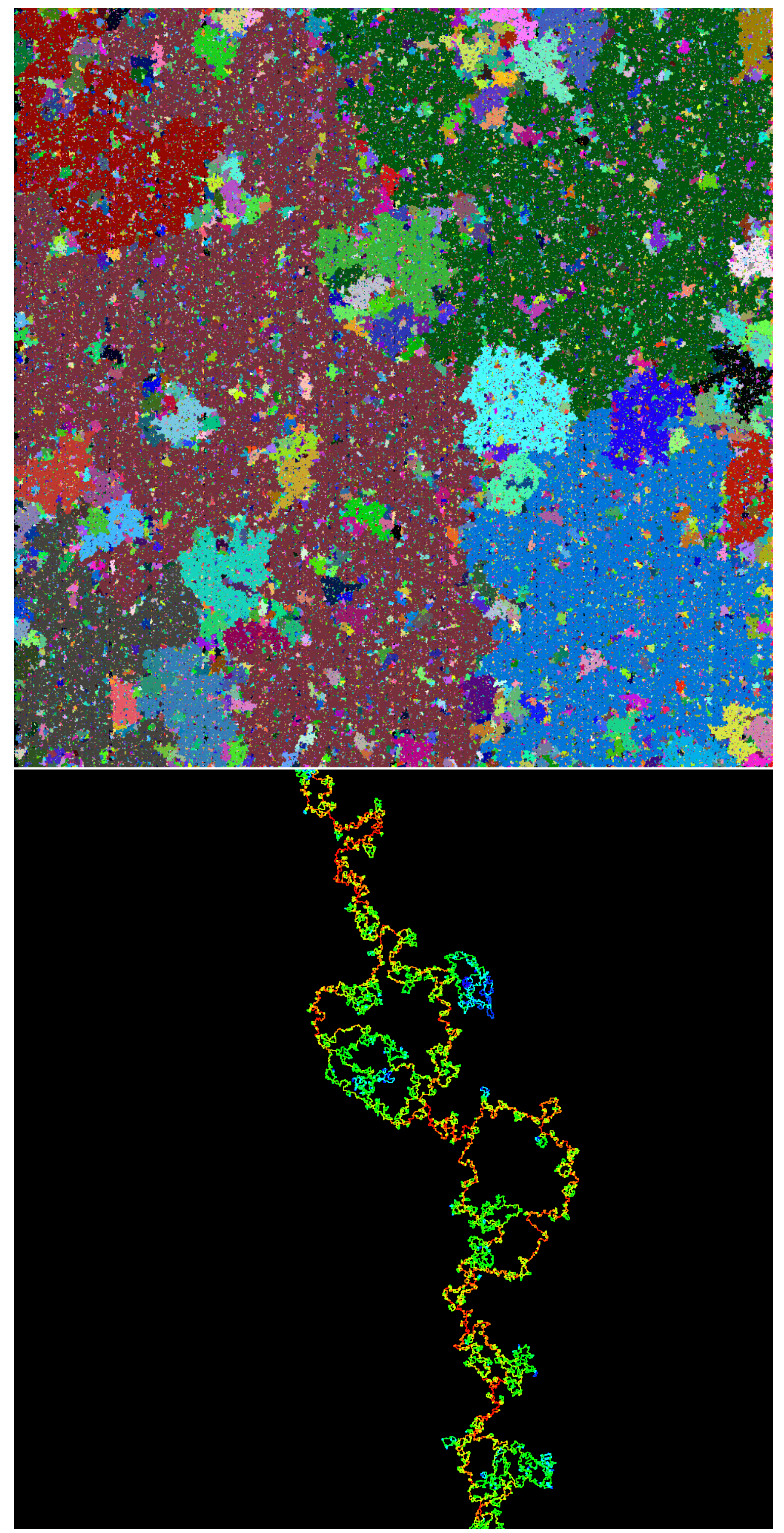

Figure 85: Continuation of figure 84. Bond percolation in $\Lambda(1024)$, $p = 0.401$. The morphology operator “dilation by a square 2” was applied to the pivots.

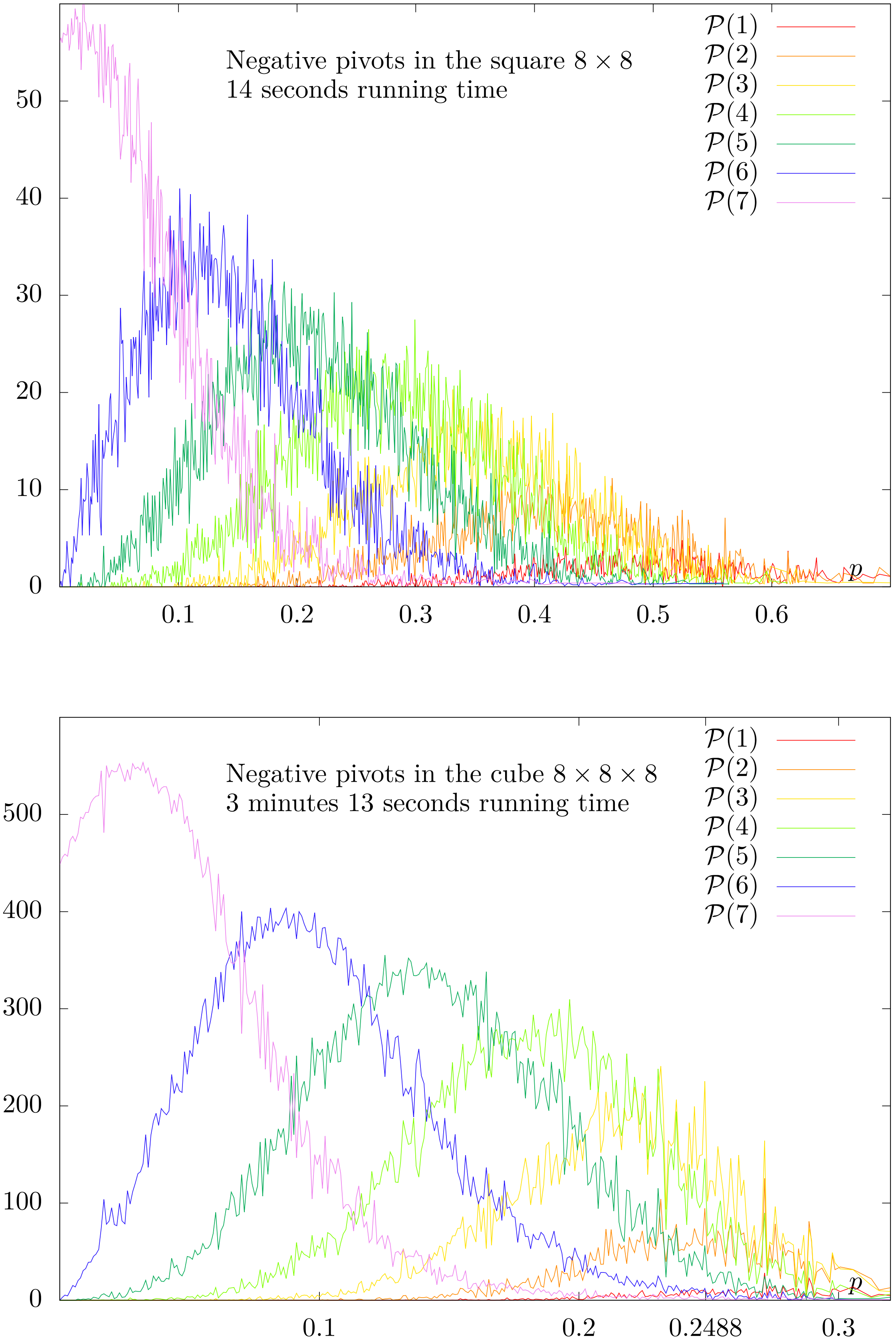

Figure 86: Negative pivots in the square and cube of side 8, as function of $p$. Each point is the average of 10 simulations. The increment is $\Delta p = 0.001$.

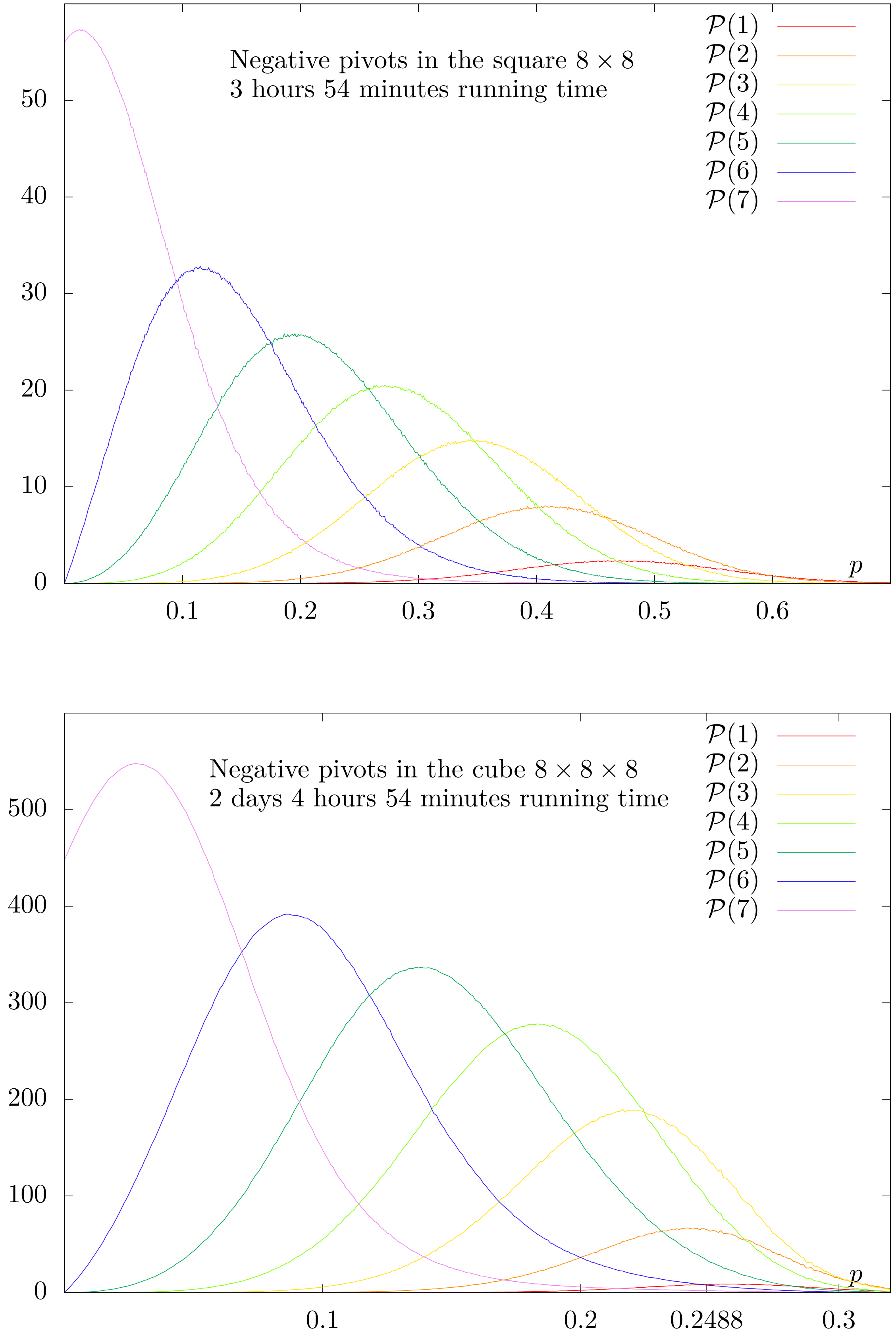

Figure 87: Negative pivots in the square and cube of side 8, as function of $p$. Each point is the average of 10000 simulations. The increment is $\Delta p = 0.001$.

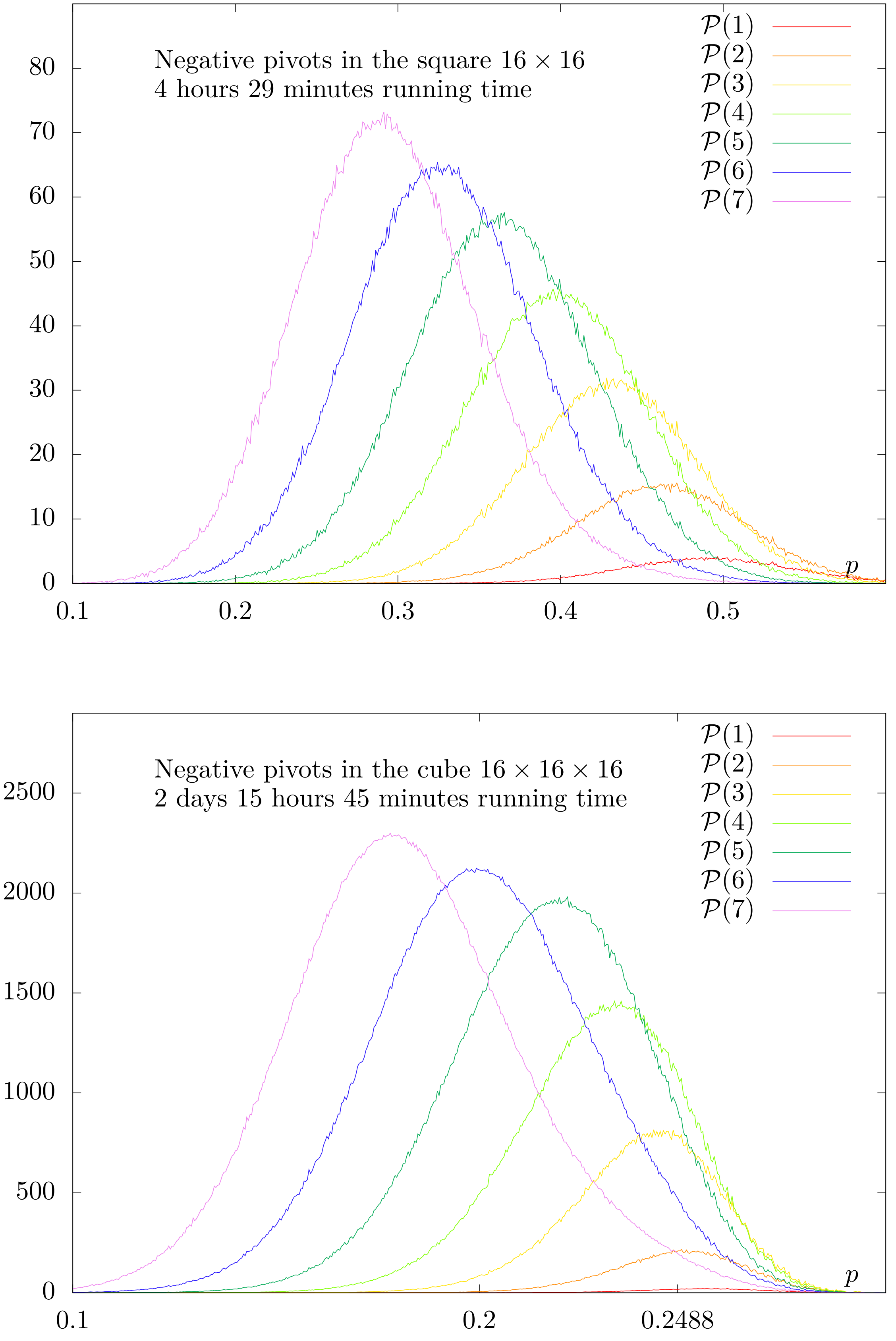

Figure 88: Negative pivots sets in square and cube of side 16, as function of $p$. Each point is the average of 1000 simulations. The increment is $\Delta p = 0.001$.

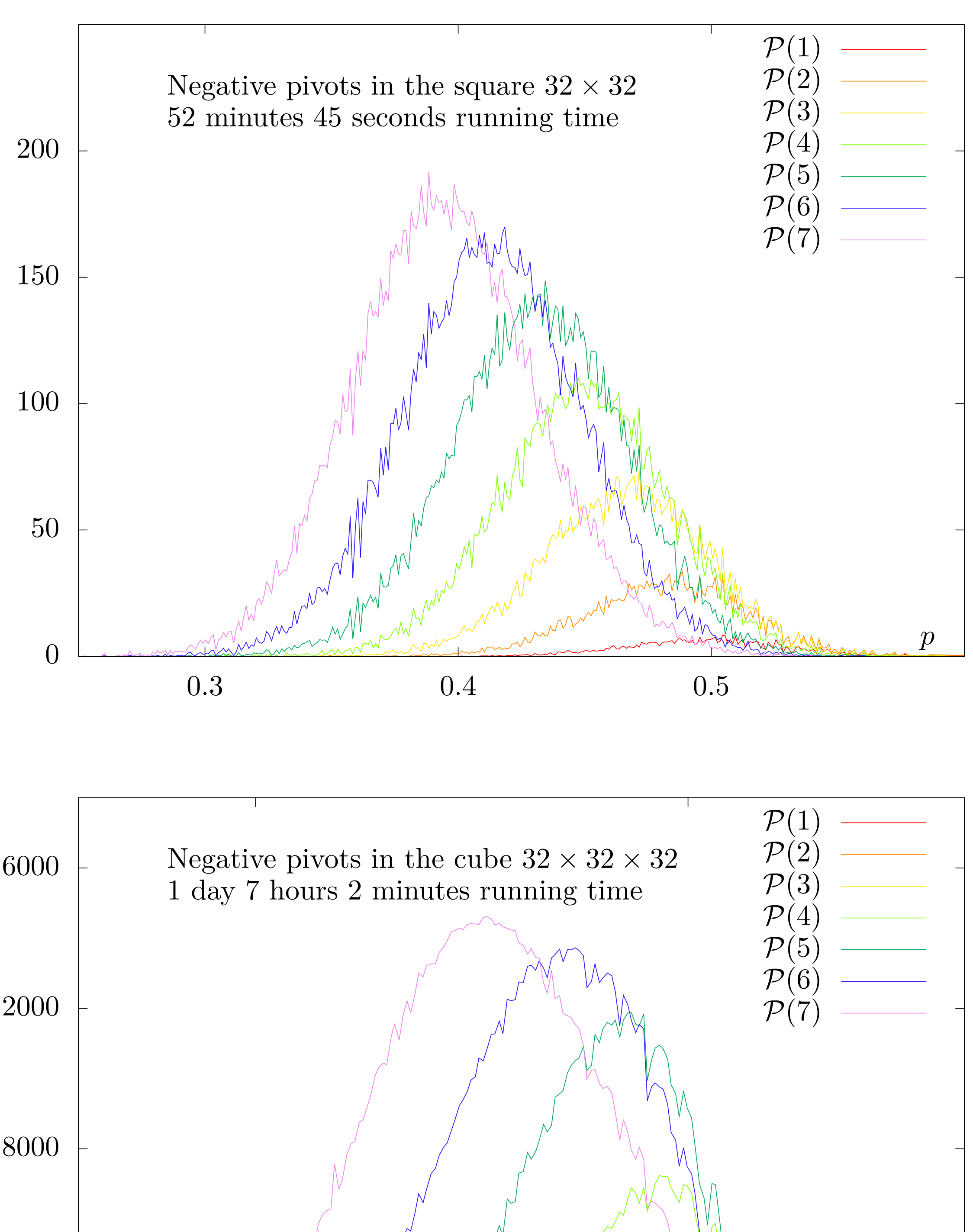

Figure 89: Negative pivots in the square and cube of side 32, as function of $p$. Each point is the average of 100 simulations. The increment is $\Delta p = 0.001$.

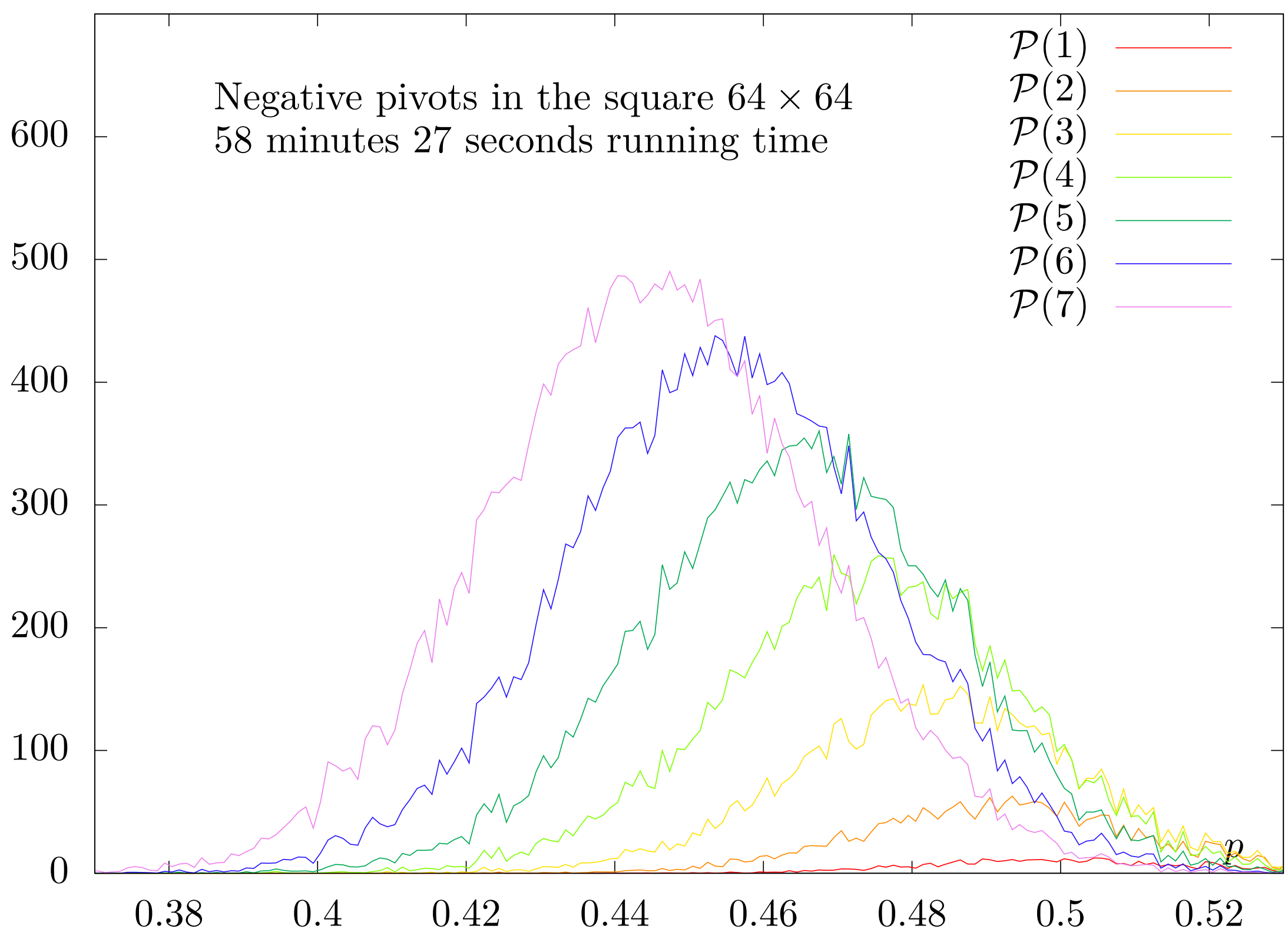

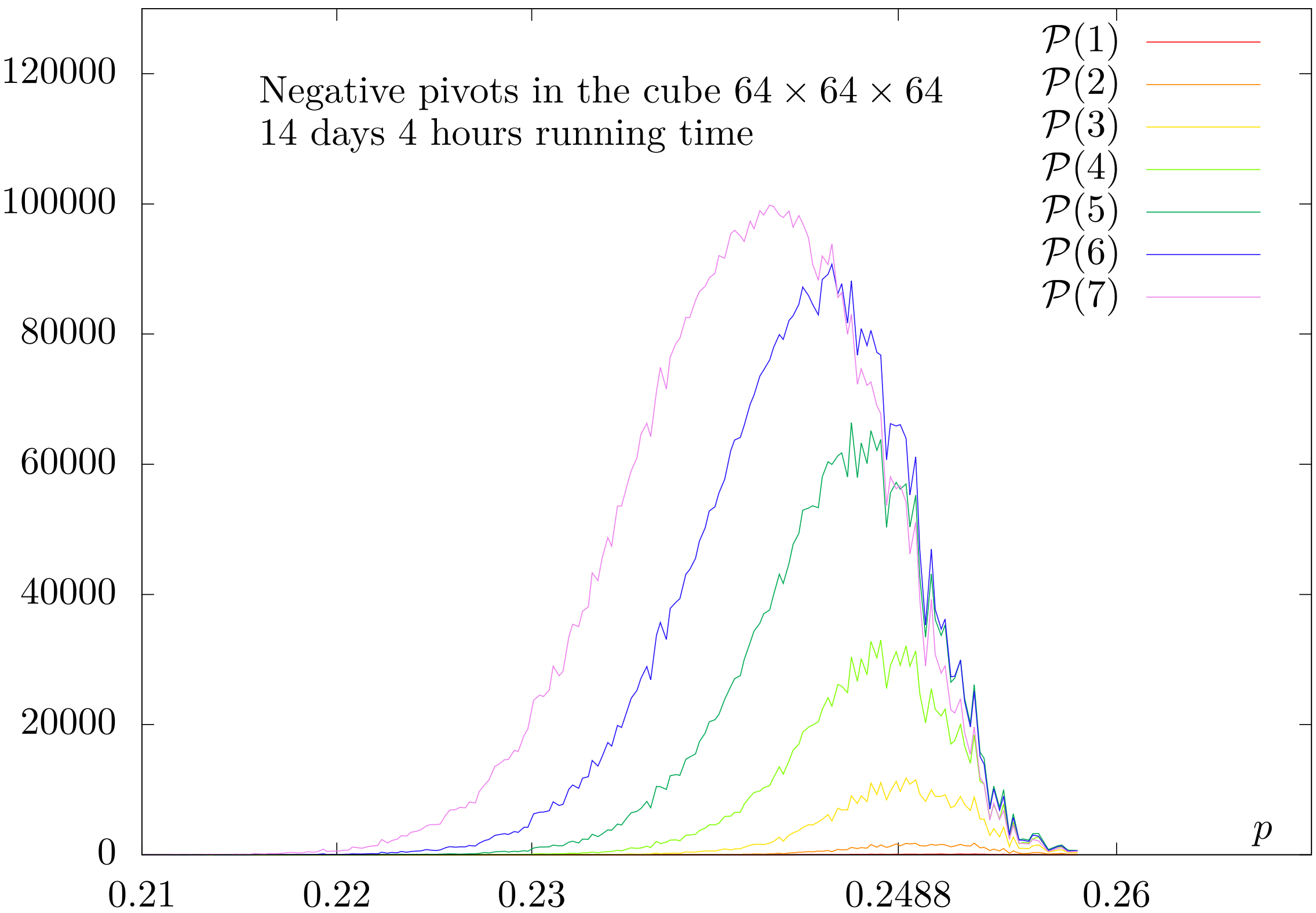

Figure 90: Negative pivots in the square and cube of side 64, as function of $p$. Each point is the average of 100 simulations. The increment is $\Delta p = 0.001$.

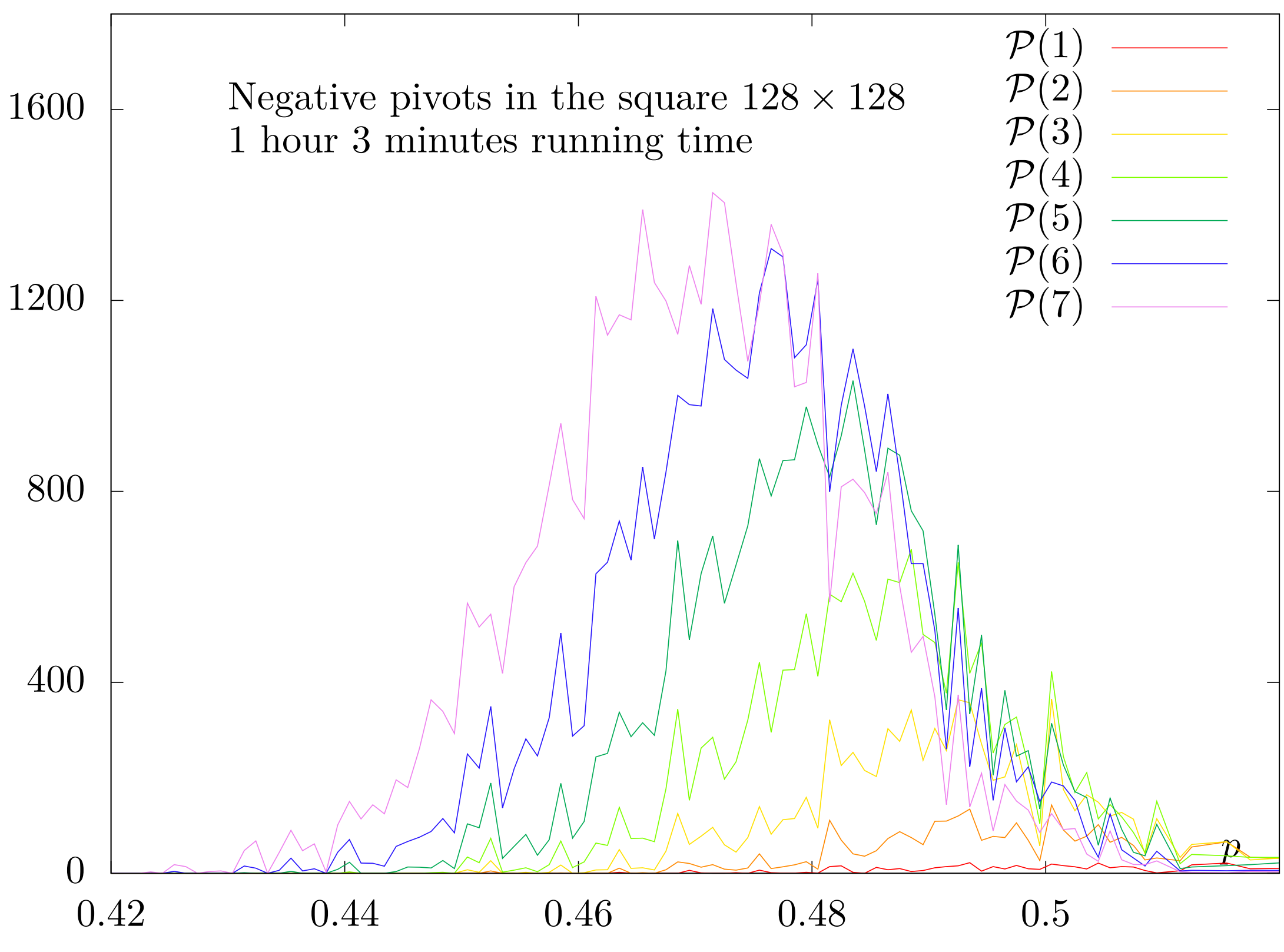

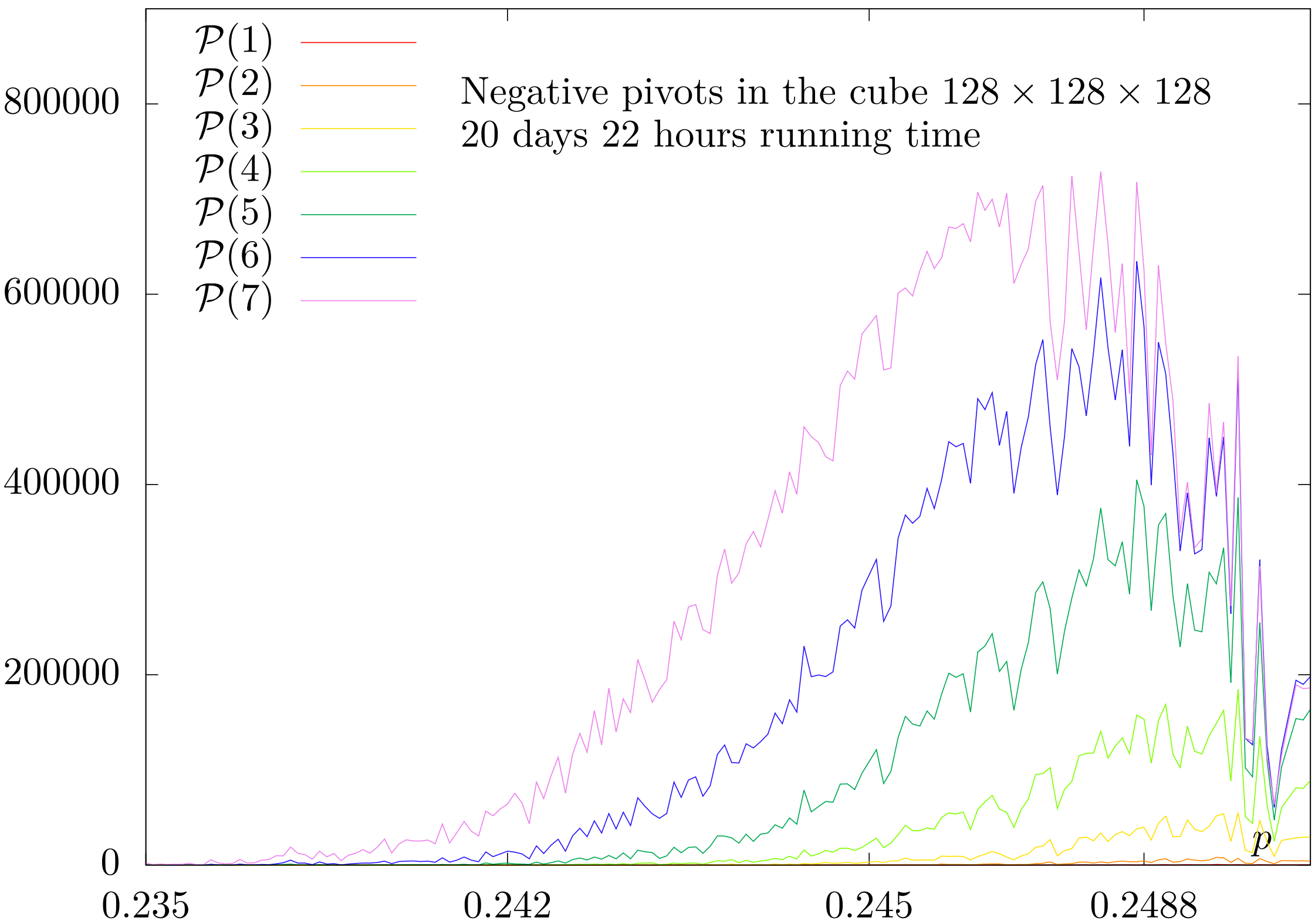

Figure 91: Negative pivots in the square and cube of side 128, as function of $p$. Each point is the average of 10 simulations. The increment is $\Delta p = 0.001$.

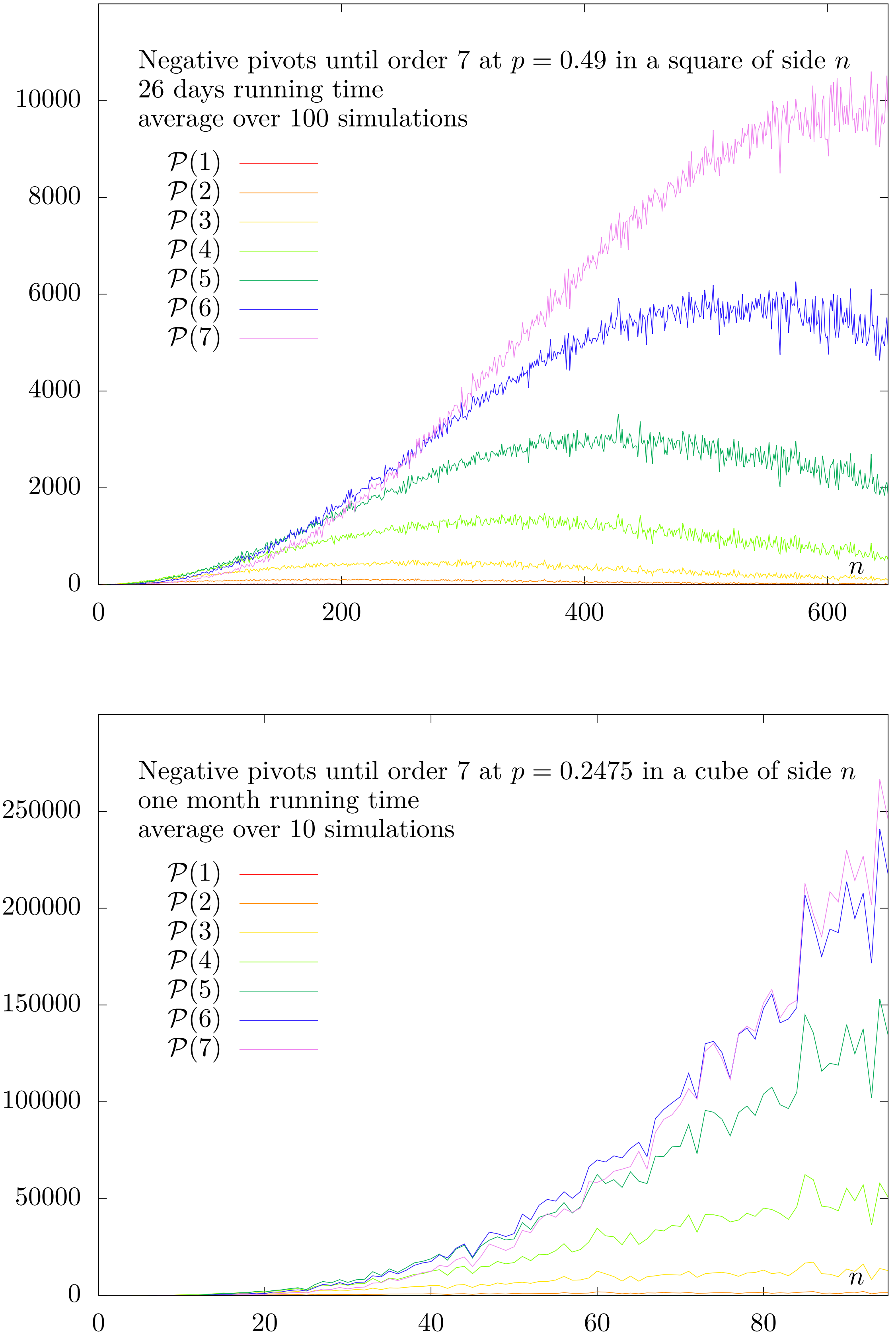

Figure 92: Negative pivots until order seven in a square of side $n$ at $p = 0.49$ and in a cube of side $n$ at $p = 0.2475$. The increment for $n$ is $\Delta n = 1$.

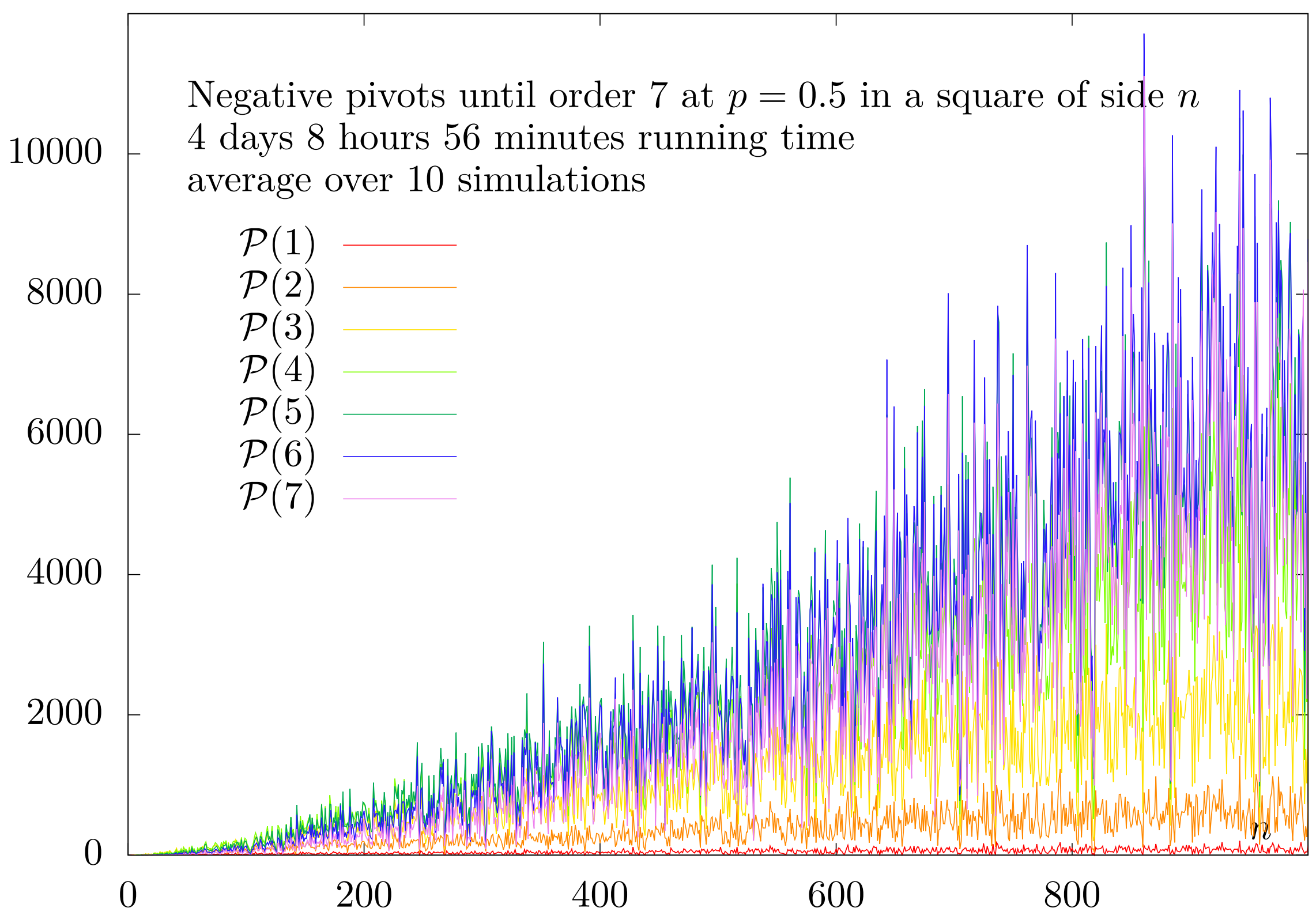

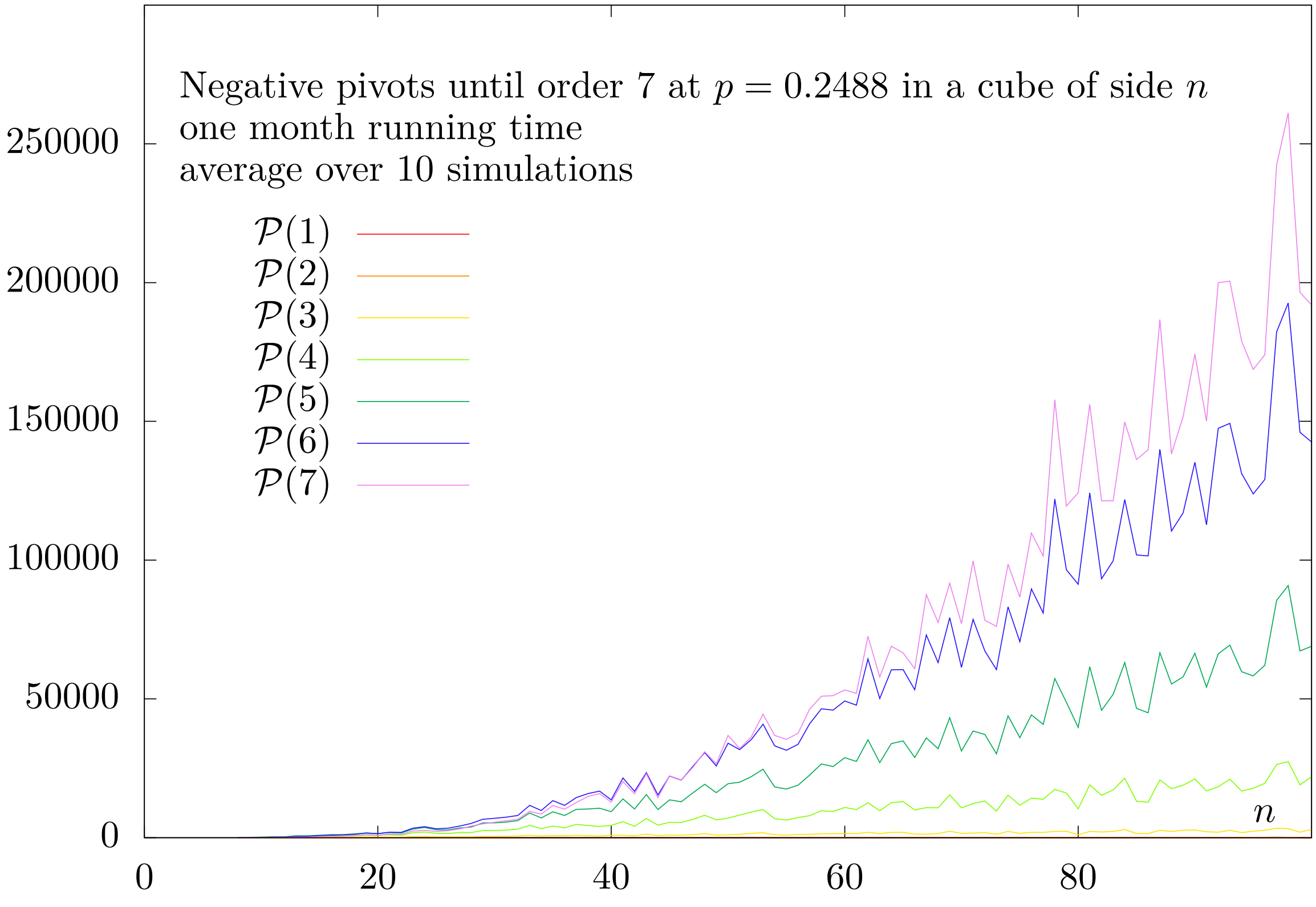

Figure 93: Negative pivots until order seven in a square of side $n$ at $p = 0.5$ and in a cube of side $n$ at $p = 0.2488$. The increment for $n$ is $\Delta n = 1$.

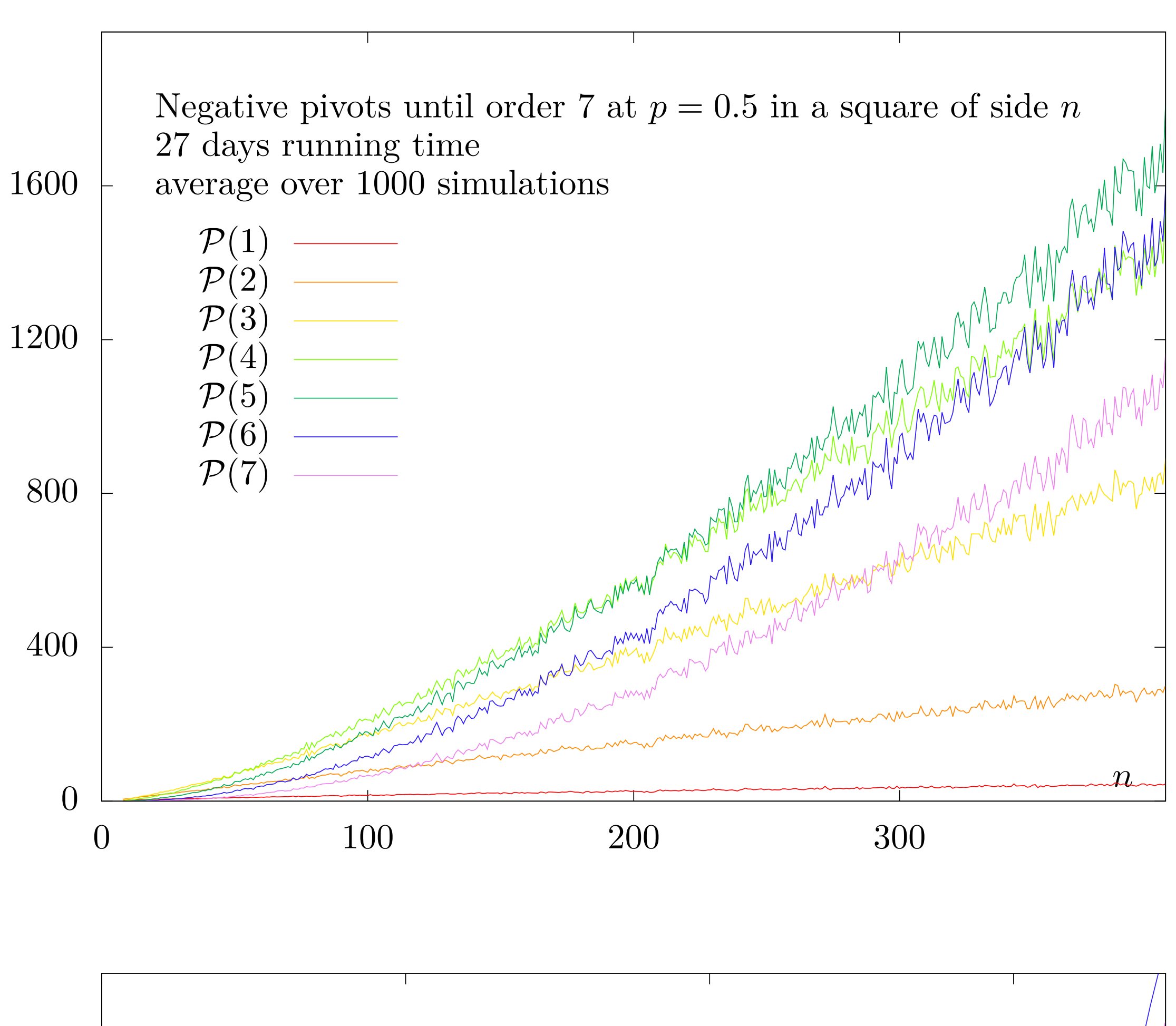

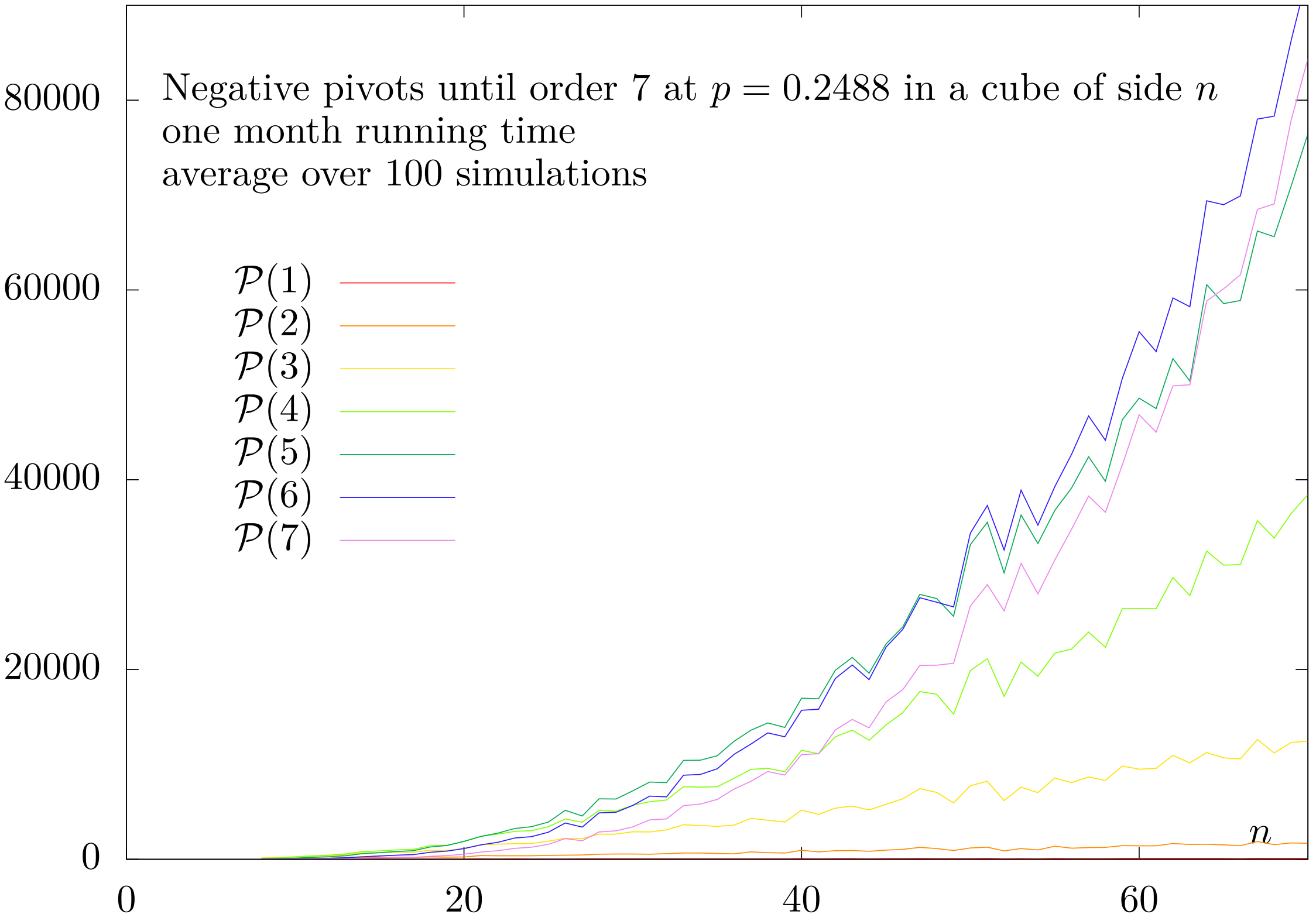

Figure 94: Negative pivots until order seven in a square of side $n$ at $p = 0.5$ and in a cube of side $n$ at $p = 0.2488$. The increment for $n$ is $\Delta n = 1$.

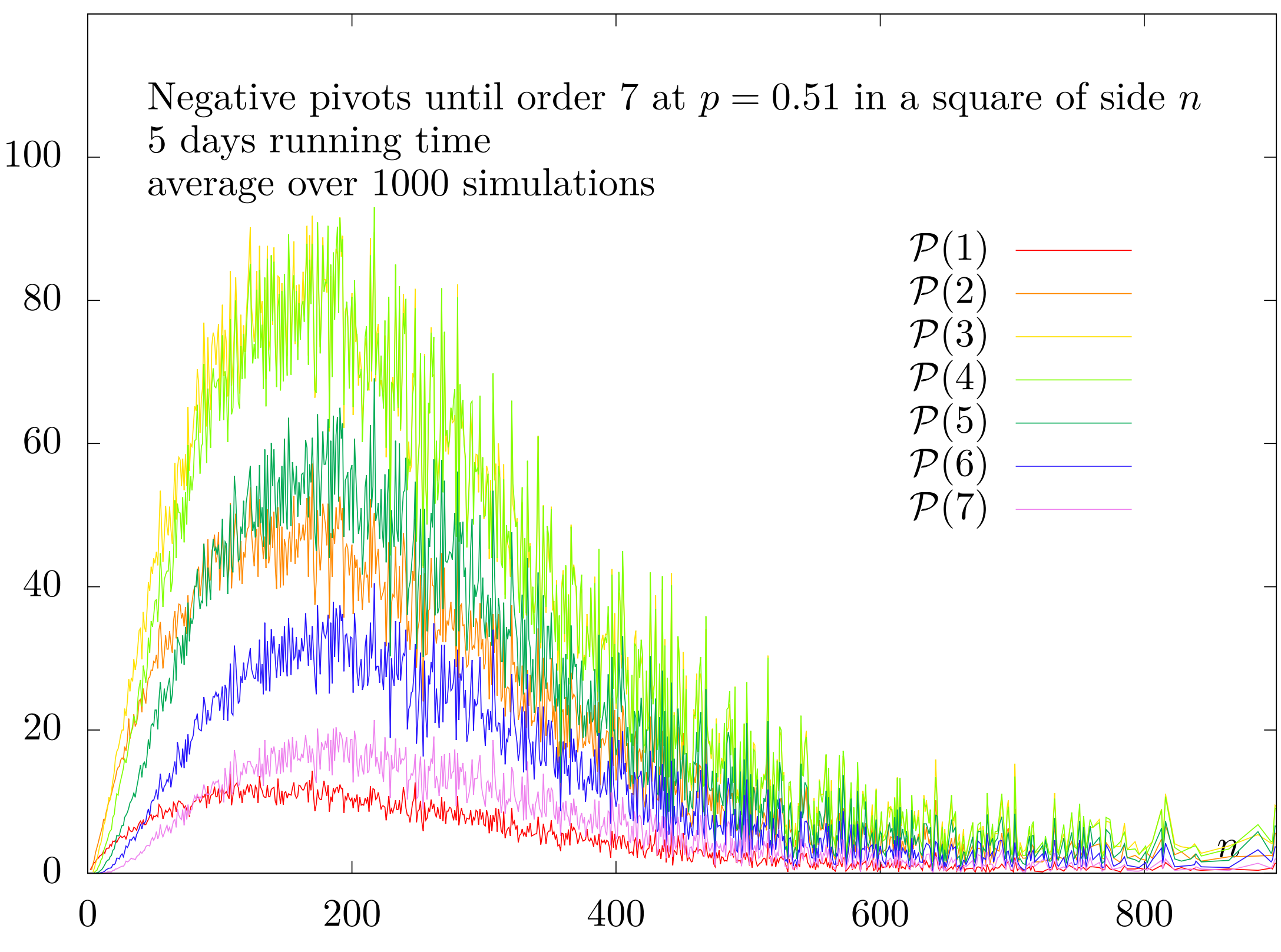

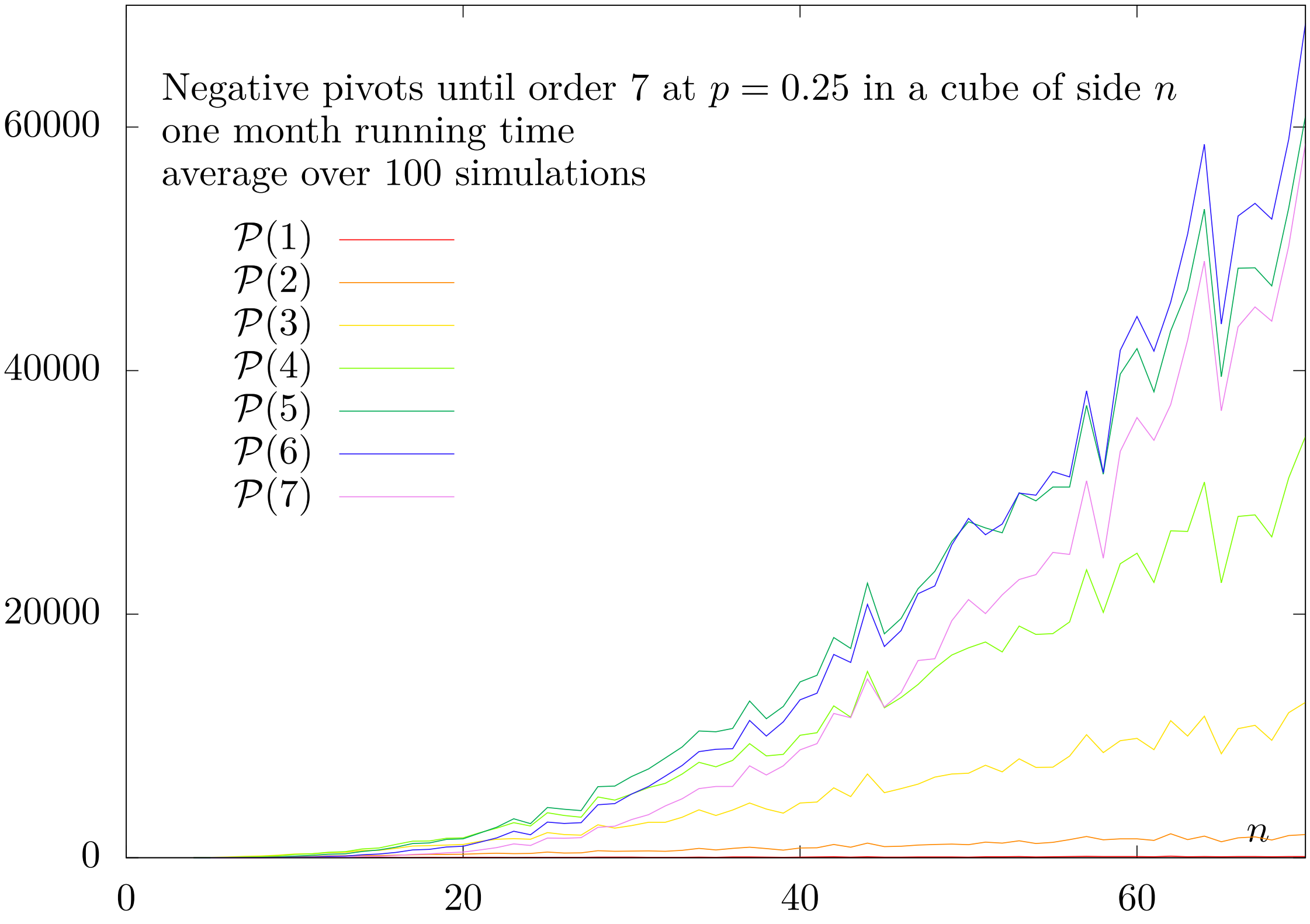

Figure 95: Negative pivots until order seven in a square of side $n$ at $p = 0.51$ and in a cube of side $n$ at $p = 0.25$. The increment for $n$ is $\Delta n = 1$.

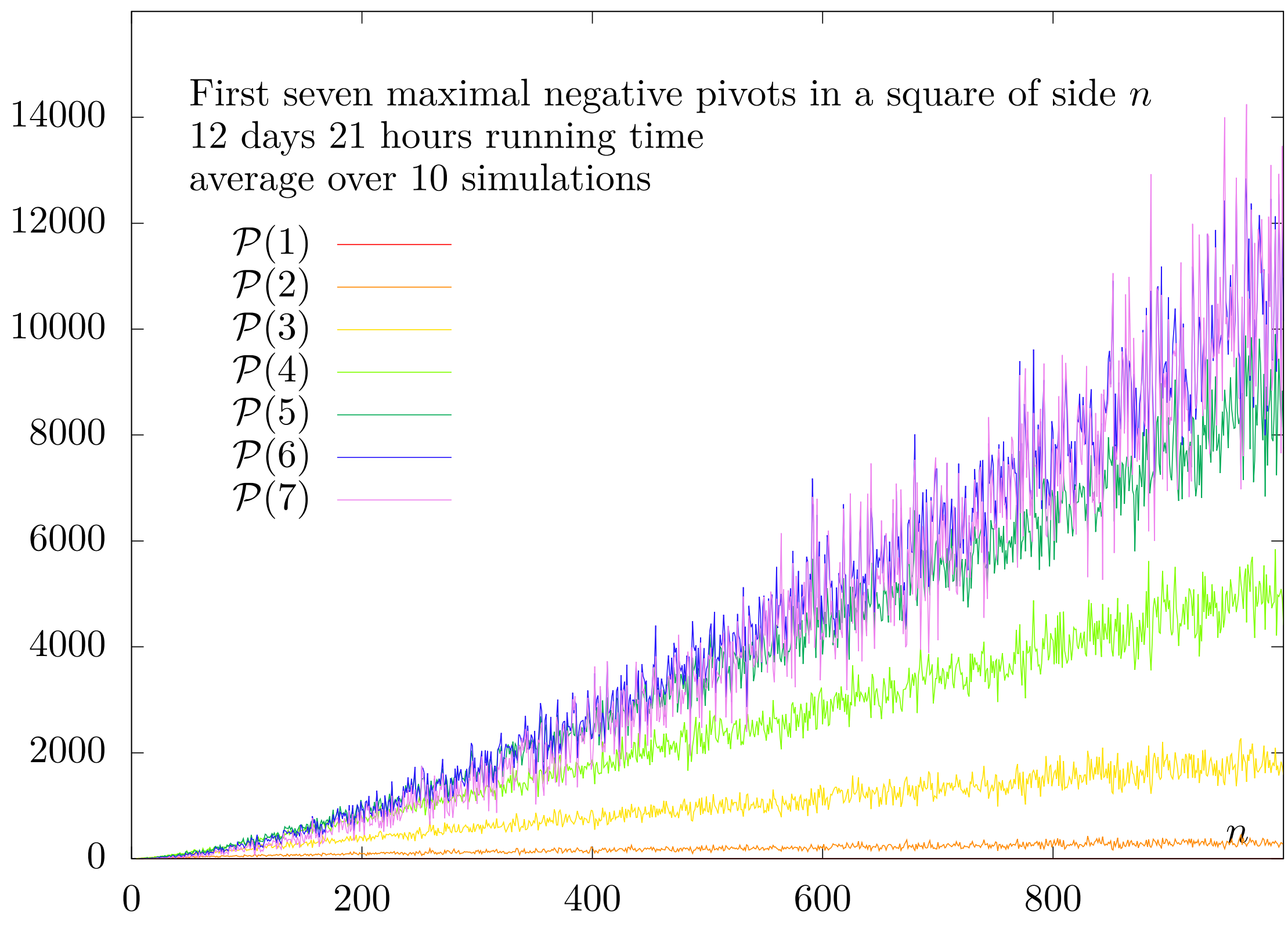

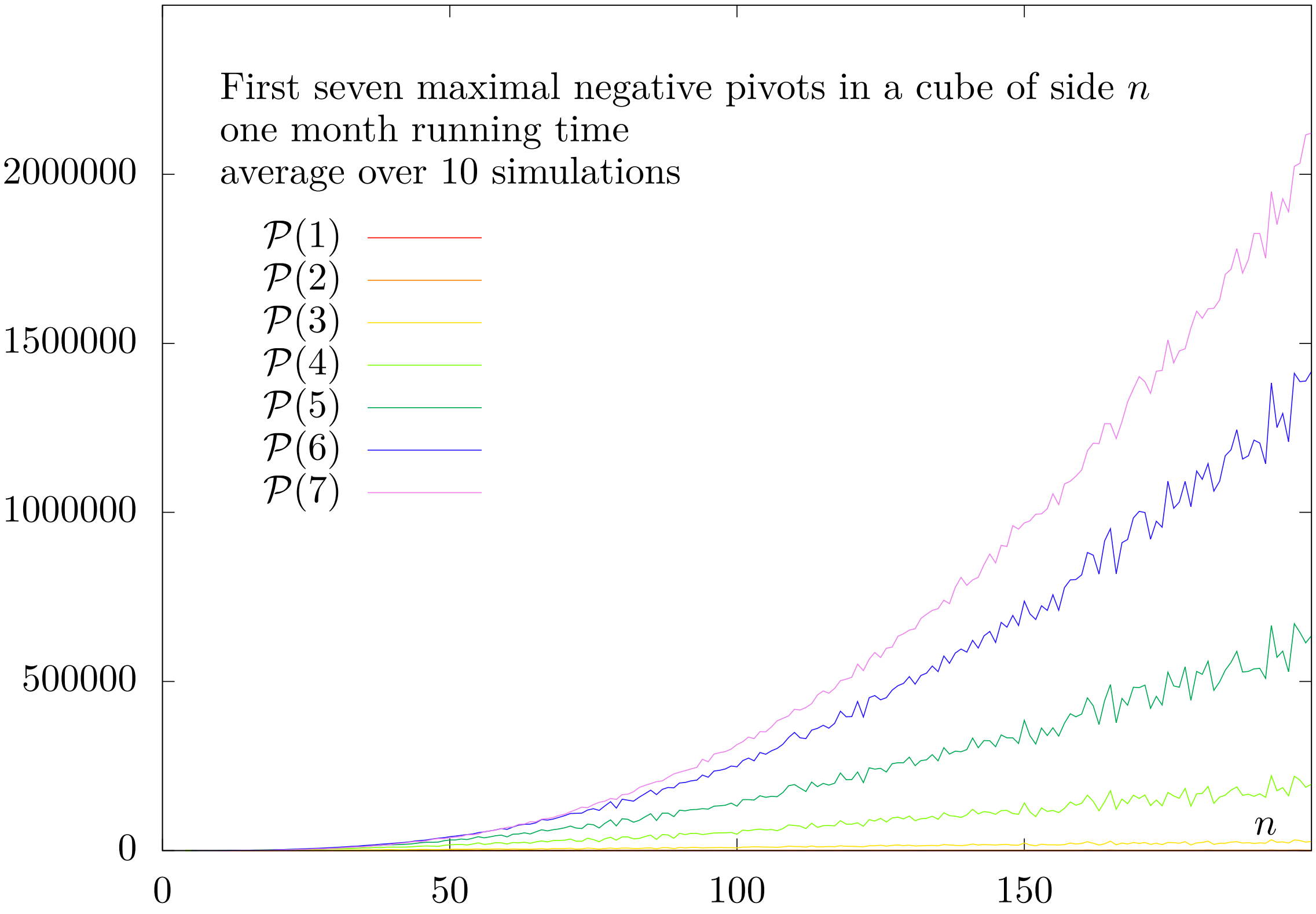

Figure 96: The first seven maximal intersection sets in a square and a cube: ultimate configuration before the connexion. The increment for $n$ is $\Delta n = 1$.

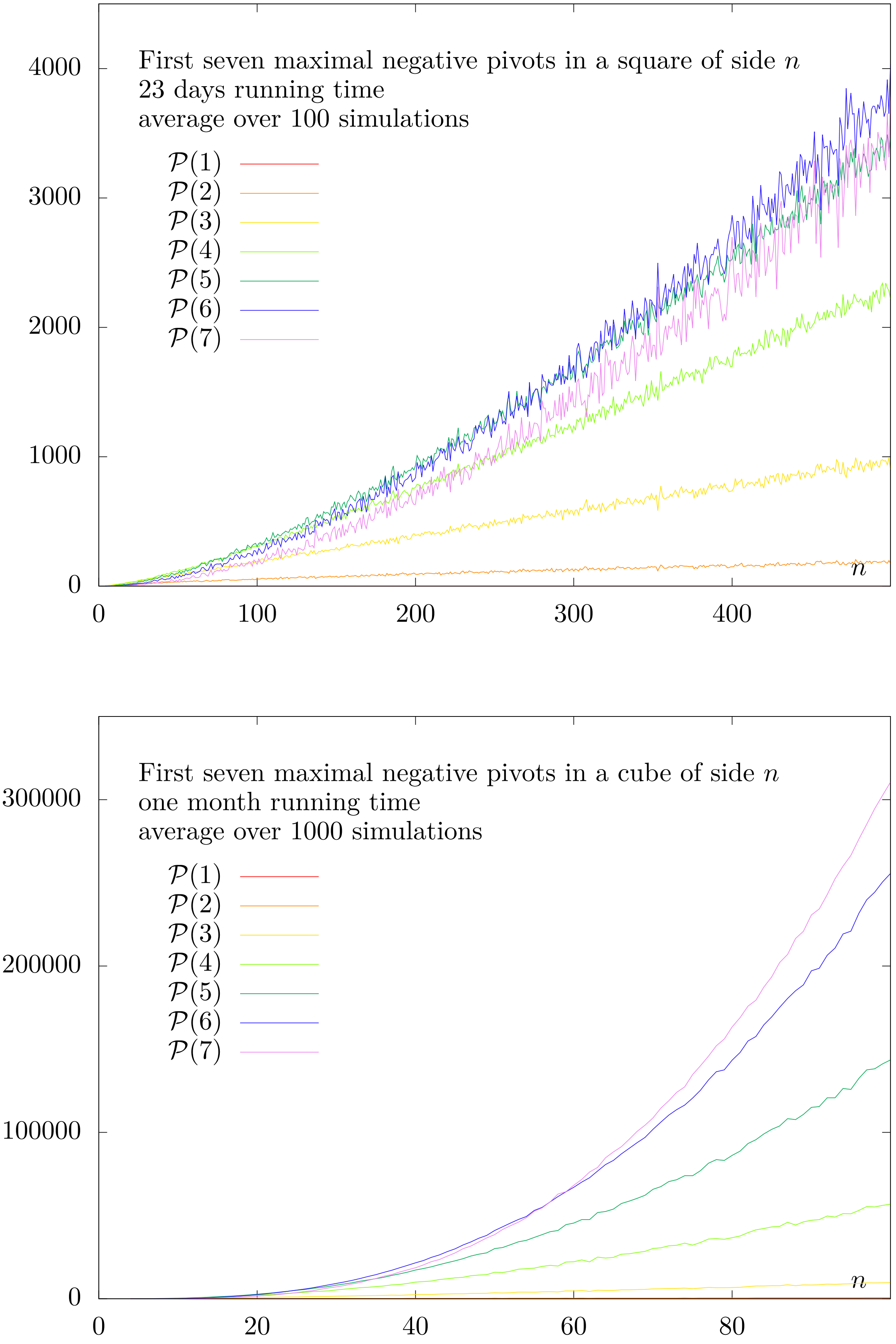

Figure 97: The first seven maximal intersection sets in a square and a cube: ultimate configuration before the connexion. The increment for $n$ is $\Delta n = 1$.

## 45.7 Second positive pivots and Tarjan's algorithm

The unsuccessful attempts of subsections 45.4, 45.5 and 45.6 had the merit of making us understand that the problem of controlling the intersection set $\mathcal{I}(2)$ or the second pivots $\mathcal{P}_2(\mathcal{D})$ are equivalent to the problem of controlling the second positive pivots $\mathcal{P}^+(\mathcal{P}(\mathcal{D}))$. Unfortunately, this question seems equally difficult than the initial one. In desperation, we can always resort to simulations. This is the goal of this subsection. We consider the event $\mathcal{D}(n)$ which is the disconnection between two opposite faces in the box $\Lambda(n)$. Moreover we consider only two-dimensional boxes, and we take for $W$ the left face of $\Lambda(n)$ and for $W'$ its right face. The definition of the left and right second positive pivotal bonds is naturally transposed in this context (simply use definition 45.9 with $\Lambda(n)$ instead of $\Pi(n)$). We performed various simulations to visualize the second positive pivots. To find the bonds of $\mathcal{P}_L^+(\mathcal{P}(\mathcal{D}))$, we will rely on a famous algorithm due to Tarjan [128]. This algorithm computes in fact the strongly connected components in a directed graph (this is more general than what we need, the graphs we consider here are undirected). Donald Knuth has written an implementation of Tarjan's algorithm in literate programming[4] (see [86], pages 512-519). An alternative algorithm is the Kosaraju-Sharir[5] algorithm [120]. In the process of computing the strongly connected components, Tarjan's algorithm identifies the bridges. These are the bonds $e$ whose endpoints are not any more connected in the graph deprived of the bond $e$. The second positive pivots are the bridges having one endpoint connected to the left or the right face, and the other endpoint connected to the endpoint of a pivotal bond for the disconnection event. Once the pivotal bonds $\mathcal{P}(\mathcal{D})$ are identified, a slight modification of Tarjan's algorithm allows to extract the relevant bridges and to find efficiently the second positive pivots $\mathcal{P}^+(\mathcal{P}(\mathcal{D}))$. In figures 98 and 99, we see the pivotal bonds for the event $\mathcal{D}$, as well as the negative and positive second pivots.

[4]On the occasion of the publication of the eBooks of The Art of Computer Programming (TAOCP), Jon Bentley asked Donald Knuth [87]: Of all the programs that you've written, what are some of which you are most proud, and why? Here are some excerpts from his reply: ... Of course I'm proud of TEX and METAFONT, because they seem to have helped to change the world, and because they led to many friendships. ... While I was preparing for Volume 4 of TAOCP in the 90s, I wrote several dozen short routines using "literate programming". My favorite is the implementation of Tarjan's beautiful algorithm for strong components, which appears on pages 512–519 of The Stanford GraphBase (1994). ...

[5]Kosaraju suggested it in 1978, Sharir discovered it independently and published it in 1981.

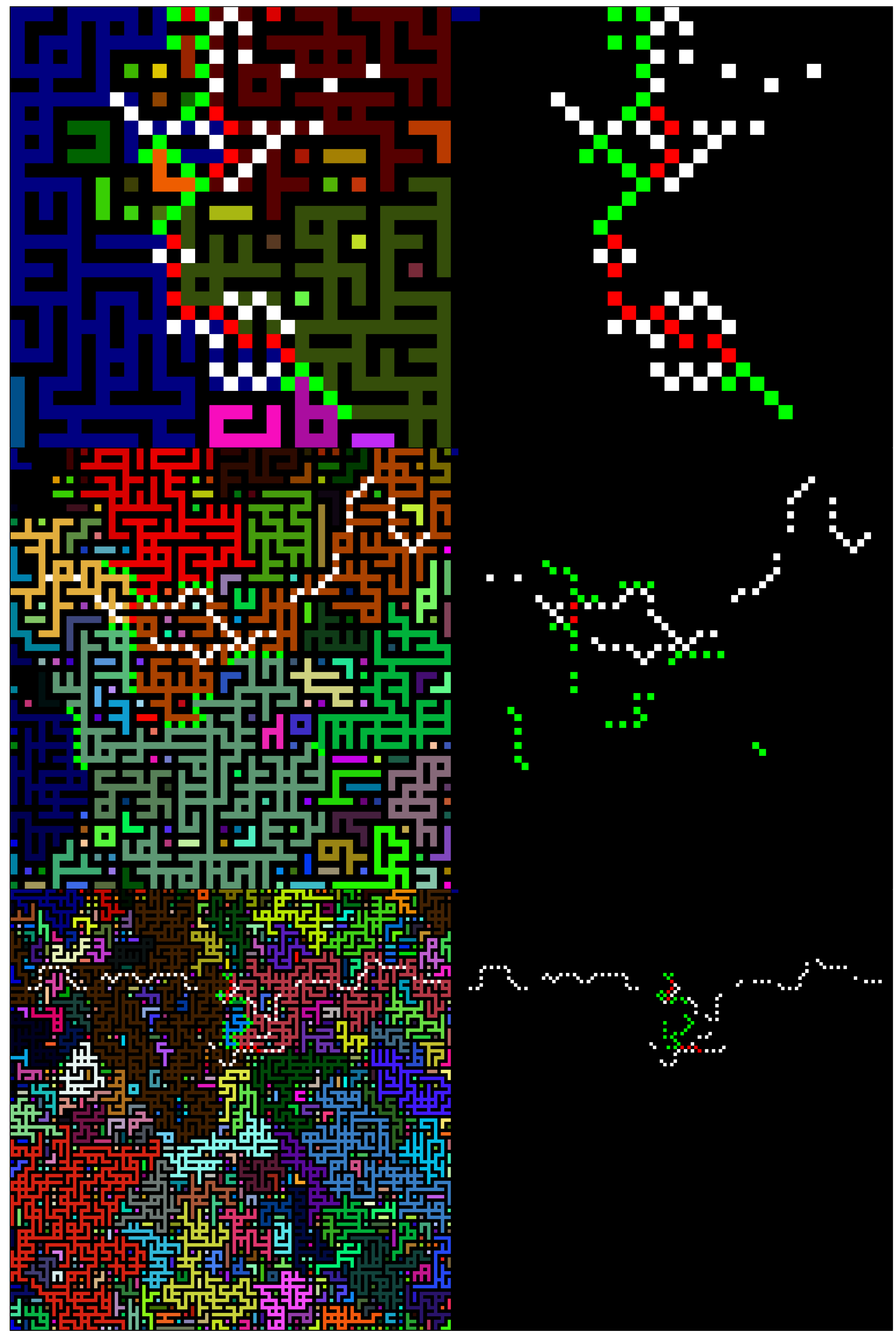

Figure 98: Pivots for the disconnection between left and right faces of a square: pivotal bonds, negative second pivotal bonds, positive second pivotal bonds. Left column: the colored clusters and the pivots. Right column: the pivots only. Boxes $\Lambda(16), \Lambda(32), \Lambda(64)$, $p = 0.443,\ 0.361,\ 0.340$ (top row to bottom row).

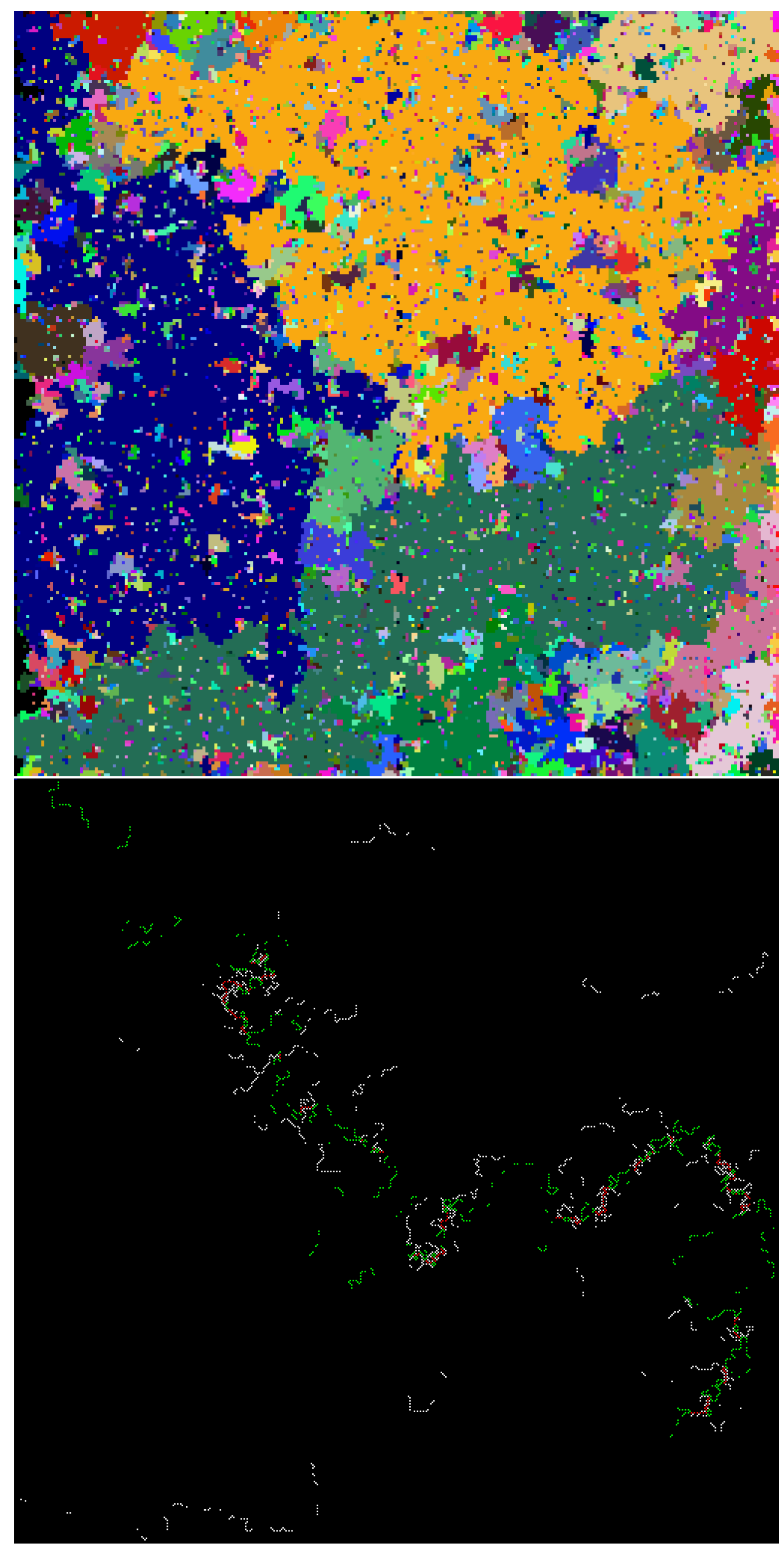

Figure 99: Continuation of figure 98. Bond percolation in $\Lambda(256)$, $p = 0.396$. With a larger box, the pivotal bonds are no longer discernible on the picture.

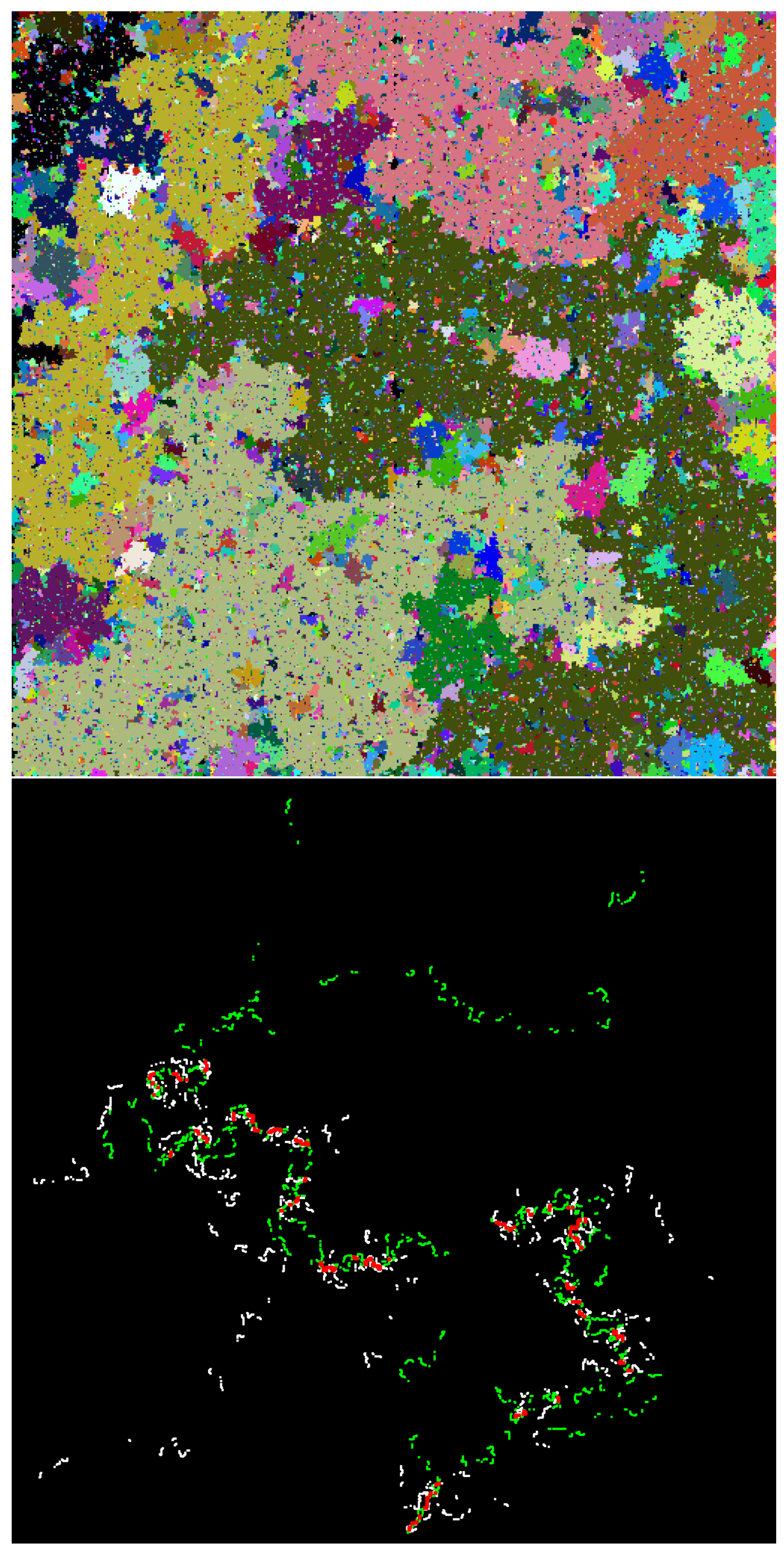

Figure 100: Continuation of figure 99. Bond percolation in $\Lambda(512)$, $p = 0.400$. Pivotal bonds are dilated by 2, negative and positive second pivotal bonds by 1.

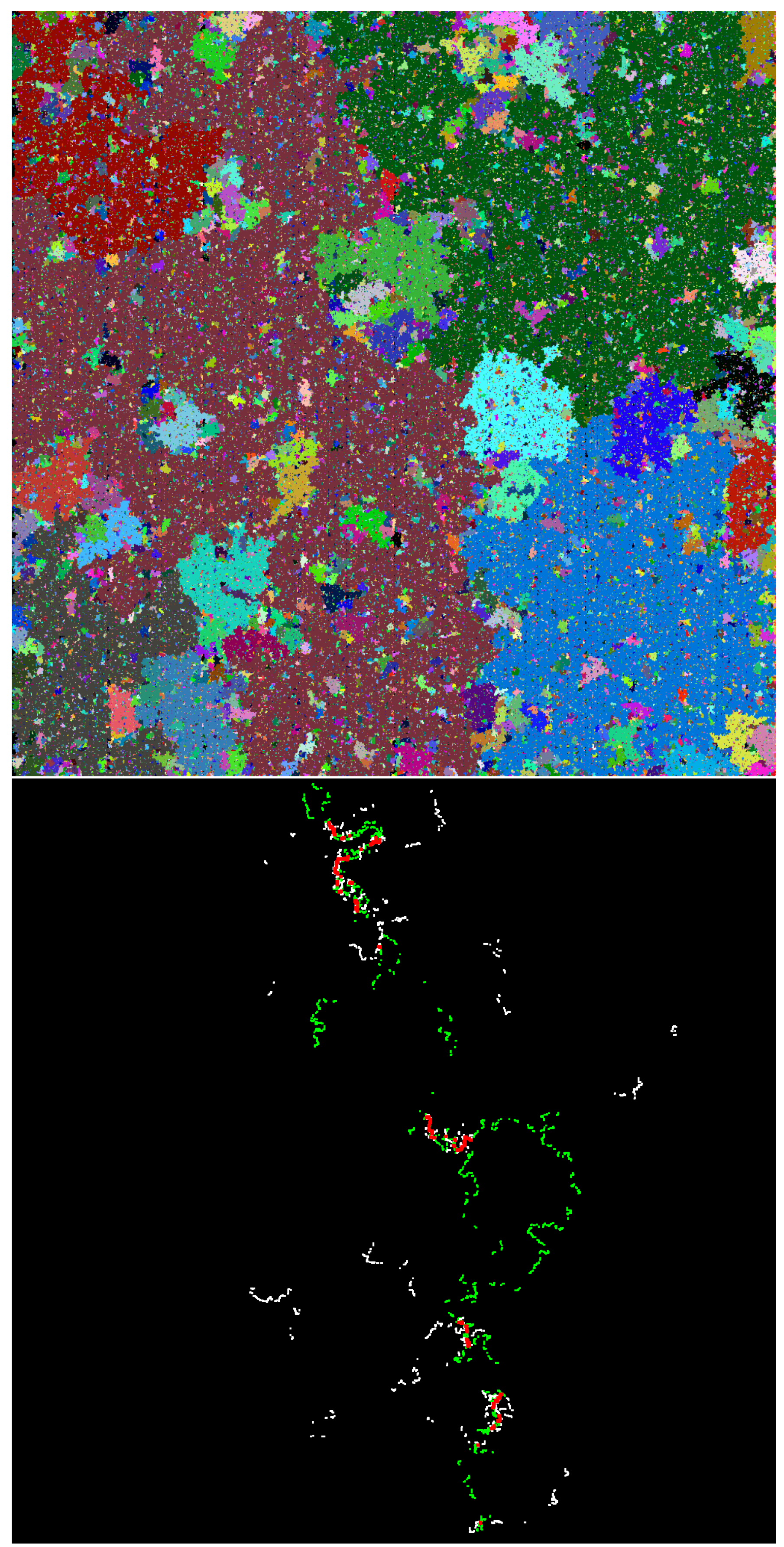

Figure 101: Continuation of figure 100. Bond percolation in $\Lambda(1024)$, $p = 0.401$. Pivotal bonds are dilated by 4, negative and positive second pivotal bonds by 2.

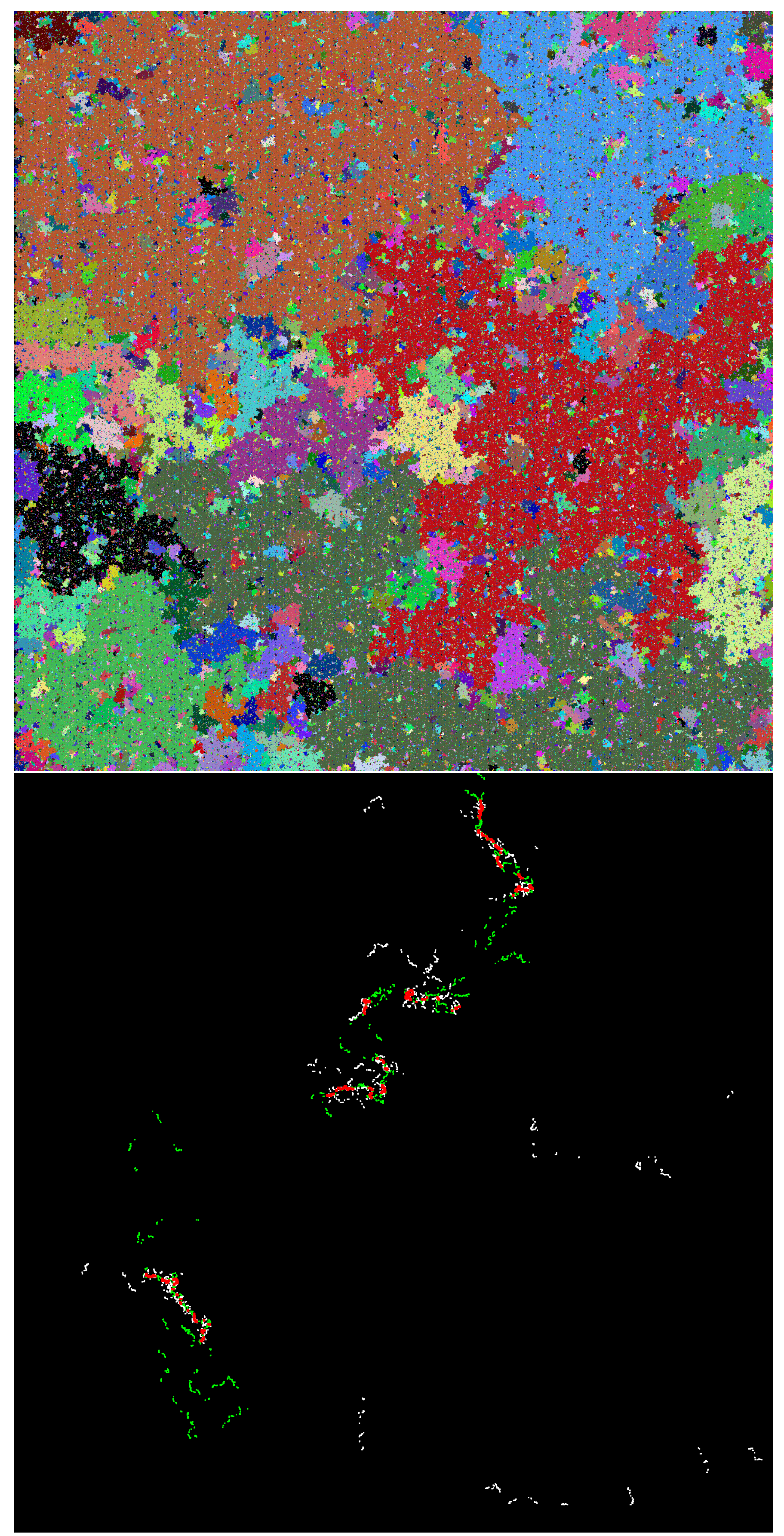

Figure 102: Continuation of figure 101. Bond percolation in $\Lambda(2048)$, $p = 0.401$. Pivotal bonds are dilated by 6, negative and positive second pivotal bonds by 3.

# 46 Two inequalities on the higher pivots ⊛

We present here a generalization of two inequalities of Talagrand, Benjamini, Kalai, Schramm, Keller and Kindler, which involve the higher pivots. Unfortunately, they do not seem to be helpful for our percolation problem. We recapitulate in the next subsection the notation used in the results and their proofs, so that the whole section is completely self-contained.

## 46.1 Basic notation

Let $n \geq 1$ be an integer and let us consider the configuration space $\Omega = \{0,1\}^n$. For $\omega \in \Omega$, we denote by $\omega(i)$ its $i$-th component, so that $\omega = \big(\omega(1), \dots, \omega(n)\big)$. Let $p \in ]0,1[$. We endow $\Omega$ with the product probability measure $P_p$ defined by

$$\forall \omega \in \Omega \qquad P_p(\omega) \,=\, \prod_{i=1}^{n} p^{\omega(i)}(1-p)^{1-\omega(i)}\,.$$

We define $n$ random variables $(Y_i, 1 \leq i \leq n)$ by setting

$$\forall i \in \{1,\dots,n\} \qquad Y_i(\omega) \,=\, \begin{cases} -\sqrt{\dfrac{p}{1-p}} & \text{if } \omega(i) = 0\,,\\ \sqrt{\dfrac{1-p}{p}} & \text{if } \omega(i) = 1\,. \end{cases}$$

For $K$ a subset of $\{1,\dots,n\}$, we define the Walsh product

$$Y_K \,=\, \prod_{i \in K} Y_i\,.$$

Let $f$ be a function from $\Omega$ to $\mathbb{R}$. For $i \in \{1,\dots,n\}$, we define the discrete derivative along the $i$-th component as

$$\forall \omega \in \Omega \qquad \delta_i f(\omega) \,=\, f(\omega^i) - f(\omega_i)\,,$$

and it can be checked through elementary computations that

$$E\big(fY_i\big) \,=\, \sqrt{p(1-p)}E\big(\delta_i f\big)\,. \tag{46.1}$$

For $K$ a subset of $\{1,\dots,n\}$, the operator $\delta_K$ is the product of the operators $\delta_i$, $i \in K$, that is

$$\delta_K \,=\, \prod_{i \in K} \delta_i\,,$$

and successive applications of the identity (46.1) yield that

$$E\big(fY_K\big) \,=\, \big(p(1-p)\big)^{|K|/2} E\big(\delta_K f\big)\,. \tag{46.2}$$

Notice that $\delta_K f$ is a function which does not depend any more on the components of $\omega$ whose index belongs to $K$.

## 46.2 The two inequalities

We are now in position to state the two main results of the section.

**Theorem 46.1.** *Let $A$ be a subset of $\{0,1\}^n$. Let $I, J$ be a partition of $\{1,\dots,n\}$. Let $d \geq 2$ and $r \geq 1$ be two integers, let $s > 0$ be a real number and let us define*

$$L = \Bigg\{ K \subset J : |K| = r, \quad \Bigg( \sum_{\substack{T \subset I \\ |T| = d-1}} \big(E(Y_T Y_K 1_A)\big)^2 \Bigg)^{1/2} \geq s\big(p(1-p)\big)^{r/2} E\big(|\delta_K 1_A|\big) \Bigg\}. \tag{46.3}$$

*There exists a positive constant $c$, which depends on $d$ and $r$, but not on $s$, such that*

$$\sum_{K \in L} \Big( E\big(|\delta_K 1_A|\big)\Big)^2 \leq \big(p(1-p)\big)^{-r} \frac{1}{c} \exp\big(-cs^{2/(d-1)}\big) P_p(A). \tag{46.4}$$

**Theorem 46.2.** *Let $f : \{0,1\}^n \to \{0,1\}$ be a Boolean function. Let $I, J$ be a partition of $\{1,\dots,n\}$, let $d \geq 2$ and $r \geq 1$ be two integers, and let us define*

$$H(J,r,f) = \big(p(1-p)\big)^r \sum_{K \subset J,\, |K| = r} E\big(|\delta_K f|\big)^2. \tag{46.5}$$

*We have the inequality*

$$\sum_{\substack{T \subset I \\ |T| = d-1}} \sum_{\substack{K \subset J \\ |K| = r}} \big(E(f Y_T Y_K)\big)^2 \leq 16^d (d!)^2 H(J,r,f) \Big(\frac{B(p)e}{d-1}\Big)^{d-1} \left(1 \vee \ln\left(\frac{E(f)}{H(J,r,f)}\right)\right)^{d-1}. \tag{46.6}$$

Let us first discuss the relationship between the previous existing works and the inequalities stated in theorems 46.1 and 46.2. The mathematicians who contributed to the development of these inequalities are, in chronological order, Talagrand [125], Benjamini, Kalai and Schramm [15], Keller and Kindler [80]. For $j$ in $\{1,\dots,n\}$, the classical notion of influence of the $j$-th coordinate on $f$ is defined as

$$I_j(f) = P_p\big(\{\, \omega \in \Omega : f(\omega^j) \neq f(\omega_j) \,\}\big),$$

and it can be rewritten with the operator $\delta_j$ as

$$I_j(f) = E\big(|\delta_j f|\big).$$

In particular, for $r = 1$, the quantity $H(J,1,f)$ becomes simply

$$H(J,1,f) = p(1-p) \sum_{j \in J} \big(I_j(f)\big)^2.$$

If we take $d = 2$, $r = 1$, $p = 1/2$ in theorems 46.1 and 46.2, and suppose in addition that $A$ is monotone, then we recover lemma 3.1 and proposition 2.3 of Talagrand [125], albeit with the additional factor $P_p(A)$ (this factor is bounded by 1 in the course of Talagrand's proof). However, we might consider events $A$ whose probability is small, and the factor $P_p(A)$ might improve considerably the inequality. The version of theorems 46.1 and 46.2 for an increasing set $A$ with $d \geq 2$, $r = 1$, $p = 1/2$ was obtained by Benjamini, Kalai and Schramm (see theorem 2.6 in [15]). It was a key inequality in the derivation of their famous noise sensitivity theorem. The version of theorems 46.1 and 46.2 for a general set $A$ (not necessarily increasing) with $d \geq 2$, $r = 1$, and with respect to the biased Bernoulli measure of parameter $p$ was obtained by Keller and Kindler (see lemma 14 in the appendix of [80] and inequality (12) in the proof of lemma 6 of [80]). Our contribution here is the generalization of these inequalities to an arbitrary positive integer $r$. Also, our inequalities include an extra factor $P_p(A)$ or $E(f)$, which is a little improvement of the previous inequalities for $r = 1$. A little final detail: Keller and Kindler had an extra condition, they required an a priori bound on $H(J, 1, f)$, of the form

$$H(J, 1, f) \leq \exp(-2(d-1)) . \tag{46.7}$$

However, when $H(J, 1, f)$ becomes large, the logarithm in (46.6) becomes small, even negative, so we can get rid of the condition (46.7) by taking the maximum between 1 and the logarithm factor.

Let us discuss next the relationship between the two theorems 46.1 and 46.2. Talagrand, Benjamini, Kalai and Schramm proved first theorem 46.1 and used it to prove theorem 46.2. Keller and Kindler [80] introduced a new scheme of proof, much more transparent than the original argument of Talagrand. They managed to prove directly theorem 46.2 for a general event $A$, thereby bypassing the estimate provided by theorem 46.1. In their appendix, they also provide a derivation of theorem 46.1 along the same lines as their proof of theorem 46.2. Thus it seems that theorem 46.1 is a stronger result than theorem 46.2. So there exist now two strategies for proving theorems 46.1 and 46.2. We can walk in the footsteps of Talagrand, or follow the streamlined arguments of Keller and Kindler. If the aim is to prove theorem 46.2 only, then the strategy of Keller and Kindler is certainly the best way. However, we were more interested in the estimate of theorem 46.1, because our goal was to gain some control on the pivots of order higher than one to solve our percolation problem. Therefore we shall examine in detail the original argument of Talagrand, and compare the outcome with the one of Keller and Kindler.

There is a slightly awkward point about the proofs. In fact, once the correct statement of the inequalities is found, the remaining work consists in checking that the existing proofs can be adapted with the necessary minor modifications. In other words, there is no notable new mathematical idea in the proofs themselves, compared to the proofs provided by Talagrand or Keller and Kindler. For instance, Benjamini, Kalai and Schramm [15] made the choice to point only the important modifications of Talagrand's proof to go from the case $d = 2$ to a general $d \geq 2$. In our case, we really wish to compare the two lines of proofs, and we do feel it necessary to provide full proofs (otherwise we would essentially stop here and leave it to the reader to check the details, but since the inequalities have undergone several generalizations since Talagrand's paper, that would not be very fair to the reader). So we do reproduce the full arguments, but with a minimum of details, and we will follow quite precisely Talagrand on the one hand, and Keller and Kindler on the other hand. Our only contributions to these proofs are limited reformulations of some passages, in particular we use at several points the conditional expectation.

**Conditional expectation.** We will mainly use the notation of Keller and Kindler throughout the proofs, with one important difference. We do not rely on the variables $\omega_T$ as Keller and Kindler, and we always use $\omega$ as the integration variable in the expectations (whereas Keller and Kindler use $x$ and $y$). The reason is that we have rewritten some passages with the help of conditional expectations. The change is cosmetic, but it gives a more probabilistic flavor to the computations. For $I$ a subset of $\{1,\dots,n\}$, we introduce the $\sigma$-field $\sigma(I)$ generated by the components of $\omega$ whose index belongs to $I$, that is,

$$\sigma(I) \,=\, \sigma\big(\omega(i), i \in I\big)\,.$$

**Large deviations inequality**. A crucial ingredient that has been used to generalize the original inequalities of Talagrand is an inequality controlling the tail of a function $f$ on $\Omega$ whose Walsh expansion has degree at most $d$. We restate here the bound used by Keller and Kindler (see theorem 10 in [80]), which comes from lemma 2.2 of [41]. We refer to the paper of Keller and Kindler for more explanations on this inequality. We say that a function $f:\{0,1\}^n \to \mathbb{R}$ has a Walsh expansion of degree at most $d$ if and only if

$$\forall T \subset \{1,\dots,n\} \qquad |T| > d \quad \Longrightarrow \quad E(fY_T) = 0\,.$$

**Theorem 46.3.** *Let $p \in ]0,1[$ and let*

$$B(p) \,=\, \frac{(1-p)-p}{2p(1-p)(\ln(1-p)-\ln p)}\,.$$

*Let $f:\{0,1\}^n \to \mathbb{R}$ be a function whose Walsh expansion has degree at most d, and assume that $E(f^2)=1$. Then for any $t \geq (2B(p)e)^{d/2}$,*

$$P_p\big(|f| \geq t\big) \leq \exp\Big(-\frac{d}{2B(p)e}t^{2/d}\Big). \tag{46.8}$$

For convenience, we state a little variant which is valid for all $t > 0$.

**Corollary 46.4.** *Let $f:\{0,1\}^n \to \mathbb{R}$ be a function whose Walsh expansion has degree at most d, and assume that $E(f^2)=1$. Then*

$$\forall t > 0 \qquad P_p\big(|f| \geq t\big) \leq \exp\Big(d-\frac{d}{2B(p)e}t^{2/d}\Big). \tag{46.9}$$

*Proof.* For $t \geq (2B(p)e)^{d/2}$, this is a consequence of (46.8). For $t < (2B(p)e)^{d/2}$, the right-hand quantity in (46.9) is larger than 1. ☐

## 46.3 Why does Talagrand suppose $A$ increasing?

A first point to clarify is to understand why Talagrand's lemma is stated for an increasing set $A$. Upon inspection, it seems that his proof uses this hypothesis only once. To explain this point precisely, let us restate Talagrand's famous lemma 3.1, in order to compare it with theorem 46.1. Given an event $A$ and an index $k \in \{1, \dots, n\}$, the event $A_k$ is the event that $A$ occurs and $k$ is pivotal for the event $A$:

$$A_k = \big\{ \omega \in A : \omega^k \notin A \text{ or } \omega_k \notin A \big\}.$$

**Lemma 46.5** (Lemma 3.1 of [125])**.** *Consider a partition $I, J$ of $\{1, \dots, n\}$. Consider an increasing set $A$ and $S > 0$. Consider*

$$L = \left\{ k \in J : \left( \sum_{i \in I} \Big( \int_{A_k} Y_i \, d\mu \Big)^2 \right)^{1/2} \geq S\mu(A_k) \right\}. \tag{46.10}$$

*There exists a universal constant $C > 0$ such that*

$$\sum_{k \in L} \mu\big(A_k\big)^2 \leq C \exp\left( -\frac{S^2}{C} \right).$$

The hypothesis that the set $A$ is increasing is used at one point in the proof. In formula (2.2) of [125], Talagrand states the following inequality, which is reused later in formulas (3.3) and (3.6) during the proof of lemma 3.1: for an increasing set $A$, we have

$$\sum_{k \in L} \mu(A_k)^2 \leq 1. \tag{46.11}$$

This inequality comes from Parseval's identity. Indeed, for an increasing set $A$, it turns out that

$$\begin{aligned} \mu(A_k) &= \mu\big(\omega(k) = 1,\, \omega^k \in A,\, \omega_k \notin A\big) \\ &= \frac{1}{2}\mu\big(\omega^k \in A,\, \omega_k \notin A\big) = \frac{1}{2}E(\delta_k 1_A) = E(Y_k 1_A) \end{aligned}$$

is a Walsh coefficient of $1_A$. For a general event $A$, we would have

$$\begin{aligned} \mu(A_k) &= \mu\big(\omega(k) = 1,\, \omega^k \in A,\, \omega_k \notin A\big) + \mu\big(\omega(k) = 0,\, \omega^k \notin A,\, \omega_k \in A\big) \\ &= \frac{1}{2}\mu\big(\omega^k \in A,\, \omega_k \notin A\big) + \frac{1}{2}\mu\big(\omega^k \notin A,\, \omega_k \in A\big) = \frac{1}{2}E\big(|\delta_k 1_A|\big), \end{aligned}$$

and nothing can guarantee that inequality (46.11) would hold. Now, if we rewrite the definition of the set $L$ in (46.3) in the particular case $d = 2$, $r = 1$, $p = 1/2$, we get

$$L = \left\{ k \in J : \left( \sum_{i \in I} \big( E(Y_i Y_k 1_A) \big)^2 \right)^{1/2} \geq \frac{1}{2} s E\big(|\delta_k 1_A|\big) \right\}. \tag{46.12}$$

There is one subtle point that makes the definitions given in (46.10) and (46.12) differ for a general set $A$. In fact, for $i \in I$ and $k \in J$, we have, thanks to formula (46.1),

$$E(Y_i Y_k 1_A) \,=\, \sqrt{p(1-p)} E(Y_i \delta_k 1_A)\,. \tag{46.13}$$

When $A$ is increasing, we have $|\delta_k 1_A| = \delta_k 1_A$ and

$$2 \int_{A_k} Y_i \, d\mu \,=\, E(Y_i \delta_k 1_A)\,, \tag{46.14}$$

thus the definitions (46.10) and (46.12) coincide. However, this is not any more the case for a general set $A$, because $\delta_k 1_A$ might take negative values. This is why the adequate way to generalize Talagrand's lemma to non-monotone events consists in replacing (46.14) by (46.13) in the definition of the set $L$. Once this is done, Talagrand's proof goes through mutatis mutandis, because the inequality that needs to be used instead of (46.11) becomes

$$\sum_{k \in L} \big(E(Y_k 1_A)\big)^2 \,\leq\, 1\,,$$

which holds thanks to Parseval's inequality. In fact, we can even use the stronger inequality

$$\sum_{k \in L} \big(E(Y_k 1_A)\big)^2 \,\leq\, P(A)\,,$$

and get a slight improvement with respect to Talagrand's lemma. The goal of Talagrand was to obtain a refinement of the Harris-Kleitman inequality, which holds only for monotone sets, this is probably why he considered only monotone events throughout his work. Let us track further where the monotonicity hypothesis is used in Talagrand's paper. Talagrand's lemma 3.1 is the crucial estimate to prove theorem 2.4 of [125], which we restate next.

**Theorem 46.6** (Theorem 2.4 of [125])**.** *There exists a universal constant $c$ such that, for any increasing subsets $A, B$ of $\{0,1\}^n$, we have*

$$\sum_{\substack{1 \leq i,k \leq n \\ i \neq k}} \Big| \int_A Y_i Y_k \, d\mu \Big| \Big| \int_B Y_i Y_k \, d\mu \Big| \,\leq\, c \sum_{1 \leq k \leq n} \mu(A_k) \mu(B_k) \ln \frac{c}{\sum_{1 \leq k \leq n} \mu(A_k)\mu(B_k)}\,. \tag{46.15}$$

In fact, in Talagrand's paper, the left-hand side of (46.15) appears as

$$\sum_{\substack{1 \leq i,k \leq n \\ i \neq k}} \Big| \int_{A_k} Y_i \, d\mu \Big| \Big| \int_{B_k} Y_i \, d\mu \Big|\,. \tag{46.16}$$

We have used the identity (46.14) to rewrite (46.16) as the left-hand side of (46.15), in order to prepare the ground for a potential generalization of the inequality. Since the generalization of lemma 3.1 to a non-monotone event involves the quantities $E(|\delta_k 1_A|)$ in place of $\mu(A_k)$, it seems that an adequate generalization of inequality (46.15) would be

$$\sum_{\substack{1\le i,k\le n\\ i\neq k}} \Big|\int_A Y_i Y_k \, d\mu\Big| \Big|\int_B Y_i Y_k \, d\mu\Big| \le$$

$$c \sum_{1\le k\le n} E(|\delta_k 1_A|) E(|\delta_k 1_B|) \ln \frac{c}{\sum_{1\le k\le n} E(|\delta_k 1_A|) E(|\delta_k 1_B|)} \,. \quad (46.17)$$

However, we run into a serious problem when trying to prove the inequality (46.17) for non-monotone sets. Indeed, at the end of step 2 of the proof of his theorem 2.4, Talagrand makes again appeal to Parseval's inequality in order to control the sums

$$\sum_{k\in L} \big(E(|\delta_k 1_A|)\big)^2 \,, \qquad \sum_{k\in L} \big(E(|\delta_k 1_B|)\big)^2 \,. \quad (46.18)$$

When the events $A, B$ are both increasing, we have

$$E(|\delta_k 1_A|) = 2E(Y_k 1_A) \,, \quad E(|\delta_k 1_B|) = 2E(Y_k 1_B) \,,$$

hence Parseval's inequality can be used, but in the general case, an additional condition is needed. Keller and Kindler extend Talagrand's theorem 2.4 to non-monotone sets, but they assume that the above sums are both bounded by 1. Finally, Talagrand used his theorem 2.4 to derive the following inequality.

**Proposition 46.7** (Proposition 2.3 of [125])**.** *There exists a universal constant $c$ such that, for any increasing subset $A$ of $\{0,1\}^n$, we have*

$$\sum_{\substack{1\le i,k\le n\\ i\neq k}} \Big(\int_A Y_i Y_k \, d\mu\Big)^2 \le c \sum_{1\le k\le n} \mu(A_k)^2 \ln \frac{c}{\sum_{1\le k\le n} \mu(A_k)^2} \,.$$

The good news is that the proper reformulation of proposition 46.7 can be proved directly with the help of theorem 46.1, without making the detour through theorem 46.6. Therefore the additional requirement needed for theorem 46.6 when dealing with non-monotone sets is not necessary for the non-monotone version of proposition 46.7.

## 46.4 The proof of Talagrand

We reproduce here the argument of Talagrand, with the necessary modifications needed to prove theorems 46.1 and 46.2, which generalize lemma 3.1 and theorem 2.4 of [125]. To facilitate the comparison with Talagrand's paper [125], we will follow quite precisely Talagrand's presentation. In particular, we use the same step numbers as Talagrand, however we mainly adopt the notation of Keller and Kindler [80], with some minor changes. First, since we will sum over subsets $K$ of $J$ instead of using a single index $k \in J$, the letter $K$ is used here to denote a subset of $J$, whereas Talagrand used it to denote a generic constant. Second, we work with the biased measure of parameter $p$, so the integer $p$ in Talagrand's proof is denoted by $m$ here. Throughout the proof, we will remove the index $p$ from $P_p$ and we simply denote by $P$ the Bernoulli product measure with parameter $p$. Table 1 sums up the correspondence between Talagrand's notation and the one used here.

| Talagrand's paper [125] | corresponding notation |
|---|---|
| uniform measure $\mu$ | Bernoulli measure $P$ |
| $x_i$ | $\omega(i)$ |
| $2r_i$ | $Y_i$ |
| $p$ | $m$ |
| $k$ an index in $L$ | $K$ a subset of $L$ |
| a constant $K$ | a constant $c$ |

Table 1: Correspondence between Talagrand's notation and ours

The proof of Talagrand is rather difficult to understand. In the words of Benjamini, Kalai and Schramm, Talagrand's argument is beautiful, but mysterious. We start with the proof of theorem 46.1 and we divide the proof into four steps exactly as Talagrand did.

We consider a subset $A$ of $\{0,1\}^n$ and a partition $I, J$ of $\{1,\dots,n\}$. We fix an integer $d \geq 2$, another integer $r \geq 1$ and a real number $s > 0$.

**Step 1.** Let us fix a set $K$ belonging to $L$. There exists a family of numbers

$$\alpha_{T,K}\,,\quad T\subset I\,,\quad |T|=d-1\,,$$

such that

$$\sum_{T\subset I,\,|T|=d-1}\alpha_{T,K}^2\;=\;1\,,$$

and

$$\Bigg(\sum_{T\subset I,\,|T|=d-1}\big(E(Y_TY_K1_A)\big)^2\Bigg)^{1/2}\;=\;\sum_{T\subset I,\,|T|=d-1}\alpha_{T,K}E(Y_TY_K1_A)\,.$$

**Step 2.** We write $\Omega$ as the product $\Omega=\Omega_1\times\Omega_2$ where

$$\Omega_1=\{0,1\}^I\,,\qquad\Omega_2=\{0,1\}^J\,.$$

Accordingly, we decompose an element $\omega$ of $\Omega$ as $\omega=(\omega_1,\omega_2)$ where $\omega_1$ is in $\Omega_1$ and $\omega_2$ is in $\Omega_2$. We define next a random variable $F_K$ by setting

$$F_K(\omega)\;=\;\sum_{T\subset I,\,|T|=d-1}\alpha_{T,K}Y_T(\omega)\,.$$

Since $Y_T(\omega)$ depends only on the components of $\omega$ whose index is in $T \subset I$, we see that $F_K(\omega) = F_K(\omega_1)$ depends only on $\omega_1$. We can thus write the expectation $E(F_K Y_K 1_A)$ as a product and use Fubini's theorem to get

$$\begin{aligned}
E(F_K Y_K 1_A) &= \int_{\Omega_1}\int_{\Omega_2} F_K(\omega_1) Y_K(\omega_2) 1_A(\omega_1,\omega_2)\, dP_1(\omega_1)\, dP_2(\omega_2) \\
&= \int_{\Omega_1} F_K(\omega_1) \bigg( \int_{\Omega_2} Y_K(\omega_2) 1_A(\omega_1,\omega_2)\, dP_2(\omega_2) \bigg)\, dP_1(\omega_1)\,.
\end{aligned}$$

Noting that, for any $\omega_1$ in $\Omega_1$,

$$\int_{\Omega_2} Y_K(\omega_2) 1_A(\omega_1,\omega_2)\, dP_2(\omega_2) \,=\, E\big(Y_K 1_A \,|\, \omega_1\big)\,,$$

we see that

$$E(F_K Y_K 1_A) \,=\, E\Big(F_K(\omega_1) E(Y_K 1_A \,|\, \omega_1)\Big)\,.$$

**Step 3.** We define

$$\Omega_1(K,0) \,=\, \Big\{ \omega_1 \in \Omega_1 : \big|F_K(\omega_1)\big| \leq 2 \Big\}\,,$$

and for $m \geq 1$,

$$\Omega_1(K,m) \,=\, \Big\{ \omega_1 \in \Omega_1 : 2^m < \big|F_K(\omega_1)\big| \leq 2^{m+1} \Big\}\,.$$

We look at the expectation of $F_K Y_K 1_A$ restricted to $\Omega_1(K,m)$:

$$\begin{aligned}
E(K,m) &= \int_{\Omega_1(K,m)} F_K(\omega_1) E\big(Y_K 1_A \,|\, \omega_1\big)\, dP_1(\omega_1) \\
&= E\Big(1_{\Omega_1(K,m)}(\omega_1) F_K(\omega_1) E\big(Y_K 1_A \,|\, \omega_1\big)\Big)\,.
\end{aligned}$$

We apply the Cauchy-Schwarz inequality, and we obtain

$$\begin{aligned}
E(K,m)^2 &\leq E\Big(1^2_{\Omega_1(K,m)}\Big) E\bigg(\Big(1_{\Omega_1(K,m)}(\omega_1) F_K(\omega_1) E\big(Y_K 1_A \,|\, \omega_1\big)\Big)^2\bigg) \\
&\leq P\big(\Omega_1(K,m)\big) 2^{2m+2} E\bigg(\Big(E\big(Y_K 1_A \,|\, \omega_1\big)\Big)^2\bigg)\,. \qquad (46.19)
\end{aligned}$$

We use next corollary 46.4 to control $P\big(\Omega_1(K,m)\big)$: for any $m \geq 1$,

$$P\big(\Omega_1(K,m)\big) \,\leq\, P\big(\big|F_K(\omega_1)\big| > 2^m\big) \,\leq\, \exp\Big(d-1-\frac{d-1}{2B(p)e} 2^{2m/(d-1)}\Big)\,. \qquad (46.20)$$

Setting

$$c_H \,=\, \frac{d-1}{2B(p)e}\,,$$

and substituting (46.20) into (46.19), we get

$$E(K,m)^2 \;\leq\; \exp\big(d-1-c_H 2^{2m/(d-1)}\big)2^{2m+2}E\Big(\Big(E\big(Y_K 1_A \,|\, \omega_1\big)\Big)^2\Big)\,. \quad (46.21)$$

We shall next sum the inequality (46.21) over $K$ in $L$. Since the elements of $L$ are included in $J$, we have by Parseval's identity

$$\sum_{K\subset L}\Big(E\big(Y_K 1_A \,|\, \omega_1\big)\Big)^2 \;\leq\; \sum_{K\subset J}\Big(E\big(Y_K 1_A \,|\, \omega_1\big)\Big)^2 \;=\; P\big(A \,|\, \omega_1\big)\,. \qquad (46.22)$$

Taking the expectation of (46.22) with respect to $\omega_1$, we have

$$E\Big(\sum_{K\subset L}\Big(E\big(Y_K 1_A \,|\, \omega_1\big)\Big)^2\Big) \;\leq\; E\Big(P\big(A \,|\, \omega_1\big)\Big) \;\leq\; P(A)\,. \qquad (46.23)$$

Summing (46.21) over $K$ in $L$ and using (46.23), we get

$$\sum_{K\subset L} E(K,m)^2 \;\leq\; \exp\big(d-1-c_H 2^{2m/(d-1)}\big)2^{2m+2}P(A)\,. \qquad (46.24)$$

**Step 4.** We shall next re-sum over $m \geq 0$. Notice that, for $K$ an element of $L$, we have

$$s\big(p(1-p)\big)^{r/2}E\big(|\delta_K 1_A|\big) \;\leq\; E(F_K Y_K 1_A) \;=\; \sum_{m\geq 0} E(K,m)\,. \qquad (46.25)$$

We shall split the sum in two parts. Let $m_0$ be an integer, to be chosen later, and let us set

$$E_0(K) \;=\; \sum_{m\leq m_0} E(K,m)\,.$$

We can bound $E_0(K)$ as follows:

$$\begin{aligned}\big|E_0(K)\big| \;&\leq\; 2^{m_0+1}E\Big(\Big|E\big(Y_K 1_A \,|\, \omega_1\big)\Big|\Big)\\ &=\; 2^{m_0+1}\big(p(1-p)\big)^{r/2}E\Big(\big|E\big(\delta_K 1_A \,|\, \omega_1\big)\big|\Big)\\ &\leq\; 2^{m_0+1}\big(p(1-p)\big)^{r/2}E\big(\big|\delta_K 1_A\big|\big)\,. \end{aligned} \qquad (46.26)$$

To alleviate the equations, we define also

$$H(K) \;=\; \big(p(1-p)\big)^{r/2}E\big(|\delta_K 1_A|\big)\,.$$

From (46.25), we have then

$$sH(K) \;\leq\; E_0(K) \;+\; \sum_{m>m_0} E(K,m)\,,$$

whence, by convexity of the square function,

$$s^2 H(K)^2 \,\leq\, 2E_0(K)^2 \,+\, 2\Big(\sum_{m>m_0} E(K,m)\Big)^2$$
$$\leq\, 2E_0(K)^2 \,+\, 2\sum_{m>m_0} 2^m E(K,m)^2 \,. \quad (46.27)$$

We next sum (46.27) over $K$ in $L$ and we use the inequalities (46.24) and (46.26):

$$s^2 \sum_{K\subset L} H(K)^2 \,\leq\, 2^{2m_0+3} \sum_{K\subset L} H(K)^2$$
$$+ \sum_{m>m_0} 2^{3m+3} \exp\big(d-1-c_H 2^{2m/(d-1)}\big) P(A) \,. \quad (46.28)$$

There exists a positive constant $c < c_H$ such that, for any $m_0 \geq 1$,

$$2\sum_{m>m_0} 2^{3m+3} \exp\big(d-1-c_H 2^{2m/(d-1)}\big) \,\leq\, \frac{1}{c} \exp\big(-c 2^{2m_0/(d-1)}\big) \,. \quad (46.29)$$

We choose now for $m_0$ the largest integer such that $2^{2m_0+3} \leq s^2/2$. We have therefore

$$2^{2m_0} \,>\, s^2/2^6 \,, \quad (46.30)$$

and it follows from (46.28) and (46.29) that

$$s^2 \sum_{K\subset L} H(K)^2 \,\leq\, \frac{1}{c} \exp\big(-c 2^{2m_0/(d-1)}\big) P(A) \,.$$

Using the inequality (46.30) on $s$, we finally obtain the desired estimate (46.4).

We proceed next to the proof of theorem 46.2. As explained in the previous section, Talagrand proved first a decoupled version of theorem 46.2, but for such a decoupled version, we would need to suppose that the sets are increasing (as Talagrand did), or to add a condition to control the sums in (46.18) (as Keller and Kindler did). Fortunately, theorem 46.2 can be deduced directly from theorem 46.1 without imposing these restrictions. Our strategy consists in adapting Talagrand's proof of the decoupled version to get the non decoupled inequality. To alleviate the expressions, we set

$$c_{r,p} \,=\, \big(p(1-p)\big)^r \,,$$

and we define the sets $L_m$, $m \geq 0$, by

$$L_0 \,=\, \left\{ K\subset J : |K| = r, \sum_{\substack{T\subset I \\ |T|=d-1}} \big(E(Y_T Y_K 1_A)\big)^2 \,\leq\, 4c_{r,p}\Big(E\big(|\delta_K 1_A|\big)\Big)^2 \right\},$$

and for $m \geq 1$,

$$L_m = \Big\{ K \subset J : |K| = r,$$

$$2^{2m} c_{r,p} \Big( E\big(|\delta_K 1_A|\big)\Big)^2 < \sum_{\substack{T\subset I\\|T|=d-1}} \big(E(Y_T Y_K 1_A)\big)^2 \leq 2^{2m+2} c_{r,p} \Big( E\big(|\delta_K 1_A|\big)\Big)^2 \Big\}.$$

It follows from theorem 46.1 that, for $m \geq 1$,

$$\sum_{K\in L_m} \Big( E\big(|\delta_K 1_A|\big)\Big)^2 \leq \frac{1}{c_{r,p} c} \exp\big(-c 2^{2m/(d-1)}\big) P(A). \tag{46.31}$$

Inequality (46.31) and the definition of $L_m$ imply furthermore that, for $m \geq 1$,

$$\sum_{K\in L_m} \sum_{\substack{T\subset I\\|T|=d-1}} \big(E(Y_T Y_K 1_A)\big)^2 \leq 2^{2m+2} \frac{1}{c} \exp\big(-c 2^{2m/(d-1)}\big) P(A).$$

Thus there exists a constant $c' > 0$ such that, for any $m_0 \geq 1$,

$$\sum_{m>m_0} \sum_{K\in L_m} \sum_{\substack{T\subset I\\|T|=d-1}} \big(E(Y_T Y_K 1_A)\big)^2 \leq \frac{1}{c'} \exp\big(-c' 2^{2m_0/(d-1)}\big) P(A).$$

Using the definition of the sets $L_m$, we have also, for any $m_0 \geq 1$,

$$\sum_{m\leq m_0} \sum_{K\in L_m} \sum_{\substack{T\subset I\\|T|=d-1}} \big(E(Y_T Y_K 1_A)\big)^2 \leq 2^{2m_0+2} \sum_{m\leq m_0} \sum_{K\in L_m} c_{r,p} \Big( E\big(|\delta_K 1_A|\big)\Big)^2.$$

Using the quantity $H(J, r, 1_A)$ introduced in (46.5), and summing the two previous inequalities, we obtain

$$\sum_{K\in L} \sum_{\substack{T\subset I\\|T|=d-1}} \big(E(Y_T Y_K 1_A)\big)^2 \leq$$

$$2^{2m_0+2} H(J, r, 1_A) + \frac{1}{c'} \exp\big(-c' 2^{2m_0/(d-1)}\big) P(A). \tag{46.32}$$

If $P(A) \leq H(J, r, 1_A)$, we take $m_0 = 1$ and we obtain that, for some constant $c'' > 0$,

$$\sum_{K\in L} \sum_{\substack{T\subset I\\|T|=d-1}} \big(E(Y_T Y_K 1_A)\big)^2 \leq c'' H(J, r, 1_A). \tag{46.33}$$

If $P(A) > H(J, r, 1_A)$, we take for $m_0$ the smallest positive integer such that

$$\exp\big(-c' 2^{2m_0/(d-1)}\big) \leq \frac{H(J, r, 1_A)}{P(A)}.$$

The integer $m_0$ satisfies

$$2^{2m_0} \geq \Big(\frac{1}{c'}\ln\Big(\frac{P(A)}{H(J,r,1_A)}\Big)\Big)^{d-1} > 2^{2m_0-2}\,, \tag{46.34}$$

and we obtain from (46.32) and (46.34) that

$$\sum_{K\in L}\sum_{\substack{T\subset I\\|T|=d-1}}\big(E(Y_TY_K1_A)\big)^2 \leq$$
$$H(J,r,1_A)\bigg(16\Big(\frac{1}{c'}\ln\Big(\frac{P(A)}{H(J,r,1_A)}\Big)\Big)^{d-1}+\frac{1}{c'}\bigg). \tag{46.35}$$

Putting together inequalities (46.33) and (46.35), we conclude that there exists a constant $\widetilde{c}>0$ such that

$$\sum_{K\in L}\sum_{\substack{T\subset I\\|T|=d-1}}\big(E(Y_TY_K1_A)\big)^2 \leq \widetilde{c}\,H(J,r,1_A)\Big(1\vee\ln\Big(\frac{P(A)}{H(J,r,1_A)}\Big)\Big)^{d-1}.$$

So far we have obtained an inequality having the form of (46.6) but with an unknown constant $\widetilde{c}$. It is difficult to extract an explicit value from the previous computations, because it would depend on the sum of some complicated series. A definitive advantage of the proof of Keller and Kindler, to be presented shortly, is that it provides an explicit value for the constant.

## 46.5 The proof of Keller and Kindler

We present here the proof of Keller and Kindler. We use essentially the same notation as theirs, but we prefer to stick to the probabilistic framework rather than using the Hilbert space formalism. Let $f$ be a function from $\Omega$ to $\mathbb{R}$. Let $I,J$ be a partition of $\{1,\dots,n\}$. Let $d\geq 2$ and $r\geq 1$ be two integers, and let us fix a subset $K$ of $J$ of cardinality $r$. We define a function $f'_K$ as

$$f'_K = \sum_{T\subset I,|T|=d-1} E(Y_TY_Kf)\,Y_T\,.$$

The function $f'_K$ depends only on the components of $\omega$ whose index belongs to $I$. The starting point of the arguments of Keller and Kindler is the following formula:

$$E(f'_KY_Kf) = \sum_{T\subset I,|T|=d-1}\big(E(Y_TY_Kf)\big)^2\,. \tag{46.36}$$

By Parseval's identity, we have

$$E({f'_K}^2) = \sum_{T\subset I,|T|=d-1}\big(E(Y_TY_Kf)\big)^2\,. \tag{46.37}$$

We normalize $f'_K$ into the function $f_K$ defined by

$$f_K \;=\; \frac{f'_K}{\sqrt{E({f'_K}^2)}}\,.$$

It follows from (46.36) and (46.37) that

$$E\big(f_K Y_K f\big) \;=\; \sqrt{\sum_{T\subset I, |T|=d-1} \big(E(Y_T Y_K f)\big)^2}\,. \tag{46.38}$$

Since $f_K$ depends only on the components of $\omega$ whose index belongs to $I$, which is disjoint from $K$, it follows from the identity (46.2) that

$$E\big(f_K Y_K f\big) \;=\; \big(p(1-p)\big)^{r/2} E\big(f_K \delta_K f\big)\,. \tag{46.39}$$

Furthermore, taking the conditional expectation with respect to $\sigma(I)$, we obtain

$$E\big(f_K \delta_K f\big) \;=\; E\Big(E\big(f_K \delta_K f \,|\, \sigma(I)\big)\Big) \;=\; E\Big(f_K E\big(\delta_K f \,|\, \sigma(I)\big)\Big)\,,$$

whence

$$\big|E\big(f_K \delta_K f\big)\big| \;\leq\; E\Big(\big|f_K\big|\big|E\big(\delta_K f \,|\, \sigma(I)\big)\big|\Big)\,. \tag{46.40}$$

We use Fubini's theorem to rewrite the right-hand quantity in (46.40) as

$$E\Big(\big|f_K\big|\big|E\big(\delta_K f \,|\, \sigma(I)\big)\big|\Big) \;=\; \int_0^\infty E\Big(1_{\{|f_K|>t\}}\big|E\big(\delta_K f \,|\, \sigma(I)\big)\big|\Big)\,dt\,. \tag{46.41}$$

Let $t_0 \geq 0$ be a parameter to be chosen later. We split the integral over $[0,+\infty[$ into two integrals over $[0,t_0]$ and $]t_0,+\infty[$, we take the square in (46.40), (46.41), and we obtain

$$\begin{aligned} \Big(E\big(f_K \delta_K f\big)\Big)^2 \;&\leq\; \bigg(\int_0^\infty E\Big(1_{\{|f_K|>t\}}\big|E\big(\delta_K f \,|\, \sigma(I)\big)\big|\Big)\,dt\bigg)^2 \\ &\leq\; \bigg(\int_0^{t_0}\cdots + \int_{t_0}^{+\infty}\cdots\bigg)^2 \;\leq\; 2\bigg(\int_0^{t_0}\cdots\bigg)^2 + 2\bigg(\int_{t_0}^{+\infty}\cdots\bigg)^2\,. \end{aligned} \tag{46.42}$$

We bound separately each integral as follows. For the integral over $[0,t_0]$, we have

$$\int_0^{t_0}\cdots \;\leq\; t_0\, E\Big(\big|E\big(\delta_K f \,|\, \sigma(I)\big)\big|\Big) \;\leq\; t_0\, E\Big(E\big(\big|\delta_K f\big| \,|\, \sigma(I)\big)\Big) \;=\; t_0\, E\big(\big|\delta_K f\big|\big)\,. \tag{46.43}$$

For the integral over $]t_0,+\infty[$, we bound the integrand with the Cauchy-Schwarz inequality:

$$E\big(1_{\{|f_K|>t\}}\big|E\big(\delta_K f \,|\, \sigma(I)\big)\big|\big) \;\leq\; \sqrt{P\big(|f_K|>t\big) E\big(\big|E\big(\delta_K f \,|\, \sigma(I)\big)\big|^2\big)}\,. \tag{46.44}$$

Substituting (46.44) into the integral, we obtain

$$\int_{t_0}^{+\infty} \cdots \leq \Big( \int_{t_0}^{+\infty} \sqrt{P\big(|f_K| > t\big)}\, dt \Big) \sqrt{E\big(\big|E\big(\delta_K f \,|\, \sigma(I)\big)\big|^2\big)} \,. \qquad (46.45)$$

Substituting (46.43) and (46.45) into (46.42), we have

$$\Big(E\big(f_K \delta_K f\big)\Big)^2 \leq 2t_0^2 \Big(E\big(\big|\delta_K f\big|\big)\Big)^2 + \\ 2\bigg( \int_{t_0}^{+\infty} \sqrt{P\big(|f_K| > t\big)}\, dt \bigg)^2 E\big(\big|E\big(\delta_K f \,|\, \sigma(I)\big)\big|^2\big) \,. \qquad (46.46)$$

In the next step, we bound $\big(E\big(\big|\delta_K f\big|\big)\big)^2$ by $E\big((\delta_K f)^2\big)$, and we use corollary 46.4 to control the integral:

$$\int_{t_0}^{+\infty} \sqrt{P\big(|f_K| > t\big)}\, dt \leq \int_{t_0}^{+\infty} \exp\Big( \frac{d-1}{2} - \frac{d-1}{4B(p)e} t^{2/(d-1)} \Big)\, dt \,. \qquad (46.47)$$

The next lemma gives a useful bound on the last integral.

**Lemma 46.8.** *Let $\alpha > 0$ and let $m$ be a positive integer. For any $T > 0$, we have*

$$\int_T^{+\infty} \exp\big(-\alpha t^{2/m}\big)\, dt \leq \frac{2m!}{\alpha^{m/2}} \exp\big(-\alpha T^{2/m}/2\big) \,. \qquad (46.48)$$

*Proof.* We first do the change of variable $u = \alpha t^{2/m}$. We have

$$\frac{du}{u} = \frac{2}{m}\frac{dt}{t}\,, \qquad t = \Big(\frac{u}{\alpha}\Big)^{m/2},$$

whence

$$\int_T^{+\infty} \exp\big(-\alpha t^{2/m}\big)\, dt = \int_{\alpha T^{2/m}}^{+\infty} \frac{m}{2}\Big(\frac{u}{\alpha}\Big)^{m/2}\Big(\frac{1}{u}\Big) \exp\big(-u\big)\, du \,.$$

We use next the inequality

$$\forall u > 0 \qquad u^{(m-1)/2} \leq (m-1)!\, \exp(u/2)\,,$$

and we make the change of variables $v = u - \alpha T^{2/m}$ to get

$$\begin{aligned} \int_T^{+\infty} \exp\big(-\alpha t^{2/m}\big)\, dt &\leq \frac{m!}{2\alpha^{m/2}} \int_{\alpha T^{2/m}}^{+\infty} \frac{1}{\sqrt{u}} \exp\big(-u/2\big)\, du \\ &\leq \frac{m!}{2\alpha^{m/2}} \exp\big(-\alpha T^{2/m}/2\big) \int_0^{+\infty} \frac{1}{\sqrt{v + \alpha T^{2/m}}} \exp\big(-v/2\big)\, dv \\ &\leq \frac{m!}{2\alpha^{m/2}} \exp\big(-\alpha T^{2/m}/2\big) \int_0^{+\infty} \frac{1}{\sqrt{v}} \exp\big(-v/2\big)\, dv \,. \qquad (46.49) \end{aligned}$$

We finally bound the last integral as follows:

$$\int_0^{+\infty} \frac{1}{\sqrt{v}} \exp\big(-v/2\big)\, dv \,\leq\, \int_0^1 \frac{1}{\sqrt{v}}\, dv + \int_1^{+\infty} \exp\big(-v/2\big)\, dv \,\leq\, 4\,. \quad (46.50)$$

Substituting (46.50) into (46.49), we obtain (46.48). ☐

We apply lemma 46.8 with

$$\alpha = \frac{d-1}{4B(p)e}\,, \quad m = d-1\,, \quad T = t_0\,,$$

and we obtain from (46.47) that

$$\int_{t_0}^{+\infty} \sqrt{P\big(|f_K| > t\big)}\, dt \,\leq\, c_1 \exp\Big(-c_2 t_0^{2/(d-1)}\Big)\,. \quad (46.51)$$

where $c_1, c_2$ are the constants defined by

$$c_1 \,=\, 2(d-1)!e^{d-1}\Big(\frac{4B(p)}{d-1}\Big)^{(d-1)/2}\,, \qquad c_2 \,=\, \frac{d-1}{8B(p)e}\,.$$

Substituting (46.51) into (46.46), we conclude that

$$\Big(E\big(f_K \delta_K f\big)\Big)^2 \,\leq\, 2t_0^2\Big(E\big(|\delta_K f|\big)\Big)^2 + \\ 2c_1^2 \exp\big(-2c_2 t_0^{2/(d-1)}\big)\, E\big(\big|E\big(\delta_K f \,|\, \sigma(I)\big)\big|^2\big)\,.$$

Coming back to equalities (46.38) and (46.39), we obtain

$$\sum_{T \subset I, |T| = d-1} \big(E(Y_T Y_K f)\big)^2 \,\leq\, 2\big(p(1-p)\big)^r \times \\ \Big(t_0^2\Big(E\big(|\delta_K f|\big)\Big)^2 \,+\, c_1^2 \exp\big(-2c_2 t_0^{2/(d-1)}\big) E\big(\big|E\big(\delta_K f \,|\, \sigma(I)\big)\big|^2\big)\Big)\,. \quad (46.52)$$

By (46.2) and Parseval's identity, we have

$$\sum_{\substack{K \subset J \\ |K| = r}} \big(p(1-p)\big)^r E\big(\big|E\big(\delta_K f \,|\, \sigma(I)\big)\big|^2\big) \,=\, \sum_{\substack{K \subset J \\ |K| = r}} E\Big(\big(E\big(Y_K f \,|\, \sigma(I)\big)\big)^2\Big) \\ \leq\, E(f^2) \,=\, E(f)\,. \quad (46.53)$$

We next sum (46.52) over $K$, we use (46.53) and the definition (46.5) of $H(J, r, f)$ to conclude that

$$\sum_{\substack{K \subset J \\ |K| = r}} \sum_{\substack{T \subset I \\ |T| = d-1}} \big(E(Y_T Y_K f)\big)^2 \,\leq \\ 2t_0^2 H(J, r, f) \,+\, 2c_1^2 \exp\big(-2c_2 t_0^{2/(d-1)}\big) E(f)\,. \quad (46.54)$$

If $E(f) \le e\,H(J,r,f)$, then we take $t_0 = 0$ (or rather send $t_0$ to 0) in (46.54) and we get that

$$\sum_{\substack{K\subset J\\|K|=r}}\ \sum_{\substack{T\subset I\\|T|=d-1}} \big(E(Y_TY_Kf)\big)^2 \;\le\; 2c_1^2\, e\, H(J,r,f)\,. \tag{46.55}$$

If $E(f) > e\,H(J,r,f)$, then we take $t_0$ such that

$$\exp\big(-2c_2t_0^{2/(d-1)}\big) \;=\; \frac{H(J,r,f)}{E(f)}\,,$$

that is,

$$t_0^2 \;=\; \left(\frac{1}{2c_2}\ln\left(\frac{E(f)}{H(J,r,f)}\right)\right)^{d-1}.$$

Plugging this value of $t_0$ in (46.54), we have

$$\begin{aligned}
\sum_{\substack{K\subset J\\|K|=r}}\ \sum_{\substack{T\subset I\\|T|=d-1}} &\big(E(Y_TY_Kf)\big)^2\\
&\le\; 2H(J,r,f)\left(\left(\frac{1}{2c_2}\ln\left(\frac{E(f)}{H(J,r,f)}\right)\right)^{d-1} +\, c_1^2\right)\\
&\le\; 2H(J,r,f)\Big(\frac{1}{(2c_2)^{d-1}} \,+\, c_1^2\Big)\left(\ln\left(\frac{E(f)}{H(J,r,f)}\right)\right)^{d-1}\\
&\le\; 16^d(d!)^2H(J,r,f)\Big(\frac{B(p)e}{d-1}\Big)^{d-1}\left(\ln\left(\frac{E(f)}{H(J,r,f)}\right)\right)^{d-1}. \qquad (46.56)
\end{aligned}$$

Putting together inequalities (46.55) and (46.56), we conclude that

$$\sum_{\substack{K\subset J\\|K|=r}}\ \sum_{\substack{T\subset I\\|T|=d-1}} \big(E(Y_TY_Kf)\big)^2 \;\le\;$$

$$16^d(d!)^2H(J,r,f)\Big(\frac{B(p)e}{d-1}\Big)^{d-1}\left(1\vee\ln\left(\frac{E(f)}{H(J,r,f)}\right)\right)^{d-1}.$$

This completes the proof of theorem 46.2.

## Part VIII

# Monotone subsets of the hypercube ⊛

The control of the pivots of arbitrary orders seem to be the missing ingredient to conclude the proof that $\theta(p_c, \mathbb{Z}^d) = 0$. Throughout sections 31, 33 and Part V, we tried to obtain a control on these pivots in the specific case of the disconnection event, but the various attempts presented there failed. Actually, the repeated failures to achieve an adequate control seem to be due to the fact that we were completely entangled in the percolation configuration, to the point that we were paralyzed by the specificities of the disconnection event. For instance, there exists a specific algorithm to discover the pivots, and most of our efforts have been devoted to take advantage of this algorithm to obtain a useful probabilistic control. This worked to some extent, but the controls were too weak to conclude. We try here to apply the Choquet principle.[6] In short, we consider the problem at a more general level, and we try to do the job for a generic monotone event.

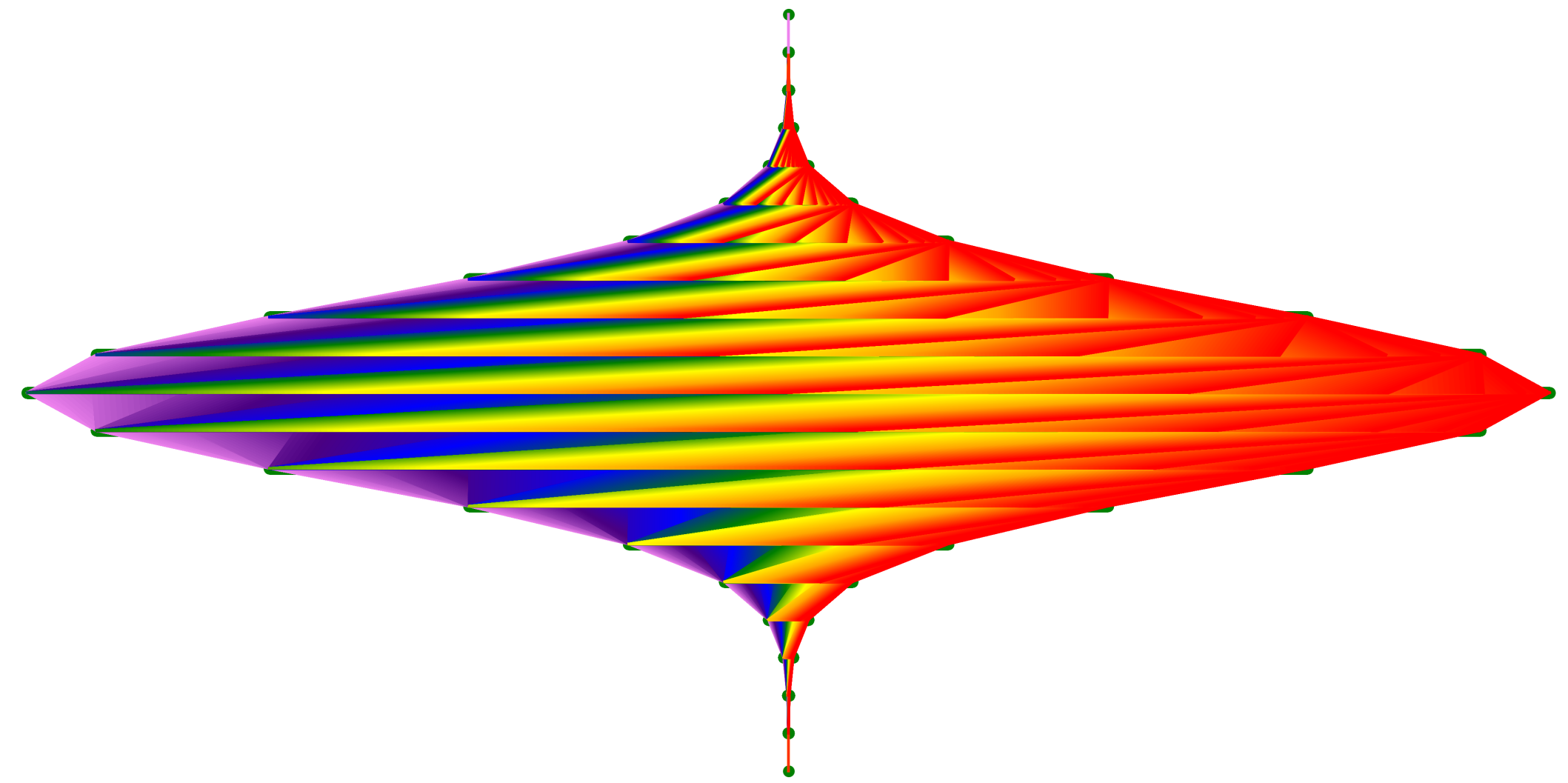

Figure 103: The hypercube graph for $n = 20$

[6] In his EMS interview [47], Michel Talagrand recounts: "One of the first things I asked professor Choquet was to please advise me on how to do research. And he gave me the most useful advice I ever got, which I called the Choquet principle, which was well known before him, but I received it from him.

> *When you attack a problem, always choose the proper setting for the problem. And the proper setting is that it should have no extra structure, it should be the minimal structure where the problem makes sense.*

The obvious reason for that is that you don't get distracted by things which are not relevant for the problem. I have applied that advice several times, with great success I would say."

## 47 The quest for the missing inequalities

Throughout this part, we have to struggle with several terminological dilemmas because we stand at the crossroads of three very rich mathematical subjects: the theory of Boolean functions, extremal set theory, and percolation. Here is a selection of three excellent books devoted to these three subjects:
• Analysis of Boolean functions by Ryan O'Donnell [110]
• Percolation by Geoffrey Grimmett [63],
• Extremal problems for finite sets by Peter Frankl and Norihide Tokushige [56].
The literature on influence and noise sensitivity is essentially written in terms of Boolean functions. Extremal set theory is formulated in terms of subsets of $\{0,1\}^n$. In percolation, the set $\{0,1\}^n$ is considered as a probability space and its elements are called configurations. Since this whole text is focused on the conjecture $\theta(p_c,\mathbb{Z}^d)=0$, we will try to adapt to the percolation terminology the definitions and concepts coming from the two other topics.

So we will mainly focus on increasing subsets of $\{0,1\}^n$. We already know two remarkable inequalities involving the pivotal set $\mathcal{P}^+(A)$. From proposition 39.5, we have

$$E\big(|\mathcal{P}^+(A)|\big)\,\leq\,\sqrt{\frac{nP_p(A)}{p(1-p)}}\,. \tag{47.1}$$

From corollary 40.3, we have furthermore

$$P_p(A)\,\leq\,\exp\Big(-\frac{1}{2n}\Big((1-p)E\big(|\mathcal{P}^+(A)|\,\big|\,A\big)\Big)^2\Big)\,. \tag{47.2}$$

Our long term goal in this part is to obtain similar inequalities for the set of the positive higher pivots defined in section 49. Our initial attempt consisted in generalizing the Margulis-Russo formula, but it ended in a resounding failure. The framework was the same as here (monotone subsets of $\{0,1\}^n$), but the question was not formulated as an extremal problem, and the technique was not adequate (differentiation with respect to the parameter $p$). Afterwards, we tried to obtain a control for the specific case of a disconnection event in percolation, and again all our attempts have been doomed to failure. So we try here a new angle of attack. In section 50, we reformulate the inequality (47.1) as an extremal problem. At several points, we will use some input from the theory of Boolean functions. In section 48, we explain the correspondence between boolean functions, subsets of $\{0,1\}^n$, antichains, and we describe a representation for the hypercube and its subsets that we will use for the numerical computations and the various examples. In section 49, we define the higher positive pivots. In sections 51 and 52, we present the maximizers for the second pivots for small values of $n$. In section 53, we formulate a conjecture on these maximizers. We discuss this conjecture, some variants and related results in the subsequent sections 54, 58, 59, 60, 61, and we present further numerical results in sections 55, 56, and 57. We conclude this part with the set of Talagrand in section 62, which helps to clarify the situation.

# 48 Increasing sets and boolean functions

Boolean functions play a prominent role in computer science and there exists a huge literature devoted to the study of boolean functions defined on $\{0,1\}^n$. We shall focus here on monotone boolean functions, and we recommend the excellent review text written by Korshunov [89]. A boolean function $f$ from $\{0,1\}^n$ to $\{0,1\}$ is monotone if

$$\forall \omega, \eta \in \{0,1\}^n \qquad \omega \leq \eta \quad \Rightarrow \quad f(\omega) \,\leq\, f(\eta)\,.$$

There is a natural one to one correspondence between monotone subsets of $\{0,1\}^n$ and monotone boolean functions: to each monotone subset $A$ of $\{0,1\}^n$, we associate the boolean function $f = 1_A$. In the other direction, to each monotone boolean function $f$ defined on $\{0,1\}^n$, we associate the monotone subset $A = f^{-1}(\{1\})$.

Let $A$ be an increasing subset of $\{0,1\}^n$. A minimal configuration of $A$ is a configuration $\sigma$ which belongs to $A$ and such that there is no configuration strictly smaller than $\sigma$ belonging to $A$, or equivalently, such that

$$\forall \eta \in A \qquad \eta \leq \sigma \quad \Rightarrow \quad \eta = \sigma\,.$$

We denote by Minimal($A$) the set of the minimal configurations of $A$. We give next the corresponding definition for a boolean function.

**Definition 48.1.** *Let $f$ be an increasing boolean function. A minimal implicant of $f$ is a minimal configuration $\sigma$ of the set $A = f^{-1}(\{1\})$. We denote by Minimal($f$) the set of the minimal implicants of $f$.*

The point is that the set Minimal($A$) (respectively Minimal($f$)) completely determines the increasing set $A$ (respectively the increasing boolean function $f$). Moreover, the sets Minimal($A$) and Minimal($f$) are antichains.

**Definition 48.2.** *An antichain in $\{0,1\}^n$ is a collection $\mathcal{A}$ of configurations of $\{0,1\}^n$ which are pairwise not comparable, i.e.,*

$$\forall \sigma, \eta \in \mathcal{A} \qquad \sigma \neq \eta \quad \Longrightarrow \quad \sigma \not\leq \eta\,, \quad \eta \not\leq \sigma\,.$$

The usual definition of an antichain in $\{0,1\}^n$ is given for subsets of $\{\,1,\dots,n\,\}$ rather than for configurations. Yet there is a one to one correspondence between configurations of $\{0,1\}^n$ and subsets of $\{\,1,\dots,n\,\}$: to each configuration $\omega$, we associate its support $\operatorname{supp}\omega$, defined by

$$\operatorname{supp}\omega \;=\; \big\{\, i \in \{\,1,\dots,n\,\} : \omega(i) = 1 \,\big\}\,.$$

Thanks to this correspondence, the above definition 48.2 of an antichain is equivalent to the usual one.

Finally, we can also give similar definitions for decreasing sets and functions. If $A$ is a decreasing set, we denote by Maximal($A$) the maximal elements of $A$; if $f$ is a decreasing boolean function, we define Maximal($f$) as Maximal$\big(f^{-1}(\{1\})\big)$. This is all the more interesting because, if $A$ is an increasing set, then its complement $A^c$ is decreasing, therefore $A$ is completely determined by Maximal($A^c$).

## 48.1 The Dedekind numbers

Of course, the structure of a monotone boolean function is much simpler than that of a general boolean function, yet there is already a lot of interesting and challenging mathematical questions dealing with monotone boolean functions. Some of them might even be more difficult than for arbitrary boolean functions. A first fundamental question is to compute the number of monotone boolean functions acting on a fixed number of variables. Dedekind stated this problem and computed the first five numbers in the nineteenth century. In his honour, the number of monotone boolean functions in $n$ variables is called the $n$-th Dedekind number and is denoted by $D(n)$. To simplify their classification, we introduce an equivalence relation on monotone boolean functions, as follows. Two monotone boolean functions $f$ and $g$ are said to be equivalent if there exists a permutation $\phi$ of $\{ 1, \dots, n \}$ such that

$$\forall x = (x_1, \dots, x_n) \in \{0,1\}^n \qquad g(x) = f\big(x_{\phi(1)}, \dots, x_{\phi(n)}\big)\,.$$

The number of equivalence classes of monotone boolean functions in $n$ variables, also called the number of unequivalent monotone boolean functions, is called the $n$-th Dedekind number of the second kind and is denoted by $R(n)$. Currently, the first nine Dedekind numbers are known, they are listed in the On-Line Encyclopedia of Integer Sequences [111, 112]. The values for $n \leq 8$ are given in table 2.

| $n$ | $D(n)$ | $R(n)$ |
|---|---|---|
| 0 | 2 | 2 |
| 1 | 3 | 3 |
| 2 | 6 | 5 |
| 3 | 20 | 10 |
| 4 | 168 | 30 |
| 5 | 7581 | 210 |
| 6 | 7828354 | 16353 |
| 7 | 2414682040998 | 490013148 |
| 8 | 56130437228687557907788 | 1392195548889993358 |

Table 2: The values of $D(n)$ and $R(n)$ for $0 \leq n \leq 8$

Only the first five values for $n \leq 4$ can reasonably be computed by hand. Each value for $n = 5, 6, 7, 8, 9$ was the achievement of a specific mathematical work assisted by a computer. Each time, the combination of some mathematical ingeniosity and the increased computer power allowed to compute the next number. The values for $n = 9$ were computed in 2023:

$$D(9) = 286386577668298411128469151667598498812366\,,$$
$$R(9) = 789204635842035040527740846300252680\,.$$

We refer to [111, 112] for the list of the works and the mathematicians who computed these numbers. Besides the exact computation of the Dedekind numbers, mathematicians have striven to obtain upper and lower bounds for the Dedekind numbers, as well as asymptotic expansions. Using the one to one correspondence between increasing boolean functions and antichains, we can interpret the $n$-th Dedekind number $D(n)$ as the number of antichains in $\{0,1\}^n$. Sperner's theorem states that the maximal cardinality of an antichain in $\{0,1\}^n$ is

$$\binom{n}{\lfloor n/2 \rfloor}\,.$$

This classical result is considered as the starting point of extremal set theory. The shortest known proof is due to Lubell [98] and it goes as follows. Let $A$ be an antichain on $\{0,1\}^n$. A chain on $\{0,1\}^n$ is a collection of configurations which is totally ordered. There exist exactly $n!$ maximal chains on $\{0,1\}^n$. For any $k \in \{\,0,\dots,n\,\}$ and any configuration $\omega \in \{0,1\}^n$ such that $|\operatorname{supp}\omega| = k$, there exist exactly $k!(n-k)!$ maximal chains containing the configuration $\omega$. Since $A$ is an antichain, a maximal chain can intersect $A$ at most once, therefore

$$\sum_{0\leq k\leq n}\ \sum_{\substack{\omega\in A\\ |\operatorname{supp}\omega|=k}} k!(n-k)!\ \leq\ n!\,.$$

Now the quantity $k!(n-k)!$ is minimal for $k = \lfloor n/2 \rfloor$, and this readily implies the upper bound of Sperner's theorem. For the lower bound, we simply observe that the collection of the configurations having exactly $\lfloor n/2 \rfloor$ components equal to 1 is an antichain having the desired cardinality.

Using Sperner's theorem, we obtain the following bounds on $D(n)$:

$$2^{\binom{n}{\lfloor n/2 \rfloor}}\ \leq\ D(n)\ \leq\ 2^{n\binom{n}{\lfloor n/2 \rfloor}}\,. \tag{48.1}$$

The lower bound corresponds to the number of antichains whose all elements have exactly $\lfloor n/2 \rfloor$ components equal to 1. The upper bound is a crude upper bound on the total number of antichains. We see from the inequalities (48.1) that the sequence $D(n)$ grows superexponentially fast. In fact, Kleitman [85] showed that

$$\log_2 D(n) \sim \binom{n}{\lfloor n/2 \rfloor}\,.$$

The reason for this result is that the vast majority of the antichains in $\{0,1\}^n$ correspond to configurations having between $\lfloor n/2 \rfloor - 3$ and $\lfloor n/2 \rfloor + 3$ components equal to 1. In a mathematical tour de force, Korshunov [89] succeeded in obtaining an equivalent of $D(n)$ as $n$ goes to $\infty$; unfortunately, the paper [89] does not seem to have been translated in english. Fortunately, the sequence $R(n)$ is more reasonable. Typically, we gain a factor of order $n!$ by considering unequivalent monotone functions.

## 48.2 Exhaustive exhausting searchs

In the course of the computation of the first Dedekind numbers, the algorithms that were used produced the exhaustive list of all the monotone boolean functions. This constitutes a precious database if one wishes to test some conjectures on monotone boolean functions. Unfortunately, the corresponding files are not published (even in electronic form) along with the mathematical papers. So we looked on the web for programs or scientific sites which would compute or contain the list of the monotone boolean functions in a small number of variables. For $n \leq 6$, Olivier Sobrie had computed the lists during his PhD thesis. These lists are available on his website `https://olivier.sobrie.be` and the Python code to generate them is available on his github site `https://github.com/oso/imbfs`. This material is used and explained in the work [52]. In order to get the files for $n = 7$, we contacted the authors of the work [124], where $D(7)$ and $R(7)$ were computed (in addition, the paper [124] provides a nicely written review of the history of the Dedekind numbers). Unfortunately, their code and the files had been lost several years ago due to some software issues. However, Tamon Stephen was a great help. He reported that their calculations had been redone at least twice subsequently. Indeed, the value of $R(7)$ is computed again in the work [20], with an apparently much faster implementation involving a hash table (unfortunately not publicly available). In his work on neural networks [24], Romain Cazé computes all the monotone boolean functions for $n \leq 7$, with the help of a python code which is available on github (`https://github.com/rcaze/PlosCB2013`). The Python code of Romain Cazé has been rewritten in C++ by Mikaël Monet in the course of his work on probabilistic databases [107], the corresponding code being available on gitlab (`https://gitlab.com/Gruyere/sperner-families-generator`). Both codes of Romain Cazé and of Mikaël Monet are fully operational. Yet the code of Mikaël Monet is faster and more efficient (thanks to its careful optimization and the power of C++), so in the end we relied on this code to recreate the full list of all unequivalent monotone boolean functions in 7 variables. Once we have this file at our disposal, we can find by a brute computation the monotone boolean functions which realize the maximum of a specific functional. Of course we need only to consider unequivalent monotone boolean functions. For $n = 6$, their number is $R(6) = 16353$, and an exhaustive search over the set of monotone boolean functions in 6 variables can be easily performed on a modern laptop in just a few seconds. For $n = 7$, their number is $R(7) = 490013148$, and things become considerably harder. Still, a non-optimized file containing all these functions has an approximate size of 22 Gigabytes, and an exhaustive search over the set of monotone boolean functions in 7 variables takes between 4 hours and 1 day, depending on its complexity. The set of functions in 8 variables is far too large to be examined exhaustively.

## 48.3 Representation of the hypercube graph

Most of the action will take place in the hypercube $\{0,1\}^n$, endowed with the graph structure in which two vertices are connected by an edge whenever they differ in exactly one coordinate. This graph is called the hypercube graph, it has $2^n$ vertices and $n2^{n-1}$ edges. There exist several ways of drawing it. Here we will study some increasing subsets of this graph, and we will be interested in the pivotal components of these sets (to be defined in section 49). To visualize these sets for small values of $n$, we choose to use a layer representation of the hypercube graph. In this representation, the vertices are organized into $n$ layers, the layer number $k$ containing the vertices having exactly $k$ components equal to 1. The edges connect vertices between two consecutive layers. In the following pictures, the vertices are arranged according to the lexicographic order on each layer. We also use a code color to distinguish the edges, which follows two rules. The downwards edges exiting a vertex all share the same color. This color is choosing according to the position of the vertex in its row, within a rainbow colormap. So the largest layer is divided in seven parts, corresponding to the seven rainbow colors violet, indigo, blue, green, yellow, orange, red. For this largest layer, the color associated to the edges of the leftmost vertex is violet, it is red for the rightmost. The other layers use the same colormap, starting with violet for the leftmost vertex. To cope with the rapidly increasing number of vertices, we use a color gradient in each of the seventh parts of each layer. It is fundamental to gain some insight into the geometry of the hypercube graph to understand the features of the monotone sets which are relevant to our problem. The hypercube graphs for $n = 5$ and $n = 20$ are depicted in figures 104 and 103. The hypercube graphs for $1 \leq n \leq 8$ are presented in figure 105 and for $9 \leq n \leq 12$, $n = 20$ in figure 106. The labels of the vertices are included only for $n \leq 5$, for larger values of $n$ they would no longer be readable. Furthermore, we will need to visualize increasing subsets of the hypercube graphs. Let $A$ be such a subset. We will mark by a green dot the vertices belonging to $A$ and by a blue cross the vertices which do not belong to $A$. Moreover, we color in red the edges in the edge boundary of $A$, that is the edges between a vertex of $A$ and a vertex of its complement. In figure 107, we present an example of a subset of $\{0, 1\}^5$. Depending on the set $A$, we might present only some interesting layers of the hypercube graphs. Later on, we will include other relevant edges with different colors.

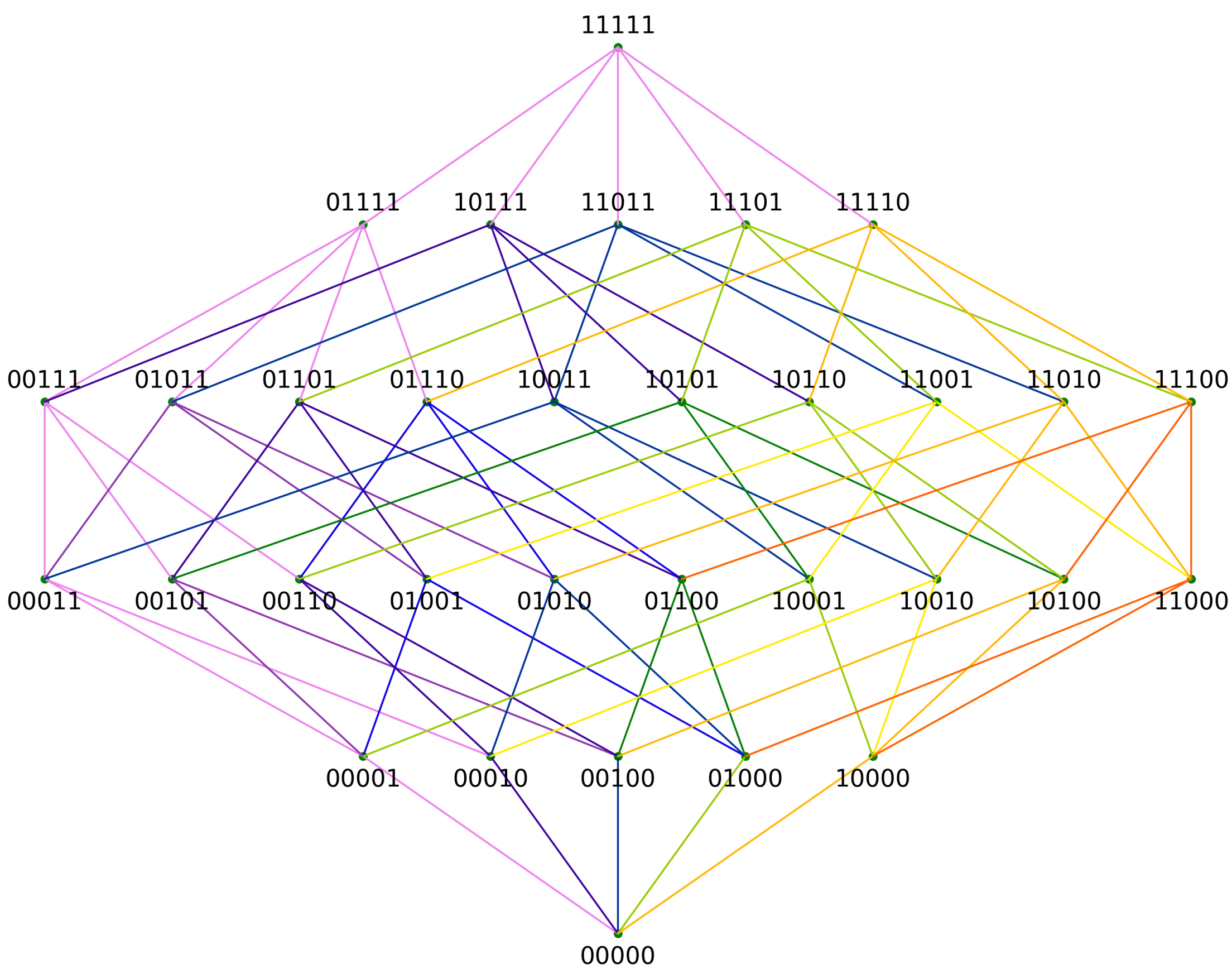

Figure 104: The hypercube graph for $n = 5$

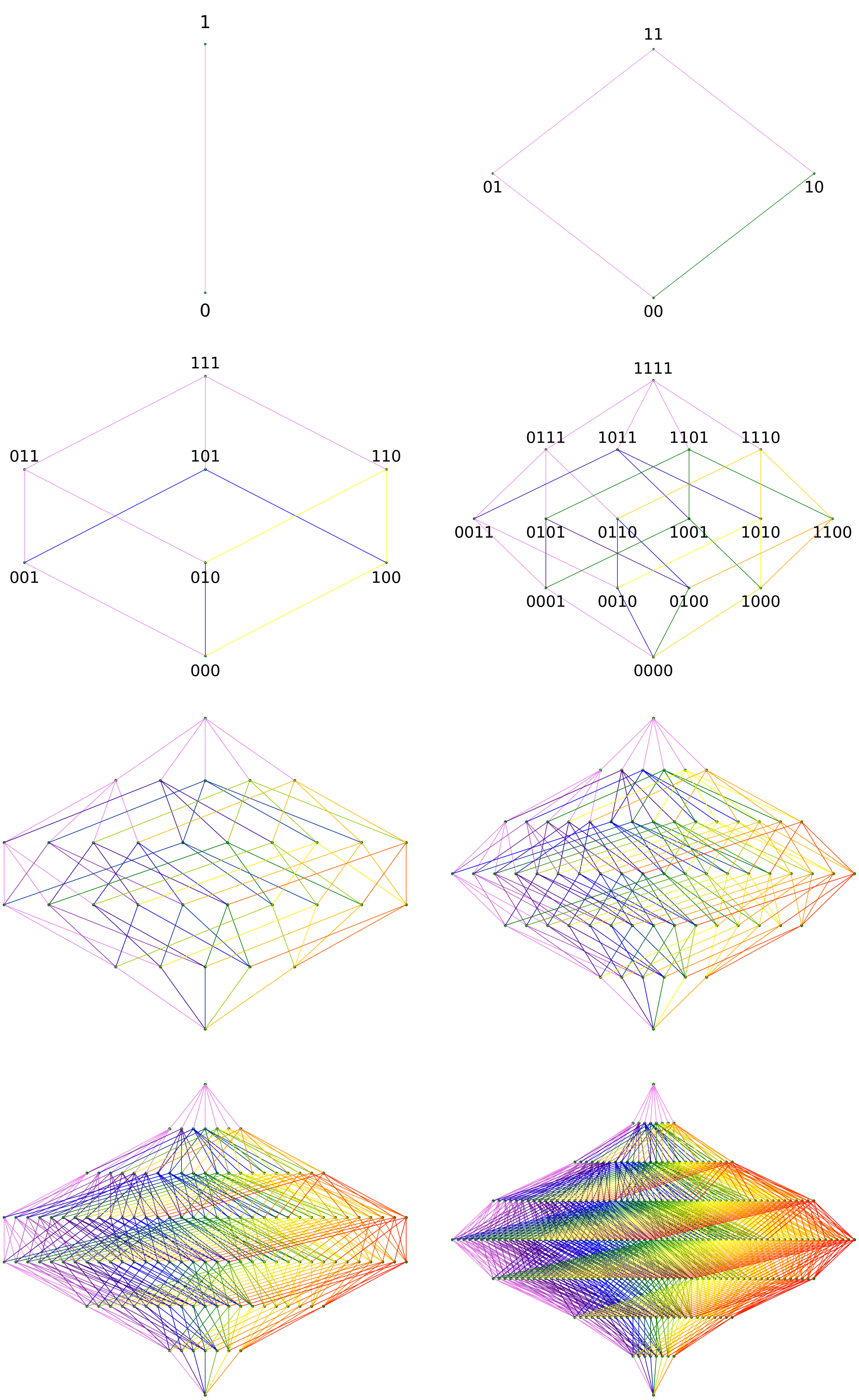

Figure 105: The hypercube graphs for $n = 1$ to $n = 8$

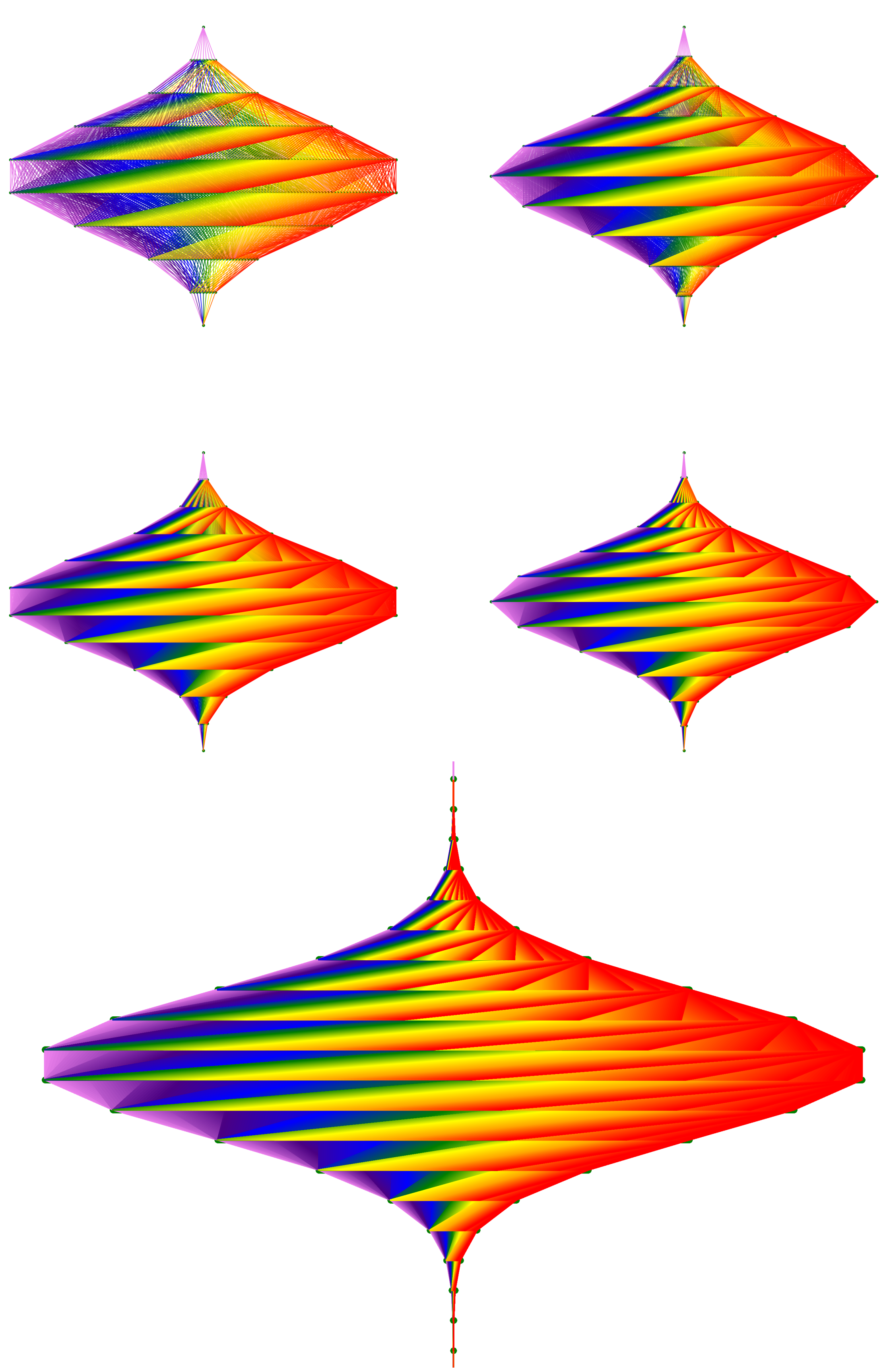

Figure 106: The hypercube graphs for $9 \leq n \leq 12$ and $n = 21$

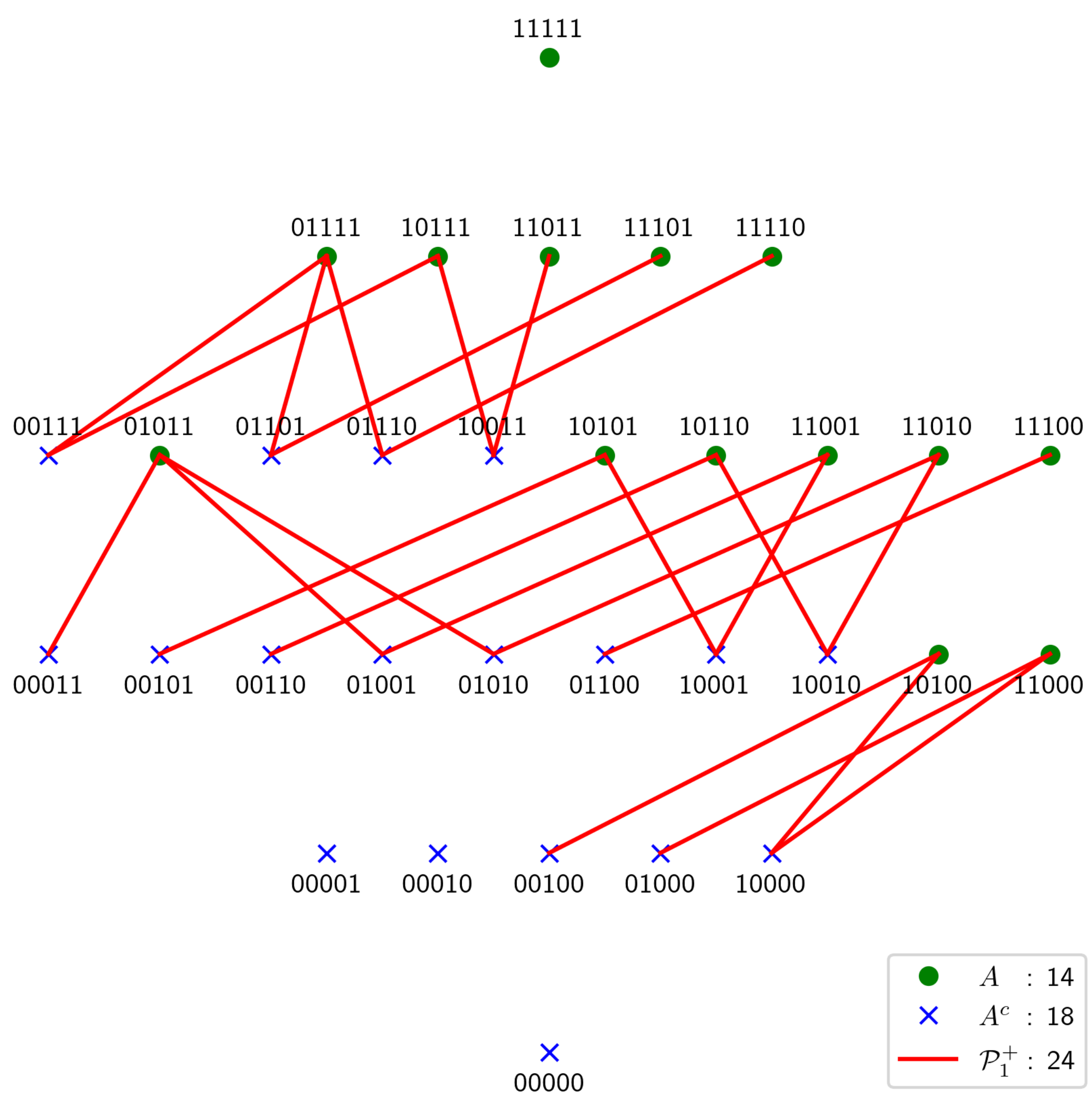

Figure 107: The set $A$ such that $\text{Minimal}(A) = \{ 10100, 11000, 01011 \}$

# 49 The higher positive pivots ⊙

As already said, our goal is to control the positive pivots of arbitrary orders, that we define next. Let $A$ be an increasing subset of $\{0,1\}^n$ and let $\omega$ be a configuration of $\{0,1\}^n$. In section 39, we have already defined the positive pivots $\mathcal{P}^+(A,\omega)$ as

$$\mathcal{P}^+(A,\omega) \,=\, \big\{\, i \in \{\,1,\dots,n\,\} : \omega_i \not\in A,\, \omega^i \in A \,\big\}\,.$$

Let $k \geq 1$. In subsection 44.1, we defined the set $\mathcal{P}_k(A,\omega)$ of the pivots of order $k$. We shall define here the set $\mathcal{P}_k^+(A,\omega)$ of the positive pivots of order $k$, which will be a specific subset of $\mathcal{P}_k(A,\omega)$.

For $I$ a subset of $\{\,1,\dots,n\,\}$, we denote by $\omega_I$ (respectively $\omega^I$) the configuration $\omega$ in which all the components whose index belongs to $I$ are set to 0 (resp. 1):

$$\forall j \in \{\,1,\dots,n\,\}$$
$$\omega_I(j) \,=\, \begin{cases} \omega(j) & \text{if } j \not\in I \\ 0 & \text{if } j \in I \end{cases}, \qquad \omega^I(j) \,=\, \begin{cases} \omega(j) & \text{if } j \not\in I \\ 1 & \text{if } j \in I \end{cases}.$$

For an integer $k \geq 0$, we define the set of the positive pivots of order $k$ as

$$\mathcal{P}_k^+(A,\omega) \,=\, \Big\{\, i \in \{\,1,\dots,n\,\} : \exists I \subset \{\,1,\dots,n\,\} \quad |I| \leq k-1, \\ \forall j \in I \quad \omega(j) = 1 \quad i \in \mathcal{P}^+\big(A,\omega_I\big) \,\Big\}\,. \tag{49.1}$$

Notice that $\mathcal{P}_1^+(A,\omega)$ coincides with the classical set $\mathcal{P}^+(A,\omega)$. However, the definition of the higher positive pivots is different from the definition of the higher pivots given in subsection 44.1. For instance, in formula (49.1), we require that all the components that are modified in $\omega$ were initially equal to 1, whereas in formula (44.1), there is no constraint on the states of the components $\omega(i)$, $i \in I$ (hence the adjective positive). However, we have the obvious inclusion

$$\forall \omega \in \{0,1\}^n \qquad \mathcal{P}_k^+(A,\omega) \,\subset\, \mathcal{P}_k(A,\omega)\,.$$

Let us list some simple immediate properties of the random set $\mathcal{P}_k^+(A,\omega)$:

- $\forall \omega \in \{0,1\}^n \quad \mathcal{P}_{n+1}^+(A,\omega) \,=\, \mathcal{P}_n^+(A,\omega)\,.$
- $\forall \omega \in \{0,1\}^n \quad \forall k \in \{\,0,\dots,n\,\} \qquad \mathcal{P}_k^+(A,\omega) \,\subset\, \mathcal{P}_{k+1}^+(A,\omega)\,.$
- $\forall k \in \{\,0,\dots,n\,\} \quad \forall \omega \leq \eta \in A \qquad \mathcal{P}_k^+(A,\omega) \,\supset\, \mathcal{P}_k^+(A,\eta)\,.$

Of course, we could define analogously the higher negative pivots, but this will not be necessary. The relevance of these definitions is that the intersection sets revealed by the intertwined explorations are indeed higher negative pivots for a decreasing disconnection event. Up to a symmetry, the main challenge now is to get a quantitative control on the mean number of positive pivots of a fixed order $k \geq 1$ associated to an increasing event.

# 50 Extremal problems for $E(|\mathcal{P}_1^+|)$

Our goal in this section is to reprove the two inequalities (47.1) and (47.2) with the help of a novel technique. Instead of relying on the Margulis-Russo formula, we will formulate two extremal problems whose solutions will yield the desired inequalities. We deal first with inequality (47.1) in subsection 50.1. For small values of $n$, that is $n \leq 7$, the solutions of this problem are obtained with the help of an exhaustive search performed by a computer, and the results are presented in subsection 50.2. The deviations inequality (47.2) is harder, the associated extremal problem is in fact a variant of the one associated to inequality (47.1) for which we introduce an additional constraint on the size of the set. This work is done in subsection 50.3.

## 50.1 Unconstrained problem

We reformulate now the inequality (47.1) as an extremal problem. To make things simpler, we consider the symmetric case $p = 1/2$. For $A$ an increasing subset of $\{0,1\}^n$, we define

$$\phi_1(A) \,=\, \sum_{\omega \in A} \left|\mathcal{P}_1^+(A,\omega)\right|.$$

We wish to compute the maximum

$$M_1(n) \,=\, \max\,\big\{\,\phi_1(A) : A \text{ increasing subset of } \{0,1\}^n\,\big\}\,, \tag{50.1}$$

and to find the configurations realizing the maximum. The solution to this problem is easy and well-known.

**Proposition 50.1.** *The value of the maximum* (50.1) *is given by*

$$M_1(n) \,=\, 2^n E\Big(\Big|\sum_{1\leq i\leq n} \omega(i) - \frac{n}{2}\Big|\Big) \,=\, \begin{cases} \displaystyle\binom{n}{\frac{n}{2}}\frac{n}{2} & \text{if } n \text{ is even}\,, \\ \displaystyle\binom{n}{\frac{n+1}{2}}\frac{n+1}{2} & \text{if } n \text{ is odd}\,. \end{cases} \tag{50.2}$$

*Any set $A$ such that*

$$\Big\{\omega : \sum_{1\leq i\leq n} \omega(i) > \frac{n}{2}\Big\} \subset A \subset \Big\{\omega : \sum_{1\leq i\leq n} \omega(i) \geq \frac{n}{2}\Big\} \tag{50.3}$$

*will realize the maximum. If $n$ is even, there are $2^{\binom{n}{n/2}}$ such sets. If $n$ is odd, there exists only one set $A^*$ satisfying* (50.3), *which is*

$$A^* \,=\, \Big\{\omega : \sum_{1\leq i\leq n} \omega(i) > \frac{n}{2}\Big\}. \tag{50.4}$$

*The boolean function $1_{A^*}$ associated to $A^*$ is called the majority function.*

*First proof.* The problem can be solved with a simple computation. Indeed, we have

$$\phi_1(A) \;=\; 2^n E\big(\big|\mathcal{P}_1^+(A)\big|1_A\big) \;=\; 2^{n-1}E\big(\big|\mathcal{P}_1^+(A)\big|\big)\,. \tag{50.5}$$

Using the Margulis-Russo formula, we have

$$E\big(\big|\mathcal{P}_1^+(A)\big|\big) \;=\; 2E\big(1_A S_n\big)\,, \tag{50.6}$$

where $S_n$ is the random sum

$$S_n(\omega) \;=\; \sum_{1\leq i\leq n}\big(2\omega(i)-1\big) \;=\; 2\Big(\sum_{1\leq i\leq n}\omega(i)-\frac{n}{2}\Big)\,. \tag{50.7}$$

Substituting the identity (50.6) in (50.5), we get

$$\phi_1(A) \;=\; 2^n E\big(1_A S_n\big)\,. \tag{50.8}$$

Obviously, we have the inequality

$$\phi_1(A) \;\leq\; 2^n E\big(1_{\{\,S_n>0\,\}}1_A S_n\big)\,.$$

Therefore the maximum (50.1) is reached with the set $A=\{\,\omega : S_n(\omega)>0\,\}$ which, in view of (50.7), coincides with the set $A^*$ defined in (50.4). The value

$$M_1(n) \;=\; 2^n E\big(1_{\{\,S_n>0\,\}}S_n\big)$$

has already been computed in formula (42.17). We offer here a different method, which consists in computing the number of pivotal components over the set $A^*$, given by

$$\phi_1(A^*) \;=\; \sum_{\omega:S_n(\omega)>0}\big|\mathcal{P}_1^+(S_n>0,\omega)\big|\,.$$

We distinguish two cases, depending on the parity of $n$.

**First case: $n=2p$ for some $p\geq 1$.** We have then

$$\phi_1(A^*) \;=\; \sum_{\omega:S_n(\omega)=2}(p+1) \;=\; \binom{2p}{p+1}(p+1) \;=\; \binom{2p}{p}p\,.$$

**Second case: $n=2p+1$ for some $p\geq 0$.** We have then

$$\phi_1(A^*) \;=\; \sum_{\omega:S_n(\omega)=1}(p+1) \;=\; \binom{2p+1}{p+1}(p+1)\,.$$

The claim concerning the maximizers follows from formulas (50.7) and (50.8). Indeed, for any set $A$, we have

$$\phi_1(A) \;=\; \phi_1(A^*)+2^n E\big(1_{A\setminus A^*}S_n\big)-2^n E\big(1_{A^*\setminus A}S_n\big)\,.$$

The expectation over $A\setminus A^*$ is non-positive, while the expectation over $A^*\setminus A$ is non-negative. Thus $\phi_1(A)=\phi_1(A^*)$ if and only if $A$ satisfies (50.3). □

Unfortunately, the previous argument does not lend itself well to generalisations. So, as a preliminary step, we will try to reprove the previous result with a more robust technique inspired from extremal set theory, a mathematical field where such questions are commonplace. The technique consists in starting with a candidate set $A$ and trying to perturb it slightly to increase its $\phi_1$ value. Of course, a set realizing the maximum of $\phi_1$ should not admit any perturbation increasing its $\phi_1$ value. If we have at our disposal an adequate collection of admissible perturbations, then we hope to be able to find and describe completely the maximizers. An important difficulty here is that an admissible perturbation of a candidate set should yield another monotone set.

*Second proof.* Let $A$ be an increasing subset of $\{0,1\}^n$. We consider first perturbations which consist in removing a configuration from $A$. The only configurations that can be removed from $A$ without destroying the monotonicity are the minimal configurations of $A$, i.e., the elements of Minimal($A$). So, let $\sigma$ belong to Minimal($A$) and let $B = A \setminus \{\,\sigma\,\}$. We compute next the variation of $\phi_1$. By definition, we have

$$\phi_1(B) - \phi_1(A) \,=\, \sum_{\omega \in B} \big|\mathcal{P}_1^+(B,\omega)\big| - \sum_{\omega \in A} \big|\mathcal{P}_1^+(A,\omega)\big|\,. \tag{50.9}$$

We recall that the Hamming distance $H(\omega,\rho)$ between two configurations $\omega,\rho$ is the number of components where they differ, i.e.,

$$H(\omega,\rho) \,=\, \big|\{\, i \in \{\,1,\dots,n\,\} : \omega(i) \neq \rho(i)\,\}\big| \,=\, \sum_{1\leq i\leq n} \big|\omega(i) - \rho(i)\big|\,.$$

It turns out that the removal of $\sigma$ will affect only the pivotal sets of the configurations which are at Hamming distance 0 or 1 from $\sigma$.

**Lemma 50.2.** *For any configuration $\omega$ such that $H(\omega,\sigma)\geq 2$, we have*

$$\mathcal{P}_1^+(A,\omega) = \mathcal{P}_1^+(B,\omega)\,.$$

*Proof.* Let $\omega \in \{0,1\}^n$ be such that $H(\omega,\sigma)\geq 2$. Let also $i$ be an index in $\{\,1,\dots,n\,\}$. Certainly we have $\omega \neq \sigma$ and $\omega_i \neq \sigma$. Thus we have the equivalences

$$\omega \in A \quad\Longleftrightarrow\quad \omega \in B\,,\qquad \omega_i \in A \quad\Longleftrightarrow\quad \omega_i \in B\,,$$

implying that

$$i \in \mathcal{P}_1^+(A,\omega) \quad\Longleftrightarrow\quad i \in \mathcal{P}_1^+(B,\omega)\,.$$

This being true for any $i \in \{\,1,\dots,n\,\}$, the claim of the lemma is proved. $\square$

It follows from lemma 50.2 that the only terms that might not vanish in the sum appearing in (50.9) correspond to configurations $\omega$ such that $H(\omega,\sigma)\leq 1$. For $\omega=\sigma$, we have $\mathcal{P}_1^+(A,\omega) = \operatorname{supp}\sigma$. Since $\sigma$ is a minimal element of $A$, the configurations $\omega$ of $A$ such that $H(\omega,\sigma)=1$ are of the form $\sigma^k$ for an index $k$ not belonging to $\operatorname{supp}\sigma$. For $\omega=\sigma^k$, with $k\notin \operatorname{supp}\sigma$, we have

$$k \in \mathcal{P}_1^+(B,\sigma^k)\,,\quad k \notin \mathcal{P}_1^+(A,\sigma^k)\,,\quad \mathcal{P}_1^+(B,\sigma^k)\setminus\{\,k\,\} \,=\, \mathcal{P}_1^+(A,\sigma^k)\,,$$

whence

$$\big|\mathcal{P}_1^+(B,\sigma^k)\big| - \big|\mathcal{P}_1^+(A,\sigma^k)\big| \,=\, 1\,.$$

We conclude that

$$\phi_1(B) - \phi_1(A) \,=\, n - 2\big|\operatorname{supp}\sigma\big| \,. \tag{50.10}$$

As a consequence, for a set $A$ which maximizes $\phi_1(A)$, there is no minimal element whose support has cardinality smaller than $n/2$. In a second step, we consider perturbations which consist in adding a configuration to the set $A$. In order to get a monotone set, the only elements that can be added to $A$ are the maximal elements of $A^c$. We perform a computation similar to the one we did for the removal of minimal elements, and we conclude that, for a set $A$ which maximizes $\phi_1(A)$, there is no maximal element of $A^c$ whose support has cardinality smaller than $n/2$. □

Undoubtedly, this second proof is considerably longer than the first one. However, its principle is completely different and we hope that it can be applied to handle the higher positive pivots. Using Stirling's formula, we obtain the asymptotic behavior of $M_1(n)$:

$$M_1(n) \,\sim\, 2^n\sqrt{\frac{n}{2\pi}} \quad \text{as} \quad n \to +\infty \,.$$

This approximation is incredibly good even for $n$ small, as illustrated in figure 108. A direct check with the help of formula (50.2) shows that $M_1(2n) = 2M_1(2n-1)$ for any $n \geq 1$.

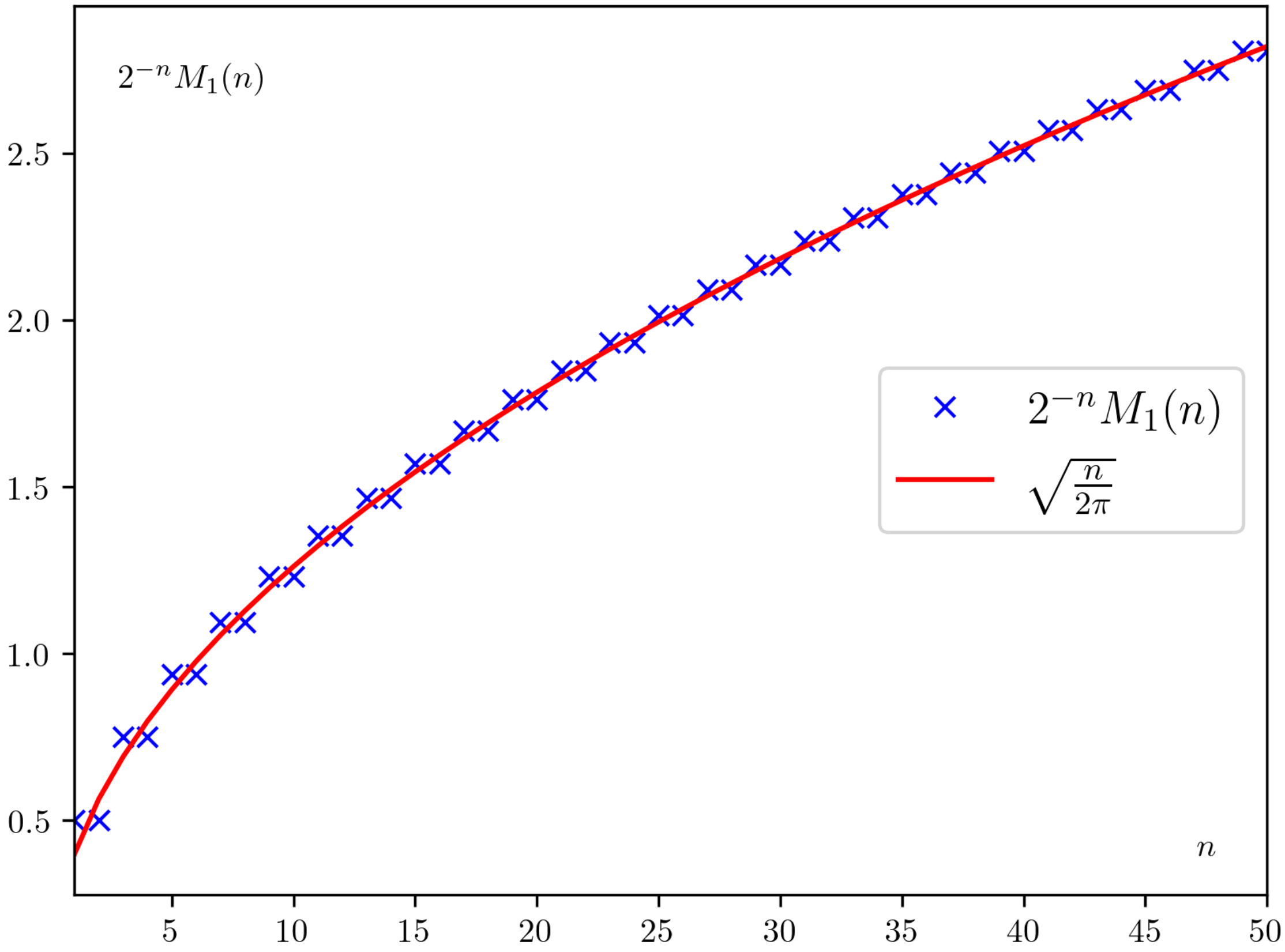

Figure 108: The curve $2^{-n}M_1(n)$ and its Stirling approximation

## 50.2 Maximizers of $E(|\mathcal{P}_1^+|)$ for $2 \le n \le 7$

We show here the results found for the maximizers of $E(|\mathcal{P}_1^+|)$ after an exhaustive search of the monotone subsets of the hypercubes for $2 \le n \le 7$. We use the representation explained in subsection 48.3, and we present the maximizers up to equivalence modulo a permutation of the components. In the figures, we focus on the edge boundary of the sets, so in some of them we remove most of the layers which do not contain a single pivot. To warm up, we start in figure 109 with the three maximizers for $n = 2$.

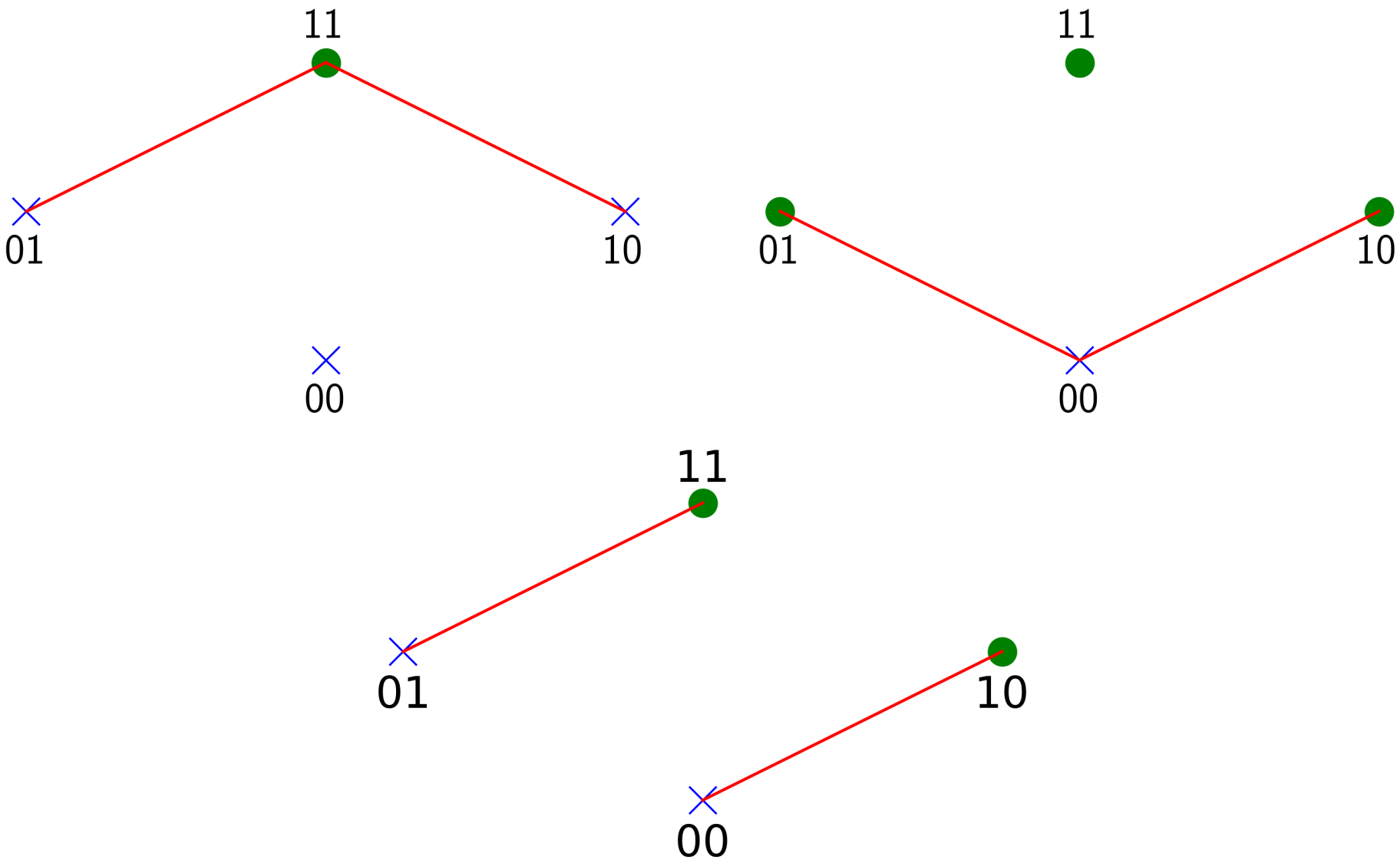

Figure 109: The three maximizers of $E(|\mathcal{P}_1^+|)$ for $n = 2$

We move on with the maximizers for $n = 4$. Up to equivalence, there are 11 maximizers which are listed in the table 3 and represented in figure 110. The number 11 should be compared with the total number of maximizers, which is equal to $2^{\binom{4}{2}} = 64$. For $n = 3, 5, 7$, there is only one maximizer, shown in figure 111. For $n = 6$, up to equivalence, there are 2136 maximizers, and $2^{\binom{6}{3}} = 2^{20} = 1048576$ in total, some of them are shown in figure 112.

| $\|A\|$ | Minimal$(A)$ |
|---|---|
| 7 | 1100 0011 |
| 7 | 1100 1010 0111 |
| 8 | 1100 1010 0101 |
| 8 | 1100 1010 0110 |
| 6 | 1100 1011 0111 |
| 9 | 1100 1010 0101 0011 |
| 9 | 1100 1010 1001 0110 |
| 5 | 1110 1101 1011 0111 |
| 8 | 1100 1010 1001 0111 |
| 10 | 1100 1010 1001 0110 0101 |
| 11 | 1100 1010 1001 0110 0101 0011 |

Table 3: The 11 maximizers of $E(|\mathcal{P}_1^+|)$ for $n = 4$

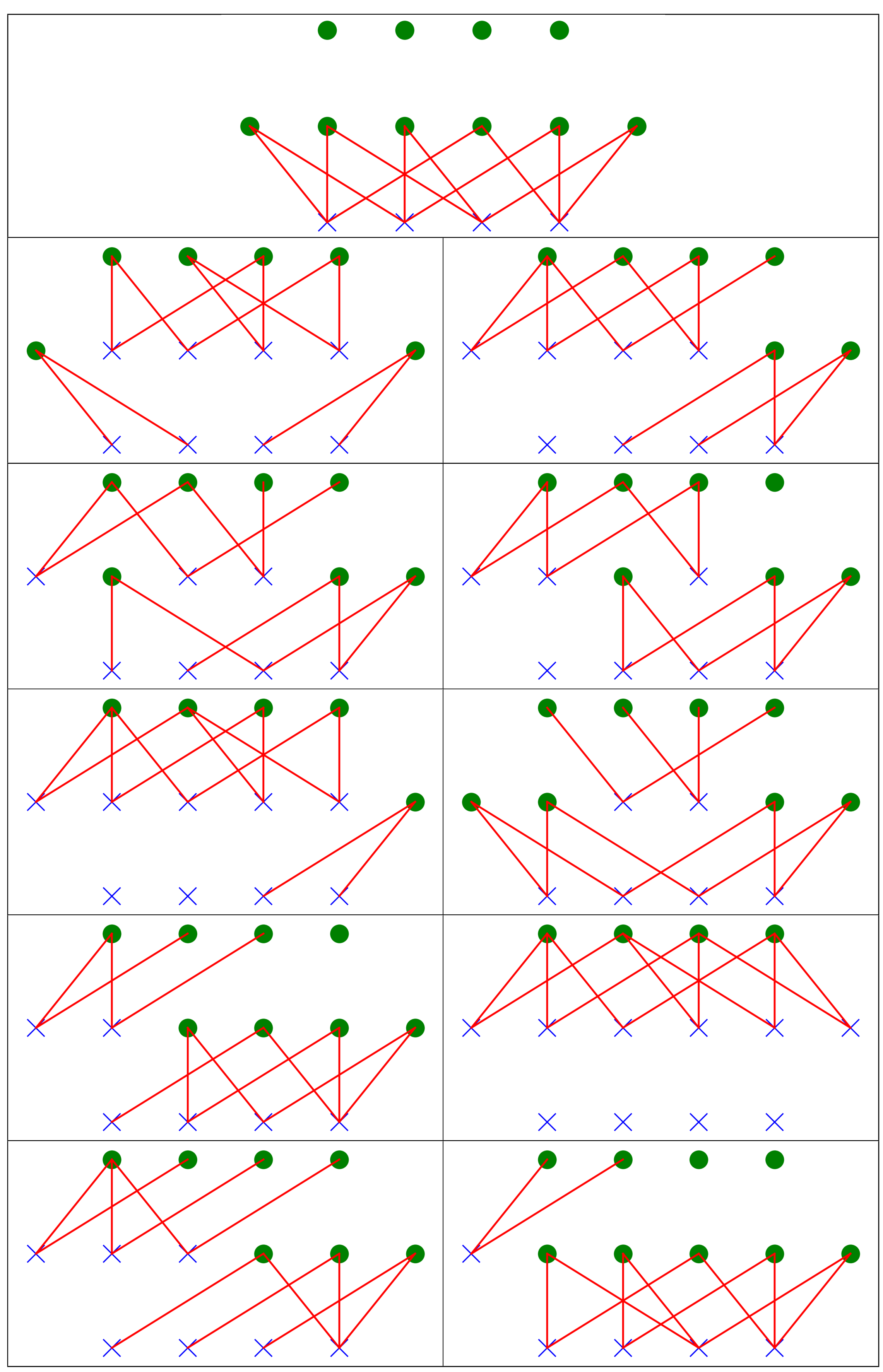

Figure 110: The 11 maximizers of $E(|\mathcal{P}_1^+|)$ for $n = 4$

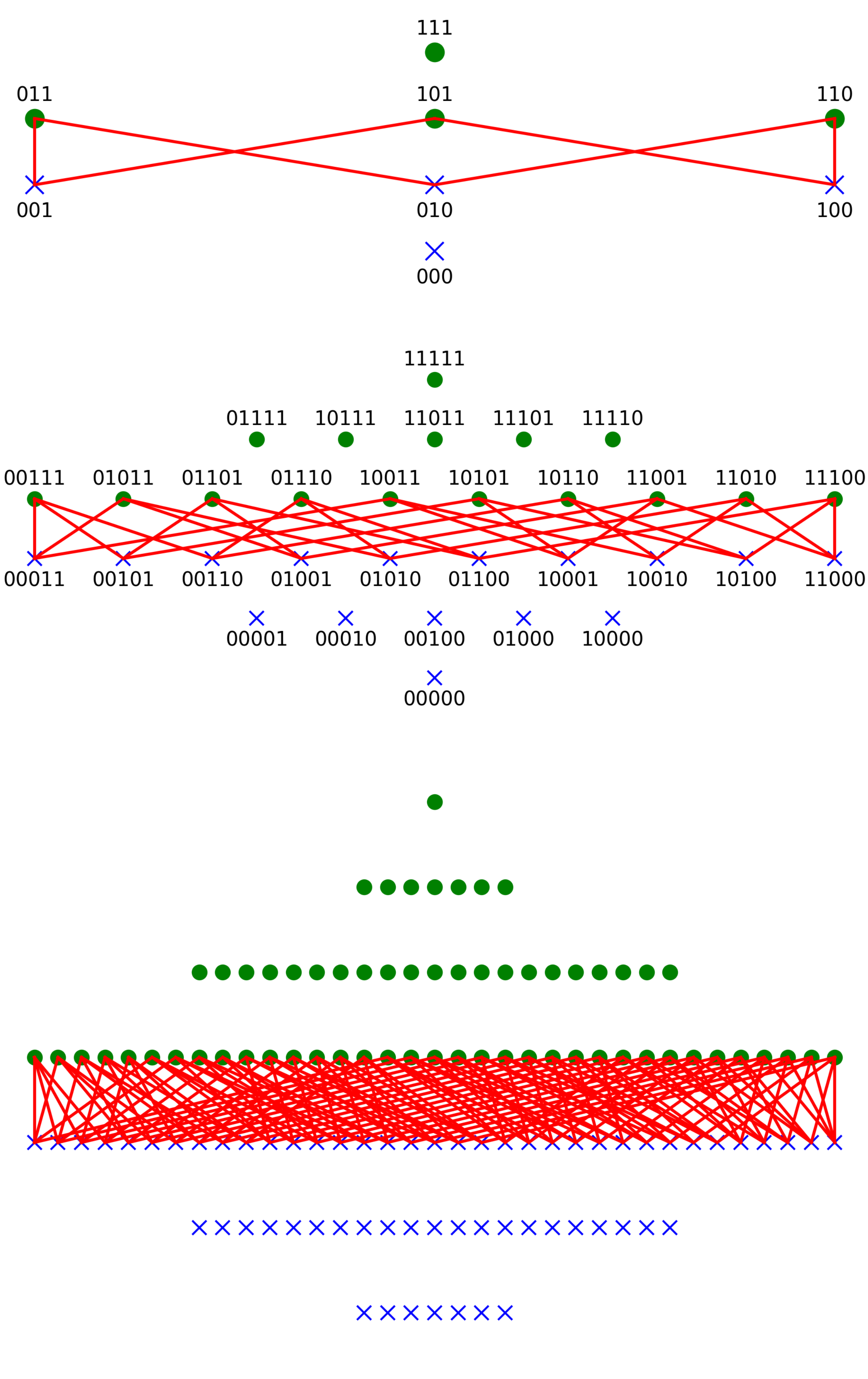

Figure 111: The unique maximizers of $E(|\mathcal{P}_1^+|)$ for $n = 3, 5, 7$

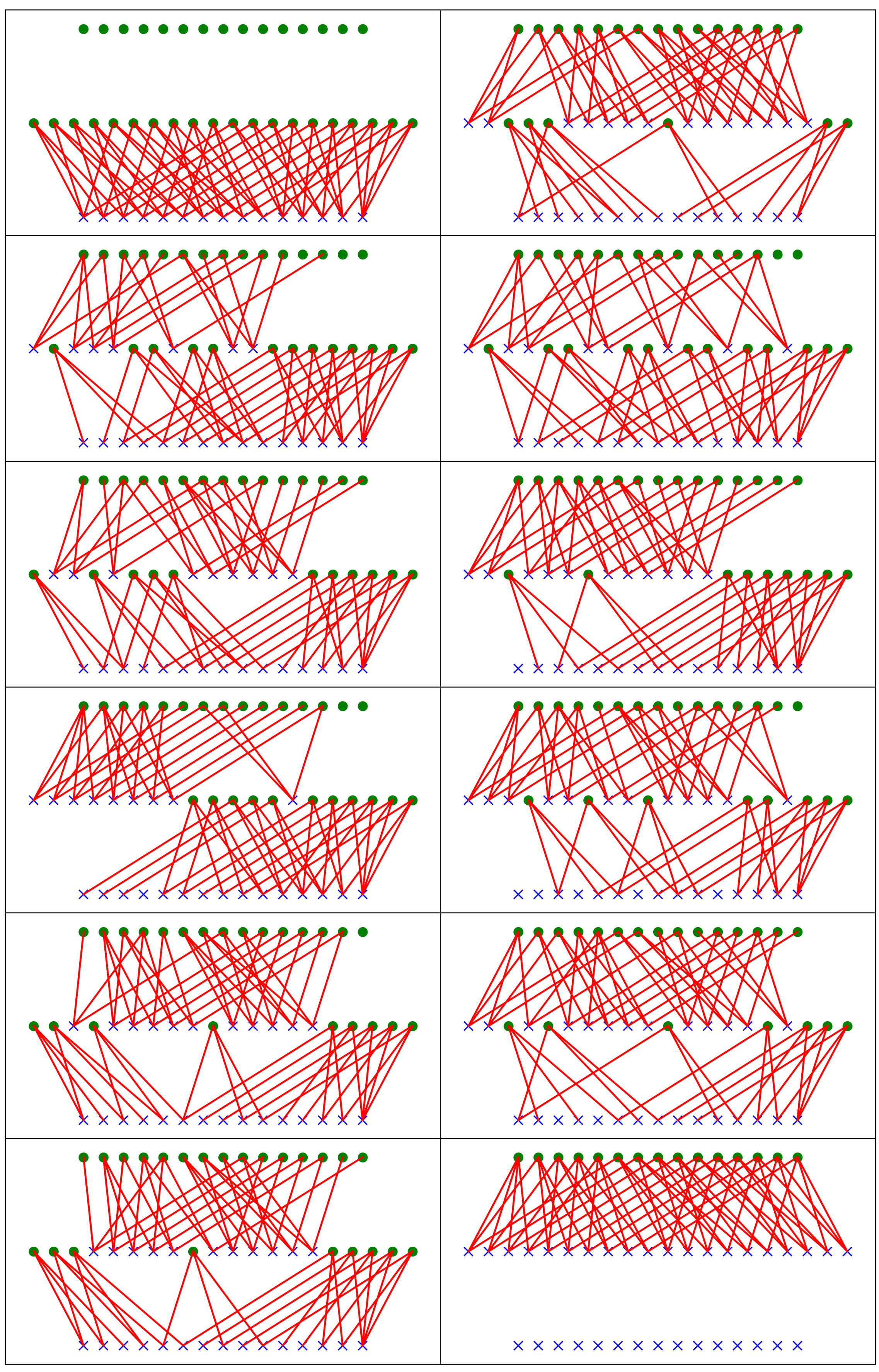

Figure 112: Some examples of the 2136 maximizers of $E(|\mathcal{P}_1^+|)$ for $n = 6$

## 50.3 Constrained problem

We try now to prove the deviations inequality (47.2) through the study of another extremal set problem. To that end, we introduce a size constraint on the sets. So, we fix an integer $N$ such that $1 \leq N < 2^n$ and we consider the maximization problem

$$M_1(n,N) \,=\, \max\big\{\,\phi_1(A) : A \text{ increasing subset of } \{0,1\}^n,\, |A|=N\,\big\}\,. \tag{50.11}$$

For $K \in \{\,0,\dots,n\,\}$, we denote by $N(K)$ the sum

$$N(K) \,=\, \sum_{K\leq k\leq n} \binom{n}{k}\,. \tag{50.12}$$

We say that a configuration $\omega$ of $\{0,1\}^n$ belongs to the level $k$ if it has exactly $k$ components equal to 1, or equivalently if its support has cardinality $k$. This terminology refers to the layer representation of the hypercube described in subsection 48.3.

**Proposition 50.3.** *Let $N$ be such that $1 \leq N < 2^n$. Let $K$ be the unique integer in $\{\,1,\dots,n\,\}$ such that*

$$N(K) \,\leq\, N < N(K-1)\,. \tag{50.13}$$

*The value of the maximum* (50.11) *is given by*

$$M_1(n,N) \,=\, 2\sum_{K\leq k\leq n}(k-K+1)\binom{n}{k} + (2K-2-n)N \tag{50.14}$$

$$= \binom{n}{K}K + (2K-2-n)\big(N-N(K)\big)\,. \tag{50.15}$$

*The maximizers are the sets $A$ such that $|A| = N$ and*

$$\Big\{\omega : \sum_{1\leq i\leq n}\omega(i)\geq K\Big\} \,\subset\, A \,\subset\, \Big\{\omega : \sum_{1\leq i\leq n}\omega(i)\geq K-1\Big\}\,. \tag{50.16}$$

*These sets are obtained by adjoining $N - N(K)$ configurations of level $K-1$ to the set $\big\{\,\omega : \sum_{i=1}^n \omega(i) \geq K\,\big\}$. The number of maximizers is equal to the binomial coefficient*

$$\binom{\binom{n}{K-1}}{N-N(K)}\,. \tag{50.17}$$

*If there is equality in* (50.13)*, i.e., if $N = N(K)$, then there is a unique maximizer, which is the set*

$$\Big\{\omega : \sum_{1\leq i\leq n}\omega(i)\geq K\Big\}\,.$$

We see from formula (50.15) that

$$\forall K \in \{\, 0, \ldots, n \,\} \qquad M_1\big(n, N(K)\big) \;=\; \binom{n}{K} K \,.$$

Let us fix $K \in \{\, 1, \ldots, n \,\}$. We claim that the map

$$N \in [N(K), N(K-1)] \mapsto M_1(n, N)$$

is simply the linear interpolation between the points

$$\left( N(K), \binom{n}{K} K \right) \qquad \text{and} \qquad \left( N(K-1), \binom{n}{K-1} (K-1) \right) . \tag{50.18}$$

Indeed, if we evaluate the formula (50.15) at the point $N = N(K-1)$, we get

$$\begin{aligned}
&\binom{n}{K} K + (2K - 2 - n)\big(N(K-1) - N(K)\big) \\
&\qquad\qquad = \binom{n}{K} K + (2K - 2 - n)\binom{n}{K-1} \\
&\qquad = K\binom{n}{K} + (K-1)\binom{n}{K-1} - (n - K + 1)\binom{n}{K-1} \\
&\qquad\qquad\qquad\qquad = (K-1)\binom{n}{K-1} .
\end{aligned}$$

Thus the formula (50.15) takes the values (50.18) at the extremities of the interval $[N(K), N(K-1)]$. Now that we understand better the structure of the map $N \mapsto M_1(n, N)$, it becomes easy to compute its maximum. Necessarily, we have

$$\max_{0 \leq N \leq 2^n} M_1(n, N) \;=\; \max_{0 \leq K \leq n} M_1\big(n, N(K)\big) \,.$$

Moreover, for $1 \leq K \leq n$, we have

$$M_1\big(n, N(K)\big) \;=\; \binom{n}{K} K \;=\; \binom{n-1}{K-1} n \,.$$

If $n-1$ is even, the last binomial coefficient is maximal for $K - 1 = (n-1)/2$, that is $K = (n+1)/2$. If $n-1$ is odd, the binomial coefficient is maximal for $K - 1 = (n - 1 \pm 1)/2$, that is $K = n/2$ or $K = n/2 + 1$. In the end, we are back to the value obtained in (50.2) when computing the maximum $M_1(n)$ of the unconstrained problem (50.1).

For $k \in \{\, 0, \ldots, n \,\}$, we define the set $A^*(k)$ as

$$A^*(k) \;=\; \Big\{\, \omega : \sum_{1 \leq i \leq n} \omega(i) \geq k \,\Big\} ,$$

so that the inclusions (50.16) can be rewritten as $A^*(K) \subset A \subset A^*(K-1)$.

*First proof.* We reuse the initial computation done in the first proof of proposition 50.1. We had obtained in formula (50.8) that

$$\phi_1(A) \,=\, 2^{n+1} E\Big(1_A(\omega)\Big(\sum_{1\leq i\leq n} \omega(i) - \frac{n}{2}\Big)\Big)\,.$$

Since we work with the uniform probability measure, the constraint $|A| = N$ can be rewritten as $E\big(1_A(\omega)\big) = N2^{-n}$, whence

$$\phi_1(A) \,=\, 2^{n+1} E\Big(1_A(\omega)\Big(\sum_{1\leq i\leq n} \omega(i) - K + 1\Big)\Big) + (2K-2-n)N\,.$$

Obviously, we have the inequalities

$$\begin{aligned}\phi_1(A) \,&\leq\, 2^{n+1} E\Big(1_A(\omega) 1_{A^*(K-1)}(\omega)\Big(\sum_{1\leq i\leq n} \omega(i) - K + 1\Big)\Big) + (2K-2-n)N \\ &\leq\, 2^{n+1} E\Big(1_{A^*(K-1)}(\omega)\Big(\sum_{1\leq i\leq n} \omega(i) - K + 1\Big)\Big) + (2K-2-n)N \\ &=\, 2 \sum_{K\leq k\leq n} (k-K+1)\binom{n}{k} + (2K-2-n)N\,. \end{aligned} \quad (50.19)$$

There is no guarantee that $A^*(K-1)$ satisfies the additional size constraint. Moreover, the first inequality in (50.19) is an equality if and only if $A \subset A^*(K-1)$, and the second inequality in (50.19) is an equality if and only if $A^*(K) \subset A$. It follows that the maximizers are exactly the sets $A$ which satisfy the inclusions (50.16) and such that $|A| = N$. These sets are obtained by adjoining $N - N(K)$ configurations of level $K-1$ to $A^*(K)$. The number of these sets is given by the binomial coefficient (50.17). In the case of equality $N = N(K)$, this binomial coefficient is equal to 1 and $A^*(K)$ is the unique maximizer. Finally, the value of the maximum is given by the last quantity in (50.19), which gives the first expression for $M_1(n,N)$ presented in formula (50.14). We shall simplify this expression by computing the number of pivotal components of a maximizer. So let $A^*$ be a set satisfying $|A^*| = N$ and $A^*(K) \subset A^* \subset A^*(K-1)$. We call additional configurations the configurations in $A^* \setminus A^*(K)$, these are configurations of level $K-1$. All the components of the additional configurations are pivotal for $A^*$. In the configurations of $A^*$ of level $K+1$ and higher, no component is pivotal. For the remaining configurations of $A^*$, which belong to level $K$, two cases occur. Either they do not contain an additional configuration, and all their components are pivotal, or they contain at least one additional configuration. In fact, any additional configuration neutralizes exactly one pivotal component in $n-K+1$ configurations of level $K$. Moreover, one potentially pivotal component of a configuration in the level $K$ can be neutralized by only one additional configuration. We have thus

$$M_1(n,N) \,=\, \phi_1(A^*) \,=\, \binom{n}{K} K + \big(N - N(K)\big)\big(K - 1 - (n-K+1)\big)\,.$$

This is the second expression for $M_1(n,N)$ presented in formula (50.15). $\square$

*Second proof.* Let $A$ be an increasing subset of $\{0,1\}^n$ such that $|A| = N$. We consider perturbations of $A$ which consist in removing a configuration and adding another one at the same time. The point is that these perturbations do not affect the constraint $|A| = N$. Moreover, we have to be careful not to destroy the monotonicity of $A$. Thus the only configurations that can be removed from $A$ are the minimal configurations of $A$, while the only elements that can be added to $A$ are the maximal elements of $A^c$. So, let $\sigma$ belong to Minimal$(A)$, let $\eta$ belong to Maximal$(A^c)$, and let $B = A \setminus \{\,\sigma\,\} \cup \{\,\eta\,\}$. We compute next the variation of $\phi_1$. By definition, we have

$$\phi_1(B) - \phi_1(A) \,=\, \sum_{\omega \in B} \big|\mathcal{P}_1^+(B,\omega)\big| - \sum_{\omega \in A} \big|\mathcal{P}_1^+(A,\omega)\big|\,. \tag{50.20}$$

It follows from lemma 50.2 that the only terms that might not vanish in the sums appearing in (50.20) correspond to configurations $\omega$ such that $H(\omega,\sigma) \leq 1$ or $H(\omega,\eta) \leq 1$. For $\omega = \sigma$, we have $\mathcal{P}_1^+(A,\omega) = \operatorname{supp}\sigma$. For $\omega = \eta$, we have $\mathcal{P}_1^+(B,\eta) = \operatorname{supp}\eta$. Since $\sigma$ is a minimal element of $A$, the configurations $\omega$ of $A \cup B$ such that $H(\omega,\sigma) = 1$ are of the form $\sigma^k$ for an index $k$ not belonging to $\operatorname{supp}\sigma$. Since $\eta$ is a maximal element of $A^c$, the configurations $\omega$ of $A \cup B$ such that $H(\omega,\eta) = 1$ are of the form $\eta^k$ for an index $k$ not belonging to $\operatorname{supp}\eta$. For $\omega = \sigma^k$, with $k \notin \operatorname{supp}\sigma$, we have

$$k \in \mathcal{P}_1^+(B,\sigma^k)\,, \quad k \notin \mathcal{P}_1^+(A,\sigma^k)\,, \quad \mathcal{P}_1^+(B,\sigma^k) \setminus \{\,k\,\} \,=\, \mathcal{P}_1^+(A,\sigma^k)\,,$$

whence

$$\big|\mathcal{P}_1^+(B,\sigma^k)\big| - \big|\mathcal{P}_1^+(A,\sigma^k)\big| \,=\, 1\,.$$

For $\omega = \eta^k$, with $k \notin \operatorname{supp}\eta$, we have

$$k \in \mathcal{P}_1^+(A,\eta^k)\,, \quad k \notin \mathcal{P}_1^+(B,\eta^k)\,, \quad \mathcal{P}_1^+(A,\eta^k) \setminus \{\,k\,\} \,=\, \mathcal{P}_1^+(B,\eta^k)\,,$$

whence

$$\big|\mathcal{P}_1^+(B,\eta^k)\big| - \big|\mathcal{P}_1^+(A,\eta^k)\big| \,=\, -1\,.$$

We conclude that

$$\begin{aligned}\phi_1(B) - \phi_1(A) \,&=\, \big|\operatorname{supp}\eta\big| - \big|\operatorname{supp}\sigma\big| + \big(n - \big|\operatorname{supp}\sigma\big|\big) - \big(n - \big|\operatorname{supp}\eta\big|\big)\\ &=\, 2\big(\big|\operatorname{supp}\eta\big| - \big|\operatorname{supp}\sigma\big|\big)\,.\end{aligned}$$

Now, if a set $A$ maximizes $\phi_1(A)$, then it does not admit any perturbation which increases the value of $\phi_1$. Therefore any maximizer $A$ of the constrained problem (50.11) is such that

$$\max\,\big\{\,|\operatorname{supp}\eta| : \eta \in A^c\,\big\} \,\leq\, \min\,\big\{\,|\operatorname{supp}\sigma| : \eta \in A\,\big\}\,. \tag{50.21}$$

However, for a monotone set $A$, we have the inequality

$$\max\,\big\{\,|\operatorname{supp}\eta| : \eta \in A^c\,\big\} \,\geq\, \min\,\big\{\,|\operatorname{supp}\sigma| : \eta \in A\,\big\} - 1\,. \tag{50.22}$$

The inequalities (50.21) and (50.22) imply that there exists an integer $K$ such that

$$\Big\{\omega:\sum_{i=1}^{n}\omega(i)\geq K\Big\}\subset A\subset\Big\{\omega:\sum_{i=1}^{n}\omega(i)\geq K-1\Big\}. \tag{50.23}$$

Necessarily, this integer $K$ has to be the unique integer defined by the inequalities (50.13). At this point, it is an easy matter to check that the value of $\phi_1$ is the same for all the sets satisfying (50.23), hence all these sets are maximizers. The explicit formulas (50.14) and (50.15) might be computed in the same way as in the first proof. The other conclusions of proposition 50.3 follow easily. □

The two expressions (50.14) and (50.15) of $M_1(n,N)$ have been obtained with very different methods. For the first formula, we produced two inequalities leading to the correct upper bound, while for the second formula we computed explicitly the number of pivotal components for an optimal set $A^*$. To be safe, let us check directly that the formulas (50.14) and (50.15) coincide. After removal of the terms appearing in both formulas, everything boils down to the curious identity stated in the next lemma.

**Lemma 50.4.** *Let $n\geq 1$ and let $K$ be such that $0\leq K\leq n$. We have*

$$\binom{n}{K}K=\sum_{K\leq k\leq n}(2k-n)\binom{n}{k}. \tag{50.24}$$

*Proof.* There is certainly a combinatorial interpretation of the identity (50.24), although we did not manage to find it in the literature. This identity is in the classical book [61], see formula (5.18) there. We shall limit ourselves here to a direct verification by induction on $K$. For $K=0$, the result is true because

$$\sum_{0\leq k\leq n}2k\binom{n}{k}=\sum_{1\leq k\leq n}2n\binom{n-1}{k-1}=2n2^{n-1}=n2^n.$$

Suppose now that the result has been proved at rank $K$ where $0\leq K<n$. In order to apply the induction hypothesis, we decompose the right-hand side of (50.24) as follows:

$$\begin{aligned}\sum_{K+1\leq k\leq n}(2k-n)\binom{n}{k}&=\sum_{K\leq k\leq n}(2k-n)\binom{n}{k}-(2K-n)\binom{n}{K}\\&=\binom{n}{K}K-(2K-n)\binom{n}{K}=(n-K)\binom{n}{K}\\&=\binom{n}{K+1}(K+1).\end{aligned}$$

This proves the identity at rank $K+1$ and the induction is completed. □

## 50.4 An associated deviations inequality

Let $A$ be an increasing subset of $\{0,1\}^n$. We would like to use the results of subsection 50.3 on the constrained extremal problem to obtain an inequality linking $P(A)$ and $E\big(|\mathcal{P}^+(A)|\,\big|\,A\big)$, similar to the one stated in proposition 40.2. Let us set

$$N = |A|\,,\quad \lambda \,=\, E\big(|\mathcal{P}^+(A)|\,\big|\,A\big)\,.$$

With the notation of subsection 50.3, we have, thanks to the fact that we work with the uniform probability,

$$\lambda \,=\, \frac{1}{N}\phi_1(A) \,\leq\, \frac{1}{N}M_1(n,N)\,. \tag{50.25}$$

For $n$ fixed, the graph of the map $N \mapsto M_1(n,N)$ is symmetric and has a maximum at $\lfloor \frac{N+1}{2} \rfloor$, as shown in figure 113. The point is that the quantity on the right-hand side of (50.25) is monotone in $N$.

**Lemma 50.5.** *The map*

$$\psi_1(N) \,=\, \frac{1}{N}M_1(n,N)$$

*is decreasing on* $\{\,1,\dots,2^n\,\}$.

*Proof.* Notice first that $\psi_1(1) = n$ and $\psi_1(2^n) = 0$. We compute next the increments of $\psi_1$. Let $N$ belong to $\{\,1,\dots,2^n\,\}$ and let $K$ be such that $N(K) \leq N < N(K-1)$. Thanks to formula (50.15), we have

$$\begin{aligned}\psi_1(N+1)-\psi_1(N) \,=\, &\frac{1}{N+1}\Bigg(\binom{n}{K}K + (2K-2-n)\big(N+1-N(K)\big)\Bigg)\\ &-\frac{1}{N}\Bigg(\binom{n}{K}K + (2K-2-n)\big(N-N(K)\big)\Bigg)\\ &= -\frac{1}{N(N+1)}\Bigg(\binom{n}{K}K - (2K-2-n)N(K)\Bigg)\,.\end{aligned}$$

We checked right after the statement of proposition 50.3 that the formula (50.15) could also be used at the point $N = N(K-1)$, thus the above expression of the increment is valid for $N = N(K-1)-1$. To evaluate this last quantity, we make appeal to the expression (50.12) of $N(K)$ and to lemma 50.4 and we obtain

$$\begin{aligned}\psi_1(N+1)-\psi_1(N) \,=\, &-\frac{1}{N(N+1)}\sum_{K\leq k\leq n}\Bigg(\big((2k-n)-(2K-2-n)\big)\binom{n}{k}\Bigg)\\ &= -\frac{1}{N(N+1)}\sum_{K\leq k\leq n}\Bigg((2k-2K+2)\binom{n}{k}\Bigg)\,.\end{aligned}$$

This last sum is clearly negative. ∎

As a consequence, a lower bound on the value of $\psi_1(N)$ yields an upper bound on $N$, given by its inverse function. More precisely, we define $\psi_1^{-1}$ by

$$\forall \lambda \in [0,n] \qquad \psi_1^{-1}(\lambda) \,=\, \max\,\big\{\, N\geq 0 : \psi_1(N)\geq \lambda \,\big\}\,.$$

We have the implication:

$$\forall N \in \,\{\,1,\dots,2^n\,\}\quad \forall \lambda\in[0,n] \qquad \psi_1(N)\geq\lambda \quad\Longrightarrow\quad N\leq \psi_1^{-1}(\lambda)\,.$$

Since we work with the uniform probability, an upper bound on $N$ can be rewritten as an upper bound on $P(A)$. In the end, the inequality (50.25) and the above implication together yield the inequality

$$P(A)\,\leq\, 2^{-n}\psi_1^{-1}\Big(E\big(|\mathcal{P}^+(A)|\,\big|\,A\big)\Big)\,. \tag{50.26}$$

This inequality has the same structure as the deviations inequality (40.5). In fact, it is an improvement of (40.5) in the uniform case $p=1/2$. We show next how an inequality similar to (40.5) can be derived from the inequality (50.26). We consider only the case where $P(A)<1/2$, or equivalently, $N<2^{n-1}$. In this case, the integer $K$ associated to $N$ through the inequalities (50.13) is such that $N(K)<2^{n-1}$. Let us consider the
• If $n$ is odd, say $n=2m-1$, then $N(m)=2^{n-1}$, thus $K>m=(n+1)/2$.
• If $n$ is even, say $n=2m$, then $N(m)>2^{n-1}$, thus $K>m=n/2$.
In both cases, we have the inequality

$$K\geq \frac{n}{2}+1 \quad\Longleftrightarrow\quad 2K-2\geq n\,,$$

and it follows from formula (50.14) that

$$\psi_1(N) \,=\, \frac{2}{N}\sum_{K\leq k\leq n}(k-K+1)\binom{n}{k}+2K-2-n\,.$$

Let $L$ be an integer such that $K\leq L\leq n$, to be chosen later. We split the sum in two and we bound separately each term as follows:

$$\begin{aligned}
\psi_1(N) \,&=\, \frac{2}{N}\Big(\sum_{K\leq k<L}\cdots+\sum_{L\leq k\leq n}\cdots\Big)+2K-2-n\\
&\leq\, \frac{2}{N}\Big((L-K)N(K)+(n-K+1)N(L)\Big)+2K-2-n\\
&\leq\, 2(L-K)+\frac{n}{N}N(L)+2K-2-n\\
&\leq\, 2L-n+\frac{n}{N}N(L)\,. \qquad (50.27)
\end{aligned}$$

So far we have supposed that $L\geq K$. However, for $L$ such that $n/2\leq L<K$, we have

$$2L-n+\frac{n}{N}N(L)\geq\, n\,.$$

Since in addition $M_1(n,N) \leq nN$ and $\psi_1(N) \leq n$, we see that the inequality (50.27) still holds. We conclude therefore that

$$\forall L \in \left\{\,\frac{n}{2}, \cdots, n\,\right\} \qquad \psi_1(N) \,\leq\, 2L - n + \frac{n}{N} N(L)\,. \tag{50.28}$$

It remains now to choose the value of $L$ in order to obtain an interesting upper bound. To that end, we rewrite $N(L)$ as the tail of a binomial random variable:

$$N(L) \,=\, \sum_{L \leq k \leq n} \binom{n}{k} \,=\, 2^n P\Big(B\big(n, \frac{1}{2}\big) \geq L\Big)\,,$$

where $B\big(n, \frac{1}{2}\big)$ is a standard binomial random variable of parameters $n, 1/2$. We shall make appeal to Hoeffding's inequality to bound $N(L)$:

$$\forall L \in \left\{\,\frac{n}{2}, \cdots, n\,\right\} \qquad N(L) \,\leq\, 2^n \exp\Big(-\frac{(2L-n)^2}{2n}\Big)\,. \tag{50.29}$$

Notice that we use the very same inequality in order to derive the deviations inequality (40.5). So we substitute (50.29) into (50.28), and we optimize over $L$. After a few tries, a reasonable choice is to take $L$ such that

$$2L - n \,=\, \sqrt{2n}\sqrt{n \ln 2 + \ln n - \ln N - \frac{1}{2} \ln\big(2n(n \ln 2 - \ln N)\big)}\,. \tag{50.30}$$

Notice that, since $N < 2^{n-1}$, then

$$2L - n \geq \sqrt{2n}\sqrt{\ln 2 + \ln n - \frac{1}{2}\ln\big(2n(n\ln 2)\big)} \,\geq\, \sqrt{\ln\Big(\frac{2}{\ln 2}\Big)} \,>\, 1\,.$$

Substituting the value of $L$ given by (50.30) into inequalities (50.29) and (50.28), we obtain

$$\psi_1(N) \,\leq\, 2\sqrt{2n}\sqrt{n \ln 2 + \ln n - \ln N}\,. \tag{50.31}$$

The inequality (50.31) yields the following upper bound on $\psi_1^{-1}$:

$$\forall \lambda \in [0,n] \qquad \psi_1^{-1}(\lambda) \,\leq\, 2^n n \exp\Big(-\frac{\lambda^2}{8n}\Big)\,. \tag{50.32}$$

Substituting (50.32) into the inequality (50.26), we get finally

$$P(A) \,\leq\, n \exp\Big(-\frac{1}{8n}\Big(E\big(|\mathcal{P}^+(A)|\,\big|\,A\big)\Big)^2\Big)\,.$$

Thus we have recovered an inequality weaker than the inequality (40.5) in the symmetric case (because of the factor $n$ in front of the exponential). Moreover we are confident that our method of proof could be adapted to treat the non symmetric case $p \neq 1/2$ as well. In any case, the proof presented here is more powerful than the one used in subsection 40.1. Indeed, we have obtained a lot of information on the optimal sets, and we believe that the inequality (50.26) is stronger than (40.5).

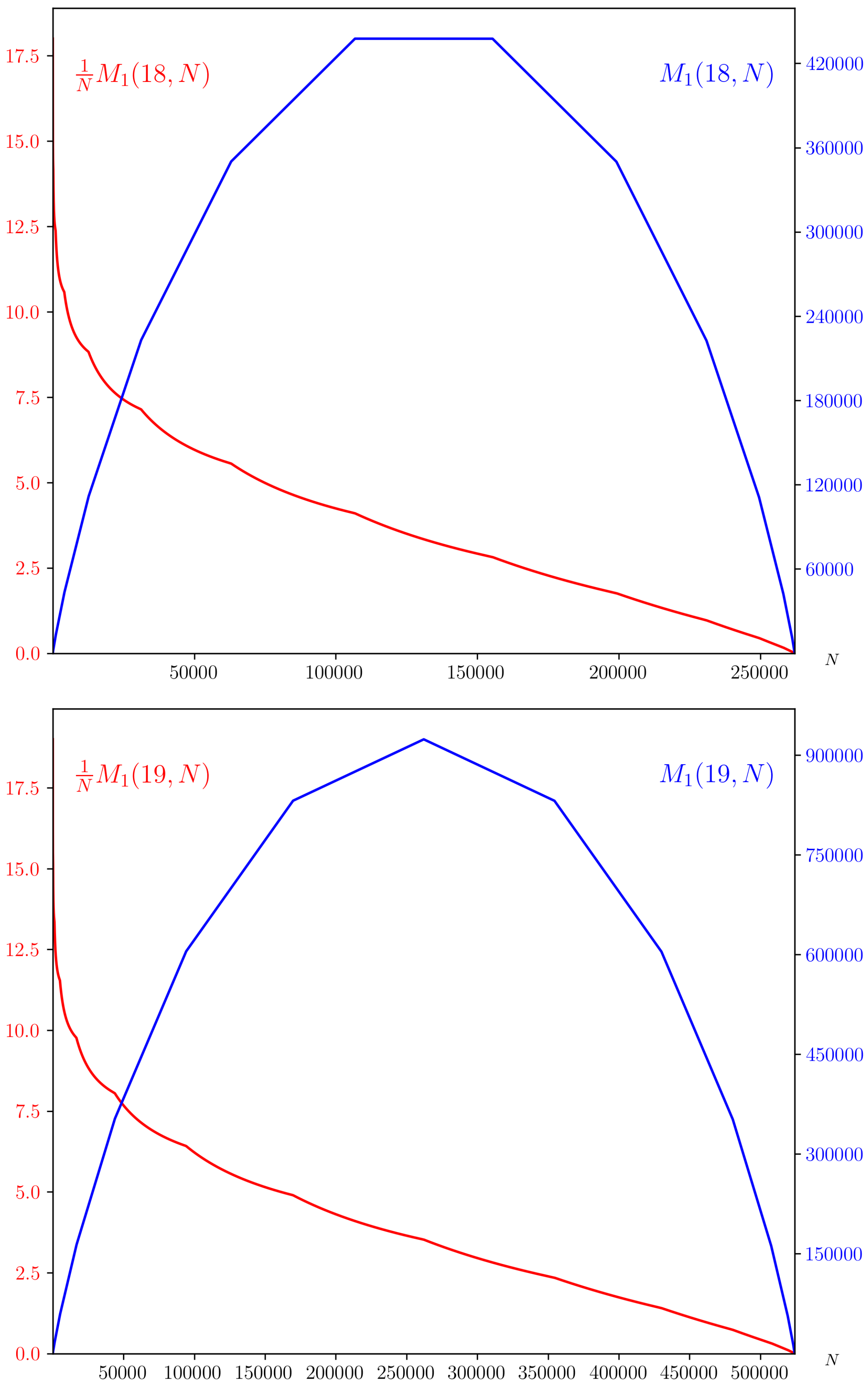

Figure 113: The graphs of $M_1(n,N)$ and $\frac{1}{N}M_1(n,N)$ for $n = 18, 19$

## 51 Maximizers of $E(|\mathcal{P}_2^+ \setminus \mathcal{P}_1^+|)$

We attack here the difficult problem of controlling the positive pivots of second order and we try to implement the same strategy as for the pivots of first order. So we first formulate the relevant extremal problem. For $A$ an increasing subset of $\{0,1\}^n$, we define

$$\phi_{2\setminus 1}(A) \;=\; \sum_{\omega \in A} \big|\mathcal{P}_2^+(A,\omega) \setminus \mathcal{P}_1^+(A,\omega)\big|\,.$$

We consider the problem

$$M_{2\setminus 1}(n) \;=\; \max\,\big\{\,\phi_{2\setminus 1}(A) : A \text{ increasing subset of } \{0,1\}^n\,\big\}\,. \tag{51.1}$$

Contrary to the problem (50.1), this new problem is extremely complicated and there is no obvious solution. In a first approach, we resort to computer-assisted computations. We performed an exhaustive analysis of all the monotone boolean functions in 7 or less variables. The table 4 presents the maximizers, which are also depicted in layer representation in figures 114 and 115. From $n=4$ onwards, the maximizer is never unique, and moreover it seems difficult to infer a general result from these first cases. So before getting involved in mathematical proofs, we consider a variant of the extremal problem (51.1).

| $n$ | $\|A\|$ | $E(\|\mathcal{P}_1^+\|)$ | $E(\|\mathcal{P}_2^+ \setminus \mathcal{P}_1^+\|)$ | Minimal$(A)$ |
|---|---|---|---|---|
| 2 | 3 | 2 | 2 | 10 01 |
| 3 | 7 | 3 | 6 | 100 010 001 |
| 4 | 13 | 8 | 12 | 1000 0100 0011 |
| | 14 | 6 | 12 | 1000 0100 0010 |
| | 9 | 12 | 12 | 1100 1010 0101 0011 |
| | 12 | 10 | 12 | 1000 0110 0101 0011 |
| | 10 | 12 | 12 | 1100 1010 1001 0110 0101 |
| | 15 | 4 | 12 | 1000 0100 0010 0001 |
| | 11 | 12 | 12 | 1100 1010 1001 0110 0101 0011 |
| 5 | 22 | 24 | 38 | 11000 10100 10010 01110 01001 00101 00011 |
| 6 | 51 | 42 | 96 | 110000 101000 100100 100011 011100 010010 010001 001010 001001 000110 000101 |
| 7 | 111 | 69 | 200 | 1100000 1010000 1001000 1000100 1000011 0111000 0110100 0101100 0100010 0100001 0011100 0010010 0010001 0001010 0001001 0000110 0000101 |

Table 4: The maximizers of $E(|\mathcal{P}_2^+ \setminus \mathcal{P}_1^+|)$ for $n=2,3,4,5,6,7$

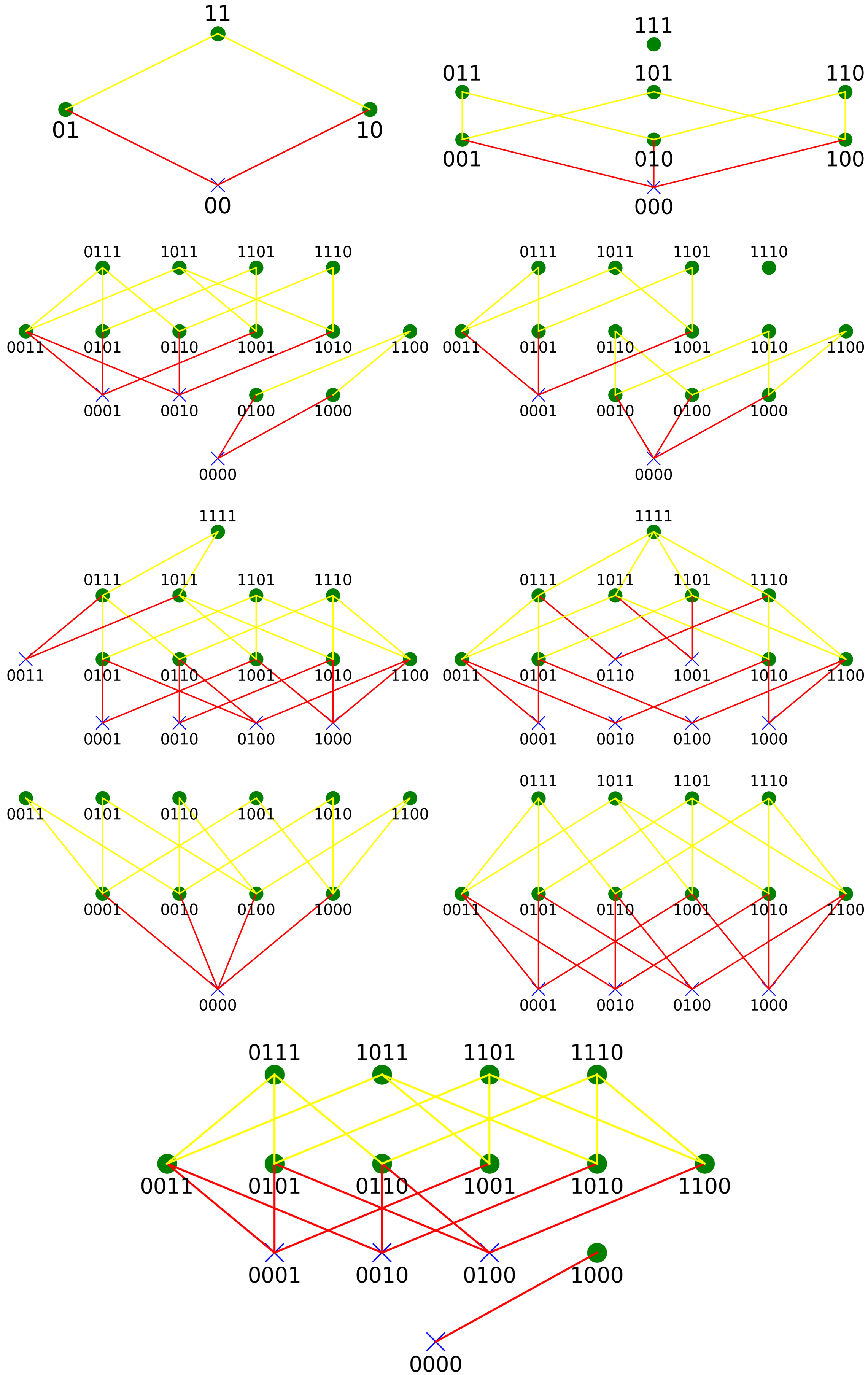

Figure 114: The maximizers of $E(|\mathcal{P}_2^+ \setminus \mathcal{P}_1^+|)$ for $n = 2, 3, 4$

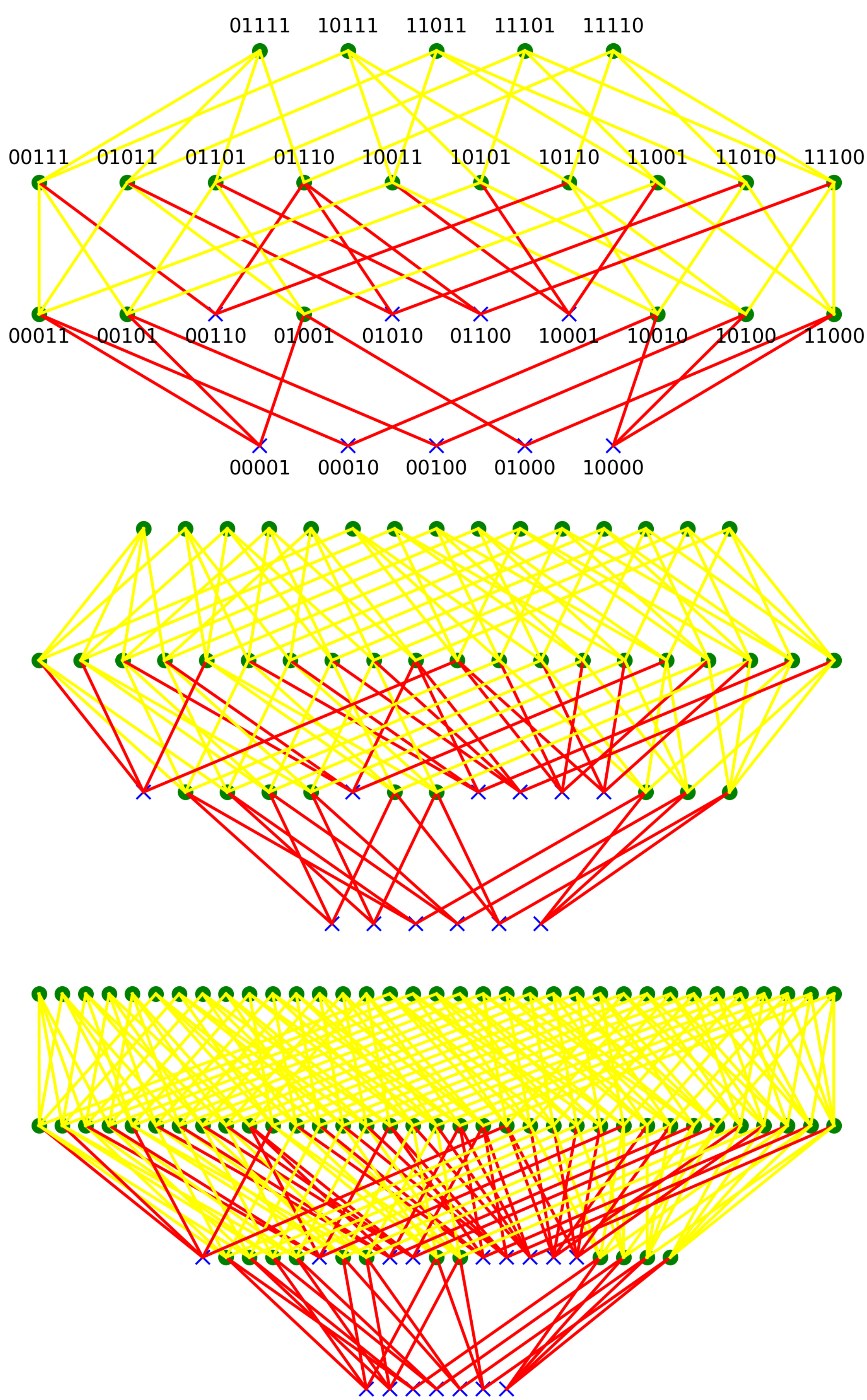

Figure 115: The unique maximizers of $E(|\mathcal{P}_2^+ \setminus \mathcal{P}_1^+|)$ for $n = 5, 6, 7$

# 52 Maximizers of $E(|\mathcal{P}_2^+|)$

Since we have a satisfactory control of the pivots of first order, we can as well try to control directly the whole set of the positive pivots of second order. So we modify the extremal problem (51.1) of the previous section in order to include as well the pivots of first order. For $A$ an increasing subset of $\{0,1\}^n$, we define

$$\phi_2(A) \;=\; \sum_{\omega \in A} \left|\mathcal{P}_2^+(A,\omega)\right|.$$

We consider the extremal problem

$$M_2(n) \;=\; \max\,\Big\{\,\phi_2(A) : A \text{ increasing subset of } \{0,1\}^n\,\Big\}. \qquad (52.1)$$

Again, we performed an exhaustive analysis of all the monotone boolean functions in 7 or less variables. The table 5 presents the maximizers, which are also depicted in layer representation in figures 116 and 117. There is a unique maximizer whenever $n \le 6$, but surprisingly, there are two maximizers for $n = 7$! Anyway, this problem seems to be better posed than the previous problem (51.1), so we shall focus on this one from now onwards. Our next goal would be to conduct a rigorous analysis leading to a useful result on the value of $M_2(n)$, and on the maximizers of the problem (52.1).

| $n$ | $\|A\|$ | $E(\|\mathcal{P}_1^+\|)$ | $E(\|\mathcal{P}_2^+\|)$ | Minimal$(A)$ |
|---|---|---|---|---|
| 2 | 3 | 2 | 4 | 10 01 |
| 3 | 4 | 6 | 9 | 110 101 011 |
| | 7 | 3 | 9 | 100 010 001 |
| 4 | 9 | 12 | 24 | 1100 1010 0101 0011 |
| | 10 | 12 | 24 | 1100 1010 1001 0110 0101 |
| | 11 | 12 | 24 | 1100 1010 1001 0110 0101 0011 |
| 5 | 22 | 24 | 62 | 11000 10100 10010 01110 01001 00101 00011 |
| 6 | 51 | 42 | 138 | 110000 101000 100100 100011 011100 010010 010001 001010 001001 000110 000101 |
| 7 | 87 | 117 | 314 | 1110000 1101000 1100100 1100010 1011000 1010100 1010010 1001100 1001010 1000101 1000011 0111100 0111010 0110001 0101001 0100110 0100101 0100011 0011001 0010110 0010101 0010011 0001110 0001101 0001011 |
| | 88 | 116 | 314 | 1110000 1101000 1100100 1100010 1011000 1010100 1010010 1001001 1000110 1000101 1000011 0111000 0110110 0110001 0101100 0101010 0100101 0100011 0011100 0011010 0010101 0010011 0001110 0001101 0001011 |

Table 5: The maximizers of $E(|\mathcal{P}_2^+|)$ for $n = 2,3,4,5,6,7$

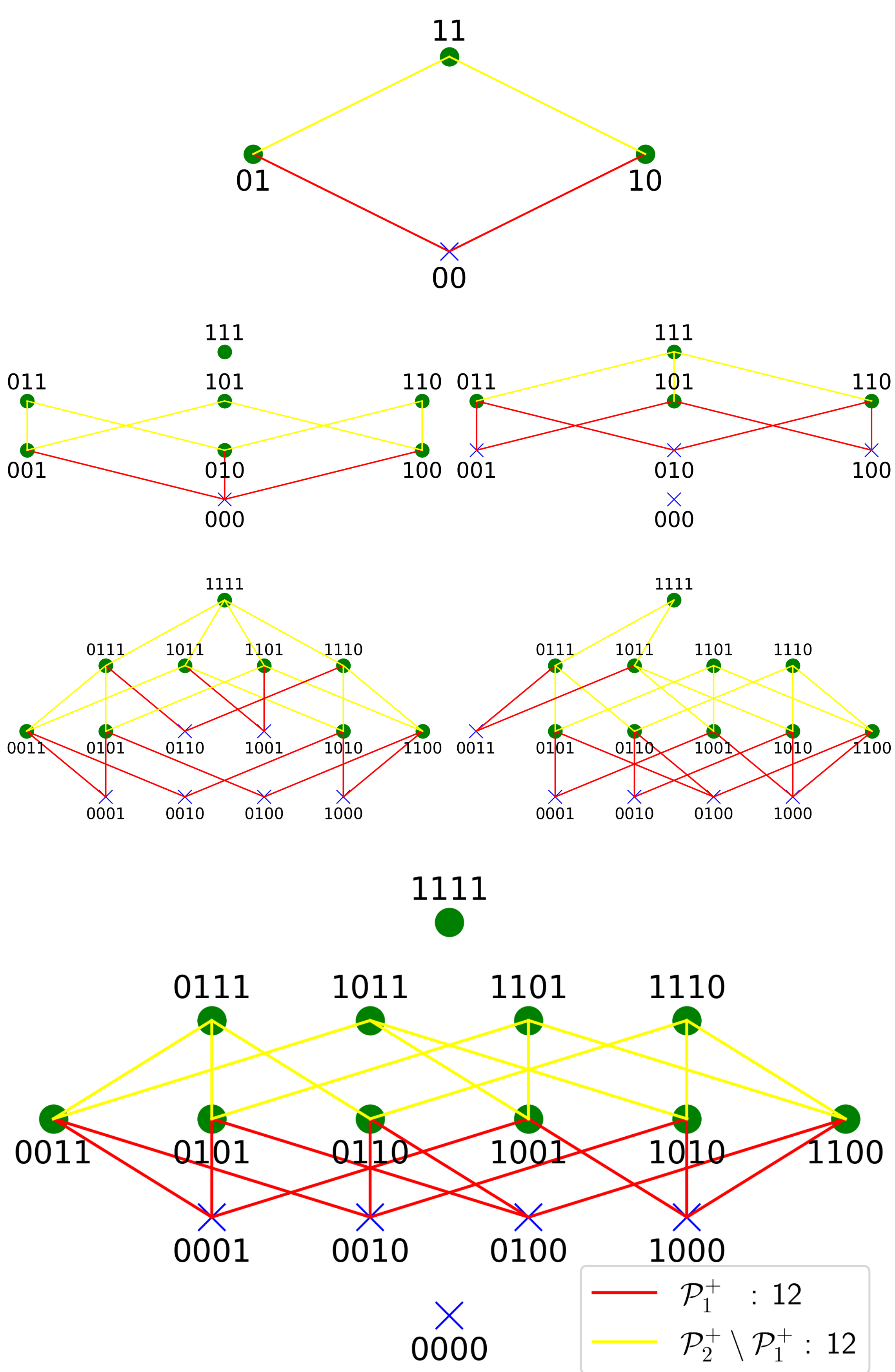

Figure 116: The maximizers of $E(|\mathcal{P}_2^+|)$ for $n = 2, 3, 4$

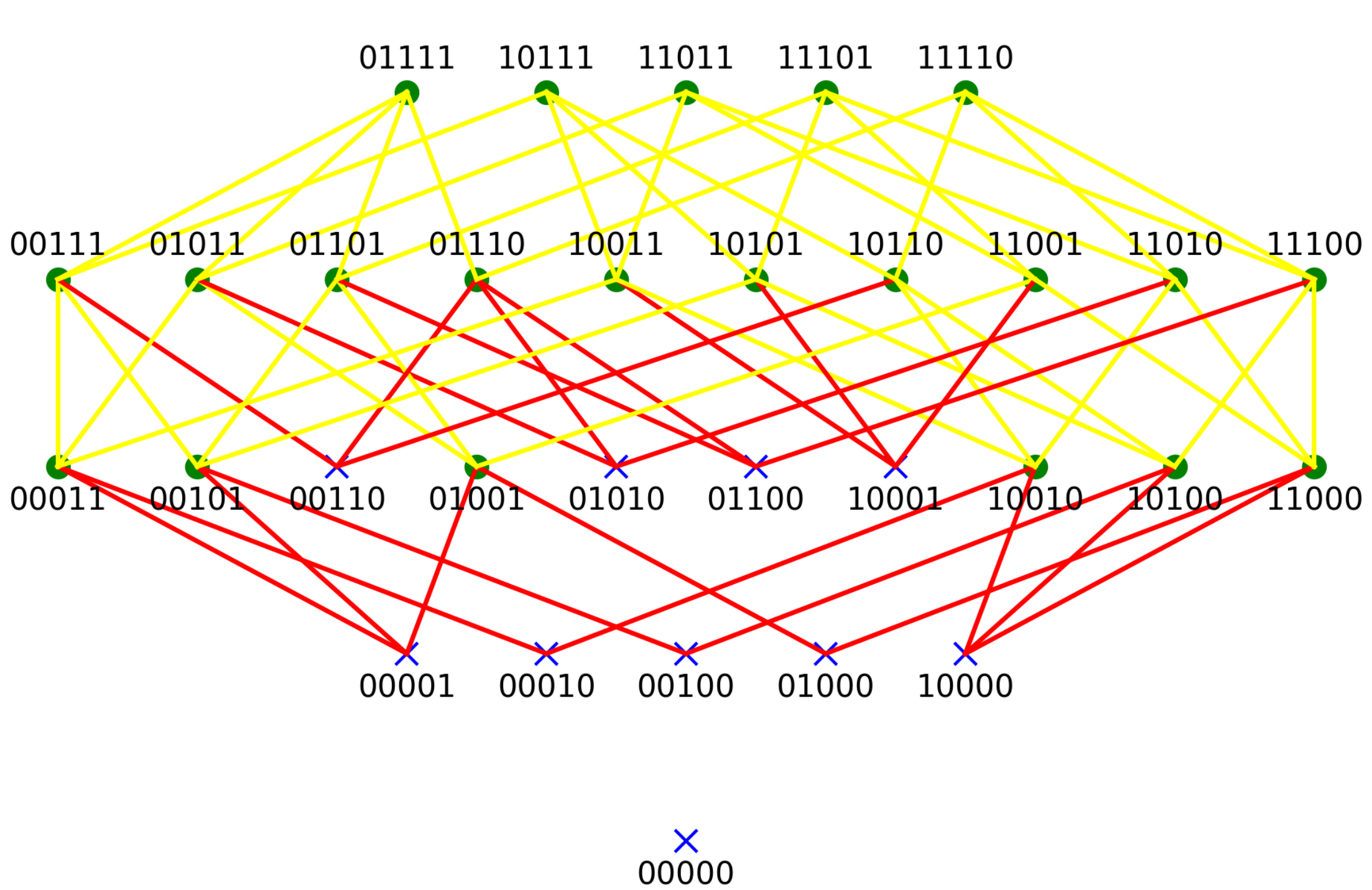

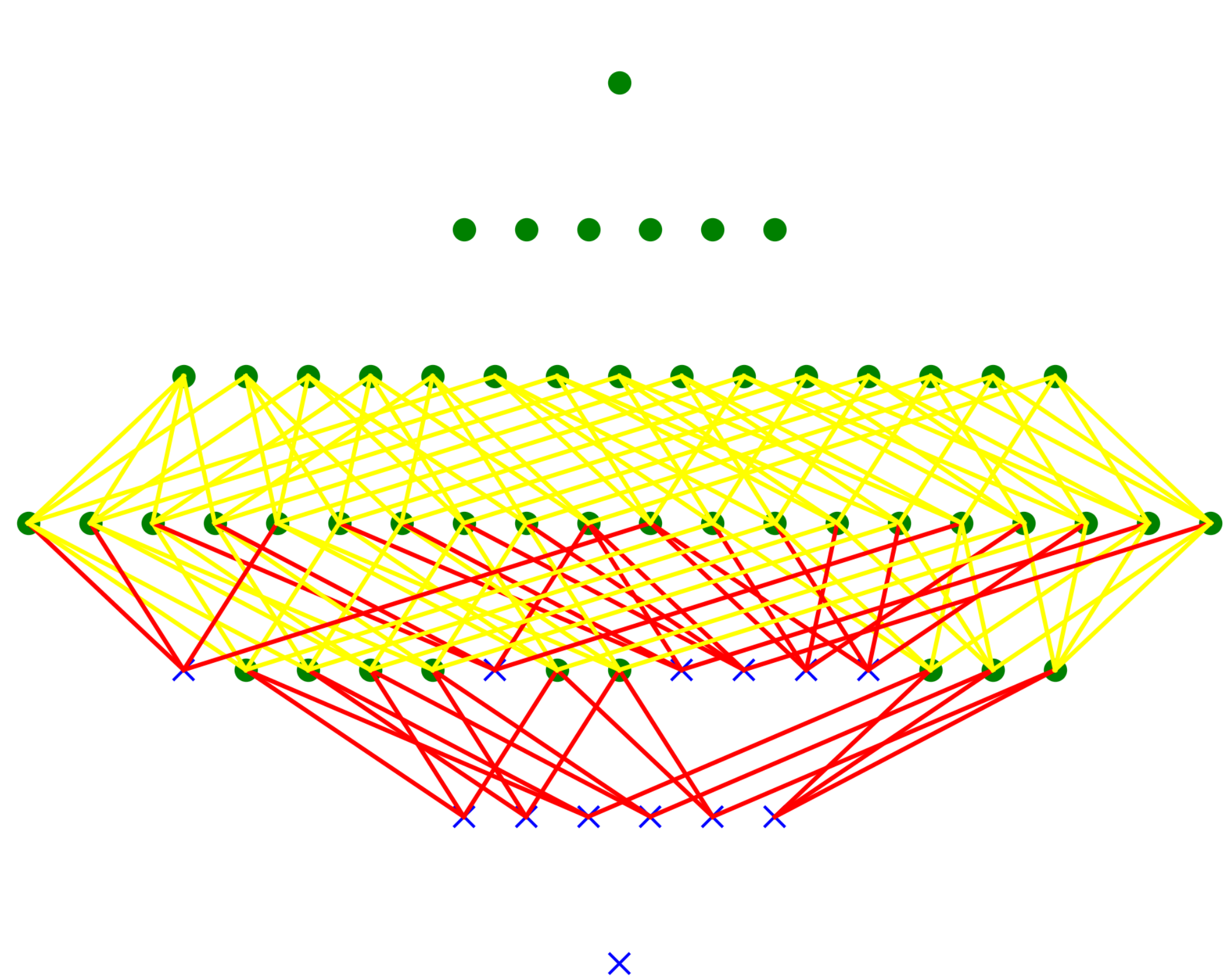

Figure 117: The unique maximizers of $E(|\mathcal{P}_2^+|)$ for $n = 5, 6$

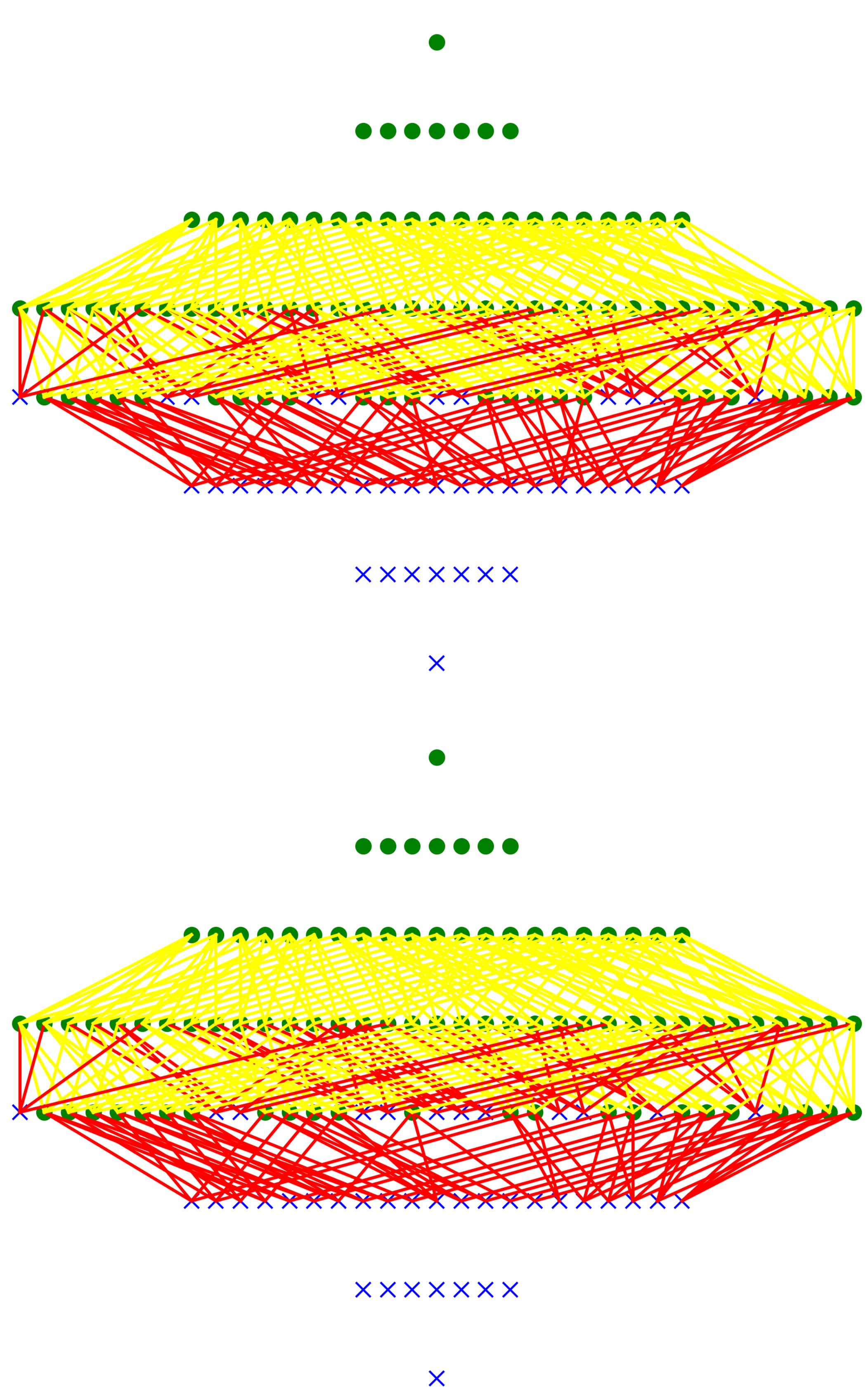

Figure 118: The two maximizers of $E(|\mathcal{P}_2^+|)$ for $n = 7$

## 53 A conjecture on the maxima of $\phi_{2\setminus 1}$ and $\phi_2$

We recall the definition of the maps $\phi_{2\setminus 1}$ and $\phi_2$ that we try to maximize. For a subset $A$ of $\{0,1\}^n$, they are given by

$$\phi_{2\setminus 1}(A) \,=\, \sum_{\omega\in A} \big|\mathcal{P}_2^+(A,\omega)\setminus\mathcal{P}_1^+(A,\omega)\big|\,,$$
$$\phi_2(A) \,=\, \sum_{\omega\in A} \big|\mathcal{P}_2^+(A,\omega)\big|\,.$$

The structure of the sets $A$ maximizing $\phi_{2\setminus 1}$ or $\phi_2$ is quite complex even for $n\leq 7$, as we saw in the two previous sections. In view of the results obtained from the exhaustive searches for $n\leq 7$, we formulate the following conjecture.

**Conjecture 53.1.** *Let $n\geq 1$ be fixed. If $A^*$ is an increasing subset of $\{0,1\}^n$ which realizes the maximum of $\phi_{2\setminus 1}$ or of $\phi_2$, then*

$$\Big\{\,\omega: \sum_{1\leq i\leq n}\omega(i)\geq\frac{n}{2}+1\Big\}\,\subset\, A^*\,\subset\,\Big\{\,\omega: \sum_{1\leq i\leq n}\omega(i)\geq\frac{n}{2}-2\Big\}\,. \tag{53.1}$$

The point of this conjecture is that it yields a valuable upper bound on the maximal values $M_{2\setminus 1}(n)$ and $M_2(n)$ of $\phi_{2\setminus 1}$ and $\phi_2$ (as defined in formulas (51.1) and (52.1)). Indeed, it would follow from the conjecture that

$$\begin{aligned}\max\big(M_{2\setminus 1}(n),M_2(n)\big)\,&\leq\, n\big(\frac{n}{2}+1\big)\Big|\Big\{\,\omega:\frac{n}{2}-2\leq\sum_{1\leq i\leq n}\omega(i)\leq\frac{n}{2}+1\Big\}\Big|\\ &\leq\,4\big(\frac{n}{2}+1\big)\binom{n}{\lfloor n/2\rfloor}\,\leq\,\frac{2\sqrt{2}(n+2)}{\sqrt{\pi n}}2^n\,\leq\,6\sqrt{n}2^n\,,\end{aligned} \tag{53.2}$$

where we have used the following classical bound on the central binomial coefficient (see [100], chapter 10, lemma 7):

$$\forall n\geq 1\qquad \binom{n}{\lfloor n/2\rfloor}\,\leq\,\sqrt{\frac{2}{\pi n}}2^n\,.$$

We would therefore conclude from (53.2) that

$$\forall A\subset\{0,1\}^n\quad A\text{ increasing}\quad\Rightarrow\quad E\big(|\mathcal{P}_2^+(A)|\big)\,\leq\,6\sqrt{n}\,.$$

The conjecture (53.1) holds for $n\leq 7$. Let us try to attack it for larger values of $n$. We had two proofs for the inequality (47.1) dealing with the first pivots. The first proof relied on the Margulis-Russo formula. Unfortunately, all our efforts to obtain a corresponding representation formula for the higher pivots have failed. A natural attempt consists in looking at the successive derivatives with respect to the parameter $p$. We have already tried very hard this strategy in section 41, with limited success. The second proof was presented in subsection 50.1, it consisted in looking at perturbations of a candidate set $A$ and their effect on the value of the functional $\phi_1$. So we try next this strategy and we consider first perturbations which consist in removing a minimal configuration of the candidate set. The computation done for $\phi_1$ shows that the removal of a minimal element $\sigma$ of a set $A$ creates a variation of $\phi_1$ equal to $\Delta\phi_1 = n-2|\operatorname{supp}\sigma|$. This variation is positive as long as $|\operatorname{supp}\sigma| < n/2$. Unfortunately, the computation of the variations of $\phi_{2\setminus 1}$ and $\phi_2$ are complicated and they did not lead to any definite conclusion. The reason is that the first pivots created by the removal of $\sigma$ might correspond to second pivots which were already present in $A$. Thus we resort to the computer. First we check throughout our database of monotone functions what happens whenever we remove a minimal element belonging to the lowest level of the set. For $n \leq 6$, the removal of such a minimal element strictly below the middle level never leads to a decrease of $\phi_2$. Unfortunately, for $n = 7$, this is not any more the case. Let us consider the increasing set $A$ whose minimal elements are

1000010 0100010 0010010 0001010 0000110 1000001 0100001
0010001 0001001 0000101 1110000 1101000 1011000 0111000
1100100 1010100 0110100 1001100 0101100 0011100 .

The removal of the configuration $\sigma = 0000101$ leads to a decrease of 7 for $\phi_2$:

$$\phi_2(A \setminus \{\,\sigma\,\}) - \phi_2(A) \;=\; -7\,, \qquad |\operatorname{supp}\sigma| \;=\; 2 \;<\; \frac{7}{2}\,.$$

This modification is illustrated on figure 119. In this figure and the following, we adopt the layer representation of the set $A$. The interesting layers of the set $A$ are presented on top, the modifications appear in the middle, and they are redeployed for clarity in the bottom. The first pivots are depicted by red edges and the second pivots by yellow edges. The modifications of the pivots induced by the removal of $\sigma$ are drawn in black in the top picture. The middle picture focuses on these modifications. The dashed edges correspond to pivots which are discarded by the modification. The continuous edges correspond to pivots which are present after the modification, but which have undergone a change of status. The orange edges correspond to pivots which were of second order before the modification and which become of first order after the modification. The configuration $\sigma$ which is removed is marked by a green disc superposed with a blue cross. In this example, with $n = 7$, the element to be removed is on the second level. So this is problematic if we wish to use this method to prove the conjecture, because we might stay stuck with a configuration containing elements of the second level. However, we are still quite close to the third and fourth layers, so maybe we should not be discouraged by this first example. Indeed, we would be very happy if the trouble occurred only when trying to remove configurations belonging to the level $\lfloor n/2 \rfloor - 1$. There is also a good news with this example. In fact, the removal of the configuration 1110000 leads to an increase of 3 for $\phi_2$. The effect of this removal is presented in figure 120. Thus we see that the removal of a minimal configuration of the lowest layer might lead to a decrease of $\phi_2$, and at the same time, some other minimal configuration might lead to an increase of $\phi_2$!

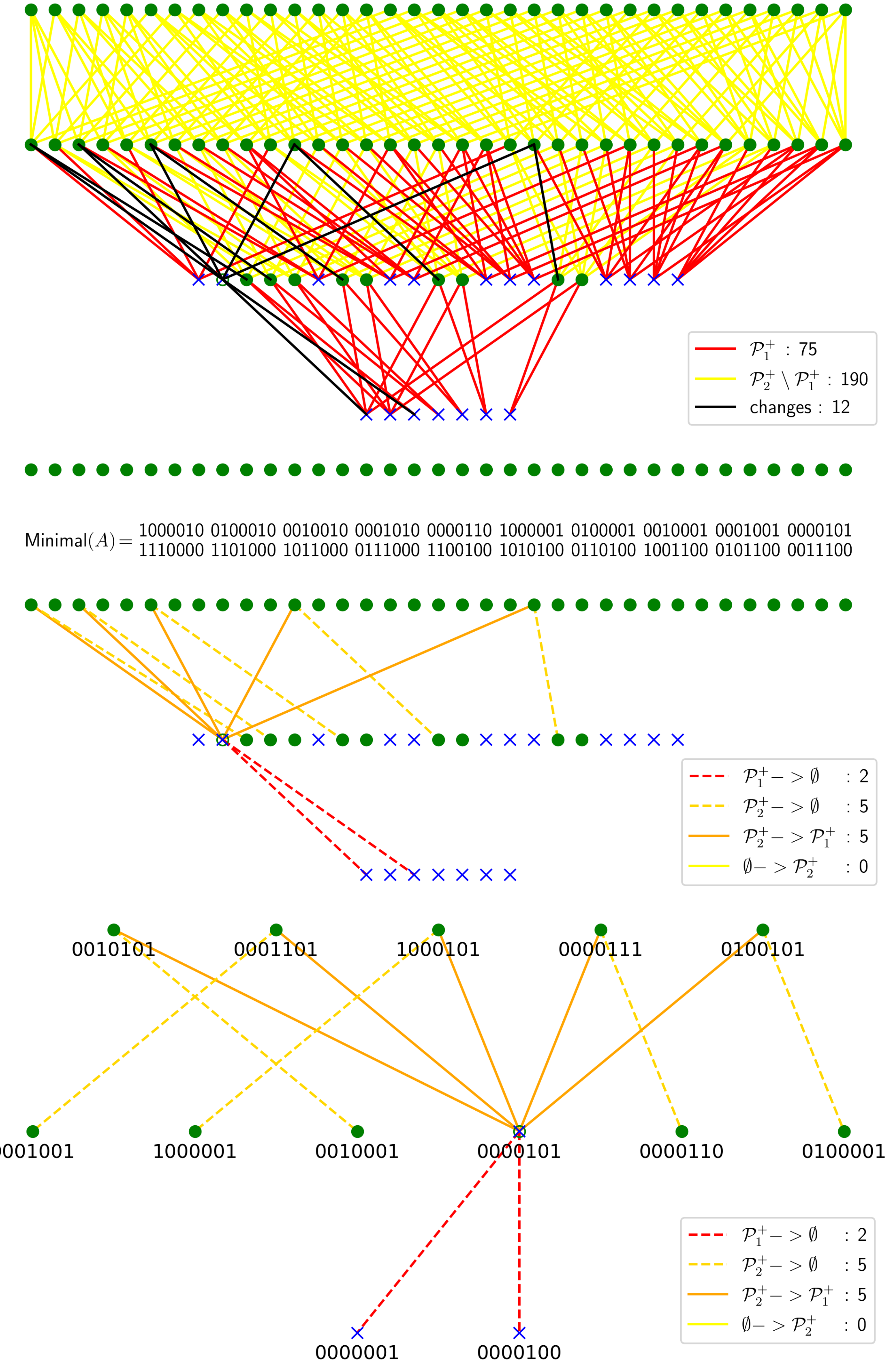

Figure 119: Removal of the configuration 0000101: $\Delta\phi_2 = -7$.
Top: the set $A$. Middle: the modifications (same representation as top).
Bottom: focus on the modifications (redeployed horizontally for clarity).

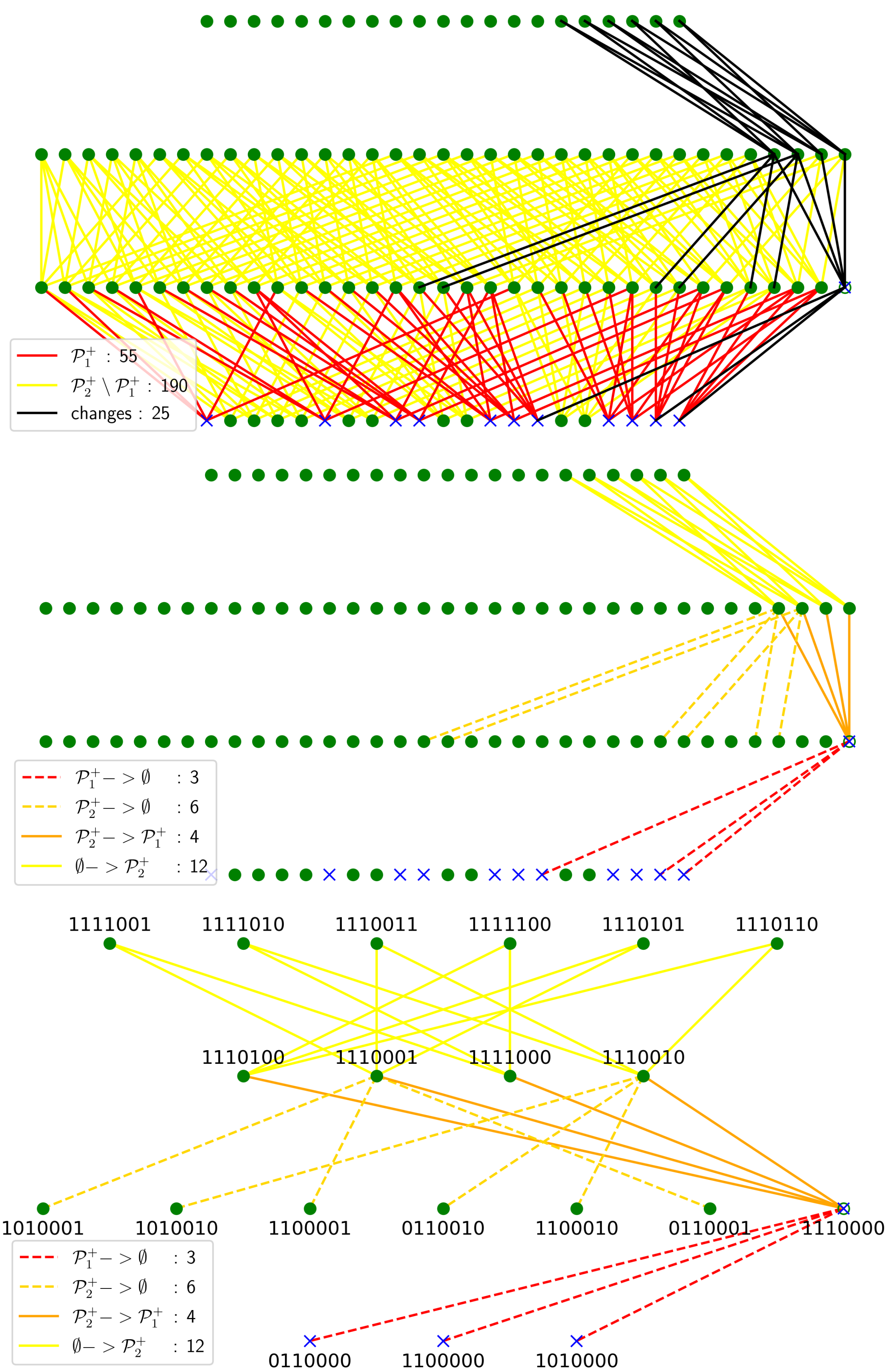

Figure 120: Removal of the configuration 1110000: $\Delta\phi_2 = +3$.
Top: the set $A$. Middle: the modifications (same representation as top).
Bottom: focus on the modifications (redeployed horizontally for clarity).

# 54 Gradient algorithm

In order to localize the maxima of $\phi_2$, we implemented a gradient algorithm on the increasing subsets of the hypercube. The gradient algorithm can be used in practice to find some local maxima through numerical computations for low values of $n$, but also for the theoretical investigation of $\phi_2$ for any $n \geq 1$. We start from the full set $A_0$, or any increasing set $A_0$. Suppose we have a mechanism for modifying an increasing set $A$. More precisely, for each increasing subset $A$, we define a set $\mathcal{V}(A)$, which we think of as its neighbours in some sense. We describe next a standard gradient algorithm for $\phi_2$ associated to this modification mechanism.

**Standard gradient algorithm**. The algorithm starts from a set $A_0$. At each iteration, the algorithm updates a candidate set $A_k$. Let $k \geq 0$, suppose that the set $A_k$ has been built and let us explain how to build the set $A_{k+1}$. We examine the neighborhood $\mathcal{V}(A_k)$ of $A_k$. If there exists $A$ in $\mathcal{V}(A_k)$ such that $\phi_2(A) > \phi_2(A_k)$, then we pick up one such set $A$ (according to some rule if there are several of them) and we set $A_{k+1} = A$. Otherwise, the algorithm terminates at step $r = k$. The output of the algorithm is a sequence $A_0, \dots, A_r$ of subsets of $\{0,1\}^n$ such that

$$\forall k \in \{\, 0 \dots, r-1 \,\} \qquad A_{k+1} \in \mathcal{V}(A_k)\,, \qquad \phi_2(A_k) \,<\, \phi_2(A_{k+1})\,,$$

and the final set $A_r$ is a local maximum, i.e.,

$$\forall A \in \mathcal{V}(A_r) \qquad \phi_2(A) \,\leq\, \phi_2(A_r)\,.$$

Since the function $\phi_2$ is strictly increasing along the sequence of sets $A_0, A_1, \dots,$ then certainly the algorithm terminates after a finite number of iterations. The pseudocode of this exploration algorithm is given in Algorithm 54.1. The crucial problem is the design of the modification mechanism. This design aims at fulfilling two conflicting goals:

• First goal: avoiding local maxima which are not global maxima. We would like the function $\phi_2$ to have as few local maxima as possible, and ideally none that are not a global maximum. This would guarantee that the output of the algorithm is always a global maximum.

• Second goal: expressing the variation of $\phi_2$. We need that the computation of the variation of $\phi_2$ between two neighbours is fast (for the numerical computations) and understandable in terms of the geometry of $A$ (to be able to derive interesting theoretical results).

A naive choice fulfilling the first goal is to take the whole collection of the increasing subsets of $\{0,1\}^n$ to be the neighborhood of each set. In this case, a single step of the gradient algorithm amounts to an exhaustive search of the whole space! More seriously, if we focus on the second goal, then we would start with mechanisms which modify the set $A$ as little as possible. The smallest possible modification of a set $A$ consists in removing or adding a configuration to $A$. For the resulting set to be still increasing, we can only remove a minimal configuration of $A$ or add a maximal configuration of $A^c$. We discuss these two operations in the next two subsections.

**Algorithm 54.1** Standard gradient algorithm

```
A(0) ← an initial increasing subset of {0,1}^n, k ← 0, Failure ← false
repeat
    Compute the maximum M_2 of φ_2 over V(A(k))
    if M_2 > φ_2(A(k)) then
        Pick up A in V(A(k)) such that φ_2(A) = M_2
        k ← k + 1, A(k) ← A
    else
        Failure ← true
    end if
until Failure
return A(k)
```

## 54.1 Removing a minimal element

We examine first the removal operation and so in our first attempt we consider the system of neighborhoods defined by

$$\mathcal{V}_{\text{remove}}(A) \,=\, \big\{\, A \setminus \{\,\sigma\,\} : \sigma \in \text{Minimal}(A) \,\big\}\,.$$

We could also consider a gradient algorithm associated to the function $\phi_{2/1}$. However, since the gradient algorithm for $\phi_1$ is successful, it is certainly better to try the gradient algorithm for $\phi_2 = \phi_1 + \phi_{2/1}$ rather than for $\phi_{2/1}$. Indeed, the removal of a minimal configuration in the lower part of the hypercube always leads to an increase of $\phi_1$, and this can only help to increase $\phi_2$ as well.

### 54.1.1 Numerical computations

We launched this gradient exploration with $n = 11$ starting from $A_0 = \{0,1\}^{11}$. For the exploration step, we try first the elements of the neighborhood which have the smallest number of ones (because we hoped that it is always possible to leave a set which is not a maximum with such a move), and we order them according to the lexicographic order. We discovered that, for the increasing set whose minimal elements are listed in table 6, the removal of the configuration 10000000010 leads to a decrease of 2 for $\phi_2$, see figure 121. In this case, we are quite far from the middle layers $5, 6$. Thus it appears that the strategy consisting of removing exclusively minimal elements belonging to the lowest level will not be successful. The good news is again that, in this example, the removal of the minimal configuration 00000011100 leads to an increase of $\phi_2$ of 47. This is shown in figure 122. Still for $n = 11$, the final configuration of the gradient algorithm is the increasing set $A$ whose minimal elements are listed in table 8. We observe with pleasure that this final configuration satisfies the conclusion (53.1) of the conjecture 53.1. In fact, we tried the gradient algorithm described previously, and a little variant in which a removal is accepted whenever it does not decrease the value of $\phi_2$ (whereas in the initial algorithm a removal was accepted when it increases the value of $\phi_2$). Let us call the first algorithm the strict gradient-remove and the second algorithm the loose gradient-remove.

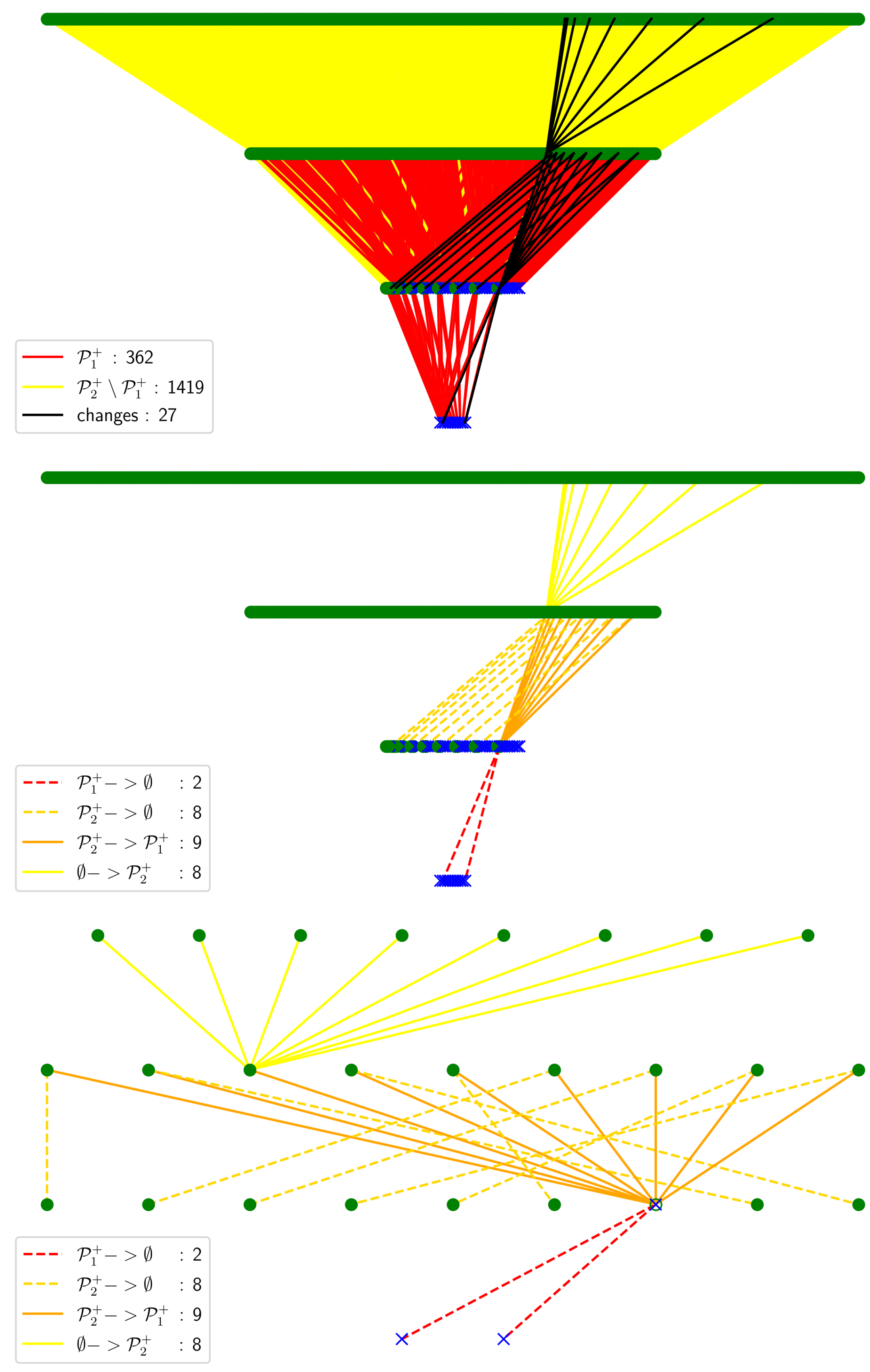

Figure 121: Removal of the configuration 10000000010: $\Delta\phi_2 = -2$.
Top: the set $A$. Middle: the modifications (same representation as top).
Bottom: focus on the modifications (redeployed horizontally for clarity).

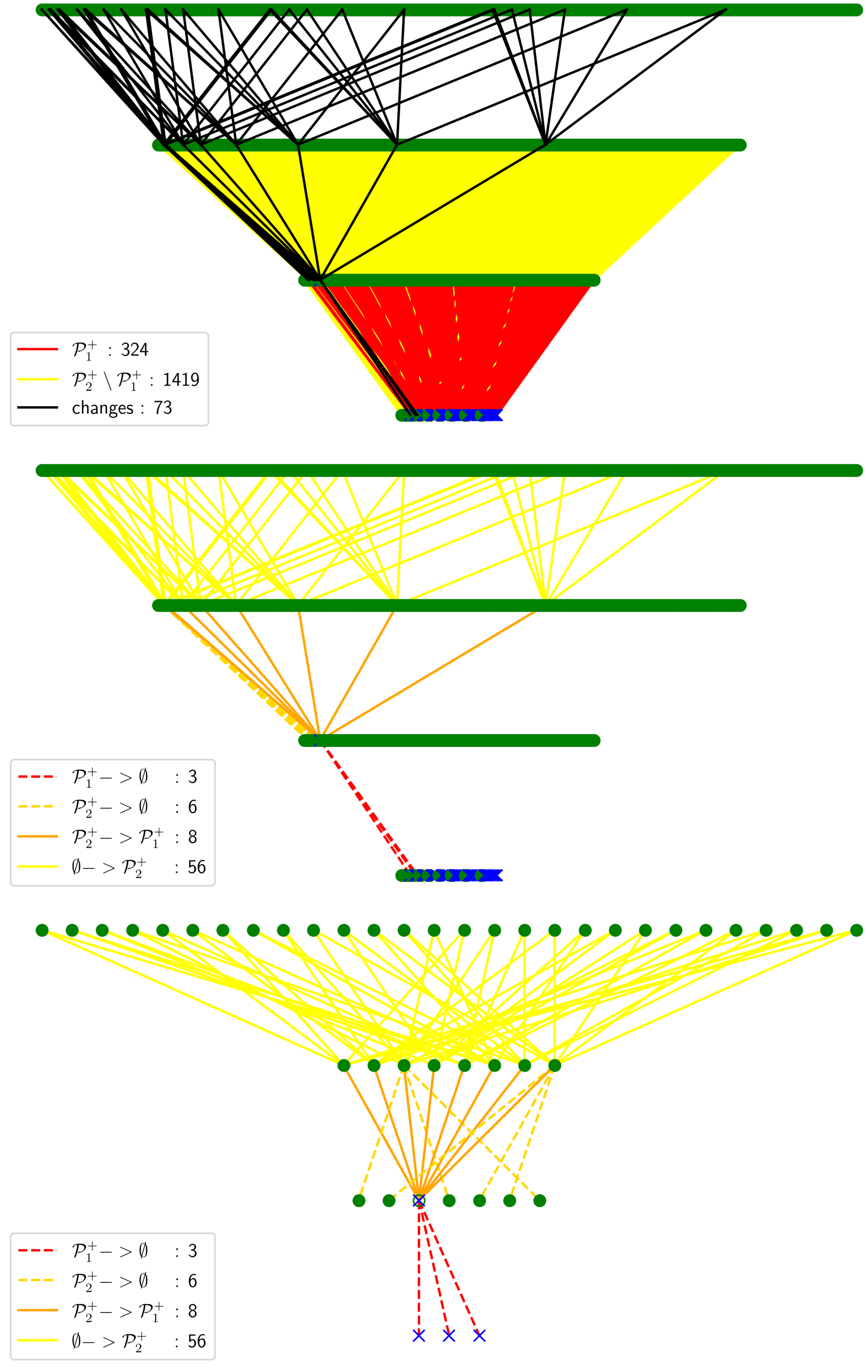

Figure 122: Removal of the configuration 00000011100: $\Delta\phi_2 = +47$.
Top: the set $A$. Middle: the modifications (same representation as top).
Bottom: focus on the modifications (redeployed horizontally for clarity).

| 10000000010 | 01000000010 | 00100000010 | 00010000010 | 00001000010 |
|---|---|---|---|---|
| 00000100010 | 00000010010 | 00000001010 | 00000000110 | 10000000001 |
| 01000000001 | 00100000001 | 00010000001 | 00001000001 | 00000100001 |
| 00000010001 | 00000001001 | 00000000101 | 00000000011 | 11100000000 |
| 11010000000 | 10110000000 | 01110000000 | 11001000000 | 10101000000 |
| 01101000000 | 10011000000 | 01011000000 | 00111000000 | 11000100000 |
| 10100100000 | 01100100000 | 10010100000 | 01010100000 | 00110100000 |
| 10001100000 | 01001100000 | 00101100000 | 00011100000 | 11000010000 |
| 10100010000 | 01100010000 | 10010010000 | 01010010000 | 00110010000 |
| 10001010000 | 01001010000 | 00101010000 | 00011010000 | 10000110000 |
| 01000110000 | 00100110000 | 00010110000 | 00001110000 | 11000001000 |
| 10100001000 | 01100001000 | 10010001000 | 01010001000 | 00110001000 |
| 10001001000 | 01001001000 | 00101001000 | 00011001000 | 10000101000 |
| 01000101000 | 00100101000 | 00010101000 | 00001101000 | 10000011000 |
| 01000011000 | 00100011000 | 00010011000 | 00001011000 | 00000111000 |
| 11000000100 | 10100000100 | 01100000100 | 10010000100 | 01010000100 |
| 00110000100 | 10001000100 | 01001000100 | 00101000100 | 00011000100 |
| 10000100100 | 01000100100 | 00100100100 | 00010100100 | 00001100100 |
| 10000010100 | 01000010100 | 00100010100 | 00010010100 | 00001010100 |
| 00000110100 | 10000001100 | 01000001100 | 00100001100 | 00010001100 |
| | 00001001100 | 00000101100 | 00000011100 | |

Table 6: The minimal elements of the increasing set of figures 121 and 122.

Since the remove operation is irreversible, the loose variant will terminate after a finite number of iterations. One can hope that the loose gradient-remove is less likely to be trapped in a local maximum than the strict gradient-remove. For $n = 7$, we examined all the sets for which the index $m$ of the lowest layer satisfies $m \leq 3$ or $m \leq 2$. The results are presented in table 7. A set is stuck if there is no better set in its neighborhood, and it is a violation if in addition the difference between the numbers max and min defined in (54.1) is strictly larger than one. The results for $m \leq 2$ are quite encouraging, yet they show that the remove operation alone will not be enough to conclude.

| gradient-remove | examined | stuck | violation | run time |
|---|---|---|---|---|
| strict, $m \leq 3$ | 474546235 | 95262779 | 7536729 | 1 d 5 h 11 m |
| loose, $m \leq 3$ | 474546235 | 26781732 | 2353773 | 23 h 32 m |
| strict, $m \leq 2$ | 62529260 | 120 | 2 | 2 h 51 m 56 s |
| loose, $m \leq 2$ | 62529260 | 36 | 2 | 2 h 44 m 52 s |

Table 7: The number of local maxima of gradient-remove for $n = 7$

| | | | | |
|---|---|---|---|---|
| 10000101001 | 00000011011 | 00000010111 | 11000000011 | 10001000011 |
| 10000001011 | 00101011000 | 00101001100 | 00101001010 | 00100101001 |
| 01000010011 | 00010101001 | 00010001110 | 00111000100 | 00011000110 |
| 00011010100 | 01010100100 | 10011000100 | 10010100100 | 10010001100 |
| 11000010001 | 00010000111 | 00001000111 | 00000001111 | 01001000011 |
| 00001011001 | 00001010101 | 10001000101 | 00101010001 | 00101000101 |
| 10010000110 | 01001010001 | 01000010101 | 00111000010 | 00110000110 |
| 00110000011 | 00101000011 | 00100001011 | 00100000111 | 01100100010 |
| 01100011000 | 01100010100 | 01100010001 | 10100101000 | 10100100010 |
| 10100100001 | 01101100000 | 01100110000 | 01100101000 | 01110000001 |
| 01100100001 | 11100010000 | 10110010000 | 10101010000 | 10100011000 |
| 10100010100 | 10100010010 | 10100010001 | 10110001000 | 00010110001 |
| 00010010101 | 10010001001 | 10010000011 | 11010000001 | 01011000001 |
| 01010001001 | 01010000101 | 01010000011 | 00100110100 | 10110000100 |
| 10110000010 | 10110000001 | 11110000000 | 01111000000 | 01110001000 |
| 01110000100 | 01110000010 | 10000110010 | 00001100101 | 10001011000 |
| 10001100010 | 00001100110 | 00001111000 | 00001101100 | 01000001110 |
| 10010010010 | 10001010010 | 11000010010 | 00010110010 | 00001110010 |
| 00010111000 | 00110110000 | 00100100101 | 01000011010 | 00100100110 |
| 01000100101 | 00011100001 | 11000101000 | 10000101100 | 01100001100 |
| 10000110100 | 11000011000 | 11000110000 | 11010010000 | 10000001101 |
| 11000001001 | 01000001101 | 00010100011 | 01000111000 | 01000101100 |
| 01001001100 | 11001001000 | 10001001100 | 10010011000 | 01001110000 |
| 01001011000 | 01011010000 | 11100000100 | 11100000010 | 10100001010 |
| 11010100000 | 01010110000 | 01010101000 | 01010100010 | 11011000000 |
| 01011100000 | 00000110110 | 01011001000 | 01011000010 | 00111100000 |
| 00101110000 | 00101100100 | 00101100010 | 00101010010 | 00001100011 |
| 01001001010 | 00100010110 | 00011101000 | 00011011000 | 00111010000 |
| 10100000101 | 01100000110 | 01010010010 | 10000011100 | 00100011100 |
| 01000110100 | 00010101010 | 10011001000 | 10011000010 | 10011100000 |
| 10000011001 | 00100011001 | 10000010101 | 01000100110 | 00111001001 |

Table 8: The 155 minimal elements of the final configuration for the gradient exploration starting from $A_0 = \{0,1\}^{11}$. The corresponding set $A$ satisfies $|A| = 1614$, $\phi_1(A) = 1874$, $\phi_{2\setminus 1}(A) = 4242$, $\phi_2(A) = 6116$,

$$\left\{\omega : \sum_{1\le i\le 11} \omega(i) \ge 6\right\} \subset A \subset \left\{\omega : \sum_{1\le i\le 11} \omega(i) \ge 4\right\}.$$

In some further experiments, we ran both algorithms for $n$ between 3 and 16. We ended up with the configurations described in table 9. Given a set $A$, the column "size" displays $|\mathrm{Minimal}(A)|$, and $\min(A)$, $\max(A^c)$ are defined as

$$\min(A) \,=\, \min\big\{\,|\operatorname{supp}\omega| : \omega\in A\,\big\}\,,\quad \max(A^c) \,=\, \max\big\{\,|\operatorname{supp}\omega| : \omega\notin A\,\big\}\,. \tag{54.1}$$

The results are annoying, because the final set found by the algorithm intersects the level 4 for $n = 12, 13, 14$ and the level 5 for $n = 16$. This does not invalidate the conjecture 53.1, however it means that the gradient-remove algorithm will not be enough to prove it, were it true.

| $n$ | size | $\min(A)$ | $\max(A^c)$ | $\|A\|$ | $\phi_1(A)$ | $\phi_2(A)$ | run time |
|---|---|---|---|---|---|---|---|
| strict gradient-remove algorithm | | | | | | | |
| 3 | 3 | 1 | 0 | 7 | 3 | 9 | 0 s |
| 4 | 6 | 2 | 1 | 11 | 12 | 24 | 0 s |
| 5 | 7 | 2 | 2 | 22 | 24 | 62 | 0 s |
| 6 | 12 | 2 | 2 | 50 | 44 | 136 | 0 s |
| 7 | 16 | 2 | 3 | 105 | 81 | 268 | 0 s |
| 8 | 30 | 3 | 4 | 182 | 230 | 653 | 0 s |
| 9 | 39 | 3 | 4 | 389 | 437 | 1383 | 0 s |
| 10 | 72 | 3 | 5 | 795 | 880 | 2942 | 1 s |
| 11 | 168 | 4 | 5 | 1605 | 1873 | 6239 | 15 s |
| 12 | 233 | 4 | 5 | 3493 | 3112 | 11364 | 42 s |
| 13 | 329 | 4 | 7 | 6549 | 7645 | 27648 | 0 h 13 m 24 s |
| 14 | 894 | 4 | 7 | 13869 | 13872 | 51478 | 1 h 10 m 6 s |
| 15 | 1438 | 5 | 7 | 28951 | 23763 | 95201 | 6 h 20 m 10 s |
| 16 | 2041 | 5 | 7 | 60250 | 37548 | 163297 | 23 h 55 m 56 s |
| loose gradient-remove algorithm | | | | | | | |
| 3 | 3 | 1 | 0 | 7 | 3 | 9 | 0 s |
| 4 | 4 | 2 | 2 | 9 | 12 | 24 | 0 s |
| 5 | 7 | 2 | 2 | 22 | 24 | 62 | 0 s |
| 6 | 12 | 2 | 2 | 50 | 44 | 136 | 0 s |
| 7 | 21 | 3 | 3 | 83 | 121 | 296 | 0 s |
| 8 | 27 | 3 | 4 | 184 | 228 | 660 | 0 s |
| 9 | 51 | 3 | 5 | 351 | 505 | 1440 | 0 s |
| 10 | 87 | 3 | 5 | 737 | 1018 | 3007 | 2 s |
| 11 | 155 | 4 | 5 | 1614 | 1874 | 6116 | 10 s |
| 12 | 197 | 4 | 6 | 3093 | 4124 | 13559 | 1 m 51 s |
| 13 | 530 | 4 | 6 | 6456 | 8130 | 27473 | 15 m 9 s |
| 14 | 898 | 5 | 7 | 13524 | 14990 | 54026 | 2 h 10 m |
| 15 | 1385 | 5 | 8 | 25087 | 35359 | 125158 | 1 d 5 h 52 m |
| 16 | 2150 | 5 | 9 | 51011 | 70646 | 256400 | 8 d 11 h 48 m |

Table 9: The characteristics of the final set of gradient-remove

### 54.1.2 The variation of $\phi_2$

We shall next try to compute the variation of $\phi_2$ induced by the removal of a minimal element $\sigma$. It turns out that the removal of $\sigma$ will affect only the second order pivots of the configurations which are at Hamming distance 0, 1 or 2 from $\sigma$. In the next proposition, we describe more precisely the effect of the removal of $\sigma$ on the pivotal sets $\mathcal{P}_1^+$ and $\mathcal{P}_2^+$.

**Proposition 54.1.** *Let $A$ be an increasing subset of $\{0,1\}^n$. Let $\sigma$ belong to Minimal($A$), and let $B = A \setminus \{\sigma\}$. For any configuration $\omega$ such that $H(\omega,\sigma) \geq 3$ or $\omega \not\geq \sigma$, we have*

$$\mathcal{P}_2^+(A,\omega) = \mathcal{P}_2^+(B,\omega)\,. \tag{54.2}$$

*For $i \notin \operatorname{supp}\sigma$, we have $i \in \mathcal{P}_1^+(B,\sigma^i)$ and moreover*

$$\begin{aligned} \mathcal{P}_1^+(A,\sigma^i)\setminus\{i\} \,&=\, \mathcal{P}_1^+(B,\sigma^i)\setminus\{i\}\,,\\ \mathcal{P}_{2\setminus 1}^+(B,\sigma^i) \,&\subset\, \mathcal{P}_{2\setminus 1}^+(A,\sigma^i)\,. \end{aligned} \tag{54.3}$$

*For $i,j \notin \operatorname{supp}\sigma$, $i<j$, we have*

$$\{i,j\} \subset \mathcal{P}_2^+(B,\sigma^{ij})\,, \tag{54.4}$$

$$\mathcal{P}_2^+(A,\sigma^{ij})\setminus\{i,j\} \,=\, \mathcal{P}_2^+(B,\sigma^{ij})\setminus\{i,j\}\,. \tag{54.5}$$

*The results (54.2) and (54.5) hold for the sets $\mathcal{P}_{2\setminus 1}^+$ as well.*

*Proof.* We do the proofs of (54.2) and (54.5) for the sets $\mathcal{P}_2^+$ only, the argument for the sets $\mathcal{P}_{2\setminus 1}^+$ is the same. Let $\omega \in \{0,1\}^n$ be such that $H(\omega,\sigma) \geq 3$ or $\omega \not\geq \sigma$. We see from the definition (49.1) that the set $\mathcal{P}_2^+(A,\omega)$ depends only on the trace of $A$ on the set of configurations

$$\{\,\omega_I : I \subset \{1,\dots,n\},\, |I| \leq 2\,\}\,.$$

The condition imposed on $\omega$ implies that the set $B$ has the same trace as $A$ on this set, therefore $\mathcal{P}_2^+(A,\omega) = \mathcal{P}_2^+(B,\omega)$. The statements concerning the first pivots are immediate. We prove next the inclusion (54.3). Let $i \notin \operatorname{supp}\sigma$ and let $j$ belong $\mathcal{P}_{2\setminus 1}^+(B,\sigma^i)$. Necessarily, the index $j$ is not equal to $i$, and by definition, there exists $k \in \{i\} \cup \operatorname{supp}\sigma \setminus \{j\}$ such that

$$\sigma^i_j \in B\,,\quad \sigma^i_k \in B\,,\quad \sigma^i_{jk} \notin B\,. \tag{54.6}$$

Necessarily, the index $k$ is not equal to $i$, and none of the configurations appearing in (54.6) can coincide with $\sigma$. Therefore we have

$$\sigma^i_j \in A\,,\quad \sigma^i_k \in A\,,\quad \sigma^i_{jk} \notin A\,, \tag{54.7}$$

and we see that $j$ belongs to $\mathcal{P}_{2\setminus 1}^+(A,\sigma^i)$. This concludes the proof of the inclusion (54.3). Unfortunately, we cannot use this argument to prove the converse inclusion. Indeed, the conditions stated in (54.7) imply those stated in (54.6) only if $k$ is distinct from $i$, and this cannot be guaranteed even if we start with an index $j \neq i$. We prove finally the inclusion (54.4) and the equality (54.5). Let $i, j \notin \operatorname{supp} \sigma$ with $i < j$. The inclusion (54.4) follows from the fact that $\sigma^i$ and $\sigma^j$ are in $B$ while $\sigma$ is not. Let $k$ belong $\mathcal{P}_2^+(A, \sigma^{ij}) \setminus \{ i, j \}$. By definition, there exists $l \in \{ i, j \} \cup \operatorname{supp} \sigma \setminus \{ k \}$ such that

$$(\sigma^{ij})_k \in A\,, \quad (\sigma^{ij})_l \in A\,, \quad (\sigma^{ij})_{kl} \notin A\,. \tag{54.8}$$

Since $k$ is different from $i$ and $j$, none of the configurations appearing in (54.8) can coincide with $\sigma$. Therefore we have

$$(\sigma^{ij})_k \in B\,, \quad (\sigma^{ij})_l \in B\,, \quad (\sigma^{ij})_{kl} \notin B\,,$$

and we see that $k$ belongs to $\mathcal{P}_2^+(B, \sigma^{ij}) \setminus \{ i, j \}$. We can moreover reverse the argument to get the converse inclusion, and this proves the equality (54.5). □

With the help of proposition 54.1, we are now able to evaluate the variation of $\phi_2$ induced by a removal.

**Lemma 54.2.** *Let $A$ be an increasing subset of $\{ 0, 1 \}^n$. Let $\sigma$ belong to Minimal($A$), and let $B = A \setminus \{ \sigma \}$. The variation of $\phi_2$ between $A$ and $B$ is given by*

$$\begin{aligned} \phi_2(B) - \phi_2(A) \,=\, n - 2|\operatorname{supp} \sigma| + \sum_{i \notin \operatorname{supp} \sigma} \Big( \big|\mathcal{P}_{2\setminus 1}^+(B, \sigma^i)\big| - \big|\mathcal{P}_{2\setminus 1}^+(A, \sigma^i)\big| \Big) \\ + \sum_{\substack{\{i,j\}: i \neq j \\ i,j \notin \operatorname{supp} \sigma}} \Big( \big|\mathcal{P}_{2\setminus 1}^+(B, \sigma^{ij})\big| - \big|\mathcal{P}_{2\setminus 1}^+(A, \sigma^{ij})\big| \Big)\,. \end{aligned} \tag{54.9}$$

*The first sum in formula* (54.9) *is non-positive, while the second sum in formula* (54.9) *is non-negative.*

*Proof.* The variation of $\phi_2$ is the sum of the variations of $\phi_1$ and $\phi_{2\setminus 1}$:

$$\phi_2(B) - \phi_2(A) \,=\, \phi_1(B) - \phi_1(A) \,+\, \phi_{2\setminus 1}(B) - \phi_{2\setminus 1}(A)\,.$$

The variation of $\phi_1$ has already been computed in formula (50.10), it is equal to $n - 2\big|\operatorname{supp} \sigma\big|$. We focus next on the variation of $\phi_{2\setminus 1}$. By definition, we have

$$\phi_{2\setminus 1}(B) - \phi_{2\setminus 1}(A) \,=\, \sum_{\omega \in B} \big|\mathcal{P}_{2\setminus 1}^+(B, \omega)\big| - \sum_{\omega \in A} \big|\mathcal{P}_{2\setminus 1}^+(A, \omega)\big|\,. \tag{54.10}$$

It follows from proposition 54.1 that the only terms that might not vanish in the sum appearing in (54.10) correspond to configurations $\omega$ such that $H(\omega, \sigma) \leq 2$ and $\omega \geq \sigma$. These are the configurations belonging to the set

$$\{ \sigma \} \,\cup\, \big\{ \sigma^i : i \notin \operatorname{supp} \sigma \big\} \,\cup\, \big\{ \sigma^{ij} : i, j \notin \operatorname{supp} \sigma,\, i < j \big\}\,.$$

Since $\sigma$ is a minimal element of $A$, the configuration $\omega = \sigma$ contains only pivots of first order, therefore the variation of $\phi_{2\setminus 1}$ reduces to the two sums appearing in formula (54.9). The inclusion (54.3) implies that the first sum in formula (54.9) is non-positive. The inclusions (54.4) and (54.5) imply that the second sum in formula (54.9) is non-negative. □

### 54.1.3 The local maxima of $\phi_2$

What can we say about the local maxima of $\phi_2$ with respect to the remove operation? To be honest, not much, unfortunately, beyond the next lemma. From now onwards, we will designate by $\mathcal{H}^n_-$ and $\mathcal{H}^n_+$ the lower and the upper part of the hypercube, i.e.,

$$\begin{aligned}\mathcal{H}^n_- &= \Big\{\, \omega\in\{\,0,1\,\}^n : |\operatorname{supp}\omega| < \frac{n}{2}\,\Big\}\,,\\ \mathcal{H}^n_+ &= \Big\{\, \omega\in\{\,0,1\,\}^n : |\operatorname{supp}\omega| > \frac{n}{2}\,\Big\}\,.\end{aligned}$$

**Lemma 54.3.** *Let $A$ be a local maximum of $\phi_2$ with respect to the remove operation. Then it satisfies the following:*

$$\forall\sigma\in \mathit{Minimal}(A)\cap\mathcal{H}^n_-\quad \exists\, i\notin\operatorname{supp}\sigma\quad \exists\, j\in\operatorname{supp}\sigma\quad \sigma^i_j\in A\,, \tag{54.11}$$

*or equivalently*

$$\forall\sigma\in \mathit{Minimal}(A)\cap\mathcal{H}^n_-\quad \exists\, i\notin\operatorname{supp}\sigma\quad i\in\mathcal{P}^+_2(A,\sigma^i)\,. \tag{54.12}$$

*Proof.* We prove first that the two conditions (54.11) and (54.12) are equivalent. Suppose that the condition (54.11) holds. Let $\sigma$ belong to $\text{Minimal}(A)\cap\mathcal{H}^n_-$ and let $i,j$ be indices associated to $\sigma$ as in (54.11). We have then

$$\sigma^i_j\in A\,,\quad \sigma=\sigma^j\in A\,,\quad \sigma_j\notin A\,, \tag{54.13}$$

and this implies that $i$ belongs to $\mathcal{P}^+_2(A,\sigma^i)$. Thus the condition (54.12) holds. Conversely, suppose that (54.12) holds. Let $\sigma$ belong to $\text{Minimal}(A)\cap\mathcal{H}^n_-$ and let $i$ be an index associated to $\sigma$ as in (54.12). Since $i$ belongs to $\mathcal{P}^+_2(A,\sigma^i)$, then there exists $j$ in $\operatorname{supp}\sigma$ such that (54.13) holds, thus the condition (54.11) holds as well.

Let next $A$ be a local maximum of $\phi_2$. Suppose that $A$ does not satisfy the condition (54.11). Let $\sigma$ be an element of $\text{Minimal}(A)\cap\mathcal{H}^n_-$ violating (54.11), i.e.,

$$\forall i\notin\operatorname{supp}\sigma\quad \forall j\in\operatorname{supp}\sigma\quad \sigma^i_j\notin A\,. \tag{54.14}$$

Let $B = A\setminus\{\,\sigma\,\}$. The variation of $\phi_2$ between $A$ and $B$ is given by formula (54.9). The condition (54.14) implies that

$$\forall i\notin\operatorname{supp}\sigma\quad \mathcal{P}^+_2(A,\sigma^i)\setminus\mathcal{P}^+_1(A,\sigma^i)=\varnothing\,,$$

and this implies furthermore that the first sum in formula (54.9) vanishes. In lemma 54.2, we proved that the second sum in formula (54.9) is non-negative, so we conclude that

$$\phi_2(B)-\phi_2(A)\;\geq\; n-2|\operatorname{supp}\sigma|\;>\;0\,.$$

The last inequality is due to the fact that $\sigma\in\mathcal{H}^n_-$, and it contradicts the fact that $A$ is a local maximum of $\phi_2$. Thus the set $A$ has to satisfy the condition (54.11). $\square$

## 54.2 Adding a maximal element

Symmetrically, we consider the add operation which consists in adding a maximal element of the complement. The corresponding system of neighborhoods is defined by

$$\mathcal{V}_{\text{add}}(A) \,=\, \big\{\, A \cup \{\,\eta\,\} : \eta \in \text{Maximal}(A^c) \,\big\}\,.$$

Fortunately, we can reuse the computations done for the remove operation. Indeed, let $A$ be an increasing subset of $\{\,0,1\,\}^n$ and let $\eta$ be an element of $\text{Maximal}(A^c)$. Let us set $C = A \cup \{\,\eta\,\}$. Then $\eta$ is an element of $\text{Minimal}(C)$, therefore the variation $\phi_2(A) - \phi_2(C)$ is given by lemma 54.2 applied to $A = C$ and $B = A$, so that, after a sign change, we have

$$\begin{aligned}\phi_2(C) - \phi_2(A) \,=\, 2|\operatorname{supp}\eta| - n + \sum_{i \notin \operatorname{supp}\eta} \Big( \big|\mathcal{P}^+_{2\setminus 1}(C,\eta^i)\big| - \big|\mathcal{P}^+_{2\setminus 1}(A,\eta^i)\big| \Big) \\ + \sum_{\substack{\{i,j\}: i\neq j \\ i,j \notin \operatorname{supp}\eta}} \Big( \big|\mathcal{P}^+_{2\setminus 1}(C,\eta^{ij})\big| - \big|\mathcal{P}^+_{2\setminus 1}(A,\eta^{ij})\big| \Big)\,. \quad (54.15)\end{aligned}$$

The first sum in formula (54.15) is non-negative, while the second sum in formula (54.15) is non-positive. In addition, we can apply the results of proposition 54.1 with $A = C$ and $B = A$. Together with formula (54.15), this yields some conditions satisfied by the local maxima of $\phi_2$ with respect to the add operation. These conditions cannot be deduced directly from those stated in lemma 54.3 for the remove operation, so we state them in the next lemma.

**Lemma 54.4.** *Let $A$ be a local maximum of $\phi_2$ with respect to the add operation. Then it satisfies the following:*

$$\begin{aligned}\forall \eta \in \mathit{Maximal}(A^c) \cap \mathcal{H}^n_+ \quad \exists i,j \notin \operatorname{supp}\eta\,, \quad i \neq j\,, \\ \forall\, k \in \operatorname{supp}\eta \quad \eta^{ij}_k \notin A \ \ \mathit{or}\ \eta^i_k \in A\,. \quad (54.16)\end{aligned}$$

*Proof.* Let $A$ be a local maximum of $\phi_2$. Suppose that $A$ does not satisfy (54.16). Let $\eta$ be an element of $\text{Maximal}(A^c) \cap \mathcal{H}^n_+$ violating (54.16), i.e.,

$$\forall i,j \notin \operatorname{supp}\eta\,, \quad i \neq j\,, \quad \exists\, k \in \operatorname{supp}\eta \quad \eta^{ij}_k \in A \text{ and } \eta^i_k \notin A\,. \qquad (54.17)$$

Let $C = A \cup \{\,\eta\,\}$. The variation of $\phi_2$ between $C$ and $A$ is given by formula (54.15). The condition (54.17) implies that

$$\forall i,j \notin \operatorname{supp}\eta\,, \quad i \neq j\,, \qquad \{\,i,j\,\} \subset \mathcal{P}^+_2(C,\eta^{ij})\,,$$

and together with the equality (54.5) applied to $A = C$ and $B = A$, this implies that the second sum in formula (54.15) vanishes. Since we know that the first sum in formula (54.15) is non-negative, we conclude that

$$\phi_2(C) - \phi_2(A) \,\geq\, 2|\operatorname{supp}\sigma| - n \,>\, 0\,.$$

The last inequality is due to the fact that $\sigma \in \mathcal{H}^n_+$, and it contradicts the fact that $A$ is a local maximum of $\phi_2$ with respect to the add operation. Thus the set $A$ has to satisfy the condition (54.16). □

We performed the same numerical computations for gradient-add as for gradient-remove. The results are presented in tables 10 and 11, they lead to similar conclusions (the content of the tables is explained before table 9).

| gradient-add | examined | stuck | violation | run time |
|---|---|---|---|---|
| strict | 474546235 | 6853 | 3512 | 18 h 50 m |
| loose | 474546235 | 390 | 177 | 18 h 27 m |

Table 10: The number of local maxima of strict gradient-add for $n = 7$

| $n$ | size | $\min(A)$ | $\max(A^c)$ | $\|A\|$ | $\phi_1(A)$ | $\phi_2(A)$ | run time |
|---|---|---|---|---|---|---|---|
| strict gradient-add algorithm | | | | | | | |
| 3 | 3 | 2 | 1 | 4 | 6 | 9 | 0 s |
| 4 | 4 | 3 | 2 | 5 | 12 | 16 | 0 s |
| 5 | 10 | 3 | 2 | 16 | 30 | 50 | 0 s |
| 6 | 15 | 4 | 3 | 22 | 60 | 90 | 0 s |
| 7 | 26 | 4 | 4 | 49 | 125 | 227 | 0 s |
| 8 | 51 | 5 | 5 | 87 | 268 | 444 | 0 s |
| 9 | 112 | 5 | 5 | 235 | 609 | 1176 | 1 s |
| 10 | 190 | 6 | 6 | 358 | 1204 | 2124 | 6 s |
| 11 | 392 | 6 | 7 | 953 | 2699 | 5799 | 55 s |
| 12 | 717 | 7 | 7 | 1466 | 5304 | 9804 | 0 h 9 m 42 s |
| 13 | 1454 | 7 | 8 | 3826 | 11728 | 26107 | 1 h 1 m 24 s |
| 14 | 2619 | 8 | 9 | 6086 | 23232 | 47760 | 7 h 4 m 12 s |
| 15 | 2983 | 8 | 9 | 12489 | 46721 | 105857 | 12 h 14 m |
| 16 | 7784 | 9 | 10 | 22298 | 94138 | 200974 | 7 d 23 h |
| loose gradient-add algorithm | | | | | | | |
| 3 | 3 | 2 | 1 | 4 | 6 | 9 | 0 s |
| 4 | 6 | 2 | 1 | 11 | 12 | 24 | 0 s |
| 5 | 10 | 3 | 2 | 16 | 30 | 50 | 0 s |
| 6 | 15 | 3 | 3 | 32 | 60 | 120 | 0 s |
| 7 | 31 | 4 | 4 | 59 | 135 | 247 | 0 s |
| 8 | 61 | 4 | 4 | 148 | 280 | 608 | 0 s |
| 9 | 112 | 5 | 5 | 235 | 609 | 1176 | 1 s |
| 10 | 210 | 5 | 5 | 596 | 1260 | 2989 | 8 s |
| 11 | 414 | 6 | 6 | 940 | 2688 | 5394 | 0 h 1 m 32 s |
| 12 | 732 | 6 | 7 | 2287 | 5500 | 13473 | 0 h 8 m 39 s |
| 13 | 1472 | 7 | 8 | 3845 | 11751 | 26159 | 1 h 3 m 56 s |
| 14 | 1352 | 7 | 9 | 7480 | 23398 | 56459 | 1 h 39 m 17 s |
| 15 | 3552 | 8 | 9 | 13160 | 47576 | 109023 | 16 h 40 m |
| 16 | 4233 | 8 | 10 | 28920 | 100066 | 243469 | 5 d 5 h 51 m |

Table 11: The characteristics of the final set of gradient-add

## 54.3 Combination of add or remove

We see from the numerical experiments that neither the remove operation nor the add operation allow to reach the global maximum of $\phi_2$. We examine next what happens if we combine the two operations. At each step, we are allowed either to remove or to add a configuration. This corresponds to the system of neighborhoods defined by

$$\mathcal{V}_{\text{add or remove}}(A) \,=\, \mathcal{V}_{\text{add}}(A) \cup \mathcal{V}_{\text{remove}}(A)\,.$$

We call gradient-aor the associated gradient algorithm. To ensure the termination of the gradient-aor algorithm, we consider only the strict variant, in which a move is performed only when it leads to an increase of $\phi_2$, otherwise the algorithm could cycle indefinitely by removing and adding a configuration which does not change the value of $\phi_2$. We performed the same numerical computations for gradient-add-or-remove as for gradient-remove and gradient-add. Already for $n = 3$, there are two sets which are stuck for gradient-add-or-remove, although they are not the global maxima. However, if we accept modifications for which the value of $\phi_2$ does not change, then we can escape from these sets and reach a global maximum. Such a sequence of moves is shown in figure 124.

| $n$ | examined | stuck | violation | run time |
|---|---|---|---|---|
| 3 | 8 | 2 | 0 | 0 s |
| 4 | 28 | 5 | 0 | 0 s |
| 5 | 208 | 3 | 0 | 0 s |
| 6 | 16351 | 39 | 0 | 1 s |
| 7 | 490013146 | 914 | 0 | 2 d 8 h 26 m |

Table 12: The number of local maxima of gradient-aor for $3 \leq n \leq 7$

The results of the exhaustive exploration for $3 \leq n \leq 7$ are presented in table 12. We are pleased to note that, although numerous local maxima do exist for $n = 6$ and $n = 7$, there is no violating set. Examples of local maxima for $n = 5, 6, 7$ are presented in figure 123 (these examples are in fact strict local maxima, every add or remove modification leads to a decrease of $\phi_2$). We are therefore tempted to believe in the following result.

**Hope 54.5.** *For any $n \geq 1$ and any increasing subset $A$ of $\{0,1\}^n$, there exists a sequence of increasing sets $(A_i, 0 \leq i \leq r)$ of $\{0,1\}^n$ such that $A_0 = A$, for any $i$ in $\{1, \dots, r\}$, the set $A_i$ is obtained from $A_{i-1}$ by the removal of a minimal configuration of $A$ or the addition of a maximal configuration of $A^c$, $\phi_2(A_{i-1}) < \phi_2(A_i)$, and moreover the final set $A^* = A_r$ of the sequence satisfies*

$$\max\big\{\,|\operatorname{supp}\omega| : \omega \notin A^*\,\big\} \,\leq\, 1 + \min\big\{\,|\operatorname{supp}\omega| : \omega \in A^*\,\big\}\,. \tag{54.18}$$

If hope 54.5 is true, then it would imply that any increasing subset $A^*$ of $\{0,1\}^n$ which realizes the maximum of $\phi_2$ would satisfy the condition (54.18), and from this we could easily obtain the desired control on the second order pivots.

Minimal($A$)={ 11100,11010,11001,10110,10101,10011,01110,01101,01011,00111 }

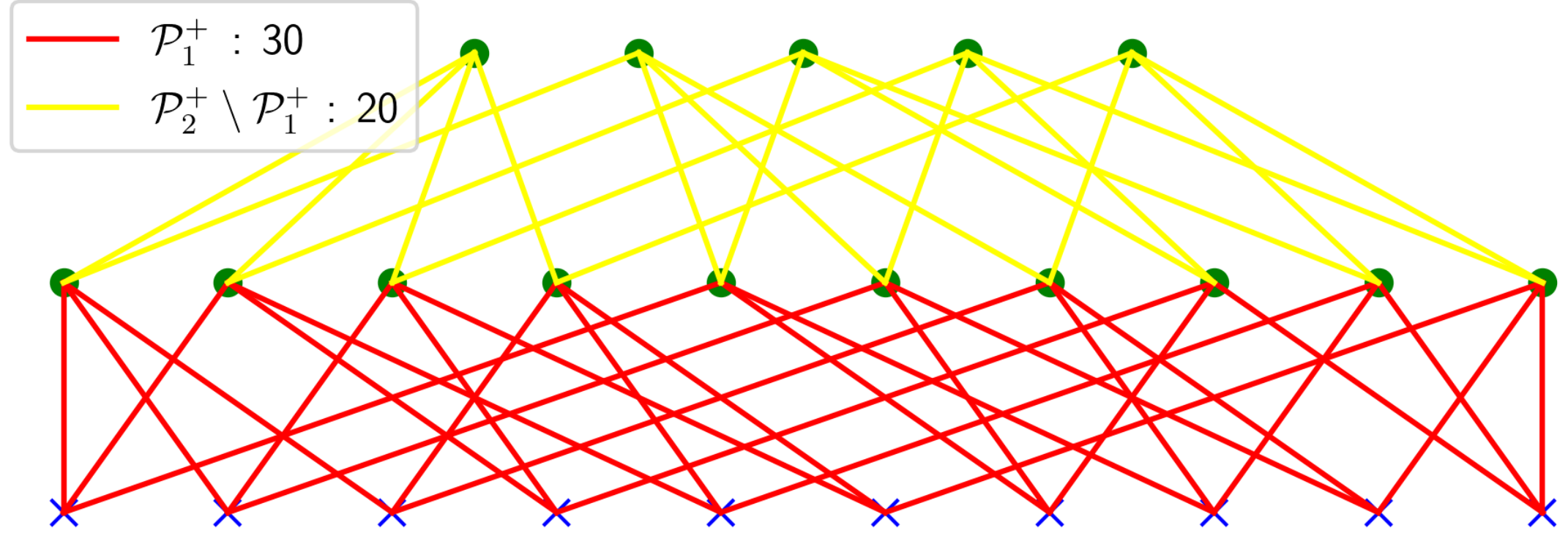

Minimal($A$)={ 111000,110100,110010,110001,101100,101010,100101,100011
011100,011001,010110,010011,001110,001101,001011,000111 }

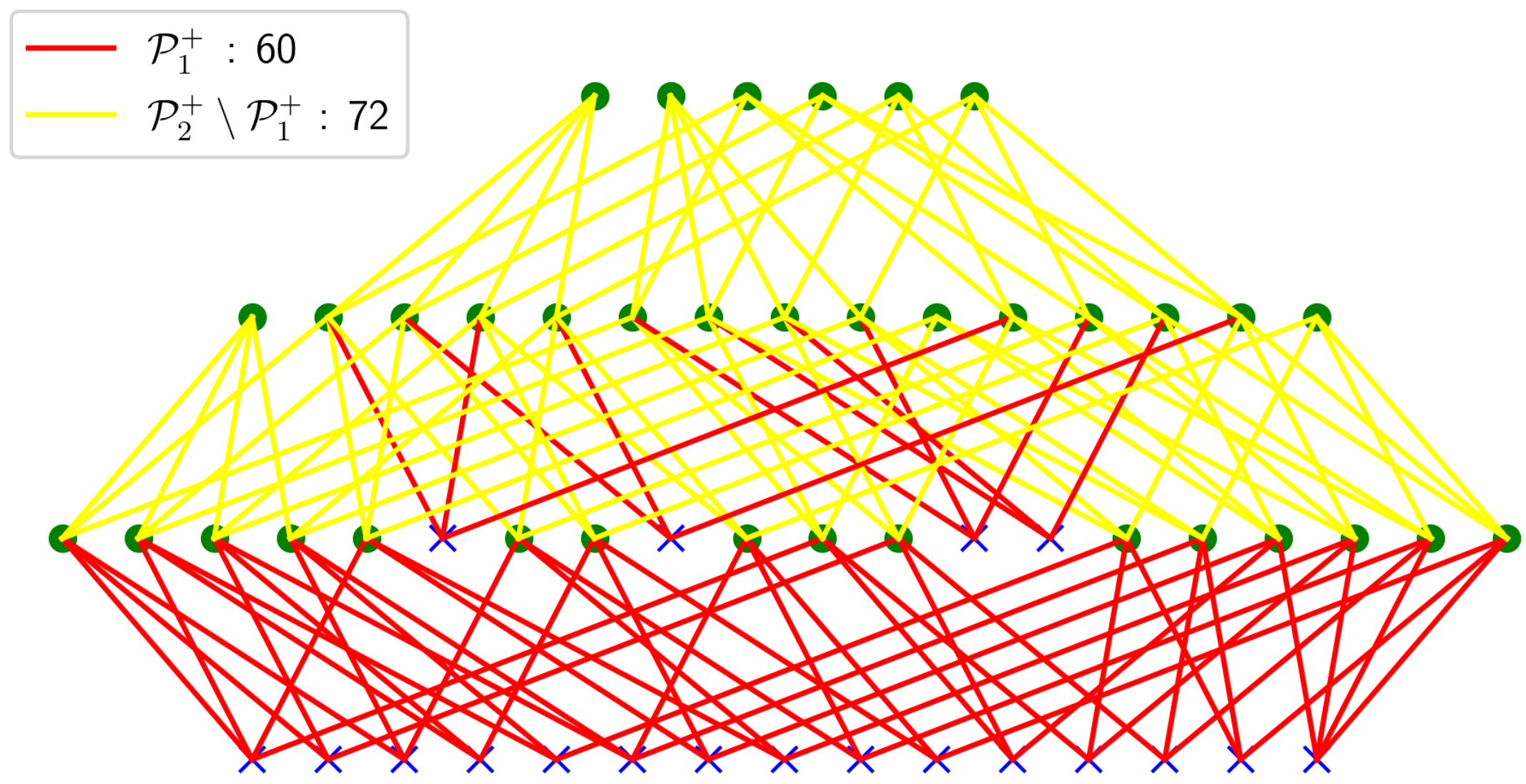

Minimal($A$)={ 1110000,1101000,1100100,1100010,1011000,1010110,1001001,1000101,
0111000,0110001,0101110,0100011,0011100,0011010,0010101,0001011,0000111 }

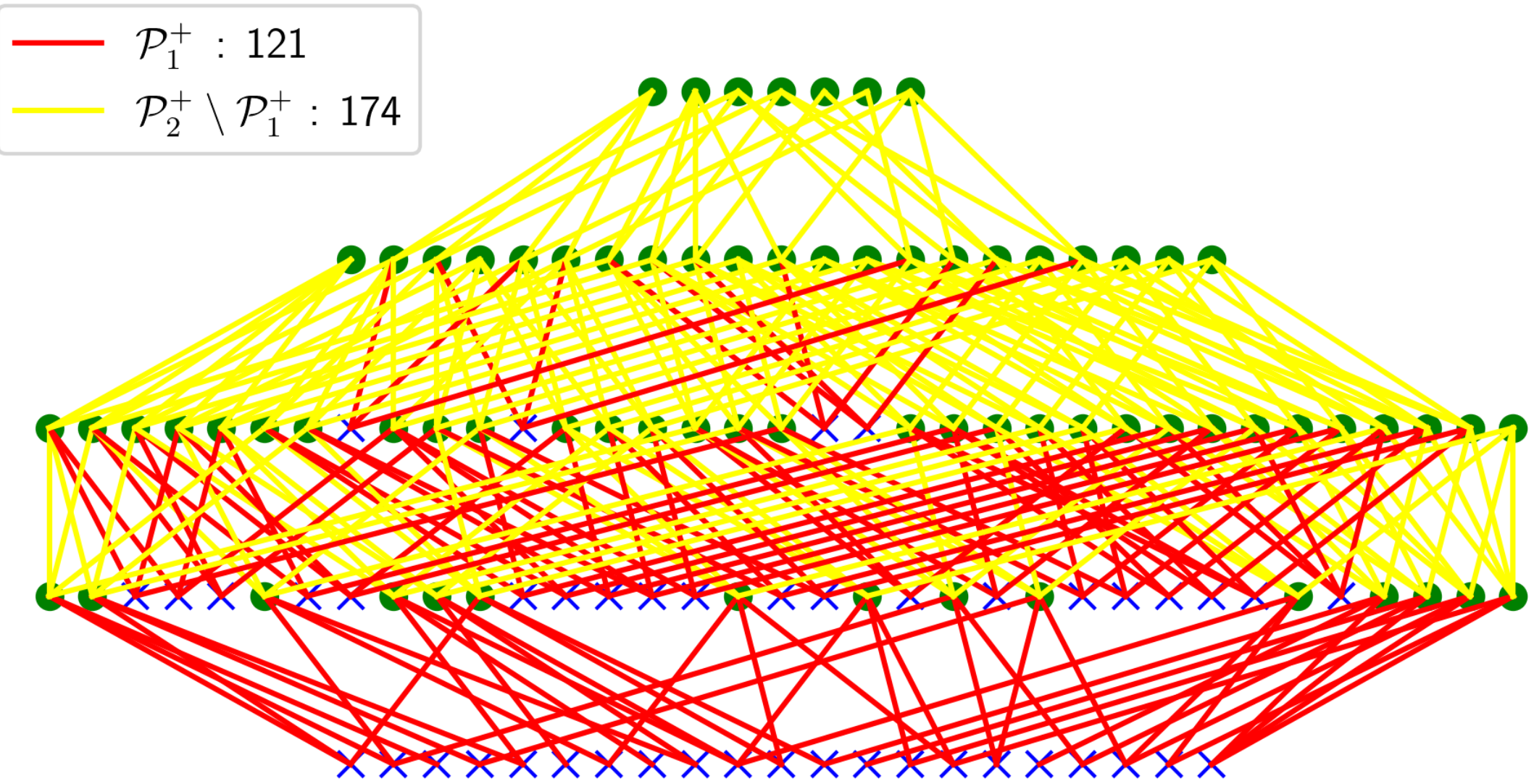

Figure 123: Strict local maxima of gradient-add-or-remove for $n = 5, 6, 7$

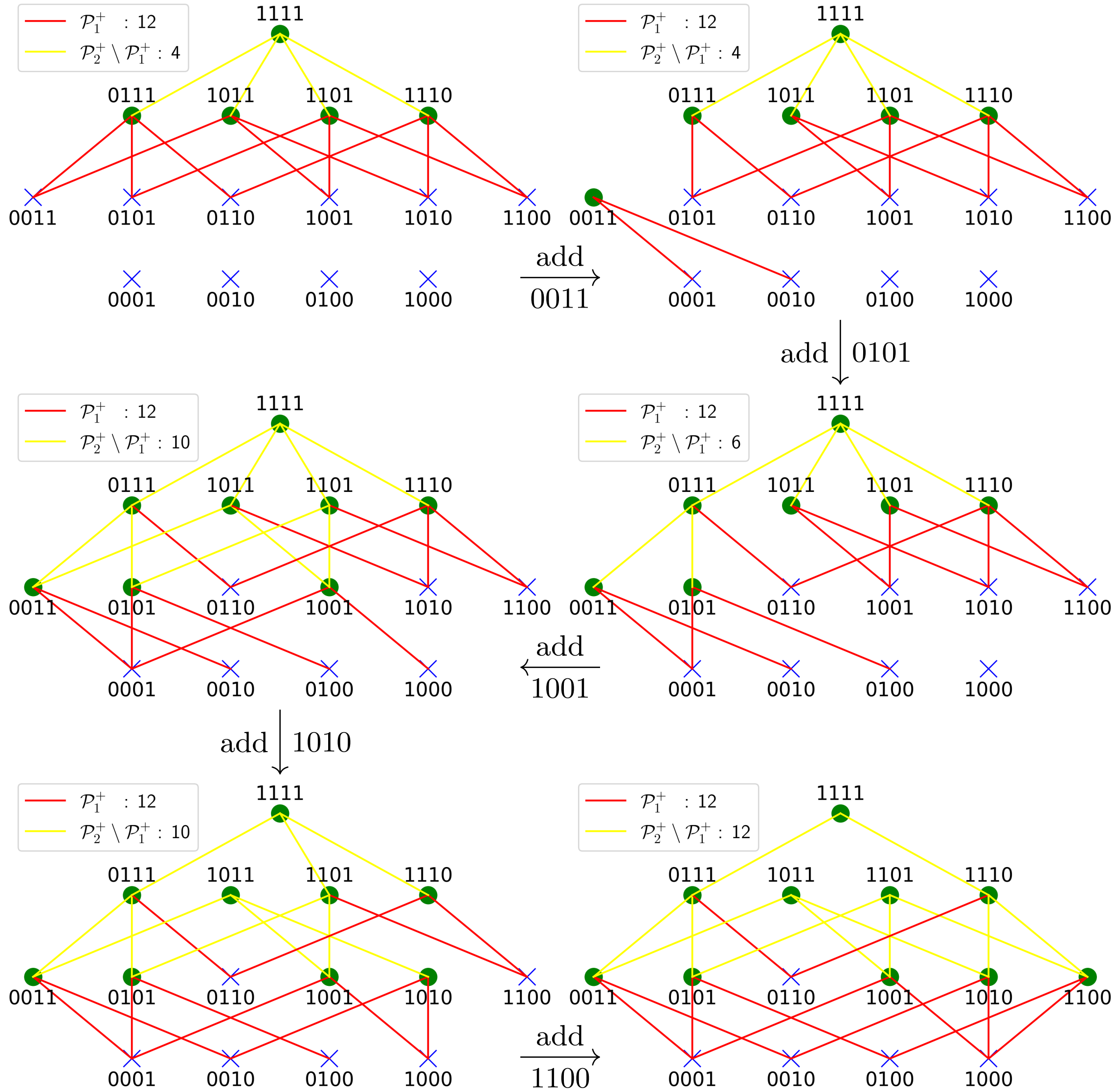

Figure 124: The set $A$ such that $\text{Minimal}(A) = \{0111, 1011, 1101, 1110\}$ is not a global maximum of $\phi_2$, and there is no add or remove modification which increases $\phi_2$. However, the sequence of the successive additions of 0011, 0101, 1001, 1010, 1100 leads to a global maximum of $\phi_2$, without ever decreasing $\phi_2$.

For $n$ larger than 7, we launched the gradient-add-or-remove algorithm starting either from the largest configuration $1\cdots1$ or from the lowest configuration $0\cdots0$. In order to escape from the non strict local maxima of $\phi_2$, we tried also a relaxed version, in which 1000 moves keeping the function $\phi_2$ unchanged are authorized. The results are presented in table 14. The results starting from $1\cdots1$ are quite satisfactory, because the final configuration $A$ reached by the algorithm always satisfies $\max(A^c) \leq 1 + \min(A)$. Unfortunately, this is not the case when starting from $0\cdots0$. Indeed, for $n \geq 12$, the algorithm is stuck in a configuration $A$ such that $\max(A^c) \geq 2 + \min(A)$. In table 13, we present the maximal elements of the complement of one of these maxima. So it seems that the formulation of hope 54.5 should be amended. However, it is difficult to grasp the geometrical structure of the local maxima of $\phi_2$. A possibility would be to weaken the inequality (54.18) of hope 54.5, for instance to require that the final configuration $A$ satisfies $\max(A^c) \leq c + \min(A)$ for some constant $c$ independent of $n$. However, we see that for $n = 16$, the constant $c$ should already be larger or equal than 4, so it is dubious that this condition would work. Another possibility would be to require only a large inequality $\phi_2(A_{i-1}) \leq \phi_2(A_i)$ along the sequence of sets of hope 54.5. However, the last example of figure 123 shows a set which is a strict local maximum, and which violates the hope 54.5. So it seems that the hope 54.5 should be abandoned.

| | | | | |
|---|---|---|---|---|
| 111110000000 | 110001100000 | 011001100000 | 011110010000 | 011001010000 |
| 110011010000 | 101011010000 | 101000110000 | 100010110000 | 010001110000 |
| 110001001000 | 100101001000 | 010011001000 | 000111001000 | 111000101000 |
| 010100101000 | 010010101000 | 000011101000 | 000110011000 | 100000111000 |
| 001000111000 | 101100000100 | 100101000100 | 010101000100 | 000111000100 |
| 100010100100 | 001010100100 | 100100010100 | 011100010100 | 011010010100 |
| 101001010100 | 100011010100 | 001100110100 | 000010110100 | 011100001100 |
| 010010001100 | 100000101100 | 001001101100 | 010100011100 | 000000111100 |
| 011110000010 | 110100100010 | 001110100010 | 001110010010 | 101101010010 |
| 010000110010 | 100100001010 | 110000101010 | 000100101010 | 010000011010 |
| 100010011010 | 110000000110 | 010011000110 | 100000100110 | 001000100110 |
| 000100100110 | 010010010110 | 001010010110 | 001001010110 | 000011010110 |
| 010000001110 | 001010001110 | 001001001110 | 001000011110 | 000001011110 |
| 001101000001 | 101000100001 | 100110100001 | 010110100001 | 001010010001 |
| 101001010001 | 010100110001 | 000010110001 | 000101110001 | 001010001001 |
| 101001001001 | 010101001001 | 000100101001 | 110000011001 | 011100011001 |
| 000010011001 | 011100000101 | 010001000101 | 100011000101 | 010000100101 |
| 000001100101 | 110000010101 | 101000010101 | 010100010101 | 101000001101 |
| 010100001101 | 000110001101 | 011000011101 | 001100011101 | 111000000011 |
| 110100000011 | 010110000011 | 100011000011 | 010011000011 | 100010100011 |
| 000011100011 | 000001110011 | 000000101011 | 001001011011 | 001000000111 |
| | 000110000111 | 000000011111 | | |

Table 13: The 107 maximal elements of the complement $A^c$ of the final configuration for gradient-add-or-remove starting from 000000000000 (the 243 minimal elements of $A$ do not fit inside a single page). The corresponding set $A$ satisfies $|A| = 3444$, $\phi_1(A) = 3210$, $\phi_2(A) = 12217$.

| $n$ | size | min | max | $\|A\|$ | $\phi_1(A)$ | $\phi_2(A)$ | run time |
|---|---|---|---|---|---|---|---|
| gradient-add-or-remove starting from $1\cdots1$ | | | | | | | |
| 8 | 47 | 5 | 5 | 83 | 260 | 448 | 0 s |
| 9 | 105 | 5 | 5 | 228 | 602 | 1225 | 3 s |
| 10 | 176 | 6 | 6 | 344 | 1176 | 2208 | 28 s |
| 11 | 382 | 6 | 6 | 943 | 2691 | 5888 | 6 m 4 s |
| 12 | 665 | 7 | 7 | 1459 | 5290 | 10570 | 44 m 35 s |
| 13 | 1191 | 7 | 8 | 3349 | 10935 | 24826 | 5 h 3 m 24 s |
| 14 | 2529 | 8 | 9 | 5991 | 23034 | 48190 | 2 d 17 h 45 m |
| 15 | 2847 | 8 | 9 | 12293 | 46443 | 105195 | 4 d 17 h 17 m |
| relaxed gradient-add-or-remove starting from $1\cdots1$ | | | | | | | |
| 8 | 54 | 4 | 4 | 141 | 280 | 640 | 0 s |
| 9 | 105 | 5 | 5 | 228 | 602 | 1225 | 3 s |
| 10 | 120 | 5 | 6 | 456 | 1176 | 2768 | 32 s |
| 11 | 382 | 6 | 6 | 943 | 2691 | 5888 | 5 m 13 s |
| 12 | 277 | 6 | 7 | 1695 | 5292 | 11839 | 45 m 43 s |
| 13 | 1191 | 7 | 8 | 3349 | 10935 | 24826 | 4 h 15 m 33 s |
| 14 | 922 | 7 | 9 | 6735 | 23096 | 52683 | 2 d 12 h 20 m |
| 15 | 2851 | 8 | 9 | 12302 | 46452 | 105270 | 6 d 16 h 17 m |
| gradient-add-or-remove starting from $0\cdots0$ | | | | | | | |
| 8 | 31 | 3 | 4 | 180 | 230 | 671 | 0 s |
| 9 | 44 | 3 | 4 | 394 | 428 | 1406 | 0 s |
| 10 | 78 | 3 | 5 | 803 | 862 | 2977 | 2 s |
| 11 | 188 | 4 | 5 | 1625 | 1819 | 6344 | 24 s |
| 12 | 243 | 4 | 6 | 3444 | 3210 | 12217 | 1 m 34 s |
| 13 | 351 | 4 | 7 | 6556 | 7628 | 27929 | 11 m 6 s |
| 14 | 1009 | 4 | 7 | 13975 | 13424 | 52404 | 1 h 23 m |
| 15 | 1508 | 5 | 7 | 28898 | 23754 | 98952 | 8 h 2 m |
| 16 | 3084 | 5 | 9 | 49571 | 74376 | 267310 | 8 d 9 h 2 m |
| relaxed gradient-add-or-remove starting from $0\cdots0$ | | | | | | | |
| 8 | 29 | 3 | 4 | 170 | 240 | 694 | 0 s |
| 9 | 43 | 3 | 4 | 392 | 432 | 1424 | 0 s |
| 10 | 80 | 3 | 5 | 801 | 868 | 2993 | 9 s |
| 11 | 185 | 4 | 5 | 1623 | 1825 | 6353 | 58 s |
| 12 | 229 | 4 | 6 | 3400 | 3344 | 12640 | 2 m 52 s |
| 13 | 356 | 4 | 7 | 6520 | 7708 | 28258 | 51 m 36 s |
| 14 | 999 | 4 | 7 | 13937 | 13562 | 52621 | 6 h 46 m |
| 15 | 1389 | 5 | 8 | 26887 | 29447 | 115925 | 1 d 12 h 39 m |
| 16 | 3064 | 5 | 9 | 49466 | 74576 | 267789 | 12 d 15 h 24 m |

Table 14: The characteristics of the final set of gradient-add-or-remove and its relaxed version for $8 \leq n \leq 16$, starting either from $1\cdots1$ or from $0\cdots0$. In the relaxed version, 1000 moves keeping the function $\phi_2$ unchanged are authorized.

## 54.4 Combination of add and remove

The previous attempts were aimed at proving the global inequality on the expected number of pivots of second order for an increasing set $A$, without any constraint on its cardinality. In particular, the modification procedures were allowed to change the cardinality of the set $A$ under consideration. In order to prove an associated deviations inequality, we should consider a modification scheme which keeps the cardinality of $A$ unchanged. Moreover, if we manage to obtain a valuable inequality with a cardinality constraint, we could use it as an intermediate step to obtain a global inequality. So we consider here an operation which removes a minimal element of $A$ and simultaneously add a maximal element of $A^c$. This corresponds to the system of neighborhoods defined by

$$\mathcal{V}_{\text{add and remove}}(A) = \big\{ A \cup \{\eta\} \setminus \{\sigma\} : \sigma \in \text{Minimal}(A), \eta \in \text{Maximal}(A^c) \big\}.$$

We call gradient-aar the associated gradient algorithm. To ensure the termination of the gradient-aar algorithm, we consider only the strict variant, in which a move is performed only when it leads to an increase of $\phi_2$, otherwise the algorithm could cycle indefinitely. The results of the exhaustive exploration for $3 \leq n \leq 7$ are presented in table 15.

| $n$ | examined | stuck | violation | run time |
|---|---|---|---|---|
| 3 | 8 | 7 | 0 | 0 s |
| 4 | 28 | 15 | 0 | 0 s |
| 5 | 208 | 36 | 0 | 0 s |
| 6 | 16351 | 221 | 0 | 4 s |
| 7 | 490013146 | 135083 | 38 | 13 d 9 h 53 |

Table 15: The number of local maxima of gradient-aar for $3 \leq n \leq 7$

A set $A$ is stuck if

$$\forall B \in \mathcal{V}_{\text{add and remove}}(A) \qquad \phi_2(B) \leq \phi_2(A),$$

and it is a violation if in addition

$$\max(A^c) > \min(A) + 1,$$

where the numbers max and min are defined as

$$\min(A) = \min\big\{ |\text{supp}\,\omega| : \omega \in A \big\}, \quad \max(A^c) = \max\big\{ |\text{supp}\,\omega| : \omega \notin A \big\}.$$

The situation is fine for $n \leq 6$, but 38 problematic sets exist for $n = 7$. There is some hope, because these 38 sets all satisfy $\max(A^c) = \min(A) + 2$. Two examples of these problematic sets are presented in figure 125. The main problem with these sets is that the gap between $\min(A)$ and $\max(A^c)$ increases. The good point is that all these violating sets have a boundary which is quite close to the middle layer of the hypercube.

$$\mathrm{Minimal}(A) = \{\,1111000, 1110100, 1110010, 1110001, 1101110,\\ 1101101, 1101011, 1010111, 1001111\,\}, \quad \min(A) = 4,\ \max(A^c) = 6$$

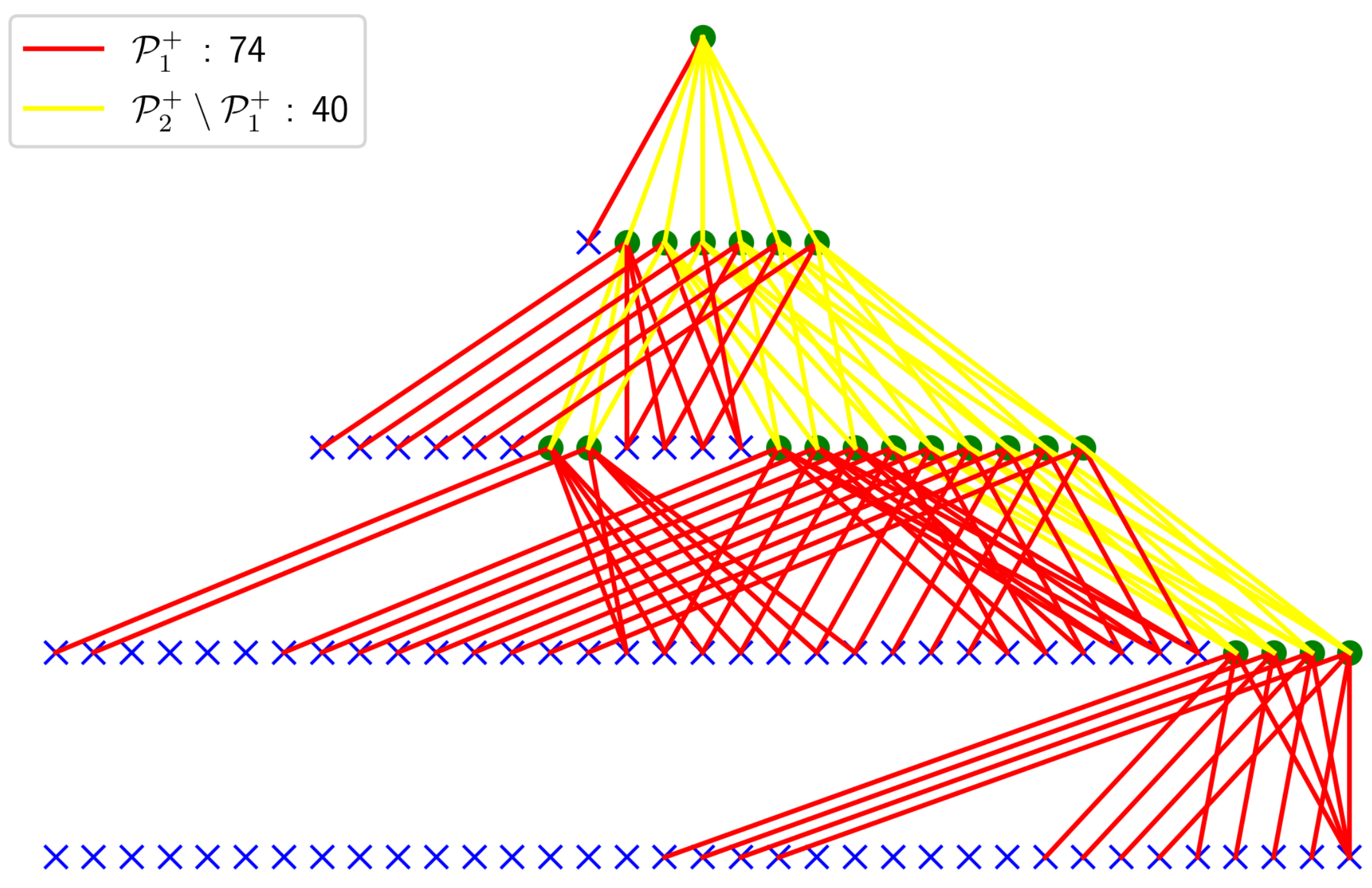

$$\mathrm{Minimal}(A) = \{\,1110000, 1101000, 1100100, 1100010, 1100001, 1011100, 1011010,\\ 1010110, 1010011, 1001101, 1001011, 0111100, 0111001, 0110101,\\ 0101110, 0101011, 0100111\,\}, \quad \min(A) = 3,\ \max(A^c) = 5$$

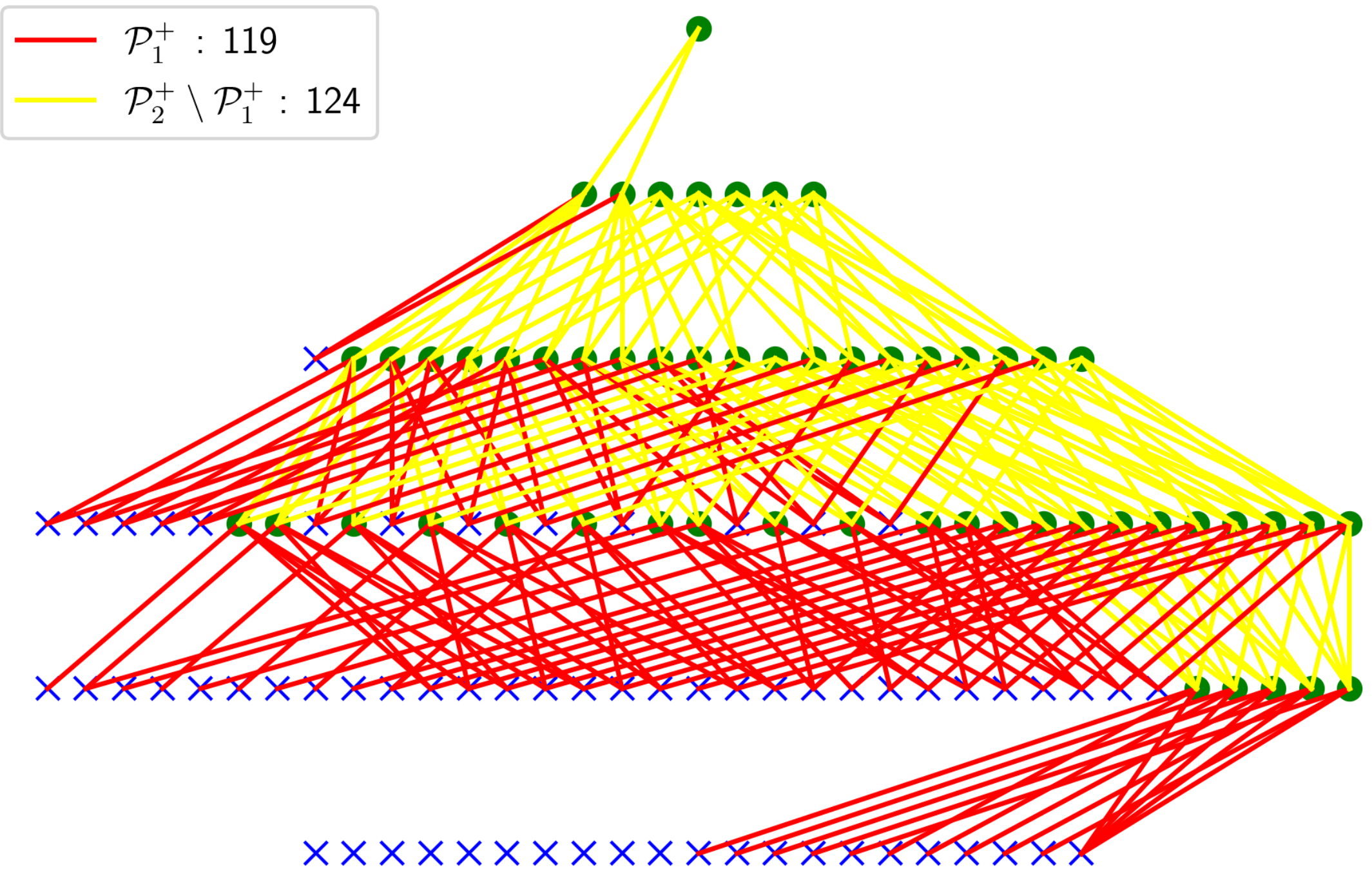

Figure 125: Two violating local maxima of gradient-add-and-remove for $n = 7$

## 54.5 Line modifications

The computation of the variation of $\phi_2$ induced by the removal or the addition of a configuration leads to complicated formulas (see (54.9) and (54.15)). So we should consider other modification schemes. Our aim and hope is to prove that the maxima of $\phi_2$ (or at least one of them) can be found among the sets $A$ satisfying $\max(A^c) - \min(A) \leq 1$ (where the numbers $\min(A)$ and $\max(A^c)$ are defined in (54.1)), or even such that

$$\frac{n}{2} - 1 \leq \min(A) , \quad \max(A^c) \leq \frac{n}{2} + 1 .$$

Therefore we consider modifications of the set $A$ that bring the set closer to these sets. A natural attempt consists in removing the lowest layer of the set $A$ if it belongs to the lower part of the hypercube $\mathcal{H}^n_-$ or adding the highest layer of $A^c$ if it belongs to the upper part of the hypercube $\mathcal{H}^n_+$. Interestingly, there exists a nice proof of Sperner's theorem which is based on such a modification scheme (see the very first theorem 1.1 in [56]).

**Lemma 54.6.** *Let $A$ be an increasing subset of $\{0,1\}^n$. Let $m$ be the index of the lowest layer intersected by $A$, i.e.,*

$$m = \min\{\, |\operatorname{supp} \sigma| : \sigma \in A \,\} ,$$

*and let $S$ be the configurations of the lowest layer of $A$, i.e.,*

$$S = \big\{\, \sigma \in A : |\operatorname{supp} \sigma| = m \,\big\} .$$

*The variation of $\phi_2$ between $A$ and $A \setminus S$ is given by*

$$\begin{aligned} \phi_2(A\setminus S) - \phi_2(A) = (n-2m)|S| + \sum_{\eta \in \{\, \sigma^i : \sigma\in S,\, i \notin \operatorname{supp}\sigma \,\}} \Big( \big|\mathcal{P}^+_{2\setminus 1}(A\setminus S, \eta)\big| - \big|\mathcal{P}^+_{2\setminus 1}(A, \eta)\big| \Big) \\ + \sum_{\eta \in \{\, \sigma^{ij} : \sigma\in S,\, i,j \notin \operatorname{supp}\sigma,\, i\neq j \,\}} \Big( \big|\mathcal{P}^+_{2\setminus 1}(A \setminus S, \eta)\big| - \big|\mathcal{P}^+_{2\setminus 1}(A, \eta)\big| \Big) . \end{aligned} \tag{54.19}$$

*The first sum in formula* (54.19) *is non-positive, while the second sum in formula* (54.19) *is non-negative.*

The proof of lemma 54.6 follows the same lines as the proof of lemma 54.2, so we do not rewrite it here. Unfortunately, the formula (54.19) is still complicated and we did not manage to use it efficiently. Even worse, there exist increasing sets $A$ such that $2m < n - 1$, and for which the removal of the lowest layer induces a decrease of the value of $\phi_2$! Figure 126 presents such an example, in fact it is the very same set which already appeared in figures 119,120. Therefore it does not seem sensible to hope to prove the conjecture 53.1 with the help of line modifications.

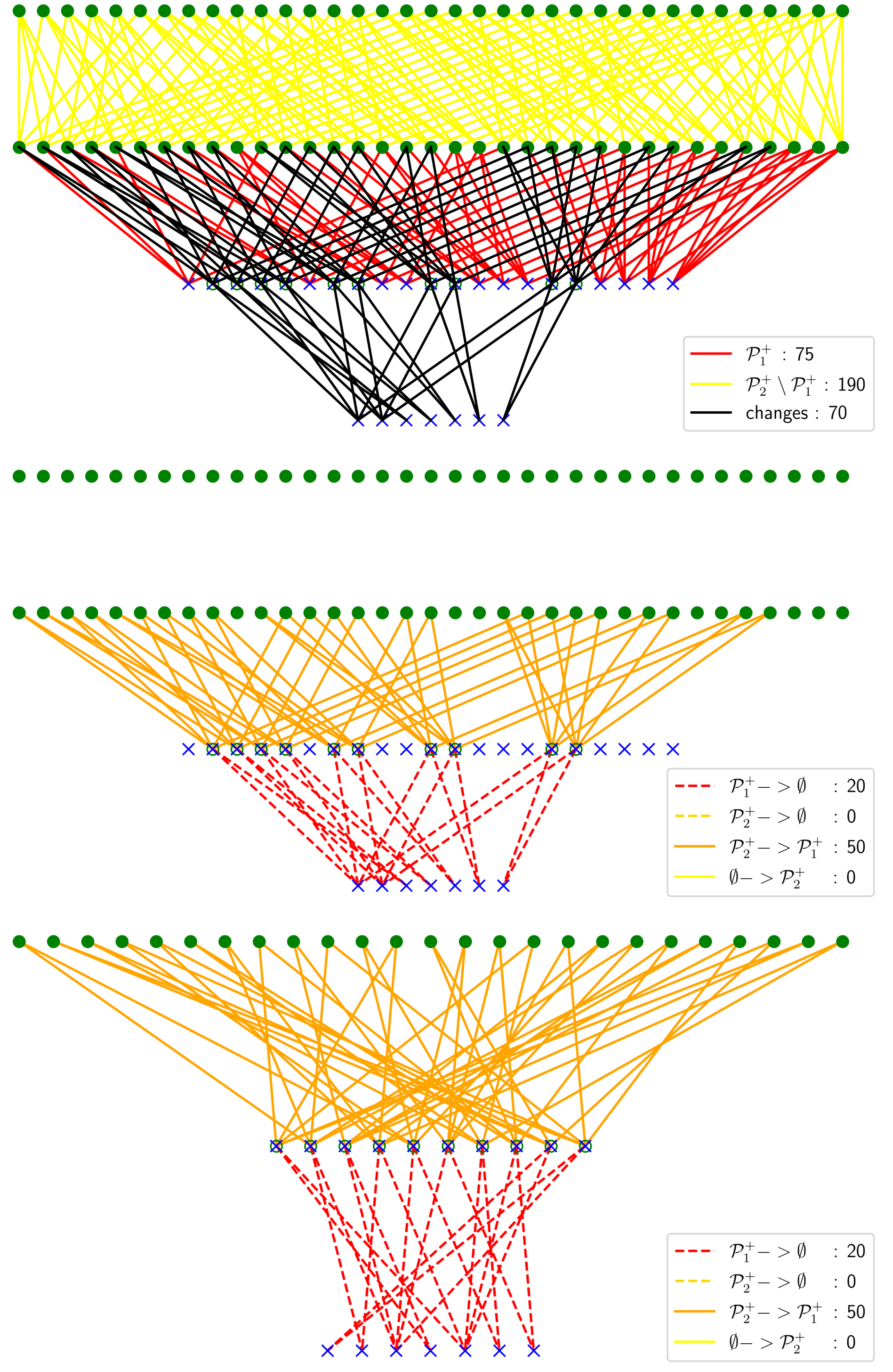

Figure 126: Same set as figures 119,120, removal of the lowest line: $\Delta\phi_2 = -30$. Top: the set $A$. Middle: the modifications (same representation as top). Bottom: focus on the modifications (redeployed horizontally for clarity).

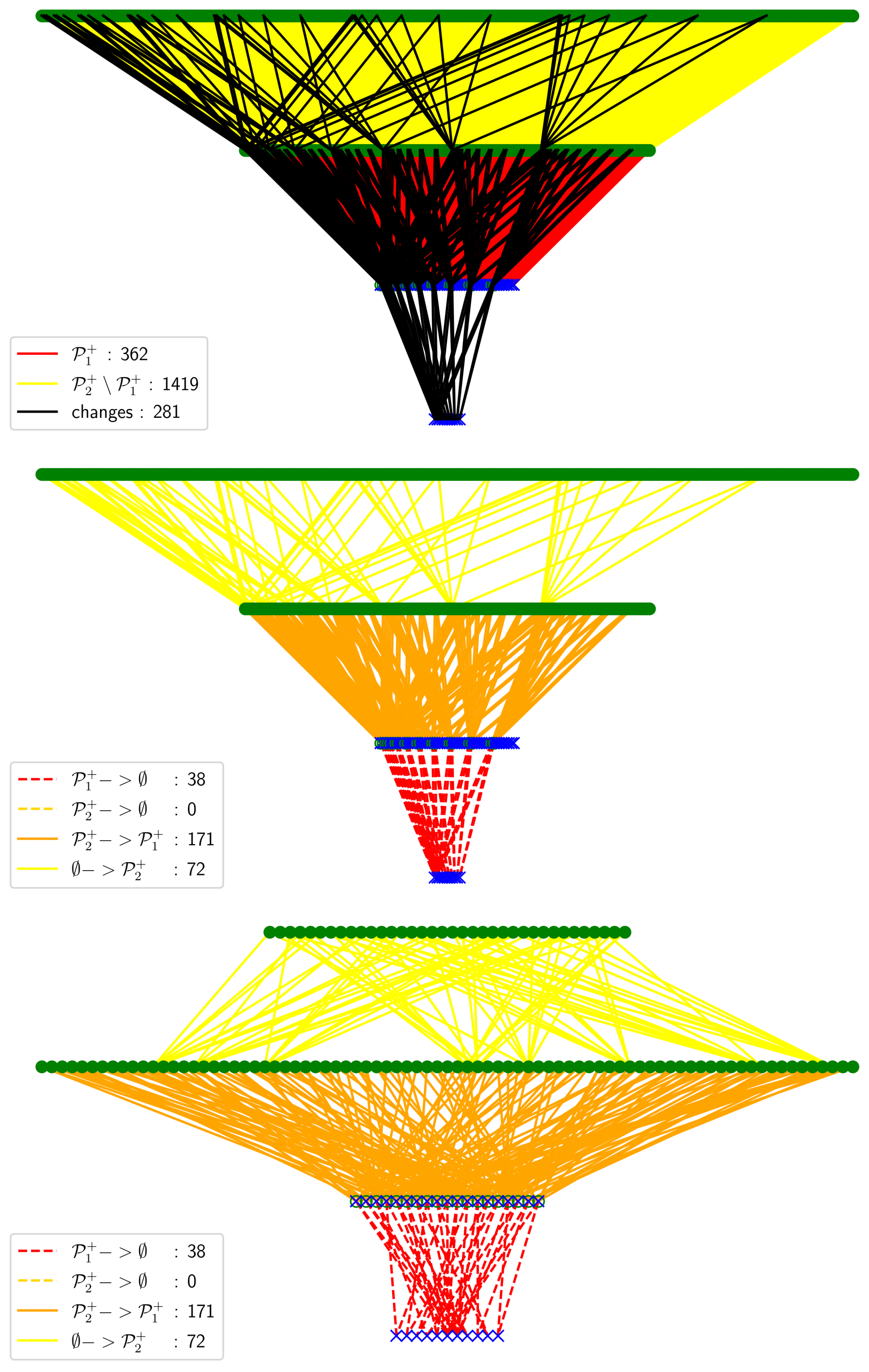

Figure 127: Same set as figures 121,122, removal of the lowest line: $\Delta\phi_2 = -33$. Top: the set $A$. Middle: the modifications (same representation as top). Bottom: focus on the modifications (redeployed horizontally for clarity).

## 54.6 Pushing modifications

The line modifications considered in the previous subsection 54.5 did not work out. We consider here a more complicated scheme of modification, which looks promising. The idea is to remove at once all the set Minimal($A$) from $A$ and to add at the same time to $A$ the set Maximal($A^c$). Of course, we should be a bit more careful if we wish that the modifications do not decrease the value of $\phi_2$, namely we should remove only configurations belonging to the lower part of the hypercube, and we should add only configurations belonging to the upper part of the hypercube. We proceed next with a more precise definition. Let $A$ be an increasing subset of $\{0,1\}^n$. We define the sets of configurations $S$ and $T$ as

$$\begin{aligned} S &= \Big\{\, \sigma \in \text{Minimal}(A) : |\operatorname{supp}\sigma| \leq n/2-1 \,\Big\}\,, \\ T &= \Big\{\, \sigma \in \text{Maximal}(A^c) : |\operatorname{supp}\sigma| \geq n/2+1 \,\Big\}\,. \end{aligned}$$

Our modification scheme would transform the set $A$ into $A\cup T\setminus S$. The hope is that this modification pushes the boundary of $A$ towards the middle layer of the hypercube in a way that the second order pivots which might have been destroyed are recreated nearby. We can try to compute the variation of $\phi_2$. We obtain a formula similar to (54.19), albeit slightly more complicated. Let us introduce the neighborhoods $M_1$ and $M_2$ of $S\cup T$ defined by

$$\begin{aligned} M_1 &= \Big\{\, \sigma^i : \sigma \in S\cup T,\, i \notin \operatorname{supp}\sigma \,\Big\}\,, \\ M_2 &= \Big\{\, \sigma^{ij} : \sigma \in S\cup T,\, i,j \notin \operatorname{supp}\sigma,\, i\neq j \,\Big\}\,. \end{aligned}$$

The formula for the variation of $\phi_2$ is then

$$\begin{aligned} \phi_2(A\cup T\setminus S)-\phi_2(A) \,=\, &\sum_{\sigma\in S}(n-2|\operatorname{supp}\sigma|)-\sum_{\sigma\in T}(n-2|\operatorname{supp}\sigma|) \\ &+\sum_{\sigma\in M_1}\Big(\big|\mathcal{P}^+_{2\setminus 1}(A\cup T\setminus S,\sigma)\big|-\big|\mathcal{P}^+_{2\setminus 1}(A,\sigma)\big|\Big) \\ &\quad+\sum_{\sigma\in M_2\setminus M_1}\Big(\big|\mathcal{P}^+_{2\setminus 1}(A\cup T\setminus S,\sigma)\big|-\big|\mathcal{P}^+_{2\setminus 1}(A,\sigma)\big|\Big)\,. \end{aligned} \tag{54.20}$$

We did not manage to exploit formula (54.20) so far, but there is still hope that it might be useful in some cases. Furthermore, it seems wise to decompose further this strategy into two slightly easier questions. We should first consider increasing sets $A$ which contain completely the upper part of the hypercube, or even a set of the form

$$\Big\{\, \omega\in\{0,1\}^n : |\operatorname{supp}\omega| > \frac{n}{2}+k \,\Big\}\,,$$

where $k$ is a fixed number. The action of the modification scheme on these sets is simpler, as it will not add any configuration, it will only remove the minimal configurations. If we managed to control the variation of $\phi_2$ under this scheme, then we could afterwards consider sets $A$ which are included in the upper part of the hypercube. In a final step, we would study the remaining sets $A$, for which the modification scheme is more complicated.

## 55 Maximizers of $E(|\mathcal{P}_2^+(A)|\,|\,A)$

We did not manage so far to derive an interesting inequality on the second order pivots. In this section, in order to gain some understanding on the structure of the functional $\phi_2$, we adopt an experimental methodology. We fix a small value of $n \leq 6$ and an integer $m$ such that $0 \leq m \leq 2^n$. We compute all the sets $A$ realizing the maximum of the problem

$$M_2(n,m) \;=\; \max\,\big\{\,\phi_2(A) : A \text{ increasing subset of } \{0,1\}^n,\, |A| = m\,\big\}\,. \quad (55.1)$$

We work here with $p = 1/2$, thus the cardinality of a set $A$ and its probability are proportional, and the problem (55.1) amounts to maximizing the conditional probability $E(|\mathcal{P}_2^+(A)|\,|\,A)$. We finally arrange the maximizers into a sequence, corresponding to increasing values of $m$. We put first the maximizers for $m = 0$, then for $m = 1$, and so on. Of course, we consider only the maximizers up to equivalence modulo a permutation of the variables. In case there exist several maximizers for some value of $m$, the sequence will contain several sets having cardinality $m$. Ideally, we would like to observe a structural property of the resulting sequence that would allow us to obtain a result valid for any $n \geq 2$. A first problem is to find an adequate representation of the sequence of sets. We use the same representation as for the visualization of the maxima of $\phi_2$, with some minor refinements. The vertices of the sets are colored in green. Among these, the vertices which cannot be removed are drawn with a circle. When a vertex can be removed and its removal would lead to an increase (respectively a decrease) of $\phi_2$, it is drawn with a square (respectively a triangle). Moreover the variation of $\phi_2$ associated with the vertex removal is written just above the vertex. The vertices of the complements of the sets are colored in blue. Among these, the vertices which cannot be added are drawn with an inclined cross. When a vertex can be added and its addition would lead to an increase (respectively a decrease) of $\phi_2$, it is drawn with a vertical cross (respectively a star). Moreover the variation of $\phi_2$ associated with the vertex addition is written just above the vertex. The first order pivots are represented by red segments. The second order pivots are represented by yellow segments.

Figure 128 presents the results for $n = 3$. This is just a warming up, in order to help checking the previous conventions. Figures 129, 130 present the results for $n = 4$. For $n = 4$, we keep all the information on the pictures, but for $n = 5$ and $n = 6$, we have to reduce the size of the numerous pictures, so we remove the information previously attached to the vertices. Figure 131 presents the results for $n = 5$. Figures 132, 133 present the results for $n = 6$. These initial results are very encouraging. Indeed, we see that the sequence of the maximizers is essentially obtained by filling progressively the hypercube in a natural manner, starting from the top and completing successively each layer, up to some minor exceptions. Unfortunately, we do not see a simple pattern emerging from the pictures, nor do we understand properly the filling mechanism. Of course, one should be very careful before drawing any definite conclusion, indeed we have only the full results for $n \leq 6$, so it would be quite dangerous to extrapolate a result valid for any value of $n$!

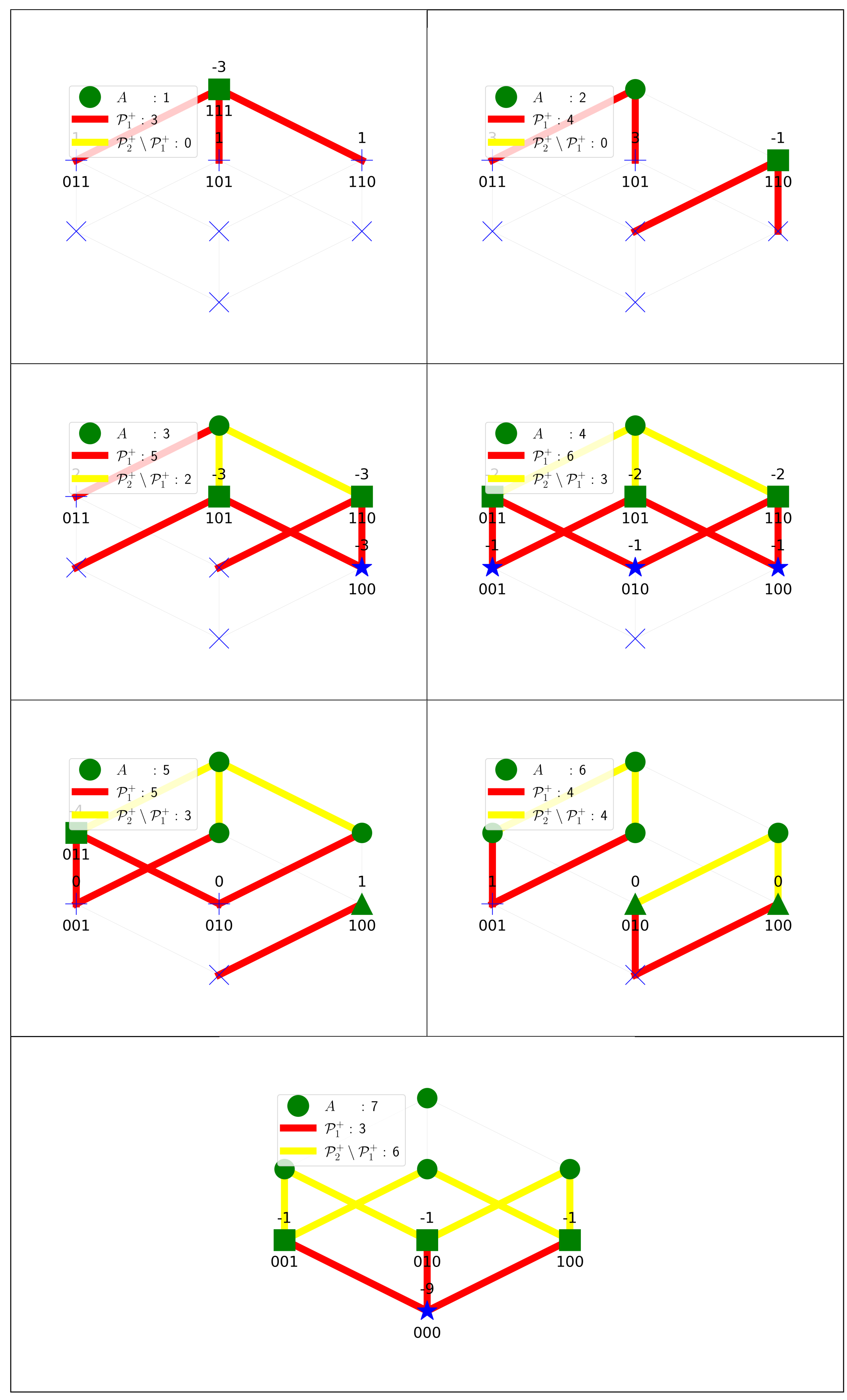

Figure 128: Maximizers of $E(|\mathcal{P}_2^+(A)|\,|\,A)$ in $\{0,1\}^3$

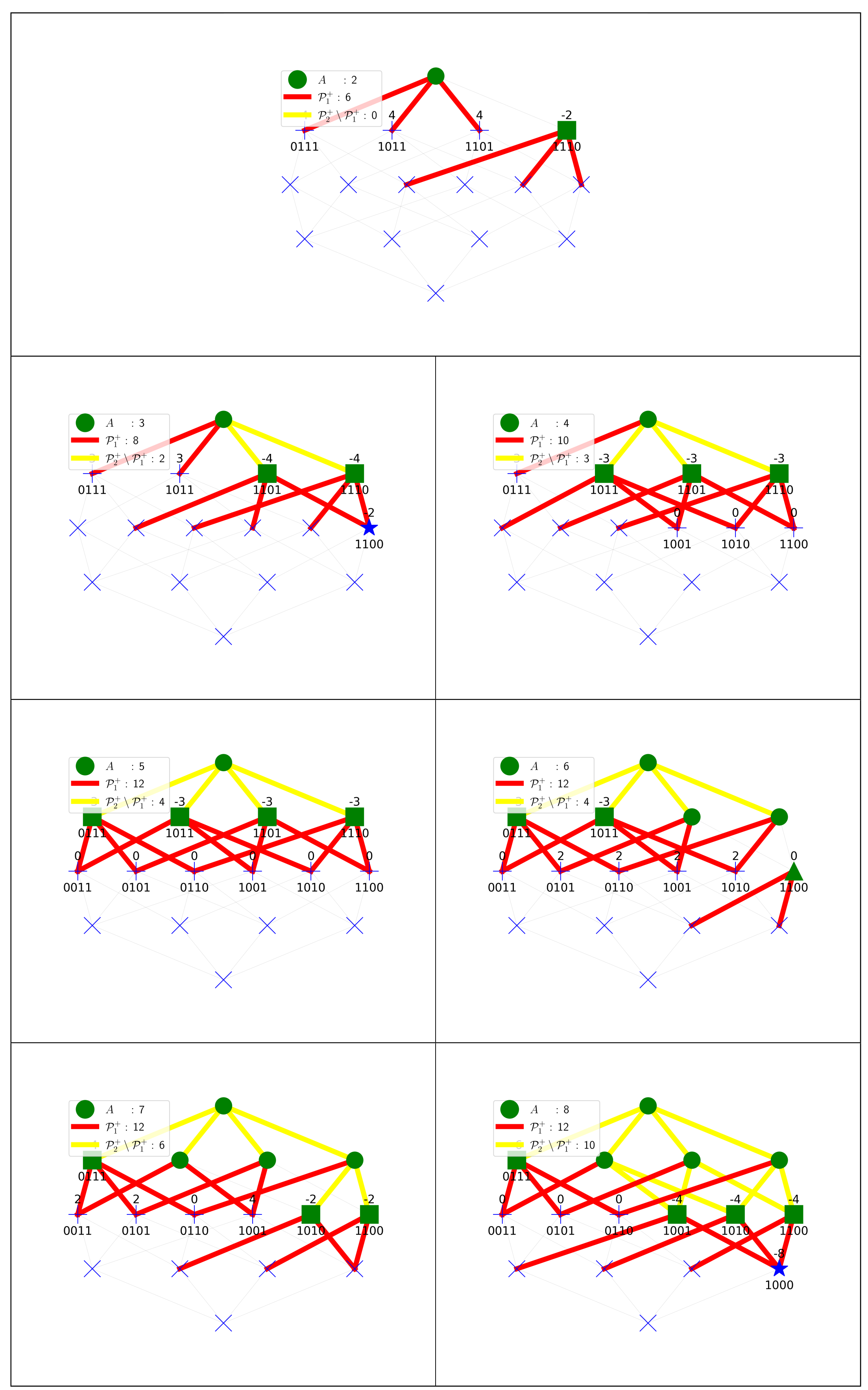

Figure 129: Maximizers of $E(|\mathcal{P}_2^+(A)|\,|\,A)$ in $\{0,1\}^4$

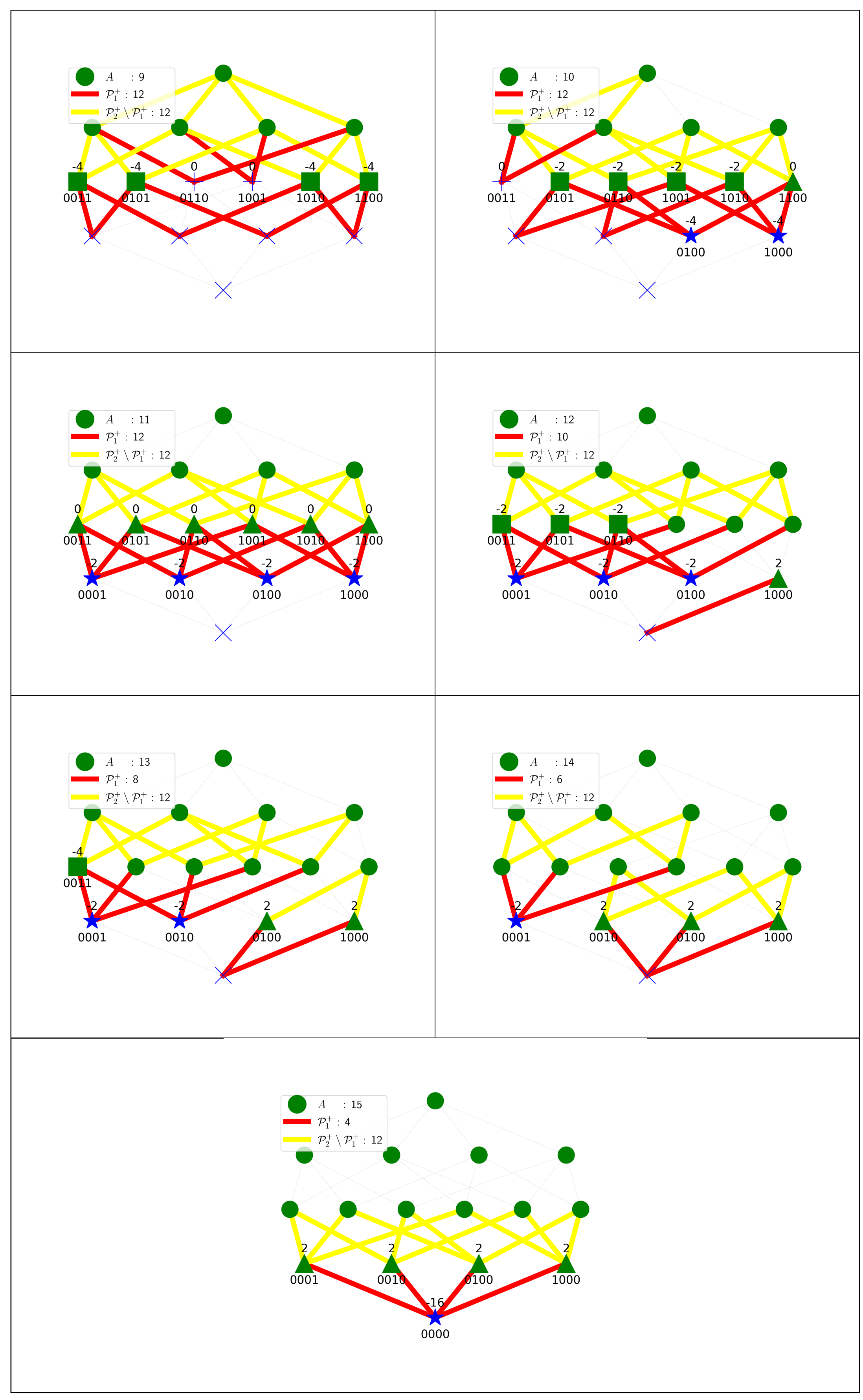

Figure 130: Continuation of figure 129

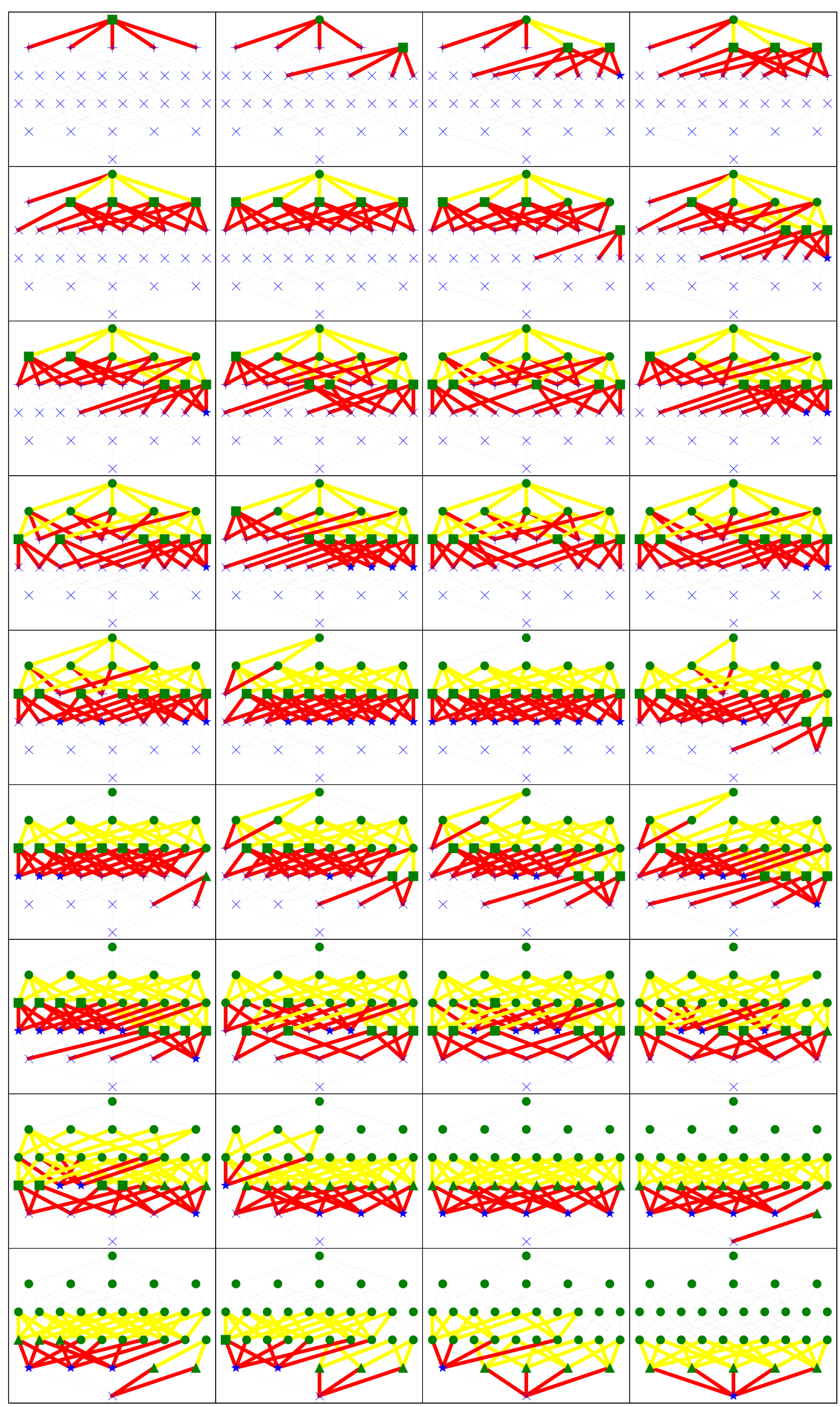

Figure 131: Maximizers of $E(|\mathcal{P}_2^+(A)|\,|\,A)$ in $\{\,0,1\,\}^5$

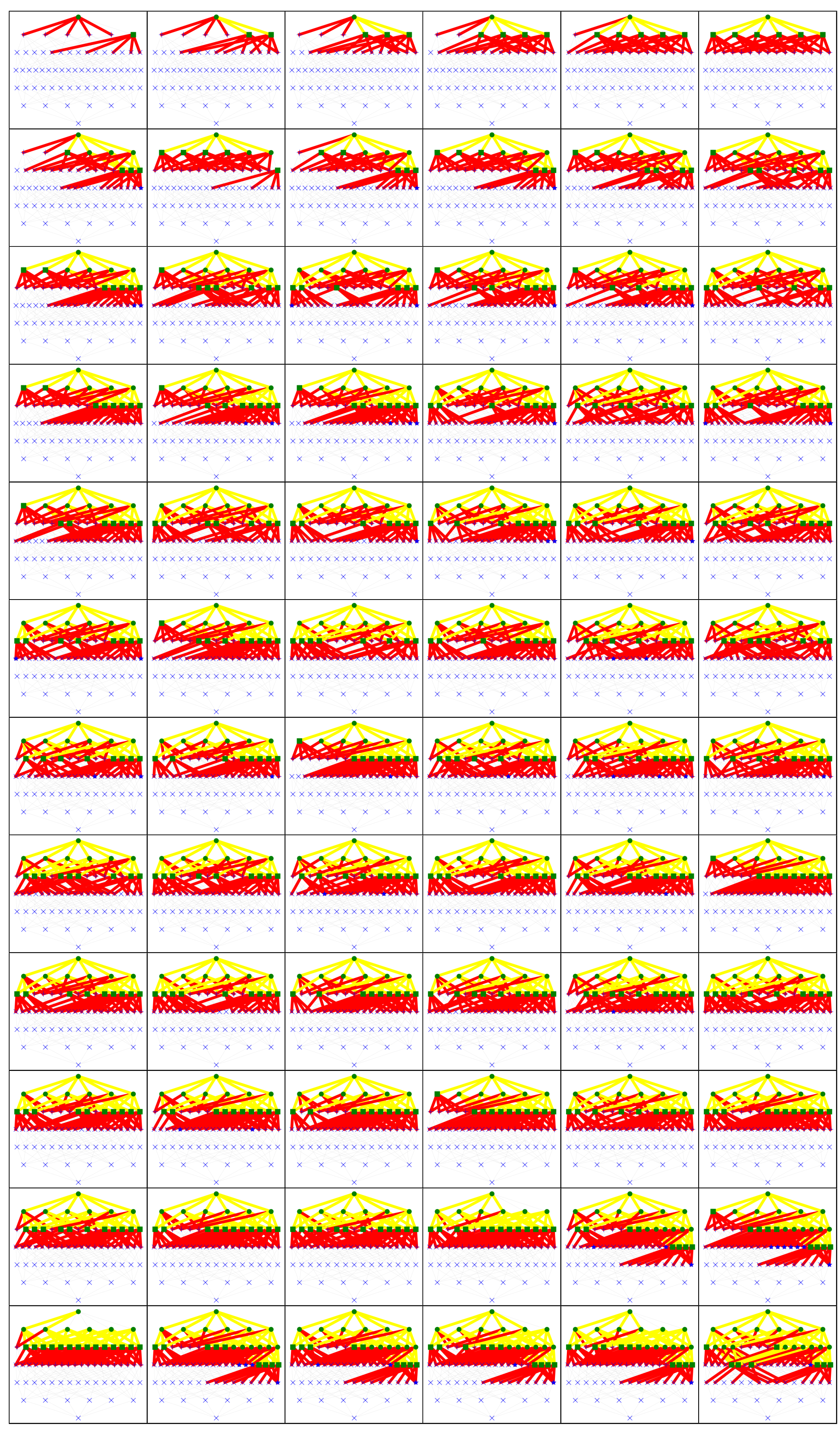

Figure 132: Maximizers of $E(|\mathcal{P}_2^+(A)|\,|\,A)$ in $\{\,0,1\,\}^6$

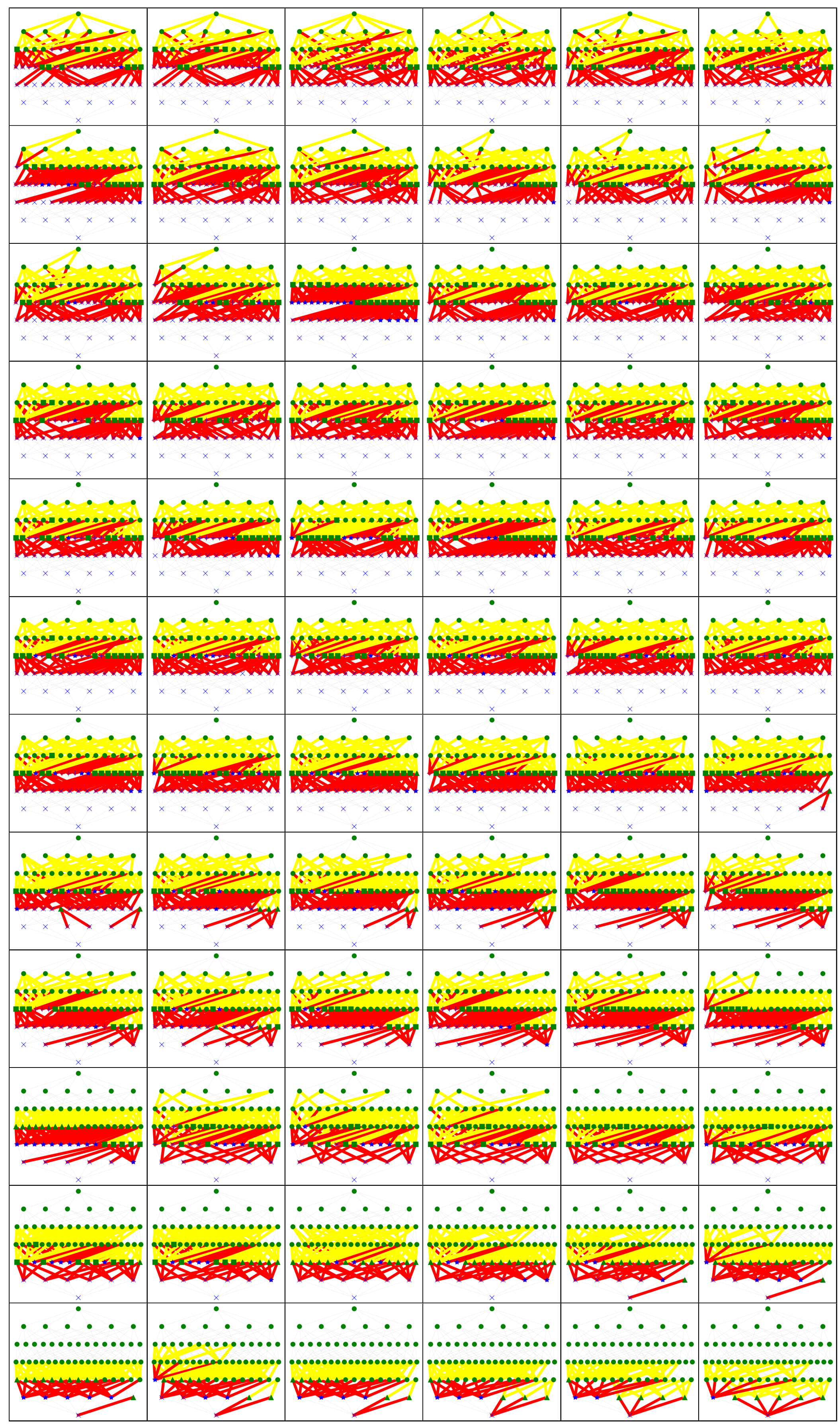

Figure 133: Continuation of figure 132

# 56 Maximizers of $E(|\mathcal{P}_3^+ \setminus \mathcal{P}_2^+|)$

We investigate here the pivots of third order with the help of numerical computations. For $A$ an increasing subset of $\{0,1\}^n$, we define

$$\phi_{3\setminus 2}(A) \,=\, \sum_{\omega\in A} \big|\mathcal{P}_3^+(A,\omega)\setminus\mathcal{P}_2^+(A,\omega)\big|\,.$$

We consider the problem

$$M_{3\setminus 2}(n) \,=\, \max\,\big\{\,\phi_{3\setminus 2}(A) : A \text{ increasing subset of } \{0,1\}^n\,\big\}\,.$$

We performed an exhaustive analysis of all the monotone boolean functions in 6 or less variables. Table 16 presents the results. The maximizers are depicted in layer representation in figures 134, 135 and 136.

| $n$ | $\lvert A\rvert$ | $\phi_2$ | $\phi_{3\setminus 2}$ | Minimal$(A)$ |
|---|---|---|---|---|
| 3 | 7 | 9 | 3 | 100 010 001 |
| 4 | 15 | 16 | 12 | 1000 0100 0010 0001 |
| 5 | 30 | 32 | 30 | 10000 01000 00100 00010 |
| | 31 | 25 | 30 | 10000 01000 00100 00010 00001 |
| 6 | 60 | 64 | 72 | 100000 010000 001000 000100 |
| | 61 | 62 | 72 | 100000 010000 001000 000100 000011 |

Table 16: The maximizers of $E(|\mathcal{P}_3^+ \setminus \mathcal{P}_2^+|)$ for $n = 3,4,5,6$

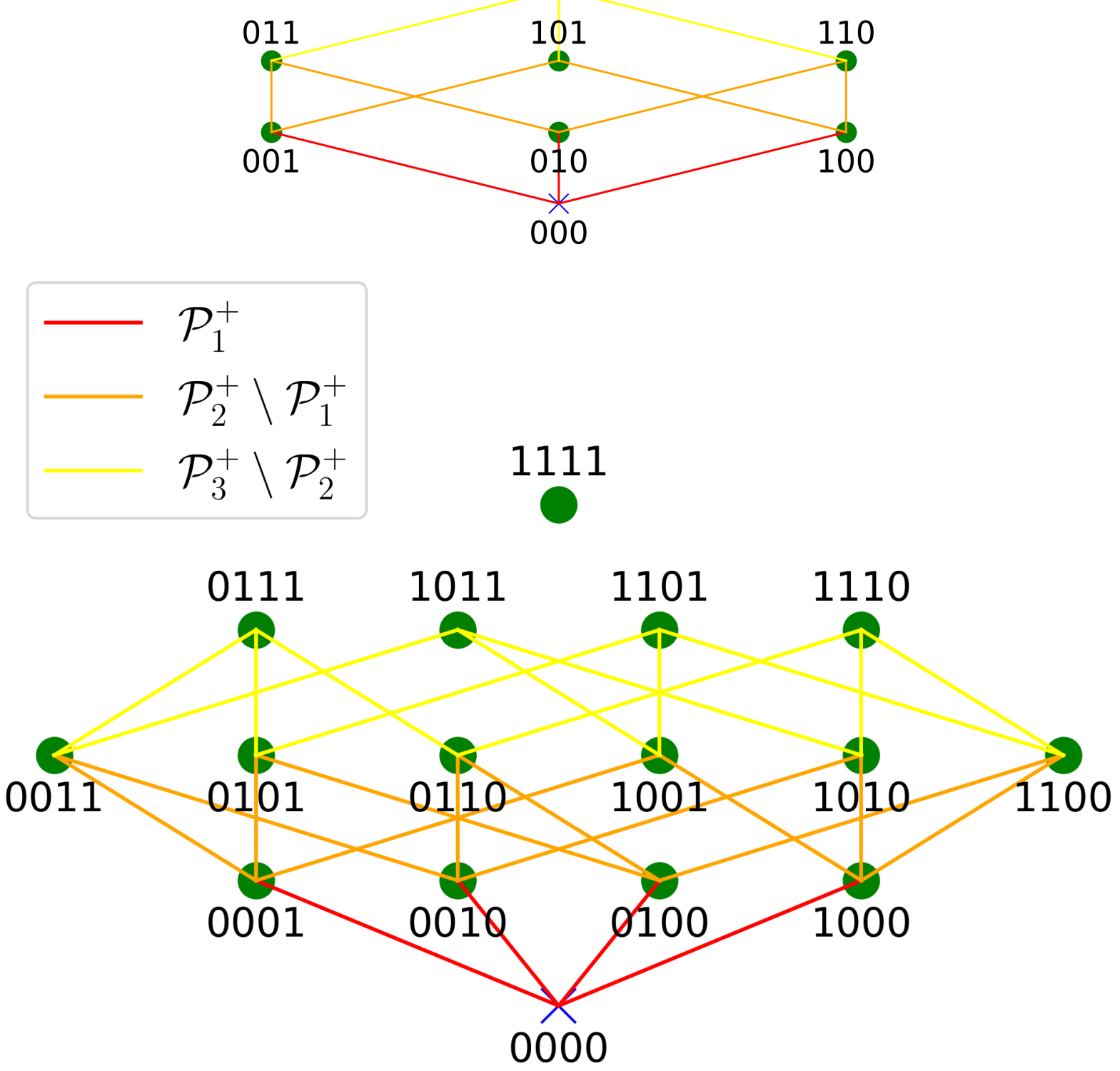

Figure 134: The maximizers of $E(|\mathcal{P}_3^+ \setminus \mathcal{P}_2^+|)$ for $n = 3,4$

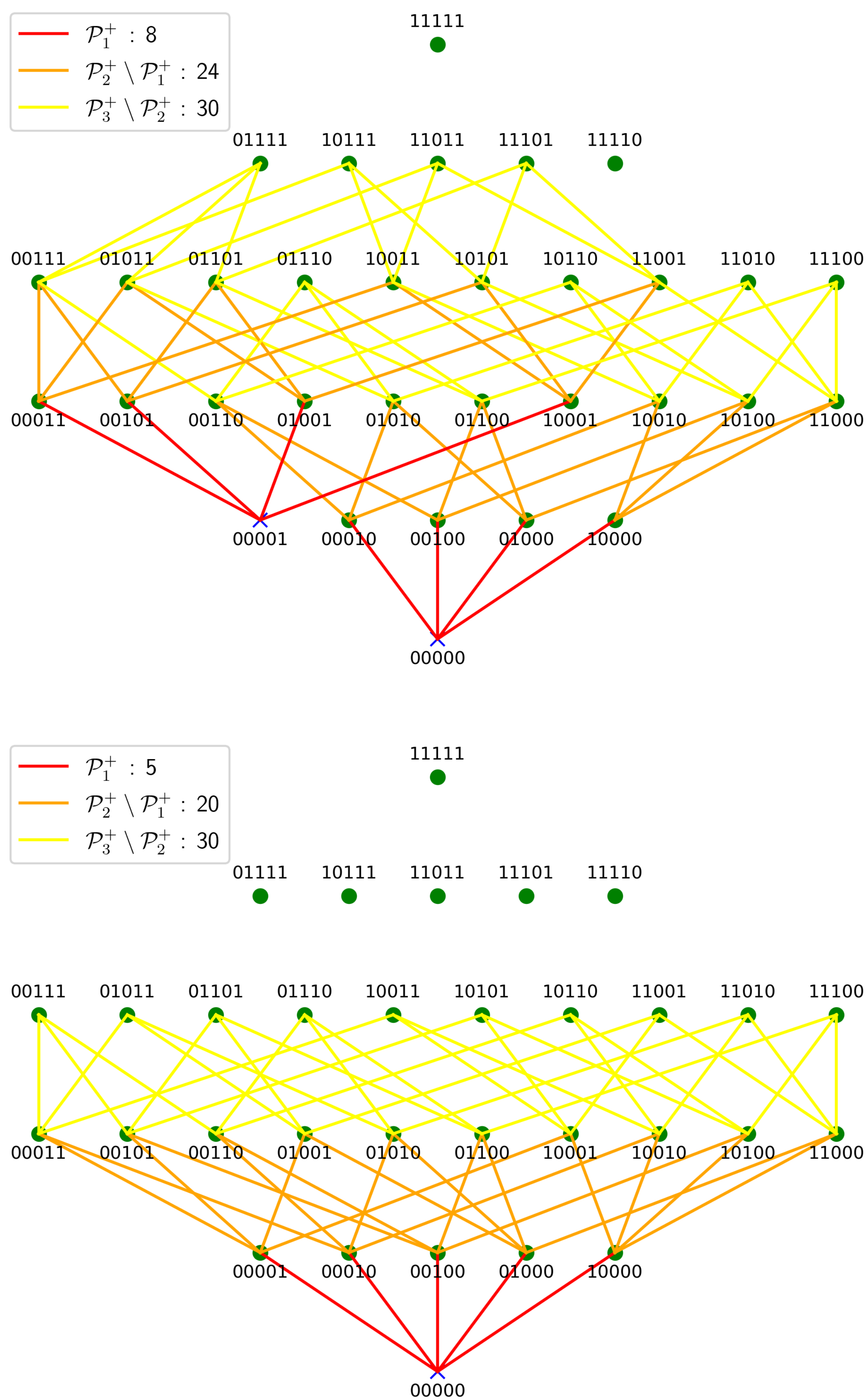

Figure 135: The maximizers of $E(|\mathcal{P}_3^+ \setminus \mathcal{P}_2^+|)$ for $n = 5$

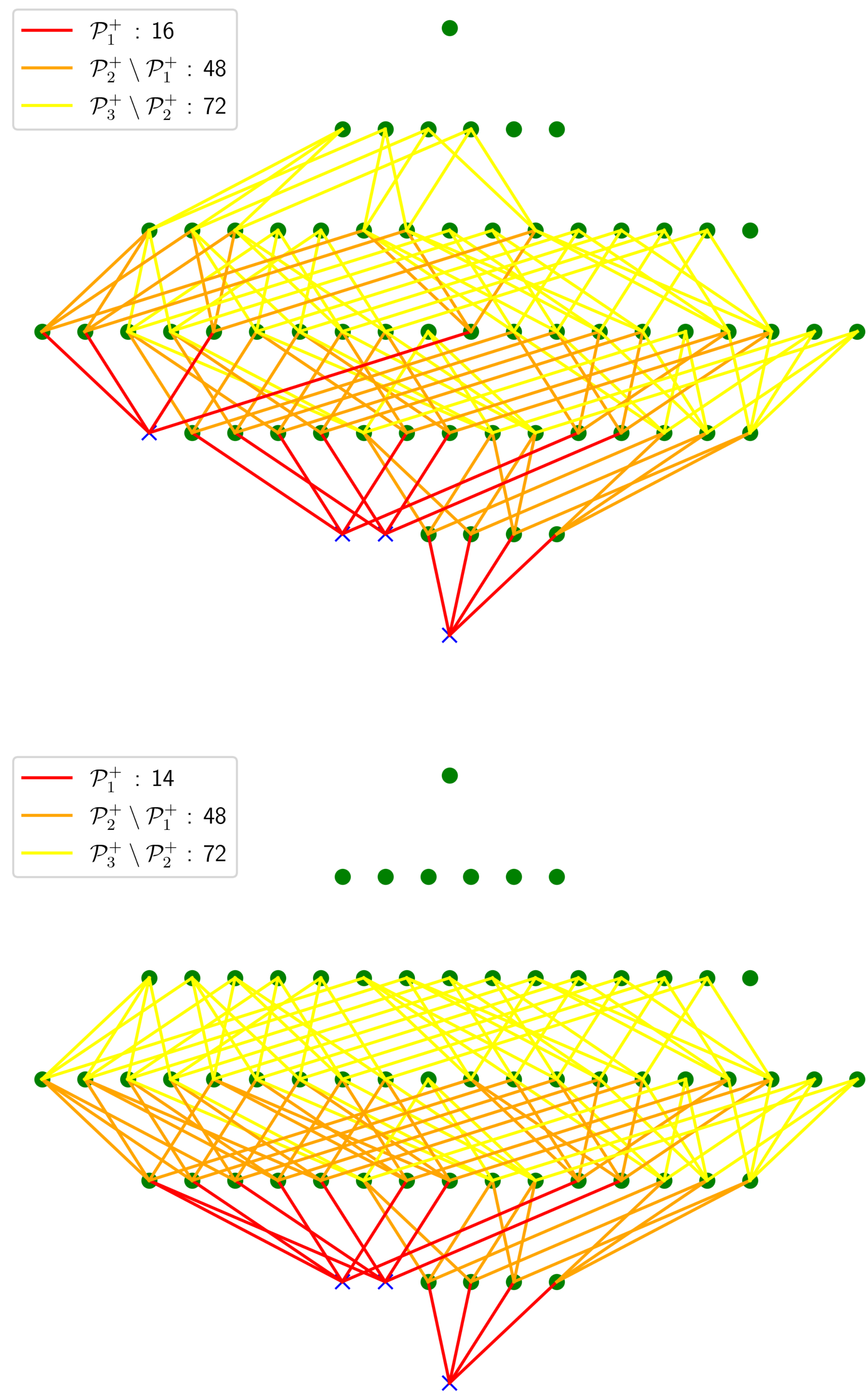

Figure 136: The maximizers of $E(|\mathcal{P}_3^+ \setminus \mathcal{P}_2^+|)$ for $n = 6$

# 57 Maximizers of $E(|\mathcal{P}_3^+|)$

We continue the numerical investigations around the pivots of third order, this time we group together all the pivots until order three. For $A$ an increasing subset of $\{0,1\}^n$, we define

$$\phi_3(A) \,=\, \sum_{\omega\in A} \big|\mathcal{P}_3^+(A,\omega)\big|\,.$$

We consider the problem

$$M_3(n) \,=\, \max\,\big\{\,\phi_3(A) : A \text{ increasing subset of } \{0,1\}^n\,\big\}\,.$$

We performed an exhaustive analysis of all the monotone boolean functions in 6 or less variables. Table 17 presents the results. The maximizers are depicted in layer representation in figures 137, 138 and 139.

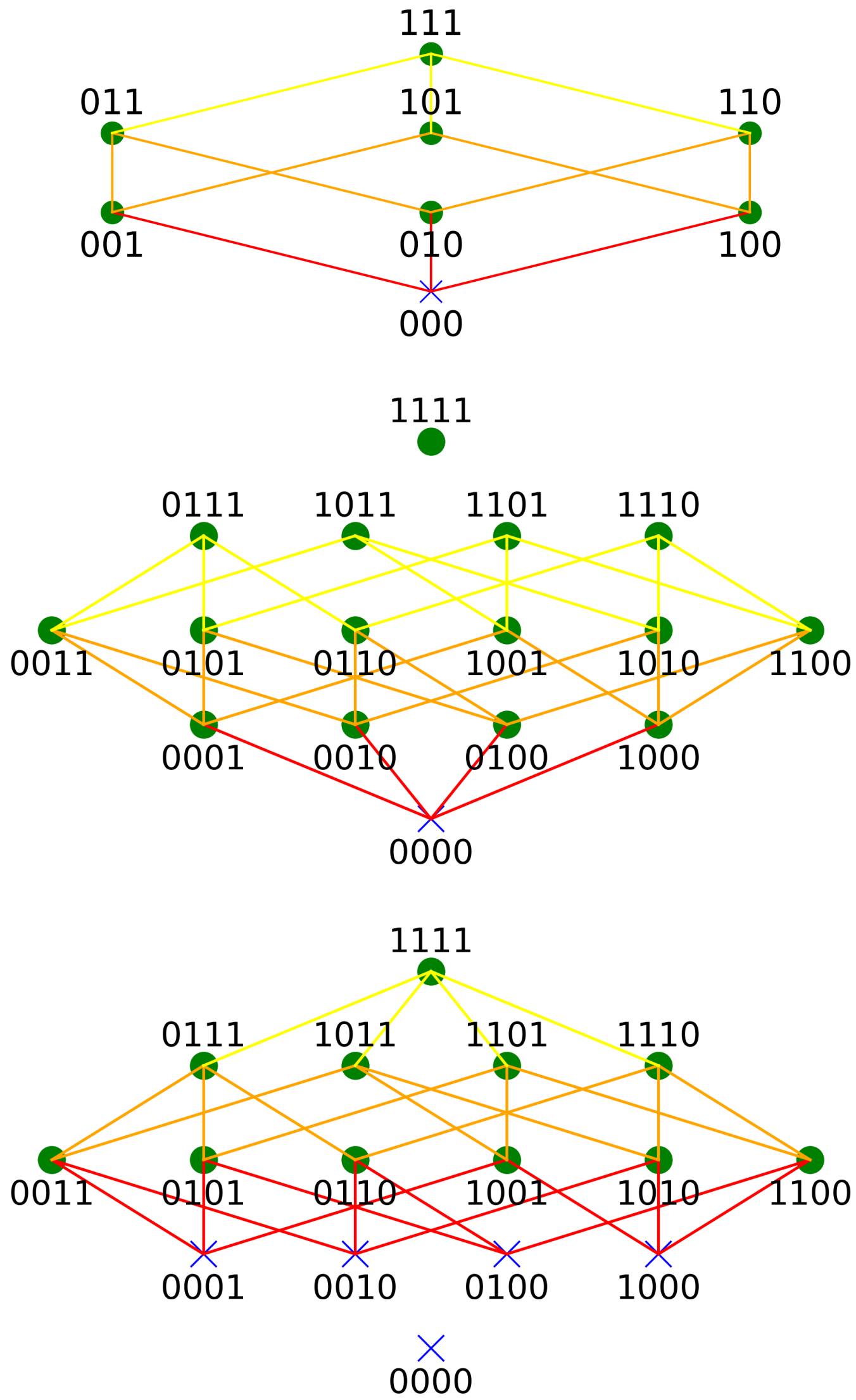

Figure 137: Maximizers of $E(|\mathcal{P}_3^+|)$ for $n=3,4$. First , second , third pivots.

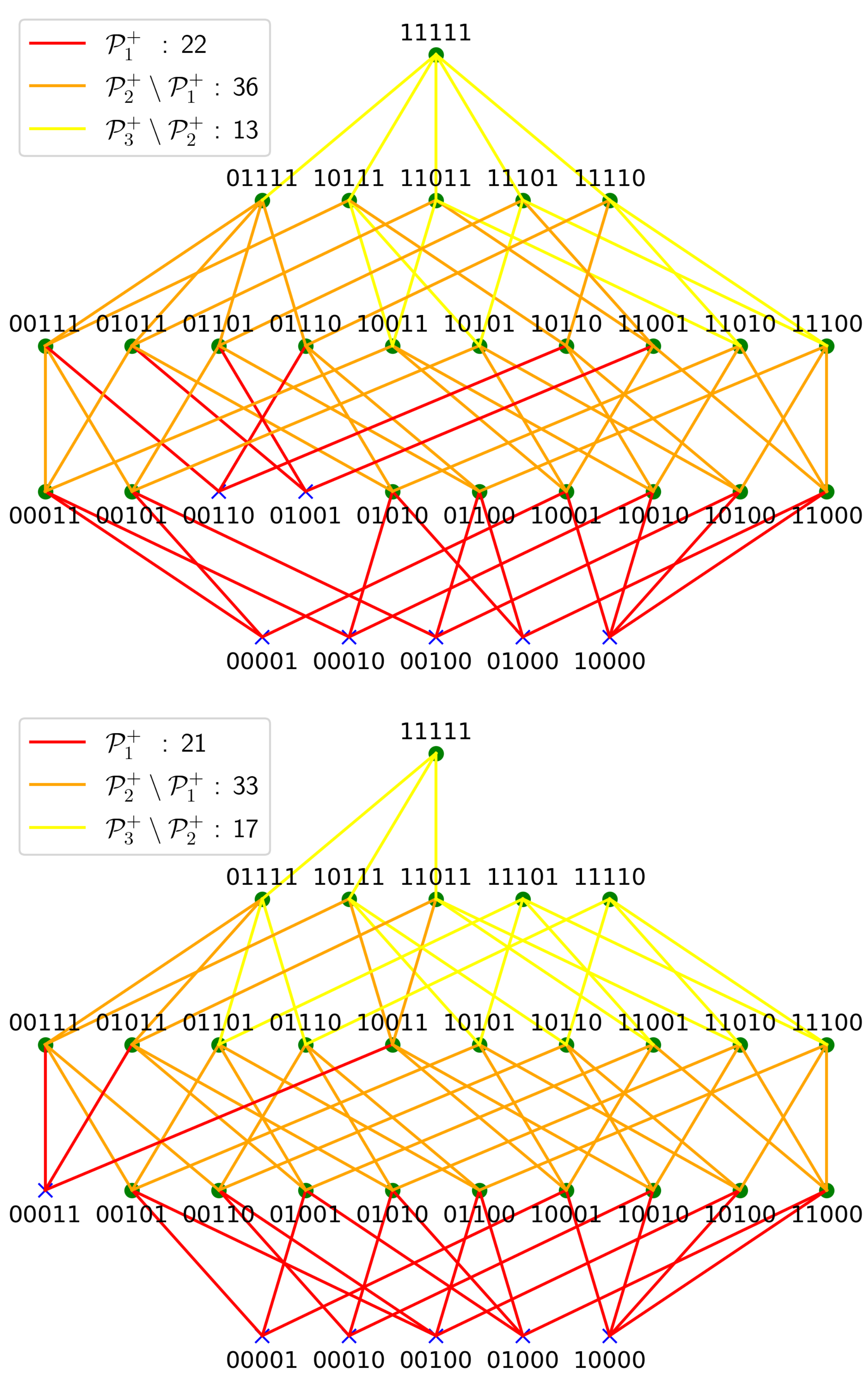

Figure 138: The maximizers of $E(|\mathcal{P}_3^+|)$ for $n = 5$

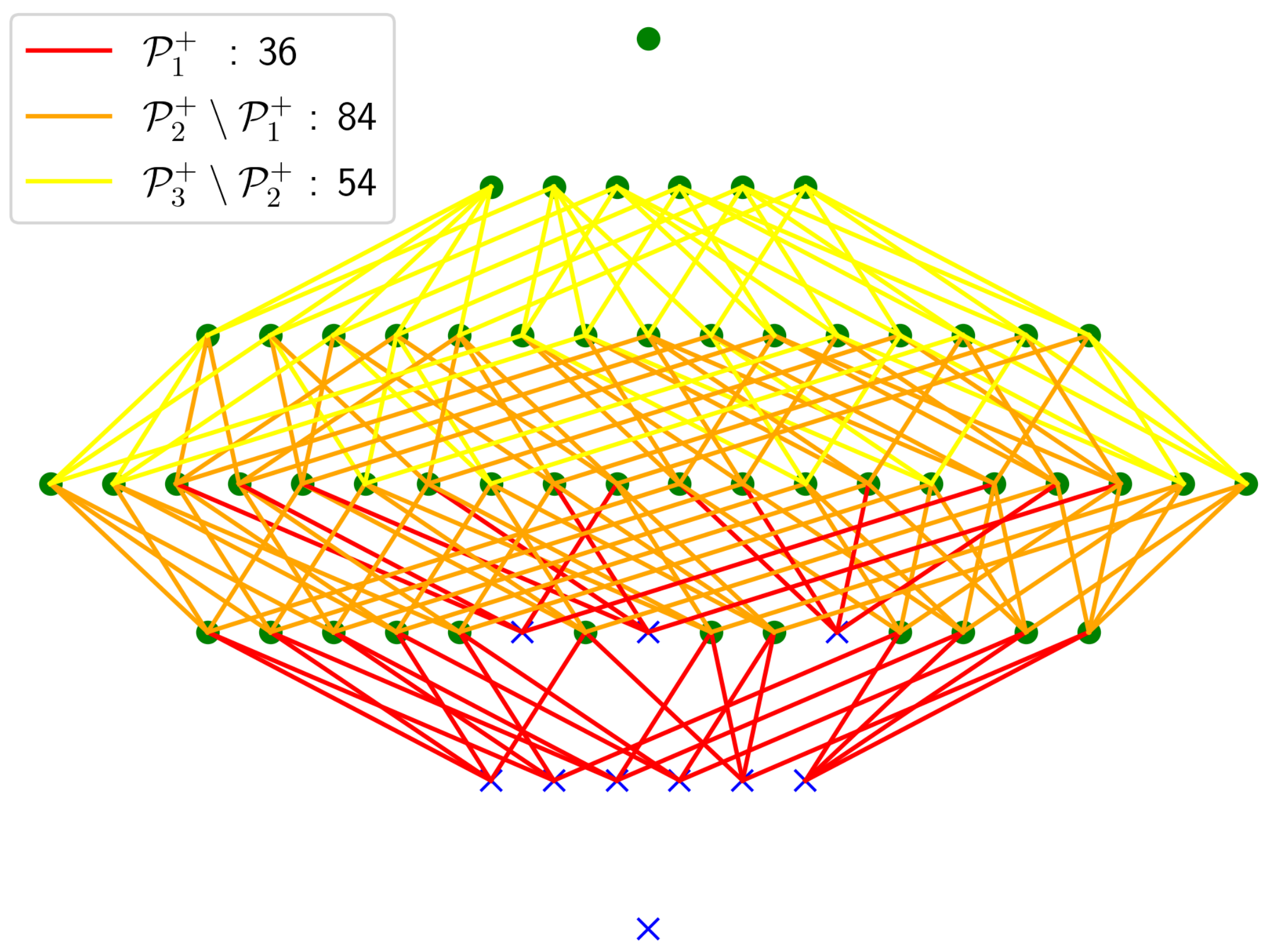

Figure 139: The maximizer of $E(|\mathcal{P}_3^+|)$ for $n = 6$

| $n$ | $\|A\|$ | $\phi_2$ | $\phi_3$ | Minimal$(A)$ |
|---|---|---|---|---|
| 3 | 7 | 9 | 12 | 100 010 001 |
| 4 | 11 | 24 | 28 | 1100 1010 1001 0110 0101 0011 |
| | 15 | 16 | 28 | 1000 0100 0010 0001 |
| 5 | 24 | 58 | 71 | 11000 10100 10010 10001 01100 01010 00101 00011 |
| | 25 | 54 | 71 | 11000 10100 10010 10001 01100 01010 01001 00110 00101 |
| 6 | 54 | 120 | 174 | 110000 101000 100100 100010 011000 010100 010001 001010 001001 000110 000101 000011 |

Table 17: The maximizers of $E(|\mathcal{P}_3^+|)$ for $n = 3, 4, 5, 6$

What conclusions can we draw from the numerical computations presented in this section and the previous one? An unpleasant fact is that there might exist several maximizers both for $\phi_{3\setminus 2}$ (for $n = 5, 6$) and $\phi_3$ (for $n = 5$). A good point is that the boundary of the maximizers seems to be close somehow to the middle layer of the hypercube. However the very structure of a maximizer, even when it is unique, appears to be quite complex. For instance the unique maximizer of $\phi_3$ for $n = 6$ shown in figure 139 has pivots of order three which are spread between four different layers.

## 58 The missing pivotal inequalities

There is one missing ingredient to complete the proof that $\theta(p_c, \mathbb{Z}^d) = 0$. We need a quantitative control on the expected number of the pivots of arbitrary fixed order. We formulate first a weak conjecture on the mean number of pivots of arbitrary order for a disconnection event. So far we did not succeed in defining properly a class of disconnection events for which the conjecture could be verified, which is why this concept is left intentionally vague.

**Conjecture 58.1** (Mean pivotal inequality)**.** *Let $n \geq 1$ and let $p \in ]0,1[$. There exists a constant $c = c(p)$ which depends on $p$ only such that, for a decreasing subset $\mathcal{D}$ of $\{0,1\}^n$ which is of the disconnection type, we have*

$$\forall k \in \{1,\dots,n\} \qquad E\big(\big|\mathcal{P}_k^+(\mathcal{D},\omega)\big|\big) \leq c\,k\,\sqrt{n}\,. \tag{58.1}$$

Conjecture 58.1, if true, would be a good sign indicating that our strategy to prove that $\theta(p_c, \mathbb{Z}^d) = 0$ might work. The main idea is the following. As always, we start with a parameter $p$ such that $\theta(p, \mathbb{Z}^d) > 0$. The presence of several disjoint large clusters in a box of side $n$ implies the absence of connections between some large subsets of the box. If these subsets have cardinality of order $n^d$, the intertwined exploration algorithms described in section 29 will reveal the presence of intersection sets of cardinality of order at least $n^{d-1}$, up to some logarithmic corrections. The exponent $d-1$ comes from the isoperimetric theorem, while the logarithmic corrections are due to the number of iterations needed for the intertwined explorations to explore fully the box. Thus, an event $\mathcal{D}$ which implies the presence of several disjoint large clusters is such that

$$\exists\,\kappa > 0 \quad \forall k \in \{1,\dots,\kappa\ln n\} \qquad E\big(\big|\mathcal{P}_k^+(\mathcal{D},\omega)\big|\,\big|\,\mathcal{D}\big) \geq \frac{1}{\kappa\ln n}\,n^{d-1}\,. \tag{58.2}$$

The event $\mathcal{D}$ being also decreasing, it would satisfy also the inequalities (58.1). Now inequalities (58.1) and (58.2) together indicate that the probability of $\mathcal{D}$ must be small. The success of our strategy depends in a crucial way on a quantitative estimate on how small $P_p(\mathcal{D})$ must be in order to ensure that inequality (58.2) holds. So we formulate next a stronger conjecture.

**Conjecture 58.2** (Deviations pivotal inequality)**.** *Let $n \geq 1$ and let $p \in ]0,1[$. There exists a constant $c = c(p)$ which depends on $p$ only such that, for a decreasing subset $\mathcal{D}$ of $\{0,1\}^n$ which is of the disconnection type, we have*

$$\forall k \in \{1,\dots,n\} \qquad P_p(\mathcal{D}) \leq \exp\Big(-\frac{c}{kn}\Big(E\big(|\mathcal{P}_k^+(\mathcal{D})|\,\big|\,\mathcal{D}\big)\Big)^2\Big)\,.$$

Conjecture 58.2, if true, would imply that $\theta(p_c, \mathbb{Z}^d) = 0$. Some weaker formulations of the conjecture would also be enough. For example, we could replace the fraction $c/k$ by a constant $c(k)$ and require a weaker form of dependence of $c(k)$ on $k$. The case $k = 1$ of these conjectures has been proved rigorously in proposition 39.5 and corollary 40.3. Despite the evidence provided by the extensive numerical computations presented in the previous sections, we were unable to find a suitable definition of the disconnection type for which we could prove the conjectures, even for $k = 2$.

# 59 Steiner symmetrization

How should we proceed to prove or disprove the conjecture 53.1? Ideally, we would like to determine exactly the set of the maximizers of $\phi_2$. Since this set looks much simpler for $\phi_2$ than for $\phi_{2\setminus 1}$, we focus here on $\phi_2$ and the extremal problem

$$M_2(n) \;=\; \max\, \big\{\, \phi_2(A) : A \text{ increasing subset of } \{0,1\}^n \,\big\}\,. \tag{59.1}$$

This extremal problem has the same flavour as a classical isoperimetric problem, were it not for the fact that we try to maximize the functional. The tricky part comes from the set of the candidates. Indeed, by requiring that the set $A$ is monotone, we drastically reduce the possibilities and we expect that the maximum value $M_2(n)$ is negligible compared to $n2^n$. The extremal problem (59.1) is very specific and does not seem to have been considered before. Our first attempts to solve this problem feel like coming up against a brick wall. Is there something in the vast literature on extremal problems that might help us? There is one generic strategy which invariably leads to the solution, whenever it can be implemented: to call in a smoothing procedure in the spirit of the Steiner symmetrization. To solve the classical isoperimetric problem in the plane, Jakob Steiner introduced an operation called the symmetrization. This operation transforms a set into another simpler set having the same area, without increasing its perimeter. For a beautiful exposition of this technique and its various generalizations, the reader can consult the book of Harper [68]. In our case, to attack the extremal problem (59.1), we are looking for a mechanism which transforms an increasing set $A$ into another increasing set $\widetilde{A}$, hopefully simpler than $A$, and such that

$$\phi_2(\widetilde{A}) \;\geq\; \phi_2(A)\,.$$

Of course, the problem is to find and design an adequate mechanism. For the case of the first pivots and the functional $\phi_1$, the second proof of proposition 50.1 relied indeed on such a mechanism, which consisted in removing from the set $A$ a minimal element belonging to a level below the middle level of the hypercube. Unfortunately, the examples of the previous section show that this mechanism is not any more adapted to the functional $\phi_2$, so we have to look for something else. Although our problem falls within the general category of extremal problems for finite graphs, it seems to be of a different nature than the classical isoperimetric problems. A first difference is that our problem is a maximization problem. A second difference is that the challenge is created by the choice of the admissible candidates, namely the increasing sets. Without this constraint, the maximization problem would not be interesting. A regretful consequence of these differences is that the classical tools used to solve isoperimetric problems on the hypercube, like the shift or compression, do not seem to be relevant. Even worse, they probably go in the wrong direction, because to increase the number of pivots of second order, it is probably better to somehow expand the minimal elements of a set $A$ rather than to compress them.

# 60 Kruskal-Katona

One of the cornerstones of extremal set theory is the Kruskal-Katona theorem (see for instance [56]), which we recall briefly next. This result relies on the so-called colex order. Let $E$ and $F$ be finite subsets of the set of the positive integers. We say that $E$ is smaller than $F$ for the colex order if

$$E \subset F \quad \text{ or } \quad \max(E \setminus F) < \max(F \setminus E) .$$

We identify an element of $\{0,1\}^n$ with a subset of $\{1,\dots,n\}$ through the map

$$\omega \in \{0,1\}^n \mapsto \operatorname{supp} \omega \subset \{1,\dots,n\} .$$

This map is one to one and onto the subsets of $\{1,\dots,n\}$, thus the colex order can be transported backwards on $\{0,1\}^n$. We consider the trace of the colex order on the configurations having exactly $k$ ones, and we denote by $A(k,m)$ the collection of the first $m$ configurations for the colex order of the set

$$\Big\{ \omega \in \{0,1\}^n : \sum_{1\le i\le n} \omega(i) = k \Big\} .$$

For instance, taking $n=6$, we have

$$\begin{aligned} A(1,7) &= \{ 100000, 010000, 001000, 000100, 000010, 000001 \} ,\\ A(2,7) &= \{ 110000, 101000, 011000, 100100, 010100, 001100, 100010 \} ,\\ A(3,7) &= \{ 111000, 110100, 101100, 011100, 111001, 110101, 110011 \} . \end{aligned}$$

Thus, to build $A(k,m)$, we consider configurations belonging to the level $k$, and we try to position the ones as far left as possible, and we take the first $m$ configurations generated this way.

**Theorem 60.1** (Kruskal-Katona)**.** *Let $k$ be an integer in $\{1,\dots,n\}$ and let $m \ge 1$. Let $A$ be a subset of $\{0,1\}^n$ consisting of $m$ configurations having exactly $k$ ones. The immediate shadow $\sigma(A)$ of $A$ is the subset*

$$\sigma(A) = \Big\{ \omega \in \{0,1\}^n : \sum_{1\le i\le n} \omega(i) = k-1, \, \exists\, \eta \in A \quad \omega \le \eta \Big\} .$$

*We have*

$$\big|\sigma(A)\big| \ge \big|\sigma(A(k,m))\big| .$$

This theorem is a kind of isoperimetric theorem on the hypercube. It has numerous applications and generalizations, and there exist several instructive proofs for it. It feels like there should exist a connection between our extremal problem and the Kruskal-Katona theorem. In their seminal paper introducing the notion of noise sensitivity [15], Benjamini, Kalai and Schramm state the inequality (50.2) (see theorem 3.3 there), and they say that it is a consequence of the Kruskal-Katona theorem. It is not clear what they had in mind, but it was certainly a proof different from the classical ones. In any case, in the example of figure 119, the minimal elements on the third layer of the hypercube constitute a typical Kruskal-Katona sequence.

# 61 The shift

One of the most efficient proof of the Kruskal-Katona theorem 60.1 relies on the shift operation, that we describe next. For $i \in \{1, \dots, n\}$, we denote by $e_i$ the configuration which has exactly one 1 at the $i$-th component and 0 elsewhere:

$$\forall j \in \{1, \dots, n\} \qquad e_i(j) = \begin{cases} 1 & \text{if} \quad i = j\,, \\ 0 & \text{if} \quad i \neq j\,. \end{cases}$$

Let $i, j$ be two distinct integers in $\{1, \dots, n\}$. For $\omega$ a configuration in $\{0,1\}^n$, we define the $(i,j)$-shifted configuration $s_{i,j}(\omega)$ as

$$s_{i,j}(\omega) = \begin{cases} \omega & \text{if} \quad (\omega(i), \omega(j)) \neq (0,1)\,, \\ \omega + e_i - e_j & \text{if} \quad (\omega(i), \omega(j)) = (0,1)\,. \end{cases} \tag{61.1}$$

We have merely rewritten the definition of the shift in our context. The shift is classically defined on subsets of $\{1, \dots, n\}$ as follows:

$$\forall F \subset \{1, \dots, n\} \qquad s_{i,j}(F) = \begin{cases} F \setminus \{j\} \cup \{i\} & \text{if} \quad i \notin F,\, j \in F\,, \\ F & \text{otherwise}\,. \end{cases} \tag{61.2}$$

Naturally, the two definitions (61.1) and (61.2) are completely equivalent, once we identify a configuration $\omega$ with its support $\operatorname{supp} \omega$. Since we struggle with a percolation problem, our choice is to state the definitions in terms of sets and configurations on $\{0,1\}^n$, rather in terms of collections of sets and subsets of $\{0,1\}^n$, as is customary in extremal set theory. The other choice, consisting in working with collections of sets, would require to rewrite the definitions related to pivots, and also to redefine the Bernoulli probability measure as a probability measure on subsets of $\{1, \dots, n\}$, and it might lead to greater misunderstandings and confusion. This being said, we will mainly stick to the notation and terminology of the beautiful book [56], with one exception. The $(i,j)$-shift $s_{i,j}(\omega)$ of a configuration $\omega$ is defined without making reference to a set of configurations, and we define separately the $(i,j)$-shift of a set, as follows. Let $A$ be a subset of $\{0,1\}^n$ and let $i, j$ be two distinct integers in $\{1, \dots, n\}$. The $(i,j)$-shift $S_{i,j}(A)$ of $A$ is defined as

$$S_{i,j}(A) = \left\{ s_{i,j}(\omega) : \omega \in A,\, s_{i,j}(\omega) \notin A \right\} \cup \left\{ \omega : \omega \in A,\, s_{i,j}(\omega) \in A \right\}.$$

This definition goes back to the classical paper of Erdös, Ko and Rado [50]. Informally, to build $S_{i,j}(A)$, we proceed as follows. We look at every configuration $\omega$ of $A$ such that $\omega(i) = 0$ and $\omega(j) = 1$. We try to replace such a configuration $\omega$ by the configuration $s_{i,j}(\omega)$ in which the components $i, j$ are exchanged, but we do the replacement only when $s_{i,j}(\omega)$ was not already present in $A$. The other configurations of $A$ are not modified. The shift operation has become a central tool in extremal set theory. We shall be interested here in its action on monotone sets.

We start with a simple observation: the shift of a monotone set is still monotone. Surprisingly, we could not find this statement in the existing literature, maybe due to a lack of knowledge. In fact, the shift operator is used in various contexts and the terminology varies from one to another. What we call here a monotone set might be called indifferently an ideal, or a downset, or an hereditary collection in extremal set theory. Moreover, there are not many works that deal with a problem similar to ours, i.e. trying to optimise a quantity over the monotone sets. A notable exception is the work of Duffus, Howard and Leader [43] on the width of downsets. In the paragraph before lemma 2 in [43], they write: "Note that the i-compression of a downset is again a downset." However their $i$-compression is different from the shift. So we include next the detailed argument for the shift.

**Proposition 61.1.** *Let $A$ be an increasing subset of $\{0,1\}^n$, and let $i, j$ be two distinct integers in $\{1, \dots, n\}$. The $(i,j)$-shift $S_{i,j}(A)$ of $A$ is still increasing.*

*Proof.* Throughout the proof, we drop the indices $i, j$ and we note simply $s, S$ the shifts $s_{i,j}$, $S_{i,j}$. Let $\sigma$ be a configuration in $S(A)$ and let $\eta$ be another configuration such that $\eta \geq \sigma$. We wish to prove that $\eta$ belongs also to $S(A)$. Notice that any configuration $\omega$ of $S(A)$ such that $(\omega(i), \omega(j)) \neq (1,0)$ belongs also to $A$. We consider several cases and subcases.

- $(\sigma(i), \sigma(j)) = (0,0)$. Then $\sigma$ is also in $A$. Since $A$ is increasing, then $\eta$ also belongs to $A$, as well as $s(\eta)$. Therefore $\eta$ is in $S(A)$.

- $(\sigma(i), \sigma(j)) = (1,1)$. Then $\sigma$ is also in $A$. Moreover $(\eta(i), \eta(j)) = (1,1)$ and $\eta$ is in both $A$ and $S(A)$.

- $(\sigma(i), \sigma(j)) = (1,0)$. We consider several subcases, depending on whether $\sigma$ is in $A$ or not.

  ⋆ $\sigma \in A$. In this case, $\eta$ is also in $A$. Since $s(\eta) = \eta$, then $\eta$ is in $S(A)$.

  ⋆ $\sigma \notin A$. In this case, $\sigma$ was created in $S(A)$ by the shift from $\widetilde{\sigma} = s_{j,i}(\sigma)$.

    ⋄ $(\eta(i), \eta(j)) = (1,1)$. We have $\eta \geq \widetilde{\sigma}$, thus $\eta$ is in $A$. Moreover $s(\eta) = \eta$, therefore $\eta$ is in $S(A)$.

    ⋄ $(\eta(i), \eta(j)) = (1,0)$. Let us set $\widetilde{\eta} = s_{j,i}(\eta)$. Since $\widetilde{\eta} \geq \widetilde{\sigma}$, then $\widetilde{\eta}$ is in $A$. Therefore $\eta = s(\widetilde{\eta})$ is in $S(A)$.

- $(\sigma(i), \sigma(j)) = (0,1)$. Since $\sigma$ is in $S(A)$, it has not been replaced by $s(\sigma)$ after the shift, and both configurations $\sigma$ and $s(\sigma)$ are present in $A$. Thus $\eta$ belongs also to $A$.

  ⋆ $(\eta(i), \eta(j)) \neq (0,1)$. We have $s(\eta) = \eta$, thus $\eta$ is also in $S(A)$.

  ⋆ $(\eta(i), \eta(j)) = (0,1)$. Since $s(\eta) \geq s(\sigma)$, then $s(\eta)$ is also in $A$, therefore both $\eta$ and $s(\eta)$ are in $S(A)$.

In all these cases, we see that $\eta$ is in $S(A)$, hence the set $S(A)$ is increasing. ☐

The good news is that the shift operator does not modify the expected size of the pivotal set.

**Proposition 61.2.** *Let $A$ be an increasing subset of $\{0,1\}^n$, and let $i,j$ be two distinct integers in $\{1,\dots,n\}$. We have*

$$E\big(|\mathcal{P}_1^+(A)|\big) \,=\, E\big(|\mathcal{P}_1^+(S_{i,j}(A))|\big)\,.$$

*Proof.* Throughout the proof, we drop the indices $i,j$ and we note simply $s,S$ the shifts $s_{i,j}$, $S_{i,j}$. We will compare $\mathcal{P}_1^+(A,\omega)$ and $\mathcal{P}_1^+(S(A),\omega)$ for $\omega$ a configuration in $A\cup S(A)$. Unfortunately, there are a lot of different cases to consider, depending on whether $\omega$ is in $A$ or $S(A)$, and the action of the shift on $\omega$. So, we fix a configuration $\omega$. We choose to organize the main cases according to the values of $\omega(i)$ and $\omega(j)$ (another possibility would be to consider first configurations belonging to $A$, and then configurations belonging to $S(A)$, but this leads to a long proof with several parts which are almost identical for both cases). Different subcases arise, because naturally the indices $i,j$ have a special role.

- $\omega(i)=\omega(j)=0$. This is the simplest case, because the shift has no effect on $\omega$ and on the configurations below $\omega$. We have $s(\omega)=\omega$, $s^{-1}(\omega)=\{\,\omega\,\}$ and

  $$\omega\in A \qquad\Longleftrightarrow\qquad \omega\in S(A)\,. \tag{61.3}$$

  So we suppose that $\omega$ is in $A\cap S(A)$. Let $k$ be an index in $\{1,\dots,n\}$. We have

  $$s^{-1}(\omega_k)=\{\,\omega_k\,\}\,,\qquad \omega_k\notin A\quad\Longleftrightarrow\quad \omega_k\notin S(A)\,, \tag{61.4}$$

  and

  $$k\in\mathcal{P}_1^+(A,\omega)\qquad\Longleftrightarrow\qquad k\in\mathcal{P}_1^+(S(A),\omega)\,. \tag{61.5}$$

  We conclude that

  $$\mathcal{P}_1^+(A,\omega)\,=\,\mathcal{P}_1^+(S(A),\omega)\,. \tag{61.6}$$

- $\omega(i)=\omega(j)=1$. This is the second simplest case, the shift has no effect on $\omega$ but it might affect some of the configurations below $\omega$. We have again $s(\omega)=\omega$, $s^{-1}(\omega)=\{\,\omega\,\}$ and the equivalence (61.3) holds. Notice also that $s^{-1}(\omega_i)=\{\,\omega_i\,\}$ and $\omega_j=s(\omega_i)$. So we suppose that $\omega$ is in $A\cap S(A)$. Let $k$ be an index in $\{1,\dots,n\}$. We consider several subcases:

  ⋆ $k\neq i,j$. In this case, the properties (61.4) and (61.5) hold. Thus

  $$\mathcal{P}_1^+(A,\omega)\setminus\{\,i,j\,\}\,=\,\mathcal{P}_1^+(S(A),\omega)\setminus\{\,i,j\,\}\,.$$

  ⋆ $k\in\{\,i,j\,\}$, $\omega_i\notin A$, $\omega_j\notin A$. In this case, $\omega_i\notin S(A)$, $\omega_j\notin S(A)$. Thus

  $$\{\,i,j\,\}\,\subset\,\mathcal{P}_1^+(A,\omega)\cap\mathcal{P}_1^+(S(A),\omega)\,.$$

  ⋆ $k\in\{\,i,j\,\}$, $\omega_i\in A$, $\omega_j\notin A$. In this case, $\omega_i\notin S(A)$, $\omega_j\in S(A)$. Thus

  $$i\in\mathcal{P}_1^+(S(A),\omega)\setminus\mathcal{P}_1^+(A,\omega)\,,\quad j\in\mathcal{P}_1^+(A,\omega)\setminus\mathcal{P}_1^+(S(A),\omega)\,.$$

From these three subcases, we conclude that

$$\mathcal{P}_1^+(S(A),\omega) \,=\, S\big(\mathcal{P}_1^+(A,\omega)\big)\,. \tag{61.7}$$

- $\omega(i) = 1$, $\omega(j) = 0$ or $\omega(i) = 0$, $\omega(j) = 1$. This case is more complicated than the previous ones, as the shift might move several relevant configurations from $A$ to $S(A)$, along with the pivots associated to them. The most efficient presentation consists in considering both configurations simultaneously. Suppose for instance that $\omega(i) = 1$, $\omega(j) = 0$ and let us set $\eta = \omega_i^j$. We have
$$s(\omega) = \omega\,,\quad s(\eta) = \omega\,,\quad s^{-1}(\omega) = \{\,\eta\,\}\,.$$
Let $k$ be an index in $\{\,1,\dots,n\,\}$. We consider several subcases, arranged in increasing complexity:

  - $\star$ $\omega \notin A\,,\eta \notin A$. In this case, $\omega \notin S(A)\,,\eta \notin S(A)$, whence
  $$\mathcal{P}_1^+(A,\omega) \,=\, \mathcal{P}_1^+(S(A),\omega) \,=\, \mathcal{P}_1^+(A,\eta) \,=\, \mathcal{P}_1^+(S(A),\eta) \,=\, \varnothing\,.$$
  - $\star$ $\omega \in A\,,\eta \notin A$. In this case, $\omega \in S(A)$ and $\eta \notin S(A)$. We notice first that
  $$\mathcal{P}_1^+(A,\eta) \,=\, \mathcal{P}_1^+(S(A),\eta) \,=\, \varnothing\,.$$
  Furthermore, since $\eta \notin A$, then $A$ contains no configuration smaller than $\eta$, and the restrictions of $A$ and $S(A)$ to the configurations smaller than $\omega$ coincide. Therefore, the properties (61.4) and (61.5) hold for $\omega$, thus
  $$\mathcal{P}_1^+(A,\omega) \,=\, \mathcal{P}_1^+(S(A),\omega)\,.$$
  - $\star$ $\omega \notin A\,,\eta \in A$. In this case, $\omega \in S(A)$ and $\eta \notin S(A)$.
    - $\diamond$ $k \notin \{\,i,j\,\}$, $\eta_k \in A$. We have $\omega_k \in S(A)$, $\eta_k \notin S(A)$, thus
    $$k \notin \mathcal{P}_1^+(A,\omega) \cup \mathcal{P}_1^+(S(A),\omega)\,,\quad k \notin \mathcal{P}_1^+(A,\eta) \cup \mathcal{P}_1^+(S(A),\eta)\,.$$
    - $\diamond$ $k \notin \{\,i,j\,\}$, $\eta_k \notin A$. We have $\omega_k \notin S(A)$, $\eta_k \notin S(A)$, thus
    $$k \in \mathcal{P}_1^+(S(A),\omega) \cap \mathcal{P}_1^+(A,\eta)\,,\quad k \notin \mathcal{P}_1^+(A,\omega) \cup \mathcal{P}_1^+(S(A),\eta)\,.$$
    - $\diamond$ $k \in \{\,i,j\,\}$. We have $\eta_j \notin A$ (because $\eta_j \leq \omega$), whence $\eta_j \notin S(A)$. Notice also that $\omega_i = \eta_j$, thus
    $$i \notin \mathcal{P}_1^+(A,\omega) \cap \mathcal{P}_1^+(S(A),\eta)\,,\quad j \in \mathcal{P}_1^+(A,\eta) \cap \mathcal{P}_1^+(S(A),\omega)\,.$$
  - $\star$ $\omega \in A\,,\eta \in A$. In this case, $\omega \in S(A)$ and $\eta \in S(A)$.
    - $\diamond$ $k \notin \{\,i,j\,\}$, $\omega_k \in A$, $\eta_k \in A$. We have $\omega_k \in S(A)$, $\eta_k \in S(A)$, thus
    $$k \notin \mathcal{P}_1^+(A,\omega) \cup \mathcal{P}_1^+(S(A),\omega)\,,\quad k \notin \mathcal{P}_1^+(A,\eta) \cup \mathcal{P}_1^+(S(A),\eta)\,.$$

$\diamond$ $k \notin \{\, i,j \,\}$, $\omega_k \notin A$, $\eta_k \notin A$. We have $\omega_k \notin S(A)$, $\eta_k \notin S(A)$, thus

$$k \in \mathcal{P}_1^+(A,\omega) \cap \mathcal{P}_1^+(S(A),\omega)\,, \quad k \in \mathcal{P}_1^+(A,\eta) \cap \mathcal{P}_1^+(S(A),\eta)\,.$$

$\diamond$ $k \notin \{\, i,j \,\}$, $\omega_k \in A$, $\eta_k \notin A$. We have $\omega_k \in S(A)$, $\eta_k \notin S(A)$, thus

$$k \notin \mathcal{P}_1^+(A,\omega) \cup \mathcal{P}_1^+(S(A),\omega)\,, \quad k \in \mathcal{P}_1^+(A,\eta) \cap \mathcal{P}_1^+(S(A),\eta)\,.$$

$\diamond$ $k \notin \{\, i,j \,\}$, $\omega_k \notin A$, $\eta_k \in A$. We have $\omega_k \in S(A)$, $\eta_k \notin S(A)$, thus

$$k \in \mathcal{P}_1^+(A,\omega) \cap \mathcal{P}_1^+(S(A),\eta)\,, \quad k \notin \mathcal{P}_1^+(S(A),\omega) \cup \mathcal{P}_1^+(A,\eta)\,.$$

$\diamond$ $k \in \{\, i,j \,\}$, $\omega_i \in A$, $\eta_j \in A$. We have $\omega_i \in S(A)$, $\eta_j \in S(A)$. Thus

$$\{\, i,j \,\} \cap \Big(\mathcal{P}_1^+(A,\omega) \cup \mathcal{P}_1^+(S(A),\omega) \cup \mathcal{P}_1^+(A,\eta) \cup \mathcal{P}_1^+(S(A),\eta)\Big) \,=\, \varnothing\,.$$

$\diamond$ $k \in \{\, i,j \,\}$, $\omega_i \notin A$, $\eta_j \notin A$. We have $\omega_i \notin S(A)$, $\eta_j \notin S(A)$. Thus

$$i \in \mathcal{P}_1^+(A,\omega) \cap \mathcal{P}_1^+(S(A),\omega)\,, \quad j \in \mathcal{P}_1^+(A,\eta) \cap \mathcal{P}_1^+(S(A),\eta)\,.$$

Notice that in fact $\omega_i = \eta_j$, so we do not have to consider the cases $\omega_i \in A$, $\eta_j \notin A$ or $\omega_i \notin A$, $\eta_j \in A$. In all the cases considered above, we check that

$$\big|\mathcal{P}_1^+(A,\omega)\big| + \big|\mathcal{P}_1^+(A,\eta)\big| \,=\, \big|\mathcal{P}_1^+(S(A),\omega)\big| + \big|\mathcal{P}_1^+(S(A),\eta)\big|\,. \tag{61.8}$$

We shall now use the previous analysis to compare the total numbers of pivots of $A$ and $S(A)$. We decompose the sum

$$\sum_{\omega\in A} \big|\mathcal{P}_1^+(A,\omega)\big| \,=\, \sum_{\omega:\omega(i)=\omega(j)=0} \cdots \;+ \sum_{\omega:\omega(i)=\omega(j)=1} \cdots \;+ \sum_{\omega:\omega(i)=1,\omega(j)=0} \cdots \;+ \sum_{\omega:\omega(i)=0,\omega(j)=1} \cdots\,. \tag{61.9}$$

From (61.6), we have

$$\sum_{\omega\in A:\omega(i)=\omega(j)=0} \big|\mathcal{P}_1^+(A,\omega)\big| \,=\, \sum_{\omega\in S(A):\omega(i)=\omega(j)=0} \big|\mathcal{P}_1^+(S(A),\omega)\big|\,. \tag{61.10}$$

From (61.7), we have

$$\sum_{\omega\in A:\omega(i)=\omega(j)=1} \big|\mathcal{P}_1^+(S(A),\omega)\big| \,=\, \sum_{\omega\in A:\omega(i)=\omega(j)=1} \big|S\big(\mathcal{P}_1^+(A,\omega)\big)\big| \,=\, \sum_{\omega\in S(A):\omega(i)=\omega(j)=1} \big|\mathcal{P}_1^+(A,\omega)\big|\,. \tag{61.11}$$

Finally, we combine the last two sums in (61.9) by summing over all the pairs $\omega, \eta$ such that $\omega(i) = 1$, $\omega(j) = 0$ and $\eta = \omega_i^j$. Decomposing further the sum into the different subcases studied previously, we obtain with the help of (61.8) that

$$\sum_{\omega,\eta} \left|\mathcal{P}_1^+(A,\omega)\right| + \left|\mathcal{P}_1^+(A,\eta)\right| = \sum_{\omega,\eta} \left|\mathcal{P}_1^+(S(A),\omega)\right| + \left|\mathcal{P}_1^+(S(A),\eta)\right|. \quad (61.12)$$

Summing the equalities (61.10), (61.11) and (61.12), we obtain

$$\sum_{\omega\in A} \left|\mathcal{P}_1^+(A,\omega)\right| = \sum_{\omega\in S(A)} \left|\mathcal{P}_1^+(S(A),\omega)\right|,$$

which is the desired result. □

With the help of proposition 61.2, we could derive an alternative proof of the results of proposition 50.1, which would run as follows. We say that a set $A$ is fully compressed if it is invariant under the application of all the shifts

$$S_{i,j}\,, \quad i,j \in \{\,1,\dots,n\,\}\,, \quad i<j\,.$$

It follows from proposition 61.2 that the quantity to be maximized is invariant by application of a shift, so we can reduce the set of the candidates to the fully compressed configurations. This is a drastic simplification, but there is still some work to be done. Indeed, the trace on the level $k$ of a fully compressed configuration can be different from the sets $A(k,m)$ appearing in the Kruskal-Katona theorem 60.1. After this first step, it remains to analyze the pivots of the fully compressed monotone sets. Maybe this was the idea of Benjamini, Kalai and Schramm, when they wrote that the inequality (50.2) is a consequence of the Kruskal-Katona theorem (see the remark after theorem 3.3 in [15]). By the way, this is also the most efficient scheme of proof of the Kruskal-Katona theorem itself: indeed, the shadow of a set never increases by applying a shift. Anyway, this is certainly not an efficient strategy to prove proposition 50.1. However, our goal is rather to find a useful generalization of this inequality to the higher pivots, and this is why we tried to examine all the different proofs of the inequality (50.2). An advantage of the above strategy is that it works also for the more difficult extremal problem of proposition 50.3, because the cardinality of a set is also invariant under the shifts.

Let us examine next the effect of the shift on the second pivots. In the simplest cases, we observe that the number of second pivots does not decrease after the application of the shift. So we could hope to reduce drastically the difficulty of the problem by considering only shifted configurations. Now, the bad news is that the shift might decrease the number of second pivots. Figure 140 shows an example of this behavior. We start with the monotone set whose minimal implicants are 1100, 1010, 0101. We apply the shift $S_{2,3}$ and we obtain the set whose minimal implicants are 1100, 1010, 0110. After applying the shift, the number of second order pivots has decreased by 2.

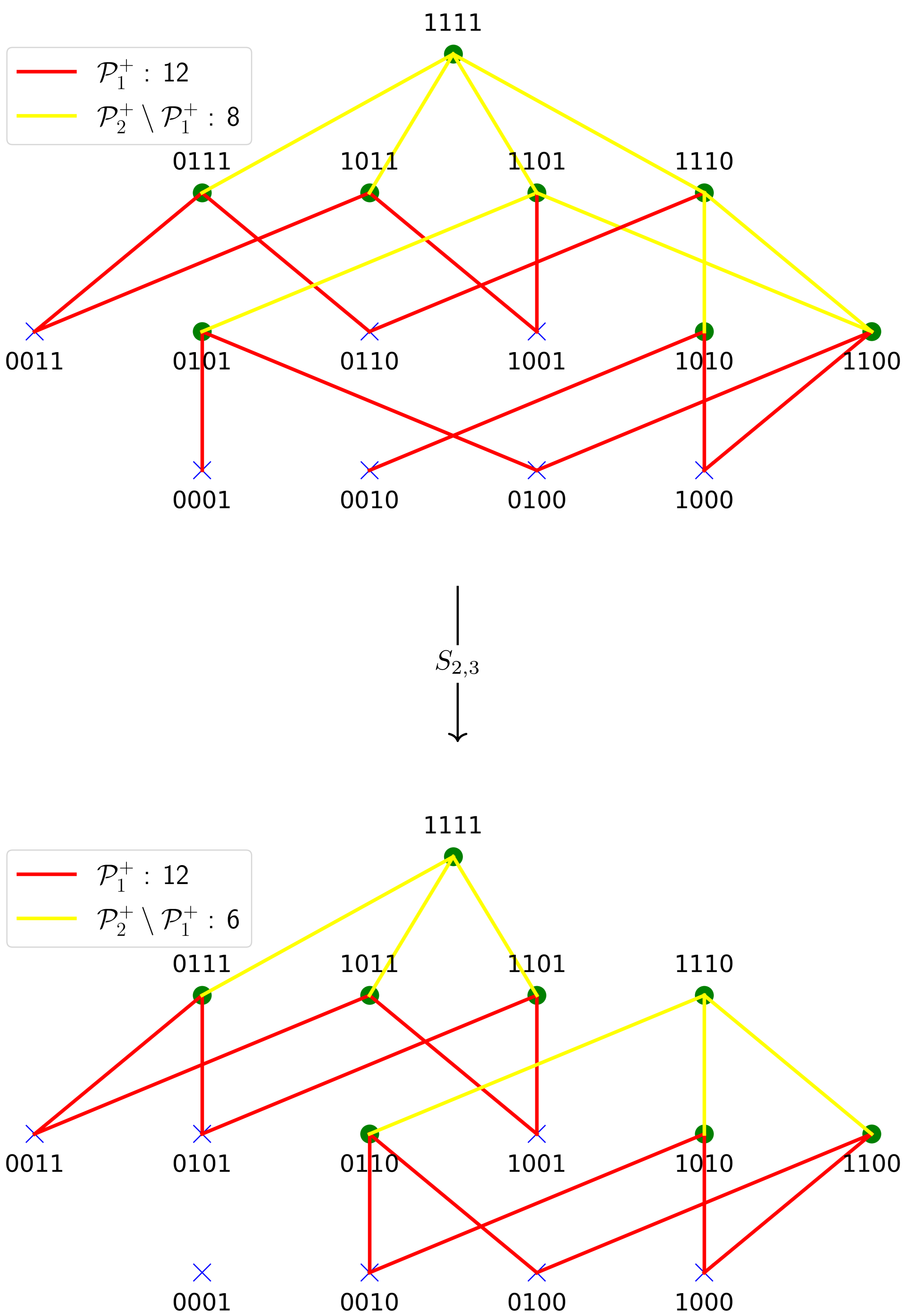

Figure 140: The shift $S_{2,3}$ causes a decrease of the second pivots.

# 62 The set of Talagrand

At the end of his paper [125], Talagrand builds an example of an increasing subset $A^\circ$ of $\{0,1\}^n$ such that, for some universal constant $K$, we have

$$P\Big(\big|\mathcal{P}_1^+(A^\circ,\omega)\big| \geq \frac{1}{K}\sqrt{n}\Big) \;\geq\; \frac{1}{K}\,. \tag{62.1}$$

For this set $A^\circ$, the expected value of $\big|\mathcal{P}_1^+(A^\circ,\omega)\big|$ is larger than $\sqrt{n}/K^2$, which is of the same order as the maximal value $M_1(n)$ given in (50.2). However, the edge boundary of Talagrand's set $A^\circ$ has a very different structure from the sets realizing the maximum $M_1(n)$. These sets were little perturbations of the set $A^*$ defined by

$$A^* \;=\; \Big\{\,\omega : \sum_{1\leq i\leq n} \omega(i) > \frac{n}{2}\,\Big\}\,.$$

What happens for the set $A^*$ is that the probability of the inner vertex boundary of $A^*$, given by

$$\partial^{\,in} A^* \;=\; \Big\{\,\omega : \sum_{1\leq i\leq n} \omega(i) = \max\Big(\Big\lceil\frac{n}{2}\Big\rceil, \Big\lceil\frac{n+1}{2}\Big\rceil\Big)\,\Big\}\,,$$

is of order $1/\sqrt{n}$, but each configuration $\omega$ in $\partial^{\,in} A^*$ has roughly $n/2$ components equal to 1, and all of them are pivotal, so that instead of (62.1) we have

$$P\Big(\big|\mathcal{P}_1^+(A^*,\omega)\big| \geq \frac{n}{2}\Big) \;\geq\; \frac{1}{K\sqrt{n}}\,.$$

Although Talagrand's set $A^\circ$ reaches the order $\sqrt{n}$ for the mean number of pivots, it is ultimately beaten by the set $A^*$ which realizes the value $M_1(n)$. For a long time, we hoped that the situation would be similar for the second order pivots, in the sense that the maximum $M_2(n)$ was reached by sets which are little perturbations of $A^*$, or at least not too different from the upper part of the hypercube. If this had been the case, then we would have been able to obtain the desired control on the pivots of order 2. Alas, Talagrand's example completely dashed this hope. We give next the details of Talagrand's construction, done at the very end of his paper [125]. We use essentially the same notation as Talagrand, and we introduce an additional event $T_j$ in order to compute an adequate lower bound on the expected number of second pivots. Talagrand's style being extremely efficient, while we tend to provide all the boring details, a few lines of Talagrand's computations become here almost three pages...

Let $\alpha$ be a fixed positive constant. Let $n\geq 1$, and let us set $p=\lfloor\alpha\sqrt{n}\rfloor$. Let

$$(Y_{i,j}, 1\leq i\leq p, 1\leq j\leq 2^p),$$

be an array of i.i.d. random variables, uniformly distributed over $\{\,1,\dots,n\,\}$. The set $A^\circ$ of Talagrand is defined as

$$A^\circ \;=\; \Big\{\,\omega\in\{0,1\}^n : \forall j\in\{\,1,\dots,2^p\,\} \quad \exists\, i\in\{\,1,\dots,p\,\} \quad \omega(Y_{i,j})=1\,\Big\}\,.$$

Let $k$ belong to $\{\,0,\dots,n\,\}$. We denote by $\mathcal{L}(k)$ the set of the configurations of $\{0,1\}^n$ having exactly $k$ components equal to 1, i.e.,

$$\mathcal{L}(k) \,=\, \Big\{\,\omega\in\{0,1\}^n : |\,\mathrm{supp}\,\omega| = k\,\Big\}\,.$$

Using the independence of the variables $Y_{i,j}$, we can compute the conditional probability

$$P\big(\omega\in A^\circ\,\big|\,\omega\in\mathcal{L}(n-k)\big) \,=\, \Big(1-\Big(\frac{k}{n}\Big)^p\Big)^{2^p}\,. \tag{62.2}$$

Now, for $k\leq n/2+p$, we have

$$\Big(1-\Big(\frac{k}{n}\Big)^p\Big)^{2^p} \,\geq\, \Big(1-\Big(\frac{1}{2}+\frac{p}{n}\Big)^p\Big)^{2^p} \,\geq\, \Big(1-\frac{1}{2^p}\Big(1+\frac{2p}{n}\Big)^p\Big)^{2^p}\,, \tag{62.3}$$

and furthermore

$$\Big(1+\frac{2p}{n}\Big)^p \,\leq\, \exp\big(\frac{2p^2}{n}\big) \,\leq\, e^{2\alpha^2}\,. \tag{62.4}$$

Putting together the three inequalities (62.2), (62.3) and (62.4), we obtain

$$\forall k\leq n/2+p \qquad P\big(\omega\in A^\circ\,\big|\,\omega\in\mathcal{L}(n-k)\big) \,\geq\, \Big(1-\frac{e^{2\alpha^2}}{2^p}\Big)^{2^p}\,. \tag{62.5}$$

The right-hand side of (62.5) does not depend on $k$ and it converges towards $\exp(-e^{2\alpha^2})$ as $n$ goes to $\infty$. Therefore, there exists a constant $c(\alpha)>0$, which depends on $\alpha$ only, such that

$$\forall n\geq 1\quad \forall k\leq\frac{n}{2}+p \qquad P\big(\omega\in A^\circ\,\big|\,\omega\in\mathcal{L}(n-k)\big) \,\geq\, c(\alpha)\,. \tag{62.6}$$

Following Talagrand, for $j\in\{\,1,\dots,2^p\,\}$, we introduce the events $U_j,V_j$ and an additional event $T_j$, defined by

$$\begin{aligned}
U_j \,&=\, \Big\{\,\omega\in\{0,1\}^n : \big|\{\,i\in\{\,1,\dots,p\,\}:\omega(Y_{i,j})=1\,\}\big| \,=\, 1\,\Big\}\,,\\
V_j \,&=\, \Big\{\,\omega\in\{0,1\}^n : \big|\{\,i\in\{\,1,\dots,p\,\}:\omega(Y_{i,j})=1\,\}\big| \,\geq\, 1\,\Big\}\,,\\
T_j \,&=\, \Big\{\,\omega\in\{0,1\}^n : \big|\{\,i\in\{\,1,\dots,p\,\}:\omega(Y_{i,j})=1\,\}\big| \,=\, 2\,\Big\}\,.
\end{aligned}$$

We have then

$$A^\circ \,=\, \bigcap_{1\leq j\leq 2^p} V_j\,,$$

and it follows from the independence of the variables $Y_{i,j}$ that

$$\begin{aligned}
P\big(U_j\,\big|\,A^\circ\cap\mathcal{L}(n-k)\big) \,&=\, P\big(U_j\,\big|\,V_j\cap\mathcal{L}(n-k)\big) \,\geq\, p\Big(1-\frac{k}{n}\Big)\Big(\frac{k}{n}\Big)^{p-1}\,,\\
P\big(T_j\,\big|\,A^\circ\cap\mathcal{L}(n-k)\big) \,&=\, P\big(T_j\,\big|\,V_j\cap\mathcal{L}(n-k)\big) \,\geq\, \frac{p(p-1)}{2}\Big(1-\frac{k}{n}\Big)^2\Big(\frac{k}{n}\Big)^{p-2}\,.
\end{aligned}$$

From now onwards, we consider only values of $k$ such that $|k - n/2| \leq p$. In this regime, the previous inequalities yield that

$$P\big(U_j \,\big|\, A^\circ \cap \mathcal{L}(n-k)\big) \;\geq\; p\Big(\frac{1}{2} - \frac{p}{n}\Big)^p\,, \tag{62.7}$$

$$P\big(T_j \,\big|\, A^\circ \cap \mathcal{L}(n-k)\big) \;\geq\; \frac{p(p-1)}{2}\Big(\frac{1}{2} - \frac{p}{n}\Big)^p\,. \tag{62.8}$$

Furthermore, since $p \leq \alpha\sqrt{n}$, we have

$$\Big(\frac{1}{2} - \frac{p}{n}\Big)^p \;\geq\; \frac{1}{2^p}\Big(1 - \frac{2p}{n}\Big)^p \;\geq\; \frac{1}{2^p}\Big(1 - \frac{2\alpha}{\sqrt{n}}\Big)^{\alpha\sqrt{n}} \;\geq\; \frac{1}{2^p}e^{-3\alpha^2}\,, \tag{62.9}$$

where the last inequality holds for $n$ large enough. Inequalities (62.7), (62.8) and (62.9) yield that

$$P\big(U_j \,\big|\, A^\circ \cap \mathcal{L}(n-k)\big) \;\geq\; \frac{p}{2^p}e^{-3\alpha^2}\,,$$

$$P\big(T_j \,\big|\, A^\circ \cap \mathcal{L}(n-k)\big) \;\geq\; \frac{p(p-1)}{2^p}e^{-3\alpha^2}\,.$$

Conditionally on $A^\circ \cap \mathcal{L}(n-k)$, the events $U_j, 1 \leq j \leq 2^p$, are independent, as well as the events $V_j, 1 \leq j \leq 2^p$. By the law of large numbers, for $p$ sufficiently large, we have

$$P\Big(\sum_{1 \leq j \leq 2^p} U_j \geq \frac{p}{2}e^{-3\alpha^2} \,\big|\, A^\circ \cap \mathcal{L}(n-k)\Big) \;\geq\; \frac{1}{2}\,,$$

$$P\Big(\sum_{1 \leq j \leq 2^p} T_j \geq \frac{p^2}{2}e^{-3\alpha^2} \,\big|\, A^\circ \cap \mathcal{L}(n-k)\Big) \;\geq\; \frac{1}{2}\,. \tag{62.10}$$

We need an ultimate conditioning before concluding. We define the random sets

$$J_1 \;=\; \Big\{\, j \in \{\,1,\dots,2^p\,\} : \omega \in U_j \,\Big\}\,,$$

$$J_2 \;=\; \Big\{\, j \in \{\,1,\dots,2^p\,\} : \omega \in T_j \,\Big\}\,.$$

For $j \in J_1$, there exists a unique $I_1(j) \in \{\,1,\dots,n\,\}$ such that $\omega(Y_{I_1(j),j}) = 1$. For $j \in J_2$, there exists a unique pair of indices $I_{21}(j), I_{22}(j) \in \{\,1,\dots,n\,\}$ such that $I_{21}(j) < I_{22}(j)$ and $\omega(Y_{I_{21}(j),j}) = \omega(Y_{I_{22}(j),j}) = 1$. It follows from the definition of $A^\circ$ that, whenever $\omega$ is in $A^\circ$, we have

$$\Big\{\, I_1(j) : 1 \leq j \leq 2^p \,\Big\} \;\subset\; \mathcal{P}_1^+(A^\circ,\omega)\,,$$

$$\Big\{\, I_{21}(j) : 1 \leq j \leq 2^p \,\Big\} \cup \Big\{\, I_{22}(j) : 1 \leq j \leq 2^p \,\Big\} \;\subset\; \mathcal{P}_2^+(A^\circ,\omega)\,. \tag{62.11}$$

The final trick is to remark that, if we condition on the set $J_1$, the indices $I_1(j), j \in J_1$ are independent and uniformly distributed over the support of $\omega$, which has cardinality $n-k$ and is of order $n/2$. Similarly, if we condition on the set $J_2$, the pair of indices $(I_{21}(j), I_{22}(j)), j \in J_2$ are independent and uniformly distributed over the pairs of distinct indices belonging to the support of $\omega$, which has cardinality $n-k$ and is of order $n/2$. For $j \in J_2$, let $I_3(j)$ be a random index obtained by choosing either $I_{21}(j)$ or $I_{22}(j)$, each with probability $1/2$. The random variables $I_3(j), j \in J_2$, are independent and uniformly distributed over the support of $\omega$, and thanks to the inclusion (62.11), they are all second order pivots. In order to obtain a lower bound on the number of first pivots (respectively the number of second pivots), we simply need to control the number of distinct values taken by the variables $I_1(j), j \in J_1$ (respectively $I_3(j), j \in J_2$).

This is a standard question, related to the classical coupon collector problem. Let us reformulate it with a simpler notation before performing the adequate computations. Let $N, m \in \{1, \dots, n\}$ and let $Y_1, \cdots, Y_m$ be $m$ independent random variables with uniform distribution over $\{1, \dots, N\}$. We would like to control

$$R_m = \big|\{ Y_k : 1 \leq k \leq m \}\big| .$$

For $i \in \{1, \dots, N\}$, let $A_i$ be the event that the value $i$ appears at least once in the sample $Y_1, \dots, Y_m$. The probability $P(A_i)$ does not depend on $i$ and is equal to

$$P(A_i) = 1 - \Big(1 - \frac{1}{N}\Big)^m .$$

With this notation, $R_m$ can be rewritten as

$$R_m = \sum_{1 \leq i \leq N} 1_{A_i} .$$

By linearity, the expectation of $R_m$ is equal to

$$E(R_m) = \sum_{1 \leq i \leq N} P(A_i) = N\Big(1 - \Big(1 - \frac{1}{N}\Big)^m\Big) . \tag{62.12}$$

To compute the second moment of $R_m$, we write

$$E\big((R_m)^2\big) = \sum_{1 \leq i,j \leq N} P(A_i \cap A_j) , \tag{62.13}$$

and we compute, for $i \neq j$,

$$\begin{aligned} P(A_i \cap A_j) &= P(A_i) + P(A_j) - P(A_i \cup A_j) \\ &= 1 - 2\Big(1 - \frac{1}{N}\Big)^m + \Big(1 - \frac{2}{N}\Big)^m , \end{aligned} \tag{62.14}$$

where the last term is the probability that neither $i$ nor $j$ appear in the sample $Y_1, \dots, Y_m$, which is the complement of the event $A_i \cup A_j$. Putting together formulas (62.12), (62.13) and (62.14), we obtain the following expression for the

variance of $R_m$:

$$\operatorname{Var}(R_m) \;=\; N(N-1)\Big(1-2\Big(1-\frac{1}{N}\Big)^m+\Big(1-\frac{2}{N}\Big)^m\Big)+E(R_m)$$
$$-\,N^2\Big(1-\Big(1-\frac{1}{N}\Big)^m\Big)^2\,. \quad (62.15)$$

Let us come back to our problem of obtaining a lower bound on the number of pivots of order one and two. Talagrand proved that there exists a positive constant $c>0$ such that

$$\forall n\geq 1 \qquad P\Big(\big|\mathcal{P}_1^+(A^\circ,\omega)\big|\geq c\sqrt{n}\Big)\;\geq\; c\,.$$

We consider next the second pivots and we try to prove a similar inequality involving $\big|\mathcal{P}_2^+(A^\circ,\omega)\big|$. To obtain the desired estimates, we will have to perform several conditioning. Let $L>0$ be a constant, and let us write

$$P\Big(\big|\mathcal{P}_2^+(A^\circ,\omega)\big|\geq Lp^2\Big)\;\geq$$
$$\sum_{\substack{k:|k-n/2|\leq p\\ j_2\geq p^2e^{-3\alpha^2}/2}} P\big(\big|\mathcal{P}_2^+(A^\circ,\omega)\big|\geq Lp^2,\,|J_2|=j_2,\,\omega\in A^\circ\cap\mathcal{L}(n-k)\big)$$
$$\geq \sum_{\substack{k:|k-n/2|\leq p\\ j_2\geq p^2e^{-3\alpha^2}/2}} P\big(\big|\mathcal{P}_2^+(A^\circ,\omega)\big|\geq Lp^2\,\big|\,|J_2|=j_2,\,\omega\in A^\circ\cap\mathcal{L}(n-k)\big)$$
$$\times P\big(|J_2|=j_2,\,\omega\in A^\circ\cap\mathcal{L}(n-k)\big)\,. \quad (62.16)$$

To estimate the conditional probability appearing in (62.16), we use the computations done on the coupon collector problem with $m=j_2$ and $N=n-k$. We have

$$P\Big(\big|\mathcal{P}_2^+(A^\circ,\omega)\big|\geq Lp^2\,\big|\,|J_2|=j_2,\,\omega\in A^\circ\cap\mathcal{L}(n-k)\Big)$$
$$\geq\; P\Big(\big|\{\,I_3(j):j\in J_2\,\}\big|\;\geq\; Lp^2\,\big|\,|J_2|=j_2,\,\omega\in A^\circ\cap\mathcal{L}(n-k)\Big)$$
$$\geq\; P\big(R_m\geq Lp^2\big)\,. \quad (62.17)$$

We will choose $L$ small enough to ensure that $Lp^2<E(R_m)$, and we use the Chebyshev inequality to control

$$P\big(R_m<Lp^2\big)\;\leq\;P\big(\big|R_m-E(R_m)\big|>E(R_m)-Lp^2\big)\;\leq\;\frac{\operatorname{Var}(R_m)}{\big(E(R_m)-Lp^2\big)^2}\,. \quad (62.18)$$

With $m=j_2\geq(p^2/2)e^{-3\alpha^2}$ and $N\geq n/2-\sqrt{n}$, we have from (62.12) that

$$E(R_m)\;\geq\;\Big(\frac{n}{2}-\sqrt{n}\Big)\Big(1-\Big(1-\frac{1}{n}\Big)^{(\alpha^2n/2)e^{-3\alpha^2}}\Big)\;\geq\;\frac{n}{3}\Big(1-\exp\big(-(\alpha^2/2)e^{-3\alpha^2}\big)\Big)\,,$$

where the last inequality holds for $n$ large enough. Thus we can choose $L$ small enough (depending on $\alpha$) to ensure that

$$E(R_m) - Lp^2 \geq \frac{1}{2} E(R_m)\,. \tag{62.19}$$

Inequalities (62.18) and (62.19) yield that

$$P\big(R_m < Lp^2\big) \leq \frac{4\,\mathrm{Var}(R_m)}{\big(E(R_m)\big)^2}\,. \tag{62.20}$$

It follows from (62.15) that $\mathrm{Var}(R_m) \leq E(R_m)$, thus we deduce from (62.20) that

$$P\big(R_m < Lp^2\big) \leq \frac{4}{E(R_m)} \leq \frac{1}{2}\,, \tag{62.21}$$

where the last inequality holds for $n$ large enough, since $E(R_m)$ goes to $\infty$ as $n$ goes to $\infty$. Substituting (62.21) into (62.17) and (62.16), we obtain

$$P\Big(\big|\mathcal{P}_2^+(A^\circ,\omega)\big| \geq Lp^2\Big) \geq \frac{1}{2} \sum_{\substack{k:|k-n/2|\leq p\\ j_2\geq p^2 e^{-3\alpha^2}/2}} P\big(|J_2| = j_2,\, \omega \in A^\circ \cap \mathcal{L}(n-k)\big)$$

$$= \frac{1}{2} \sum_{k:|k-n/2|\leq p} P\Big(|J_2| \geq \frac{p^2}{2} e^{-3\alpha^2} \,\Big|\, \omega \in A^\circ \cap \mathcal{L}(n-k)\Big) P\big(\omega \in A^\circ \cap \mathcal{L}(n-k)\big)\,. \tag{62.22}$$

We finally use the inequalities (62.10) and (62.6) to conclude from (62.22) that

$$P\Big(\big|\mathcal{P}_2^+(A^\circ,\omega)\big| \geq Lp^2\Big) \geq \frac{1}{4} \sum_{k:|k-n/2|\leq p} P\big(\omega \in A^\circ \cap \mathcal{L}(n-k)\big)$$

$$= \frac{1}{4} \sum_{k:|k-n/2|\leq p} P\big(\omega \in A^\circ \,\big|\, \omega \in \mathcal{L}(n-k)\big) P\big(\omega \in \mathcal{L}(n-k)\big)$$

$$\geq \frac{c(\alpha)}{4} \sum_{k:|k-n/2|\leq p} P\big(\omega \in \mathcal{L}(n-k)\big)\,.$$

The central limit theorem ensures that this last sum is larger than a fixed positive constant. So we conclude that there exists a positive constant $c$ such that, for $n$ large enough,

$$P\Big(\big|\mathcal{P}_2^+(A^\circ,\omega)\big| \geq cn\Big) \geq c\,.$$

This inequality furnishes a lower bound on the maximum $M_2(n)$.

**Proposition 62.1.** *There exists a positive constant $c$ such that, for any $n \geq 1$,*

$$M_2(n) = \max\big\{\, \phi_2(A) : A \text{ increasing subset of } \{0,1\}^n \,\big\} \geq cn\,.$$

Furthermore, we could proceed in the same way for the pivots of higher order, and we would obtain a similar result.

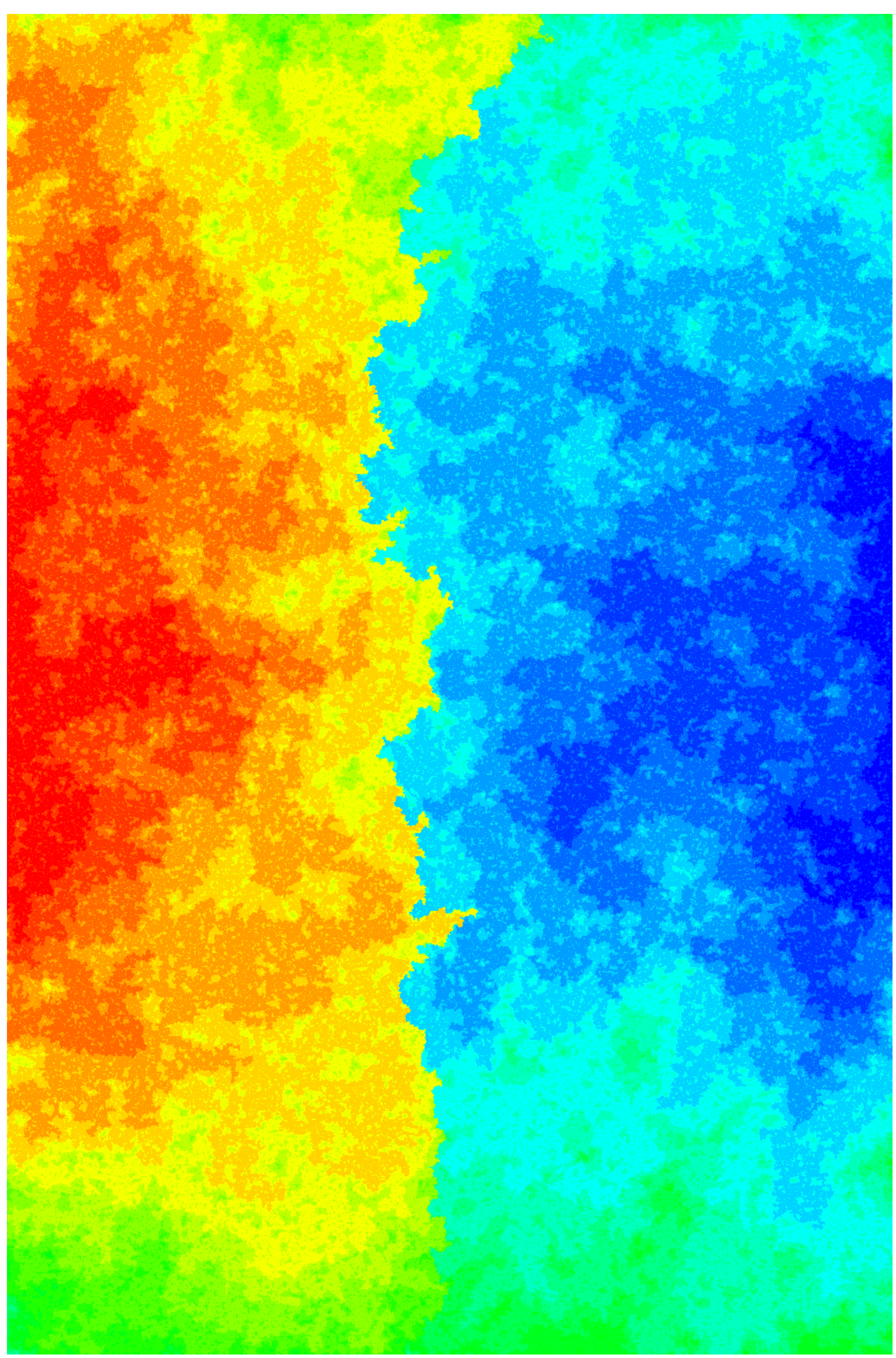

Figure 141: The heart of the matter: disconnection between the middle thirds of left and right sides in the rectangle $1550 \times 1024$, bond percolation, $p = 0.49$.

# Part IX
# The $\mathcal{V}$-local pivots ⊛

The example of Talagrand presented in section 62 has shattered all hope to obtain a bound on the pivots of order $k \geq 2$ for a general monotone event. The relevant bound, if it exists, should be sought in a more restricted framework. One possibility would be to come back directly to the percolation model and to study exclusively a disconnection event like those presented in the three challenges. However, this avenue has already been thoroughly explored. Even for the second pivots, we were unable to reach a satisfactory conclusion. In this part, we try an intermediate approach. While we stay at the general level of the product space $\{0,1\}^n$, we aim at controlling specific subsets of the pivots of order $k$, which we call the $\mathcal{V}$-local pivots. In section 63, we explain why we need to localize the search of the pivots if we want to get an upper bound of the correct order. We define then the $\mathcal{V}$-local pivots in section 64, and we give the formula which is the starting point to control them. However, there is still a long way to go before obtaining an operational inequality. In section 65, we reorganize the starting formula to prepare the ground for the ensuing computations. This reorganization is done with the help of a new quantity called the pivotality. The links between the signature and the pivotality are explored in section 66. The next goal is to prove a recursive inequality to control the pivots of order $k$. In fact, this recursive inequality will concern the $k$ pivotal sets of cardinality $k$, which are introduced in section 67. After proving an important bound on the signature in section 68, the recursive inequality is finally obtained in section 69. Here is what is looks like:

$$E\Big(\Big|\Big\{\, I\in\mathcal{I}(\mathcal{V},k):\text{pivotality}(I,\omega)=k\,\Big\}\Big|1_A(\omega)\Big)\ \leq\ \Big|E\big(L_k^{\mathcal{V}}1_A\big)\Big| +$$
$$\sum_{\substack{1\leq\ell\leq k-1\\0\leq h\leq k-\ell}}\binom{k}{h}\frac{E\Big(\Big|\Big\{\, J\in\mathcal{I}(\mathcal{V},\ell):\text{pivotality}(J,\omega)=\ell\,\Big\}\Big|1_A(\omega)\Big)}{(1-p)^{k-\ell}}\,M(n)\binom{V(n)}{k-\ell}\,.$$

Afterwards, we show how to exploit the recursive inequality in order to obtain quantitative bounds. In section 70, we perform an induction leading to a nice upper bound on the expected cardinality of the $k$ pivotal sets of cardinality $k$. This upper bound is converted into an upper bound on the expected cardinality of the $\mathcal{V}$-local pivots of order $k$ in section 71. Although this is a crucial step, it will not be enough to make our strategy work in the percolation model. In fact, we will need an inequality controlling the upper deviations from the expected mean. In section 72, we prove a general deviations inequality involving a random variable whose tail obeys a stretched exponential behavior. We apply then this inequality in section 73. The main input there is a classical large deviations inequality for the Bernoulli chaos, stemming from hypercontractivity. In the final section of this part, section 74, we present an exponential inequality which gives a control on the probability of an event in terms of the conditional expectation of the cardinality of the $\mathcal{V}$-local pivots of order $k$ associated to it.

# 63 Localization of the exploration

The example of Talagrand presented in the previous section 62 shows that, unfortunately, the cardinality of the set $\mathcal{P}_k^+(A,\omega)$ of the positive pivots of any order $k\geq 2$ of a monotone subset $A$ of $\{0,1\}^n$ can potentially reach the order $n$. Yet we know that $\mathcal{P}_k^+(A,\omega)\subset\mathcal{P}_k(A,\omega)$, and we proved in corollary 44.2 that

$$\forall k\geq 1\qquad \mathcal{P}_k(A)\,=\,\bigcup_{1\leq\ell\leq k}\pi\Big(\mathcal{S}^{\ell,-1}(A)\cup\mathcal{S}^{\ell,+1}(A)\Big)\,,\tag{63.1}$$

where the sets $\mathcal{S}^{\ell,s}(A,\omega)$ were defined in (44.2), they are subsets of $\{\,1,\dots,n\,\}^\ell$, and $\pi$ is the projection operator defined as follows: for any collection $\mathcal{I}$ of subsets of $\{\,1,\dots,n\,\}$,

$$\pi(\mathcal{I})\,=\,\Big\{\,i\in\{\,1,\dots,n\,\}:\exists I\in\mathcal{I}\,,\quad i\in I\,\Big\}\,=\,\bigcup_{I\in\mathcal{I}}I\,.\tag{63.2}$$

On top of that, the generalized Margulis-Russo formula stated in theorem 41.2, together with the expression (43.4) of the $\ell$-th derivative of $P_p(A)$ involving the Bernoulli chaos of order $\ell$, yield that

$$\sum_{-2^{\ell-1}\leq s\leq 2^{\ell-1}}s\,E\big(\big|\mathcal{S}^{\ell,s}(A)\big|\big)\,=\,\sum_{\substack{1\leq i_1,\dots,i_\ell\leq n\\ \text{pairwise distinct}}}E\Big(X_{i_1}(\omega)\,\cdots\,X_{i_\ell}(\omega)1_A(\omega)\Big)\,.\tag{63.3}$$

Applying the Cauchy–Schwarz inequality to (63.3), and using proposition 43.2, we obtain that

$$\Big|\sum_{-2^{\ell-1}\leq s\leq 2^{\ell-1}}s\,E\big(\big|\mathcal{S}^{\ell,s}(A)\big|\big)\Big|\,\leq\,\ell!\sqrt{P_p(A)\binom{n}{\ell}\frac{1}{p^\ell(1-p)^\ell}}\,.\tag{63.4}$$

The chain of formulas (63.1), (63.3), (63.4) generate a constraint on the law of the random set $\mathcal{P}_k^+(A)$. Unfortunately, this constraint is not strong enough to get a useful bound on $\big|\mathcal{P}_k^+(A)\big|$. The main weakness of the previous inequality (63.4) seems to be the exponent $\ell/2$ of $n$ in the upper bound. To control successfully the disconnection event in percolation, we would have liked to have an upper bound of order $\sqrt{n}$. The power $\ell/2$ stems from the fact that formula (63.3) involves the full Bernoulli chaos of order $\ell$, which is a homogeneous polynomial of degree $\ell$ with $\binom{n}{\ell}$ terms. When dealing with a generic monotone event on $\{0,1\}^n$, there is limited hope of improving the previous inequalities. Indeed, to detect all the elements of the sets $\mathcal{S}^{\ell,s}(A,\omega)$, it seems necessary to try all the subsets of $\{\,1,\dots,n\,\}$ of cardinality $\ell$, and this is the mission accomplished by the Bernoulli chaos of order $\ell$. If we want to use this strategy and to end up with an upper bound of order $\sqrt{n}$, then the only hope is to work with a homogeneous polynomial having significantly less terms than the Bernoulli chaos. This amounts to localize the search to a restricted subset of the pivots. Inevitably, we might miss some pivots of order $\ell$. Later on, we will introduce a condition guaranteeing that the relevant pivots are correctly detected by a restricted homogeneous polynomial.

## 64 The starting point

A covering $\mathcal{V}$ of the index set $\{1,\dots,n\}$ is a collection of subsets of $\{1,\dots,n\}$ whose union is $\{1,\dots,n\}$, i.e.,

$$\{1,\dots,n\} = \bigcup_{V\in\mathcal{V}} V\,.$$

The constructions that follow depend on a covering $\mathcal{V}$ of $\{1,\dots,n\}$. When the time comes to apply the results to the percolation model, we will make a specific choice for $\mathcal{V}$. For the time being, we leave $\mathcal{V}$ as a superscript in the objects that depend on $\mathcal{V}$. For $i\in\{1,\dots,n\}$, we denote $\mathcal{V}(i)$ the sets of the covering $\mathcal{V}$ which contain $i$, i.e.,

$$\mathcal{V}(i) = \big\{\,V : V\in\mathcal{V},\, i\in V\,\big\}\,.$$

We denote by $M(n)$ the maximal multiplicity of the covering, i.e.,

$$M(n) = \max\Big\{\big|\mathcal{V}(i)\big| : i\in\{1,\dots,n\}\Big\}\,, \tag{64.1}$$

and by $V(n)$ the maximal cardinality of a set belonging to the covering, i.e.,

$$V(n) = \max\big\{\,|V| : V\in\mathcal{V}\,\big\}\,. \tag{64.2}$$

In the future applications, we will choose the covering so that $M(n)$ and $V(n)$ are negligible compared to $n$. We define next what we mean by a $\mathcal{V}$-local pivot. In short, a $\mathcal{V}$-local pivot $i$ is a pivot in the sense of the definition given in subsection 44.1, with the additional requirement that the associated subset $I$ and the pivot $i$ must be included in an element of the covering $\mathcal{V}$. For the sake of completeness, we rewrite next the full definition of a $\mathcal{V}$-local pivot.

Let $A\subset\Omega$ be an event. Let $\omega$ be a configuration of $\Omega$ and let $k\geq 1$. We say that a component $i$ in $\{1,\dots,n\}$ is a $\mathcal{V}$-local pivot of order $k$ for the event $A$ in the configuration $\omega$ if there exist $V$ in $\mathcal{V}(i)$, a subset $I$ of $V\setminus\{i\}$ of cardinality strictly less than $k$, and a configuration $\sigma$ in $\Omega(I)$ such that, if we transform the configuration $\omega$ by forcing the components $I$ to be in the same states as in $\sigma$, then the component $i$ becomes pivotal for $A$. We denote by $\mathcal{P}_k^{\mathcal{V}}(A,\omega)$ the set of the $k$-pivotal components for $A$ in the configuration $\omega$, hence

$$\begin{multline}\mathcal{P}_k^{\mathcal{V}}(A,\omega) = \Big\{ i\in\{1,\dots,n\} : \exists\, V\in\mathcal{V}(i) \quad \exists\, I\subset V\setminus\{i\} \quad |I|\leq k-1 \\ \exists\,\sigma\in\Omega(I)\quad i\in\mathcal{P}\big(A,\omega(I\leftarrow\sigma)\big)\Big\}\,, \end{multline} \tag{64.3}$$

where the notation $\omega(I\leftarrow\sigma)$ was defined in (41.2). We make the convention that $\mathcal{P}_0^{\mathcal{V}}(A,\omega)=\varnothing$. For $k=1$, the set $\mathcal{P}_1^{\mathcal{V}}(A,\omega)$ coincides with the classical set $\mathcal{P}(A,\omega)$ of the pivots for $A$ in $\omega$.

Let us fix $k\geq 1$. We design next an adequate mechanism to detect the $\mathcal{V}$-local pivots of order $k$. With the help of the covering $\mathcal{V}$, we define a specific deterministic subset $\mathcal{I}(\mathcal{V},k)$ of the subsets of $\{1,\dots,n\}$ having cardinality $k$ by setting

$$\mathcal{I}(\mathcal{V},k) \,=\, \bigcup_{V\in\mathcal{V}} \bigcup_{\substack{i_1,\dots,i_k\in V\\ \text{pairwise distinct}}} \Big\{ i_1,\dots,i_k \Big\}.$$

The number of elements in $\mathcal{V}$ is bounded by $nM(n)$, hence

$$\forall k\in\{1,\dots,n\} \qquad \big|\mathcal{I}(\mathcal{V},k)\big| \,\leq\, n\,M(n)\,V(n)^k\,. \tag{64.4}$$

The key instrument to detect the $\mathcal{V}$-local pivots of order $k$ is the homogeneous polynomial $L_k^{\mathcal{V}}(\omega)$ of order $k$ given by

$$L_k^{\mathcal{V}}(\omega) \,=\, \sum_{I\in\mathcal{I}(\mathcal{V},k)} \prod_{i\in I} X_i(\omega)\,.$$

The detection of the $\mathcal{V}$-local pivots of order $k$ of the event $A$ is simply achieved by taking the expectation of the product of $L_k^{\mathcal{V}}$ and $1_A$. To see this, we shall redo some computations that were done in section 41, but in a different order. In section 41, we started by differentiating several times the probability $P_p(A)$. When trying to generalize the Russo-Margulis formula, we obtained an expression involving the Bernoulli chaos. Now, we wish to obtain analogous formulas for a restricted homogeneous polynomial, hence we redo the relevant steps of the computation starting with this polynomial. First, we have

$$E\big(L_k^{\mathcal{V}}1_A\big) \,=\, \sum_{I\in\mathcal{I}(\mathcal{V},k)} E\Big(1_A\prod_{i\in I} X_i\Big)\,. \tag{64.5}$$

Let $I$ be a fixed element of $\mathcal{I}(\mathcal{V},k)$. Applying $k$ times formula (39.2), we obtain

$$E\Big(1_A\prod_{i\in I} X_i\Big) \,=\, E\Big(\prod_{i\in I}\delta_i 1_A\Big)\,. \tag{64.6}$$

We recall that the signature of a Boolean function $f$ defined on $\Omega(I)$ is

$$\mathcal{S}(f) \,=\, \sum_{\sigma\in\Omega(I)} (-1)^{|I|-o(\sigma)} f(\sigma)\,. \tag{64.7}$$

Denoting by $1_A\big(\omega(I\leftarrow\cdot)\big)$ the Boolean function $\sigma\in\Omega(I)\mapsto 1_A\big(\omega(I\leftarrow\sigma)\big)$ (where $\omega(I\leftarrow\sigma)$ was defined in (41.2)), we have, using formula (41.6)

$$\prod_{i\in I}\delta_i\,1_A(\omega) \,=\, \mathcal{S}\big(1_A(\omega(I\leftarrow\cdot))\big)\,. \tag{64.8}$$

Substituting the identities (64.6) and (64.8) into (64.5), we obtain

$$E\big(L_k^{\mathcal{V}}1_A\big) \,=\, \sum_{I\in\mathcal{I}(\mathcal{V},k)} E\Big(\mathcal{S}\big(1_A(\omega(I\leftarrow\cdot))\big)\Big)\,. \tag{64.9}$$

Formula (64.9) is the starting point to obtain a control on the $\mathcal{V}$-local pivots of order $k$.

## 65 Reorganization of the sum

We estimate first the left-hand side of formula (64.9). Using the fact that the random variables $X_i$ are i.i.d. and centered, we have, thanks to the bound (64.4),

$$E\big((L_k^{\mathcal{V}} 1_A)^2\big) \,\leq\, \frac{\big|\mathcal{I}(\mathcal{V},k)\big|}{p^k(1-p)^k} \,\leq\, \frac{nM(n)V(n)^k}{p^k(1-p)^k}\,. \tag{65.1}$$

We examine next the right-hand side of formula (64.9). In a first step, we move back the summation inside the expectation

$$E\big(L_k^{\mathcal{V}} 1_A\big) \,=\, E\bigg(\sum_{I\in\mathcal{I}(\mathcal{V},k)} \mathcal{S}\big(1_A(\omega(I\leftarrow\cdot))\big)\bigg)\,. \tag{65.2}$$

Let us fix a configuration $\omega$ and a set $I$ in $\mathcal{I}(\mathcal{V},k)$. The map $1_A\big(\omega(I\leftarrow\cdot)\big)$ is a Boolean function defined on $\Omega(I)$. The set $I$ is a subset of $\{\,1,\dots,n\,\}$ having cardinality $k$. We can list the elements of $I$ in ascending order, and identify $I$ with $\{\,1,\dots,k\,\}$. Once this identification is done, a configuration $\sigma$ in $\Omega(I)$ can be seen as an element of $\{0,1\}^k$, and the map $1_A(\omega(I\leftarrow\cdot))$ as a Boolean function defined on $\{0,1\}^k$. We denote this Boolean function by $\mathcal{F}(A,\omega,I)$. Let us denote by $\mathcal{F}(k)$ the set of all the Boolean functions defined on $\{0,1\}^k$. We rewrite the sum in (65.2) as

$$\sum_{I\in\mathcal{I}(\mathcal{V},k)} \mathcal{S}\big(1_A(\omega(I\leftarrow\cdot))\big) \,=\, \sum_{f\in\mathcal{F}(k)} \sum_{\substack{I\in\mathcal{I}(\mathcal{V},k)\\ \mathcal{F}(A,\omega,I)=f}} \mathcal{S}(f)\,, \tag{65.3}$$

where the signature $\mathcal{S}(f)$ of a function $f$ in $\mathcal{F}(k)$ is now defined by

$$\mathcal{S}(f) \,=\, \sum_{\sigma\in\{0,1\}^k} (-1)^{k-o(\sigma)} f(\sigma)\,.$$

Only the functions $f$ such that $\mathcal{S}(f)\neq 0$ contribute to the above sum. Let us denote by $\mathcal{F}^*(k)$ the subset of $\mathcal{F}(k)$ consisting of the functions having a non-vanishing signature:

$$\mathcal{F}^*(k) \,=\, \big\{\, f\in\mathcal{F}(k) : \mathcal{S}(f)\neq 0\,\big\}\,.$$

Substituting (65.3) into (65.2), and taking the sums over $f, I$ outside the expectation, we obtain

$$E\big(L_k^{\mathcal{V}} 1_A\big) \,=\, \sum_{f\in\mathcal{F}^*(k)} \mathcal{S}(f) \sum_{I\in\mathcal{I}(\mathcal{V},k)} P\big(\mathcal{F}(A,\omega,I)=f\big)\,. \tag{65.4}$$

Let us fix a function $f$ in $\mathcal{F}(k)$ and let us denote by $\operatorname{supp} f$ its support, defined as

$$\operatorname{supp} f \,=\, \big\{\,\sigma\in\{0,1\}^k : f(\sigma)=1\,\big\}\,.$$

The event $\{\,\mathcal{F}(A,\omega,I)=f\,\}$ does not depend on the restriction of $\omega$ to $I$. Let us endow the set $\{0,1\}^k$ with the Bernoulli probability measure of parameter $p$. We have then

$$P\big(A\,\big|\,\mathcal{F}(A,\omega,I)=f\big)\;=\;P\big(\operatorname{supp} f\big)\,,$$

whence

$$\begin{aligned}P\big(A\cap\mathcal{F}(A,\omega,I)=f\big)\;&=\;P\big(A\,\big|\,\mathcal{F}(A,\omega,I)=f\big)P\big(\mathcal{F}(A,\omega,I)=f\big)\\&=\;P\big(\operatorname{supp} f\big)P\big(\mathcal{F}(A,\omega,I)=f\big)\,.\end{aligned}\tag{65.5}$$

If in addition $f$ is in $\mathcal{F}^*(k)$, then $f$ is not constant, so that $0<P\big(\operatorname{supp} f\big)<1$ and we can rewrite (65.5) as

$$P\big(\mathcal{F}(A,\omega,I)=f\big)\;=\;\frac{1}{P\big(\operatorname{supp} f\big)}P\big(A\cap\mathcal{F}(A,\omega,I)=f\big)\,.\tag{65.6}$$

Let us define the renormalized signature $\overline{\mathcal{S}}$ by setting

$$\forall f\in\mathcal{F}^*(k)\qquad\overline{\mathcal{S}}(f)\;=\;\frac{\mathcal{S}(f)}{P\big(\operatorname{supp} f\big)}\,.$$

Substituting (65.6) into (65.4), we get

$$E\big(L_k^{\mathcal{V}}1_A\big)\;=\;\sum_{f\in\mathcal{F}^*(k)}\overline{\mathcal{S}}(f)\sum_{I\in\mathcal{I}(\mathcal{V},k)}P\big(A\cap\mathcal{F}(A,\omega,I)=f\big)\,.$$

In the next step, we write the probability as an expectation, and we put back the summations inside the expectation:

$$\begin{aligned}E\big(L_k^{\mathcal{V}}1_A\big)\;&=\;\sum_{f\in\mathcal{F}^*(k)}\overline{\mathcal{S}}(f)\sum_{I\in\mathcal{I}(\mathcal{V},k)}E\Big(1_A(\omega)1_{\{\,\mathcal{F}(A,\omega,I)=f\,\}}\Big)\\&=E\Bigg(\sum_{f\in\mathcal{F}^*(k)}\sum_{I\in\mathcal{I}(\mathcal{V},k)}\overline{\mathcal{S}}(f)1_A(\omega)1_{\{\,\mathcal{F}(A,\omega,I)=f\,\}}\Bigg)\\&=E\Bigg(\sum_{I\in\mathcal{I}(\mathcal{V},k)}\sum_{f\in\mathcal{F}^*(k)}\overline{\mathcal{S}}(f)1_A(\omega)1_{\{\,\mathcal{F}(A,\omega,I)=f\,\}}\Bigg)\\&=E\Bigg(\sum_{I\in\mathcal{I}(\mathcal{V},k)}\overline{\mathcal{S}}\big(1_A(\omega(I\leftarrow\cdot))\big)1_A(\omega)\Bigg)\,,\end{aligned}\tag{65.7}$$

where in the last line we use the renormalized signature $\overline{\mathcal{S}}$ of a Boolean function $f$ on $\Omega(I)$, which is naturally defined as

$$\overline{\mathcal{S}}(f)\;=\;\frac{1}{P\big(\operatorname{supp} f\big)}\sum_{\sigma\in\Omega(I)}(-1)^{|I|-o(\sigma)}f(\sigma)\;=\;\frac{1}{P\big(\operatorname{supp} f\big)}\mathcal{S}(f)\,.\tag{65.8}$$

The formula (65.8) does not make sense for the null function. One possibility is to declare that $\overline{\mathcal{S}}$ is 0 for the null function. Another one is to restrict ourselves to the functions having a non-vanishing signature, which are the only ones which contribute to the sums. So, we introduce the set

$$\mathcal{S}^*\big(k,\omega\big) \,=\, \Big\{\, I\in\mathcal{I}(\mathcal{V},k) : \mathcal{S}\big(1_A(\omega(I\leftarrow\cdot))\big)\neq 0 \,\Big\}\,,$$

and we rewrite (65.7) as

$$E\big(L_k^{\mathcal{V}}1_A\big) \,=\, E\Bigg(\sum_{I\in\mathcal{S}^*(k,\omega)}\overline{\mathcal{S}}\big(1_A(\omega(I\leftarrow\cdot))\big)1_A(\omega)\Bigg)\,. \tag{65.9}$$

The advantage of formula (65.9) over formula (65.2) is that it lends itself well to conditioning by the event $A$, as we will do later. For the time being, we wish to obtain a recursive formula in order to control the $\mathcal{V}$-local pivots of order $k$. To that end, we decompose further the sum appearing in the expectation in (65.9) with the help of the pivotality function, that we introduce next.

Let $\ell\in\{\,1,\dots,n\,\}$. For $I$ a subset of $\{\,1,\dots,n\,\}$, we denote by $\pi_\ell(I)$ the collection of the subsets of $I$ having cardinality $\ell$, i.e.,

$$\pi_\ell(I) \,=\, \Big\{\,\{\,i_1,\dots,i_\ell\,\} : i_1,\dots,i_\ell \text{ pairwise distinct elements of } I\,\Big\}\,.$$

We extend next the projection $\pi_\ell$ to collections of subsets: for any collection $\mathcal{I}$ of subsets of $\{\,1,\dots,n\,\}$, we define

$$\pi_\ell(\mathcal{I}) \,=\, \bigcup_{I\in\mathcal{I}}\pi_\ell(I)\,.$$

Given a configuration $\omega$, we define the pivotality of a subset $I$ of $\{\,1,\dots,n\,\}$ as

$$\text{pivotality}(I,\omega) \,=\, \min\,\Big\{\,\ell\in\{\,1,\dots,n\,\} : \pi_\ell(I)\cap\mathcal{S}^*(\ell,\omega)\neq\varnothing\,\Big\}\,. \tag{65.10}$$

We make the usual convention that $\min\varnothing=+\infty$, thus the pivotality of $I$ might be infinite; this means that, in the configuration $\omega$, the occurrence of $A$ is not decided by the components in $I$. In principle, we should write pivotality$(A,I,\omega)$, but since the event $A$ is fixed throughout the computations, we omit it when possible. The set $\mathcal{S}^*(k,\omega)$ can then be written as the disjoint union

$$\mathcal{S}^*(k,\omega) \,=\, \bigcup_{1\leq\ell\leq k}\,\big\{\,I\in\mathcal{S}^*(k,\omega) : \text{pivotality}(I) \,=\, \ell\,\big\}\,. \tag{65.11}$$

In a second step, we take advantage of the decomposition (65.11) to rewrite the sum in (65.2) as

$$\sum_{I\in\mathcal{S}^*(k,\omega)}\cdots \quad=\sum_{1\leq\ell\leq k}\;\sum_{\substack{I\in\mathcal{S}^*(k,\omega)\\ \text{pivotality}(I,\omega)=\ell}}\cdots\quad. \tag{65.12}$$

Reporting (65.12) into (65.9), and taking the sum over $\ell$ outside the expectation, we obtain

$$E\big(L_k^{\mathcal{V}}1_A\big) \,=\, \sum_{1\leq\ell\leq k}E\Bigg(\sum_{\substack{I\in\mathcal{S}^*(k,\omega)\\ \text{pivotality}(I,\omega)=\ell}}\overline{\mathcal{S}}\big(1_A(\omega(I\leftarrow\cdot))\big)1_A(\omega)\Bigg)\,. \tag{65.13}$$

## 66 Signature and pivotal sets

We prove here two results on the signature $\mathcal{S}\big(1_A(\omega(I \leftarrow \cdot))\big)$ for $I \in \mathcal{I}(\mathcal{V}, k)$. To simplify the exposition, we introduce a specific definition for a pivotal set. It is far from obvious to extend to a set the notion of a pivotal component. There exist various possible generalizations which have been employed in the theory of boolean functions. Biswas and Sarkar [16] discuss this problem and present various applications. As one might expect, different definitions serve different purposes. We choose here a simple one, that is suited only for monotone sets.

For $\omega$ a configuration in $\Omega$ and $I$ a subset of $\{0,1\}^n$, we denote by $\omega^I$ (respectively $\omega_I$) the configuration $\omega$ where the components whose index is in $I$ are set to 1 (respectively to 0), that is,

$$\forall j \in \{\,1,\dots,n\,\}$$
$$\omega^I(j) \,=\, \begin{cases} \omega(j) & \text{if } j \notin I \\ \ \ 1 & \text{if } j \in I \end{cases}, \qquad \omega_I(j) \,=\, \begin{cases} \omega(j) & \text{if } j \notin I \\ \ \ 0 & \text{if } j \in I \end{cases}.$$

**Definition 66.1.** *Let $A$ be a monotone subset of $\{0,1\}^n$ and let $\omega$ be a configuration in $\{0,1\}^n$. A subset $I$ of $\{\,1,\dots,n\,\}$ is said to be pivotal for $A$ in the configuration $\omega$ if $1_A(\omega^I) \neq 1_A(\omega_I)$.*

In short, a subset $I$ is pivotal in a configuration $\omega$ if the occurrence of $A$ depends on the status of the components of $\omega$ whose index belongs to $I$. The next lemma formalizes the obvious fact that, for a monotone event, when a subset $I$ has a non-vanishing signature, it is pivotal.

**Lemma 66.2.** *Let $A$ be a monotone subset of $\{0,1\}^n$. Let $\omega$ be a configuration in $\{0,1\}^n$ and let $I$ be a subset of $\{\,1,\dots,n\,\}$ such that $\mathcal{S}\big(1_A(\omega(I \leftarrow \cdot))\big) \neq 0$. Then $I$ is pivotal for $A$ in $\omega$.*

*Proof.* Let $A, \omega$ and $I$ be as in the statement of the lemma. Recalling the definition (64.7) of the signature, we see that there exist $\sigma_0, \sigma_1$ in $\Omega(I)$ such that

$$1_A\big(\omega(I \leftarrow \sigma_0)\big) \,=\, 0\,, \qquad 1_A\big(\omega(I \leftarrow \sigma_1)\big) \,=\, 1\,. \tag{66.1}$$

Obviously, we have

$$\omega_I \,\leq\, \omega(I \leftarrow \sigma_0) \leq \omega^I\,, \qquad \omega_I \,\leq\, \omega(I \leftarrow \sigma_1) \leq \omega^I\,. \tag{66.2}$$

Since the set $A$ is monotone, the equalities (66.1) and the inequalities (66.2) readily imply that $1_A(\omega^I) \neq 1_A(\omega_I)$. ☐

The next lemma is a sort of weak converse of lemma 66.2, it will help to detect the existence of subsets having a non-vanishing signature.

**Lemma 66.3.** *Let $A$ be a monotone subset of $\{0,1\}^n$. Let $\omega$ be a configuration in $\{0,1\}^n$ and let $I$ be a subset of $\{\,1,\dots,n\,\}$ which is pivotal for $A$ in $\omega$. Then there exists a subset $J$ of $I$ such that $\mathcal{S}\big(1_A(\omega(J \leftarrow \cdot))\big) \in \{-1,+1\}$.*

*Proof.* Let $A, \omega, I$ be as in the statement of the lemma. We do the proof only for an increasing set $A$ (the arguments for a decreasing set are very similar). Let us consider the collection of sets

$$\mathcal{I}(\omega, I) \;=\; \Big\{\, J \subset I : \omega_J \notin A,\, \omega^J \in A \,\Big\}.$$

It follows from the hypothesis that $I$ is in $\mathcal{I}$, therefore $\mathcal{I}$ is not empty. Let $J^*$ be an element of $\mathcal{I}(\omega, I)$ which is minimal for the set inclusion. We consider two cases, depending on whether $\omega$ is in $A$ or not.

• $\omega \in A$. Suppose that there exists $j$ in $J^*$ such that $\omega(j) = 0$. Then the set $K = J^* \setminus \{\, j \,\}$ would satisfy

$$\omega_K = \omega_J \notin A\,, \quad \omega^K \geq \omega \in A\,.$$

Since $A$ is increasing, we see that $\omega^K \in A$, thus $K$ belongs to $\mathcal{I}(\omega, I)$, but this would contradict the minimality of $J^*$. Therefore $J^*$ is included into $\operatorname{supp}\omega$. Furthermore, it follows from the minimality of $J^*$ that

$$\forall K \subset J\,, \quad K \neq J\,, \qquad \omega_K \in A\,.$$

We are now in position to compute the signature of $1_A(\omega(J^* \leftarrow \cdot))$:

$$\begin{aligned}
\mathcal{S}\big(1_A(\omega(J^* \leftarrow \cdot))\big) \;&=\; \sum_{\sigma \in \Omega(J^*)} (-1)^{|J^*| - o(\sigma)} 1_A(\omega(J^* \leftarrow \sigma)) \\
&=\; \sum_{0 \leq k \leq |J^*|} \sum_{\substack{\sigma \in \Omega(J^*) \\ o(\sigma) = k}} (-1)^{|J^*| - o(\sigma)} 1_A(\omega(J^* \leftarrow \sigma)) \\
&=\; \sum_{1 \leq k \leq |J^*|} \sum_{\substack{\sigma \in \Omega(J^*) \\ o(\sigma) = k}} (-1)^{|J^*| - k} \;=\; \sum_{1 \leq k \leq |J^*|} \binom{|J^*|}{k} (-1)^{|J^*| - k} \;=\; -1\,.
\end{aligned}$$

This completes the case where $\omega \in A$.

• $\omega \notin A$. The arguments are similar to the previous case, so we only outline the differences. We have that $J^*$ is included into $\{\, 1, \dots, n \,\} \setminus \operatorname{supp}\omega$ and

$$\forall K \subset J\,, \quad K \neq J\,, \qquad \omega^K \notin A\,.$$

It follows that the only term which does not vanish in the expression of the signature is the one associated to the configuration $\sigma$ having all its components equal to 1, thus $\mathcal{S}\big(1_A(\omega(J^* \leftarrow \cdot))\big) = 1$. ☐

With the help of the two previous lemmas, we can reformulate the definition of the quantity $\text{pivotality}(I, \omega)$ associated to a monotone event $A$.

**Corollary 66.4.** *Let $A$ be a monotone subset of $\{0,1\}^n$ and let $I$ be a subset of $\{\, 1, \dots, n \,\}$. The pivotality of $I$ in a configuration $\omega$ of $\{0,1\}^n$ is equal to*

$$\textit{pivotality}(I, \omega) \;=\; \min\, \Big\{\, |J| : J \subset I,\;\; J \textit{ pivotal set for } A \textit{ in } \omega \,\Big\}.$$

*Proof.* Lemmas 66.2 and 66.3 together imply that, for a monotone set $A$, the above definition of the pivotality coincides with the one given in (65.10). ☐

## 67 The $\ell$ pivotal sets of cardinality $k$

Notice that, for $k \geq 2$ and $I \in \mathcal{I}(\mathcal{V}, k)$, the event $\{ \text{pivotality}(I, \omega) = k \}$ depends also on the components of $\omega$ whose index belongs to $I$. In this sense, the event $\{ \text{pivotality}(I, \omega) = k \}$ differs notably from the event $\{ I \text{ is a pivotal set} \}$, as the next lemma witnesses.

**Lemma 67.1.** *Let $A$ be a decreasing subset of $\{0,1\}^n$. A set $I$ in $\mathcal{I}(\mathcal{V}, 1)$ is a singleton $I = \{ i \}$ and we have the equivalence*

$$\text{pivotality}(\{ i \}, \omega) = 1 \quad \Longleftrightarrow \quad i \in \mathcal{P}^-(A, \omega) \quad \Longleftrightarrow \quad \mathcal{S}\big(1_A(\omega(I \leftarrow \cdot))\big) = 1 \,.$$

*Let $k \in \{ 2, \dots, n \}$ and let $I$ belong to $\mathcal{I}(\mathcal{V}, k)$. Let $\omega$ be a configuration such that pivotality$(I, \omega) = k$. We have*

$$\begin{aligned} \omega \in A \quad &\Longrightarrow \quad \forall i \in I \quad \omega(i) = 0 \,, \quad \mathcal{S}\big(1_A(\omega(I \leftarrow \cdot))\big) = -1 \,, \\ \omega \notin A \quad &\Longrightarrow \quad \forall i \in I \quad \omega(i) = 1 \,, \quad \mathcal{S}\big(1_A(\omega(I \leftarrow \cdot))\big) = (-1)^k \,. \end{aligned}$$

*Proof.* The case $k = 1$ is straightforward. Let $k \geq 2$ and let $\omega$ be a configuration such that pivotality$(I, \omega) = k$. Necessarily, the set $I$ is pivotal for $A$ in $\omega$. Moreover, it follows from lemma 66.3 that no strict subset of $I$ is pivotal for $A$ in $\omega$. We consider two cases:

• $\omega \in A$. Let $J = I \setminus \operatorname{supp} \omega$. We have

$$\omega^J = \omega^I \notin A, \quad \omega_J = \omega \in A \,,$$

thus $J$ is not empty and it is a pivotal set for $A$ in $\omega$. Therefore $J = I$ and $I \cap \operatorname{supp} \omega = \varnothing$. It follows that $\omega(i) = 0$ for $i \in I$ and moreover

$$\forall J \subset I \quad J \neq I \quad \Longrightarrow \quad \omega^J \in A \,.$$

We have all the information required to compute the signature:

$$\begin{aligned} \mathcal{S}\big(1_A(\omega(I \leftarrow \cdot))\big) \;&=\; \sum_{\sigma \in \Omega(I)} (-1)^{|I| - o(\sigma)} 1_A(\omega(I \leftarrow \sigma)) \\ &= \sum_{0 \leq k \leq |I|-1} \; \sum_{\substack{\sigma \in \Omega(I) \\ |\operatorname{supp} \sigma| = k}} (-1)^{|I| - o(\sigma)} \;=\; \sum_{0 \leq k \leq |I|-1} \binom{|I|}{k} (-1)^{|I|-k} \;=\; -1 \,. \end{aligned}$$

• $\omega \notin A$. Let $J = I \cap \operatorname{supp} \omega$. We have

$$\omega^J = \omega \notin A, \quad \omega_J = \omega_I \in A \,,$$

thus $J$ is not empty and it is a pivotal set for $A$ in $\omega$. Therefore $J = I$ and $I \subset \operatorname{supp} \omega$. It follows that $\omega(i) = 1$ for $i \in I$ and moreover $\omega(I \leftarrow \sigma)$ is not in $A$ unless $\operatorname{supp} \sigma$ is empty. Thus

$$\mathcal{S}\big(1_A(\omega(I \leftarrow \cdot))\big) \;=\; \sum_{\substack{\sigma \in \Omega(I) \\ |\operatorname{supp} \sigma| = 0}} (-1)^{|I| - o(\sigma)} \;=\; (-1)^k \,.$$

This completes the proof of the lemma. $\square$

We consider next the more delicate case of the $\ell$ pivotal sets of cardinality $k$. We will need the lower bound stated in the next lemma.

**Lemma 67.2.** *Let $A$ be a decreasing subset of $\{0,1\}^n$. Let $k \in \{2,\dots,n\}$, let $I$ belong to $\mathcal{I}(\mathcal{V},k)$ and let $\ell \in \{1,\dots,k\}$. Let $\omega$ be a configuration belonging to $A$ and such that pivotality$(I,\omega)=\ell$. We have*

$$P\Big(\operatorname{supp}\big(1_A(\omega(I\leftarrow\cdot))\big)\Big) \,\geq\, (1-p^\ell)(1-p)^{k-\ell}\,. \tag{67.1}$$

*Proof.* Let $k,\ell$, $A$ and $\omega$ be as in the hypothesis of the lemma. By corollary 66.4, there exists a subset $J$ of $I$ of cardinality $\ell$, such that $J$ is pivotal for $A$ in $\omega$, but no strict subset of $J$ is pivotal for $A$ in $\omega$. Let $K=J\setminus\operatorname{supp}\omega$. We have

$$\omega^K=\omega^J\notin A,\quad \omega_K=\omega\in A\,,$$

thus $K$ is not empty and it is a pivotal set for $A$ in $\omega$. Therefore $K=J$ and $I\cap\operatorname{supp}\omega=\varnothing$. It follows that $\omega(i)=0$ for $i\in J$ and moreover

$$\forall K\subset J\quad K\neq J\quad\Longrightarrow\quad \omega^K\in A\,. \tag{67.2}$$

Since the set $A$ is decreasing, the property (67.2) implies that

$$\forall K\subset J\quad K\neq J\quad\Longrightarrow\quad \big(\omega^K\big)_{I\setminus J}\in A\,.$$

The above property yields that

$$\begin{aligned}P\Big(\operatorname{supp}\big(1_A(\omega(I\leftarrow\cdot))\big)\Big) &\geq \sum_{K\subset J,K\neq J} P\Big(\big(\omega^K\big)_{I\setminus J}\Big) = \sum_{K\subset J,K\neq J} p^{|K|}(1-p)^{|I|-|K|}\\ &=(1-p)^{k-\ell}\sum_{K\subset J,K\neq J} p^{|K|}(1-p)^{|J|-|K|} \,=\, (1-p)^{k-\ell}(1-p^\ell)\,.\end{aligned}$$

This lower bound is uniform with respect to the configuration $\omega$ and it is the desired inequality (67.1). ☐

# 68 Bounding the signature

We try here to derive quantitative bounds on the signature of a subset in terms of its pivotality. The next proposition presents an important refinement of the second result of proposition 44.1. This result will be crucial to exploit the identity (65.13) and to get a control on the $\mathcal{V}$-local pivots.

**Proposition 68.1.** *Let $A$ be a monotone subset of $\{0,1\}^n$. For any $\omega$ in $\{0,1\}^n$, we have the implication*

$$\begin{aligned}&\forall k\in\{1,\dots,n\}\quad \forall I\in\mathcal{I}(\mathcal{V},k)\quad \forall\ell\in\{1,\dots,k\}\\ &\qquad \textit{pivotality}(I,\omega)=\ell\quad\Longrightarrow\quad \big|\mathcal{S}\big(1_A(\omega(I\leftarrow\cdot))\big)\big|\leq\sum_{0\leq h\leq k-\ell}\binom{k}{h}.\end{aligned} \tag{68.1}$$

*Proof.* Let us fix $k \geq 2$ and a configuration $\omega \in \Omega$. Let $I$ be an element of $\mathcal{I}(\mathcal{V}, k)$ such that pivotality$(I, \omega) = \ell$. It follows from the definition (41.4) and the formula (41.3) that

$$\mathcal{S}\big(1_A(\omega(I \leftarrow \cdot))\big) \,=\, \prod_{i \in I} \delta_i \, 1_A(\omega)\,. \tag{68.2}$$

Let $i^*$ be an index in $I$ and let us set $J = I \setminus \{\, i^* \,\}$. The operators $\delta_i$ commute. We put $\delta_{i^*}$ at the beginning of the product and we develop the remaining product over $J$ with the help of formula (41.3):

$$\begin{aligned}\prod_{i \in I} \delta_i \, 1_A(\omega) \,&=\, \delta_{i^*}\Big(\prod_{i \in J} \delta_i\Big)\, 1_A(\omega) \\ &=\, \delta_{i^*}\Big(\sum_{\eta \in \Omega(J)} (-1)^{|J| - o(\eta)} 1_A\big(\omega(J \leftarrow \eta)\big)\Big) \\ &=\, \sum_{\eta \in \Omega(J)} (-1)^{|J| - o(\eta)} \delta_{i^*}\Big(1_A\big(\omega(J \leftarrow \eta)\big)\Big)\,. \end{aligned} \tag{68.3}$$

For $K$ a subset of $J$, we denote by $\eta(K)$ the configuration of $\Omega(J)$ obtained by modifying all the components of $\omega$ in $J \setminus K$ while keeping the components of $\omega$ in $K$, i.e.,

$$\forall i \in J \qquad \eta(K)(i) \,=\, \begin{cases} \omega(i) & \text{if } i \in K\,, \\ 1 - \omega(i) & \text{if } i \in J \setminus K\,. \end{cases}$$

With this notation, we can rewrite the last sum of formula (68.3) as a sum over the subsets $K$ of $J$, and we obtain from (68.2) and (68.3) that

$$\mathcal{S}\big(1_A(\omega(I \leftarrow \cdot))\big) \,=\, \sum_{K \subset J} (-1)^{|J| - o(\eta(K))} \delta_{i^*}\Big(1_A\big(\omega(J \leftarrow \eta(K))\big)\Big)\,. \tag{68.4}$$

Let $K$ be a subset of $J$. The configuration $\eta(K)$ coincides with $\omega$ on $K$, therefore

$$1_A\big(\omega(J \leftarrow \eta(K))\big) \,=\, 1_A\big(\omega(J \setminus K \leftarrow \eta(K))\big)\,.$$

Suppose furthermore that $K$ is such that

$$\delta_{i^*}\Big(1_A\big(\omega(J \setminus K \leftarrow \eta(K))\big)\Big) \neq 0\,. \tag{68.5}$$

Recalling that $\{\, i^* \,\} \cup J \setminus K = I \setminus K$, we have

$$\omega_{I \setminus K} \,\leq\, \omega\big(J \setminus K \leftarrow \eta(K)\big)_{i^*} \,\leq\, \omega\big(J \setminus K \leftarrow \eta(K)\big)^{i^*} \,\leq\, \omega^{I \setminus K}\,. \tag{68.6}$$

Since $A$ is monotone, the inequalities (68.6) and the condition (68.5) imply that

$$1_A\big(\omega_{I \setminus K}\big) \,\neq\, 1_A\big(\omega^{I \setminus K}\big)\,.$$

Applying lemma 66.3, we see that

$$\text{pivotality}(I,\omega) \;\leq\; |I \setminus K|\,.$$

We conclude that the non-vanishing terms in the sum (68.4) correspond to subsets $K$ such that $|K| \leq k-\ell$, hence

$$\mathcal{S}\big(1_A(\omega(I\leftarrow\cdot))\big) \;=\; \sum_{\substack{K\subset J\\ |K|\leq k-\ell}} (-1)^{|J|-o(\eta(K))}\delta_{i^*}\Big(1_A\big(\omega(J\leftarrow\eta(K))\big)\Big)\,. \tag{68.7}$$

Taking the absolute values in (68.7), we obtain

$$\big|\mathcal{S}\big(1_A(\omega(I\leftarrow\cdot))\big)\big| \;\leq\; \sum_{0\leq h\leq k-\ell}\binom{k}{h}\,.$$

This proves the implication (68.1). □

# 69 A recursive inequality

From now onwards, we suppose that the event $A$ is decreasing. To alleviate the formulas, we introduce the notation

$$\mathcal{T}(I,\omega) \;=\; \overline{\mathcal{S}}\big(1_A(\omega(I\leftarrow\cdot))\big)1_A(\omega)\,.$$

We restart from formula (65.13), and we isolate the term corresponding to $\ell=k$:

$$E\Bigg(\sum_{\substack{I\in\mathcal{S}^*(k,\omega)\\ \text{pivotality}(I,\omega)=k}}\mathcal{T}(I,\omega)\Bigg) \;=\; E\big(L_k^{\mathcal{V}}1_A\big) \;-\; \sum_{1\leq\ell\leq k-1} E\Bigg(\sum_{\substack{I\in\mathcal{S}^*(k,\omega)\\ \text{pivotality}(I,\omega)=\ell}}\mathcal{T}(I,\omega)\Bigg)\,. \tag{69.1}$$

We deal first with the left-hand quantity, that we rewrite as follows:

$$E\Bigg(\sum_{\substack{I\in\mathcal{S}^*(k,\omega)\\ \text{pivotality}(I,\omega)=k}}\mathcal{T}(I,\omega)\Bigg) \;=\; \sum_{I\in\mathcal{I}(\mathcal{V},k)} E\Big(1_{\left\{\substack{I\in\mathcal{S}^*(k,\omega)\\ \text{pivotality}(I,\omega)=k}\right\}}\mathcal{T}(I,\omega)\Big)\,. \tag{69.2}$$

Let us fix $I$ in $\mathcal{I}(\mathcal{V},k)$. It follows from lemma 67.1 that

$$\omega\in A\,,\quad \text{pivotality}(I,\omega)=k \quad\Longrightarrow\quad \overline{\mathcal{S}}\big(1_A(\omega(I\leftarrow\cdot))\big) = -\frac{1}{1-p^k}\,.$$

Therefore

$$1_{\left\{\substack{I\in\mathcal{S}^*(k,\omega)\\ \text{pivotality}(I,\omega)=k}\right\}}\mathcal{T}(I,\omega) \;=\; -\frac{1}{1-p^k}\,1_{\{\,\text{pivotality}(I,\omega)=k\,\}}1_A(\omega)\,. \tag{69.3}$$

Substituting (69.3) into (69.2), and re-introducing the sum over $I$ inside the expectation, we get

$$E\Bigg(\sum_{\substack{I\in\mathcal{S}^*(k,\omega)\\ \text{pivotality}(I,\omega)=k}}\mathcal{T}(I,\omega)\Bigg) = \\ -\frac{1}{1-p^k}\,E\Big(\Big|\Big\{\,I\in\mathcal{I}(\mathcal{V},k):\text{pivotality}(I,\omega)=k\,\Big\}\Big|1_A(\omega)\Big)\,. \quad (69.4)$$

Substituting (69.4) into the identity (69.1), and taking the absolute value, we get

$$\frac{1}{1-p^k}E\Big(\Big|\Big\{\,I\in\mathcal{I}(\mathcal{V},k):\text{pivotality}(I,\omega)=k\,\Big\}\Big|1_A(\omega)\Big) \\ \leq \Big|E\big(L_k^{\mathcal{V}}1_A\big)\Big| + \sum_{1\leq\ell\leq k-1}\Bigg|E\Bigg(\sum_{\substack{I\in\mathcal{S}^*(k,\omega)\\ \text{pivotality}(I,\omega)=\ell}}\mathcal{T}(I,\omega)\Bigg)\Bigg| \\ \leq \Big|E\big(L_k^{\mathcal{V}}1_A\big)\Big| + \sum_{1\leq\ell\leq k-1}E\Bigg(\sum_{\substack{I\in\mathcal{S}^*(k,\omega)\\ \text{pivotality}(I,\omega)=\ell}}\big|\mathcal{T}(I,\omega)\big|\Bigg)\,. \quad (69.5)$$

Let us fix $I$ in $I(\mathcal{V},k)$. Let $f$ be a decreasing Boolean function on $\Omega(I)$ such that $\mathcal{S}(f)\neq 0$. Necessarily, the function $f$ is not constant and its support is not void. Moreover, we have

$$\big|\overline{\mathcal{S}}(f)\big| = \frac{1}{P\big(\text{supp}\,f\big)}\big|\mathcal{S}(f)\big|\,. \quad (69.6)$$

Let us fix a configuration $\omega$. For $I\in\mathcal{S}^*(k,\omega)$, the function $f=1_A(\omega(I\leftarrow\cdot))$ is decreasing on $\Omega(I)$ and has a non-vanishing signature, so we can apply the inequality (69.6) and we get

$$\big|\mathcal{T}(I,\omega)\big| = \big|\overline{\mathcal{S}}\big(1_A(\omega(I\leftarrow\cdot))\big)\big|1_A(\omega) \leq \frac{\big|\mathcal{S}\big(1_A(\omega(I\leftarrow\cdot))\big)\big|}{P\big(\text{supp}\,1_A(\omega(I\leftarrow\cdot))\big)}1_A(\omega)\,. \quad (69.7)$$

We control the denominator in the right-hand side of (69.7) with the help of (67.1). We use successively the inequalities (69.7) and (68.1) to control the sum inside the last expectation of (69.5), and we obtain

$$\sum_{\substack{I\in\mathcal{S}^*(k,\omega)\\ \text{pivotality}(I,\omega)=\ell}}\big|\mathcal{T}(I,\omega)\big| \leq \sum_{\substack{I\in\mathcal{S}^*(k,\omega)\\ \text{pivotality}(I,\omega)=\ell}}\frac{\big|\mathcal{S}\big(1_A(\omega(I\leftarrow\cdot))\big)\big|}{(1-p^\ell)(1-p)^{k-\ell}}1_A(\omega) \\ \leq \sum_{0\leq h\leq k-\ell}\binom{k}{h}\frac{\Big|\Big\{\,I\in\mathcal{S}^*(k,\omega):\text{pivotality}(I,\omega)=\ell\,\Big\}\Big|}{(1-p^\ell)(1-p)^{k-\ell}}1_A(\omega)\,. \quad (69.8)$$

We focus next on the last cardinality in (69.8). We have the obvious inequality

$$\Big|\Big\{ I\in\mathcal{S}^*(k,\omega) : \text{pivotality}(I,\omega)=\ell \Big\}\Big| \\ \leq \Big|\Big\{ I\in\mathcal{I}(\mathcal{V},k) : \text{pivotality}(I,\omega)=\ell \Big\}\Big|. \quad (69.9)$$

To bound the right-hand quantity in (69.9), we rely on the following lemma.

**Lemma 69.1.** *Let $A$ be a decreasing subset of $\{0,1\}^n$. For any configuration $\omega$ in $A$ and any $k,\ell$ such that $1\leq\ell\leq k$, we have*

$$\Big|\Big\{ I\in\mathcal{I}(\mathcal{V},k) : \textit{pivotality}(I,\omega)=\ell \Big\}\Big| \leq \\ \Big|\Big\{ J\in\mathcal{I}(\mathcal{V},\ell) : \textit{pivotality}(J,\omega)=\ell \Big\}\Big| M(n)\binom{V(n)}{k-\ell}. \quad (69.10)$$

*Proof.* Let $I$ be an element of $\mathcal{I}(\mathcal{V},k)$ such that $\text{pivotality}(I,\omega)=\ell$. By definition, there exist $V$ in $\mathcal{V}$ such that $I\subset V$ and a subset $J$ of $I$ of cardinality $\ell$ such that $J$ is in $\mathcal{S}^*(\ell,\omega)$. Since $\text{pivotality}(I,\omega)=\ell$, then $\text{pivotality}(J,\omega)=\ell$. Furthermore, the remaining elements of $I$ are included in $V$, i.e., $I\setminus J\subset V$. We deduce from these considerations that

$$\Big\{ I\in\mathcal{I}(\mathcal{V},k) : \text{pivotality}(I,\omega)=\ell \Big\} \subset \bigcup_{\substack{J\in\mathcal{S}^*(\ell,\omega)\\ \text{pivotality}(J,\omega)=\ell}} \bigcup_{\substack{V\in\mathcal{V}\\ J\subset V}} J\cup\binom{V\setminus J}{k-\ell}, \quad (69.11)$$

where $\binom{V\setminus J}{k-\ell}$ denotes the collection of the subsets of $V\setminus J$ having $k-\ell$ elements. This is a classical extension of the notation for the binomial coefficients, in particular, we have

$$\forall J\subset V \qquad \Big|\binom{V\setminus J}{k-\ell}\Big| \leq \binom{|V\setminus J|}{k-\ell} \leq \binom{V(n)}{k-\ell}. \quad (69.12)$$

The number of elements $V$ of $\mathcal{V}$ containing a fixed index, and a fortiori a fixed set $J$, is bounded from above by $M(n)$. The desired inequality (69.10) follows readily from the inequality (69.12) and the inclusion (69.11). □

Using the chain of inequalities (69.5), (69.8), (69.9), and (69.10) we obtain the much sought-after recursive inequality

$$E\bigg(\Big|\Big\{ I\in\mathcal{I}(\mathcal{V},k) : \text{pivotality}(I,\omega)=k \Big\}\Big|1_A(\omega)\bigg) \leq \Big|E\big(L_k^{\mathcal{V}}1_A\big)\Big| + \\ \sum_{\substack{1\leq\ell\leq k-1\\ 0\leq h\leq k-\ell}} \binom{k}{h} \frac{E\Big(\Big|\Big\{ J\in\mathcal{I}(\mathcal{V},\ell) : \text{pivotality}(J,\omega)=\ell \Big\}\Big|1_A(\omega)\Big)}{(1-p)^{k-\ell}} M(n)\binom{V(n)}{k-\ell}. \quad (69.13)$$

# 70 Induction

In order to alleviate the formulas, we first introduce the following notation:

$$\forall\omega\in\Omega\quad\forall k\in\{\,1,\dots,n\,\}\quad \mathcal{SP}(k,\omega)\,=\,\Big\{\,I\in\mathcal{I}(\mathcal{V},k):\text{pivotality}(I,\omega)=k\,\Big\}\,.$$

We will iterate inequality (69.13), or rather we will use it to prove by induction an interesting bound on $\big|\mathcal{SP}(k,\omega)\big|$. Before doing so, we modify it into a weaker inequality, but one that is better suited to this objective. Let us fix $k\geq 1$. Thanks to the bound (65.1) on $\big|E\big(L_k^{\mathcal{V}}1_A\big)\big|$, we have (recall that $M(n)\geq 1$ and $V(n)\geq 1$)

$$\Big|E\big(L_k^{\mathcal{V}}1_A\big)\Big|\,\leq\,\sqrt{\frac{nM(n)V(n)^k}{p^k(1-p)^k}}\,\leq\,\sqrt{n}\bigg(\frac{M(n)V(n)}{p(1-p)}\bigg)^k\,. \tag{70.1}$$

We have also the rough inequalities

$$\forall\ell\in\{\,1,\dots,k-1\,\}\qquad\sum_{0\leq h\leq k-\ell}\binom{k}{h}\binom{V(n)}{k-\ell}\,\leq\,k\big(kV(n)\big)^{k-\ell}\,. \tag{70.2}$$

Using the inequalities (70.1) and (70.2), we transform inequality (69.13) into

$$E\big(\big|\mathcal{SP}(k,\omega)\big|1_A(\omega)\big)\,\leq\,
\sqrt{n}\bigg(\frac{M(n)V(n)}{p(1-p)}\bigg)^k\,+\,k\,M(n)\sum_{1\leq\ell\leq k-1}E\big(\big|\mathcal{SP}(\ell,\omega)\big|1_A(\omega)\big)\Big(\frac{kV(n)}{1-p}\Big)^{k-\ell}\,. \tag{70.3}$$

We are now ready to prove the crucial inequality on $E\big(\big|\mathcal{SP}(k,\omega)\big|\big)$.

**Proposition 70.1.** *Let $A$ be a decreasing subset of $\{0,1\}^n$. We have*

$$\forall k\in\{\,1,\dots,n\,\}\qquad E\big(\big|\mathcal{SP}(k,\omega)\big|1_A(\omega)\big)\,\leq\,\sqrt{n}\bigg(\frac{2M(n)V(n)}{p(1-p)}\bigg)^k k^{3k}\,. \tag{70.4}$$

*Proof.* We proceed by induction on $k$. For $k=1$, we have $\mathcal{SP}(1,\omega)\,=\,\mathcal{P}(A,\omega)$. As $A$ is decreasing, we have

$$E\big(\big|\mathcal{P}(A,\omega)\big|1_A(\omega)\big)\,=\,(1-p)E\big(\big|\mathcal{P}(A,\omega)\big|\big)\,.$$

Applying the inequality of proposition 39.5, we obtain

$$E\big(\big|\mathcal{P}(A,\omega)\big|\big)\,\leq\,\sqrt{\frac{nP_p(A)}{p(1-p)}}\,\leq\,\sqrt{n}\frac{2M(n)V(n)}{p(1-p)}\,,$$

We conclude that

$$E\big(\big|\mathcal{SP}(1,\omega)1_A(\omega)\big|\big)\,\leq\,\sqrt{n}\frac{2M(n)V(n)}{p}\,.$$

Notice that the inequality (70.3) with $k = 1$ gives an even better inequality, with no factor 2. Let now $k \geq 2$ and suppose that we have proved the result until the rank $k-1$:

$$\forall \ell \in \{\,1,\dots,k-1\,\} \qquad E\big(\big|\mathcal{SP}(\ell,\omega)\big|1_A(\omega)\big) \,\leq\, \sqrt{n}\bigg(\frac{2M(n)V(n)}{p(1-p)}\bigg)^{\ell}\ell^{3\ell}\,.$$

Substituting all these inequalities into the inequality (70.3), we obtain

$$\begin{aligned} &E\big(\big|\mathcal{SP}(k,\omega)\big|1_A(\omega)\big) \,\leq \\ &\sqrt{n}\bigg(\frac{M(n)V(n)}{p(1-p)}\bigg)^{k} + k\,M(n)\sum_{1\leq\ell\leq k-1}\sqrt{n}\bigg(\frac{2M(n)V(n)}{p(1-p)}\bigg)^{\ell}\ell^{3\ell}\Big(\frac{kV(n)}{1-p}\Big)^{k-\ell} \\ &\qquad\qquad \leq\, \sqrt{n}\bigg(\frac{M(n)V(n)}{p(1-p)}\bigg)^{k}\bigg(1+2^{k-1}k\sum_{1\leq\ell\leq k-1}\ell^{3\ell}k^{k-\ell}\bigg)\,. \qquad (70.5) \end{aligned}$$

Noticing that

$$k\sum_{1\leq\ell\leq k-1}\ell^{3\ell}k^{k-\ell}\,\leq\, k\sum_{1\leq\ell\leq k-1}k^{3\ell}k^{k-\ell}\,\leq\, k^2k^{3k-2}\,=\,k^{3k}\,,$$

we bound the last parenthesis in (70.5) by $2^kk^{3k}$ and we conclude that

$$E\big(\big|\mathcal{SP}(k,\omega)\big|1_A(\omega)\big) \,\leq\, \sqrt{n}\bigg(\frac{2M(n)V(n)}{p(1-p)}\bigg)^{k}k^{3k}\,.$$

This conclude the induction step and the proof of proposition 70.1. □

Since $\mathcal{SP}(k,\omega)$ is a random subset of $\{\,1,\dots,n\,\}$ having cardinality $k$, we have the obvious inequality:

$$\forall k \in \{\,1,\dots,n\,\} \qquad \big|\mathcal{SP}(k,\omega)\big| \,\leq\, \binom{n}{k} \,\leq\, n^k\,.$$

Therefore the inequality (70.4) is completely useless whenever $k$ or the product $M(n)V(n)$ is of order $n$. Later on, we will apply the inequality with a covering $\mathcal{V}$ such that $M(n)V(n)$ is bounded by a power of $\ln n$ and $k$ is negligible compared to $\ln n$. The whole significance of inequality (70.4) will appear in this context.

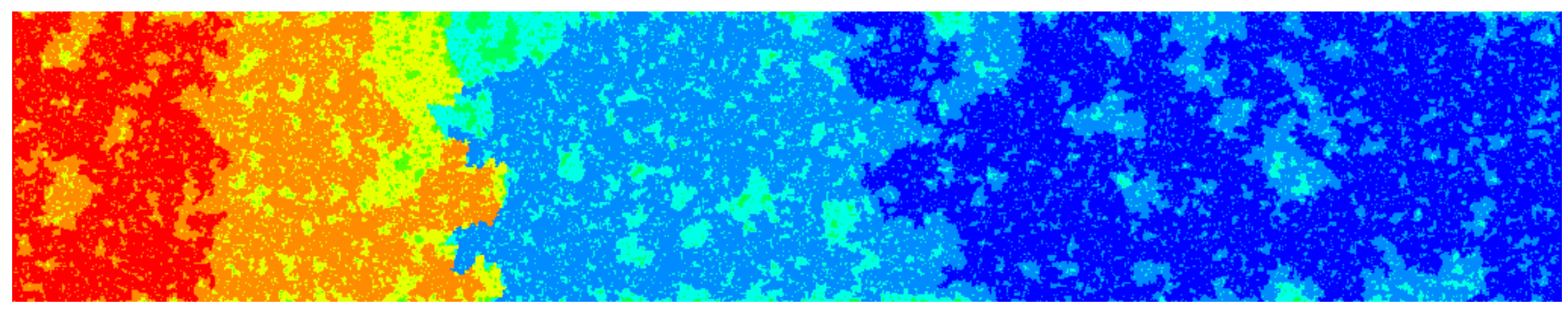

## 71 Control of the $\mathcal{V}$-local pivots

For the time being, we show how the control of the set $\mathcal{SP}(k,\omega)$ yields a control on the local pivots of order $k$. We will not derive a bound on the set $\mathcal{P}_k^{\mathcal{V}}(A,\omega)$ defined in (64.3), but rather on a specific subset of it, that we define next. The set $\mathcal{P}_k^{-,\mathcal{V}}(A,\omega)$ of the negative $\mathcal{V}$-local pivots of order $k$ for $A$ in $\omega$ is

$$\mathcal{P}_k^{-,\mathcal{V}}(A,\omega) = \Big\{ i\in\{1,\dots,n\} : \exists V\in\mathcal{V}(i) \quad \exists I\subset V\setminus\{i\} \quad |I|\leq k-1 \\ i\in\mathcal{P}^-\big(A,\omega^I\big)\,,\quad \forall j\in I\cup\{i\}\quad \omega(j)=0\Big\}\,.$$

The notation $\omega^I$ was defined at the beginning of section 66, the set $\mathcal{SP}(k,\omega)$ at the beginning of section 70, and the projection $\pi$ in (63.2).

**Lemma 71.1.** *Let $A$ be a decreasing subset of $\{0,1\}^n$. We have*

$$\forall\omega\in A\quad\forall k\in\{1,\dots,n\}\qquad \mathcal{P}_k^{-,\mathcal{V}}(A,\omega)\setminus\mathcal{P}_{k-1}^{-,\mathcal{V}}(A,\omega)\,\subset\,\pi\big(\mathcal{SP}(k,\omega)\big)\,.$$

*Proof.* Let $\omega\in A$, $k\in\{1,\dots,n\}$ and let $i$ belong to $\mathcal{P}_k^{-,\mathcal{V}}(A,\omega)\setminus\mathcal{P}_{k-1}^{-,\mathcal{V}}(A,\omega)$. It follows from the very definition of $\mathcal{P}_k^{-,\mathcal{V}}(A,\omega)$ that $\omega(i)=0$ and

$$\exists V\in\mathcal{V}(i)\quad\exists I\subset V\setminus\{i\}\quad |I|=k-1\quad i\in\mathcal{P}^-\big(A,\omega^I\big)\,,\quad\forall j\in I\quad\omega(j)=0\,.$$

We have $|I|=k-1$ because $i$ is not in $\mathcal{P}_{k-1}^{\mathcal{V}}(A,\omega)$. Let $J=I\cup\{i\}$. Since $V$ is in $\mathcal{V}(i)$, then $J$ is in $\mathcal{I}(\mathcal{V},k)$. The condition $i\in\mathcal{P}^-\big(A,\omega^I\big)$ means that

$$\omega^I_i\in A\,,\qquad (\omega^I)^i=\omega^J\notin A\,.$$

We compute next pivotality$(J,\omega)$. Notice that $\omega_J=\omega\in A$ while $\omega^J\notin A$. Lemma 66.3 tells us that there exists a subset $K$ of $J$ such that

$$\mathcal{S}\big(1_A(\omega(K\leftarrow\cdot))\big)\in\{-1,+1\}\,.$$

Therefore pivotality$(J,\omega)$ is not infinite. Suppose that pivotality$(J,\omega)=\ell<k$. By definition, there exists a subset $K$ of $J$ such that $\mathcal{S}\big(1_A(\omega(K\leftarrow\cdot))\big)\neq0$ and $|K|=\ell$. Applying lemma 66.2, we would have $\omega_K\in A$, $\omega^K\notin A$. As $\omega(i)=0$ and $K\setminus\{i\}\subset I$, we have $\omega^{K\setminus\{i\}}\leq\omega^I_i\in A$, therefore $\omega^{K\setminus\{i\}}$ is in $A$ and $\omega^{K\setminus\{i\}}\neq\omega^K$. This proves that $i$ must be present in $K$, and we have

$$i\in\mathcal{P}^-\big(A,\omega^{K\setminus\{i\}}\big)\,,\qquad |K\setminus\{i\}|\leq\ell-1\leq k-2\,.$$

This stands in contradiction with the hypothesis that $i\notin\mathcal{P}_{k-1}^{-,\mathcal{V}}(A,\omega)$. Therefore we have $\ell=k$ and pivotality$(J,\omega)=k$. This implies that $J\in\mathcal{SP}(k,\omega)$ and we can conclude that $i\in\pi\big(\mathcal{SP}(k,\omega)\big)$. ☐

Lemma 71.1 and proposition 70.1 together yield the following inequality on the negative $\mathcal{V}$-local pivots of order $k$.

**Corollary 71.2.** *Let $A$ be a decreasing subset of $\{0,1\}^n$. We have*

$$\forall\omega\in A\quad\forall k\in\{1,\dots,n\}\\ E\Big(\big|\mathcal{P}_k^{-,\mathcal{V}}(A,\omega)\setminus\mathcal{P}_{k-1}^{-,\mathcal{V}}(A,\omega)\big|1_A(\omega)\Big)\;\leq\;\sqrt{n}\bigg(\frac{2M(n)V(n)}{p(1-p)}\bigg)^k k^{3k+1}\,.\quad(71.1)$$

# 72 A general deviations inequality

A major goal was achieved when we obtained the inequality (71.1), which gives a bound on the expected number of the negative $\mathcal{V}$-local pivots of order $k$. This inequality can already be used to obtain a quantitative estimate on the probability of a decreasing event $A$. Let us explain briefly how. Suppose that we have a decreasing event $A$, and that the conditional expectation of the negative $\mathcal{V}$-local pivots of order $k$ is larger than $n^\alpha$, for some exponent $\alpha > 1$:

$$E\Big(\big|\mathcal{P}_k^{-,\mathcal{V}}(A) \setminus \mathcal{P}_{k-1}^{-,\mathcal{V}}(A)\big| \,\Big|\, A\Big) \,\geq\, n^\alpha\,. \tag{72.1}$$

Dividing the inequality (71.1) by $P(A)$, we have

$$E\Big(\big|\mathcal{P}_k^{-,\mathcal{V}}(A) \setminus \mathcal{P}_{k-1}^{-,\mathcal{V}}(A)\big| \,\Big|\, A\Big) \,\leq\, \sqrt{n}\bigg(\frac{2M(n)V(n)}{p(1-p)}\bigg)^k k^{3k+1}\frac{1}{P(A)}\,. \tag{72.2}$$

Combining together (72.1) and (72.2), we conclude that

$$P(A) \,\leq\, \frac{1}{n^{\alpha-1/2}}\bigg(\frac{2M(n)V(n)}{p(1-p)}\bigg)^k k^{3k+1}\,. \tag{72.3}$$

Of course, inequality (72.3) is useful only if $M(n), V(n)$ and $k$ are not too large. In our future applications, the quantities $M(n)$ and $V(n)$ will be bounded by a power of $\ln n$, and $k$ will be negligible compared to $(\ln n)^\beta$ for some $\beta < 1$. This will ensure that the power of $n$ in inequality (72.3) is not altered.

However, to implement our strategy in the percolation model, we will need a stronger exponential inequality controlling the upper deviations of the negative $\mathcal{V}$-local pivots of order $k$. The technique to derive this exponential inequality is the same as the one which was employed to derive the exponential inequality associated to the Margulis-Russo formula in subsection 40.1. However, as was to be expected, there are some complications. One of them is that the random variable at the heart of the estimates is not any more the simple i.i.d. sum $S_n(\omega)$, but the random homogeneous polynomial $L_k^{\mathcal{V}}(\omega)$. The key probabilistic input leading to the exponential inequality was the subgaussian behavior of $S_n(\omega)$. For $k \geq 2$, the random variable $L_k^{\mathcal{V}}(\omega)$ is not any more subgaussian. Fortunately, its tail still obeys a stretched exponential inequality. So, to warm-up before dealing with the local pivots, we prove here a general deviations inequality associated to a random variable whose tail has a stretched exponential behavior.

**Proposition 72.1.** *Let $a \geq e$ and let $b, \alpha$ be positive constants. Let $(L_n)_{n\geq 1}$ be a sequence of random variables which satisfies*

$$\forall n \geq 1 \qquad \forall u > 0 \qquad P\big(|L_n| \geq u\big) \,\leq\, a\exp\Big(-\frac{bu^{2\alpha}}{n^\alpha}\Big)\,. \tag{72.4}$$

*For any event $A \subset \{0,1\}^n$, we have*

$$P(A) \,\leq\, a\exp\bigg(-\frac{b}{n^\alpha}c(\alpha)\Big(\big|E\big(L_n \,\big|\, A\big)\big|\Big)^{2\alpha}\bigg)\,, \tag{72.5}$$

*where the constant $c(\alpha)$ is given by*

$$c(\alpha) \;=\; \begin{cases} \dfrac{1}{4^\alpha} & \text{if} \quad \alpha \geq 1/2\,, \\ \dfrac{1}{4}\Big(\dfrac{\alpha}{\Gamma\big(\frac{1}{2\alpha}\big)}\Big)^{2\alpha} & \text{if} \quad \alpha < 1/2\,, \end{cases}$$

*and $\Gamma$ is the classical function defined by*

$$\forall x > 0 \qquad \Gamma(x) \;=\; \int_0^{+\infty} t^{x-1} e^{-t}\, dt\,.$$

*The constant $c(\alpha)$ satisfies*

$$\forall \alpha \in ]0,1] \qquad c(\alpha) \;\geq\; \frac{1}{4}\alpha^{2\alpha+1}\,. \tag{72.6}$$

*Proof.* The proof is an adaptation of the proof of proposition 40.2. We start with

$$E\big(|L_n|\,\big|\,A\big) \;=\; \frac{1}{P(A)}\int_A |L_n|\, dP\,. \tag{72.7}$$

Using Fubini's theorem, we rewrite the integral in (72.7) as follows:

$$\int_A |L_n|\, dP \;=\; \int_A \Big(\int_0^{+\infty} 1_{u\leq |L_n|}\, du\Big)\, dP \;=\; \int_0^{+\infty} P\big(A, |L_n| \geq u\big)\, du\,.$$

Let $t > 0$ and let us split this integral in two:

$$\int_0^{+\infty} P\big(A, |L_n| \geq u\big)\, du \;=\; \int_0^t \cdots + \int_t^{+\infty} \cdots \;\leq\; tP(A) + \int_t^{+\infty} P\big(|L_n| \geq u\big)\, du\,. \tag{72.8}$$

Reporting (72.8) in (72.7), we have

$$E\big(|L_n|\,\big|\,A\big) \;\leq\; \frac{1}{P(A)}\int_t^{+\infty} P\big(|L_n| \geq u\big)\, du \,+\, t\,. \tag{72.9}$$

Reporting (72.4) into (72.9), we obtain

$$E\big(|L_n|\,\big|\,A\big) \;\leq\; \frac{a}{P(A)}\int_t^{+\infty} \exp\Big(-\frac{bu^{2\alpha}}{n^\alpha}\Big)\, du \,+\, t\,. \tag{72.10}$$

We make the change of variable $cu = \sqrt{n}v$, with $c = b^{1/(2\alpha)}$, so that the previous inequality (72.10) becomes

$$E\big(|L_n|\,\big|\,A\big) \;\leq\; \frac{a\sqrt{n}}{cP(A)}\int_{ct/\sqrt{n}}^{+\infty} \exp\big(-v^{2\alpha}\big)\, dv \,+\, t\,.$$

This inequality holds for any $t > 0$, so we replace $t$ by $\frac{\sqrt{n}t}{c}$, and we get

$$E\big(|L_n|\,\big|\,A\big) \;\leq\; \frac{\sqrt{n}}{c}\bigg(\frac{a}{P(A)}\int_t^{+\infty} \exp\big(-v^{2\alpha}\big)\, dv \,+\, t\bigg)\,. \tag{72.11}$$

The right-hand quantity admits a unique global minimum at

$$t^* \;=\; \Big(\ln \frac{a}{P(A)}\Big)^{1/(2\alpha)} \;\geq\; 1\,.$$

We try next to compute a simple upper bound of the right-hand side of (72.11) by looking either at $t = t^*$ or at a value close to it. We distinguish several cases, according to the value of $\alpha$:

• $\alpha = 1/2$. In this case, we can compute the integral, and we get

$$\frac{a}{P(A)} \int_{t^*}^{+\infty} \exp\big(-v\big)\, dv \;=\; \frac{a}{P(A)} \exp\big(-t^*\big) \;=\; 1 \;\leq\; t^*\,. \tag{72.12}$$

Substituting (72.12) into (72.11), we get

$$E\big(|L_n|\,\big|\,A\big) \;\leq\; \frac{\sqrt{n}}{c} 2t^* \;=\; 2\frac{\sqrt{n}}{c} \ln \frac{a}{P(A)}\,.$$

• $\alpha > 1/2$. We perform an integration by parts to bound the integral as follows: for any $t > 0$,

$$\begin{aligned}
\int_t^{+\infty} \exp\big(-v^{2\alpha}\big)\, dv \;&=\; \int_t^{+\infty} \frac{1}{2\alpha v^{2\alpha-1}} 2\alpha v^{2\alpha-1} \exp\big(-v^{2\alpha}\big)\, dv \\
&= \Big[-\frac{1}{2\alpha v^{2\alpha-1}} \exp\big(-v^{2\alpha}\big)\Big]_t^{+\infty} - \int_t^{+\infty} \frac{2\alpha-1}{2\alpha v^{2\alpha}} \exp\big(-v^{2\alpha}\big)\, dv \\
&= \frac{1}{2\alpha t^{2\alpha-1}} \exp\big(-t^{2\alpha}\big) - \frac{2\alpha-1}{2\alpha} \int_t^{+\infty} \frac{1}{v^{2\alpha}} \exp\big(-v^{2\alpha}\big)\, dv\,.
\end{aligned} \tag{72.13}$$

Taking advantage of the fact that $2\alpha - 1 > 0$, we deduce from (72.13) that

$$\int_t^{+\infty} \exp\big(-v^{2\alpha}\big)\, dv \;\leq\; \frac{1}{2\alpha t^{2\alpha-1}} \exp\big(-t^{2\alpha}\big)\,. \tag{72.14}$$

Using (72.14) at $t = t^*$, we have

$$\frac{a}{P(A)} \int_{t^*}^{+\infty} \exp\big(-v^{2\alpha}\big)\, dv \;\leq\; \frac{1}{2\alpha (t^*)^{2\alpha-1}} \;\leq\; t^*\,. \tag{72.15}$$

Substituting (72.15) into (72.11), we get

$$E\big(|L_n|\,\big|\,A\big) \;\leq\; \frac{\sqrt{n}}{c} 2t^* \;=\; 2\sqrt{n}\Big(\frac{1}{b}\ln \frac{a}{P(A)}\Big)^{1/(2\alpha)}\,. \tag{72.16}$$

• $0 < \alpha < 1/2$. This case is more complicated, because the integration by parts does not lead to a simple inequality. Indeed, the last term in (72.13) is positive. So we proceed differently. Let us fix $t > 0$. We make the change of variable $u = v^{2\alpha}$ in the integral. We have

$$\frac{du}{u} \;=\; 2\alpha \frac{dv}{v}\,, \qquad v = u^{1/(2\alpha)}\,,$$

whence

$$\int_t^{+\infty} e^{-v^{2\alpha}}\,dv \;=\; \int_{t^{2\alpha}}^{+\infty} \frac{1}{2\alpha} u^{1/(2\alpha)-1} e^{-u}\,du\,. \tag{72.17}$$

We bound next the integral with the help of the $\Gamma$ function as follows:

$$\begin{aligned}
\int_{t^{2\alpha}}^{+\infty} \frac{1}{2\alpha} u^{1/(2\alpha)-1} e^{-u}\,du \;&\le\; e^{-t^{2\alpha}/2} \int_{t^{2\alpha}}^{+\infty} \frac{1}{2\alpha} u^{1/(2\alpha)-1} e^{-u/2}\,du \\
&\le\; e^{-t^{2\alpha}/2} \frac{1}{2\alpha} \int_0^{+\infty} u^{1/(2\alpha)-1} e^{-u/2}\,du \\
&=\; e^{-t^{2\alpha}/2} \frac{2^{1/(2\alpha)}}{2\alpha} \int_0^{+\infty} w^{1/(2\alpha)-1} e^{-w}\,dw \\
&=\; e^{-t^{2\alpha}/2} \frac{2^{1/(2\alpha)}}{2\alpha} \Gamma\Big(\frac{1}{2\alpha}\Big)\,.
\end{aligned} \tag{72.18}$$

Substituting (72.17) and (72.18) into (72.11), we get

$$E\big(|L_n| \,\big|\, A\big) \;\le\; \frac{\sqrt{n}}{c}\bigg(\frac{a}{P(A)} e^{-t^{2\alpha}/2} d(\alpha) \,+\, t\bigg)\,,$$

where $d(\alpha) \ge 1$ is given by

$$d(\alpha) \;=\; \frac{2^{1/(2\alpha)}}{2\alpha}\Gamma\Big(\frac{1}{2\alpha}\Big) \;=\; 2^{1/(2\alpha)}\Gamma\Big(1+\frac{1}{2\alpha}\Big) \;\ge\; 1\,.$$

Taking finally $t = 2^{1/(2\alpha)} t^*$, we conclude that

$$\begin{aligned}
E\big(|L_n| \,\big|\, A\big) \;&\le\; \frac{\sqrt{n}}{c}\bigg(d(\alpha) \,+\, \Big(2\ln\frac{a}{P(A)}\Big)^{1/(2\alpha)}\bigg) \\
&\le\; \frac{\sqrt{n}}{c} 2d(\alpha)\Big(2\ln\frac{a}{P(A)}\Big)^{1/(2\alpha)} \le\; \frac{\sqrt{n}}{\alpha}\,\Gamma\Big(\frac{1}{2\alpha}\Big)\Big(\frac{4}{b}\ln\frac{a}{P(A)}\Big)^{1/(2\alpha)}.
\end{aligned} \tag{72.19}$$

In each case $\alpha \ge 1/2$ and $\alpha < 1/2$, we use finally the fact that

$$\Big|E\big(L_n \,\big|\, A\big)\Big| \;\le\; E\big(|L_n| \,\big|\, A\big)\,,$$

and we deduce the inequality (72.5) from the inequalities (72.16) and (72.19). It remains to prove the lower bound (72.6). To that end, we will use a simple upper bound on the function $\Gamma$. Let $x \ge 1$. The functional equation satisfied by $\Gamma$ implies that

$$\Gamma(x) \;=\; (x-1)\cdots(x-\lfloor x\rfloor+1)\,\Gamma(x-\lfloor x\rfloor+1)\,. \tag{72.20}$$

Now $x-\lfloor x\rfloor+1$ is in the interval $[1,2]$. Moreover $\Gamma(1)=\Gamma(2)=1$, and since the function $\Gamma$ is logarithmically convex on $]0,+\infty[$, we have $\Gamma(x-\lfloor x\rfloor+1)\le 1$. It follows then from (72.20) that $\Gamma(x) \;\le\; x^{\lfloor x\rfloor-1} \;\le\; x^x$. Substituting this inequality in the expression of $c(\alpha)$, we conclude that, for $\alpha < 1/2$, we have

$$c(\alpha) \;=\; \frac{1}{4}\bigg(\frac{\alpha}{\Gamma\big(\frac{1}{2\alpha}\big)}\bigg)^{2\alpha} \;\ge\; \frac{1}{4}\alpha^{2\alpha+1}\,.$$

This inequality holds also for $\alpha$ in $[1/2,1]$. $\square$

# 73 Application to $L_k^{\mathcal{V}}$

In this section, we will apply the general deviations inequality of proposition 72.1 to the random variable $L_k^{\mathcal{V}}$. To that end, we need first a control on the tail behavior of $L_k^{\mathcal{V}}$. The adequate inequality is stated in the next proposition, it comes from a general inequality on random polynomials, which can be proved via hypercontractivity.

**Proposition 73.1.** *Let $k \in \{1, \dots, n\}$. The random variable $L_k^{\mathcal{V}}$, defined by*

$$L_k^{\mathcal{V}}(\omega) \;=\; \sum_{I \in \mathcal{I}(\mathcal{V},k)} \prod_{i \in I} X_i(\omega)\,,$$

*obeys the following deviations inequality:*

$$\forall t > 0 \qquad P_p\big(|L_k^{\mathcal{V}}| \geq t\big) \;\leq\; \exp\left(k - \frac{k\, t^{2/k}}{C(p)\, M(n) V(n)\, n^{1/k}}\right), \tag{73.1}$$

*where*

$$C(p) \;=\; 2B(p)e \;=\; \frac{e(1-2p)}{p(1-p)\ln(1/p-1)}\,. \tag{73.2}$$

*Proof.* The function $L_k^{\mathcal{V}} : \{0,1\}^n \to \mathbb{R}$ is a homogeneous polynomial of order $k$. Applying theorem 46.3, we obtain

$$\forall t > 0 \qquad P_p\Big(|L_k^{\mathcal{V}}| \geq t\sqrt{E\big(\big(L_k^{\mathcal{V}}\big)^2\big)}\Big) \;\leq\; \exp\Big(k - \frac{k}{2B(p)e} t^{2/k}\Big)\,.$$

We make the change of variable $t \leftrightarrow t/\sqrt{E\big(\big(L_k^{\mathcal{V}}\big)^2\big)}$ and we obtain

$$\forall t > 0 \qquad P_p\big(|L_k^{\mathcal{V}}| \geq t\big) \;\leq\; \exp\left(k - \frac{kt^{2/k}}{2B(p)e\big(E\big(\big(L_k^{\mathcal{V}}\big)^2\big)\big)^{1/k}}\right).$$

We use finally the upper bound on $E\big(\big(L_k^{\mathcal{V}}\big)^2\big)$ computed in (65.1), and we get the desired inequality (73.1). $\square$

It remains only to make the synthesis between the two propositions 72.1 and 73.1 to obtain the following corollary.

**Corollary 73.2.** *Let $k \in \{1, \dots, n\}$. For any event $A \subset \Omega$, we have*

$$P(A) \;\leq\; \exp\left(k - \frac{\Big(\big|E\big(L_k^{\mathcal{V}} \,|\, A\big)\big|\Big)^{2/k}}{4k^{2/k} C(p)\, M(n) V(n) n^{1/k}}\right),$$

*where the constant $C(p)$ was defined in* (73.2).

*Proof.* In view of the inequality of proposition 73.1 on the tail behavior of $L_k^{\mathcal{V}}$, we can apply the general deviations inequality of proposition 72.1 with

$$a = e^k\,, \quad b = \frac{k}{C(p)\, M(n) V(n)}\,, \quad \alpha = \frac{1}{k}\,,$$

and we use also the inequality (72.6) to obtain the desired result. $\square$

## 74 The exponential inequalities

We divide the recursive inequality (69.13) by $P(A)$ and we obtain the following conditional version:

$$E\Big(\big|\{\, I\in\mathcal{I}(\mathcal{V},k) : \text{pivotality}(I,\omega)=k\,\}\big|\,\Big|\,A\Big)\,\leq\,\big|E\big(L_k^{\mathcal{V}}\,|\,A\big)\big|+$$

$$\sum_{\substack{1\leq\ell\leq k-1\\0\leq h\leq k-\ell}}\binom{k}{h}\frac{E\Big(\big|\{\, J\in\mathcal{I}(\mathcal{V},\ell) : \text{pivotality}(J,\omega)=\ell\,\}\big|\,\Big|\,A\Big)}{(1-p)^{k-\ell}}\,M(n)\binom{V(n)}{k-\ell}.$$

To alleviate the formulas, we introduce the notation, for $\ell\in\{\,1,\dots,k\,\}$,

$$E^{\mathcal{V}}(\ell,A)\,=\,E\Big(\big|\{\, J\in\mathcal{I}(\mathcal{V},\ell) : \text{pivotality}(J,\omega)=\ell\,\}\big|\,\Big|\,A\Big).\tag{74.1}$$

We use the rough inequalities (70.2) to obtain the weaker inequality

$$E^{\mathcal{V}}(k,A)\,\leq\,\big|E\big(L_k^{\mathcal{V}}\,|\,A\big)\big|+\sum_{1\leq\ell\leq k-1}E^{\mathcal{V}}(\ell,A)\,M(n)\,k\Big(\frac{kV(n)}{1-p}\Big)^{k-\ell}.$$

We rewrite this inequality in order to get a lower bound on $\big|E\big(L_k^{\mathcal{V}}\,|\,A\big)\big|$:

$$\big|E\big(L_k^{\mathcal{V}}\,|\,A\big)\big|\,\geq\,\Big(E^{\mathcal{V}}(k,A)-\sum_{1\leq\ell\leq k-1}E^{\mathcal{V}}(\ell,A)\,M(n)\,k\Big(\frac{kV(n)}{1-p}\Big)^{k-\ell}\Big)^+,\tag{74.2}$$

where $(x)^+=\max(x,0)$ is the positive part of $x$. Recall that the recursive inequality was established only for a decreasing event. Substituting inequality (74.2) into the inequality of corollary 73.2, we obtain the following result.

**Proposition 74.1.** *Let $A$ be a decreasing subset of $\{0,1\}^n$. We have*

$$\forall k\in\{\,1,\dots,n\,\}\tag{74.3}$$

$$P(A)\,\leq\,\exp\left(k-\frac{\left(\Big(E^{\mathcal{V}}(k,A)-\sum_{1\leq\ell\leq k-1}E^{\mathcal{V}}(\ell,A)\,M(n)\,k\Big(\frac{kV(n)}{1-p}\Big)^{k-\ell}\Big)^+\right)^{2/k}}{4k^{2/k}C(p)\,M(n)\frac{V(n)}{1-p}n^{1/k}}\right)$$

*where $E^{\mathcal{V}}(\ell,A)$, $1\leq\ell\leq k$, are the conditional expectations defined in* (74.1), *and the constant $C(p)$ was defined in* (73.2).

For $k=1$, we have $E^{\mathcal{V}}(1,A)=E\big(\mathcal{P}^-(A)\,|\,A\big)$ and inequality (74.3) reads

$$P(A)\,\leq\,\exp\left(1-\frac{(1-p)\Big(E\big(\mathcal{P}^-(A)\,|\,A\big)\Big)^2}{4C(p)\,M(n)V(n)n}\right),$$

which is essentially the exponential inequality (40.5) we already knew for monotone events. There is an apparent deterioration because of the presence of the quantities $M(n),V(n)$; however, these can be taken equal to 1, because the negative pivots of order 1 can be detected without exploring any neighborhood. Needless to say, these inequalities look quite awkward for $k\geq 2$. We will see very soon how they can be applied to our percolation problems.

Part X

# Partial answers to the challenges

So, where do we stand now? To get interesting answers to the three challenges raised in section 30, we need to obtain quantitative estimates on the number of pivots of order $k$, and this for any fixed value $k \geq 2$. This is the main obstacle we have been facing for more than fifteen years. In part V, we tried to work out this question for a general monotone set, but this attempt was a resounding failure. Indeed, the example of Talagrand, developed in section 62, shows that the number of pivots of any order $k \geq 2$ of a monotone subset of $\{0,1\}^n$ can reach the order $n$. However, in our first brutal attempt of control in section 31, we could show that, for a disconnection event in percolation, the number of higher intersection sets satisfies the following bound:

$$\forall k \geq 1 \qquad |\mathcal{I}(k)| \;\leq\; c(k,p,d)\big(n^d \ln n\big)^{1-2^{-k}}\,, \tag{VI.4}$$

where $c(k)$ is a constant depending on $k, p, d$ (see (31.80) in subsection 31.6). The intersection set of order $k$ is a subset of the pivots of order $k$, yet we believe that the argument used to prove (VI.4) could be modified to obtain a similar bound on the whole set of the closed pivots of order $k$, but still for the disconnection event. We can draw two main conclusions from this bound. First, the example of Talagrand falls into a different category than the disconnection event in percolation. Second, we must absolutely take into account some specificities of the disconnection event under consideration, beyond the mere fact that it is decreasing. The main point that led to the inequality (VI.4) is the existence of a special exploration algorithm for finding the higher intersections. Moreover, it seems that the inequalities (VI.4) could be considerably improved if we could develop a better probabilistic control of the algorithm. We have been trying to do so throughout sections 33 and 45, focusing only on the second pivots, but without success. Due to our inability to adequately control the conditional distribution of the percolation configuration knowing the position of the first pivots, there is a missing factor in the estimates. This leads to an inevitable deterioration of the result. This is the crux of the matter. If only we could derive an estimate which would be in accordance with the numerical simulations, then the way would be open to resolving the conjecture $\theta(p_c, \mathbb{Z}^d) = 0$.

In order to move forward and to go around this obstacle, we work with the additional hypothesis 75.1 of stretched exponential decay. With the help of this hypothesis, we will rule out the existence of pivots of order $k$ which are far away from each other. It remains then to control the local pivots of order $k$. This is achieved with the results on the $\mathcal{V}$-local pivots developed in part VII. Although the results obtained in this way are not entirely satisfactory, they contribute to exclude certain specific scenarios.

# 75 The strategy

Throughout part VIII, we shall assume the following hypothesis.

**Hypothesis 75.1** (Stretched exponential decay)**.** *We suppose that $p$ is such that*

$$\liminf_{n\to\infty}\,\frac{1}{\ln n}\ln\ln\Big(\frac{1}{P_p\big(n\leq|C(0)|<+\infty\big)}\Big)\,>\,0\,.$$

Hypothesis 75.1 implies the following estimate on the tail of the finite clusters: there exist two positive constants $a,\alpha$ such that

$$\forall n\geq 1\qquad P_p\big(n\leq|C(0)|<+\infty\big)\,\leq\,a\exp(-n^\alpha)\,.\tag{75.1}$$

In each of the three challenging questions raised in section 30, we considered a disconnection event $\mathcal{D}$ occurring in a finite domain $\Lambda$. There were slight variations in the notation, because these events depended on additional parameters. Let us recapitulate:

• Disconnection of a sub-box (subsection 30.1): The domain $\Lambda$ is a cubic box. For $\Gamma$ a cubic sub-box of $\Lambda$, we consider the disconnection event

$$\mathcal{D}(\Gamma,\Lambda)\,=\,\big\{\,\Gamma\not\longleftrightarrow\partial^{\,in}\Lambda\,\big\}\,,$$

• The largest cluster in a box (subsection 30.2): The domain $\Lambda$ is a cubic box and the event is

$$\mathcal{D}(x,m)\,=\,\big\{\,x\not\longleftrightarrow\partial^{\,in}\Lambda(n)\,,|C(x)|\geq m^d\,\big\}\,.$$

• Disconnection in a cuboid (subsection 30.3): The domain is a prism $\Pi(n)$. Denoting by $L$ and $R$ the middle hypersquares of the two opposite largest faces of $\Pi(n)$, the event is

$$\mathcal{D}(n)\,=\,\big\{\,L\not\longleftrightarrow R\text{ in }\Pi(n)\,\big\}\,.$$

We explain next how we will use the stretched exponential hypothesis to reach interesting conclusions for these three events. The strategy is the same in each case. From now onwards, we denote by $\mathcal{D}$ one of these three events, and by $\Lambda$ the associated domain. Let $M\geq 1$ be an additional integer parameter. Let $\mathcal{E}$ be the event that the percolation configuration restricted to $\Lambda$ does not contain a finite cluster of cardinality strictly larger than $M$. Recall that a finite cluster of $\Lambda$ is a cluster which does not contain a vertex belonging to the inner vertex boundary of $\Lambda$, hence

$$\mathcal{E}\,=\,\Big\{\forall x\in\Lambda\qquad x\not\longleftrightarrow\partial^{\,in}\Lambda\quad\Longrightarrow\quad|C(x)|\leq M\Big\}\,.$$

An easy upper bound on $\mathcal{E}^c$ is given by the union bound:

$$\begin{aligned}P(\mathcal{E}^c)\,&\leq\,\sum_{x\in\Lambda}P\big(x\not\longleftrightarrow\partial^{\,in}\Lambda,\,|C(x)|>M\big)\\&\leq\,\sum_{x\in\Lambda}P\big(M<|C(x)|<+\infty\big)\,=\,|\Lambda|P\big(M<|C(0)|<+\infty\big)\,.\end{aligned}\tag{75.2}$$

Inequality (75.2), together with the hypothesis (75.1), yield that

$$P(\mathcal{E}^c)\,\leq\,|\Lambda|\,a\exp(-M^\alpha)\,.$$

Let $r\geq 1$. We choose

$$M\,=\,\Big((r+1)\ln|\Lambda|\Big)^{1/\alpha}\,,$$

and we obtain that

$$P(\mathcal{E}^c)\,=\,O\Big(\frac{1}{|\Lambda|^r}\Big)\,=\,O\Big(\frac{1}{n^{dr}}\Big)\,.$$

This estimate will allow us to get rid of the configurations containing finite clusters of intermediate sizes. So we bound $P(\mathcal{D})$ as follows:

$$P(\mathcal{D})\,\leq\,P(\mathcal{D}\cap\mathcal{E})\,+\,P(\mathcal{E}^c)\,\leq\,P(\mathcal{D}\cap\mathcal{E})\,+\,O\Big(\frac{1}{n^{dr}}\Big)\,.\tag{75.3}$$

We use then the method employed in subsection 30.1 to bound $P(\mathcal{D}\cap\mathcal{E})$, namely, we define a Min-cut problem, we consider the associated intertwined explorations and their final iteration number. The only difference is that we deal all the time with the event $\mathcal{D}\cap\mathcal{E}$ instead of $\mathcal{D}$. The next step should in principle be adapted to each specific disconnection event. However, to ease the reading, we will adopt a slightly relaxed presentation, in which we do not care for the constants appearing in the inequalities involving the intersection sets $\mathcal{I}(t)$. We present only the general scheme, which works for the three challenges. The first two challenges depend on a parameter $m\geq 1$ (which can be taken equal to $n$ when dealing with the third challenge). The three challenges occur in a domain whose side length is typically $n$. This domain is either $\Lambda(2n)$, or $\Lambda(n)$ or the prism $\Pi(n)$. In the sequel, we denote it by $\Lambda$. The only thing that matters is that its cardinality scales as $n^d$ as $n$ goes to $\infty$. Proceeding exactly as in subsection 30.1, we will obtain an inequality of the form

$$P(\mathcal{D}\cap\mathcal{E})\,\leq\,(\ln n)\max_{1\,\leq\,t\,\leq\,\frac{\ln n}{\beta(n)}}P\Big(\,|\mathcal{I}(t)|\,\geq\,\frac{m^{d-1}\beta(n)}{\ln n},\,\mathcal{D}\cap\mathcal{E}\,\Big)+n^{d-\beta(n)}\,,\tag{75.4}$$

where $\beta(n)$ is a function such that $\lim_{n\to+\infty}\beta(n)=+\infty$. The intersection set $\mathcal{I}(t)$ is in fact a subset of the negative pivots of order $t$ for the disconnection event $\mathcal{D}$. Thus everything boils down to estimating

$$\max_{1\,\leq\,t\,\leq\,\frac{\ln n}{\beta(n)}}P\Big(\,\big|\mathcal{P}_t^-(\mathcal{D})\big|\,\geq\,\frac{m^{d-1}\beta(n)}{\ln n},\,\mathcal{D}\cap\mathcal{E}\,\Big)\,,$$

but this seems still to be very difficult. In fact, we do not have at our disposal a general method for estimating the tail of the number of pivots. However, we have a deviations inequality which provides a control on the probability of the event once we have an estimate on the conditional expectation of the number of pivots. So our strategy is to obtain a lower bound on this last quantity.

# 76 Locality of the pivots

We will make appeal to the quantitative estimates on the $\mathcal{V}$-local pivots developed in part VII. For the covering $\mathcal{V}$, we take the collection of the cubic sub-boxes of $\Lambda$ having side length equal to

$$L\,=\,\left\lceil 3M\frac{\ln n}{\beta(n)}\right\rceil\,=\,\left\lceil 3\Big(2\ln|\Lambda|\Big)^{1/\alpha}\frac{\ln n}{\beta(n)}\right\rceil.$$

More precisely, the covering $\mathcal{V}$ is defined as

$$\mathcal{V}\,=\,\Big\{\,V\subset\Lambda:\exists\,x\in\Lambda\quad V=x+\Big]-\frac{L}{2},\frac{L}{2}\Big]^d\,\Big\}.$$

With this choice for the covering, the associated quantities $M(n)$ and $V(n)$ defined in (64.1) and (64.2) are both bounded from above by $L^d$, and there exist positive constants $c,b$ such that

$$\max\big(M(n),V(n)\big)\,\leq\,c(\ln n)^b\,. \tag{76.1}$$

Let us fix $t$ such that $2\leq t\leq(\ln n)/\beta(n)$. It was proved in lemma 29.5 that $\mathcal{I}(t)$ is a subset of $\mathcal{P}_t^-(\mathcal{D})$. We claim that, when the event $\mathcal{D}\cap\mathcal{E}$ occurs, the set $\mathcal{I}(t)$ is in fact a subset of $\mathcal{P}_t^{-,\mathcal{V}}(\mathcal{D})$. Let us re-examine quickly the proof of lemma 29.5. We consider only the case of the even indices, so we suppose that $t=2k$. Let $\omega$ be a configuration in $\mathcal{D}\cap\mathcal{E}$, and let us consider a bond $e$ belonging to $\mathcal{I}(2k)$. By the definition of $\mathcal{I}(2k)$, we can write $e=\langle x,x'\rangle$, where $x$ is in $\mathcal{A}(2k)$ and $x'$ in $\mathcal{A}(2k-1)$. In the proof of lemma 29.5, it was shown that there exist a path $x_0,e_0,x_1,e_1,\cdots,e_{\ell-1},x_\ell$ in $\mathcal{A}(2k)$ joining a vertex $x_0$ of $W$ to $x_\ell=x$ such that at most $k$ bonds among $e_0,\cdots,e_{\ell-1}$ are closed, and another path $x'_0,e'_0,x'_1,e'_1,\cdots,e'_{\ell'-1},x'_{\ell'}$ in $\mathcal{A}(2k-1)$ joining a vertex $x'_0$ of $W'$ to $x'_{\ell'}=x'$ such that at most $k-1$ bonds among $e'_0,\cdots,e'_{\ell'-1}$ are closed (the sets $W$ and $W'$ are the sets which are disconnected when $\mathcal{D}$ occurs). We denote by $E$ (respectively $E'$) the set of the closed bonds belonging to the first path (respectively the second path). We have thus $|E|\leq k$, $|E'|\leq k-1$, and moreover, upon opening the bonds of $E\cup E'$, we get a configuration belonging to $\mathcal{D}$ in which $e$ is pivotal. This shows that $e$ is in $\mathcal{P}_{2k}^-(\mathcal{D},\omega)$. To prove that $e$ is in $\mathcal{P}_{2k}^{-,\mathcal{V}}(\mathcal{D},\omega)$, we need to check in addition that the set of bonds $E\cup E'\cup\{\,e\,\}$ is included in a box belonging to the covering $\mathcal{V}$. For simplicity, we consider only the first challenge, the disconnection of a sub-box $W'=\Gamma$ from the boundary of $W=\Lambda(2n)$; for the other two challenges, small adaptations of the next argument are needed. By construction, the second path meets only finite clusters of $\Lambda(2n)$, therefore the number of vertices visited between two consecutive closed bonds of the path is less than or equal to $M$. For the first path, if we define

$$h\,=\,\max\big\{\,i\in\{\,0,\dots,\ell\,\}:x_i\longleftrightarrow W\,\big\}\,,$$

then the same reasoning applies for the subpath $x_h,e_h,x_{h+1},e_{h+1},\cdots,e_{\ell-1},x_\ell$. The definition of $h$ ensures also that all the bonds $e_0,\dots,e_{h-1}$ are open, thus the set $E$ is necessarily included in the bonds of the former sub-path. Therefore, we can enumerate the bonds of $E\cup E'\cup\{\,e\,\}$ in a sequence such that the distance between two consecutive bonds is less than or equal to $M$. This implies that $E\cup E'\cup\{\,e\,\}$ is included in a box $V$ belonging to $\mathcal{V}$ and we conclude that $e$ is indeed in $\mathcal{P}_{2k}^{\mathcal{V}-}(\mathcal{D},\omega)$. So, we have the inclusion

$$\forall\omega\in\mathcal{D}\cap\mathcal{E}\qquad \mathcal{I}(t,\omega)\,\subset\,\mathcal{P}_t^{-,\mathcal{V}}(\mathcal{D},\omega)\,. \tag{76.2}$$

Furthermore, the event $\mathcal{E}$ is decreasing and

$$\forall\omega\in\mathcal{D}\cap\mathcal{E}\qquad \mathcal{P}_t^{-,\mathcal{V}}(\mathcal{D},\omega)\,\subset\,\mathcal{P}_t^{-,\mathcal{V}}(\mathcal{D}\cap\mathcal{E},\omega)\,. \tag{76.3}$$

It follows from (76.2) and (76.3) that

$$P\Big(\big|\mathcal{I}(t)\big|\,\geq\,\frac{m^{d-1}\beta(n)}{\ln n},\,\mathcal{D}\cap\mathcal{E}\Big)\,\leq\,P\Big(\big|\mathcal{P}_t^{-,\mathcal{V}}(\mathcal{D}\cap\mathcal{E})\big|\,\geq\,\frac{m^{d-1}\beta(n)}{\ln n},\,\mathcal{D}\cap\mathcal{E}\Big)\,. \tag{76.4}$$

Moreover, the good old Markov inequality gives

$$P\Big(\big|\mathcal{P}_t^{-,\mathcal{V}}(\mathcal{D}\cap\mathcal{E})\big|\geq\frac{m^{d-1}\beta(n)}{\ln n},\,\mathcal{D}\cap\mathcal{E}\Big)\,\leq\,\frac{\ln n}{m^{d-1}\beta(n)}E\Big(\big|\mathcal{P}_t^{-,\mathcal{V}}(\mathcal{D}\cap\mathcal{E})\big|1_{\mathcal{D}\cap\mathcal{E}}\Big)\,. \tag{76.5}$$

Substituting successively (76.5) in (76.4) and in (75.4), we get

$$P(\mathcal{D}\cap\mathcal{E})\,\leq\,\frac{(\ln n)^2}{m^{d-1}\beta(n)}\max_{1\,\leq\,t\,\leq\,\frac{\ln n}{\beta(n)}}E\Big(\,\big|\mathcal{P}_t^{-,\mathcal{V}}(\mathcal{D}\cap\mathcal{E})\big|1_{\mathcal{D}\cap\mathcal{E}}\,\Big)+n^{d-\beta(n)}\,. \tag{76.6}$$

We consider two cases. If $P(\mathcal{D}\cap\mathcal{E})\leq 2n^{d-\beta(n)}$, then we have a very satisfactory upper bound on $P(\mathcal{D}\cap\mathcal{E})$. Otherwise, suppose that

$$P(\mathcal{D}\cap\mathcal{E})\,>\,2n^{d-\beta(n)}\,. \tag{76.7}$$

Using (76.7), we rewrite (76.6) as

$$\max_{1\,\leq\,t\,\leq\,\frac{\ln n}{\beta(n)}}E\Big(\,\big|\mathcal{P}_t^{-,\mathcal{V}}(\mathcal{D}\cap\mathcal{E})\big|\,\Big|\,\mathcal{D}\cap\mathcal{E}\,\Big)\,\geq\,\frac{m^{d-1}\beta(n)}{2(\ln n)^2}\,.$$

Thus there exists $t^*$ such that $1\leq t^*\leq\frac{\ln n}{\beta(n)}$ and

$$E\Big(\,\big|\mathcal{P}_{t^*}^{-,\mathcal{V}}(\mathcal{D}\cap\mathcal{E})\big|\,\Big|\,\mathcal{D}\cap\mathcal{E}\,\Big)\,\geq\,\frac{m^{d-1}\beta(n)}{2(\ln n)^2}\,. \tag{76.8}$$

It follows from lemma 71.1 that

$$\begin{aligned}E\Big(\,\big|\mathcal{P}_{t^*}^{-,\mathcal{V}}(\mathcal{D}\cap\mathcal{E})\big|\,\Big|\,\mathcal{D}\cap\mathcal{E}\,\Big)\,&\leq\,E\Big(\,\big|\pi\big(\mathcal{SP}(t^*,\omega)\big)\big|\,\Big|\,\mathcal{D}\cap\mathcal{E}\,\Big)\\&\leq\,t^*E\Big(\,\big|\mathcal{SP}(t^*,\omega)\big|\,\Big|\,\mathcal{D}\cap\mathcal{E}\,\Big)\,,\end{aligned} \tag{76.9}$$

where the set $\mathcal{SP}(t^*,\omega)$ was defined at the beginning of section 70. Using the quantity $E^{\mathcal{V}}(\cdot,\cdot)$ defined in (74.1), we conclude from (76.8) and (76.9) that

$$E^{\mathcal{V}}(t^*,\mathcal{D}\cap\mathcal{E})\,\geq\,\frac{m^{d-1}\beta(n)}{2t^*(\ln n)^2}\,\geq\,\frac{m^{d-1}\big(\beta(n)\big)^2}{2(\ln n)^3}\,. \tag{76.10}$$

# 77 A convolution equation

We make a pause here to analyze the recurrence relation associated to the recursive inequality (69.13). Let $M, V$ be two positive integers (which play the role of the quantities $M(n), V(n)$ appearing in (69.13)). We study here the sequence $(u_k)_{k\geq 0}$ defined by a non-negative initial value $u_0$ and the recurrence relation

$$\forall k \geq 1 \qquad u_k \,=\, \sum_{0\leq \ell\leq k-1} u_\ell \, M\, k\big(kV\big)^{k-\ell}\,. \tag{77.1}$$

To be precise, the value $V$ will be taken to be $V(n)/(1-p)$, rather than $V(n)$.

**Proposition 77.1.** *The sequence $(u_k)_{k\geq 0}$ is non-negative, non-decreasing and satisfies*

$$\forall k \geq 1 \qquad u_0\big(MV\big)^k (k!)^2 \,\leq\, u_k \,\leq\, u_0\big(2MV\big)^k k^{3k}\,. \tag{77.2}$$

*Proof.* A simple induction shows that the sequence $(u_k)_{k\geq 0}$ is non-negative. The fact that it is non-decreasing is due to the parameters $M, V$ being larger or equal than one. Indeed, for any $k \geq 1$, looking at the term corresponding to $\ell = k-1$ in the convolution formula (77.1), we have

$$u_k \,\geq\, u_{k-1}\, M\, k\big(kV\big)^{k-(k-1)} \,=\, u_{k-1}\, MV\, k^2\,. \tag{77.3}$$

Inequality (77.3) readily implies that the sequence is non-decreasing. Moreover, by iterating this inequality, we obtain the lower bound on $u_k$ stated in (77.2). We prove next by induction on $k$ the upper bound on $u_k$ stated in (77.2). For $k = 1$, we have

$$u_1 \,=\, u_0\, M\, V \,\leq\, u_0\big(2MV\big)\,.$$

Let now $k \geq 2$ and suppose that we have proved the result until the rank $k-1$:

$$\forall \ell \in \{\,1,\dots,k-1\,\} \qquad u_\ell \,\leq\, \big(2MV\big)^\ell \ell^{3\ell} u_0\,.$$

Substituting all these inequalities into the convolution equality (77.1), we obtain

$$\begin{aligned} u_k \,&\leq\, \sum_{0\leq \ell\leq k-1} \big(2MV\big)^\ell \ell^{3\ell} u_0\, M\, k\big(kV\big)^{k-\ell} \\ &\leq\, u_0\big(2MV\big)^k k \sum_{0\leq \ell\leq k-1} \ell^{3\ell} k^{k-\ell}\,. \end{aligned} \tag{77.4}$$

We have next

$$k \sum_{0\leq \ell\leq k-1} \ell^{3\ell} k^{k-\ell} \,\leq\, k \sum_{0\leq \ell\leq k-1} k^{3\ell} k^{k-\ell} \,\leq\, k^2 k^{3k-2} \,=\, k^{3k}\,. \tag{77.5}$$

Substituting (77.5) into (77.4), we get

$$u_k \,\leq\, u_0\big(2MV\big)^k k^{3k}\,.$$

This proves the result at rank $k$ and it concludes the induction step. ☐

Our next goal is to compare an arbitrary sequence $(v_k)_{k\geq 1}$ with a sequence $(u_k)_{k\geq 0}$ which satisfies the convolution equality (77.1). The trick consists in starting with a value $u_0$ which is larger or equal than the sum of the upper deviations of $(v_k)_{k\geq 1}$ with respect to the convolution identity.

**Proposition 77.2.** *Let $(v_k)_{k\geq 1}$ be an arbitrary sequence of non-negative real numbers. Let us fix $K \geq 1$ and let us define*

$$\Delta(K) \,=\, \sum_{1\leq k\leq K} \Big(v_k - \sum_{1\leq \ell\leq k-1} v_\ell \, M\, k\big(kV\big)^{k-\ell}\Big)^+ . \tag{77.6}$$

*A sequence $(u_k)_{k\geq 0}$ satisfying the recurrence relation* (77.1) *which starts from a value $u_0$ such that $u_0 \geq \Delta(K)$ satisfies $u_K \geq v_K$.*

*Proof.* We proceed by induction on $K$. For $K = 1$, we have

$$\Delta(1) = (v_1)^+ = v_1\,, \quad u_1 = u_0 MV \geq v_1\,.$$

Let now $K \geq 2$ and suppose that we have proved the result until the rank $K-1$. We consider the sequence $(u_k)_{k\geq 0}$ satisfying the recurrence relation (77.1) which starts from $u_0 = \Delta(K)$. Since

$$\forall \ell \in \{\,1,\dots,K-1\,\} \qquad \Delta(K) \,\geq\, \Delta(\ell)\,,$$

it follows from the induction hypothesis that

$$\forall \ell \in \{\,1,\dots,K-1\,\} \qquad u_\ell \,\geq\, v_\ell\,.$$

Let us define

$$\delta(K) \,=\, \Big(v_K - \sum_{1\leq \ell\leq K-1} v_\ell \, M\, K\big(KV\big)^{K-\ell}\Big)^+ .$$

We have then

$$v_K \,\leq\, \delta(K) + \sum_{1\leq \ell\leq K-1} v_\ell \, M\, K\big(KV\big)^{K-\ell} \,\leq\, u_0 + \sum_{1\leq \ell\leq K-1} u_\ell \, M\, K\big(KV\big)^{K-\ell} \,\leq\, u_K\,.$$

This concludes the induction step and the proof of the proposition. ☐

We finally combine the results of propositions 77.1 and 77.2 to obtain a very general inequality, which will be used in conjunction with the exponential inequalities of section 74.

**Corollary 77.3.** *Let $(v_k)_{k\geq 1}$ be an arbitrary sequence of non-negative real numbers. For any $K \geq 1$, we have*

$$v_K \,\leq\, \Bigg(\sum_{1\leq k\leq K} \Big(v_k - \sum_{1\leq \ell\leq k-1} v_\ell \, M\, k\big(kV\big)^{k-\ell}\Big)^+\Bigg)\big(2MV\big)^K K^{3K}\,. \tag{77.7}$$

*Proof.* Let us fix $K \geq 1$. We consider the sequence $(u_k)_{k\geq 0}$ satisfying (77.1) which starts from $u_0 = \Delta(K)$, where $\Delta(K)$ is defined in (77.6). We simply combine the upper bound (77.2) on $u_K$ given by proposition 77.1 and proposition 77.2 to obtain the inequality (77.7). ☐

# 78 Estimate on $P(\mathcal{D})$

We apply corollary 77.3 to the sequence $v_k = E^{\mathcal{V}}(k, \mathcal{D}\cap\mathcal{E})$, $1 \leq k \leq n$ at the index $K = t^*$. For the parameters $M, V$, we take the quantities $M = M(n)$, $V = V(n)$ associated to the covering $\mathcal{V}$. In fact, we will need an extra factor $1-p$, so we rather use $V = V(n)/(1-p)$. We use also the inequality (76.10) to bound from below $E^{\mathcal{V}}(t^*, \mathcal{D}\cap\mathcal{E})$ and we obtain

$$\frac{m^{d-1}\big(\beta(n)\big)^2}{2(\ln n)^3} \leq E^{\mathcal{V}}(t^*, \mathcal{D}\cap\mathcal{E}) \leq \\ \Bigg(\sum_{1\leq k\leq t^*}\Big(E^{\mathcal{V}}(k, \mathcal{D}\cap\mathcal{E}) - \sum_{1\leq \ell\leq k-1} E^{\mathcal{V}}(\ell, \mathcal{D}\cap\mathcal{E})\, M(n)\, k\Big(\frac{kV(n)}{1-p}\Big)^{k-\ell}\Big)^+\Bigg) \\ \times\Big(2M(n)\frac{V(n)}{1-p}\Big)^{t^*}(t^*)^{3t^*}\,.$$

Therefore there exists $k^*$ such that $1 \leq k^* \leq t^* \leq \frac{\ln n}{\beta(n)}$ and

$$\Big(E^{\mathcal{V}}(k^*, \mathcal{D}\cap\mathcal{E}) - \sum_{1\leq \ell\leq k^*-1} E^{\mathcal{V}}(\ell, \mathcal{D}\cap\mathcal{E})\, M(n)\, k^*\Big(\frac{k^*V(n)}{1-p}\Big)^{k^*-\ell}\Big)^+ \\ \geq \frac{1}{\Big(2M(n)\frac{V(n)}{1-p}\Big)^{t^*}(t^*)^{3t^*+1}}\frac{m^{d-1}\big(\beta(n)\big)^2}{2(\ln n)^3} \\ \geq \frac{1}{\Big(2M(n)\frac{V(n)}{1-p}\Big)^{\frac{\ln n}{\beta(n)}}\Big(\frac{\ln n}{\beta(n)}\Big)^{\frac{3\ln n}{\beta(n)}+1}}\frac{m^{d-1}\big(\beta(n)\big)^2}{2(\ln n)^3}\,. \tag{78.1}$$

The time has come to use a quantitative estimate on $\beta(n)$. The estimate of the termination time was done in proposition 30.1. The function $\beta(n)$ there was the one of proposition 16.1 (actually, the function $\beta(n)$ there should be replaced by $\beta(n/2)$). This function $\beta(n)$ was defined in (16.3) in the course of the proof of proposition 16.1, by setting

$$\forall n\geq 3 \qquad \beta(n) = \min\Big(\sqrt{\ln n}, \sqrt{\alpha\big(\sqrt{\ln n}\big)}\Big)\,. \tag{78.2}$$

The function $\alpha(t)$ was introduced in corollary 14.5, where it was also proved that, for $t$ larger than a value $t_0(\theta(p), p, d)$ depending on $\theta(p), p, d$, we have

$$\alpha(t) \geq \ln\Big(\frac{\theta(p)}{2}\Big) - \frac{1}{2}\ln\Bigg(-\phi\Big(P_p\Big(\Big(\frac{t}{4d}\Big)^{d/(d+1)} < |C(0)| < +\infty\Big)\Big)\Bigg). \tag{78.3}$$

Inequality (78.3) and the hypothesis 75.1 imply that there exist $\alpha > 0$ such that

$$\forall t\geq t_0(\theta(p), p, d) \qquad \alpha(t) \geq t^{\frac{\alpha d}{3(d+1)}}\,. \tag{78.4}$$

Substituting inequality (78.4) into (78.2), we conclude that there exist $\gamma$ in $]0,1/2[$ and a value $n_0(\theta(p),p,d)$ depending on $\theta(p),p,d$ such that

$$\forall n\geq n_0(\theta(p),p,d)\qquad \beta(n)\,\geq\,(\ln n)^\gamma\,. \tag{78.5}$$

This implies furthermore that

$$\forall n\geq n_0(\theta(p),p,d)\qquad k^*\,\leq\,(\ln n)^{1-\gamma}\,. \tag{78.6}$$

Substituting (78.5) into (78.1), we get that, for $n\geq n_0(\theta(p),p,d)$,

$$\Big(E^{\mathcal{V}}(k^*,\mathcal{D}\cap\mathcal{E})-\sum_{1\leq\ell\leq k^*-1}E^{\mathcal{V}}(\ell,\mathcal{D}\cap\mathcal{E})\,M(n)\,k^*\Big(\frac{k^*V(n)}{1-p}\Big)^{k^*-\ell}\Big)^+$$
$$\geq\,\frac{m^{d-1}}{\Big(2M(n)\dfrac{V(n)}{1-p}\ln n\Big)^{5(\ln n)^{1-\gamma}}}\,. \tag{78.7}$$

It remains to apply the exponential inequality of proposition 74.1 to $\mathcal{D}\cap\mathcal{E}$ and the index $k^*$. Suppose that $m=n^\delta$, for some $\delta$ such that $(2d-2)\delta>1$. We use the inequalities (76.1) and (78.7) to bound the fraction in the exponential appearing in (74.3), and we obtain that

$$P(\mathcal{D}\cap\mathcal{E})\,\leq\,\exp\left(k^*-\frac{\left(\dfrac{n^{(d-1)\delta}}{\Big(2M(n)\frac{V(n)}{1-p}\ln n\Big)^{5(\ln n)^{1-\gamma}}}\right)^{2/k^*}}{4(k^*)^{2/k^*}C(p)\,M(n)\frac{V(n)}{1-p}n^{1/k^*}}\right)$$
$$\leq\,\exp\left(k^*-\frac{\left(\dfrac{n^{2(d-1)\delta-1}}{\big(2\,c^2(\ln n)^{2b+1}\big)^{10(\ln n)^{1-\gamma}}}\right)^{1/k^*}}{4\,C(p)\,c^2(\ln n)^{2b+2}}\right). \tag{78.8}$$

Thanks to the condition $2(d-1)\delta-1>0$, for $n$ large enough, we will have

$$\frac{n^{2(d-1)\delta-1}}{\big(2\,c^2(\ln n)^{2b+1}\big)^{10(\ln n)^{1-\gamma}}}\,>\,1\,. \tag{78.9}$$

Using once more the upper bound (78.6) on $k^*$, we conclude from inequalities (78.8) and (78.9) that

$$P(\mathcal{D}\cap\mathcal{E})\,\leq\,\exp\left((\ln n)^{1-\gamma}-\frac{1}{2^{12}\,C(p)\,c^{22}(\ln n)^{22b+12}}\Big(n^{2(d-1)-\delta}\Big)^{(\ln n)^{\gamma-1}}\right).$$

Coming back to inequality (75.3), we conclude that

$$P(\mathcal{D})\,=\,O\Big(\frac{1}{n^{dr}}\Big)\,.$$

We have finally succeeded in getting a quantitative estimate on the probability of the disconnection event $\mathcal{D}$. Notice that we can take any value for the parameter $r$. The initial assumptions were that $p$ is such that $\theta(p,\mathbb{Z}^d)>0$ and the hypothesis 75.1.

# 79 Conclusions

In the final section of this part, we revisit the three challenging questions introduced in section 30, and we present the results we can achieve if we assume that the parameter $p$ is such that the hypothesis 75.1 is in force and $\theta(p,\mathbb{Z}^d)>0$.

## 79.1 Disconnection of a sub-box

Let $m,n$ be two integers with $1\leq m\leq n$. We work inside $\Lambda=\Lambda(2n)$ and we fix a sub-box of $\Lambda(2n)$ which is a translate of $\Lambda(m)$ centered at a vertex $x$ of $\Lambda(n)$. The disconnection event $\mathcal{D}$ is

$$\mathcal{D}(x,m)\,=\,\big\{\,x+\Lambda(m)\not\longleftrightarrow\partial^{\,in}\Lambda(2n)\,\big\}\,.$$

Applying the previous strategy, we conclude that if $m=n^\delta$, for some $\delta$ such that $(2d-2)\delta>1$, then

$$\forall r\geq 1\qquad P(\mathcal{D})\,=\,O\Big(\frac{1}{n^{dr}}\Big)\,. \tag{79.1}$$

We sum this estimate over all possible values of $x$ in $\Lambda(n)$ and all values of $m$ in the interval $[n^{1/(2d-2)+\varepsilon},2n]$, and we conclude that, for any $\varepsilon>0$,

$$P\Big(\begin{matrix}\text{there exists } x\in\Lambda(n) \text{ and } m\geq n^{1/(2d-2)+\varepsilon}\\ \text{such that } x+\Lambda(m)\not\longleftrightarrow\partial^{\,in}\Lambda(2n)\end{matrix}\Big)\,=\,O\Big(\frac{1}{n^{d(r-1)-1}}\Big)\,.$$

## 79.2 The largest cluster in a box

Let $m,n$ be two integers with $1\leq m\leq n$. We work inside $\Lambda=\Lambda(n)$ and we fix a vertex $x$ of $\Lambda(n)$. We consider the event $\mathcal{D}(x,m)$ defined by

$$\mathcal{D}(x,m)\,=\,\big\{\,x\not\longleftrightarrow\partial^{\,in}\Lambda(n)\,,|C(x)|\geq m^d\,\big\}\,.$$

Our strategy yields the very same result (79.1) as for the previous question. Summing (79.1) over $x$ in $\Lambda(n)$ and over $m$, we conclude that, for any $\varepsilon>0$,

$$P\Big(\begin{matrix}\text{there exists a vertex } x \text{ in } \Lambda(n) \text{ such that}\\ x\not\longleftrightarrow\partial^{\,in}\Lambda(n) \text{ and } |C(x)|\geq n^{d/(2d-2)+\varepsilon}\end{matrix}\Big)\,=\,O\Big(\frac{1}{n^{d(r-1)-1}}\Big)\,.$$

Thus we have succeeded in showing that, for any $\varepsilon>0$, with high probability, the largest finite cluster in $\Lambda(n)$ has cardinality less than $n^{d/(2d-2)+\varepsilon}$.

## 79.3 Disconnection in a cuboid

This question would need some little additional work to fit in our strategy. Namely, some further arguments are required to show that, when the event $\mathcal{E}$ occurs, the pivots of order $k$ are indeed local. We do not give the details here, but in the end, using the notations of subsection 30.3, we would prove that

$$\forall r>0\qquad P\big(L\not\longleftrightarrow R\text{ in }\Pi(n)\big)\,=\,O\Big(\frac{1}{n^r}\Big)\,.$$

This is a precious quantitative estimate on a disconnection in a finite domain.

## Part XI

# Proof for the site model ⊙ ∘

This part contains the final steps of the proof of theorem 9.2. The central tool used in the proof is the multiple intertwined explorations process, it is constructed in section 80. The proof is ignited in section 85. It contains two main steps. In the first step, carried out in section 87, we derive a control on the number of large clusters present in a cubical shell. In a second step, carried out in section 93.3, we use the isoperimetric theorem to show that, with overwhelming probability, the largest cluster in a finite box is of volume order. The central issue is to obtain a quantitative control on the intersections revealed by the multiple intertwined explorations. Unfortunately, this program could be completed only locally with the inequalities developed in subsection 89.1.

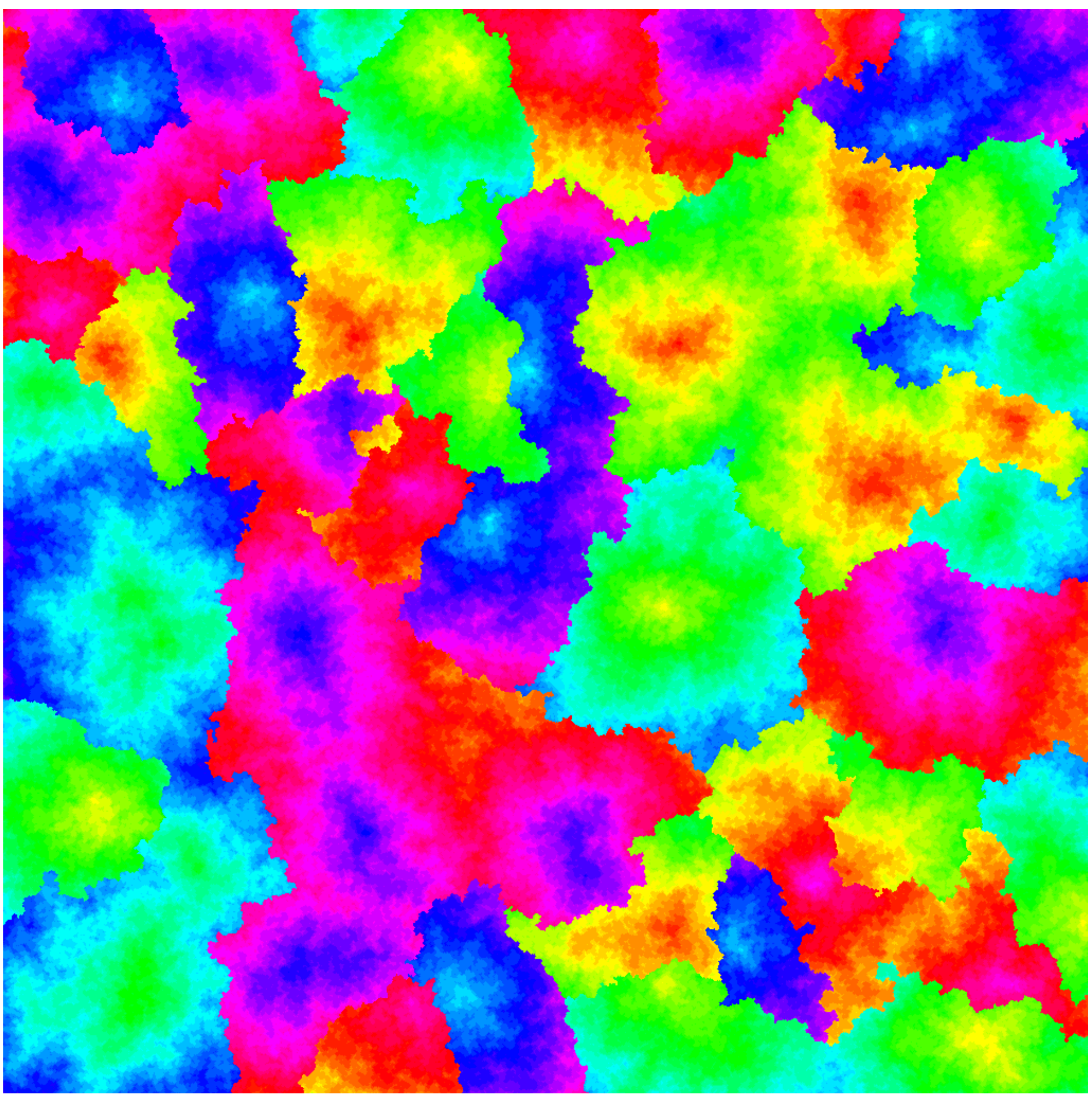

Figure 142: Fifty intertwined explorations, $\Lambda(1024)$, site percolation, $p = 0.5$.

# 80 Intertwined site explorations

In this section, we will build the counterpart for the site percolation model of the objects and algorithms defined in section 26 for the bond model. Although we follow essentially the same guidelines to design the taboo exploration, the two models present several little subtle differences, due to the nature of the blocking mechanism at work. For instance, for the bond model, the waiting set was a part of the outer boundary of the aggregate, whereas for the site model, it will be a part of the inner boundary. This is why we have to write completely the definitions and the proofs for each model.

## 80.1 Waiting site aggregates

We recall first some definitions that were introduced in subsections 12.1 and 12.3. As usual, we use $|\cdot|$ to denote the Euclidean norm. The set $\mathcal{N}(A)$ of the neighbours of a subset $A$ of $\mathbb{Z}^d$ is defined as

$$\mathcal{N}(A) \,=\, \big\{\, y\in\mathbb{Z}^d\setminus A : \exists\, x\in A \qquad |x-y|=1 \,\big\}\,.$$

In case $A$ is a singleton $A=\{\,x\,\}$, we write simply $\mathcal{N}(x)$ instead of $\mathcal{N}(\{\,x\,\})$. Let $D$ be a subset of $\mathbb{Z}^d$. For $A$ a subset of $D$, we define its internal boundary $\partial^{in}_D A$ in $D$ and its external boundary $\partial^{out}_D A$ in $D$ by

$$\begin{aligned}\partial^{in}_D A \,&=\, \big\{\, x\in A : \exists y\in D\setminus A \quad |x-y|=1 \,\big\}\,,\\ \partial^{out}_D A \,&=\, \big\{\, x\in D\setminus A : \exists y\in A \quad |x-y|=1 \,\big\}\,.\end{aligned}$$

For $x$ in $D$, the open cluster of $x$ in $D$ is the set $C(x,D)$ of the sites which are connected to $x$ by an open path, i.e.,

$$C(x,D) \,=\, \big\{\, y\in D : x\longleftrightarrow y \text{ in } D \,\big\}\,.$$

We define also the set $\overline{C}(x,D)$ of the sites whose travel distance in $D$ to $x$ is null, i.e.,

$$\forall x\in D \qquad \overline{C}(x,D) \,=\, C(x,D)\cup\partial^{out}_D C(x,D)\,.$$

There remains a slight ambiguity to be resolved in the case of a closed site. We make the convention that, for a closed site $x$ of $D$, we have $C(x,D)=\varnothing$ and $\partial^{out}_D C(x,D)=\{\,x\,\}$, so that $\overline{C}(x,D)=\{\,x\,\}$. Another possible convention would be to say that $C(x,D)=\{\,x\,\}$ whenever the site $x$ is closed. However, our choice is in accordance with another choice we did when defining the site percolation model. Indeed, we considered the random graph obtained by removing all the closed sites and all the edges incident to them. An alternative choice would be to keep all the vertices of $\mathbb{Z}^d$ and to remove all the edges having one endpoint closed. With this choice, we would naturally decide that $C(x,D)=\{\,x\,\}$ whenever the site $x$ is closed. Another consequence of this last choice would be that, in the pictures, all closed sites should be given a label and a color, and the visual effect would be worse than the simple blurring created by the presence of the black pixels. Unfortunately, this little decision is not innocent, and it creates some constraints for the mathematical definitions that are to follow.

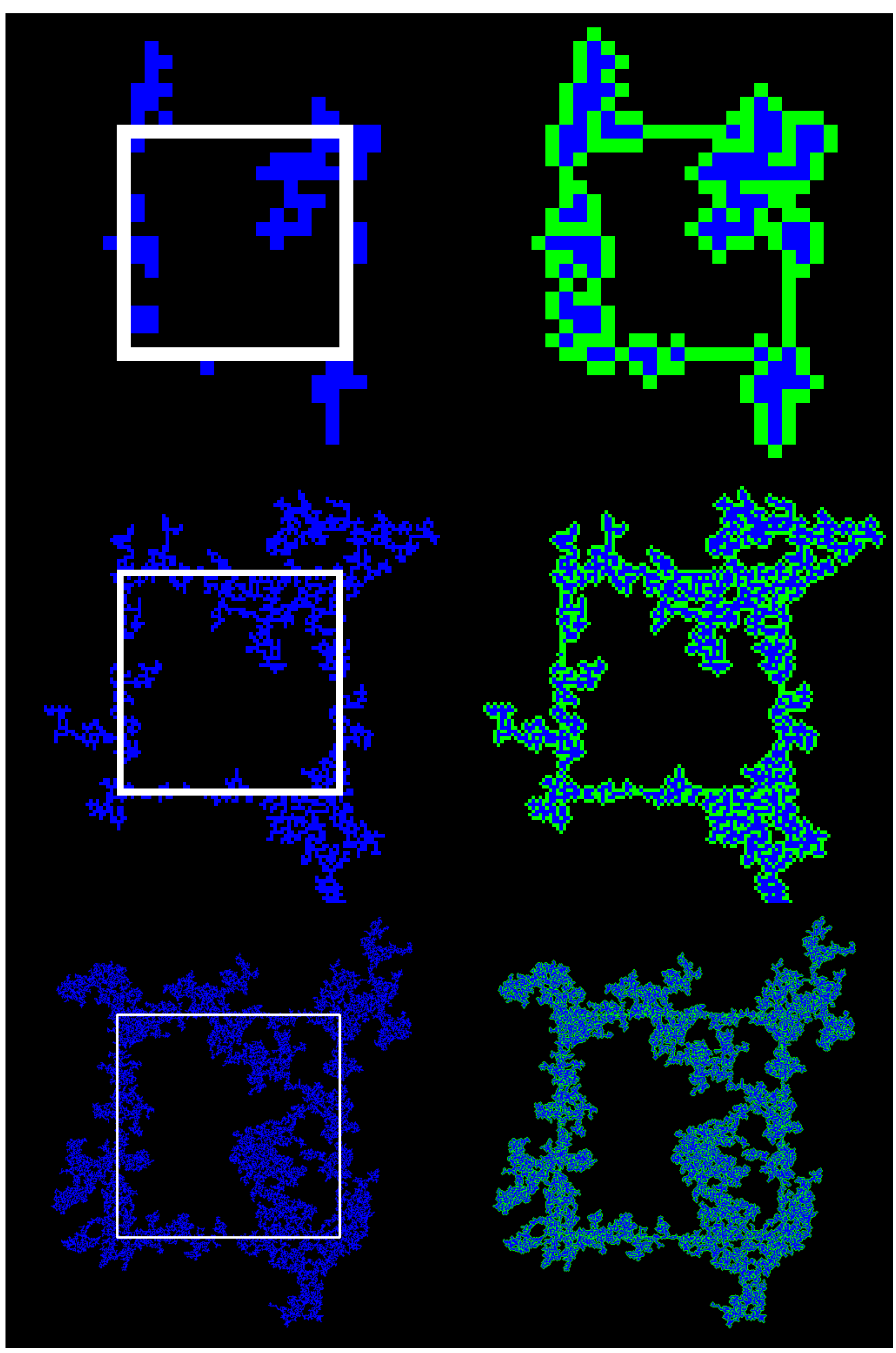

Figure 143: Three examples of $C(A,D)$ and $\overline{C}(A,D)$: $A$ is the white square, top: $\Lambda(32)$, $p = 0.4$, middle: $\Lambda(128)$, $p = 0.53$, bottom: $\Lambda(512)$, $p = 0.57$. $C(A,D)$ is in blue, $\overline{C}(A,D) \setminus C(A,D)$ is in green.

For $A$ a subset of $\mathbb{Z}^d$ (not necessarily included in $D$), we define

$$C(A,D) \,=\, \bigcup_{x\in A\cap D} C(x,D)\,,\qquad \overline{C}(A,D) \,=\, \bigcup_{x\in A\cap D} \overline{C}(x,D)\,.$$

**Definition 80.1.** *A subset $\mathcal{A}$ of $D$ is an aggregate of clusters in $D$, or simply an aggregate, if*

$$\forall x\in\mathcal{A}\qquad C(x,D)\subset\mathcal{A}\,.$$

*A set $\mathcal{W}$ is a waiting set associated to an aggregate $\mathcal{A}$ if it is a subset of its inner boundary in $D$, i.e., $\mathcal{W}\subset\partial_D^{in}\mathcal{A}$. An aggregate $\mathcal{A}$ endowed with a waiting set $\mathcal{W}$ is called a waiting aggregate.*

Let us mention a few simple properties of aggregates. If $\mathcal{A}$ is an aggregate and $x$ is not in $\mathcal{A}$, then $C(x,D)\subset D\setminus\mathcal{A}$. Whether or not a set $\mathcal{A}$ is an aggregate depends solely on the set of the sites in $\mathcal{A}$ that are open. In particular, any set obtained by adjoining closed sites to an aggregate is still an aggregate. For instance, given a subset $A$ of $\mathbb{Z}^d$, any set $C$ such that

$$C(A,D)\,\subset\, C\,\subset\,\overline{C}(A,D)$$

is also an aggregate. Finally, the union of two aggregates is also an aggregate.

At this stage, the reader may wonder why the waiting aggregates have to be defined differently in the bond and the site models. It is in fact a consequence of the nature of the random blocking mechanism at work in each model. The purpose of the waiting aggregate is to record the portion of the configuration which has been examined by an explorer, as well as the part of the configuration from where the exploration will restart upon reactivation. This is why the interpretation of the sets $\mathcal{A},\mathcal{W}$ of the waiting aggregates are quite different in the bond and the site models.

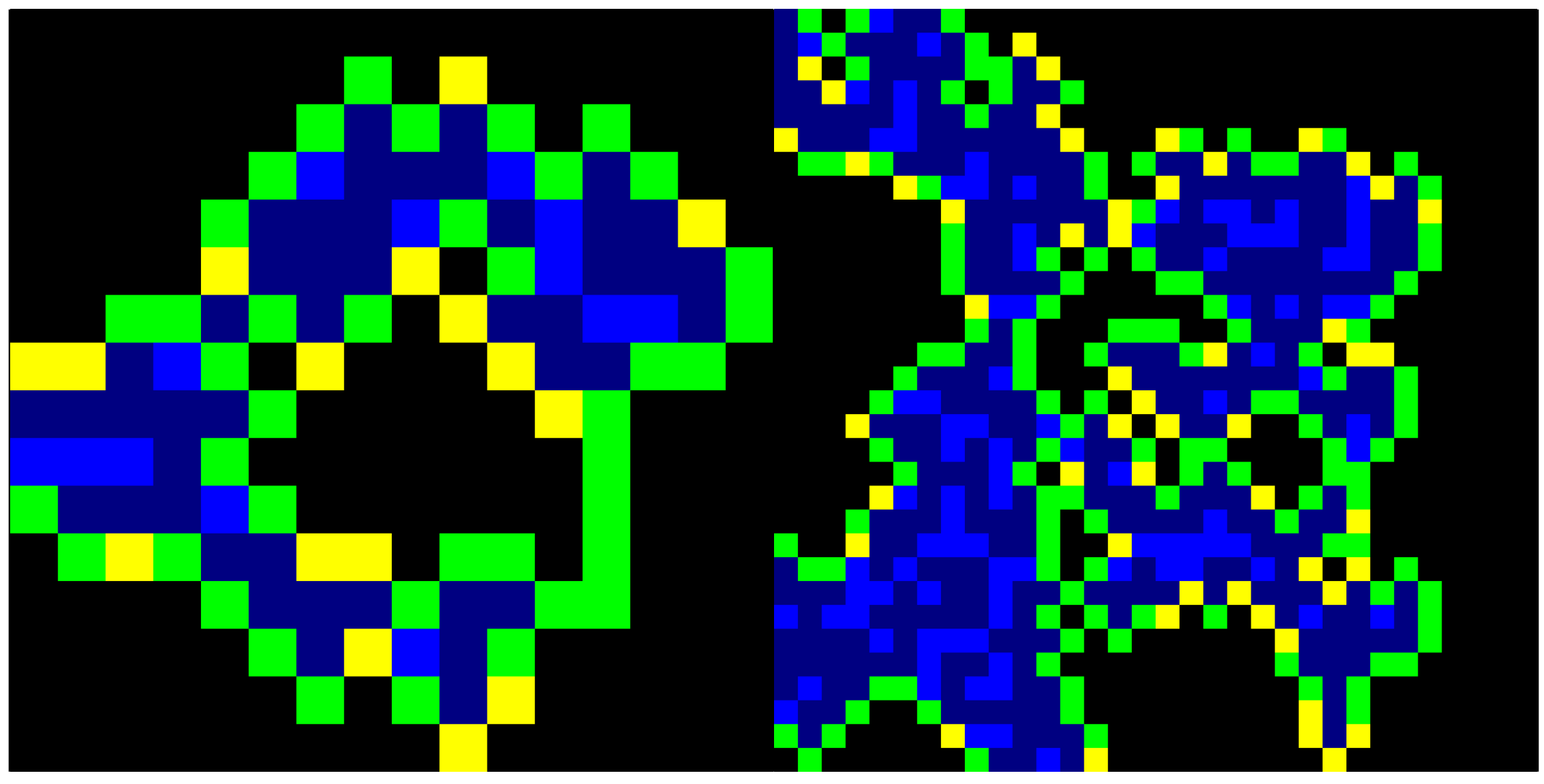

Figure 144: Two waiting aggregates: $\mathcal{A}$ in blue, $\partial_D^{in}\mathcal{A}$ in green. $\mathcal{W}$ in yellow.
Left: $\Lambda(16)$, $p=0.4$ (left). Right: $\Lambda(32)$, $p=0.5$ (right).

## 80.2 Taboo site exploration algorithm

We start from the standard exploration algorithm described in section 5, and we introduce a taboo set $\mathcal{T}$, which the algorithm is not permitted to visit. The pseudocode of the taboo exploration algorithm is presented as algorithm 80.1. We give next the corresponding mathematical definition. For $A$ a subset of $D$, we define

$$\overline{A}^D \;=\; A \cup \partial_D^{out} A\,.$$

Let $A$ be a subset of $\mathbb{Z}^d$ and let $\mathcal{T}$ be a subset of $D$. Upon termination, the taboo exploration algorithm starting from $A$ with taboo set $\mathcal{T}$ will have explored the set $\overline{\text{Clusters}}\,(A, D, \mathcal{T})$ defined as

$$\overline{\text{Clusters}}\,(A, D, \mathcal{T}) \;=\; \big(D \cap A \setminus \mathcal{T}\big) \cup \overline{C\big(A, D \setminus \mathcal{T}\big)}^{D\setminus\mathcal{T}}\,. \tag{80.1}$$

This set consists of the sites which can be reached in $D \setminus \mathcal{T}$ from $A \cap D \setminus \mathcal{T}$ with a null travel time. The definition of the set $\text{Clusters}\,(A, D, \mathcal{T})$ for the bond model in subsection 29.2 was essentially the same, although the definition of the travel time is different for the site model. So, if one wishes to give a formulation which works for both models, this seems to be the correct choice. As for the bond model, the map $\overline{\text{Clusters}}\,(A, D, \mathcal{T})$ is non-increasing with respect to $\mathcal{T}$. Another consequence of the definition (80.1) is that $\overline{\text{Clusters}}\,(A, D, \mathcal{T})$ is a subset of $D \setminus \mathcal{T}$, even if $A$ is not. In fact, the starting set of the taboo exploration is taken to be $A \cap D \setminus \mathcal{T}$, and the exploration is confined to $D \setminus \mathcal{T}$. During the exploration, the algorithm records a set of waiting sites, which is

$$\text{Waiting}\,(A, D, \mathcal{T}) \;=\; \partial_D^{in}\overline{\text{Clusters}}\,(A, D, \mathcal{T})\,. \tag{80.2}$$

The next round of exploration will start from the set $\text{Waiting}\,(A, D, \mathcal{T})$.

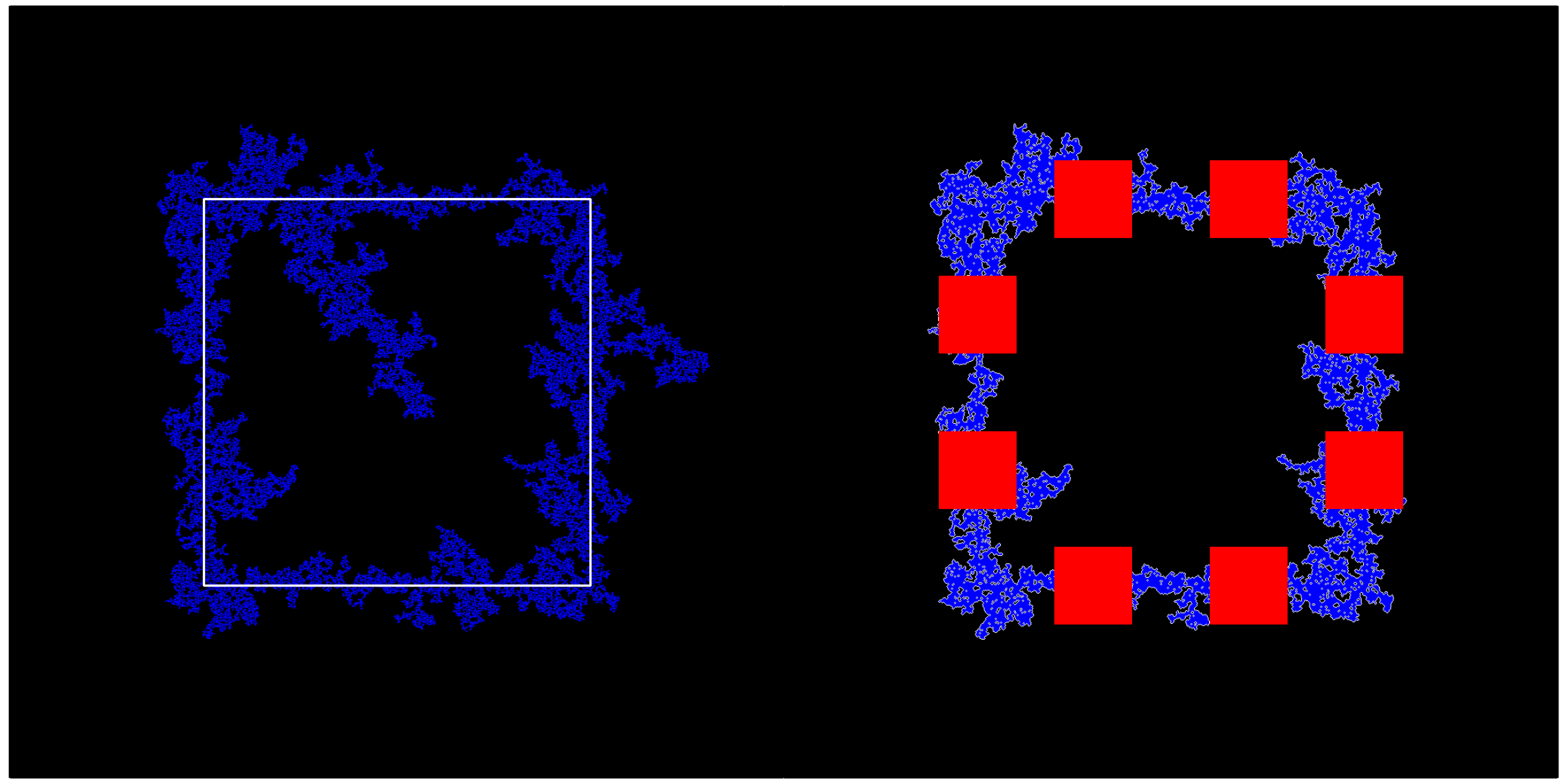

Figure 145: Taboo site exploration, $D = \Lambda(1024)$, $p = 0.577$, $A = \partial\Lambda(256)$. Left: $C(A, D)$. Right: $\text{Clusters}\,(A, D, \mathcal{T})$, $\text{Waiting}\,(A, D, \mathcal{T})$, $\mathcal{T} = 8$ squares. The clusters of $C(A, D)$ are visible in figure 146.

**Algorithm 80.1** Taboo_Site_Explore_Clusters $(A, D, \mathcal{T})$

**Require:** a non-empty starting set $A$, an authorized domain $D$, a taboo set $\mathcal{T}$
**Ensure:** $\overline{C}$ is $\overline{\text{Clusters}}\,(A, D, \mathcal{T})$

$\overline{C} \leftarrow \varnothing$
**repeat**
  Pick up $x \in A$ and remove it from $A$ ▷ $A$ is the set of the active sites
  **if** $x$ is in $D \setminus \mathcal{T}$ **then**
    $\overline{C} \leftarrow \overline{C} \cup \{x\}$
    **if** $x$ is open **then**
      **for** $y \in \mathcal{N}(x) \cap D \setminus \mathcal{T}$ **do**
        **if** $y \notin \overline{C}$ **then** ▷ Check that $y$ has not been examined
          $A \leftarrow A \cup \{y\}$ ▷ Put $y$ in the set of active sites
        **end if**
      **end for**
    **end if**
  **end if**
**until** $A = \varnothing$
**return** $\overline{C}$

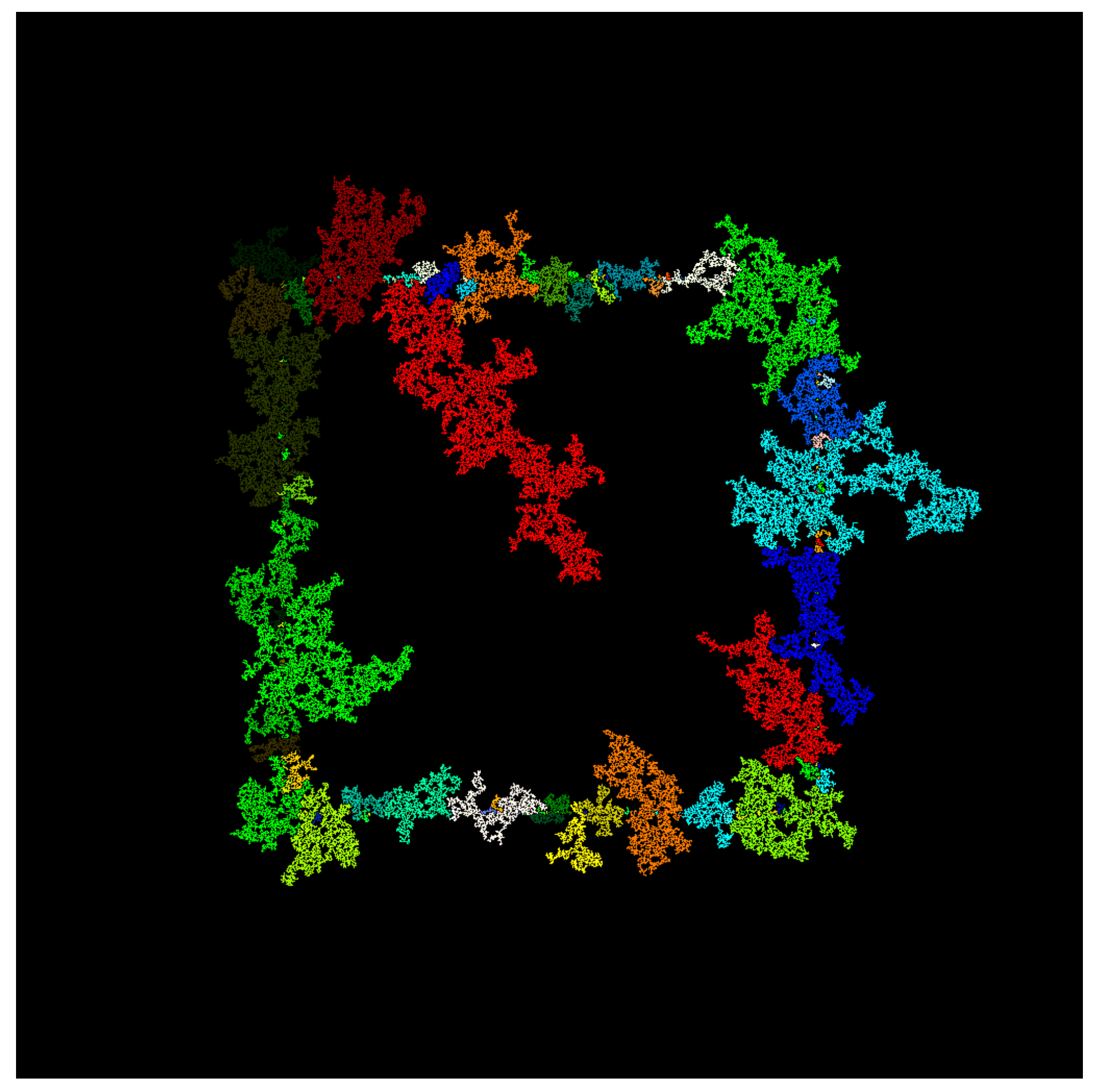

Figure 146: Clusters of $C\big(\partial\Lambda(256), \Lambda(1024)\big)$, site percolation at $p = 0.577$.

## 80.3 Coherence site exploration

The exploration processes which are to be defined shortly will involve pairs of aggregates satisfying certain conditions. We call such pairs coherent, and we define them next.

**Definition 80.2.** *Two waiting site aggregates $(\mathcal{A},\mathcal{W})$, $(\mathcal{A}',\mathcal{W}')$ are coherent if*
- *the sets $\mathcal{A}$ and $\mathcal{A}'$ are disjoint, i.e., $\mathcal{A}\cap\mathcal{A}'=\varnothing$;*
- *the inner boundary of $\mathcal{A}\cup\mathcal{A}'$ in $D$ is covered by $\mathcal{W}$ and $\mathcal{W}'$, i.e.,*

$$\partial_D^{in}\big(\mathcal{A}\cup\mathcal{A}'\big)\,\subset\,\mathcal{W}\cup\mathcal{W}'\,.$$

We proceed now to the definition of the map Coherent_Site_Explore. This map encodes the step at the heart of the intertwined explorations, that we shall build later.

**Definition 80.3.** *Let $(\mathcal{A},\mathcal{W})$, $(\mathcal{A}',\mathcal{W}')$ be two coherent waiting site aggregates. We define a map Coherent_Site_Explore on the pairs of coherent waiting aggregates by setting*

$$\begin{multline*}\mathit{Coherent_Site_Explore}\Big((\mathcal{A},\mathcal{W}),(\mathcal{A}',\mathcal{W}')\Big)\,=\\ \Big(\mathcal{A}\cup\overline{\mathit{Clusters}}\big(\mathcal{N}(\mathcal{W}),D\setminus\mathcal{A},\mathcal{A}'\big),\,\mathit{Waiting}\big(\mathcal{N}(\mathcal{W}),D\setminus\mathcal{A},\mathcal{A}'\big)\Big).\quad(80.3)\end{multline*}$$

The second waiting set $\mathcal{W}'$ plays no role in the previous definition, and it could be omitted at this point. However, for reasons that will appear later, we prefer to have a map defined on pairs of waiting aggregates, so we keep $\mathcal{W}'$. The definition 80.3 is more intricate than its counterpart for the bond model, given in definition 29.2, but that seems to be the order of things. The only apparent difference between the two definitions is that we use the set $\overline{\text{Clusters}}$ here instead of Clusters there, but this is not purely cosmetic. In fact, we did not manage to find a working definition for the site model directly, we had to resolve first the situation for the bond model in order to understand properly how the coherent map should be defined. Let us introduce the notation

$$(\mathcal{A}'',\mathcal{W}'')\,=\,\text{Coherent_Site_Explore}\big((\mathcal{A},\mathcal{W}),(\mathcal{A}',\mathcal{W}')\big)\,.$$

The definition of the sets $\mathcal{A}''$, $\mathcal{W}''$ can be rewritten as follows:

$$\mathcal{A}''\,=\,\mathcal{A}\cup\big(D\cap\mathcal{N}(\mathcal{W})\setminus(\mathcal{A}\cup\mathcal{A}')\big)\cup\overline{C\big(\mathcal{N}(\mathcal{W}),D\setminus(\mathcal{A}\cup\mathcal{A}')\big)}^{D\setminus(\mathcal{A}\cup\mathcal{A}')},\quad(80.4)$$

$$\mathcal{W}''\,=\,\partial_{D\setminus\mathcal{A}}^{in}\Big(\big(D\cap\mathcal{N}(\mathcal{W})\setminus(\mathcal{A}\cup\mathcal{A}')\big)\cup\overline{C\big(\mathcal{N}(\mathcal{W}),D\setminus(\mathcal{A}\cup\mathcal{A}')\big)}^{D\setminus(\mathcal{A}\cup\mathcal{A}')}\Big).\quad(80.5)$$

Formula (80.4) results from a direct application of the definitions (80.1) and (80.3), while formula (80.5) comes directly from the definitions (80.2) and (80.3). Notice that the maps Coherent_Bond_Explore and Coherent_Site_Explore are quite different. For instance, for the site model, the exploration does not restart exactly from the waiting set $\mathcal{W}$, but from its neighborhood $\mathcal{N}(\mathcal{W})$. In the bond model, it suffices to look at the bonds which have one endpoint in the waiting set, whereas in the site model, the sites of the waiting set have already been examined during the previous round, so we have to look at their neighbours when restarting. Of course, we could have taken this set of neighbours as the waiting set, in which case the waiting set would again have been a subset of the outer boundary of the aggregate. In fact, we tried this approach, but it has some disadvantages. Computationally, it leads to useless operations, because a lot of sites might be put in the waiting set and discarded later on. During the exploration of a cluster, it happens very frequently that a neighbour of a closed site is reachable through an open path that will be discovered after the visit to the closed site (see figure 147). Thus it is not a good idea to compute during the exploration all the neighbours of the closed sites and to put them on a waiting list. The best strategy seems to record only the closed sites, and look at their neighbours once the exploration of the clusters is completed, as we do here.

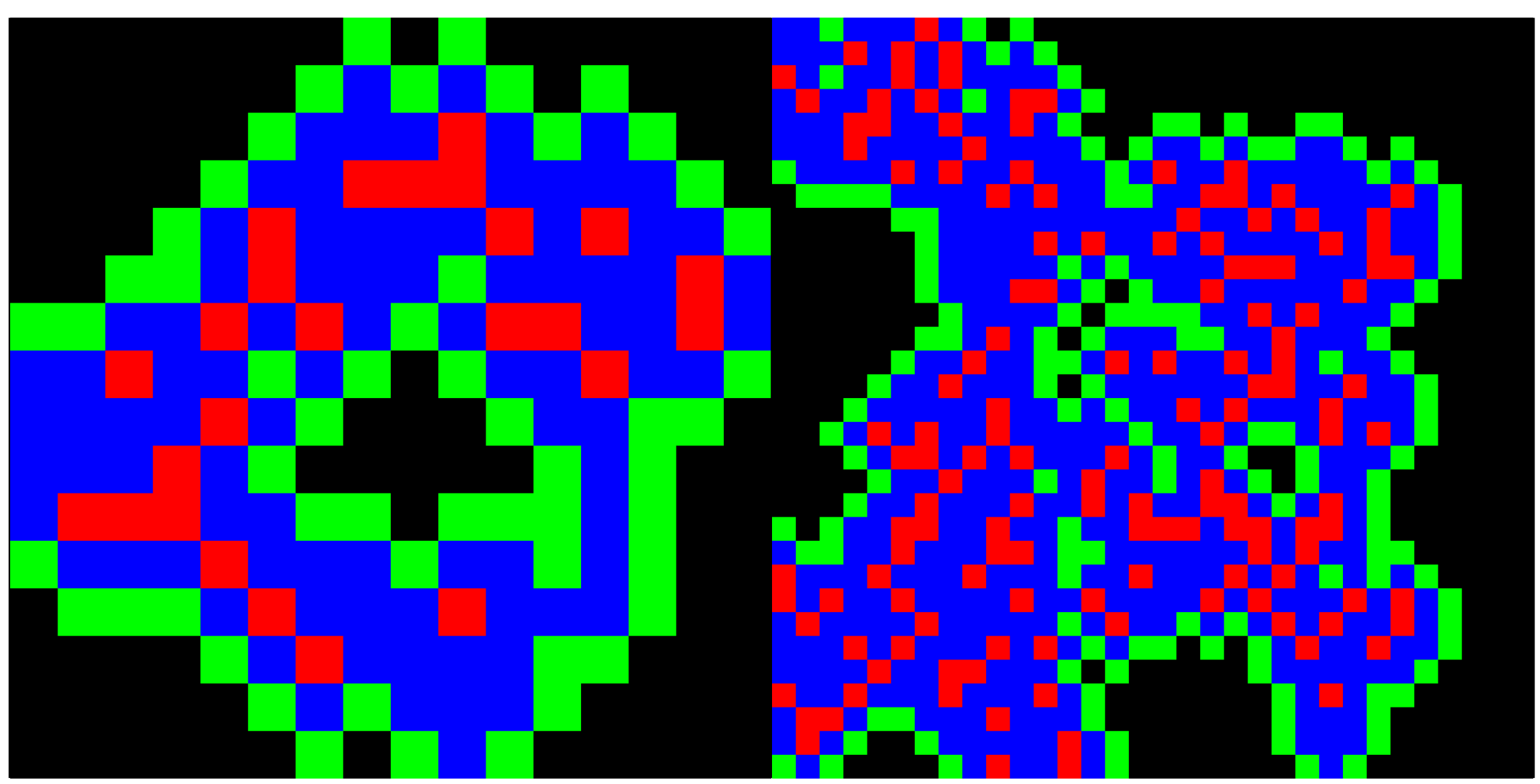

Figure 147: Same configurations as in figure 144. $\mathcal{A}$ in blue, $\partial_D^{out}\mathcal{A}$ in green. The sites which were put in the waiting set and then removed are in red.

The point of keeping the two waiting sets $\mathcal{W}, \mathcal{W}'$ when defining the map Coherent_Site_Explore is expressed through the next result, which shows that the coherence is maintained by the map.

**Proposition 80.4.** *Suppose that* $(\mathcal{A}, \mathcal{W})$, $(\mathcal{A}', \mathcal{W}')$ *are two coherent waiting site aggregates in* $D$*. Then* $(\mathcal{A}'', \mathcal{W}'')$ *is a waiting site aggregate in* $D$ *and moreover* $(\mathcal{A}', \mathcal{W}')$ *and* $(\mathcal{A}'', \mathcal{W}'')$ *are coherent.*

*Proof.* Let us prove first that $\mathcal{A}''$ is an aggregate in $D$. Let $x$ be an open site belonging to $\mathcal{A}''$. If $x$ belongs to $\mathcal{A}$, then $C(x, D)$ is included in $\mathcal{A}$ because $\mathcal{A}$ is an aggregate in $D$. Suppose that $x$ is in $\mathcal{A}'' \setminus \mathcal{A}$. We see from formula (80.4) that

$$\mathcal{A}'' \setminus \mathcal{A} = \big(D \cap \mathcal{N}(\mathcal{W}) \setminus (\mathcal{A} \cup \mathcal{A}')\big) \cup \overline{C\big(\mathcal{N}(\mathcal{W}), D \setminus (\mathcal{A} \cup \mathcal{A}')\big)}^{D \setminus (\mathcal{A} \cup \mathcal{A}')}$$

Thus there exists $y$ in $\mathcal{N}(\mathcal{W}) \cap D \setminus (\mathcal{A} \cup \mathcal{A}')$ such that

$$x \in \overline{C\big(y, D \setminus (\mathcal{A} \cup \mathcal{A}')\big)}^{D\setminus(\mathcal{A}\cup\mathcal{A}')}.$$

By definition, there exists a path $z_0, z_1, \dots, z_r$ in $D \setminus (\mathcal{A} \cup \mathcal{A}')$ joining $z_0 = y$ to $z_r = x$ such that all the sites $z_0, \dots, z_{r-1}$ are open. Since $x$ is open as well, we conclude that $x$ and $y$ are connected by an open path in $D \setminus (\mathcal{A} \cup \mathcal{A}')$, whence

$$C\big(x, D \setminus (\mathcal{A} \cup \mathcal{A}')\big) \;=\; C\big(y, D \setminus (\mathcal{A} \cup \mathcal{A}')\big) . \tag{80.6}$$

Since $\mathcal{A}$ and $\mathcal{A}'$ are two aggregates in $D$ and $y \notin \mathcal{A} \cup \mathcal{A}'$, then it follows that

$$C(y, D) \;\subset\; D \setminus (\mathcal{A} \cup \mathcal{A}') ,$$

and this implies furthermore that

$$C(y, D) \;=\; C\big(y, D \setminus (\mathcal{A} \cup \mathcal{A}')\big) . \tag{80.7}$$

The equalities (80.6) and (80.7) imply that $C(x, D) = C(y, D) \subset \mathcal{A}''$. We have thus proved that $\mathcal{A}''$ is an aggregate in $D$. From the definitions (80.4) and (80.5), we see that

$$\mathcal{W}'' \;=\; \partial^{in}_{D\setminus\mathcal{A}}\big(\mathcal{A}'' \setminus \mathcal{A}\big) . \tag{80.8}$$

Let $x$ be a site of $\mathcal{W}''$. From formula (80.8), we see that the site $x$ is in $\mathcal{A}'' \setminus \mathcal{A}$ and there exists $y$ in $D \setminus (\mathcal{A} \cup \mathcal{A}' \cup \mathcal{A}'')$ such that $|x - y| = 1$. This implies that $x$ is in $\partial^{in}_D \mathcal{A}''$. Thus $\mathcal{W}'' \subset \partial^{in}_D \mathcal{A}''$ and $(\mathcal{A}'', \mathcal{W}'')$ is a waiting aggregate in $D$. We prove finally that $(\mathcal{A}', \mathcal{W}')$ and $(\mathcal{A}'', \mathcal{W}'')$ are coherent. Since $(\mathcal{A}, \mathcal{W})$ and $(\mathcal{A}', \mathcal{W}')$ are coherent, then the two aggregates $\mathcal{A}$ and $\mathcal{A}'$ are disjoint. Moreover, we have

$$\overline{\text{Clusters}}\,\big(\mathcal{N}(\mathcal{W}), D \setminus \mathcal{A}, \mathcal{A}'\big) \;\subset\; D \setminus (\mathcal{A} \cup \mathcal{A}') ,$$

whence

$$\mathcal{A}' \cap \overline{\text{Clusters}}\,\big(\mathcal{N}(\mathcal{W}), D \setminus \mathcal{A}, \mathcal{A}'\big) \;=\; \varnothing .$$

It follows directly from the definition (80.3) (or formula (80.5)) that $\mathcal{A}' \cap \mathcal{A}'' = \varnothing$. It remains to prove that

$$\partial^{in}_D\big(\mathcal{A}' \cup \mathcal{A}''\big) \;\subset\; \mathcal{W}' \cup \mathcal{W}'' . \tag{80.9}$$

Let $x$ be a site in $\partial^{in}_D\big(\mathcal{A}' \cup \mathcal{A}''\big)$. By definition, there exists $y$ in $D \setminus (\mathcal{A}' \cup \mathcal{A}'')$ such that $|x - y| = 1$. We consider several cases, and we use the same colors as in picture 148:

- $x \in \mathcal{A} \cup \mathcal{A}'$ . As $y$ is not in $\mathcal{A} \cup \mathcal{A}'$ , we see that $x$ belongs to $\partial^{in}_D\big(\mathcal{A} \cup \mathcal{A}'\ \big)$, which is itself included in $\mathcal{W} \cup \mathcal{W}'$ because $(\mathcal{A}, \mathcal{W})$ and $(\mathcal{A}', \mathcal{W}')$ are coherent.
  - ⋆ $x \in \mathcal{A}$. Since $\mathcal{W}' \cap \mathcal{A} = \varnothing$, then $x$ is in $\mathcal{W}$. This would imply that $y$ is in

$$\mathcal{N}(\mathcal{W})\setminus(\mathcal{A}'\cup\mathcal{A}'') \;\subset\; \mathcal{N}(\mathcal{W})\setminus(\mathcal{A}\cup\mathcal{A}') \;\subset\; \overline{\text{Clusters}}\,\big(\mathcal{N}(\mathcal{W}), D\setminus\mathcal{A}, \mathcal{A}'\big) \;\subset\; \mathcal{A}'' ,$$

    which is absurd. So this case cannot occur.
  - ⋆ $x \in \mathcal{A}'$. Since $\mathcal{W} \cap \mathcal{A}' = \varnothing$, then $x$ is in $\mathcal{W}'$.
- $x \in \mathcal{A}'' \setminus \mathcal{A} = \overline{\text{Clusters}}\,\big(\mathcal{N}(\mathcal{W}), D \setminus \mathcal{A}, \mathcal{A}'\big)$. The site $y$ is a neighbor of $x$ in $D \setminus (\mathcal{A}' \cup \mathcal{A}'')$, therefore $x$ is in $\mathcal{W}'' = \text{Waiting}\,(\mathcal{N}(\mathcal{W}), D \setminus \mathcal{A}, \mathcal{A}')$.

This completes the proof of the inclusion (80.9). □

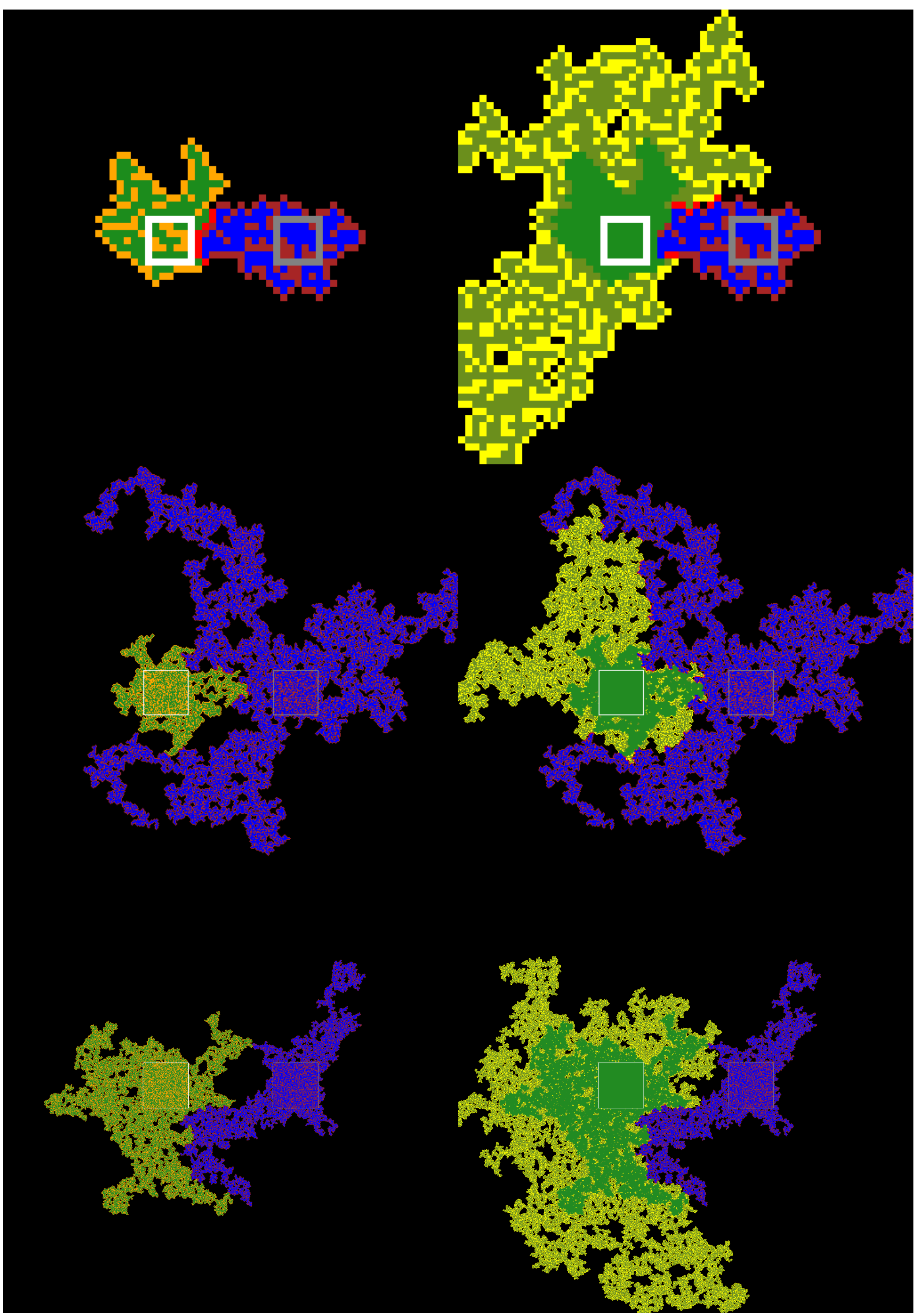

Figure 148: Three examples of the map Coherent_Site_Explore.
$\mathcal{A} = C(\text{white } \square, D)$, $\mathcal{W} = \partial^{in}\mathcal{A}$, $\mathcal{A}' = C(\text{gray } \square, D)$, $\mathcal{W}' = \partial^{in}\mathcal{A}'$,
$\mathcal{A}'' \setminus \mathcal{A} = \overline{\text{Clusters}}\,(\mathcal{W}, D \setminus \mathcal{A}, \mathcal{A}')$, $\mathcal{W}'' = \text{Waiting}\,(\mathcal{W}, D \setminus \mathcal{A}, \mathcal{A}')$.
Top picture: $D = \Lambda(64)$, $p = 0,58$. Middle picture: $D = \Lambda(512)$, $p = 0,583$.
Bottom picture: $D = \Lambda(1024)$, $p = 0,587$. $\mathcal{T}^+(0)$ on left, $\mathcal{T}_1^{-0}(1)$ on right.

## 80.4 Multiple intertwined explorations

The goal of this subsection is to define multiple intertwined explorations. In the case of the bond model, we defined only two intertwined explorations, and this was enough to attack the three challenging questions. However, to prove theorem 9.2, we need to extend the construction to an arbitrary number of intertwined explorations. We give first an adequate definition of coherence for $N \geq 2$ waiting aggregates.

**Definition 80.5.** *Let $N \geq 2$ be an integer. A $N$-tuple of waiting site aggregates*

$$(\mathcal{A}_i, \mathcal{W}_i, 1 \leq i \leq N)$$

*is said to be coherent if*

• *the aggregates are pairwise disjoint, i.e.,*

$$\forall i,j \in \{1,\dots,N\} \qquad i \neq j \quad \Longrightarrow \quad \mathcal{A}_i \cap \mathcal{A}_j = \varnothing\,;$$

• *the union of the waiting sets covers the inner boundary of the union of the aggregates, i.e.,*

$$\partial_D^{in}\Big(\bigcup_{1\leq i\leq N} \mathcal{A}_i\Big) \subset \bigcup_{1\leq i\leq N} \mathcal{W}_i\,.$$

To make some formulas more readable, we adopt a vertical notation: the waiting aggregate $(\mathcal{A}, \mathcal{W})$ might be presented in vector format, as $\binom{\mathcal{A}}{\mathcal{W}}$. We introduce a specific union operation for waiting aggregates, denoted by $\sqcup$, and defined as follows: for $N$ waiting aggregates $(\mathcal{A}_i, \mathcal{W}_i, 1 \leq i \leq N)$, we set

$$\bigsqcup_{1\leq i\leq N} \binom{\mathcal{A}_i}{\mathcal{W}_i} = \begin{pmatrix} \displaystyle\bigcup_{1\leq i\leq N} \mathcal{A}_i \\ \partial_D^{in}\Big(\displaystyle\bigcup_{1\leq i\leq N} \mathcal{A}_i\Big) \cap \Big(\displaystyle\bigcup_{1\leq i\leq N} \mathcal{W}_i\Big) \end{pmatrix}.$$

We take the intersection with the inner boundary for the waiting set in order to fulfil the requirement that the latter set is included in the former. A consequence of the definitions is the following stability property.

**Stability of coherence under $\bigsqcup$.** If $(\mathcal{A}_i, \mathcal{W}_i, 1 \leq i \leq N)$ is a coherent $N$-tuple of waiting site aggregates, and if $I_1, \dots, I_M$ is a partition of $\{1,\dots,N\}$, then the $M$-tuple of waiting site aggregates

$$\left(\bigsqcup_{i\in I_1} \binom{\mathcal{A}_i}{\mathcal{W}_i}, \dots, \bigsqcup_{i\in I_M} \binom{\mathcal{A}_i}{\mathcal{W}_i}\right)$$

is also coherent.

We shall take advantage of this stability property to construct the multiple intertwined explorations. A possibility would have been to extend the definition of the map Coherent_Site_Explore to $N$-tuples of waiting aggregates. However, we can bring ourselves back to the case of two explorers. We state next a preliminary result before proceeding to the construction.

**Proposition 80.6.** *Let $N \geq 2$ be an integer. Suppose that $(\mathcal{A}_i, \mathcal{W}_i, 1 \leq i \leq N)$ is a coherent $N$-tuple of waiting site aggregates. Let $i$ in $\{1, \dots, N\}$ and let*

$$\begin{pmatrix} \mathcal{A}'_i \\ \mathcal{W}'_i \end{pmatrix} = \mathit{Coherent_Site_Explore}\left( \begin{pmatrix} \mathcal{A}_i \\ \mathcal{W}_i \end{pmatrix}, \bigsqcup_{\substack{1 \leq j \leq N \\ j \neq i}} \begin{pmatrix} \mathcal{A}_j \\ \mathcal{W}_j \end{pmatrix} \right). \tag{80.10}$$

*The $N$-tuple of waiting site aggregates*

$$\begin{pmatrix} \mathcal{A}_1 \\ \mathcal{W}_1 \end{pmatrix}, \dots, \begin{pmatrix} \mathcal{A}_{i-1} \\ \mathcal{W}_{i-1} \end{pmatrix}, \begin{pmatrix} \mathcal{A}'_i \\ \mathcal{W}'_i \end{pmatrix}, \begin{pmatrix} \mathcal{A}_{i+1} \\ \mathcal{W}_{i+1} \end{pmatrix}, \dots, \begin{pmatrix} \mathcal{A}_n \\ \mathcal{W}_n \end{pmatrix}$$

*is a coherent $N$-tuple of waiting site aggregates.*

*Proof.* Thanks to the previous stability property, the two arguments of the map Coherent_Site_Explore in (80.10) form a pair of coherent waiting site aggregates. In addition, proposition 80.4 tells us that the two waiting site aggregates

$$\begin{pmatrix} \mathcal{A}'_i \\ \mathcal{W}'_i \end{pmatrix}, \qquad \bigsqcup_{\substack{1 \leq j \leq N \\ j \neq i}} \begin{pmatrix} \mathcal{A}_j \\ \mathcal{W}_j \end{pmatrix} \tag{80.11}$$

are coherent. This implies that

$$\mathcal{A}'_i \cap \Big( \bigcup_{1 \leq j \leq N,\, j \neq i} \mathcal{A}_j \Big) = \varnothing .$$

From the coherence of the $N$ waiting aggregates $(\mathcal{A}_j, 1 \leq j \leq N, j \neq i)$, we know that they are pairwise disjoint, and we conclude that the first condition for coherence is satisfied. We know also that

$$\partial_D^{in} \Big( \bigcup_{1 \leq j \leq N,\, j \neq i} \mathcal{A}_j \Big) \subset \bigcup_{1 \leq j \leq N,\, j \neq i} \mathcal{W}_j . \tag{80.12}$$

The coherence of the two waiting aggregates in (80.11) yields that

$$\partial_D^{in} \Big( \mathcal{A}'_i \cup \bigcup_{\substack{1 \leq j \leq N \\ j \neq i}} \mathcal{A}_j \Big) \subset \mathcal{W}'_i \cup \left( \partial_D^{in} \Big( \bigcup_{\substack{1 \leq j \leq N \\ j \neq i}} \mathcal{A}_j \Big) \cap \bigcup_{\substack{1 \leq j \leq N \\ j \neq i}} \mathcal{W}_j \right). \tag{80.13}$$

The second condition for the coherence follows from (80.12) and (80.13). □

We combine now several taboo explorations to create the multiple intertwined explorations. In the case of two starting sets, the construction is the same as in subsection 29.4 for the bond model. The only difference is that we use $\overline{\text{Clusters}}(\cdot, \cdot, \cdot)$ instead of $\text{Clusters}(\cdot, \cdot, \cdot)$ and Coherent_Site_Explore instead of Coherent_Bond_Explore. All the comments of the construction for the bond model apply as well for the site model. The two constructions might certainly be unified, indeed the proofs differ at some minor points. Yet these constructions are cornerstones of the proofs. To improve readability, we prefer to write the details completely for each model, at the expense of some additional length. For the bond model, we focused on disconnection events involving only two sets (typically opposite faces of a box or a prism), and we needed only to consider two intertwined explorations. This is not any more the case when we try to prove the theorem 9.2, we will consider disconnection events involving several large clusters, and therefore we need to construct multiple intertwined explorations. There is a further difference concerning the termination time. Whereas in the bond model the explorations continued until there was nothing left to explore, we consider here a different termination time, and we will stop the explorations when a specific region has been completely explored. Naturally these extensions require new notation, which we introduce next. We define first what we call the standard framework for the intertwined explorations.

**Definition 80.7** (Standard framework.)**.** *Let $D$ be a finite connected subset of $\mathbb{Z}^d$, Let $N \geq 2$ and let $w_1, \dots, w_N$ be $N$ distinct sites of $D$. Let $G$ be a subset of $D$. We consider a site percolation configuration in $D$ such that the sites $w_1, \cdots, w_N$ are pairwise disconnected.*

The aim is to start from the sites $w_1, \cdots, w_N$ and to explore the entire region $G$. Each site will be assigned to one of $N$ explorers, who will explore the percolation configuration together according to a precise set of rules. There is a lot of freedom for the choice of these rules. We present here the simplest scheme:

• there is fixed deterministic order on the $N$ explorers;
• the explorers are called in turn, according to this order, during an inner loop;
• the inner loop is repeated until the region $G$ has been fully explored.

When called upon, an explorer restarts its exploration from its waiting set, with taboo set the territory already explored by all the other explorers. Mathematically, we shall construct $N$ sequences of random waiting aggregates

$$\big(\mathcal{A}_i(t), \mathcal{W}_i(t)\big), \quad 1 \leq i \leq N\,, \quad 0 \leq t \leq T\,,$$

whose meaning is the following: for $i$ in $\{\,1, \dots, N\,\}$,

• $\mathcal{A}_i(t)$ is the set of the sites explored by the explorer $i$ after iteration $t$;
• $\mathcal{W}_i(t)$ is the waiting set of the explorer $i$ after iteration $t$.

We give next the full mathematical description of the $N$ sequences of random waiting aggregates.

**Initialization.** We define first the $N$ initial sets at $t = 0$. We set

$$\mathcal{A}_1(0) \,=\, \overline{\text{Clusters}}\big(\{\,w_1\,\}, D, \varnothing\big)\,, \quad \mathcal{W}_1(0) \,=\, \partial_D^{in} \mathcal{A}_1(0)\,.$$

Let $i$ be an index in $\{\,2, \dots, N\,\}$ and suppose that $\mathcal{A}_j(0)$, $\mathcal{W}_j(0)$, $1 \leq j < i$, have been defined. The sets $\mathcal{A}_i(0)$, $\mathcal{W}_i(0)$ are defined as

$$\mathcal{A}_i(0) \,=\, \overline{\text{Clusters}}\Big(\{\,w_i\,\}, D, \bigcup_{1 \leq j < i} \mathcal{A}_j(0)\Big)\,,$$

$$\mathcal{W}_i(0) \,=\, \partial_D^{in} \mathcal{A}_i(0) \setminus \bigcup_{1 \leq j < i} \mathcal{A}_j(0)\,. \tag{80.14}$$

We then define the next terms of the sequences using two nested iterations: the outer iteration is applied to times $t \geq 0$, and the inner iteration is applied to indices $i$ in $\{1, \dots, N\}$.

**Outer iteration on $t$.** Let $t \geq 0$ be fixed and suppose that the sets $\mathcal{A}_i(t), \mathcal{W}_i(t)$, $1 \leq i \leq N$, have been defined. We define the sets at time $t+1$ through the following inner iteration.

**Inner iteration on $i$.** Let $i$ be an index in $\{1, \dots, N-1\}$ and suppose that the sets $\mathcal{A}_j(t+1)$, $\mathcal{W}_j(t+1)$, $1 \leq j < i$, have been defined. We define the sets $\mathcal{A}_i(t+1)$, $\mathcal{W}_i(t+1)$ by setting

$$\begin{pmatrix} \mathcal{A}_i(t+1) \\ \mathcal{W}_i(t+1) \end{pmatrix} = \text{Coherent_Site_Explore}\left( \begin{pmatrix} \mathcal{A}_i(t) \\ \mathcal{W}_i(t) \end{pmatrix}, \bigsqcup_{1\leq j<i} \begin{pmatrix} \mathcal{A}_j(t+1) \\ \mathcal{W}_j(t+1) \end{pmatrix} \sqcup \bigsqcup_{i<k\leq N} \begin{pmatrix} \mathcal{A}_k(t) \\ \mathcal{W}_k(t) \end{pmatrix} \right). \tag{80.15}$$

**Termination.** The inner iteration always terminates at $N$. The outer iteration terminates at the random time $T$, one iteration after the region $G$ has been fully explored, i.e.,

$$T = \min\left\{ t \geq 1 : G \subset \bigcup_{1\leq i\leq N} \mathcal{A}_i(t) \right\} + 1 .$$

The $+1$ after the minimum ensures that, upon termination, there are no more waiting sites in $G$, i.e.,

$$G \cap \bigcup_{1\leq i\leq N} \mathcal{W}_i(T) = \varnothing .$$

This completes the construction of the multiple intertwined explorations. Naturally, these objects can equivalently be defined with the help of an algorithm. The corresponding pseudocode is presented in algorithm 80.2. However, for the subsequent analysis, we will rely on the previous mathematical construction.

The mechanism at the heart of the intertwined explorations is encoded in formula (80.15). Let us try to reformulate it by returning to the definition of Coherent_Site_Explore. We define first the taboo set

$$\mathcal{T}_i(t+1) = \mathcal{A}_1(t+1) \cup \cdots \cup \mathcal{A}_{i-1}(t+1) \cup \mathcal{A}_{i+1}(t) \cup \cdots \cup \mathcal{A}_N(t) . \tag{80.16}$$

The sets $\mathcal{A}_i(t+1)$ and $\mathcal{W}_i(t+1)$ are then given by the following formulas:

$$\begin{aligned} \mathcal{A}_i(t+1) &= \mathcal{A}_i(t) \cup \overline{\text{Clusters}}\big(\mathcal{N}(\mathcal{W}_i(t)), D \setminus \mathcal{A}_i(t), \mathcal{T}_i(t+1)\big) , \qquad (80.17) \\ \mathcal{W}_i(t+1) &= \text{Waiting}\big(\mathcal{N}(\mathcal{W}_i(t)), D \setminus \mathcal{A}_i(t), \mathcal{T}_i(t+1)\big) . \end{aligned}$$

Using the definition (80.2) of the map Waiting and the formula (80.17), the waiting set $\mathcal{W}_i(t+1)$ can alternatively be defined directly from $\mathcal{A}_i(t+1)$ as expressed in the following formula:

$$\mathcal{W}_i(t+1) = \partial^{in}_{D\setminus\mathcal{A}_i(t)} \mathcal{A}_i(t+1) . \tag{80.18}$$

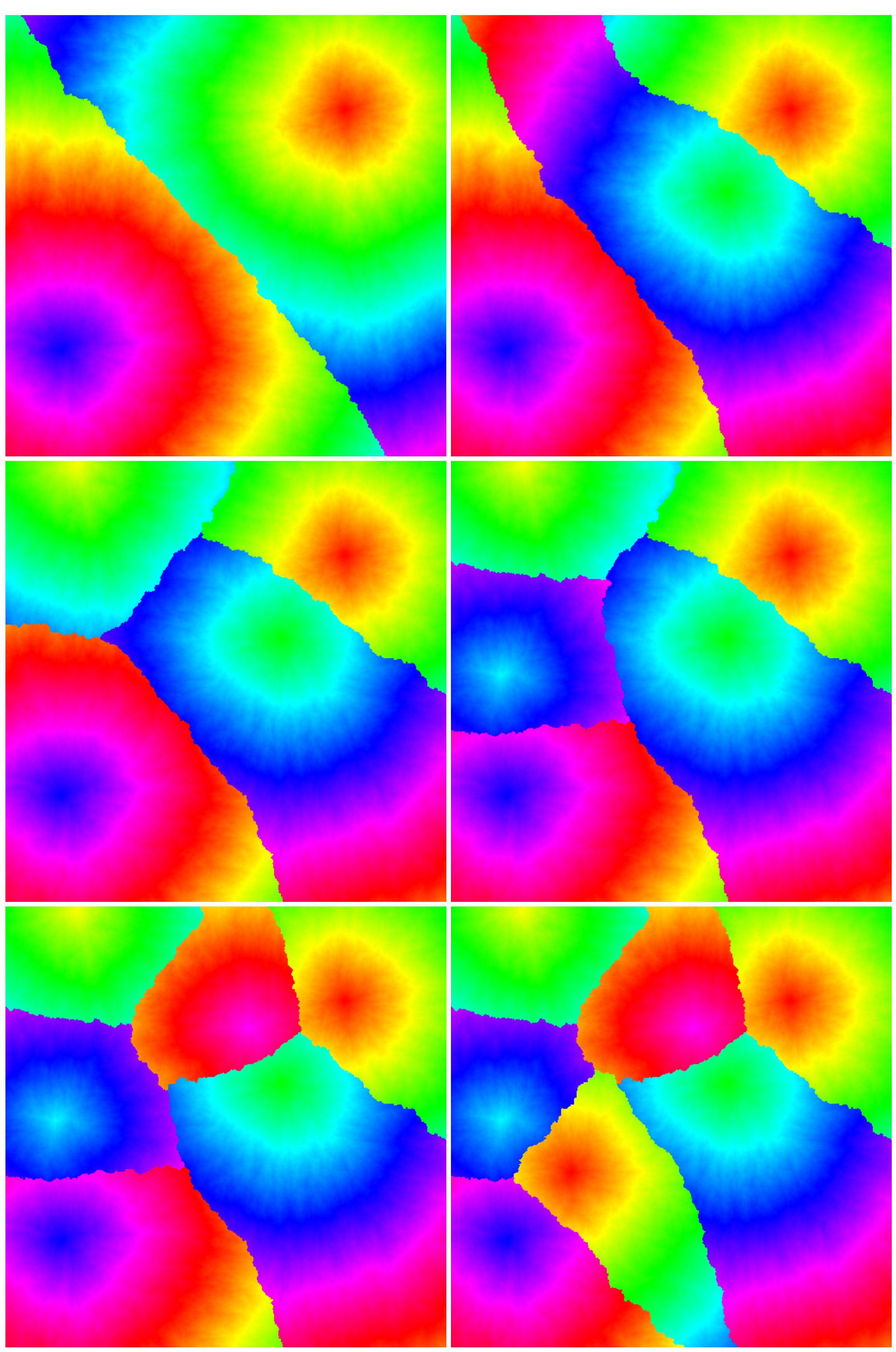

Figure 149: Examples of intertwined explorations in $\Lambda(512)$, $p = 0.2$.
Same configuration, from two explorers (top left) to seven (bottom right).

**Algorithm 80.2** Intertwined_Site_Exploration $(w_1, \dots, w_N, D, G)$

**Require:** $N$ non-empty starting sites $w_1, \dots, w_N$,
an authorized domain $D$, a goal region $G$

**Ensure:** $T$ is the termination time of the multiple intertwined explorations, $\mathcal{A}_1(T), \cdots, \mathcal{A}_N(T)$ are the sets visited by the $N$ explorers

**for all** $x \in D$ **do** ▷ Initialization
  visited$[x] \leftarrow$ false
**end for**
$t \leftarrow -1$
**repeat** ▷ Outer loop
  $t \leftarrow t+1$
  **for** $i = 1$ **to** $N$ **do**
    **if** $t = 0$ **then**
      $\mathcal{T} \leftarrow \mathcal{A}_1(0) \cup \cdots \cup \mathcal{A}_{i-1}(0)$
      $A \leftarrow \{ w_i \}$
    **else**
      $\mathcal{T} \leftarrow \mathcal{A}_1(t) \cup \cdots \cup \mathcal{A}_{i-1}(t) \cup \mathcal{A}_{i+1}(t-1) \cup \cdots \cup \mathcal{A}_N(t-1)$
      $A \leftarrow \mathcal{N}(\mathcal{W}_i(t-1)) \cap D \setminus \mathcal{T}$
    **end if**
    $\overline{C} \leftarrow \varnothing, W \leftarrow \varnothing$
    **repeat** ▷ Inner loop
      Pick up $x \in A$
      $A \leftarrow A \setminus \{x\}$
      **if** not visited$[x]$ **then**
        visited$[x] \leftarrow$ true
        $\overline{C} \leftarrow \overline{C} \cup \{ x \}$
        **if** $x$ is closed **then**
          $W \leftarrow W \cup \{ x \}$
        **else if** $x$ is open **then**
          $A \leftarrow A \cup \big(\mathcal{N}(x) \cap D \setminus \mathcal{T}\big)$
        **end if**
      **end if**
    **until** $A = \varnothing$
    $\mathcal{A}_i(t) \leftarrow \mathcal{A}_i(t-1) \cup \overline{C}$
    $\mathcal{W}_i(t) \leftarrow W$ ▷ See the remark below
  **end for**
**until** $G \subset \mathcal{A}_1(t) \cup \cdots \cup \mathcal{A}_N(t)$ ▷ $G$ is fully explored
$T \leftarrow t$
**return** $T+1$, $\mathcal{A}_1(T), \cdots, \mathcal{A}_N(T)$

**Remark on $\mathcal{W}_i(t)$:** In the pseudocode, the waiting set $\mathcal{W}_i(t)$ is defined as $W$, which is the set of the closed sites of $\mathcal{A}_i(t) \setminus \mathcal{A}_i(t-1)$. If we stick to the mathematical definition (80.18), we should have set $\mathcal{W}_i(t) \leftarrow W \cap \partial^{in}_{D \setminus \mathcal{A}_i(t-1)} \mathcal{A}_i(t)$. However, in the inner loop, the sites of $W \cap \big(\mathcal{A}_i(t-1) \cup \mathcal{T}\big)$ are ignored as they are recorded in the array visited$[\cdot]$, and the neighbours of the remaining sites will be duly examined. This choice bypasses the computation of $\partial^{in}_{D \setminus \mathcal{A}_i(t-1)} \mathcal{A}_i(t)$.

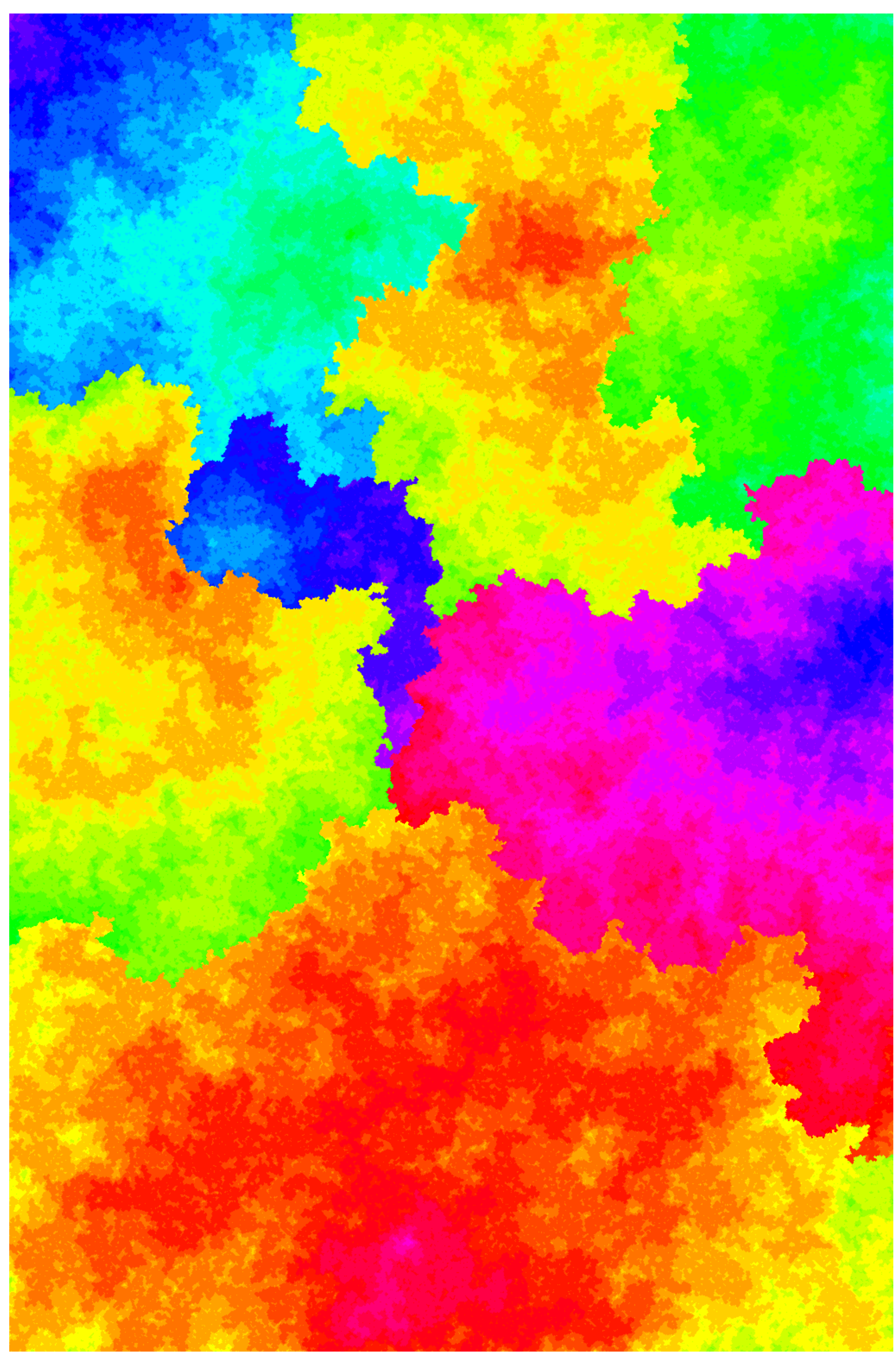

Figure 150: 7 intertwined explorations, $1550 \times 1024$, site percolation, $p = 0.57$.

## 80.5 Termination

We prove here the counterpart for the site model of proposition 29.6. The role of the taboo sets used in the multiple intertwined explorations is to prevent each explorer from entering into the set of the sites already discovered by the other explorers. Although each explorer forbids the others to penetrate his own territory, the union of the sets discovered by all the explorers coincides with the set that a single non-taboo exploration would visit if it were started from the union of the initial sites of the explorers. Our next goal is to give a precise formulation of this crucial result. We derive then two important consequences. The first one is that the multiple intertwined explorations terminate after a finite number of iterations. The second one is an explicit upper bound on the termination time, expressed in terms of a maximum of travel times within $D$. Throughout this subsection, we consider the standard framework of definition 80.7 and the associated $N$ intertwined explorations $(\mathcal{A}_i(t),\mathcal{W}_i(t),1\leq i\leq N,0\leq t\leq T)$. For convenience, we recall the definitions (12.6), (12.7) of the ball $\mathcal{B}(A,D,t)$ and the shell $\text{Shell}\,(A,D,t)$:

$$\begin{aligned}\forall t\geq 0\quad\forall A\subset D\qquad \mathcal{B}(A,D,t)\,&=\,\big\{\,x\in D:T_D(A,x)\leq t\,\big\}\,,\\ \text{Shell}\,(A,D,t)\,&=\,\big\{\,x\in D:T_D(A,x)=t\,\big\}\,.\end{aligned}$$

To alleviate the formulas, we adopt the following convention. When a map takes a set as an argument, and this set is reduced to a singleton $\{\,w\,\}$, we write simply $w$ instead of $\{\,w\,\}$. For instance, we write $\mathcal{B}(w,D,t)$, $\text{Shell}\,(w,D,t)$. This will help a lot, because each explorer will start from a single site. We begin with the following simple lemma.

**Lemma 80.8.** *For any $i$ in $\{\,1,\dots,N\,\}$ and $t$ in $\{\,0,\dots,T\,\}$, the travel time within $\mathcal{A}_i(t)$ between $w_i$ and any site in $\mathcal{A}_i(t)$ is less than or equal to $t$, so that*

$$\mathcal{A}_i(t)\,=\,\mathcal{B}\big(w_i,\mathcal{A}_i(t),t\big)\,.$$

*Proof.* We prove the result by induction on $t$. For $t=0$, it is an immediate consequence of the definition (80.14) of the $N$ initial sets. Let $t\geq 0$ and suppose that the result has been proved at rank $t$. Let $i$ belong to $\{\,1,\dots,N\,\}$. It follows from formulas (80.17) and (80.1) that

$$\begin{aligned}\mathcal{A}_i(t+1)\setminus\mathcal{A}_i(t)\,&=\,\overline{\text{Clusters}}\,\big(\mathcal{N}(\mathcal{W}_i(t)),D\setminus\mathcal{A}_i(t),\mathcal{T}_i(t+1)\big)\\ &\subset\,\overline{\text{Clusters}}\,\big(\mathcal{N}(\mathcal{W}_i(t)),\mathcal{A}_i(t+1),\varnothing\big)\subset\,\text{Shell}\,\big(\mathcal{A}_i(t),\mathcal{A}_i(t+1),1\big)\,.\end{aligned}\tag{80.19}$$

The induction hypothesis at rank $t$ yields that

$$\mathcal{A}_i(t)\,=\,\mathcal{B}\big(w_i,\mathcal{A}_i(t),t\big)\,\subset\,\mathcal{B}\big(w_i,\mathcal{A}_i(t+1),t\big)\,.\tag{80.20}$$

It follows from the two inclusions (80.19) and (80.20) that

$$\mathcal{A}_i(t{+}1)\subset\mathcal{A}_i(t)\cup\text{Shell}\,\big(\mathcal{B}(w_i,\mathcal{A}_i(t{+}1),t),\mathcal{A}_i(t{+}1),1\big)\subset\mathcal{B}\big(w_i,\mathcal{A}_i(t{+}1),t{+}1\big).$$

The converse inclusion is obvious, and the result holds at rank $t+1$. □

To shorten the formulas, we define the set $W$ of the initial sites as

$$W = \{ w_1, \dots, w_N \} .$$

The main result of this subsection is stated in the next proposition.

**Proposition 80.9.** *We consider the standard framework of definition 80.7 and the $N$ intertwined explorations $(\mathcal{A}_i(t), \mathcal{W}_i(t), 1 \leq i \leq N, 0 \leq t \leq T)$. We have*

$$\forall t \in \{ 0, \dots, T \} \qquad \bigcup_{1 \leq i \leq N} \mathcal{A}_i(t) = \mathcal{B}(W, D, t) . \tag{80.21}$$

Since the successive shells $\text{Shell}\,(W, D, s)$, $0 \leq s \leq t$, form a partition of the ball $\mathcal{B}(W, D, t)$, the result (80.21) can equivalently be reformulated as

$$\forall t \in \{ 0, \dots, T \} \qquad \text{Shell}\,(W, D, t) = \bigcup_{1 \leq i \leq N} \mathcal{A}_i(t) \setminus \bigcup_{1 \leq i \leq N} \mathcal{A}_i(t-1) . \tag{80.22}$$

In the above formula, we make the convention that $\mathcal{A}_i(-1) = \varnothing$ for $1 \leq i \leq N$.

*Proof.* For $t = 0$, the definition (80.14) of the sets $\mathcal{A}_i(0)$, $1 \leq i \leq N$, yields that

$$\begin{aligned} \bigcup_{1 \leq i \leq N} \mathcal{A}_i(0) &= \bigcup_{1 \leq i \leq N} \overline{\text{Clusters}} \Big( \{ w_i \}, D, \bigcup_{1 \leq j < i} \mathcal{A}_j(0) \Big) \\ &\subset \bigcup_{1 \leq i \leq N} \overline{C}\big(w_i, D\big) = \overline{C}\big(W, D\big) = \text{Shell}\,(W, D, 0) . \end{aligned}$$

We prove next the converse inclusion. Let $x$ belong to the set $\text{Shell}\,(W, D, 0)$. By definition, there exists a path $z_0, z_1, \dots, z_r$ in $D$ joining a site $z_0 = w_i$ of $W$ to $z_r = x$ such that all the sites $z_0, \dots, z_{r-1}$ are open. If $x$ is also open, then $x$ is in $C\big(w_i, D\big) \subset \mathcal{A}_i(0)$. If $x$ is closed, then either $x$ is in $\bigcup_{1 \leq j < i} \mathcal{A}_j(0)$ or in

$$\mathcal{A}_i(0) = \overline{\text{Clusters}} \Big( \{ w_i \}, D, \bigcup_{1 \leq j < i} \mathcal{A}_j(0) \Big) .$$

We proceed next by induction on $t$ to prove the formula (80.22). Suppose that the result has been proved until rank $t$ for some $t \geq 0$. By formula (80.17), for $i$ in $\{ 1, \dots, N \}$, we have the inclusion

$$\begin{aligned} \mathcal{A}_i(t+1) \setminus \mathcal{A}_i(t) &\subset \overline{\text{Clusters}} \big( \mathcal{N}(\mathcal{W}_i(t)), D \setminus \mathcal{A}_i(t), \mathcal{T}_i(t+1) \big) \\ &\subset \overline{\text{Clusters}} \big( \mathcal{N}(\mathcal{W}_i(t)), D \setminus \mathcal{A}_i(t), \varnothing \big) \subset \text{Shell}\,\big( \mathcal{A}_i(t), D, 1 \big) . \end{aligned} \tag{80.23}$$

In the sequel of the proof, we denote by $\mathcal{A}(t)$ the set $\bigcup_{1 \leq i \leq N} \mathcal{A}_i(t)$. Using the induction hypothesis, we have

$$\mathcal{A}(t) = \mathcal{B}(W, D, t) . \tag{80.24}$$

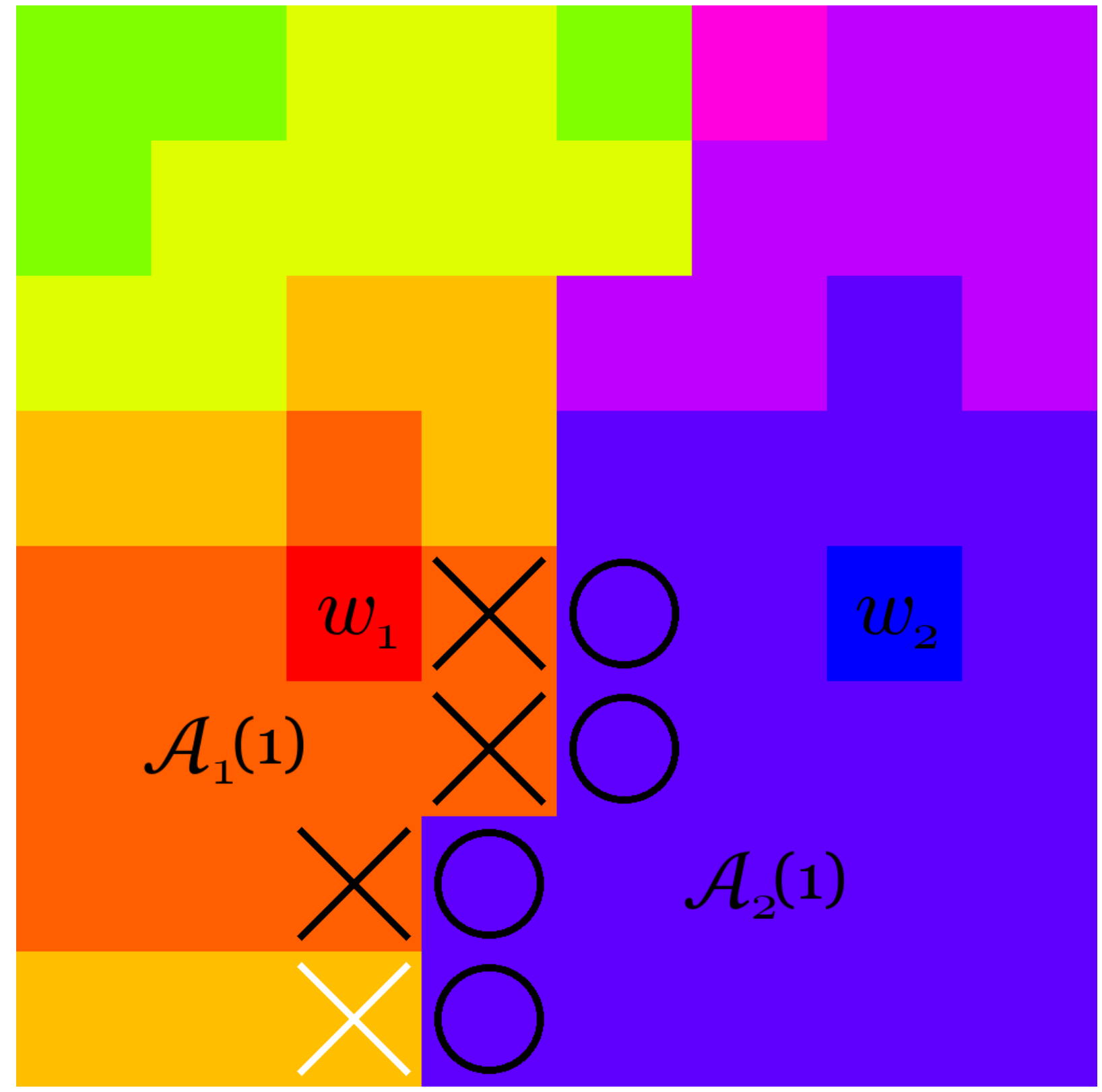

Figure 151: The explorations started from $w_1$, $w_2$ in the box $\Lambda(8)$ with $p = 0.3$: $\mathcal{A}_1(0)$, $\mathcal{A}_1(1) \setminus \mathcal{A}_1(0)$, $\mathcal{A}_1(2) \setminus \mathcal{A}_1(1)$, $\mathcal{A}_1(3) \setminus \mathcal{A}_1(2)$, $\mathcal{A}_1(4) \setminus \mathcal{A}_1(3)$, $\mathcal{A}_2(0)$, $\mathcal{A}_2(1) \setminus \mathcal{A}_2(0)$, $\mathcal{A}_2(2) \setminus \mathcal{A}_2(1)$, $\mathcal{A}_2(3) \setminus \mathcal{A}_2(2)$.
Crosses, $w_1$, $w_2$ are closed, circles are open. The white cross is in $\overline{\text{Clusters}}\big(\mathcal{N}(\mathcal{W}_i(1)), D, \mathcal{A}(1)\big)$ for $i = 1, 2$, but it is not in $\overline{\text{Clusters}}\big(\mathcal{N}(\mathcal{W}_2(1)), D \setminus \mathcal{A}_2(1), \mathcal{T}_2(2)\big)$.

The equality (80.24) and the inclusion (80.23) together yield that

$$\begin{aligned} \bigcup_{1 \leq i \leq N} \mathcal{A}_i(t+1) \ &\subset \ \mathcal{A}(t) \cup \bigcup_{1 \leq i \leq N} \text{Shell}\big(\mathcal{A}_i(t), D, 1\big) \\ &\subset \bigcup_{0 \leq s \leq t+1} \text{Shell}\,(W, D, s) \ = \ \mathcal{B}(W, D, t+1)\,. \end{aligned}$$

We prove next the converse inclusion. Using again (80.24), we have

$$\begin{aligned} \text{Shell}\,(W, D, t+1) \ &= \ \text{Shell}\big(\mathcal{A}(t), D, 1\big) \ = \ \overline{C}\Big(\mathcal{N}\big(\partial_D^{in}\mathcal{A}(t)\big), D \setminus \mathcal{A}(t)\Big) \\ &\subset \overline{C}\bigg(\mathcal{N}\Big(\bigcup_{1 \leq i \leq N} \mathcal{W}_i(t)\Big), D \setminus \mathcal{A}(t)\bigg) \ \subset \ \bigcup_{1 \leq i \leq N} \overline{C}\Big(\mathcal{N}\big(\mathcal{W}_i(t)\big), D \setminus \mathcal{A}(t)\Big)\,. \end{aligned}$$

The penultimate inclusion above stems from the fact that, thanks to proposition 80.6, the $N$ waiting aggregates of the sequence at time $t$ are coherent. Suppose that $x$ is in $\text{Shell}\,(W, D, t+1)$. Using the above inclusion, we see that there exists $i$ in $\{ 1, \dots, N \}$ such that

$$x \in \overline{C}\big(\mathcal{N}\big(\mathcal{W}_i(t)\big), D \setminus \mathcal{A}(t)\big) \;=\; \overline{\text{Clusters}}\,\big(\mathcal{N}(\mathcal{W}_i(t)), D, \mathcal{A}(t)\big)\,. \tag{80.25}$$

Notice that there might exist several indices $i$ satisfying (80.25) (see figure 151). Let $i^*$ be the smallest index $i$ in $\{ 1, \dots, N \}$ such that (80.25) holds. By definition, there exists a path $z_0, \dots, z_r$ in $D \setminus \mathcal{A}(t)$ joining a site $z_0 = y$ of $\mathcal{N}\big(\mathcal{W}_{i^*}(t)\big)$ to $z_r = x$ such that the sites $z_0, \dots, z_{r-1}$ are open. Since the path $z_0, \dots, z_r$ is included in $D \setminus \mathcal{A}(t)$, then

$$\{ z_0, \dots, z_r \} \;\subset\; \overline{\text{Clusters}}\,\big(\mathcal{N}(\mathcal{W}_{i^*}(t)), D, \mathcal{A}(t)\big)\,.$$

We claim that $x$ is not in

$$\mathcal{T}_{i^*}(t+1) \;=\; \mathcal{A}_1(t+1) \cup \cdots \cup \mathcal{A}_{i^*-1}(t+1) \cup \mathcal{A}_{i^*+1}(t) \cup \cdots \cup \mathcal{A}_N(t)\,.$$

We already know that $x$ is not in $\mathcal{A}(t)$. Suppose that $x$ is in $\mathcal{A}_j(t+1) \setminus \mathcal{A}(t)$ for some $j < i^*$. By definition, there would exist a path $z'_0, \dots, z'_s$ in

$$D \setminus \big(\mathcal{A}_j(t) \cup \mathcal{T}_j(t+1)\big)$$

joining a site $z'_0$ of $\mathcal{N}\big(\mathcal{W}_j(t)\big)$ to $z'_s = x$ such that the sites $z'_0, \dots, z'_{s-1}$ are open. Thus we would have

$$\{ z'_0, \dots, z'_s \} \;\subset\; \overline{\text{Clusters}}\,\big(\mathcal{N}(\mathcal{W}_j(t)), D, \mathcal{A}(t)\big)\,,$$

and this would contradict the minimality of $i^*$. Furthermore, we claim that no site of the path $z_0, \dots, z_{r-1}$ is in $\mathcal{T}_{i^*}(t+1)$. Let us fix $j$ in $\{ 1, \dots, i^* - 1 \}$ and suppose that one site among $\{ z_0, \dots, z_{r-1} \}$ belongs to $\mathcal{A}_j(t+1)$. Since $\mathcal{A}_j(t+1)$ is an aggregate, and all the sites of the path $\{ z_0, \dots, z_{r-1} \}$ are open, then they all belong to $\mathcal{A}_j(t+1)$. Moreover, they are not in $\mathcal{A}(t)$, which contains $\mathcal{A}_j(t)$, thus they are in the increment

$$\mathcal{A}_{j+1}(t) \setminus \mathcal{A}_j(t) \;=\; \overline{\text{Clusters}}\,\big(\mathcal{N}(\mathcal{W}_j(t)), D \setminus \mathcal{A}_j(t), \mathcal{T}_j(t+1)\big)\,. \tag{80.26}$$

It follows from the definition (80.16) of the set $\mathcal{T}_j(t+1)$ that

$$\overline{\text{Clusters}}\,\big(\mathcal{N}(\mathcal{W}_j(t)), D \setminus \mathcal{A}_j(t), \mathcal{T}_j(t+1)\big) \subset \overline{\text{Clusters}}\,\big(\mathcal{N}(\mathcal{W}_j(t)), D, \mathcal{A}(t)\big). \tag{80.27}$$

Since $x$ is not in $\mathcal{A}(t) \cup \mathcal{T}_{i^*}(t+1)$, which contains $\mathcal{A}(t) \cup \mathcal{T}_j(t+1)$, then $x$, being a neighbour of $z_{r-1}$, is also in the set (80.26), and it would follow from (80.27) that

$$\{ z_0, \dots, z_r \} \;\subset\; \overline{\text{Clusters}}\,\big(\mathcal{N}(\mathcal{W}_j(t)), D, \mathcal{A}(t)\big)\,.$$

This would again contradict the minimality of $i^*$. So we have proved that none of the sites among $\{ z_0, \dots, z_{r-1} \}$ is in $\mathcal{A}_j(t+1)$. This is true for any $j < i^*$, hence no site of the path $z_0, \dots, z_r$ is in $\mathcal{T}_{i^*}(t+1)$. Thus the site $x$ is in fact in $\overline{\text{Clusters}}\,\big(\mathcal{N}(\mathcal{W}_{i^*}(t)), D \setminus \mathcal{A}_{i^*}(t), \mathcal{T}_{i^*}(t+1)\big)$, which is included in $\mathcal{A}_{i^*}(t+1)$. We have proved that $\text{Shell}\,(W, D, t+1)$ is included in $\bigcup_{1 \le i \le N} \mathcal{A}_i(t+1)$. This concludes the induction step and the proof of (80.21). $\square$

With the help of proposition 80.9, we can obtain a control on the termination time $T$ in a finite connected domain $D$. The identity (80.21) implies that, for any $t \geq 0$, the union $\mathcal{A}_1(t) \cup \cdots \cup \mathcal{A}_N(t)$ contains all the sites that are within travel distance from $W$ in $D$ less than or equal to $t$. This yields the following upper bound on $T$.

**Corollary 80.10.** *The termination time $T$ is less than or equal to one plus the maximal travel distance between the set $W$ and a site of $G$, i.e.,*

$$T \;\leq\; 1 + \max\,\big\{\, T_D(W,x) : x \in G \,\big\}\,.$$

We state next another interesting consequence of proposition 80.9 concerning the structure of the sets $\mathcal{A}_1(t), \dots, \mathcal{A}_N(t)$, $0 \leq t \leq T$.

**Corollary 80.11.** *Let us fix $i$ in $\{\,1,\dots,N\,\}$, and let us define*

$$\widehat{\mathcal{A}}_i(T) \;=\; \bigcup_{1\leq j\leq N,\, j\neq i} \mathcal{A}_j(T)\,. \tag{80.28}$$

*For any $t$ in $\{\,0,\dots,T\,\}$, we have*

$$\mathcal{A}_i(t) \;=\; \mathcal{B}\big(w_i, D\setminus \widehat{\mathcal{A}}_i(T), t\big)\,. \tag{80.29}$$

*Proof.* Let us fix $i$ in $\{\,1,\dots,N\,\}$. For $t=0$, the formula (80.29) can be checked directly, and for $t=T$, it is a consequence of lemma 80.8. We consider next an integer $t$ in $\{\,1,\dots,T-1\,\}$. As $\mathcal{A}_i(t)$ is included in $D\setminus\widehat{\mathcal{A}}_i(T)$, it follows from lemma 80.8 that

$$\mathcal{A}_i(t) \;=\; \mathcal{B}\big(w_i, \mathcal{A}_i(t), t\big) \;\subset\; \mathcal{B}\big(w_i, D\setminus \widehat{\mathcal{A}}_i(T), t\big)\,. \tag{80.30}$$

By proposition 80.9, for any $s$ in $\{\,1,\dots,T\,\}$, we have

$$\mathcal{A}_i(s)\setminus\mathcal{A}_i(s-1) \;=\; \mathcal{A}_i(s)\;\setminus \bigcup_{1\leq j\leq N} \mathcal{A}_j(s-1) \;\subset\; \text{Shell}\,(W,D,s)\,. \tag{80.31}$$

We decompose $D\setminus\mathcal{A}_i(t)$ as

$$D\setminus\mathcal{A}_i(t) \;=\; \Big(\bigcup_{t<s\leq T}\big(\mathcal{A}_i(s)\setminus\mathcal{A}_i(s-1)\big)\Big)\cap\big(D\setminus\mathcal{A}_i(T)\big)\,,$$

and we write

$$\begin{aligned}
\mathcal{B}\big(w_i, D\setminus \widehat{\mathcal{A}}_i(T), t\big)&\setminus\mathcal{A}_i(t)\\
&= \bigcup_{t<s\leq T}\Big(\mathcal{B}\big(w_i, D\setminus \widehat{\mathcal{A}}_i(T), t\big)\cap\big(\mathcal{A}_i(s)\setminus\mathcal{A}_i(s-1)\big)\Big)\\
&\qquad\qquad \cup\Big(\mathcal{B}\big(w_i, D\setminus \widehat{\mathcal{A}}_i(T), t\big)\cap\big(D\setminus\mathcal{A}_i(T)\big)\Big)\,.
\end{aligned} \tag{80.32}$$

To control the term inside the big union in (80.32), we use the inclusion (80.31): for $s$ such that $t<s\leq T$, we have

$$\begin{aligned}
\mathcal{B}\big(w_i, D\setminus \widehat{\mathcal{A}}_i(T), t\big)&\cap\big(\mathcal{A}_i(s)\setminus\mathcal{A}_i(s-1)\big)\\
&\subset\; \mathcal{B}\big(w_i, \mathcal{A}_i(T), t\big)\;\cap\;\text{Shell}\,(W,D,s)\;=\;\varnothing\,.
\end{aligned}$$

For the last term in (80.32), we have

$$\mathcal{B}\big(w_i, D\setminus\widehat{\mathcal{A}}_i(T), t\big)\cap\big(D\setminus\mathcal{A}_i(T)\big)\ \subset\ \mathcal{B}\big(W, D, t\big)\cap\Big(D\setminus\big(\widehat{\mathcal{A}}_i(T)\cup\mathcal{A}_i(T)\big)\Big)\,. \tag{80.33}$$

By proposition 80.9, we have also

$$\widehat{\mathcal{A}}_i(T)\cup\mathcal{A}_i(T)\ =\ \bigcup_{1\leq i\leq N}\mathcal{A}_i(T)\ =\ \mathcal{B}\big(W, D, T\big)\,. \tag{80.34}$$

Together, formulas (80.32), (80.33) and (80.34) imply that

$$\mathcal{B}\big(w_i, D\setminus\widehat{\mathcal{A}}_i(T), t\big)\setminus\mathcal{A}_i(t)\ =\ \varnothing\,. \tag{80.35}$$

The equality (80.29) follows from formulas (80.30) and (80.35). □

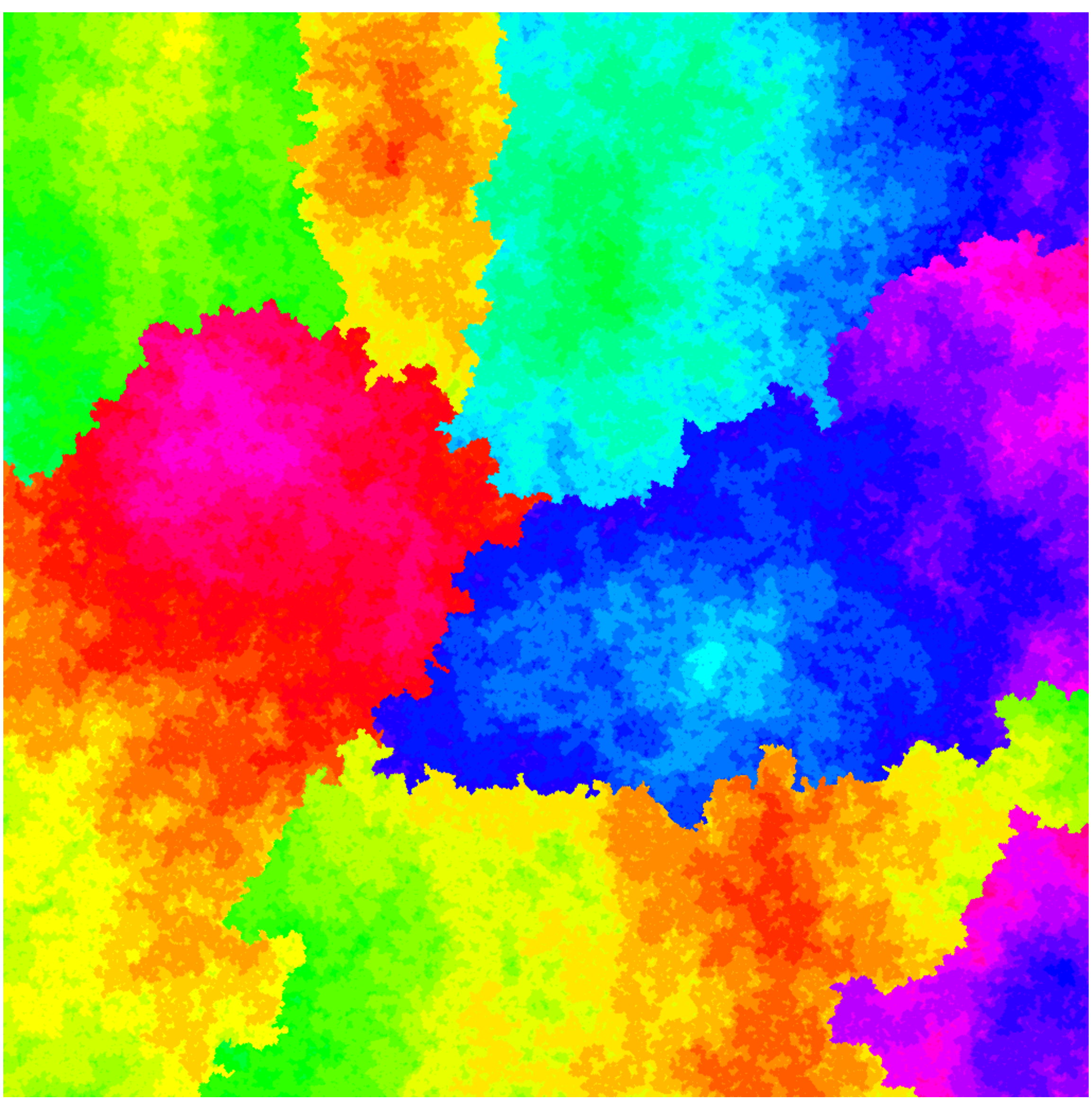

Figure 152: 7 intertwined explorations, 1024 × 1024, site percolation, $p = 0.57$.

## 80.6 Intersections in the site model

Our next aim is to define properly the intersection sets associated to the multiple intertwined explorations. We recall that the edge boundary $\Delta_D A$ of a subset $A$ of $D$ is the set of the edges of $D$ which have one endpoint in $A$ and the other in $D \setminus A$, i.e.,

$$\Delta_D A \,=\, \big\{\, e = \langle x, y\rangle \in \mathbb{E}^d : x \in A, y \in D \setminus A \,\big\} \,.$$

**Definition 80.12.** *We consider the standard framework of definition 80.7 and the associated $N$ intertwined explorations $(\mathcal{A}_i(t), \mathcal{W}_i(t), 1 \leq i \leq N, 0 \leq t \leq T)$. The first intersection set $\mathcal{I}(0)$ is*

$$\mathcal{I}(0) \,=\, \bigcup_{1\leq j<k\leq N} \big(\Delta_D \mathcal{A}_j(0) \cap \Delta_D \mathcal{A}_k(0)\big) \,.$$

*We define the successive intersection sets $(\mathcal{I}(t), 1 \leq t \leq T)$ by setting*

$$\forall t \geq 1 \qquad \mathcal{I}(t) \,=\, \bigcup_{1\leq j<k\leq N} \big(\Delta_D \mathcal{A}_j(t) \cap \Delta_D \mathcal{A}_k(t)\big) \setminus \bigcup_{0\leq s\leq t-1} \mathcal{I}(s) \,. \qquad (80.36)$$

*We say that an edge is an intersection edge if it belongs to an intersection set.*

A quick comparison with the intersection sets for the bond model defined in (29.15) would lead to the false impression that the definition (80.36) is a straightforward generalization to the site model and to multiple explorations. Once again, the reality is more complex. In the bond model, the intersection sets were closed edges preventing the connection between two aggregates. In the site model, the intersection sets are still edges preventing the connection between two aggregates, but in the sense that at least one endpoint of the edge has to be closed. In the case one endpoint is open, the other endpoint, which is closed, is pivotal for the connection between the two aggregates. From the construction of the intertwined explorations, we see that the closed endpoint was acquired first by an explorer, and the open endpoint later on by another explorer. This closed endpoint was a taboo site for the second explorer, which played a significant role: if it had not been in the taboo set, it would have been acquired by the second explorer. We call such a site a relevant taboo site. Furthermore, the open site is also a relevant taboo site for the first explorer. If both endpoints of the edge are closed, the endpoints are not necessarily pivotal, but they are both relevant taboo sites, albeit at a later iteration. We shall make precise the notion of relevant taboo sites in the next subsection.

The intersection sets and the relevant taboo sites will play a crucial role in the proof of theorem 9.2. We will focus on the intersection edges whose both endpoints are in the goal region $G$. We decided that the termination time is one step after the region $G$ has been fully explored, and this ensures that no edge of $\mathcal{I}(T)$ has both its endpoints in $G$.

## 80.7 The relevant taboo sites

We are now ready to define precisely the set of the relevant taboo sites.

**Definition 80.13.** *Let $i$ belong to $\{ 1,\dots,N \}$ and $t$ to $\{ 0,\dots,T \}$. A site $x$ is said to be a relevant taboo site for the explorer $i$ at iteration $t$ if*

• *by the time the explorer $i$ is called at iteration $t$, the site $x$ has already been visited by another explorer;*

• *if the previous condition was not met, the site $x$ would be visited by the explorer $i$ at iteration $t$.*

*We denote by $\mathcal{T}_i^*(t)$ the set of the relevant taboo sites for the explorer $i$ at iteration $t$.*

We will next characterize the relevant taboo sites with the help of the waiting aggregates of the intertwined explorations. Let $i$ be an integer in $\{ 1,\dots,N \}$. For $t=0$, we have

$$\mathcal{T}_i^*(0) \,=\, \overline{\text{Clusters}}\,\big(\{ w_i \}, D, \varnothing\big) \cap \Big( \bigcup_{1\leq j<i} \mathcal{A}_j(0)\Big)\,, \tag{80.37}$$

and the sites of $\mathcal{T}_i^*(0)$ are in fact the pivotal sites that prevent $w_i$ from being connected to one site among $w_1,\dots,w_{i-1}$. Let next $t\geq 1$. We recall that the taboo set for the explorer $i$ at iteration $t$ is

$$\mathcal{T}_i(t) \,=\, \mathcal{A}_1(t)\cup\cdots\cup\mathcal{A}_{i-1}(t)\cup\mathcal{A}_{i+1}(t-1)\cup\cdots\cup\mathcal{A}_N(t-1)\,,$$

and that $\mathcal{A}_i(t)$ is obtained by adjoining to $\mathcal{A}_i(t-1)$ the regions visited by the taboo exploration starting from the neighbours of the waiting set $\mathcal{W}_i(t-1)$, i.e.,

$$\mathcal{A}_i(t) \,=\, \mathcal{A}_i(t-1)\cup\overline{\text{Clusters}}\,\big(\mathcal{N}(\mathcal{W}_i(t-1)), D\setminus\mathcal{A}_i(t-1), \mathcal{T}_i(t)\big)\,.$$

**Lemma 80.14.** *Let us fix $i$ in $\{ 1,\dots,N \}$ and $t$ in $\{ 1,\dots,T \}$. A site $x$ belongs to $\mathcal{T}_i^*(t)$ if and only if it belongs to $\mathcal{T}_i(t)$ and to*

$$\overline{\textit{Clusters}}\,\big(\mathcal{N}(\mathcal{W}_i(t-1)), D\setminus\mathcal{A}_i(t-1), \mathcal{T}_i(t)\setminus\{ x \}\big)\,. \tag{80.38}$$

*Proof.* The first condition in definition 80.13 amounts to saying that the site $x$ belongs to the taboo set $\mathcal{T}_i(t)$. The second condition says that, if we remove $x$ from the taboo set of the explorer $i$ at time $t$, then the site $x$ will be visited by the explorer $i$. This means that $x$ belongs to the set (80.38). □

The formulation proposed in lemma 80.14 is somewhat cryptic, but it seems difficult to give an explicit definition of $\mathcal{T}_i^*(t)$. Notice that we ask that $x$ belongs to the set (80.38), and the definition of this set involves $x$ itself!

The next definitions will help to analyze further the intersection edges and the relevant taboo sites.

**Definition 80.15** (Visiting iteration and visitor index)**.** *For a site $x$ in $D$, we define the visiting iteration $T(x)$ of $x$ as the iteration number at which $x$ is visited by an explorer:*

$$T(x) \,=\, \min\big\{\, t\geq 0: \exists\, i\in\{\,1,\dots,N\,\}\quad x\in\mathcal{A}_i(t)\,\big\}\,.$$

*In case the site $x$ is never visited by an explorer, the visiting iteration $T(x)$ is infinite. Whenever $T(x)$ is finite, we define also the visitor index $I(x)$ as the index of the explorer which visited $x$:*

$$\forall x\in D\qquad T(x)<+\infty\quad\Longrightarrow\quad \exists\,!I(x)\in\{\,1,\dots,N\,\}\quad x\in\mathcal{A}_{I(x)}\big(T(x)\big)\,.$$

From the construction of the multiple intertwined explorations, we observe that the visiting iteration and the visitor index must satisfy the following constraints.

**Lemma 80.16.** *For any $x,y$ in $D$ such that $|x-y|=1$, we have*

$$\begin{gathered}T(x),T(y)<+\infty\quad\Longrightarrow\quad \big|T(x)-T(y)\big|\,\leq\,1\,,\\ I(x)\neq I(y)\quad\Longrightarrow\quad \{x,y\}\in\mathcal{I}\big(\max\big(T(x),T(y)\big)\big)\,.\end{gathered}$$

*Proof.* Suppose for instance that $T(x)=s$ and $T(y)\geq s+2$. At iteration $s$, the explorer $I(x)$ visits $x$, and $y$ is not visited by any other explorer. This implies that $x$ is closed, and that $y$ is in the neighborhood of the waiting set of the explorer $I(x)$. Necessarily, at iteration $s+1$, the site $y$ will be visited, either by the explorer $I(x)$ or by an explorer whose index is strictly smaller than $I(x)$. But this contradicts the fact that $T(y)\geq s+2$. Thus $T(y)\leq s+1=T(x)+1$. A similar argument shows that $T(x)\leq T(y)+1$. The second implication is a direct consequence of the definition of the intersection sets. ☐

A bit more can be said for the endpoints of an intersection edge, as stated in the next lemma.

**Lemma 80.17.** *Let $t$ belong to $\{\,0,\dots,T\,\}$ and let $e=\{x,y\}$ be an edge in $\mathcal{I}(t)$. We have*

$$\min\big\{\,T(x),T(y)\,\big\}\,\geq\,t-1\,,\qquad \max\big\{\,T(x),T(y)\,\big\}\,=\,t\,.$$

*Proof.* We already know from lemma 80.16 that $|T(x)-T(y)|\leq 1$. By the definition of $\mathcal{I}(t)$, there exist two indices $j<k$ in $\{\,1,\dots,N\,\}$ such that $x$ is in $\mathcal{A}_j(t)$ and $y$ in $\mathcal{A}_k(t)$. This implies readily that

$$\max\big\{\,T(x),T(y)\,\big\}\,\leq\,t\,.$$

We claim also that $\max\big\{T(x),T(y)\big\}\,\geq\,t$. Suppose that $T(x)<t$ and $T(y)<t$. In this case, the edge $e$ would belong to $\Delta_D\mathcal{A}_j(t-1)\cap\Delta_D\mathcal{A}_k(t-1)$, however

$$\Delta_D\mathcal{A}_j(t-1)\cap\Delta_D\mathcal{A}_k(t-1)\,\subset\,\bigcup_{1\leq s\leq t-1}\mathcal{I}(s)\,,$$

and this would contradict the fact that $e$ is in $\mathcal{I}(t)$. ☐

With the help of the definition 80.15, we will analyze more precisely what happens at the endpoints of the intersection edges. We prove next that both endpoints of an intersection edge are relevant taboo sites. We define $\mathcal{T}^*(t)$ as the set of all the relevant taboo sites at time $t$, i.e.,

$$\forall t \in \{\,0,\dots,T\,\} \qquad \mathcal{T}^*(t) \,=\, \bigcup_{1\leq i\leq N} \mathcal{T}_i^*(t)\,.$$

**Proposition 80.18.** *Let $t$ be an integer in $\{\,0,\dots,T-1\,\}$ and let $e=\{x,y\}$ be an edge of $G$ in $\mathcal{I}(t)$. Then at least one site among $x,y$ is closed and*

$$\{x,y\}\,\subset\,\mathcal{T}^*(t)\cup\mathcal{T}^*(t+1)\,.$$

*Proof.* Let $t$ be an integer in $\{\,0,\dots,T\,\}$, let $x,y$ be two neighbouring sites in $G$ and suppose that the edge $e=\{x,y\}$ is in $\mathcal{I}(t)$. Because of the definition of $T$, we must have $t<T$. We denote by $i=I(x)$ (respectively $j=I(y)$) the visitor index of $x$ (respectively $y$). We consider several cases, depending on the status of $x$, $y$, their visiting iterations and their visitor indices:

• $x$ and $y$ are open. This is not compatible with the fact that $e$ is in $\mathcal{I}(t)$.

• $x$ is closed, $y$ is open, $T(x)=t-1, T(y)=t$. The site $x$ was visited at iteration $t-1$ by the explorer $i$. At iteration $t$, the site $y$ is visited by the explorer $j$. If $x$ had not been already visited, then it would be visited by the explorer $j$, thus $x$ is a relevant taboo site for the explorer $j$ at iteration $t$. Moreover, at iteration $t$, the site $y$ would be visited by the explorer $i$ if it had not been previously visited by the explorer $j$. Thus we must have $j<i$. We conclude that the site $x$ is in $\mathcal{T}_j^*(t)$ and the site $y$ in $\mathcal{T}_i^*(t)$.

• $x$ is closed, $y$ is open, $T(x)=t, T(y)=t$. Both sites $x$ and $y$ are visited at iteration $t$. The site visited first must be $x$, otherwise the explorer $j$ would have the opportunity to visit $x$ before the explorer $i$. Therefore we have $i<j$, and the site $x$ is in $\mathcal{T}_j^*(t)$. Moreover, the explorer $i$ will try to visit $y$ at iteration $t+1$, thus the site $y$ is in $\mathcal{T}_i^*(t+1)$.

• $x$ and $y$ are closed, $T(x)=t-1, T(y)=t$. The site $x$ was visited at iteration $t-1$ by the explorer $i$. At iteration $t$, the site $y$ is visited by the explorer $j$. This must occur before the explorer $i$ has a chance to visit $y$, therefore we have $i>j$. We conclude that the site $x$ is in $\mathcal{T}_j^*(t+1)$ and the site $y$ in $\mathcal{T}_i^*(t)$.

• $x$ and $y$ are closed, $T(x)=t, T(y)=t$. Both sites $x,y$ are visited at iteration $t$. Both sites will be relevant taboo sites, but at the subsequent iteration. Thus the site $x$ is in $\mathcal{T}_j^*(t+1)$ and the site $y$ in $\mathcal{T}_i^*(t+1)$.

We see that, in each case, the edge $e=\{x,y\}$ is in $\mathcal{I}(t)\cup\mathcal{I}(t+1)$. ☐

Let us define the set $\mathcal{T}^{*0}(t)$ of the closed relevant taboo sites at time $t$:

$$\forall t \in \{\,0,\dots,T\,\} \qquad \mathcal{T}^{*0}(t) \,=\, \big\{\,x\in\mathcal{T}^*(t):\omega(x)=0\,\big\}\,.$$

Proposition 80.18 yields an inequality between the cardinalities of the intersection edges and the closed relevant taboo sites, that we state in the next corollary.

**Corollary 80.19.** *For $t$ in $\{0,\dots,T-1\}$, we have*

$$\big|\mathcal{I}(t)\big| \,\leq\, 2d\,\big|\mathcal{T}^{*0}(t)\big| + 2d\,\big|\mathcal{T}^{*0}(t+1)\big|\,.$$

We prove finally a partial converse to proposition 80.18: the relevant taboo sites are endpoints of intersection edges.

**Lemma 80.20.** *Let $t$ be an integer in $\{0,\dots,T\}$ and let $x$ belong to $\mathcal{T}^*(t)$. There exists a neighbour $y$ of $x$ such that the edge $e=\{x,y\}$ is in $\mathcal{I}(t-1)\cup\mathcal{I}(t)$, and moreover $t-2\leq T(x)\leq t$.*

*Proof.* Let $t$ belong to $\{0,\dots,T\}$ and let $x$ in $\mathcal{T}^*(t)$. We denote by $i=I(x)$ the visitor index of $x$ and we suppose that $x$ is in $\mathcal{T}_j^*(t)$ for some $j\neq i$ in $\{1,\dots,N\}$. Since $x$ is a relevant taboo site for the explorer $j$ at iteration $t$, then the site $x$ has been visited by the explorer $i$ before the explorer $j$ has a chance to visit it, thus $T(x)\leq t$. The second condition in the definition 80.13 implies that there exists also a neighbour $y$ of $x$ in $D$ which would allow the explorer $j$ to visit $x$ at time $t$ if it had not already been visited by the explorer $i$. This neighbour $y$ is such that $I(y)=j$ and either ($y$ is open and $T(y)=t$) or ($y$ is closed and $T(y)=t-1$). We examine separately these two cases, along with the subcases $i<j$ and $i>j$:

- $y$ is open and $T(y)=t$.
  - ⋆ $i<j$. Necessarily $T(x)=t$, otherwise the explorer $i$ would visit $y$.
  - ⋆ $i>j$. Necessarily $T(x)=t-1$, otherwise the explorer $j$ would visit $x$.

  In both subcases, the edge $\{x,y\}$ is in $\mathcal{I}(t)$.
- $y$ is closed and $T(y)=t-1$.
  - ⋆ $i<j$. Either $T(x)=t-1$ and $x$ is closed, or $T(x)=t$. Then the edge $\{x,y\}$ is in $\mathcal{I}(t-1)$ in the first case and in $\mathcal{I}(t)$ in the second case.
  - ⋆ $i>j$. We must have $T(x)=t-1$, otherwise the explorer $j$ would visit $x$. Then the edge $\{x,y\}$ is in $\mathcal{I}(t-1)$.

In all the cases, we see that the edge $\{x,y\}$ is in $\mathcal{I}(t-1)\cup\mathcal{I}(t)$. □

Let $i$ in $\{1,\dots,N\}$. We partition further the relevant taboo sites $\mathcal{T}_i^*(t)$ for the explorer $i$ at time $t$ into two subsets according to the visiting time, as follows. We define $\mathcal{T}_i^-(0)=\varnothing$, $\mathcal{T}_i^+(0)=\mathcal{T}_i^*(0)$, and for $t$ in $\{1,\dots,T\}$,

$$\begin{aligned}
\mathcal{T}_i^-(t) &= \big\{\,x\in\mathcal{T}_i^*(t): T(x)=t-1 \text{ or } T(x)=t-2\,\big\}\,,\\
\mathcal{T}_i^+(t) &= \big\{\,x\in\mathcal{T}_i^*(t): T(x)=t\,\big\}\,.
\end{aligned}$$

It follows from the classification of the intersection edges done in the course of the proof of proposition 80.18 that

$$\forall t\in\{0,\dots,T\}\qquad \mathcal{T}_i^*(t)\,=\,\mathcal{T}_i^-(t)\,\cup\,\mathcal{T}_i^+(t)\,.$$

We define also, for $t$ in $\{0,\dots,T\}$,

$$\mathcal{T}^-(t)\,=\,\bigcup_{1\leq i\leq N}\mathcal{T}_i^-(t)\,,\qquad \mathcal{T}^+(t)\,=\,\bigcup_{1\leq i\leq N}\mathcal{T}_i^+(t)\,. \tag{80.39}$$

In the case of the bond model, the intersection edges were pivotal edges, albeit of order larger than one (see lemma 29.5). The situation is more complex in the site model, we only have the following result.

**Lemma 80.21.** *The sites of $\mathcal{T}^+(0)$ are the pivotal sites for the event $\mathcal{D}$. For $t$ in $\{1,\dots,T\}$, we have the inclusions*

$$\mathcal{T}^-(t) \subset \mathcal{P}_{2t}(\mathcal{D}) \cup \Big\{ x\in D : T_D\big(x,\mathcal{P}_{2t-1}(\mathcal{D})\big) \leq t \Big\},$$
$$\mathcal{T}^+(t) \subset \mathcal{P}_{2t+1}(\mathcal{D}) \cup \Big\{ x\in D : T_D\big(x,\mathcal{P}_{2t}(\mathcal{D})\big) \leq t \Big\}.$$

*Proof.* Notice that $\mathcal{T}^-(0)=\varnothing$, so that $\mathcal{T}^*(0)=\mathcal{T}^+(0)$. The sites of $\mathcal{T}^*(0)$ have already been examined, see formula (80.37). Let $t$ in $\{1,\dots,T\}$ and let $x$ be a site belonging to $\mathcal{T}^-(t)$. By the definition of $\mathcal{T}^-(t)$, there exists $j$ in $\{1,\dots,N\}$ such that the site $x$ is in $\mathcal{T}_j^-(t)$, whence $T_D(w_j,x)\leq t$. In particular, there exists a path $y_0,\cdots,y_\ell$ in $D$ joining $y_0=w_j$ to $y_\ell=x$ such that at most $t$ sites among $y_0,\cdots,y_{\ell-1}$ are closed. We denote by $E$ the set of these sites. Setting $i=I(x)$, we have also that $x$ is in $\mathcal{A}_i(t-1)$, whence $T_{\mathcal{A}_i(t-1)}(w_i,x)\leq t-1$ by lemma 80.8. Thus there exists a path $z_0,\cdots,z_m$ in $\mathcal{A}_i(t-1)$ joining $z_0=w_i$ to $z_m=x$ such that at most $t-1$ sites among $z_0,\cdots,z_{m-1}$ are closed. We denote by $F$ the set of these sites. Let $\widetilde{\omega}$ be the configuration obtained from the configuration $\omega$ by opening all the sites of $E\cup F$. We consider two cases:

• $\widetilde{\omega}$ is not in $\mathcal{D}$. In this case, one of the sites in $E\cup F$ belongs to $\mathcal{P}_{2t-1}(\mathcal{D})$, but all the sites in $E\cup F$ are at travel distance less than or equal to $t$ from $x$.

• $\widetilde{\omega}$ is in $\mathcal{D}$. In this case the site $x$ is pivotal for $\mathcal{D}$ in $\widetilde{\omega}$, and $x$ is in $\mathcal{P}_{2t}(\mathcal{D},\omega)$.

The case of $\mathcal{T}^+(t)$ is a slight variation of the previous argument, essentially $\mathcal{A}_i(t-1)$ has to be replaced by $\mathcal{A}_i(t)$ and the set $F$ contains at most $t$ sites. ☐

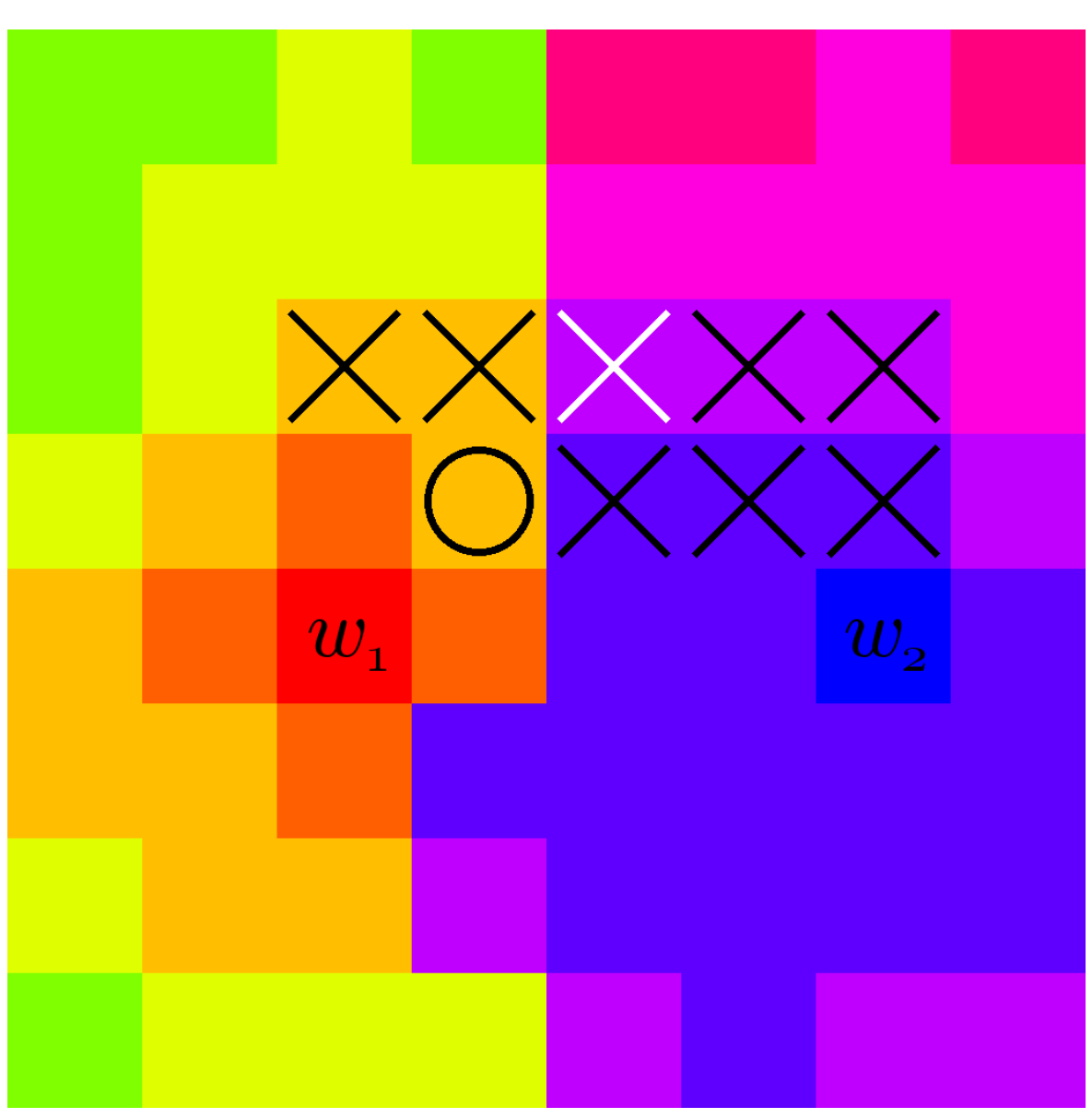

Figure 153: The explorations started from $w_1$, $w_2$ in $\Lambda(8)$ with $p=0.2$: $\mathcal{A}_1(0)$, $\mathcal{A}_1(1)\setminus\mathcal{A}_1(0)$, $\mathcal{A}_1(2)\setminus\mathcal{A}_1(1)$, $\mathcal{A}_2(0)$, $\mathcal{A}_2(1)\setminus\mathcal{A}_2(0)$, $\mathcal{A}_2(2)\setminus\mathcal{A}_2(1)$, $\times$, $w_1$, $w_2$ are closed, $\circ$ is open. The site with white $\times$ is in $\mathcal{T}^-(3)\setminus\mathcal{P}_6(D)$.

# 81 Reconstruction

The interest of the relevant taboo sites is that they allow to describe precisely the interaction between the explorers. More precisely, during the iteration $t$, the explorer $i$ perceives the other explorers only through his set of relevant taboo sites $\mathcal{T}_i^*(t)$. We formalize next this point. Let us fix an index $i$ in $\{1,\dots,N\}$. Corollary 80.11 shows that the sequence $(A_i(t), 0\leq t\leq T)$ can be computed by performing a taboo exploration starting from $w_i$ with the taboo set $\widehat{\mathcal{A}}_i(T)$ defined in (80.28). This is the content of formula (80.29), which we recall next:

$$\forall t\in\{0,\dots,T\}\qquad \mathcal{A}_i(t)\,=\,\mathcal{B}\big(w_i, D\setminus\widehat{\mathcal{A}}_i(T), t\big)\,. \tag{81.1}$$

The taboo set $\widehat{\mathcal{A}}_i(T)$, defined in (80.28), consists of all the sites visited by the other explorers throughout the exploration. In the next lemma, we show that a much smaller taboo set can be used to get the same result.

**Proposition 81.1.** *Let us fix $i$ in $\{1,\dots,N\}$, and let us define*

$$\forall t\in\{0,\dots,T\}\qquad \overline{\mathcal{T}}_i(t)\,=\,\bigcup_{0\leq s\leq t}\mathcal{T}_i^*(s)\,.$$

*For any $t$ in $\{0,\dots,T\}$, we have*

$$\mathcal{A}_i(t)\,=\,\mathcal{B}\big(w_i, D\setminus\overline{\mathcal{T}}_i(t), t\big)\,. \tag{81.2}$$

*Moreover, the set $\overline{\mathcal{T}}_i(t)$ is the smallest possible taboo set for which the identity* (81.2) *holds, and we have*

$$\forall x\in\ \overline{\mathcal{T}}_i(t)\quad x\in\mathcal{B}\Big(w_i, D\setminus\big(\overline{\mathcal{T}}_i(t)\setminus\{x\}\big), t\Big)\,. \tag{81.3}$$

*Proof.* Let $i$ in $\{1,\dots,N\}$ and $t$ in $\{0,\dots,T\}$. By construction, the set $\overline{\mathcal{T}}_i(t)$ is included in $\widehat{\mathcal{A}}_i(T)$, therefore, using formula (81.1), we have

$$\mathcal{A}_i(t)\,=\,\mathcal{B}\big(w_i, D\setminus\widehat{\mathcal{A}}_i(T), t\big)\,\subset\,\mathcal{B}\big(w_i, D\setminus\overline{\mathcal{T}}_i(t), t\big)\,.$$

To prove the converse inclusion, we proceed by induction on $t$. Let us start with $t=0$. The set $\overline{\mathcal{T}}_i(0)=\mathcal{T}_i^*(0)$ consists of the sites which are pivotal for the connection between $w_i$ and a site $w_j$ with an index $j<i$. Recall that the set $\mathcal{A}_i(0)$ is obtained through the following taboo exploration (see definition (80.14)):

$$\mathcal{A}_i(0)\,=\,\overline{\text{Clusters}}\Big(\{w_i\}, D, \bigcup_{1\leq j<i}\mathcal{A}_j(0)\Big)\,. \tag{81.4}$$

The relevant taboo sites for this exploration are the sites in the set

$$\Big\{x\in\bigcup_{1\leq j<i}\mathcal{A}_j(0) : \exists y\in\mathcal{A}_i(0)\cap\mathcal{N}(x)\quad \omega(y)=1\Big\}\,. \tag{81.5}$$

Yet the above set coincides with $\mathcal{T}_i^*(0)$, thus (81.4) can be rewritten as

$$\mathcal{A}_i(0)\,=\,\mathcal{B}\big(w_i, D\setminus\overline{\mathcal{T}}_i(0), 0\big)\,. \tag{81.6}$$

We notice finally that, for any $x$ in the set (81.5), we have $T_D(w_i,x)=0$. On the event $\mathcal{D}$, the sites in $\overline{\mathcal{T}_i}(0)$ are closed, and a path joining $w_i$ to $x$ in $D$ with null travel time can visit only open sites, apart from $x$. Thus

$$T_{D\setminus(\overline{\mathcal{T}_i}(0)\setminus\{\,x\,\})}(w_i,x)=0\,,$$

and therefore

$$x\in\mathcal{B}\big(w_i,D\setminus(\overline{\mathcal{T}_i}(0)\setminus\{\,x\,\}),0\big)\,.$$

This proves (81.3) for $t=0$ and that $\overline{\mathcal{T}_i}(0)$ is the smallest set satisfying (81.6).

We perform next the induction step. Let $t$ in $\{\,0,\dots,T-1\,\}$ and suppose that the equality (81.2) has been proved until the integer $t$. We write

$$\begin{aligned}\mathcal{B}\big(w_i,D\setminus\overline{\mathcal{T}_i}(t+1),t+1\big)\;=\;&\mathcal{B}\big(w_i,D\setminus\overline{\mathcal{T}_i}(t+1),t\big)\\&\cup\text{Shell}\,\big(w_i,D\setminus\overline{\mathcal{T}_i}(t+1),t+1\big)\,.\end{aligned}\tag{81.7}$$

Since $\overline{\mathcal{T}_i}(t)$ is included in $\overline{\mathcal{T}_i}(t+1)$, then

$$\mathcal{B}\big(w_i,D\setminus\overline{\mathcal{T}_i}(t+1),t\big)\;\subset\;\mathcal{B}\big(w_i,D\setminus\overline{\mathcal{T}_i}(t),t\big)\;=\;\mathcal{A}_i(t)\,,\tag{81.8}$$

where we used the induction hypothesis to get the last equality. We deal next with the last term in (81.7). Let $x$ be a site in $\text{Shell}\,\big(w_i,D\setminus\overline{\mathcal{T}_i}(t+1),t+1\big)$. By definition, there exists a path $y_0,y_1,\cdots,y_\ell$ in $D\setminus\overline{\mathcal{T}_i}(t+1)$ joining $y_0=w_i$ to $y_\ell=x$ such that exactly $t+1$ sites among $y_0,\cdots,y_{\ell-1}$ are closed. Let $k$ be the index of the $(t+1)$-th closed site along the path, i.e., $k$ is the unique index such that

$$\big|\big\{\,j:0\leq j\leq k,\ \omega(y_j)=0\,\big\}\big|\;=\;t+1\,.$$

Let us set $z=y_{k+1}$. We have then

$$\begin{gathered}y_k\in\text{Shell}\,\big(w_i,D\setminus\overline{\mathcal{T}_i}(t+1),t\big)\,,\quad|y_k-z|=1\,,\\ z\in\text{Shell}\,\big(w_i,D\setminus\overline{\mathcal{T}_i}(t+1),t+1\big)\,,\quad T_{D\setminus\overline{\mathcal{T}_i}(t+1)}(z,x)=0\,.\end{gathered}\tag{81.9}$$

Using the induction hypothesis and (81.9), we see that $y_k$ belongs to

$$\begin{aligned}\text{Shell}\,\big(w_i,D\setminus\overline{\mathcal{T}_i}(t),t\big)\;&=\;\mathcal{A}_i(t)\setminus\mathcal{A}_i(t-1)\\&=\;\overline{\text{Clusters}}\,\big(\mathcal{N}(\mathcal{W}_i(t-1)),D\setminus\mathcal{A}_i(t-1),\mathcal{T}_i(t)\big)\,.\end{aligned}$$

Furthermore, the site $z$ is in

$$\text{Shell}\,\big(w_i,D\setminus\overline{\mathcal{T}_i}(t),t+1\big)\;\subset\;D\setminus\mathcal{B}\big(w_i,D\setminus\overline{\mathcal{T}_i}(t),t\big)\;\subset\;D\setminus\mathcal{A}_i(t)\,.$$

Suppose that

$$z\in\mathcal{T}_i(t+1)\;=\;\mathcal{A}_1(t+1)\cup\cdots\cup\mathcal{A}_{i-1}(t+1)\cup\mathcal{A}_{i+1}(t)\cup\cdots\cup\mathcal{A}_N(t)\,.\tag{81.10}$$

This would further imply that $z$ would belong to $\overline{\mathcal{T}_i}(t+1)$, which is absurd. Therefore (81.10) does not hold, and $z$ is not in $\mathcal{T}_i(t+1)$. In particular, the site $z$ is in $D\setminus\big(\mathcal{A}_i(t)\cup\mathcal{T}_i(t)\big)$ and it is a neighbour of $y_k$, therefore $y_k$ is in

$$\begin{aligned}\mathcal{W}_i(t)\;&=\;\text{Waiting}\,\big(\mathcal{N}(\mathcal{W}_i(t-1)),D\setminus\mathcal{A}_i(t-1),\mathcal{T}_i(t)\big)\\&=\;\partial^{\,in}_{D\setminus\mathcal{A}_i(t-1)}\big(\mathcal{A}_i(t)\setminus\mathcal{A}_i(t-1)\big)\,,\end{aligned}$$

and $z$ is in

$$\mathcal{N}(\mathcal{W}_i(t))\cap D\setminus\big(\mathcal{A}_i(t)\cup\mathcal{T}_i(t+1)\big)\ \subset\ \mathcal{A}_i(t+1)\,.$$

If $z=x$, then we are done. Otherwise, suppose that $z\neq x$, or equivalently that $k+1<\ell$. Since the sites $y_{k+1},\dots,y_{\ell-1}$ must be open, and $\mathcal{A}_i(t+1)$ is an aggregate, then $y_{\ell-1}$ is also in $\mathcal{A}_i(t+1)$. Suppose that $x$ is in $\mathcal{T}_i(t+1)$. Since $x$ is the neighbour of $y_{\ell-1}$, then $x$ would belong to $\overline{\mathcal{T}_i}(t+1)$, which is absurd. Thus $x$ is not in $\mathcal{T}_i(t+1)$, and it belongs also to

$$\overline{\text{Clusters}}\,\big(\mathcal{N}(\mathcal{W}_i(t)),D\setminus\mathcal{A}_i(t),\mathcal{T}_i(t+1)\big)\ \subset\ \mathcal{A}_i(t+1)\,.$$

We have proved that

$$\text{Shell}\,\big(w_i,D\setminus\overline{\mathcal{T}_i}(t+1),t+1\big)\ \subset\ \mathcal{A}_i(t+1)\,,$$

and together with (81.7) and (81.8), this proves the identity (81.2) at rank $t+1$. To complete the induction step, it remains to prove the statement on the minimality of $\overline{\mathcal{T}_i}(t+1)$. Let $x$ be a site in $\overline{\mathcal{T}_i}(t+1)$. We will show that $x$ is also in

$$\mathcal{B}\big(w_i,D\setminus(\overline{\mathcal{T}_i}(t+1)\setminus\{\,x\,\}),t+1\big)\,. \tag{81.11}$$

Let $s$ be the smallest integer in $\{\,0,\dots,t+1\,\}$ such that $x$ is in $\mathcal{T}_i^*(s)$. If $s=0$, then $x$ belongs to (80.37). Notice that all the sites in $\mathcal{T}_i^*(0)$ must be closed. The travel time between $w_i$ and any site of $\mathcal{T}_i^*(0)$ is null. In particular, we have $T_D(w_i,x)=0$, and there exists a path $y_0,y_1,\cdots,y_\ell$ in $D$ joining $y_0=w_i$ to $y_\ell=x$ such that $y_0,\cdots,y_{\ell-1}$ are all open; thus these sites belong to $\mathcal{A}_i(0)$, and they are not in $\mathcal{T}_i(t+1)$. This implies that $x$ is also in

$$\begin{aligned}\overline{\text{Clusters}}\,\Big(\{\,w_i\,\},D\setminus\big(\mathcal{T}_i(t+1)\setminus\{\,x\,\}\big),\varnothing\Big)\ &=\ \mathcal{B}\big(w_i,D\setminus(\mathcal{T}_i(t+1)\setminus\{\,x\,\}),0\big)\\ &\subset\ \mathcal{B}\big(w_i,D\setminus(\overline{\mathcal{T}_i}(t+1)\setminus\{\,x\,\}),t+1\big)\,.\end{aligned}$$

Suppose next that $1\leq s\leq t+1$. According to lemma 80.14, we have

$$x\in\mathcal{T}_i(s)\cap\overline{\text{Clusters}}\,\big(\mathcal{N}(\mathcal{W}_i(s-1)),D\setminus\mathcal{A}_i(s-1),\mathcal{T}_i(s)\setminus\{\,x\,\}\big)\,,$$

hence there exists a path $y_0,\cdots,y_\ell$ in $D\setminus\big(\mathcal{A}_i(s-1)\cup(\mathcal{T}_i(s)\setminus\{\,x\,\})\big)$ joining a site $y_0$ of $\mathcal{N}(\mathcal{W}_i(s-1))$ to $y_\ell=x$ such that $y_0,\cdots,y_{\ell-1}$ are all open. Among such paths, we pick one having minimal length, so that $y_0,\cdots,y_{\ell-1}$ are all distinct from $x$. We consider first the case $\ell\geq 1$. In this case, the site $y=y_{\ell-1}$ is a neighbour of $x$, and

$$y\in\overline{\text{Clusters}}\,\big(\mathcal{N}(\mathcal{W}_i(s-1)),D\setminus\mathcal{A}_i(s-1),\mathcal{T}_i(s)\big)\,. \tag{81.12}$$

We show next that the taboo set $\mathcal{T}_i(s)$ in (81.12) can be replaced by $\overline{\mathcal{T}_i}(t+1)$, without altering the result. It follows from the definitions of $\mathcal{A}_i(s)$ and $\overline{\mathcal{T}_i}(s)$ that

$$\begin{aligned}\mathcal{A}_i(s)\setminus\mathcal{A}_i(s-1)\ &=\ \overline{\text{Clusters}}\,\big(\mathcal{N}(\mathcal{W}_i(s-1)),D\setminus\mathcal{A}_i(s-1),\mathcal{T}_i(s)\big)\\ &=\ \overline{\text{Clusters}}\,\big(\mathcal{N}(\mathcal{W}_i(s-1)),D\setminus\mathcal{A}_i(s-1),\overline{\mathcal{T}_i}(s)\big)\\ &=\ \overline{\text{Clusters}}\,\big(\mathcal{N}(\mathcal{W}_i(s-1)),D\setminus\mathcal{A}_i(s-1),\overline{\mathcal{T}_i}(t+1)\big)\,.\end{aligned} \tag{81.13}$$

The last equality stems from the fact that

$$\overline{\mathcal{T}}_i(s) \subset \overline{\mathcal{T}}_i(t+1) \,, \qquad \overline{\mathcal{T}}_i(t+1) \cap \mathcal{A}_i(s) \,=\, \varnothing \,.$$

We conclude from (81.12) and (81.13) that

$$y \in \overline{\text{Clusters}}\,\big(\mathcal{N}(\mathcal{W}_i(s-1)), D \setminus \mathcal{A}_i(s-1), \overline{\mathcal{T}}_i(t+1)\big) \,. \tag{81.14}$$

Since $y$ is a neighbour of $x$ and it is open, then (81.14) implies that

$$x \in \overline{\text{Clusters}}\,\big(\mathcal{N}(\mathcal{W}_i(s-1)), D \setminus \mathcal{A}_i(s-1), \overline{\mathcal{T}}_i(t+1) \setminus \{\, x \,\}\big) \,. \tag{81.15}$$

In the case $\ell = 0$, the site $x$ belongs to $\mathcal{N}(\mathcal{W}_i(s-1)) \setminus \mathcal{A}_i(s-1)$, and (81.15) holds as well. Furthermore $\mathcal{W}_i(s-1) \subset \mathcal{A}_i(s-1)$, whence

$$\begin{aligned} &\overline{\text{Clusters}}\,\big(\mathcal{N}(\mathcal{W}_i(s-1)), D \setminus \mathcal{A}_i(s-1), \overline{\mathcal{T}}_i(t+1) \setminus \{\, x \,\}\big) \\ &\qquad\qquad \subset\ \text{Shell}\,\big(\mathcal{A}_i(s-1), D \setminus (\overline{\mathcal{T}}_i(t+1) \setminus \{\, x \,\}), 1\big) \,. \end{aligned} \tag{81.16}$$

We prove next the following inclusion:

$$\text{Shell}\,\big(\mathcal{A}_i(s-1), D \setminus (\overline{\mathcal{T}}_i(t+1) \setminus \{\, x \,\}), 1\big) \,\subset\, \mathcal{B}\big(w_i, D \setminus (\overline{\mathcal{T}}_i(t+1) \setminus \{\, x \,\}), t+1\big) \,.$$

We have the obvious inclusion

$$\mathcal{B}\big(w_i, D \setminus (\overline{\mathcal{T}}_i(t+1) \setminus \{\, x \,\}), s\big) \,\subset\, \mathcal{B}\big(w_i, D \setminus (\overline{\mathcal{T}}_i(t+1) \setminus \{\, x \,\}), t+1\big) \,. \tag{81.17}$$

We decompose the first set in (81.17) as

$$\begin{aligned} \mathcal{B}\big(w_i, D \setminus (\overline{\mathcal{T}}_i(t+1) \setminus \{\, x \,\}), s\big) \,&=\, \mathcal{B}\big(w_i, D \setminus (\overline{\mathcal{T}}_i(t+1) \setminus \{\, x \,\}), s-1\big) \\ &\cup \text{Shell}\,\Big(\mathcal{B}\big(w_i, D \setminus (\overline{\mathcal{T}}_i(t+1) \setminus \{\, x \,\}), s-1\big), D \setminus (\overline{\mathcal{T}}_i(t+1) \setminus \{\, x \,\}), 1\Big) \,. \end{aligned} \tag{81.18}$$

Thanks to the induction hypothesis at rank $s-1$, we have

$$\mathcal{A}_i(s-1) \,=\, \mathcal{B}\big(w_i, D \setminus \overline{\mathcal{T}}_i(s-1), s-1\big) \,. \tag{81.19}$$

As $x$ belongs to $\mathcal{T}_i^*(s) \,=\, \overline{\mathcal{T}}_i(s) \setminus \overline{\mathcal{T}}_i(s-1)$, and the set $\overline{\mathcal{T}}_i(t+1)$ is disjoint from $\mathcal{A}_i(s-1)$, we deduce from (81.19) that

$$\mathcal{A}_i(s-1) \,=\, \mathcal{B}\big(w_i, D \setminus (\overline{\mathcal{T}}_i(t+1) \setminus \{\, x \,\}), s-1\big) \,. \tag{81.20}$$

Using the equality (81.20), formula (81.18) can be rewritten as

$$\begin{aligned} \mathcal{B}\big(w_i, D \setminus (\overline{\mathcal{T}}_i(t+1) \setminus \{\, x \,\}), s\big) \,&=\, \mathcal{A}_i(s-1) \\ &\qquad \cup \text{Shell}\,\big(\mathcal{A}_i(s-1), D \setminus (\overline{\mathcal{T}}_i(t+1) \setminus \{\, x \,\}), 1\big) \,. \end{aligned} \tag{81.21}$$

As $x$ belongs to the first set in the inclusion (81.16), it belongs to the second one, and thanks to (81.21), also to the first set in (81.17), which is included in (81.11). This completes the induction step. $\square$

For $i$ in $\{1,\dots,N\}$ and $t$ in $\{1,\dots,T\}$, we define the sets $\mathcal{T}_i^{*0}(t)$, $\mathcal{T}_i^{*1}(t)$ of the closed and open relevant taboo sites for the explorer $i$ at time $t$ as

$$\mathcal{T}_i^{*0}(t) \,=\, \big\{\, x\in\mathcal{T}_i^*(t):\omega(x)=0\,\big\}\,,\quad \mathcal{T}_i^{*1}(t) \,=\, \big\{\, x\in\mathcal{T}_i^*(t):\omega(x)=1\,\big\}\,.$$

The next proposition shows how the closed relevant taboo sites at time $t$ can be reached with the help of an adequate taboo exploration, once the relevant taboo sites until time $t-1$ and the open relevant taboo sites at time $t$ are known.

**Proposition 81.2.** *Let us fix $i$ in $\{1,\dots,N\}$. For any $t$ in $\{1,\dots,T\}$, we have*

$$\mathcal{T}_i^{*0}(t)\setminus\overline{\mathcal{T}}_i(t-1) \,=\, \mathcal{T}_i(t)\cap\mathcal{B}\Big(w_i,D\setminus\big(\overline{\mathcal{T}}_i(t-1)\cup\mathcal{T}_i^{*1}(t)\big),t\Big)\,. \tag{81.22}$$

*Proof.* Using the result (81.3) of proposition 81.1, we have

$$\forall x\in\mathcal{T}_i^{*0}(t)\setminus\overline{\mathcal{T}}_i(t-1)$$
$$x\in\mathcal{B}\big(w_i,D\setminus(\overline{\mathcal{T}}_i(t)\setminus\{\,x\,\}),t\big)\,\subset\,\mathcal{B}\Big(w_i,D\setminus\big(\overline{\mathcal{T}}_i(t-1)\cup\mathcal{T}_i^{*1}(t)\big),t\Big)\,.$$

Moreover $\mathcal{T}_i^{*0}(t)\subset\mathcal{T}_i(t)$, thus

$$\mathcal{T}_i^{*0}(t)\setminus\overline{\mathcal{T}}_i(t-1) \,\subset\, \mathcal{T}_i(t)\cap\mathcal{B}\Big(w_i,D\setminus\big(\overline{\mathcal{T}}_i(t-1)\cup\mathcal{T}_i^{*1}(t)\big),t\Big)\,. \tag{81.23}$$

We prove next the converse inclusion. Since the relevant taboo sites until time $t-1$ for the explorer $i$ are all contained in the set $\overline{\mathcal{T}}_i(t-1)$, and the sites of $\mathcal{T}_i^{*1}(t)$ are not in $\mathcal{A}_i(t-1)$, we have, thanks to proposition 81.1,

$$\mathcal{B}\Big(w_i,D\setminus\big(\overline{\mathcal{T}}_i(t-1)\cup\mathcal{T}_i^{*1}(t)\big),t-1\Big)=\mathcal{B}\Big(w_i,D\setminus\overline{\mathcal{T}}_i(t-1),t-1\Big)=\mathcal{A}_i(t-1). \tag{81.24}$$

Using the identity (81.24), we decompose the last set of formula (81.23) as

$$\mathcal{B}\Big(w_i,D\setminus\big(\overline{\mathcal{T}}_i(t-1)\cup\mathcal{T}_i^{*1}(t)\big),t\Big)\;=$$
$$\mathcal{A}_i(t-1)\cup\text{Shell}\,\Big(\mathcal{A}_i(t-1),D\setminus\big(\overline{\mathcal{T}}_i(t-1)\cup\mathcal{T}_i^{*1}(t)\big),1\Big)\,.$$

We know that $\mathcal{T}_i(t)$ and $\mathcal{A}_i(t-1)$ are disjoint, so it remains to study

$$\mathcal{T}_i(t)\cap\text{Shell}\,\Big(\mathcal{A}_i(t-1),D\setminus\big(\overline{\mathcal{T}}_i(t-1)\cup\mathcal{T}_i^{*1}(t)\big),1\Big)\,. \tag{81.25}$$

Let $x$ belong to the set (81.25). By definition, there exists a path $y_0,y_1,\cdots,y_\ell$ in $D\setminus\big(\overline{\mathcal{T}}_i(t-1)\cup\mathcal{T}_i^{*1}(t)\big)$ joining a site $y_0$ of $\mathcal{A}_i(t-1)$ to $y_\ell=x$ such that exactly one site among $y_0,\cdots,y_{\ell-1}$ is closed, say $y_k$, where $0\leq k\leq\ell-1$. We consider two cases:

• $y_k$ is not in $\mathcal{A}_i(t-1)$. Since $y_0,\dots,y_{k-1}$ are open and are in $\mathcal{A}_i(t-1)$, then $y_k$ must belong to $\overline{\mathcal{T}}_i(t-1)$, otherwise it would be also in $\mathcal{A}_i(t-1)$. However the whole path $y_0,\cdots,y_\ell$ is in $D\setminus\overline{\mathcal{T}}_i(t-1)$, and we have a contradiction.

• $y_k$ is in $\mathcal{A}_i(t-1)$. Since $y_{k+1},\dots,y_{\ell-1}$ are open, they are not in $\mathcal{T}_i^{*0}(t)$, hence they are in $D\setminus\big(\overline{\mathcal{T}}_i(t-1)\cup\mathcal{T}_i^*(t)\big)$ and $y_{\ell-1}$ is in

$$\text{Shell}\Big(\mathcal{A}_i(t-1), D\setminus\big(\overline{\mathcal{T}}_i(t-1)\cup\mathcal{T}_i^*(t)\big),1\Big)$$
$$=\ \text{Shell}\big(\mathcal{A}_i(t-1), D\setminus\overline{\mathcal{T}}_i(t),1\big)\ =\ \mathcal{A}_i(t)\setminus\mathcal{A}_i(t-1)\,.$$

The site $y_{\ell-1}$ is an open site of the aggregate $\mathcal{A}_i(t)$, and $x$ is a neighbour of $y_{\ell-1}$, but $x$ is not in $\mathcal{A}_i(t)$. Therefore $x$ is closed, and it belongs to $\mathcal{T}_i^*(t)$.

From these two cases, we conclude that $x$ is in $\mathcal{T}_i^{*0}(t)$, and that the set (81.25) is included in $\mathcal{T}_i^{*0}(t)\setminus\overline{\mathcal{T}}_i(t-1)$. This completes the proof of the converse inclusion of (81.23) and the proof of (81.22). □

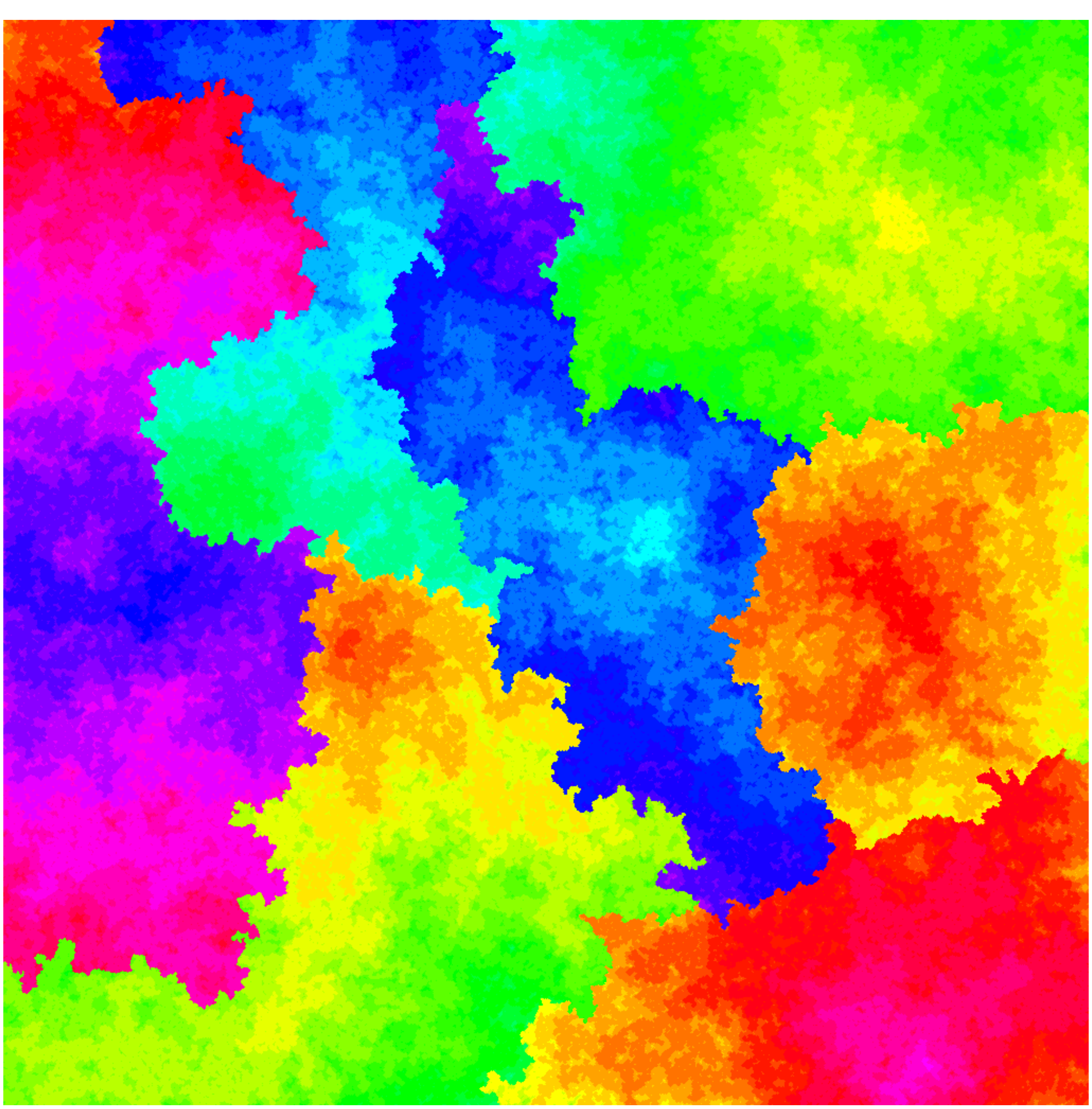

Figure 154: 7 intertwined explorations, $1024\times1024$, site percolation, $p=0.57$.

# 82 Differences in the bond and site models

Let us pause to discuss the differences that arise in the bond and site models with regard to the intertwined explorations. We denote generically by $(\mathcal{A}(t), \mathcal{W}(t))$ the waiting aggregate of an explorer at step $t$, and by $\mathcal{T}(t)$ its taboo set, in both the bond and site models. Let us sum up how we recover the relevant information from the waiting aggregate of the explorers in the bond and the site models:

• **Bond model.** The set $\mathcal{A}(t)$ is a bond aggregate and $\mathcal{W}(t)$ is a subset of the outer boundary of $\mathcal{A}(t)$ (in particular $\mathcal{W}(t)$ is not included in $\mathcal{A}(t)$). The explored edges are those having one endpoint in $\mathcal{A}(t)$. When called upon to continue the exploration, the explorer will look at all the sites of $\mathcal{W}(t)$ which have not been previously visited. He will start exploring the cluster of each of these sites, using the set $\mathcal{T}(t)$ as taboo set.

• **Site model.** The set $\mathcal{A}(t)$ is a site aggregate and $\mathcal{W}(t)$ is a subset of the inner boundary of $\mathcal{A}(t)$ (in particular $\mathcal{W}(t)$ is included in $\mathcal{A}(t)$). The explored sites are those in $\mathcal{A}(t)$. When called upon to continue the exploration, the explorer will look at all the sites which are neighbours to a site of $\mathcal{W}(t)$ and which have not been previously visited. He will start exploring the cluster of each of these sites, using the set $\mathcal{T}(t)$ as taboo set.

We see that there are fundamental differences in the interpretation of the waiting aggregates for the bond and the site models. As was discussed already at the beginning of section 80, the waiting aggregates must be defined differently according to the model, The taboo algorithm must also be adapted to each model, even if the guiding principle for its conception is the same. This being said, there are still a lot of possible choices. For instance, in the site model, we could define an aggregate as a set $\mathcal{A}$ having the following property:

$$\forall x \in \mathcal{A} \qquad \overline{C}(x, D) \subset \mathcal{A}\,.$$

With this definition, an aggregate $\mathcal{A}$ would be closed for the travel time $T_D$ in the sense that, as soon as it contains a site $x$, it contains all the sites which are at null travel time from $x$ in $D$, i.e.,

$$\forall x \in \mathcal{A} \quad \forall y \in D \qquad \big[\quad T_D(x, y) = 0 \quad \Longrightarrow \quad y \in \mathcal{A} \quad\big]\,. \tag{82.1}$$

This solution has some advantages (the definition of the intersection sets would be simpler), but also some drawbacks (the aggregates would not be pairwise disjoint, as they are in the bond model). The main problem is that we would need an additional structure to remember all the sites which have been visited by the explorer. A neat possibility would be to associate three sets to each explorer: the sets of the open sites and the closed sites discovered by the explorer, and the waiting set of the sites that the explorer would start visiting on its next call. For the bond model, we could do the same, by working exclusively with sets of bonds instead of sites. The price we would have to pay is that we would need to work with three sets all the time: the maps Coherent_Site_Explore, Coherent_Bond_Explore would have three arguments instead of two.

## 83 General intertwinings ⊛

We present here a more general scheme of intertwined explorations. Whereas in the scheme of subsection 80.4 the explorers were called in a fixed order one after another, we allow here for a more general calling protocol. This protocol is encoded in a calling function denoted by $\Phi$ and the intertwined explorations associated to it are denoted by

$$\mathcal{A}_i(k),\,\mathcal{W}_i(k),\quad 1\leq i\leq N\,,\quad 0\leq k\leq K\,.$$

We change the index $t$ into $k$, because there is no systematic iteration through the explorers, and $k$ represents now the number of explorer calls performed by the intertwined explorations. The calling function $\Phi$ is defined on $\mathbb{N}$ with values in $\{\,1,\dots,N\,\}$, its meaning is the following: the index of the explorer called at step $k$ is $\Phi(k)$. So the first explorer to be called is the explorer number $\Phi(0)$ and we set

$$\begin{gathered}\mathcal{A}_{\Phi(0)}(0)\,=\,\overline{\text{Clusters}}\big(\{\,w_{\Phi(0)}\,\},D,\varnothing\big)\,,\quad \mathcal{W}_{\Phi(0)}(0)\,=\,\partial_D^{in}\mathcal{A}_{\Phi(0)}(0)\,,\\ \forall i\in\{\,1,\dots,N\,\}\setminus\{\,\Phi(0)\,\}\qquad \mathcal{A}_i(0)\,=\,\mathcal{W}_i(0)\,=\,\varnothing\,.\end{gathered}$$

Let $k\geq 1$ and suppose that $\mathcal{A}_i(j)$, $\mathcal{W}_i(j)$, $1\leq i\leq N$, $0\leq j<k$, have been defined. The explorer called at step $k$ is $\Phi(k)$. The other explorers are stalled during step $k$ and we define

$$\forall i\in\{\,1,\dots,N\,\}\setminus\{\,\Phi(k)\,\}\qquad \mathcal{A}_i(k)\,=\,\mathcal{A}_i(k-1)\,,\quad \mathcal{W}_i(k)\,=\,\mathcal{W}_i(k-1)\,.$$

We distinguish next two cases, depending on whether the explorer $\Phi(k)$ is called for the first time or not at step $k$.

• **First call of explorer number $\Phi(k)$:** we suppose that

$$\forall \ell\in\{\,0,\dots,k-1\,\}\qquad \Phi(\ell)\neq\Phi(k)\,.$$

We define

$$\begin{gathered}\mathcal{A}_{\Phi(k)}(k)\,=\,\overline{\text{Clusters}}\Big(\{\,w_{\Phi(k)}\},D,\bigcup_{\substack{1\leq j\leq N\\ j\neq\Phi(k)}}\mathcal{A}_j(k-1)\Big)\,,\\ \mathcal{W}_{\Phi(k)}(k)\,=\,\partial^{in}_{D\,\setminus\,\mathcal{A}_{\Phi(k)}(k-1)}\mathcal{A}_{\Phi(k)}(k)\,.\end{gathered}$$

• **Subsequent calls of explorer number $\Phi(k)$:** we suppose that

$$\exists\,\ell\in\{\,0,\dots,k-1\,\}\qquad \Phi(\ell)=\Phi(k)\,.$$

We define

$$\begin{pmatrix}\mathcal{A}_{\Phi(k)}(k)\\ \mathcal{W}_{\Phi(k)}(k)\end{pmatrix}\,=\,\\ \text{Coherent_Site_Explore}\left(\begin{pmatrix}\mathcal{A}_{\Phi(k)}(k-1)\\ \mathcal{W}_{\Phi(k)}(k-1)\end{pmatrix},\bigsqcup_{\substack{1\leq j\leq N\\ j\neq\Phi(k)}}\begin{pmatrix}\mathcal{A}_j(k-1)\\ \mathcal{W}_j(k-1)\end{pmatrix}\right).$$

The construction's logic is as follows. When an explorer is called upon, he explores the configuration starting from his waiting set. Naturally, he has already explored the sites in his waiting set, so he starts from the neighbours of these sites that have not yet been explored. During exploration, the set of sites explored by the other explorers is used as the taboo set. Once he has finished exploring, the inner boundary of the explored region becomes his new waiting set. Moreover, to compute this inner boundary, the sites belonging to the taboo set are excluded. As in the other scheme, the process terminates at the random index $K(\Phi)$, when the goal region $G$ has been fully explored, i.e.,

$$K(\Phi)\,=\,\min\Big\{\,k\geq 1:G\subset\bigcup_{1\leq i\leq N}\mathcal{A}_i(k)\,\Big\}\,.$$

This completes the construction of the general multiple intertwined explorations associated to a calling function $\Phi$. In this way, to each calling function $\Phi$, we can associate a sequence of intersection sets $\big(\mathcal{I}^\Phi(k), 0\leq k\leq K(\Phi)\big)$. A possible strategy to vindicate the conjecture would be to study the whole collection of these intersection sets, for all reasonable calling functions $\Phi$. By reasonable, we mean a calling function such that each explorer is called after a controlled interval, so that the termination time of the process is not dramatically affected. To control this termination time, we would use the following adaptation of proposition 80.9.

**Proposition 83.1.** *Let $D$ be a finite domain in $\mathbb{Z}^d$. Let $N\geq 2$ and let $w_1,\dots,w_N$ be distinct sites of $D$. Let $\Phi:\mathbb{N}\to\{\,1,\dots,N\,\}$ be a calling function. Suppose that we are given a site percolation configuration in which $w_1,\dots,w_N$ are pairwise disconnected in $D$. Let*

$$\mathcal{A}_i(k),\,\mathcal{W}_i(k),\quad 1\leq i\leq N\,,\quad 0\leq k\leq K(\Phi)\,,$$

*be the $N$ intertwined explorations starting from $w_1,\dots,w_N$ associated to $\Phi$. For $k\geq 0$, we define $s(k)$ as the largest integer such that, until step $k$ included, each of the $N$ explorers has been called at least $s(k)$ times, i.e.,*

$$s(k)\,=\,\max\Big\{\,s\geq 0:\forall i\in\{\,1,\dots,N\,\}\quad\big|\{\,\ell:0\leq\ell\leq k,\,\Phi(\ell)=i\,\}\big|\geq s\,\Big\}\,.$$

*We have then*

$$\forall k\in\{\,0,\dots,K\,\}\qquad\bigcup_{0\leq s\leq s(k)}\mathit{Shell}\Big(\{\,w_1,\dots,w_N\,\},D,s\Big)\,\subset\,\bigcup_{1\leq i\leq N}\mathcal{A}_i(k)\,.$$

The idea would then be to consider all the intersection sets revealed by the general multiple intertwined explorations associated to a reasonable collection of calling functions. We would try to rule out some specific disconnection events that create too large intersection sets. However, we did not succeed so far along this way. Another further extension would be to consider random calling protocols, which would take into account the current state of the exploration.

# 84 A max-cut lemma

A central tool that will be needed in the proof of theorem 9.2 is a specific max-cut result, that we present in this section. This result is quite general and classical, we have written a version that is tailored for our needs.

**Set collections with bounded multiplicity.** A collection $(E_j, j \in J)$ of subsets of $\mathbb{Z}^d$ has a bounded multiplicity $m \geq 1$ if

$$\forall x \in \mathbb{Z}^d \qquad \Big|\{\, j \in J : x \in E_j \,\}\Big| \,\leq\, m\,.$$

We will use the following inequality, which is derived from a classical double counting argument.

**Lemma 84.1.** *Let $(E_j, j \in J)$ be a collection of finite subsets of $\mathbb{Z}^d$ having bounded multiplicity $m \geq 1$. We have*

$$\sum_{j \in J} |E_j| \,\leq\, m \Big| \bigcup_{j \in J} E_j \Big|\,. \tag{84.1}$$

*Proof.* The hypothesis implies that

$$\forall x \in \mathbb{Z}^d \qquad \sum_{j \in J} 1_{E_j}(x) \,\leq\, m\, 1_{\bigcup_{j\in J} E_j}(x)\,.$$

We sum this inequality over $x$ in $\mathbb{Z}^d$ and we get

$$\sum_{x \in \mathbb{Z}^d} \sum_{j \in J} 1_{E_j}(x) \,\leq\, m \sum_{x \in \mathbb{Z}^d} 1_{\bigcup_{j\in J} E_j}(x) \,=\, m \Big| \bigcup_{j \in J} E_j \Big|\,.$$

We interchange the first two sums and we obtain the inequality (84.1). □

**The max directed cut problem.** This is a classical problem of optimization. A good reference is the book [136] of Vijay Vazirani, see page 23, problem 2.15 in section 2.4 there. We consider the specific case of the complete graph, that we describe next. Let $I$ be a finite set and let $w$ be a weight function defined on the oriented edges of the complete graph over the index set $I$, i.e., the function $w$ is defined on

$$\big\{\, (k, \ell) : k, \ell \in I,\, k \neq \ell \,\big\}$$

and it takes its values in $\mathbb{R}^+$. A partition of $I$ in two sets $K$ and $L$ is called a cut of $I$, and to each cut $K, L$ we associate its weight defined by

$$w(K, L) \,=\, \sum_{k \in K} \sum_{\ell \in L} w(k, \ell)\,.$$

The max directed cut problem consists in computing the maximum

$$\max \big\{\, w(K, L) : K, L \text{ cut of } I \,\big\}\,, \tag{84.2}$$

and more generally, in finding and studying the cuts $K, L$ which realize the maximum of the function $w(K, L)$. As long as the set $I$ is finite, the maximum is well-defined and optimal cuts do exist. We define the total weight $w(I)$ of the graph as

$$w(I) = \sum_{k,\ell \in I, k \neq \ell} w(k, \ell) \, .$$

The maximum (84.2) is larger than or equal to one quarter of the total weight $w(I)$ of the graph. For completeness, we state and prove this classical result in the next lemma.

**Lemma 84.2.** *If* $(K^*, L^*)$ *is a cut realizing the maximum* (84.2)*, then*

$$w(K^*, L^*) \geq \frac{1}{4} w(I) \, . \tag{84.3}$$

*Proof.* We present two arguments. The first one corresponds to the outcome of a local search algorithm (see [136], page 23, exercise 2.2 of section 2.4). Let $(K^*, L^*)$ be a cut realizing the maximum (84.2). We cannot increase its weight, so any cut obtained from $(K^*, L^*)$ by moving one vertex from one set of the cut to the other will have a weight less than or equal to $w(K^*, L^*)$. Let $k \in K^*$. We have thus

$$w\big(K^* \setminus \{ k \}, L^* \cup \{ k \}\big) \leq w(K^*, L^*) \, . \tag{84.4}$$

Inequality (84.4) can be rewritten as

$$\sum_{\ell \in L^*} w(k, \ell) \geq \sum_{h \in K^* \setminus \{ k \}} w(h, k) \, . \tag{84.5}$$

Summing inequality (84.5) over $k$ in $K^*$, we get

$$w(K^*, L^*) \geq \sum_{k \in K^*} \sum_{h \in K^* \setminus \{ k \}} w(h, k) = \sum_{k,h \in K^*, \, k \neq h} w(h, k) \, . \tag{84.6}$$

By exchanging the roles of $K^*$ and $L^*$, or by performing the same computation starting with an element of $L^*$, we obtain the symmetric inequality

$$w(K^*, L^*) \geq \sum_{h,\ell \in L^*, \, h \neq \ell} w(\ell, h) \, . \tag{84.7}$$

We finally sum the two inequalities (84.6) and (84.7), and we use the fact that $w(K^*, L^*) \geq w(L^*, K^*)$ to conclude that

$$4w(K^*, L^*) \geq \sum_{k,h \in K^*, \, k \neq h} w(h, k) + \sum_{h,\ell \in L^*, \, h \neq \ell} w(\ell, h) + w(K^*, L^*) + w(L^*, K^*) \, .$$

The right-hand quantity is precisely equal to $w(I)$ and we obtain the desired inequality (84.3).

The second argument is a simple probabilistic argument (see [136], page 138, exercise 16.6 of section 16.5). We create a random cut $K, L$ by assigning independently each element of $I$ to one of the two sets $K$ or $L$ with equal probability 1/2. We compute then the expectation of $w(K,L)$:

$$E\big(w(K,L)\big) \;=\; \sum_{k,\ell\in I,k\neq\ell} w(k,\ell)P\big(k\in K,\ell\in L\big) \;=\; \frac{1}{4}w(I)\,.$$

Since the expectation of $w(K,L)$ is equal to $w(I)/4$, certainly there exists a choice for the sets $K,L$ such that $w(K,L)$ is larger than or equal to $w(I)/4$, hence the inequality (84.3). □

We are now ready to state the max-cut result that we will rely on. Notice that we already made appeal to this result in the proof of lemma 33.2 (which was a side result, that does not intervene in the proof of theorem 9.2).

**Lemma 84.3.** *Let $D$ be a finite subset of $\mathbb{Z}^d$ and let $(E_i, i\in I)$, $(F_i, i\in I)$ be two collection of subsets of $D$ such that the collection*

$$E_i\cap F_j\,,\quad i,j\in I\,,\quad i\neq j\,,$$

*has bounded multiplicity $m\geq 1$. There exist two disjoint subsets $K,L$ of $I$ such that*

$$\Big|\Big(\bigcup_{k\in K} E_k\Big)\cap\Big(\bigcup_{\ell\in L} F_\ell\Big)\Big| \;\geq\; \frac{1}{4m}\Big|\bigcup_{k,\ell\in I,\,k\neq\ell} E_k\cap F_\ell\Big|\,. \tag{84.8}$$

*Proof.* The problem of finding the subsets $K$ and $L$ of $I$ can be recast as a max directed cut problem. We consider the complete graph on the index set $I$, and we define a weight function $w$ on the edges of this graph as follows:

$$\forall k,\ell\in I\qquad w(k,\ell)\;=\;\big|E_k\cap F_\ell\big|\,. \tag{84.9}$$

By lemma 84.2, there exist two disjoint subsets $K,L$ of $I$ such that

$$w(K,L)\;\geq\;\frac{1}{4}w(I)\,. \tag{84.10}$$

Coming back to the definition (84.9) of the weight $w(k,\ell)$, the inequality (84.10) can be rewritten as

$$\sum_{k\in K}\sum_{\ell\in L}\big|E_k\cap F_\ell\big|\;\geq\;\frac{1}{4}\sum_{k,\ell\in I,k\neq\ell}\big|E_k\cap F_\ell\big|\,. \tag{84.11}$$

By the simple union bound, we have

$$\Big|\bigcup_{k,\ell\in I,\,k\neq\ell} E_k\cap F_\ell\Big|\;\leq\;\sum_{k,\ell\in I,k\neq\ell}\big|E_k\cap F_\ell\big|\,. \tag{84.12}$$

We apply lemma 84.1 to the collection of sets $E_i\cap F_j$, $i,j\in I$, and we obtain

$$\sum_{k\in K}\sum_{\ell\in L}\big|E_k\cap F_\ell\big|\;\leq\;m\Big|\Big(\bigcup_{k\in K} E_k\Big)\cap\Big(\bigcup_{\ell\in L} F_\ell\Big)\Big|\,. \tag{84.13}$$

Combining inequalities (84.11), (84.12) and (84.13), we obtain (84.8). □

# 85 Ignition of the proof

We start here the proof of theorem 9.2. The proof involves several results proved in the previous sections. The central object that we analyze is a certain collection of boundary clusters. Let us recall briefly the relevant definitions and results.

For convenience, we will work within the box $\Lambda(2n)$ and we will study the trace on $\Lambda(2n)$ of the clusters of the larger box $\Lambda(4n)$. The advantage of this choice is that the distance between $\Lambda(2n)$ and the boundary $\partial^{\,in}\Lambda(4n)$ of $\Lambda(4n)$ is equal to $n$ (rather than $n/2$ if we were working with $\Lambda(n)$ and $\Lambda(2n)$). This not only saves the factor 2 in the denominator of several formulas, but it also avoids irrelevant discussions about the parity of $n$. Needless to say, everything could be written for $\Lambda(n)$ and $\Lambda(2n)$. Throughout the proof, we will look at the percolation configuration restricted to the box $\Lambda(4n)$. To alleviate the formulas, in some places we use $D$ to denote $\Lambda(4n)$. We would also like to point out that several results are valid and stated for an arbitrary finite domain $D$ of $\mathbb{Z}^d$.

The strategy of the proof is to study the collection of the clusters which intersect both $\Lambda(2n)$ and $\partial^{\,in}\Lambda(4n)$, and to show that there are not too many of them. Let us first introduce a notation for the clusters that intersect the boundary of a box, or rather a general finite domain $D$.

**Definition 85.1.** *Let $D$ be a finite subset of $\mathbb{Z}^d$. A cluster of the percolation configuration restricted to $D$ is a boundary cluster if it intersects $\partial^{\,in}D$. We denote by $\mathcal{C}_{bd}(D)$ the collection of the boundary clusters of the percolation configuration restricted to $D$.*

To alleviate the formulas, for $n\geq 1$, we write simply $\mathcal{C}_{\mathrm{bd}}(n)$ instead of $\mathcal{C}_{\mathrm{bd}}\big(\Lambda(n)\big)$. Let us fix $\varepsilon>0$. We consider the event $\mathcal{F}_{\mathrm{bd}}(n,\varepsilon)$ defined by

$$\mathcal{F}_{\mathrm{bd}}(n,\varepsilon)\,=\,\Big\{\sum_{C\in\mathcal{C}_{\mathrm{bd}}(4n)}\big|C\cap\Lambda(2n)\big|\,\geq\,\big(\theta(p)-\varepsilon\big)\big|\Lambda(2n)\big|\Big\}\,. \tag{85.1}$$

With the help of the spatial ergodic theorem, it was proved in sections 8.5 and 10 that

$$\lim_{n\to\infty}P\big(\mathcal{F}_{\mathrm{bd}}(n,\varepsilon)\big)\,=\,1\,.$$

This is the starting point of the proof. The argument of the proof requires also a control on the sizes of the finite clusters of $\Lambda(4n)$. These clusters were loosely defined in subsection 33.6 for bond percolation. Let us give a formal definition, still for a general finite domain $D$.

**Definition 85.2.** *Let $D$ be a finite subset of $\mathbb{Z}^d$. A cluster of the percolation configuration restricted to $D$ is called finite if it does not intersect $\partial^{\,in}D$. We denote by $\mathcal{C}_{fi}(D)$ the collection of the finite clusters of the percolation configuration restricted to $D$.*

The point is that the finite clusters are left unchanged if we look at the whole percolation configuration in $\mathbb{Z}^d$ instead of the configuration restricted to $D$. To alleviate the formulas, for $n\geq 1$, we write simply $\mathcal{C}_{\mathrm{fi}}(n)$ instead of $\mathcal{C}_{\mathrm{fi}}\big(\Lambda(n)\big)$.

Let $d \geq 3$ and let $p$ in $]0,1[$ be a parameter such that $\theta(p,\mathbb{Z}^d) > 0$. Let us fix an integer $M > 0$. We consider the event $\mathcal{F}_{\text{fi}}(n,M)$ defined by

$$\mathcal{F}_{\text{fi}}(n,M) \,=\, \big\{\, \forall C \in \mathcal{C}_{\text{fi}}(4n) \quad |C| < M \,\big\}\,.$$

It follows from the hypothesis (9.16) that

$$\exists\, \alpha > 0 \quad \exists M_0 \geq 0 \quad \forall M \geq M_0 \qquad P_p\big(M \leq |C(0)| < +\infty\big) \,\leq\, \exp(-M^\alpha)\,. \tag{85.2}$$

We use this estimate and a standard union bound to control the probability of the event $\mathcal{F}_{\text{fi}}(n,M)$:

$$\forall n \geq 1 \quad \forall M \geq M_0 \qquad P\big(\mathcal{F}_{\text{fi}}(n,M)^c\big) \,\leq\, \big|\Lambda(4n)\big| \exp(-M^\alpha)\,. \tag{85.3}$$

Using the simple bound $|\Lambda(4n)| \leq (4n+1)^d \leq n^{3d}$, valid for $n \geq 3$, and setting

$$M(n) \,=\, \big\lceil (4d \ln n)^{1/\alpha} \big\rceil\,, \tag{85.4}$$

we conclude from (85.3) that there exists an integer $n_0(M_0,\alpha,p,d)$ depending on $M_0,\alpha,p,d$ such that

$$\forall n \geq n_0(M_0,\alpha,p,d) \qquad P\big(\mathcal{F}_{\text{fi}}(n,M(n))^c\big) \,\leq\, \frac{1}{n}\,. \tag{85.5}$$

The next important result that comes into play is the one concerning the travel time to the boundary obtained in section 16. Proposition 16.1 states that there exists a non-decreasing function $\beta(n)$ defined for $n \geq 3$ satisfying $\beta(n) \leq \sqrt{\ln n}$ for $n \geq 3$, $\lim_{n\to+\infty} \beta(n) = +\infty$, and such that, if we define the event $\mathcal{G}(n)$ as

$$\mathcal{G}(n) \,=\, \Big\{\, \forall x \in \Lambda(4n) \quad T_{\Lambda(4n)}\big(x, \partial^{\,in}\Lambda(4n)\big) \,\leq\, \frac{\ln n}{\beta(n)} \,\Big\}\,, \tag{85.6}$$

then we have

$$\forall n \geq 3 \qquad P\big(\mathcal{G}(n)\big) \,\geq\, 1 - n^{3d-\beta(n)}\,. \tag{85.7}$$

We use next the hypothesis (9.16) in order to obtain a quantitative estimate on $\beta(n)$. The function $\beta(n)$ was defined in (16.3) in the course of the proof of proposition 16.1, by setting

$$\forall n \geq 3 \qquad \beta(n) \,=\, \min\Big(\sqrt{\ln n}, \sqrt{\alpha\big(\sqrt{\ln n}\big)}\Big)\,. \tag{85.8}$$

In turn, the function $\alpha(t)$ was introduced in corollary 14.5, where it was proved that, for $t$ larger than a value $t_0\big(\theta(p),p,d\big)$ depending on $\theta(p),p,d$, we have

$$\alpha(t) \,\geq\, \ln\Big(\frac{\theta(p)}{2}\Big) - \frac{1}{2}\ln\Bigg( -\,\phi\Big(P_p\Big(\Big(\frac{t}{4d}\Big)^{d/(d+1)} \!<\! |C(0)| \!<\! +\infty\Big)\Big)\Bigg)\,. \tag{85.9}$$

Inequality (85.9) and the tail estimate (85.2) imply that there exists $\alpha > 0$ such that

$$\forall t \geq t_0\big(\theta(p),p,d\big) \qquad \alpha(t) \,\geq\, t^{\frac{\alpha d}{3(d+1)}}\,. \tag{85.10}$$

Substituting inequality (85.10) into (85.8), we conclude that there exist $\gamma$ in $]0,1/2[$ and a value $n_1\big(\theta(p),\alpha,p,d\big)$ depending on $\theta(p),\alpha,p,d$ such that

$$\forall n \geq n_1\big(\theta(p),\alpha,p,d\big) \qquad \beta(n) \,\geq\, (\ln n)^{\gamma}\,. \tag{85.11}$$

This lower bound on $\beta(n)$ will play a crucial role in the proof. Putting together the estimates (85.7) and (85.11), we conclude that, for $n \geq n_1\big(\theta(p),\alpha,p,d\big)$,

$$P\Big(\forall x \in \Lambda(4n) \quad T_{\Lambda(4n)}\big(x,\partial^{\,in}\Lambda(4n)\big) \,\leq\, (\ln n)^{1-\gamma}\Big) \,\geq\, 1 - n^{3d-(\ln n)^{\gamma}}\,.$$

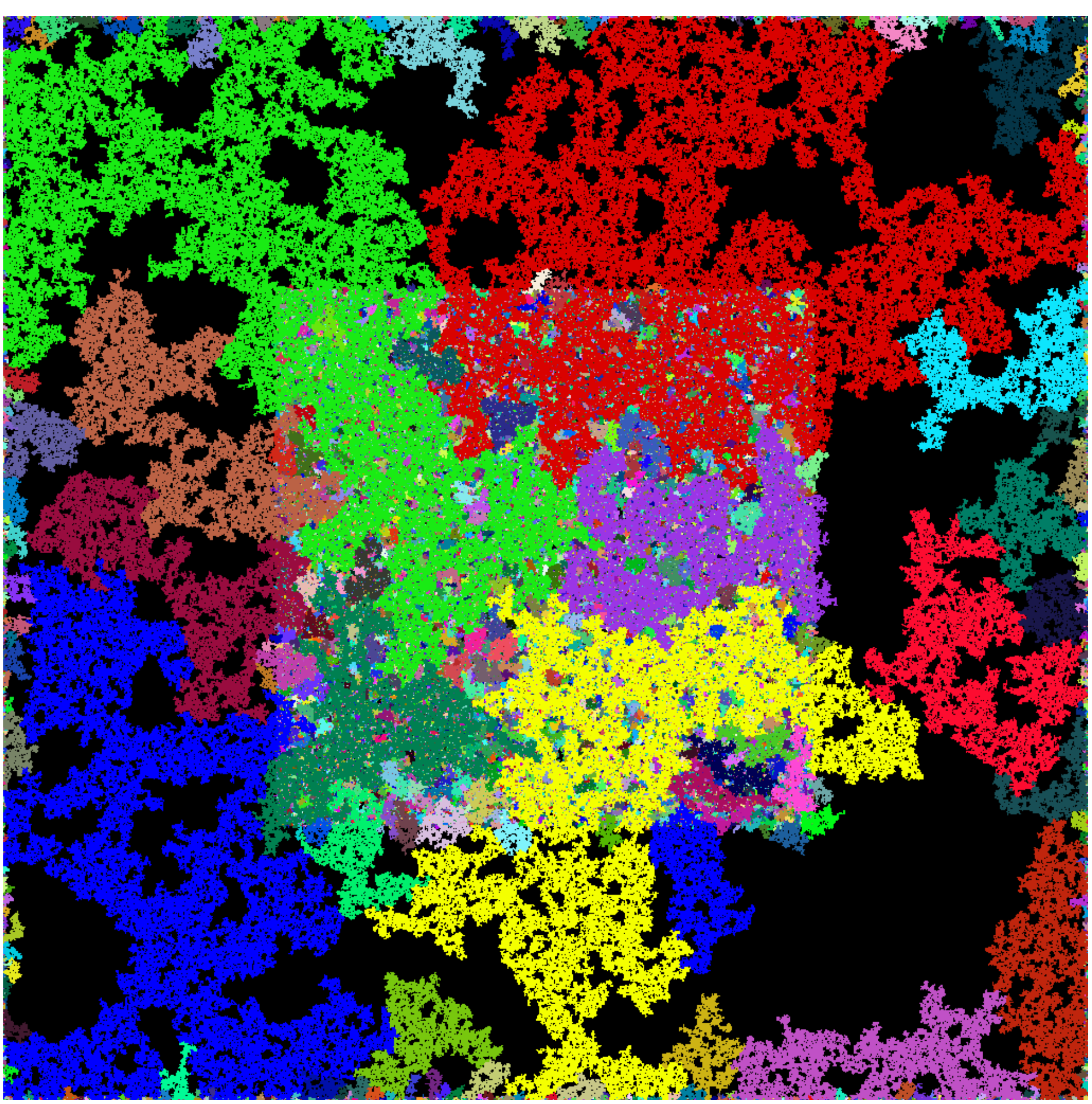

Figure 155: Bond percolation, p=0.496, 1024x1024 box, crossing clusters. The clusters touching $\Lambda(512)$ or the boundary of $\Lambda(1024)$ are colored.

# 86 Outline of the proof

Let us give a loose description of the core of our argument to prove theorem 9.2. We suppose that $p$ is a parameter such that $\theta(p, \mathbb{Z}^d) > 0$, and we fix $\varepsilon > 0$. With this hypothesis, we can prove that, for $n$ large, with high probability, the collection of the boundary clusters of the box $\Lambda(4n)$ fills a fraction larger than $\theta(p) - \varepsilon$ of the volume of the half box $\Lambda(2n)$. We wish to prove that, with high probability, at least one cluster of this collection has a trace on $\Lambda(2n)$ which has a cardinality larger than $(\theta(p) - 2\varepsilon)|\Lambda(2n)|$. To do so, we proceed in two steps. In the first step, we prove that the number of the boundary clusters of $\Lambda(4n)$ which reach $\Lambda(2n)$ is at most of order $n^{d-2}$. In the second step, we prove that it is unlikely that the collection of the boundary clusters of $\Lambda(4n)$ can be partitioned into two subcollections which both cover a positive fraction of the volume of $\Lambda(2n)$. In both steps, we show that the scenarios that differ from those described above imply the existence of an excess of closed sites. This excess of closed sites is identified as the set of the closed relevant taboo sites of an adequate process of multiple intertwined explorations. The ultimate mechanism to control the size of this set is quite rudimentary and rests on Hoeffding's inequality. In fact, we will show that this set is located in the extremities of the intersection edges of the multiple intertwined explorations starting from the boundary clusters. By modifying slightly the taboo explorations associated to this process, we prove that this set is included in the set discovered by two distinct taboo explorations. The trouble is that the sequence of the states of the sites discovered by a taboo exploration is not a naive sequence of i.i.d. Bernoulli random variables. It would be the case if, say, the taboo sets of each exploration were fixed deterministic sets. However, in the multiple intertwined explorations, each explorer uses as taboo set the set of the sites discovered so far by all the other explorers. Naturally, this creates a huge problem, and it is this problem that we never managed to solve. As soon as we provide an explorer with some partial information on the sites discovered by the other explorers, this creates a bias, and it affects the distribution of the remaining sites. When an explorer is called, we would like to tell him which sites in the neighborhood of its waiting set are forbidden. Unfortunately, this apparent innocent information simultaneous tells that the other sites, which are allowed for the exploration, are not connected to the aggregates of the other explorers. This is negative information, which is positively correlated with the existence of large closed interfaces, and it is very hard to control its impact. In the framework of the bond model, we tried with the help of a brutal combinatorial control on the relevant taboo sites. This brutal control led to some non-trivial results, but they were not strong enough to reach the desired conclusions.

We try here a radically different approach. Instead of preventing the explorers to enter into the territory discovered by the others, we will implement a stopping rule to avoid entering too much inside this forbidden territory. Here is how it works. Let $M \geq 1$ be a fixed integer. Imagine that there are $N$ boundary clusters, and that these clusters are the only clusters of $\Lambda(4n)$ having cardinality larger than $M$. We pick up $N$ sites $w_1, \dots, w_N$, such that exactly one site belongs to each of the $N$ boundary clusters. We start the multiple intertwined explorations from $w_1, \dots, w_N$. During the first iteration of the multiple intertwined explorations, the $N$ boundary clusters will be explored. In the subsequent iterations, the remaining clusters which are discovered will have a cardinality strictly smaller than $M$. From the second iteration onwards, we forget the taboo sets, instead we use the following stopping rule: as soon as an explorer sees a connected set of open sites of cardinality $M$, he will abandon its exploration, but he will not forget what he has discovered. The benefits of this mechanism are twofold. On one hand, the rule is a genuine stopping rule, and it does not impact the distribution of the sites that are being explored. On the other hand, it permits to avoid a deep penetration into the territories that should have stayed forbidden. The parameter $M$ will be chosen carefully, so that the deterioration of the estimates due to the repeated incursions into the taboo sets remains limited. The drawback is that, in order to be implemented successfully, this strategy requires an adequate control over the size of the finite clusters. The precise details are the object of the subsequent sections.

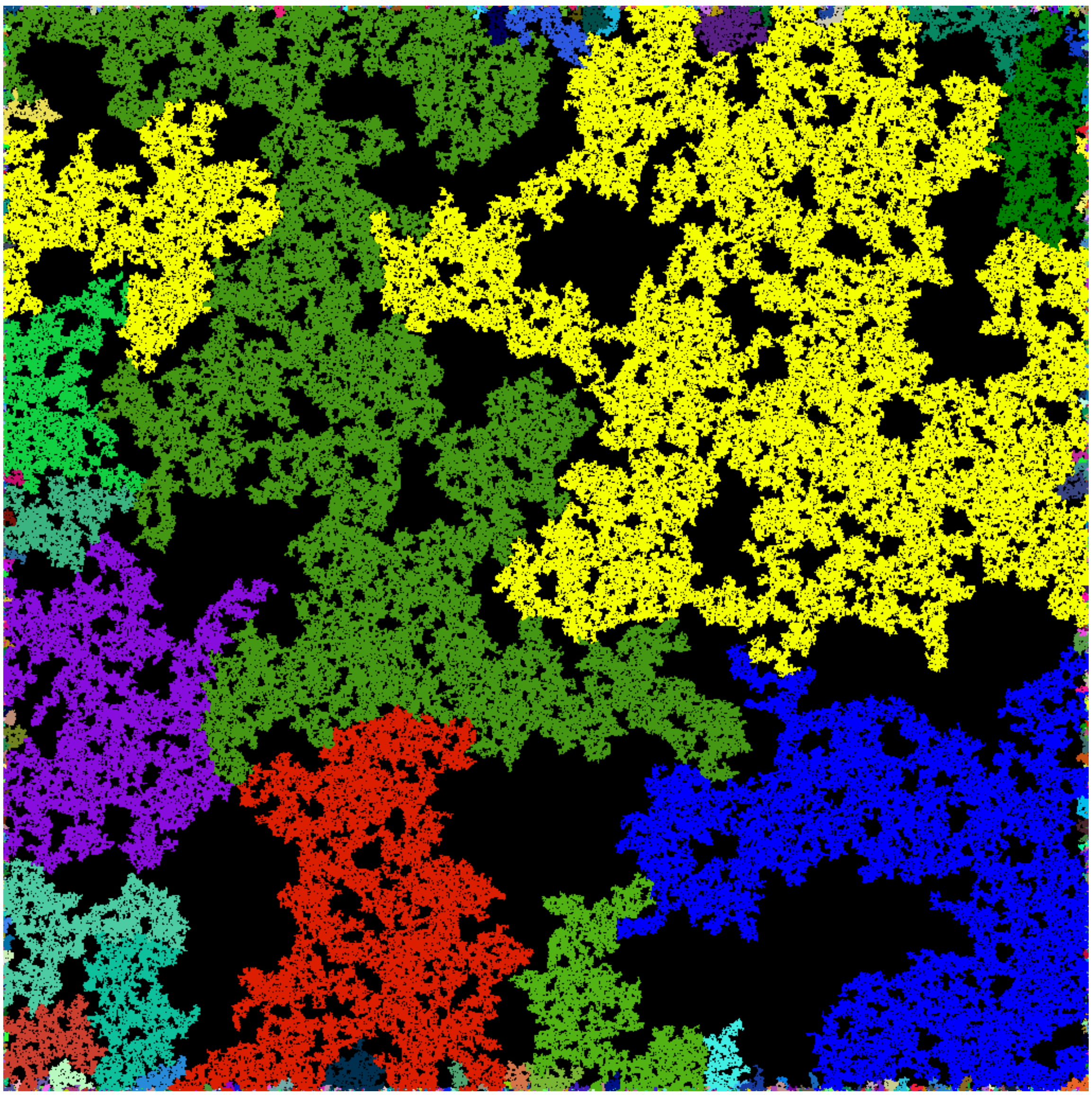

Figure 156: Bond percolation, p=0.498, 1024x1024 box, the boundary clusters.

# 87 The number of annulus-crossing clusters

Let us denote by $\mathcal{N}^*$ the number of boundary clusters in $\Lambda(4n)$ which intersect the box $\Lambda(2n)$, i.e.,

$$\mathcal{N}^* = \left|\left\{ C \in \mathcal{C}_{\mathrm{bd}}(4n) : C \cap \Lambda(2n) \neq \varnothing \right\}\right|.$$

Our main endeavour will be to obtain a quantitative control on $\mathcal{N}^*$. Obviously, we have

$$\mathcal{N}^* \leq \left|\partial^{in}\Lambda(4n)\right| \leq 2d(4n+1)^{d-1},$$

but we would like to obtain a better polynomial bound. We are in fact returning to the crucial issue explained in section 8.6.

## 87.1 The strategy of the proof

Let us try to explain our strategy before starting the precise proof. Suppose that $\mathcal{N}^*$ is large. This means that there are a lot of disjoint clusters in the box $\Lambda(4n)$ which cross the annulus $\Lambda(4n) \setminus \Lambda(2n)$. As $\theta(p)$ is positive, it is not difficult to have a crossing, yet to have disjoint crossing clusters requires that they are not connected together. In their influential paper [84], Kesten and Zhang obtained a very efficient control on the probability of having two disjoint crossing clusters with the help of the slab technology. There is no hope of proving a similar result at the critical point, or in its vicinity. However, if $\mathcal{N}^*$ is large, say of polynomial order, then the presence of $\mathcal{N}^*$ disjoint crossing clusters entails the existence of many closed sites. Namely, each crossing cluster must be surrounded by a set of closed sites which prevents him from being connected to the other crossing clusters. We are going to estimate the probability of this scenario. Of course, the closed surrounding sets might be quite far from each other, and there might be a lot of different choices for these sets. This makes our task very complicated. Our technique consists in picking one vertex in each annulus-crossing cluster, and launching multiple intertwined explorations from these vertices. Upon termination, the intertwined explorations will have revealed a large set of intersection edges, having cardinality larger than $n\mathcal{N}^*$. As we saw in proposition 80.18, each intersection edge has at least one endpoint which is a closed relevant taboo site, thus we will observe at least $n\mathcal{N}^*/(2d)$ relevant taboo sites.

The proof is divided intro three main stages. The first stage is carried out in section 88. We select an iteration $t$ during which there is a big increase in the number of relevant taboo sites. With the help of a max-cut result, we show that the relevant taboo sites for the iteration $t$ can be localized by a pair of taboo explorations. In the second stage, we define and study what we call the genuine truncated explorations. We show in section 89 that the taboo explorations of the first stage can be approximated by specific genuine truncated explorations. The final third stage is carried out in section 90, where we compute quantitative estimates on the number of the relevant taboo sites. The whole argument is completed and recapitulated in subsection 91.1.

## 87.2 Ignition of the first stage of the proof

We wish to obtain a quantitative control on $\mathcal{N}^*$ on the event $\mathcal{F}_{\text{fi}}(n,M)\cap\mathcal{G}(n)$. Although only the events $\mathcal{F}_{\text{fi}}(n,M)$ and $\mathcal{G}(n)$ are relevant to obtain the control on $\mathcal{N}^*$, we will require that the event $\mathcal{F}_{\text{bd}}(n,\varepsilon)$ occurs as well, this will allow us to conclude something on the maximal cluster in $\Lambda(4n)$. To alleviate the formulas, we set

$$\mathcal{E} \,=\, \mathcal{F}_{\text{fi}}(n,M)\cap\mathcal{G}(n)\cap\mathcal{F}_{\text{bd}}(n,\varepsilon)\,. \tag{87.1}$$

We try here to put into action the strategy described in subsection 87.1 and we start now with the detailed proof. We denote by $\mathcal{N}_{\text{bd}}$ the number of boundary clusters in $\Lambda(4n)$ and by $\mathcal{N}^*$ the number of clusters intersecting simultaneously $\partial^{\,in}\Lambda(4n)$ and $\Lambda(2n)$, i.e.,

$$\begin{aligned}\mathcal{N}_{\text{bd}} \,&=\, \big|\mathcal{C}_{\text{bd}}(4n)\big|\,,\\ \mathcal{N}^* \,&=\, \Big|\big\{\,C\in\mathcal{C}_{\text{bd}}(4n): C\cap\Lambda(2n)\neq\varnothing\,\big\}\Big|\,.\end{aligned}$$

Let us fix an integer $N^*\geq 2$. We wish to estimate $P\big(\mathcal{N}^*\geq N^*,\mathcal{E}\big)$. We start by conditioning on the value of $\mathcal{N}_{\text{bd}}$:

$$P\big(\mathcal{N}^*\geq N^*,\mathcal{E}\big) \,=\, \sum_{N^*\leq N\leq|\partial^{\,in}\Lambda(4n)|} P\big(\mathcal{N}^*\geq N^*,\mathcal{E},\mathcal{N}_{\text{bd}}=N\big)\,. \tag{87.2}$$

Let us fix next $N$ in $\{\,N^*,\dots,|\partial^{\,in}\Lambda(4n)|\,\}$ and let us suppose that $\mathcal{N}_{\text{bd}}=N$. We pick up $N$ sites $w_1,\dots,w_N$ in $\partial^{\,in}\Lambda(4n)$, such that exactly one site belongs to each of the $N$ boundary clusters. We decompose further the probability in the sum (87.2) according to the possible choices of $w_1,\dots,w_N$, and we use a standard union bound to get

$$\begin{aligned}P\big(\mathcal{N}^*\geq N^*,\mathcal{E},\mathcal{N}_{\text{bd}}=N\big) \,&=\, P\begin{pmatrix}\mathcal{N}^*\geq N^*,\mathcal{E},\exists\, w_1,\dots,w_N\in\partial^{\,in}\Lambda(4n)\\ \mathcal{C}_{\text{bd}}(4n)=\big\{\,C\big(w_i,\Lambda(4n)\big),1\leq i\leq N\,\big\}\end{pmatrix}\\ \leq \sum_{w_1,\dots,w_N\in\partial^{\,in}\Lambda(4n)} & P\Big(\mathcal{N}^*\geq N^*,\mathcal{E},\mathcal{C}_{\text{bd}}(4n)=\big\{\,C\big(w_i,\Lambda(4n)\big),1\leq i\leq N\,\big\}\Big)\,.\end{aligned} \tag{87.3}$$

From now onwards, we fix the $N$ sites $w_1,\dots,w_N$ in $\partial^{\,in}\Lambda(4n)$, and we focus on the event appearing in the last sum above, which we denote by $\mathcal{H}$:

$$\mathcal{H} \,=\, \mathcal{E}\cap\big\{\,\mathcal{N}^*\geq N^*,\mathcal{C}_{\text{bd}}(4n)=\big\{\,C\big(w_i,\Lambda(4n)\big),1\leq i\leq N\,\big\}\,\big\}\,. \tag{87.4}$$

We consider a configuration of percolation belonging to $\mathcal{H}$ and we launch the multiple intertwined explorations starting from $w_1,\dots,w_N$. The goal set $G$ is taken to be the whole box $\Lambda(4n)$. The event $\mathcal{G}(n)$ and the inequality (85.11) guarantee that the termination time $T$ is at most $(\ln n)^{1-\gamma}+1$, hence smaller than $2(\ln n)^{1-\gamma}$ for $n\geq 3$.

## 87.3 From $\mathcal{N}^*$ to $\overline{\mathcal{T}}^0(T)$

Let us consider the final aggregates $(\mathcal{A}_i(T), 1 \leq i \leq N)$ produced by the intertwined explorations. On one hand, they form a partition of $\Lambda(4n)$. On the other hand, each aggregate contains exactly one boundary cluster. We know also that there are at least $N^*$ boundary clusters that intersect $\Lambda(2n)$. Let us consider a boundary cluster $C$ in $\Lambda(4n)$ that intersects the box $\Lambda(2n)$. This boundary cluster meets one of the $2d$ faces of $\Lambda(4n)$, say the face $F(\Lambda(4n), v)$ where $v$ or $-v$ is one of the vectors of the canonical basis of $\mathbb{R}^d$,

$$u_i = (0, \dots, 0, 1, 0, \dots, 0)\,, \quad 1 \leq i \leq d$$

(the notation $F(\Lambda(4n), v)$ was defined in (18.1)). Furthermore, the cluster $C$ contains an open path inside $\Lambda(4n)$ leading to a site of $\Lambda(2n)$. Suppose for instance that $v = u_i$ for some $i$ in $\{ 1, \dots, d \}$. Then this path has to cross the parallelepiped

$$\begin{aligned}
\Gamma_i(2n) \;&=\; \Lambda(4n) \cap \big(3nu_i + \Lambda(4n)\big) \\
&=\; \big\{ \big(x(1), x(2), \cdots, x(d)\big) \in \Lambda(4n) : n \leq x(i) \leq 2n \,\big\}
\end{aligned}$$

in the direction of the $i$-th axis, meaning that it contains a site $x$ such that $x(i) = 2n$ and a site $y$ such that $y(i) = n$. Notice that the parallelepiped $\Gamma_i(2n)$ was introduced in section 9 to define the famous good event, albeit for a different purpose. Of course, the face of $\Lambda(4n)$ varies from one boundary cluster $C$ to another. Since there are at least $N^*$ such boundary clusters but only $2d$ faces, then there exists a face and a parallelepiped shared by at least $N^*/(2d)$ clusters. Suppose that this common face corresponds to $u_1$, and let us denote by $\mathcal{C}_1^*(n)$ the subcollection of the boundary clusters which intersect both $F(\Lambda(4n), u_1)$ and $\Lambda(2n)$. For $s$ in $\mathbb{Z}^d$, we denote by $H(s)$ the hyperplane $\{ s \} \times \mathbb{Z}^{d-1}$, i.e.,

$$\forall s \in \mathbb{Z} \qquad H(s) \,=\, \big\{ \big(s, x(2), \cdots, x(d)\big) : x(2), \dots, x(d) \in \mathbb{Z} \,\big\}\,.$$

Let us fix $s$ in $\{ n, \dots, 2n \}$. The collection $\mathcal{C}_1^*(n)$ is such that

$$\forall C \in \mathcal{C}_1^*(n) \qquad C \cap H(s) \,\neq\, \varnothing\,.$$

In particular, the traces of the aggregates $(\mathcal{A}_i(T), 1 \leq i \leq N)$ on $H(s)$ form a partition of $\Lambda(4n) \cap H(s)$ containing at least $N^*/(2d)$ non-empty sets. Based on this information, the next lemma provides a lower bound on the number of intersection edges included in $H(s)$.

**Lemma 87.1.** *Let $\Lambda$ be a $(d-1)$-dimensional cubic box, and let $\pi = (\pi_i, 1 \leq i \leq r)$ be a partition of $\Lambda$ in $r$ non-empty sets, where $r \geq 2$. We have*

$$\Big| \bigcup_{1 \leq j < k \leq r} \big(\Delta_\Lambda \pi_j \cap \Delta_\Lambda \pi_k\big) \Big| \,\geq\, \frac{r}{2}\,.$$

*Proof.* Let $i$ belong to $\{ 1, \dots, r \}$. Since $r \geq 2$, the set $\pi_i$ does not contain the whole box $\Lambda$, therefore $\partial^{\,out} \pi_i \cap \Lambda$ is not empty, and thus

$$\exists \ell \in \{ 1, \dots, r \} \setminus \{ i \} \quad \exists x \in \pi_i \quad \exists y \in \pi_\ell \qquad |x - y| = 1\,.$$

Therefore the edge $\{x, y\}$ is in $\Delta_\Lambda \pi_i \cap \Delta_\Lambda \pi_\ell$. Of course, there might be many choices for the index $\ell$. We select one of them, for instance the smallest, and we denote it by $\phi(i)$. We have then

$$\Big| \bigcup_{1\leq j<k\leq r} \big(\Delta_\Lambda \pi_j \cap \Delta_\Lambda \pi_k\big)\Big| \,\geq\, \frac{1}{2}\Big| \bigcup_{1\leq i\leq r} \big(\Delta_\Lambda \pi_i \cap \Delta_\Lambda \pi_{\phi(i)}\big)\Big| \,\geq\, \frac{r}{2}\,.$$

We have to divide by 2 because an edge $\{x, y\}$ can belong to two sets of the collection $\Delta_\Lambda \pi_i \cap \Delta_\Lambda \pi_{\phi(i)}$, $1 \leq i \leq r$. □

We apply lemma 87.1 to the $(d-1)$-dimensional box $\Lambda(4n) \cap H(s)$ and the partition induced by the traces of the aggregates $(\mathcal{A}_i(T), 1 \leq i \leq N)$ on $H(s)$, and we obtain

$$\Big| \bigcup_{1\leq j<k\leq N} \big(\Delta_{\Lambda(4n)\cap H(s)} \mathcal{A}_j(T) \cap \Delta_{\Lambda(4n)\cap H(s)} \mathcal{A}_k(T)\big)\Big| \,\geq\, \frac{N^*}{4d}\,. \tag{87.5}$$

The hyperplanes $H(s)$, $n \leq s \leq 2n$, are pairwise disjoint. Summing the inequality (87.5) over $s$ in $\{n, \dots, 2n\}$, we conclude that

$$\Big| \bigcup_{0\leq t\leq T-1} \mathcal{I}(t)\Big| \,\geq\, \Big| \bigcup_{1\leq j<k\leq N} \big(\Delta_{\Lambda(4n)} \mathcal{A}_j(T) \cap \Delta_{\Lambda(4n)} \mathcal{A}_k(T)\big)\Big| \,\geq\, \frac{N^* n}{4d}\,. \tag{87.6}$$

For $t$ in $\{0, \dots, T\}$, we shall denote by $\overline{\mathcal{T}}^0(t)$ the set of all the closed relevant taboo sites until time $t$, i.e.,

$$\forall t \in \{0, \dots, T\} \qquad \overline{\mathcal{T}}^0(t) \,=\, \bigcup_{0\leq s\leq t} \mathcal{T}^{*0}(s)\,.$$

Using corollary 80.19, we deduce from (87.6) that, on the event $\mathcal{H}$,

$$\big|\overline{\mathcal{T}}^0(T)\big| \,\geq\, \frac{N^* n}{16 d^2}\,. \tag{87.7}$$

We have proved that any configuration of the event $\mathcal{H}$ defined in (87.4) satisfies the inequality (87.7), whence

$$\begin{aligned} P(\mathcal{H}) \,&=\, P\left(\begin{matrix} \mathcal{N}^* \geq N^*, \mathcal{E}, \big|\overline{\mathcal{T}}^0(T)\big| \,\geq\, \dfrac{N^* n}{16d^2}\,, \\ \mathcal{C}_{\mathrm{bd}}(4n) = \big\{\, C\big(w_i, \Lambda(4n)\big), 1 \leq i \leq N \,\big\} \end{matrix}\right) \\ &\leq\, P\left(\mathcal{D}(w_1, \cdots, w_N),\, \mathcal{E},\, T \leq 2(\ln n)^{1-\gamma},\, \big|\overline{\mathcal{T}}^0(T)\big| \,\geq\, \frac{N^* n}{16d^2}\right), \end{aligned} \tag{87.8}$$

where $\mathcal{D}(w_1, \cdots, w_N)$ is the event that $w_1, \dots, w_N$ are pairwise disconnected in $\Lambda(4n)$. As usual, in the sequel, we will write simply $\mathcal{D}$ instead of $\mathcal{D}(w_1, \cdots, w_N)$.

## 87.4 A leap in the relevant taboo sites

We need to do one more operation before entering the second stage of the proof. We examine the sequence $(\left|\overline{\mathcal{T}}^0(t)\right|, 0 \leq t \leq T)$ in order to detect an index $t$ such that there is a big gap between $\left|\overline{\mathcal{T}}^0(t-1)\right|$ and $\left|\overline{\mathcal{T}}^0(t)\right|$ and, at the same time, $\left|\overline{\mathcal{T}}^0(t)\right|$ is close to the final value $\left|\overline{\mathcal{T}}^0(T)\right|$. To do that, we will use the following little lemma.

**Lemma 87.2.** *Let $k \geq 1$ and let $u_1, \dots, u_k$ be $k$ positive integers such that $u_k > u_1$. Let $c$ be a real number such that*

$$1 \leq c^k < \frac{u_k}{u_1}.$$

*There exists an index $i^*$ in $\{2, \dots, k\}$ such that*

$$u_{i^*} \geq \frac{1}{c^k} u_k, \qquad u_{i^*-1} \leq \frac{1}{c} u_{i^*}. \tag{87.9}$$

*Proof.* Let us consider the set of indices $I$ defined by

$$I = \Big\{ i \in \{2, \dots, k\} : u_i \geq c\, u_{i-1} \Big\}.$$

Suppose that the set $I$ is empty. We would then have

$$\forall i \in \{2, \dots, k\} \qquad u_i < c\, u_{i-1}.$$

Multiplying together all these inequalities would lead to $u_k < c^{k-1} u_1$, which is absurd. Thus $I$ is not empty. Let now $i^*$ be the maximum element of $I$. By definition, we have

$$\forall i \in \{i^*+1, \dots, k\} \qquad u_i < c\, u_{i-1}.$$

In the case where $i^* < k$, multiplying together these inequalities, we get,

$$u_k < c^{k-i^*} u_{i^*} \leq c^k u_{i^*}.$$

In the case where $i^* = k$, we have obviously $u_k \leq c^k u_k$. In both cases, we see that the index $i^*$ satisfies the two required inequalities (87.9). □

Let us consider a configuration realizing the event appearing in the last probability of formula (87.8). In particular, we have

$$\left|\overline{\mathcal{T}}^0(T)\right| \geq \frac{N^* n}{16 d^2}. \tag{87.10}$$

Notice that the sequence $(\left|\overline{\mathcal{T}}^0(t)\right|, 0 \leq t \leq T)$ is non-decreasing and that we have $T \geq 1$. Let $T_0$ be the index of the first positive element in this sequence. Suppose that

$$1 \leq \left|\overline{\mathcal{T}}^0(T_0)\right| \leq \frac{N^* n}{16 d^2 \psi(n)}, \tag{87.11}$$

where $\psi(n) > 1$ is a function that we will choose later. In view of (87.10), this implies also that $T_0 < T$. We apply lemma 87.2 to the sequence

$$u_i \,=\, \big|\overline{\mathcal{T}}^0(T_0-1+i)\big|\,, \qquad 1\leq i\leq T-T_0+1\,,$$

and with a parameter $c$ such that $1 \leq c^{T-T_0+1} \leq \psi(n)$. To get rid of $T_0$, we will impose the stronger condition

$$1\,\leq\, c^{T+1}\,\leq\,\psi(n)\,. \tag{87.12}$$

The exact value of $c$ will be adjusted later, and it will depend on $n$. Lemma 87.2 provides us with an index $i^*$ such that the time $t^* = T_0 - 1 + i^*$ satisfies

$$\big|\overline{\mathcal{T}}^0(t^*)\big|\,\geq\,\frac{N^*n}{16d^2c^{T-T_0+1}}\,,\qquad \big|\overline{\mathcal{T}}^0(t^*-1)\big|\,\leq\,\frac{1}{c}\big|\overline{\mathcal{T}}^0(t^*)\big|\,. \tag{87.13}$$

The occurrence of the event $\mathcal{G}(n)$ defined in (85.6) and the lower bound (85.11) on $\beta(n)$ imply that the termination time $T$ is less than or equal to $\bar{t}(n)$, where

$$\bar{t}(n)\,=\,2(\ln n)^{1-\gamma}\,. \tag{87.14}$$

So we replace the first inequality in (87.13) by the weaker inequality

$$\big|\overline{\mathcal{T}}^0(t^*)\big|\,\geq\,\frac{N^*n}{16d^2c^{\bar{t}(n)}}\,.$$

If the condition (87.11) does not hold, then we have $\big|\overline{\mathcal{T}}^0(T_0)\big| > N^*n/(16d^2\psi(n))$. We distinguish also two cases: either $T_0 = 0$ or $\big[T_0 \geq 1$ and $\overline{\mathcal{T}}^0(T_0-1) = \varnothing\big]$. We decompose the probability in (87.8) according to these different cases and the values of $T_0$ and $t^*$:

$$P(\mathcal{H})\,\leq\,P\Big(\mathcal{D},\,\mathcal{E},\,\big|\overline{\mathcal{T}}^0(0)\big|\,>\,\frac{N^*n}{16d^2\psi(n)}\Big) \tag{87.15}$$

$$+\sum_{1\leq t\leq \bar{t}(n)}P\begin{pmatrix}\mathcal{D},\,\mathcal{E},\,\overline{\mathcal{T}}^0(t-1)\,=\,\varnothing,\\ \big|\overline{\mathcal{T}}^0(t)\big|\,>\,\dfrac{N^*n}{16d^2\psi(n)}\end{pmatrix} \tag{87.16}$$

$$+\sum_{1\leq t^*\leq \bar{t}(n)}P\begin{pmatrix}\mathcal{D},\,\mathcal{E},\,\big|\overline{\mathcal{T}}^0(t^*)\big|\,>\,\dfrac{N^*n}{16d^2c^{\bar{t}(n)}}\,,\\ \big|\overline{\mathcal{T}}^0(t^*-1)\big|\,\leq\,\dfrac{1}{c}\big|\overline{\mathcal{T}}^0(t^*)\big|\end{pmatrix}. \tag{87.17}$$

From now onwards, we fix $t$ such that $1 \leq t \leq \bar{t}(n)$, a function $\rho : \mathbb{N} \to [0,1]$, and a function $\phi : \mathbb{N} \to [1,+\infty[$ and we focus on the probability

$$P\begin{pmatrix}\mathcal{D},\,\mathcal{E},\,\big|\overline{\mathcal{T}}^0(t^*)\big|\,>\,\dfrac{N^*n}{\phi(n)}\,,\\ \big|\overline{\mathcal{T}}^0(t^*-1)\big|\,\leq\,\rho(n)\,\big|\overline{\mathcal{T}}^0(t^*)\big|\end{pmatrix}.$$

We consider two choices for $\phi(n)$ and $\rho(n)$:
- $\phi(n) = 16d^2\psi(n)$ and $\rho(n) = 0$, this corresponds to the probability in (87.16);
- $\phi(n) = 16d^2c^{\bar{t}(n)}$ and $\rho(n) = 1/c$ where $c$ is some value satisfying (87.12), this corresponds to the probability in (87.17).

# 88 Localization of the relevant taboo sites

Ultimately, we wish to design genuine explorations that manage to localize precisely the relevant taboo sites. Unfortunately, this does not seem to be feasible. In this section, we localize the relevant taboo sites with the help of pairs of taboo explorations. In the end, we will show that a positive fraction of the relevant taboo sites is included in the intersection of only two taboo explorations.

The constructions carried out here do not depend on the specific geometry of $\Lambda(4n)$, so we present them in our generic situation. Let $D$ be a finite domain in $\mathbb{Z}^d$ and let $G$ be a subset of $D$. Let $N \geq 2$ and let $w_1, \dots, w_N$ be $N$ distinct sites of $D$. To shorten the formulas, we define

$$W \,=\, \{\, w_1, \dots, w_N \,\}\,.$$

Let $\mathcal{D}$ be the event that $w_1, \dots, w_N$ are pairwise disconnected in $D$, i.e.,

$$\mathcal{D} \,=\, \bigcap_{1 \leq j < k \leq N} \{\, w_j \not\longleftrightarrow w_k \text{ in } D \,\}\,.$$

We consider a site percolation configuration $\omega$ belonging to $\mathcal{D}$ and the $N$ intertwined explorations with goal set $G$ starting from $w_1, \dots, w_N$:

$$\mathcal{A}_i(t),\, \mathcal{W}_i(t), \quad 1 \leq i \leq N\,, \quad 0 \leq t \leq T\,.$$

Our aim is to localize and control the relevant taboo sets $\mathcal{T}^*(t)$, $0 \leq t \leq T$.

## 88.1 The first relevant taboo set

We focus here on $\mathcal{T}^{*0}(0)$. The sites of $\mathcal{T}^{*0}(0)$ are exactly the sites of $D$ that are visited by two cluster explorations starting from two different sites of $W$. In fact, on the event $\mathcal{D}$, the $N$ sets $C(w_i, D)$, $1 \leq i \leq N$, are pairwise disjoint, and we have

$$\mathcal{T}^{*0}(0) \,=\, \bigcup_{1 \leq k < \ell \leq N} \partial^{\,out} C(w_k, D) \cap \partial^{\,out} C(w_\ell, D)\,.$$

We apply next lemma 84.3 to the family of sets

$$E_i \,=\, F_i \,=\, \partial^{\,out} C(w_i, D)\,, \quad 1 \leq i \leq N\,.$$

This is a special case where the two collections $(E_i, i \in I)$, $(F_i, i \in I)$ are the same. For a site $x$ to belong to $\partial^{\,out} C(w_i, D)$, it must have a neighbour in $C(w_i, D)$. On the event $\mathcal{D}$, the clusters $C(w_i, D)$ are pairwise disjoint, thus a site $x$ can belong to at most $2d$ sets among $\partial^{\,out} C(w_i, D)$, $1 \leq i \leq N$. Therefore the collection

$$\partial^{\,out} C(w_i, D) \cap \partial^{\,out} C(w_j, D)\,, \quad 1 \leq i, j \leq N\,, i \neq j\,,$$

has multiplicity less than $4d^2$. By lemma 84.3, there exist two disjoint subsets $K, L$ of $\{\, 1, \dots, N \,\}$ such that

$$\Big| \Big( \bigcup_{k \in K} \partial^{\,out} C(w_k, D) \Big) \cap \Big( \bigcup_{\ell \in L} \partial^{\,out} C(w_\ell, D) \Big) \Big| \,\geq\, \frac{1}{16d^2} \big| \mathcal{T}^{*0}(0) \big|\,. \tag{88.1}$$

## 88.2 The successive relevant taboo sets

In this subsection and the next three, we fix $t$ in $\{1,\dots,T\}$. Let us explain briefly the strategy behind the computations performed in the next subsections. The main idea is that the relevant taboo sites can be reached by two distinct taboo explorations starting from two different sites of $W$. On one hand, they are obviously reached by the first explorer that visits them. On the other hand, another explorer manages to visit a neighboring site during the same iteration, or the iteration before. By tuning adequately the taboo sets associated to these two explorers, we highlight two distinct taboo explorations that reach simultaneously a given relevant taboo site. The taboo explorations will be used only to localize the closed relevant taboo sites. In subsection 88.3, we will prove that an open relevant taboo site possesses a neighbour which is a closed relevant taboo site. Thus, the control of the sets of the closed relevant taboo sites yields automatically a control of the sets of the open relevant taboo sites. There is an extra difficulty linked with the visiting time of the relevant taboo sites and their states. For the closed relevant taboo sites, we develop first a control on the sites visited at time $t-1$. This will yield a control on the open relevant taboo sites involved in the localization of the closed relevant taboo sites visited at time $t$. In turn, this last result will give a control on all the open relevant taboo sites at time $t$.

In the definition (80.39), we introduced the partition of $\mathcal{T}^*(t)$ into the two subsets $\mathcal{T}^-(t)$ and $\mathcal{T}^+(t)$. We will focus on the closed sites contained in these sets, namely $\mathcal{T}^{-0}(t)$ and $\mathcal{T}^{+0}(t)$. We classify the sites of these two sets according to the indices of the explorers for which they are taboo. For $i$ in $\{1,\dots,N\}$, we define

$$\begin{aligned}
\mathcal{T}_i^{-0}(t) \;&=\; \big\{\, x\in\mathcal{T}_i^*(t) : T(x)=t-1 \text{ or } T(x)=t-2,\, \omega(x)=0 \,\big\}\,,\\
\mathcal{T}_i^{+0}(t) \;&=\; \big\{\, x\in\mathcal{T}_i^*(t) : T(x)=t,\, \omega(x)=0 \,\big\}\,.
\end{aligned}$$

Notice that these sets are not necessarily disjoint. Indeed, a site can be a relevant taboo site for several explorers. Nevertheless, we have

$$\begin{aligned}
\mathcal{T}^{-0}(t) \;&=\; \bigcup_{1\leq i\leq N} \mathcal{T}_i^{-0}(t) \;=\; \mathcal{T}^-(t)\cap\mathcal{T}^{*0}(t)\,,\\
\mathcal{T}^{+0}(t) \;&=\; \bigcup_{1\leq i\leq N} \mathcal{T}_i^{+0}(t) \;=\; \mathcal{T}^+(t)\cap\mathcal{T}^{*0}(t)\,.
\end{aligned}$$

Furthermore, in order to control $\big|\mathcal{T}^*(t)\big|$, we must not forget the open relevant taboo sites, that is the set $\mathcal{T}^{*1}(t)$. In the next three subsections 88.3, 88.4, 88.5, we show how to control successively $\big|\mathcal{T}^{*1}(t)\big|$, $\big|\mathcal{T}^{-0}(t)\big|$, $\big|\mathcal{T}^{+0}(t)\big|$. The outcome of all the computations is presented in subsection 88.7. For $t$ in $\{0,\dots,T\}$, we shall denote by $\overline{\mathcal{T}}^1(t)$ the set of all the open relevant taboo sites until time $t$, i.e.,

$$\forall t\in\{0,\dots,T\}\qquad \overline{\mathcal{T}}^1(t) \;=\; \bigcup_{0\leq s\leq t}\mathcal{T}^{*1}(s)\,.$$

Table 18 lists the various taboo sets we have introduced and their meaning.

| name | interpretation |
|---|---|
| $\mathcal{T}^*(t)$ | relevant taboo sites (rts in short) at time $t$ |
| $\mathcal{T}^{*0}(t)$ | closed rts at time $t$ |
| $\mathcal{T}^{*1}(t)$ | open rts at time $t$ |
| $\mathcal{T}^-(t)$ | rts at time $t$, visited at time $t-1$ or $t-2$ |
| $\mathcal{T}^+(t)$ | rts at time $t$ which are visited at time $t$ |
| $\mathcal{T}_i^*(t)$ | rts at time $t$ for the explorer $i$ |
| $\mathcal{T}_i^-(t)$ | rts at time $t$ which were visited at time $t-1$ or $t-2$ and are taboo for the explorer $i$ |
| $\mathcal{T}_i^+(t)$ | rts at time $t$ which are visited at time $t$ and are taboo for the explorer $i$ |
| $\mathcal{T}_i^{-0}(t)$ | closed rts at time $t$, visited at time $t-1$ or $t-2$, and taboo for the explorer $i$ |
| $\mathcal{T}_i^{+0}(t)$ | closed rts at time $t$, which are visited at time $t$ and are taboo for the explorer $i$ |
| $\overline{\mathcal{T}}_i(t)$ | the union of the sets $\mathcal{T}_i(s)$, $0 \le s \le t$ |
| $\overline{\mathcal{T}}(t)$ | the union of the sets $\mathcal{T}^*(s)$, $0 \le s \le t$ |
| $\overline{\mathcal{T}}^0(t)$ | the union of the sets $\mathcal{T}^0(s)$, $0 \le s \le t$ |
| $\overline{\mathcal{T}}^1(t)$ | the union of the sets $\mathcal{T}^1(s)$, $0 \le s \le t$ |

Table 18: The various sets of relevant taboo sites (rts in short)

Here is a brief summary of the system of notation for the relevant taboo sites:

- The subscript $i$ refers to the explorer $i$;
- The argument $t$ in the parentheses refers to the iteration time;
- The superscript $*$ means that the sites can be either closed or open;
- The superscript 0 means that the sites are closed;
- The superscript 1 means that the sites are open;
- The superscript $-$ means that the visiting time is $t-1$ or $t-2$;
- The superscript $+$ means that the visiting time is $t$;
- Overlining means that we take the union of all previous sets of the same type.

## 88.3 Localization of $\mathcal{T}^{*1}(t)$

We focus here on the set $\mathcal{T}^{*1}(t)$ of the open relevant taboo sites at time $t$.

**Lemma 88.1.** *For $t$ in $\{1,\dots,T\}$, any site in $\mathcal{T}^{-1}(t)$ has a neighbour belonging to $\mathcal{T}^{+0}(t-1)$, any site in $\mathcal{T}^{+1}(t)$ has a neighbour belonging to $\mathcal{T}^{-0}(t)$.*

*Proof.* Let $t$ in $\{1,\dots,T\}$ and let $x$ be a site belonging to $\mathcal{T}^{*1}(t)$. Let $i$ be the visitor index of $x$ and let $j$ be such that $x$ is in $\mathcal{T}_j^{*1}(t)$. By definition, the site $x$ would be visited by the explorer $j$ at time $t$ if it had not been already visited by the explorer $i$. To put it another way, during the iteration $t$, the explorer $j$ tries to visit the site $x$, but he realizes that $x$ has already been visited. This attempt of visit originates from a neighbour $y$ of $x$. Necessarily $y$ must be closed, otherwise it would belong to the aggregate of the explorer $i$. Since the attempt to visit $x$ starting from $y$ occurs during the iteration $t$ and $y$ is closed, then the visiting time of $y$ must be equal to $t-1$. To sum up, we have

$$x \text{ is open}\,,\quad T(x)=t-1 \text{ or } T(x)=t\,,\quad y \text{ is closed}\,,\quad T(y)=t-1\,.$$

At the precise moment when the explorer $i$ visits $x$, he tries to visit the neighbours of $x$. So he attempts to visit $y$, but he fails because $y$ was already visited by $j$ (otherwise $y$ would not have $j$ for visitor index!). This means that $y$ is a closed relevant taboo site for the explorer $i$ at time $T(x)$. We conclude that $y$ belongs to $\mathcal{T}^{+0}(t-1)$ if $T(x)=t-1$ and to $\mathcal{T}^{-0}(t)$ if $T(x)=t$. ☐

An example of two intertwined explorations and their relevant taboo sites is presented in figure 157, and again in figure 158 with more details on the first taboo sets. An immediate consequence of lemma 88.1 is the following inequalities:

$$\forall t\in\{1,\dots,T\}\qquad \big|\mathcal{T}^{-1}(t)\big|\,\le\,2d\,\big|\mathcal{T}^{+0}(t-1)\big|\,,\quad \big|\mathcal{T}^{+1}(t)\big|\,\le\,2d\,\big|\mathcal{T}^{-0}(t)\big|\,. \tag{88.2}$$

Therefore the control of $\big|\mathcal{T}^{*1}(t)\big|$ is guaranteed once we have an adequate control on $\big|\mathcal{T}^{*0}(t-1)\big|$ and $\big|\mathcal{T}^{-0}(t)\big|$. The next goal is to control these last sets. As we will see, the control of $\big|\mathcal{T}^{-0}(t)\big|$ will rest on the control of $\big|\mathcal{T}^{-1}(t)\big|$!

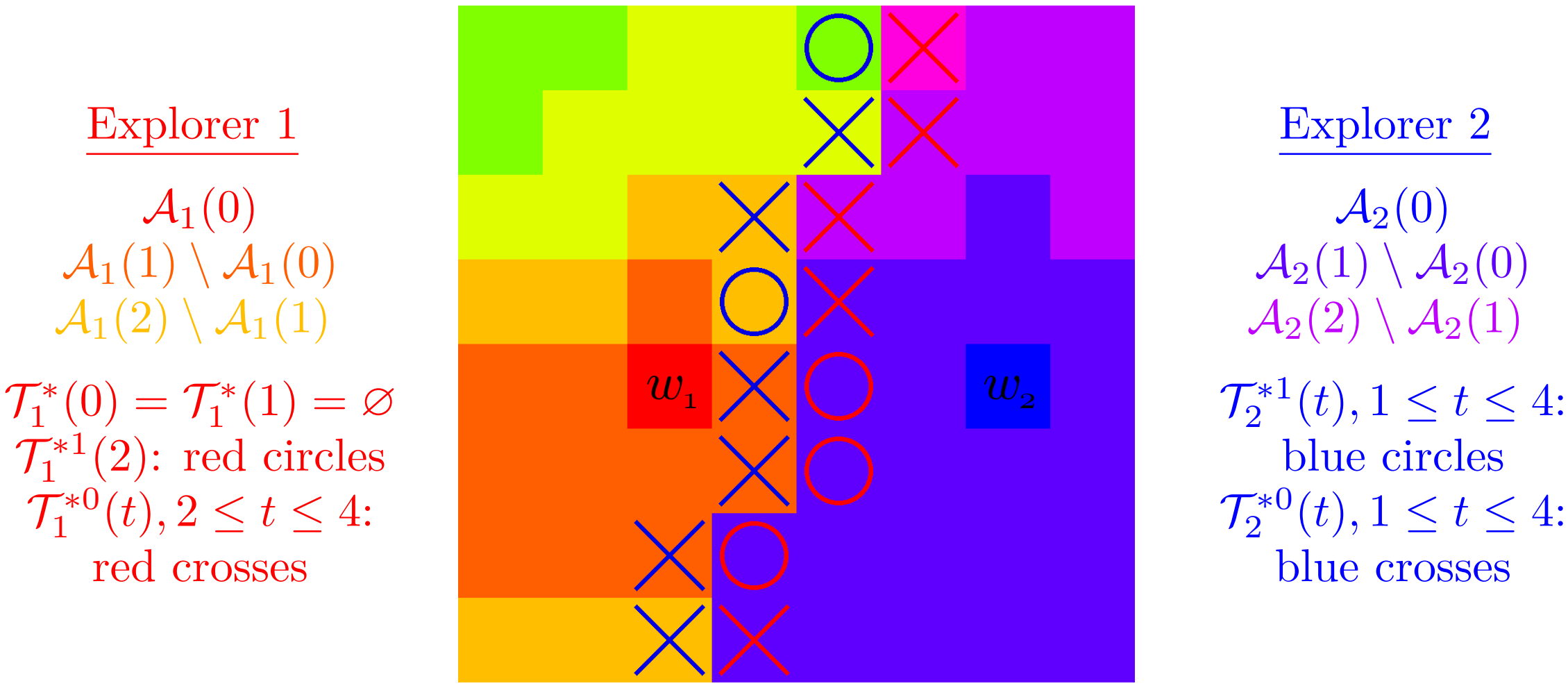

Figure 157: A site in $\mathcal{T}^{*1}(t)$ always has a neighbour in $\mathcal{T}^{*0}(t-1)\cup\mathcal{T}^{*0}(t)$.

## 88.4 Localization of $\mathcal{T}^{-0}(t)$

We focus here on the set $\mathcal{T}^{-0}(t)$. We suppose that we have already gained some control on $\overline{\mathcal{T}}(t-1)$, and we try to control $\mathcal{T}^{-0}(t)$ with the help of inequalities that involve $\overline{\mathcal{T}}(t-1)$ and adequate taboo explorations. By the definition (80.39), a taboo site in $\mathcal{T}^{-}(t)$ is visited at time $t-1$ or $t-2$, therefore

$$\mathcal{T}^{-0}(t) \,=\, \bigcup_{1\leq i,j\leq N, i\neq j} \mathcal{T}_i^{-0}(t)\cap\mathcal{A}_j(t-1)\,.$$

Let us consider a site $x$ belonging to $\mathcal{T}_i^{-0}(t)$. By definition, the explorer $i$ attempted to visit $x$ at time $t$. This attempt must originate from a neighbor of $x$, thus one of the neighbors of $x$ must belong to $\mathcal{A}_i(t)$. In addition, the aggregates $(\mathcal{A}_j(t-1), 1\leq j\leq N)$ are pairwise disjoint, as well as $(\mathcal{A}_i(t), 1\leq i\leq N)$, thus a site $x$ can belong to at most $2d$ sets of the family

$$\mathcal{T}_i^{-0}(t)\cap\mathcal{A}_j(t-1)\,,\quad 1\leq i,j\leq N,\, i\neq j\,.$$

By lemma 84.3 applied to this family with $m = 2d$, there exist two disjoint subsets $K^-, L^-$ of $\{\,1,\dots,N\,\}$ such that

$$\frac{1}{8d}\,|\mathcal{T}^{-0}(t)| \,\leq\, \Big|\Big(\bigcup_{k\in K^-}\mathcal{T}_k^{-0}(t)\Big)\cap\Big(\bigcup_{\ell\in L^-}\mathcal{A}_\ell(t-1)\Big)\Big|\,. \tag{88.3}$$

It can happen that a site belongs simultaneously to the sets $\mathcal{T}^{-0}(t)$ and $\overline{\mathcal{T}}(t-1)$. Imagine for instance that the site $x$ is closed and visited by explorer $i$ at time $t-1$, that one of its neighbors is open and visited by explorer $j$ at time $t-1$, with $j>i$, and finally another of its neighbors is open and visited by explorer $k$ at time $t$, with $k<i$. Then $x$ is in both taboo sets $\mathcal{T}_j^{+0}(t-1)$ and $\mathcal{T}_k^{-0}(t)$. As we supposed that the relevant taboo sites at time $t-1$ have already been controlled, we need only to control the relevant taboo sites at time $t$ which are not relevant at time $t-1$. We have the following inclusions:

$$\forall k\in K^-\qquad \mathcal{T}_k^{-0}(t)\,\subset\,\big(\mathcal{T}_k^{-0}(t)\setminus\overline{\mathcal{T}}_k(t-1)\big)\cup\overline{\mathcal{T}}_k(t-1)\,,$$

which, together with (88.3), yield

$$\frac{1}{8d}\,|\mathcal{T}^{-0}(t)| \,\leq\, \Big|\Big(\bigcup_{k\in K^-}\mathcal{T}_k^{-0}(t)\setminus\overline{\mathcal{T}}_k(t-1)\Big)\cap\Big(\bigcup_{\ell\in L^-}\mathcal{A}_\ell(t-1)\Big)\Big|+\big|\overline{\mathcal{T}}(t-1)\big|\,. \tag{88.4}$$

Let us define

$$\forall k\in K^-\qquad \widehat{\mathcal{T}}_k^-(t)\,=\,\overline{\mathcal{T}}_k(t-1)\cup\mathcal{T}_k^{-1}(t)\,. \tag{88.5}$$

It follows from proposition 81.2 that, for any $k$ in $K^-$,

$$\begin{aligned}\mathcal{T}_k^{-0}(t)\setminus\overline{\mathcal{T}}(t-1)\,&\subset\,\mathcal{T}_k^{*0}(t)\setminus\overline{\mathcal{T}}_k(t-1)\\ &\subset\mathcal{B}\Big(w_k, D\setminus\big(\overline{\mathcal{T}}_k(t-1)\cup\mathcal{T}_k^{*1}(t)\big),t\Big)\,\subset\,\mathcal{B}\big(w_k, D\setminus\widehat{\mathcal{T}}_k^-(t),t\big)\,.\end{aligned} \tag{88.6}$$

Taking the union over $k$ in $K^-$, we deduce from the inclusion (88.6) that

$$\bigcup_{k\in K^-}\mathcal{T}_k^{-0}(t)\setminus\overline{\mathcal{T}}_k(t-1)\,\subset\,\bigcup_{k\in K^-}\mathcal{B}\big(w_k, D\setminus\widehat{\mathcal{T}}_k^-(t),t\big)\,. \tag{88.7}$$

The inclusion (88.7) readily implies that

$$\Big(\bigcup_{k\in K^-}\mathcal{T}_k^{-0}(t)\setminus\overline{\mathcal{T}}_k(t-1)\Big)\cap\Big(\bigcup_{\ell\in L^-}\mathcal{A}_\ell(t-1)\Big)\subset\\ \Big(\bigcup_{k\in K^-}\mathcal{B}\big(w_k,D\setminus\widehat{\mathcal{T}}_k^-(t),t\big)\Big)\cap\Big(\bigcup_{\ell\in L^-}\mathcal{A}_\ell(t-1)\Big).\quad(88.8)$$

We deduce from (88.4) and (88.8) that

$$\frac{1}{8d}\left|\mathcal{T}^{-0}(t)\right|\leq\left|\Big(\bigcup_{k\in K^-}\mathcal{B}\big(w_k,D\setminus\widehat{\mathcal{T}}_k^-(t),t\big)\Big)\cap\Big(\bigcup_{\ell\in L^-}\mathcal{A}_\ell(t-1)\Big)\right|+\left|\overline{\mathcal{T}}(t-1)\right|.\quad(88.9)$$

By proposition 81.1, we have

$$\forall\ell\in L^-\qquad\mathcal{A}_\ell(t-1)\,=\,\mathcal{B}\big(w_\ell,D\setminus\overline{\mathcal{T}}_\ell(t-1),t-1\big).\quad(88.10)$$

We conclude from (88.9) and (88.10) that

$$\frac{1}{8d}\left|\mathcal{T}^{-0}(t)\right|\leq\left|\Big(\bigcup_{k\in K^-}\mathcal{B}\big(w_k,D\setminus\widehat{\mathcal{T}}_k^-(t),t\big)\Big)\cap\Big(\bigcup_{\ell\in L^-}\mathcal{B}\big(w_\ell,D\setminus\overline{\mathcal{T}}_\ell(t-1),t-1\big)\Big)\right|\\ +\left|\overline{\mathcal{T}}(t-1)\right|.\quad(88.11)$$

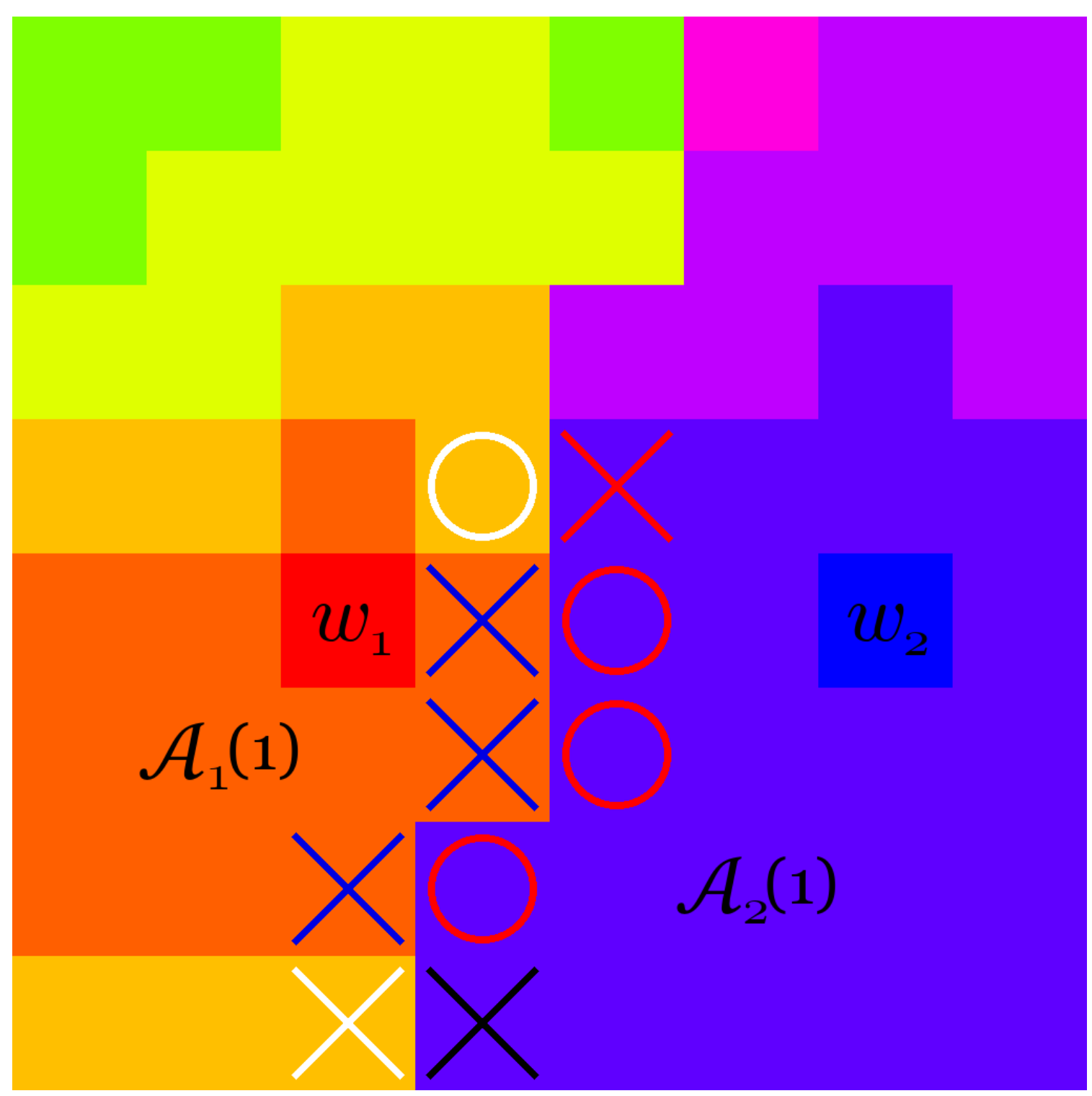

Figure 158: The explorations started from $w_1$, $w_2$ in the box $\Lambda(8)$ with $p=0.3$: $\mathcal{A}_1(0)$, $\mathcal{A}_1(1)\setminus\mathcal{A}_1(0)$, $\mathcal{A}_1(2)\setminus\mathcal{A}_1(1)$, $\mathcal{A}_2(0)$, $\mathcal{A}_2(1)\setminus\mathcal{A}_2(0)$, $\mathcal{A}_2(2)\setminus\mathcal{A}_2(1)$, crosses, $w_1$, $w_2$ are closed, circles are open. red circles are in $\mathcal{T}_1^{-1}(2)$, blue crosses are in $\mathcal{T}_2^{-0}(2)$, red cross is in $\mathcal{T}_1^{-0}(3)$, white circle is in $\mathcal{T}_2^{+1}(3)$, white cross is in $\mathcal{T}_2^{+0}(3)$.

## 88.5 Localization of $\mathcal{T}^{+0}(t)$

We focus here on the set $\mathcal{T}^{+0}(t)$. The technique is globally similar to the one employed for the set $\mathcal{T}^{-0}(t)$, but with some crucial differences. For instance, the taboo set used to localize the set $\mathcal{T}^{+0}(t)$ is not the same as for the set $\mathcal{T}^{-0}(t)$. Thus we cannot simply rely on the inequalities of the previous subsection, and we have to derive from scratch the required inequalities. It follows from the definition (80.39) of $\mathcal{T}^{+}(t)$ that

$$\mathcal{T}^{+0}(t) \,=\, \bigcup_{1\leq i,j\leq N, i\neq j} \mathcal{T}_i^{+0}(t)\cap \mathcal{A}_j(t)\,.$$

Using the very same argument as in the beginning of the previous subsection, we obtain that there exist two disjoint subsets $K^+, L^+$ of $\{\,1,\dots,N\,\}$ such that

$$\frac{1}{8d}\,|\mathcal{T}^{+0}(t)| \,\leq\, \Big|\Big(\bigcup_{k\in K^+}\mathcal{T}_k^{+0}(t)\Big)\cap\Big(\bigcup_{\ell\in L^+}\mathcal{A}_\ell(t)\Big)\Big|\,. \tag{88.12}$$

Contrary to the case of the set $\mathcal{T}^{-0}(t)$ studied in subsection 88.4, the set $\mathcal{T}^{+0}(t)$ is disjoint from $\overline{\mathcal{T}}(t-1)$. Indeed, the visiting time of the sites in $\mathcal{T}^{+0}(t)$ is equal to $t$, hence they cannot belong to $\overline{\mathcal{T}}(t-1)$ as well. We have thus

$$\forall k\in K^+\qquad \mathcal{T}_k^{+0}(t)\,\subset\,\mathcal{T}_k^{*0}(t)\setminus\overline{\mathcal{T}}_k(t-1)\,. \tag{88.13}$$

For $k$ in $K^+$, we define

$$\widehat{\mathcal{T}}_k^{*}(t)\,=\,\overline{\mathcal{T}}_k(t-1)\cup\mathcal{T}_k^{-0}(t)\cup\mathcal{T}_k^{*1}(t)\,. \tag{88.14}$$

It follows from proposition 81.2 that

$$\forall k\in K^+\qquad \mathcal{T}_k^{+0}(t)\setminus\overline{\mathcal{T}}_k(t-1)\,\subset\,\mathcal{B}\big(w_k,D\setminus\widehat{\mathcal{T}}_k^{*}(t),t\big)\,. \tag{88.15}$$

We deduce from the two inclusions (88.13) and (88.15) that

$$\bigcup_{k\in K^+}\mathcal{T}_k^{+0}(t)\,\subset\,\bigcup_{k\in K^+}\mathcal{B}\big(w_k,D\setminus\widehat{\mathcal{T}}_k^{*}(t),t\big)\,. \tag{88.16}$$

The inclusion (88.16) readily implies that

$$\Big(\bigcup_{k\in K^+}\mathcal{T}_k^{+0}(t)\Big)\cap\Big(\bigcup_{\ell\in L^+}\mathcal{A}_\ell(t)\Big)\,\subset\,\Big(\bigcup_{k\in K^+}\mathcal{B}\big(w_k,D\setminus\widehat{\mathcal{T}}_k^{*}(t),t\big)\Big)\cap\Big(\bigcup_{\ell\in L^+}\mathcal{A}_\ell(t)\Big)\,. \tag{88.17}$$

We deduce from (88.12) and (88.17) that

$$\frac{1}{8d}\,|\mathcal{T}^{+0}(t)| \,\leq\, \Big|\Big(\bigcup_{k\in K^+}\mathcal{B}\big(w_k,D\setminus\widehat{\mathcal{T}}_k^{*}(t),t\big)\Big)\cap\Big(\bigcup_{\ell\in L^+}\mathcal{A}_\ell(t)\Big)\Big|\,. \tag{88.18}$$

We still need to rework the union $\bigcup_{\ell\in L^+}\mathcal{A}_\ell(t)$. By proposition 81.1, we have

$$\forall \ell\in L^+\qquad \mathcal{A}_\ell(t)\,=\,\mathcal{B}\big(w_\ell,D\setminus\overline{\mathcal{T}}_\ell(t),t\big)\,. \tag{88.19}$$

Since our goal is to control the set $\mathcal{T}_\ell^{*0}(t)$, we try to get rid of it in (88.19). In the next lemma, we examine what happens if we replace $\overline{\mathcal{T}}_\ell(t)$ by $\widehat{\mathcal{T}}_\ell^{*}(t)$ in (88.19).

**Lemma 88.2.** *For any $t$ in $\{1,\dots,T\}$ and for any $\ell$ in $L^+$, we have*

$$\mathcal{B}\big(w_\ell, D\setminus\widehat{\mathcal{T}}_\ell^*(t), t\big) \,=\, \mathcal{A}_\ell(t)\,\cup\,\mathcal{T}_\ell^{+0}(t)\,. \tag{88.20}$$

*Proof.* Let us fix $\ell$ in $L^+$ and $t$ in $\{1,\dots,T\}$. Since $\widehat{\mathcal{T}}_\ell^*(t)$ is included in $\overline{\mathcal{T}}_\ell(t)$, then

$$\mathcal{A}_\ell(t) \,=\, \mathcal{B}\big(w_\ell, D\setminus\overline{\mathcal{T}}_\ell(t), t\big) \,\subset\, \mathcal{B}\big(w_\ell, D\setminus\widehat{\mathcal{T}}_\ell^*(t), t\big)\,. \tag{88.21}$$

In addition, it follows from proposition 81.1 that

$$\forall x\in\mathcal{T}_\ell^{+0}(t)\quad x\in\mathcal{B}\Big(w_\ell, D\setminus\big(\overline{\mathcal{T}}_\ell(t)\setminus\{x\}\big), t\Big) \,\subset\, \mathcal{B}\big(w_\ell, D\setminus\widehat{\mathcal{T}}_\ell^*(t), t\big)\,. \tag{88.22}$$

We conclude from (88.21) and (88.22) that

$$\mathcal{A}_\ell(t)\,\cup\,\mathcal{T}_\ell^{+0}(t) \,\subset\, \mathcal{B}\big(w_\ell, D\setminus\widehat{\mathcal{T}}_\ell^*(t), t\big)\,. \tag{88.23}$$

We prove next the converse inclusion. The set $\overline{\mathcal{T}}_\ell(t-1)$ is included in $\widehat{\mathcal{T}}_\ell^*(t)$, so that

$$\mathcal{B}\big(w_\ell, D\setminus\widehat{\mathcal{T}}_\ell^*(t), t-1\big) \,=\, \mathcal{B}\big(w_\ell, D\setminus\overline{\mathcal{T}}_\ell(t), t-1\big) \,=\, \mathcal{A}_\ell(t-1)\,,$$

whence

$$\begin{aligned}\mathcal{B}\big(w_\ell, D\setminus\widehat{\mathcal{T}}_\ell^*(t), t\big) \,&=\, \mathcal{B}\Big(\mathcal{B}\big(w_\ell, D\setminus\widehat{\mathcal{T}}_\ell^*(t), t-1\big), D\setminus\widehat{\mathcal{T}}_\ell^*(t), 1\Big)\\ &=\, \mathcal{B}\big(\mathcal{A}_\ell(t-1), D\setminus\widehat{\mathcal{T}}_\ell^*(t), 1\big)\,,\end{aligned}$$

$$\mathcal{B}\big(w_\ell, D\setminus\widehat{\mathcal{T}}_\ell^*(t), t\big)\,\setminus\,\mathcal{A}_\ell(t-1) \,=\, \mathcal{B}\big(\mathcal{A}_\ell(t-1), D\setminus\widehat{\mathcal{T}}_\ell^*(t), 1\big)\,\setminus\,\mathcal{A}_\ell(t-1)\,.$$

Let $x$ be a site in $\mathcal{B}\big(w_\ell, D\setminus\widehat{\mathcal{T}}_\ell^*(t), t\big)\setminus\mathcal{A}_\ell(t-1)$. From the above identity, there exists a path $y_0, y_1, \cdots, y_r$ in $D\setminus\widehat{\mathcal{T}}_\ell^*(t)$ joining a site $y_0$ of $\mathcal{A}_\ell(t-1)$ to $y_r = x$ such that exactly one site among $y_0,\cdots,y_{r-1}$ is closed, say $y_s$, where $0\le s\le r-1$. Suppose that $s\ge 1$. Since $y_0,\cdots,y_{s-1}$ are open and $\mathcal{A}_\ell(t-1)$ is an aggregate, then $y_{s-1}$ is also in $\mathcal{A}_\ell(t-1)$. We know also that $y_s$ is not in $\widehat{\mathcal{T}}_\ell^*(t)$, which contains $\overline{\mathcal{T}}_\ell(t-1)$, thus $y_s$ is not a relevant taboo site for the explorer $\ell$ at time $t-1$. Therefore $y_s$ belongs also to $\mathcal{A}_\ell(t-1)$. This is also true if $s=0$. The sites $y_{s+1},\cdots,y_{r-1}$ are open and they are not in $\widehat{\mathcal{T}}_\ell^*(t)$, which contains $\mathcal{T}_\ell^{*1}(t)$, thus they are not relevant taboo sites for the explorer $\ell$ at time $t$. Therefore the whole subpath $y_{s+1},\cdots,y_{r-1}$ is included in $\mathcal{A}_\ell(t)$. We consider two cases:

- The site $x=y_r$ is open. Since $x$ is not in $\mathcal{T}_\ell^{*1}(t)$, then $x$ is in $\mathcal{A}_\ell(t)$.
- The site $x=y_r$ is closed. If $x$ is not in $\mathcal{A}_\ell(t)$, then $x$ must be a relevant taboo site for the explorer $\ell$ at time $t$. Thus $x$ belongs to $\overline{\mathcal{T}}_\ell(t)\setminus\widehat{\mathcal{T}}_\ell^*(t) \,=\, \mathcal{T}_\ell^{+0}(t)$.

We conclude that $x$ is in $\mathcal{A}_\ell(t)\cup\mathcal{T}_\ell^{+0}(t)$, thus proving that

$$\mathcal{B}\big(w_\ell, D\setminus\widehat{\mathcal{T}}_\ell^*(t), t\big)\setminus\mathcal{A}_\ell(t-1) \,\subset\, \mathcal{A}_\ell(t)\,\cup\,\mathcal{T}_\ell^{+0}(t)\,.$$

This inclusion, together with (88.23), completes the proof of (88.20). $\square$

Lemma 88.2 and the inclusion (88.18) together imply that

$$\frac{1}{8d}\,\big|\mathcal{T}^{+0}(t)\big| \,\le\, \Big|\Big(\bigcup_{k\in K^+}\mathcal{B}\big(w_k, D\setminus\widehat{\mathcal{T}}_k^*(t), t\big)\Big)\cap\Big(\bigcup_{\ell\in L^+}\mathcal{B}\big(w_\ell, D\setminus\widehat{\mathcal{T}}_\ell^*(t), t\big)\Big)\Big|\,. \tag{88.24}$$

## 88.6 The key point

We prove here the following crucial fact: all the sites belonging to the sets appearing as intersections in the right-hand side of the inequalities (88.11) and (88.24) are closed. Why is that so important?

Our strategy is to show that the presence of numerous crossing clusters implies the existence of an excess of closed sites. We have already proved that a large number of crossing clusters implies that one of the intersection sets must be large, which in turn implies that one of the sets of the closed relevant taboo sites $\mathcal{T}^{*0}(t)$ must be large. The inequalities (88.11) and (88.24) imply then that the intersection of the outcome of two well-chosen taboo explorations is large. We prove precisely here that these taboo explorations can meet only on closed sites. Thus this intersection will exhibit an excess of closed sites!

**Lemma 88.3.** *For any pair $k,\ell$ of distinct indices in $\{1,\dots,N\}$, and any $t$ in $\{1,\dots,T\}$, the sites belonging to*

$$\mathcal{B}\big(w_k, D\setminus \widehat{\mathcal{T}}_k^-(t),t\big)\cap \mathcal{B}\big(w_\ell, D\setminus \overline{\mathcal{T}}_\ell(t-1),t-1\big)\,, \tag{88.25}$$

$$\mathcal{B}\big(w_k, D\setminus \widehat{\mathcal{T}}_k^*(t),t\big)\cap \mathcal{B}\big(w_\ell, D\setminus \widehat{\mathcal{T}}_\ell^*(t),t\big) \tag{88.26}$$

*are closed.*

*Proof.* The proof is very similar for the two sets, however there are some details which differ. We start by presenting some elements that are the same for both sets. Afterwards, we deal separately with each set. Let $k,\ell,t$ be as in the statement of the lemma. By proposition 81.1, we have

$$\mathcal{A}_k(t-1)\;=\;\mathcal{B}\big(w_k, D\setminus \overline{\mathcal{T}}_k(t-1),t-1\big)\,. \tag{88.27}$$

Using this equality, we see that $\mathcal{A}_k(t-1)$ is $T_{D\setminus\overline{\mathcal{T}}_k(t-1)}$-closed in the sense defined in (82.1). The equality (88.27) implies also that

$$\mathcal{B}\big(w_k, D\setminus \widehat{\mathcal{T}}_k^*(t),t-1\big)\;\subset\;\mathcal{B}\big(w_k, D\setminus \widehat{\mathcal{T}}_k^-(t),t-1\big)\;\subset\;\mathcal{A}_k(t-1)\,, \tag{88.28}$$

$$\mathcal{B}\big(w_k, D\setminus \overline{\mathcal{T}}_k(t-1),t-1\big)\cap\mathcal{B}\big(w_\ell, D\setminus \overline{\mathcal{T}}_\ell(t-1),t-1\big)\;=\;\varnothing\,. \tag{88.29}$$

Using the identity (12.9) on the balls and the inclusion (88.28), we have

$$\begin{aligned}\mathcal{B}\big(w_k, D\setminus \widehat{\mathcal{T}}_k^-(t),t\big)\;&=\;\mathcal{B}\Big(\mathcal{B}\big(w_k, D\setminus \widehat{\mathcal{T}}_k^-(t),t-1\big), D\setminus \widehat{\mathcal{T}}_k^-(t),1\Big)\\ &\qquad\subset\;\mathcal{B}\big(\mathcal{A}_k(t-1), D\setminus \widehat{\mathcal{T}}_k^-(t),1\big)\,,\end{aligned}$$

and similarly

$$\mathcal{B}\big(w_k, D\setminus \widehat{\mathcal{T}}_k^*(t),t\big)\;\subset\;\mathcal{B}\big(\mathcal{A}_k(t-1), D\setminus \widehat{\mathcal{T}}_k^*(t),1\big)\,.$$

We treat next separately the two sets (88.25) and (88.26). There are a lot of similarities in the two arguments, however there are a few unavoidable differences, so it seems best to keep them separated.

We start with the set

$$\mathcal{B}\big(w_k, D \setminus \widehat{\mathcal{T}}_k^-(t), t\big) \cap \mathcal{B}\big(w_\ell, D \setminus \overline{\mathcal{T}}_\ell(t-1), t-1\big)\,. \tag{88.30}$$

Let $x$ be a site belonging to the set (88.30). It follows from (88.29) that this set is disjoint from $\mathcal{B}\big(w_k, D \setminus \overline{\mathcal{T}}_k(t-1), t-1\big)$, thus $x$ is in

$$\begin{aligned}&\mathcal{B}\big(w_k, D \setminus \widehat{\mathcal{T}}_k^-(t), t\big) \setminus \mathcal{B}\big(w_k, D \setminus \overline{\mathcal{T}}_k(t-1), t-1\big) \,\subset\\ &\mathcal{B}\Big(\mathcal{B}\big(w_k, D \setminus \widehat{\mathcal{T}}_k^-(t), t-1\big), D \setminus \widehat{\mathcal{T}}_k^-(t), 1\Big) \setminus \mathcal{B}\big(w_k, D \setminus \widehat{\mathcal{T}}_k^-(t), t-1\big)\,.\end{aligned} \tag{88.31}$$

By the very definition (12.7) of the ball $\mathcal{B}(\cdot,\cdot,\cdot)$, there exists a path $y_0, y_1, \cdots, y_r$ in $D \setminus \widehat{\mathcal{T}}_k^-(t)$ joining a site $y_0$ of $\mathcal{B}\big(w_k, D \setminus \widehat{\mathcal{T}}_k^-(t), t-1\big)$ to $y_r = x$ such that exactly one site among $y_0, \cdots, y_{r-1}$ is closed, say $y_s$, where $0 \leq s \leq r-1$. Since $x$ belongs to the set (88.31), then $x$ and $y_0$ are distinct and also $r \geq 1$. We claim that $y_s$ belongs to $\mathcal{A}_k(t-1)$. The site $y_0$ is in $\mathcal{A}_k(t-1)$ thanks to the inclusion (88.28). Suppose that $s \geq 1$. Since $\mathcal{A}_k(t-1)$ is an aggregate and the sites $y_0, \cdots, y_{s-1}$ are open, then they all belong to $\mathcal{A}_k(t-1)$. Furthermore, the site $y_s$ is in $D \setminus \widehat{\mathcal{T}}_k^-(t) \subset D \setminus \overline{\mathcal{T}}_k(t-1)$ and within null travel distance from $y_{s-1}$. Using the remark right after (88.27), we conclude that $y_s$ is indeed in $\mathcal{A}_k(t-1)$. Recall that our goal is to prove that the site $x$ is closed. Suppose that $x$ is open. Since $y_s$ is the only closed site among $y_0, \cdots, y_{r-1}$ and we assumed that $y_r = x$ is open, then all the sites $y_{s+1}, \cdots, y_r$ are open. Yet the site $x$ belongs also to $\mathcal{B}\big(w_\ell, D \setminus \overline{\mathcal{T}}_\ell(t-1), t-1\big)$, which is equal to the aggregate $\mathcal{A}_\ell(t-1)$, therefore all the sites $y_{s+1}, \cdots, y_r$ must also belong to $\mathcal{A}_\ell(t-1)$. To sum up, we have

$$y_s \in \mathcal{A}_k(t-1)\,, \quad \omega(y_s) = 0\,, \quad y_{s+1} \in \mathcal{A}_\ell(t-1)\,, \quad \omega(y_{s+1}) = 1\,.$$

This corresponds to the following scenario: the site $y_{s+1}$ is visited at iteration $t-1$ or before by the explorer $\ell$ and it is open, but the site $y_s$ was previously visited by the explorer $k$. The explorer $k$ will try to visit $y_{s+1}$ either at the time of visit of $y_s$, or at the subsequent iteration, thus $y_{s+1}$ is a relevant taboo site for the explorer $k$. If this attempt of visit of $y_{s+1}$ occurs at time $t-1$ or before, then $y_{s+1}$ is in $\overline{\mathcal{T}}_k(t-1)$, if it occurs at time $t$, then $y_{s+1}$ is in $\mathcal{T}_k^{-1}(t)$, because it was visited at time $t-1$ or before. Thus $y_{s+1}$ is in $\widehat{\mathcal{T}}_k^-(t)$. This contradicts the fact that the path $y_0, y_1, \cdots, y_r$ is in $D \setminus \widehat{\mathcal{T}}_k^-(t)$. We conclude that the site $x$ must be closed.

We treat now the set

$$\mathcal{B}\big(w_k, D \setminus \widehat{\mathcal{T}}_k^*(t), t\big) \cap \mathcal{B}\big(w_\ell, D \setminus \widehat{\mathcal{T}}_\ell^*(t), t\big)\,. \tag{88.32}$$

Let $x$ be a site belonging to the set (88.32). If $x$ belongs to

$$\mathcal{B}\big(w_k, D \setminus \widehat{\mathcal{T}}_k^*(t), t-1\big) \,\subset\, \mathcal{B}\big(w_k, D \setminus \overline{\mathcal{T}}_k(t-1), t-1\big)\,,$$

then, using also the inclusion (88.28), we are back to the previous case. From now onwards, we suppose that $x$ is not in $\mathcal{B}\big(w_k, D \setminus \widehat{\mathcal{T}}_k^*(t), t-1\big)$. By the very definition (12.7) of the ball $\mathcal{B}(\cdot,\cdot,\cdot)$, there exists a path $y_0, y_1, \cdots, y_r$ in $D \setminus \widehat{\mathcal{T}}_k^*(t)$ joining a site $y_0$ of $\mathcal{B}\big(w_k, D \setminus \widehat{\mathcal{T}}_k^*(t), t-1\big)$ to $y_r = x$ such that exactly one site among $y_0, \cdots, y_{r-1}$ is closed, say $y_s$, where $0 \leq s \leq r-1$. We claim that $y_s$ belongs to $\mathcal{A}_k(t-1)$. The site $y_0$ is in $\mathcal{A}_k(t-1)$ thanks to the inclusion (88.28). Suppose that $s \geq 1$. Since $\mathcal{A}_k(t-1)$ is an aggregate and the sites $y_0, \cdots, y_{s-1}$ are open, then they all belong to $\mathcal{A}_k(t-1)$. Furthermore, the site $y_s$ is in $D \setminus \widehat{\mathcal{T}}_k^*(t) \subset D \setminus \overline{\mathcal{T}}_k(t-1)$ and within null travel distance from $y_{s-1}$. Using the remark after (88.27), we conclude that $y_s$ is indeed in $\mathcal{A}_k(t-1)$. Recall that our goal is to prove that the site $x$ is closed. Suppose that $x$ is open. Since $y_s$ is the only closed site among $y_0, \cdots, y_{r-1}$ and we assumed that $y_r = x$ is open, then all the sites $y_{s+1}, \cdots, y_r$ are open. Yet the site $x$ belongs to $\mathcal{A}_\ell(t)$, which is an aggregate, therefore all the sites $y_{s+1}, \cdots, y_r$ must also belong to $\mathcal{A}_\ell(t)$. To sum up, we have

$$y_s \in \mathcal{A}_k(t-1)\,, \quad \omega(y_s) = 0\,, \quad y_{s+1} \in \mathcal{A}_\ell(t)\,, \quad \omega(y_{s+1}) = 1\,.$$

This corresponds to the following scenario: the site $y_{s+1}$ is visited at iteration $t$ or before by the explorer $\ell$ and it is open, but the site $y_s$ was already visited by the explorer $k$ at time $t-1$ or before. The explorer $k$ will try to visit $y_{s+1}$ either at the time of visit of $y_s$, or at the subsequent iteration, thus $y_{s+1}$ is a relevant taboo site for the explorer $k$. If this attempt of visit of $y_{s+1}$ occurs at time $t-1$ or before, then $y_{s+1}$ is in $\overline{\mathcal{T}}_k(t-1)$, if it occurs at time $t$, then $y_{s+1}$ is in $\mathcal{T}_k^{*1}(t)$. This contradicts the fact that the path $y_0, y_1, \cdots, y_r$ is in $D \setminus \widehat{\mathcal{T}}_k^*(t)$. We conclude that the site $x$ must be closed. □

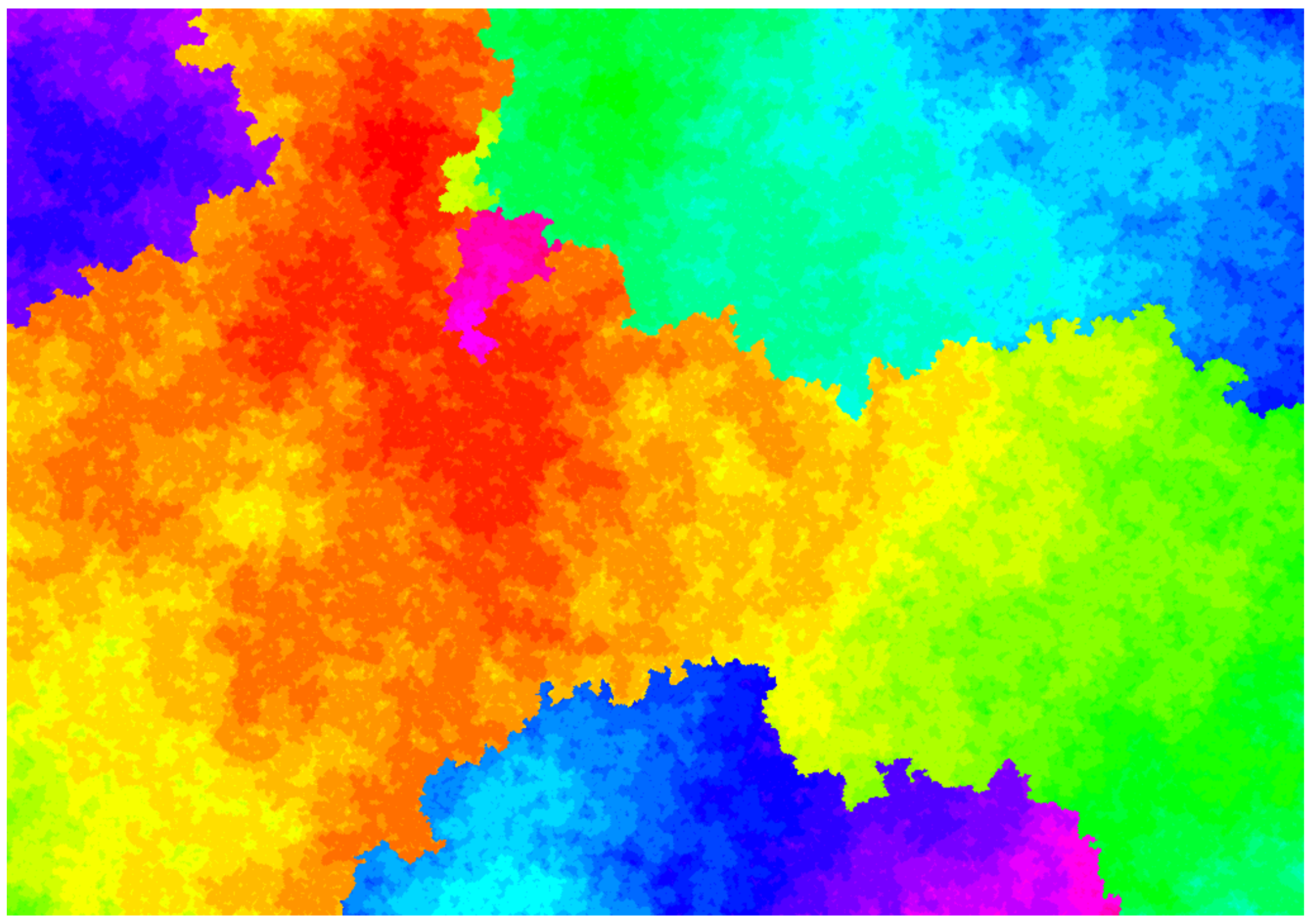

Figure 159: 7 intertwined explorations, $1024 \times 716$, site percolation, $p = 0.57$.

## 88.7 Synthesis of the inequalities

For any $t$ in $\{1,\dots,T\}$, the set $\mathcal{T}^*(t)$ is decomposed as

$$\mathcal{T}^*(t) \,=\, \mathcal{T}^{-1}(t)\cup\mathcal{T}^{+1}(t)\cup\mathcal{T}^{-0}(t)\cup\mathcal{T}^{+0}(t)\,.$$

In the previous three subsections, we have obtained inequalities on the cardinalities of the sets of the right-hand side. Let us rewrite these inequalities in the order they will be used. We fix an integer $t$ in $\{1,\dots,T\}$ and we imagine that we have a proper control on $\overline{\mathcal{T}}(t-1)$. The first inequality in (88.2) gives

$$\big|\mathcal{T}^{-1}(t)\big| \,\leq\, 2d\,\big|\mathcal{T}^{+0}(t-1)\big|\,. \tag{88.33}$$

We will then use the inequality (88.11), which says:

$$\begin{gathered}
\exists\, K^-,L^-\subset\{1,\dots,N\}\qquad K^-\cap L^-=\varnothing\,,\\
\frac{1}{8d}\,\big|\mathcal{T}^{-0}(t)\big|\,\leq\\
\Big|\Big(\bigcup_{k\in K^-}\mathcal{B}\big(w_k,D\setminus\widehat{\mathcal{T}}_k^-(t),t\big)\Big)\cap\Big(\bigcup_{\ell\in L^-}\mathcal{B}\big(w_\ell,D\setminus\overline{\mathcal{T}}_\ell(t-1),t-1\big)\Big)\Big|+\big|\overline{\mathcal{T}}(t-1)\big|\,.
\end{gathered} \tag{88.34}$$

The sets $\widehat{\mathcal{T}}_k^-(t)$ were defined in (88.5). The point is that the sets $\widehat{\mathcal{T}}_k^-(t)$, $k\in K^-$, are all included in $\overline{\mathcal{T}}(t-1)\cup\mathcal{T}^{-1}(t)$, and $\big|\mathcal{T}^{-1}(t)\big|$ has been controlled by the inequality (88.33). We move on and we use the second inequality in (88.2):

$$\big|\mathcal{T}^{+1}(t)\big| \,\leq\, 2d\,\big|\mathcal{T}^{-0}(t)\big|\,. \tag{88.35}$$

As $\big|\mathcal{T}^{-0}(t)\big|$ has been controlled by the inequality (88.34), this yields a control on $\big|\mathcal{T}^{+1}(t)\big|$. We are now ready to call the inequality (88.24), which tells:

$$\begin{gathered}
\exists\, K^+,L^+\subset\{1,\dots,N\}\qquad K^+\cap L^+=\varnothing\,,\\
\frac{1}{8d}\,\big|\mathcal{T}^{+0}(t)\big|\,\leq\,\Big|\Big(\bigcup_{k\in K^+}\mathcal{B}\big(w_k,D\setminus\widehat{\mathcal{T}}_k^*(t),t\big)\Big)\cap\Big(\bigcup_{\ell\in L^+}\mathcal{B}\big(w_\ell,D\setminus\widehat{\mathcal{T}}_\ell^*(t),t\big)\Big)\Big|\,.
\end{gathered} \tag{88.36}$$

The sets $\widehat{\mathcal{T}}_k^*(t)$ were defined in (88.14). The point is that the sets $\widehat{\mathcal{T}}_k^*(t)$, $k\in K^+$, $\widehat{\mathcal{T}}_\ell^*(t)$, $\ell\in L^+$, are all included in $\overline{\mathcal{T}}(t-1)\cup\mathcal{T}^{-0}(t)\cup\mathcal{T}^{+1}(t)$, and $\big|\mathcal{T}^{-0}(t)\big|$, $\big|\mathcal{T}^{+1}(t)\big|$ have been controlled in inequalities (88.34), (88.35). The four inequalities (88.33), (88.34), (88.35), (88.36) achieve the objective we set ourselves at the end of subsection 88.2. Indeed, the sets of sites involved in the right-hand side of the inequalities are the outcome of taboo explorations whose taboo sets do not contain any closed relevant taboo site at step $t$, apart from those already belonging to $\overline{\mathcal{T}}(t-1)$. Moreover, the open relevant taboo sites at step $t$ are controlled thanks to the two inequalities (88.33) and (88.35). We are thus in position to use the control on $|\overline{\mathcal{T}}(t-1)|$ do derive a control on $\big|\mathcal{T}^*(t)\big|$. There is still some work to be done to complete this program, we must approximate the sets appearing in the right-hand side of (88.34) and (88.36) with those visited by genuine explorations.

## 88.8 Ignition of the second stage of the proof

We use now the inequalities obtained in subsection 88.7 to conclude the first stage of the proof and to ignite the second stage. We come back to the probability we wish to estimate in (87.15). We use the inequalities (88.34), (88.36) and the standard union bound to obtain the following result. For any $t$ such that $1 \leq t \leq \overline{t}(n) = 2(\ln n)^{1-\gamma}$, we have

$$P\begin{pmatrix} \mathcal{D},\, \mathcal{E},\, \big|\overline{\mathcal{T}}^0(t)\big| > \dfrac{N^* n}{\phi(n)}\,, \\ \big|\overline{\mathcal{T}}^0(t-1)\big| \,\leq\, \rho(n)\,\big|\overline{\mathcal{T}}^0(t)\big| \end{pmatrix} \leq$$

$$P\begin{pmatrix} \mathcal{D},\, \mathcal{E}, \quad \big|\overline{\mathcal{T}}^0(t)\big| > \dfrac{N^* n}{\phi(n)}\,, \quad \big|\overline{\mathcal{T}}^0(t-1)\big| \leq \rho(n)\,\big|\overline{\mathcal{T}}^0(t)\big|\,, \\ \exists\, K^-, L^-, K^+, L^+ \subset \{\,1,\dots,N\,\} \quad K^- \cap L^- = \varnothing,\, K^+ \cap L^+ = \varnothing, \\ \dfrac{1}{8d}\,\big|\mathcal{T}^{-0}(t)\big| \,\leq\, \big|\overline{\mathcal{T}}(t-1)\big| + \\ \Big|\Big(\displaystyle\bigcup_{k\in K^-} \mathcal{B}\big(w_k, D\setminus \widehat{\mathcal{T}}_k^-(t), t\big)\Big) \cap \Big(\bigcup_{\ell\in L^-} \mathcal{B}\big(w_\ell, D\setminus \overline{\mathcal{T}}_\ell(t-1), t-1\big)\Big)\Big| \\ \dfrac{1}{8d}\,\big|\mathcal{T}^{+0}(t)\big| \,\leq \\ \Big|\Big(\displaystyle\bigcup_{k\in K^+} \mathcal{B}\big(w_k, D\setminus \widehat{\mathcal{T}}_k^*(t), t\big)\Big) \cap \Big(\bigcup_{\ell\in L^+} \mathcal{B}\big(w_\ell, D\setminus \widehat{\mathcal{T}}_\ell^*(t), t\big)\Big)\Big| \end{pmatrix}$$

$$\leq \sum_{\substack{K^-, L^-, K^+, L^+ \subset \{\,1,\dots,N\,\} \\ K^-\cap L^- = \varnothing,\, K^+\cap L^+ = \varnothing}} \tag{88.37}$$

$$P\begin{pmatrix} \mathcal{D},\, \mathcal{E}, \quad \big|\overline{\mathcal{T}}^0(t)\big| > \dfrac{N^* n}{\phi(n)}\,, \quad \big|\overline{\mathcal{T}}^0(t-1)\big| \leq \rho(n)\,\big|\overline{\mathcal{T}}^0(t)\big|\,, \\ \dfrac{1}{8d}\,\big|\mathcal{T}^{-0}(t)\big| \,\leq\, \big|\overline{\mathcal{T}}(t-1)\big| + \\ \Big|\Big(\displaystyle\bigcup_{k\in K^-} \mathcal{B}\big(w_k, D\setminus \widehat{\mathcal{T}}_k^-(t), t\big)\Big) \cap \Big(\bigcup_{\ell\in L^-} \mathcal{B}\big(w_\ell, D\setminus \overline{\mathcal{T}}_\ell(t-1), t-1\big)\Big)\Big| \\ \dfrac{1}{8d}\,\big|\mathcal{T}^{+0}(t)\big| \,\leq \\ \Big|\Big(\displaystyle\bigcup_{k\in K^+} \mathcal{B}\big(w_k, D\setminus \widehat{\mathcal{T}}_k^*(t), t\big)\Big) \cap \Big(\bigcup_{\ell\in L^+} \mathcal{B}\big(w_\ell, D\setminus \widehat{\mathcal{T}}_\ell^*(t), t\big)\Big)\Big| \end{pmatrix}.$$

This inequality concludes the first stage of the proof. We fix next four subsets $K^-, L^-, K^+, L^+$ of $\{\,1,\dots,N\,\}$ such that $K^-\cap L^- = \varnothing$, $K^+\cap L^+ = \varnothing$, and we focus on the probability appearing inside the sum in (88.37). The goal of the second stage of the proof is to approximate the taboo explorations appearing inside the probability by genuine explorations.

# 89 From taboo to genuine

In the previous section 88, we managed to localize efficiently a positive fraction of the closed relevant taboo sets with the help of two taboo explorations. We wish next to approximate these two taboo explorations by two genuine explorations. The approximations rest on a variant of the standard exploration algorithm, that we call the truncated exploration. We introduce it in subsection 89.1. We define precisely genuine explorations and sets in subsection 89.2, and we prove an important stability result in subsection 89.3. In subsection 89.4, we discuss self-avoiding explorations. We define and study the truncated balls in subsection 89.5. Based on these notions, we build in subsection 89.6 a variant of the multiple intertwined explorations, which we naturally call genuine truncated explorations, in reference to the truncation. We show then how the multiple intertwined explorations can be approximated by the genuine truncated explorations. The third stage of the proof is ignited in subsection 89.7.

## 89.1 Spatial truncation

We introduce an additional parameter $M$ in order to truncate spatially the cluster explorations. This modification is significant, because it affects the basic exploration process which is at the heart of several constructions and proofs. The truncated exploration has two distinguished features. It is a genuine exploration starting from a single site, and it is local, in the sense that it explores only sites at distance less than or equal to $M$ from the starting site. Here the distance is the one associated to the norm $|\cdot|_1$. We give next the precise definitions.

For $x$ in $\mathbb{Z}^d$ and $t \geq 0$, we denote by $B_1(x,t)$ the closed ball centered at $x$ of radius $t$ for the norm $|\cdot|_1$, it was already defined in (14.19). We extend the definition by setting, for $A$ a subset of $\mathbb{Z}^d$,

$$\forall t \geq 0 \qquad B_1(A,t) = \left\{ y \in \mathbb{Z}^d : \exists\, x \in A \quad |x-y|_1 \leq t \right\}.$$

Upon termination, the truncated exploration starting from $x$ will have explored the set $\overline{\text{Cluster}}\,(x,M,D)$ defined by

$$\overline{\text{Cluster}}\,(x,M,D) = \overline{C}\big(x, D \cap B_1(x,M)\big) \tag{89.1}$$

(the set $\overline{C}(x,D)$ was defined in (12.12)). The set $\overline{\text{Cluster}}\,(x,M,D)$ coincides with $\overline{C}(x,D)$ in case the open cluster $C(x,D)$ of $x$ in $D$ has strictly less than $M$ open sites, but in general it is a proper subset of it. For instance, the truncated exploration will not visit the sites of $B_1(x,M)$ which belong to $C(x,D)$ if they are not connected to $x$ by an open path included in $B_1(x,M)$. The point is that the truncated exploration will completely explore all the finite clusters, provided that the event $\mathcal{E}$ occurs.

We extend next the definition (89.1) of the truncated exploration by considering a whole set of starting sites. Let $A$ be a finite subset of $\mathbb{Z}^d$. We define

$$\overline{\text{Clusters}}\,(A,M,D) = \bigcup_{x \in A \cap D} \overline{\text{Cluster}}\,(x,M,D)\,. \tag{89.2}$$

An immediate consequence of this definition is that the truncated exploration of a set is local, in the sense that, for any finite subset $A$ of $\mathbb{Z}^d$,

$$\overline{\text{Clusters}}\,(A, M, D) \;\subset\; B_1(A, M)\,. \tag{89.3}$$

Now the question that comes up is: given a subset $A$ of $\mathbb{Z}^d$, possibly random, can we also realize the set $\overline{\text{Clusters}}\,(A, M, D)$ as the outcome of a genuine exploration algorithm?

In order to answer positively to this question, we have to build explicitly an exploration algorithm which computes the set $\overline{\text{Clusters}}\,(A, M, D)$. Here is one possibility. We consider a deterministic subset $A$ of $\mathbb{Z}^d$ and a finite domain $D$ included in $\mathbb{Z}^d$. We order the elements of $A \cap D$ into a sequence $a_1, \dots, a_m$ according to some deterministic rule. For the first step, we launch a truncated exploration in the domain $D$ starting from $a_1$. We denote by $\mathcal{E}_1$ the outcome of this exploration, i.e., $\mathcal{E}_1 = \overline{\text{Cluster}}\,(a_1, M, D)$. Let $k$ be an integer such that $2 \le k \le m$ and suppose that the first $k-1$ sites $a_1, \dots, a_{k-1}$ have been processed. We define $\mathcal{E}_k = \overline{\text{Cluster}}\,\big(a_k, M, D\big)$. The algorithm terminates at step $k = m$. The set of sites explored by this algorithm is

$$\mathcal{E} \;=\; \bigcup_{1 \le k \le m} \mathcal{E}_k\,.$$

The truncated exploration starting from one site is genuine, but it is not obvious that this exploration algorithm is still genuine!

To answer this question, we could construct directly a genuine algorithm whose outcome is the random set $\mathcal{E}$. We certainly could do that in the specific case of the truncated explorations by mimicking the shell exploration algorithm. In case the starting set $A$ is random, the question is even more delicate. The problem is that a random set might bring with itself some information on the configuration, which might destroy the genuineness of the algorithm. Furthermore, as we will soon iterate the previous construction, we are concerned with this case. Instead of relying on a specific construction, we will prove a more general stability result on genuine exploration algorithms in the next subsection.

Finally, we mention that there is another possibility for truncating the explorations, which may seem more natural than the above one. As we wish that the truncated explorations completely explore the clusters having strictly less than $M$ open sites, we could simply decide to stop the explorations once they have discovered $M$ open sites. The advantage is that the approximation results would in principle be better. Yet, one problem with this truncation rule is that there is no natural direct mathematical definition for the output of the exploration, and we would need to give an iterative definition. This is still possible, with the help of a variant of the formalism employed in section 13 for genuine exploration algorithms. We essentially need to introduce an adequate stopping time with respect to the sequence of the random variables driving the algorithm. We would then have a definition of $\overline{\text{Cluster}}\,(x, M, D)$ which would depend not only on $M$, but also on all the decision rules involved in the algorithm. This would cause trouble much later on in the argument, and this is why we have to work with the spatial truncation defined previously.

## 89.2 Genuine explorations and sets

Genuine exploration algorithms have been defined rigorously in section 13. The algorithms considered there were built upon a sequence of i.i.d. Bernoulli random variables $(X_k)_{k\geq 1}$ with parameter $p$. The main point was the possibility of using Hoeffding's inequality on the set of sites visited by the algorithm. We will change slightly the point of view here, by focusing more on the sequence of the sites explored by the algorithm, while the randomness driving the exploration will come directly from the percolation configuration. We highlight two essential properties of the random set of sites visited by a genuine exploration algorithm. First, this set is a stopping set, a classical concept already studied in the literature on random sets. Second, this set is a genuine set, according to the new definition 89.4. As we will discuss several mathematical aspects of genuine exploration algorithms and genuine sets, we cannot dispense from giving precise definitions, even if the formalism seems unnecessarily heavy.

Throughout the rest of this section, we fix a finite subset $D$ of $\mathbb{Z}^d$. We introduce also a cemetery state $\dagger$ that does not belong to $D$, which serves as a termination signal, and we set

$$D_\dagger \,=\, D \cup \{\,\dagger\,\}\,.$$

**Definition 89.1.** *A genuine exploration algorithm $\Phi$ in $D$ consists of a finite sequence of decision rules $\Phi = (\phi_k, 1 \leq k \leq \ell)$, where, for $1 \leq k \leq \ell$, $\phi_k$ is a map defined on $D^k \times \{\,0,1\,\}^k$ with values in $D_\dagger$.*

Until now, we have only considered exploration algorithms which never come back to a site that was already visited. This is the case if, for any $k$ in $\{\,1,\dots,\ell\,\}$, the map $\phi_k$ always proposes a new site, i.e.,

$$\forall x_1,\dots,x_k \in D \quad \forall \varepsilon_1,\dots,\varepsilon_k \in \{0,1\} \qquad \phi_k\big(x_1,\dots,x_k,\varepsilon_1,\dots,\varepsilon_k\big) \notin \{\,x_1,\dots,x_k\,\}\,. \tag{89.4}$$

This condition implies that the termination time $T$ is less or equal than $|D|$. Moreover, under (89.4), the values of the decision rules $\phi_k(x_1,\dots,x_k,\varepsilon_1,\dots,\varepsilon_k)$ are irrelevant whenever the sites $x_1,\dots,x_k$ are not pairwise distinct, as this never happens when running such an exploration. From now onwards, we relax the constraint (89.4) and we consider more general schemes of exploration, which are allowed to pass several times through a site. This is why the length $\ell$ of the sequence of decision rules might be much larger than $|D|$.

Suppose that we are given a genuine exploration algorithm $\Phi$ in $D$ and a percolation configuration in $D$. Given a starting point, we construct a random sequence by iterating the maps $(\phi_k, 1 \leq k \leq \ell)$ until we land into the cemetery state $\dagger$. To ensure the termination of the algorithm, we assume that the ultimate map $\phi_\ell$ is identically equal to $\dagger$, i.e.,

$$\forall x_1,\dots,x_\ell \in D \quad \forall \varepsilon_1,\dots,\varepsilon_\ell \in \{0,1\} \qquad \phi_\ell\big(x_1,\dots,x_\ell,\varepsilon_1,\dots,\varepsilon_\ell\big) \,=\, \dagger\,. \tag{89.5}$$

We give next the precise mathematical construction. Let $x$ be a starting point in $D$. The exploration algorithm starting from $x$ produces a random sequence of sites $X_1, X_2, \dots, X_T$ defined by the following iterative procedure:

• **Initialization:** We set $X_1 = x$.

Let $k \geq 1$ and suppose that the sites $X_1, \dots, X_k$ have been defined.

• **Termination:** If $\phi_k\big(X_1, \dots, X_k, \omega(X_1), \dots, \omega(X_k)\big)$ is equal to †, then the termination time $T$ is set to $k$, and the exploration terminates.

• **Iteration step $k$:** If $\phi_k\big(X_1, \dots, X_k, \omega(X_1), \dots, \omega(X_k)\big)$ is not equal to †, then we define $X_{k+1}$ as

$$X_{k+1} \,=\, \phi_k\big(X_1, \dots, X_k, \omega(X_1), \dots, \omega(X_k)\big)\,. \tag{89.6}$$

The condition (89.5) guarantees that the termination time $T$ is less than or equal to $\ell$. Upon termination, the algorithm returns the random sequence $X_1, \dots, X_T$. It follows from the construction that the random time $T$ is a stopping time with respect to the sequence of the observed states $\omega(X_1), \dots, \omega(X_T)$, i.e.,

$$\forall t \in \{\,1, \dots, \ell\,\} \qquad \{\,T = t\,\} \in \sigma\big(\omega(X_1), \dots, \omega(X_t)\big)\,. \tag{89.7}$$

A formula like $\sigma\big(\omega(X_1), \dots, \omega(X_t)\big)$ might create some confusion, because of the simultaneous presence of $\omega$ and the random variables $X_1, \dots, X_t$, which depend themselves on $\omega$! To dissipate any misunderstanding, we introduce the random variables $(V_i, 1 \leq i \leq \ell)$ defined by

$$\forall i \in \{\,1, \dots, \ell\,\} \qquad V_i \,=\, \begin{cases} \omega(X_i) & \text{if} \quad i \leq T\,, \\ 0 \quad \text{if} & T < i \leq \ell\,. \end{cases} \tag{89.8}$$

With this notation, the awkward-looking formula (89.7) can be rewritten safely as

$$\forall t \in \{\,1, \dots, \ell\,\} \qquad \{\,T = t\,\} \in \sigma\big(V_1, \dots, V_t\big)\,.$$

The $\sigma$-field of the events observable up to the stopping time $T$ is classically defined as

$$\begin{aligned} \mathcal{F}_T \,&=\, \sigma\big(V_t, 1 \leq t \leq T\big) \\ &=\, \Big\{\, A \subset \{0,1\}^D : \forall t \in \{\,1, \dots, \ell\,\} \quad A \cap \{\,T = t\,\} \in \sigma\big(V_1, \dots, V_t\big)\,\Big\}\,. \end{aligned} \tag{89.9}$$

Let us consider the random sequence $X_1, \dots, X_T$. The first site is deterministic, equal to $x$. A straightforward induction shows that, for any $k \geq 2$, on the event $\{\,T \geq k\,\}$, the site $X_k$ is a deterministic function of $V_1, \dots, V_{k-1}$. In light of the definition (89.9), we see that the random sequence $X_1, \dots, X_T$ is measurable with respect to the $\sigma$-field $\mathcal{F}_T$. We record this result in the next proposition, for future use.

**Proposition 89.2.** *The random sequence $X_1, \dots, X_T$ is a deterministic function of the sequence $\omega(X_1), \dots, \omega(X_T)$.*

As the exploration might visit several times a given site, the sites $X_1, \dots, X_T$ are not necessarily distinct. We form the random set containing all these sites, we call it the outcome of the genuine exploration starting from $x$, and we denote it by $\Phi(x)$:

$$\forall x \in D \qquad \Phi(x) \,=\, \big\{\, X_1, \dots, X_T \,\big\}\,.$$

These random sets are stopping sets, a concept which generalizes to random sets the stopping time property for stochastic processes indexed by one parameter. Here is the definition.

**Definition 89.3.** *A random subset $\mathcal{E}$ of $D$ is said to be a stopping set if it is measurable with respect to the states of the sites it contains, i.e.,*

$$\forall E \subset D \qquad \{\, \mathcal{E} = E \,\} \in \sigma\big(\omega(x), x \in E\big)\,.$$

Stopping sets have naturally been studied in the general theory of random sets (see definition 2.27 in chapter 5 of [106]). The recent paper [93] contains additional important references on the theory of stopping sets. The definition is more technical in a continuous setting. Fortunately, we deal here only with finite sets, for which the definition is natural and straightforward. Just to see how it works, let us prove that the outcome $\Phi(x)$ of a genuine exploration in $D$ starting from $x$ is a stopping set. Let $E$ be any subset of $D$ containing $x$, let us denote by $s$ its cardinality and by $x_1, \dots, x_s$ its elements, ordered in such a way that $x_1 = x$. Denoting by $\mathcal{S}(t,s)$ the set of the surjective maps from $\{\, 1, \dots, t \,\}$ to $\{\, 1, \dots, s \,\}$, we have

$$\{\, \Phi(x) = E \,\} \,=\, \bigcup_{s \le t \le \ell} \ \bigcup_{\substack{\varphi \in \mathcal{S}(t,s)\\ \varphi(1)=1}} \Big\{\, T = t,\, \forall k \in \{\, 1, \dots, t \,\} \quad X_k = x_{\varphi(k)} \,\Big\}$$

$$= \bigcup_{s \le t \le \ell} \ \bigcup_{\substack{\varphi \in \mathcal{S}(t,s)\\ \varphi(1)=1}} \left\{ \begin{matrix} \phi_k\big(X_1, \dots, X_k, \omega(X_1), \dots, \omega(X_k)\big) = x_{\varphi(k+1)},\ 1 \le k < t \\ \phi_t\big(X_1, \dots, X_t, \omega(X_1), \dots, \omega(X_t)\big) = \dagger \end{matrix} \right\}.$$

We see on this last expression that the event $\{\, \Phi(x) = E \,\}$ is a finite union of events that depend only on the random variables $\big(\omega(y), y \in E\big)$. This proves that $\Phi(x)$ is indeed a stopping set. Besides being a stopping set, the outcome of a genuine exploration enjoys another fundamental property, which will play a crucial role in the sequel. We call this property genuineness and we define it next.

**Definition 89.4.** *A random subset $\mathcal{E}$ of $D$ is said to be a genuine set if there exists a sequence $(U_i)_{1 \le i \le \ell}$ of i.i.d. Bernoulli random variables with parameter $p$, and a random time $S$, which is a stopping time with respect to the sequence $(U_i)_{1 \le i \le \ell}$, such that the two random vectors*

$$\Big(|\mathcal{E}|, \sum_{x \in \mathcal{E}} \omega(x)\Big) \quad \textit{and} \quad \Big(S, \sum_{1 \le i \le S} U_i\Big)$$

*have the same joint distribution.*

Notice that the definition concerns only the distribution of $\big(|\mathcal{E}|, \sum_{x\in\mathcal{E}} \omega(x)\big)$: the random variables $S$ and $(U_i)_{1\leq i\leq \ell}$ do not need to be defined on the same probability space as the percolation process. Oddly enough (or maybe not?), we already ran into a variant of this condition in the subsubsection 45.3.1 when we tried to extend the Margulis-Russo formula to a random domain, see the hypothesis 45.6 there. At that time, we claimed that the outcome of a genuine exploration satisfies this hypothesis 45.6, but without proof. Yet this result will play a central role later on, so now is the time to provide a full proof.

In the subsubsection 45.3.1, we dealt with decision rules which were required to return a new site at each step (see the condition (45.30)). The corresponding condition in our context would be the following:

$$\begin{gathered}\forall k \in \{\, 1, \dots, |D| \,\}\\ \forall x_1, \dots, x_k \text{ pairwise distinct elements of } D \quad \forall \varepsilon_1, \dots, \varepsilon_k \in \{0,1\}\\ \phi_k\big(x_1, \dots, x_k, \varepsilon_1, \dots, \varepsilon_k\big) \notin \{\, x_1, \dots, x_k \,\}\,. \end{gathered} \tag{89.10}$$

This condition implies automatically that the termination time $T$ is less or equal than $|D|$. In a first attempt, we tried to work with genuine explorations satisfying the additional condition (89.10). Everything worked quite fine for the results we are discussing here, but unfortunately we ran into trouble when trying to work out the details of the last lemma of next section (lemma 89.13). The goal of the next section is to prove that we can compose genuine explorations and still obtain a genuine set. The argument requires at some point that we exhibit a system of decision rules producing the same set as the composition of two genuine explorations. However, when working with rules satisfying (89.10), it seems that additional complicated compatibility conditions would have to be imposed in order for the proof to work. This is why in the end we relaxed the definition of a genuine algorithm and we removed the condition (89.10).

The goal of the rest of this section is to show that the outcome of a genuine exploration is a genuine set. We first gather some useful constructions and intermediate results on genuine explorations. We consider a genuine exploration algorithm $\Phi$ in $D$ and a site $x$ of $D$. We denote by $T$ the termination time of the exploration starting from $x$. Let $X_1, X_2, \dots, X_T$ be the sequence of sites produced by the exploration starting from $x$. We define the cleansing $\widehat{X}_1, \dots, \widehat{X}_{\widehat{T}}$ of the sequence $X_1, \dots, X_T$ to be the sequence obtained by keeping only the first occurrence of each site in $X_1, \dots, X_T$. It is formally defined by the following iterative procedure:

• **Initialization:** We set $i_1 = 1$ and $\widehat{X}_1 = X_1 = x$.

Let $k \geq 1$ and suppose that $i_k$ and $\widehat{X}_1, \dots, \widehat{X}_k$ have been defined.

• **Termination:** If $\{\, X_{i_k+1}, \dots, X_T \,\}$ is included in $\{\, \widehat{X}_1, \dots, \widehat{X}_k \,\}$, then we set $\widehat{T} = k$ and the construction is finished.

• **Iteration step $k$:** Otherwise, we define

$$i_{k+1} \,=\, \min\big\{\, i > i_k : X_i \notin \{\, \widehat{X}_1, \dots, \widehat{X}_k \,\}\,\big\}\,, \qquad \widehat{X}_{k+1} \,=\, X_{i_{k+1}}\,.$$

An immediate consequence of the construction is that

$$\forall k \in \{\,1,\dots,\widehat{T}\,\} \qquad \big|\{\,\widehat{X}_1,\dots,\widehat{X}_k\,\}\big| \,=\, k\,. \tag{89.11}$$

We explain next how to compute the position of an element in the cleansing. We define

$$\forall k \in \{\,1,\dots,T\,\} \qquad \varphi(k) \,=\, \min\Big\{\,\big|\{\,X_1,\dots,X_i\,\}\big| : i \in \{\,1,\dots,k\,\}, X_k = X_i\,\Big\}. \tag{89.12}$$

**Lemma 89.5.** *For any $k$ in $\{\,1,\dots,T\,\}$, we have $X_k = \widehat{X}_{\varphi(k)}$.*

*Proof.* We prove the result by induction. We have $X_1 = \widehat{X}_{\varphi(1)}$. Let $k$ be such that $1 \le k < T$, and suppose that we have proved that

$$\forall i \in \{\,1,\dots,k\,\} \qquad X_i = \widehat{X}_{\varphi(i)}\,.$$

We consider two cases:

$\star$ $X_{k+1} \in \{\,X_1,\dots,X_k\,\}$. Let $i$ be the index of the first occurrence of $X_{k+1}$, i.e.,

$$i \,=\, \min\,\big\{\,j \in \{\,1,\dots,k\,\} : X_j = X_{k+1}\,\big\}\,.$$

It follows from formula (89.12) that $\varphi(k+1) = \varphi(i)$ and from the induction hypothesis that

$$X_{k+1} = X_i = \widehat{X}_{\varphi(i)} = \widehat{X}_{\varphi(k+1)}\,.$$

$\star$ $X_{k+1} \not\in \{\,X_1,\dots,X_k\,\}$. Formulas (89.11) and (89.12) give

$$\big|\{\,\widehat{X}_1,\dots,\widehat{X}_{\varphi(k+1)}\,\}\big| \,=\, \varphi(k+1) \,=\, \big|\{\,X_1,\dots,X_{k+1}\,\}\big| \,=\, \big|\{\,X_1,\dots,X_k\,\}\big| + 1\,. \tag{89.13}$$

By the induction hypothesis, we have

$$\big\{\,\widehat{X}_1,\dots,\widehat{X}_{\varphi(k)}\,\big\} \,=\, \big\{\,X_1,\dots,X_k\,\big\}\,. \tag{89.14}$$

Since $\varphi(k+1) \le k+1$, then $\widehat{X}_{\varphi(k+1)}$ belongs to $\big\{\,X_1,\dots,X_{k+1}\,\big\}$. The equalities (89.13) and (89.14) imply that $X_{k+1} = \widehat{X}_{\varphi(k+1)}$. This completes the induction step and the proof. □

In general, we cannot recover a sequence from its cleansing. However, in the particular case of a sequence produced by an exploration algorithm, we have the following useful result.

**Proposition 89.6.** *The two sequences $X_1,\dots,X_T$ and $\omega(X_1),\dots,\omega(X_T)$ are deterministic functions of the sequence of states $\omega(\widehat{X}_1),\dots,\omega(\widehat{X}_{\widehat{T}})$.*

Here is an explanation for this result. The exploration algorithm which builds the sequence $X_1,\dots,X_T$ starts from $x$ and examine successively the states of the sites of this sequence. At each step, it gathers the information discovered so far to decide the next site to be examined. In total, it will have examined exactly the states of the sites $\widehat{X}_1,\dots,\widehat{X}_{\widehat{T}}$. The stopping criterion is also checked with the help of this information. This argument is very convincing, so why should we write a proof of this result? In fact, the writing of this kind of proof, however uninteresting it might be in the end, is important to fix little details that might become important at some unexpected stage. In several instances, the demand for a full proof ultimately lead to non-negligible modifications of the initial constructions.

*Proof.* We shall use again the auxiliary variables $(V_i, 1 \leq i \leq \ell)$ introduced in (89.8). As we discussed previously, $T$ is a stopping time with respect to the sequence $(V_i, 1 \leq i \leq \ell)$, and the random sequence $X_1, \dots, X_T$ is measurable with respect to the $\sigma$-field $\mathcal{F}_T$. This means that $X_1, \dots, X_T$ can be obtained through a deterministic procedure once we know the values of the random variables $V_1, \dots, V_T$ (notice that $T$ is a random time!). To make the proof more transparent, we introduce the variables $(\widehat{V}_i, 1 \leq i \leq \ell)$ given by

$$\forall i \in \{\, 1, \dots, \ell \,\} \qquad \widehat{V}_i \,=\, \begin{cases} \omega(\widehat{X}_i) & \text{if} \quad i \leq \widehat{T}\,, \\ 0 \quad \text{if} & \widehat{T} < i \leq \ell\,. \end{cases} \tag{89.15}$$

The sequence $(\widehat{V}_i, 1 \leq i \leq \ell)$ is built by completing $(\omega(\widehat{X}_i), 1 \leq i \leq \widehat{T})$ with zeroes until we obtain a sequence of length $\ell$. We show next that the sequence $V_1, \dots, V_T$ can be recovered once we know the sequence $(\widehat{V}_1, \dots, \widehat{V}_\ell)$. As the sequence $X_1, \dots, X_T$ is defined iteratively, it seems that the natural way to obtain a rigorous proof is to define another iterative procedure that builds $X_1, \dots, X_T$, which is exclusively based on the information contained in $(\widehat{V}_1, \dots, \widehat{V}_\ell)$. This is the purpose of the next construction. We define iteratively a sequence $(\widetilde{X}_k, 1 \leq k \leq T^*)$ of sites in $D$ and a map $\varphi$ from $\{\, 1, \dots, T^* \,\}$ to $\{\, 1, \dots, \widehat{T} \,\}$ as follows:

• **Initialization:** We set $\widetilde{X}_1 = x$ and $\varphi(1) = 1$.

Let $k \geq 1$ and suppose that $\widetilde{X}_1, \dots, \widetilde{X}_k$, $\varphi(1), \dots, \varphi(k)$ have been defined.

• **Termination:** If $\phi_k\big(\widetilde{X}_1, \dots, \widetilde{X}_k, \widehat{V}_{\varphi(1)}, \dots, \widehat{V}_{\varphi(k)}\big)$ is equal to †, then we set $T^* = k$ and the construction is finished.

• **Iteration step $k$:** Suppose that $\phi_k\big(\widetilde{X}_1, \dots, \widetilde{X}_k, \widehat{V}_{\varphi(1)}, \dots, \widehat{V}_{\varphi(k)}\big)$ is not equal to †. We define $\widetilde{X}_{k+1}$ as

$$\widetilde{X}_{k+1} \,=\, \phi_k\big(\widetilde{X}_1, \dots, \widetilde{X}_k, \widehat{V}_{\varphi(1)}, \dots, \widehat{V}_{\varphi(k)}\big)\,. \tag{89.16}$$

To define $\varphi(k+1)$, we consider two cases:

$\star$ $\widetilde{X}_{k+1} \in \{\, \widetilde{X}_1, \dots, \widetilde{X}_k \,\}$. We define $\varphi(k+1)$ as

$$\varphi(k+1) \,=\, \min \Big\{ \big|\{\, \widetilde{X}_1, \dots, \widetilde{X}_i \,\}\big| : i \in \{\, 1, \dots, k \,\}, \widetilde{X}_{k+1} = \widetilde{X}_i \Big\}\,.$$

$\star$ $\widetilde{X}_{k+1} \notin \{\, \widetilde{X}_1, \dots, \widetilde{X}_k \,\}$. We define $\varphi(k+1)$ as

$$\varphi(k+1) \,=\, \big|\{\, \widetilde{X}_1, \dots, \widetilde{X}_k \,\}\big| + 1\,.$$

From the construction, we see that the whole sequence $\widetilde{X}_1, \dots, \widetilde{X}_{T^*}$ and the map $\varphi$ are determined by the sequence $\widehat{V}_1, \dots, \widehat{V}_\ell$. We show next that $\widetilde{X}_1, \dots, \widetilde{X}_{T^*}$ coincides with the sequence we were trying to reconstruct.

**Lemma 89.7.** *For any index $k$ such that $1 \leq k \leq \min(T, T^*)$, we have $\widetilde{X}_k = X_k$ and $\widetilde{X}_k = \widehat{X}_{\varphi(k)}$.*

*Proof.* We prove the lemma by induction. We have $\widetilde{X}_1 = X_{\varphi(1)} = X_1 = x$. Let $k$ be such that $1 \leq k < \min(T, T^*)$, and suppose that we have proved that

$$\forall i \in \{ 1, \dots, k \} \qquad \widetilde{X}_i = X_i = \widehat{X}_{\varphi(i)} \,.$$

The construction rule (89.16), together with the definition (89.15) of the variables $(\widehat{V}_i, 1 \leq i \leq \ell)$, yields that

$$\begin{aligned} \widetilde{X}_{k+1} &= \phi_k\big(\widetilde{X}_1, \dots, \widetilde{X}_k, \widehat{V}_{\varphi(1)}, \dots, \widehat{V}_{\varphi(k)}\big) \\ &= \phi_k\big(\widetilde{X}_1, \dots, \widetilde{X}_k, \omega(\widehat{X}_{\varphi(1)}), \dots, \omega(\widehat{X}_{\varphi(k)})\big) \\ &= \phi_k\big(\widetilde{X}_1, \dots, \widetilde{X}_k, \omega(\widetilde{X}_1), \dots, \omega(\widetilde{X}_k)\big) \\ &= \phi_k\big(X_1, \dots, X_k, \omega(X_1), \dots, \omega(X_k)\big) \,=\, X_{k+1} \,. \end{aligned}$$

It remains to prove that $\widetilde{X}_{k+1} = \widehat{X}_{\varphi(k+1)}$. Following the definition of $\varphi(k+1)$, we consider two cases:

$\star$ $\widetilde{X}_{k+1} \in \{ \widetilde{X}_1, \dots, \widetilde{X}_k \}$. Let $i$ be the smallest index in $\{ 1, \dots, k \}$ such that $\widetilde{X}_{k+1} = \widetilde{X}_i$. We have then $\varphi(k+1) = \varphi(i)$ and

$$\widetilde{X}_{k+1} = \widetilde{X}_i = \widehat{X}_{\varphi(i)} = \widehat{X}_{\varphi(k+1)} \,.$$

$\star$ $\widetilde{X}_{k+1} \notin \{ \widetilde{X}_1, \dots, \widetilde{X}_k \}$. By the induction hypothesis, we have

$$\{ \widehat{X}_1, \dots, \widehat{X}_{\varphi(k)} \} \,=\, \{ X_1, \dots, X_k \} \,=\, \{ \widetilde{X}_1, \dots, \widetilde{X}_k \} \,, \tag{89.17}$$

and we have also proved above that $X_{k+1} = \widetilde{X}_{k+1}$. Using formulas (89.11) and (89.12), we obtain then

$$\begin{aligned} \big|\{ \widehat{X}_1, \dots, \widehat{X}_{\varphi(k+1)} \}\big| \,&=\, \varphi(k+1) \,=\, \big|\{ X_1, \dots, X_{k+1} \}\big| \\ &= \big|\{ \widetilde{X}_1, \dots, \widetilde{X}_{k+1} \}\big| \,=\, \big|\{ \widetilde{X}_1, \dots, \widetilde{X}_k \}\big| + 1 \,. \end{aligned} \tag{89.18}$$

Since $\varphi(k+1) \leq k+1$, then $\widehat{X}_{\varphi(k+1)}$ belongs to $\{ X_1, \dots, X_{k+1} \}$. The equalities (89.17) and (89.18) imply that $\widehat{X}_{\varphi(k+1)} = X_{k+1} = \widetilde{X}_{k+1}$. This completes the induction step and the proof of the lemma. □

In view of lemma 89.7, the only point that remains to be proved is that $T = T^*$. Suppose that $T < T^*$. Let us consider the situation at time $T$. From the definition of the stopping time $T$, we know that

$$\phi_T\big(X_1, \dots, X_T, \omega(X_1), \dots, \omega(X_T)\big)$$

is then equal to †. This would imply that

$$\begin{aligned}\phi_T\big(\widetilde{X}_1,\dots,\widetilde{X}_T,\widehat{V}_{\varphi(1)},\dots,\widehat{V}_{\varphi(T)}\big) &= \phi_T\big(\widetilde{X}_1,\dots,\widetilde{X}_T,\omega(\widetilde{X}_1),\dots,\omega(\widetilde{X}_T)\big)\\ &= \phi_T\big(X_1,\dots,X_T,\omega(X_1),\dots,\omega(X_T)\big)\end{aligned}$$

is also equal to †, so that $T^*$ should also be equal to $T$. Symmetrically, suppose that $T^* < T$. By the definition of $T^*$,

$$\phi_{T^*}\big(\widetilde{X}_1,\dots,\widetilde{X}_{T^*},\omega(\widetilde{X}_1),\dots,\omega(\widetilde{X}_{T^*})\big) \tag{89.19}$$

would be equal to †. By lemma 89.7, the quantity (89.19) would be equal to

$$\phi_{T^*}\big(X_1,\dots,X_{T^*},\omega(X_1),\dots,\omega(X_{T^*})\big)\,,$$

so that $T$ should be equal to $T^*$, and we arrive again at a contradiction. Therefore $T^* = T$. As the sequence $\widetilde{X}_1,\dots,\widetilde{X}_{T^*}$ is completely determined by the starting site $x$ and $\widehat{V}_1,\dots,\widehat{V}_{\widehat{T}}$, we have reached the desired conclusion. ☐

We are now ready to state and prove the main result of this section.

**Proposition 89.8.** *Let $\Phi$ be a genuine exploration algorithm in $D$. For any $x$ in $D$, its outcome $\Phi(x)$ starting from $x$ is a genuine random subset of $D$.*

*Proof.* We consider $\Phi$, $x$ as in the statement of the proposition. According to the definition 89.4 of a genuine set, we have to build a sequence $(U_i)_{1\le i\le |D|}$ of i.i.d. Bernoulli random variables with parameter $p$ and a random index $S$, which is a stopping time with respect to the sequence $(U_i)_{1\le i\le|D|}$, so that

$$\Big(|\Phi(x)|,\sum_{y\in\Phi(x)}\omega(y)\Big)\quad\text{and}\quad\Big(S,\sum_{1\le i\le S}U_i\Big)$$

have the same joint distribution. Let $X_1,\dots,X_T$ be the sequence of sites produced by the exploration $\Phi$ starting from $x$, and let $\widehat{X}_1,\dots,\widehat{X}_{\widehat{T}}$ be its cleansing. Let $(V_i)_{1\le i\le|D|}$ be a sequence of i.i.d. Bernoulli random variables with parameter $p$, which is completely independent of the percolation process. The relevant sequence $(U_i)_{1\le i\le|D|}$ is defined as

$$\forall i\in\{\,1,\dots,|D|\,\}\qquad U_i\,=\,\begin{cases}\omega(\widehat{X}_i) & \text{if}\quad i\le\widehat{T}\,,\\ V_i & \text{if}\quad i>\widehat{T}\,.\end{cases} \tag{89.20}$$

For the random index $S$, we take simply $S=\widehat{T}$. Let us check that this construction enjoys the desired properties. By the definition of $\Phi(x)$ and the construction of the cleansing sequence, we have

$$\big|\Phi(x)\big|\,=\,\big|\{\,X_1,\dots,X_T\,\}\big|\,=\,\big|\{\,\widehat{X}_1,\dots,\widehat{X}_{\widehat{T}}\,\}\big|\,=\,\widehat{T}\,=\,S\,. \tag{89.21}$$

In addition, we have the equality

$$\sum_{y\in\Phi(x)}\omega(y)\,=\,\sum_{1\le i\le\widehat{T}}\omega(\widehat{X}_i)\,=\,\sum_{1\le i\le S}U_i\,.$$

We prove next that $\widehat{T}$ is a stopping time with respect to the sequence $(U_i)_{1\leq i\leq |D|}$. We see from formula (89.21) that $\widehat{T}$ is measurable with respect to the sequence $X_1,\dots,X_T$. In turn, proposition 89.6 tells us that the sequence $X_1,\dots,X_T$ is a deterministic function of the sequence of states $\omega(\widehat{X}_1),\dots,\omega(\widehat{X}_{\widehat{T}})$. From the previous discussion, and the definition (89.20) of the sequence $(U_i)_{1\leq i\leq |D|}$, we conclude that

$$\forall \widehat{t}\in\{\,1,\dots,|D|\,\}\qquad \{\,\widehat{T}=\widehat{t}\,\}\in\sigma\big(\omega(\widehat{X}_1),\dots,\omega(\widehat{X}_{\widehat{t}})\big)\;=\;\sigma\big(U_1,\dots,U_{\widehat{t}}\big)\,.$$

It remains to prove that the sequence $(U_i)_{1\leq i\leq |D|}$ is i.i.d.. Let $(\delta_i, 1\leq i\leq |D|)$ be a fixed sequence of elements of $\{0,1\}$ such that

$$P\big(U_i=\delta_i,\,1\leq i\leq |D|\big)\;>\;0\,. \tag{89.22}$$

We wish to compute the above probability. To that end, we condition on $\widehat{T}$ and we write

$$P\big(U_i=\delta_i,\,1\leq i\leq |D|\big)\;=\;\sum_{1\leq \widehat{t}\leq |D|}P\big(U_i=\delta_i,\,1\leq i\leq |D|,\,\widehat{T}=\widehat{t}\,\big)\,. \tag{89.23}$$

It turns out that at most one term in the previous sum is non-zero, as we prove in the next lemma. We could give an argument specific to our current construction, however the underlying reason is the stopping time structure. So we present next the general result, valid for any stopping time for a random binary sequence.

**Lemma 89.9.** *Let $(B_i, 1\leq i\leq \ell)$ be a sequence of random variables defined on a probability space $(\Omega,\mathcal{F},\mathcal{P})$ with values in $\{0,1\}$, and let $T:\Omega\to\{\,1,\dots,\ell\,\}$ be a stopping time for $(B_i, 1\leq i\leq \ell)$. Let $(b_i, 1\leq i\leq \ell)$ be a fixed sequence of elements of $\{0,1\}$ such that the event $\big\{\,B_1=b_1,\dots,B_\ell=b_\ell\big\}$ is not empty. There exists at most one index $t^*$ in $\{\,1,\dots,\ell\,\}$ such that*

$$\Big\{\,B_1=b_1,\dots,B_T=b_T,T=t^*\,\Big\}\neq\varnothing\,. \tag{89.24}$$

*Suppose that* (89.24) *holds for the index $t^*$. We have in addition*

$$\Big\{\,B_1=b_1,\dots,B_{t^*}=b_{t^*}\,\Big\}\;\subset\;\{\,T=t^*\,\}\,. \tag{89.25}$$

*Proof.* The stopping time property means that

$$\forall t\in\{\,1,\dots,\ell\,\}\qquad \{\,T=t\,\}\in\sigma\big(B_1,\dots,B_t\big)\,.$$

Thus, for any $t$ in $\{\,1,\dots,\ell\,\}$, there exists a subset $A(t)$ of $\{0,1\}^t$ such that

$$\{\,T=t\,\}\;=\;\Big\{\omega\in\Omega:\big(B_1,\dots,B_t\big)\in A(t)\Big\}\,. \tag{89.26}$$

As the sets $\{\,T=t\,\}$, $1\leq t\leq \ell$, are pairwise disjoint, we have also

$$\forall t\in\{\,1,\dots,\ell\,\}\quad \forall s\in\{\,1,\dots,t-1\,\}$$
$$(b_1,\dots,b_t)\in A(t)\quad\Longrightarrow\quad (b_1,\dots,b_s)\notin A(s)\,. \tag{89.27}$$

Suppose that there exist two values $t_1, t_2$ such that $1 \leq t_1 < t_2 \leq \ell$ and the two events

$$\mathcal{E}_i \,=\, \Big\{\, B_1 = b_1, \dots, B_T = b_T, T = t_i \,\Big\}\,, \qquad i = 1, 2\,, \tag{89.28}$$

are not empty. Let us pick up a configuration $\omega_1$ in $\mathcal{E}_1$ and a configuration $\omega_2$ in $\mathcal{E}_2$. Let us fix $i$ in $\{1, 2\}$. The configuration $\omega_i$ is such that

$$B_1(\omega_i) = b_1, \dots, B_{t_i}(\omega_i) = b_{t_i}, \quad \text{and} \quad T = t_i\,.$$

The equality (89.26) implies that $(b_1, \dots, b_{t_i})$ is in $A(t_i)$. Thus we have

$$t_1 < t_2\,, \quad (b_1, \dots, b_{t_1}) \in A(t_1)\,, \quad (b_1, \dots, b_{t_2}) \in A(t_2)\,.$$

This stands in contradiction with the condition (89.27). Therefore we cannot have two events $\mathcal{E}_1, \mathcal{E}_2$ as in (89.28) which are non-empty.

Suppose next that (89.24) holds for the index $t^*$. This implies that

$$(b_1, \dots, b_{t^*}) \in A(t^*)\,. \tag{89.29}$$

The inclusion (89.25) is a consequence of (89.29). □

We apply lemma 89.9 to the sequence $(U_i)_{1 \leq i \leq |D|}$ and the stopping time $\widehat{T}$. The sum in formula (89.23) reduces to a unique term. To alleviate the next formulas, we denote by $s$ the value $t^*$ given by lemma 89.9. We have

$$\begin{aligned} P\big(U_i = \delta_i,\, 1 \leq i \leq |D|\big) \,&=\, P\big(U_i = \delta_i,\, 1 \leq i \leq |D|, \widehat{T} = s\big) \\ &=\, P\big(\omega(\widehat{X}_1) = \delta_1, \dots, \omega(\widehat{X}_s) = \delta_s, V_{s+1} = \delta_{s+1}, \dots, V_{|D|} = \delta_{|D|}, \widehat{T} = s\big)\,. \end{aligned} \tag{89.30}$$

As the sequence $(V_i)_{1 \leq i \leq |D|}$ is independent of the percolation process, we have

$$\begin{aligned} &P\big(\omega(\widehat{X}_1) = \delta_1, \dots, \omega(\widehat{X}_s) = \delta_s, V_{s+1} = \delta_{s+1}, \dots, V_{|D|} = \delta_{|D|}, \widehat{T} = s\big) \,= \\ &\quad P\big(\omega(\widehat{X}_1) = \delta_1, \dots, \omega(\widehat{X}_s) = \delta_s, \widehat{T} = s\big) P\big(V_{s+1} = \delta_{s+1}, \dots, V_{|D|} = \delta_{|D|}\big)\,, \end{aligned} \tag{89.31}$$

and we focus now on the probability $P\big(\omega(\widehat{X}_1) = \delta_1, \dots, \omega(\widehat{X}_s) = \delta_s, \widehat{T} = s\big)$. We know from the condition (89.22) that this probability is not equal to 0. We decompose it according to the values of the sequence $X_1, \dots, X_T$ and its states:

$$P\big(\omega(\widehat{X}_1) = \delta_1, \dots, \omega(\widehat{X}_s) = \delta_s, \widehat{T} = s\big) \,= \sum_{\substack{1 \leq t \leq \ell \\ x_1, \dots, x_t \in D \\ \widehat{x}_1, \dots, \widehat{x}_s \in D \\ \varepsilon_1, \dots, \varepsilon_t \in \{0,1\}}} P \begin{pmatrix} X_1 = x_1, \dots, X_t = x_t, T = t, \\ \omega(X_1) = \varepsilon_1, \dots, \omega(X_t) = \varepsilon_t, \\ \widehat{X}_1 = \widehat{x}_1, \dots, \widehat{X}_s = \widehat{x}_s, \widehat{T} = s, \\ \omega(\widehat{X}_1) = \delta_1, \dots, \omega(\widehat{X}_s) = \delta_s \end{pmatrix}. \tag{89.32}$$

By proposition 89.6, the two sequences $X_1, \dots, X_T$ and $\omega(X_1), \dots, \omega(X_T)$ are deterministic functions of the sequence of states $\omega(\widehat{X}_1), \dots, \omega(\widehat{X}_{\widehat{T}})$. In addition, the sequence $\widehat{X}_1, \dots, \widehat{X}_{\widehat{T}}$ is a deterministic function of the sequence $X_1, \dots, X_T$. Therefore, once we fix the values $s$ and $\delta_1, \dots, \delta_s$, there is a unique choice for the value of $t \geq s$, the sequences of sites $(x_1, \dots, x_t)$, $(\widehat{x}_1, \dots, \widehat{x}_s)$ in $D$ and states $\varepsilon_1, \dots, \varepsilon_t$ in $\{0, 1\}$, for which the probability

$$P\begin{pmatrix} X_1 = x_1, \dots, X_t = x_t, T = t, \\ \omega(X_1) = \varepsilon_1, \dots, \omega(X_t) = \varepsilon_t, \\ \widehat{X}_1 = \widehat{x}_1, \dots, \widehat{X}_s = \widehat{x}_s, \widehat{T} = s, \\ \omega(\widehat{X}_1) = \delta_1, \dots, \omega(\widehat{X}_s) = \delta_s \end{pmatrix}$$

might be non-zero. Because of the construction of the sequence $X_1, \dots, X_T$, and the definition of the termination time $T$, the two sequences $x_1, \dots, x_t$ and $\varepsilon_1, \dots, \varepsilon_t$ are such that $x_1 = x$,

$$\forall k \in \{1, \dots, t-1\} \qquad x_{k+1} = \phi_k\big(x_1, \dots, x_k, \varepsilon_1, \dots, \varepsilon_k\big), \tag{89.33}$$

$$\phi_t\big(x_1, \dots, x_t, \varepsilon_1, \dots, \varepsilon_t\big) = \dagger. \tag{89.34}$$

Furthermore, the sequence $\widehat{x}_1, \dots, \widehat{x}_s$ has to be the cleansing of the sequence $x_1, \dots, x_t$, so that

$$\{x_1, \dots, x_t\} = \{\widehat{x}_1, \dots, \widehat{x}_s\}. \tag{89.35}$$

Finally, the states $\varepsilon_1, \dots, \varepsilon_t$ must be compatible with the values $\delta_1, \dots, \delta_s$, in the sense that

$$\forall k \in \{1, \dots, t\} \quad \forall h \in \{1, \dots, s\} \qquad x_k = \widehat{x}_h \implies \varepsilon_k = \delta_h. \tag{89.36}$$

In any case, the sum in (89.32) is reduced to a single term, the term corresponding to the choices for $t$, $(x_1, \dots, x_t)$, $(\widehat{x}_1, \dots, \widehat{x}_s)$, and $\varepsilon_1, \dots, \varepsilon_t$ described above, i.e.,

$$P\big(\omega(\widehat{X}_1) = \delta_1, \dots, \omega(\widehat{X}_s) = \delta_s, \widehat{T} = s\big) = P\begin{pmatrix} X_1 = x_1, \dots, X_t = x_t, T = t, \\ \omega(X_1) = \varepsilon_1, \dots, \omega(X_t) = \varepsilon_t, \\ \widehat{X}_1 = \widehat{x}_1, \dots, \widehat{X}_s = \widehat{x}_s, \widehat{T} = s, \\ \omega(\widehat{X}_1) = \delta_1, \dots, \omega(\widehat{X}_s) = \delta_s \end{pmatrix}. \tag{89.37}$$

The cleansed sequence $\widehat{X}_1, \dots, \widehat{X}_{\widehat{T}}$ is a deterministic function of the sequence $X_1, \dots, X_T$, thus the values $s$ and $\delta_1, \dots, \delta_s$ can be recovered from the value $t$ and the two sequences $x_1, \dots, x_t$, $\varepsilon_1, \dots, \varepsilon_t$, whence

$$P\begin{pmatrix} X_1 = x_1, \dots, X_t = x_t, T = t, \\ \omega(X_1) = \varepsilon_1, \dots, \omega(X_t) = \varepsilon_t, \\ \widehat{X}_1 = \widehat{x}_1, \dots, \widehat{X}_s = \widehat{x}_s, \widehat{T} = s, \\ \omega(\widehat{X}_1) = \delta_1, \dots, \omega(\widehat{X}_s) = \delta_s \end{pmatrix} = P\begin{pmatrix} X_1 = x_1, \dots, X_t = x_t, T = t, \\ \omega(X_1) = \varepsilon_1, \dots, \omega(X_t) = \varepsilon_t \end{pmatrix}. \tag{89.38}$$

We claim that

$$\left\{ \begin{matrix} X_1 = x_1, \dots, X_t = x_t, T = t, \\ \omega(X_1) = \varepsilon_1, \dots, \omega(X_t) = \varepsilon_t \end{matrix} \right\} = \Big\{ \omega(x_1) = \varepsilon_1, \dots, \omega(x_t) = \varepsilon_t \Big\} . \quad (89.39)$$

Obviously, the event on the left is included in the event on the right. We prove the converse inequality in the next lemma.

**Lemma 89.10.** *Suppose that the event* $\Big\{ \omega(x_1) = \varepsilon_1, \dots, \omega(x_t) = \varepsilon_t \Big\}$ *occurs. We have then* $T = t$ *and*

$$\forall k \in \{ 1, \dots, t \} \qquad X_k = x_k \, .$$

*Proof.* We prove first by induction that

$$\forall k \in \{ 1, \dots, \min(t, T) \} \qquad X_k = x_k \, . \quad (89.40)$$

We have $X_1 = x_1 = x$. Let $k$ be such that $1 \leq k < \min(t, T)$, and suppose that we have proved that

$$\forall i \in \{ 1, \dots, k \} \qquad X_i = x_i \, .$$

The construction rules (89.6) and (89.33), together with the induction hypothesis, yield that

$$\begin{aligned} X_{k+1} \,&=\, \phi_k\big(X_1, \dots, X_k, \omega(X_1), \dots, \omega(X_k)\big) \\ &=\, \phi_k\big(x_1, \dots, x_k, \varepsilon_1, \dots, \varepsilon_k\big) \,=\, x_{k+1} \, . \end{aligned}$$

This completes the induction step and the proof of (89.40). It remains to prove that $T = t$. Suppose that $T < t$. The definition of $T$ implies that

$$\phi_T\big(X_1, \dots, X_T, \omega(X_1), \dots, \omega(X_T)\big) \quad (89.41)$$

is equal to †. With the help of (89.40), we see that the quantity (89.41) is equal to

$$\phi_T\big(x_1, \dots, x_T, \varepsilon_1, \dots, \varepsilon_T\big) \, , \quad (89.42)$$

but it follows from (89.33) that the quantity (89.42) is equal to the site $x_{T+1}$ of $D$, which is absurd. Suppose that $T > t$. The condition (89.34) implies that

$$\phi_t\big(x_1, \dots, x_t, \varepsilon_1, \dots, \varepsilon_t\big) \quad (89.43)$$

is equal to †. With the help of (89.40), we see that the quantity (89.43) is equal to

$$\phi_t\big(X_1, \dots, X_t, \omega(X_1), \dots, \omega(X_t)\big) \, , \quad (89.44)$$

but it follows from the definition of $T$ that the quantity (89.44) is equal to the site $X_{t+1}$ of $D$, which is absurd. The only possibility is that $T = t$. ☐

Another possibility to prove the equality (89.39) is to make appeal to the proposition 89.2, which tells that

$$\left\{\begin{matrix} X_1 = x_1, \dots, X_t = x_t, T = t, \\ \omega(X_1) = \varepsilon_1, \dots, \omega(X_t) = \varepsilon_t \end{matrix}\right\} = \Big\{ \omega(X_1) = \varepsilon_1, \dots, \omega(X_t) = \varepsilon_t, T = t \Big\} . \tag{89.45}$$

However, to show that the events on the right in formulas (89.39) and (89.45) are the same, it seems that we would need to argue essentially as in lemma 89.10.

Using the chain of equalities (89.37), (89.38), (89.39), we obtain finally that

$$P\big(\omega(\widehat{X}_1) = \delta_1, \dots, \omega(\widehat{X}_s) = \delta_s, \widehat{T} = s\big) = P\big(\omega(x_1) = \varepsilon_1, \dots, \omega(x_t) = \varepsilon_t\big) . \tag{89.46}$$

The point of this reformulation is to show that the probability ultimately depends on the states of the sites $x_1, \dots, x_t$. It follows from the formulas (89.33), (89.34), (89.35) and (89.36) that the sites $x_1, \dots, x_t$ and the values $\varepsilon_1, \dots, \varepsilon_t$ are deterministic functions of $\delta_1, \dots, \delta_s$. Furthermore, using again the formulas (89.35) and (89.36), we have

$$P\big(\omega(x_1) = \varepsilon_1, \dots, \omega(x_t) = \varepsilon_t\big) = P\big(\omega(\widehat{x}_1) = \delta_1, \dots, \omega(\widehat{x}_s) = \delta_s\big) . \tag{89.47}$$

Substituting successively (89.47) in (89.46) and then in (89.31) and (89.30), we obtain finally the much-sought-after formula

$$\begin{aligned} &P\big(U_i = \delta_i,\, 1 \le i \le |D|\big) = \\ &\qquad P\big(\omega(\widehat{x}_1) = \delta_1, \dots, \omega(\widehat{x}_s) = \delta_s\big) P\big(V_{s+1} = \delta_{s+1}, \dots, V_{|D|} = \delta_{|D|}\big) . \end{aligned} \tag{89.48}$$

As the sites $\widehat{x}_1, \dots, \widehat{x}_s$ are distinct, we have

$$P\big(\omega(\widehat{x}_1) = \delta_1, \dots, \omega(\widehat{x}_s) = \delta_s\big) = \prod_{1 \le i \le s} p^{\delta_i} (1-p)^{1-\delta_i} . \tag{89.49}$$

Recall that the sequence $(V_i)_{1 \le i \le |D|}$ is i.i.d. Bernoulli with parameter $p$. The equalities (89.48) and (89.49) yield

$$P\big(U_i = \delta_i,\, 1 \le i \le |D|\big) = \prod_{1 \le i \le |D|} p^{\delta_i} (1-p)^{1-\delta_i} . \tag{89.50}$$

A final remark is needed to conclude the proof. In fact, we started the previous computation with a sequence of values $(\delta_i, 1 \le i \le |D|)$ such that

$$P\big(U_i = \delta_i,\, 1 \le i \le |D|\big) > 0 . \tag{89.51}$$

Thus we have proved the formula (89.50) for any sequence $(\delta_i, 1 \le i \le |D|)$ satisfying (89.51). As the sum of both sides of (89.50) over all the sequences $(\delta_i, 1 \le i \le |D|)$, without constraint, is equal to 1, we conclude that (89.50) holds in fact for all the sequences $(\delta_i, 1 \le i \le |D|)$. This proves that the sequence $(U_i)_{1 \le i \le |D|}$ is i.i.d. and it completes the proof. □

## 89.3 Composing genuine explorations

Now that we have defined precisely the relevant concepts, we state in the next proposition an important result on the stability of genuine algorithms.

**Proposition 89.11.** *Let $\Phi = (\phi_k, 1 \le k \le \ell)$ and $\Psi = (\psi_k, 1 \le k \le m)$ be two genuine exploration algorithms in $D$. Let $x$ be a site of $D$ and let $X$ be a random site in $D$ which is measurable with respect to the outcome of $\Phi$ starting from $x$, i.e., there exists a deterministic map $f$ defined on the collection of the subsets of $D$ with values in $D$ such that $X = f\big(\Phi(x)\big)$. The random set $\Phi(x) \cup \Psi(X)$ is the outcome of a genuine exploration in $D$ starting from $x$.*

*Proof.* Let $\Phi, \Psi, x, X, f$ be as in the hypothesis of the proposition. We will define decision rules $\Gamma = (\gamma_k, 1 \le k \le \ell + m)$ so that the outcome of the genuine exploration $\Gamma$ in $D$ starting from $x$ is precisely $\Phi(x) \cup \Psi(X)$. The principle of the construction is the following:
• as long as the exploration $\Phi$ is not finished, we use the decision rules of $\Phi$;
• once the exploration $\Phi$ has terminated, we continue with $\Psi$, starting from $X$.
Unfortunately, the corresponding formalism is a bit heavy. For $1 \le k \le \ell$, we define a termination time map $T_k$ associated to $\Phi$, as the map from $D^k \times \{\,0, 1\,\}^k$ to $\{\,0, \dots, k\,\}$ given by

$$\begin{aligned}
&\forall x_1, \dots, x_k \in D \qquad \forall \varepsilon_1, \dots, \varepsilon_k \in \{0, 1\}\\
&T_k\big(x_1, \dots, x_k, \varepsilon_1, \dots, \varepsilon_k\big) \,=\\
&\qquad\begin{cases} 0 \quad \text{if} \quad \forall j \in \{\,1, \dots, k\,\} \quad \phi_j\big(x_1, \dots, x_j, \varepsilon_1, \dots, \varepsilon_j\big) \in D\,,\\ \min\Big\{\, j \in \{\,1, \dots, k\,\} : \phi_j\big(x_1, \dots, x_j, \varepsilon_1, \dots, \varepsilon_j\big) = \dagger \,\Big\} \quad \text{otherwise}\,. \end{cases}
\end{aligned} \tag{89.52}$$

The meaning of this definition is the following. If $T_k\big(x_1, \dots, x_k, \varepsilon_1, \dots, \varepsilon_k\big)$ is equal to 0, then the exploration $\Phi$ starting at $x_1$ will not have terminated after having seen the data $\big(x_1, \dots, x_k, \varepsilon_1, \dots, \varepsilon_k\big)$; if this quantity is equal to a positive integer $j$, then $j$ is the termination time of the same exploration processing the previous data. From step $j$ onwards, the exploration will restart from the site $f\big(\{\,x_1, \dots, x_j\,\}\big)$, and it will be driven by the decision rules of $\Psi$. We define next the decision rules of the exploration $\Gamma$. Let $k$ be such that $1 \le k \le \ell + m$, let $x_1, \dots, x_k$ belong to $D$ and let $\varepsilon_1, \dots, \varepsilon_k$ belong to $\{0, 1\}$. We set $j = T_k\big(x_1, \dots, x_k, \varepsilon_1, \dots, \varepsilon_k\big)$ and we consider the following four cases:
• $j = 0$. In this case, the exploration $\Phi$ has not terminated, so we define

$$\gamma_k\big(x_1, \dots, x_k, \varepsilon_1, \dots, \varepsilon_k\big) \,=\, \phi_k\big(x_1, \dots, x_k, \varepsilon_1, \dots, \varepsilon_k\big)\,.$$

• $j = k$. In this case, the exploration $\Phi$ has just terminated, so we define

$$\gamma_k\big(x_1, \dots, x_k, \varepsilon_1, \dots, \varepsilon_k\big) \,=\, f\big(\{\,x_1, \dots, x_k\,\}\big)\,.$$

• $1 \le j < k$ and $\psi_{k-j}\big(x_{j+1}, \dots, x_k, \varepsilon_{j+1}, \dots, \varepsilon_k\big) = \dagger$. In this case, both explorations $\Phi$ and $\Psi$ have terminated, so we define

$$\gamma_k\big(x_1, \dots, x_k, \varepsilon_1, \dots, \varepsilon_k\big) \,=\, \dagger\,.$$

• $1 \le j < k$ and $\psi_{k-j}\big(x_{j+1},\dots,x_k,\varepsilon_{j+1},\dots,\varepsilon_k\big) \neq \dagger$. In this case, the exploration $\Phi$ has terminated, but not $\Psi$. We define

$$\gamma_k\big(x_1,\dots,x_k,\varepsilon_1,\dots,\varepsilon_k\big) \,=\, \psi_{k-j}\big(x_{j+1},\dots,x_k,\varepsilon_{j+1},\dots,\varepsilon_k\big)\,.$$

The condition (89.5) ensures that $\phi_\ell(\cdot)$ is identically equal to $\dagger$, therefore the integer $j$ in the above construction is less than or equal to $\ell$. Furthermore, the map $\psi_m$ returns only the cemetery $\dagger$, therefore the last map $\gamma_{\ell+m}$ is constant equal to $\dagger$. Thus the decision rules $\Gamma = (\gamma_k, 1 \le k \le \ell + m)$ define a genuine exploration algorithm in $D$. Let us now compare the sequence of sites visited by the three explorations $\Phi(x)$, $\Psi(X)$ and $\Gamma(x)$. To that end, we introduce a specific notation for each sequence:

• $X_1,\dots,X_T$ is the sequence produced by $\Phi$ starting from $x$;

• $Y_1,\dots,Y_U$ is the sequence produced by $\Psi$ starting from $X = f\big(\Phi(x)\big)$;

• $Z_1,\dots,Z_V$ is the sequence produced by $\Gamma$ starting from $x$.

Naturally, these sequences depend on the percolation configuration and they are random, as well as their lengths $T, U, V$. The proof of the proposition will be completed by proving that $V = T + U$ and

$$(Z_1,\dots,Z_T) \,=\, (X_1,\dots,X_T)\,,\quad (Z_{T+1},\dots,Z_V) \,=\, (Y_1,\dots,Y_U)\,. \tag{89.53}$$

There are two parts. The first part deals with the sequence $X_1,\dots,X_T$, it is the object of lemma 89.12. The second part deals with the sequence $Y_1,\dots,Y_U$, it is the object of lemma 89.13. At first glance, these lemmas look boring and uninteresting, for two reasons. The first one is that the results they state seem to be immediate consequences of the definition of $\Gamma$. The second one is that the proofs are indeed mechanical verifications of the results. Our first intention was to omit them. However, we tried initially to work with explorations satisfying the condition (89.4). When writing the proofs, it appeared that some complicated compatibility condition was needed in order for the analogous results to hold. This matter is discussed further in subsection 89.4. This is why we have chosen to remove completely the condition (89.4) from the definition of genuine explorations. In the end, for the sake of completeness, we decided to include the proofs, tedious as they may be.

**Lemma 89.12.** *We have $T < V$ and $X_k = Z_k$ for $1 \le k \le T$.*

*Proof.* We start by proving that

$$\forall k \in \{\,1,\dots,\min(T,V)\,\}\qquad X_k = Z_k\,. \tag{89.54}$$

Due to the iterative constructions of the sequences, it seems that, once more, the lemma has to be proved by induction. As this induction contains no difficulty, let us do it quickly. We have $X_1 = Z_1 = x$. Let $k$ be such that $1 \le k < \min(T,V)$, and suppose that we have proved that

$$\forall i \in \{\,1,\dots,k\,\}\qquad X_i = Z_i\,.$$

We apply the construction rule (89.6) for $Z_{k+1}$. The definition of $Z_{k+1}$ depends on

$$\gamma_k\big(Z_1,\dots,Z_k,\omega(Z_1),\dots,\omega(Z_k)\big)\,,$$

which himself depends on

$$T_k\big(Z_1,\dots,Z_k,\omega(Z_1),\dots,\omega(Z_k)\big)\,.$$

So we start by computing this quantity. We assumed that $k<T$, meaning that

$$\forall i\in\{\,1,\dots,k\,\}\qquad \phi_i\big(X_1,\dots,X_i,\omega(X_1),\dots,\omega(X_i)\big)\neq\dagger\,. \tag{89.55}$$

Using the induction hypothesis, we rewrite (89.55) as

$$\forall i\in\{\,1,\dots,k\,\}\qquad \phi_i\big(Z_1,\dots,Z_i,\omega(Z_1),\dots,\omega(Z_i)\big)\neq\dagger\,. \tag{89.56}$$

Looking at the definition (89.52) of $T_k(\cdot)$, we see that (89.56) implies that

$$T_k\big(Z_1,\dots,Z_k,\omega(Z_1),\dots,\omega(Z_k)\big)\;=\;0\,.$$

Therefore we are in the first case in the definition of $\gamma_k$, where the exploration $\Phi$ has not terminated, and we conclude that

$$\begin{aligned}
Z_{k+1}\;&=\;\gamma_k\big(Z_1,\dots,Z_k,\omega(Z_1),\dots,\omega(Z_k)\big)\\
&=\;\phi_k\big(Z_1,\dots,Z_k,\omega(Z_1),\dots,\omega(Z_k)\big)\\
&=\;\phi_k\big(X_1,\dots,X_k,\omega(X_1),\dots,\omega(X_k)\big)\;=\;X_{k+1}\,.
\end{aligned}$$

This completes the induction step and the proof of (89.54). In the course of the previous proof, we have also shown that

$$\forall k\in\{\,1,\dots,\min(T,V)-1\,\}\qquad T_k\big(Z_1,\dots,Z_k,\omega(Z_1),\dots,\omega(Z_k)\big)\;=\;0\,. \tag{89.57}$$

By construction, the termination time of $\Gamma$ is strictly larger than the termination time of $\Phi$. The result (89.57) implies that $V\geq\min(T,V)+1$, whence $V>T$. ☐

**Lemma 89.13.** *We have $V=T+U$ and $Z_k=Y_{k-T}$ for $T+1\leq k\leq T+U$.*

*Proof.* We start by proving that

$$\forall k\in\{\,T+1,\dots,\min(T+U,V)\,\}\qquad Z_k=Y_{k-T}\,. \tag{89.58}$$

As always, the proof has to be done by induction. The very definition of $T$ and the previous lemma 89.12 together yield that

$$\begin{aligned}
\dagger\;&=\;\phi_T\big(X_1,\dots,X_T,\omega(X_1),\dots,\omega(X_T)\big)\\
&=\;\phi_T\big(Z_1,\dots,Z_T,\omega(Z_1),\dots,\omega(Z_T)\big)\,,
\end{aligned}$$

whence

$$T_T\big(Z_1,\dots,Z_T,\omega(Z_1),\dots,\omega(Z_T)\big)\;=\;T\,. \tag{89.59}$$

The awkward-looking notation $T_T$ means that we substitute $k$ by $T$ in the definition (89.52) of $T_k$. Let now $k$ be an index such that $k \geq T+1$. It follows from (89.59) that

$$T_k\big(Z_1, \dots, Z_k, \omega(Z_1), \dots, \omega(Z_k)\big) \;=\; T\,. \tag{89.60}$$

Let us compute $Z_{k+1}$. The identity (89.60) indicates that we are in the third or fourth case for the definition of $\gamma_k$, where the exploration $\Phi$ has terminated. We use then the exploration $\Psi$ starting from $X = f(\Phi(x))$. From the construction of the sequence $Z_1, \dots, Z_V$, we have $Z_{T+1} = Y_1 = X$. This is the desired result for $k = T+1$. Next, we carry out the induction step. Let $k$ be such that $T < k < \min(T+U, V)$, and suppose that the result has been proved until rank $k$. Since $k < V$ and $k - T < U$, then

$$\begin{aligned}\psi_{k-T}\big(Z_{T+1}, \dots, Z_k, \omega(Z_{T+1}), \dots, \omega(Z_k)\big) \;&\neq\; \dagger\,,\\ \psi_{k-T}\big(Y_1, \dots, Y_{k-T}, \omega(Y_1), \dots, \omega(Y_{k-T})\big) \;&\neq\; \dagger\,.\end{aligned}$$

Using the induction hypothesis, we have

$$\begin{aligned}Z_{k+1} \;&=\; \psi_{k-T}\big(Z_{T+1}, \dots, Z_k, \omega(Z_{T+1}), \dots, \omega(Z_k)\big)\\ &=\; \psi_{k-T}\big(Y_1, \dots, Y_{k-T}, \omega(Y_1), \dots, \omega(Y_{k-T})\big) \;=\; Y_{k+1-T}\,.\end{aligned}$$

This proves the result at rank $k+1$ and the induction step is completed. Thus we have proved (89.58). It remains to prove that $V = T + U$. We know from lemma 89.12 that $V > T$. Suppose that $V < T + U$. Using (89.58) for $T + 1 \leq k \leq V$, we have

$$\begin{aligned}&\psi_{V-T}\big(Z_{T+1}, \dots, Z_V, \omega(Z_{T+1}), \dots, \omega(Z_V)\big)\\ &\qquad\qquad =\; \psi_{V-T}\big(Y_1, \dots, Y_{V-T}, \omega(Y_1), \dots, \omega(Y_{V-T})\big)\,.\end{aligned} \tag{89.61}$$

As $V - T < U$, it follows from the definition of $U$ that the quantity of the right-hand side of (89.61) is not equal to $\dagger$. This implies that

$$\gamma_V\big(Z_1, \dots, Z_V, \omega(Z_1), \dots, \omega(Z_V)\big)$$

is not equal to $\dagger$, but this stands in contradiction with the definition of $V$. Suppose that $V > T + U$. Using again (89.58) for $T + 1 \leq k \leq T + U$, we have

$$\begin{aligned}&\psi_U\big(Z_{T+1}, \dots, Z_{T+U}, \omega(Z_{T+1}), \dots, \omega(Z_{T+U})\big)\\ &\qquad\qquad =\; \psi_U\big(Y_1, \dots, Y_U, \omega(Y_1), \dots, \omega(Y_U)\big)\,.\end{aligned} \tag{89.62}$$

It follows from the definition of $U$ that the quantity of the right-hand side of (89.62) is equal to $\dagger$. This implies that

$$\gamma_{T+U}\big(Z_1, \dots, Z_{T+U}, \omega(Z_1), \dots, \omega(Z_{T+U})\big)$$

is also equal to $\dagger$, but this contradicts the fact that $V > T + U$. The only possible case is that $V = T + U$. $\square$

Lemmas 89.12 and 89.13 yield the two equalities (89.53), which in turn gives that $\Gamma(x) = \Phi(x) \cup \Psi(X)$, and we are done. $\square$

## 89.4 Self-avoiding explorations

A sequence of decision rules $(\phi_k, 1 \leq k \leq |D|)$ is said to be self-avoiding if it satisfies the following condition: for any $k$ in $\{ 1, \dots, |D| \}$,

$$\forall x_1, \dots, x_k \text{ pairwise distinct elements of } D \quad \forall \varepsilon_1, \dots, \varepsilon_k \in \{0,1\}$$
$$\phi_k\big(x_1, \dots, x_k, \varepsilon_1, \dots, \varepsilon_k\big) \notin \{ x_1, \dots, x_k \} . \quad (89.63)$$

Until the subsection 89.2, we had been working with self-avoiding explorations. There are some advantages in imposing the constraint (89.63). First, it goes well with the idea of an exploration mechanism, where at each step a new site will be examined. Second, it automatically ensures that the termination time is at most equal to $|D|$. Third, self-avoiding explorations go faster and are more efficient at exploring than explorations that might visit several times the sites. So, why did we relax this constraint in the formal definition 89.1 of a genuine exploration?

In fact, all the results of subsection 89.2 hold for self-avoiding explorations as well. In addition, several arguments are greatly simplified when dealing with self-avoiding explorations. For example, the cleansing operation is no longer necessary, as the sequence of sites produced by a self-avoiding exploration is automatically constituted of pairwise distinct sites. In fact, the subsections 89.2 and 89.3 were initially written for self-avoiding explorations, but serious trouble appeared for the last lemma 89.13!

To be honest, some additional complications started to arise earlier in the proof of proposition 89.11. Indeed, if one defines a genuine set to be a set which is the outcome of a self-avoiding genuine exploration, then the exploration $\Gamma$ constructed by combining the explorations $\Phi$ and $\Psi$ should also be self-avoiding. Let us consider the framework of the proof of proposition 89.11. Here is a natural way to define $\Gamma$ so that it is self-avoiding.

**Alternative definition of** $\Gamma$. Let $k$ be such that $1 \leq k \leq |D|$, let $x_1, \dots, x_k$ be pairwise distinct elements of $D$ and let $\varepsilon_1, \dots, \varepsilon_k$ belong to $\{0, 1\}$. By the way, it should be required that the map $f$ introduced in the statement of proposition 89.11 returns a site not belonging to its argument:

$$\forall A \subset D \qquad f(A) \notin A .$$

The above condition is problematic for $A = D$, so we should also allow that $f$ takes the value †. We set $j = T_k\big(x_1, \dots, x_k, \varepsilon_1, \dots, \varepsilon_k\big)$ and we consider the following four cases:

• $j = 0$. In this case, the exploration $\Phi$ has not terminated, so we define

$$\gamma_k\big(x_1, \dots, x_k, \varepsilon_1, \dots, \varepsilon_k\big) \;=\; \phi_k\big(x_1, \dots, x_k, \varepsilon_1, \dots, \varepsilon_k\big) .$$

• $j = k$. In this case, the exploration $\Phi$ has just terminated, so we define

$$\gamma_k\big(x_1, \dots, x_k, \varepsilon_1, \dots, \varepsilon_k\big) \;=\; f\big(\{ x_1, \dots, x_k \}\big) .$$

• $1 \le j < k$ and $\psi_{k-j}\big(x_{j+1},\dots,x_k,\varepsilon_{j+1},\dots,\varepsilon_k\big) = \dagger$. In this case, the exploration $\Phi$ and $\Psi$ have both terminated, so we define

$$\gamma_k\big(x_1,\dots,x_k,\varepsilon_1,\dots,\varepsilon_k\big) \,=\, \dagger\,.$$

• $1 \le j < k$ and $\psi_{k-j}\big(x_{j+1},\dots,x_k,\varepsilon_{j+1},\dots,\varepsilon_k\big) \neq \dagger$. In this case, the exploration $\Phi$ has terminated, but not $\Psi$. We cannot simply define

$$\gamma_k\big(x_1,\dots,x_k,\varepsilon_1,\dots,\varepsilon_k\big) \,=\, \psi_{k-j}\big(x_{j+1},\dots,x_k,\varepsilon_{j+1},\dots,\varepsilon_k\big)\,,$$

because it could happen that $\psi_{k-j}\big(x_{j+1},\dots,x_k,\varepsilon_{j+1},\dots,\varepsilon_k\big)$ is equal to one of the sites $x_1,\dots,x_j$, and the resulting decision rule would not satisfy the constraint (89.63). To avoid this problem, we iterate forward the map $\psi$ until we find a site that is not already present in the sequence. So we define a sequence of sites $y_1,y_2,\dots,y_\ell$ and a sequence of states $\delta_1,\dots,\delta_{\ell-1}$ in $\{0,1\}$ through the following iterative procedure:

• **Initialization:** We set $y_1 = f\big(\{\,x_1,\dots,x_j\,\}\big)$.

• **Termination:** Let $\ell \ge 1$ and suppose that the sites $y_1,\dots,y_\ell$ have been defined. If $y_\ell$ does not belong to the set $\{\,x_1,\dots,x_j\,\}$, the construction stops at step $\ell$ and we set $\ell^* = \ell$.

• **Step $\ell$:** Let $\ell \ge 1$, suppose that the sites $y_1,\dots,y_\ell$ have been defined, as well as the states $\delta_1,\dots,\delta_{\ell-1}$, and that the termination criterion is not met. This means that the site $y_\ell$ is already present in the set $\{\,x_1,\dots,x_j\,\}$. Let $i$ be the unique index such that $y_\ell = x_i$. We define $\delta_\ell = \varepsilon_i$ and

$$y_{\ell+1} \,=\, \psi_{k+\ell-j}\big(x_{j+1},\dots,x_k,y_1,\dots,y_\ell,\varepsilon_{j+1},\dots,\varepsilon_k,\delta_1,\dots,\delta_\ell\big)\,. \qquad (89.64)$$

We claim that this procedure terminates. Indeed, suppose that it reaches the step $\ell \ge 1$ and that it produces the sequence $y_1,\dots,y_\ell$. As the exploration $\Psi$ satisfies (89.63), in view of (89.64), the sequence $y_1,\dots,y_\ell$ is such that

$$\forall h \in \{\,1,\dots,\ell\,\} \qquad y_h \notin \big\{\,x_{j+1},\dots,x_k,y_1,\dots,y_{h-1}\,\big\}\,. \qquad (89.65)$$

If the termination criterion is not met until step $\ell$ included, then we have also

$$\forall h \in \{\,1,\dots,\ell\,\} \qquad y_h \in \big\{\,x_1,\dots,x_j\,\big\}\,. \qquad (89.66)$$

The condition (89.65) implies that the sites $y_1,\dots,y_\ell$ are pairwise distinct, hence by the condition (89.66) we must also have $\ell \le j$. Therefore the procedure terminates at some step $\ell^* \le j+1$. We finally define

$$\gamma_k\big(x_1,\dots,x_k,\varepsilon_1,\dots,\varepsilon_k\big) \,=\, y_{\ell^*}\,. \qquad (89.67)$$

The termination criterion and the condition (89.65) guarantee that $y_{\ell^*}$ is neither in $\big\{\,x_1,\dots,x_j\,\big\}$ nor in $\big\{\,x_{j+1},\dots,x_k\,\big\}$. Notice that $y_{\ell^*}$ can be equal to $\dagger$. In any case, the formula (89.67) gives a value which is not in $\big\{\,x_1,\dots,x_k\,\big\}$, the map $\gamma_k$ satisfies (89.63), and the decision rules $\Gamma = (\gamma_k, 1 \le k \le |D|)$ define a genuine exploration algorithm in $D$, which is self-avoiding. So far, everything is fine. Lemma 89.12 is still valid and the same proof applies. The trouble is with lemma 89.13, which is no longer valid. The statement of lemma 89.13 should be replaced by a weaker one, for instance we can prove that

$$\{ Z_{T+1},\dots,Z_V \} = \{ Y_1,\dots,Y_U \} \setminus \{ X_1,\dots,X_T \} . \tag{89.68}$$

The equality (89.68) would be enough to conclude that the outcome of the exploration $\Gamma$ is $\Phi(x)\cup\Psi(X)$. However, there are little troubles with a statement like (89.68). First, the proof of (89.68) is rather cumbersome. Indeed, the strategy would be the same as for lemma 89.13, based on an induction and the iterative definition of $\Gamma$. Yet, as the definition of $\Gamma$ is much more involved, so is the proof. We would have to go through the same steps, hence inside the induction step, we would need to iterate the exploration $\Psi$ until we find a site not belonging to $\{ X_1,\dots,X_T \}$. Second, a statement like (89.68) induces a change of perspective, compared to lemma 89.12. It indicates that it is more appropriate to work with random sets rather than random sequences of sites. In the same vein, the exploration $\Gamma$ built above enjoys strong properties of symmetry: in the case where $T_k\big(x_1,\dots,x_k,\varepsilon_1,\dots,\varepsilon_k\big)=j\geq 1$, the rule depends on the map $\Psi$ and the set $\{ x_1,\dots,x_j \}$ (and not on the order of the sites in the sequence $x_1,\dots,x_j$). This would prompt us to introduce a stronger reasonable symmetry property on the decision rules. Denoting by $\mathfrak{S}_k$ the group of permutations of $\{ 1,\dots,k \}$, we would work with explorations $\Phi=(\phi_k, 1\leq k\leq |D|)$ satisfying

$$\begin{gathered}\forall k\in\{ 1,\dots,|D| \}\quad \forall x_1,\dots,x_k\in D\quad \forall\varepsilon_1,\dots,\varepsilon_k\in\{0,1\}\quad \forall\sigma\in\mathfrak{S}_k\\ \phi_k\big(x_1,\dots,x_k,\varepsilon_1,\dots,\varepsilon_k\big) = \phi_k\big(x_{\sigma(1)},\dots,x_{\sigma(k)},\varepsilon_{\sigma(1)},\dots,\varepsilon_{\sigma(k)}\big) .\end{gathered} \tag{89.69}$$

However, the conditions (89.69) mean essentially that the decision rules should be a function of the set of the visited sites, independently of the order of their appearance. This indicates that the decision rule $\phi_k$ should take as input the set of the visited sites and their states. This subcategory of explorations, which was already mentioned in section 13 when we introduced genuine explorations, seems to be the best suited for our current needs. This is probably the one that should be used to shorten the proofs and improve the effectiveness of the presentation. For instance, an equality like (89.68) would be easier to formulate and prove with decision rules involving sets rather than sequences. But we must be careful, because we risk losing out in another way. Indeed, all the probabilistic estimates are ultimately based on sequences of random variables, and so it seems inevitable that we will have to reorganize at some point the set of the sites discovered during the explorations into an adequate sequence. Furthermore, the more general scheme of decision rules based on sequences might be useful for other purposes, there are some delicate problems for which it is essential to take into account the order in which sites are visited when choosing the next one. In the end, we will stick to our current framework and we will continue to work with sequences of sites, rather than sets of sites. In addition, throughout the sequel, we do not require any more that the genuine explorations are self-avoiding.

## 89.5 The truncated closed balls

We fix an integer $M \geq 1$. In this subsection, we define and study the truncated closed balls of radius $t$ associated to the truncated explorations with parameter $M$, which are the counterpart of the balls $\mathcal{B}(A, D, t)$ defined in (12.7). The difference is that, for exploring the successive shells around the initial set, we employ the truncated exploration instead of the standard cluster exploration. For $A$ a subset of $\mathbb{Z}^d$, we define $\widehat{\mathcal{B}}_M(A, D, 0) = A \cap D$ and

$$\forall t \geq 0 \qquad \widehat{\mathcal{B}}_M(A, D, t+1) \,=\, \overline{\text{Clusters}}\Big(\widehat{\mathcal{B}}_M(A, D, t) \cup \mathcal{N}\big(\widehat{\mathcal{B}}_M(A, D, t)\big), M, D\Big)\,. \tag{89.70}$$

This definition ensures that $\widehat{\mathcal{B}}_M(A, D, t)$ is increasing in $A$ and $D$ with respect to the set inclusion: for any $t \geq 0$ and any finite subsets $A, B, D, E$ of $\mathbb{Z}^d$, we have

$$A \subset B,\ D \subset E \quad\Longrightarrow\quad \widehat{\mathcal{B}}_M(A, D, t) \subset \widehat{\mathcal{B}}_M(B, E, t)\,. \tag{89.71}$$

More interestingly, for any $t \geq 0$, the map $A \mapsto \widehat{\mathcal{B}}_M(A, D, t)$ is compatible with the union, as stated in the next lemma.

**Lemma 89.14.** *For any $t \geq 0$, we have*

$$\forall A, B \subset \mathbb{Z}^d \qquad \widehat{\mathcal{B}}_M(A \cup B, D, t) \,=\, \widehat{\mathcal{B}}_M(A, D, t) \cup \widehat{\mathcal{B}}_M(B, D, t)\,.$$

*Proof.* Let $t \geq 0$ and let $A, B$ be two subsets of $\mathbb{Z}^d$. It follows immediately from (89.71) that

$$\widehat{\mathcal{B}}_M(A, D, t) \cup \widehat{\mathcal{B}}_M(B, D, t) \,\subset\, \widehat{\mathcal{B}}_M(A \cup B, D, t)\,.$$

Let us prove the converse inclusion by induction on $t$. For $t = 0$, we have

$$\begin{aligned}\widehat{\mathcal{B}}_M(A \cup B, D, 0) \,&=\, (A \cup B) \cap D \\ &=\, (A \cap D) \cup (B \cap D) \,=\, \widehat{\mathcal{B}}_M(A, D, 0) \cup \widehat{\mathcal{B}}_M(B, D, 0)\,.\end{aligned}$$

Let next $t \geq 0$ and suppose that the result has been proved at rank $t$. We use the definition (89.70) to write

$$\widehat{\mathcal{B}}_M(A\cup B, D, t+1) \,=\, \overline{\text{Clusters}}\Big(\widehat{\mathcal{B}}_M(A\cup B, D, t) \cup \mathcal{N}\big(\widehat{\mathcal{B}}_M(A\cup B, D, t)\big), M, D\Big)\,. \tag{89.72}$$

From the induction hypothesis, we have

$$\widehat{\mathcal{B}}_M(A \cup B, D, t) \,=\, \widehat{\mathcal{B}}_M(A, D, t) \cup \widehat{\mathcal{B}}_M(B, D, t)\,,$$

whence also

$$\begin{gathered}\widehat{\mathcal{B}}_M(A \cup B, D, t) \cup \mathcal{N}\big(\widehat{\mathcal{B}}_M(A \cup B, D, t)\big) \,\subset \\ \widehat{\mathcal{B}}_M(A, D, t) \cup \widehat{\mathcal{B}}_M(B, D, t) \cup \mathcal{N}\big(\widehat{\mathcal{B}}_M(A, D, t)\big) \cup \mathcal{N}\big(\widehat{\mathcal{B}}_M(B, D, t)\big)\,.\end{gathered} \tag{89.73}$$

It follows from the definition (89.2) that the map $A \mapsto \overline{\text{Clusters}}\,(A, M, D)$ is increasing with respect to the set inclusion and compatible with the union. Therefore the formulas (89.72) and (89.73) yield that

$$\begin{aligned}\widehat{\mathcal{B}}_M(A\cup B, D, t+1) \,\subset\, &\overline{\text{Clusters}}\Big(\widehat{\mathcal{B}}_M(A,D,t)\,\cup\,\mathcal{N}\big(\widehat{\mathcal{B}}_M(A,D,t)\big), M, D\Big)\\ &\cup\overline{\text{Clusters}}\Big(\widehat{\mathcal{B}}_M(B,D,t)\,\cup\,\mathcal{N}\big(\widehat{\mathcal{B}}_M(B,D,t)\big), M, D\Big)\\ &\qquad = \widehat{\mathcal{B}}_M(A,D,t+1)\cup\widehat{\mathcal{B}}_M(B,D,t+1)\,.\end{aligned}$$

This completes the induction step and the proof of the lemma. □

We have also the semigroup property:

$$\forall s\geq 0\quad \forall t\geq 0\qquad \widehat{\mathcal{B}}_M(A,D,s+t) \,=\, \widehat{\mathcal{B}}_M\big(\widehat{\mathcal{B}}_M(A,D,s),D,t\big)\,. \tag{89.74}$$

The truncation implies that these sets are smaller than the balls $\mathcal{B}(A, D, t)$:

$$\forall t\geq 0\qquad \widehat{\mathcal{B}}_M(A,D,t)\,\subset\,\mathcal{B}(A,D,t)\,.$$

An immediate consequence of the definitions is that the truncated balls can be compared with the deterministic $B_1$ balls, as shown by the next lemma.

**Lemma 89.15.** *For any subset $A$ of $D$, we have the inclusions*

$$\forall t\geq 0\qquad \widehat{\mathcal{B}}_M\big(A,D,t\big)\,\subset\, B_1\big(A\cap D, t(M+1)\big)\,.$$

*Proof.* We do the proof by induction on $t$. For $t = 0$, we have indeed

$$\widehat{\mathcal{B}}_M\big(A,D,0\big) \,=\, A\cap D \,=\, B_1\big(A\cap D,0\big)\,.$$

Let $t \geq 0$ and suppose that the result has been proved at rank $t$. By definition, we have

$$\widehat{\mathcal{B}}_M(A,D,t+1) \,=\, \overline{\text{Clusters}}\Big(\widehat{\mathcal{B}}_M(A,D,t)\,\cup\,\mathcal{N}\big(\widehat{\mathcal{B}}_M(A,D,t)\big), M, D\Big)\,. \tag{89.75}$$

The induction hypothesis yields that

$$\widehat{\mathcal{B}}_M\big(A,D,t\big)\,\subset\, B_1\big(A\cap D, t(M+1)\big)\,,$$

whence also

$$\mathcal{N}\big(\widehat{\mathcal{B}}_M(A,D,t)\big)\,\subset\, B_1\big(A\cap D, t(M+1)+1\big)\,.$$

We substitute the two previous inclusions in (89.75) and we obtain

$$\widehat{\mathcal{B}}_M(A,D,t+1)\,\subset\,\overline{\text{Clusters}}\Big(B_1\big(A\cap D, t(M+1)+1\big), M, D\Big)\,.$$

We use finally the property (89.3), and we conclude that

$$\widehat{\mathcal{B}}_M(A,D,t+1)\,\subset\, B_1\big(A\cap D, t(M+1)+1+M\big) \,=\, B_1\big(A\cap D, (t+1)(M+1)\big)\,.$$

This completes the induction step and the proof. □

We shall next prove that the truncated balls are genuine. The key fact is stated in the next lemma.

**Lemma 89.16.** *If $\mathcal{E}$ is the outcome of a genuine exploration algorithm, then the same is true for $\overline{\text{Clusters}}\big(\mathcal{E} \cup \mathcal{N}(\mathcal{E}), M, D\big)$ as well.*

*Proof.* The proof rests on repeated applications of proposition 89.11. We consider an arbitrary ordering $x_1, \dots, x_{|D|}$ of the sites of $D$. Let $A$ be a subset of $D$ and let $k$ belong to $\{1, \dots, |D|\}$. We define

$$I_k(A) = \max\Big\{\, i \in \{1, \dots, |D|\} : x_i \in A,\ \big|\{\, x_j : 1 \leq j < i,\ x_j \in A \,\}\big| \leq k-1 \,\Big\},$$

$$f_k(A) \,=\, x_{I_k(A)}\,.$$

The integer $I_k(A)$ is the index $i$ of the site in the ordered sequence $x_1, \dots, x_{|D|}$ such that $x_i$ is in $A$, and there are at most $k-1$ sites before $x_i$ in the sequence which are in $A$. The site $f_k(A)$ is simply the site whose index is $I_k(A)$. Let $\mathcal{E}$ be a random genuine set. For any site $x$ in $D$, the set $\overline{\text{Cluster}}\,(x, M, D)$ is the outcome of a genuine exploration algorithm. By proposition 89.11, the same is true for the sets

$$\mathcal{E} \cup \overline{\text{Clusters}}\,\Big(f_k\big(\mathcal{E} \cup \mathcal{N}(\mathcal{E})\big), M, D\Big), \quad 1 \leq k \leq |D|\,.$$

Using the previous definitions, we have the representation

$$\overline{\text{Clusters}}\,\big(\mathcal{E} \cup \mathcal{N}(\mathcal{E}), M, D\big) \,=\, \bigcup_{1 \leq k \leq |D|} \mathcal{E} \cup \overline{\text{Clusters}}\,\Big(f_k\big(\mathcal{E} \cup \mathcal{N}(\mathcal{E})\big), M, D\Big)\,. \tag{89.76}$$

A simple induction, or repeated applications of proposition 89.11 shows that the set (89.76) is the outcome of a genuine exploration algorithm. □

Once we have lemma 89.16, an immediate induction shows that the set (89.78) is genuine. We record this result in the next corollary.

**Corollary 89.17.** *For any subset $W$ of $D$, any $t \geq 0$, the random set*

$$\widehat{\mathcal{B}}_M\big(\mathit{Shell}\,(W, D), D, t\big)$$

*is genuine.*

*Proof.* For $t = 0$, the set $\widehat{\mathcal{B}}_M\big(\text{Shell}\,(W, D), D, 0\big)$ coincides with $\text{Shell}\,\big(W, D\big)$, and it is one of our first examples of a genuine set. It was defined in (12.5) and it is the outcome of the standard exploration algorithm (see algorithm 5.1), initialized with the set $W$ for the set of the active sites. For $t \geq 0$, we use the definition (89.70) to write

$$\widehat{\mathcal{B}}_M\big(\text{Shell}\,(W, D), D, t+1\big) \,=\,$$
$$\overline{\text{Clusters}}\,\bigg(\widehat{\mathcal{B}}_M\big(\text{Shell}\,\big(W, D\big), D, t\big) \,\cup\, \mathcal{N}\Big(\widehat{\mathcal{B}}_M\big(\text{Shell}\,(W, D), D, t\big)\Big), M, D\bigg)\,.$$

We conclude by a simple induction, or repeated applications of lemma 89.16. □

## 89.6 Genuine truncated approximations

For the sake of the discussion, we replace the four sets $K^-, L^-, K^+, L^+$ appearing in the inequalities (88.34), (88.36) by a generic set. Let us fix a subset $I$ of $\{1, \dots, N\}$. Our goal here is to approximate the sets of the form

$$\bigcup_{i\in I} \mathcal{B}\big(w_i, D\setminus \overline{\mathcal{T}}_i(t-1), t-1\big)\,, \quad \bigcup_{i\in I} \mathcal{B}\big(w_i, D\setminus \widehat{\mathcal{T}}_i^-(t), t\big)\,, \quad \bigcup_{i\in I} \mathcal{B}\big(w_i, D\setminus \widehat{\mathcal{T}}_i^*(t), t\big) \tag{89.77}$$

by sets which are the outcome of genuine explorations. The difficulty arises from the presence of the taboo sets and the random nature of these sets. Because they are random, they convey information about the status of sites outside of them, thereby creating a nasty disturbance in the distribution of the states of the sites belonging to the set visited by the taboo explorations. It is a recurring issue that we keep encountering. In part V, we tried to perform a conditioning on the taboo sets, but the ensuing estimates were too weak. Here, we are going to get round this problem with the help of the hypothesis (9.16) and the truncated explorations. Our first try would be to simply remove the taboo sets in the expressions (89.77), i.e., to compare them with the two sets

$$\bigcup_{i\in I} \mathcal{B}\big(w_i, D, t-1\big)\,, \quad \bigcup_{i\in I} \mathcal{B}\big(w_i, D, t\big)\,.$$

However, this is doomed to fail. If a taboo set was preventing the exploration to enter into a large cluster, then this large cluster will create a huge discrepancy between the two sets. In addition, we cannot restrict ourselves to the sole exploration of small clusters, our intention is precisely to start the intertwined explorations from vertices which belong to distinct large clusters of the configuration. Thus we need a mechanism preventing that the non-taboo explorations differ too much from the taboo explorations of the sets (89.77). This is really a delicate matter, because it is vital that the mechanism does not destroy the genuine behavior of the explorations. For instance, that was the problem when we tried to condition on the taboo sets: if we provide any information whatsoever on the size of the clusters of the taboo sites, or on their connections to some other subset, the distribution of the states of the vertices of these clusters is not any more a Bernoulli product measure, and the probabilistic estimates are ruined. Even worse, the same applies to non-taboo sites, because being non-taboo means not having been visited by the other explorers. We adopt here a radically different strategy, namely we implement a decision rule which is measurable with respect to the information gathered by the genuine exploration. In the probabilistic language, we introduce a stopping time associated to the exploration process, and we stop the exploration when the explored cluster is too large. This is where the truncated explorations come into play.

Our second approximation try is to remove the taboo sets in the expressions (89.77), to use the standard exploration scheme when $t = 0$, but to use truncated balls instead of balls for the subsequent times. Following the previous ideas, the last set in the expressions (89.77) will be compared to the set

$$\bigcup_{i\in I}\widehat{\mathcal{B}}_M\big(\text{Shell}\,(w_i),D,t\big)\;=\;\widehat{\mathcal{B}}_M\Big(\text{Shell}\,\big(\{\,w_i:i\in I\,\},D\big),D,t\Big)\,. \tag{89.78}$$

The above two sets are equal thanks to lemma 89.14. There is an additional subtlety, which is essential for this strategy to succeed: we will condition on the event that the finite clusters have cardinality smaller than $M$. On this event, all the clusters of the percolation configuration restricted to $D$ which do not meet the boundary $\partial^{in}D$ will be included into finite clusters of $D$, and will also have cardinality strictly smaller than $M$. In particular, from the second iteration onwards, whenever an explorer finds himself visiting a cluster with cardinality larger than or equal to $M$, it means that he is on the wrong track because he went through a site that was taboo, hence we can safely truncate his exploration.

Upon termination of the intertwined explorations, all the intersections occurring inside the set $G$ will have been localized. Using the results of subsection 88.7, we can recover the relevant taboo sites of the intertwined explorations until time $T$ as the intersections of adequately chosen taboo explorations. To detect the first set of closed relevant taboo sites $\mathcal{T}^{*0}(0)$, we can use explorations having an empty taboo set (this was done in subsection 88.1). As these explorations are genuine, we can use directly the functional $S$ and Hoeffding's inequality to obtain a quantitative control on $\big|\mathcal{T}^{*0}(0)\big|$. It is much more difficult to control the relevant taboo sets $(\mathcal{T}^{*t}(0),\ t\geq 1)$ occurring during the subsequent iterations. Indeed, for $t\geq 1$, the relevant taboo sites $\mathcal{T}^{*}(t-1)$ for the first $t-1$ iterations play an essential role to localize the set $\mathcal{T}^{*0}(t)$, and we cannot get rid of $\mathcal{T}^{*}(t-1)$ in order to localize precisely the relevant taboo sites for the iteration $t$. What is worse, as the set $\mathcal{T}^{*}(t-1)$ is random, we have no control on the distribution of the states of the sites visited by a taboo exploration that uses $\mathcal{T}^{*}(t-1)$ as taboo set, or a subset of it. In particular, we do not have any probabilistic control for the functional $S$ computed on the set visited by this taboo exploration.

Our strategy to get around the previous problem consists in approximating the taboo explorations by genuine truncated explorations. The idea is that, after the first iteration and the exploration of the clusters of $w_1,\dots,w_N$, the intersections can be localized by exploring only the finite clusters. When the event $\mathcal{E}$ occurs, the finite clusters have cardinality smaller than $M$. Therefore, from the second iteration onwards, whenever an explorer finds himself visiting a cluster of size $M$ or larger, he should stop the exploration. This is a nice simple rule, which has the big advantage of not destroying the genuine character of the exploration. For a technical reason, we will in fact confine the exploration to the sites that are at distance less than or equal to $M$ from the starting site. To sum up, instead of performing taboo explorations, after the first iteration, we will employ explorers who perform a genuine truncated exploration. Of course, the sets visited by the taboo explorations and the genuine truncated explorations differ, and our next task is to develop a control of these differences.

At first sight, it seems that a genuine truncated exploration always explores a larger set of sites than any taboo exploration starting from the same initial set. This would certainly be the case if the unexplored region contained only clusters of size smaller than $M$: in this case, the genuine truncated exploration becomes a genuine exploration, and it will explore systematically the whole shells surrounding the starting set, contrary to the taboo exploration, which is submitted to some restrictions. However, things are more complicated when the unexplored region contains also clusters of diameter larger than $M$: in this case, it might happen that the taboo exploration is allowed to visit completely such a cluster, while the genuine truncated exploration remains confined to a ball of diameter $M$. As we assume that the finite clusters have a size smaller than $M$, one way to avoid this last problem is to start the exploration from a set of sites containing exactly one vertex of each boundary cluster. This is what we did in subsection 87.

We will next develop approximation results between the three sets (89.77) and the sets (89.78). These approximations are essential for the proof of theorem 9.2. Let us sum up the idea behind the approximations. As we start the intertwined explorations with an initial set of vertices $W$ which contains exactly one vertex of each boundary cluster, all the boundary clusters are explored during the first iteration. During the subsequent iterations, the explorers can discover only finite clusters, and thanks to the occurrence of the event $\mathcal{H}$, all these finite clusters have cardinality smaller than $M$. Thus it is natural to compare the aggregate of the explorer $i$ with $\text{Shell}\,(w_i, D)$ for the first iteration, and with $\widehat{\mathcal{B}}_M\big(\text{Shell}\,(w_i, D), D, 1\big)$ for the second one, and so on so forth. Typically, these approximating sets will be larger than the original ones, the difficult part consists in controlling the sites that have been added with respect to the taboo explorations. The easiest and most natural approximation is for the first set in (89.77), which uses $\overline{\mathcal{T}_i}(t-1)$ as taboo set, and which corresponds exactly to the aggregate $\mathcal{A}_i(t-1)$ of the intertwined explorations. It is carried out in the subsubsection 89.6.1. The more difficult cases of the last two sets in (89.77), which use $\widehat{\mathcal{T}}_i^{\,-}(t)$, $\widehat{\mathcal{T}}_i^{\,*}(t)$ as taboo sets are carried out simultaneously in the subsubsections 89.6.3, 89.6.2. Naturally, the proofs of the three approximations rest on similar arguments. We do simultaneously the approximations for the last two sets $\widehat{\mathcal{T}}_i^{\,-}(t)$, $\widehat{\mathcal{T}}_i^{\,*}(t)$, and one might wonder whether we could do the three approximations simultaneously. In fact, we tried to unify them and to present a general approximation result that would imply the three of them at once. The drawback with such a presentation is that in the hypothesis of the general result, we would need to include a list of complicated conditions satisfied by the sequence of the taboo sets so that the approximation works. This list becomes quickly very obscure and we feel that this presentation would not facilitate the reading at all. In the end, we opt for the more lengthy but more readable presentation: first we do the easiest case of $\overline{\mathcal{T}_i}(t-1)$, and then the more subtle cases of $\widehat{\mathcal{T}}_i^{\,-}(t)$, $\widehat{\mathcal{T}}_i^{\,*}(t)$. By the way, we will need again this kind of results in section 92 in a more complicated setting.

### 89.6.1 The case of the aggregates

In this subsubsection, we carry out the approximation of $\mathcal{B}\big(w_i, D\setminus\overline{\mathcal{T}_i}(t-1), t-1\big)$. This set is obtained with an exploration using $\overline{\mathcal{T}_i}(t-1)$ as taboo set. The identity (81.2) in proposition 81.1 shows that the set $\mathcal{B}\big(w_i, D\setminus\overline{\mathcal{T}_i}(t-1), t-1\big)$ is precisely the aggregate $\mathcal{A}_i(t-1)$ of the explorer $i$. This is why the approximation is simpler for this set than for the others, and we present it separately. The approximations will involve error terms, and these terms are adequate enlargements of the sets of the relevant taboo sites, that we introduce now. For $t$ in $\{0,\dots,T\}$ and $s\geq 0$, we define

$$\begin{aligned} B^{t,0}_{1,M}(s) &= B_1\big(\mathcal{T}^{*0}(t), s(M+1)\big)\,,\\ B^{t,1}_{1,M}(s) &= B_1\big(\mathcal{T}^{*1}(t), s(M+1)+M-1\big)\,. \end{aligned} \tag{89.79}$$

The next proposition will be applied for the value $t-1$, but for its statement and its proof it is more convenient to use $t$ instead of $t-1$.

**Proposition 89.18.** *We consider a configuration realizing the event $\mathcal{H}$ defined in* (87.4) *and a subset $I$ of $\{1,\dots,N\}$. Setting $W(I)=\{w_i: i\in I\}$, we have, for any $t$ in $\{0,\dots,T\}$,*

$$\begin{gathered} \bigcup_{i\in I}\mathcal{B}\big(w_i, D\setminus\overline{\mathcal{T}_i}(t), t\big) \subset \widehat{\mathcal{B}}_M\Big(\mathit{Shell}\big(W(I),D\big),D,t\Big)\subset\\ \bigcup_{i\in I}\mathcal{B}\big(w_i, D\setminus\overline{\mathcal{T}_i}(t), t\big)\cup\bigcup_{0\leq s\leq t}\big(B^{s,0}_{1,M}(t-s)\cup B^{s,1}_{1,M}(t+1-s)\big)\,. \end{gathered} \tag{89.80}$$

*Proof.* We shall prove the proposition by induction on $t$ and we start with $t=0$. Let us fix $i$ in $I$. By definition, $\widehat{\mathcal{B}}_M\big(\mathrm{Shell}\,(w_i,D),D,0\big) = \mathrm{Shell}\,(w_i,D)$. On the event $\mathcal{H}$, we have $\overline{\mathcal{T}_i}(0)=\mathcal{T}^{*0}_i(0)$ and $\mathcal{T}^{*1}_i(0)=\varnothing$, thus

$$\mathcal{B}\big(w_i, D\setminus\overline{\mathcal{T}_i}(0), 0\big) = \mathrm{Shell}\,(w_i,D)\setminus\overline{\mathcal{T}_i}(0)\,, \tag{89.81}$$

$$\mathcal{B}\big(w_i, D\setminus\overline{\mathcal{T}_i}(0), 0\big)\cup B^{0,0}_{1,M}(0)\cup B^{0,1}_{1,M}(1) = \mathrm{Shell}\,(w_i,D)\cup\mathcal{T}^{*0}(0)\,. \tag{89.82}$$

Taking the union over $i$ in $I$ of the equality (89.81), we get

$$\begin{aligned} \bigcup_{i\in I}\mathcal{B}\big(w_i, D\setminus\overline{\mathcal{T}_i}(0), 0\big) &= \bigcup_{i\in I}\mathrm{Shell}\,(w_i,D)\setminus\overline{\mathcal{T}_i}(0)\subset\bigcup_{i\in I}\mathrm{Shell}\,(w_i,D)\\ &= \mathrm{Shell}\,\big(W(I),D\big) = \widehat{\mathcal{B}}_M\big(\mathrm{Shell}\,\big(W(I),D\big),D,0\big)\,. \end{aligned}$$

This proves the first inclusion in (89.80) for $t=0$. For the second inclusion, we take the union over $i$ in $I$ of the equality (89.82), and we get

$$\begin{aligned} \mathrm{Shell}\,\big(\{w_i: i\in I\},D\big) &\subset \bigcup_{i\in I}\mathrm{Shell}\,(w_i,D)\cup\mathcal{T}^{*0}(0)\\ &\subset\bigcup_{i\in I}\mathcal{B}\big(w_i, D\setminus\overline{\mathcal{T}_i}(0), 0\big)\cup B^{0,0}_{1,M}(0)\cup B^{0,1}_{1,M}(1)\,. \end{aligned}$$

Thus the inclusions (89.80) hold for $t = 0$.

We fix next $t \geq 1$ and we prove the first inclusion in (89.80) for $t$. Let us fix again $i$ in $I$. The set $\overline{\mathcal{T}_i}(t)$ prevents the explorer $i$ to explore the clusters of the vertices $w_1, \dots, w_N$, thus

$$\mathcal{B}\big(\text{Shell}\,(w_i, D), D \setminus \overline{\mathcal{T}_i}(t), t\big) \setminus \text{Shell}\,(w_i, D) \ \subset\ D \setminus \bigcup_{1 \leq j \leq N} C(w_j, D)\,.$$

The occurrence of the event $\mathcal{H}$ guarantees that all the clusters in the above set are finite clusters of $\Lambda(4n)$, and also that all these finite clusters have cardinality strictly less than $M$. This implies that

$$\begin{aligned}\mathcal{B}\big(w_i, D \setminus \overline{\mathcal{T}_i}(t), t\big) \ &=\ \mathcal{B}\big(\text{Shell}\,(w_i, D), D \setminus \overline{\mathcal{T}_i}(t), t\big)\\ &\qquad\subset\ \widehat{\mathcal{B}}_M\Big(\text{Shell}\,\big(W(I), D\big), D, t\Big)\,. \qquad (89.83)\end{aligned}$$

Taking the union over $i$ in $I$ in (89.83), we get

$$\bigcup_{i \in I} \mathcal{B}\big(w_i, D \setminus \overline{\mathcal{T}_i}(t), t\big) \ \subset\ \widehat{\mathcal{B}}_M\Big(\text{Shell}\,\big(W(I), D\big), D, t\Big)\,,$$

which is precisely the first inclusion in (89.80) for $t$.

The proof of the second inclusion in (89.80) is more delicate. We examine the case $t = 1$ before attacking the general case. Let $x$ be a site belonging to

$$\widehat{\mathcal{B}}_M\Big(\text{Shell}\,\big(W(I), D\big), D, 1\Big) \ \setminus\ \bigcup_{i \in I} \mathcal{B}\big(w_i, D \setminus \overline{\mathcal{T}_i}(1), 1\big)\,. \qquad (89.84)$$

By definition of $\widehat{\mathcal{B}}_M\big(\text{Shell}\,\big(W(I), D\big), D, 1\big)$, there exists a path $y_0, y_1, \cdots, y_r$ in $D$ joining a site $y_0$ of Shell $\big(W(I), D\big)$ to $y_r = x$ such that the sites $y_1, \cdots, y_{r-1}$ are in $D \setminus$Shell $\big(W(I), D\big)$, they are all open, and they are visited before $x$ by the truncated exploration starting from $y_1$. The truncation imposes that $y_2, \dots, y_r$ are in $B_1(y_1, M)$. Moreover, there exists a path $z_0, z_1, \cdots, z_s$ in $D$ joining a site $z_0$ of $W(I)$ to $z_s = y_0$ such that the sites $z_0, \cdots, z_{s-1}$ are in $D$ and they are all open. Let $i$ in $I$ be such that $z_0 = w_i$. Since $x$ is not in $\mathcal{B}\big(w_i, D \setminus \overline{\mathcal{T}_i}(1), 1\big)$, then the path $z_0, \cdots, z_{s-1}, y_0, \cdots, y_r$ must go through $\overline{\mathcal{T}_i}(1)$. Recall that

$$\overline{\mathcal{T}_i}(1) \ =\ \mathcal{T}_i^{*0}(0) \cup \mathcal{T}_i^{*0}(1) \cup \mathcal{T}_i^{*1}(1)\,.$$

The sites $z_0, \cdots, z_{s-1}$ are all open and connected to $w_i$, so they cannot belong to $\overline{\mathcal{T}_i}(1)$. Thus there exists $q$ such that $y_q$ is in $\overline{\mathcal{T}_i}(1)$. If $y_q$ is open, then it is in $\overline{\mathcal{T}_i}^1(1) = \mathcal{T}_i^{*1}(1)$, $q \geq 1$ and since both $y_q$ and $x$ are in $B_1(y_1, M)$, then $|x - y_q|_1 \leq 2M$, and $x$ is in $B_1\big(\mathcal{T}_i^{*1}(1), 2M\big)$. If $y_q$ is closed, then either $q = 0$ or $q = r$. If $q = r$, then $x$ is in $\mathcal{T}_i^{*0}(0) \cup \mathcal{T}_i^{*0}(1)$. If $q = 0$, then $y_0$ is in $\mathcal{T}_i^{*0}(0) \cup \mathcal{T}_i^{*0}(1)$. Yet $y_0$ belongs to Shell $(w_i, D)$, thus it is either in $\mathcal{A}_i(0)$ or in $\mathcal{T}_i^{*0}(0)$, but it cannot be in $\mathcal{T}_i^{*0}(1) \setminus \mathcal{T}_i^{*0}(0)$. Therefore $y_0$ is in $\mathcal{T}_i^{*0}(0)$, and $x$ is in $B_1\big(\overline{\mathcal{T}_i}^0(0), M + 1\big)$.

We conclude that the set (89.84) is included in

$$B_1\big(\overline{\mathcal{T}}_i^0(0),M+1\big)\cup\overline{\mathcal{T}}_i^0(1)\cup B_1\big(\overline{\mathcal{T}}_i^1(1),2M\big)\ \subset\ B^{0,0}_{1,M}(1)\cup B^{1,0}_{1,M}(0)\cup B^{1,1}_{1,M}(1)\,.$$

This completes the case $t=1$.

Let now $t\geq 2$ and suppose that the inclusion (89.80) has been proved until rank $t-1$. With the help of the semigroup property (89.74), we rewrite $\widehat{\mathcal{B}}_M\big(\text{Shell}\,(W(I),D),D,t\big)$ as

$$\widehat{\mathcal{B}}_M\big(\text{Shell}\,(W(I),D),D,t\big)\ =\ \widehat{\mathcal{B}}_M\Big(\widehat{\mathcal{B}}_M\big(\text{Shell}\,(W(I),D),D,t-1\big),D,1\Big)\,. \tag{89.85}$$

We use next the induction hypothesis. The identity (89.85) and the inclusion (89.80) at rank $t-1$ together yield that

$$\begin{aligned}
&\widehat{\mathcal{B}}_M\big(\text{Shell}\,(W(I),D),D,t\big)\ \subset\\
&\widehat{\mathcal{B}}_M\Big(\bigcup_{i\in I}\mathcal{B}\big(w_i,D\setminus\overline{\mathcal{T}}_i(t-1),t-1\big)\cup\bigcup_{0\leq s\leq t-1}\big(B^{s,0}_{1,M}(t-1-s)\,\cup\,B^{s,1}_{1,M}(t-s)\big),D,1\Big)\\
&\qquad\subset\ \widehat{\mathcal{B}}_M\Big(\bigcup_{i\in I}\mathcal{B}\big(w_i,D\setminus\overline{\mathcal{T}}_i(t-1),t-1\big),D,1\Big)\cup\\
&\qquad\qquad\qquad\widehat{\mathcal{B}}_M\Big(\bigcup_{0\leq s\leq t-1}\big(B^{s,0}_{1,M}(t-1-s)\,\cup\,B^{s,1}_{1,M}(t-s)\big),D,1\Big)\,.
\end{aligned} \tag{89.86}$$

We study separately the two sets on the right-hand side of (89.86). Let us start with the first set. It follows from (81.2) in proposition 81.1 and the definition of the sets $\mathcal{T}_i^*(t)$, $i\in I$, that

$$\forall i\in I\qquad \mathcal{A}_i(t-1)\ =\ \mathcal{B}\big(w_i,D\setminus\overline{\mathcal{T}}_i(t-1),t-1\big)\ =\ \mathcal{B}\big(w_i,D\setminus\overline{\mathcal{T}}_i(t),t-1\big)\,,$$

whence

$$\widehat{\mathcal{B}}_M\Big(\bigcup_{i\in I}\mathcal{B}\big(w_i,D\setminus\overline{\mathcal{T}}_i(t-1),t-1\big),D,1\Big)\ =\ \widehat{\mathcal{B}}_M\Big(\bigcup_{i\in I}\mathcal{B}\big(w_i,D\setminus\overline{\mathcal{T}}_i(t),t-1\big),D,1\Big)\,. \tag{89.87}$$

We would like to compare the set of the right-hand side with the set where $\overline{\mathcal{T}}(t)$ is used as a taboo set when performing the truncated exploration, i.e., the set

$$\widehat{\mathcal{B}}_M\Big(\bigcup_{i\in I}\mathcal{B}\big(w_i,D\setminus\overline{\mathcal{T}}_i(t),t-1\big),D\setminus\overline{\mathcal{T}}(t),1\Big)\,. \tag{89.88}$$

The point is that the previous set (89.88) is further included in the set

$$\bigcup_{i\in I}\mathcal{B}\big(w_i,D\setminus\overline{\mathcal{T}}_i(t),t\big)\,. \tag{89.89}$$

Thus, in the next lemma, we compare the set of the right-hand side of (89.87) with the set (89.89).

**Lemma 89.19.** *For any $t$ in $\{1,\dots,T\}$, we have the inclusion*

$$\widehat{\mathcal{B}}_M\Big(\bigcup_{i\in I}\mathcal{B}\big(w_i,D\setminus\overline{\mathcal{T}}_i(t),t-1\big),D,1\Big)\subset \bigcup_{i\in I}\mathcal{B}\big(w_i,D\setminus\overline{\mathcal{T}}_i(t),t\big)\cup\overline{\mathcal{T}}^0(t)\cup B_1\big(\overline{\mathcal{T}}^1(t),2M\big)\,. \quad (89.90)$$

*Proof.* Let $x$ belong to

$$\bigcup_{i\in I}\widehat{\mathcal{B}}_M\Big(\mathcal{B}\big(w_i,D\setminus\overline{\mathcal{T}}_i(t),t-1\big),D,1\Big)\,.$$

By definition, there exist $i$ in $I$ and a path $y_0,y_1,\cdots,y_r$ in $D$ joining a site $y_0$ of $\mathcal{B}\big(w_i,D\setminus\overline{\mathcal{T}}_i(t),t-1\big)$ to $y_r=x$ such that the sites $y_1,\cdots,y_{r-1}$ are in $D\setminus\mathcal{B}\big(w_i,D\setminus\overline{\mathcal{T}}_i(t),t-1\big)$, they are all open, and they are visited before $x$ by the truncated exploration starting from $y_1$. The truncation imposes that $y_2,\dots,y_r$ are in $B_1(y_1,M)$. If $r=0$ and the path is reduced to $y_0=x$, or if the whole path $y_0,y_1,\cdots,y_r$ does not visit $\overline{\mathcal{T}}_i(t)$, then $x$ is in

$$\mathcal{B}\Big(\mathcal{B}\big(w_i,D\setminus\overline{\mathcal{T}}_i(t),t-1\big),D\setminus\overline{\mathcal{T}}_i(t),1\Big)\;=\;\mathcal{B}\big(w_i,D\setminus\overline{\mathcal{T}}_i(t),t\big)\,.$$

Otherwise, let $s$ be the smallest index such that $y_s$ is in $\overline{\mathcal{T}}_i(t)$. As $y_0$ belongs to $\mathcal{B}\big(w_i,D\setminus\overline{\mathcal{T}}_i(t),t-1\big)$, which is disjoint from $\overline{\mathcal{T}}_i(t)$, we must have $s\geq 1$. Since both $y_s$ and $x$ are in $B_1(y_1,M)$, then $|x-y_s|_1\leq 2M$. We consider two cases:

- $y_s$ is closed. Then $s$ is equal to $r$, and $y_s=x$ is in $\overline{\mathcal{T}}^0(t)$.
- $y_s$ is open. Then $y_s$ is in $\overline{\mathcal{T}}_i^1(t)$ and $x$ is in $B_1\big(\overline{\mathcal{T}}^1(t),2M\big)$.

In each case, the site $x$ belongs to one of the sets appearing in the right-hand side of (89.90). This concludes the proof of the inclusion stated in the lemma. □

For the second set and last set of (89.86), we use lemma 89.15 to obtain

$$\begin{aligned}\widehat{\mathcal{B}}_M\Big(\bigcup_{0\leq s\leq t-1}\big(B^{s,0}_{1,M}(t-1-s)\,\cup\,B^{s,1}_{1,M}(t-s)\big),D,1\Big)\\ \subset B_1\Big(\bigcup_{0\leq s\leq t-1}\big(B^{s,0}_{1,M}(t-1-s)\,\cup\,B^{s,1}_{1,M}(t-s)\big),M+1\Big)\\ \subset\bigcup_{0\leq s\leq t-1}\big(B^{s,0}_{1,M}(t-s)\,\cup\,B^{s,1}_{1,M}(t+1-s)\big)\,.\end{aligned}\quad (89.91)$$

Using the identity (89.87), together with the inclusions (89.90) and (89.91), we deduce from (89.86) that

$$\begin{aligned}\widehat{\mathcal{B}}_M\Big(\text{Shell}\,\big(W(I),D\big),D,t\Big)\subset\bigcup_{i\in I}\mathcal{B}\big(w_i,D\setminus\overline{\mathcal{T}}_i(t),t\big)\cup\overline{\mathcal{T}}^0(t)\\ \cup B_1\big(\overline{\mathcal{T}}^1(t),2M\big)\,\cup\bigcup_{0\leq s\leq t-1}\big(B^{s,0}_{1,M}(t-s)\,\cup\,B^{s,1}_{1,M}(t+1-s)\big)\,.\end{aligned}$$

This is the inclusion (89.80) at rank $t$. □

### 89.6.2 Upper bounds

In this subsubsection, we fix an index $i$ in $\{1,\dots,N\}$ and we shall show that both subsets $\mathcal{B}\big(w_i, D\setminus\widehat{\mathcal{T}}_i^*(t),t\big)$ and $\mathcal{B}\big(w_i, D\setminus\widehat{\mathcal{T}}_i^-(t),t\big)$ are included in $\widehat{\mathcal{B}}_M\big(\text{Shell}\,(w_i,D),D,t\big)$. These inclusions yield an upper bound for the approximations. For convenience, we recall that, for $t$ in $\{1,\dots,T\}$, we have

$$\widehat{\mathcal{T}}_i^-(t)\,=\,\overline{\mathcal{T}}_i(t-1)\cup\mathcal{T}_i^{-1}(t)\,,\qquad \widehat{\mathcal{T}}_i^*(t)\,=\,\overline{\mathcal{T}}_i(t-1)\cup\mathcal{T}_i^{-0}(t)\cup\mathcal{T}_i^{*1}(t)\,.$$

**Proposition 89.20.** *We consider a configuration realizing the event $\mathcal{H}$ defined in* (87.4) *and an index $i$ in $\{1,\dots,N\}$. For any $t$ in $\{1,\dots,T\}$, we have*

$$\mathcal{B}\big(w_i, D\setminus\widehat{\mathcal{T}}_i^*(t),t\big)\,\subset\,\mathcal{B}\big(w_i, D\setminus\widehat{\mathcal{T}}_i^-(t),t\big)\,\subset\,\widehat{\mathcal{B}}_M\big(\mathit{Shell}\,(w_i,D),D,t\big)\,.\quad(89.92)$$

*Proof.* The first inclusion in (89.92) is straightforward, because $\widehat{\mathcal{T}}_i^-(t)$ is included in $\widehat{\mathcal{T}}_i^*(t)$. We prove next the second inclusion in (89.92). Let us fix $t$ in $\{1,\dots,T\}$. We claim that

$$\mathcal{B}\big(w_i, D\setminus\widehat{\mathcal{T}}_i^-(t),t\big)\,\subset\,\widehat{\mathcal{B}}_M\Big(\mathcal{B}\big(w_i,D\setminus\overline{\mathcal{T}}_i(t-1),t-1\big),D\setminus\widehat{\mathcal{T}}_i^-(t),1\Big)\,.\quad(89.93)$$

In order to prove this inclusion, we rewrite first the left-hand side of (89.93) with the help of the usual semigroup property (12.9):

$$\mathcal{B}\big(w_i, D\setminus\widehat{\mathcal{T}}_i^-(t),t\big)\,=\,\mathcal{B}\Big(\mathcal{B}\big(w_i,D\setminus\widehat{\mathcal{T}}_i^-(t),t-1\big),D\setminus\widehat{\mathcal{T}}_i^-(t),1\Big)\,.\qquad(89.94)$$

As $\widehat{\mathcal{T}}_i^-(t)$ contains $\overline{\mathcal{T}}_i(t-1)$, we deduce from (89.94) that

$$\mathcal{B}\big(w_i, D\setminus\widehat{\mathcal{T}}_i^-(t),t\big)\,\subset\,\mathcal{B}\Big(\mathcal{B}\big(w_i,D\setminus\overline{\mathcal{T}}_i(t-1),t-1\big),D\setminus\widehat{\mathcal{T}}_i^-(t),1\Big)\,.\quad(89.95)$$

Using the identity (81.2) of proposition 81.1, the inclusion (89.95) becomes

$$\mathcal{B}\big(w_i, D\setminus\widehat{\mathcal{T}}_i^-(t),t\big)\,\subset\,\mathcal{B}\big(\mathcal{A}_i(t-1),D\setminus\widehat{\mathcal{T}}_i^-(t),1\big)\,.$$

We consider then two sites $x,y$ such that

$$x\in\mathcal{A}_i(t-1)\,,\quad y\notin\mathcal{A}_i(t-1)\,,\quad y\in D\setminus\widehat{\mathcal{T}}_i^-(t)\,,\quad |x-y|_1=1\,,\qquad(89.96)$$

and we analyze the contribution of $y$ to the two sets appearing in (89.93). During the iteration $t$, when the explorer $i$ becomes active, it will visit the site $y$ unless that site has been visited before. We analyze the various possible cases for $y$:

• The site $y$ is in a boundary cluster. In this case, it is visited during the first iteration by an explorer other than $i$. Thus $y$ is in $\mathcal{T}_i^*(t-1)$. As $\mathcal{T}_i^*(t-1)$ is included in $\widehat{\mathcal{T}}_i^-(t)$ for $t\geq 1$, this case cannot occur.

• The site $y$ is not in $\mathcal{T}_i^*(t)$. In this case, the site $y$ is not in a boundary cluster so that its cluster is finite.

• The site $y$ is in $\mathcal{T}_i^*(t)$. In this case, the site $y$ is in

$$\mathcal{T}_i^*(t)\setminus\widehat{\mathcal{T}}_i^-(t)\,\subset\,\mathcal{T}_i^{*0}(t)\cup\mathcal{T}_i^{+1}(t)\,.$$

We consider two subcases:

$\star$ If $y$ is in $\mathcal{T}_i^{*0}(t)$, then it is closed.

$\star$ If $y$ is in $\mathcal{T}_i^{+1}(t)$, then it was visited by an explorer during the iteration $t \geq 1$, thus it is not in a boundary cluster.

We see that, in all the above cases, either $y$ is closed or the cluster of $y$ is finite, whence

$$\overline{\text{Cluster}}\,(y, M, D) \;=\; \overline{C}(y, D)\,.$$

Therefore the contribution of $y$ to the two sets appearing in (89.93) is the same. This is true for any pair of sites $x, y$ satisfying (89.96), and this completes the proof of the inclusion (89.93). This inclusion holds for any $t$ in $\{\,1, \dots, T\,\}$ and we will rely on it to prove by induction the second inclusion in (89.92). We start with the case $t = 1$. The inclusion (89.93) for $t = 1$ yields

$$\begin{aligned}\mathcal{B}\big(w_i, D \setminus \widehat{\mathcal{T}}_i^{-}(1), 1\big) \,&\subset\, \widehat{\mathcal{B}}_M\Big(\mathcal{B}\big(w_i, D \setminus \overline{\mathcal{T}}_i(0), 0\big), D \setminus \widehat{\mathcal{T}}_i^{-}(1), 1\Big)\\ &\subset\, \widehat{\mathcal{B}}_M\big(\mathcal{B}(w_i, D, 0), D, 1\big) \;=\; \widehat{\mathcal{B}}_M\big(\text{Shell}\,(w_i, D), D, 1\big)\,.\end{aligned}$$

Thus the inclusion (89.92) holds for $t = 1$. Let $t \geq 2$ and suppose that the inclusion (89.92) has been proved at rank $t-1$. The induction hypothesis at rank $t-1$ gives that

$$\begin{aligned}&\mathcal{B}\big(w_i, D \setminus \overline{\mathcal{T}}_i(t-1), t-1\big)\\ &\qquad \subset\, \mathcal{B}\big(w_i, D \setminus \widehat{\mathcal{T}}_i^{-}(t-1), t-1\big) \,\subset\, \widehat{\mathcal{B}}_M\big(\text{Shell}\,(w_i, D), D, t-1\big)\,.\end{aligned} \tag{89.97}$$

Substituting (89.97) into (89.93), we obtain

$$\begin{aligned}\mathcal{B}\big(w_i, D \setminus \widehat{\mathcal{T}}_i^{-}(t), t\big) \,&\subset\, \widehat{\mathcal{B}}_M\Big(\widehat{\mathcal{B}}_M\big(\text{Shell}\,(w_i, D), D, t-1\big), D \setminus \widehat{\mathcal{T}}_i^{-}(t), 1\Big)\\ &\subset\, \widehat{\mathcal{B}}_M\big(\text{Shell}\,(w_i, D), D, t\big)\,.\end{aligned}$$

This completes the induction step and the proof. $\square$

### 89.6.3 Lower bounds

We need to have an approximation result for the two sets

$$\bigcup_{i \in I} \mathcal{B}\big(w_i, D \setminus \widehat{\mathcal{T}}_i^{-}(t), t\big)\,, \qquad \bigcup_{i \in I} \mathcal{B}\big(w_i, D \setminus \widehat{\mathcal{T}}_i^{*}(t), t\big)\,, \tag{89.98}$$

which is the counterpart of proposition 89.18. Throughout this subsubsection, we fix a subset $I$ of $\{\,1, \dots, N\,\}$. We show first that $\widehat{\mathcal{B}}_M\big(\text{Shell}\,\big(W(I), D\big), D, t\big)$ is included in the union of $\bigcup_{i \in I} \mathcal{B}\big(w_i, D \setminus \widehat{\mathcal{T}}_i^{-}(t), t\big)$ and an adequate enlargement of the sets of the relevant taboo sites. We show next a similar result between $\mathcal{B}\big(w_i, D \setminus \widehat{\mathcal{T}}_i^{-}(t), t\big)$ and $\mathcal{B}\big(w_i, D \setminus \widehat{\mathcal{T}}_i^{*}(t), t\big)$, valid for any $i$ in $I$. These two results yield a lower bound for the approximation of the two sets (89.98) by $\widehat{\mathcal{B}}_M\big(\text{Shell}\,\big(W(I),D\big),D,t\big)$. For convenience, we recall that for $t$ in $\{\,1,\dots,T\,\}$, we have

$$\begin{aligned}\widehat{\mathcal{T}}_i^{-}(t) \,&=\, \overline{\mathcal{T}}_i(t-1)\cup\mathcal{T}_i^{-1}(t)\,,\\ \widehat{\mathcal{T}}_i^{*}(t) \,&=\, \overline{\mathcal{T}}_i(t-1)\cup\mathcal{T}_i^{-0}(t)\cup\mathcal{T}_i^{*1}(t)\,.\end{aligned}$$

As $\widehat{\mathcal{T}}_i^{-}(t)\subset\widehat{\mathcal{T}}_i^{*}(t)\subset\overline{\mathcal{T}}_i(t)$, we have the inclusions

$$\mathcal{B}\big(w_i,D\setminus\overline{\mathcal{T}}_i(t),t\big)\,\subset\,\mathcal{B}\big(w_i,D\setminus\widehat{\mathcal{T}}_i^{*}(t),t\big)\,\subset\,\mathcal{B}\big(w_i,D\setminus\widehat{\mathcal{T}}_i^{-}(t),t\big)\,,$$

and we can use directly proposition 89.18 to conclude that

$$\widehat{\mathcal{B}}_M\Big(\text{Shell}\,\big(W(I),D\big),D,t\Big)\,\subset \tag{89.99}$$

$$\bigcup_{i\in I}\mathcal{B}\big(w_i,D\setminus\widehat{\mathcal{T}}_i^{*}(t),t\big)\,\cup\bigcup_{0\le s\le t}\big(B_{1,M}^{s,0}(t-s)\,\cup\,B_{1,M}^{s,1}(t-s+1)\big) \tag{89.100}$$

$$\subset\,\bigcup_{i\in I}\mathcal{B}\big(w_i,D\setminus\widehat{\mathcal{T}}_i^{-}(t),t\big)\,\cup\bigcup_{0\le s\le t}\big(B_{1,M}^{s,0}(t-s)\,\cup\,B_{1,M}^{s,1}(t-s+1)\big)\,. \tag{89.101}$$

The trouble with this approximation is that the error term involves the set $\mathcal{T}^{*0}(t)$ of the closed relevant taboo sites at time $t$. We are trying to control this set with the help of the inequalities (88.34), (88.36) and we will fail if the error created when approximating the terms on the right-hand side of (88.34) and (88.36) depend on $\mathcal{T}^{*0}(t)$. The point is that the set $\widehat{\mathcal{T}}_i^{-}(t)$ does not contain $\mathcal{T}_i^{*0}(t)$, thus we should be able to remove it from the last union in (89.101). Yet this would not be enough. Indeed, the goal of the first inclusion (89.99) is to obtain a control on $\mathcal{T}^{-0}(t)$, and the control on $\mathcal{T}^{-0}(t)$ will in turn be used to get a control on $\mathcal{T}^{+1}(t)$, so it is also essential that this last set does not appear in the last union in (89.100). We will still use the set $B_{1,M}^{t,0}(s)$, defined in (89.79), but instead of the set $B_{1,M}^{t,1}(s)$, we shall use the set $\widehat{B}_{1,M}^{-,t,1}(s)$ given by

$$\forall t\in\{\,1,\dots,T\,\}\quad\forall s\geq 0\qquad \widehat{B}_{1,M}^{-,t,1}(s)\,=\,B_1\big(\overline{\mathcal{T}}^{1}(t-1)\cup\mathcal{T}^{-1}(t),s(M+1)+M-1\big)\,.$$

In the next proposition, we manage to prove a refined lower bound on the first set in (89.98) that suits our purpose.

**Proposition 89.21.** *We consider a configuration realizing the event $\mathcal{H}$ defined in* (87.4) *and a subset $I$ of $\{\,1,\dots,N\,\}$. For any $t$ in $\{\,1,\dots,T\,\}$, we have*

$$\widehat{\mathcal{B}}_M\Big(\mathit{Shell}\big(W(I),D\big),D,t\Big)\,\subset \tag{89.102}$$

$$\bigcup_{i\in I}\mathcal{B}\big(w_i,D\setminus\widehat{\mathcal{T}}_i^{-}(t),t\big)\,\cup\bigcup_{0\le s\le t-1}\Big(B_{1,M}^{s,0}(t-s)\,\cup\,\widehat{B}_{1,M}^{-,s+1,1}(t-s)\Big)\,.$$

*Proof.* We shall prove the proposition by induction on $t$ and we start with $t = 1$. Let $i$ belong to $I$. As the set $\widehat{\mathcal{T}}_i^-(1)$ contains $\mathcal{T}_i(0)$, we have

$$\mathcal{B}\big(\text{Shell}\,(w_i, D), D \setminus \widehat{\mathcal{T}}_i^-(1), 1\big) \setminus \text{Shell}\,\big(w_i, D\big) \;\subset\; D \setminus \bigcup_{1 \le j \le N} C(w_j, D)\,.$$

The occurrence of the event $\mathcal{H}$ guarantees that all the clusters in the above set are finite clusters of $\Lambda(4n)$, and also that all these finite clusters have cardinality strictly less than $M$. This implies that

$$\begin{aligned} \mathcal{B}\big(w_i, D \setminus \widehat{\mathcal{T}}_i^-(1), 1\big) \;&=\; \mathcal{B}\big(\text{Shell}\,(w_i, D), D \setminus \widehat{\mathcal{T}}_i^-(1), 1\big) \\ &= \widehat{\mathcal{B}}_M\big(\text{Shell}\,(w_i, D), D \setminus \widehat{\mathcal{T}}_i^-(1), 1\big) \;\subset\; \widehat{\mathcal{B}}_M\big(\text{Shell}\,(w_i, D), D, 1\big)\,. \end{aligned} \tag{89.103}$$

Taking the union of (89.103) over $i$ in $I$, we get

$$\begin{aligned} \bigcup_{i \in I} \mathcal{B}\big(w_i, D \setminus \widehat{\mathcal{T}}_i^-(1), 1\big) \;&=\; \bigcup_{i \in I} \widehat{\mathcal{B}}_M\big(\text{Shell}\,(w_i, D), D \setminus \widehat{\mathcal{T}}_i^-(1), 1\big) \\ &\subset \bigcup_{i \in I} \widehat{\mathcal{B}}_M\big(\text{Shell}\,(w_i, D), D, 1\big) \;=\; \widehat{\mathcal{B}}_M\Big(\text{Shell}\,\big(W(I), D\big), D, 1\Big)\,. \end{aligned}$$

Let $x$ be a site belonging to

$$\widehat{\mathcal{B}}_M\Big(\text{Shell}\,\big(W(I), D\big), D, 1\Big) \;\setminus\; \bigcup_{i \in I} \widehat{\mathcal{B}}_M\big(\text{Shell}\,(w_i, D), D \setminus \widehat{\mathcal{T}}_i^-(1), 1\big)\,. \tag{89.104}$$

By definition of $\widehat{\mathcal{B}}_M\big(\text{Shell}\,\big(W(I), D\big), D, 1\big)$, there exists a path $y_0, y_1, \cdots, y_r$ in $D$ joining a site $y_0$ of Shell $\big(W(I), D\big)$ to $y_r = x$ such that the sites $y_1, \cdots, y_{r-1}$ are in $D \setminus \text{Shell}\,\big(W(I), D\big)$, they are all open, and they are visited before $x$ by the truncated exploration starting from $y_1$. The truncation imposes that $y_2, \ldots, y_r$ are in $B_1(y_1, M)$. (if $x$ is in Shell $\big(W(I), D\big)$, then the sequence is reduced to $x = y_0$ and we have $r = 0$). We have also

$$\text{Shell}\,\big(W(I), D\big) \;=\; \bigcup_{i \in I} \text{Shell}\,(w_i, D)\,,$$

so that there exists an index $i$ in $I$ such that $y_0$ is in Shell $(w_i, D)$. Since $x$ is not in $\widehat{\mathcal{B}}_M\big(\text{Shell}\,(w_i, D), D \setminus \widehat{\mathcal{T}}_i^-(1), 1\big)$, then the path $y_0, \cdots, y_r$ must go through $\widehat{\mathcal{T}}_i^-(1)$. Thus there exists $s$ such that $y_s$ is in $\widehat{\mathcal{T}}_i^-(1)$. For $s \ge 1$, both $y_s$ and $x$ are in $B_1(y_1, M)$, whence $|x - y_s|_1 \;\le\; 2M$. If $s = 0$, then $y_s = y_0$ is in

$$\text{Shell}\,(w_i, D) \cap \widehat{\mathcal{T}}_i^-(1) \;=\; \text{Shell}\,(w_i, D) \cap \mathcal{T}_i^*(0)\,,$$

so that $x$ is in $B_1(\mathcal{T}_i^{*0}(0), M + 1)$. If $0 < s < r$, then $y_s$ is in $\mathcal{T}_i^{-1}(1)$, so that $x$ is in $B_1(\mathcal{T}_i^{-1}(1), 2M)$. If $s = r$, then $x$ is in $\widehat{\mathcal{T}}_i^-(1) = \mathcal{T}_i(0) \cup \mathcal{T}_i^{-1}(1)$. We conclude that the set (89.104) is included in $B_1(\mathcal{T}_i^{*0}(0), M + 1) \cup B_1(\mathcal{T}_i^{-1}(1), 2M)$. This completes the case $t = 1$.

We perform now the induction step. Let $t \geq 2$ and suppose that the inclusion (89.102) has been proved until rank $t-1$. Using the semigroup property (89.74), we write $\widehat{\mathcal{B}}_M\big(\text{Shell}\,\big(W(I),D\big),D,t\big)$ as

$$\widehat{\mathcal{B}}_M\big(\text{Shell}\,\big(W(I),D\big),D,t\big) \;=\; \widehat{\mathcal{B}}_M\Big(\widehat{\mathcal{B}}_M\big(\text{Shell}\,\big(W(I),D\big),D,t-1\big),D,1\Big)\,. \tag{89.105}$$

We use next the induction hypothesis. The identity (89.105) and the inclusion (89.102) at rank $t-1$ together yield that

$$\begin{aligned}\widehat{\mathcal{B}}_M\big(\text{Shell}\,\big(W(I),D\big),D,t\big) \;\subset\; &\widehat{\mathcal{B}}_M\Bigg(\bigcup_{i\in I}\mathcal{B}\big(w_i,D\setminus\widehat{\mathcal{T}}_i^-(t-1),t-1\big)\\ &\cup\bigcup_{0\leq s\leq t-2}\Big(B_{1,M}^{s,0}(t-1-s)\,\cup\,\widehat{B}_{1,M}^{-,s+1,1}(t-1-s)\Big),D,1\Bigg)\\ \subset\;&\bigcup_{i\in I}\widehat{\mathcal{B}}_M\Big(\mathcal{B}\big(w_i,D\setminus\widehat{\mathcal{T}}_i^-(t-1),t-1\big),D,1\Big)\cup\\ &\widehat{\mathcal{B}}_M\Big(\bigcup_{0\leq s\leq t-2}\Big(B_{1,M}^{s,0}(t-1-s)\,\cup\,\widehat{B}_{1,M}^{-,s+1,1}(t-1-s)\Big),D,1\Big)\,.\end{aligned} \tag{89.106}$$

We study separately the two sets on the right-hand side of (89.106). Let us start with the first set.

**Lemma 89.22.** *For $t$ in $\{\,2,\dots,T\,\}$ and any $i$ in $I$, we have the inclusion*

$$\begin{aligned}&\widehat{\mathcal{B}}_M\Big(\mathcal{B}\big(w_i,D\setminus\widehat{\mathcal{T}}_i^-(t-1),t-1\big),D,1\Big)\;\subset\;\mathcal{B}\big(w_i,D\setminus\widehat{\mathcal{T}}_i^-(t),t\big)\cup\\ &\overline{\mathcal{T}}_i^0(t-2)\cup B_1\big(\mathcal{T}_i^{*0}(t-1),M+1\big)\cup B_1\big(\overline{\mathcal{T}}_i^1(t-1)\cup\mathcal{T}_i^{-1}(t),2M\big)\,.\end{aligned} \tag{89.107}$$

Before doing the proof, let us explain the origin of the various terms appearing in (89.107). These terms correspond to relevant taboo sites of $\widehat{\mathcal{T}}_i^-(t)$ which are neighbours of a site of $\mathcal{B}\big(w_i,D\setminus\widehat{\mathcal{T}}_i^-(t-1),t-1\big)$, or belong to this set. The truncated exploration $\widehat{\mathcal{B}}_M(\cdot,D,1)$ will explore these sites and this will create a discrepancy compared to the taboo exploration $\mathcal{B}\big(w_i,D\setminus\widehat{\mathcal{T}}_i^-(t),t\big)$. More precisely, let $x$ be a site which has a neighbour $y$ in $\mathcal{B}\big(w_i,D\setminus\widehat{\mathcal{T}}_i^-(t-1),t-1\big)$. Several cases can occur:

- The neighbour $y$ is in fact in $\widehat{\mathcal{T}}_i^-(t)\setminus\widehat{\mathcal{T}}_i^-(t-1)$, hence in $\mathcal{T}_i^*(t-1)\cup\mathcal{T}_i^{-1}(t)$. Then the discrepancy will be included in $B_1(x,M)\subset B_1(y,M+1)$.
- The site $x$ is closed. Then $\overline{\text{Cluster}}\,(x,M,D)$ is reduced to $x$. It will create a discrepancy only if $x$ is in $\widehat{\mathcal{T}}_i^-(t)$. As $x$ is closed, it is then in $\overline{\mathcal{T}}_i^0(t-1)$.
- The site $x$ is open. Then $\overline{\text{Cluster}}\,(x,M,D)$ is included in $B_1(x,M)$. There will be a discrepancy if the truncated exploration starting from $x$ has to go through sites of $\widehat{\mathcal{T}}_i^-(t)$. If $B_1(x,M)$ contains open relevant taboo sites of $\widehat{\mathcal{T}}_i^-(t)$, then the discrepancies will be within distance $2M$ of these sites. For the closed relevant taboo sites, the discrepancies are limited to $\overline{\mathcal{T}}_i^0(t-1)\cap B_1(x,M)$.

All the possible causes of discrepancies listed above are taken into account in the formula (89.107). We are now ready to start the rigorous proof of the lemma.

*Proof.* Let $x$ belong to $\widehat{\mathcal{B}}_M\big(\mathcal{B}(w_i, D\setminus\widehat{\mathcal{T}}_i^-(t-1), t-1), D, 1\big)$. By definition, there exists a path $y_0,\cdots,y_r$ in $D$ joining a site $y_0$ of $\mathcal{B}\big(w_i, D\setminus\widehat{\mathcal{T}}_i^-(t-1), t-1\big)$ to $y_r = x$ such that the sites $y_1,\cdots,y_{r-1}$ are all open, and they are visited before $x$ by the truncated exploration starting from $y_1$ (we can very well have $r=0$ and $x=y_0$). The truncation imposes that $y_2,\dots,y_r$ are in $B_1(y_1, M)$. Furthermore, there exists a path $z_0,\cdots,z_s$ in $D\setminus\widehat{\mathcal{T}}_i^-(t-1)$ from $z_0=w_i$ to $z_s=y_0$ such that the number of closed sites among $z_0,\cdots,z_{s-1}$ is less than or equal to $t-1$. Notice that $y_0$ might be closed, but in any case the path $z_0,\cdots,z_s,y_1,\cdots,y_{r-1}$ contains at most $t$ closed sites. Now we have

$$\widehat{\mathcal{T}}_i^-(t) \,=\, \widehat{\mathcal{T}}_i^-(t-1)\cup\mathcal{T}_i^*(t-1)\cup\mathcal{T}_i^{-1}(t)\,.$$

We consider several cases:

• No site among $z_0,\cdots,z_s,y_1,\cdots,y_r$ is in $\widehat{\mathcal{T}}_i^-(t)$. Then $x$ is in $\mathcal{B}\big(w_i, D\setminus\widehat{\mathcal{T}}_i^-(t), t\big)$.

• The path $y_1,\cdots,y_{r-1}$ visits $\widehat{\mathcal{T}}_i^-(t)$. Since these sites are open, then $x$ is in

$$B_1\big(\overline{\mathcal{T}}_i^1(t-1)\cup\mathcal{T}_i^{-1}(t), 2M\big)\,.$$

• The site $y_0$ is in $\widehat{\mathcal{T}}_i^-(t)$. As $y_0$ is not in $\widehat{\mathcal{T}}_i^-(t-1)$, it is in

$$\widehat{\mathcal{T}}_i^-(t)\setminus\widehat{\mathcal{T}}_i^-(t-1)\,\subset\,\mathcal{T}_i^*(t-1)\cup\mathcal{T}_i^{-1}(t)\,.$$

Then $x$ is in $B_1\big(\mathcal{T}_i^*(t-1)\cup\mathcal{T}_i^{-1}(t), M+1\big)$.

• The site $y_r$ is in $\widehat{\mathcal{T}}_i^-(t)$. Then $x=y_r$ is in $\widehat{\mathcal{T}}_i^-(t)$.

• The path $z_0,\cdots,z_{s-1}$ visits $\widehat{\mathcal{T}}_i^-(t)$. We know also that this path does not visit $\widehat{\mathcal{T}}_i^-(t-1)$. Let $q<s$ be the smallest index such that $z_q$ is in $\widehat{\mathcal{T}}_i^-(t)$. The path $z_0,\cdots,z_{q-1}$ is a path in $D\setminus\widehat{\mathcal{T}}_i^-(t)$ from $z_0=w_i$ to $z_{q-1}$ containing at most $t-1$ closed sites. In particular, this path does not contain any site belonging to $\overline{\mathcal{T}}_i(t-1)$, hence the whole path is included in $\mathcal{A}_i(t-1)$. Suppose that $z_q$ is not in $\overline{\mathcal{T}}_i(t-1)$. In this case, the path $z_0,\cdots,z_q$ is a path in $D\setminus\overline{\mathcal{T}}_i(t-1)$ showing that $T_{D\setminus\overline{\mathcal{T}}_i(t-1)}(w_i,z_q)\leq t-1$ and the site $z_q$ would be in

$$\mathcal{B}\big(w_i, D\setminus\overline{\mathcal{T}}_i(t-1), t-1\big)\,=\,\mathcal{A}_i(t-1)\,.$$

However, this is not compatible with the fact that $z_q$ is also in $\widehat{\mathcal{T}}_i^-(t)$. Thus this case cannot occur. We conclude that the site $z_q$ is in

$$\begin{gathered}\overline{\mathcal{T}}_i(t-1)\setminus\widehat{\mathcal{T}}_i^-(t-1)\,\subset\,\overline{\mathcal{T}}_i(t-1)\setminus\big(\overline{\mathcal{T}}_i(t-2)\cup\mathcal{T}_i^{-1}(t-1)\big)\\ \subset\,\big(\mathcal{T}_i^{*0}(t-1)\cup\mathcal{T}_i^{*1}(t-1)\big)\setminus\mathcal{T}_i^{-1}(t-1)\,\subset\,\mathcal{T}_i^{*0}(t-1)\cup\mathcal{T}_i^{+1}(t-1)\,.\end{gathered}\tag{89.108}$$

Suppose that $z_q$ is open. Then the visiting time of $z_q$ is equal to $t-1\geq 1$, so that $z_q$ was not visited during the first iteration at $t=0$. Therefore the cluster of $z_q$ is finite and its cardinality is strictly less than $M$. If $z_q$ is closed, its cluster is obviously finite. Furthermore, it follows from the definition of the relevant taboo sets that the last set in (89.108) is disjoint from

$$\mathcal{A}_i(t-2)\,=\,\mathcal{B}\big(w_i, D\setminus\overline{\mathcal{T}}_i(t-2), t-2\big)\,\supset\,\mathcal{B}\big(w_i, D\setminus\widehat{\mathcal{T}}_i^-(t-1), t-2\big)\,.$$

Therefore $z_q$ is not in $\mathcal{B}\big(w_i, D \setminus \widehat{\mathcal{T}}_i^-(t-1), t-2\big)$. This implies that the path $z_0, \cdots, z_{q-1}$ contains $t-1$ closed sites. However, the path $z_0, \cdots, z_{s-1}$ contains at most $t-1$ closed sites, thus the sites $z_q, \dots, z_{s-1}$ are all open. As the cluster of $z_q$ is finite and contains these sites, we must have $s-1-q+1 < M$. Since the site $z_q$ is open, it is not in $\mathcal{T}_i^{*0}(t-1)$, thus it follows from (89.108) that $z_q$ belongs to $\mathcal{T}_i^{+1}(t-1)$. The path $z_q, \cdots, z_{s-1}, z_s, y_1, \cdots, y_r$ starts at a site in $\mathcal{T}_i^{+1}(t-1)$, ends at $y_r = x$ and its displacement is less than or equal to

$$|y_r - z_q|_1 \,\leq\, |y_r - y_1|_1 + 1 + |z_s - z_q|_1 \,\leq\, M + 1 + s - q \,\leq\, 2M\,.$$

Therefore $x$ is in $B_1\big(\mathcal{T}_i^{+1}(t-1), 2M\big)$.

Based on the inequalities derived in the previous cases, we conclude that

$$\widehat{\mathcal{B}}_M\Big(\mathcal{B}\big(w_i, D \setminus \widehat{\mathcal{T}}_i^-(t-1), t-1\big), D, 1\Big) \,\subset\, \mathcal{B}\big(w_i, D \setminus \widehat{\mathcal{T}}_i^-(t), t\big) \cup \\ \overline{\mathcal{T}}_i^0(t-2) \cup B_1\big(\mathcal{T}_i^{*0}(t-1), M+1\big) \cup B_1\big(\overline{\mathcal{T}}_i^1(t-1) \cup \mathcal{T}_i^{-1}(t), 2M\big)\,.$$

This concludes the proof of the inclusion stated in the lemma. □

Why is there an enlargement by $2M$ around the set $\overline{\mathcal{T}}_i^1(t-1) \cup \mathcal{T}_i^{-1}(t)$ in (89.107)? It seems more natural to have an enlargement by $M+1$, as is the case for the set $\overline{\mathcal{T}}_i^0(t-1)$. Indeed, when we add a taboo site $y$ to the exploration, the potential consequence is that the truncated exploration will not look at the cluster of $y$. As the exploration is truncated, the exploration starting from $y$ will stay confined to the ball $B_1(y, M)$. However, the distance $2M$ is necessary in order to cover the following scenario. Imagine that $y$ is a taboo site in $\mathcal{T}_i^{+1}(t-1)$. If we remove $y$ from the set of the taboo sites at time $t-1$, then the cluster of $y$ will be explored. As $y$ is visited at time $t-1 \geq 1$, it belongs to a finite cluster. Thus the exploration of this cluster is confined to $B_1(y, M)$, but it might visit a closed site on the boundary of the cluster of $y$, which itself is a neighbour of a site $z$ belonging to an unexplored cluster of size $M-1$. The cluster of $z$ will be explored at time $t$, and this can lead to visit sites which are at distance $M-1$ from $z$, thus at distance $2M$ from the site $y$ which was removed from the taboo set at time $t-1$. The sites in $\mathcal{T}_i^{-1}(t) \setminus \mathcal{T}_i^{+1}(t-1)$ are not relevant taboo sites at time $t-1$ for the explorer $i$, hence we can simply forget them during the iteration $t-1$. We stress also that lemma 89.22 holds only under various hypotheses that have been done previously. Specifically, the occurrence of the event $\mathcal{H}$ implies that the starting set $\{\, w_1, \dots, w_N \,\}$ contains exactly one vertex of each boundary cluster. This is crucial for the argument to work, because we rely on the fact that all the boundary clusters are explored during the first iteration at $t = 0$, and that all the remaining clusters, being finite clusters, have cardinality strictly smaller than $M$. This is why the influence of the taboo sites can be localized in balls of radius $M$.

The proof of the previous lemma is a bit cumbersome. In fact, there is trick that we use for several results, so it seems worth recording it in the following separate lemma.

**Lemma 89.23.** *Let $A, \mathcal{T}, \mathcal{S}$ be three subsets of $D$. If Shell$(A, D \setminus \mathcal{T}) \cap \mathcal{S} = \varnothing$, then*

$$\mathcal{B}\big(A, D\setminus\mathcal{T}, 1\big) \,\subset\, \mathcal{B}\big(A, D\setminus(\mathcal{T}\cup\mathcal{S}), 1\big) \cup \mathit{Shell}(\mathcal{S}, D)\,. \tag{89.109}$$

*Proof.* Let $x$ be a site in $\mathcal{B}\big(A, D\setminus\mathcal{T}, 1\big)$. By definition, there exists a path $y_0, \cdots, y_r$ in $D\setminus\mathcal{T}$ joining a site $y_0$ of $A\cap D\setminus\mathcal{T}$ to $y_r = x$ such that the path $y_0, \cdots, y_{r-1}$ contains at most one closed site. We consider two cases.

• The path $y_0, \cdots, y_r$ does not visit $\mathcal{S}$. In this case, the site $x$ is in

$$\mathcal{B}\big(A, D\setminus(\mathcal{T}\cup\mathcal{S}), 1\big)\,.$$

• The path $y_0, \cdots, y_r$ visits $\mathcal{S}$. We know in addition that this path does not visit $\mathcal{T}$. Let $q$ be the smallest index such that $y_q$ is in $\mathcal{S}$. As $\mathcal{S}$ is disjoint from Shell$\,(A, D\setminus\mathcal{T})$, the site $y_q$ is not in Shell$\,(A, D\setminus\mathcal{T})$. The path $y_0, \cdots, y_q$ goes from $A$ to $y_q$ in $D\setminus\mathcal{T}$, thus it must go through a closed site before reaching $y_q$. Therefore the sites $y_q, \cdots, y_{r-1}$ are all open, and the site $x$ is in Shell$\,(\mathcal{S}, D)$. We conclude from the two previous cases that $x$ is either in $\mathcal{B}\big(A, D\setminus(\mathcal{T}\cup\mathcal{S}), 1\big)$ or in Shell$\,(\mathcal{S}, D)$. This yields the desired inclusion (89.109). ☐

In order to check the hypothesis Shell$\,(A, D\setminus\mathcal{T})\cap\mathcal{S} = \varnothing$ of lemma 89.23, we will rely on the following result.

**Lemma 89.24.** *Let $i$ belong to $I$ and let $t$ in $\{\,0,\dots,T\,\}$. For any subsets $A, \mathcal{T}$ of $D$, we have*

$$A\subset\mathcal{A}_i(t),\,\overline{\mathcal{T}_i}(t)\subset\mathcal{T} \quad\Longrightarrow\quad \mathit{Shell}\big(A, D\setminus\mathcal{T}\big) \,\subset\, \mathcal{A}_i(t)\,.$$

*Proof.* Let $A, \mathcal{T}$ be two subsets of $D$ such that $A\subset\mathcal{A}_i(t)$ and $\overline{\mathcal{T}_i}(t)\subset\mathcal{T}$. We have

$$\text{Shell}\,\big(A, D\setminus\mathcal{T}\big) \,\subset\, \text{Shell}\,\big(\mathcal{A}_i(t), D\setminus\overline{\mathcal{T}_i}(t)\big)\,. \tag{89.110}$$

By proposition 81.1, we have $\mathcal{A}_i(t) = \mathcal{B}\big(w_i, D\setminus\overline{\mathcal{T}_i}(t), t\big)$. Using this equality, we see that $\mathcal{A}_i(t)$ is $T_{D\setminus\overline{\mathcal{T}_i}(t)}$-closed in the sense defined in (82.1). Therefore

$$\text{Shell}\,\big(\mathcal{A}_i(t), D\setminus\overline{\mathcal{T}_i}(t)\big) \;=\; \mathcal{A}_i(t)\,. \tag{89.111}$$

The claim of the lemma follows from the inclusion (89.110) and the equality (89.111). ☐

With the help of lemmas 89.23 and 89.24, we can prove the interesting intermediate result stated in the next lemma, which is easier to grasp than the result of lemma 89.22.

**Lemma 89.25.** *For $t$ in $\{\,2,\dots,T\,\}$ and any $i$ in $I$, we have the inclusion*

$$\begin{aligned}\mathcal{B}\big(w_i, D\setminus\widehat{\mathcal{T}}_i^{-}(t-1), t-1\big) \,\subset&\\ \mathcal{B}\big(w_i, D\setminus\widehat{\mathcal{T}}_i^{-}(t), t-1\big) \,\cup\, &\mathcal{T}_i^{*0}(t-1) \,\cup\, \mathit{Shell}\big(\mathcal{T}_i^{+1}(t-1), D\big)\,.\end{aligned} \tag{89.112}$$

*Proof.* In the proof of lemma 89.22, we relied on the definition of the balls $\mathcal{B}(w_i, \cdot, t)$ and we examined the structure of the paths connecting $w_i$ to a site in the ball. We use here a different strategy, namely we use the semigroup property of the balls, the key identity presented in proposition 81.1 and the previous lemma 89.23. To do so, we go one step backwards and we take the aggregate $\mathcal{A}_i(t-2)$ as the starting point (this is legitimate, recall that the result holds only for $t \geq 2$). Let us begin. We fix $t$ in $\{2, \dots, T\}$ and $i$ in $I$. Using the semigroup property (12.9), we have

$$\mathcal{B}\big(w_i, D \setminus \widehat{\mathcal{T}}_i^-(t-1), t-1\big) \,=\, \mathcal{B}\Big(\mathcal{B}\big(w_i, D \setminus \widehat{\mathcal{T}}_i^-(t-1), t-2\big), D \setminus \widehat{\mathcal{T}}_i^-(t-1), 1\Big)\,. \tag{89.113}$$

The set $\widehat{\mathcal{T}}_i^-(t-1)$ contains the set $\overline{\mathcal{T}}_i(t-2)$ of all the taboo sites which are relevant to build the aggregate $\mathcal{A}_i(t-2)$. Using the identity (81.2), we have

$$\mathcal{B}\big(w_i, D \setminus \widehat{\mathcal{T}}_i^-(t-1), t-2\big) \,=\, \mathcal{B}\big(w_i, D \setminus \overline{\mathcal{T}}_i(t-2), t-2\big) \,=\, \mathcal{A}_i(t-2)\,. \tag{89.114}$$

Substituting (89.114) into (89.113), we obtain that

$$\mathcal{B}\big(w_i, D \setminus \widehat{\mathcal{T}}_i^-(t-1), t-1\big) \,=\, \mathcal{B}\big(\mathcal{A}_i(t-2), D \setminus \widehat{\mathcal{T}}_i^-(t-1), 1\big)\,. \tag{89.115}$$

We wish to control how the above set is modified if we use $\overline{\mathcal{T}}_i(t-1)$ for the taboo set, instead of $\widehat{\mathcal{T}}_i^-(t-1)$. This is where lemma 89.23 comes into play. So we apply lemma 89.23 with the choices

$$A \,=\, \mathcal{A}_i(t-2)\,, \quad \mathcal{T} \,=\, \widehat{\mathcal{T}}_i^-(t-1)\,, \quad \mathcal{S} \,=\, \overline{\mathcal{T}}_i(t-1) \setminus \widehat{\mathcal{T}}_i^-(t-1)\,.$$

Recall that $\widehat{\mathcal{T}}_i^-(t-1)$ contains $\overline{\mathcal{T}}_i(t-2)$. Thanks to lemma 89.24, we have

$$\text{Shell}\,\big(\mathcal{A}_i(t-2), D \setminus \widehat{\mathcal{T}}_i^-(t-1)\big) \,=\, \mathcal{A}_i(t-2)\,.$$

We can now check the hypothesis of lemma 89.23:

$$\text{Shell}\,\big(\mathcal{A}_i(t-2), D \setminus \widehat{\mathcal{T}}_i^-(t-1)\big) \cap \mathcal{S} \,=\, \mathcal{A}_i(t-2) \cap \big(\overline{\mathcal{T}}_i(t-1) \setminus \widehat{\mathcal{T}}_i^-(t-1)\big) \,=\, \varnothing\,.$$

We apply lemma 89.23, and we get

$$\begin{aligned} &\mathcal{B}\big(\mathcal{A}_i(t-2), D \setminus \widehat{\mathcal{T}}_i^-(t-1), 1\big) \,\subset \\ &\qquad \mathcal{B}\big(\mathcal{A}_i(t-2), D \setminus \overline{\mathcal{T}}_i(t-1), 1\big) \,\cup\, \text{Shell}\,\big(\overline{\mathcal{T}}_i(t-1) \setminus \widehat{\mathcal{T}}_i^-(t-1), D\big)\,. \end{aligned} \tag{89.116}$$

Using again the identity (81.2), we have

$$\mathcal{A}_i(t-2) \,=\, \mathcal{B}\big(w_i, D \setminus \overline{\mathcal{T}}_i(t-2), t-2\big) \,=\, \mathcal{B}\big(w_i, D \setminus \overline{\mathcal{T}}_i(t-1), t-2\big)\,,$$

and this allows us to rewrite the first set in (89.116) as

$$\begin{aligned} \mathcal{B}\big(\mathcal{A}_i(t-2), D \setminus \overline{\mathcal{T}}_i(t-1), 1\big) \,&=\, \mathcal{B}\Big(\mathcal{B}\big(w_i, D \setminus \overline{\mathcal{T}}_i(t-1), t-2\big), D \setminus \overline{\mathcal{T}}_i(t-1), 1\Big) \\ &=\, \mathcal{B}\big(w_i, D \setminus \overline{\mathcal{T}}_i(t-1), t-1\big)\,. \end{aligned} \tag{89.117}$$

Furthermore, we have

$$\overline{\mathcal{T}}_i(t-1) \setminus \widehat{\mathcal{T}}_i^-(t-1) \subset \overline{\mathcal{T}}_i(t-1) \setminus \big(\overline{\mathcal{T}}_i(t-2) \cup \mathcal{T}_i^{-1}(t-1)\big)$$
$$\subset \mathcal{T}_i^{+0}(t-1) \cup \mathcal{T}_i^{+1}(t-1) \cup \big(\mathcal{T}_i^{-0}(t-1) \setminus \overline{\mathcal{T}}_i(t-2)\big)\,,$$

from which we deduce that

$$\text{Shell}\,\big(\overline{\mathcal{T}}_i(t-1) \setminus \widehat{\mathcal{T}}_i^-(t-1), D\big) \subset$$
$$\mathcal{T}_i^{+0}(t-1) \cup \mathcal{T}_i^{-0}(t-1) \cup \text{Shell}\,\big(\mathcal{T}_i^{+1}(t-1), D\big)\,. \quad (89.118)$$

The identity (89.115) and the three inclusions (89.116), (89.117), (89.118) yield

$$\mathcal{B}\big(w_i, D \setminus \widehat{\mathcal{T}}_i^-(t-1), t-1\big) \subset$$
$$\mathcal{B}\big(w_i, D \setminus \overline{\mathcal{T}}_i(t-1), t-1\big) \,\cup\, \mathcal{T}_i^{*0}(t-1) \,\cup\, \text{Shell}\,\big(\mathcal{T}_i^{+1}(t-1), D\big)\,. \quad (89.119)$$

The set $\widehat{\mathcal{T}}_i^-(t)$ contains the set $\overline{\mathcal{T}}_i(t-1)$ of all the taboo sites which are relevant to build the aggregate $\mathcal{A}_i(t-1)$. By proposition 81.1, we have

$$\mathcal{B}\big(w_i, D \setminus \overline{\mathcal{T}}_i(t-1), t-1\big) \,=\, \mathcal{A}_i(t-1) \,=\, \mathcal{B}\big(w_i, D \setminus \widehat{\mathcal{T}}_i^-(t), t-1\big)\,. \quad (89.120)$$

The desired result (89.112) follows from (89.119) and (89.120). □

We need one more lemma to present an alternative proof of lemma 89.22.

**Lemma 89.26.** *Let $A, \mathcal{T}$ be two subsets of $D$. Letting*

$$\mathcal{T}^0 \,=\, \big\{\, x \in \mathcal{T} : x \text{ is closed} \,\big\}\,, \qquad \mathcal{T}^1 \,=\, \big\{\, x \in \mathcal{T} : x \text{ is open} \,\big\}\,,$$

*we have*

$$\widehat{\mathcal{B}}_M\big(A, D, 1\big) \subset \mathcal{B}\big(A, D \setminus \mathcal{T}, 1\big) \cup \big(\mathcal{T}^0 \setminus A\big) \cup B_1\big(\mathcal{T}^0 \cap A, M+1\big) \cup B_1\big(\mathcal{T}^1, 2M\big)\,.$$

*Proof.* From the definition (89.70), we have

$$\widehat{\mathcal{B}}_M(A, D, 1) \,=\, \overline{\text{Clusters}}\,\big(A \cup \mathcal{N}(A), M, D\big) \,=\, \bigcup_{x \in A \,\cup\, \mathcal{N}(A)} \overline{\text{Cluster}}\,\big(x, M, D\big)\,. \quad (89.121)$$

The set $\widehat{\mathcal{B}}_M(A, D, 1)$ is always included in $\mathcal{B}\big(A, D, 1\big)$. If a site is in $\widehat{\mathcal{B}}_M(A, D, 1)$ but not in $\mathcal{B}\big(A, D \setminus \mathcal{T}, 1\big)$, it is because of the effect of the taboo set $\mathcal{T}$. We consider a site $x$ of $\mathcal{T}$ and we discuss several cases:

• $x \in \mathcal{T}^0 \setminus A$. As $x$ is closed, the truncated exploration might visit $x$, but it will then go on exploring other sites. For the taboo exploration, the site $x$ does not prevent the exploration of other sites than $x$ itself.

• $x \in \mathcal{T} \cap A$. The site $x$ might prevent the taboo exploration from starting from a neighbour $y$ of $x$ which is outside $A$. The sites that might be concerned are included in $B_1(x, M+1)$.

• $x \in \mathcal{T}^1$. When exploring a $B_1$ ball of radius $M$, the truncated exploration might go through $x$, and then reach a site at distance $2M$ from $x$. Such sites might be missed by the taboo exploration.

We see that all the damages created by the taboo set $\mathcal{T}$ are included in the sets of the right-hand side of (89.121). □

With the help of the last four lemmas, we can propose a streamlined alternative proof of lemma 89.22.

*Alternative proof of lemma 89.22.* We fix $t$ in $\{2, \dots, T\}$ and $i$ in $I$. By lemma 89.25, we have

$$\mathcal{B}\big(w_i, D \setminus \widehat{\mathcal{T}}_i^-(t-1), t-1\big) \subset \\ \mathcal{B}\big(w_i, D \setminus \widehat{\mathcal{T}}_i^-(t), t-1\big) \cup \mathcal{T}_i^{*0}(t-1) \cup \text{Shell}\big(\mathcal{T}_i^{+1}(t-1), D\big) . \quad (89.122)$$

We remark next that the sites in the set $\mathcal{T}_i^{+1}(t-1)$ are visited during the iteration $t-1 \geq 1$. However, all the boundary clusters are visited during the first step at $t=0$, therefore all the remaining sites belong to finite clusters. As the event $\mathcal{H}$ occurs, the finite clusters have strictly less than $M$ open sites, and we have the inclusion

$$\text{Shell}\big(\mathcal{T}_i^{+1}(t-1), D\big) \subset B_1\big(\mathcal{T}_i^{+1}(t-1), M-1\big) . \quad (89.123)$$

Let us precise why we can put $M-1$ in (89.123). If $x$ is in $\text{Shell}\big(\mathcal{T}_i^{+1}(t-1), D\big)$, then there exists a path $y_0, \cdots, y_r$ in $D$ joining a site $y_0$ of $\mathcal{T}_i^{+1}(t-1)$ to $x$ such that the sites $y_0, \cdots, y_{r-1}$ are all open. Since the cluster of $y_0$ is finite, then $r < M$ and $|x-y_0|_1 \leq r \leq M-1$.

Putting together the two inclusions (89.122) and (89.123), we conclude that

$$\mathcal{B}\big(w_i, D \setminus \widehat{\mathcal{T}}_i^-(t-1), t-1\big) \subset \\ \mathcal{B}\big(w_i, D \setminus \widehat{\mathcal{T}}_i^-(t), t-1\big) \cup \mathcal{T}_i^{*0}(t-1) \cup B_1\big(\mathcal{T}_i^{+1}(t-1), M-1\big) . \quad (89.124)$$

Using the monotonicity of $\widehat{\mathcal{B}}_M(\cdot,\cdot,\cdot)$, we deduce from (89.124) that

$$\widehat{\mathcal{B}}_M\Big(\mathcal{B}\big(w_i, D \setminus \widehat{\mathcal{T}}_i^-(t-1), t-1\big), D, 1\Big) \subset \widehat{\mathcal{B}}_M\Big(\mathcal{B}\big(w_i, D \setminus \widehat{\mathcal{T}}_i^-(t), t-1\big), D, 1\Big) \\ \cup \widehat{\mathcal{B}}_M\big(\mathcal{T}_i^{*0}(t-1), D, 1\big) \cup \widehat{\mathcal{B}}_M\Big(B_1\big(\mathcal{T}_i^{+1}(t-1), M-1\big), D, 1\Big) .$$

We use lemma 89.15 to control the last two terms in the above formula:

$$\widehat{\mathcal{B}}_M\Big(\mathcal{B}\big(w_i, D \setminus \widehat{\mathcal{T}}_i^-(t-1), t-1\big), D, 1\Big) \subset \widehat{\mathcal{B}}_M\Big(\mathcal{B}\big(w_i, D \setminus \widehat{\mathcal{T}}_i^-(t), t-1\big), D, 1\Big) \\ \cup B_1\big(\mathcal{T}_i^{*0}(t-1), M+1\big) \cup B_1\big(\mathcal{T}_i^{+1}(t-1), 2M\big) . \quad (89.125)$$

We apply next lemma 89.26 with $A = \mathcal{B}\big(w_i, D \setminus \widehat{\mathcal{T}}_i^-(t), t-1\big)$ and $\mathcal{T} = \widehat{\mathcal{T}}_i^-(t)$. Noticing that

$$\mathcal{T}^0 \setminus A \subset \overline{\mathcal{T}}_i^0(t-1) , \qquad \mathcal{T}^0 \cap A = \varnothing , \qquad \mathcal{T}^1 \subset \overline{\mathcal{T}}_i^1(t-1) \cup \mathcal{T}_i^{-1}(t) ,$$

we obtain

$$\widehat{\mathcal{B}}_M\Big(\mathcal{B}\big(w_i, D \setminus \widehat{\mathcal{T}}_i^-(t), t-1\big), D, 1\Big) \subset \\ \mathcal{B}\big(w_i, D \setminus \widehat{\mathcal{T}}_i^-(t), t\big) \cup \overline{\mathcal{T}}_i^0(t-1) \cup B_1\big(\overline{\mathcal{T}}_i^1(t-1) \cup \mathcal{T}_i^{-1}(t), 2M\big) . \quad (89.126)$$

The two inclusions (89.125) and (89.126) together imply the result (89.107) stated in lemma 89.22. ☐

If we take into account the three intermediate lemmas 89.23, 89.25 and 89.26, then the alternative proof of lemma 89.22 presented above is not considerably shorter than the first proof. Its advantage is that it does not make appeal to the paths connecting $w_i$ to a site at travel distance $t$, and it is more automatic. However, each proof has its own advantage, so in the end we decided to keep both presentations. We are now ready to complete the induction step that we started to prove the inclusion (89.102). We had arrived at formula (89.106). We apply lemma 89.22 to each set $\widehat{\mathcal{B}}_M\big(\mathcal{B}\big(w_i, D\setminus\widehat{\mathcal{T}}_i^-(t-1), t-1\big), D, 1\big)$ and we take the union of (89.107) over $i$ in $I$. We obtain

$$\begin{aligned}\bigcup_{i\in I}\widehat{\mathcal{B}}_M\Big(\mathcal{B}\big(w_i, D\setminus\widehat{\mathcal{T}}_i^-(t-1), t-1\big), D, 1\Big) \subset \\ \bigcup_{i\in I}\Big(\mathcal{B}\big(w_i, D\setminus\widehat{\mathcal{T}}_i^-(t), t\big)\cup\overline{\mathcal{T}}_i^0(t-2)\cup B_1\big(\mathcal{T}_i^{*0}(t-1), M+1\big) \\ \cup B_1\big(\overline{\mathcal{T}}_i^1(t-1)\cup\mathcal{T}_i^{-1}(t), 2M\big)\Big) \\ \subset \Big(\bigcup_{i\in I}\mathcal{B}\big(w_i, D\setminus\widehat{\mathcal{T}}_i^-(t), t\big)\Big)\cup\overline{\mathcal{T}}_i^0(t-2)\cup B_{1,M}^{t-1,0}(1)\cup\widehat{B}_{1,M}^{-,t,1}(1)\,.\end{aligned} \quad (89.127)$$

For the second set of the right-hand side of (89.106), we use lemma 89.15 to obtain

$$\begin{aligned}\widehat{\mathcal{B}}_M\left(\bigcup_{0\leq s\leq t-2}\Big(B_{1,M}^{s,0}(t-1-s)\,\cup\,\widehat{B}_{1,M}^{-,s+1,1}(t-1-s)\Big), D, 1\right) \\ \subset B_1\left(\bigcup_{0\leq s\leq t-2}\Big(B_{1,M}^{s,0}(t-1-s)\,\cup\,\widehat{B}_{1,M}^{-,s+1,1}(t-1-s)\Big), M+1\right) \\ \subset \bigcup_{0\leq s\leq t-2}\Big(B_{1,M}^{s,0}(t-s)\,\cup\,\widehat{B}_{1,M}^{-,s+1,1}(t-s)\Big)\,.\end{aligned} \quad (89.128)$$

Using the inclusions (89.127) and (89.128), we deduce from (89.106) that

$$\begin{aligned}\widehat{\mathcal{B}}_M\Big(\mathrm{Shell}\,\big(W(I), D\big), D, t\Big) \subset \\ \Big(\bigcup_{i\in I}\mathcal{B}\big(w_i, D\setminus\widehat{\mathcal{T}}_i^-(t), t\big)\Big)\,\cup\,\overline{\mathcal{T}}^0(t-2)\cup B_{1,M}^{t-1,0}(1)\cup\widehat{B}_{1,M}^{-,t,1}(1) \\ \cup\bigcup_{0\leq s\leq t-2}\Big(B_{1,M}^{s,0}(t-s)\,\cup\,\widehat{B}_{1,M}^{-,s+1,1}(t-s)\Big)\,.\end{aligned} \quad (89.129)$$

We check next that

$$\begin{aligned}\overline{\mathcal{T}}^0(t-2)\cup B_{1,M}^{t-1,0}(1)\cup\widehat{B}_{1,M}^{-,t,1}(1)\cup\bigcup_{0\leq s\leq t-2}\Big(B_{1,M}^{s,0}(t-s)\,\cup\,\widehat{B}_{1,M}^{-,s+1,1}(t-s)\Big) \\ = \bigcup_{0\leq s\leq t-1}\Big(B_{1,M}^{s,0}(t-s)\,\cup\,\widehat{B}_{1,M}^{-,s+1,1}(t-s)\Big)\,.\end{aligned} \quad (89.130)$$

The inclusion (89.102) at rank $t$ follows from (89.129) and (89.130). □

We develop next an approximation for the second set appearing in (89.98). This approximation will be based on the following lemma.

**Lemma 89.27.** *For any $t$ in $\{1,\dots,T\}$ and any $i$ in $I$, we have the inclusion*

$$\mathcal{B}\big(w_i, D\backslash\widehat{\mathcal{T}}_i^{\,-}(t), t\big) \subset \mathcal{B}\big(w_i, D\backslash\widehat{\mathcal{T}}_i^{\,*}(t), t\big) \cup \mathcal{T}_i^{\,-0}(t) \cup B_1\big(\mathcal{T}_i^{\,+1}(t), M-1\big)\,. \quad (89.131)$$

In practice, when using (89.131), we might replace $M-1$ by $M$ to alleviate the formulas. We shall give two proofs of lemma 89.27. The first one is based on lemma 89.23, the second one is a direct argument based on the analysis of paths.

*First proof of lemma 89.27.* Let us fix $t$ in $\{1,\dots,T\}$ and $i$ in $I$. By the definition of the balls $\mathcal{B}(w_i,\cdot,t)$, we have

$$\mathcal{B}\big(w_i, D\setminus\widehat{\mathcal{T}}_i^{\,-}(t), t\big) \,=\, \mathcal{B}\Big(\mathcal{B}\big(w_i, D\setminus\widehat{\mathcal{T}}_i^{\,-}(t), t-1\big), D\setminus\widehat{\mathcal{T}}_i^{\,-}(t), 1\Big)\,. \quad (89.132)$$

Let us examine the set $\mathcal{B}\big(w_i, D\setminus\widehat{\mathcal{T}}_i^{\,-}(t), t-1\big)$. We have the inclusions

$$\overline{\mathcal{T}}_i(t-1) \,\subset\, \widehat{\mathcal{T}}_i^{\,-}(t) \,\subset\, \widehat{\mathcal{T}}_i^{\,*}(t) \,\subset\, \widehat{\mathcal{A}}_i(T)\,, \quad (89.133)$$

where the taboo set $\widehat{\mathcal{A}}_i(T)$ was defined in (80.28). From the inclusions (89.133), we deduce that

$$\begin{aligned}\mathcal{B}\big(w_i, D\setminus\widehat{\mathcal{A}}_i(T), t-1\big) \,&\subset\, \mathcal{B}\big(w_i, D\setminus\widehat{\mathcal{T}}_i^{\,*}(t), t-1\big)\\ &\subset\, \mathcal{B}\big(w_i, D\setminus\widehat{\mathcal{T}}_i^{\,-}(t), t-1\big) \,\subset\, \mathcal{B}\big(w_i, D\setminus\overline{\mathcal{T}}_i(t-1), t-1\big)\,. \quad (89.134)\end{aligned}$$

However, using the identity (81.1) and proposition 81.1, we obtain that

$$\mathcal{A}_i(t-1) \,=\, \mathcal{B}\big(w_i, D\setminus\widehat{\mathcal{A}}_i(T), t-1\big) \,=\, \mathcal{B}\big(w_i, D\setminus\overline{\mathcal{T}}_i(t-1), t-1\big)\,,$$

thus the four sets appearing in (89.134) are equal to $\mathcal{A}_i(t-1)$. Coming back to (89.132), we can replace $\widehat{\mathcal{T}}_i^{\,-}(t)$ by $\widehat{\mathcal{T}}_i^{\,*}(t)$ for the innermost taboo set, and we get

$$\mathcal{B}\big(w_i, D\setminus\widehat{\mathcal{T}}_i^{\,-}(t), t\big) \,=\, \mathcal{B}\Big(\mathcal{B}\big(w_i, D\setminus\widehat{\mathcal{T}}_i^{\,*}(t), t-1\big), D\setminus\widehat{\mathcal{T}}_i^{\,-}(t), 1\Big)\,.$$

In the next step, we would like to replace also $\widehat{\mathcal{T}}_i^{\,-}(t)$ by $\widehat{\mathcal{T}}_i^{\,*}(t)$ for the outermost taboo set. We have the inclusion

$$\begin{aligned}\widehat{\mathcal{T}}_i^{\,*}(t)\setminus\widehat{\mathcal{T}}_i^{\,-}(t) \,&\subset\, \big(\mathcal{T}_i^{\,-0}(t)\cup\mathcal{T}_i^{\,*1}(t)\big)\setminus\big(\mathcal{T}_i^{\,-1}(t)\cup\overline{\mathcal{T}}_i(t-1)\big)\\ &\subset\, \big(\mathcal{T}_i^{\,-0}(t)\cup\mathcal{T}_i^{\,+1}(t)\big)\setminus\overline{\mathcal{T}}_i(t-1)\,. \quad (89.135)\end{aligned}$$

We are going to apply lemma 89.23 with the choices

$$A \,=\, \mathcal{B}\big(w_i, D\setminus\widehat{\mathcal{T}}_i^{\,*}(t), t-1\big)\,, \quad \mathcal{T} \,=\, \widehat{\mathcal{T}}_i^{\,-}(t)\,, \quad \mathcal{S} \,=\, \widehat{\mathcal{T}}_i^{\,*}(t)\setminus\widehat{\mathcal{T}}_i^{\,-}(t)\,.$$

Recall that $\widehat{\mathcal{T}}_i^{\,-}(t)$ contains $\overline{\mathcal{T}}_i(t-1)$. Thanks to lemma 89.24 and the inclusions (89.134), we have

$$\text{Shell}\Big(\mathcal{B}\big(w_i, D\setminus\widehat{\mathcal{T}}_i^{\,*}(t), t-1\big), D\setminus\widehat{\mathcal{T}}_i^{\,-}(t)\Big) \,\subset\, \mathcal{A}_i(t-1)\,. \quad (89.136)$$

We check the hypothesis of lemma 89.23. With the help of (89.136), we have

$$\text{Shell}\,\Big(\mathcal{B}\big(w_i, D\setminus \widehat{\mathcal{T}}_i^*(t), t-1\big), D\setminus \widehat{\mathcal{T}}_i^-(t)\Big)\cap\Big(\widehat{\mathcal{T}}_i^*(t)\setminus \widehat{\mathcal{T}}_i^-(t)\Big)$$
$$\subset \mathcal{A}_i(t-1)\cap\Big(\big(\mathcal{T}_i^{-0}(t)\cup \mathcal{T}_i^{+1}(t)\big)\setminus \overline{\mathcal{T}}_i(t-1)\Big) = \varnothing\,.$$

We are thus in position to use lemma 89.23 to obtain the further inclusion

$$\mathcal{B}\big(w_i, D\setminus \widehat{\mathcal{T}}_i^-(t), t\big) \subset$$
$$\mathcal{B}\Big(\mathcal{B}\big(w_i, D\setminus \widehat{\mathcal{T}}_i^*(t), t-1\big), D\setminus \widehat{\mathcal{T}}_i^*(t), 1\Big)\cup \text{Shell}\,\big(\widehat{\mathcal{T}}_i^*(t)\setminus \widehat{\mathcal{T}}_i^-(t), D\big)$$
$$\subset \mathcal{B}\big(w_i, D\setminus \widehat{\mathcal{T}}_i^*(t), t\big)\cup \text{Shell}\,\big(\mathcal{T}_i^{-0}(t)\cup \mathcal{T}_i^{+1}(t), D\big) \quad (89.137)$$

(we have used again the inclusion (89.135) in the last step). It remains to study the last set of the previous formula. As the sites of $\mathcal{T}_i^{-0}(t)$ are closed, we have

$$\text{Shell}\,\big(\mathcal{T}_i^{-0}(t), D\big) = \mathcal{T}_i^{-0}(t)\,. \quad (89.138)$$

Furthermore, the sites of $\mathcal{T}_i^{+1}(t)$ are visited at time $t\geq 1$, therefore they belong to finite clusters. As the event $\mathcal{H}$ occurs, the finite clusters have strictly less than $M$ open sites, and we have the inclusion

$$\text{Shell}\,\big(\mathcal{T}_i^{+1}(t), D\big) \subset B_1\big(\mathcal{T}_i^{+1}(t), M-1\big) \quad (89.139)$$

(see the explanation after (89.123) for the $M-1$). Putting together the inclusions (89.137), (89.138), (89.139), we obtain the inclusion (89.131). ☐

*Second proof of lemma 89.27.* Let $x$ belong to $\mathcal{B}\big(w_i, D\setminus\widehat{\mathcal{T}}_i^-(t), t\big)$. By definition, there exists a path $y_0,\cdots,y_r$ in $D\setminus \widehat{\mathcal{T}}_i^-(t)$ from $y_0=w_i$ to $y_r=x$ such that the number of closed sites among $y_0,\cdots,y_{r-1}$ is less than or equal to $t$. We consider several cases:

• No site among $y_0,\cdots,y_r$ is in $\widehat{\mathcal{T}}_i^*(t)$. Then $x$ is in $\mathcal{B}\big(w_i, D\setminus \widehat{\mathcal{T}}_i^*(t), t\big)$.

• The path $y_0,\cdots,y_r$ visits $\widehat{\mathcal{T}}_i^*(t)$. We know also that this path does not visit $\widehat{\mathcal{T}}_i^-(t)$. Let $s$ be the smallest index such that $y_s$ is in $\widehat{\mathcal{T}}_i^*(t)$. The site $y_s$ is in

$$\widehat{\mathcal{T}}_i^*(t)\setminus \widehat{\mathcal{T}}_i^-(t) \subset \big(\mathcal{T}_i^{-0}(t)\cup \mathcal{T}_i^{+1}(t)\big)\setminus \overline{\mathcal{T}}_i(t-1)\,.$$

Yet the site $y_s$ is also in $\mathcal{B}\big(w_i, D\setminus \widehat{\mathcal{T}}_i^-(t), t\big)$. By proposition 81.2, we have

$$\big(\mathcal{T}_i^*(t)\setminus\overline{\mathcal{T}}_i(t-1)\big)\cap\mathcal{B}\big(w_i, D\setminus \widehat{\mathcal{T}}_i^-(t), t-1\big)\subset \big(\mathcal{T}_i^*(t)\setminus\overline{\mathcal{T}}_i(t-1)\big)\cap\mathcal{A}_i(t-1) = \varnothing.$$

Thus the site $y_s$ is not in $\mathcal{B}\big(w_i, D\setminus \widehat{\mathcal{T}}_i^-(t), t-1\big)$. This implies that there are at least $t$ closed sites along the subpath $y_0,\cdots,y_{s-1}$. We consider two subcases:

⋆ $s=r$. In this case, we have $y_s=x$, thus $x$ belongs to $\mathcal{T}_i^{-0}(t)\cup \mathcal{T}_i^{+1}(t)$.

⋆ $s<r$. As the path $y_0,\cdots,y_{r-1}$ contains at most $t$ closed sites, the site $y_s$ must be open and it belongs to $\mathcal{T}_i^{+1}(t)$. Furthermore, the sites $y_s, y_{s+1},\dots,y_{r-1}$ are all open. The site $y_s$ is visited at time $t\geq 1$, therefore it belongs to a finite

cluster, and the open path $y_s, y_{s+1}, \ldots, y_{r-1}$ contains strictly less than $M$ open sites, so that $r - s < M$. The subpath $y_s, \cdots, y_r$ is a path testifying that $x$ is in $B_1\big(\mathcal{T}_i^{+1}(t), M-1\big)$.

In each case, the site $x$ belongs to one of the sets appearing in the right-hand side of (89.131). This concludes the proof of the inclusion stated in the lemma. $\square$

In order to obtain a lower bound for the second set appearing in (89.98), we apply lemma 89.27 to each set $\mathcal{B}\big(w_i, D \setminus \widehat{\mathcal{T}}_i^-(t), t\big)$ and we take the union of (89.131) over $i$ in $I$. This yields immediately the following result.

**Corollary 89.28.** *We consider a configuration realizing the event $\mathcal{H}$ defined in* (87.4) *and a subset $I$ of $\{1, \ldots, N\}$. For any $t$ in $\{1, \ldots, T\}$, we have*

$$\bigcup_{i \in I} \mathcal{B}\big(w_i, D \backslash \widehat{\mathcal{T}}_i^-(t), t\big) \subset \Big(\bigcup_{i \in I} \mathcal{B}\big(w_i, D \backslash \widehat{\mathcal{T}}_i^*(t), t\big)\Big) \cup \mathcal{T}^{-0}(t) \cup B_1\big(\mathcal{T}^{+1}(t), M\big)\,.$$

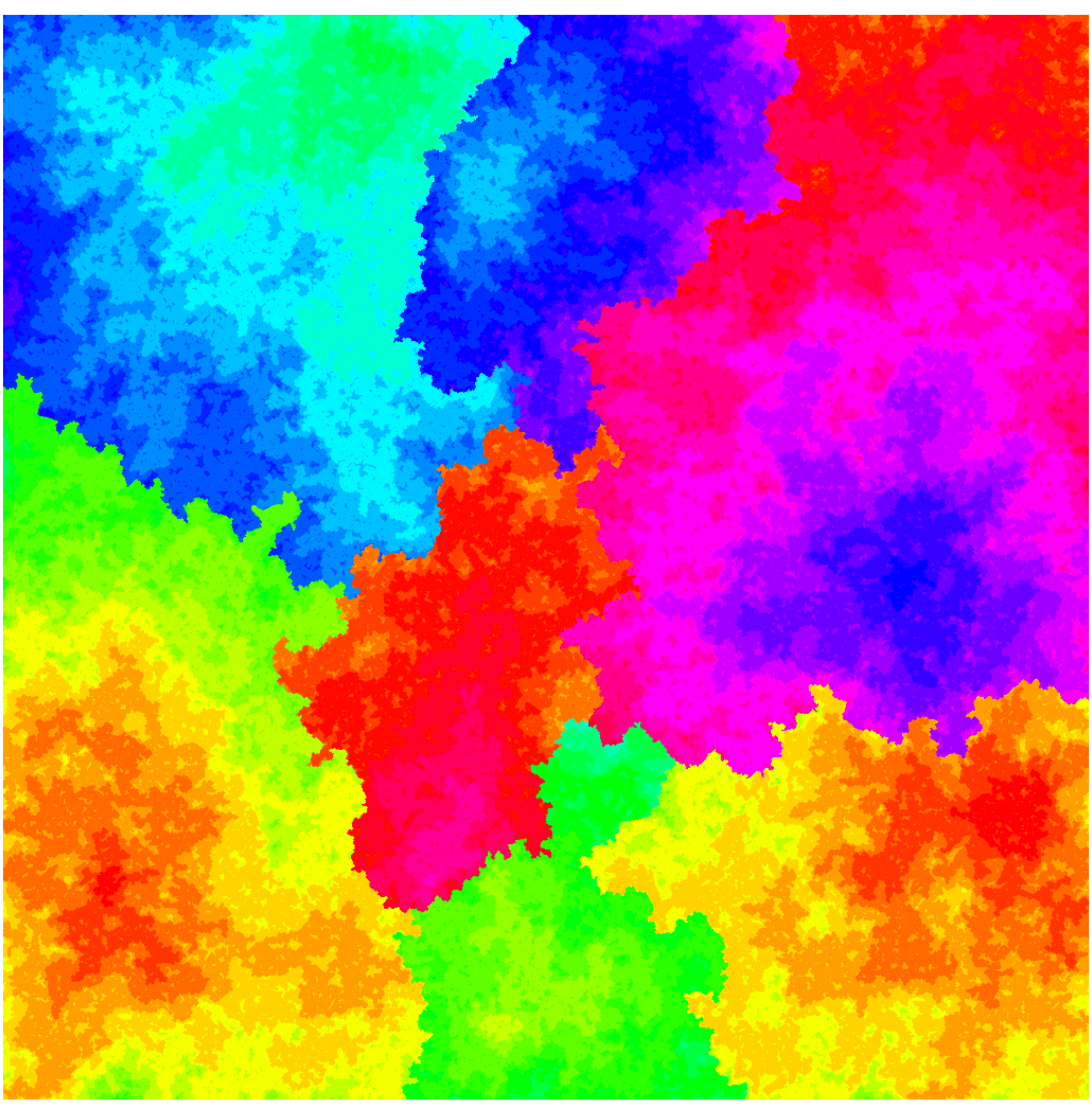

Figure 160: 7 intertwined explorations, $1024 \times 1024$, site percolation, $p = 0.57$.

## 89.7 Ignition of the third stage of the proof

We use now the inclusions obtained in the previous subsections to conclude the second stage of the proof and to ignite the third stage. We come back to the probability we wish to estimate, i.e., the probability appearing inside the sum in (88.37). So we fix four subsets $K^-, L^-, K^+, L^+$ of $\{1, \dots, N\}$ such that $K^- \cap L^- = \varnothing$, $K^+ \cap L^+ = \varnothing$, and we are concerned with estimating the probability

$$P\left(\begin{array}{l} \qquad \mathcal{D},\, \mathcal{E}, \quad \big|\overline{\mathcal{T}}^0(t)\big| > \dfrac{N^*n}{\phi(n)}, \quad \big|\overline{\mathcal{T}}^0(t-1)\big| \leq \rho(n)\,\big|\overline{\mathcal{T}}^0(t)\big|, \\ \dfrac{1}{8d}\,\big|\mathcal{T}^{-0}(t)\big| \leq \big|\overline{\mathcal{T}}(t-1)\big| + \\ \Big|\Big(\bigcup\limits_{k\in K^-} \mathcal{B}\big(w_k, D\setminus\widehat{\mathcal{T}}_k^-(t), t\big)\Big) \cap \Big(\bigcup\limits_{\ell\in L^-} \mathcal{B}\big(w_\ell, D\setminus\overline{\mathcal{T}}_\ell(t-1), t-1\big)\Big)\Big| \\ \dfrac{1}{8d}\,\big|\mathcal{T}^{+0}(t)\big| \leq \\ \quad \Big|\Big(\bigcup_{k\in K^+} \mathcal{B}\big(w_k, D\setminus\widehat{\mathcal{T}}_k^*(t), t\big)\Big) \cap \Big(\bigcup_{\ell\in L^+} \mathcal{B}\big(w_\ell, D\setminus\widehat{\mathcal{T}}_\ell^*(t), t\big)\Big)\Big| \end{array}\right). \tag{89.140}$$

We consider a configuration realizing the event in the probability above and we apply the approximation results to each set appearing inside the probability. We start with the set $\overline{\mathcal{L}}^-$ defined by

$$\overline{\mathcal{L}}^- = \bigcup_{\ell\in L^-} \mathcal{B}\big(w_\ell, D\setminus\overline{\mathcal{T}}_\ell(t-1), t-1\big)\,.$$

Its approximation by truncated genuine explorations is the set $\widehat{\mathcal{L}}_M^-$ defined by

$$\widehat{\mathcal{L}}_M^- = \bigcup_{\ell\in L^-} \widehat{\mathcal{B}}_M\big(\text{Shell}\,(w_\ell, D), D, t-1\big)\,.$$

We apply proposition 89.18 at iteration $t-1$ and we obtain that

$$\overline{\mathcal{L}}^- \subset \widehat{\mathcal{L}}_M^-\,, \qquad \widehat{\mathcal{L}}_M^- \setminus \overline{\mathcal{L}}^- \subset \overline{\text{Error}}(M,t)\,, \tag{89.141}$$

where the error term $\overline{\text{Error}}(M,t)$ is given by

$$\overline{\text{Error}}(M,t) = \bigcup_{0\leq s\leq t-1} \big(B_{1,M}^{s,0}(t-1-s) \cup B_{1,M}^{s,1}(t-s)\big)\,.$$

The important point is that this first error term involves only the relevant taboo sites up to time $t-1$. We continue with the set $\widehat{\mathcal{K}}^-$ defined by

$$\widehat{\mathcal{K}}^- = \bigcup_{k\in K^-} \mathcal{B}\big(w_k, D\setminus\widehat{\mathcal{T}}_k^-(t), t\big)\,.$$

Its approximation by truncated genuine explorations is the set $\widehat{\mathcal{K}}_M^-$ defined by

$$\widehat{\mathcal{K}}_M^- = \bigcup_{k\in K^-} \widehat{\mathcal{B}}_M\big(\text{Shell}\,(w_k, D), D, t\big)\,.$$

With the help of the propositions 89.20 and 89.21, we obtain that

$$\widehat{\mathcal{K}}^- \subset \widehat{\mathcal{K}}^-_M \,, \qquad \widehat{\mathcal{K}}^-_M \setminus \widehat{\mathcal{K}}^- \subset \widehat{\text{Error}}^-(M,t)\,, \tag{89.142}$$

where the error term $\widehat{\text{Error}}^-(M,t)$ is given by

$$\widehat{\text{Error}}^-(M,t) \,=\, \bigcup_{0\leq s\leq t-1} \big(B^{s,0}_{1,M}(t-s) \cup \widehat{\mathcal{B}}^{-,s+1,1}_{1,M}(t-s)\big)\,.$$

We notice that this second error term involves the relevant taboo sites up to time $t-1$, but also the set $\mathcal{T}^{-1}(t)$. The point is that this last set is also controlled by $\overline{\mathcal{T}}^{+0}(t-1)$. We terminate with the two sets $\widehat{\mathcal{K}}^+$ and $\widehat{\mathcal{L}}^+$ defined by

$$\widehat{\mathcal{K}}^+ \,=\, \bigcup_{k\in K^+} \mathcal{B}\big(w_k, D\setminus \widehat{\mathcal{T}}^*_k(t), t\big)\,, \quad \widehat{\mathcal{L}}^+ \,=\, \bigcup_{\ell\in L^+} \mathcal{B}\big(w_\ell, D\setminus \widehat{\mathcal{T}}^*_\ell(t), t\big)\,.$$

Their respective approximations by truncated genuine explorations are the sets $\widehat{\mathcal{K}}^+_M$ and $\widehat{\mathcal{L}}^+_M$ defined by

$$\widehat{\mathcal{K}}^+_M \,=\, \bigcup_{k\in K^+} \widehat{\mathcal{B}}_M\big(\text{Shell}\,(w_k,D), D, t\big)\,,$$
$$\widehat{\mathcal{L}}^+_M \,=\, \bigcup_{\ell\in L^+} \widehat{\mathcal{B}}_M\big(\text{Shell}\,(w_\ell,D), D, t\big)\,.$$

With the help of propositions 89.20, 89.21 and corollary 89.28, we obtain that

$$\widehat{\mathcal{K}}^+ \subset \widehat{\mathcal{K}}^+_M \,, \qquad \widehat{\mathcal{K}}^+_M \setminus \widehat{\mathcal{K}}^+ \subset \widehat{\text{Error}}^+(M,t)\,, \tag{89.143}$$
$$\widehat{\mathcal{L}}^+ \subset \widehat{\mathcal{L}}^+_M \,, \qquad \widehat{\mathcal{L}}^+_M \setminus \widehat{\mathcal{L}}^+ \subset \widehat{\text{Error}}^+(M,t)\,, \tag{89.144}$$

where the error term $\widehat{\text{Error}}^+(M,t)$ is given by

$$\widehat{\text{Error}}^+(M,t) \,=\, \widehat{\text{Error}}^-(M,t) \cup \mathcal{T}^{-0}(t) \cup B_1\big(\mathcal{T}^{+1}(t), M\big)\,. \tag{89.145}$$

We notice that this third error term involves the sets $\mathcal{T}^{-0}(t)$ and $\mathcal{T}^{+1}(t)$, in addition to the previous error term $\widehat{\text{Error}}^-(M,t)$. The point is that $\mathcal{T}^{-0}(t)$ will be controlled before we need to control $\widehat{\text{Error}}^+(M,t)$, and $\mathcal{T}^{+1}(t)$ will be controlled once $\mathcal{T}^{-0}(t)$ is.

Let us try to compute useful upper bounds on the cardinalities of the error sets $\overline{\text{Error}}(M,t)$, $\widehat{\text{Error}}^-(M,t)$, $\widehat{\text{Error}}^+(M,t)$. We recall that, on the event $\mathcal{E}$, the termination time is bounded by $\overline{t}(n) = 2(\ln n)^{1-\gamma}$, so we need only to consider values of $t$ smaller than or equal to $\overline{t}(n)$. We rely on the following crude bounds. We control the cardinality of the balls $B_1$ by

$$\forall x\in\mathbb{Z}^d \quad \forall r\in\mathbb{N} \qquad \big|B_1(x,r)\big| \,\leq\, \Big|\{\,-r,\dots,r\,\}^d\Big| \,\leq\, (2r+1)^d\,,$$

and we apply systematically the previous inequality in the definitions of the sets involved in the error terms. We obtain the following inequalities: for any $s,t$ such that $0 \leq s < t \leq \bar{t}(n)$,

$$\left|B^{s,0}_{1,M}(t-s-1)\right| \leq \left|\mathcal{T}^{*0}(s)\right|\big(2(t-s-1)(M+1)+1\big)^d \leq \left|\mathcal{T}^{*0}(s)\right|\big(5\bar{t}(n)M\big)^d, \tag{89.146}$$

$$\left|B^{s,1}_{1,M}(t-s)\right| \leq \left|\mathcal{T}^{*1}(s)\right|\big(2(t-s)(M+1)+M\big)^d \leq \left|\mathcal{T}^{*1}(s)\right|\big(5\bar{t}(n)M\big)^d. \tag{89.147}$$

The inequality (89.146) yields that

$$\Big|\bigcup_{0\leq s\leq t-1} B^{s,0}_{1,M}(t-s-1)\Big| \leq \big(5\bar{t}(n)M\big)^d \sum_{0\leq s\leq t-1} \left|\mathcal{T}^{*0}(s)\right|. \tag{89.148}$$

Similarly, we deduce from the inequality (89.147) that

$$\Big|\bigcup_{0\leq s\leq t-1} B^{s,1}_{1,M}(t-s)\Big| \leq \big(5\bar{t}(n)M\big)^d \sum_{0\leq s\leq t-1} \left|\mathcal{T}^{*1}(s)\right|. \tag{89.149}$$

Adding together the inequalities (89.148) and (89.149), we have

$$\left|\overline{\mathrm{Error}}(M,t)\right| \leq \big(5\bar{t}(n)M\big)^d \left|\overline{\mathcal{T}}(t-1)\right|. \tag{89.150}$$

We focus next on the second error term $\widehat{\mathrm{Error}}^-(M,t)$. For any $s,t$ such that $0 \leq s < t \leq \bar{t}(n)$, we have

$$\begin{aligned}\left|\widehat{\mathcal{B}}^{-,s+1,1}_{1,M}(t-s)\right| &\leq \left|B_1\big(\overline{\mathcal{T}}^1(s)\cup\mathcal{T}^{-1}(s+1),(t-s)(M+1)+M-1\big)\right| \\ &\quad \Big(\left|\overline{\mathcal{T}}^1(s)\right| + \left|\mathcal{T}^{-1}(s+1)\right|\Big)\big(2(t-s)(M+1)+2M-1\big)^d \\ &\leq \Big(\left|\overline{\mathcal{T}}^1(s)\right| + \left|\mathcal{T}^{-1}(s+1)\right|\Big)\big(9\bar{t}(n)M\big)^d. \end{aligned} \tag{89.151}$$

The inequality (89.151) yields that

$$\Big|\bigcup_{0\leq s\leq t-1} \widehat{\mathcal{B}}^{-,s+1,1}_{1,M}(t-s)\Big| \leq \big(9\bar{t}(n)M\big)^d \sum_{0\leq s\leq t-1} \Big(\left|\overline{\mathcal{T}}^1(s)\right| + \left|\mathcal{T}^{-1}(s+1)\right|\Big). \tag{89.152}$$

We conclude from (89.148) and (89.152) that

$$\begin{aligned}\left|\widehat{\mathrm{Error}}^-(M,t)\right| &\leq \big(9\bar{t}(n)M\big)^d \Big(\left|\overline{\mathcal{T}}(t-1)\right| + \left|\overline{\mathcal{T}}^{-1}(t)\right|\Big) \\ &\leq \big(9\bar{t}(n)M\big)^d \Big(2\left|\overline{\mathcal{T}}(t-1)\right| + \left|\mathcal{T}^{-1}(t)\right|\Big). \end{aligned} \tag{89.153}$$

We finally consider the third error term $\widehat{\mathrm{Error}}^+(M,t)$. It follows immediately from (89.145) that

$$\left|\widehat{\mathrm{Error}}^+(M,t)\right| \leq \left|\widehat{\mathrm{Error}}^-(M,t)\right| + \left|\mathcal{T}^{-0}(t)\right| + \left|\mathcal{T}^{+1}(t)\right|\big(2M+1\big)^d. \tag{89.154}$$

We rewrite next the last two inequalities appearing in the probability (89.140) with the help of our new notation, they become

$$\frac{1}{8d}\left|\mathcal{T}^{-0}(t)\right| \leq \left|\overline{\mathcal{T}}(t-1)\right| + \left|\widehat{\mathcal{K}}^{-}\cap\overline{\mathcal{L}}^{-}\right|, \tag{89.155}$$

$$\frac{1}{8d}\left|\mathcal{T}^{+0}(t)\right| \leq \left|\widehat{\mathcal{K}}^{+}\cap\widehat{\mathcal{L}}^{+}\right|. \tag{89.156}$$

The two sets $\widehat{\mathcal{K}}^{-}\cap\overline{\mathcal{L}}^{-}$ and $\widehat{\mathcal{K}}^{+}\cap\widehat{\mathcal{L}}^{+}$ are obtained as the intersections of two taboo explorations. We try to replace them by their approximations by genuine truncated explorations. Using (89.141) and (89.142), we have

$$\widehat{\mathcal{K}}^{-}\cap\overline{\mathcal{L}}^{-} \subset \widehat{\mathcal{K}}^{-}_M\cap\widehat{\mathcal{L}}^{-}_M, \tag{89.157}$$

$$\big(\widehat{\mathcal{K}}^{-}_M\cap\widehat{\mathcal{L}}^{-}_M\big)\setminus\big(\widehat{\mathcal{K}}^{-}\cap\overline{\mathcal{L}}^{-}\big) \subset \overline{\text{Error}}(M,t)\cup\widehat{\text{Error}}^{-}(M,t). \tag{89.158}$$

Similarly, using (89.143) and (89.144), we have

$$\widehat{\mathcal{K}}^{+}\cap\widehat{\mathcal{L}}^{+} \subset \widehat{\mathcal{K}}^{+}_M\cap\widehat{\mathcal{L}}^{+}_M, \tag{89.159}$$

$$\big(\widehat{\mathcal{K}}^{+}_M\cap\widehat{\mathcal{L}}^{+}_M\big)\setminus\big(\widehat{\mathcal{K}}^{+}\cap\widehat{\mathcal{L}}^{+}\big) \subset \widehat{\text{Error}}^{+}(M,t). \tag{89.160}$$

The inclusions (89.157), (89.159) and the inequalities (89.155), (89.156) yield that

$$\frac{1}{8d}\left|\mathcal{T}^{-0}(t)\right| \leq \left|\overline{\mathcal{T}}(t-1)\right| + \left|\widehat{\mathcal{K}}^{-}_M\cap\widehat{\mathcal{L}}^{-}_M\right|, \tag{89.161}$$

$$\frac{1}{8d}\left|\mathcal{T}^{+0}(t)\right| \leq \left|\widehat{\mathcal{K}}^{+}_M\cap\widehat{\mathcal{L}}^{+}_M\right|. \tag{89.162}$$

Before proceeding further, we recall a crucial point. It was proved in lemma 88.3 that all the sites in the sets $\widehat{\mathcal{K}}^{-}\cap\overline{\mathcal{L}}^{-}$, $\widehat{\mathcal{K}}^{+}\cap\widehat{\mathcal{L}}^{+}$ are closed. Using the inclusions (89.158), (89.160), we can thus bound from above the numbers of open sites in $\widehat{\mathcal{K}}^{-}_M\cap\widehat{\mathcal{L}}^{-}_M$, $\widehat{\mathcal{K}}^{+}_M\cap\widehat{\mathcal{L}}^{+}_M$ as follows:

$$\begin{aligned}\Big|\big\{\, x\in\widehat{\mathcal{K}}^{-}_M\cap\widehat{\mathcal{L}}^{-}_M : x \text{ is open}\,\big\}\Big| &\leq \Big|\big(\widehat{\mathcal{K}}^{-}_M\cap\widehat{\mathcal{L}}^{-}_M\big)\setminus\big(\widehat{\mathcal{K}}^{-}\cap\overline{\mathcal{L}}^{-}\big)\Big| \\ &\leq \left|\overline{\text{Error}}(M,t)\cup\widehat{\text{Error}}^{-}(M,t)\right|,\end{aligned} \tag{89.163}$$

$$\begin{aligned}\Big|\big\{\, x\in\widehat{\mathcal{K}}^{+}_M\cap\widehat{\mathcal{L}}^{+}_M : x \text{ is open}\,\big\}\Big| &\leq \Big|\big(\widehat{\mathcal{K}}^{+}_M\cap\widehat{\mathcal{L}}^{+}_M\big)\setminus\big(\widehat{\mathcal{K}}^{+}\cap\widehat{\mathcal{L}}^{+}\big)\Big| \\ &\leq \left|\widehat{\text{Error}}^{+}(M,t)\right|.\end{aligned} \tag{89.164}$$

We use next the controls on the error terms (89.150), (89.153), (89.154) and we obtain

$$\left|\overline{\text{Error}}(M,t)\cup\widehat{\text{Error}}^{-}(M,t)\right| \leq \big(9\bar{t}(n)M\big)^d\Big(3\left|\overline{\mathcal{T}}(t-1)\right| + \left|\mathcal{T}^{-1}(t)\right|\Big), \tag{89.165}$$

$$\begin{aligned}&\left|\widehat{\text{Error}}^{+}(M,t)\right| \leq \\ &\big(9\bar{t}(n)M\big)^d\Big(2\left|\overline{\mathcal{T}}(t-1)\right| + \left|\mathcal{T}^{-1}(t)\right|\Big) + \left|\mathcal{T}^{-0}(t)\right| + \left|\mathcal{T}^{+1}(t)\right|(3M)^d.\end{aligned} \tag{89.166}$$

We want to obtain inequalities which involve only the sets $\overline{\mathcal{T}}^0(t-1)$, $\mathcal{T}^{-0}(t)$ and $\mathcal{T}^{+0}(t)$. In order to rework the inequalities (89.165) and (89.166), we make appeal to the inequalities (88.2), which we recall for convenience, but with a variable $s$ instead of $t$:

$$\forall s\in\{\,1,\dots,t\,\}\qquad \big|\mathcal{T}^{-1}(s)\big|\,\leq\,2d\,\big|\mathcal{T}^{+0}(s-1)\big|\,,\tag{89.167}$$

$$\forall s\in\{\,1,\dots,t\,\}\qquad \big|\mathcal{T}^{+1}(s)\big|\,\leq\,2d\,\big|\mathcal{T}^{-0}(s)\big|\,.\tag{89.168}$$

We get rid of the term $\big|\mathcal{T}^{-1}(t)\big|$ in (89.165) by applying the inequality (89.167) with $s=t$, and we obtain

$$\big|\mathcal{T}^{-1}(t)\big|\,\leq\,2d\,\big|\mathcal{T}^{+0}(t-1)\big|\,\leq\,2d\,\big|\overline{\mathcal{T}}(t-1)\big|\,.\tag{89.169}$$

For the term $\big|\mathcal{T}^{+1}(t)\big|$, we apply the inequality (89.168) with $s=t$. Notice that the effect of each inequality is different. While the term $\big|\mathcal{T}^{-1}(t)\big|$ will merge with $\big|\overline{\mathcal{T}}(t-1)\big|$ in (89.165), the term $\big|\mathcal{T}^{-0}(t)\big|$ remains present in (89.166). We deal next with $\big|\overline{\mathcal{T}}(t-1)\big|$. For $t=1$, on the event $\mathcal{D}$, we have $\overline{\mathcal{T}}(0)=\overline{\mathcal{T}}^0(0)$ and $\overline{\mathcal{T}}^1(0)=\varnothing$. For $t\geq 2$, summing (89.167) and (89.168) from 1 to $t-1$, we obtain

$$\begin{aligned}\big|\overline{\mathcal{T}}^1(t-1)\big|\,&\leq\,\sum_{1\leq s\leq t-1}\Big(\big|\mathcal{T}^{-1}(s)\big|+\big|\mathcal{T}^{+1}(s)\big|\Big)\\&\leq\,2d\sum_{1\leq s\leq t-1}\Big(\big|\mathcal{T}^{+0}(s-1)\big|+\big|\mathcal{T}^{-0}(s)\big|\Big)\,\leq\,4d\,\big|\overline{\mathcal{T}}^0(t-1)\big|\,,\end{aligned}$$

whence

$$\big|\overline{\mathcal{T}}(t-1)\big|\,=\,\big|\overline{\mathcal{T}}^0(t-1)\big|+\big|\overline{\mathcal{T}}^1(t-1)\big|\,\leq\,(4d+1)\,\big|\overline{\mathcal{T}}^0(t-1)\big|\,.\tag{89.170}$$

Using the inequalities (89.169), (89.168) with $s=t$, and (89.170), we deduce from the two inequalities (89.165), (89.166) that

$$\big|\overline{\mathrm{Error}}(M,t)\cup\widehat{\mathrm{Error}}^-(M,t)\big|\,\leq\,28d^2\big(9\overline{t}(n)M\big)^d\big|\overline{\mathcal{T}}^0(t-1)\big|\,,\tag{89.171}$$

$$\big|\widehat{\mathrm{Error}}^+(M,t)\big|\,\leq\,30d^2\big(9\overline{t}(n)M\big)^d\big(\big|\overline{\mathcal{T}}^0(t-1)\big|+\big|\mathcal{T}^{-0}(t)\big|\big)\,.\tag{89.172}$$

We use next all the inequalities (89.161), (89.162), (89.163), (89.164), (89.171), (89.172) and we conclude that

$$\Big|\big\{\,x\in\widehat{\mathcal{K}}_M^-\cap\widehat{\mathcal{L}}_M^-:x\text{ is open}\,\big\}\Big|\,\leq\,30d^2\big(9\overline{t}(n)M\big)^d\big|\overline{\mathcal{T}}^0(t-1)\big|\,,\tag{89.173}$$

$$\Big|\big\{\,x\in\widehat{\mathcal{K}}_M^+\cap\widehat{\mathcal{L}}_M^+:x\text{ is open}\,\big\}\Big|\,\leq\,30d^2\big(9\overline{t}(n)M\big)^d\big(\big|\overline{\mathcal{T}}^0(t-1)\big|+\big|\mathcal{T}^{-0}(t)\big|\big)\,.\tag{89.174}$$

With the help of the inequalities (89.161), (89.162), (89.170), (89.173), (89.174), we bound the probability (89.140) as follows:

$$P\left(\begin{array}{c}\mathcal{D},\,\mathcal{E},\quad \left|\overline{\mathcal{T}}^0(t)\right| > \dfrac{N^*n}{\phi(n)}, \quad \left|\overline{\mathcal{T}}^0(t-1)\right| \leq \rho(n)\left|\overline{\mathcal{T}}^0(t)\right|, \\ \dfrac{1}{8d}\left|\mathcal{T}^{-0}(t)\right| \leq \left|\overline{\mathcal{T}}(t-1)\right| + \\ \left|\Big(\bigcup_{k\in K^-}\mathcal{B}\big(w_k, D\setminus\widehat{\mathcal{T}}_k^-(t),t\big)\Big)\cap\Big(\bigcup_{\ell\in L^-}\mathcal{B}\big(w_\ell, D\setminus\overline{\mathcal{T}}_\ell(t-1),t-1\big)\Big)\right| \\ \dfrac{1}{8d}\left|\mathcal{T}^{+0}(t)\right| \leq \\ \left|\Big(\bigcup_{k\in K^+}\mathcal{B}\big(w_k, D\setminus\widehat{\mathcal{T}}_k^*(t),t\big)\Big)\cap\Big(\bigcup_{\ell\in L^+}\mathcal{B}\big(w_\ell, D\setminus\widehat{\mathcal{T}}_\ell^*(t),t\big)\Big)\right|\end{array}\right)$$

$$\leq \tag{89.175}$$

$$P\left(\begin{array}{c}\mathcal{D},\,\mathcal{E},\quad \left|\overline{\mathcal{T}}^0(t)\right| > \dfrac{N^*n}{\phi(n)}, \quad \left|\overline{\mathcal{T}}^0(t-1)\right| \leq \rho(n)\left|\overline{\mathcal{T}}^0(t)\right|, \\ \dfrac{1}{8d}\left|\mathcal{T}^{-0}(t)\right| \leq 5d\left|\overline{\mathcal{T}}^0(t-1)\right| + \left|\widehat{\mathcal{K}}_M^-\cap\widehat{\mathcal{L}}_M^-\right|, \\ \dfrac{1}{8d}\left|\mathcal{T}^{+0}(t)\right| \leq \left|\widehat{\mathcal{K}}_M^+\cap\widehat{\mathcal{L}}_M^+\right|, \\ \left|\left\{\,x\in\widehat{\mathcal{K}}_M^-\cap\widehat{\mathcal{L}}_M^- : x \text{ is open}\,\right\}\right| \leq 30d^2\big(9\bar{t}(n)M\big)^d\left|\overline{\mathcal{T}}^0(t-1)\right|, \\ \left|\left\{\,x\in\widehat{\mathcal{K}}_M^+\cap\widehat{\mathcal{L}}_M^+ : x \text{ is open}\,\right\}\right| \leq \\ 30d^2\big(9\bar{t}(n)M\big)^d\big(\left|\overline{\mathcal{T}}^0(t-1)\right| + \left|\mathcal{T}^{-0}(t)\right|\big)\end{array}\right).$$

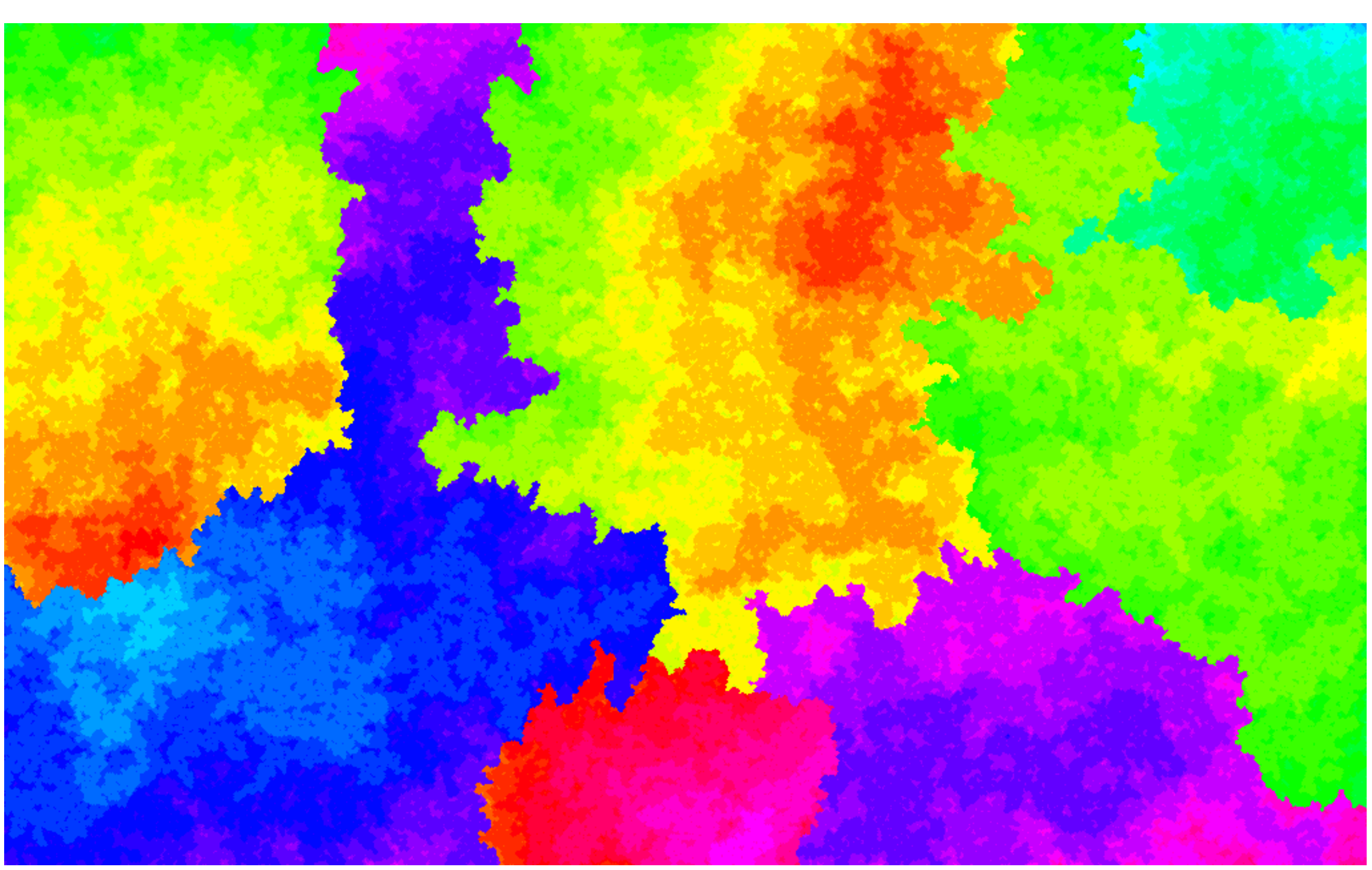

Figure 161: 7 intertwined explorations, 1024 × 633, site percolation, $p = 0.57$.

# 90 Control of the relevant taboo sites

In this section, we shall derive probabilistic estimates on the number of the relevant taboo sites. These estimates rely on the approximations of the set of the relevant taboo sites by genuine explorations carried out in section 89. For genuine explorations, we employ the usual probabilistic control. In subsection 90.1, we deal with the first taboo set $\mathcal{T}^{*0}(0)$, which plays a special role. The higher taboo sets $(\overline{\mathcal{T}}^0(t), t \geq 1)$ are treated in subsections 90.2 and 90.3. We explain how these sets can be approximately localized with the help of two genuine explorations in subsection 90.2, and we derive then a quantitative control in subsection 90.3.

## 90.1 Control of $\mathcal{T}^{*0}(0)$

The set $\mathcal{T}^{*0}(0)$ coincides with $\overline{\mathcal{T}}(0) = \mathcal{T}^*(0)$, in fact it consists precisely of the pivotal sites for $\mathcal{D}$. Since the event $\mathcal{D}$ is decreasing, we can use the general results on pivots. Rewriting the equality (40.2) in our context, we get

$$P\big(\, x \in \mathcal{P}^-(\mathcal{D})\,\big) \,=\, \frac{1}{1-p} P\big(\{\, x \in \mathcal{P}^-(\mathcal{D})\,\} \cap \mathcal{D}\big)\,. \tag{90.1}$$

We sum next (90.1) over all the sites of $D$, and we apply the general inequality (39.10) to the decreasing event $\mathcal{D}$ to conclude that

$$E\big(|\mathcal{T}^{*0}(0)|; \mathcal{D}\big) \,=\, E\big(|\mathcal{P}^-(\mathcal{D})|; \mathcal{D}\big) \,\leq\, \sqrt{\frac{1-p}{p}|D|P(\mathcal{D})} \,\leq\, \frac{1}{p}\sqrt{|D|P(\mathcal{D})}\,. \tag{90.2}$$

This is a good start. Indeed, in the case of a box $\Lambda(n)$, the term $\sqrt{|D|}$ will be of order $n^{d/2}$, thus any event forcing the presence of a number of pivots of a higher order of magnitude is deemed to be unlikely. Rather than the inequality (90.2) on the expectations, we will need a deviations inequality in order to quantify how unlikely the events become when the intersection sets are too large. The first natural attempt is to try with the exponential inequality of proposition 40.2, and to proceed in a way similar to what we did in part VI for the three challenging questions raised in section 30. Unfortunately, we face an additional difficulty, because we are trying to control events which are not purely monotone, but which are the intersection of increasing and decreasing events. Indeed, the event $\mathcal{D} \cap \mathcal{E}$ is

$$\mathcal{D} \cap \mathcal{E} \,=\, \mathcal{D} \cap \mathcal{F}_{\text{fi}}(n, M) \cap \mathcal{G}(n) \cap \mathcal{F}_{\text{bd}}(n, \varepsilon)\,.$$

In the above intersection, the first two events $\mathcal{D}$, $\mathcal{F}_{\text{fi}}(n, M)$ are decreasing, while the last two events $\mathcal{G}(n)$ and $\mathcal{F}_{\text{bd}}(n, \varepsilon)$ are increasing. In order to derive a control on $P(\mathcal{D} \cap \mathcal{E})$ from an upper bound on $E\big(|\mathcal{T}^{*0}(0)|; \mathcal{D} \cap \mathcal{E}\big)$, we would need exponential inequalities that are much more complex than those developed throughout part VI. Since these inequalities were ultimately based on Hoeffding's inequality, we prefer to use a more elementary approach, which, however, has the disadvantage of only working for the specific event $\mathcal{D}\cap\mathcal{E}$. This approach is entirely based on the essential fact that the intersection sets can be revealed with the help of the multiple intertwined explorations.

To alleviate the formulas, we define

$$\forall I\subset\{\,1,\dots,N\,\}\qquad \mathcal{E}(I)\,=\,\overline{\text{Clusters}}\Big(\{\,w_i,i\in I\,\},D\Big)\,=\,\bigcup_{i\in I}\overline{C}(w_i,D)\,.\tag{90.3}$$

For any deterministic set $I$, the set $\mathcal{E}(I)$ is the set of the sites visited by a genuine exploration algorithm starting from $\{\,w_i,i\in I\,\}$. On the event $\mathcal{D}$, we have that, for any pair $K,L$ of disjoint subsets of $\{\,1,\dots,N\,\}$,

$$\mathcal{E}(K)\cap\mathcal{E}(L)\,=\,\Big(\bigcup_{k\in K}\partial^{\,out}C(w_k,D)\Big)\cap\Big(\bigcup_{\ell\in L}\partial^{\,out}C(w_\ell,D)\Big)\,.\tag{90.4}$$

Let $i_0$ be a positive integer. Putting together (88.1) and (90.4), we obtain

$$\begin{aligned}P\big(\big|\mathcal{T}^{*0}(0)\big|\geq i_0,\mathcal{D},\mathcal{E}\big)\,&\leq\,P\begin{pmatrix}\exists K,L\subset\{\,1,\dots,N\,\}\quad K\cap L=\varnothing\,,\\ \big|\mathcal{E}(K)\cap\mathcal{E}(L)\big|\geq i_0/16d^2\,,\ \mathcal{D}\,,\ \mathcal{E}\end{pmatrix}\\ &\leq\sum_{\substack{K,L\subset\{\,1,\dots,N\,\}\\ K\cap L=\varnothing}}P\Big(\big|\mathcal{E}(K)\cap\mathcal{E}(L)\big|\geq\frac{i_0}{16d^2},\,\mathcal{D}\Big)\,.\end{aligned}\tag{90.5}$$

Beware that there is an unfortunate collision with the notations: the event $\mathcal{E}$ was defined in (87.1), while the set $\mathcal{E}(I)$ is defined in (90.3). Naturally, we have

$$\mathcal{E}(K)\cup\mathcal{E}(L)\,=\,\mathcal{E}(K\cup L)\,.\tag{90.6}$$

As usual, the estimates are based on the functional $S$, originally introduced in (13.4), that we recall for convenience: for $A$ a finite subset of $D$, we define

$$S(A)\,=\,\frac{1}{p}\Big|\{\,x\in A:x\text{ is open}\,\}\Big|-\frac{1}{1-p}\Big|\{\,x\in A:x\text{ is closed}\,\}\Big|\,.\tag{90.7}$$

The function $S$ is additive, thus we have

$$S\big(\mathcal{E}(K\cup L)\big)+S\big(\mathcal{E}(K)\cap\mathcal{E}(L)\big)\,=\,S\big(\mathcal{E}(K)\big)+S\big(\mathcal{E}(L)\big)\,.\tag{90.8}$$

On the event $\mathcal{D}$, the sites of $\mathcal{E}(K)\cap\mathcal{E}(L)$ are all closed, therefore

$$S\big(\mathcal{E}(K)\cap\mathcal{E}(L)\big)\,=\,-\frac{1}{1-p}\big|\mathcal{E}(K)\cap\mathcal{E}(L)\big|\,.\tag{90.9}$$

Using the probabilistic control (13.5) on the three explorations starting from $K$, $L$ and $K\cup L$, the equalities (90.6), (90.8) and (90.9), we obtain that

$$P\Big(\big|\mathcal{E}(K)\cap\mathcal{E}(L)\big|\geq\frac{i_0}{16d^2},\,\mathcal{D}\Big)\,\leq\,6|\Lambda(4n)|\exp\Big(-\frac{(pi_0)^2}{1152d^4|\Lambda(4n)|}\Big)\,.\tag{90.10}$$

We substitute (90.10) into (90.5). This yields

$$\forall i_0\geq 1\qquad P\big(\big|\mathcal{T}^{*0}(0)\big|\geq i_0,\mathcal{D},\mathcal{E}\big)\,\leq\,6\,2^{2N}\,|\Lambda(4n)|\exp\Big(-\frac{(pi_0)^2}{1152d^4|\Lambda(4n)|}\Big)\,.\tag{90.11}$$

## 90.2 Reduction to two genuine explorations

To alleviate the formulas, we introduce the notation

$$\tau(n) \;=\; 30d^2\big(9\bar{t}(n)M\big)^d\,.$$

Using the definitions (87.14) and (85.4) of $\bar{t}(n)$ and $M$, we have

$$\begin{aligned}\tau(n) \;&=\; 30d^2\Big(18(\ln n)^{1-\gamma}(4d\ln n)^{1/\alpha}\Big)^d\\ &=\; 30d^2 18^d(4d)^{d/\alpha}(\ln n)^{d-d\gamma+d/\alpha}\,.\end{aligned}\tag{90.12}$$

We increase slightly the constants in the inequalities involved in the event appearing in the second probability of (89.175). More precisely, we replace each occurrence of a term of the form constant $\times\, d^2\big(9\bar{t}(n)M\big)^d$ by the function $\tau(n)$. The event appearing in the second probability of (89.175) is included in the event

$$\left\{\begin{matrix}\mathcal{D},\,\mathcal{E},\quad \big|\overline{\mathcal{T}}^0(t)\big| > \dfrac{N^*n}{\phi(n)},\quad \big|\overline{\mathcal{T}}^0(t-1)\big| \,\leq\, \rho(n)\,\big|\overline{\mathcal{T}}^0(t)\big|\,,\\ \dfrac{1}{8d}\,\big|\mathcal{T}^{-0}(t)\big| \,\leq\, \big|\widehat{\mathcal{K}}^-_M\cap\widehat{\mathcal{L}}^-_M\big| + 5d\big|\overline{\mathcal{T}}^0(t-1)\big|\,,\\ \dfrac{1}{8d}\,\big|\mathcal{T}^{+0}(t)\big| \,\leq\, \big|\widehat{\mathcal{K}}^+_M\cap\widehat{\mathcal{L}}^+_M\big|\,,\\ \Big|\big\{\,x\in\widehat{\mathcal{K}}^-_M\cap\widehat{\mathcal{L}}^-_M : x\text{ is open}\,\big\}\Big| \leq \tau(n)\big|\overline{\mathcal{T}}^0(t-1)\big|\,,\\ \Big|\big\{\,x\in\widehat{\mathcal{K}}^+_M\cap\widehat{\mathcal{L}}^+_M : x\text{ is open}\,\big\}\Big| \leq \tau(n)\big(\big|\overline{\mathcal{T}}^0(t-1)\big| + \big|\mathcal{T}^{-0}(t)\big|\big)\end{matrix}\right\}.$$

Decomposing $\overline{\mathcal{T}}^0(t)$ as

$$\overline{\mathcal{T}}^0(t) \;=\; \mathcal{T}^{-0}(t)\cup\mathcal{T}^{+0}(t)\cup\overline{\mathcal{T}}^0(t-1)\,,$$

we obtain the inequality

$$\big|\overline{\mathcal{T}}^0(t)\big| \;\leq\; \big|\mathcal{T}^{-0}(t)\big| + \big|\mathcal{T}^{+0}(t)\big| + \big|\overline{\mathcal{T}}^0(t-1)\big|\,.\tag{90.13}$$

Recall that $\rho(n)\leq 1/2$. We suppose also that $n\geq 3$, so that $\bar{t}(n)\geq 2$ and

$$\tau(n) \;\geq\; 30d^2 9^d\,.\tag{90.14}$$

Using (90.13), we have the following inclusions of events:

$$\begin{aligned}&\Big\{\,\big|\overline{\mathcal{T}}^0(t-1)\big| \leq \rho(n)\,\big|\overline{\mathcal{T}}^0(t)\big|\,\Big\}\subset\\ &\qquad\qquad\Big\{\,\big|\mathcal{T}^{-0}(t)\big| + \big|\mathcal{T}^{+0}(t)\big| \geq \frac{1}{2}\,\big|\overline{\mathcal{T}}^0(t)\big|\,\Big\}\subset\\ &\Big\{\,\big|\mathcal{T}^{-0}(t)\big| > \frac{1}{\tau(n)^2}\,\big|\overline{\mathcal{T}}^0(t)\big|\,\Big\}\cup\left\{\begin{matrix}\big|\mathcal{T}^{-0}(t)\big| \leq \dfrac{1}{\tau(n)^2}\,\big|\overline{\mathcal{T}}^0(t)\big|\,,\\ \big|\mathcal{T}^{+0}(t)\big| > \dfrac{1}{4}\,\big|\overline{\mathcal{T}}^0(t)\big|\end{matrix}\right\}.\end{aligned}\tag{90.15}$$

With the help of the inclusion (90.15), we split in two parts the last probability in (89.175), each one involving only two genuine explorations, as follows:

$$
P\begin{pmatrix}
\mathcal{D},\,\mathcal{E},\quad \big|\overline{\mathcal{T}}^0(t)\big| > \dfrac{N^*n}{\phi(n)},\quad \big|\overline{\mathcal{T}}^0(t-1)\big| \le \rho(n)\,\big|\overline{\mathcal{T}}^0(t)\big|\,, \\
\dfrac{1}{8d}\,\big|\mathcal{T}^{-0}(t)\big| \le \big|\widehat{\mathcal{K}}_M^- \cap \widehat{\mathcal{L}}_M^-\big| + 5d\big|\overline{\mathcal{T}}^0(t-1)\big|\,, \\
\dfrac{1}{8d}\,\big|\mathcal{T}^{+0}(t)\big| \le \big|\widehat{\mathcal{K}}_M^+ \cap \widehat{\mathcal{L}}_M^+\big|\,, \\
\Big|\big\{\, x\in \widehat{\mathcal{K}}_M^- \cap \widehat{\mathcal{L}}_M^- : x \text{ is open}\,\big\}\Big| \le \tau(n)\big|\overline{\mathcal{T}}^0(t-1)\big|\,, \\
\Big|\big\{\, x\in \widehat{\mathcal{K}}_M^+ \cap \widehat{\mathcal{L}}_M^+ : x \text{ is open}\,\big\}\Big| \le \tau(n)\big(\big|\overline{\mathcal{T}}^0(t-1)\big| + \big|\mathcal{T}^{-0}(t)\big|\big)
\end{pmatrix} \tag{90.16}
$$

$$
\le
$$

$$
P\begin{pmatrix}
\mathcal{D},\,\mathcal{E},\quad \big|\overline{\mathcal{T}}^0(t)\big| > \dfrac{N^*n}{\phi(n)},\quad \big|\overline{\mathcal{T}}^0(t-1)\big| \le \rho(n)\,\big|\overline{\mathcal{T}}^0(t)\big|\,, \\
\big|\mathcal{T}^{-0}(t)\big| > \dfrac{1}{\tau(n)^2}\,\big|\overline{\mathcal{T}}^0(t)\big|\,, \\
\dfrac{1}{8d}\,\big|\mathcal{T}^{-0}(t)\big| \le \big|\widehat{\mathcal{K}}_M^- \cap \widehat{\mathcal{L}}_M^-\big| + 5d\big|\overline{\mathcal{T}}^0(t-1)\big|\,, \\
\Big|\big\{\, x\in \widehat{\mathcal{K}}_M^- \cap \widehat{\mathcal{L}}_M^- : x \text{ is open}\,\big\}\Big| \le \tau(n)\big|\overline{\mathcal{T}}^0(t-1)\big|
\end{pmatrix}
$$

$$
+
$$

$$
P\begin{pmatrix}
\mathcal{D},\,\mathcal{E},\, \big|\overline{\mathcal{T}}^0(t)\big| > \dfrac{N^*n}{\phi(n)},\, \big|\overline{\mathcal{T}}^0(t-1)\big| \le \rho(n)\,\big|\overline{\mathcal{T}}^0(t)\big|\,, \\
\big|\mathcal{T}^{-0}(t)\big| \le \dfrac{1}{\tau(n)^2}\,\big|\overline{\mathcal{T}}^0(t)\big|,\, \big|\mathcal{T}^{+0}(t)\big| > \dfrac{1}{4}\big|\overline{\mathcal{T}}^0(t)\big|\,, \\
\dfrac{1}{8d}\,\big|\mathcal{T}^{+0}(t)\big| \le \big|\widehat{\mathcal{K}}_M^+ \cap \widehat{\mathcal{L}}_M^+\big|\,, \\
\Big|\big\{\, x\in \widehat{\mathcal{K}}_M^+ \cap \widehat{\mathcal{L}}_M^+ : x \text{ is open}\,\big\}\Big| \le \tau(n)\big(\big|\overline{\mathcal{T}}^0(t-1)\big| + \big|\mathcal{T}^{-0}(t)\big|\big)
\end{pmatrix}
$$

$$
\le
$$

$$
P\begin{pmatrix}
\mathcal{D},\,\mathcal{E},\quad \big|\overline{\mathcal{T}}^0(t)\big| > \dfrac{N^*n}{\phi(n)}\,, \\
\dfrac{1}{8d\tau(n)^2}\,\big|\overline{\mathcal{T}}^0(t)\big| \le \big|\widehat{\mathcal{K}}_M^- \cap \widehat{\mathcal{L}}_M^-\big| + 5d\rho(n)\big|\overline{\mathcal{T}}^0(t)\big|\,, \\
\Big|\big\{\, x\in \widehat{\mathcal{K}}_M^- \cap \widehat{\mathcal{L}}_M^- : x \text{ is open}\,\big\}\Big| \le \tau(n)\rho(n)\big|\overline{\mathcal{T}}^0(t)\big|
\end{pmatrix} \tag{90.17}
$$

$$
+
$$

$$
P\begin{pmatrix}
\mathcal{D},\,\mathcal{E},\, \big|\overline{\mathcal{T}}^0(t)\big| > \dfrac{N^*n}{\phi(n)}\,, \\
\dfrac{1}{32d}\big|\overline{\mathcal{T}}^0(t)\big| \le \big|\widehat{\mathcal{K}}_M^+ \cap \widehat{\mathcal{L}}_M^+\big|\,, \\
\Big|\big\{\, x\in \widehat{\mathcal{K}}_M^+ \cap \widehat{\mathcal{L}}_M^+ : x \text{ is open}\,\big\}\Big| \le \Big(\tau(n)\rho(n) + \dfrac{1}{\tau(n)}\Big)\big|\overline{\mathcal{T}}^0(t)\big|
\end{pmatrix}. \tag{90.18}
$$

To get the last inequality, we replaced the inequalities inside the probabilities by weaker ones, which have the advantage of involving only the set $\mathcal{T}^0(t)$. The events in the probabilities (90.17) and (90.18) provide a lower bound on the cardinalities of the two random sets $\widehat{\mathcal{K}}_M^- \cap \widehat{\mathcal{L}}_M^-$, $\widehat{\mathcal{K}}_M^+ \cap \widehat{\mathcal{L}}_M^+$, as well as an upper bound on the number of open sites they contain. We will convert this information into a bound on the value of the functional $S$ on these sets. Instead of doing twice the same computation, we record it into the next lemma.

**Lemma 90.1.** *For any subset $A$ of $D$, we have*

$$S(A) \;=\; \frac{1}{p(1-p)}\Big(\big|\{\,x\in A: x \textit{ is open}\,\}\big| - p|A|\Big)\,. \tag{90.19}$$

*Proof.* We simply use the definition (90.7) of $S$, and we compute $S$ as follows:

$$\begin{aligned} S(A) \;&=\; \frac{1}{p}\big|\{\,x\in A: x \text{ is open}\,\}\big| - \frac{1}{1-p}\big|\{\,x\in A: x \text{ is closed}\,\}\big| \\ &= \frac{1}{p}\big|\{\,x\in A: x \text{ is open}\,\}\big| - \frac{1}{1-p}\Big(|A| - \big|\{\,x\in A: x \text{ is open}\,\}\big|\Big) \\ &= \frac{1}{p(1-p)}\Big(\big|\{\,x\in A: x \text{ is open}\,\}\big| - p|A|\Big)\,. \end{aligned}$$

This is precisely the formula (90.19). □

We use next the inequalities on the cardinalities of the sets and the number of open sites to bound from above the probabilities (90.17) and (90.18) by the probabilities of two events involving the functional $S$. The probability (90.17) is less or equal than

$$P\left(\begin{gathered} \mathcal{D},\,\mathcal{E},\quad \big|\overline{\mathcal{T}}^0(t)\big| \;>\; \frac{N^*n}{\phi(n)}\,, \\ S\big(\widehat{\mathcal{K}}_M^- \cap \widehat{\mathcal{L}}_M^-\big) \;\leq\; \frac{1}{p(1-p)}\left(\tau(n)\rho(n) - p\Big(\frac{1}{8d\tau(n)^2} - 5d\rho(n)\Big)\right)\big|\overline{\mathcal{T}}^0(t)\big| \end{gathered}\right), \tag{90.20}$$

and the probability (90.18) is less or equal than

$$P\left(\begin{gathered} \mathcal{D},\,\mathcal{E},\, \big|\overline{\mathcal{T}}^0(t)\big| > \frac{N^*n}{\phi(n)}\,, \\ S\big(\widehat{\mathcal{K}}_M^+ \cap \widehat{\mathcal{L}}_M^+\big) \;\leq\; \frac{1}{p(1-p)}\Big(\tau(n)\rho(n) + \frac{1}{\tau(n)} - \frac{p}{32d}\Big)\big|\overline{\mathcal{T}}^0(t)\big| \end{gathered}\right). \tag{90.21}$$

Both probabilities (90.20) and (90.21) involve two genuine explorations, and at first sight they look similar. There is a tiny difference in fact, because among the explorations, we have one which requires $t-1$ iterations, $\widehat{\mathcal{L}}_M^-$, while the three others, $\widehat{\mathcal{L}}_M^+$, $\widehat{\mathcal{K}}_M^-$, $\widehat{\mathcal{K}}_M^+$ perform $t$ iterations. However, this detail does not affect the estimates to be made in any way, so we treat both probabilities simultaneously. We simply remove the sign $-$ or $+$ in the superscript, and we consider two random sets $\widehat{\mathcal{K}}_M$ and $\widehat{\mathcal{L}}_M$ that are obtained as the output of two genuine truncated explorations having two disjoint starting sets, i.e., of the form

$$\widehat{\mathcal{K}}_M \;=\; \bigcup_{k\in K}\widehat{\mathcal{B}}_M\big(\text{Shell}\,(w_k,D),D,t\big)\,,\quad \widehat{\mathcal{L}}_M \;=\; \bigcup_{\ell\in L}\widehat{\mathcal{B}}_M\big(\text{Shell}\,(w_\ell,D),D,s\big)\,,\tag{90.22}$$

where $K,L$ are two deterministic disjoint subsets of $\{1,\dots,N\}$, and $s$ is equal to $t-1$ or $t$. We derive next a probabilistic control on $S\big(\widehat{\mathcal{K}}_M\cap\widehat{\mathcal{L}}_M\big)$. We start by rewriting the definition of the two sets $\widehat{\mathcal{K}}_M,\widehat{\mathcal{L}}_M$ as

$$\begin{aligned}\widehat{\mathcal{K}}_M \;&=\; \widehat{\mathcal{B}}_M\Big(\text{Shell}\,\big(\{\,w_k,k\in K\,\},D\big),D,t\Big)\,,\\ \widehat{\mathcal{L}}_M \;&=\; \widehat{\mathcal{B}}_M\Big(\text{Shell}\,\big(\{\,w_\ell,\ell\in L\,\},D\big),D,s\Big)\,.\end{aligned}$$

It follows from corollary 89.17 that the two sets above are the output of genuine truncated explorations. In case $s=t$, the third set $\widehat{\mathcal{K}}_M\cup\widehat{\mathcal{L}}_M$ is equal to

$$\widehat{\mathcal{K}}_M\cup\widehat{\mathcal{L}}_M \;=\; \widehat{\mathcal{B}}_M\Big(\text{Shell}\,\big(\{\,w_j,j\in K\cup L\,\},D\big),D,t\Big)\,,$$

hence it is also the output of genuine truncated explorations. In case $s=t-1$, we make appeal to the proposition 89.11 to assert that $\widehat{\mathcal{K}}_M\cup\widehat{\mathcal{L}}_M$, being the union of two genuine sets, is also a genuine set. Therefore, the probabilistic control (13.5) can be applied to the three sets $\widehat{\mathcal{K}}_M,\widehat{\mathcal{L}}_M,\widehat{\mathcal{K}}_M\cup\widehat{\mathcal{L}}_M$. We evaluate next the functional $S$ on these sets. Thanks to the additivity of $S$, we have

$$S\big(\widehat{\mathcal{K}}_M\cup\widehat{\mathcal{L}}_M\big)+S\big(\widehat{\mathcal{K}}_M\cap\widehat{\mathcal{L}}_M\big) \;=\; S\big(\widehat{\mathcal{K}}_M\big)+S\big(\widehat{\mathcal{L}}_M\big)\,.\tag{90.23}$$

It follows from the equality (90.23) that, for any $\lambda\geq 0$,

$$\begin{gathered}P\Big(S\big(\widehat{\mathcal{K}}_M\cap\widehat{\mathcal{L}}_M\big)\leq-\lambda\Big)\;\leq\\ P\Big(S\big(\widehat{\mathcal{K}}_M\big)\leq-\frac{\lambda}{3}\Big)+P\Big(S\big(\widehat{\mathcal{L}}_M\big)\leq-\frac{\lambda}{3}\Big)+P\Big(S\big(\widehat{\mathcal{K}}_M\cup\widehat{\mathcal{L}}_M\big)\geq\frac{\lambda}{3}\Big)\,.\end{gathered}$$

Using the probabilistic control (13.5) on the genuine sets $\widehat{\mathcal{K}}_M$, $\widehat{\mathcal{L}}_M$, $\widehat{\mathcal{K}}_M\cup\widehat{\mathcal{L}}_M$, we conclude that

$$\forall\lambda\geq 0\qquad P\Big(S\big(\widehat{\mathcal{K}}_M\cap\widehat{\mathcal{L}}_M\big)\leq-\lambda\Big)\;\leq\;6|\Lambda(4n)|\exp\Big(-\frac{2(p(1-p)\lambda)^2}{9|\Lambda(4n)|}\Big)\,.\tag{90.24}$$

We are going to use the estimate (90.24) to control the probabilities (90.20) and (90.21). To do so, we must ensure that the upper bounds on $S$ given there are negative. We proceed separately for each probability, and we start with the probability (90.20). We will require that $\rho(n)$ is such that

$$\rho(n)\;\leq\;\frac{p}{32d^2\tau(n)^3}\,.\tag{90.25}$$

This inequality, together with the lower bound (90.14) on $\tau(n)$, yields that

$$\tau(n)\rho(n)-p\Big(\frac{1}{8d\tau(n)^2}-5d\rho(n)\Big)\;\leq\;-\frac{p}{16d\tau(n)^2}\;<\;0\,.\tag{90.26}$$

Using the inequality (90.26), and the lower bound on $\left|\overline{\mathcal{T}}^0(t)\right|$, we obtain that the probability (90.20) is less or equal than

$$P\left(S\big(\widehat{\mathcal{K}}_M^-\cap\widehat{\mathcal{L}}_M^-\big)\,\leq\,-\frac{N^*n}{16d(1-p)\tau(n)^2\phi(n)}\right).\tag{90.27}$$

We consider next the probability (90.21). We will require that $\tau(n)$ is such that

$$\tau(n)\,\geq\,\frac{64d}{p}\,.\tag{90.28}$$

The inequalities (90.25), (90.28), together with the lower bound (90.14) on $\tau(n)$, yield that

$$\tau(n)\rho(n)+\frac{1}{\tau(n)}-\frac{p}{32d}\,\leq\,\frac{p}{32d^2\tau(n)^2}+\frac{1}{\tau(n)}-\frac{p}{32d}\,\leq\,-\frac{p}{128d}\,\leq\,-\frac{p}{16d\tau(n)^2}\,.\tag{90.29}$$

Using the inequality (90.29), and the lower bound on $\left|\overline{\mathcal{T}}^0(t)\right|$, we obtain that the probability (90.21) is less or equal than

$$P\left(S\big(\widehat{\mathcal{K}}_M^+\cap\widehat{\mathcal{L}}_M^+\big)\,\leq\,-\frac{N^*n}{16d(1-p)\tau(n)^2\phi(n)}\right).\tag{90.30}$$

As the pairs of sets $\big(\widehat{\mathcal{K}}_M^-,\widehat{\mathcal{L}}_M^-\big)$, $\big(\widehat{\mathcal{K}}_M^+,\widehat{\mathcal{L}}_M^+\big)$ are both of the form (90.22), we can employ the estimate (90.24) to control the two probabilities (90.27), (90.30). Thanks to our judicious conditions (90.25) and (90.28), we obtain the same upper bound for both, namely

$$6|\Lambda(4n)|\exp\Big(-\frac{(pN^*n)^2}{1152d^2|\Lambda(4n)|\tau(n)^4\phi(n)^2}\Big)\,.$$

This gives an upper bound for the two probabilities (90.17) and (90.18), which in turn yields an upper bound for the probability (90.16). We substitute this upper bound in the inequality (89.175), and we obtain the following upper bound for the probability (89.140):

$$P\left(\begin{array}{c}\mathcal{D},\,\mathcal{E},\quad\left|\overline{\mathcal{T}}^0(t)\right|\,>\,\dfrac{N^*n}{\phi(n)},\quad\left|\overline{\mathcal{T}}^0(t-1)\right|\,\leq\,\rho(n)\left|\overline{\mathcal{T}}^0(t)\right|,\\ \dfrac{1}{8d}\left|\mathcal{T}^{-0}(t)\right|\,\leq\,\left|\overline{\mathcal{T}}(t-1)\right|+\\ \left|\Big(\bigcup\limits_{k\in K^-}\mathcal{B}\big(w_k,D\setminus\widehat{\mathcal{T}}_k^-(t),t\big)\Big)\cap\Big(\bigcup\limits_{\ell\in L^-}\mathcal{B}\big(w_\ell,D\setminus\overline{\mathcal{T}}_\ell(t-1),t-1\big)\Big)\right|\\ \dfrac{1}{8d}\left|\mathcal{T}^{+0}(t)\right|\,\leq\\ \left|\Big(\bigcup_{k\in K^+}\mathcal{B}\big(w_k,D\setminus\widehat{\mathcal{T}}_k^*(t),t\big)\Big)\cap\Big(\bigcup_{\ell\in L^+}\mathcal{B}\big(w_\ell,D\setminus\widehat{\mathcal{T}}_\ell^*(t),t\big)\Big)\right|\end{array}\right)$$

$$\leq\,12|\Lambda(4n)|\exp\Big(-\frac{(pN^*n)^2}{1152d^2|\Lambda(4n)|\tau(n)^4\phi(n)^2}\Big)\,.\tag{90.31}$$

## 90.3 Control of $\overline{\mathcal{T}}^0(t)$

We introduced the functions $\phi(n)$ and $\rho(n)$ at the end of section 87, because we had to deal with two different scenarios for the sequence $(\left|\overline{\mathcal{T}}^0(t)\right|, 0 \leq t \leq T)$. In order to proceed with the inequality (90.31), we have to use the specific expressions of the functions $\phi(n)$ and $\rho(n)$ in each scenario. The first scenario is the case where $\overline{\mathcal{T}}^0(t-1)$ is empty. This case is the simplest and it corresponds to the probability appearing in the sum (87.16) and the choices

$$\phi(n) \,=\, 16d^2\psi(n)\,, \qquad \rho(n) \,=\, 0\,.$$

With these values, and the expression (90.12) of $\tau(n)$, the upper bound obtained in (90.31) becomes

$$12|\Lambda(4n)|\exp\Big(-\frac{(pN^*n)^2}{10^{12}\,d^{14}18^{4d}(4d)^{4d/\alpha}|\Lambda(4n)|(\ln n)^{4d-4d\gamma+4d/\alpha}\psi(n)^2}\Big)\,. \tag{90.32}$$

We have replaced the numerical constant $1152\cdot 16^2\cdot 30^4$ by the larger value $10^{12}$. The second scenario is the case where $\overline{\mathcal{T}}^0(t-1)$ is not empty. This case is quite delicate, it corresponds to the probability appearing in the sum (87.17) and the choices

$$\phi(n) \,=\, 16d^2c^{\bar{t}(n)}\,, \quad \rho(n) \,=\, \frac{1}{c}\,.$$

The value $c$ was introduced in subsection 87.4. The time has now come to make an adequate choice for $c$. The point is that we must be able to apply lemma 87.2 to the sequence $\big(\left|\overline{\mathcal{T}}^0(T_0-1+i)\right|, 1 \leq i \leq T-T_0+1\big)$. In order to fulfill the hypothesis of lemma 87.2, we had imposed the condition (87.12), which we recall:

$$1 \,\leq\, c^{T+1} \,\leq\, \psi(n)\,. \tag{90.33}$$

On the event $\mathcal{E}$, we know that $T \leq (2\ln n)^{1-\gamma}$. Therefore the condition (90.33) is implied by the stronger condition

$$1 \,\leq\, c^{(2\ln n)^{1-\gamma}+1} \,\leq\, \psi(n)\,. \tag{90.34}$$

We can work with any value $c$ as long as it satisfies the inequalities (90.34). On the event in the probability (90.31), we have

$$\left|\overline{\mathcal{T}}^0(t)\right| \,>\, \frac{N^*n}{16d^2c^{\bar{t}(n)}}\,. \tag{90.35}$$

We wish to choose $c \geq 1$ in such a way that the polynomial order of $n$ in inequality (90.35) is not affected. Moreover, we wish to use the inequality (90.31), but we have to remember that this inequality was derived under the hypothesis that the condition (90.25) is satisfied, and we must therefore choose $c$ accordingly. The condition (90.25), rewritten in terms of $c$, $\bar{t}(n)$ and $M$, instead of $\rho(n)$ and $\tau(n)$, reads:

$$\frac{1}{c} \,\leq\, \frac{p}{864000\,d^8\big(9\bar{t}(n)M\big)^{3d}}\,.$$

For instance, we would be happy if there existed $n_2 \geq 1$ such that the following two conditions were satisfied:

$$\forall n \geq n_2 \qquad c \,\geq\, \frac{1}{p} 864000\, d^8 \big(9\bar{t}(n)M\big)^{3d}\,, \tag{90.36}$$

$$\exists\,\delta \in ]0,1[ \quad \forall n \geq n_2 \qquad c^{\bar{t}(n)} \,\leq\, \exp\big((\ln n)^{\delta}\big)\,. \tag{90.37}$$

We use the expressions of $M(n)$ and $\bar{t}(n)$ given in (85.4) and (87.14) in order to rewrite the inequalities (90.36) and (90.37):

$$\forall n \geq n_2 \qquad c \,\geq\, \frac{1}{p} 864000\, d^8 \Big(18(\ln n)^{1-\gamma} \big\lceil (4d \ln n)^{1/\alpha} \big\rceil\Big)^{3d}\,, \tag{90.38}$$

$$\exists\,\delta \in ]0,1[ \quad \forall n \geq n_2 \qquad \ln c \,\leq\, \frac{1}{2}(\ln n)^{\gamma+\delta-1}\,. \tag{90.39}$$

Recall that the exponent $\gamma$ is in the interval $]0,1/2[$. We take $\delta = 1-\gamma/2$, we remove $\gamma$ from the condition (90.38) and we obtain the following inequalities, which are stronger and simpler than (90.38) and (90.39):

$$\forall n \geq n_2 \qquad c \,\geq\, \frac{1}{p} 864000\, d^8 18^{3d} (4d)^{3d/\alpha} (\ln n)^{3d+3d/\alpha}\,, \tag{90.40}$$

$$\forall n \geq n_2 \qquad \ln c \,\leq\, \frac{1}{2}(\ln n)^{\gamma/2}\,. \tag{90.41}$$

What we mean is that the inequalities (90.40) and (90.41) imply the desired inequalities (90.36) and (90.37). So we try the value $c$ given by the expression appearing in (90.40), i.e.,

$$c \,=\, \frac{1}{p} 864000\, d^8 18^{3d} (4d)^{3d/\alpha} (\ln n)^{3d+3d/\alpha}\,. \tag{90.42}$$

With this choice, as $n$ goes to $\infty$, $\ln c$ is negligible compared to $(\ln n)^{\gamma/2}$, hence there exists a value $n_2(\alpha,\gamma,p,d)$ depending on $\alpha,\gamma,p,d$ such that

$$\forall n \geq n_2(\alpha,\gamma,p,d) \quad \ln\Big(\frac{1}{p} 864000\, d^8 18^{3d} (4d)^{3d/\alpha} (\ln n)^{3d+3d/\alpha}\Big) \,<\, \frac{1}{2}(\ln n)^{\gamma/2}\,. \tag{90.43}$$

The choice (90.42) for $c$ and the inequality (90.43) guarantee that (90.40) and (90.41) are satisfied with $n_2 = n_2(\alpha,\gamma,p,d)$. We must also satisfy the constraint (90.34) on $c$. However, we had not chosen the function $\psi(n)$ so far. To avoid any interference between this constraint and (90.41), we take

$$\psi(n) \,=\, \exp\Big(\frac{1}{2}(\ln n)^{1-\gamma/4}\Big)\,. \tag{90.44}$$

The function $\psi(n)$ is sub-polynomial, which is crucial. We check next that, with these choices, the condition (90.34) is satisfied. We have

$$c^{(2\ln n)^{\gamma}+1} \,\leq\, \exp\Big(\big((2\ln n)^{1-\gamma}+1\big)(\ln n)^{\gamma/2}\Big) \,\leq\, \exp\Big(3(\ln n)^{1-\gamma/2}\Big)\,.$$

As $1-\gamma/2 < 1-\gamma/4$, there exists a value $n_3(\alpha,\gamma,p,d)$ depending on $\alpha,\gamma,p,d$ such that

$$\forall n \geq n_3(\alpha,\gamma,p,d) \qquad c^{(2\ln n)^{1-\gamma}+1} \leq \psi(n)\,. \tag{90.45}$$

We must also fulfill the condition (90.28), but this is pretty easy: there exists a value $n_4(\alpha,p,d)$ depending on $\alpha,p,d$ such that

$$\forall n \geq n_4(\alpha,p,d) \qquad M^d \,=\, \big\lceil (4d\ln n)^{1/\alpha}\big\rceil^d \,\geq\, \frac{64d}{p}\,,$$

and this ensures that

$$\forall n \geq n_4(\alpha,p,d) \qquad \tau(n) \,=\, 30d^2\big(9\bar{t}(n)M\big)^d \,\geq\, \frac{64d}{p}\,.$$

With these values, the upper bound obtained in (90.31) becomes

$$12|\Lambda(4n)|\exp\Big(-\frac{(pN^*n)^2}{10^{12}\,d^{14}18^{4d}(4d)^{4d/\alpha}|\Lambda(4n)|(\ln n)^{4d-4d\gamma+4d/\alpha}c^{2\bar{t}(n)}}\Big)\,. \tag{90.46}$$

In fact, this quantity is similar to the other upper bound (90.32), the difference being that the function $\psi(n)$ has been replaced by $c^{\bar{t}(n)}$. We control this factor with the help of (90.41), we have the bound

$$c^{\bar{t}(n)} \,\leq\, \exp\Big(\frac{1}{2}(\ln n)^{\gamma/2}2(\ln n)^{1-\gamma}\Big) \,=\, \exp\Big((\ln n)^{1-\gamma/2}\Big)\,. \tag{90.47}$$

Substituting the inequality (90.47) in (90.46), we obtain the following upper bound for (90.31):

$$12|\Lambda(4n)|\exp\Big(-\frac{(pN^*n)^2e^{-2(\ln n)^{1-\gamma/2}}}{10^{12}\,d^{14}18^{4d}(4d)^{4d/\alpha}|\Lambda(4n)|(\ln n)^{4d-4d\gamma+4d/\alpha}}\Big)\,. \tag{90.48}$$

We see on this last quantity that the presence of the positive exponent $\gamma$ is vital. It ensures that the term $\exp\big(-2(\ln n)^{1-\gamma/2}\big)$ is sub-polynomial in $n$, and that it does not interfere with the term $N^*n$ on the polynomial scale. This exponent $\gamma$ has come a long way. It originates from the inequality (85.11), stated in the ignition of the proof, which itself comes from corollary 14.5, proposition 16.1 and the hypothesis (9.16). Notice that, for inequality (85.11) to hold, the integer $n$ has to be larger than a value $n_1(\theta(p),\alpha,p,d)$ depending on $\theta(p),\alpha,p,d$. The quantity (90.48) should be put in relationship with the other upper bound (90.32). Now that the function $\psi(n)$ has finally been chosen in (90.44), we can also substitute its value in this other upper bound (90.32), which becomes

$$12|\Lambda(4n)|\exp\Big(-\frac{(pN^*n)^2e^{-(\ln n)^{3/4}}}{10^{12}\,d^{14}18^{4d}(4d)^{4d/\alpha}|\Lambda(4n)|(\ln n)^{4d-4d\gamma+4d/\alpha}}\Big)\,. \tag{90.49}$$

We observe that the two bounds (90.48) and (90.49) look almost the same, the only difference between them are the two factors $e^{-2(\ln n)^{1-\gamma/2}}$ and $e^{-(\ln n)^{3/4}}$ appearing in the numerator of the fraction in the exponentials. The point is that they are both sub-polynomial. To reach this result, we had to tune adequately the choice of the function $\psi(n)$. We keep the largest of these two upper bounds, that is (90.48), and we report it successively in the inequalities (90.31) and (89.140), and we conclude the following:

$$\forall n \geq \max\Big(n_0(M_0,\alpha,p,d), n_1(\theta(p),\alpha,p,d),\\ n_2(\alpha,\gamma,p,d), n_3(\alpha,\gamma,p,d), n_4(\alpha,p,d)\Big)$$

$$P\left(\begin{array}{c} \mathcal{D},\,\mathcal{E},\quad \big|\overline{\mathcal{T}}^0(t)\big| > \dfrac{N^*n}{\phi(n)},\quad \big|\overline{\mathcal{T}}^0(t-1)\big| \leq \rho(n)\,\big|\overline{\mathcal{T}}^0(t)\big|, \\ \dfrac{1}{8d}\,\big|\mathcal{T}^{-0}(t)\big| \leq \big|\overline{\mathcal{T}}(t-1)\big| + \\ \Big|\Big(\displaystyle\bigcup_{k\in K^-}\mathcal{B}\big(w_k, D\setminus\widehat{\mathcal{T}}_k^-(t),t\big)\Big)\cap\Big(\bigcup_{\ell\in L^-}\mathcal{B}\big(w_\ell, D\setminus\overline{\mathcal{T}}_\ell(t-1),t-1\big)\Big)\Big| \\ \dfrac{1}{8d}\,\big|\mathcal{T}^{+0}(t)\big| \leq \\ \Big|\Big(\bigcup_{k\in K^+}\mathcal{B}\big(w_k, D\setminus\widehat{\mathcal{T}}_k^*(t),t\big)\Big)\cap\Big(\bigcup_{\ell\in L^+}\mathcal{B}\big(w_\ell, D\setminus\widehat{\mathcal{T}}_\ell^*(t),t\big)\Big)\Big| \end{array}\right)$$

$$\leq \qquad (90.50)$$

$$24|\Lambda(4n)|\exp\Big(-\frac{(pN^*n)^2e^{-2(\ln n)^{1-\gamma/2}}}{10^{12}\,d^{14}18^{4d}(4d)^{4d/\alpha}|\Lambda(4n)|(\ln n)^{4d-4d\gamma+4d/\alpha}}\Big)\,.$$

Of course, we could be more relaxed about the dependence of the values $n_0$, $n_1$, $n_2$, $n_3$, $n_4$ on the other parameters. However, we have to keep track of the conditions imposed on $n$ in the course of the argument. Indeed, the ultimate goal is to prove that some event occurring in $\Lambda(4n)$ is typical, meaning that the limit of its probability goes to 1 as $n$ goes to $\infty$. Everything would fall apart if, by some misfortune, we were to impose unworkable conditions on $n$.

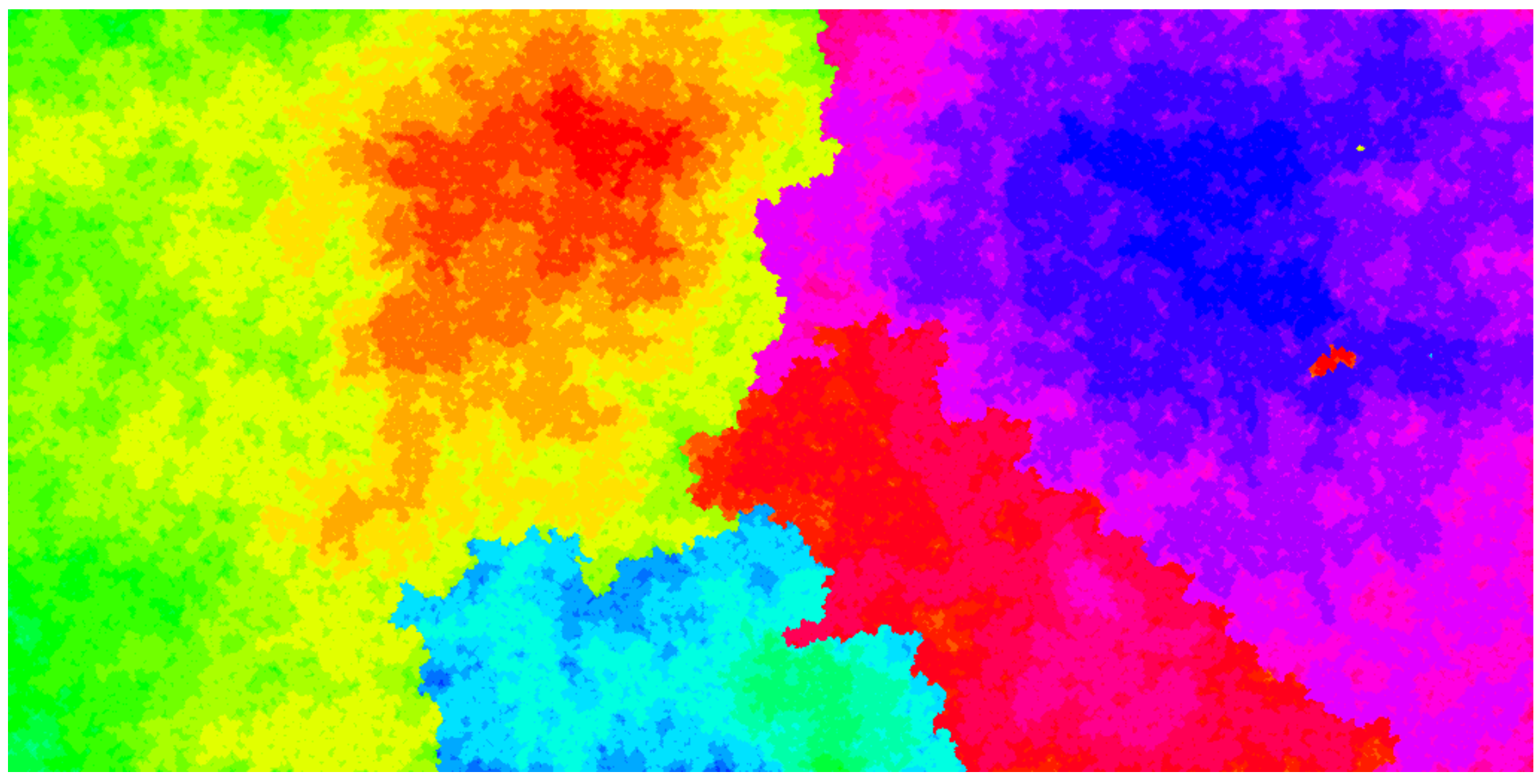

Figure 162: 7 intertwined explorations, $1024 \times 512$, site percolation, $p = 0.57$.

# 91 The exponent $3/2$

It remains to put together all the estimates we have done so far in order to reap the rewards of our work. This is achieved in subsection 91.1. An overview of the global flow of the inequalities and the three stages of the proof is presented in subsection 91.2.

## 91.1 Recapitulation

We shall put in action the inequality (90.50) that we have obtained in section 90. We rewind the argument presented in the previous sections and we play it back again, but this time we present the conditions imposed on $n$ at the very beginning. In the ignition of the proof, we introduced the integers $M_0$, $n_0(M_0,\alpha,p,d)$, $n_1(\theta(p),\alpha,p,d)$. Let us recall briefly the properties defining these integers (for the details, see section 85). We suppose that

$$\limsup_{n\to\infty}\frac{1}{\ln n}\ln\ln\Big(\frac{1}{P_p\big(n\leq|C(0)|<+\infty\big)}\Big)\;>\;0\,.$$

Thanks to this hypothesis, we obtain the existence of a positive exponent $\alpha$ and an integer $M_0$ as follows:

$$\exists\,\alpha>0\quad\exists M_0\geq 0\quad\forall M\geq M_0\qquad P_p\big(M\leq|C(0)|<+\infty\big)\;\leq\;\exp(-M^\alpha)\,. \tag{91.1}$$

The integer $n_0(M_0,\alpha,p,d)$ was defined in formula (85.5):

$$\exists\,n_0(M_0,\alpha,p,d)\geq 1\qquad\forall n\geq n_0(M_0,\alpha,p,d)\qquad P\big(\mathcal{F}_{\text{fi}}(n,M(n))^c\big)\;\leq\;\frac{1}{n}\,.$$

The above result rested on the property (91.1), the choice $M(n)=\big\lceil(4d\ln n)^{1/\alpha}\big\rceil$ and a simple union bound. The integer $n_1(\theta(p),\alpha,p,d)$ was defined in (85.11):

$$\exists\,\gamma\in]0,1/2[\quad\exists\,n_1(\theta(p),\alpha,p,d)\geq 1\quad\forall n\geq n_1\big(\theta(p),\alpha,p,d\big)\quad\beta(n)\;\geq\;(\ln n)^\gamma\,. \tag{91.2}$$

This result is much more sophisticated, it rests on proposition 16.1 and corollary 14.5, where the function $\beta(n)$ is introduced and a lower bound is given on $\beta(n)$ in terms of another function $\alpha(t)$. In turn, the hypothesis (9.16) yields an explicit lower bound on $\alpha(t)$, which permits to prove (91.2). Naturally, the above results hold only for a parameter $p$ such that $\theta(p)$ is positive. Thus the exponent $\gamma$, and the integers $M_0$, $n_0(M_0,\alpha,p,d)$, $n_1(\theta(p),\alpha,p,d)$ are fixed in the ignition of the proof. We define next $\overline{t}(n)=2(\ln n)^{1-\gamma}$, we introduce the value $c=c(n)$ given in formula (90.42) by

$$c\;=\;c(n)\;=\;\frac{1}{p}864000\,d^8 18^{3d}(4d)^{3d/\alpha}(\ln n)^{3d+3d/\alpha}\,, \tag{91.3}$$

and we choose

$$\psi(n)\;=\;\exp\Big(\frac{1}{2}(\ln n)^{1-\gamma/4}\Big)\,.$$

The integer $n_2(\alpha,\gamma,p,d)$ was defined towards the end of the argument in formula (90.43), but logically it should also be introduced right at the beginning. It is defined by the following property:

$$\exists\, n_2(\alpha,\gamma,p,d)\geq 1 \quad \forall n\geq n_2(\alpha,\gamma,p,d)$$
$$\ln\Big(\frac{1}{p}864000\, d^8 18^{3d}(4d)^{3d/\alpha}(\ln n)^{3d+3d/\alpha}\Big)\,<\,\frac{1}{2}(\ln n)^{\gamma/2}\,.$$

The integer $n_3(\alpha,\gamma,p,d)$ was introduced to satisfy the condition (90.45). Since $\gamma<1/2$, there exists a value $n_3(\alpha,\gamma,p,d)\geq n_2(\alpha,\gamma,p,d)$ depending on $\alpha,\gamma,p,d$ such that

$$\forall n\geq n_3(\alpha,\gamma,p,d)\qquad c(n)^{(2\ln n)^{1-\gamma}+1}\,\leq\,\exp\Big(3(\ln n)^{1-\gamma/2}\Big)\,\leq\,\psi(n)\,.$$

Finally, the integer $n_4(\alpha,p,d)$ was introduced to satisfy the condition (90.28): there exists a value $n_4(\alpha,p,d)$ depending on $\alpha,p,d$ such that

$$\forall n\geq n_4(\alpha,p,d)\qquad M^d\,=\,\big\lceil(4d\ln n)^{1/\alpha}\big\rceil^d\,\geq\,\frac{64d}{p}\,.$$

So we start now the whole story, and we prove then the inequality (90.50), that we recall next:

$$\forall n\geq \max\Big(n_0(M_0,\alpha,p,d),n_1(\theta(p),\alpha,p,d),$$
$$n_2(\alpha,\gamma,p,d),n_3(\alpha,\gamma,p,d),n_4(\alpha,p,d)\Big)$$

$$P\left(\begin{array}{c}\mathcal{D},\,\mathcal{E},\quad\big|\overline{\mathcal{T}}^0(t)\big|\,>\,\dfrac{N^*n}{\phi(n)}\,,\quad\big|\overline{\mathcal{T}}^0(t-1)\big|\,\leq\,\rho(n)\,\big|\overline{\mathcal{T}}^0(t)\big|\,,\\ \dfrac{1}{8d}\,\big|\mathcal{T}^{-0}(t)\big|\,\leq\,\big|\overline{\mathcal{T}}(t-1)\big|+\\ \Big|\Big(\displaystyle\bigcup_{k\in K^-}\mathcal{B}\big(w_k,D\setminus\widehat{\mathcal{T}}_k^-(t),t\big)\Big)\cap\Big(\displaystyle\bigcup_{\ell\in L^-}\mathcal{B}\big(w_\ell,D\setminus\overline{\mathcal{T}}_\ell(t-1),t-1\big)\Big)\Big|\\ \dfrac{1}{8d}\,\big|\mathcal{T}^{+0}(t)\big|\,\leq\\ \Big|\Big(\bigcup_{k\in K^+}\mathcal{B}\big(w_k,D\setminus\widehat{\mathcal{T}}_k^*(t),t\big)\Big)\cap\Big(\bigcup_{\ell\in L^+}\mathcal{B}\big(w_\ell,D\setminus\widehat{\mathcal{T}}_\ell^*(t),t\big)\Big)\Big|\end{array}\right)$$
$$\leq \tag{91.4}$$
$$24|\Lambda(4n)|\exp\Big(-\frac{(pN^*n)^2e^{-2(\ln n)^{1-\gamma/2}}}{10^{12}\,d^{14}18^{4d}(4d)^{4d/\alpha}|\Lambda(4n)|(\ln n)^{4d-4d\gamma+4d/\alpha}}\Big)\,.$$

Notice that this inequality holds for both choices:

- $\phi(n)=16d^2\psi(n)$ and $\rho(n)=0$;
- $\phi(n)=16d^2c^{\bar{t}(n)}$ and $\rho(n)=1/c(n)$ where $c(n)$ is given by (91.3).

To alleviate the next formulas, we regroup all the sub-polynomial factors into the notation

$$\nu(n)\,=\,\frac{p^2e^{-2(\ln n)^{1-\gamma/2}}}{10^{12}\,d^{14}18^{4d}(4d)^{4d/\alpha}(\ln n)^{4d-4d\gamma+4d/\alpha}}\,.$$

We start from the initial inequality (87.2), and we derive the chain of inequalities leading until formula (88.37). We substitute the inequality (91.4) into the sum in formula (88.37), and we obtain

$$P\begin{pmatrix} \mathcal{D},\,\mathcal{E},\,\big|\overline{\mathcal{T}}^0(t)\big|\,>\,\dfrac{N^*n}{\phi(n)}\,, \\ \big|\overline{\mathcal{T}}^0(t-1)\big|\,\leq\,\rho(n)\,\big|\overline{\mathcal{T}}^0(t)\big| \end{pmatrix}\leq\sum_{\substack{K^-,L^-,K^+,L^+\subset\{\,1,\dots,N\,\}\\ K^-\cap L^-=\varnothing,\,K^+\cap L^+=\varnothing}} 24|\Lambda(4n)|e^{-\frac{\nu(n)(N^*n)^2}{|\Lambda(4n)|}}$$

$$\leq\,2^{4N}24|\Lambda(4n)|e^{-\frac{\nu(n)(N^*n)^2}{|\Lambda(4n)|}}\,.\qquad(91.5)$$

We then return to the inequality (87.15). There were three different terms on the right-hand side. We use the inequality (90.11) to control the first term (87.15), and the inequality (91.5) to control the second term (87.16) and the third term (87.17), and we get

$$P(\mathcal{H})\,\leq\qquad(91.6)$$

$$6\,2^{2N}|\Lambda(4n)|\exp\Big(-\frac{(pN^*n)^2e^{-(\ln n)^{3/4}}}{1152\cdot16^2\,d^8|\Lambda(4n)|}\Big)\,+\,2\overline{t}(n)\,2^{4N}24|\Lambda(4n)|e^{-\frac{\nu(n)(N^*n)^2}{|\Lambda(4n)|}}$$

$$\leq\,97(\ln n)\,2^{4N}|\Lambda(4n)|e^{-\frac{\nu(n)(N^*n)^2}{|\Lambda(4n)|}}\,.\qquad(91.7)$$

To obtain the final term (91.7), we have used the fact that the first term in the right-hand side of (91.6) is smaller than the second one. We substitute the above inequality into the inequality (87.3). This gives

$$P\big(\mathcal{N}^*\geq N^*,\mathcal{E},\mathcal{N}_{\mathrm{bd}}=N\big)\,\leq\sum_{w_1,\dots,w_N\in\partial^{\,in}\Lambda(4n)}97(\ln n)\,2^{4N}|\Lambda(4n)|e^{-\frac{\nu(n)(N^*n)^2}{|\Lambda(4n)|}}$$

$$\leq\,97|\Lambda(4n)|^{2N+2}e^{-\frac{\nu(n)(N^*n)^2}{|\Lambda(4n)|}}\,.$$

We finally arrive at the initial inequality (87.2). Setting

$$c(d,\alpha)\,=\,1152\cdot16^2\cdot30^4\,d^{14}18^{4d}(4d)^{4d/\alpha}\,,$$

we obtain that

$$P\big(\mathcal{N}^*\geq N^*,\mathcal{E}\big)\,\leq\sum_{N^*\leq N\leq|\partial^{\,in}\Lambda(4n)|}97|\Lambda(4n)|^{2N+2}e^{-\frac{\nu(n)(N^*n)^2}{|\Lambda(4n)|}}$$

$$\leq\,97\exp\left(\big(2|\partial^{\,in}\Lambda(4n)|+3\big)\,\ln\big(|\Lambda(4n)|\big)-\frac{(pN^*n)^2e^{-2(\ln n)^{1-\gamma/2}}}{c(d,\alpha)|\Lambda(4n)|(\ln n)^{4d-4d\gamma+4d/\alpha}}\right).\qquad(91.8)$$

Let us pause and contemplate this inequality. Was it worth the work? What can we conclude from it?

Of course, the inequality becomes useful only when the right-hand side goes to 0. Let us take $N^* = n^{d-3/2} e^{3(\ln n)^{1-\gamma/2}}$. We see that

$$\lim_{n\to\infty} P\Big(\mathcal{N}^* \geq n^{d-3/2} e^{3(\ln n)^{1-\gamma/2}}, \mathcal{E}\Big) \,=\, 0\,.$$

This implies that, when $\mathcal{E}$ occurs, with probability going to 1 as $n$ goes to $\infty$,

$$\mathcal{N}^* \,\leq\, n^{d-3/2+o(1)}\,. \tag{91.9}$$

Furthermore, on the event $\mathcal{E}$, the event $\mathcal{F}_{\mathrm{bd}}(n,\varepsilon)$ defined in (85.1) occurs as well, whence

$$\sum_{C\in\mathcal{C}_{\mathrm{bd}}(4n)} \big|C\cap\Lambda(2n)\big| \,\geq\, \big(\theta(p)-\varepsilon\big)\big|\Lambda(2n)\big|\,. \tag{91.10}$$

The inequality (91.9) tells that the number of terms in the sum on the left-hand side of (91.10) is at most $n^{d-3/2+o(1)}$. From there, we deduce that, with probability going to 1 as $n$ goes to $\infty$, the box $\Lambda(4n)$ contains a boundary cluster $C$ such that

$$\big|C\cap\Lambda(2n)\big| \,\geq\, \big(\theta(p)-\varepsilon\big)\big|\Lambda(2n)\big|n^{3/2-d+o(1)} \,=\, n^{3/2+o(1)}\,. \tag{91.11}$$

This is still not what we wish, but we have obtained an explicit power of $n$ much larger than 1, and this is a huge step forward!

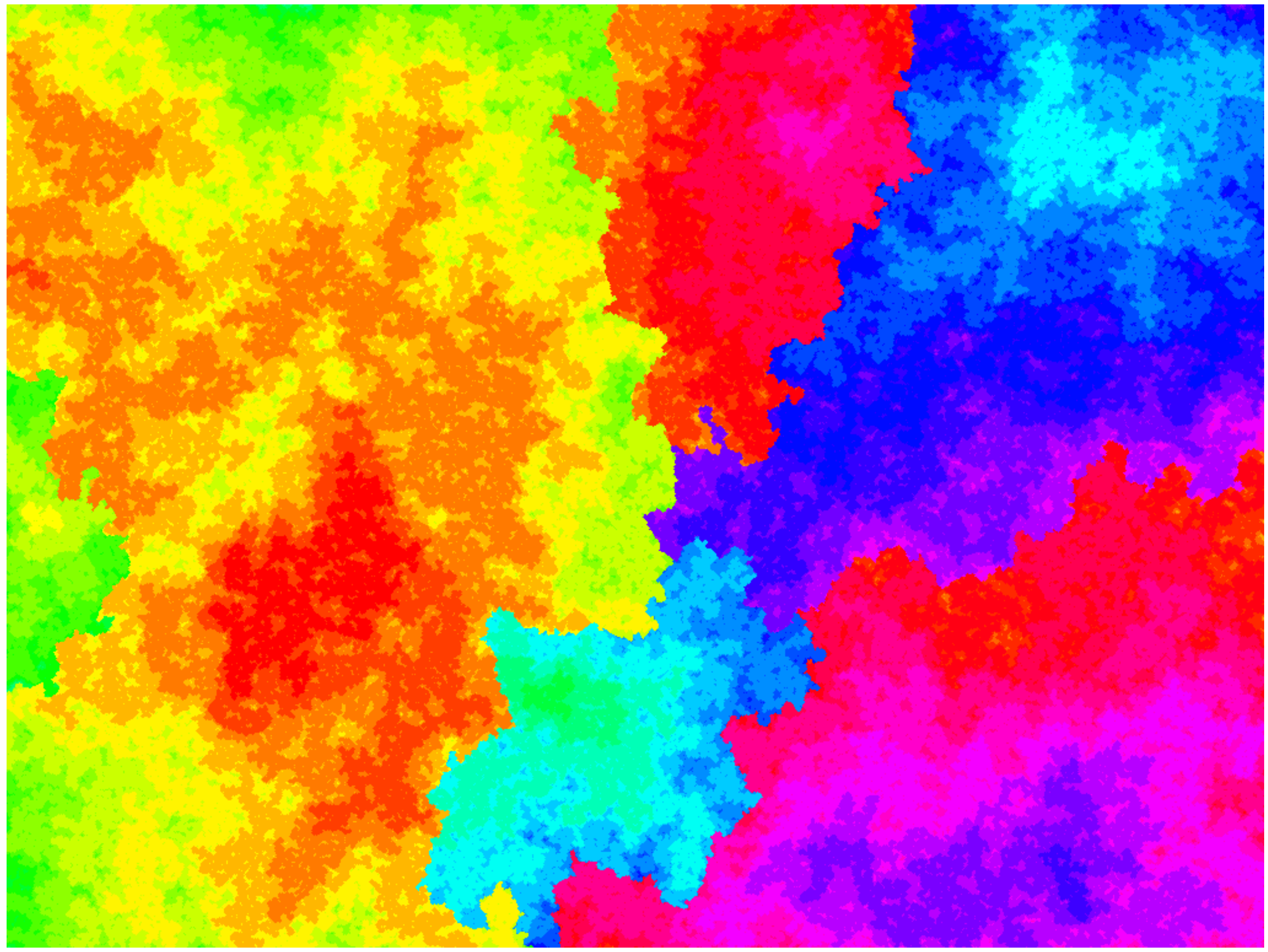

Figure 163: 7 intertwined explorations, $1024\times 768$, site percolation, $p = 0.57$.

## 91.2 Overview of the three stages

In order to get an overview of the flow of the sums and the inequalities, we recapitulate here the three stages of the argument, but to simplify, we leave the auxiliary functions $\phi(n)$, $\psi(n)$, $\rho(n)$, and we set also $a_n = N^*n/(16d^2\psi(n))$:

$$P\big(\mathcal{N}^* \geq N^*, \mathcal{E}\big) \leq \sum_{N^* \leq N \leq |\partial^{\,in}\Lambda(4n)|} P\Big(\mathcal{N}^* \geq N^*, \mathcal{E}, \exists\, w_1, \dots, w_N \in \partial^{\,in}\Lambda(4n) \quad \mathcal{C}_{\mathrm{bd}}(4n) = \big\{\, C\big(w_i, \Lambda(4n)\big), 1 \leq i \leq N\,\big\}\Big)$$

$$\leq \sum_{\substack{w_1,\dots,w_N \in \partial^{\,in}\Lambda(4n)\\ N^* \leq N \leq |\partial^{\,in}\Lambda(4n)|}} P\Big(\mathcal{N}^* \geq N^*, \mathcal{E}, \mathcal{C}_{\mathrm{bd}}(4n) = \big\{\, C\big(w_i, \Lambda(4n)\big), 1 \leq i \leq N\,\big\}\Big)$$

Ignition of the first stage

$$\leq \sum_{\substack{w_1,\dots,w_N \in \partial^{\,in}\Lambda(4n)\\ N^* \leq N \leq |\partial^{\,in}\Lambda(4n)|}} P\Big(\mathcal{D},\, \mathcal{E},\, \big|\overline{\mathcal{T}}^0(0)\big| > \frac{N^*n}{16d^2\psi(n)}\Big) + \sum_{\substack{w_1,\dots,w_N \in \partial^{\,in}\Lambda(4n)\\ N^* \leq N \leq |\partial^{\,in}\Lambda(4n)|\\ 1 \leq t \leq \overline{t}(n)}} 2P \left(\begin{matrix} \mathcal{D},\, \mathcal{E},\, \big|\overline{\mathcal{T}}^0(t)\big| > \dfrac{N^*n}{\phi(n)}, \\ \big|\overline{\mathcal{T}}^0(t-1)\big| \leq \rho(n)\,\big|\overline{\mathcal{T}}^0(t)\big| \end{matrix}\right)$$

$$\leq \sum_{\substack{w_1,\dots,w_N \in \partial^{\,in}\Lambda(4n)\\ N^* \leq N \leq |\partial^{\,in}\Lambda(4n)|}} P\left(\begin{matrix}\mathcal{D},\, \mathcal{E}\\ \big|\overline{\mathcal{T}}^0(0)\big| > a_n\end{matrix}\right) + \sum_{\substack{w_1,\dots,w_N \in \partial^{\,in}\Lambda(4n)\\ N^* \leq N \leq |\partial^{\,in}\Lambda(4n)|\\ 1 \leq t \leq \overline{t}(n)}} 2P \left(\begin{matrix} \mathcal{D},\, \mathcal{E}, \quad \big|\overline{\mathcal{T}}^0(t)\big| > \dfrac{N^*n}{\phi(n)}, \quad \big|\overline{\mathcal{T}}^0(t-1)\big| \leq \rho(n)\,\big|\overline{\mathcal{T}}^0(t)\big|, \\ \exists\, K^-, L^-, K^+, L^+ \subset \{\, 1, \dots, N\,\} \quad K^- \cap L^- = \varnothing,\, K^+ \cap L^+ = \varnothing, \\ \dfrac{1}{4d}\,\big|\mathcal{T}^{-0}(t)\big| \leq \big|\overline{\mathcal{T}}(t-1)\big| + \\ \Big|\Big(\bigcup\limits_{k \in K^-} \mathcal{B}\big(w_k, D \setminus \widehat{\mathcal{T}}_k^-(t), t\big)\Big) \cap \Big(\bigcup\limits_{\ell \in L^-} \mathcal{B}\big(w_\ell, D \setminus \overline{\mathcal{T}}_\ell(t-1), t-1\big)\Big)\Big| \\ \dfrac{1}{4d}\,\big|\mathcal{T}^{+0}(t)\big| \leq \\ \Big|\Big(\bigcup\limits_{k \in K^+} \mathcal{B}\big(w_k, D \setminus \widehat{\mathcal{T}}_k^*(t), t\big)\Big) \cap \Big(\bigcup\limits_{\ell \in L^+} \mathcal{B}\big(w_\ell, D \setminus \widehat{\mathcal{T}}_\ell^*(t), t\big)\Big)\Big| \end{matrix}\right)$$

$$a_n = \frac{N^*n}{16d^2\psi(n)}$$

End of the first stage and

Ignition of the second stage

$$
\leq \sum_{\substack{w_1,\dots,w_N\in\partial^{\,in}\Lambda(4n)\\ N^*\leq N\leq|\partial^{\,in}\Lambda(4n)|}} P\left(\begin{matrix}\mathcal{D},\,\mathcal{E}\\ \big|\overline{\mathcal{T}}^0(0)\big|>a_n\end{matrix}\right)+
\sum_{\substack{w_1,\dots,w_N\in\partial^{\,in}\Lambda(4n)\\ N^*\leq N\leq|\partial^{\,in}\Lambda(4n)|\\ 1\leq t\leq\bar{t}(n)\\ K^-,L^-,K^+,L^+\subset\{\,1,\dots,N\,\}\\ K^-\cap L^-=\varnothing,\,K^+\cap L^+=\varnothing}} 2P\left(\begin{matrix}
\mathcal{D},\,\mathcal{E},\quad \big|\overline{\mathcal{T}}^0(t)\big|>\dfrac{N^*n}{\phi(n)},\quad \big|\overline{\mathcal{T}}^0(t-1)\big|\leq\rho(n)\,\big|\overline{\mathcal{T}}^0(t)\big|,\\
\dfrac{1}{4d}\,\big|\mathcal{T}^{-0}(t)\big|\leq\big|\overline{\mathcal{T}}(t-1)\big|+\\
\Big|\Big(\bigcup_{k\in K^-}\mathcal{B}\big(w_k,D\setminus\widehat{\mathcal{T}}_k^-(t),t\big)\Big)\cap\Big(\bigcup_{\ell\in L^-}\mathcal{B}\big(w_\ell,D\setminus\overline{\mathcal{T}}_\ell(t-1),t-1\big)\Big)\Big|\\
\dfrac{1}{4d}\,\big|\mathcal{T}^{+0}(t)\big|\leq\\
\Big|\Big(\bigcup_{k\in K^+}\mathcal{B}\big(w_k,D\setminus\widehat{\mathcal{T}}_k^*(t),t\big)\Big)\cap\Big(\bigcup_{\ell\in L^+}\mathcal{B}\big(w_\ell,D\setminus\widehat{\mathcal{T}}_\ell^*(t),t\big)\Big)\Big|
\end{matrix}\right)
$$

End of the second stage and
Ignition of the third stage

$$
\leq \sum_{\substack{w_1,\dots,w_N\in\partial^{\,in}\Lambda(4n)\\ N^*\leq N\leq|\partial^{\,in}\Lambda(4n)|}} P\left(\begin{matrix}\mathcal{D},\,\mathcal{E}\\ \big|\overline{\mathcal{T}}^0(0)\big|>a_n\end{matrix}\right)+
\sum_{\substack{w_1,\dots,w_N\in\partial^{\,in}\Lambda(4n)\\ N^*\leq N\leq|\partial^{\,in}\Lambda(4n)|\\ 1\leq t\leq\bar{t}(n)\\ K^-,L^-,K^+,L^+\subset\{\,1,\dots,N\,\}\\ K^-\cap L^-=\varnothing,\,K^+\cap L^+=\varnothing}} 2P\left(\begin{matrix}
\mathcal{D},\,\mathcal{E},\quad \big|\overline{\mathcal{T}}^0(t)\big|>\dfrac{N^*n}{\phi(n)},\quad \big|\overline{\mathcal{T}}^0(t-1)\big|\leq\rho(n)\,\big|\overline{\mathcal{T}}^0(t)\big|,\\
\dfrac{1}{4d}\,\big|\mathcal{T}^{-0}(t)\big|\leq\big|\widehat{\mathcal{K}}_M^-\cap\widehat{\mathcal{L}}_M^-\big|+5d\big|\overline{\mathcal{T}}^0(t-1)\big|,\\
\dfrac{1}{4d}\,\big|\mathcal{T}^{+0}(t)\big|\leq\big|\widehat{\mathcal{K}}_M^+\cap\widehat{\mathcal{L}}_M^+\big|,\\
\Big|\big\{\,x\in\widehat{\mathcal{K}}_M^-\cap\widehat{\mathcal{L}}_M^-:x\text{ is open}\,\big\}\Big|\leq\tau(n)\big|\overline{\mathcal{T}}^0(t-1)\big|,\\
\Big|\big\{\,x\in\widehat{\mathcal{K}}_M^+\cap\widehat{\mathcal{L}}_M^+:x\text{ is open}\big\}\Big|\leq\tau(n)\big(\big|\overline{\mathcal{T}}^0(t-1)\big|+\big|\mathcal{T}^{-0}(t)\big|\big)
\end{matrix}\right)
$$

$$
i_n=\frac{N^*n}{8d(1-p)\tau(n)^2\phi(n)}
$$

$$
\leq \sum_{\substack{w_1,\dots,w_N\in\partial^{\,in}\Lambda(4n)\\ N^*\leq N\leq|\partial^{\,in}\Lambda(4n)|}} P\left(\begin{matrix}\mathcal{D},\,\mathcal{E}\\ \big|\overline{\mathcal{T}}^0(0)\big|>a_n\end{matrix}\right)+\sum_{\substack{w_1,\dots,w_N\in\partial^{\,in}\Lambda(4n)\\ N^*\leq N\leq|\partial^{\,in}\Lambda(4n)|\\ 1\leq t\leq\bar{t}(n)}}\ \sum_{\substack{K^-,L^-,K^+,L^+\subset\{\,1,\dots,N\,\}\\ K^-\cap L^-=\varnothing,\,K^+\cap L^+=\varnothing}} 4P\Big(S\big(\widehat{\mathcal{K}}_M\cap\widehat{\mathcal{L}}_M\big)\leq -i_n\Big)
$$

End of the third stage

$$
\leq\ 97\exp\left(\big(2|\partial^{\,in}\Lambda(4n)|+3\big)\ln\big(|\Lambda(4n)|\big)-\frac{(pN^*n)^2e^{-2(\ln n)^{1-\gamma/2}}}{10^{12}\,d^{14}18^{4d}(4d)^{4d/\alpha}|\Lambda(4n)|(\ln n)^{4d-4d\gamma+4d/\alpha}}\right).
$$

In the real computations, the terms in the last inequalities appear twice with different choices for $\phi(n)$ and $\rho(n)$.

# 92 Improving the exponent from $3/2$ to $2$

So far, we have proved that, when the event $\mathcal{E}$ occurs, with probability going to 1 as $n$ goes to $\infty$, the box $\Lambda(4n)$ contains a cluster of size larger than $n^{3/2+o(1)}$. As good as that may be, this is still not enough to achieve our goal. We shall improve the previous approach in order to reach the exponent 2 instead of 3/2. The exponent 3/2 was obtained thanks to the last inequality (91.8) of subsection 91.1, it came from the expression inside the exponential:

$$\big(2|\partial^{\,in}\Lambda(4n)|+3\big)\ln\big(|\Lambda(4n)|\big)-\frac{(pN^*n)^2e^{-(\ln n)^{1-\gamma/2}}}{c(d,\alpha)|\Lambda(4n)|(\ln n)^{4d-4d\gamma+4d/\alpha}}\,. \tag{92.1}$$

We would like to choose $N^*$ as small as possible, but at the same time we need that the expression (92.1) goes to $-\infty$, this is why we end up with $N^*=n^{3/2+o(1)}$. We will try here to obtain a better inequality than (91.8), more precisely we will try to gain on the boundary term $|\partial^{\,in}\Lambda(4n)|$ in (92.1). This term appeared when we took an upper bound of the quantity $|\Lambda(4n)|^{2N+2}$, where the integer $N$ was the number of boundary clusters of the percolation configuration in $\Lambda(4n)$. In order to improve this term, we will launch the intertwined explorations from a specific subset of the boundary clusters, at the expense of some additional complications. As we want to control the number $\mathcal{N}^*$ of the boundary clusters which intersect the box $\Lambda(2n)$, the natural candidate is precisely this collection of clusters. However, if the intertwined explorations are launched only from the boundary clusters in $\Lambda(4n)$ that intersect $\Lambda(2n)$, then some important points of our argument do not hold any more, namely:

• The termination time. The occurrence of $\mathcal{G}(n)$ is not enough to guarantee that the termination time is less than $\overline{t}(n)=2(\ln n)^{1-\gamma}$.

• Approximations of the taboo explorations. After the first iteration, the taboo explorations will not have explored all the boundary clusters. Thus there will remain large boundary clusters which might be explored afterwards, and this will damage the approximation of the taboo explorations by genuine ones.

To compensate for these problems, we introduce in subsection 92.2 a first sequence $\big(\Lambda(h,n),\,0\leq h\leq\overline{h}(n)\big)$ of increasing concentric boxes included in $\Lambda(4n)$. For each $h\geq 0$, the width of the annulus $\Lambda(h+1,n)\setminus\Lambda(h,n)$ is $\lambda(n)$, which is of order $n/\ln n$. This allows us to choose an index $h$ such that the number $\mathcal{N}^*(h+1,n)$ of boundary clusters intersecting $\Lambda(h+1,n)$ is less than 8 times the number of boundary clusters intersecting $\Lambda(h-2,n)$. We take then the inner box $\Lambda(h,n)$ as the goal set, and we start the intertwined explorations from the collection of the boundary clusters that intersect $\Lambda(h+1,n)$. We show in subsection 92.3 that the existence of $\mathcal{N}^*(h+1,n)/8$ distinct boundary clusters crossing the annulus $\Lambda(h+1,n)\setminus\Lambda(h-2,n)$ implies the existence of many closed relevant taboo sites, in fact a number of order $\mathcal{N}^*\lambda(n)/16d$. The hard part of the work is to develop a probabilistic control on this number. The new difficulty that arises compared to the previous argument is that, as we do not start the intertwined explorations from the whole collection of the boundary clusters, some of these boundary clusters will be explored after the first iteration. This will ruin the approximation of the taboo explorations by genuine truncated explorations, at least in the whole box $\Lambda(4n)$. To go around this problem, we rely on the event $\mathcal{F}_{\text{fi}}(n, M)$. The occurrence of this event guarantees that the finite clusters have a diameter at most $M$, thus the potential damage created by the union of the boundary clusters can travel a distance at most $M$ during each iteration step. To stay away from the potentially dangerous region, where the approximations cannot be carried out, we introduce in subsection 92.5 a second sequence of boxes, this time decreasing, $\Gamma = \big(\Gamma(t), 0 \leq t \leq T\big)$, included in $\Lambda(h, n) \setminus \Lambda(h-1, n)$. For each $t \geq 0$, the width of the annulus $\Gamma(t) \setminus \Gamma(t+1)$ is an adequate multiple of $M$, and this prevents the unwanted boundary clusters from exerting a nasty influence in $\Gamma(t)$ during iteration $t$. Instead of controlling the relevant taboo sites in the whole box, we will control them along this sequence of boxes. Each iteration creates a new error term, which propagates spatially towards the interior of the boxes. We will develop an adequate sequence of inequalities in order to achieve a satisfactory control.

We shall redo the whole argument presented throughout the sections 85, 87, 88, 90, 91, 89. This argument is called the 3/2-argument in the sequel, while the new argument will be called the 2-argument. Our first intention was simply to point out the few places where the 3/2-argument has to be modified. Although the ignition of the proof is essentially unchanged, the subsequent steps have to undergo numerous changes, that pervade the whole text, and we ended up with a presentation that was extremely difficult to follow. Thus we decided to present an autonomous proof, in which we go faster for the passages already done in the 3/2-argument, and we reuse the few lemmas that do not need any adaptation at all. The drawback is that it makes the text longer, the advantage is that the new 2-argument can be read independently of the previous 3/2-argument. Another possibility would have been to include the required changes from the beginning, and to present only the 2-argument, but this would have added one layer of complexity, and it would have considerably obscured the presentation.

The proof is presented in the sequel of the section. Although the general strategy is the same as in the 3/2-argument, several proofs have to be completely reworked with respect to their 3/2-counterpart. The presentation of the new 2-argument has been completely overhauled compared with the 3/2-argument. The next figures 164 and 165 are diagrams presenting the organization of the 3/2-argument and the 2-argument. As evidenced by the intersecting arrows, the 3/2-argument is more entangled than the 2-argument. In the 3/2-argument, we have grouped into separate sections the localization and approximation results, as well as the probabilistic controls of the relevant taboo sites. The drawback is that they need to be used several times for the various relevant taboo sets, hence the entanglement. The new presentation may be less efficient, but it is also more intuitive. Essentially, the various inequalities and estimates are computed precisely when they are needed, in the order of appearance of the relevant taboo sites. In the 3/2-argument, they were grouped according to their nature and they were assembled at the end. The differences between these two presentations are clearly visible in the lower halves of figures 164 and 165: the left and right columns have essentially been swapped.

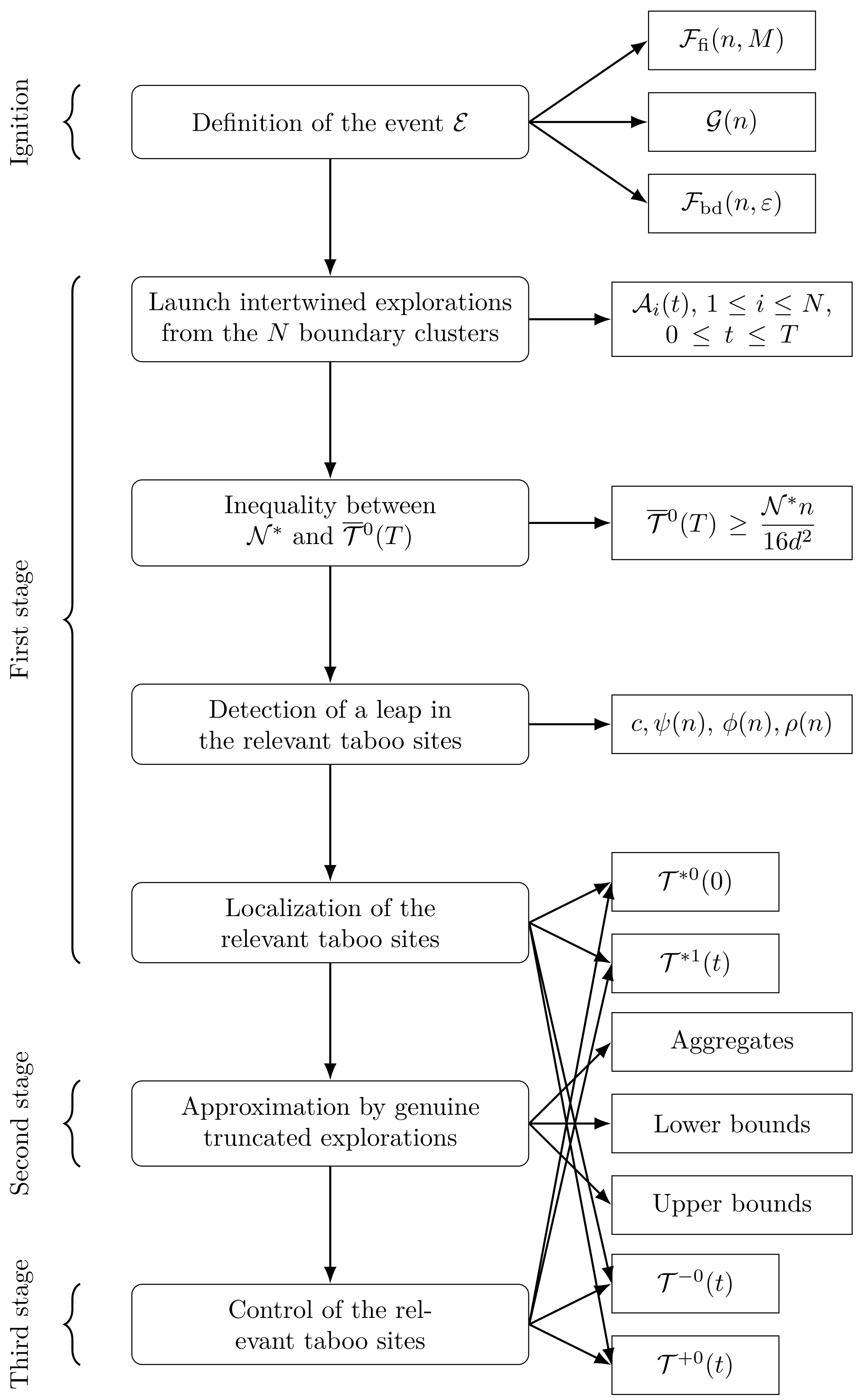

Figure 164: Scheme of the 3/2-argument

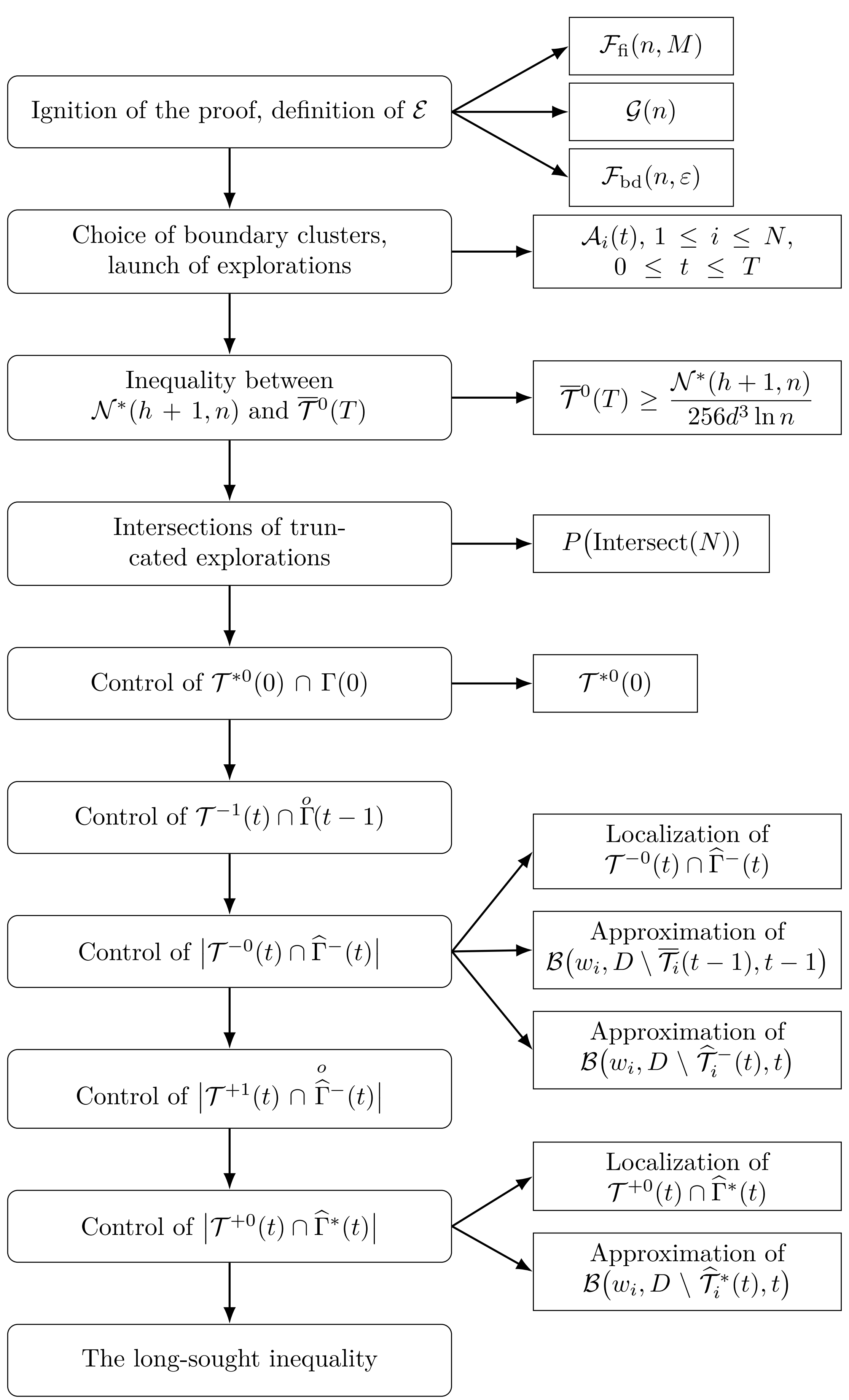

Figure 165: Scheme of the 2-argument

## 92.1 Ignition of the proof

The starting point is exactly the same as for the proof of the exponent 3/2, the ignition of the proof explained in section 85 is essentially unchanged. We reproduce quickly the main points.

Let us fix $\varepsilon > 0$. We consider the event $\mathcal{F}_{\text{bd}}(n,\varepsilon)$ defined by

$$\mathcal{F}_{\text{bd}}(n,\varepsilon) \,=\, \Big\{ \sum_{C\in\mathcal{C}_{\text{bd}}(4n)} \big|C\cap\Lambda(2n)\big| \,\geq\, \big(\theta(p)-\varepsilon\big)\big|\Lambda(2n)\big| \Big\}\,. \tag{92.2}$$

We proved in sections 8.5 and 10 that

$$\lim_{n\to\infty} P\big(\mathcal{F}_{\text{bd}}(n,\varepsilon)\big) \,=\, 1\,. \tag{92.3}$$

Let us fix an integer $M>0$. We consider the event $\mathcal{F}_{\text{fi}}(n,M)$ defined by

$$\mathcal{F}_{\text{fi}}(n,M) \,=\, \Big\{ \forall C\in\mathcal{C}_{\text{fi}}(4n) \quad |C|<M \Big\}\,.$$

It follows from the hypothesis (9.16) that

$$\exists\,\alpha>0 \quad \exists M_0\geq 0 \quad \forall M\geq M_0 \qquad P_p\big(M\leq |C(0)|<+\infty\big) \,\leq\, \exp(-M^\alpha)\,. \tag{92.4}$$

This estimate and a standard union bound yield that

$$\forall n\geq 1 \quad \forall M\geq M_0 \qquad P\big(\mathcal{F}_{\text{fi}}(n,M)^c\big) \,\leq\, \big|\Lambda(4n)\big|\exp(-M^\alpha)\,. \tag{92.5}$$

We set

$$M(n) \,=\, \big\lceil (4d\ln n)^{1/\alpha}\big\rceil\,. \tag{92.6}$$

Using (92.5), we see that there exists $n_0(M_0,\alpha,p,d)$ such that

$$\forall n\geq n_0(M_0,\alpha,p,d) \qquad P\big(\mathcal{F}_{\text{fi}}(n,M(n))^c\big) \,\leq\, \frac{1}{n}\,. \tag{92.7}$$

Proposition 16.1 states that there exists a non-decreasing function $\beta(n)$ defined for $n\geq 3$ satisfying $\beta(n)\leq\sqrt{\ln n}$ for $n\geq 3$, $\lim_{n\to+\infty}\beta(n)=+\infty$, and such that

$$\forall n\geq 3 \qquad P\Big(\forall x\in\Lambda(4n) \quad T_{\Lambda(4n)}\big(x,\partial^{\,in}\Lambda(4n)\big) \,\leq\, \frac{\ln n}{\beta(n)}\Big) \,\geq\, 1-n^{3d-\beta(n)}\,. \tag{92.8}$$

We use next the hypothesis (9.16) to obtain a quantitative estimate on $\beta(n)$. The function $\beta(n)$ was defined in (16.3) by setting

$$\forall n\geq 3 \qquad \beta(n) \,=\, \min\Big(\sqrt{\ln n}, \sqrt{\alpha\big(\sqrt{\ln n}\big)}\Big)\,. \tag{92.9}$$

The function $\alpha(t)$ was introduced in corollary 14.5, where it was also proved that, for $t$ larger than a value $t_0\big(\theta(p),p,d\big)$ depending on $\theta(p),p,d$, we have

$$\alpha(t) \,\geq\, \ln\Big(\frac{\theta(p)}{2}\Big) - \frac{1}{2}\ln\bigg(-\phi\Big(P_p\Big(\Big(\frac{t}{4d}\Big)^{d/(d+1)}<|C(0)|<+\infty\Big)\Big)\bigg)\,. \tag{92.10}$$

Inequality (92.10) and the tail estimate (92.4) imply that

$$\exists \alpha > 0 \quad \forall t \geq t_0\big(\theta(p), p, d\big) \qquad \alpha(t) \,\geq\, t^{\frac{\alpha d}{3(d+1)}}\,. \tag{92.11}$$

Substituting inequality (92.11) into (92.9), we conclude that there exist $\gamma$ in $]0, 1/2[$ and a value $n_1\big(\theta(p), p, d\big)$ depending on $\theta(p), p, d$ such that

$$\forall n \geq n_1\big(\theta(p), p, d\big) \qquad \beta(n) \,\geq\, (\ln n)^{\gamma}\,. \tag{92.12}$$

We define the event $\mathcal{G}(n)$ as

$$\mathcal{G}(n) \,=\, \Big\{\, \forall x \in \Lambda(4n) \quad T_{\Lambda(4n)}\big(x, \partial^{\,in}\Lambda(4n)\big) \,\leq\, (\ln n)^{1-\gamma} \,\Big\}\,. \tag{92.13}$$

Putting together the estimates (92.8) and (92.12), we conclude that

$$\forall n \geq n_1\big(\theta(p), p, d\big) \qquad P\big(\mathcal{G}(n)\big) \,\geq\, 1 - n^{3d-(\ln n)^{\gamma}}\,. \tag{92.14}$$

To alleviate the formulas, we set

$$\mathcal{E} \,=\, \mathcal{F}_{\text{fi}}(n, M) \cap \mathcal{G}(n) \cap \mathcal{F}_{\text{bd}}(n, \varepsilon)\,.$$

It follows from (92.3), (92.7), (92.14) that

$$\lim_{n\to\infty} P\big(\mathcal{E}\big) \,=\, 1\,.$$

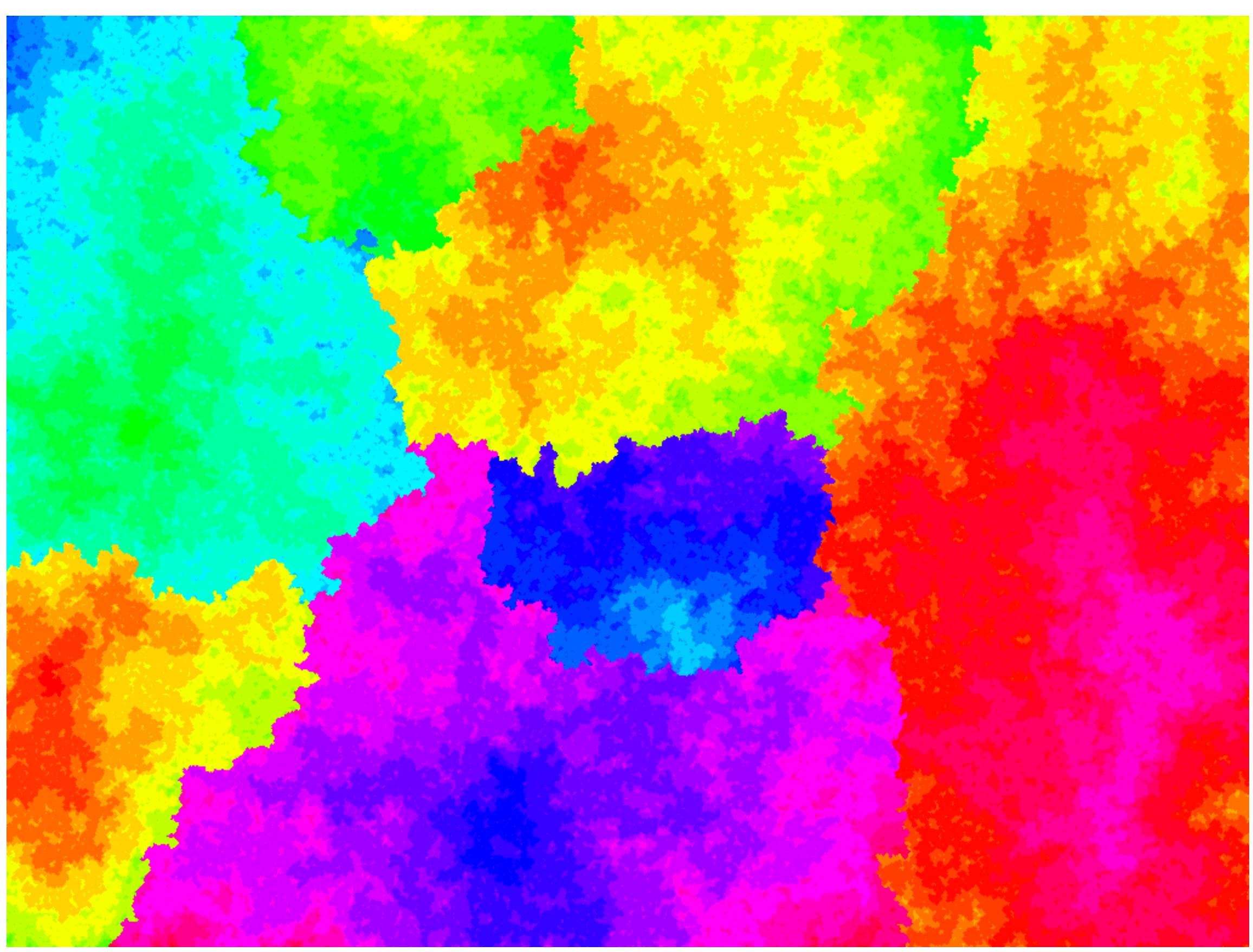

Figure 166: 7 intertwined explorations, $1024 \times 768$, site percolation, $p = 0.57$.

## 92.2 The choice of annulus-crossing clusters

Here comes a main change compared to the 3/2-argument: we will choose more carefully the boundary clusters that are used to launch the multiple intertwined explorations. We define an intermediate scale $\lambda(n)$ and an associated number $\overline{h}(n)$ by setting

$$\forall n \geq 1 \qquad \lambda(n) \,=\, \Big\lfloor \frac{n \ln 2}{d \ln n} \Big\rfloor\,, \quad \overline{h}(n) \,=\, \Big\lfloor \frac{n}{\lambda(n)} \Big\rfloor\,. \tag{92.15}$$

There exists an integer $n_3 = n_3(d,\alpha) \geq 16^d$ such that

$$\forall n \geq n_3(d,\alpha) \qquad M(n) \ln n \,<\, \frac{n}{2d \ln n} \,\leq\, \Big\lfloor \frac{n \ln 2}{d \ln n} \Big\rfloor\,. \tag{92.16}$$

From now onwards, we suppose that

$$n \,\geq\, n_3(d,\alpha)\,.$$

As $\gamma$ is positive, the condition (92.16) implies in particular that

$$M(n)(\ln n)^{1-\gamma} \,\leq\, \lambda(n)\,. \tag{92.17}$$

We introduce a sequence of concentric boxes included in $\Lambda(4n)$ by setting

$$\Lambda(h,n) \,=\, \Lambda\big(2n + 2h\lambda(n)\big)\,, \quad 0 \leq h \leq \overline{h}(n)\,. \tag{92.18}$$

For $h$ in $\{\,0,\dots,\overline{h}(n)\,\}$, we define the collection $\mathcal{C}^*(h,n)$ of the boundary clusters in $\Lambda(4n)$ which intersect the box $\Lambda(h,n)$:

$$\mathcal{C}^*(h,n) \,=\, \big\{\, C \in \mathcal{C}_{\text{bd}}\big(\Lambda(4n)\big) : C \cap \Lambda(h,n) \neq \varnothing \,\big\}\,.$$

We denote by $\mathcal{N}^*(h,n)$ their number, i.e.,

$$\forall h \in \{\,0,\dots,\overline{h}(n)\,\} \qquad \mathcal{N}^*(h,n) \,=\, \big|\mathcal{C}^*(h,n)\big|\,.$$

Obviously, we have the inclusions

$$\forall h \in \{\,1,\dots,\overline{h}(n)\,\} \qquad \mathcal{C}^*(h-1,n) \,\subset\, \mathcal{C}^*(h,n)\,,$$

and the sequence $\mathcal{N}^*(h,n)$, $0 \leq h \leq \overline{h}(n)$, is non-decreasing. As the number of boundary clusters in $\Lambda(4n)$ is bounded by the cardinality of $\partial^{\,in}\Lambda(4n)$, we have also the upper bound:

$$\mathcal{N}^*\big(\overline{h}(n),n\big) \,\leq\, \big|\mathcal{C}_{\text{bd}}\big(\Lambda(4n)\big)\big| \,\leq\, \big|\partial^{\,in}\Lambda(4n)\big|\,. \tag{92.19}$$

In the 3/2-argument, we showed that, when the number of boundary clusters in $\Lambda(4n)$ reaching the box $\Lambda(2n)$ was too large, there was an excess of closed sites in the annulus $\Lambda(4n) \setminus \Lambda(2n)$. In the new strategy, we will consider a smaller annulus. Namely, we will look for a value of $h$ such that $\mathcal{N}^*(h-2,n)$ and $\mathcal{N}^*(h+1,n)$ are comparable, and we will show that, when the number of boundary clusters in $\Lambda(4n)$ reaching the box $\Lambda(h-2,n)$ is too large, there is an excess of closed sites in the annulus $\Lambda(h-1,n) \setminus \Lambda(h-2,n)$. We start this program by choosing adequately the value of $h$, with the help of the following lemma.

**Lemma 92.1.** *Suppose that $n \geq 16^d$. If $\mathcal{N}^*(0,n)$ is not equal to 0, then there exists $H$ such that*

$$2 \,\leq\, H \,\leq\, \overline{h}(n) - 1\,, \qquad \mathcal{N}^*(H-2,n) \,\geq\, \frac{1}{8}\mathcal{N}^*(H+1,n)\,.$$

*Proof.* Suppose that the result is not true. We would then have

$$\forall h \in \{\,2,\dots,\overline{h}(n)-1\,\} \qquad \mathcal{N}^*(h-2,n) \,<\, \frac{1}{8}\mathcal{N}^*(h+1,n)\,. \tag{92.20}$$

Let $h^*$ be the largest integer such that $3h^* \leq \overline{h}(n)$. Iterating the inequality in (92.20), we would get

$$\mathcal{N}^*(0,n) \,<\, \frac{1}{8^{h^*}}\mathcal{N}^*(3h^*,n)\,. \tag{92.21}$$

It follows from (92.21) and the hypothesis $\mathcal{N}^*(0,n) \geq 1$ that

$$\mathcal{N}^*(3h^*,n) \,\geq\, 8^{h^*} \,=\, 2^{3h^*} \,\geq\, 2^{\overline{h}(n)-3} \,\geq\, 2^{\frac{n}{\lambda(n)}-4}\,. \tag{92.22}$$

In the last steps, we have used the inequalities

$$3h^* \,\leq\, \overline{h}(n) \,<\, 3h^*+3\,, \qquad \overline{h}(n) \,\leq\, \frac{n}{\lambda(n)} \,<\, \overline{h}(n)+1\,.$$

Putting together the inequalities (92.19) and (92.22), we would obtain

$$\big|\partial^{\,in}\Lambda(4n)\big| \,\geq\, 2^{\frac{n}{\lambda(n)}-4}\,. \tag{92.23}$$

On one hand, we have

$$\big|\partial^{\,in}\Lambda(4n)\big| \,\leq\, 2d(4n+1)^{d-1} \,\leq\, 6^d n^{d-1}\,. \tag{92.24}$$

On the other hand, the choice (92.15) of $\lambda(n)$ implies that

$$2^{\frac{n}{\lambda(n)}-4} \,=\, \frac{1}{16}\exp\Big(n\ln 2\Big/\Big\lfloor\frac{n\ln 2}{d\ln n}\Big\rfloor\Big) \,\geq\, \frac{1}{16}n^d\,. \tag{92.25}$$

We see that the three inequalities (92.23), (92.24), (92.25) cannot hold simultaneously for $n \geq 16^d$. Thus the result stated in the lemma is true. □

In the 3/2-argument, we have obtained directly a control on $\mathcal{N}^* = \mathcal{N}^*(0,n)$. We will try instead to control $\mathcal{N}^*(H+1,n)$, where $H$ is the integer provided by lemma 92.1. As $\mathcal{N}^*(H+1,n) \geq \mathcal{N}^*(0,n) = \mathcal{N}^*$, any improved estimate on $\mathcal{N}^*(H+1,n)$ will yield a better result for $\mathcal{N}^*$ as well. Let us rewrite the ignition of the first stage of the proof, as it should be for this new strategy. The passage below replaces the computation presented in formulas (87.2) and (87.3) of subsection 87.2. Let us fix an integer $N^* \geq 2$. We wish to estimate $P\big(\mathcal{N}^*(H+1,n)\geq N^*,\mathcal{E}\big)$. This time, we start by conditioning on the value of $H$ and then on the value of $\mathcal{N}^*(H+1,n)$:

$$
\begin{aligned}
P\big(\mathcal{N}^*(H+1,n)\geq N^*,\mathcal{E}\big) \,&=\, \sum_{0\leq h<\overline{h}(n)} P\big(\mathcal{N}^*(H+1,n)\geq N^*,\mathcal{E},H=h\big)\\
&= \sum_{0\leq h<\overline{h}(n)}\ \sum_{N^*\leq N\leq|\partial^{\,in}\Lambda(4n)|} P\big(\mathcal{E},H=h,\mathcal{N}^*(h+1,n)=N\big)\,. \quad (92.26)
\end{aligned}
$$

We fix next $h$ in $\{\,0,\dots,\overline{h}(n)-1\,\}$, $N$ in $\{\,N^*,\dots,|\partial^{\,in}\Lambda(4n)|\,\}$, and we suppose that $H=h$, $\mathcal{N}^*(h+1,n)=N$. Given $N$ sites $w_1,\dots,w_N$ in $\partial^{\,in}\Lambda(4n)$, we define the event $\mathcal{H}'(h,N,w_1,\dots,w_N)$ as

$$
\mathcal{H}'(h,N,w_1,\dots,w_N)\,=\,\mathcal{E}\cap\left\{\begin{matrix} H=h,\mathcal{N}^*(h+1,n)=N,\\ \mathcal{C}^*(h+1,n)=\big\{\,C\big(w_i,\Lambda(4n)\big),1\leq i\leq N\,\big\}\end{matrix}\right\}. \quad (92.27)
$$

The event $\mathcal{H}'(h,N,w_1,\dots,w_N)$ depends on $h,N$ and $w_1,\cdots,w_N$, but, later on, we will denote it simply $\mathcal{H}'$ to alleviate the formulas. We decompose and bound the probability in the sum (92.26) as follows:

$$
\begin{aligned}
P\left(\begin{matrix}\mathcal{E},H=h,\\ \mathcal{N}^*(h+1,n)=N\end{matrix}\right) \,&=\, P\Big(\bigcup_{w_1,\dots,w_N\in\partial^{\,in}\Lambda(4n)}\mathcal{H}'(h,N,w_1,\dots,w_N)\Big)\\
&\leq \sum_{w_1,\dots,w_N\in\partial^{\,in}\Lambda(4n)} P\big(\mathcal{H}'(h,N,w_1,\dots,w_N)\big)\,. \quad (92.28)
\end{aligned}
$$

We fix then $N$ elements $w_1,\dots,w_N$ in $\partial^{\,in}\Lambda(4n)$, and we analyze the configurations belonging to the event $\mathcal{H}'(h,N,w_1,\dots,w_N)$. We consider a configuration of percolation belonging to $\mathcal{H}'$ and we launch the multiple intertwined explorations starting from $w_1,\dots,w_N$. The goal set $G$ is taken to be the box $\Lambda(h,n)$. We denote by $\mathcal{D}$ the event that $w_1,\dots,w_N$ are pairwise disconnected in $\Lambda(4n)$. We derive next an upper bound on the termination time.

**Lemma 92.2.** *Suppose that the event $\mathcal{E}\cap\{\,H=h\,\}$ occurs. Then each site of $\Lambda(h,n)$ is at travel distance in $\Lambda(h+1,n)$ less than or equal to $(\ln n)^{1-\gamma}$ from a cluster of $\mathcal{C}^*(h+1,n)$, i.e.,*

$$
\forall x\in\Lambda(h,n)\quad\exists\, C\in\mathcal{C}^*(h+1,n)\quad T_{\Lambda(h+1,n)}(x,C)\,\leq\,(\ln n)^{1-\gamma}\,.
$$

*Proof.* Let $x$ be a site in $\Lambda(h,n)$. Since the event $\mathcal{G}(n)$ occurs, then

$$
T_{\Lambda(4n)}\big(x,\partial^{\,in}\Lambda(4n)\big)\,\leq\,(\ln n)^{1-\gamma}\,.
$$

Let $x_0,x_1,\cdots,x_r$ be a path in $\Lambda(4n)$ realizing the travel time between $x$ and $\partial^{\,in}\Lambda(4n)$. We have thus $x_0=x$, $x_r\in\partial^{\,in}\Lambda(4n)$, and at most $(\ln n)^{1-\gamma}$ sites among $x_0,x_1,\cdots,x_{r-1}$ are closed. Let $b$ be the index of the first site of the path that belongs to a boundary cluster of $\Lambda(4n)$:

$$
b\,=\,\min\Big\{\,i\in\{\,0,\dots,r\,\}:x_i\longleftrightarrow\partial^{\,in}\Lambda(4n)\,\Big\}\,.
$$

By the very definition of $b$, the sites $x_0, \dots, x_{b-1}$ are not in a boundary cluster of $\Lambda(4n)$, hence they all belong to finite clusters of $\Lambda(4n)$. On the event $\mathcal{F}_{\text{fi}}(n, M)$, the finite clusters have cardinalities strictly less than $M$. Thus the sequence $x_0, \dots, x_{b-1}$ cannot contain a run of $M$ consecutive open sites, whence

$$\Big|\big\{\, i \in \{\, 0, \dots, b-1 \,\} : x_i \text{ is closed} \,\big\}\Big| \;\geq\; \frac{b}{M}\,.$$

Yet the left-hand quantity is bounded by $(\ln n)^{1-\gamma}$, therefore

$$|x - x_b|_1 \;\leq\; b \;\leq\; M(\ln n)^{1-\gamma}\,. \tag{92.29}$$

As $x$ is in $\Lambda(h, n)$, it follows from (92.29), (92.17) that $x_b$, and in fact the whole sequence $x_0, \dots, x_{b-1}$, is in $\Lambda(h+1, n)$. This implies the desired result. $\square$

Lemma 92.2 guarantees that the termination time $T$ is at most $(\ln n)^{1-\gamma} + 1$, hence smaller than $\bar{t}(n) = 2(\ln n)^{1-\gamma}$ for $n \geq 3$.

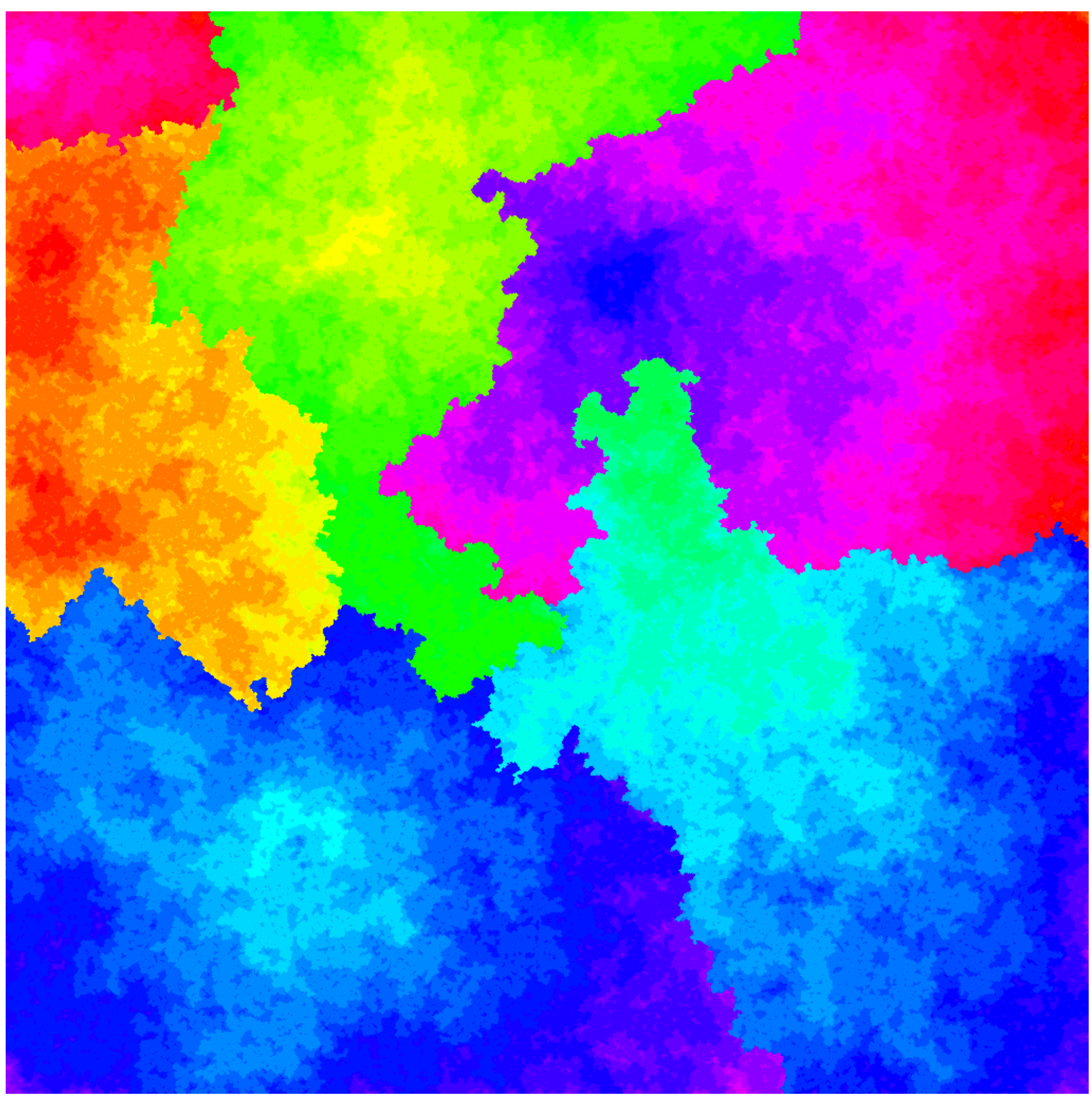

Figure 167: 7 intertwined explorations, $1024 \times 1024$, site percolation, $p = 0.57$.

## 92.3 From $\mathcal{N}^*(h+1,n)$ to $\overline{\mathcal{T}}^0(T) \cap \Lambda(h-1,n)$

The computations presented in this subsection are essentially the same as in subsection 87.3, except for one detail. Because the goal region is the box $\Lambda(h,n)$ instead of $\Lambda(4n)$, the final aggregates $(\mathcal{A}_i(T), 1 \leq i \leq N)$ produced by the intertwined explorations do not form a partition of $\Lambda(4n)$, instead they constitute a collection of pairwise disjoint sets that cover the box $\Lambda(h,n)$. As a consequence, the summation of the inequality (92.30) should not be carried over $s$ in $\{\, n,\dots,2n \,\}$, but rather over $s$ such that $n \leq s \leq n + h\lambda(n)$. Unfortunately, another problem will come later on, due to the spatial propagation of the error term, and we will not be able to control the relevant taboo sets in the box $\Lambda(h,n)$. We will need an extra buffer zone, so we will rather carry the summation over the indices $s$ in the set

$$S(h,n) \,=\, \big\{\, s : n + (h-2)\lambda(n) \,<\, s \,\leq\, n + (h-1)\lambda(n) \,\big\}\,.$$

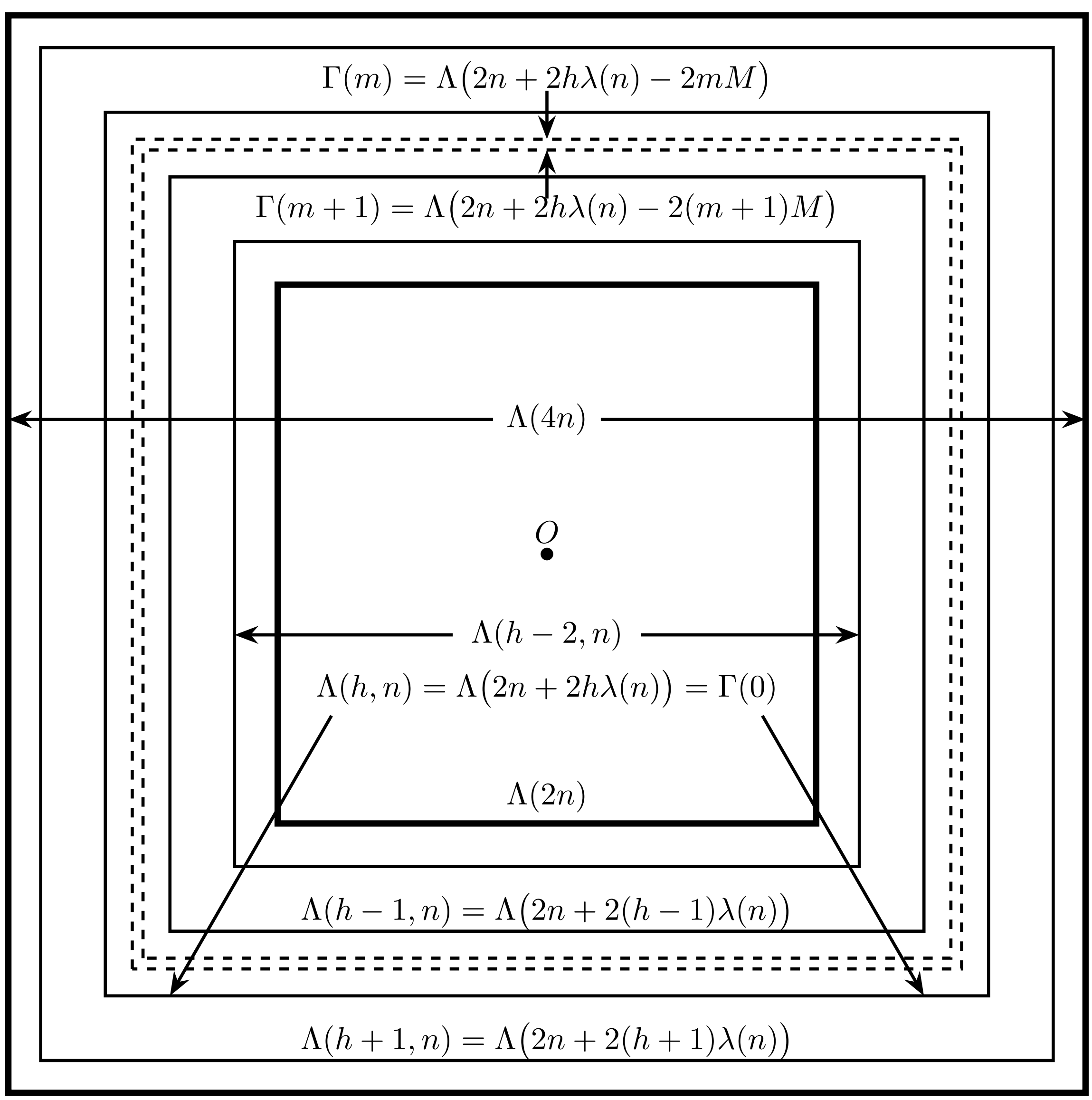

Figure 168: The boxes $\Lambda(h-2,n),\dots,\Lambda(h+1,n)$, $\Gamma(m)$, $\Gamma(m+1)$

Let us consider the final aggregates $(\mathcal{A}_i(T), 1 \leq i \leq N)$ produced by the intertwined explorations. They are pairwise distinct and they cover the box $\Lambda(h,n)$, because this box was taken to be the goal set of the intertwined explorations. In addition, each aggregate contains exactly one boundary cluster of the collection $\mathcal{C}^*(h+1,n)$. We know also that there are at least $N/8$ boundary clusters that intersect $\Lambda(h-2,n)$. Let us consider a boundary cluster $C$ in $\Lambda(4n)$ that intersects the box $\Lambda(h-2,n)$. The cluster $C$ contains an open path $\gamma$ inside $\Lambda(4n)$ leading from $\partial^{in}\Lambda(4n)$ to a site of $\Lambda(h-2,n)$. Let us consider the subpath $\overline{\gamma}$ of $\gamma$ which starts at its last visit to $\partial^{in}\Lambda(h-1,n)$ before it enters in $\Lambda(h-2,n)$. The subpath $\overline{\gamma}$ is confined to the box $\Lambda(h-1,n)$, it starts at one of the $2d$ faces of $\Lambda(h-1,n)$ and later on it enters $\Lambda(h-2,n)$. Suppose for instance that it starts at the face $F(\Lambda(h-1,n),u_1)$. It will then cross the slab

$$\left\{ \big(s, x(2), \cdots, x(d)\big) : s \in S(h,n), x(2), \dots, x(d) \in \mathbb{Z} \right\}$$

in the direction of the first axis. The face of $\Lambda(h-1,n)$ varies from one boundary cluster $C$ to another. Since there are at least $N/8$ such boundary clusters but only $2d$ faces, then there exists a face shared by at least $N/(16d)$ clusters. Suppose again that this common face corresponds to $u_1$, and let us denote by $\mathcal{C}_1^*(n)$ the corresponding subcollection of boundary clusters. For $s$ in $\mathbb{Z}^d$, we denote by $H_1(s)$ the hyperplane $\{ s \} \times \mathbb{Z}^{d-1}$, i.e.,

$$\forall s \in \mathbb{Z} \qquad H_1(s) \,=\, \left\{ \big(s, x(2), \cdots, x(d)\big) : x(2), \dots, x(d) \in \mathbb{Z} \right\}.$$

Let us fix $s$ in $S(h,n)$. All the clusters of the collection $\mathcal{C}_1^*(n)$ intersect $H_1(s)$. Thus the traces of the aggregates $(\mathcal{A}_i(T), 1 \leq i \leq N)$ on $H_1(s)$ form a partition of $\Lambda(h-1,n) \cap H_1(s)$ containing at least $N/(16d)$ non-empty sets. We apply lemma 87.1 to the $(d-1)$-dimensional box $\Lambda(h-1,n) \cap H_1(s)$ and the previous partition, and we obtain

$$\Big| \bigcup_{1 \leq j < k \leq N} \big( \Delta_{\Lambda(h-1,n)\cap H_1(s)} \mathcal{A}_j(T) \cap \Delta_{\Lambda(h-1,n)\cap H_1(s)} \mathcal{A}_k(T) \big) \Big| \,\geq\, \frac{N}{32d}\,. \quad (92.30)$$

The hyperplanes $\big(H_1(s), s \in S(h,n)\big)$ are pairwise disjoint. Summing the inequality (92.30) over $s$ in $S(h,n)$, we conclude that

$$\Big| \bigcup_{0 \leq t \leq T-1} \mathcal{I}(t) \cap \Lambda(h-1,n) \Big| \,\geq\, \frac{N\lambda(n)}{32d}\,. \qquad (92.31)$$

By proposition 80.18, at least one extremity of each edge in an intersection set is a closed relevant taboo site. Therefore

$$\big|\overline{\mathcal{T}}^0(T) \cap \Lambda(h-1,n)\big| \,\geq\, \frac{1}{4d} \Big| \bigcup_{0 \leq t \leq T-1} \mathcal{I}(t) \cap \Lambda(h-1,n) \Big|\,. \qquad (92.32)$$

We deduce from (92.31), (92.32) and (92.16) that

$$\big|\overline{\mathcal{T}}^0(T) \cap \Lambda(h-1,n)\big| \,\geq\, \frac{Nn}{256 d^3 \ln n}\,. \qquad (92.33)$$

Thus any configuration of the event $\mathcal{H}'$ satisfies the inequality (92.33).

## 92.4 Spatial propagation

The construction of the genuine truncated approximations is unchanged, and the material from subsections 89.1 to 89.5 can be reused without modification. This is not the case for the subsection 89.6, the various approximation results must be completely reformulated. The cause of the new trouble is that the intertwined explorations are stopped once the box $\Lambda(h,n)$ is fully explored, thus they will not detect all the relevant taboo sites outside of this box. Whereas the approximation results presented in subsection 89.6 were controlled with the help of all the relevant taboo sites, we cannot afford this luxury here and the error term should involve only the relevant taboo sites belonging to $\Lambda(h,n)$. This significantly complicates the task, as errors now propagate both over time and space. In order to explain more precisely the problem, let us discuss the case of the taboo sets of the first two steps.

We start with the first step. For $t=0$, the set $\overline{\mathcal{T}}(0)=\mathcal{T}^*(0)$ is the set of the pivotal sites for $\mathcal{D}$, it can be controlled in the same way as usual, the technique used in subsection 90.3 applies without modification. Furthermore, for any subset $I$ of $\{1,\dots,N\}$, we have

$$\bigcup_{i\in I}\mathcal{B}\big(w_i,D\setminus\overline{\mathcal{T}}_i(0),0\big)\,\subset\,\widehat{\mathcal{B}}_M\Big(\text{Shell}\,\big(W(I),D\big),D,0\Big)$$
$$=\,\text{Shell}\,\big(W(I),D\big)\,\subset\,\bigcup_{i\in I}\mathcal{B}\big(w_i,D\setminus\overline{\mathcal{T}}_i(0),0\big)\cup\overline{\mathcal{T}}(0)\,.\quad(92.34)$$

We examine now what happens for the second step. The next goal is to control the set $\overline{\mathcal{T}}(1)$. We decompose it as

$$\overline{\mathcal{T}}(1)\,=\,\mathcal{T}^0(0)\cup\mathcal{T}^{-0}(1)\cup\mathcal{T}^{+0}(1)\cup\mathcal{T}^{-1}(1)\cup\mathcal{T}^{+1}(1)\,.$$

The set $\mathcal{T}^0(0)$ has already been controlled. The easiest set to control among the remaining sets is $\mathcal{T}^{-1}(1)$. Indeed, lemma 88.1 tells that any site in $\mathcal{T}^{-1}(1)$ has a neighbour which is in $\mathcal{T}^{+0}(0)$. Thus we have

$$\big|\mathcal{T}^{-1}(1)\big|\,\leq\,2d\,\big|\mathcal{T}^{+0}(0)\big|\,.$$

At this point, we have gained control of $\mathcal{T}^0(0)$ and $\mathcal{T}^{-1}(1)$. We aim next at controlling the set $\mathcal{T}^{-0}(1)$. Suppose that we proceed as in section 88.4. We would first localize the set $\mathcal{T}^{-1}(0)$ with the help of two taboo explorations. Exactly as in (88.11), we would obtain that there exist two disjoint subsets $K^-,L^-$ of $\{1,\dots,N\}$ such that

$$\frac{1}{8d}\,|\mathcal{T}^{-0}(1)|\leq\Big|\Big(\bigcup_{k\in K^-}\mathcal{B}\big(w_k,D\setminus\widehat{\mathcal{T}}_k^-(1),1\big)\Big)\cap\Big(\bigcup_{\ell\in L^-}\mathcal{B}\big(w_\ell,D\setminus\overline{\mathcal{T}}_\ell(0),0\big)\Big)\Big|$$
$$+\,\big|\overline{\mathcal{T}}(0)\big|\,,\quad(92.35)$$

where the set $\widehat{\mathcal{T}}_k^-(1)$ is defined as

$$\forall k\in K^-\qquad\widehat{\mathcal{T}}_k^-(1)\,=\,\overline{\mathcal{T}}_k(0)\cup\mathcal{T}_k^{-1}(1)\,.$$

In order to exploit the inequality (92.35), we have to approximate the two taboo explorations by two genuine ones. The approximation of the second taboo exploration, associated to the taboo sets $\overline{\mathcal{T}}_\ell(0)$, is provided by the two inclusions (92.34). Let us focus on the first exploration in (92.35). We would like to control the discrepancy between the two sets

$$\bigcup_{k\in K^-}\mathcal{B}\big(w_k,D\setminus\widehat{\mathcal{T}}_k^-(1),1\big)\,,\qquad \widehat{\mathcal{B}}_M\Big(\mathrm{Shell}\,\big(W(K^-),D\big),D,1\Big)\,.\tag{92.36}$$

A discrepancy between the two sets can have two origins:

• it might come from the first set exploring a cluster which contains more than $M$ open sites. On the event $\mathcal{F}_{\mathrm{fi}}(n,M)$, any such cluster has to be a boundary cluster. Yet the initial set $W(K^-)$ contains a site of each cluster of the collection $\mathcal{C}^*(h+1,n)$, hence the boundary clusters which might create a discrepancy do not intersect $\Lambda(h+1,n)$.

• or it might come from the second set exploring the cluster of a relevant taboo site at time 1. As these explorations are truncated, these sites must be either closed taboo sites at time 1, or at lattice distance at most $M$ from the set of the open taboo sites at time 1, hence in the set

$$\bigcup_{k\in K^-}\Big(\overline{\mathcal{T}}_k(0)\cup B_1\big(\mathcal{T}_k^{-1}(1),M\big)\Big)\,\subset\,\overline{\mathcal{T}}(0)\cup B_1\big(\mathcal{T}^{-1}(1),M\big)\,.$$

We have no mechanism to control the size of the boundary clusters which do not belong to $\mathcal{C}^*(h+1,n)$, so to avoid dealing with the first source of discrepancy, we will focus on the traces of the sets (92.36) in $\Lambda(h+1,n)$. Instead of two inclusions in the whole box $\Lambda(4n)$ as we had for $t=0$ in (92.34), we have only two inclusions for the traces of the sets in the box $\Lambda(h+1,n)$, i.e.,

$$\begin{multlined}\bigcup_{k\in K^-}\mathcal{B}\big(w_k,D\setminus\widehat{\mathcal{T}}_k^-(1),1\big)\cap\Lambda(h+1,n)\\ \subset\,\widehat{\mathcal{B}}_M\Big(\mathrm{Shell}\,\big(W(K^-),D\big),D,1\Big)\cap\Lambda(h+1,n)\,\subset\\ \Bigg(\bigcup_{k\in K^-}\mathcal{B}\big(w_k,D\setminus\widehat{\mathcal{T}}_k^-(1),1\big)\cup\overline{\mathcal{T}}^0(0)\cup B_1\big(\mathcal{T}^{-1}(1),M\big)\Bigg)\cap\Lambda(h+1,n)\,.\end{multlined}\tag{92.37}$$

Since we can only control the traces of the sets on $\Lambda(h+1,n)$, we cannot exploit the localization inequality (92.35). Fortunately, we have an analog of (92.35) for the traces of the sets on $\Lambda(h+1,n)$, which is

$$\begin{multlined}\frac{1}{8d}\,\big|\mathcal{T}^{-0}(1)\cap\Lambda(h+1,n)\big|\,\leq\,\big|\overline{\mathcal{T}}(0)\cap\Lambda(h+1,n)\big|+\\ \Big|\Big(\bigcup_{k\in K^-}\mathcal{B}\big(w_k,D\setminus\widehat{\mathcal{T}}_k^-(1),1\big)\Big)\cap\Big(\bigcup_{\ell\in L^-}\mathcal{B}\big(w_\ell,D\setminus\overline{\mathcal{T}}_\ell(0),0\big)\Big)\cap\Lambda(h+1,n)\Big|\,.\end{multlined}\tag{92.38}$$

We use now the inclusions (92.34) and (92.37) in order to approximate the two taboo explorations of the right-hand side of (92.38). We obtain

$$\frac{1}{8d}\left|\mathcal{T}^{-0}(1)\cap\Lambda(h+1,n)\right| \leq \left|\overline{\mathcal{T}}(0)\cap\Lambda(h+1,n)\right|+ \\ \left|\widehat{\mathcal{B}}_M\Big(\text{Shell}\big(W(K^-),D\big),D,1\Big)\cap\text{Shell}\big(W(L^-),D\big)\cap\Lambda(h+1,n)\right| \\ +\left|\Big(\overline{\mathcal{T}}^0(0)\cup B_1\big(\mathcal{T}^{-1}(1),M\big)\Big)\cap\Lambda(h+1,n)\right|. \quad (92.39)$$

The first and third term of the right-hand side of (92.39) have been controlled. The second term is the intersection of two genuine explorations, hence we can control it as usual. Thus the inequality (92.39) provides the desired control on the set $\mathcal{T}^{-0}(1)\cap\Lambda(h+1,n)$. While we can control $\mathcal{T}^0(0)$ and $\mathcal{T}^{-1}(1)$ in the whole box $\Lambda(4n)$, we can only manage to control the trace of $\mathcal{T}^{-0}(1)$ on $\Lambda(h+1,n)$. This is a major difference with the scheme of section 90. The bad news is that things will get worse for the subsequent taboo sets. We move on and we try to control the set $\mathcal{T}^{+1}(1)$. Lemma 88.1 tells that any site in $\mathcal{T}^{+1}(1)$ has a neighbour which is in $\mathcal{T}^{-0}(1)$. However, we have only a control on the trace of $\mathcal{T}^{-0}(1)$ on $\Lambda(h+1,n)$. Hence the application of the lemma yields only a control on the trace of $\mathcal{T}^{+1}(1)$ on the interior of $\Lambda(h+1,n)$, as follows:

$$\left|\mathcal{T}^{+1}(1)\cap\overset{\circ}{\widehat{\Lambda(h+1,n)}}\right| \leq 2d\left|\mathcal{T}^{-0}(1)\cap\Lambda(h+1,n)\right|.$$

At this point, we have gained control of the following sets:

$$\mathcal{T}^0(0)\,,\quad \mathcal{T}^{-1}(1)\,,\quad \mathcal{T}^{-0}(1)\cap\Lambda(h+1,n)\,,\quad \mathcal{T}^{+1}(1)\cap\overset{\circ}{\widehat{\Lambda(h+1,n)}}\,. \quad (92.40)$$

To complete the control of the relevant taboo sets of the second step, it remains to control $\mathcal{T}^{+0}(1)$. Suppose again that we proceed as in section 88.4. We would localize this set with the help of two taboo explorations. Exactly as in (88.24), we would obtain that there exist two disjoint subsets $K^+,L^+$ of $\{1,\dots,N\}$ such that

$$\frac{1}{8d}\left|\mathcal{T}^{+0}(t)\right| \leq \left|\Big(\bigcup_{k\in K^+}\mathcal{B}\big(w_k,D\setminus\widehat{\mathcal{T}}_k^*(t),t\big)\Big)\cap\Big(\bigcup_{\ell\in L^+}\mathcal{B}\big(w_\ell,D\setminus\widehat{\mathcal{T}}_\ell^*(t),t\big)\Big)\right|, \quad (92.41)$$

where the set $\widehat{\mathcal{T}}_k^*(t)$ is defined as

$$\widehat{\mathcal{T}}_k^*(t) = \overline{\mathcal{T}}_k(t-1)\cup\mathcal{T}_k^{-0}(t)\cup\mathcal{T}_k^{*1}(t)\,.$$

In order to exploit the inequality (92.41), we have to approximate the two taboo explorations by two genuine ones. These two taboo explorations are of the same nature. Since $\widehat{\mathcal{T}}_k^*(1)$ contains $\widehat{\mathcal{T}}_k^-(1)$, we have directly the inclusions

$$\forall k\in K^+\quad \mathcal{B}\big(w_k,D\setminus\widehat{\mathcal{T}}_k^*(1),1\big)\subset\mathcal{B}\big(w_k,D\setminus\widehat{\mathcal{T}}_k^-(1),1\big)\,. \quad (92.42)$$

To carry out the approximation, we would rely on lemma 89.27, which gives that, for any $k$ in $K^+$,

$$\mathcal{B}\big(w_k,D\setminus\widehat{\mathcal{T}}_k^-(1),1\big)\subset\mathcal{B}\big(w_k,D\setminus\widehat{\mathcal{T}}_k^*(1),1\big)\cup\mathcal{T}_k^{-0}(1)\cup B_1\big(\mathcal{T}_k^{+1}(1),M\big)\,. \quad (92.43)$$

The inclusions (92.42) and (92.43) show that, in order to control the discrepancy between the sets $\mathcal{B}\big(w_k, D\setminus\widehat{\mathcal{T}}_k^-(1),1\big)$ and $\mathcal{B}\big(w_k, D\setminus\widehat{\mathcal{T}}_k^*(1),1\big)$, we need to control the taboo sets $\mathcal{T}_k^{-0}(1)$ and $\mathcal{T}_k^{+1}(1)$. We run here again into a serious problem, because we could only control the sets listed in (92.40). So the only hope at this point is to try to carry out the approximation in a restricted box. This will require a variant of lemma 89.27. Recalling the definition (92.18), and setting

$$\ell \;=\; 2n+2(h+1)\lambda(n)\,,$$

we have

$$\Lambda(h+1,n)=\Lambda(\ell)\,,\qquad \overset{\circ}{\widehat{\Lambda(h+1,n)}} \;=\; \Lambda\big(\ell-2\big)\,.$$

Instead of the inclusion (92.43), we will show that, for any $k$ in $K^+$,

$$\begin{gathered}\mathcal{B}\big(w_k, D\setminus\widehat{\mathcal{T}}_k^-(1),1\big)\cap\Lambda\big(\ell-2-2M\big)\;\subset\\ \Big(\mathcal{B}\big(w_k, D\setminus\widehat{\mathcal{T}}_k^*(1),1\big)\cup\mathcal{T}_k^{-0}(1)\cup B_1\big(\mathcal{T}_k^{+1}(1),M\big)\Big)\cap\Lambda\big(\ell-2-2M\big)\,.\end{gathered}\tag{92.44}$$

The point of (92.44) is that

$$\mathcal{T}_k^{-0}(1)\cap\Lambda\big(\ell-2-2M\big)\;\subset\;\mathcal{T}_k^{-0}(1)\cap\Lambda(h+1,n)\,,\tag{92.45}$$

$$B_1\big(\mathcal{T}_k^{+1}(1),M\big)\cap\Lambda\big(\ell-2-2M\big)\;\subset\;B_1\big(\mathcal{T}_k^{+1}(1)\cap\Lambda(\ell-2),M\big)\,.\tag{92.46}$$

Now the quantities on the right-hand side of (92.45) and (92.46) have been controlled, hence we can manage to control the discrepancy between the two sets

$$\mathcal{B}\big(w_k, D\setminus\widehat{\mathcal{T}}_k^-(1),1\big)\cap\Lambda\big(\ell-2-2M\big)\,,\quad \mathcal{B}\big(w_k, D\setminus\widehat{\mathcal{T}}_k^*(1),1\big)\cap\Lambda\big(\ell-2-2M\big)\,.$$

Furthermore, we can reuse the approximation of the first of the two sets above provided by the two inclusions (92.37). As we can only approximate the traces of the taboo explorations on the box $\Lambda\big(\ell-2-2M\big)$, instead of the inequality (92.41), we will rely on the following restricted variant: there exist two disjoint subsets $K^+, L^+$ of $\{1,\dots,N\}$ such that

$$\begin{gathered}\frac{1}{8d}\big|\mathcal{T}^{+0}(1)\cap\Lambda\big(\ell-2-2M\big)\big|\;\le\\ \Big|\Big(\bigcup_{k\in K^+}\mathcal{B}\big(w_k, D\setminus\widehat{\mathcal{T}}_k^*(1),1\big)\Big)\cap\Big(\bigcup_{\ell\in L^+}\mathcal{B}\big(w_\ell, D\setminus\widehat{\mathcal{T}}_\ell^*(1),1\big)\Big)\cap\Lambda\big(\ell-2-2M\big)\Big|\,.\end{gathered}$$

In the end, we will manage to control the set $\mathcal{T}^{+0}(1)\cap\Lambda\big(\ell-2-2M\big)$. We observe that, during each step, we have to reduce the size of the box where the control of the relevant taboo sets can be carried out. The reason is that, the errors in the approximations must be propagated in time as well as in space. The next sections are devoted to the precise implementation of an iterative approximation scheme, which provides a control on all the sets of the relevant taboo sites.

## 92.5 The successive rings of relevant taboo sites

Heuristically, we infer from the observation of the second step $t = 1$ that, once we have controlled the relevant taboo sets at step $t$ inside a box $\Lambda(\ell)$, we can control the relevant taboo sets at step $t+1$ inside a smaller box, like $\Lambda(\ell - 2 - 2M)$. Notice also that the step $t = 0$ is a bit special, because we have a control on $\mathcal{T}(0)$ and $\mathcal{T}^{-1}(1)$ in the full box $\Lambda(4n)$, so it is to be expected that each step requires to reduce even more the control box. The next challenge is to build the general scheme for estimating all the subsequent taboo sets until the termination time. To do so, we introduce adequate boxes and sets of taboo sites, precisely tuned to match the approximations carried out during each step. We set

$$\ell(n) \;=\; 2n + 2(h-1)\lambda(n) + 8T(M+2)\,,$$

and we introduce a sequence of boxes $\Gamma = \big(\Gamma(t), 0 \leq t \leq T\big)$, defined by

$$\forall t \in \{\,0,\dots,T\,\} \qquad \Gamma(t) \;=\; \Lambda\big(\ell(n) - 8t(M+2)\big)\,.$$

This sequence of boxes is tailored so that $\Gamma(T)$ coincides with $\Lambda(h-1,n)$. Furthermore, we have $\Gamma(0) = \Lambda\big(\ell(n)\big)$. We adjust the parameters to ensure that $\Gamma(0)$ is included in $\Lambda(h,n)$, i.e., we require that

$$4T(M+2) \;\leq\; \lambda(n)\,. \tag{92.47}$$

The occurrence of the event $\mathcal{G}(n)$ defined in (92.13) implies that

$$T \;\leq\; \bar{t}(n) \;=\; 2(\ln n)^{1-\gamma}\,. \tag{92.48}$$

The value $M$ was defined in (92.6). Taking into account (92.48) and (92.6), we see that a sufficient condition ensuring (92.47) is

$$8(\ln n)^{1-\gamma}\Big((4d\ln n)^{1/\alpha} + 2\Big) \;\leq\; \lambda(n)\,.$$

For $\Lambda$ a cubic box, we denote by $\overset{o}{\Lambda}$ its interior, i.e., $\overset{o}{\Lambda} = \Lambda \setminus \partial^{in}\Lambda$. Taking the interior of a box amounts to reducing by two its side length, i.e.,

$$\forall n \geq 1 \qquad \overset{\circ}{\widehat{\Lambda(2n)}} \;=\; \Lambda(2n-2)\,.$$

We define intermediate boxes, by setting, for $0 \leq t \leq T$,

$$\begin{aligned}
\widehat{\Gamma}^-(t) \;&=\; \Lambda\big(\ell(n) - 8(t-1)(M+2) - 4(M+2)\big)\,,\\
\widehat{\Gamma}^*(t) \;&=\; \Lambda\big(\ell(n) - 8(t-1)(M+2) - 6(M+2)\big)\,.
\end{aligned}$$

We have the following sequence of inclusions:

$$\Gamma(0) \supset \overset{\circ}{\widehat{\Gamma(0)}} \supset \widehat{\Gamma}^-(1) \supset \overset{\circ}{\widehat{\widehat{\Gamma}^-(1)}} \supset \cdots \supset \widehat{\Gamma}^-(T) \supset \overset{\circ}{\widehat{\widehat{\Gamma}^-(T)}} \supset \widehat{\Gamma}^*(T) \supset \Gamma(T)\,.$$

We will develop iterative inequalities to control the relevant taboo sets in the boxes of the sequence $\Gamma = \big(\Gamma(t), 0 \leq t \leq T\big)$. The two main constraints that are needed in the argument are encountered in (92.115) and (92.166), they are

$$B_1\big(\widehat{\Gamma}^-(t), 2M\big) \;\subset\; \overset{\circ}{\widehat{\Gamma(t-1)}}\,, \qquad B_1\big(\widehat{\Gamma}^*(t), M-1\big) \;\subset\; \overset{\circ}{\widehat{\widehat{\Gamma}^-(t)}}\,.$$

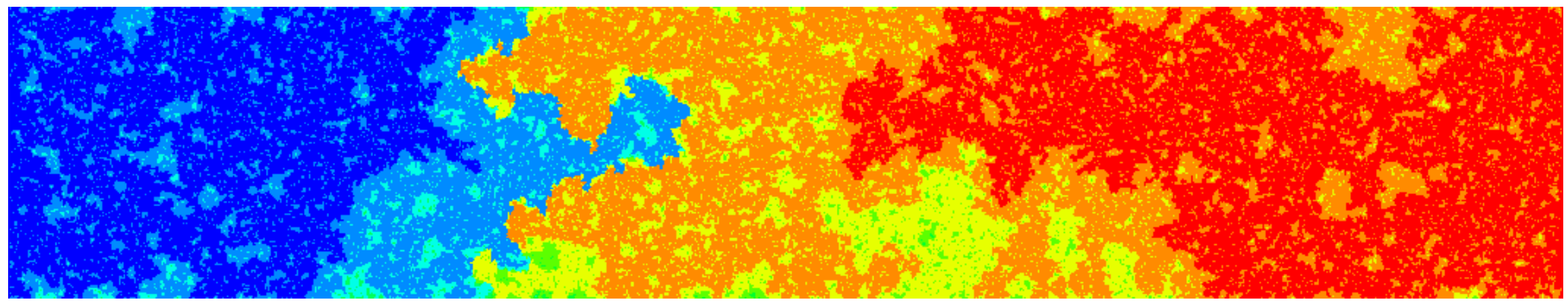

$\Gamma(t-1) = \Lambda\big(2n + 2(h-1)\lambda(n) + 8(T-t+1)M\big)$

$\widehat{\Gamma}^-(t)$

$2M+3$

$\widehat{\Gamma}^*(t)$

$M+1$

$\Gamma(t)$

$2M+4$ $M+2$ $M+2$ $O$

$\overset{\circ}{\overparen{\widehat{\Gamma}^-(t)}}$

$\overset{\circ}{\overparen{\Gamma(t-1)}}$

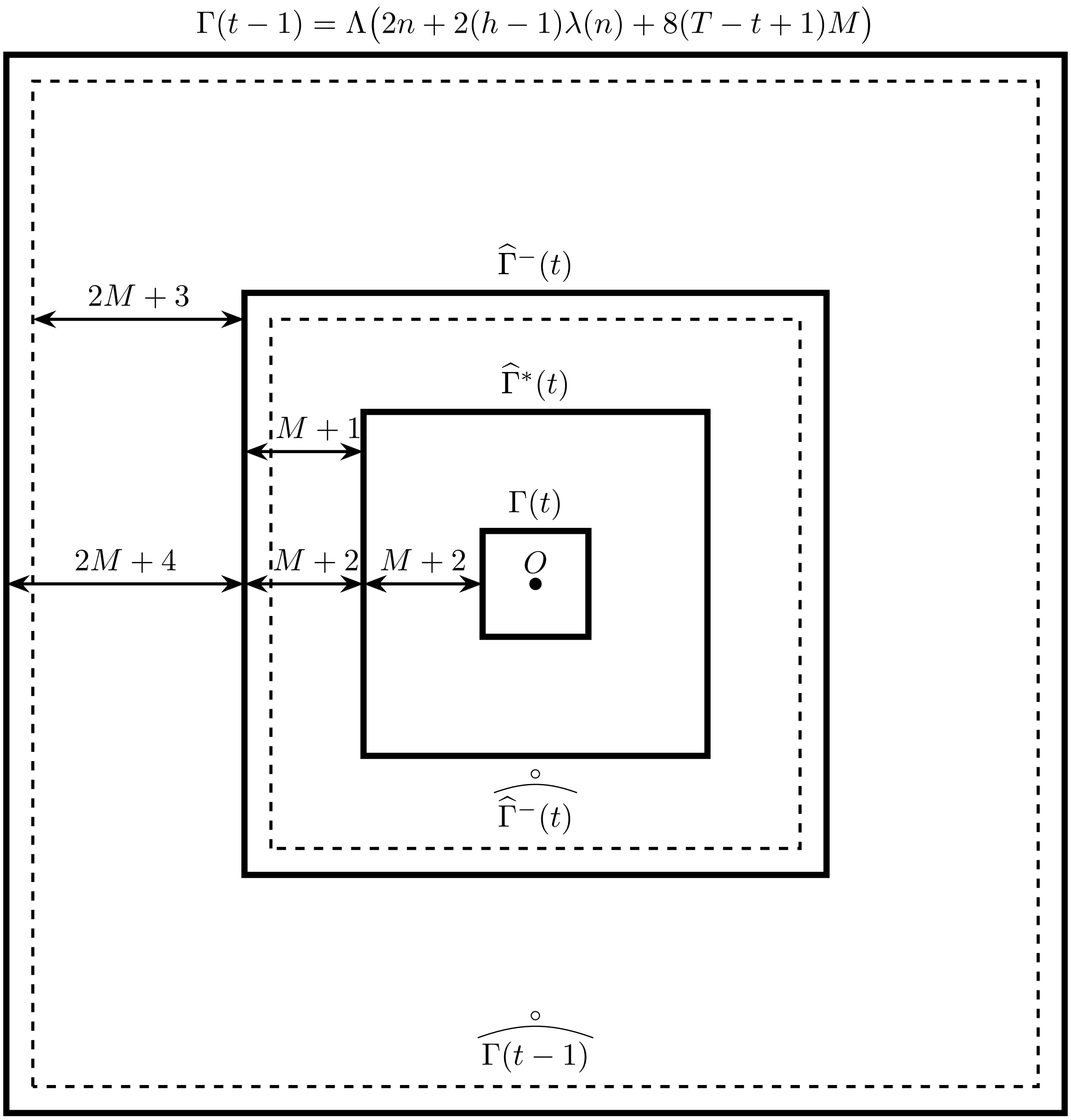

Figure 169: The boxes $\Gamma(t-1)$, $\widehat{\Gamma}^-(t)$, $\widehat{\Gamma}^*(t)$, $\Gamma(t)$, $\overset{\circ}{\overparen{\Gamma(t-1)}}$, $\overset{\circ}{\overparen{\widehat{\Gamma}^-(t)}}$.
Top and bottom: Intertwined explorations in $1024 \times 192$, $p = 0.505$.

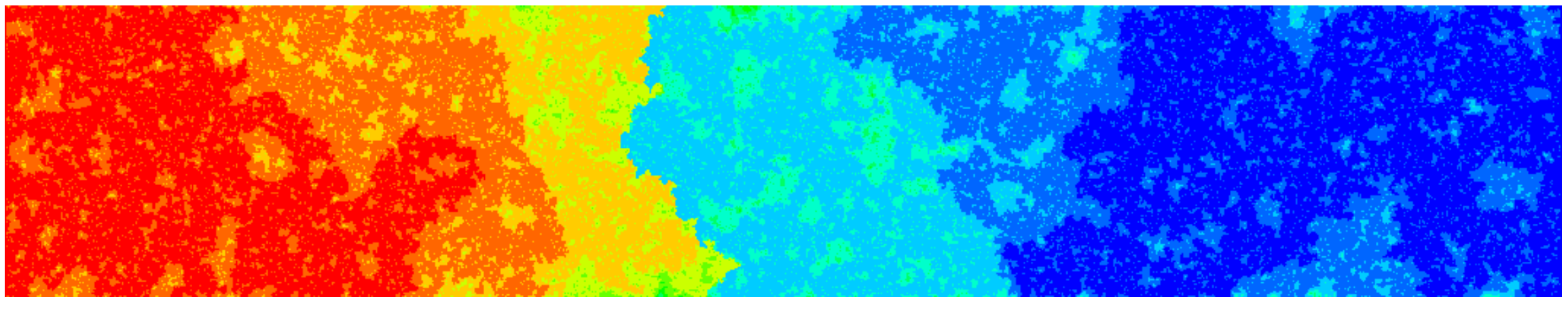

## 92.6 Intersections of truncated explorations

In the 3/2-argument, the strategy consisted in detecting a leap in the sequence of the relevant taboo sets, and in identifying two genuine truncated explorations whose intersection contains the excess of closed sites responsible for this leap. The crucial control was obtained by applying Hoeffding's inequality to these two genuine truncated explorations. These final steps were carried out in the subsection 90.2. The detection of the leap in subsection 87.4 was a bit technical, it involved the introduction of the constant $c$ and the auxiliary function $\psi(n)$, which are adequately tuned after several lengthy computations.

We shall proceed differently here. We introduce a new typical event right before starting the computations, which provides a control on all the genuine truncated explorations that might ever be used in the argument. On this event, and with this control, we obtain a chain of deterministic inequalities between the successive sets of the relevant taboo sites $(\big|\overline{\mathcal{T}}(t)\big|, 0\leq t\leq T)$ in the sequence of the nested boxes $\Gamma = \big(\Gamma(t), 0\leq t\leq T\big)$. In a nutshell, instead of applying Hoeffding's inequality right at the end of a sequence of complicated inclusions and inequalities, we applied it from the beginning to all the explorations that might potentially be involved in the localization of the relevant taboo sites. The argument is conceptually the same, but we hope that this different presentation sheds new light on the mechanism at work.

The new event depends only on the parameter $N$, which is the number of boundary clusters chosen to start the intertwined explorations. It involves the genuine truncated explorations introduced in subsection 89.2. For $A$ a subset of $D$ and $t\geq 0$, we define

$$\mathcal{H}(A,t) \;=\; \bigcup_{x\in A}\widehat{\mathcal{B}}_M\big(\text{Shell}\,(x,D),D,t\big) \;=\; \widehat{\mathcal{B}}_M\Big(\text{Shell}\,\big(A,D\big),D,t\Big)\,. \qquad (92.49)$$

We introduce now the event Intersect$(N)$, that provides a uniform control on the functional $S$ for the intersections of two genuine truncated explorations whose starting sets have a cardinality less than or equal to $N$. It is defined as

$$\begin{aligned}\text{Intersect}(N) \;=\; \Big\{\,&\forall A,B\subset\Lambda(4n)\,,\quad |A|\leq N\,,\;|B|\leq N\,,\\ &\forall m\in\{\,1,\dots,4n\,\}\,,\quad \forall s,t\in\{\,0,\dots,2\ln n\,\}\\ &\qquad S\big(\mathcal{H}(A,s)\cap\mathcal{H}(B,t)\cap\Lambda(m)\big)\;>\;-c_H n^{d/2}\sqrt{N\ln n}\,\Big\}\,,\end{aligned}$$

where $c_H$ is the constant given by

$$c_H \;=\; c_H(d,p) \;=\; \frac{5^{d+1}}{p(1-p)}\,. \qquad (92.50)$$

To be honest, the event Intersect$(N)$ depends also on $n$ and $p$. However, we use the simple notation Intersect$(N)$ to emphasize the dependence on $N$. Later on, we will use it with an adequately chosen value of $N$, that depends on the number of boundary clusters we wish to control. Let us compute the probability of the complementary event $\text{Intersect}(N)^c$. By a standard union bound, we have

$$P\Big(\text{Intersect}(N)^c\Big) \leq \sum_{\substack{A,B\subset\Lambda(4n)\\ |A|\leq N\,,|B|\leq N}} \sum_{\substack{1\leq m\leq 4n\\ 0\leq s,t\leq 2\ln n}} P\Big(S\Big(\mathcal{H}(A,s)\cap\mathcal{H}(B,t)\cap\Lambda(m)\Big) \leq -c_H n^{d/2}\sqrt{N\ln n}\Big)\,. \tag{92.51}$$

Let us fix some elements $A, B, m, s, t$ appearing in the above sum. It follows from corollary 89.17 that the set $\mathcal{H}(A, s)$ is the output of genuine truncated explorations. The same is true for the set $\mathcal{H}(B, t)$. By proposition 89.11, the set $\mathcal{H}(A, s) \cup \mathcal{H}(B, t)$, being the union of two genuine sets, is also a genuine set. Yet we are interested in their traces on $\Lambda(m)$.

**Lemma 92.3.** *Let $\mathcal{E}$ be the outcome of a genuine exploration algorithm in $D$. For any subset $A$ of $D$, the trace $\mathcal{E} \cap A$ of $\mathcal{E}$ on $A$ enjoys the probabilistic control* (13.5)*.*

*Proof.* We use the definition 89.1 of a genuine algorithm: there exists finite sequence of decision rules $\Phi = (\phi_k, 1 \leq k \leq \ell)$, where, for $1 \leq k \leq \ell$, $\phi_k$ is a map defined on $D^k \times \{0, 1\}^k$ with values in $D_\dagger$. Let $x$ be the starting point of the algorithm. We define then a sequence of random variables $(X_k, 0 \leq k \leq T)$ by setting $X_1 = x$, and for $k \geq 1$,

$$X_{k+1} = \phi_k\big(X_1, \dots, X_k, \omega(X_1), \dots, \omega(X_k)\big)\,.$$

The algorithm stops at the termination $T$ defined by

$$T = \min\big\{\, k \geq 0 : X_{k+1} = \dagger \,\big\}\,.$$

The set $\mathcal{E}$ is then equal to

$$\mathcal{E} = \big\{\, X_1, \dots, X_T \,\big\} = \big\{\, \widehat{X}_1, \dots, \widehat{X}_T \,\big\}\,,$$

where $\widehat{X}_1, \dots, \widehat{X}_{\widehat{T}}$ is the cleansing of the sequence $X_1, \dots, X_T$. We remove from the sequence $\widehat{X}_1, \dots, \widehat{X}_{\widehat{T}}$ the sites that do not belong to $A$, and we obtain a subsequence $Y_1, \dots, Y_S$. The trace $\mathcal{E} \cap A$ of $\mathcal{E}$ on $A$ is then

$$\mathcal{E} \cap A = \big\{\, \widehat{X}_1, \dots, \widehat{X}_T \,\big\} \cap A = \big\{\, Y_1, \dots, Y_S \,\big\}\,.$$

With a standard recurrence, we can prove that, for any $i \leq S$, the variables $\omega(Y_1), \dots, \omega(Y_i)$ are i.i.d. Bernoulli with parameter $p$. Indeed, each time the exploration algorithm visits a new site $Y_i$ in $A$, it has no information on its status $\omega(Y_i)$, which is distributed according to an independent Bernoulli random variable. Ideally, we would like to say that $\mathcal{E} \cap A$ is a genuine set. However, the random time $S$ is not a stopping time, as required by the definition 89.4. $\square$

By lemma 92.3, the probabilistic control (13.5) can be applied to the trace on $\Lambda(m)$ of $\mathcal{H}(A, s)$, $\mathcal{H}(B, t)$, $\mathcal{H}(A, s) \cup \mathcal{H}(B, t)$. We evaluate next the functional

$S$ on these sets. Thanks to the additivity of $S$, we have

$$S\Big(\big(\mathcal{H}(A,s)\cup\mathcal{H}(B,t)\big)\cap\Lambda(m)\Big)+S\Big(\big(\mathcal{H}(A,s)\cap\mathcal{H}(B,t)\big)\cap\Lambda(m)\Big)$$
$$= S\big(\mathcal{H}(A,s)\cap\Lambda(m)\big)+S\big(\mathcal{H}(B,t)\cap\Lambda(m)\big)\,. \quad (92.52)$$

It follows from the equality (92.52) that, for any $\lambda\geq 0$,

$$P\Big(S\Big(\big(\mathcal{H}(A,s)\cap\mathcal{H}(B,t)\big)\cap\Lambda(m)\Big)\leq -\lambda\Big)\;\leq\; P\Big(S\big(\mathcal{H}(A,s)\cap\Lambda(m)\big)\leq -\frac{\lambda}{3}\Big)$$
$$+P\Big(S\big(\mathcal{H}(B,t)\cap\Lambda(m)\big)\leq -\frac{\lambda}{3}\Big)+P\Big(S\Big(\big(\mathcal{H}(A,s)\cup\mathcal{H}(B,t)\big)\cap\Lambda(m)\Big)\geq\frac{\lambda}{3}\Big)\,.$$

Using the probabilistic control (13.5), we have, for any $\lambda\geq 0$,

$$P\Big(S\Big(\big(\mathcal{H}(A,s)\cap\mathcal{H}(B,t)\big)\cap\Lambda(m)\Big)\leq -\lambda\Big)\;\leq\; 6|\Lambda(4n)|\exp\Big(-\frac{2\big(p(1-p)\lambda\big)^2}{9|\Lambda(4n)|}\Big)\,. \quad (92.53)$$

We apply the inequality (92.53) with the value

$$\lambda(n)\;=\;c_H n^{d/2}\sqrt{N\ln n}\;=\;\frac{5^{d+1}}{p(1-p)}n^{d/2}\sqrt{N\ln n}\,.$$

Using the simple bound $|\Lambda(4n)|\;\leq\;5^d n^d$, we see that

$$\frac{2\big(p(1-p)\lambda(n)\big)^2}{9|\Lambda(4n)|}\;=\;\frac{2\,5^{2d+2}}{9|\Lambda(4n)|}n^d N\ln n\;\geq\;5^{d+1}N\ln n\,,$$

from which we deduce that

$$P\Big(S\Big(\big(\mathcal{H}(A,s)\cap\mathcal{H}(B,t)\big)\cap\Lambda(m)\Big)\leq -c_H n^{d/2}\sqrt{N\ln n}\Big)$$
$$\leq\;6|\Lambda(4n)|\exp\big(-5^{d+1}N\ln n\big)\,. \quad (92.54)$$

This estimate is uniform with respect to the sets $A,B$ and the times $s,t$. From now onwards, we suppose that $n\geq 5$. We substitute (92.54) in (92.51), and we get

$$P\big(\mathrm{Intersect}(N)^c\big)\;\leq\;24n(2\ln n+1)^2|\Lambda(4n)|^{2N+1}\exp\big(-5^{d+1}N\ln n\big)\,.$$

As $n\geq 5$, we have

$$24n(2\ln n+1)^2\;\leq\;|\Lambda(4n)|=(4n+1)^d\,,$$
$$(2N+2)\ln\big(|\Lambda(4n)|\big)\;\leq\;d(2N+2)\ln(4n+1)\;\leq\;8dN\ln n\,,$$

whence, recalling that we work with $d\geq 3$,

$$P\big(\mathrm{Intersect}(N)^c\big)\;\leq\;\exp\Big((2N+2)\ln\big(|\Lambda(4n)|\big)-5^{d+1}N\ln n\Big)$$
$$\leq\;\exp\big((8d-5^{d+1})N\ln n\big)\;\leq\;\exp\big(-4^dN\ln n\big)\,. \quad (92.55)$$

## 92.7 Control of $\mathcal{T}^{*0}(0) \cap \Gamma(0)$

Recall that $\mathcal{D}$ is the event that $w_1, \dots, w_N$ are pairwise disconnected in $\Lambda(4n)$. We denote also the box $\Lambda(4n)$ by $D$. On the event $\mathcal{D}$, we have

$$\mathcal{T}^{*0}(0) \,=\, \bigcup_{1 \le k < \ell \le N} \partial^{\,out} C(w_k, D) \cap \partial^{\,out} C(w_\ell, D)\,.$$

We could control the whole set $\mathcal{T}^{*0}(0)$ in $D$. However, for the scheme we will implement, it is enough to control the set $\mathcal{T}^{*0}(0) \cap \Gamma(0)$. To do so, we apply next lemma 84.3 to the family of sets

$$E_i \,=\, F_i \,=\, \partial^{\,out} C(w_i, D) \cap \Gamma(0)\,, \quad 1 \le i \le N\,.$$

The collection

$$\partial^{\,out} C(w_i, D) \cap \partial^{\,out} C(w_j, D) \cap \Gamma(0)\,, \quad 1 \le i, j \le N\,, i \ne j\,,$$

has multiplicity less than $4d^2$. By lemma 84.3, there exist two disjoint subsets $K(0), L(0)$ of $\{\,1, \dots, N\,\}$ such that

$$\Big|\Big(\bigcup_{k \in K(0)} \partial^{\,out} C(w_k, D)\Big) \cap \Big(\bigcup_{\ell \in L(0)} \partial^{\,out} C(w_\ell, D)\Big) \cap \Gamma(0)\Big| \ge \frac{1}{16d^2} \big|\mathcal{T}^{*0}(0) \cap \Gamma(0)\big|. \tag{92.56}$$

The set of the left-hand side of (92.56) is the intersection of two genuine explorations, hence its cardinality can be controlled as usual. However, in order to streamline the proof, we will do it with the help of the event Intersect$(N)$. In fact, we simply remark that the previous sets can be represented as follows. We define the sets $A, B$ by

$$A \,=\, \big\{\, w_k : k \in K(0)\,\big\}\,, \quad B \,=\, \big\{\, w_\ell : \ell \in L(0)\,\big\}\,,$$

and we have then

$$\Big(\bigcup_{k \in K(0)} \partial^{\,out} C(w_k, D)\Big) \cap \Big(\bigcup_{\ell \in L(0)} \partial^{\,out} C(w_\ell, D)\Big) \cap \Gamma(0)$$
$$= \text{Shell}\,(A, D) \cap \text{Shell}\,(B, D) \cap \Gamma(0) \,=\, \mathcal{H}\big(A, 0\big) \cap \mathcal{H}\big(B, 0\big) \cap \Gamma(0)$$

(the notation $\mathcal{H}(\cdot, 0)$ was introduced in (92.49)). On the event Intersect$(N)$, we have automatically

$$S\big(\mathcal{H}(A, 0) \cap \mathcal{H}(B, 0)\big) \,>\, -c_H n^{d/2} \sqrt{N \ln n}\,. \tag{92.57}$$

As we have discussed several times, the sites in $\mathcal{H}(A, 0) \cap \mathcal{H}(B, 0)$ are pivotal sites for the event $\mathcal{D}$, they are all closed, so that

$$S\big(\mathcal{H}(A, 0) \cap \mathcal{H}(B, 0)\big) \,=\, -\frac{1}{1-p} \big|\mathcal{H}(A, 0) \cap \mathcal{H}(B, 0)\big|\,. \tag{92.58}$$

We deduce from (92.56), (92.57) and (92.58) that

$$\big|\mathcal{T}^{*0}(0) \cap \Gamma(0)\big| \,\le\, 16d^2 c_H n^{d/2} \sqrt{N \ln n}\,. \tag{92.59}$$

This provides the desired control on $\mathcal{T}^{*0}(0) \cap \Gamma(0)$.

## 92.8 Control of $\mathcal{T}^{-1}(t) \cap \overset{\circ}{\widehat{\Gamma(t-1)}}$

We use lemma 88.1 in order to control the set $\mathcal{T}^{-1}(t)$. There is a little complication compared to subsection 88.3. The lemma guarantees that an open relevant taboo site has a neighbour which is a closed relevant taboo site. Once we have controlled the trace of the set $\mathcal{T}^{+0}(t-1)$ in the box $\Gamma(t-1)$, we can only control the trace of the set $\mathcal{T}^{-1}(t)$ on the interior of $\Gamma(t-1)$. More precisely, we obtain the following inequality:

$$\left|\mathcal{T}^{-1}(t) \cap \overset{\circ}{\widehat{\Gamma(t-1)}}\right| \leq 2d\left|\mathcal{T}^{+0}(t-1) \cap \Gamma(t-1)\right|. \tag{92.60}$$

## 92.9 Localization of $\mathcal{T}^{-0}(t) \cap \widehat{\Gamma}^{-}(t)$

We suppose that all the relevant taboo sites at time $t-1$ have already been controlled in the box $\widehat{\Gamma}^{-}(t)$, and we try here to localize the set $\mathcal{T}^{-0}(t) \cap \widehat{\Gamma}^{-}(t)$. By the definition (80.39), a taboo site in $\mathcal{T}^{-0}(t)$ is visited at time $t-1$ or $t-2$, therefore

$$\mathcal{T}^{-0}(t) \cap \widehat{\Gamma}^{-}(t) = \bigcup_{1\leq i,j\leq N, i\neq j} \mathcal{T}_i^{-0}(t) \cap \mathcal{A}_j(t-1) \cap \widehat{\Gamma}^{-}(t).$$

Let us consider a site $x$ belonging to $\mathcal{T}_i^{-0}(t) \cap \widehat{\Gamma}^{-}(t)$. By definition, the explorer $i$ attempted to visit $x$ at time $t$. This attempt must originate from a neighbor of $x$, thus one of the neighbors of $x$ must belong to $\mathcal{A}_i(t)$. In addition, the aggregates $(\mathcal{A}_j(t-1), 1 \leq j \leq N)$ are pairwise disjoint, as well as $(\mathcal{A}_i(t), 1 \leq i \leq N)$, thus a site $x$ can belong to at most $2d$ sets of the family

$$\mathcal{T}_i^{-0}(t) \cap \mathcal{A}_j(t-1) \cap \widehat{\Gamma}^{-}(t), \quad 1 \leq i, j \leq N,\, i \neq j.$$

We apply lemma 84.3 to the above family: there exist two disjoint subsets $K^{-}(t), L^{-}(t)$ of $\{1, \dots, N\}$ such that

$$\frac{1}{8d}\left|\mathcal{T}^{-0}(t) \cap \widehat{\Gamma}^{-}(t)\right| \leq \left|\Big(\bigcup_{k\in K^{-}(t)} \mathcal{T}_k^{-0}(t)\Big) \cap \Big(\bigcup_{\ell\in L^{-}(t)} \mathcal{A}_\ell(t-1)\Big) \cap \widehat{\Gamma}^{-}(t)\right|. \tag{92.61}$$

As we supposed that the relevant taboo sites at time $t-1$ have already been controlled in the box $\widehat{\Gamma}^{-}(t)$, we need only to control the relevant taboo sites at time $t$ which are not relevant at time $t-1$. We have the following inclusions:

$$\forall k \in K^{-}(t) \qquad \mathcal{T}_k^{-0}(t) \subset \big(\mathcal{T}_k^{-0}(t) \setminus \overline{\mathcal{T}}_k^{0}(t-1)\big) \cup \overline{\mathcal{T}}_k^{0}(t-1),$$

which, together with (92.61), yield

$$\begin{aligned} \frac{1}{8d}\left|\mathcal{T}^{-0}(t) \cap \widehat{\Gamma}^{-}(t)\right| \leq \left|\overline{\mathcal{T}}^{0}(t-1) \cap \widehat{\Gamma}^{-}(t)\right| + \\ \left|\Big(\bigcup_{k\in K^{-}(t)} \mathcal{T}_k^{-0}(t) \setminus \overline{\mathcal{T}}_k^{0}(t-1)\Big) \cap \Big(\bigcup_{\ell\in L^{-}(t)} \mathcal{A}_\ell(t-1)\Big) \cap \widehat{\Gamma}^{-}(t)\right|. \end{aligned} \tag{92.62}$$

Let us define

$$\forall k \in K^-(t) \qquad \widehat{\mathcal{T}}_k^-(t) \,=\, \overline{\mathcal{T}}_k(t-1) \cup \mathcal{T}_k^{-1}(t)\,.$$

It follows from proposition 81.2 that, for any $k$ in $K^-(t)$,

$$\begin{aligned} \mathcal{T}_k^{-0}(t) \setminus \overline{\mathcal{T}}_k(t-1) \,&\subset\, \mathcal{T}_k^{*0}(t) \setminus \overline{\mathcal{T}}_k(t-1) \\ &\subset\, \mathcal{B}\Big(w_k, D \setminus \big(\overline{\mathcal{T}}_k(t-1) \cup \mathcal{T}_k^{*1}(t)\big), t\Big) \,\subset\, \mathcal{B}\big(w_k, D \setminus \widehat{\mathcal{T}}_k^-(t), t\big)\,. \end{aligned} \quad (92.63)$$

Starting from the inclusion (92.63), we take the trace on $\widehat{\Gamma}^-(t)$ and the union over $k$ in $K^-(t)$, and we get

$$\Big(\bigcup_{k \in K^-(t)} \mathcal{T}_k^{-0}(t) \setminus \overline{\mathcal{T}}_k(t-1)\Big) \cap \widehat{\Gamma}^-(t) \,\subset\, \bigcup_{k \in K^-(t)} \mathcal{B}\big(w_k, D \setminus \widehat{\mathcal{T}}_k^-(t), t\big) \cap \widehat{\Gamma}^-(t)\,. \quad (92.64)$$

The inclusion (92.64) readily implies that

$$\begin{aligned} &\Big(\bigcup_{k \in K^-(t)} \mathcal{T}_k^{-0}(t) \setminus \overline{\mathcal{T}}_k(t-1)\Big) \cap \Big(\bigcup_{\ell \in L^-(t)} \mathcal{A}_\ell(t-1)\Big) \cap \widehat{\Gamma}^-(t) \,\subset \\ &\qquad \Big(\bigcup_{k \in K^-(t)} \mathcal{B}\big(w_k, D \setminus \widehat{\mathcal{T}}_k^-(t), t\big)\Big) \cap \Big(\bigcup_{\ell \in L^-(t)} \mathcal{A}_\ell(t-1)\Big) \cap \widehat{\Gamma}^-(t)\,. \end{aligned} \quad (92.65)$$

We deduce from (92.62) and (92.65) that

$$\begin{aligned} \frac{1}{8d}\big|\mathcal{T}^{-0}(t) \cap \widehat{\Gamma}^-(t)\big| \,&\leq\, \big|\overline{\mathcal{T}}^0(t-1) \cap \widehat{\Gamma}^-(t)\big| + \\ &\Big|\Big(\bigcup_{k \in K^-(t)} \mathcal{B}\big(w_k, D \setminus \widehat{\mathcal{T}}_k^-(t), t\big)\Big) \cap \Big(\bigcup_{\ell \in L^-(t)} \mathcal{A}_\ell(t-1)\Big) \cap \widehat{\Gamma}^-(t)\Big|\,. \end{aligned} \quad (92.66)$$

By proposition 81.1, we have

$$\forall \ell \in L^-(t) \qquad \mathcal{A}_\ell(t-1) \,=\, \mathcal{B}\big(w_\ell, D \setminus \overline{\mathcal{T}}_\ell(t-1), t-1\big)\,. \quad (92.67)$$

Using (92.66), (92.67), we obtain that

$$\begin{aligned} &\frac{1}{8d}\big|\mathcal{T}^{-0}(t) \cap \widehat{\Gamma}^-(t)\big| \,\leq\, \big|\overline{\mathcal{T}}^0(t-1) \cap \widehat{\Gamma}^-(t)\big| + \\ &\Big|\Big(\bigcup_{k \in K^-(t)} \mathcal{B}\big(w_k, D \setminus \widehat{\mathcal{T}}_k^-(t), t\big)\Big) \cap \Big(\bigcup_{\ell \in L^-(t)} \mathcal{B}\big(w_\ell, D \setminus \overline{\mathcal{T}}_\ell(t-1), t-1\big)\Big) \cap \widehat{\Gamma}^-(t)\Big|\,. \end{aligned} \quad (92.68)$$

We conclude that the set $\mathcal{T}^{-0}(t) \cap \widehat{\Gamma}^-(t)$ can be localized with the help of the two taboo explorations

$$\bigcup_{k \in K^-(t)} \mathcal{B}\big(w_k, D \setminus \widehat{\mathcal{T}}_k^-(t), t\big)\,, \qquad \bigcup_{\ell \in L^-(t)} \mathcal{B}\big(w_\ell, D \setminus \overline{\mathcal{T}}_\ell(t-1), t-1\big)\,.$$

The next difficult problem will be to control the size of the intersection of these two sets.

## 92.10 Approximation of $\mathcal{B}\big(w_i, D\setminus\overline{\mathcal{T}}_i(t-1), t-1\big)$

Let us fix a subset $I$ of $\{1,\dots,N\}$. In this subsection, we try to approximate the taboo exploration

$$\bigcup_{i\in I}\mathcal{B}\big(w_i, D\setminus\overline{\mathcal{T}}_i(t-1), t-1\big)$$

by a genuine truncated taboo exploration. We obtained such a result in the whole box in subsubsection 89.6.1. Unfortunately, the result of proposition 89.18 cannot be used here, because the starting set for the intertwined explorations does not contain a vertex of each boundary cluster. As a consequence, the taboo explorations might very well visit a large boundary cluster after the first step, thereby creating a huge discrepancy with the approximating genuine truncated exploration. We have no control on the size of these boundary clusters, thus there is no hope of proving a result as strong as proposition 89.18. However, the starting set contains a vertex of each boundary cluster which intersects the box $\Lambda(h+1,n)$. We know also that the finite clusters have cardinality strictly less than $M$, and this will help to prove a weaker local version of the approximation result. For convenience, we state the proposition for the sets at time $t$ instead of $t-1$. To start with, the error terms should be defined around the relevant taboo sites included in a specific ring associated to the sequence of boxes $\Gamma$. Instead of the sets $B^{t,0}_{1,M}(s)$, $B^{t,1}_{1,M}(s)$ defined in (89.79), we will use the sets $B^{t,0}_{1,\Gamma}(s)$, $B^{t,1}_{1,\Gamma}(s)$, $B^{t,-}_{1,\Gamma}(s)$ defined as follows. For $t$ in $\{0,\dots,T\}$ and $s\geq 0$, we set

$$B^{t,0}_{1,\Gamma}(s) \;=\; B_1\big(\mathcal{T}^{*0}(t)\cap\Gamma(t), s(M+1)\big)\,, \tag{92.69}$$

$$B^{t,1}_{1,\Gamma}(s) \;=\; B_1\big(\mathcal{T}^{*1}(t)\cap\Gamma(t), s(M+1)+M-1\big)\,. \tag{92.70}$$

For $t$ in $\{1,\dots,T\}$ and $s\geq 0$, we set

$$B^{t,-}_{1,\Gamma}(s) \;=\; B_1\big(\overline{\mathcal{T}}(t-1)\cap\Gamma(t), s(M+1)+M-1\big)\,. \tag{92.71}$$

For convenience, we extend the definition (31.31) of the balls $\mathcal{B}(A,D,t)$ to the case where $A$ is not necessarily a subset of $D$: for $A$ any subset of $\mathbb{Z}^d$, we define

$$\mathcal{B}(A,D,t) \;=\; \mathcal{B}(A\cap D,D,t) \;=\; \big\{\,x\in D: T_D(A\cap D,x)\leq t\,\big\}\,.$$

Here is the proper reformulation of proposition 89.18.

**Proposition 92.4.** *We consider a configuration realizing the event $\mathcal{H}'$ defined in (92.27) and a subset $I$ of $\{1,\dots,N\}$. Setting $W(I)=\{\,w_i: i\in I\,\}$, we have, for any $t$ in $\{0,\dots,T-1\}$,*

$$\begin{aligned}\bigcup_{i\in I}\mathcal{B}\big(w_i, D\setminus\overline{\mathcal{T}}_i(t), t\big)\cap\widehat{\Gamma}^-(t+1) \;\subset\; \widehat{\mathcal{B}}_M\Big(\mathit{Shell}\big(W(I),D\big),D,t\Big)\cap\widehat{\Gamma}^-(t+1)\\ \subset\;\Big(\bigcup_{i\in I}\mathcal{B}\big(w_i, D\setminus\overline{\mathcal{T}}_i(t), t\big)\cap\widehat{\Gamma}^-(t+1)\Big)\cup\\ \bigcup_{1\leq s\leq t} B^{s,-}_{1,\Gamma}(t-s+1)\,\cup \bigcup_{0\leq s\leq t}\big(B^{s,0}_{1,\Gamma}(t-s)\,\cup\, B^{s,1}_{1,\Gamma}(t-s+1)\big)\,.\end{aligned} \tag{92.72}$$

*Proof.* We shall prove the proposition by induction on $t$ and we start with $t = 0$. Let us fix $i$ in $I$. We have

$$\widehat{\mathcal{B}}_M\big(\text{Shell}\,(w_i, D), D, 0\big) \;=\; \text{Shell}\,(w_i, D)\,.$$

On the event $\mathcal{H}'$, we have $\overline{\mathcal{T}_i}(0) = \mathcal{T}_i^{*0}(0)$ and $\mathcal{T}_i^{*1}(0) = \varnothing$, thus

$$B^{0,0}_{1,\Gamma}(0) \;=\; \mathcal{T}^{*0}(0) \cap \Gamma(0)\,, \quad B^{0,1}_{1,\Gamma}(0) \;=\; \varnothing\,,$$

$$\mathcal{B}\big(w_i, D \setminus \overline{\mathcal{T}_i}(0), 0\big) \;=\; \text{Shell}\,(w_i, D) \setminus \overline{\mathcal{T}_i}(0)\,, \tag{92.73}$$

$$\text{Shell}\,(w_i, D) \cap \widehat{\Gamma}^-(1) \;\subset\; \Big(\mathcal{B}\big(w_i, D \setminus \overline{\mathcal{T}_i}(0), 0\big) \cup B^{0,0}_{1,\Gamma}(0) \cup B^{0,1}_{1,\Gamma}(0)\Big) \cap \widehat{\Gamma}^-(1)\,. \tag{92.74}$$

Taking the union over $i$ in $I$ of the equality (92.73), we get

$$\begin{aligned}\bigcup_{i\in I} \mathcal{B}\big(w_i, D \setminus \overline{\mathcal{T}_i}(0), 0\big) \;&=\; \bigcup_{i\in I} \text{Shell}\,(w_i, D) \setminus \overline{\mathcal{T}_i}(0) \;\subset\; \bigcup_{i\in I} \text{Shell}\,(w_i, D) \\ &=\; \text{Shell}\,\big(W(I), D\big) \;=\; \widehat{\mathcal{B}}_M\Big(\text{Shell}\,\big(W(I), D\big), D, 0\Big)\,.\end{aligned}$$

Taking the trace of the sets on $\widehat{\Gamma}^-(1)$, we obtain the first inclusion in (92.72) for $t = 0$. For the second inclusion, we take the union over $i$ in $I$ of the inclusion (92.74), and we get

$$\begin{aligned}\text{Shell}\,\big(\{\, w_i : i \in I \,\}, D\big) \cap \widehat{\Gamma}^-(1) \;&\subset\; \bigcup_{i\in I} \text{Shell}\,(w_i, D) \cap \widehat{\Gamma}^-(1) \\ &\subset\; \Big(\bigcup_{i\in I} \mathcal{B}\big(w_i, D \setminus \overline{\mathcal{T}_i}(0), 0\big) \cup B^{0,0}_{1,\Gamma}(0) \cup B^{0,1}_{1,\Gamma}(0)\Big) \cap \widehat{\Gamma}^-(1)\,.\end{aligned}$$

Thus the inclusions (92.72) hold for $t = 0$.

Let now $t \geq 1$ and suppose that the two inclusions (92.72) have been proved until rank $t - 1$. Let us fix $i$ in $I$. Applying proposition 81.1, we have that

$$\begin{aligned}\mathcal{A}_i(t) \;=\; \mathcal{B}\big(w_i, D \setminus \overline{\mathcal{T}_i}(t), t\big) \;&=\; \mathcal{B}\big(w_i, D \setminus \overline{\mathcal{T}}(t), t\big) \\ &=\; \mathcal{B}\Big(\mathcal{B}\big(w_i, D \setminus \overline{\mathcal{T}}(t), t-1\big), D \setminus \overline{\mathcal{T}}(t), 1\Big)\,.\end{aligned} \tag{92.75}$$

Since $\overline{\mathcal{T}_i}(t-1)$ is contained in $\overline{\mathcal{T}}(t)$, then

$$\mathcal{B}\big(w_i, D \setminus \overline{\mathcal{T}}(t), t-1\big) \;\subset\; \mathcal{B}\big(w_i, D \setminus \overline{\mathcal{T}_i}(t-1), t-1\big)\,. \tag{92.76}$$

It follows from (92.75) and (92.76) that

$$\mathcal{B}\big(w_i, D \setminus \overline{\mathcal{T}_i}(t), t\big) \;\subset\; \mathcal{B}\Big(\mathcal{B}\big(w_i, D \setminus \overline{\mathcal{T}_i}(t-1), t-1\big), D \setminus \overline{\mathcal{T}}(t), 1\Big)\,. \tag{92.77}$$

We take the trace of the inclusion (92.77) on the box $\widehat{\Gamma}^-(t+1)$:

$$\begin{aligned}&\mathcal{B}\big(w_i, D \setminus \overline{\mathcal{T}_i}(t), t\big) \cap \widehat{\Gamma}^-(t+1) \subset \\ &\qquad \mathcal{B}\Big(\mathcal{B}\big(w_i, D \setminus \overline{\mathcal{T}_i}(t-1), t-1\big), D \setminus \overline{\mathcal{T}}(t), 1\Big) \cap \widehat{\Gamma}^-(t+1)\,.\end{aligned} \tag{92.78}$$

We claim that

$$\mathcal{B}\Big(\mathcal{B}\big(w_i, D\setminus \overline{\mathcal{T}}_i(t-1), t-1\big), D\setminus \overline{\mathcal{T}}(t), 1\Big) \cap \widehat{\Gamma}^-(t+1) \,=\, \mathcal{B}\Big(\mathcal{B}\big(w_i, D\setminus \overline{\mathcal{T}}_i(t-1), t-1\big) \cap \Gamma(t), D\setminus \overline{\mathcal{T}}(t), 1\Big) \cap \widehat{\Gamma}^-(t+1)\,. \quad (92.79)$$

Obviously, the set of the right-hand side of (92.79) is included in the set of the left-hand side, and we need only to prove the converse inclusion. Let $x$ be a site belonging to

$$\mathcal{B}\Big(\mathcal{B}\big(w_i, D\setminus \overline{\mathcal{T}}_i(t-1), t-1\big), D\setminus \overline{\mathcal{T}}(t), 1\Big) \cap \widehat{\Gamma}^-(t+1)\,.$$

We consider two cases:

• $x$ is in $\mathcal{B}\big(w_i, D\setminus \overline{\mathcal{T}}_i(t-1), t-1\big)$. We know also that $x$ is in $\widehat{\Gamma}^-(t+1)$, which is included in $\Gamma(t)$, thus $x$ is in

$$\Big(\mathcal{B}\big(w_i, D\setminus \overline{\mathcal{T}}_i(t-1), t-1\big) \cap \Gamma(t)\Big) \cap \widehat{\Gamma}^-(t+1) \,\subset\, \mathcal{B}\Big(\mathcal{B}\big(w_i, D\setminus \overline{\mathcal{T}}_i(t-1), t-1\big) \cap \Gamma(t), D\setminus \overline{\mathcal{T}}(t), 1\Big) \cap \widehat{\Gamma}^-(t+1)\,.$$

• $x$ is not in $\mathcal{B}\big(w_i, D\setminus \overline{\mathcal{T}}_i(t-1), t-1\big)$. By definition, there exists a path $y_0, \cdots, y_r$ in $D\setminus \overline{\mathcal{T}}(t)$ joining a site $y_0$ of $\mathcal{B}\big(w_i, D\setminus \overline{\mathcal{T}}_i(t-1), t-1\big)$ to $y_r = x$ such that the sites $y_1, \cdots, y_{r-1}$ are in $D \setminus \mathcal{B}\big(w_i, D\setminus \overline{\mathcal{T}}_i(t-1), t-1\big)$ and they are all open. If $r=1$, then the path is reduced to $y_0, x$, and $y_0$ is in $B_1\big(\widehat{\Gamma}^-(t+1), 1\big)$, which is included in $\Gamma(t)$. Suppose next that $r \geq 2$. Let $C$ be the open cluster of $y_1$ in $\Lambda(4n)$. The taboo set $\overline{\mathcal{T}}(t)$ prevents the explorer $i$ from exploring the clusters of the vertices $w_1, \dots, w_N$, thus $C$ does not belong to the collection $\mathcal{C}^*(h+1, n)$ (recall the definition (92.27) of the event $\mathcal{H}'$). This means that $C$ is not a boundary cluster which intersects the box $\Lambda(h+1, n)$. Yet the site $y_{r-1}$ is a neighbour of $x$, and $x$ is in the box $\Gamma(t)$, included in $\Gamma(0)$, itself included in $\Lambda(h, n)$. Thus $C$ intersects $\Lambda(h+1, n)$. Therefore $C$ is not a boundary cluster and it is a finite cluster. On the event $\mathcal{H}'$, the finite clusters have cardinalities strictly less than $M$, thus $r-1 < M$ and $|y_0 - x|_1 \leq M$. Since $x$ is in $\widehat{\Gamma}^-(t+1)$, then the site $y_0$ is in $\Gamma(t)$. Thus, in both cases $r=1$ and $r \geq 2$, the site $x$ is in

$$\mathcal{B}\big(y_0, D\setminus \overline{\mathcal{T}}(t), 1\big) \cap \widehat{\Gamma}^-(t+1) \,\subset\, \mathcal{B}\Big(\mathcal{B}\big(w_i, D\setminus \overline{\mathcal{T}}_i(t-1), t-1\big) \cap \Gamma(t), D\setminus \overline{\mathcal{T}}(t), 1\Big) \cap \widehat{\Gamma}^-(t+1)\,.$$

This completes the proof of the converse inclusion and of the equality (92.79). With the help of (92.79), we rewrite (92.78) as

$$\mathcal{B}\big(w_i, D\setminus \overline{\mathcal{T}}_i(t), t\big) \cap \widehat{\Gamma}^-(t+1) \,\subset\, \mathcal{B}\Big(\mathcal{B}\big(w_i, D\setminus \overline{\mathcal{T}}_i(t-1), t-1\big) \cap \Gamma(t), D\setminus \overline{\mathcal{T}}(t), 1\Big) \cap \widehat{\Gamma}^-(t+1)\,. \quad (92.80)$$

Taking the union over $i$ in $I$ of the inclusion (92.80), we get

$$\begin{aligned}
&\bigcup_{i\in I} \mathcal{B}\big(w_i, D\setminus \overline{\mathcal{T}}_i(t), t\big)\cap \widehat{\Gamma}^-(t+1)\\
&\qquad\subset \bigcup_{i\in I}\mathcal{B}\Big(\mathcal{B}\big(w_i, D\setminus \overline{\mathcal{T}}_i(t-1), t-1\big)\cap\Gamma(t), D\setminus\overline{\mathcal{T}}(t),1\Big)\cap\widehat{\Gamma}^-(t+1)\\
&\qquad= \mathcal{B}\Big(\bigcup_{i\in I}\mathcal{B}\big(w_i, D\setminus \overline{\mathcal{T}}_i(t-1), t-1\big)\cap\Gamma(t), D\setminus\overline{\mathcal{T}}(t),1\Big)\cap\widehat{\Gamma}^-(t+1)\\
&\qquad\subset \mathcal{B}\Big(\bigcup_{i\in I}\mathcal{B}\big(w_i, D\setminus \overline{\mathcal{T}}_i(t-1), t-1\big)\cap\widehat{\Gamma}^-(t), D\setminus\overline{\mathcal{T}}(t),1\Big)\cap\widehat{\Gamma}^-(t+1)\,.
\end{aligned}\tag{92.81}$$

To get the last inclusion, we have used the fact that $\Gamma(t)$ is included in $\widehat{\Gamma}^-(t)$. We are now in position to use the induction hypothesis, that is the first inclusion in (92.72) at rank $t-1$. Substituting this inclusion in the last term of (92.81), and replacing the taboo set $\overline{\mathcal{T}}(t)$ by the smaller taboo set $\overline{\mathcal{T}}(t-1)$, we obtain

$$\begin{aligned}
&\bigcup_{i\in I} \mathcal{B}\big(w_i, D\setminus \overline{\mathcal{T}}_i(t), t\big)\cap \widehat{\Gamma}^-(t+1)\subset\\
&\mathcal{B}\bigg(\widehat{\mathcal{B}}_M\Big(\text{Shell}\,\big(W(I),D\big),D,t-1\Big)\cap\widehat{\Gamma}^-(t), D\setminus\overline{\mathcal{T}}(t-1),1\bigg)\cap\widehat{\Gamma}^-(t+1)\,.
\end{aligned}$$

To conclude the proof, we will show that

$$\begin{aligned}
\mathcal{B}\bigg(\widehat{\mathcal{B}}_M\Big(\text{Shell}\,\big(W(I),D\big),D,t-1\Big)\cap\widehat{\Gamma}^-(t), D\setminus\overline{\mathcal{T}}(t-1),1\bigg)\cap\widehat{\Gamma}^-(t+1)&\\
\subset\ \widehat{\mathcal{B}}_M\Big(\text{Shell}\,\big(W(I),D\big),D,t\Big)\cap\widehat{\Gamma}^-(t+1)&
\end{aligned}\tag{92.82}$$

In order to prove this inclusion, we introduce the notation

$$E_M\ =\ \widehat{\mathcal{B}}_M\Big(\text{Shell}\,\big(W(I),D\big),D,t-1\Big)\cap\widehat{\Gamma}^-(t)\,.$$

Let $x$ be a site belonging to $\mathcal{B}\big(E_M, D\setminus\overline{\mathcal{T}}(t-1),1\big)\cap\widehat{\Gamma}^-(t+1)$. We consider two cases:

• $x$ is in $E_M$. Then $x$ is in

$$E_M\cap\widehat{\Gamma}^-(t+1)\ \subset\ \widehat{\mathcal{B}}_M\Big(\text{Shell}\,\big(W(I),D\big),D,t\Big)\cap\widehat{\Gamma}^-(t+1)\,.$$

• $x$ is not in $E_M$. By definition, there exists a path $y_0,\cdots,y_r$ in $D\setminus\overline{\mathcal{T}}(t-1)$ joining a site $y_0$ of $E_M$ to $y_r=x$ such that the sites $y_1,\cdots,y_{r-1}$ are in $D\setminus E_M$ and they are all open. If $r=1$, then the path is reduced to $y_0,x$, and $x$ is in

$$\widehat{\mathcal{B}}_M\big(E_M,D,1\big)\cap\widehat{\Gamma}^-(t+1)\ \subset\ \widehat{\mathcal{B}}_M\Big(\text{Shell}\,\big(W(I),D\big),D,t\Big)\cap\widehat{\Gamma}^-(t+1)\,.$$

Suppose next that $r \geq 2$. Let $C$ be the open cluster of $y_1$ in $\Lambda(4n)$. The site $x$ is in $\widehat{\Gamma}^-(t+1)$, thus the site $y_{r-1}$ is in $B_1\big(\widehat{\Gamma}^-(t+1),1\big)$ and $C$ intersects $\Lambda(h+1,n)$. If $C$ was a boundary cluster, then it would be included in Shell $\big(W(I),D\big)$, so that $y_{r-1}$ would belong to $E_M$, which is absurd. Therefore $C$ is a finite cluster. On the event $\mathcal{H}'$, the finite clusters have cardinalities strictly less than $M$, thus $r-1<M$ and

$$x \,\in\, \widehat{\mathcal{B}}_M\big(y_0, D\setminus\overline{\mathcal{T}}(t-1),1\big)\,\subset\,\widehat{\mathcal{B}}_M\big(E_M, D\setminus\overline{\mathcal{T}}(t-1),1\big)\,.$$

We conclude that $x$ is in

$$\begin{aligned}\widehat{\mathcal{B}}_M\big(E_M, D\setminus\overline{\mathcal{T}}(t-1),1\big)\cap\widehat{\Gamma}^-(t+1)\,&\subset\,\widehat{\mathcal{B}}_M\big(E_M,D,1\big)\cap\widehat{\Gamma}^-(t+1)\\ &=\,\widehat{\mathcal{B}}_M\Big(\text{Shell}\,\big(W(I),D\big),D,t\Big)\cap\widehat{\Gamma}^-(t+1)\,.\end{aligned}$$

This completes the proof of the inclusion (92.82), as well as the induction step for proving the first inclusion in (92.72).

We carry out now the induction step to prove the second inclusion in (92.72). Using the semigroup property (89.74), we write $\widehat{\mathcal{B}}_M\big(\text{Shell}\,\big(W(I),D\big),D,t\big)$ as

$$\widehat{\mathcal{B}}_M\big(\text{Shell}\,\big(W(I),D\big),D,t\big)\,=\,\widehat{\mathcal{B}}_M\Big(\widehat{\mathcal{B}}_M\big(\text{Shell}\,\big(W(I),D\big),D,t-1\big),D,1\Big)\,. \tag{92.83}$$

We take the trace of the identity (92.83) on $\widehat{\Gamma}^-(t+1)$ and we get

$$\begin{aligned}&\widehat{\mathcal{B}}_M\big(\text{Shell}\,\big(W(I),D\big),D,t\big)\cap\widehat{\Gamma}^-(t+1)\,=\\ &\qquad\widehat{\mathcal{B}}_M\Big(\widehat{\mathcal{B}}_M\big(\text{Shell}\,\big(W(I),D\big),D,t-1\big),D,1\Big)\cap\widehat{\Gamma}^-(t+1)\,.\end{aligned} \tag{92.84}$$

In order to rework the right-hand side of (92.84), we need the following result.

**Lemma 92.5.** *For any subset $A$ of $D$, and any box $\Lambda$ included in $D$, we have*

$$\widehat{\mathcal{B}}_M\big(A,D,1\big)\cap\Lambda\,=\,\widehat{\mathcal{B}}_M\Big(A\cap B_1(\Lambda,M+1),D\cap B_1(\Lambda,2M),1\Big)\cap\Lambda\,. \tag{92.85}$$

*Proof.* The set of the right-hand side of (92.85) is obviously included in the set of the left-hand side. We have to prove the converse inclusion. Let $x$ be a site in $\widehat{\mathcal{B}}_M\big(A,D,1\big)\cap\Lambda$. If $x$ is in $A\cap\Lambda$, then $x$ is also in

$$\widehat{\mathcal{B}}_M\Big(A\cap B_1(\Lambda,M+1),D\cap B_1(\Lambda,2M),1\Big)\cap\Lambda\,. \tag{92.86}$$

Otherwise, it follows from the definition of $\widehat{\mathcal{B}}_M(\cdot,\cdot,\cdot)$ that $x$ belongs to

$$\overline{\text{Clusters}}\,(A,M,D)\cup\overline{\text{Clusters}}\,\big(\mathcal{N}(A),M,D\big)\,=\,\bigcup_{y\in A\cup\mathcal{N}(A)}\overline{\text{Cluster}}\,(y,M,D)\,.$$

Thus there exists a site $y$ in $A\cup\mathcal{N}(A)$ such that $x$ belongs to $\overline{\text{Cluster}}\,(y,M,D)$. It follows from the definition of $\overline{\text{Cluster}}\,(y,M,D)$ that

$$\overline{\text{Cluster}}\,(y,M,D)\,\subset\,B_1(y,M)\,.$$

Since $x$ is in $B_1(y,M)\cap\Lambda$, then $y$ is in $B_1(\Lambda,M)$ and

$$\overline{\text{Cluster}}\,(y,M,D)\,\subset\,B_1(\Lambda,2M)\,.$$

If $y$ is in $A$, then $x$ is in the set of the right-hand side of (92.85). Suppose next that $y$ is in $\mathcal{N}(A)\setminus A$. In this case, the site $y$ has a neighbour $z$ which belongs to $A$. This site $z$ is in $A\cap B_1(\Lambda,M+1)$ and moreover,

$$\begin{aligned}&\overline{\text{Cluster}}\,(y,M,D)\,\subset\,\overline{\text{Cluster}}\,\big(y,M,D\cap B_1(\Lambda,2M)\big)\\ &\subset\,\overline{\text{Clusters}}\,\big(\mathcal{N}(z),M,D\cap B_1(\Lambda,2M)\big)\,\subset\,\widehat{\mathcal{B}}_M\big(z,D\cap B_1(\Lambda,2M),1\big)\,.\end{aligned}\tag{92.87}$$

We conclude from (92.87) that $x$ belongs to the set (92.86). This proves the converse inclusion and it completes the proof of the lemma. □

Notice that we have tuned the boxes $\widehat{\Gamma}^-(t+1)$ and $\Gamma(t)$ so that

$$B_1\big(\widehat{\Gamma}^-(t+1),M+1\big)\,\subset\,\Gamma(t)\,.\tag{92.88}$$

Using lemma 92.5, and the inclusion (92.88), we have

$$\begin{aligned}&\widehat{\mathcal{B}}_M\Big(\widehat{\mathcal{B}}_M\big(\text{Shell}\,\big(W(I),D\big),D,t-1\big),D,1\Big)\cap\widehat{\Gamma}^-(t+1)\\ &\subset\,\widehat{\mathcal{B}}_M\Big(\widehat{\mathcal{B}}_M\big(\text{Shell}\,\big(W(I),D\big),D,t-1\big)\cap B_1\big(\widehat{\Gamma}^-(t+1),M+1\big),D,1\Big)\cap\widehat{\Gamma}^-(t+1)\\ &\qquad\subset\,\widehat{\mathcal{B}}_M\Big(\widehat{\mathcal{B}}_M\big(\text{Shell}\,\big(W(I),D\big),D,t-1\big)\cap\Gamma(t),D,1\Big)\cap\widehat{\Gamma}^-(t+1)\\ &\quad\subset\,\widehat{\mathcal{B}}_M\Big(\widehat{\mathcal{B}}_M\big(\text{Shell}\,\big(W(I),D\big),D,t-1\big)\cap\widehat{\Gamma}^-(t),D,1\Big)\cap\widehat{\Gamma}^-(t+1)\,.\end{aligned}\tag{92.89}$$

We use next the induction hypothesis. The identity (92.84), the inclusion (92.89) and the inclusion (92.72) at rank $t-1$ together yield that

$$\begin{aligned}&\widehat{\mathcal{B}}_M\big(\text{Shell}\,\big(W(I),D\big),D,t\big)\cap\widehat{\Gamma}^-(t+1)\,\subset\\ &\widehat{\mathcal{B}}_M\Bigg(\Big(\bigcup_{i\in I}\mathcal{B}\big(w_i,D\setminus\overline{\mathcal{T}}_i(t-1),t-1\big)\cap\widehat{\Gamma}^-(t)\Big)\cup\bigcup_{1\le s\le t-1}B^{s,-}_{1,\Gamma}(t-s)\,\cup\\ &\qquad\bigcup_{0\le s\le t-1}\big(B^{s,0}_{1,\Gamma}(t-1-s)\,\cup\,B^{s,1}_{1,\Gamma}(t-s)\big),D,1\Bigg)\cap\widehat{\Gamma}^-(t+1)\\ &\quad\subset\,\Bigg(\widehat{\mathcal{B}}_M\Big(\bigcup_{i\in I}\mathcal{B}\big(w_i,D\setminus\overline{\mathcal{T}}_i(t-1),t-1\big),D,1\Big)\cap\widehat{\Gamma}^-(t+1)\Bigg)\cup\\ &\qquad\widehat{\mathcal{B}}_M\Big(\bigcup_{1\le s\le t-1}B^{s,-}_{1,\Gamma}(t-s)\,\cup\bigcup_{0\le s\le t-1}\big(B^{s,0}_{1,\Gamma}(t-1-s)\,\cup\,B^{s,1}_{1,\Gamma}(t-s)\big),D,1\Big)\,.\end{aligned}\tag{92.90}$$

We study separately the two sets on the right-hand side of (92.90). Let us start with the first set. It follows from (81.2) in proposition 81.1 and the definition of the sets $\mathcal{T}_i^*(t)$, $i \in I$, that

$$\forall i \in I \qquad \mathcal{A}_i(t-1) \,=\, \mathcal{B}\big(w_i, D \setminus \overline{\mathcal{T}}_i(t-1), t-1\big) \,=\, \mathcal{B}\big(w_i, D \setminus \overline{\mathcal{T}}_i(t), t-1\big)\,,$$

whence

$$\widehat{\mathcal{B}}_M\Big(\bigcup_{i\in I}\mathcal{B}\big(w_i, D \setminus \overline{\mathcal{T}}_i(t-1), t-1\big), D, 1\Big) \cap \widehat{\Gamma}^-(t+1) \,=\, \widehat{\mathcal{B}}_M\Big(\bigcup_{i\in I}\mathcal{B}\big(w_i, D \setminus \overline{\mathcal{T}}_i(t), t-1\big), D, 1\Big) \cap \widehat{\Gamma}^-(t+1)\,. \quad (92.91)$$

We would like to compare the set of the right-hand side with the set obtained by non-truncated explorations, but with taboo set $\overline{\mathcal{T}}(t)$, that is the set

$$\mathcal{B}\Big(\bigcup_{i\in I}\mathcal{B}\big(w_i, D \setminus \overline{\mathcal{T}}_i(t), t-1\big), D \setminus \overline{\mathcal{T}}(t), 1\Big) \cap \widehat{\Gamma}^-(t+1) \\ = \bigcup_{i\in I}\mathcal{B}\big(w_i, D \setminus \overline{\mathcal{T}}_i(t), t\big) \cap \widehat{\Gamma}^-(t+1)\,.$$

The next lemma is a localized variant of lemma 89.19.

**Lemma 92.6.** *For any $t$ in $\{1, \dots, T\}$, we have the inclusion*

$$\widehat{\mathcal{B}}_M\Big(\bigcup_{i\in I}\mathcal{B}\big(w_i, D \setminus \overline{\mathcal{T}}_i(t), t-1\big), D, 1\Big) \cap \widehat{\Gamma}^-(t+1) \,\subset\, \Big(\bigcup_{i\in I}\mathcal{B}\big(w_i, D \setminus \overline{\mathcal{T}}_i(t), t\big) \cap \widehat{\Gamma}^-(t+1)\Big) \cup B_1\big(\overline{\mathcal{T}}(t-1) \cap \Gamma(t), 2M\big) \\ \cup \big(\mathcal{T}^{*0}(t) \cap \widehat{\Gamma}^-(t+1)\big) \cup B_1\big(\mathcal{T}^{*1}(t) \cap \Gamma(t), 2M\big)\,. \quad (92.92)$$

*Proof.* Obviously, we have the inclusion

$$\bigcup_{i\in I}\mathcal{B}\big(w_i, D \setminus \overline{\mathcal{T}}_i(t), t-1\big) \,\subset\, \bigcup_{i\in I}\mathcal{B}\big(w_i, D \setminus \overline{\mathcal{T}}_i(t), t\big)\,.$$

Let $x$ belong to

$$\bigcup_{i\in I}\widehat{\mathcal{B}}_M\Big(\mathcal{B}\big(w_i, D \setminus \overline{\mathcal{T}}_i(t), t-1\big), D, 1\Big) \cap \widehat{\Gamma}^-(t+1) \setminus \bigcup_{i\in I}\mathcal{B}\big(w_i, D \setminus \overline{\mathcal{T}}_i(t), t-1\big)\,.$$

By definition, there exist $i$ in $I$ and a path $y_0, y_1, \cdots, y_r$ in $D$ joining a site $y_0$ of $\mathcal{B}\big(w_i, D \setminus \overline{\mathcal{T}}_i(t), t-1\big)$ to $y_r = x$ such that $r \geq 1$, the sites $y_1, \cdots, y_{r-1}$ are in $D \setminus \mathcal{B}\big(w_i, D \setminus \overline{\mathcal{T}}_i(t), t-1\big)$, they are all open, and they are visited before $x$ by the truncated exploration starting from $y_1$. The truncation imposes that $y_2, \dots, y_r$ are in $B_1(y_1, M)$. If the whole path $y_0, y_1, \cdots, y_r$ does not visit $\overline{\mathcal{T}}_i(t)$, then $x$ is in

$$\mathcal{B}\Big(\mathcal{B}\big(w_i, D\setminus\overline{\mathcal{T}}_i(t), t-1\big), D\setminus\overline{\mathcal{T}}_i(t), 1\Big)\cap\widehat{\Gamma}^-(t+1) = \mathcal{B}\big(w_i, D\setminus\overline{\mathcal{T}}_i(t), t\big)\cap\widehat{\Gamma}^-(t+1).$$

Otherwise, let $s$ be the smallest index such that $y_s$ is in $\overline{\mathcal{T}}_i(t)$. As $y_0$ belongs to $\mathcal{B}\big(w_i, D\setminus\overline{\mathcal{T}}_i(t), t-1\big)$, which is disjoint from $\overline{\mathcal{T}}_i(t)$, we must have $s\geq 1$. Since both $y_s$ and $x$ are in $B_1(y_1,M)$, then $|x-y_s|_1\leq 2M$. Moreover $x$ is in $\widehat{\Gamma}^-(t+1)$, thus $y_s$ is in $\Gamma(t)$. Suppose that $y_s$ is in $\overline{\mathcal{T}}_i(t-1)$. Then $y_s$ is in $\overline{\mathcal{T}}(t-1)\cap\Gamma(t)$, and the path $y_{s+1},\cdots,y_r$ testifies that $x$ is in $B_1\big(\overline{\mathcal{T}}(t-1)\cap\Gamma(t),2M\big)$. Suppose next that $y_s$ is in $\overline{\mathcal{T}}_i(t)\setminus\overline{\mathcal{T}}_i(t-1)\subset\mathcal{T}_i^*(t)$. We consider two cases:

• The site $y_s$ is closed. Then the index $s$ must be equal to $r$, and $y_s=x$ is in $\mathcal{T}_i^{*0}(t)\cap\widehat{\Gamma}^-(t+1)$.

• The site $y_s$ is open. Then $y_s$ is in $\mathcal{T}_i^{*1}(t)\cap\Gamma(t)$. In addition, the sites $y_s,\cdots,y_r$ form a path testifying that $x$ is in $B_1\big(\mathcal{T}_i^{*1}(t)\cap\Gamma(t),2M\big)$.

In each case, the site $x$ belongs to one of the sets appearing in the right-hand side of (92.92). This concludes the proof of the inclusion stated in the lemma. □

For the second set on the right-hand side of (92.90), we use lemma 89.15 to obtain

$$\begin{aligned}
\widehat{\mathcal{B}}_M\Big(&\bigcup_{1\leq s\leq t-1} B^{s,-}_{1,\Gamma}(t-s)\,\cup\bigcup_{0\leq s\leq t-1}\big(B^{s,0}_{1,\Gamma}(t-1-s)\,\cup\,B^{s,1}_{1,\Gamma}(t-s)\big),D,1\Big)\\
&\subset B_1\Big(\bigcup_{1\leq s\leq t-1} B^{s,-}_{1,\Gamma}(t-s)\,\cup\bigcup_{0\leq s\leq t-1}\big(B^{s,0}_{1,\Gamma}(t-1-s)\,\cup\,B^{s,1}_{1,\Gamma}(t-s)\big),M+1\Big)\\
&\quad\subset\Big(\bigcup_{1\leq s\leq t-1} B^{s,-}_{1,\Gamma}(t+1-s)\,\cup\bigcup_{0\leq s\leq t-1}\big(B^{s,0}_{1,\Gamma}(t-s)\,\cup\,B^{s,1}_{1,\Gamma}(t+1-s)\big)\Big)\,.
\end{aligned}\tag{92.93}$$

Using the identity (92.91), and the two inclusions (92.92), (92.93), we deduce from (92.90), that

$$\begin{aligned}
\widehat{\mathcal{B}}_M\Big(\text{Shell}\,\big(W(I),D\big),D,t\Big)\cap\widehat{\Gamma}^-(t+1)\,\subset\,\Big(\bigcup_{i\in I}\mathcal{B}\big(w_i,D\setminus\overline{\mathcal{T}}_i(t),t\big)\cap\widehat{\Gamma}^-(t+1)\Big)\\
\cup\,B_1\big(\overline{\mathcal{T}}(t-1)\cap\Gamma(t),2M\big)\cup\big(\mathcal{T}^{*0}(t)\cap\widehat{\Gamma}^-(t+1)\big)\cup B_1\big(\mathcal{T}^{*1}(t)\cap\Gamma(t),2M\big)\\
\cup\bigcup_{1\leq s\leq t-1} B^{s,-}_{1,\Gamma}(t+1-s)\,\cup\bigcup_{0\leq s\leq t-1}\big(B^{s,0}_{1,\Gamma}(t-s)\,\cup\,B^{s,1}_{1,\Gamma}(t+1-s)\big)\,.
\end{aligned}\tag{92.94}$$

Furthermore, we check that

$$B_1\big(\overline{\mathcal{T}}(t-1)\cap\Gamma(t),2M\big)\,=\,B^{t,-}_{1,\Gamma}(1)\,,\tag{92.95}$$

$$\mathcal{T}^{*0}(t)\cap\widehat{\Gamma}^-(t+1)\,\subset\mathcal{T}^{*0}(t)\cap\Gamma(t)\,=\,B^{t,0}_{1,\Gamma}(0)\,,\tag{92.96}$$

$$B_1\big(\mathcal{T}^{*1}(t)\cap\Gamma(t),2M\big)\,=\,B^{t,1}_{1,\Gamma}(1)\,.\tag{92.97}$$

Substituting (92.95), (92.96), (92.97) into (92.94), we obtain the second inclusion (92.72) at rank $t$. This completes the induction step for proving the second inclusion in (92.72), as well as the proof of the proposition. □

## 92.11 Approximation of $\mathcal{B}\big(w_i, D \setminus \widehat{\mathcal{T}}_i^-(t), t\big)$

In this subsection, we try to approximate the taboo exploration

$$\bigcup_{i\in I} \mathcal{B}\big(w_i, D \setminus \widehat{\mathcal{T}}_i^-(t), t\big)\,.$$

The problematic for approximating this set is similar to the one we encountered in the previous subsection when dealing with the set $\mathcal{B}\big(w_i, D \setminus \overline{\mathcal{T}}_i(t), t\big)$. In propositions 89.20, 89.21, we already carried out an approximation of this set for the intertwined explorations starting from all the boundary clusters of $\Lambda(4n)$. Naturally, this approximation is not valid any more when we start with a restricted subset of the boundary clusters. Therefore, we have to develop an adequate localized variant of this approximation.

Instead of the sets $B^{t,0}_{1,\Gamma}(s)$, $B^{t,1}_{1,\Gamma}(s)$, $B^{t,-}_{1,\Gamma}(s)$ defined in (92.69), (92.70), (92.71), we shall use the set $\overline{B}^t_{1,\Gamma}(s)$ defined as

$$\forall t \in \{\,0,\dots,T\,\}\quad \forall s\geq 0 \qquad \overline{B}^t_{1,\Gamma}(s) \,=\, B_1\big(\overline{\mathcal{T}}(t)\cap\Gamma(t), s(M+1)+M-1\big)\,, \tag{92.98}$$

and the set $\widehat{B}^{t,1}_{1,\Gamma}(s)$ given by

$$\begin{aligned}\forall t \in \{\,1,\dots,T\,\}\quad \forall s\geq 0 &\\ \widehat{B}^{t,1}_{1,\Gamma}(s) \,=\, B_1\Big(&\mathcal{T}^{-1}(t)\cap \overset{\circ}{\widehat{\Gamma(t-1)}}, s(M+1)+M-1\Big)\\ &\cup B_1\big(\mathcal{T}^{+1}(t-1)\cap\Gamma(t-1), s(M+1)+M-1\big)\,. \end{aligned}\tag{92.99}$$

Here is the proper reformulation for the set $\mathcal{B}\big(w_i, D \setminus \widehat{\mathcal{T}}_i^-(t), t\big)$ of the results stated in propositions 89.20, 89.21.

**Proposition 92.7.** *We consider a configuration realizing the event $\mathcal{H}'$ defined in* (92.27) *and a subset $I$ of $\{\,1,\dots,N\,\}$. For any $t$ in $\{\,1,\dots,T\,\}$, we have*

$$\begin{aligned}\bigcup_{i\in I} \mathcal{B}\big(w_i, D \setminus \widehat{\mathcal{T}}_i^-(t), t\big)\cap\widehat{\Gamma}^-(t) \,\subset\,&\\ \widehat{\mathcal{B}}_M\Big(\mathit{Shell}\big(W(I),D\big),D,t\Big)\cap\widehat{\Gamma}^-(t) \,\subset\,& \end{aligned}\tag{92.100}$$

$$\Big(\bigcup_{i\in I} \mathcal{B}\big(w_i, D \setminus \widehat{\mathcal{T}}_i^-(t), t\big)\cap\widehat{\Gamma}^-(t)\Big)\cup\bigcup_{1\leq s\leq t}\Big(\overline{B}^{s-1}_{1,\Gamma}(t+1-s)\cup\widehat{B}^{s,1}_{1,\Gamma}(t+1-s)\Big)\,.$$

*Proof.* We shall prove the proposition by induction on $t$ and we start with $t=1$. Let $i$ belong to $I$. As the set $\widehat{\mathcal{T}}_i^-(1)$ contains $\mathcal{T}_i(0)$, we have

$$\widehat{\Gamma}^-(1)\cap\mathcal{B}\big(\mathrm{Shell}\,(w_i,D), D\setminus\widehat{\mathcal{T}}_i^-(1),1\big)\setminus\mathrm{Shell}\,\big(w_i,D\big) \,\subset\, \widehat{\Gamma}^-(1)\setminus\bigcup_{1\leq j\leq N} C(w_j,D)\,.$$

Since the event $\mathcal{H}'$ occurs, then the only boundary clusters which intersect $\Lambda(h+1,n)$, and a fortiori $\widehat{\Gamma}^-(1)$, are the $N$ clusters $C(w_j,D)$, $1\leq j\leq N$. Thus all the clusters that intersect the above set are finite clusters of $\Lambda(4n)$. In addition, the event $\mathcal{H}'$ guarantees that all these finite clusters have cardinality strictly less than $M$. This implies that

$$\begin{aligned}\mathcal{B}\big(w_i, D\setminus\widehat{\mathcal{T}}_i^-(1),1\big)\cap\widehat{\Gamma}^-(1) &= \mathcal{B}\big(\text{Shell}\,(w_i,D),D\setminus\widehat{\mathcal{T}}_i^-(1),1\big)\cap\widehat{\Gamma}^-(1)\\ &= \widehat{\mathcal{B}}_M\big(\text{Shell}\,(w_i,D),D\setminus\widehat{\mathcal{T}}_i^-(1),1\big)\cap\widehat{\Gamma}^-(1)\\ &\subset \widehat{\mathcal{B}}_M\big(\text{Shell}\,(w_i,D),D,1\big)\cap\widehat{\Gamma}^-(1)\,. \end{aligned}\tag{92.101}$$

Taking the union of (92.101) over $i$ in $I$, we get

$$\bigcup_{i\in I}\mathcal{B}\big(w_i, D\setminus\widehat{\mathcal{T}}_i^-(1),1\big)\cap\widehat{\Gamma}^-(1) = \bigcup_{i\in I}\widehat{\mathcal{B}}_M\big(\text{Shell}\,(w_i,D),D\setminus\widehat{\mathcal{T}}_i^-(1),1\big)\cap\widehat{\Gamma}^-(1)$$
$$\subset \bigcup_{i\in I}\widehat{\mathcal{B}}_M\big(\text{Shell}\,(w_i,D),D,1\big)\cap\widehat{\Gamma}^-(1) = \widehat{\mathcal{B}}_M\Big(\text{Shell}\,\big(W(I),D\big),D,1\Big)\cap\widehat{\Gamma}^-(1)\,.$$

Let $x$ be a site belonging to

$$\widehat{\mathcal{B}}_M\Big(\text{Shell}\,\big(W(I),D\big),D,1\Big)\cap\widehat{\Gamma}^-(1)\ \setminus\ \bigcup_{i\in I}\widehat{\mathcal{B}}_M\big(\text{Shell}\,(w_i,D),D\setminus\widehat{\mathcal{T}}_i^-(1),1\big)\,.\tag{92.102}$$

By definition of $\widehat{\mathcal{B}}_M\big(\text{Shell}\,\big(W(I),D\big),D,1\big)$, there exists a path $y_0,y_1,\cdots,y_r$ in $D$ joining a site $y_0$ of $\text{Shell}\,\big(W(I),D\big)$ to $y_r=x$ such that the sites $y_1,\cdots,y_{r-1}$ are in $D\setminus\text{Shell}\,\big(W(I),D\big)$, they are all open, and they are visited before $x$ by the truncated exploration starting from $y_1$. The truncation imposes that $y_2,\dots,y_r$ are in $B_1(y_1,M)$ (if $x$ is in $\text{Shell}\,\big(W(I),D\big)$, then the sequence is reduced to $x=y_0$ and we have $r=0$). Since $x$ is in $\widehat{\Gamma}^-(1)$, then the whole path $y_0,\dots,y_r$ is included in $\Gamma(0)$. We have also

$$\text{Shell}\,\big(W(I),D\big)\ =\ \bigcup_{i\in I}\text{Shell}\,(w_i,D)\,,$$

so that there exists an index $i$ in $I$ such that $y_0$ is in $\text{Shell}\,(w_i,D)$. Since $x$ is not in $\widehat{\mathcal{B}}_M\big(\text{Shell}\,(w_i,D),D\setminus\widehat{\mathcal{T}}_i^-(1),1\big)$, then the path $y_0,\cdots,y_r$ must go through $\widehat{\mathcal{T}}_i^-(1)$. Thus there exists $s$ such that $y_s$ is in $\widehat{\mathcal{T}}_i^-(1)$. If $s=0$, then $y_s=y_0$ is in

$$\text{Shell}\,(w_i,D)\cap\widehat{\mathcal{T}}_i^-(1)\cap\Gamma(0)\ =\ \text{Shell}\,(w_i,D)\cap\mathcal{T}_i^*(0)\cap\Gamma(0)\,,$$

so that $x$ is in

$$B_1(\mathcal{T}_i^{*0}(0)\cap\Gamma(0),M+1)\ \subset\ \overline{B}^0_{1,\Gamma}(1)\,.$$

Suppose now that $s>0$. Since both $y_s$ and $x$ are in $B_1(y_1,M)$, then $|x-y_s|_1\leq 2M$. Therefore, the site $y_s$ is in

$$B_1(\widehat{\Gamma}^-(1),2M)\ \subset\ \overset{\circ}{\widehat{\Gamma(0)}}\,.$$

If $0<s<r$, then $y_s$ is in $\mathcal{T}_i^{-1}(1)\cap\overset{\circ}{\widehat{\Gamma(0)}}$, so that $x$ is in $B_1(\mathcal{T}_i^{-1}(1)\cap\overset{\circ}{\widehat{\Gamma(0)}},M)$, which is included in $\widehat{B}^{1,1}_{1,\Gamma}(1)$. If $s=r$, then $x$ is in

$$\widehat{\mathcal{T}}_i^-(1)\cap\widehat{\Gamma}^-(1)\ =\ \Big(\mathcal{T}_i^{*0}(0)\cup\mathcal{T}_i^{-1}(1)\Big)\cap\overset{\circ}{\widehat{\Gamma(0)}}\ \subset\ \overline{B}^0_{1,\Gamma}(1)\cup\widehat{B}^{1,1}_{1,\Gamma}(1)\,.$$

We conclude that the set (92.102) is included in $\overline{B}^{0}_{1,\Gamma}(1)\cup\widehat{B}^{1,1}_{1,\Gamma}(1)$. This completes the case $t=1$.

We perform now the induction step. Let $t\geq 2$ and suppose that the two inclusions (92.100) has been proved until rank $t-1$. Let us fix $i$ in $I$. The usual semigroup property (12.9) yields that

$$\begin{aligned}\mathcal{B}\big(w_i,D\setminus\widehat{\mathcal{T}}_i^-(t),t\big)\,\subset\,&\mathcal{B}\Big(\mathcal{B}\big(w_i,D\setminus\widehat{\mathcal{T}}_i^-(t),t-1\big),D\setminus\widehat{\mathcal{T}}_i^-(t),1\Big)\\&\subset\,\mathcal{B}\Big(\mathcal{B}\big(w_i,D\setminus\overline{\mathcal{T}}_i(t-1),t-1\big),D\setminus\widehat{\mathcal{T}}_i^-(t),1\Big)\,.\end{aligned}$$

We claim that

$$\begin{aligned}&\mathcal{B}\big(w_i,D\setminus\widehat{\mathcal{T}}_i^-(t),t\big)\cap\widehat{\Gamma}^-(t)\,\subset\\&\quad\widehat{\mathcal{B}}_M\Big(\mathcal{B}\big(w_i,D\setminus\overline{\mathcal{T}}_i(t-1),t-1\big)\cap\Gamma(t-1),D\setminus\widehat{\mathcal{T}}_i^-(t),1\Big)\cap\widehat{\Gamma}^-(t)\,.\end{aligned}\tag{92.103}$$

In order to prove this inclusion, with the help of the identity (81.2) of proposition 81.1, we rewrite first the left-hand side of (92.103) as

$$\begin{aligned}&\mathcal{B}\big(w_i,D\setminus\widehat{\mathcal{T}}_i^-(t),t\big)\cap\widehat{\Gamma}^-(t)\\&\qquad=\mathcal{B}\Big(\mathcal{B}\big(w_i,D\setminus\overline{\mathcal{T}}_i(t-1),t-1\big),D\setminus\widehat{\mathcal{T}}_i^-(t),1\Big)\cap\widehat{\Gamma}^-(t)\\&\qquad\qquad=\mathcal{B}\big(\mathcal{A}_i(t-1),D\setminus\widehat{\mathcal{T}}_i^-(t),1\big)\cap\widehat{\Gamma}^-(t)\,.\end{aligned}$$

We consider then two sites $x,y$ such that

$$x\in\mathcal{A}_i(t-1)\,,\quad y\notin\mathcal{A}_i(t-1)\,,\quad y\in D\setminus\widehat{\mathcal{T}}_i^-(t)\,,\quad |x-y|_1=1\,,\tag{92.104}$$

and we analyze the contribution of $y$ to the two sets appearing in (92.103). During the iteration $t$, when the explorer $i$ becomes active, it will visit the site $y$ unless that site has been visited before. We analyze the various possible cases for $y$:

• The site $y$ is open and it is in a boundary cluster. If the exploration starting from $y$ contributes to $\mathcal{B}\big(w_i,D\setminus\widehat{\mathcal{T}}_i^-(t),t\big)\cap\widehat{\Gamma}^-(t)$, then the cluster of $y$ must intersect the box $\Gamma(t-1)$. However, the boundary clusters that intersect $\Gamma(t-1)$ belong to the family $C(w_j,D)$, $1\leq j\leq N$. Therefore the site $y$ was visited during the first iteration by an explorer other than $i$, and $y$ is in $\mathcal{T}_i^{-1}(t)$. As $\mathcal{T}_i^{-1}(t)$ is included in $\widehat{\mathcal{T}}_i^-(t)$ for $t\geq 1$, this case cannot occur.

• The site $y$ is not in $\mathcal{T}_i^*(t)$. In this case, the site $y$ is not in a boundary cluster that intersects the box $\Gamma(0)$. If $y$ is in a boundary cluster, then the exploration of this cluster and its boundary does not reach the box $\widehat{\Gamma}^-(t)$. Otherwise, the cluster of $y$ is finite.

• The site $y$ is in $\mathcal{T}_i^*(t)$. In this case, the site $y$ is in

$$\mathcal{T}_i^*(t)\setminus\widehat{\mathcal{T}}_i^-(t)\,=\,\mathcal{T}_i^{*0}(t)\cup\mathcal{T}_i^{+1}(t)\,.$$

We consider two subcases:

$\star$ If $y$ is in $\mathcal{T}_i^{*0}(t)$, then it is closed.
$\star$ If $y$ is in $\mathcal{T}_i^{+1}(t)$, then it was visited by an explorer during the iteration $t \geq 1$, thus it is not in a boundary cluster that intersects the box $\Gamma(0)$. If $y$ is in a boundary cluster, then the exploration of this cluster and its boundary does not reach the box $\widehat{\Gamma}^-(t)$. Otherwise, the cluster of $y$ is finite.

We see that, in all the above cases, either the cluster of $y$ is a boundary cluster that does not intersect the box $\Gamma(0)$, or the cluster of $y$ is finite (this includes the case where $y$ is closed). The only case when the exploration of the cluster of $y$ visits a site of $\widehat{\Gamma}^-(t)$ is when it is finite. Let us suppose that it is the case, i.e., the cluster of $y$ is finite and

$$\overline{C}(y,D) \cap \widehat{\Gamma}^-(t) \,\neq\, \varnothing\,. \tag{92.105}$$

As the event $\mathcal{H}'$ occurs, we have also

$$\overline{\text{Cluster}}\,(y,M,D) \,=\, \overline{C}(y,D)\,, \tag{92.106}$$

and it follows from (92.105) and (92.106) that $\overline{C}(y,D)$ is included in $\Gamma(t-1)$. In particular, the site $y$ is in $\Gamma(t-1)$, and its contribution to the two sets appearing in (92.103) is the same. This is true for any pair of sites $x,y$ satisfying (92.104), and this completes the proof of the inclusion (92.103). This inclusion holds for any $t$ in $\{\,2,\dots,T\,\}$. We complete now the proof of the first inclusion in (92.100). Let us fix $t\geq 2$. We take the union over $i$ in $I$ of the inclusion (92.103) and we get

$$\begin{aligned}
&\bigcup_{i\in I}\mathcal{B}\big(w_i, D\setminus\widehat{\mathcal{T}}_i^-(t),t\big)\cap\widehat{\Gamma}^-(t)\\
&\quad\subset\bigcup_{i\in I}\widehat{\mathcal{B}}_M\Big(\mathcal{B}\big(w_i,D\setminus\overline{\mathcal{T}}_i(t-1),t-1\big)\cap\Gamma(t-1),D\setminus\widehat{\mathcal{T}}_i^-(t),1\Big)\cap\widehat{\Gamma}^-(t)\\
&\qquad\subset\bigcup_{i\in I}\widehat{\mathcal{B}}_M\Big(\mathcal{B}\big(w_i,D\setminus\overline{\mathcal{T}}_i(t-1),t-1\big)\cap\Gamma(t-1),D,1\Big)\cap\widehat{\Gamma}^-(t)\\
&\quad\subset\widehat{\mathcal{B}}_M\Bigg(\bigcup_{i\in I}\mathcal{B}\big(w_i,D\setminus\overline{\mathcal{T}}_i(t-1),t-1\big)\cap\Gamma(t-1),D,1\Bigg)\cap\widehat{\Gamma}^-(t)\,.
\end{aligned}\tag{92.107}$$

Notice that, to get the second inclusion in (92.107), we have used the monotonicity of the truncated exploration $\widehat{\mathcal{B}}_M(\cdot,\cdot,\cdot)$ with respect to its second argument. In order to control the last union in (92.107), we use the first inclusion in (92.100) at iteration $t-1$ (more precisely the trace on $\Gamma(t-1)\subset\widehat{\Gamma}^-(t-1)$ of this inclusion) and we obtain that

$$\begin{aligned}
&\bigcup_{i\in I}\mathcal{B}\big(w_i, D\setminus\widehat{\mathcal{T}}_i^-(t),t\big)\cap\widehat{\Gamma}^-(t)\\
&\qquad\subset\widehat{\mathcal{B}}_M\Bigg(\widehat{\mathcal{B}}_M\Big(\text{Shell}\,\big(W(I),D\big),D,t-1\Big)\cap\Gamma(t-1),D,1\Bigg)\cap\widehat{\Gamma}^-(t)\\
&\qquad\qquad\subset\widehat{\mathcal{B}}_M\Big(\text{Shell}\,\big(W(I),D\big),D,t\Big)\cap\widehat{\Gamma}^-(t)\,.
\end{aligned}$$

This completes the proof of the first inclusion in (92.100). We prove next the second inclusion of (92.100) at rank $t$. Using the semigroup property (89.74), we write $\widehat{\mathcal{B}}_M\big(\text{Shell}\,\big(W(I),D\big),D,t\big)$ as

$$\widehat{\mathcal{B}}_M\big(\text{Shell}\,\big(W(I),D\big),D,t\big) \;=\; \widehat{\mathcal{B}}_M\Big(\widehat{\mathcal{B}}_M\big(\text{Shell}\,\big(W(I),D\big),D,t-1\big),D,1\Big)\,. \tag{92.108}$$

We use next the induction hypothesis. The identity (92.108) and the inclusion (92.100) at rank $t-1$ together yield that

$$\begin{aligned}
&\widehat{\mathcal{B}}_M\big(\text{Shell}\,\big(W(I),D\big),D,t\big)\cap\widehat{\Gamma}^-(t)\;\subset\\
&\quad\widehat{\mathcal{B}}_M\Bigg(\Big(\bigcup_{i\in I}\mathcal{B}\big(w_i,D\setminus\widehat{\mathcal{T}}_i^-(t-1),t-1\big)\cap\widehat{\Gamma}^-(t-1)\Big)\\
&\qquad\qquad\cup\bigcup_{1\leq s\leq t-1}\Big(\overline{B}^{s-1}_{1,\Gamma}(t-s)\cup\widehat{B}^{s,1}_{1,\Gamma}(t-s)\Big),D,1\Bigg)\cap\widehat{\Gamma}^-(t)\\
&\subset\;\Bigg(\bigcup_{i\in I}\widehat{\mathcal{B}}_M\Big(\mathcal{B}\big(w_i,D\setminus\widehat{\mathcal{T}}_i^-(t-1),t-1\big),D,1\Big)\cap\widehat{\Gamma}^-(t)\Bigg)\cup\\
&\qquad\qquad\widehat{\mathcal{B}}_M\Bigg(\bigcup_{1\leq s\leq t-1}\Big(\overline{B}^{s-1}_{1,\Gamma}(t-s)\cup\widehat{B}^{s,1}_{1,\Gamma}(t-s)\Big),D,1\Bigg)\,.
\end{aligned} \tag{92.109}$$

We study separately the two sets on the right-hand side of (92.109). Let us start with the first set. The next lemma is a localized variant of lemma 89.22. We need this localized variant, because, unfortunately, we cannot apply directly lemma 89.22. Indeed, there is one hypothesis which is missing. For lemma 89.22, we assumed that the intertwined explorations explore all the boundary clusters of $\Lambda(4n)$ during the first iteration, and this is not any more the case here. Instead, during the first iteration, we explore all the boundary clusters which intersect the box $\Gamma(0)$. This first problem could be dealt with by restricting all the inclusions to $\Gamma(0)$, but this would not be enough. The second problem is that we never manage to control all the relevant taboo sites in $\Gamma(0)$. For $t\geq 1$, we can only control the relevant taboo sites at time $t$ in the box $\Gamma(t-1)$, and we have no hope of doing better. Indeed, as we ignore the boundary clusters which do not intersect $\Gamma(0)$ when choosing the starting set for the intertwined explorations, we must at all costs avoid their influence. As the finite clusters have diameter less than $M$, this influence travels a distance at most $M$ from one iteration to the next. Therefore, the box within which we can achieve a control of the relevant taboo sets must be reduced by a distance at least $M$ at each step. For our strategy to work, we need therefore a version of lemma 89.22 that is compatible with our chain of inequalities.

**Lemma 92.8.** *For $t$ in $\{2,\dots,T\}$ and any $i$ in $I$, we have the inclusion*

$$\widehat{\mathcal{B}}_M\Big(\mathcal{B}\big(w_i, D\setminus\widehat{\mathcal{T}}_i^-(t-1), t-1\big), D, 1\Big)\cap\widehat{\Gamma}^-(t)\ \subset$$

$$\Big(\mathcal{B}\big(w_i, D\setminus\widehat{\mathcal{T}}_i^-(t), t\big)\cap\widehat{\Gamma}^-(t)\Big)\cup B_1\big(\overline{\mathcal{T}_i}(t-1)\cap\Gamma(t-1), 2M\big)$$
$$\cup\, B_1\big(\mathcal{T}_i^{-1}(t)\cap\overset{\circ}{\overbrace{\Gamma(t-1)}}, 2M\big)\cup B_1\big(\mathcal{T}_i^{+1}(t-1)\cap\Gamma(t-1), 2M\big)\,. \quad (92.110)$$

Unfortunately, lemma 92.8 is not a consequence of lemma 89.22. However, the proof follows the same lines as the proof of lemma 89.22. We had in fact two proofs available for lemma 89.22. We choose here to adapt the alternative streamlined proof, which is conceptually more transparent. This proof relied on the lemmas 89.23, 89.24, 89.25, 89.26, which are still valid here (these lemmas do not require that all the boundary clusters of $\Lambda(4n)$ are explored during the first iteration). One could also hope that a few limited modifications of the proof of lemma 89.22 would yield the proof of its localized version. Unfortunately, the need to take into account the spatial propagation of the estimates requires to reorganize seriously the chain of inclusions. In the end, the best option seems to rewrite completely the proof.

*Proof.* We fix $t$ in $\{2,\dots,T\}$ and $i$ in $I$. To alleviate the next formulas, we introduce the notation

$$\widetilde{\Gamma}(t)\ =\ B_1\big(\widehat{\Gamma}^-(t), M+1\big)\,,\qquad \widetilde{\widetilde{\Gamma}}(t)\ =\ B_1\big(\widehat{\Gamma}^-(t), 2M\big)\,.$$

We apply lemma 92.5 to the set $A=\mathcal{B}\big(w_i, D\setminus\widehat{\mathcal{T}}_i^-(t-1), t-1\big)$ and the box $\widehat{\Gamma}^-(t)$, and we get

$$\widehat{\mathcal{B}}_M\Big(\mathcal{B}\big(w_i, D\setminus\widehat{\mathcal{T}}_i^-(t-1), t-1\big), D, 1\Big)\cap\widehat{\Gamma}^-(t)\ \subset$$
$$\widehat{\mathcal{B}}_M\Big(\mathcal{B}\big(w_i, D\setminus\widehat{\mathcal{T}}_i^-(t-1), t-1\big)\cap\widetilde{\Gamma}(t), D\cap\widetilde{\widetilde{\Gamma}}(t), 1\Big)\cap\widehat{\Gamma}^-(t)\,. \quad (92.111)$$

By lemma 89.25, we have

$$\mathcal{B}\big(w_i, D\setminus\widehat{\mathcal{T}}_i^-(t-1), t-1\big)\ \subset$$
$$\mathcal{B}\big(w_i, D\setminus\widehat{\mathcal{T}}_i^-(t), t-1\big)\ \cup\ \mathcal{T}_i^{*0}(t-1)\ \cup\ \text{Shell}\,\big(\mathcal{T}_i^{+1}(t-1), D\big)\,. \quad (92.112)$$

We substitute the inclusion (92.112) in the right-hand side of (92.111) and we obtain, thanks to the monotonicity of $\widehat{\mathcal{B}}_M(\cdot,\cdot,\cdot)$,

$$\widehat{\mathcal{B}}_M\Big(\mathcal{B}\big(w_i, D\setminus\widehat{\mathcal{T}}_i^-(t-1), t-1\big), D, 1\Big)\cap\widehat{\Gamma}^-(t)\ \subset$$
$$\Bigg(\widehat{\mathcal{B}}_M\Big(\mathcal{B}\big(w_i, D\setminus\widehat{\mathcal{T}}_i^-(t), t-1\big)\cap\widetilde{\Gamma}(t), D\cap\widetilde{\widetilde{\Gamma}}(t), 1\Big)$$
$$\cup\,\widehat{\mathcal{B}}_M\Big(\mathcal{T}_i^{*0}(t-1)\cap\widetilde{\Gamma}(t), D\cap\widetilde{\widetilde{\Gamma}}(t), 1\Big)$$
$$\cup\,\widehat{\mathcal{B}}_M\Big(\text{Shell}\,\big(\mathcal{T}_i^{+1}(t-1), D\big)\cap\widetilde{\Gamma}(t), D\cap\widetilde{\widetilde{\Gamma}}(t), 1\Big)\Bigg)\cap\widehat{\Gamma}^-(t)\,. \quad (92.113)$$

We examine next the three terms in the right-hand side of (92.113). For the first term, we apply lemma 89.26 with

$$A = \mathcal{B}\big(w_i, D\setminus \widehat{\mathcal{T}}_i^-(t), t-1\big)\cap\widetilde{\Gamma}(t)\,,\quad \mathcal{T} = \widehat{\mathcal{T}}_i^-(t)\quad\text{in}\quad D\cap\widetilde{\widetilde{\Gamma}}(t)\,.$$

Noticing that

$$\mathcal{T}^0\setminus A\subset\overline{\mathcal{T}}_i^0(t-1)\,,\qquad \mathcal{T}^0\cap A=\varnothing\,,\qquad \mathcal{T}^1\subset\overline{\mathcal{T}}_i^1(t-1)\cup\mathcal{T}_i^{-1}(t)\,,$$

we obtain

$$\begin{aligned}\widehat{\mathcal{B}}_M\Big(\mathcal{B}\big(w_i, D\setminus\widehat{\mathcal{T}}_i^-(t), t-1\big)\cap\widetilde{\Gamma}(t), D\cap\widetilde{\widetilde{\Gamma}}(t), 1\Big)&\subset\mathcal{B}\big(w_i, D\setminus\widehat{\mathcal{T}}_i^-(t), t\big)\\ \cup\big(\overline{\mathcal{T}}_i^0(t-1)\cap\widetilde{\widetilde{\Gamma}}(t)\big)\cup B_1\Big(\big(\overline{\mathcal{T}}_i^1(t-1)\cup\mathcal{T}_i^{-1}(t)\big)&\cap\widetilde{\widetilde{\Gamma}}(t), 2M\Big)\,.\end{aligned}\tag{92.114}$$

Using the inclusions

$$\widetilde{\widetilde{\Gamma}}(t)\subset\overset{\circ}{\widehat{\Gamma(t-1)}}\subset\Gamma(t-1)\,,\tag{92.115}$$

we have

$$\begin{aligned}B_1&\Big(\big(\overline{\mathcal{T}}_i^1(t-1)\cup\mathcal{T}_i^{-1}(t)\big)\cap\widetilde{\widetilde{\Gamma}}(t), 2M\Big)\\ &\qquad\subset B_1\big(\overline{\mathcal{T}}_i(t-1)\cap\widetilde{\widetilde{\Gamma}}(t), 2M\big)\cup B_1\big(\mathcal{T}_i^{-1}(t)\cap\widetilde{\widetilde{\Gamma}}(t), 2M\big)\\ &\quad\subset B_1\big(\overline{\mathcal{T}}_i(t-1)\cap\Gamma(t-1), 2M\big)\cup B_1\big(\mathcal{T}_i^{-1}(t)\cap\overset{\circ}{\widehat{\Gamma(t-1)}}, 2M\big)\,.\end{aligned}\tag{92.116}$$

Substituting (92.116) into (92.114), we obtain

$$\begin{aligned}\widehat{\mathcal{B}}_M\Big(\mathcal{B}\big(w_i, D\setminus\widehat{\mathcal{T}}_i^-(t), t-1\big)\cap\widetilde{\Gamma}(t), D\cap\widetilde{\widetilde{\Gamma}}(t), 1\Big)&\subset\mathcal{B}\big(w_i, D\setminus\widehat{\mathcal{T}}_i^-(t), t\big)\\ \cup B_1\big(\overline{\mathcal{T}}_i(t-1)\cap\Gamma(t-1), 2M\big)\cup B_1\big(\mathcal{T}_i^{-1}(t)&\cap\overset{\circ}{\widehat{\Gamma(t-1)}}, 2M\big)\,.\end{aligned}\tag{92.117}$$

For the second and third terms of (92.113), we apply lemma 89.15, and we get

$$\widehat{\mathcal{B}}_M\Big(\mathcal{T}_i^{*0}(t-1)\cap\widetilde{\Gamma}(t), D\cap\widetilde{\widetilde{\Gamma}}(t), 1\Big)\subset B_1\Big(\mathcal{T}_i^{*0}(t-1)\cap\widetilde{\Gamma}(t), M+1\Big)\,,\tag{92.118}$$

$$\begin{aligned}\widehat{\mathcal{B}}_M\Big(\text{Shell}\,\big(\mathcal{T}_i^{+1}(t-1), D\big)\cap\widetilde{\Gamma}(t), D\cap\widetilde{\widetilde{\Gamma}}(t), 1\Big)&\subset\\ B_1\Big(\text{Shell}\,\big(\mathcal{T}_i^{+1}(t-1), D\big)&\cap\widetilde{\Gamma}(t), M+1\Big)\,.\end{aligned}\tag{92.119}$$

It remains to study $\text{Shell}\,\big(\mathcal{T}_i^{+1}(t-1), D\big)\cap\widetilde{\Gamma}(t)$. The sites in the set $\mathcal{T}_i^{+1}(t-1)$ are visited during the iteration $t-1\geq 1$. However, all the boundary clusters that intersect the box $\Gamma(0)$ are visited during the first step at $t=0$, therefore all the remaining sites of $\Gamma(0)$ belong to finite clusters. As the event $\mathcal{H}'$ occurs, the finite clusters have strictly less that $M$ open sites, and we have the inclusion

$$\text{Shell}\,\big(\mathcal{T}_i^{+1}(t-1), D\big)\cap\widetilde{\Gamma}(t)\subset B_1\big(\mathcal{T}_i^{+1}(t-1), M-1\big)\cap\widetilde{\Gamma}(t)\,.\tag{92.120}$$

Let us precise why we can put $M-1$ in (92.120). If $x$ is in Shell $\big(\mathcal{T}_i^{+1}(t-1), D\big)$, then there exists a path $y_0, \cdots, y_r$ in $D$ joining a site $y_0$ of $\mathcal{T}_i^{+1}(t-1)$ to $x$ such that the sites $y_0, \cdots, y_{r-1}$ are all open. Since the cluster of $y_0$ is finite, then $r < M$ and $|x - y_0|_1 \leq r \leq M-1$. We substitute the inclusion (92.120) into (92.119), and we obtain

$$\widehat{\mathcal{B}}_M\Big(\text{Shell}\,\big(\mathcal{T}_i^{+1}(t-1), D\big) \cap \widetilde{\Gamma}(t), D \cap \widetilde{\widetilde{\Gamma}}(t), 1\Big) \subset \\ B_1\Big(B_1\big(\mathcal{T}_i^{+1}(t-1), M-1\big) \cap \widetilde{\Gamma}(t), M+1\Big)\,. \quad (92.121)$$

The following little result will help to rework the right-hand side of (92.121).

**Lemma 92.9.** *For any subset $A$ of $D$, and any box $\Lambda$ included in $D$, we have*

$$\forall \ell \geq 1 \qquad B_1\big(A, \ell\big) \cap \Lambda \;=\; B_1\big(A \cap B_1(\Lambda, \ell), \ell\big) \cap \Lambda\,. \quad (92.122)$$

*Proof.* The set of the right-hand side of (92.122) is obviously included in the set of the left-hand side. We have to prove the converse inclusion. Let $x$ be a site in $B_1\big(A, \ell\big) \cap \Lambda$. If $x$ is in $A \cap \Lambda$, then $x$ is also in

$$B_1\big(A \cap B_1(\Lambda, \ell), \ell\big)\,. \quad (92.123)$$

Otherwise, there exists a site $y$ in $A$ such that $|y-x|_1 \leq \ell$. Since $x$ is in $\Lambda$, then $y$ is in $A \cap B_1(\Lambda, \ell)$, and we conclude that $x$ is again in the set (92.123). This proves the converse inclusion and it completes the proof of the lemma. ☐

We apply lemma 92.9 with $A = \mathcal{T}_i^{+1}(t-1)$, $\Lambda = \widetilde{\Gamma}(t)$, $\ell = M-1$, and we get

$$B_1\big(\mathcal{T}_i^{+1}(t-1), M-1\big) \cap \widetilde{\Gamma}(t) \;=\; B_1\Big(\mathcal{T}_i^{+1}(t-1) \cap B_1\big(\widetilde{\Gamma}(t), M-1\big), M-1\Big) \cap \widetilde{\Gamma}(t) \\ \subset\; B_1\big(\mathcal{T}_i^{+1}(t-1) \cap \Gamma(t-1), M-1\big)\,. \quad (92.124)$$

We substitute the inclusion (92.124) in (92.121), and we obtain

$$\widehat{\mathcal{B}}_M\Big(\text{Shell}\,\big(\mathcal{T}_i^{+1}(t-1), D\big) \cap \widetilde{\Gamma}(t), D \cap \widetilde{\widetilde{\Gamma}}(t), 1\Big) \\ \subset\; B_1\Big(B_1\big(\mathcal{T}_i^{+1}(t-1) \cap \Gamma(t-1), M-1\big), M+1\Big) \\ \subset\; B_1\big(\mathcal{T}_i^{+1}(t-1) \cap \Gamma(t-1), 2M\big)\,. \quad (92.125)$$

The inclusions (92.113), (92.117), (92.118), (92.125) yield the result (92.110) stated in lemma 92.8. ☐

We take the union over $i$ in $I$ of the inclusion (92.110), and we get

$$\bigcup_{i \in I} \widehat{\mathcal{B}}_M\Big(\mathcal{B}\big(w_i, D \setminus \widehat{\mathcal{T}}_i^-(t-1), t-1\big), D, 1\Big) \cap \widehat{\Gamma}^-(t) \;\subset \\ \bigcup_{i \in I}\Bigg(\Big(\mathcal{B}\big(w_i, D \setminus \widehat{\mathcal{T}}_i^-(t), t\big) \cap \widehat{\Gamma}^-(t)\Big) \cup B_1\big(\overline{\mathcal{T}}_i(t-1) \cap \Gamma(t-1), 2M\big) \\ \cup\, B_1\big(\mathcal{T}_i^{-1}(t) \cap \overset{\circ}{\widehat{\Gamma(t-1)}}, 2M\big) \cup B_1\big(\mathcal{T}_i^{+1}(t-1) \cap \Gamma(t-1), 2M\big)\Bigg)\,. \quad (92.126)$$

For any $i$ in $I$, we have

$$B_1\big(\overline{\mathcal{T}}_i(t-1)\cap\Gamma(t-1),2M\big)\,\subset\,\overline{B}^{t-1}_{1,\Gamma}(1)\,,\tag{92.127}$$

$$B_1\big(\mathcal{T}_i^{-1}(t)\cap\overset{\circ}{\widehat{\Gamma(t-1)}},2M\big)\cup B_1\big(\mathcal{T}_i^{+1}(t-1)\cap\Gamma(t-1),2M\big)\,\subset\,\widehat{B}^{t,1}_{1,\Gamma}(1)\,.\tag{92.128}$$

Using the inclusions (92.127) and (92.128), we can rewrite (92.126) as

$$\begin{multline}\bigcup_{i\in I}\widehat{\mathcal{B}}_M\Big(\mathcal{B}\big(w_i,D\setminus\widehat{\mathcal{T}}_i^-(t-1),t-1\big),D,1\Big)\cap\widehat{\Gamma}^-(t)\,\subset\\ \Big(\bigcup_{i\in I}\mathcal{B}\big(w_i,D\setminus\widehat{\mathcal{T}}_i^-(t),t\big)\cap\widehat{\Gamma}^-(t)\Big)\,\cup\,\overline{B}^{t-1}_{1,\Gamma}(1)\,\cup\,\widehat{B}^{t,1}_{1,\Gamma}(1)\,.\end{multline}\tag{92.129}$$

For the second set of the right-hand side of (92.109), we use lemma 89.15 to obtain

$$\begin{multline}\widehat{\mathcal{B}}_M\Bigg(\bigcup_{1\le s\le t-1}\Big(\overline{B}^{s-1}_{1,\Gamma}(t-s)\cup\widehat{B}^{s,1}_{1,\Gamma}(t-s)\Big),D,1\Bigg)\\ \subset\,B_1\Big(\bigcup_{1\le s\le t-1}\Big(\overline{B}^{s-1}_{1,\Gamma}(t-s)\cup\widehat{B}^{s,1}_{1,\Gamma}(t-s)\Big),M+1\Big)\\ \subset\bigcup_{1\le s\le t-1}\Big(\overline{B}^{s-1}_{1,\Gamma}(t+1-s)\cup\widehat{B}^{s,1}_{1,\Gamma}(t+1-s)\Big)\,.\end{multline}\tag{92.130}$$

Using the inclusions (92.129) and (92.130), we deduce from (92.109) that

$$\begin{multline}\widehat{\mathcal{B}}_M\Big(\text{Shell}\,\big(W(I),D\big),D,t\Big)\cap\widehat{\Gamma}^-(t)\,\subset\\ \Big(\bigcup_{i\in I}\mathcal{B}\big(w_i,D\setminus\widehat{\mathcal{T}}_i^-(t),t\big)\cap\widehat{\Gamma}^-(t)\Big)\,\cup\,\overline{B}^{t-1}_{1,\Gamma}(1)\,\cup\,\widehat{B}^{t,1}_{1,\Gamma}(1)\\ \cup\bigcup_{1\le s\le t-1}\Big(\overline{B}^{s-1}_{1,\Gamma}(t+1-s)\cup\widehat{B}^{s,1}_{1,\Gamma}(t+1-s)\Big)\,.\end{multline}\tag{92.131}$$

We check next that

$$\begin{multline}\overline{B}^{t-1}_{1,\Gamma}(1)\,\cup\,\widehat{B}^{t,1}_{1,\Gamma}(1)\,\cup\bigcup_{1\le s\le t-1}\Big(\overline{B}^{s-1}_{1,\Gamma}(t+1-s)\cup\widehat{B}^{s,1}_{1,\Gamma}(t+1-s)\Big)\\ =\bigcup_{1\le s\le t}\Big(\overline{B}^{s-1}_{1,\Gamma}(t+1-s)\cup\widehat{B}^{s,1}_{1,\Gamma}(t+1-s)\Big)\,.\end{multline}\tag{92.132}$$

The inclusion (92.100) at rank $t$ follows from formulas (92.131) and (92.132). $\square$

## 92.12 Control of $\left|\mathcal{T}^{-0}(t) \cap \widehat{\Gamma}^{-}(t)\right|$ and $\left|\mathcal{T}^{+1}(t) \cap \widehat{\overset{\circ}{\widehat{\Gamma}^{-}(t)}}\right|$

The starting point to control $\left|\mathcal{T}^{-0}(t) \cap \widehat{\Gamma}^{-}(t)\right|$ is the inequality (92.68): there exist two disjoint subsets $K^{-}(t), L^{-}(t)$ of $\{1, \dots, N\}$ such that, setting

$$\begin{aligned} \overline{\mathcal{L}}^{-}(t) &= \bigcup_{\ell \in L^{-}(t)} \mathcal{B}\big(w_\ell, D \setminus \overline{\mathcal{T}}_\ell(t-1), t-1\big) \cap \widehat{\Gamma}^{-}(t)\,, \\ \widehat{\mathcal{K}}^{-}(t) &= \bigcup_{k \in K^{-}(t)} \mathcal{B}\big(w_k, D \setminus \widehat{\mathcal{T}}_k^{-}(t), t\big) \cap \widehat{\Gamma}^{-}(t)\,, \end{aligned}$$

we have

$$\frac{1}{8d}\left|\mathcal{T}^{-0}(t) \cap \widehat{\Gamma}^{-}(t)\right| \leq \left|\overline{\mathcal{T}}^{0}(t-1) \cap \widehat{\Gamma}^{-}(t)\right| + \left|\widehat{\mathcal{K}}^{-}(t) \cap \overline{\mathcal{L}}^{-}(t)\right|. \tag{92.133}$$

We use next the approximations of the two unions in (92.133) given in propositions 92.4 and 92.7. The set $\overline{\mathcal{L}}^{-}(t)$ is approximated by the truncated genuine explorations $\overline{\mathcal{L}}_M^{-}(t)$ defined by

$$\overline{\mathcal{L}}_M^{-}(t) = \bigcup_{\ell \in L^{-}(t)} \widehat{\mathcal{B}}_M\big(\mathrm{Shell}\,(w_\ell, D), D, t-1\big) \cap \widehat{\Gamma}^{-}(t)\,.$$

We apply proposition 92.4 at iteration $t-1$ and we take the trace on the box $\widehat{\Gamma}^{-}(t)$ of the inclusion (92.72). This yields

$$\overline{\mathcal{L}}^{-}(t) \subset \overline{\mathcal{L}}_M^{-}(t)\,, \qquad \overline{\mathcal{L}}_M^{-}(t) \setminus \overline{\mathcal{L}}^{-}(t) \subset \overline{\mathrm{Error}}_{\widehat{\Gamma}^{-}}(M,t)\,, \tag{92.134}$$

where the error term $\overline{\mathrm{Error}}_{\widehat{\Gamma}^{-}}(M,t)$ is given by

$$\overline{\mathrm{Error}}_{\widehat{\Gamma}^{-}}(M,t) = \bigcup_{1 \leq s \leq t-1} B_{1,\Gamma}^{s,-}(t-s) \cup \bigcup_{0 \leq s \leq t-1} \big(B_{1,\Gamma}^{s,0}(t-s-1) \cup B_{1,\Gamma}^{s,1}(t-s)\big)\,.$$

The important point is that this first error term involves only the relevant taboo sites up to time $t-1$ (it is the superscript $s$ in the sets $B_{1,\Gamma}^{s,\cdot}(\cdot)$ that indicates the iteration step). The set $\widehat{\mathcal{K}}^{-}(t)$ is approximated by the truncated genuine explorations $\widehat{\mathcal{K}}_M^{-}(t)$ defined by

$$\widehat{\mathcal{K}}_M^{-}(t) = \bigcup_{k \in K^{-}(t)} \widehat{\mathcal{B}}_M\big(\mathrm{Shell}\,(w_k, D), D, t\big) \cap \widehat{\Gamma}^{-}(t)\,.$$

We apply proposition 92.7 at iteration $t$. This yields

$$\widehat{\mathcal{K}}^{-}(t) \subset \widehat{\mathcal{K}}_M^{-}(t)\,, \qquad \widehat{\mathcal{K}}_M^{-}(t) \setminus \widehat{\mathcal{K}}^{-}(t) \subset \widehat{\mathrm{Error}}_{\widehat{\Gamma}^{-}}^{-}(M,t)\,, \tag{92.135}$$

where the error term $\widehat{\mathrm{Error}}_{\widehat{\Gamma}^{-}}^{-}(M,t)$ is given by

$$\widehat{\mathrm{Error}}_{\widehat{\Gamma}^{-}}^{-}(M,t) = \bigcup_{1 \leq s \leq t} \Big(\overline{B}_{1,\Gamma}^{s-1}(t+1-s) \cup \widehat{B}_{1,\Gamma}^{s,1}(t+1-s)\Big)\,.$$

We notice that this second error term involves the relevant taboo sites up to time $t-1$, but also the set $\mathcal{T}^{-1}(t)$. The point is that this last set has been controlled in subsection 92.8 by $\overline{\mathcal{T}}^{+0}(t-1)$. We recall here a key point. Lemma 88.3 tells us that all the sites in the set $\widehat{\mathcal{K}}^-(t)\cap\overline{\mathcal{L}}^-(t)$ are closed. As the genuine approximations $\overline{\mathcal{L}}_M^-(t)$, $\widehat{\mathcal{K}}_M^-(t)$ are larger than the sets $\overline{\mathcal{L}}^-(t)$, $\widehat{\mathcal{K}}^-(t)$, we have

$$\big\{\,x\in\widehat{\mathcal{K}}_M^-(t)\cap\overline{\mathcal{L}}_M^-(t):x\text{ is open}\,\big\}\,\subset\,\big(\widehat{\mathcal{K}}_M^-(t)\cap\overline{\mathcal{L}}_M^-(t)\big)\setminus\big(\widehat{\mathcal{K}}^-(t)\cap\overline{\mathcal{L}}^-(t)\big)\,,\tag{92.136}$$

and the inequality (92.133) implies that

$$\frac{1}{8d}\,\big|\mathcal{T}^{-0}(t)\cap\widehat{\Gamma}^-(t)\big|\,\leq\,\big|\overline{\mathcal{T}}^0(t-1)\cap\widehat{\Gamma}^-(t)\big|+\Big|\big\{\,x\in\widehat{\mathcal{K}}_M^-(t)\cap\overline{\mathcal{L}}_M^-(t):x\text{ is closed}\,\big\}\Big|\,.\tag{92.137}$$

We try next to control the right-hand set in (92.136). It follows from (92.134) and (92.135) that

$$\big(\widehat{\mathcal{K}}_M^-(t)\cap\overline{\mathcal{L}}_M^-(t)\big)\setminus\big(\widehat{\mathcal{K}}^-(t)\cap\overline{\mathcal{L}}^-(t)\big)\,\subset\,\widehat{\text{Error}}_{\widehat{\Gamma}^-}^-(M,t)\cup\overline{\text{Error}}_{\widehat{\Gamma}^-}(M,t)\,.\tag{92.138}$$

The two inclusions (92.136) and (92.138) imply that

$$\Big|\big\{\,x\in\widehat{\mathcal{K}}_M^-(t)\cap\overline{\mathcal{L}}_M^-(t):x\text{ is open}\,\big\}\Big|\,\leq\,\Big|\widehat{\text{Error}}_{\widehat{\Gamma}^-}^-(M,t)\cup\overline{\text{Error}}_{\widehat{\Gamma}^-}(M,t)\Big|\,.\tag{92.139}$$

Let us try to compute upper bounds on the cardinalities of the error sets $\widehat{\text{Error}}_{\widehat{\Gamma}^-}^-(M,t)$, $\overline{\text{Error}}_{\widehat{\Gamma}^-}(M,t)$. We recall that, on the event $\mathcal{E}$, the termination time is bounded by $\bar{t}(n)=2(\ln n)^{1-\gamma}$, so we need only to consider values of $t$ smaller than or equal to $\bar{t}(n)$. We wish to get upper bounds on the cardinalities of the sets

$$B_{1,\Gamma}^{t,-}(s)\,,\quad B_{1,\Gamma}^{t,0}(s)\,,\quad B_{1,\Gamma}^{t,1}(s)\,,\quad \overline{B}_{1,\Gamma}^{t}(s)\,,\quad \widehat{B}_{1,\Gamma}^{t,1}(s)\,,$$

defined respectively in (92.71), (92.69), (92.70), (92.98), (92.99). We use the following bounds on the balls $B_1(\cdot,\cdot)$:

$$\forall x\in\mathbb{Z}^d\quad\forall r\geq 1\qquad\big|B_1(x,r)\big|\,\leq\,(2r+1)^d\,\leq\,(3r)^d\,.$$

To alleviate the formulas, we set

$$D(n)\,=\,\big(4\bar{t}(n)M\big)^d\,.$$

It follows immediately from the definitions of the sets that, for any $t$ such that $0\leq t\leq\bar{t}(n)$, we have

$$\begin{aligned}
\big|B_{1,\Gamma}^{s,-}(t-s)\big|\,&\leq\,\big|\overline{\mathcal{T}}(s-1)\cap\Gamma(s)\big|\,D(n)\,,\quad 1\leq s\leq t-1\,,\\
\big|B_{1,\Gamma}^{s,0}(t-s-1)\big|\,&\leq\,\big|\mathcal{T}^{*0}(s)\cap\Gamma(s)\big|\,D(n)\,,\quad 0\leq s\leq t-1\,,\\
\big|B_{1,\Gamma}^{s,1}(t-s)\big|\,&\leq\,\big|\mathcal{T}^{*1}(s)\cap\Gamma(s)\big|\,D(n)\,,\quad 0\leq s\leq t-1\,,\\
\big|\overline{B}_{1,\Gamma}^{s-1}(t+1-s)\big|\,&\leq\,\big|\overline{\mathcal{T}}(s-1)\cap\Gamma(s-1)\big|\,D(n)\,,\quad 1\leq s\leq t\,,
\end{aligned}$$

and for the last set, for $1 \leq s \leq t$,

$$\left|\widehat{B}^{s,1}_{1,\Gamma}(t+1-s)\right| \;\leq\; \left|\mathcal{T}^{-1}(s)\cap \overset{\circ}{\widehat{\Gamma(s-1)}}\right| D(n) + \left|\mathcal{T}^{+1}(s-1)\cap\Gamma(s-1)\right| D(n)\,. \tag{92.140}$$

Using the standard union bound, we deduce from the previous inequalities that, for any $t$ such that $0 \leq t \leq \bar{t}(n)$, we have

$$\Big|\bigcup_{1\leq s\leq t-1} B^{s,-}_{1,\Gamma}(t-s)\Big| \;\leq\; D(n)\sum_{1\leq s\leq t-1}\left|\overline{\mathcal{T}}(s-1)\cap\Gamma(s)\right|, \tag{92.141}$$

$$\Big|\bigcup_{0\leq s\leq t-1} B^{s,0}_{1,\Gamma}(t-s-1)\Big| \;\leq\; D(n)\sum_{0\leq s\leq t-1}\left|\mathcal{T}^{*0}(s)\cap\Gamma(s)\right|, \tag{92.142}$$

$$\Big|\bigcup_{0\leq s\leq t-1} B^{s,1}_{1,\Gamma}(t-s)\Big| \;\leq\; D(n)\sum_{0\leq s\leq t-1}\left|\mathcal{T}^{*1}(s)\cap\Gamma(s)\right|, \tag{92.143}$$

$$\Big|\bigcup_{1\leq s\leq t} \overline{B}^{s-1}_{1,\Gamma}(t+1-s)\Big| \;\leq\; D(n)\sum_{1\leq s\leq t}\left|\overline{\mathcal{T}}(s-1)\cap\Gamma(s-1)\right|. \tag{92.144}$$

We deal now with the last set. From (92.140), we deduce that

$$\Big|\bigcup_{1\leq s\leq t} \widehat{B}^{s,1}_{1,\Gamma}(t+1-s)\Big| \;\leq\; D(n)\sum_{1\leq s\leq t}\left|\mathcal{T}^{-1}(s)\cap \overset{\circ}{\widehat{\Gamma(s-1)}}\right| + D(n)\sum_{1\leq s\leq t}\left|\mathcal{T}^{+1}(s-1)\cap\Gamma(s-1)\right|. \tag{92.145}$$

The sets $\mathcal{T}^{-1}(s)$, $1 \leq s \leq t$, were controlled in the subsection 92.8. Thanks to the inequality (92.60), we have

$$\sum_{1\leq s\leq t}\left|\mathcal{T}^{-1}(s)\cap \overset{\circ}{\widehat{\Gamma(s-1)}}\right| \;\leq\; 2d\sum_{1\leq s\leq t}\left|\mathcal{T}^{+0}(s-1)\cap\Gamma(s-1)\right|. \tag{92.146}$$

Putting together the inequalities (92.145) and (92.146), we obtain

$$\begin{aligned}\Big|\bigcup_{1\leq s\leq t} \widehat{B}^{s,1}_{1,\Gamma}(t+1-s)\Big| \;\leq\;& \\ 2dD(n)\sum_{1\leq s\leq t}\left|\mathcal{T}^{+0}(s-1)\cap\Gamma(s-1)\right| &+ D(n)\sum_{1\leq s\leq t}\left|\mathcal{T}^{+1}(s-1)\cap\Gamma(s-1)\right| \\ &\leq\; 3dD(n)\sum_{1\leq s\leq t}\left|\overline{\mathcal{T}}(s-1)\cap\Gamma(s-1)\right|. \end{aligned} \tag{92.147}$$

Adding together the inequalities (92.141), (92.142), (92.143), (92.144), (92.147), we get

$$\Big|\widehat{\mathrm{Error}}_{\widehat{\Gamma}^-}(M,t)\cup\overline{\mathrm{Error}}_{\widehat{\Gamma}^-}(M,t)\Big| \,\leq\, (2+3d)D(n)\sum_{1\leq s\leq t}\big|\overline{\mathcal{T}}(s-1)\cap\Gamma(s-1)\big|$$
$$+\,D(n)\sum_{0\leq s\leq t-1}\big|\big(\mathcal{T}^{*0}(s)\cup\mathcal{T}^{*1}(s)\big)\cap\Gamma(s)\big|$$
$$\leq\,4dD(n)\sum_{1\leq s\leq t}\big|\overline{\mathcal{T}}(s-1)\cap\Gamma(s-1)\big|\,. \quad (92.148)$$

Substituting (92.148) in (92.139), we get

$$\Big|\big\{\,x\in\widehat{\mathcal{K}}^-_M(t)\cap\overline{\mathcal{L}}^-_M(t):x\text{ is open}\,\big\}\Big| \,\leq\, 4dD(n)\sum_{1\leq s\leq t}\big|\overline{\mathcal{T}}(s-1)\cap\Gamma(s-1)\big|\,. \quad (92.149)$$

We finally use the fact the event Intersect$(N)$ occurs. The sets $\widehat{\mathcal{K}}^-_M(t)$, $\overline{\mathcal{L}}^-_M(t)$ are among those concerned by the event Intersect$(N)$, whence

$$S\big(\widehat{\mathcal{K}}^-_M(t)\cap\overline{\mathcal{L}}^-_M(t)\big)\,>\,-c_Hn^{d/2}\sqrt{N\ln n}\,, \quad (92.150)$$

where the constant $c_H$ is defined in (92.50). The inequalities (92.149), (92.150) imply that

$$\Big|\big\{\,x\in\widehat{\mathcal{K}}^-_M(t)\cap\overline{\mathcal{L}}^-_M(t):x\text{ is closed}\,\big\}\Big|\,\leq$$
$$\frac{1}{p}4dD(n)\sum_{1\leq s\leq t}\big|\overline{\mathcal{T}}(s-1)\cap\Gamma(s-1)\big|+c_Hn^{d/2}\sqrt{N\ln n}\,. \quad (92.151)$$

We substitute (92.151) in (92.137), we multiply by $8d$, and this gives

$$\big|\mathcal{T}^{-0}(t)\cap\widehat{\Gamma}^-(t)\big|\,\leq\,8d\big|\overline{\mathcal{T}}^0(t-1)\cap\widehat{\Gamma}^-(t)\big|+$$
$$32\frac{d^2}{p}D(n)\sum_{1\leq s\leq t}\big|\overline{\mathcal{T}}(s-1)\cap\Gamma(s-1)\big|+16dc_Hn^{d/2}\sqrt{N\ln n}\,. \quad (92.152)$$

Recalling that $\widehat{\Gamma}^-(t)$ is a sub-box of $\Gamma(t-1)$, the first term of the right-hand side of (92.152) can be absorbed in the second one, and we conclude that

$$\big|\mathcal{T}^{-0}(t)\cap\widehat{\Gamma}^-(t)\big|\,\leq\,33\frac{d^2}{p}D(n)\sum_{0\leq s\leq t-1}\big|\overline{\mathcal{T}}(s)\cap\Gamma(s)\big|+16dc_Hn^{d/2}\sqrt{N\ln n}\,. \quad (92.153)$$

As in the subsection 92.8, we use lemma 88.1 in order to control the set $\mathcal{T}^{+1}(t)$, and we obtain the following inequality:

$$\big|\mathcal{T}^{+1}(t)\cap\overset{\circ}{\widehat{\widehat{\Gamma}^-(t)}}\big|\,\leq\,2d\,\big|\mathcal{T}^{-0}(t)\cap\widehat{\Gamma}^-(t)\big|\,. \quad (92.154)$$

## 92.13 Localization of $\mathcal{T}^{+0}(t) \cap \widehat{\Gamma}^*(t)$

We focus here on the set $\mathcal{T}^{+0}(t) \cap \widehat{\Gamma}^*(t)$. By the definition (80.39), a taboo site in $\mathcal{T}^{+0}(t)$ is visited at time $t$, therefore

$$\mathcal{T}^{+0}(t) \cap \widehat{\Gamma}^*(t) \,=\, \bigcup_{1\leq i,j\leq N, i\neq j} \mathcal{T}_i^{+0}(t) \cap \mathcal{A}_j(t) \cap \widehat{\Gamma}^*(t)\,.$$

A site $x$ can belong to at most $2d$ sets of the family

$$\mathcal{T}_i^{+0}(t) \cap \mathcal{A}_j(t) \cap \widehat{\Gamma}^*(t)\,, \quad 1 \leq i,j \leq N,\, i \neq j\,.$$

By lemma 84.3 applied to this family with $m = 2d$, there exist two disjoint subsets $K^+(t), L^+(t)$ of $\{\,1,\dots,N\,\}$ such that

$$\frac{1}{8d}\left|\mathcal{T}^{+0}(t) \cap \widehat{\Gamma}^*(t)\right| \,\leq\, \bigg|\Big(\bigcup_{k\in K^+(t)} \mathcal{T}_k^{+0}(t)\Big) \cap \Big(\bigcup_{\ell\in L^+(t)} \mathcal{A}_\ell(t)\Big) \cap \widehat{\Gamma}^*(t)\bigg|\,. \quad (92.155)$$

Contrary to the case of the set $\mathcal{T}^{-0}(t)$ studied in subsection 92.9, the set $\mathcal{T}^{+0}(t)$ is disjoint from $\overline{\mathcal{T}}(t-1)$. Indeed, the visiting time of the sites in $\mathcal{T}^{+0}(t)$ is equal to $t$, hence they cannot belong to $\overline{\mathcal{T}}(t-1)$ as well. We have thus

$$\forall k \in K^+(t) \qquad \mathcal{T}_k^{+0}(t) \,\subset\, \mathcal{T}_k^{*0}(t) \setminus \overline{\mathcal{T}}_k(t-1)\,. \qquad (92.156)$$

For $i$ in $\{\,1,\dots,N\,\}$, we define

$$\widehat{\mathcal{T}}_i^*(t) \,=\, \overline{\mathcal{T}}_i(t-1) \cup \mathcal{T}_i^{-0}(t) \cup \mathcal{T}_i^{*1}(t)\,.$$

It follows from proposition 81.2 that

$$\forall k \in K^+(t) \qquad \mathcal{T}_k^{*0}(t) \setminus \overline{\mathcal{T}}_k(t-1) \,\subset\, \mathcal{B}\big(w_k, D \setminus \widehat{\mathcal{T}}_k^*(t), t\big)\,. \qquad (92.157)$$

We deduce from the two inclusions (92.156) and (92.157) that

$$\bigcup_{k\in K^+(t)} \mathcal{T}_k^{+0}(t) \,\subset\, \bigcup_{k\in K^+(t)} \mathcal{B}\big(w_k, D \setminus \widehat{\mathcal{T}}_k^*(t), t\big)\,. \qquad (92.158)$$

The inclusion (92.158) readily implies that

$$\Big(\bigcup_{k\in K^+(t)} \mathcal{T}_k^{+0}(t)\Big) \cap \Big(\bigcup_{\ell\in L^+(t)} \mathcal{A}_\ell(t)\Big) \,\subset\, \Big(\bigcup_{k\in K^+(t)} \mathcal{B}\big(w_k, D \setminus \widehat{\mathcal{T}}_k^*(t), t\big)\Big) \cap \Big(\bigcup_{\ell\in L^+(t)} \mathcal{A}_\ell(t)\Big)\,. \quad (92.159)$$

We deduce from (92.155) and (92.159) that

$$\frac{1}{8d}\left|\mathcal{T}^{+0}(t) \cap \widehat{\Gamma}^*(t)\right| \,\leq\, \bigg|\Big(\bigcup_{k\in K^+(t)} \mathcal{B}\big(w_k, D \setminus \widehat{\mathcal{T}}_k^*(t), t\big)\Big) \cap \Big(\bigcup_{\ell\in L^+(t)} \mathcal{A}_\ell(t)\Big) \cap \widehat{\Gamma}^*(t)\bigg|\,. \quad (92.160)$$

We still need to rework a bit the union $\bigcup_{\ell\in L^+(t)} \mathcal{A}_\ell(t)$. Firstly, by proposition 81.1, we have

$$\forall \ell \in L^+(t) \qquad \mathcal{A}_\ell(t) \,=\, \mathcal{B}\big(w_\ell, D\setminus \overline{\mathcal{T}}_\ell(t), t\big)\,. \tag{92.161}$$

Secondly, since our goal is to control the set $\mathcal{T}_\ell^{+0}(t)$, we try to get rid of it in (92.161). Thanks to lemma 88.2, we have

$$\mathcal{B}\big(w_\ell, D\setminus \widehat{\mathcal{T}}_\ell^*(t), t\big) \,=\, \mathcal{A}_\ell(t)\,\cup\,\mathcal{T}_\ell^{+0}(t)\,. \tag{92.162}$$

The equality (92.162) and the inclusion (92.160) together imply that

$$\frac{1}{8d}\left|\mathcal{T}^{+0}(t)\cap\widehat{\Gamma}^*(t)\right| \leq \\ \Bigg|\Big(\bigcup_{k\in K^+(t)}\mathcal{B}\big(w_k, D\setminus \widehat{\mathcal{T}}_k^*(t), t\big)\Big)\cap\Big(\bigcup_{\ell\in L^+(t)}\mathcal{B}\big(w_\ell, D\setminus \widehat{\mathcal{T}}_\ell^*(t), t\big)\Big)\cap\widehat{\Gamma}^*(t)\Bigg|\,. \tag{92.163}$$

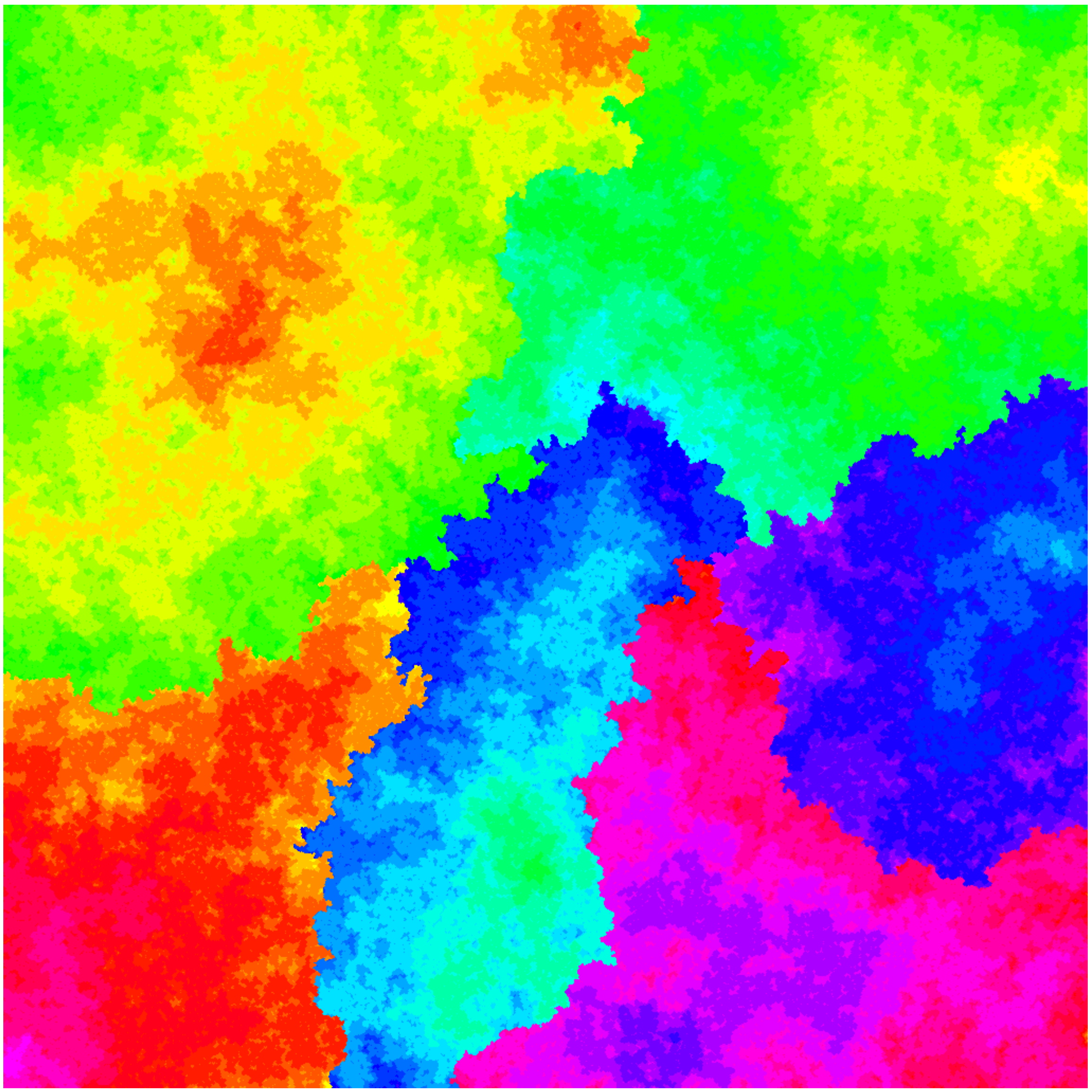

Figure 170: 7 intertwined explorations, $1024\times 1024$, site percolation, $p=0.57$.

## 92.14 Approximation of $\mathcal{B}\big(w_i, D\setminus\widehat{\mathcal{T}}_i^*(t),t\big)$

In this subsection, we try to approximate the taboo exploration

$$\bigcup_{i\in I}\mathcal{B}\big(w_i, D\setminus\widehat{\mathcal{T}}_i^*(t),t\big)\,.$$

Once more, we have to develop a local version of the approximation given in propositions 89.21, 89.20. Fortunately, we can rely on the work done previously for the set $\mathcal{B}\big(w_i, D\setminus\widehat{\mathcal{T}}_i^-(t),t\big)$. So, we will reuse the sets $\overline{B}^{\,t}_{1,\Gamma}(s)$, $\widehat{B}^{t,1}_{1,\Gamma}(s)$ defined in (92.98), (92.99), that we recall for convenience:

$$\forall t\in\{\,0,\dots,T\,\}\quad\forall s\geq 0\qquad \overline{B}^{\,t}_{1,\Gamma}(s)\;=\;B_1\big(\overline{\mathcal{T}}(t)\cap\Gamma(t),s(M+1)+M-1\big)\,,$$

$$\forall t\in\{\,1,\dots,T\,\}\quad\forall s\geq 0$$

$$\begin{aligned}\widehat{B}^{t,1}_{1,\Gamma}(s)\;=\;&B_1\Big(\mathcal{T}^{-1}(t)\cap\overbrace{\Gamma(t-1)}^{\circ},s(M+1)+M-1\Big)\\&\cup B_1\big(\mathcal{T}^{+1}(t-1)\cap\Gamma(t-1),s(M+1)+M-1\big)\,.\end{aligned}$$

We will base our approximation result on that of proposition 92.7. Although a derivation from scratch would probably lead to slightly better bounds, the following version will suffice.

**Proposition 92.10.** *We consider a configuration realizing the event $\mathcal{H}'$ defined in (92.27) and a subset $I$ of $\{\,1,\dots,N\,\}$. Setting $W(I)=\{\,w_i:i\in I\,\}$, we have, for any $t$ in $\{\,1,\dots,T\,\}$,*

$$\begin{aligned}\bigcup_{i\in I}\mathcal{B}\big(w_i, D\setminus\widehat{\mathcal{T}}_i^*(t),t\big)\cap\widehat{\Gamma}^*(t)\;\subset\;&\bigcup_{i\in I}\mathcal{B}\big(w_i, D\setminus\widehat{\mathcal{T}}_i^-(t),t\big)\cap\widehat{\Gamma}^*(t)\\&\subset\;\widehat{\mathcal{B}}_M\Big(\mathit{Shell}\big(W(I),D\big),D,t\Big)\cap\widehat{\Gamma}^*(t)\;\subset\end{aligned}\tag{92.164}$$

$$\begin{aligned}&\Big(\bigcup_{i\in I}\mathcal{B}\big(w_i, D\setminus\widehat{\mathcal{T}}_i^-(t),t\big)\cap\widehat{\Gamma}^*(t)\Big)\cup\bigcup_{1\leq s\leq t}\Big(\overline{B}^{\,s-1}_{1,\Gamma}(t+1-s)\cup\widehat{B}^{s,1}_{1,\Gamma}(t+1-s)\Big)\\&\subset\;\Big(\bigcup_{i\in I}\mathcal{B}\big(w_i, D\setminus\widehat{\mathcal{T}}_i^*(t),t\big)\cap\widehat{\Gamma}^*(t)\Big)\cup\bigcup_{1\leq s\leq t}\Big(\overline{B}^{\,s-1}_{1,\Gamma}(t+1-s)\cup\widehat{B}^{s,1}_{1,\Gamma}(t+1-s)\Big)\\&\qquad\cup\;\big(\mathcal{T}^{-0}(t)\cap\widehat{\Gamma}^*(t)\big)\;\cup\;B_1\big(\mathcal{T}^{+1}(t)\cap\overbrace{\widehat{\Gamma}^-(t)}^{\circ},M-1\big)\,.\end{aligned}$$

*Proof.* Let us fix $t$ in $\{\,1,\dots,T\,\}$. For any $i$ in $I$, the set $\widehat{\mathcal{T}}_i^*(t)$ contains $\widehat{\mathcal{T}}_i^-(t)$, whence

$$\mathcal{B}\big(w_i, D\setminus\widehat{\mathcal{T}}_i^*(t),t\big)\cap\widehat{\Gamma}^*(t)\;\subset\;\mathcal{B}\big(w_i, D\setminus\widehat{\mathcal{T}}_i^-(t),t\big)\cap\widehat{\Gamma}^*(t)\,.$$

Taking the union over $i$ in $I$, we obtain the first inclusion in (92.164). The second and third inclusions in (92.164) are immediate consequences of proposition 92.7 and the fact that the box $\widehat{\Gamma}^*(t)$ is included in the box $\widehat{\Gamma}^-(t)$. It remains to prove the fourth inclusion in (92.164). To that end, we will use the next lemma, which is a local version of lemma 89.27.

**Lemma 92.11.** *For any $t$ in $\{1,\dots,T\}$ and any $i$ in $I$, we have the inclusion*

$$\mathcal{B}\big(w_i, D\setminus\widehat{\mathcal{T}}_i^-(t),t\big)\cap\widehat{\Gamma}^*(t)\subset\Big(\mathcal{B}\big(w_i, D\setminus\widehat{\mathcal{T}}_i^*(t),t\big)\cap\widehat{\Gamma}^*(t)\Big)$$
$$\cup\big(\mathcal{T}_i^{-0}(t)\cap\widehat{\Gamma}^*(t)\big)\cup B_1\big(\mathcal{T}_i^{+1}(t)\cap\overset{\circ}{\widehat{\widehat{\Gamma}^-(t)}},M-1\big)\,. \quad (92.165)$$

*Proof.* There were two proofs for lemma 89.27. Although the first one is more automatic, the second one is shorter. We choose here to adapt the second proof. Let $x$ belong to $\mathcal{B}\big(w_i, D\setminus\widehat{\mathcal{T}}_i^-(t),t\big)\cap\widehat{\Gamma}^*(t)$. By definition, there exists a path $y_0,\cdots,y_r$ in $D\setminus\widehat{\mathcal{T}}_i^-(t)$ from $y_0=w_i$ to $y_r=x$ such that the number of closed sites among $y_0,\cdots,y_{r-1}$ is less than or equal to $t$. We consider several cases:

• No site among $y_0,\cdots,y_r$ is in $\widehat{\mathcal{T}}_i^*(t)$. Then $x$ is in $\mathcal{B}\big(w_i,D\setminus\widehat{\mathcal{T}}_i^*(t),t\big)\cap\widehat{\Gamma}^*(t)$.

• The path $y_0,\cdots,y_r$ visits $\widehat{\mathcal{T}}_i^*(t)$. We know also that this path does not visit $\widehat{\mathcal{T}}_i^-(t)$. Let $s$ be the smallest index such that $y_s$ is in $\widehat{\mathcal{T}}_i^*(t)$. The site $y_s$ is in

$$\widehat{\mathcal{T}}_i^*(t)\setminus\widehat{\mathcal{T}}_i^-(t)\;\subset\;\big(\mathcal{T}_i^{-0}(t)\cup\mathcal{T}_i^{+1}(t)\big)\setminus\overline{\mathcal{T}}_i(t-1)\,.$$

Yet the site $y_s$ is also in $\mathcal{B}\big(w_i,D\setminus\widehat{\mathcal{T}}_i^-(t),t\big)$. By proposition 81.2, we have

$$\big(\mathcal{T}_i^*(t)\setminus\overline{\mathcal{T}}_i(t-1)\big)\cap\mathcal{B}\big(w_i,D\setminus\widehat{\mathcal{T}}_i^-(t),t-1\big)\subset\big(\mathcal{T}_i^*(t)\setminus\overline{\mathcal{T}}_i(t-1)\big)\cap\mathcal{A}_i(t-1)=\varnothing.$$

Thus the site $y_s$ is not in $\mathcal{B}\big(w_i,D\setminus\widehat{\mathcal{T}}_i^-(t),t-1\big)$. This implies that there are at least $t$ closed sites along the subpath $y_0,\cdots,y_{s-1}$. We consider two subcases:

⋆ $s=r$. In this case, we have $y_s=x$, and $x$ is in $\big(\mathcal{T}_i^{-0}(t)\cup\mathcal{T}_i^{+1}(t)\big)\cap\widehat{\Gamma}^*(t)$.

⋆ $s<r$. As the path $y_0,\cdots,y_{r-1}$ contains at most $t$ closed sites, the site $y_s$ must be open and it belongs to $\mathcal{T}_i^{+1}(t)$. Furthermore, the sites $y_s,y_{s+1},\dots,y_{r-1}$ are all open. The open path $y_s,\dots,y_{r-1}$ intersects $B_1(\widehat{\Gamma}^*(t),1)$, hence also the box $\Lambda(h-1,n)$. The site $y_s$ is visited at time $t\geq 1$, but all the boundary clusters intersecting $\Lambda(h-1,n)$ are visited at time 0. Therefore the whole path $y_s,\dots,y_{r-1}$ belongs to a finite cluster, and it contains strictly less than $M$ open sites, so that $r-s<M$. Since $y_r=x$ is in $\widehat{\Gamma}^*(t)$, then $y_s$ is in

$$B_1\big(x,M-1\big)\;\subset\;B_1\big(\widehat{\Gamma}^*(t),M-1\big)\;\subset\;\overset{\circ}{\widehat{\widehat{\Gamma}^-(t)}}\,. \quad (92.166)$$

The subpath $y_s,\cdots,y_r$ testifies that the site $x$ is in $B_1\big(\mathcal{T}_i^{+1}(t)\cap\overset{\circ}{\widehat{\widehat{\Gamma}^-(t)}},M-1\big)$.

In each case, the site $x$ belongs to one of the sets appearing in the right-hand side of (92.165). This concludes the proof of the inclusion stated in the lemma. □

We apply lemma 92.11 to each $i$ in $I$, and we take the union over $i$ in $I$. We obtain that

$$\bigcup_{i\in I}\mathcal{B}\big(w_i,D\setminus\widehat{\mathcal{T}}_i^-(t),t\big)\cap\widehat{\Gamma}^*(t)\;\subset\;\Big(\bigcup_{i\in I}\mathcal{B}\big(w_i,D\setminus\widehat{\mathcal{T}}_i^*(t),t\big)\cap\widehat{\Gamma}^*(t)\Big)$$
$$\cup\;\big(\mathcal{T}^{-0}(t)\cap\widehat{\Gamma}^*(t)\big)\;\cup\;B_1\big(\mathcal{T}^{+1}(t)\cap\overset{\circ}{\widehat{\widehat{\Gamma}^-(t)}},M-1\big)\,. \quad (92.167)$$

The inclusion (92.167) implies the fourth inclusion in (92.164). □

## 92.15 Control of $\left|\mathcal{T}^{+0}(t)\cap\widehat{\Gamma}^*(t)\right|$

The starting point to control $\left|\mathcal{T}^{+0}(t)\cap\widehat{\Gamma}^*(t)\right|$ is the inequality (92.163): there exist two disjoint subsets $K^+(t), L^+(t)$ of $\{1,\dots,N\}$ such that, setting

$$\begin{aligned}\widehat{\mathcal{K}}^+(t) &= \bigcup_{k\in K^+(t)} \mathcal{B}\big(w_k, D\setminus\widehat{\mathcal{T}}_k^*(t), t\big)\cap\widehat{\Gamma}^*(t)\,,\\ \widehat{\mathcal{L}}^+(t) &= \bigcup_{\ell\in L^+(t)} \mathcal{B}\big(w_\ell, D\setminus\widehat{\mathcal{T}}_\ell^*(t), t\big)\cap\widehat{\Gamma}^*(t)\,,\end{aligned}$$

we have

$$\frac{1}{8d}\left|\mathcal{T}^{+0}(t)\cap\widehat{\Gamma}^*(t)\right| \leq \left|\widehat{\mathcal{K}}^+(t)\cap\widehat{\mathcal{L}}^+(t)\right|. \tag{92.168}$$

The sets $\widehat{\mathcal{K}}^+(t)$, $\widehat{\mathcal{L}}^+(t)$ are respectively approximated by the truncated genuine explorations $\widehat{\mathcal{K}}_M^+(t)$, $\widehat{\mathcal{L}}_M^+(t)$ defined by

$$\begin{aligned}\widehat{\mathcal{K}}_M^+(t) &= \bigcup_{k\in K^+(t)} \widehat{\mathcal{B}}_M\big(\mathrm{Shell}\,(w_k,D),D,t\big)\cap\widehat{\Gamma}^*(t)\,,\\ \widehat{\mathcal{L}}_M^+(t) &= \bigcup_{\ell\in L^+(t)} \widehat{\mathcal{B}}_M\big(\mathrm{Shell}\,(w_\ell,D),D,t\big)\cap\widehat{\Gamma}^*(t)\,.\end{aligned}$$

We apply proposition 92.10 at iteration $t$. This yields

$$\widehat{\mathcal{K}}^+(t)\subset\widehat{\mathcal{K}}_M^+(t)\,,\qquad \widehat{\mathcal{K}}_M^+(t)\setminus\widehat{\mathcal{K}}^+(t)\subset\widehat{\mathrm{Error}}^*_{\widehat{\Gamma}^*}(M,t)\,, \tag{92.169}$$

$$\widehat{\mathcal{L}}^+(t)\subset\widehat{\mathcal{L}}_M^+(t)\,,\qquad \widehat{\mathcal{L}}_M^+(t)\setminus\widehat{\mathcal{L}}^+(t)\subset\widehat{\mathrm{Error}}^*_{\widehat{\Gamma}^*}(M,t)\,, \tag{92.170}$$

where the error term $\widehat{\mathrm{Error}}^*_{\widehat{\Gamma}^*}(M,t)$ is given by

$$\begin{aligned}\widehat{\mathrm{Error}}^*_{\widehat{\Gamma}^*}(M,t) &= \bigcup_{1\leq s\leq t}\Big(\overline{B}^{s-1}_{1,\Gamma}(t+1-s)\cup\widehat{B}^{s,1}_{1,\Gamma}(t+1-s)\Big)\\ &\qquad\cup\big(\mathcal{T}^{-0}(t)\cap\widehat{\Gamma}^*(t)\big)\cup B_1\big(\mathcal{T}^{+1}(t)\cap\overset{\circ}{\widehat{\widehat{\Gamma}^-(t)}},M-1\big)\\ &= \widehat{\mathrm{Error}}^-_{\widehat{\Gamma}^-}(M,t)\cup\big(\mathcal{T}^{-0}(t)\cap\widehat{\Gamma}^*(t)\big)\cup B_1\big(\mathcal{T}^{+1}(t)\cap\overset{\circ}{\widehat{\widehat{\Gamma}^-(t)}},M-1\big)\,.\end{aligned} \tag{92.171}$$

The term $\widehat{\mathrm{Error}}^-_{\widehat{\Gamma}^-}(M,t)$ has already been controlled in (92.148). The second term $\mathcal{T}^{-0}(t)\cap\widehat{\Gamma}^*(t)$ has been controlled in subsection 92.12 by the inequality (92.137). The set appearing in the third term $\mathcal{T}^{+1}(t)\cap\hat{\Gamma}^-(t)$ has also been controlled in subsection 92.12 by the inequality (92.154). Notice that some subsets of the relevant taboo sites at time $t$ do appear in the error term. The point is that these subsets have already been controlled. Once more, lemma 88.3 tells us that all the sites in the set $\widehat{\mathcal{K}}^+(t)\cap\widehat{\mathcal{L}}^+(t)$ are closed. As the genuine approximations $\widehat{\mathcal{K}}_M^+(t)$, $\widehat{\mathcal{L}}_M^+(t)$ are larger than the sets $\widehat{\mathcal{K}}^+(t)$, $\widehat{\mathcal{L}}^+(t)$, we have

$$\big\{\,x\in\widehat{\mathcal{K}}_M^+(t)\cap\widehat{\mathcal{L}}_M^+(t): x\text{ is open}\,\big\}\subset\big(\widehat{\mathcal{K}}_M^+(t)\cap\widehat{\mathcal{L}}_M^+(t)\big)\setminus\big(\widehat{\mathcal{K}}^+(t)\cap\widehat{\mathcal{L}}^+(t)\big)\,, \tag{92.172}$$

and the inequality (92.168) implies that

$$\frac{1}{8d}\left|\mathcal{T}^{+0}(t)\cap\widehat{\Gamma}^*(t)\right| \;\leq\; \Big|\big\{\, x\in\widehat{\mathcal{K}}^+_M(t)\cap\widehat{\mathcal{L}}^+_M(t) : x \text{ is closed}\,\big\}\Big|\,. \qquad (92.173)$$

It follows from (92.169) and (92.170) that

$$\big(\widehat{\mathcal{K}}^+_M(t)\cap\widehat{\mathcal{L}}^+_M(t)\big)\setminus\big(\widehat{\mathcal{K}}^+(t)\cap\widehat{\mathcal{L}}^+(t)\big) \;\subset\; \widehat{\mathrm{Error}}^*_{\widehat{\Gamma}^*}(M,t)\,. \qquad (92.174)$$

The two inclusions (92.172) and (92.174) imply that

$$\Big|\big\{\, x\in\widehat{\mathcal{K}}^+_M(t)\cap\widehat{\mathcal{L}}^+_M(t) : x \text{ is open}\,\big\}\Big| \;\leq\; \Big|\widehat{\mathrm{Error}}^*_{\widehat{\Gamma}^*}(M,t)\Big|\,. \qquad (92.175)$$

Using (92.171), we bound the cardinality of $\widehat{\mathrm{Error}}^*_{\widehat{\Gamma}^*}(M,t)$ as follows:

$$\Big|\widehat{\mathrm{Error}}^*_{\widehat{\Gamma}^*}(M,t)\Big| \;\leq\; \Big|\widehat{\mathrm{Error}}^-_{\widehat{\Gamma}^-}(M,t)\Big| + \Big|\mathcal{T}^{-0}(t)\cap\widehat{\Gamma}^*(t)\Big| \;+\; \Big|B_1\big(\mathcal{T}^{+1}(t)\cap\widehat{\overset{\circ}{\widehat{\Gamma}^-(t)}}\,,\,M-1\big)\Big|\,. \qquad (92.176)$$

Substituting (92.148), (92.153) and (92.154) in (92.176), and using the inclusion $\widehat{\Gamma}^*(t)\subset\widehat{\Gamma}^-(t)$, we obtain

$$\Big|\widehat{\mathrm{Error}}^*_{\widehat{\Gamma}^*}(M,t)\Big| \;\leq\; \Big(4d+33\frac{d^2}{p}\big(1+2d(2M)^d\big)\Big)D(n)\sum_{1\leq s\leq t}\left|\overline{\mathcal{T}}(s-1)\cap\Gamma(s-1)\right| + 16d\big(1+2d(2M)^d\big)c_H n^{d/2}\sqrt{N\ln n}\,. \qquad (92.177)$$

We finally use the same technique as for $\left|\mathcal{T}^{-0}(t)\cap\widehat{\Gamma}^-(t)\right|$, i.e., we make appeal to the event Intersect$(N)$. The occurrence of Intersect$(N)$ guarantees that

$$S\big(\widehat{\mathcal{K}}^+_M(t)\cap\widehat{\mathcal{L}}^+_M(t)\big) \;>\; -c_H n^{d/2}\sqrt{N\ln n}\,, \qquad (92.178)$$

where $c_H$ is the constant introduced in (92.50). Furthermore, we have

$$\Big|\big\{\, x\in\widehat{\mathcal{K}}^+_M(t)\cap\widehat{\mathcal{L}}^+_M(t) : x \text{ is closed}\,\big\}\Big| \;\leq\; \frac{1}{p}\Big|\big\{\, x\in\widehat{\mathcal{K}}^+_M(t)\cap\widehat{\mathcal{L}}^+_M(t) : x \text{ is open}\,\big\}\Big| - S\big(\widehat{\mathcal{K}}^+_M(t)\cap\widehat{\mathcal{L}}^+_M(t)\big)\,. \qquad (92.179)$$

Using the inequalities (92.175), (92.177), (92.178), we deduce from (92.179) that

$$\Big|\big\{\, x\in\widehat{\mathcal{K}}^+_M(t)\cap\widehat{\mathcal{L}}^+_M(t) : x \text{ is closed}\,\big\}\Big| \;\leq\; 99\frac{d^3}{p^2}(2M)^d D(n)\sum_{1\leq s\leq t}\left|\overline{\mathcal{T}}(s-1)\cap\Gamma(s-1)\right| + 48\frac{d^2}{p}(2M)^d c_H n^{d/2}\sqrt{N\ln n}\,. \qquad (92.180)$$

We substitute (92.180) in (92.173), we multiply by $8d$, and this gives

$$\left|\mathcal{T}^{+0}(t)\cap\widehat{\Gamma}^*(t)\right| \;\leq\; 792\frac{d^4}{p^2}(2M)^d D(n)\sum_{1\leq s\leq t}\left|\overline{\mathcal{T}}(s-1)\cap\Gamma(s-1)\right| + 384\frac{d^3}{p}(2M)^d c_H n^{d/2}\sqrt{N\ln n}\,. \qquad (92.181)$$

## 92.16 The long-sought inequality

For any $t$ in $\{1,\dots,T\}$, the set $\mathcal{T}^*(t)$ is decomposed as

$$\mathcal{T}^*(t) \,=\, \mathcal{T}^{-1}(t)\cup\mathcal{T}^{+1}(t)\cup\mathcal{T}^{-0}(t)\cup\mathcal{T}^{+0}(t)\,.$$

In the previous three subsections, we have obtained inequalities on the cardinalities of the sets of the right-hand side, more precisely, on their traces on different sub-boxes of the sequence $\Gamma$. Let us recall these inequalities to understand how they can be used together. To simplify them, we define

$$a(d,p,n) \,=\, 792\frac{d^4}{p^2}(2M)^dD(n)\,,\qquad b(d,p,n) \,=\, 384\frac{d^3}{p}(2M)^dc_Hn^{d/2}\sqrt{N\ln n}\,.$$

These quantities depend on $d,p,n,N$ , but we will write simply $a$ and $b$ in the sequel. We fix an integer $t$ in $\{1,\dots,T\}$. We rewrite below, in order, the inequalities (92.60), (92.153), (92.154), (92.181):

$$\big|\mathcal{T}^{-1}(t)\cap\overset{\circ}{\widehat{\Gamma(t-1)}}\big| \,\leq\, 2d\,\big|\mathcal{T}^{+0}(t-1)\cap\Gamma(t-1)\big|\,,$$

$$\big|\mathcal{T}^{-0}(t)\cap\widehat{\Gamma}^-(t)\big| \,\leq\, a\sum_{0\leq s\leq t-1}\big|\overline{\mathcal{T}}(s)\cap\Gamma(s)\big|+b\,,\tag{92.182}$$

$$\big|\mathcal{T}^{+1}(t)\cap\overset{\circ}{\widehat{\widehat{\Gamma}^-(t)}}\big| \,\leq\, 2d\,\big|\mathcal{T}^{-0}(t)\cap\widehat{\Gamma}^-(t)\big|\,,\tag{92.183}$$

$$\big|\mathcal{T}^{+0}(t)\cap\widehat{\Gamma}^*(t)\big| \,\leq\, a\sum_{0\leq s\leq t-1}\big|\overline{\mathcal{T}}(s)\cap\Gamma(s)\big|+b\,.$$

We do next two operations. We substitute the second inequality (92.182) in the third one (92.183), and we use the inclusions

$$\Gamma(t-1)\,\supset\,\overset{\circ}{\widehat{\Gamma(t-1)}}\,\supset\,\widehat{\Gamma}^-(t)\,\supset\,\overset{\circ}{\widehat{\widehat{\Gamma}^-(t)}}\,\supset\,\widehat{\Gamma}^*(t)\,\supset\,\Gamma(t)$$

to shrink the traces of the sets of the left-hand sides to $\Gamma(t)$. After doing this, we deduce from the four previous inequalities the following weaker versions:

$$\big|\mathcal{T}^{-1}(t)\cap\Gamma(t)\big| \,\leq\, 2d\,\big|\mathcal{T}^{+0}(t-1)\cap\Gamma(t-1)\big|\,,$$

$$\big|\mathcal{T}^{-0}(t)\cap\Gamma(t)\big| \,\leq\, a\sum_{0\leq s\leq t-1}\big|\overline{\mathcal{T}}(s)\cap\Gamma(s)\big|+b\,,$$

$$\big|\mathcal{T}^{+1}(t)\cap\Gamma(t)\big| \,\leq\, 2da\sum_{0\leq s\leq t-1}\big|\overline{\mathcal{T}}(s)\cap\Gamma(s)\big|+2db$$

$$\big|\mathcal{T}^{+0}(t)\cap\Gamma(t)\big| \,\leq\, a\sum_{0\leq s\leq t-1}\big|\overline{\mathcal{T}}(s)\cap\Gamma(s)\big|+b\,.$$

Adding these four inequalities, we get the following long-sought inequality:

$$\big|\mathcal{T}^*(t)\cap\Gamma(t)\big| \,\leq\, 3da\sum_{0\leq s\leq t-1}\big|\overline{\mathcal{T}}(s)\cap\Gamma(s)\big|+3db\,.$$

Adding the sum of the right-hand side to both sides, we have furthermore

$$\sum_{0\leq s\leq t}\big|\overline{\mathcal{T}}(s)\cap\Gamma(s)\big| \,\leq\, 4da\sum_{0\leq s\leq t-1}\big|\overline{\mathcal{T}}(s)\cap\Gamma(s)\big|+3db\,.$$

This is nothing more than a simple affine recurrence inequality. It demonstrates that the sequence

$$\sum_{0\leq s\leq t}\big|\overline{\mathcal{T}}(s)\cap\Gamma(s)\big|,\quad 0\leq t\leq T\,,$$

is below the arithmetic-geometric sequence $(u_t, 0\leq t\leq T)$ defined by

$$\begin{aligned} u_0 &= \big|\overline{\mathcal{T}}(0)\cap\Gamma(0)\big|\,,\\ \forall t\in\{\,1,\dots,T\,\}\qquad u_t &= 4dau_{t-1}+3db\,. \end{aligned}$$

Solving the recurrence, we find that

$$\forall t\in\{\,1,\dots,T\,\}\qquad u_t \,=\, (4da)^t u_0+3db\frac{(4da)^t-1}{4da-1}\,. \tag{92.184}$$

Using the fact that $4da-1\geq 3d$, we have

$$3db\frac{(4da)^T-1}{4da-1}\,\leq\, b(4da)^T\,. \tag{92.185}$$

Substituting (92.185) in (92.184), we conclude that

$$\sum_{0\leq s\leq T}\big|\overline{\mathcal{T}}(s)\cap\Gamma(s)\big| \,\leq\, (4da)^T\Big(\big|\overline{\mathcal{T}}(0)\cap\Gamma(0)\big|+b\Big)\,.$$

We use finally the control (92.59) on $\overline{\mathcal{T}}^0(0)\cap\Gamma(0)$, and the bound $T\leq\bar{t}(n)$ on the termination time. We replace the constants $a,b$ by their values, as well as $D(n)$, and setting

$$\nu(n)=769\frac{d^3}{p}2^d\big\lceil(4d\ln n)^{\frac{d}{\alpha}+\frac{1}{2}}\big\rceil c_H\bigg(3168\frac{d^5}{p^2}\Big(16(\ln n)^{1-\gamma}\big\lceil(4d\ln n)^{2/\alpha}\big\rceil\Big)^d\bigg)^{2(\ln n)^{1-\gamma}}$$

we obtain

$$\sum_{0\leq s\leq T}\big|\overline{\mathcal{T}}(s)\cap\Gamma(s)\big| \,\leq\, \nu(n)\,n^{d/2}\sqrt{N}\,. \tag{92.186}$$

The prefactor $\nu(n)$ is sub-polynomial, in fact we have

$$\begin{aligned} \nu(n) \,&\leq\, \exp\Big(2(\ln n)^{1-\gamma}\big(d+\frac{2d}{\alpha}\big)\ln\ln n+o\big((\ln n)^{1-\gamma}\ln\ln n\big)\Big)\\ &\leq\, \exp\Big(o\big((\ln n)^{1-\gamma/2}\big)\Big)\,. \end{aligned} \tag{92.187}$$

As we see, the role of the exponent $\gamma>0$ is absolutely crucial!

## 92.17 The exponent 2

We finally gather the results obtained so far. In fact, from the subsection 92.3 until the subsection 92.16, we have been analyzing what happens in a configuration of the event $\mathcal{H}'$ defined in (92.27), that we recall for convenience:

$$\mathcal{H}' = \mathcal{E} \cap \Big\{ H = h, \mathcal{N}^*(h+1,n) = N, \mathcal{C}^*(h+1,n) = \big\{ C\big(w_i, \Lambda(4n)\big), 1 \leq i \leq N \big\} \Big\}.$$

What we did is to launch the multiple intertwined explorations starting from $w_1, \dots, w_N$, with goal set $G$ the box $\Lambda(h,n)$. In subsection 92.3, inequality (92.33), we showed that, for any configuration in $\mathcal{H}'$, we have the inequality:

$$\big|\overline{\mathcal{T}}^0(T) \cap \Lambda(h-1,n)\big| \,\geq\, \frac{Nn}{256 d^3 \ln n}\,. \tag{92.188}$$

Afterwards, we have been developing inequalities on the relevant taboo sets of these intertwined explorations, more specifically on their traces in the sequence of nested sub-boxes $\Gamma = \big(\Gamma(t), 0 \leq t \leq T\big)$. We introduced a further event in subsection 92.6, named $\text{Intersect}(N)$, and we showed in inequality (92.186) of subsection 92.16 that the configurations belonging to $\mathcal{H}' \cap \text{Intersect}(N)$ satisfy

$$\sum_{0 \leq s \leq T} \big|\overline{\mathcal{T}}(s) \cap \Gamma(s)\big| \,\leq\, \nu(n)\, n^{d/2} \sqrt{N}\,, \tag{92.189}$$

where the prefactor $\nu(n)$, given in (92.187), is sub-polynomial. Since the boxes $\Gamma = \big(\Gamma(t), 0 \leq t \leq T\big)$ all contain $\Lambda(h-1,n)$, the inequality (92.189) implies that

$$\big|\overline{\mathcal{T}}^0(T) \cap \Lambda(h-1,n)\big| \,\leq\, \nu(n)\, n^{d/2} \sqrt{N}\,.$$

For the values of $N$ such that

$$\nu(n)\, n^{d/2} \sqrt{N} \,<\, \frac{Nn}{256 d^3 \ln n}\,, \tag{92.190}$$

the two inequalities (92.188) and (92.189) are not compatible, therefore there are no configurations in $\mathcal{H}' \cap \text{Intersect}(N)$. The magic point is that $N$ appears on both sides of (92.190), albeit with a different power! We divide the inequality (92.190) by $\sqrt{N}$, we take the squares and we rewrite it as

$$N \,>\, \big(256 d^3 \nu(n) \ln n\big)^2\, n^{d-2}\,. \tag{92.191}$$

From the previous discussion, we conclude that, for any integer $N$ satisfying (92.191), the event $\mathcal{H}'$ is included in the complement of $\text{Intersect}(N)$, so that, by (92.55),

$$P(\mathcal{H}') \,\leq\, P\big(\text{Intersect}(N)^c\big) \,\leq\, \exp\big(-4^d N \ln n\big)\,. \tag{92.192}$$

Recall that the event $\mathcal{H}'$, introduced in (92.27), depends on $h, N$, and the vertices $w_1, \dots, w_N$. We choose now

$$N^{**} \,=\, \big(257 d^3 \nu(n) \ln n\big)^2\, n^{d-2}\,, \tag{92.193}$$

and we derive an upper bound on the probability $P\big(\mathcal{N}^*(H+1,n)\geq N^{**},\mathcal{E}\big)$. As all the integers $N$ intervening in the sum (92.26) are larger than or equal to $N^{**}$, they satisfy the inequality (92.191), and we can use the uniform upper bound (92.192). So, substituting (92.192) in the sum (92.28), we get that, for $N\geq N^{**}$,

$$P\begin{pmatrix}\mathcal{E},H=h,\\ \mathcal{N}^*(h+1,n)=N\end{pmatrix}\ \leq\sum_{w_1,\dots,w_N\in\partial^{\,in}\Lambda(4n)}\exp\big(-4^dN\ln n\big)$$
$$\leq\ \exp\Big(N\ln\big(|\Lambda(4n)|\big)-4^dN\ln n\Big)\ \leq\ \exp\big(-dN\ln n\big)\,. \quad (92.194)$$

We substitute next (92.194) into the previous sum (92.26), and we obtain

$$P\big(\mathcal{N}^*(H+1,n)\geq N^{**},\mathcal{E}\big)\ \leq\sum_{0\leq h<\overline{h}(n)}\ \sum_{N^{**}\leq N\leq|\partial^{\,in}\Lambda(4n)|}\exp\big(-dN\ln n\big)$$
$$\leq\ \overline{h}(n)\big|\partial^{\,in}\Lambda(4n)\big|\exp\big(-dN^{**}\ln n\big)$$
$$\leq\ n(4n+1)^d\exp\big(-257^2d^9\ln n\big)\,. \quad (92.195)$$

Let us denote by $\mathcal{N}^*$ the number of boundary clusters in $\Lambda(4n)$ which intersect the box $\Lambda(2n)$, i.e.,

$$\mathcal{N}^*\ =\ \big|\big\{\,C\in\mathcal{C}_{\mathrm{bd}}(4n):C\cap\Lambda(2n)\neq\varnothing\,\big\}\big|\,.$$

Any boundary cluster intersecting $\Lambda(2n)$ intersects also $\Lambda(H+1,n)$, therefore

$$\mathcal{N}^*(H+1,n)\ \geq\ \mathcal{N}^*\,. \quad (92.196)$$

Putting together (92.195) and (92.196), and replacing $N^{**}$ by its value (92.193), we conclude that

$$\lim_{n\to\infty}P\Big(\mathcal{N}^*\geq\big(257d^3\nu(n)\ln n\big)^2n^{d-2},\mathcal{E}\Big)\ =\ 0\,.$$

This implies that, when $\mathcal{E}$ occurs, with probability going to 1 as $n$ goes to $\infty$,

$$\mathcal{N}^*\ \leq\ n^{d-2+o(1)}\,. \quad (92.197)$$

On the event $\mathcal{E}$, the event $\mathcal{F}_{\mathrm{bd}}(n,\varepsilon)$ defined in (92.2) occurs as well, whence

$$\sum_{C\in\mathcal{C}_{\mathrm{bd}}(4n)}\big|C\cap\Lambda(2n)\big|\ \geq\ \big(\theta(p)-\varepsilon\big)\big|\Lambda(2n)\big|\,. \quad (92.198)$$

The inequality (92.197) says that the number of terms in the sum on the left-hand side of (92.198) is at most $n^{d-2+o(1)}$. From there, we deduce that, with probability going to 1 as $n$ goes to $\infty$, the box $\Lambda(4n)$ contains a boundary cluster $C$ such that

$$\big|C\cap\Lambda(2n)\big|\ \geq\ \big(\theta(p)-\varepsilon\big)\big|\Lambda(2n)\big|n^{2-d+o(1)}\ =\ n^{2+o(1)}\,.$$

Thanks to the construction carried out in this whole section, we have improved the exponent 3/2 obtained in (91.11) to the exponent 2. We have not reached the ultimate goal, but this is still a big step forward!

# 93 Isoperimetry comes into play

We already discussed extensively the discrete isoperimetric inequality in section 11. It is very likely that this inequality will play a major role in the resolution of the conjecture $\theta(p_c, \mathbb{Z}^d) = 0$. The time has now come for the isoperimetric inequality to enter the scene. However, as explained in section 11, the isoperimetric inequality is a very rich topic in itself, and there exist various inequalities involving different notions of boundary. We will state in the next subsection the one inequality that is to be used in the proof of theorem 1.2. We put it in action in subsection 93.2, but only in $(d-1)$-dimensional boxes. This leads to a slight improvement of the exponent 2 to $2+\varepsilon$. This is a necessary step before carrying the ultimate parts of the proof, in subsection 93.3, with the help of the full inequality in dimension $d$.

## 93.1 The relative isoperimetric inequality

As our whole story takes place in the box $\Lambda(4n)$, we certainly need a statement dealing with subsets of $\Lambda(4n)$. Béla Bollobás and Imre Leader have obtained an essentially sharp edge-isoperimetric inequality in this context. We restate next their result with our notation. For $\Lambda$ a cubic box in $\mathbb{Z}^d$ and $A$ a subset of $\mathbb{Z}^d$, we define the edge boundary $\partial_\Lambda^{edge} A$ of $A$ in $\Lambda$ as

$$\partial_\Lambda^{edge} A \,=\, \big\{\, e = \langle x, y\rangle \in \mathbb{E}^d : x \in A \cap \Lambda,\, y \in \Lambda \setminus A \,\big\}\,.$$

**Theorem 93.1** (Theorem 3 of [18])**.** *Let $d \geq 1$, $m \geq 1$ and let $A$ be a subset of $\Lambda(m)$ such that $1 \leq |A| \leq \big|\Lambda(m)\big|/2$. Then*

$$\Big|\partial_{\Lambda(m)}^{edge} A\Big| \,\geq\, \min\Big\{\, r|A|^{1-1/r}\big|\Lambda(m)\big|^{1/r-1/d} : 1 \leq r \leq d \,\Big\}\,. \tag{93.1}$$

Let us explain precisely the correspondence between our notation and [18]. Let $m \geq 1$ and let us set $\ell = \lfloor m/2 \rfloor$. With our definitions, we have

$$\Lambda(m) \cap \mathbb{Z}^d \,=\, \{\, -\ell, \dots, \ell \,\}^d\,, \qquad \big|\Lambda(m) \cap \mathbb{Z}^d\big| \,=\, (2\ell+1)^d\,.$$

The set above is isometric to $\{\, 0, \dots, 2\ell \,\}^d$, and theorem 93.1 corresponds to the specific case of theorem 3 of [18] for the values $k = 2\ell + 1$ and $n = d$ (notice that in theorem 3 of [18], the set $[k]^n$ is $\{\, 0, \dots, k-1 \,\}^n$ and $k^n = \big|[k]^n\big|$).

The inequality (93.1) is sharp in the sense that the equality can be realized. Béla Bollobás and Imre Leader [18] even provide examples of sets for which there is equality. The inequality (93.1) is in fact stronger than what we will need really. We will essentially rely on the discrete relative isoperimetric inequality, that we state next.

**Theorem 93.2** (discrete relative isoperimetric inequality)**.** *Let $d \geq 2$ be fixed. There exists a constant $c_{iso}$ depending only upon the dimension $d$ such that, for any subset $A$ of $\mathbb{Z}^d$, for any cubic box $\Lambda$ whose sides are parallel to the axes,*

$$\min\big(|A \cap \Lambda|, |\Lambda \setminus A|\big) \,\leq\, c_{iso}\,\big|\partial_\Lambda^{edge} A\big|^{\frac{d}{d-1}}\,. \tag{93.2}$$

*In addition, the constant $c_{iso}$ is less than or equal to $1/2$.*

This is certainly a classical well-known result. However, we did not find the exact statement above in the literature. So, we show next how it follows from theorem 93.1.

*Proof.* Let $A$ be a subset of $\mathbb{Z}^d$ and let us set

$$B\,=\,A\cap\Lambda\,,\qquad C\,=\,\Lambda\setminus A\,.$$

We suppose that neither $B$ nor $C$ is empty, otherwise there is nothing to prove. As $|B|+|C|=|\Lambda|$, we have either $|B|\leq|\Lambda|/2\leq|C|$ or $|C|\leq|\Lambda|/2\leq|B|$. Suppose for instance that $|B|\leq|\Lambda|/2$. We apply theorem 93.1 to the box $\Lambda$ and the set $B$, and we bound from below the term in the minimum (93.1) as follows: for any $r$ in $\{1,\dots,d\}$, we have $r\geq 2^{1-1/r}$, whence

$$r|B|^{1-1/r}\big|\Lambda\big|^{1/r-1/d}\,\geq\,\frac{|B|}{\big|\Lambda\big|^{1/d}}\;2\bigg(\frac{|\Lambda|}{2|B|}\bigg)^{1/r}\,. \tag{93.3}$$

Since $\big|\Lambda\big|\geq 2|B|$, then

$$\forall r\in\{1,\dots,d\}\qquad\bigg(\frac{|\Lambda|}{2|B|}\bigg)^{1/r}\,\geq\,\bigg(\frac{|\Lambda|}{2|B|}\bigg)^{1/d}\,. \tag{93.4}$$

We substitute (93.4) in (93.3), and we obtain that

$$r|B|^{1-1/r}\big|\Lambda\big|^{1/r-1/d}\,\geq\,\big(2|B|\big)^{1-1/d}\,. \tag{93.5}$$

Taking the minimum over $r$ in $\{1,\dots,d\}$ in (93.5), we conclude from (93.1) that

$$\big|\partial^{edge}_{\Lambda}B\big|\,\geq\,\big(2|B|\big)^{1-1/d}\,=\,\Big(2\min\big(|A\cap\Lambda|,|\Lambda\setminus A|\big)\Big)^{1-1/d}\,.$$

We finally observe that $\big|\partial^{edge}_{\Lambda}B\big|=\big|\partial^{edge}_{\Lambda}A\big|$, and we obtain the desired inequality (93.2) with a constant $c_{\text{iso}}=1/2$. In the case where $|C|\leq|B|$, we apply the previous argument to $C$ instead of $B$. □

We call the inequality (93.2) the discrete relative isoperimetric inequality in reference to its continuous counterpart, which is more classical and well-known. However, it is more technical as it requires a proper definition of the perimeter in the continuous setting. We just give next the minimal necessary definitions for stating and discussing the relative isoperimetric inequality in a ball. For a modern account of the theory of Caccioppoli sets, see one of the books [9, 59, 139] (this is absolutely not necessary to read the rest of this section). The perimeter of a Borel set $E$ of $\mathbb{R}^d$ in an open set $O$ is defined as

$$\mathcal{P}(E,O)\,=\,\sup\,\Big\{\,\int_E\operatorname{div}f(x)\,dx:f\in C^\infty_c(O,B(0,1))\,\Big\}\,,$$

where $C^\infty_c(O,B(0,1))$ is the set of the $C^\infty$ vector functions from $\mathbb{R}^d$ to the unit Euclidean ball $B(0,1)$ having a compact support included in $O$ and div is the usual divergence operator, defined for a $C^1$ vector function $f$ with scalar components $(f_1, \dots, f_d)$ as

$$\operatorname{div} f \,=\, \frac{\partial f_1}{\partial x_1} + \cdots + \frac{\partial f_d}{\partial x_d}\,.$$

The set $E$ is said to have finite perimeter in an open set $O$ if $\mathcal{P}(E, O)$ is finite. The set $E$ is said to be of locally finite perimeter or to be a Caccioppoli set if $\mathcal{P}(E, O)$ is finite for every bounded open set $O$ of $\mathbb{R}^d$. We are now ready to state the result we would like to discuss.

**Theorem 93.3** (relative isoperimetric inequality in a ball)**.** *Let $d \geq 2$. There exists a constant $c_{iso}$ depending only upon the dimension $d$ such that, for any Caccioppoli set $E$, any Euclidean ball $B$ of $\mathbb{R}^d$,*

$$\min\big(\mathcal{L}^d(E \cap B), \mathcal{L}^d(B \setminus E)\big) \,\leq\, c_{iso}\, \mathcal{P}\big(E, \overset{o}{B}\big)^{d/(d-1)}\,.$$

This result can be found for instance in [9], formula 3.37 in remark 3.45, or in [139], theorem 5.4.3. It seems that theorem 93.3 belongs to another mathematical world, much more technical. So, why did we bother presenting it here?

First, it is interesting to note that theorem 93.2 is a discrete counterpart of theorem 93.3. Second, and even more interestingly, theorem 93.2 can be rigorously deduced from such a continuous statement. The technique is standard. To a discrete subset $A$ of $\mathbb{Z}^d$, we associate a continuous set $\widetilde{A}$ by putting a unit cube centered at each vertex of $A$, and we apply the continuous isoperimetric inequality to $\widetilde{A}$. Naturally, the volume of $\widetilde{A}$ is equal to the cardinality of $A$, and the perimeter of $\widetilde{A}$ is equal to the cardinality of its edge boundary. There is a remaining difficulty, because theorem 93.3 deals with balls, and we would like to have a cubic box $\Lambda$ instead of an Euclidean ball $B$. This is absolutely possible, the relative isoperimetric inequality holds in domains more general than balls, for instance in bounded Lipschitz domains (see for instance chapter 3, exercise 3.13 of [9], or theorem 5.11.1 of [139]). However, the proofs of this kind of results involve a lot of difficult techniques pertaining to continuous spaces, so it seems that the most reasonable approach is to use the isoperimetric theorem 93.1.

We discuss finally the discrete isoperimetric inequality of Deuschel and Pisztora, stated in proposition 2.2 of [37]. This inequality is well-suited to handle several percolation questions, as it was developed along the coarse-graining estimates of Pisztora [114]. It looks like that the inequality in proposition 2.2 of [37] is weaker than the isoperimetric inequality (93.1). In fact, it is obtained by applying the weak discrete isoperimetric inequality (11.2) to each connected component of the complement of the set $A$ (hence the sum over the connected components in formula 2.2 of [37]). This weaker inequality would certainly suffice for the next subsection, where we improve the exponent 2 to $2+\varepsilon$. However, for the subsequent final section, we absolutely need a discrete isoperimetric inequality inside a cube of $\mathbb{Z}^d$ which furnishes the optimal exponent $1 - 1/d$.

## 93.2 Improving the exponent from 2 to $2+\alpha\big(1-1/(d-1)\big)$

In this subsection, we shall revisit the argument of section 92 leading to the exponent 2, and with the help of the relative isoperimetric inequality, we shall strengthen some of the estimates. This will lead to an improvement in the exponent, from 2 to $2+\varepsilon$.

The key improvement will be for the inequality (92.30). In order to improve this inequality, we will show that, typically, the boundary clusters are thick far from their extremities, in the sense that their traces on the boundaries of boxes contain a polynomial number of sites. The idea is quite simple, let us try to explain it before delving into the mathematical details. Let $m,\ell$ be two integers such that $m<\ell<n$. Consider a boundary cluster $C$ in $\Lambda(4n)$ that intersects the box $\Lambda(m)$. Suppose that the trace of $C$ on the boundary of the box $\Lambda(\ell)$ contains $R$ sites. If we close by force these $R$ sites, then the cluster $C$ will split up into several clusters. Among these clusters, those which intersect $\Lambda(m)$ are confined to the box $\Lambda(\ell)$, thus they become finite clusters of the percolation configuration in $\Lambda(4n)$. Furthermore, these clusters must also contain an open path from $\Lambda(m)$ to the boundary of $\Lambda(\ell)$, thus they are finite clusters of cardinality of order $(\ell-m)/2$. Yet, thanks to the hypothesis (9.16) and the tail estimate (92.4), the probability of observing such a large finite cluster in $\Lambda(4n)$ is smaller than $\big|\Lambda(4n)\big|\exp\big(-(\ell-m)^\alpha/2^\alpha\big)$. Now, there is also a combinatorial factor that must be taken into account, which corresponds to the closing of the $R$ sites. This factor is of order $\big|\Lambda(4n)\big|^R$, and we see that the previous scenario is very unlikely as long as $R$ is much smaller than $(\ell-m)^\alpha$. We can ensure that this is the case over an adequately chosen range of parameters.

We proceed now to the mathematical implementation of the idea explained above. We will essentially do a modification of the argument leading to the exponent 2, and we reuse all the notation introduced throughout section 92. As we approach the end of the proof, we are not trying to obtain the best possible estimate, but only the one that is truly useful for drawing a conclusion. We set

$$\forall n\geq 2\qquad \psi(n)\,=\,\frac{1}{4}\left\lceil\frac{n\ln 2}{d\ln n}\right\rceil.$$

We will use the inequalities

$$\forall n\geq 2\qquad \frac{n}{8d\ln n}\,\leq\,\psi(n)\,\leq\,\frac{1}{4}\big(\lambda(n)+1\big)\,. \tag{93.6}$$

We introduce the event Thick$(n)$ defined as

$$\text{Thick}(n)\,=\,\Big\{\forall C\in\mathcal{C}_{\text{bd}}(4n)\quad\forall m\in\{\,1,\dots,4n\,\}$$
$$C\cap\Lambda(m)\neq\varnothing\quad\Longrightarrow\quad\forall\ell\in\big\{\,m+\psi(n),\dots,4n\,\big\}\qquad\big|C\cap\partial^{\,in}\Lambda(\ell)\big|\geq n^{\alpha/2}\Big\}.$$

**Proposition 93.4.** *There exists an integer $n_5(\alpha,d,p)$ depending only on $d,\alpha,p$ such that*

$$\forall n\geq n_5(\alpha,d,p)\qquad P\big(\mathit{Thick}(n)^c\big)\,\leq\,\exp\left(-\Big(\frac{n}{34d\ln n}\Big)^\alpha\right). \tag{93.7}$$

*Proof.* We estimate the probability of the complementary event:

$$
\begin{aligned}
P(\text{Thick}(n)^c) &= \\
&P\Big(\exists C \in \mathcal{C}_{\text{bd}}(4n) \quad \exists m \in \{\,1,\dots,4n\,\} \quad C \cap \Lambda(m) \neq \varnothing, \\
&\qquad \exists \ell \in \{\,m+\psi(n),\dots,4n\,\} \qquad \big|C \cap \partial^{in}\Lambda(\ell)\big| < n^{\alpha/2}\Big) \\
&\leq \sum_{\substack{1\leq m\leq 4n\\ m+\psi(n)\leq \ell\leq 4n}} P\begin{pmatrix}\exists C \in \mathcal{C}_{\text{bd}}(4n) \quad C\cap\Lambda(m)\neq\varnothing,\\ \big|C\cap\partial^{in}\Lambda(\ell)\big| < n^{\alpha/2}\end{pmatrix} \\
&\leq \sum_{\substack{1\leq m\leq 4n\\ m+\psi(n)\leq \ell\leq 4n}} \sum_{x\in\Lambda(m)} P\begin{pmatrix} x\longleftrightarrow \partial^{in}\Lambda(4n),\\ \big|C(x)\cap\partial^{in}\Lambda(\ell)\big| < n^{\alpha/2}\end{pmatrix}. \qquad (93.8)
\end{aligned}
$$

We fix next two values of $m, \ell$ such that $1 \leq m \leq 4n$, $m+\psi(n) \leq \ell \leq 4n$, as well as a site $x$ of $\Lambda(m)$. We take the probability inside the sum and we try to estimate it. We denote by $R$ the cardinality of $C(x) \cap \partial^{in}\Lambda(\ell)$. We start by conditioning on the value of $R$ and the set $C(x)\cap\partial^{in}\Lambda(\ell)$:

$$
P\begin{pmatrix} x\longleftrightarrow \partial^{in}\Lambda(4n),\\ \big|C(x)\cap\partial^{in}\Lambda(\ell)\big| < n^{\alpha/2}\end{pmatrix} \leq \sum_{1\leq r<n^{\alpha/2}} \sum_{\substack{X\subset\partial^{in}\Lambda(\ell)\\ |X|=r}} P\begin{pmatrix} x\longleftrightarrow \partial^{in}\Lambda(4n),\ R=r,\\ C(x)\cap\partial^{in}\Lambda(\ell) = X\end{pmatrix}. \qquad (93.9)
$$

We fix $r$ such that $1 \leq r < n^{\alpha/2}$ and a subset $X$ of $\partial^{in}\Lambda(\ell)$ such that $|X| = r$. We denote by $\mathcal{H}(x,\ell,r,X)$ the event

$$
\mathcal{H} \,=\, \mathcal{H}(x,\ell,r,X) \,=\, \begin{Bmatrix} x\longleftrightarrow \partial^{in}\Lambda(4n),\ R=r,\\ C(x)\cap\partial^{in}\Lambda(\ell) = X\end{Bmatrix}.
$$

Let $\phi : \mathcal{H} \to \{0,1\}^{\Lambda(4n)}$ be the map defined by

$$
\forall \omega\in\mathcal{H} \quad \forall y \in \Lambda(4n) \qquad \phi(\omega)(y) \,=\, \begin{cases} \omega(y) & \text{if } y\notin X\,,\\ 0 & \text{if } y \in X\,.\end{cases}
$$

The map $\phi$ closes by force the sites of $X$. Notice that, by definition, these sites were open in $\omega$. The map $\phi$ leaves the other sites unchanged. We have thus

$$
\forall \omega \in \mathcal{H} \qquad P(\omega) \,=\, \Big(\frac{p}{1-p}\Big)^r P\big(\phi(\omega)\big)\,. \qquad (93.10)
$$

In addition, after the closure of the sites of $X$, the cluster of $x$ does not intersect any more the boundary $\partial^{in}\Lambda(\ell)$, hence it is confined to $\Lambda(\ell)$. We know also that, before the closure of the sites, this cluster was a boundary cluster of $\Lambda(4n)$. Necessarily, the site $x$ is connected in $\Lambda(\ell)\setminus\partial^{in}\Lambda(\ell)$ to a site $y$ which has itself a neighbour in $\partial^{\,in}\Lambda(\ell)$. This connection is not affected by the closure operation, as it is realized without using any site in $\partial^{\,in}\Lambda(\ell)$. With formulas, we have $\mathcal{H}(x,\ell,r,X)\subset\mathcal{F}(x,\ell)$, where $\mathcal{F}(x,\ell)$ is the event

$$\mathcal{F}\,=\,\mathcal{F}(x,\ell)\,=\,\Big\{\,x\longleftrightarrow\partial^{\,in}\big(\Lambda(\ell)\setminus\partial^{\,in}\Lambda(\ell)\big)\text{ in }\Lambda(\ell)\setminus\partial^{\,in}\Lambda(\ell)\,\Big\}\,.$$

As the event $\mathcal{F}(x,\ell)$ is measurable with respect to the sites in $\Lambda(\ell)\setminus\partial^{\,in}\Lambda(\ell)$, the image $\phi(\mathcal{H})$ of $\mathcal{H}$ through $\phi$ is still included in $\mathcal{F}$, so that

$$\phi(\mathcal{H})\,\subset\,\mathcal{F}\cap\big\{\,x\not\longleftrightarrow\partial^{\,in}\Lambda(4n)\,\big\}\,.\tag{93.11}$$

Summing (93.10) over $\omega$ in $\mathcal{H}$, we obtain

$$P(\mathcal{H})\,=\,\Big(\frac{p}{1-p}\Big)^r\sum_{\omega\in\mathcal{H}}P\big(\phi(\omega)\big)\,.\tag{93.12}$$

The map $\phi$ is one to one, thus

$$\sum_{\omega\in\mathcal{H}}P\big(\phi(\omega)\big)\,=\,\sum_{\omega\in\phi(\mathcal{H})}P(\omega)\,=\,P\big(\phi(\mathcal{H})\big)\,.\tag{93.13}$$

The two equalities (93.12), (93.13) and the inclusion (93.11) yield that

$$P(\mathcal{H})\,=\,\Big(\frac{p}{1-p}\Big)^rP\big(\phi(\mathcal{H})\big)\,\leq\,\Big(\frac{1}{1-p}\Big)^{n^{\alpha/2}}P\Big(\mathcal{F}\cap\big\{\,x\not\longleftrightarrow\partial^{\,in}\Lambda(4n)\,\big\}\Big)\,.\tag{93.14}$$

We estimate next this last probability. On the event $\mathcal{F}\cap\big\{\,x\not\longleftrightarrow\partial^{\,in}\Lambda(4n)\,\big\}$, the cluster of $x$ is a finite cluster of $\Lambda(4n)$ having cardinality larger than or equal to

$$\frac{1}{2}(\ell-m-2)\,\geq\,\frac{1}{2}\psi(n)-1\,\geq\,\frac{1}{4}\psi(n)$$

(we take $n$ large enough to ensure that $\psi(n)\geq 4$). It follows from the tail estimate (92.4) that

$$\begin{aligned}P\Big(\mathcal{F}\cap\big\{\,x\not\longleftrightarrow\partial^{\,in}\Lambda(4n)\,\big\}\Big)&\\ &\leq\,P\Big(\frac{1}{4}\psi(n)\leq|C(0)|<+\infty\Big)\,\leq\,\exp\Big(-\frac{1}{4^\alpha}\psi(n)^\alpha\Big)\,.\end{aligned}$$

We substitute this estimate in (93.14), and then in the sums in (93.9), and we get

$$\begin{aligned}P\left(\begin{matrix}x\longleftrightarrow\partial^{\,in}\Lambda(4n),\\ \big|C(x)\cap\partial^{\,in}\Lambda(\ell)\big|<n^{\alpha/2}\end{matrix}\right)\,&\leq\,\sum_{\substack{1\leq r<n^{\alpha/2}\\ X\subset\partial^{\,in}\Lambda(\ell),|X|=r}}\Big(\frac{1}{1-p}\Big)^{n^{\alpha/2}}\exp\Big(-\frac{1}{4^\alpha}\psi(n)^\alpha\Big)\\ &\leq\,n^{\alpha/2}\big|\partial^{\,in}\Lambda(\ell)\big|^{n^{\alpha/2}}\Big(\frac{1}{1-p}\Big)^{n^{\alpha/2}}\exp\Big(-\frac{1}{4^\alpha}\psi(n)^\alpha\Big)\\ &\leq\,\exp\left(\frac{\alpha}{2}\ln n+dn^{\alpha/2}\ln(4n+1)-n^{\alpha/2}\ln(1-p)-\frac{1}{4^\alpha}\Big(\frac{n}{8d\ln n}\Big)^\alpha\right)\,.\end{aligned}\tag{93.15}$$

We have used the inequality (93.6) in the last step. The dominant term inside the exponential in (93.15) is the last one, so that

$$\exists\,\widetilde{n}(\alpha,d,p)\geq 1\quad\forall n\geq\widetilde{n}(\alpha,d,p)$$
$$\frac{\alpha}{2}\ln n+dn^{\alpha/2}\ln(4n+1)-n^{\alpha/2}\ln(1-p)-\Big(\frac{n}{32d\ln n}\Big)^{\alpha}\,\leq\,-\Big(\frac{n}{33d\ln n}\Big)^{\alpha}\,. \tag{93.16}$$

Let us fix next an integer $n$ larger than $\widetilde{n}(\alpha,d,p)$. Using (93.16), we conclude from (93.15) that

$$P\left(\begin{matrix}x\longleftrightarrow\partial^{\,in}\Lambda(4n),\\ \big|C(x)\cap\partial^{\,in}\Lambda(\ell)\big|<n^{\alpha/2}\end{matrix}\right)\,\leq\,\exp\left(-\Big(\frac{n}{33d\ln n}\Big)^{\alpha}\right). \tag{93.17}$$

This estimate is uniform with respect to $\ell,m$ and $x$ in $\Lambda(m)$. We can thus substitute (93.17) in the initial sums (93.8), and we get

$$\begin{aligned}P\big(\mathrm{Thick}(n)^c\big)\,&\leq\,\sum_{\substack{1\leq m\leq 4n\\ m+\psi(n)\leq\ell\leq 4n}}\ \sum_{x\in\Lambda(m)}\exp\left(-\Big(\frac{n}{33d\ln n}\Big)^{\alpha}\right)\\ &\leq\,16n^2\big|\Lambda(4n)\big|\exp\left(-\Big(\frac{n}{33d\ln n}\Big)^{\alpha}\right).\end{aligned} \tag{93.18}$$

We can find an integer $n_5(\alpha,d,p)\,\geq\,\widetilde{n}(\alpha,d,p)$ such that the last quantity in (93.18) is smaller than the upper bound in (93.7). ☐

Proposition 93.4 shows that the event $\mathrm{Thick}(n)$ is typical. We restart now the argument of section 92, and we introduce the following changes. We add the event $\mathrm{Thick}(n)$ in the event $\mathcal{E}$, which becomes

$$\mathcal{E}\,=\,\mathcal{F}_{\mathrm{fi}}(n,M)\cap\mathcal{G}(n)\cap\mathcal{F}_{\mathrm{bd}}(n,\varepsilon)\cap\mathrm{Thick}(n)\,.$$

Thanks to proposition 93.4, we still have

$$\lim_{n\to\infty}P\big(\mathcal{E}\big)\,=\,1\,. \tag{93.19}$$

We proceed then exactly as in subsection 92.2, we define the boxes $\Lambda(h,n)$, and with the help of lemma 92.1, we show that there exists a random index $H$ such that

$$\mathcal{N}^*(H-2,n)\,\geq\,\frac{1}{8}\mathcal{N}^*(H+1,n)\,.$$

We condition then on the value of $H,\mathcal{N}^*(H+1,n)$, i.e., we suppose that $H=h$ and $\mathcal{N}^*(H+1,n)=N$. We arrive at subsection 92.3, and we try to improve the inequality between $\mathcal{N}^*(h+1,n)$ and $\overline{\mathcal{T}}^0(T)\cap\Lambda(h-1,n)$. Instead of looking at the traces of the aggregates on hyperplanes, we will look at their traces on the boundaries of intermediate boxes. The idea is the following. Thanks to the event Thick($n$), the trace of each boundary cluster of $\mathcal{C}^*(h-2,n)$ on these boundaries contains a polynomial number of sites. The same is true for the aggregates. The relative isoperimetric inequality for $(d-1)$-dimensional boxes implies that, around each aggregate, there is a polynomial number of relevant taboo sites in the boundary of each intermediate box. We will obtain a uniform lower bound on this polynomial number, which will result in an improvement of the exponent of $n$ in the final inequality.

We explain now how to implement rigorously the previous idea, and we rewrite the corresponding passages of the subsection 92.3. We consider the final aggregates $(\mathcal{A}_i(T), 1 \leq i \leq N)$ produced by the intertwined explorations, and we recall briefly the essential facts:

- They are pairwise disjoint and they cover the box $\Lambda(h,n)$.
- Each aggregate contains exactly one boundary cluster of $\mathcal{C}^*(h+1,n)$.
- There are at least $N/8$ boundary clusters that intersect $\Lambda(h-2,n)$.

Let us consider a boundary cluster $C$ in $\Lambda(4n)$ that intersects the box $\Lambda(h-2,n)$. As the event Thick($n$) occurs, we have

$$\forall \ell \in \big\{\, n+(h-2)\lambda(n)+\psi(n),\dots,4n \,\big\} \qquad \big|C\cap \partial^{\,in}\Lambda(\ell)\big| \geq n^{\alpha/2}\,. \quad (93.20)$$

Let us define

$$\overline{\mathcal{S}}(h,n) \,=\, \{\, n+(h-2)\lambda(n)+\psi(n),\dots,n+(h-2)\lambda(n)+2\psi(n)\,\}\,, \quad (93.21)$$

and let us fix $s$ in $\overline{\mathcal{S}}(h,n)$. It follows from (93.20) that

$$\big|C\cap \partial^{\,in}\Lambda(s)\big| \,\geq\, n^{\alpha/2}\,.$$

If an aggregate $\mathcal{A}_i(T)$ contains the boundary cluster $C$, then it will also satisfy

$$\big|\mathcal{A}_i(T)\cap \partial^{\,in}\Lambda(s)\big| \,\geq\, n^{\alpha/2}\,.$$

In addition, we know that there are at least $N/8$ such aggregates. Thus the traces of the aggregates $(\mathcal{A}_i(T), 1 \leq i \leq N)$ on $\partial^{\,in}\Lambda(s)$ form a partition of $\partial^{\,in}\Lambda(s)$ containing at least $N/8$ sets having cardinality larger than or equal to $n^{\alpha/2}$. We shall next use the following lemma to bound from below the cardinality of the intersection sets. For $\pi$ a subset of $\partial^{\,in}\Lambda$, the set $\Delta_{\partial^{\,in}\Lambda}\pi$ is defined as

$$\Delta_{\partial^{\,in}\Lambda}\pi \,=\, \big\{\, e=\langle x,y\rangle \in \mathbb{E}^d : x\in\pi, y\in \partial^{\,in}\Lambda\setminus\pi \,\big\}\,.$$

In fact, the two notations $\Delta_{\partial^{\,in}\Lambda}\pi$ and $\partial^{edge}_{\partial^{\,in}\Lambda}\pi$ designate the same set of edges. To keep in line with the intersection sets of the intertwined explorations, we will use the lighter notation $\Delta_{\partial^{\,in}\Lambda}\pi$.

**Lemma 93.5.** *Let $\Lambda$ be a cubic box in $\mathbb{Z}^d$, and let $\pi = (\pi_i, 1 \leq i \leq r)$ be a partition of $\partial^{\,in}\Lambda$ in $r$ sets. We suppose that $r \geq 32d$ and that at least $r/8$ sets of the partition $\pi$ have a cardinality larger than or equal to $v$, where $v \geq 1$. We have then*

$$\Big|\bigcup_{1\leq j<k\leq r}\big(\Delta_{\partial^{\,in}\Lambda}\pi_j \cap \Delta_{\partial^{\,in}\Lambda}\pi_k\big)\Big| \,\geq\, \frac{1}{64d^2}rv^{1-1/(d-1)}\,.$$

*Proof.* We decompose $\partial^{in}\Lambda$ as the union of $2d$ faces, $F_1, \dots, F_{2d}$. Each face is isometric to the same $(d-1)$-dimensional box, that we denote by $\Lambda_{d-1}$. Let $i$ be an index in $\{1, \dots, r\}$ such that $|\pi_i| \geq v$. There exists a face $F_j$ such that $\big|\pi_i \cap F_j\big| \geq v/2d$. Naturally, the face $F_j$ depends on the set $\pi_i$. However, there are only $2d$ faces, and at least $r/8$ sets $\pi_i$. By the pigeonhole principle, there exists a face which is associated to at least $r/(16d)$ sets $\pi_i$. Suppose that this face is $F_1$. Then there exists a subset $I$ of $\{1, \dots, r\}$ such that $|I| \geq r/(16d)$ and

$$\forall i \in I \qquad \big|\pi_i \cap F_1\big| \geq \frac{v}{2d}.$$

We shift our attention to the face $F_1$, which is isometric to $\Lambda_{d-1}$. Let $i$ belong to $I$. Since $r \geq 32d$, then $|I| \geq 2$ and there exist at least two distinct sets $\pi_j$, $\pi_k$ such that

$$\big|\pi_j \cap F_1\big| \geq \frac{v}{2d}, \qquad \big|\pi_k \cap F_1\big| \geq \frac{v}{2d}.$$

This implies that

$$\forall i \in I \qquad \min\Big(\big|\pi_i \cap F_1\big|, \big|F_1 \setminus \pi_i\big|\Big) \geq \frac{v}{2d}. \tag{93.22}$$

We apply the discrete relative isoperimetric inequality in the $(d-1)$-dimensional box $F_1$. Theorem 93.2 and the lower bounds (93.22) yield that

$$\forall i \in I \qquad \big|\Delta_{F_1}\pi_i\big| = \big|\Delta_{F_1}(\pi_i \cap F_1)\big| \geq \Big(\frac{v}{2d}\Big)^{1-1/(d-1)}. \tag{93.23}$$

Let us fix $i$ in $I$. Let $e$ be an edge in $\Delta_{F_1}\pi_i$. One extremity of $e$ is in $\pi_i$, and the other belongs to another set of the partition, say $\pi_j$. The edge $e$ belongs to exactly one set of the collection $\big(\Delta_{\partial^{in}\Lambda}\pi_j \cap \Delta_{\partial^{in}\Lambda}\pi_k,\, 1 \leq j < k \leq r\big)$, but possibly to two sets among $\Delta_{\partial^{in}\Lambda}\pi_i$, $i \in I$. The previous remarks and the inequalities (93.23) readily imply that

$$\Big|\bigcup_{1\leq j<k\leq r}\big(\Delta_{\partial^{in}\Lambda}\pi_j \cap \Delta_{\partial^{in}\Lambda}\pi_k\big)\Big| \geq \Big|\bigcup_{i\in I}\Delta_{\partial^{in}\Lambda}\pi_i\Big|$$
$$\geq \frac{1}{2}\sum_{i\in I}\big|\Delta_{\partial^{in}\Lambda}\pi_i\big| \geq \frac{r}{32d}\Big(\frac{v}{2d}\Big)^{1-1/(d-1)} \geq \frac{1}{64d^2}rv^{1-1/(d-1)}.$$

This is the desired inequality. □

We apply lemma 93.5 to the box $\Lambda(s)$ and the partition induced by the traces of the aggregates $(\mathcal{A}_i(T), 1 \leq i \leq N)$ on $\partial^{in}\Lambda(s)$, with $r = N$ and $v = n^{\alpha/2}$. Setting

$$\delta = \frac{\alpha(d-2)}{2(d-1)},$$

we obtain

$$\Big|\bigcup_{1\leq j<k\leq N}\big(\Delta_{\partial^{in}\Lambda(s)}\mathcal{A}_j(T) \cap \Delta_{\partial^{in}\Lambda(s)}\mathcal{A}_k(T)\big)\Big| \geq \frac{1}{64d^2}Nn^{\delta}. \tag{93.24}$$

The previous inequality holds for any $s$ in the set of indices $\overline{\mathcal{S}}(h,n)$ that was defined in (93.21). The boundaries $\big(\partial^{\,in}\Lambda(s), s\in\overline{\mathcal{S}}(h,n), s \text{ even}\big)$, are pairwise disjoint. Summing the inequality (93.24) over the even integers $s$ in $\overline{\mathcal{S}}(h,n)$, and using (92.15), we conclude that

$$\Big|\bigcup_{0\leq t\leq T-1}\mathcal{I}(t)\cap\Lambda(h-1,n)\Big|\geq$$
$$\Big|\bigcup_{\substack{s\in\overline{\mathcal{S}}(h,n), s\text{ even}\\ 1\leq j<k\leq N}}\big(\Delta_{\partial^{\,in}\Lambda(s)}\mathcal{A}_j(T)\cap\Delta_{\partial^{\,in}\Lambda(s)}\mathcal{A}_k(T)\big)\Big|$$
$$\geq\frac{1}{512d^2}Nn^{\delta}\lambda(n)\,\geq\,\frac{1}{1024d^3\ln n}Nn^{1+\delta}\,. \tag{93.25}$$

Using corollary 80.19, we deduce from (93.25) that, on the event $\mathcal{H}'$,

$$\big|\overline{\mathcal{T}}^0(T)\cap\Lambda(h-1,n)\big|\,\geq\,\frac{Nn^{1+\delta}}{4096d^4\ln n}\,. \tag{93.26}$$

We have proved that any configuration of the event $\mathcal{H}'$ defined in (92.27) satisfies the inequality (93.26). This inequality is better than the one obtained in subsection 92.3, see formula (92.33). Indeed, the exponent of $n$ is strictly larger here, and this will make all the difference for the final argument!

With the inequality (93.26) in hand, we continue the argument, from the subsection 92.4 until the subsection 92.16 included. There is no change at all in this whole part. We arrive at the final subsection 92.17, where the exponent two was obtained. Throughout this final subsection, we use our improved inequality (93.26) instead of the previous one (92.188). In practice, from formula (92.188) onwards, the following change should be done:

$$\text{replace the value}\qquad\frac{Nn}{64d^3\ln n}\qquad\text{by}\qquad\frac{Nn^{1+\delta}}{4096d^4\ln n}\,.$$

The condition (92.190) becomes

$$\nu(n)\,n^{d/2}\sqrt{N}\,<\,\frac{Nn^{1+\delta}}{4096d^4\ln n}\,.$$

We go on as in subsection 92.17 and we conclude that

$$\lim_{n\to\infty}P\Big(\mathcal{N}^*\geq\big(4097d^4\nu(n)\ln n\big)^2\,n^{d-2-2\delta},\mathcal{E}\Big)\,=\,0\,.$$

This implies that, when $\mathcal{E}$ occurs, with probability going to 1 as $n$ goes to $\infty$,

$$\mathcal{N}^*\,\leq\,n^{d-2-\alpha(d-2)/(d-1)+o(1)}\,. \tag{93.27}$$

We conclude also that, with probability going to 1 as $n$ goes to $\infty$, the box $\Lambda(4n)$ contains a boundary cluster $C$ such that

$$\big|C\cap\Lambda(2n)\big|\,\geq\,n^{2+\alpha(d-2)/(d-1)+o(1)}\,.$$

Thanks to the isoperimetric inequality in dimension $d-1$, we have gained a little on the exponent 2. However, as we will see shortly, what really matters for the sequel is the inequality (93.27) on the number of boundary clusters.

## 93.3 Improving the exponent from $2+\alpha(d-2)/(d-1)$ to $d$

We are now ready to present the final part of the argument. We will use the full power of the discrete relative isoperimetric inequality in dimension $d$. We restart again the argument of section 92 (for the last time!), with some additional changes. We proceed first exactly as in the previous section, until formula (93.26), which is recalled next:

$$\left|\overline{\mathcal{T}}^0(T)\cap\Lambda(h-1,n)\right| \;\geq\; \frac{Nn^{1+\frac{\alpha(d-2)}{2(d-1)}}}{4096d^4\ln n}\,. \tag{93.28}$$

At this point, we retake the details of the subsection 92.17. Thus, we rewrite the final part of the argument of subsection 92.17, starting right after formula (92.188), which is replaced by the better inequality (93.28). The following steps until formula (93.31) are the same as in subsection 92.17, the only difference being that the right-hand side of (92.188) is replaced by the one of (93.28). Still, we prefer to rewrite them with the adequate values, otherwise the reader would have to jump back and forth between the two sections, and at the same time make the replacement of values alluded previously. To sum up, the argument should be read as follows:

• add the event Thick($n$) in the definition of $\mathcal{E}$, and take the beginning of the subsection 92.17 until formula (92.188);

• at this point, use the computations of subsection 93.2 to derive the inequality (93.26), which is the improvement of (92.188) obtained by applying the relative isoperimetric inequality in the boundary of the intermediate boxes;

• once the inequality (93.26) (recalled above in (93.28)) is proved, continue from here.

We recall the definition (92.27) of the event $\mathcal{H}'$:

$$\mathcal{H}' \;=\; \mathcal{E}\cap\left\{\begin{array}{c} H=h,\mathcal{N}^*(h+1,n)=N,\\ \mathcal{C}^*(h+1,n)=\big\{\,C\big(w_i,\Lambda(4n)\big),1\leq i\leq N\,\big\}\end{array}\right\}.$$

Notice that the event $\mathcal{H}'$ depends on $h,N$ and the sites $w_1,\dots,w_N$. We shall discuss further the inequalities satisfied by the configurations in $\mathcal{H}'$. In subsection 92.6, we introduced the event Intersect($N$), and we showed in inequality (92.186) that the configurations belonging to $\mathcal{H}'\cap\text{Intersect}(N)$ satisfy

$$\sum_{0\leq s\leq T}\left|\overline{\mathcal{T}}(s)\cap\Gamma(s)\right| \;\leq\; \nu(n)\,n^{d/2}\sqrt{N}\,, \tag{93.29}$$

where the prefactor $\nu(n)$, given in (92.187), is sub-polynomial. Since the boxes $\Gamma=\big(\Gamma(t),0\leq t\leq T\big)$ all contain $\Lambda(h-1,n)$, the inequality (93.29) implies that

$$\left|\overline{\mathcal{T}}^0(T)\cap\Lambda(h-1,n)\right| \;\leq\; \nu(n)\,n^{d/2}\sqrt{N}\,. \tag{93.30}$$

For the values of $N$ such that

$$\nu(n)\,n^{d/2}\sqrt{N} \;<\; \frac{Nn^{1+\frac{\alpha(d-2)}{2(d-1)}}}{4096d^4\ln n}\,, \tag{93.31}$$

the two inequalities (93.28) and (93.30) are not compatible, this means that there are no configurations in $\mathcal{H}'\cap\text{Intersect}(N)$. We divide the inequality (93.31) by $\sqrt{N}$, we take the squares and we rewrite it as

$$N\,>\,\big(4096d^4\nu(n)\ln n\big)^2\,n^{d-2-\alpha(d-2)/(d-1)}\,.$$

We choose now

$$N^{**}\,=\,\big(4097d^4\nu(n)\ln n\big)^2\,n^{d-2-\alpha(d-2)/(d-1)}\,.$$

We study next the configurations belonging to $\mathcal{H}'$, and we consider separately the cases where $N\geq N^{**}$ and $N<N^{**}$. From the previous discussion, we see that, for any integer $N\geq N^{**}$, the event $\mathcal{H}'$ is included in the complement of $\text{Intersect}(N)$, so that, by (92.55),

$$\forall N\geq N^{**}\qquad P(\mathcal{H}')\,\leq\,P\big(\text{Intersect}(N)^c\big)\,\leq\,\exp\big(-4^dN\ln n\big)\,.\qquad(93.32)$$

We consider now the more delicate case $N<N^{**}$. Let us fix $N<N^{**}$ and a configuration belonging to $\mathcal{H}'\cap\text{Intersect}(N)$. According to inequality (93.30), we have

$$\big|\overline{\mathcal{T}}^0(T)\cap\Lambda(h-1,n)\,\leq\,\nu(n)\,n^{d/2}\sqrt{N^{**}}\,=\,\mu(n)\,n^\gamma\,,\qquad(93.33)$$

where we have set

$$\mu(n)\,=\,4097d^4\nu(n)^2(\ln n)\,,\qquad\gamma\,=\,d-1-\frac{\alpha(d-2)}{2(d-1)}\,.$$

The prefactor $\nu(n)$, given in (92.187), is sub-polynomial, and so is $\mu(n)$. Later on, we will sum the inequalities over $h$, so we would like to get rid of the dependence in $h$ in (93.33). In addition, we are ultimately interested in the box $\Lambda(2n)$. As $\Lambda(2n)$ is contained in $\Lambda(h-1,n)$ for any $h$ in $\{\,0,\dots,\overline{h}(n)\,\}$, the inequality (93.33) implies that

$$\big|\overline{\mathcal{T}}^0(T)\cap\Lambda(2n)\big|\,\leq\,\mu(n)\,n^\gamma\,.\qquad(93.34)$$

Let us stress that the random set $\overline{\mathcal{T}}^0(T)$ is decided by the intertwined explorations, which themselves depend on the choice of the starting sites $w_1,\dots,w_N$. We would also like to get rid of that dependence. To that end, we try to deduce from the inequality (93.34) a weaker inequality that involves directly the boundary clusters of $\Lambda(4n)$ which intersect $\Lambda(2n)$. Proposition 80.18 yields that

$$\Big|\bigcup_{0\leq t\leq T-1}\mathcal{I}(t)\cap\Lambda(2n)\Big|\,\leq\,4d\big|\overline{\mathcal{T}}^0(T)\cap\Lambda(2n)\big|\,.\qquad(93.35)$$

Coming back to the definition 80.12 of the intersection sets, we have furthermore

$$\bigcup_{1\leq j<k\leq N}\big(\Delta_{\Lambda(2n)}\mathcal{A}_j(T)\cap\Delta_{\Lambda(2n)}\mathcal{A}_k(T)\big)\,\subset\,\bigcup_{0\leq t\leq T-1}\mathcal{I}(t)\cap\Lambda(2n)\,.\qquad(93.36)$$

If an edge $e$ is in $\Delta_{\Lambda(2n)}\mathcal{A}_j(T)\cap\Delta_{\Lambda(2n)}\mathcal{A}_k(T)$ for some $j<k$, then it does not belong to $\Delta_{\Lambda(2n)}\mathcal{A}_i(T)$ for $i$ in $\{\,1,\dots,N\,\}\setminus\{\,j,k\,\}$. We have thus

$$\Big|\bigcup_{1\le j<k\le N}\big(\Delta_{\Lambda(2n)}\mathcal{A}_j(T)\cap\Delta_{\Lambda(2n)}\mathcal{A}_k(T)\big)\Big| \;=\; \frac{1}{2}\sum_{1\le i\le N}\big|\Delta_{\Lambda(2n)}\mathcal{A}_i(T)\big|\,. \tag{93.37}$$

The two inequalities (93.34), (93.35), the inclusion (93.36), and the equality (93.37) together yield that

$$\sum_{1\le i\le N}\big|\Delta_{\Lambda(2n)}\mathcal{A}_i(T)\big| \;\le\; 8d\mu(n)\,n^{\gamma}\,. \tag{93.38}$$

Let $i^*$ be an index in $\{\,1,\dots,N\,\}$ such that

$$\big|\mathcal{A}_{i^*}(T)\cap\Lambda(2n)\big| \;=\; \max\,\big\{\,\big|\mathcal{A}_i(T)\cap\Lambda(2n)\big| : 1\le i\le N\,\big\}\,.$$

We apply next the discrete relative isoperimetric inequality in $\Lambda(2n)$ to $\mathcal{A}_{i^*}(T)$. Theorem 93.2 yields that

$$\min\Big(\big|\mathcal{A}_{i^*}(T)\cap\Lambda(2n)\big|,\big|\Lambda(2n)\setminus\mathcal{A}_{i^*}(T)\big|\Big) \;\le\; \big|\Delta_{\Lambda(2n)}\mathcal{A}_{i^*}(T)\big|^{d/(d-1)}\,. \tag{93.39}$$

Using the bound (93.38), we deduce from (93.39) that

$$\min\Big(\big|\mathcal{A}_{i^*}(T)\cap\Lambda(2n)\big|,\big|\Lambda(2n)\setminus\mathcal{A}_{i^*}(T)\big|\Big) \;\le\; \big(8d\mu(n)\,n^{\gamma}\big)^{d/(d-1)} \;\le\; \rho(n)n^{d-\eta}\,, \tag{93.40}$$

where

$$\rho(n) \;=\; \big(8d\mu(n)\big)^{d}\,,\qquad \eta \;=\; \alpha\frac{d(d-2)}{2(d-1)^2}\,.$$

Suppose that

$$\big|\mathcal{A}_{i^*}(T)\cap\Lambda(2n)\big| \;\le\; \rho(n)n^{d-\eta}\,. \tag{93.41}$$

We have then

$$\forall i\in\{\,1,\dots,N\,\}\qquad \big|\mathcal{A}_i(T)\cap\Lambda(2n)\big| \;\le\; \rho(n)n^{d-\eta}\,. \tag{93.42}$$

We define a sequence of sets $(\mathcal{B}_k, 1\le k\le N)$ by setting

$$\forall k\in\{\,1,\dots,N\,\}\qquad \mathcal{B}_k \;=\; \bigcup_{1\le i\le k}\mathcal{A}_i(T)\cap\Lambda(2n)\,. \tag{93.43}$$

Let us fix $k$ in $\{\,1,\dots,N\,\}$. We have then, using (93.38),

$$\big|\Delta_{\Lambda(2n)}\mathcal{B}_k\big| \;\le\; \Big|\bigcup_{1\le i\le N}\Delta_{\Lambda(2n)}\mathcal{A}_i(T)\Big| \;\le\; \sum_{1\le i\le N}\big|\Delta_{\Lambda(2n)}\mathcal{A}_i(T)\big| \;\le\; 8d\mu(n)\,n^{\gamma}\,.$$

We apply the discrete relative isoperimetric inequality in $\Lambda(2n)$ to the set $\mathcal{B}_k$, and we get

$$\min\Big(\big|\mathcal{B}_k\big|,\big|\Lambda(2n)\setminus\mathcal{B}_k\big|\Big) \;\le\; \big(8d\mu(n)\,n^{\gamma}\big)^{d/(d-1)} \;\le\; \rho(n)n^{d-\eta}\,. \tag{93.44}$$

Furthermore, the sequence $(\mathcal{B}_k, 1 \leq k \leq N)$ is increasing. Let us consider the index $j^*$ defined by

$$j^* \,=\, \max\Big\{\, i : 1 \leq i \leq N,\, \big|\mathcal{B}_i\big| \leq \rho(n) n^{d-\eta} \,\Big\}\,.$$

The inequalities (93.42) imply that the set above contains 1, thus $j^*$ is well-defined and belongs to $\{\,1,\dots,N\,\}$. As $\mathcal{B}_N$ is equal to the whole box $\Lambda(2n)$, the index $j^*$ is strictly smaller than $N$. Let us set $k^* = j^* + 1$. On one hand, we have

$$\big|\mathcal{B}_{j^*}\big| \,\leq\, \rho(n) n^{d-\eta} \,<\, \big|\mathcal{B}_{k^*}\big|\,. \tag{93.45}$$

On the other hand, it follows from the definition (93.43) and the bound (93.42) that

$$\big|\mathcal{B}_{k^*}\big| \,\leq\, \big|\mathcal{B}_{j^*}\big| + \big|\mathcal{A}_{k^*}(T) \cap \Lambda(2n)\big| \,\leq\, 2\rho(n) n^{d-\eta}\,. \tag{93.46}$$

Since the prefactor $\rho(n)$ is sub-polynomial, and the exponent $\eta$ is positive, then

$$\begin{aligned} &\exists\, n_6(\alpha,d,p) \geq 1 \quad \forall n \geq n_6(\alpha,d,p) \\ &\qquad 2\rho(n) n^{d-\eta} \,<\, \min\Big(\frac{1}{2}\big|\Lambda(2n)\big|, \big(\theta(p)-\varepsilon\big)\big|\Lambda(2n)\big|, \varepsilon\big|\Lambda(2n)\big|\Big)\,. \end{aligned} \tag{93.47}$$

The parameter $\varepsilon$ was chosen at the very beginning of the proof in subsection 92.1. From now onwards, we suppose that $n$ is larger than $n_6(\alpha,d,p)$. The inequalities (93.46) and (93.47) imply that

$$\big|\Lambda(2n) \setminus \mathcal{B}_{k^*}\big| \,\geq\, \big|\Lambda(2n)\big| - \big|\mathcal{B}_{k^*}\big| \,\geq\, \big|\Lambda(2n)\big| - 2\rho(n) n^{d-\eta} > \frac{1}{2}\big|\Lambda(2n)\big|\,.$$

The isoperimetric inequality (93.44) would then give that $\big|\mathcal{B}_{k^*}\big| \leq \rho(n) n^{d-\eta}$, thereby contradicting (93.45). We conclude that the inequality (93.41) does not hold. The isoperimetric inequality (93.40) leads then to the other case, where

$$\big|\Lambda(2n) \setminus \mathcal{A}_{i^*}(T)\big| \,\leq\, \rho(n) n^{d-\eta}\,. \tag{93.48}$$

With this information in hand, we finally come back to the boundary clusters of $\Lambda(4n)$ which intersect $\Lambda(2n)$, that is the collection

$$\mathcal{C}^*(n) \,=\, \big\{\, C \in \mathcal{C}_{\mathrm{bd}}(4n) : C \cap \Lambda(2n) \neq \varnothing \,\big\}\,.$$

The collection $\mathcal{C}^*(n)$ is a sub-collection of $\mathcal{C}^*(h+1,n)$. As each aggregate contains exactly one cluster of $\mathcal{C}^*(h+1,n)$, each aggregate contains at most one cluster of $\mathcal{C}^*(n)$. Suppose that the aggregate $\mathcal{A}_{i^*}(T)$ does not contain any cluster of $\mathcal{C}^*(n)$. We would then have

$$\bigcup_{C \in \mathcal{C}^*(n)} C \cap \Lambda(2n) \,\subset\, \Lambda(2n) \setminus \mathcal{A}_{i^*}(T)\,,$$

and it would follow from (93.48) and (93.47) that

$$\sum_{C \in \mathcal{C}_{\mathrm{bd}}(4n)} \big|C \cap \Lambda(2n)\big| \,\leq\, \rho(n) n^{d-\eta} \,<\, \big(\theta(p)-\varepsilon\big)\big|\Lambda(2n)\big|\,.$$

This is not compatible with the occurrence of the event $\mathcal{F}_{\text{bd}}(n,\varepsilon)$. Thus this case is not possible, and the aggregate $\mathcal{A}_{i^*}(T)$ must contain a cluster of $\mathcal{C}^*(n)$. We denote this cluster by $C^*$. We have then

$$\bigcup_{C\in\mathcal{C}^*(n),C\neq C^*} C\cap\Lambda(2n) \,\subset\, \Lambda(2n)\setminus\mathcal{A}_{i^*}(T)\,,$$

whence, using (93.48) and (93.47),

$$\Big|\bigcup_{C\in\mathcal{C}^*(n),C\neq C^*} C\cap\Lambda(2n)\Big| \,\leq\, \rho(n)n^{d-\eta} \,<\, \varepsilon\big|\Lambda(2n)\big|\,. \tag{93.49}$$

Since the event $\mathcal{F}_{\text{bd}}(n,\varepsilon)$ occurs, we have

$$\sum_{C\in\mathcal{C}_{\text{bd}}(4n)} \big|C\cap\Lambda(2n)\big| \,\geq\, \big(\theta(p)-\varepsilon\big)\big|\Lambda(2n)\big|\,. \tag{93.50}$$

The two inequalities (93.49) and (93.50) yield the inequality we have been striving for:

$$\big|C^*\cap\Lambda(2n)\big| \,\geq\, \big(\theta(p)-2\varepsilon\big)\big|\Lambda(2n)\big|\,. \tag{93.51}$$

At this point, we have proved that, for $N<N^{**}$, any configuration belonging to $\mathcal{H}'\cap\text{Intersect}(N)$ contains a boundary cluster $C^*$ satisfying (93.51).

We will now gather the previous considerations on the event $\mathcal{H}'$ in order to describe the typical scenario occurring along with the event $\mathcal{E}$. We decompose $P\big(\mathcal{E},H=h\big)$ according to the value of $\mathcal{N}^*(h+1,n)$:

$$P\big(\mathcal{E},H=h\big) \,=\, \sum_{1\leq N\leq|\partial^{\,in}\Lambda(4n)|} P\begin{pmatrix}\mathcal{E},H=h,\\ \mathcal{N}^*(h+1,n)=N\end{pmatrix}\,. \tag{93.52}$$

We decompose further the event inside the last probability with the help of the event $\mathcal{H}'$, as follows:

$$P\begin{pmatrix}\mathcal{E},H=h,\\ \mathcal{N}^*(h+1,n)=N\end{pmatrix} \,=\, P\Big(\bigcup_{w_1,\dots,w_N\in\partial^{\,in}\Lambda(4n)} \mathcal{H}'(h,N,w_1,\dots,w_N)\Big)\,. \tag{93.53}$$

We will bound from above the probability (93.53), and we distinguish the cases $N<N^{**}$ and $N\geq N^{**}$. We examine first the easier case $N\geq N^{**}$. We use simply the standard union bound to write

$$P\begin{pmatrix}\mathcal{E},H=h,\\ \mathcal{N}^*(h+1,n)=N\end{pmatrix} \,\leq\, \sum_{w_1,\dots,w_N\in\partial^{\,in}\Lambda(4n)} P\big(\mathcal{H}'(h,N,w_1,\dots,w_N)\big)\,. \tag{93.54}$$

We substitute the uniform upper bound (93.32) in the sum (93.54), and we get

$$\begin{aligned}P\begin{pmatrix}\mathcal{E},H=h,\\ \mathcal{N}^*(h+1,n)=N\end{pmatrix} \,&\leq\, \sum_{w_1,\dots,w_N\in\partial^{\,in}\Lambda(4n)} \exp\big(-4^dN\ln n\big)\\ &\leq\, \exp\Big(N\ln\big(|\Lambda(4n)|\big)-4^dN\ln n\Big) \,\leq\, \exp\big(-dN\ln n\big)\,. \end{aligned}\tag{93.55}$$

We examine now the case $N < N^{**}$. We define the event Ultimate to be the event we want to show is typical:

$$\text{Ultimate} = \Big\{ \exists C^* \in \mathcal{C}_{\text{bd}}(4n) \qquad \big|C^* \cap \Lambda(2n)\big| \geq \big(\theta(p) - 2\varepsilon\big)\big|\Lambda(2n)\big| \Big\} .$$

We have shown previously in (93.51) that, for $N < N^{**}$, we have the inclusion

$$\mathcal{H}'(h, N, w_1, \dots, w_N) \cap \text{Intersect}(N) \subset \text{Ultimate} \cap \text{Intersect}(N) . \qquad (93.56)$$

The good point is that the event in the right-hand side of (93.56) does not depend any more on $w_1, \dots, w_N$. We take the union of (93.56) over the possible choices of $w_1, \dots, w_N$, and we obtain

$$\bigcup_{w_1,\dots,w_N \in \partial^{\,in}\Lambda(4n)} \mathcal{H}'(h, N, w_1, \dots, w_N) \cap \text{Intersect}(N) \subset \text{Ultimate} \cap \text{Intersect}(N) . \qquad (93.57)$$

Coming back to the equality (93.53), we use the inclusion (93.57) to conclude that

$$P\begin{pmatrix} \mathcal{E}, H = h, \\ \mathcal{N}^*(h+1, n) = N \end{pmatrix} \leq P\begin{pmatrix} \mathcal{E}, H = h, \\ \mathcal{N}^*(h+1, n) = N, \\ \text{Ultimate} \end{pmatrix} + P\big(\text{Intersect}(N)^c\big) . \qquad (93.58)$$

We substitute finally the inequalities (93.55) and (93.58) into the sum (93.52), and we obtain, using once more the upper bound (92.55),

$$\begin{aligned} P\big(\mathcal{E}, H = h\big) \leq \sum_{N^{**} \leq N \leq \left|\partial^{\,in}\Lambda(4n)\right|} & \exp\big(-dN \ln n\big) + \\ \sum_{1 \leq N < N^{**}} & \left( P\begin{pmatrix} \mathcal{E}, H = h, \\ \mathcal{N}^*(h+1, n) = N, \\ \text{Ultimate} \end{pmatrix} + P\big(\text{Intersect}(N)^c\big) \right) \\ & \leq P\begin{pmatrix} \mathcal{E}, H = h, \\ \text{Ultimate} \end{pmatrix} + \big|\partial^{\,in}\Lambda(4n)\big| \exp\big(-d \ln n\big) . \end{aligned} \qquad (93.59)$$

The trick is that we were able to bring the summation over $N$ back into the probability. Notice that the event Ultimate depends on $n, \varepsilon$, but not on $N$. We finally sum (93.59) over $h$ in $\{\, 0, \dots, \overline{h}(n) \,\}$ and we get

$$\begin{aligned} P(\mathcal{E}) \leq \sum_{0 \leq h < \overline{h}(n)} & \left( P\begin{pmatrix} \mathcal{E}, H = h, \\ \text{Ultimate} \end{pmatrix} + \big|\partial^{\,in}\Lambda(4n)\big| \exp\big(-d \ln n\big) \right) \\ & \leq P\big(\mathcal{E}, \text{Ultimate}\big) + \overline{h}(n) \,\big|\partial^{\,in}\Lambda(4n)\big| \exp\big(-d \ln n\big) . \end{aligned}$$

Thanks to (93.19) and the fact that

$$\lim_{n\to\infty}\ \overline{h}(n)\left|\partial^{\,in}\Lambda(4n)\right|\exp\big(-d\ln n\big)\ =\ 0\,,$$

we conclude that

$$\lim_{n\to\infty} P\big(\text{Ultimate}\big)\ =\ 1\,.$$

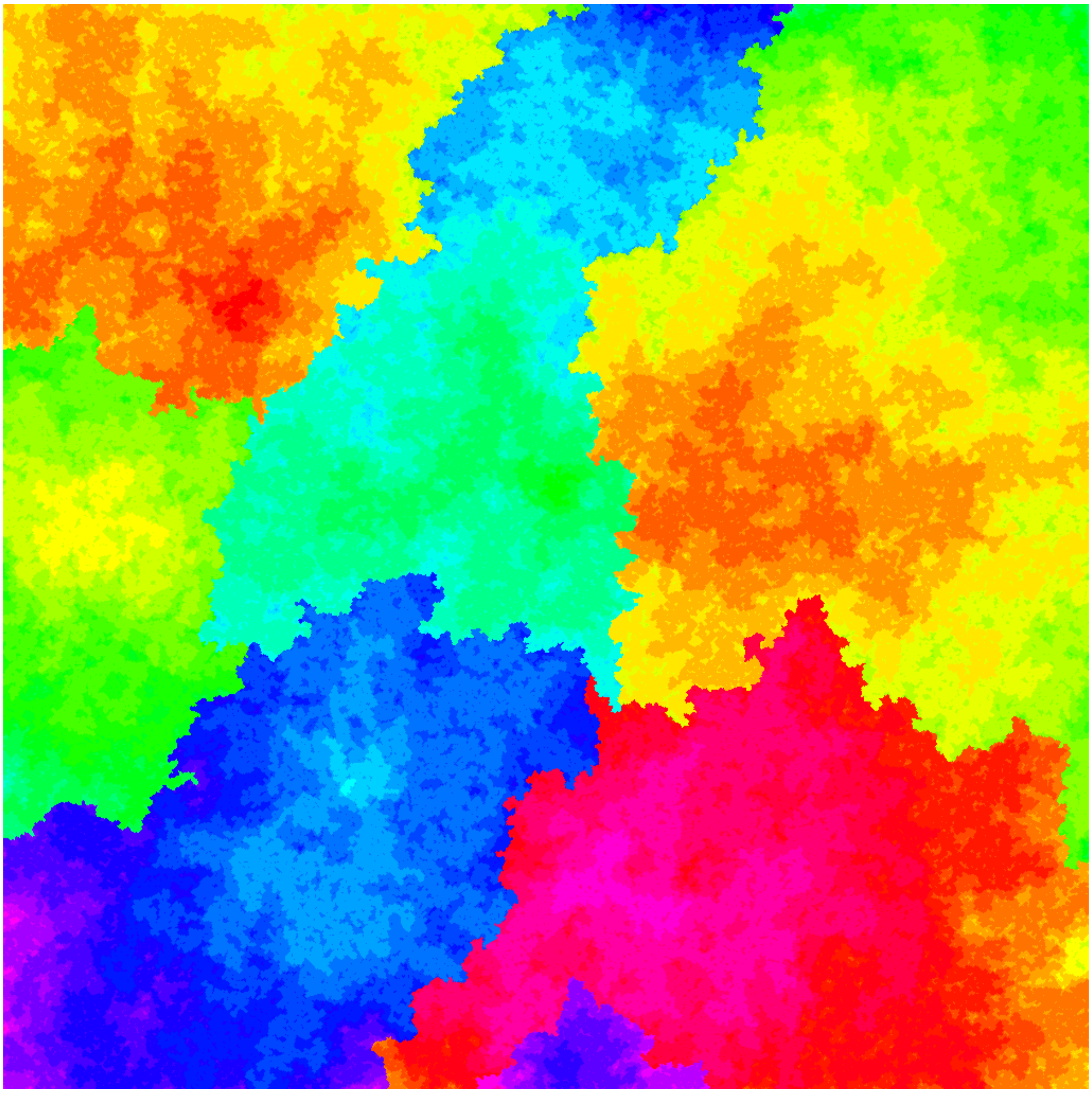

Figure 171: 7 intertwined explorations, $1024\times 1024$, site percolation, $p=0.57$.

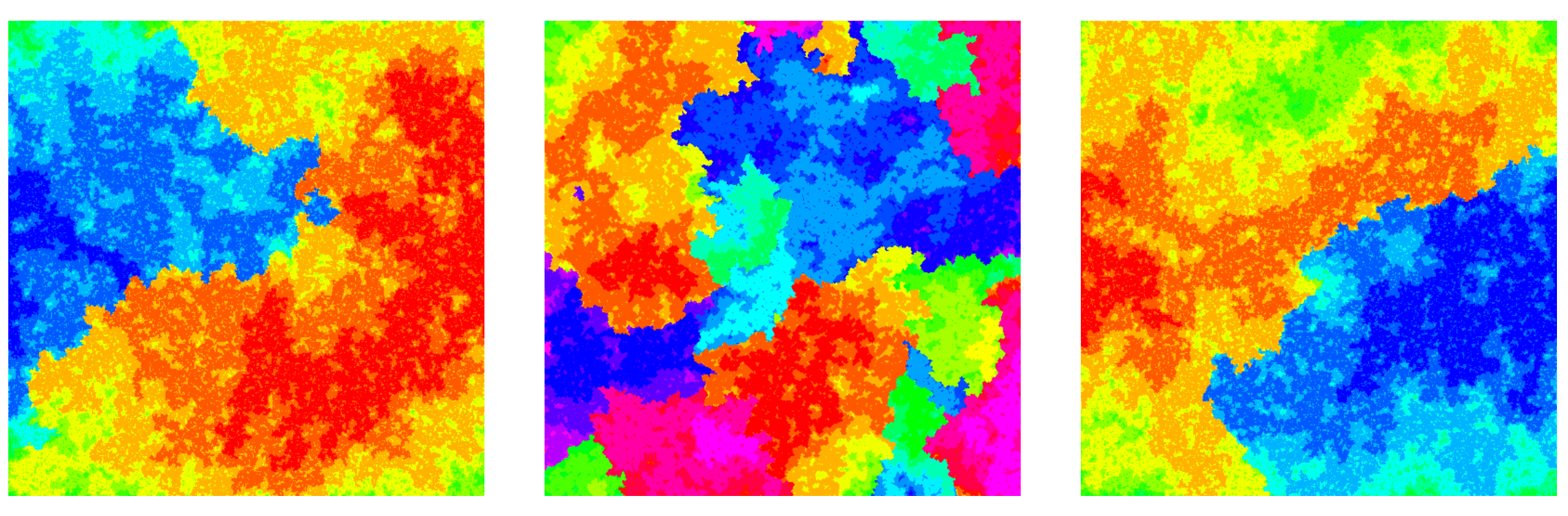

# A Appendix: A useful lemma

The next lemma is an adaptation to the case of integer-valued random variables of a basic lemma of Talagrand (lemma 2.2.3 of [126]).

**Lemma A.1.** *Consider a random variable $N$ with values in $\mathbb{N}$ which satisfies*

$$\forall i \geq 1 \qquad P(N \geq i) \;\leq\; A \exp\Big(-\frac{i^2}{B^2}\Big)\,, \tag{A.1}$$

*for certain numbers $A \geq 3$ and $B \geq 1$. Then*

$$E(N) \leq 2B\sqrt{\ln(A)}\,. \tag{A.2}$$

*Proof.* Let $i_0 \geq 1$ be a fixed integer. We decompose the expectation of $N$ as follows:

$$E(N) \,=\, \sum_{i \geq 1} P(N \geq i) \,=\, \sum_{1 \leq i \leq i_0} P(N \geq i) + \sum_{i > i_0} P(N \geq i)\,. \tag{A.3}$$

In the last sum, we bound $P(N \geq i)$ by an integral with the help of (A.1):

$$P(N \geq i) \;\leq\; A \exp\Big(-\frac{i^2}{B^2}\Big) \;\leq\; \int_{i-1}^{i} A \exp\Big(-\frac{u^2}{B^2}\Big)\,du\,. \tag{A.4}$$

We substitute (A.4) into (A.3), and we obtain

$$E(N) \;\leq\; i_0 + \int_{i_0}^{+\infty} A \exp\Big(-\frac{u^2}{B^2}\Big)\,du\,.$$

We continue the computation as Talagrand:

$$\int_{i_0}^{+\infty} A \exp\Big(-\frac{u^2}{B^2}\Big)\,du \;\leq\; \frac{1}{i_0}\int_{i_0}^{+\infty} Au \exp\Big(-\frac{u^2}{B^2}\Big)\,du \;\leq\; \frac{AB^2}{2i_0}\exp\Big(-\frac{i_0^2}{B^2}\Big)\,,$$

whence

$$E(N) \;\leq\; i_0 + \frac{AB^2}{2i_0}\exp\Big(-\frac{i_0^2}{B^2}\Big)\,.$$

We choose $i_0 \;=\; B\sqrt{\ln(A)}$. The hypothesis on $A$ and $B$ ensures that $i_0 \geq 1$ and we get the desired inequality (A.2). $\square$

# Index